\documentclass[a4paper,10pt]{extarticle}
 
\usepackage[utf8]{inputenc}
\usepackage[T1]{fontenc}
\usepackage[english]{babel}
\usepackage{lmodern}
\usepackage{geometry}
\usepackage{url,xspace}
\usepackage{amsmath,amssymb,stmaryrd,spalign,bbm,mathtools,mathrsfs,dsfont,thmtools}
\usepackage{xfrac}
\usepackage[thmmarks,amsmath,framed,thref,hyperref]{ntheorem}
\usepackage{comment}
\usepackage[colorlinks=true,citecolor=teal,linkcolor=blue]{hyperref}
\usepackage{upgreek}
\usepackage{array}
\usepackage{enumitem}
\usepackage{changepage}
\usepackage[nottoc]{tocbibind} 
\usepackage{euscript}
\usepackage{cite}
\usepackage{graphicx}
\usepackage[toc,page]{appendix}
\usepackage{etoolbox}
\usepackage{pgfplots}
\pgfplotsset{width=10cm,compat=1.9}

\usepackage{tikz, xcolor, float}
\usetikzlibrary{arrows.meta,patterns,calc}

\usepackage{makecell}
\usepackage{tikz-cd}

\usepgfplotslibrary{external}

\newtheorem{theorem}{Theorem}[section]
    \newtheorem{corollary}[theorem]{Corollary}
    \newtheorem{lemma}[theorem]{Lemma}
    \newtheorem{proposition}[theorem]{Proposition}
    
    \newtheorem{assumption}[theorem]{Assumption}
    \newtheorem{example}[theorem]{Example}
    \theorembodyfont{\upshape}
    \newtheorem{definition}[theorem]{Definition}
    \newtheorem{notation}[theorem]{Notation}
    \newtheorem{remark}[theorem]{Remark}
    \newtheorem{conjecture}[theorem]{Conjecture}

\newtheorem{question}{Question}

\theoremstyle{nonumberplain}
\theoremheaderfont{\itshape}
\theorembodyfont{\normalfont}
\theoremseparator{.\,—}
\theoremsymbol{}
\newtheorem{proof-wo}{Proof}
\theoremsymbol{\ensuremath{\blacksquare}}
\newtheorem{proof}{Proof}

\newcommand{\eus}{\EuScript}

\newcommand{\B}{\mathrm{B}}
\newcommand{\C}{\mathrm{C}}     \newcommand{\Ccinfty}{\mathrm{C}^\infty_c}     \newcommand{\Ccinftydiv}{\mathrm{C}^\infty_{c,\sigma}}
\newcommand{\D}{\mathrm{D}}

\newcommand{\E}{\mathrm{E}}
\renewcommand{\d}{\mathrm{d}}
\renewcommand{\H}{\mathrm{H}}
\newcommand{\I}{\mathrm{I}}
\renewcommand{\L}{\mathrm{L}}
\renewcommand{\P}{\mathrm{P}}
\newcommand{\N}{\mathrm{N}}
\newcommand{\T}{\mathrm{T}}
\newcommand{\K}{\mathrm{K}}
\newcommand{\R}{\mathrm{R}}
\newcommand{\W}{\mathrm{W}}
\newcommand{\X}{\mathrm{X}}
\newcommand{\Y}{\mathrm{Y}}
\newcommand{\Z}{\mathrm{Z}}

\newcommand{\mfree}{\text{\footnotesize{$\circ$}}}

\DeclareMathAlphabet\gothic{U}{euf}{m}{n}

\renewcommand{\AA}{\mathbb{A}}
\newcommand{\RR}{\mathbb{R}}

\newcommand{\car}{\mathbbm{1}}
\newcommand{\NN}{\mathbb{N}}
\renewcommand{\SS}{\mathbb{S}}
\newcommand{\QQ}{\mathbb{Q}}
\newcommand{\CC}{\mathbb{C}}
\newcommand{\ZZ}{\mathbb{Z}}

\newcommand{\PP}{\mathbb{P}}

\newcommand{\e}{\varepsilon}
\newcommand{\less}{\lesssim}
\newcommand{\llb}{\llbracket}
\newcommand{\rrb}{\rrbracket}

\newcommand{\ff}{\mathbf{f}}\newcommand{\FF}{\mathbf{F}}
\renewcommand{\gg}{\mathbf{g}}
\newcommand{\hh}{\mathbf{h}}
\newcommand{\uu}{\mathbf{u}}
\newcommand{\vv}{\mathbf{v}}
\newcommand{\ww}{\mathbf{w}}

\renewcommand{\cal}{\mathcal}
		
      \newcommand{\BesSmo}{{\cal B}}		
		
      \newcommand{\Distrib}{{\cal D'}} 
      \newcommand{\F}{{\cal F}}
      \newcommand{\ExtDiv}{{{\cal E}_\sigma}}\newcommand{\Ext}{{{\cal E}}}

\newcommand{\cp}{{\cal P}}
      \renewcommand{\S}{{\cal S}}

\let\div\undefined
\DeclareMathOperator{\div}{div\,}

\DeclareMathOperator{\curl}{curl\,}
\newcommand{\iprod}{\mathbin{\lrcorner}}

\DeclareMathOperator{\supp}{supp\,}
\DeclareMathOperator{\sgn}{sgn\,}

\counterwithin{equation}{section}

\setlist[enumerate,1]{label={\textit{(\roman*)}}}

\title{Optimal regularity results for the Stokes--Dirichlet problem\thanks{
\textbf{keywords:} \textit{Evolutionary Stokes problem, Stokes resolvent problem, Maximal regularity, Rough bounded domains, Endpoint (homogeneous) Besov spaces}}}
\author{
 Dominic \textsc{Breit}\thanks{Universität Duisburg-Essen, Fakultät für Mathematik, Thea-Leymann-Str. 9, 45127 Essen, Germany\\ \textbf{e-mail:} dominic.breit@uni-due.de} 
  \and
  Anatole \textsc{Gaudin}\thanks{Universität Duisburg-Essen, Fakultät für Mathematik, Thea-Leymann-Str. 9, 45127 Essen, Germany\\ \textbf{e-mail:} anatole.gaudin@uni-due.de}
}
\date{}

\begin{document}

\maketitle

 \begin{abstract}
We develop a sharp maximal regularity theory for the resolvent and evolution Stokes equations with no-slip boundary conditions, focusing on bounded domains of low regularity. Our framework covers the full scales of Besov and Sobolev spaces, $\B^s_{p,q}$ and $\H^{s,p}$, including endpoint cases such as $\L^\infty$. In particular, for the evolution problem, we show that in many situations the spatial regularity of $\partial_t \uu$, $-\Delta \uu +\nabla\mathfrak{p}$ in $\B^s_{p,q}$ can be controlled by the forcing term $\ff$ in the same space with suitable time integrability, but also provide the exact amount of regularity for the boundary that allows to control separately  $\nabla^2 \uu$ and $\nabla \mathfrak{p}$. Our approach also allows extending the classical $\L^p$-theory for $1\leqslant p\leqslant\infty$, giving a complete picture that includes both Bessel potential spaces $\H^{s,p}$ and Besov spaces $\B^s_{p,q}$, $p,q\in[1,\infty]$.

Our first main result establishes resolvent estimates in the half-space encompassing endpoint function spaces, while the second addresses bounded domains of minimal boundary regularity. In both cases we derive resolvent bounds, prove boundedness of the $\mathbf{H}^\infty$-functional calculus for the Stokes--Dirichlet operator, and characterize precisely the domains of its fractional powers.

In the half space setting, we work with homogeneous Sobolev and Besov spaces following the notion due to Bahouri, Chemin and Danchin, further refined by the second author. These spaces possess well-defined product laws, making them suitable for applications to the Navier--Stokes equations and other boundary value problems.  Yet, we are due to  analyze solenoidal subspaces and properties of the plain Laplacian in this framework. The analysis of solenoidal function spaces provides here a complete toolkit for the study of incompressible fluid flows. The resolvent estimate is then derived via Uka\"i-type formulas, which yield new insights even in the classical half-space case. Nevertheless, new arguments are needed for the endpoint spaces, motivating the approach by homogeneous function spaces. As a consequence of our analysis, we obtain an explicit description for the Stokes--Dirichlet operator on $\L^\infty(\RR^n_+)$, which seems completely new.

For bounded domains, our method combines standard localization and flattening with the half-space estimates. We obtain sharp results for a wide class of rough domains under minimal assumptions on boundary regularity. To this end, we rely on Sobolev multiplier theory. The assumptions coincide with those of Maz'ya--Shaposhnikova, already shown to be optimal in the case of the Laplace equation with Dirichlet boundary conditions. Specifically, we prove resolvent estimates in $\B^s_{p,q}$ and $\H^{s,p}$ provided the boundary charts have a small multiplier norm on the trace space $\B^{s+2-1/p}_{p,q}$. As a special case of our general strategy, we establish a full $\L^p$--theory for bounded $\C^{1,\alpha}$-domains for all $p\in(1,\infty)$, $\alpha>0$. In this case, solutions lie nearly in $\H^{1+\alpha+1/p,p}(\Omega)$. Moreover, for such domains, the resolvent problem in $\L^\infty(\Omega)$ produces solutions in $\C^{1,\alpha}(\overline{\Omega})$. Finally, we also obtain an $\L^1$-type theory in the Sobolev spaces $\W^{s,1}$.
\end{abstract}

\tableofcontents

\section{Introduction}
The incompressible Stokes equations
\begin{align}\label{stokes}
\partial_t\uu -\Delta\uu +\nabla\mathfrak p=\ff ,\quad\div \uu =0,
\end{align}
 are undoubtedly one of the most fundamental problems in fluid mechanics. Although their applicability to real-world problems is limited (as nonlinear effects are neglected), understanding \eqref{stokes} is the first step in the analysis of the more realistic Navier--Stokes equations.
Depending on the topology of the underlying domain and the functional analytic set-up
several important aspects remain open as of today.

\subsection{Historical overview and first motivations}\label{sec:intro1}
Suppose that \eqref{stokes} is posed in a domain $\Omega\subset\RR^n$ (with $n\geqslant 2$) with boundary $\partial\Omega$, supplemented with no-slip boundary conditions
\begin{align}\label{stokesboundary}
\uu_{|_{\partial\Omega}}=0
\end{align}
and initial datum $\uu (0)=0$. The first step is to consider the case of the half-space\footnote{We ignore the case $\Omega=\RR^n$ as the problem can be reduced to the heat equation.}
\begin{align*}
\RR_+^n:=\{x=(x_1,\dots,x_n)\in\RR^n:\,\,x_n>0\}.
\end{align*}
Hereby
\begin{align*}
\partial\RR_+^n=\{x=(x_1,\dots,x_n)\in\RR^n:\,\,x_n=0\}
\end{align*}
is the simplest possible boundary, being flat and arbitrarily smooth. However, it is unbounded and thus requires the need to work with homogeneous function spaces (we comment in details on this in Section~\ref{sec:intro2} below).
In that case, it is possible to derive a representation formula for the solution of \eqref{stokes}--\eqref{stokesboundary} by mean of Fourier transform (with respect to the first $n-1$ variables) and prove, whenever  $\uu(0)=0$,
\begin{align}\label{est1}
\|(\partial_t\uu,\nabla^2\uu,\nabla\mathfrak p)\|_{\L^q(\RR_+,\L^p(\RR_+^n))}\lesssim_{p,q,n}\,\|\ff\|_{\L^q(\RR_+,\L^p(\RR_+^n))}
\end{align}
for $p,q\in(1,\infty)$ with the help of some tools from harmonic analysis (continuity of singular integrals of Calder{\`o}n-Zygmund type). Such an estimate \eqref{est1} implies
\begin{align*}
\|(\partial_t\uu,-\Delta \uu +\nabla\mathfrak p)\|_{\L^q(\RR_+,\L^p(\RR_+^n))}\lesssim_{p,q,n}\,\|\ff\|_{\L^q(\RR_+,\L^p(\RR_+^n))}
\end{align*}
for the Stokes--Dirichlet problem, so for $A\uu := \AA_\mathcal{D} \uu := \PP_{\RR^{n}_+}(-\Delta \uu) = -\Delta \uu + \nabla\mathfrak{p}$,  with $\PP_{\RR^{n}_+}$ to be the Leray projection onto solenoidal vector fields, for $\X=\L^p(\RR^n_{+})$,
\begin{align}\label{est1-2}
\|(\partial_t\uu,A\uu)\|_{\L^q(\RR_+,\X)}\lesssim_{q,\X,A}\,\|\ff\|_{\L^q(\RR_+,\X)}.
\end{align}
The estimate \eqref{est1} was first proven by Solonnikov \cite[Theorem 3.1]{Solonnikov1973}, when $p=q$.\footnote{A similar strategy to the one exhibited by Solonnikov can be applied in the easier elliptic case. First results are due to Cattabriga \cite{Cattabriga1961}.
A comprehensive presentation can be found in Galdi \cite[Chapter IV]{bookGaldi2011}.}
A simplified approach can be derived from the work by Uka\"{i} \cite{Ukai1987} (we will also benefit from Uka\"{i}'s formula in Chapter~\ref{Sec:StokesHalfSpace}). It has been obtained that the solution can be written as a decomposition of the heat kernel (the solution operator to the heat equation in the half space) and Riesz transforms. Since the latter are well-known to be continuous on $\L^p$, $1<p<\infty$, the analysis is thereby essentially reduced to the heat equation.  An estimate such as \textbf{\eqref{est1-2}} for the (elliptic) operator $A$ is called a $\L^q(\X)$\textbf{--maximal regularity estimate}\footnote{ For $\X=\L^p$, sometimes, one writes  $\L^q(\L^p)$, $\L^q_t(\L^p_x)$ or $\L^q_t(\L^p(\Omega))$, when one wants to make clear the integrability exponents and the domain for the time and space variables.} (of parabolic type), and $\AA_\mathcal{D}$ is the said Stokes--Dirichlet operator.

At the beginning of the century, a similar idea to recover a full $\L^p$-theory by Desch, Hieber and Pr{\"u}ss \cite{DeshHieberPruss2001} has been successfully achieved. They provided a similar but different representation formula, allowing even to prove uniform boundedness and analyticity of the corresponding semigroup $(e^{-t\AA_\mathcal{D}})_{t>0}$ on $\L^\infty(\RR^n_+)$, \cite[Section~4]{DeshHieberPruss2001}. However, no precise description of the generator on $\L^\infty(\RR^n_+)$, \textit{i.e.}, neither an explicit description for the domain of $\AA_\mathcal{D}$, nor an explicit formula for $\mathbb{A}_\mathcal{D}$ is available, in particular, there is no given precise sense for which the underlying equation is actually solved.

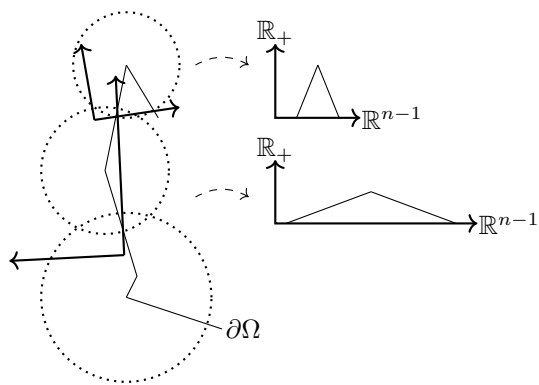
\begin{figure}[H]
\begin{center}
\begin{tikzpicture}[scale=1.4]
  \begin{scope}
    \draw plot [ tension=0.5] coordinates { (0.1,-1.2) (0.2,-1)
       (-.1,0) (0.1,1)};
       \draw (0.1,1) -- (0.4,0.5);
        \draw (0.1,-1.2) -- (1,-1.5);
      \draw[thick,dotted] (0.1,-1.2) circle (0.8);
       \draw[thick,dotted] (-.1,0) circle (0.6);
        \draw[thick,dotted] (0.1,1) circle (0.5);
        \draw [thick,<-] (0.6,0.6) -- (-0.2,0.48);
        \draw [thick,->] (-0.2,0.48) -- (-0.325,1.2);
        \draw [thick, <->] (1.5,1.2) -- (1.5,0.5) -- (2.3,0.5); 
        \draw [thick, <->] (1.5,0.1) -- (1.5,-0.5) -- (3.4,-0.5);
        \node at (1.5,1.3) {$\RR_+$};
        \node at (1.5,0.2) {$\RR_+$};
        \node at (3.7,-0.5) {$\RR^{n-1}$};
        \node at (2.6,0.5) {$\RR^{n-1}$};
        \draw (1.6,-0.5) -- (2.4,-0.2); 
        \draw (2.4,-0.2) -- (3.2,-0.5);
        \draw (1.7,0.5) -- (1.9,1);
         \draw (1.9,1) -- (2.1,0.5);
        \draw[thick,<-] (0,0.9) -- (0.08,-0.8);         
        \draw[thick,->] (0.08,-0.8) -- (-1,-0.855);
         \draw [<-,dashed] (1.25,1) to [out=150,in=30] (0.75,1);
          \draw [<-,dashed] (1.25,-0.25) to [out=150,in=30] (0.75,-0.25);
          \node at (1.2,-1.5) {$\partial\Omega$};
  \end{scope}
\end{tikzpicture}
\caption{This figure shows the local Lipschitz parametrisations of the boundary $\partial\Omega$ of some Lipschitz domain $\Omega$. The boundary can be covered neighbourhoods such that their intersection with $\partial\Omega$ arises -- after translating and rotating -- as the graph of a Lipschitz function.}\label{fig:1}
\end{center}
\end{figure}

With an estimate for the half-space at hand, one can aim for an analogous result such as \eqref{est1-2} for bounded domains $\Omega$ with curved boundaries, instead of $\RR^n_+$. This is typically based on introducing local coordinates
and flattening of the boundary, as can be seen in Figure~\ref{fig:1}. This procedure crucially hinges on the regularity of the charts
(the regularity of the underlying domain is defined via that of the local charts, see Definition~\ref{def:besovboundary}). As it turns out the most basic requirement is a Lipschitz boundary (and even under the assumption of a $\C^1$-boundary, there are counterexamples to regularity in the elliptic case, see \cite[Theorems~A~\&~1.2]{JerisonKenig1995} for the Laplace equation subject to Dirichlet boundary condition). 

 A classical result due to Solonnikov, \textit{cf.}~\cite[Theorem 4.1]{Solonnikov1973}, states that 
\begin{align}\label{est2}
\|(\partial_t\uu ,\nabla^2\uu ,\nabla\mathfrak p)\|_{\L^p(\RR_+\times\Omega)}\lesssim_{p,n,\Omega}\,\|\ff \|_{\L^p(\RR_+\times\Omega)}
\end{align}
for $p\in(1,\infty)$ provided $\Omega$ is a bounded domain with boundary $\partial\Omega$ of class $\C^2$. In this case, one can easily pass from the ``homogeneous'' to ``non-homogeneous'', and \textit{vice versa}, for expressions on the left-hand side. Indeed, by Poincar\'e-type inequalities one can easily control $\uu $ and $\nabla\uu $ (in view of \eqref{stokes}) as well as $\mathfrak p$ (under the normalisation condition $\int_\Omega \mathfrak p\,\d x=0$).

Estimates like \eqref{est2} for linearized Navier-Stokes-type systems  have been applied successfully in many different contexts to derive regularity, existence and uniqueness\footnote{having either that the time of existence is small, or that the initial datum is sufficiently small for a suitable norm.} of a wide variety of such non-linear systems on smooth domains. See for instance, and this is far from being exhaustive, \cite{GigaSohr1991,MuloneSolonnikov1995,Monniaux1999,Saal2007,FarwigKozonoSohr2008,MitreaMonniaux2009-1,Monniaux2013,Monniaux2021,HieberMonniaux2013,HieberNesensohnPrussSchade2016,Hieber2020}.

\medbreak

We refer to the work of Hieber and Saal \cite{HieberSaal2018} for an exhaustive review of the Stokes--Dirichlet problem and its related functional analytic properties in the $\L^p$-framework, $1\leqslant p\leqslant \infty$, on smooth enough domains $\Omega$.

\medbreak

In the following we enlighten two central points to motivate our approach by Besov spaces. We take the Dirichlet--Laplace problem as a reference for comparison:
\begin{itemize}
    \item \textbf{A direct study in the framework of Besov spaces is relevant:} In several situations, $\L^1$-in-time maximal regularity seems to be necessary in order to preserve certain quantities. This applies, in particular, to free boundary problems, see for instance the central works of Danchin and Mucha, Ogawa and Shimizu, or  Danchin, Hieber, Mucha and Tolksdorf \cite{DanchinMucha2009,DanchinMucha2015,OgawaShimizu2016,OgawaShimizu2021,OgawaShimizu2022,OgawaShimizu2024,DanchinHieberMuchaTolk2020}. The natural framework requires then to deal with $\L^q(\B^{s}_{p,q})$-maximal regularity for $q\in[1,\infty]$ and $s\in\RR$. Having a Besov space to describe the spatial regularity is necessary whenever $q=1$, see \cite[Theorems~4.5~\&~4.7]{BechtelBuiKunstmann2024}, \cite[Theorems~4.7~\&~4.15]{DaPratoGrisvard1975},  and \cite[Theorem~2.20]{DanchinHieberMuchaTolk2020} (notice the real interpolation identity $(\L^p,\W^{m,p})_{\theta,q} = \B^{m\theta}_{p,q}$, $m\in\NN$, $\theta\in(0,1)$, $p,q\in[1,\infty]$).
    
    Furthermore, when the spatial integrability is $p=1$ or $p=\infty$, one should directly study elliptic regularity on  $\L^1$ or $\L^\infty$. This appears to be an extremely tedious task, even for the Laplace problem subject to homogeneous Dirichlet boundary conditions on a smooth domain. To give a precise description for the regularity of the solution is already for $\Omega=\RR^n_+$ far from being straightforward. In this case the solution does not belong to $\W^{2,1}$ or $\W^{2,\infty}$, see, for instance, the standard counterexamples in \cite[Chapter~7,~Section~7.1.3]{GiaquintaMartinazzi}. Instead, one could directly study elliptic regularity at the level of Besov spaces $\B^{s}_{1,q}$ or $\B^{s}_{\infty,q}$, which are fitting better for the study of elliptic regularity. In this case (still having in mind the Dirichlet Laplace problem) a gain of exact two derivatives$^{5}$ for the solution (\textit{i.e.} membership to $\B^{s+2}_{1,q}$ or $\B^{s+2}_{\infty,q}$). Furthermore, such kind of results might not be directly reachable by real interpolation from the corresponding theories on $\L^1$ or $\L^\infty$. This particularly applies when the boundary of the underlying domain is not sufficiently smooth.
    
    Indeed, if we consider the Dirichlet Laplacian on $\L^\infty$ and $\Omega$ is a bounded $\C^{1,\alpha}$-domain, $\alpha\in(0,1)$, there is no information other than ``$-\Delta u =f$ and $f\in\L^\infty(\Omega)$ implies $u\in\C^{1,\alpha}(\overline{\Omega})$'', see \textit{e.g.} \cite[Theorem~5.21]{GiaquintaMartinazzi}. This sole information with the gap of $(1-\alpha)$ derivatives in the usual sense  makes it extremely difficult to compute directly the real interpolation space  $(\L^\infty(\Omega), \D_{\infty}(\Delta_\mathcal{D}))_{\theta,q}$ explicitly in terms of standard Besov spaces. Here $\D_{p}$ denotes the domain of an operator on $\L^p$, where $p\in[1,\infty]$.
    This real interpolation space should morally coincide with a closed subspace of $\B^{2\theta}_{\infty,q}(\Omega)$, whenever $2\theta < 1+\alpha$, but this cannot be deduced from the sole knowledge $\D_{\infty}(\Delta_\mathcal{D})\subset \C^{1,\alpha}(\overline{\Omega})$. A direct analysis of the regularity properties in Besov spaces, possibly of negative regularity, could avoide these specific issues. 

    \item \textbf{The $\L^p$-theory on domains not allowing for a gain of two full derivatives for the solution is of interest.}\footnote{By the gain of two derivatives we mean that the solution to an elliptic PDE enjoys two derivatives more than the right-hand side, for instance a forcing in $\L^2$ implies a solution in $\W^{2,2}$.} The standard approach for $\L^q_t(\L^p(\Omega))$--maximal regularity estimates -- firstly implemented in a systematic way  by Lady\v{z}enskaya, Solonnikov and Ural'tseva for generic parabolic problems of second order \cite{LadyhenskayaSolonnikovUraltseva68} -- is not suitable for domains with a boundary of low regularity. Indeed, this method requires a full gain of two derivatives on the solution with respect to the ground space $\L^p$. The approach based on a decomposition of the boundary requires at least two full derivatives on the boundary charts that describes $\partial\Omega$ (so morally speaking, dealing only with bounded $\C^{2}$ or $\C^{1,1}$-domains).

    \medbreak

    However, it has been proved for the Dirichlet--Laplace operator $\Delta_\mathcal{D}$ that 
    \begin{itemize}
        \item there exists a bounded $\C^{1}$-domain $\Omega$ such that for all $1< p<\infty$, one has for any $\varepsilon>0$, $\D_p(\Delta_\mathcal{D})\not\subset \W^{1+\frac{1}{p}+\varepsilon,p}(\Omega)$ while
        $\D_p(\Delta_\mathcal{D})\subset \W^{1+\frac{1}{p},p}(\Omega)$, for all $1<p<\infty$, see the work of Costabel \cite[Theorem~1.2]{Costabel2019};
        \item for $1<p<\infty$, there exists a bounded $\C^{1,1-\frac{1}{p}}$-domain such that $\D_p(\Delta_\mathcal{D})\not\subset \W^{2,p}(\Omega)$, see  the paper by Filonov \cite[Theorem~2.7,~p.247~\&~Remark,~p.255]{Filonov1997};
        \item for any $p,q\in(1,\infty)$, any bounded Lipschitz domain $\Omega$, the Dirichlet--Heat problem
        \begin{align*}
            \partial_t u- \Delta u = f,\quad u_{|_{\partial\Omega}}=0,\quad u(0)=0,
        \end{align*}
        admits a unique solution with the maximal regularity estimate
        \begin{align*}
            \lVert (\partial_t u,\Delta u)\rVert_{\L^q(\RR_+,\L^p(\Omega))} \lesssim_{p,q,\Omega}  \lVert f\rVert_{\L^q(\RR_+,\L^p(\Omega))},
        \end{align*}
        see \cite[Theorem~1.1]{Lamberton87} and \cite[Proposition~4.1]{CoulhonLamberton1986}, in combination with, \textit{e.g.}, \cite[Theorem~4.28]{bookOuhabaz2005}.
    \end{itemize}
    So the $\L^q_t(\L^p(\Omega))$--maximal regularity, $1<p,q<\infty$, is still available for the Dirichlet Laplacian although the boundary is too irregular to allow for a gain of two derivatives for the elliptic problem. At first glance, it seems reasonable to expect a similar behavior for the Stokes--Dirichlet problem on bounded domains of low regularity
    One should recall the strategy originally introduced by Solonnikov, Lady\v{z}enskaya and Ural'tseva for the Stokes problem in this context.

    \medbreak

    The reader might question the legitimacy of such dogged determination to reduce the regularity of the boundary of the domain, beyond its purely theoretical interest. The answer is straightforward: it ties to the description of concrete problems from fluid dynamics. One could, for instance, consider the depth of the sea or the roughness of the cliffs upon which the ocean breaks. Both are far from being described locally by smooth functions. A different and more precise class of problems can be found in the class of free-boundary problems, or more precisely fluids with moving boundaries. For this kind of problems, the functions that locally describe the graph of the boundary of the fluid domain are part of the unknowns of the problem. A particular class of equations  are the ones from  Fluid-Structure interactions, which aim at studying fluid flows around or inside an object (a ``body'') subject to deformations due to the dynamics of the fluid flow. Such equations from Fluid-Structure interactions are of particular interest, and have been widely studied in the literature. A well-studied topic is the motion of a rigid body immersed in a viscous incompressible fluid, see \cite{RB1,RB2,RB3,RB4,RB5} and the references therein. In this case, the position of the body describing (parts of) the boundary of the fluid domain is unknown and moves in time. It is, however, smooth. More delicate is the interaction with elastic structures located at the boundary of the fluid domain.\footnote{Such models have been suggested, for instance, in the study of blood vessels \cite{Quart1}.} The displacement of the structure determining the regularity of the boundary evolves according to some evolutionary problem and is of limited smoothness. Typically, weak solutions (and thus the boundary of the fluid domain) only belong to $\W^{2,2}$.
The global regularity of the velocity field becomes critical in this setting, even in the two-dimensional case, \textit{cf.} \cite{Breit2024FSI2D,GrHi2016,Shkoller2010}. Thus, in the scope of understanding such involved non-linear problems, the study of the linear (non-homogeneous) Stokes equations and related problems within a framework that goes below the standard $\C^2$-boundary regularity becomes crucial. 
\end{itemize}

\subsection{Homogeneous function spaces}\label{sec:intro2}

We now comment on homogeneous function spaces
suitable for the study of PDEs in unbounded domains following the discussion in \cite[Introduction]{Gaudin2023Lip}.

If partial differential equations are considered in unbounded domains, such as $\RR^n$ or $\RR^n_+$, where generally an inequality of Poincar\'e-type is not available, one can only control certain derivatives of the solution but not the solution itself.
Homogeneous Sobolev and Besov spaces can naturally accommodate such a behavior. 
Consider the Laplace equation in the whole space, \textit{i.e.}, 
\begin{equation*}\tag{$\mathcal{L}$}
-\Delta u =f \text{ in } \RR^n.\label{eq:LapEqRn}
\end{equation*}
One cannot expect $u\in\L^2(\RR^n)$ for arbitrary $f\in\L^2(\RR^n)$\footnote{This also happens on $\L^p(\RR^n)$, $p\in(1,\infty)$.  This phenomenon is due to the Laplacian not being invertible on $\L^2(\RR^n)$. On $\L^2(\RR^n)$, one can check it by the lack of boundedness of the Fourier symbol $\xi\mapsto|\xi|^{-2}$. On $\L^p(\RR^n)$, it can be checked by contradiction using the Closed Graph Theorem and a dilation argument.  Invertibilty of the elliptic operator here is the key point leading to the use of homogeneous function spaces:  on bounded domains, morally, Poincaré-types inequalities provide invertibility if one is able to control the gradient of the solution.}. Indeed, the only control we obtain is for the (semi-)norm defined by $$\lVert u\rVert_{\dot{\H}^{2,2}(\RR^n)} := \lVert\nabla^2u\rVert_{\L^2(\RR^n)} = \lVert f \rVert_{\L^2(\RR^n)}.$$
Thus, it appears that we need function spaces whose norms precisely control the regularity and integrability arising from the equations.

In the context of non-linear partial differential equations (even on the whole space), one has to be clear about the definition of homogeneous function spaces. Typically, elements of homogeneous function spaces, such as homogeneous Sobolev and Besov spaces, are defined as equivalence classes of tempered distributions modulo polynomials, denoted by ${\S'(\RR^n)}/{\mathcal{P}(\RR^n)}$, as in \cite[Chapter~6,~Section~6.3]{BerghLofstrom1976}, \cite[Chapter~5]{bookTriebel1983}, or \cite[Chapter~2]{bookSawano2018}.
This construction, modulo polynomials, proves unsuitable for nonlinear partial differential equations for several reasons.
In particular, according to \cite{Bourdaud1988, Bourdaud2013}, \cite[Chapter~2,~Section~2.4]{Triebel2015}, and \cite[Chapter~2,~Section~2.4.3]{bookSawano2018}, it is not clear that one can canonically choose the polynomial part to obtain an element of $\S'(\RR^n)$ in a way independent of the Sobolev index of the function spaces considered.
Furthermore,

\begin{itemize}
    \item There are issues with the definition of (para-)products:
    Given two equivalence classes $[u],[v]\in{\S'(\RR^n)}/{\mathcal{P}(\RR^n)}$ with representatives $u+P, v+Q, u+\Tilde{P}, v+\Tilde{Q} \in\S'(\RR^n)$,, we have
\begin{align*}
(u+P) (v+Q) - (u+\Tilde{P}) (v+\Tilde{Q}) =& (P- \Tilde{P})v + (Q-\Tilde{Q})u+ P Q - \Tilde{P}\Tilde{Q}.
\end{align*}
Even with meaningful bilinear (para)product estimates, $(P-\Tilde{P})v + (Q-\Tilde{Q})u$ is not a polynomial in general. This means the product ultimately depends on the choice of representatives. An issue which persists even when $P,Q,\Tilde{P},\Tilde{Q}$ are constants.
\item There are issues with the definition on domains by restriction:
    Defining the restriction of an element that belongs to ${\S'(\RR^n)}/{\mathcal{P}(\RR^n)}$ to a domain $\Omega$ seems ambiguous. For instance, consider $[f]\in {\S'(\RR)}/{\mathcal{P}(\RR)}$ such that has at least one representative, $f\in\L^1_{\text{loc}}(\RR)\cap\S'(\RR)$
    Then $x\mapsto f(x)-x^2$ and $x\mapsto f(x)-x^3$ are two distinct representatives that admit different restrictions in $\mathcal{D}'((a,b))$ for any $-\infty<a<b\leqslant \infty$. Thus, the restriction in the distributional sense seems meaningless with regard to the quotient structure. In particular, it appears unclear on how to define the support, having in mind the quotient structure.

    \item There are issues with the composition with a diffeomorphism:
    Assuming that $u+P, u+Q\in\S'(\RR^n)$ are two representatives of $[u]\in {\S'(\RR^n)}/{\mathcal{P}(\RR^n)}$, and $\Psi$ is a smooth diffeomorphism of $\RR^n$, the meaning of the expression
\begin{align*}
    u\circ \Psi + P\circ\Psi- (u\circ \Psi + Q\circ\Psi) = P\circ\Psi - Q\circ\Psi
\end{align*}
is ambiguous. Even if we assume that $u\circ \Psi\in\S'(\RR^n)$, the objects $P\circ\Psi$ and $Q\circ\Psi$ might not even qualify as tempered distributions. If it was, it should not depend on the choice of $P$ and $Q$ which is again unclear. This creates a significant challenge in transferring properties of homogeneous function spaces from the whole space or half-space to a bent half-space via a change of coordinates, particularly in defining traces on the boundary.

    \item There are issues concerning distribution theory and the completion of function spaces:
    The completion of any function space should be canonically identified as subspace of the distributions $\mathcal{D}'(\RR^n)$. However, it can be verified that
    \begin{align*}
    \mathrm{C}^\infty_c(\RR^n)\not\subset \dot{\H}^{-{n}/{2},2}(\RR^n),
    \end{align*}
    by considering a smooth compactly supported and non-negative bump function $\varphi\in\Ccinfty(\RR^n)$, such that $\int_{\RR^n} \varphi =1$.
    Consequently, the completion of Schwartz functions for the $\dot{\H}^{{n}/{2},2}(\RR^n)$-norm cannot be canonically embedded in $\mathcal{D}'(\RR^n)$, as the (pre-)dual space lacks test functions.
    This implies that if one takes any abstract completion to define these spaces, one would encounter elements that cannot be canonically identified with actual distributions.

    \item There are issues concerning localization, and stability by multiplication by smooth functions.
    A consequence of the previous points is that as soon as $s\leqslant - \sfrac{n}{p'}$, one loses stability under multiplication by smooth and compactly supported functions, that is we have
    \begin{align*}
    \mathrm{C}^\infty_c(\RR^n)\cdot \dot{\H}^{s,p}(\RR^n) \not\subset \dot{\H}^{s,p}(\RR^n).
    \end{align*}
    Indeed, let us give a counterexample in the simplified case $p=n=2$, $s=-1$. For $\varphi\in\Ccinfty(\RR^2)$, supported in the unit ball, as in the previous point, one has $\varphi\notin\dot{\H}^{-1,2}(\RR^2)$. However it can be checked, provided one has $a\in\RR^2$ with large Euclidean norm, that $$\varphi(\cdot - a) - \varphi(\cdot + a) \in\dot{\H}^{-1,2}(\RR^2).$$ Consequently multiplication by a bump function $\Psi$ with value $1$ on $\B_1(a)$, yields
    \begin{align*}
        \Psi[\varphi(\cdot - a) - \varphi(\cdot + a)] = \varphi(\cdot - a)\notin \dot{\H}^{-1,2}(\RR^2).
    \end{align*}
    Therefore, the naive localisation procedure proves to be inapplicable as soon as homogeneous functions of too negative regularity are involved.
\end{itemize}
In particular, the realization of homogeneous function spaces on special Lipschitz domains provided by Costabel, M${}^\text{c}$Intosh and Taggart \cite{CostabelTaggartMcIntosh2013}, built on ${\S'(\RR^n)}/{\mathcal{P}(\RR^n)}$, appears to be inapplicable for linear problems with boundary values and unsuitable for non-linear problems.

To circumvent the issues mentioned above related to a choice of a representative, Bahouri, Chemin, and Danchin proposed in \cite[Chapter~2]{bookBahouriCheminDanchin} to consider a subspace of $\S'(\RR^n)$ consisting solely of tempered distributions without any (non-zero) polynomial part, see \cite[Examples,~p.23]{bookBahouriCheminDanchin}. This subspace, denoted $\S'_h(\RR^n)$, is defined at the beginning of Subsection \ref{subsec:SobolevBesovRn}. Using $\S'_h(\RR^n)$ as an ambient space, Bahouri, Chemin, and Danchin constructed homogeneous Besov spaces $\dot{\mathrm{B}}^{s}_{p,q}(\RR^n)$. With their construction, the homogeneous Besov spaces are complete whenever $(s,p,q)\in\RR\times[1,\infty]^2$ satisfies
\begin{align*}\tag{$\mathcal{C}_{s,p,q}$}
    \left[ s<\frac{n}{p} \right]\text{ or }\left[q=1\text{ and } s\leqslant\frac{n}{p} \right]\text{. }
\end{align*}
The approach from \cite[Chapter~2]{bookBahouriCheminDanchin} has turned out fruitful in the analysis of (nonlinear) partial differential equations in unbounded domains. The lack of completeness in some cases is nowadays
an accepted sacrifice for restoring the definition of para-products, and a natural distribution theory. Unfortunately, it is not possible to construct homogeneous function space with all three properties (completeness, definition of para-products and distribution theory).

The work of Bahouri, Chemin and Danchin has been further extended by Danchin and Mucha to homogeneous Besov spaces on $\RR^n_+$, \textit{c.f.}~\cite{DanchinMucha2009,DanchinMucha2015}. Danchin, Hieber, Mucha, and Tolksdorf examined briefly homogeneous Sobolev spaces $\dot{\H}^{m,p}$ on $\RR^n$ and $\RR^n_+$ for $m\in\mathbb{N}$, $p\in(1,\infty)$, in \cite[Chapter~2]{DanchinHieberMuchaTolk2020} and an additional analysis of homogeneous Besov spaces on $\RR^n_+$. 
The second author extended this construction in a previous work \cite{Gaudin2022}, encompassing the entire scale of homogeneous Sobolev spaces, in the reflexive range, on the half-space $\RR^n_+$. He also extended widely this analysis to special Lipschitz domains $\Omega$ in \cite{Gaudin2023Lip}, encompassing the endpoint cases $p=1,\infty$, and removing the requirements of completeness for the central results to hold. These extensions were achieved by investigating interpolation and density properties, as well as the meaning of traces on the boundary.

\subsection{Domains with minimal regularity}\label{sec:intro3}

At the end of Section~\ref{sec:intro1}, it has been mentioned that one would like to lower the regularity assumptions on the boundary of the underlying domain. In order, to solve uniquely the corresponding evolution problem,  one needs first to investigate the regularity property of the corresponding purely elliptic problem. In this regard, it has been mentionned that even if one wants to reach a gain of two exact derivatives,  minimal regularity requirements are needed and $\C^{1,1-\sfrac{1}{p}}$-regularity for the boundary is not enough to reach $\W^{2,p}$-regularity on $\L^p$ for the Dirichlet Laplacian. What would be the minimal requirements for such a statement ? This will lead us to the necessary use of domains $\Omega$ within an appropriate class described by the multiplier theory by Maz'ya and Shaposhnikova \cite{MazyaShaposhnikova2009}, which turns to be a necessary condition in the case of the Dirichlet Laplacian. We highlight below several insights that indicate that the membership to the multiplier class is necessary for the regularity of the Dirichlet Laplacian. We then mention several results obtained from this framework for the Stokes--Dirichlet problem.

\subsubsection{Illustration by the Dirichlet Laplacian}\label{sec:intro3DirLap}
As alluded to above at the end of Section~\ref{sec:intro1}, the question of maximal (elliptic) regularity in bounded domains becomes more involved when domains of minimal regularity are considered. Corresponding results for the Laplace equation
\begin{align}\label{laplace}
-\Delta u=f,\quad u_{|_{\partial\Omega}}=0,
\end{align}
 are classical and it has been shown that the solution satisfies
\begin{align}\label{laplace2}
\|u\|_{\W^{2,p}(\Omega)}\lesssim_{p,n,\Omega} \|f\|_{\L^p(\Omega)}
\end{align}
provided $\Omega$ is a Lipschitz domain and the local boundary charts belong to the class of Sobolev multipliers $\mathcal{M}_{\W}^{2-1/p,p}(\RR^{n-1})$ on the trace space $\W^{2-1/p,p}(\partial\Omega)$ of $\W^{2,p}(\Omega)$ (analogous results hold for higher order and fractional differentiabilities) and have a small norm.
Roughly speaking the space $\mathcal{M}_{\W}^{s,p}(\RR^{n})$ consists of all functions $\varphi\,:\,\RR^{n}\longrightarrow\CC$
such that the mapping
\begin{align*}
\W^{s-1,p}(\RR^{n})\ni v\longmapsto v \nabla\varphi\in \W^{s-1,p}(\RR^{n})^n
\end{align*}
is bounded, see Section~\ref{sec:SM} for a proper introduction.
If $p$ is sufficiently large, e.g. $p>n$ such that  $\W^{1-1/p,p}(\RR^{n-1})$ is an algebra, it holds
$\mathcal{M}_{\W}^{2-1/p,p}(\RR^{n-1})\simeq \W^{2-1/p,p}_{\mathrm{unif}}(\RR^{n-1})$ where 
\begin{align}
    \lVert v\rVert_{\X_{\mathrm{unif}}(\RR^{n})} := \sup_{ x_0\in\RR^n} \lVert \zeta(\cdot-x_0) v \rVert_{\X(\RR^{n})}, \label{eq:introUnifSpaces}
\end{align}
for a given $\zeta\in\Ccinfty(\RR^n)$, such that $0\leqslant \zeta\leqslant 1$, $\zeta_{|_{\B_1(0)}}=1$ and $\zeta_{|_{\B_2(0)^c}}=0$.

A comprehensive picture is given in \cite[Chapters~3,~4,~6,~8,~9~\&~14]{MazyaShaposhnikova2009}. The sharpness of this assumption for \eqref{laplace2}, giving the $\W^{2,p}$-regularity on $\L^p(\Omega)$ for the problem \eqref{laplace}, is given in \cite[Chapter~14,~p.563]{MazyaShaposhnikova2009}.
 If $p>n$ the argument is elementary: We make the particular choice
\begin{align*}
\Omega=\{(x,y)\in\RR^2:\,y >\phi(x)\}\cap\B_1(0) \subset\RR^2
\end{align*}
for a given Lipschitz function $\phi:\RR\rightarrow\RR$. We suppose that $\phi$ is smooth outside a neighborhood of the origin say outside $\B_{1/2}(0)$. By the Riemannian mapping theorem, there exists a holomorphic function\footnote{Here we denote $\B_1^+(0):=\B_1(0)\cap \RR_+^n$.}
\begin{align*}
\mathcal{G}:\Omega\rightarrow \RR^2_+\cap \B_1^+(0)
\end{align*}
with fixed point 0.
Setting $w:=\mathrm{Im}(\mathcal{G})$, we have $w=0$ on $\{(x,y)\in\RR^2:\,y =\phi(x)\}\cap\B_1(0)$. Consider $\eta\in\Ccinfty(\B_1(0))$, such that $\eta_{|_{\B_{1/2}(0)}}=1$ and set $u:=\eta w$.
Suppose now that $u\in \W^{2,p}(\Omega)$ for some $p>1$, which implies that $w\in \W^{2,p}(\Omega\cap\B_{1/2}(0))$. then it holds that $u$ solves \eqref{laplace} with 
\begin{align*}
    f:=-\Delta \eta w -2 \nabla \eta\cdot\nabla w .
\end{align*}
Since $f$ is supported outside $\B_{1/2}(0)$, where the boundary of $\Omega$ smooth, and $w$ is harmonic with $w=0$ on $\{y =\phi(x)\}$ the function $f$ is continuous. In particular, we have $f\in \L^p(\Omega)$ for any $p<\infty$.

Since $w$ has zero boundary values by construction, the tangential derivative vanishes as well (see \cite[Theorem~9.5.1,~eq.~(9.5.3)]{MazyaShaposhnikova2009} for a rigorous proof), \textit{i.e.},
\begin{align*}
\mathcal{T}_{\partial\Omega}\partial_x  (w\circ \boldsymbol{\Phi})+\mathcal{T}_{\partial\Omega}\,\partial_y (w\circ \boldsymbol{\Phi})\partial_x \phi=0\quad\text{on}\quad\{(x,y)\in\RR^2:\,y =\phi(x)\}\cap\B_1(0)\subset\partial\Omega.
\end{align*}
Here $\boldsymbol{\Phi}$ is an extension of $\phi$ (see \eqref{eq:Phi} below for details) and $\mathcal{T}_{\partial\Omega}$ the trace operator related to $\partial\Omega$.
Noticing that $\mathcal{T}_{\partial\Omega}\,\partial_y (w\circ \boldsymbol{\Phi})$ is strictly positive (by Hopf's maximum principle) this is equivalent to
\begin{align*}
\partial_x \phi=-\frac{\mathcal{T}_{\partial\Omega}\,\partial_x  (w\circ \boldsymbol{\Phi})}{\mathcal{T}_{\partial\Omega}\,\partial_y (w\circ \boldsymbol{\Phi})}\quad\text{on}\quad \{(x,y)\in\RR^2:\,y =\phi(x)\}\cap\B_1(0).
\end{align*}
For $n=2$ and $p>2$, we have that $\W^{1-1/p,p}(\partial\Omega)$ is an algebra which implies $\phi\in\W^{2-1/p,p}_{\mathrm{loc}}(\RR)$ as desired, note that this implies a small multiplier norm, \textit{i.e.},
\begin{align*}
    \lVert \phi\rVert_{ \mathcal{M}_{\W}^{2-1/p,p}(\RR)} = \lVert \phi'\rVert_{ \mathcal{M}_{\W}^{1-1/p,p}(\RR)}\sim_{p} \lVert \phi'\rVert_{ {\W}^{1-1/p,p}_{\mathrm{unif}}(\RR)}<\epsilon
\end{align*}
for some small $\epsilon>0$, see \cite[Corollary~4.3.8~\&~Section~14.6.3]{MazyaShaposhnikova2009}. A similar argument can be applied to other function spaces which form a multiplication algebra and allow for a corresponding trace embedding.

\subsubsection{Towards the Stokes--Dirichlet problem}

Going back to the (stationary and non-stationary) Stokes--Dirichlet  problem, a result in the spirit of \eqref{laplace2} (and under the same assumptions on $\partial\Omega$) has been shown very recently for the steady Stokes problem by the first author, see \cite[Section 3]{Breit2024FSI2D}. 
This is not only of theoretical interest. Fluid dynamical problems in domains with irregular
boundary appear quite naturally in applications. In fact, the analysis in \cite[Section~3]{Breit2024FSI2D} is motivated by problems in fluid-structure interaction already mentioned above. In \cite{Breit2024FSI2D} the boundary of the fluid domain describes a flexible elastic shell.

\medbreak

While the estimate in the steady case from \cite[Section 3]{Breit2024FSI2D} is proved in a more general framework including fractional Sobolev(--Soblodeckij) spaces of non-negative order, the corresponding result only deals with $\L^p$-spaces, $1<p<\infty$, and always requires a full gain of $2$ spatial derivatives for the strategy to apply. To be precise the first author proved in \cite[Section 3]{BreitSchauder} that \eqref{est2} holds provided
$\Omega$ is a Lipschitz domain and the local boundary charts belong to the class of Sobolev multipliers $\mathcal{M}_{\W}^{2-1/p,p}$ -- the same assumption which is already sharp for the Laplace equation. Additionally, this approach does not carry over analytic properties, such as, for instance, the $\mathbf{H}^\infty$-calculus, or a possible description for the domain of the fractional powers of the Stokes--Dirichlet operator, see the following section. 

\subsection{Precise framing of the main questions and the related state of the art}

We highlight in this section the main questions — Questions~\ref{quest:mainQ1}, \ref{quest:mainQ2}, and \ref{quest:mainQ3} below — that this work aims to answer, providing a contextual presentation and a review of the corresponding results achieved in the current literature.

\medbreak

The estimate \eqref{est1-2} implies necessarily that $(\D(A),A)$ is closable on $\X$ with the necessary estimate
\begin{align*}
    \lVert \lambda(\lambda\I+A)^{-1} f\rVert_{\X} \lesssim_{\X,A,\mu}\lVert f \rVert_{\X},\quad \forall f\in\X,\quad\forall \lambda\in\Sigma_{\mu}=\{z\in\CC^\ast\,:\,|\arg(z)|<\mu\},
\end{align*}
for some $\mu\in(\sfrac{\pi}{2},\pi)$. This
means that the operator $(\D(A),A)$ is $(\pi-\mu)$--\textbf{sectorial} on $\X$. This implies that one can construct the holomorphic and uniformly bounded semigroup $(e^{-tA})_{t\geqslant 0}$ on $\X$, with  $i\RR^\ast\subset\rho(A)$, see for instance \cite[Chapters~2~\&~3]{bookHaase2006} and \cite[Beginning of Section~2]{CoulhonLamberton1986} for more details. A summary of the theory of sectorial operators and maximal regularity is given in Section~\ref{sec:MaxRegIntro}. 

\medbreak

Thus for $A \uu=\AA_\mathcal{D}\uu=-\Delta\uu+\nabla\mathfrak{p}$, this naturally forces us to begin with the following Stokes--Dirichlet resolvent problem: For $\lambda\in\CC\setminus(-\infty,0]$ and a $\Omega$ domain with boundary at least Lipschitz continuous we consider:
\begin{equation*}\tag{DS${_\lambda}$}\label{eq:IntroSystStokes}
    \left\{ \begin{array}{rllr}
         \lambda \uu - \Delta \uu +\nabla \mathfrak{p} &= \ff \text{, }&&\text{ in } \Omega\text{,}\\
        \div \, \uu &= 0\text{, } &&\text{ in } \Omega\text{,}\\
        \uu_{|_{\partial\Omega}} &=0\text{, } &&\text{ on } \partial\Omega\text{.}
    \end{array}
    \right.
\end{equation*}
Given a function space $\X(\Omega)$\footnote{Typically, here $\X$ is an $\L^p$ space, a Sobolev space $\H^{s,p}$ or $\W^{s,p}$, or a Besov space $\B^{s}_{p,q}$.} and $\ff\in\X(\Omega)$ divergence--free with $\ff\cdot \nu =0$ on $\partial\Omega$, we ask for existence and uniqueness of a solution $(\uu,\nabla \mathfrak{p})$, with $\uu\in\X(\Omega)$, and $-\Delta \uu + \nabla \mathfrak{p}\in\X(\Omega)$\footnote{This implicitly means that $\nabla \mathfrak{p}$ should be uniquely determined from $\uu$ which should be uniquely determined by $\ff$.}. More precisely, we seek $\mu\in(\sfrac{\pi}{2},\pi)$ such that we have well-posedness for \eqref{eq:IntroSystStokes} for all $\lambda\in\Sigma_\mu$ and all divergence--free $\ff\in\X(\Omega)$ with vanishing normal component on the boundary. Moreover, one of the following estimates should hold:
\begin{itemize}
    \item the pure resolvent estimate
    \begin{align}
        |\lambda|\lVert \uu\rVert_{\X(\Omega)} \lesssim_{\X,\Omega,\mu} \lVert \ff\rVert_{\X(\Omega)}\label{eq:IntroPureResolvEst};
    \end{align}
    \item the gradient resolvent estimate
    \begin{align}
        |\lambda|\lVert \uu\rVert_{\X(\Omega)} + |\lambda|^\frac{1}{2}\lVert \nabla \uu\rVert_{\X(\Omega)} \lesssim_{\X,\Omega,\mu} \lVert \ff\rVert_{\X(\Omega)}\label{eq:IntroGradResolvEst};
    \end{align}
    \item the full regularity resolvent estimate
    \begin{align}
        |\lambda|\lVert \uu\rVert_{\X(\Omega)} + |\lambda|^\frac{1}{2}\lVert \nabla \uu\rVert_{\X(\Omega)}  + \lVert (\nabla^2 \uu,\nabla \mathfrak{p})\rVert_{\X(\Omega)} \lesssim_{\X,\Omega,\mu} \lVert \ff\rVert_{\X(\Omega)}\label{eq:IntroFullRegResolvEst}.
    \end{align}
\end{itemize}

\medbreak

\noindent In the case of bounded domains, the best known results  so far are as follows.
\begin{itemize}
    \item Estimate \eqref{eq:IntroGradResolvEst} holds for bounded Lipschitz domains $\Omega$: the case $\X=\H^{s,2}$, $|s|<\frac{1}{2}$ is studied by Mitrea and Monniaux \cite{MitreaMonniaux2008}, and the case $\X=\L^p$, for all $\frac{3}{2}-\varepsilon_\Omega<p<3+\varepsilon_\Omega$ for some small $\varepsilon_\Omega>0$, by Shen \cite{Shen2012}.
    \item It was shown by Solonnikov \cite{Solonnikov1977} that estimate \eqref{eq:IntroFullRegResolvEst} holds if $\Omega$ is a bounded $\C^2$/$\C^{1,1}$-domain and $\X=\L^p$, where $1<p<\infty$ is arbitrary.
    \item Estimate \eqref{eq:IntroGradResolvEst} holds for $\X=\L^p$, for all $p\in(1,\infty)$, and $\Omega$ is a bounded $\C^{1}$-domain. This is proved in the very recent work of Geng and Shen \cite{GengShen2024}.
    \item If $\X=\L^\infty$ and $\Omega$ is a bounded  $\C^{2}$-domain estimate \eqref{eq:IntroGradResolvEst} is obtained by   Abe, Giga and Hieber \cite{AbeGigaHieber2015}. Very recently it has been obtained by Geng and Shen \cite{GengShen2025} for bounded $\C^{1,\alpha}$-domains, where $\alpha>0$.
    \item The case $\X=\L^\infty$ with a bounded $\C^1$-domain $\Omega$, has recently been considered by Geng and Shen \cite{GengShen2025}. They proved \eqref{eq:IntroPureResolvEst}, which is sharp in the sense that  the estimate \eqref{eq:IntroGradResolvEst} is known to fail for the plain Dirichlet--Laplacian in this context.
\end{itemize}
For more details and historical background, we refer to the review by Hieber and Saal \cite[Sections~2.8~\&~4]{HieberSaal2018}.

For domains of intermediate regularity no precise regularity result seems to be known: If $\Omega$ is say a bounded $\C^{1,\alpha}$-domain, $\alpha\in(0,1)$ fixed but arbitrary, then \textit{a priori} \eqref{eq:IntroGradResolvEst} holds while \eqref{eq:IntroFullRegResolvEst} does not. Hence, one may wonder if we can provide any intermediate statement. For instance, the result by Giaquinta and Modica \cite[Part~II,~Theorem~1.3,~b)]{GiaquintaModica1982} shows  for $\ff\in\L^\infty(\Omega,\CC^n)$ (divergence-free with vanishing normal component on the boundary) that the solution to \eqref{eq:IntroSystStokes} satisfies $\uu\in\C^{1,\alpha}(\overline{\Omega},\CC^n)=\B^{1+\alpha}_{\infty,\infty}(\Omega,\CC^n)$. However, it is not known if we have the estimate
\begin{align*}
    |\lambda|^\frac{1-\alpha}{2}\lVert \uu \rVert_{\B^{1+\alpha}_{\infty,\infty}(\Omega)} \lesssim_{\mu,\alpha,n}^{\Omega} \lVert \ff \rVert_{\L^{\infty}(\Omega)}.
\end{align*}
Similarly, for $\X=\L^p$, $1<p<\infty$, we ask in complement to \eqref{eq:IntroGradResolvEst} for the estimate
\begin{align}
    |\lambda|^{1-\frac{s}{2}}\lVert \uu \rVert_{\H^{s,p}(\Omega)} \lesssim_{p,s,n}^{\mu,\alpha,\Omega} \lVert \ff \rVert_{\L^{p}(\Omega)}\qquad \text{or}\qquad|\lambda|^{1-\frac{s}{2}}\lVert \uu \rVert_{\W^{s,p}(\Omega)} \lesssim_{p,s,n}^{\mu,\alpha,\Omega} \lVert \ff \rVert_{\L^{p}(\Omega)}\label{eq:IntroFractIntermResolvEstStokes}
\end{align}
for some $s\in(1,2]$ depending possibly on $p$, $\alpha$ and $\Omega$.

The interest in \eqref{eq:IntroFractIntermResolvEstStokes} arises from the fact that it translates to the following semigroup decay estimate of heat-smoothing type:
\begin{align*}
    \lVert e^{-t\AA_\mathcal{D}}\rVert_{\L^{p}(\Omega)\rightarrow\H^{s,p}(\Omega)},\,\lVert e^{-t\AA_\mathcal{D}}\rVert_{\L^{p}(\Omega)\rightarrow\W^{s,p}(\Omega)} \lesssim_{p,s,n}^\Omega t^{-\frac{s}{2}},\quad t>0.
\end{align*}
This leads to the following questions.

\begin{question}\label{quest:mainQ1} For which function spaces $\X$ depending on $\Omega$ or for which regularity on $\Omega$ depending on $\X$ do we have either \eqref{eq:IntroGradResolvEst} or \eqref{eq:IntroFullRegResolvEst} for \eqref{eq:IntroSystStokes} on $\X(\Omega)$?

\noindent If \eqref{eq:IntroGradResolvEst} holds but \eqref{eq:IntroFullRegResolvEst} does not, can we still provide more precise and quantitative regularity results on $\uu$ ?
\end{question}

\begin{question}\label{quest:mainQ2} Let $\Omega$ be a bounded domain with a boundary regularity strictly rougher than $\C^2$. For either the semigroup or the resolvent operator arising from \eqref{eq:IntroSystStokes}, 
can we still provide precise results in the endpoint cases $\X=\L^{1},\L^\infty,\B^{s}_{1,q}$ or $\B^{s}_{\infty,q}$ ?
\end{question}

On account of the counterexample by Filonov \cite[Theorem~2.7,~p.247~\&~Remark,~p.255]{Filonov1997} for the Dirichlet Laplacian (as exhibited near the end of Section~\ref{sec:intro1}), for domains of class $\C^{1,1-\frac{1}{p}}$ or rougher, one cannot expect \eqref{eq:IntroFractIntermResolvEstStokes} to hold for $s=2$. Note that by \cite[Proposition~6.6.4]{bookHaase2006}, one has for all $s\in[1,2]$,
\begin{align*}
    |\lambda|^{1-\frac{s}{2}}\lVert \AA_\mathcal{D}^{\sfrac{s}{2}} \uu \rVert_{\L^{p}(\Omega)} \lesssim_{p,s,n}^{\mu,\alpha,\Omega} \lVert \ff \rVert_{\L^{p}(\Omega)}.
\end{align*}
Here, it would be sufficient to have the following description of the fractional powers on $\L^p(\Omega,\CC^n)$:
\begin{align}
    &\D_p(\AA_\mathcal{D}^{\sfrac{s}{2}} )=\{\,\vv\in\H^{s,p}(\Omega,\CC^n)\,:\,\div \vv =0,\quad \text{ and } \vv_{|_{\partial\Omega}}=0 \,\},\label{eq:IntroFracPowStokesDescrip}\\
    \quad\text{ with }\quad&\lVert\AA_\mathcal{D}^{\sfrac{s}{2}}\ww\rVert_{\L^{p}(\Omega)}\sim_{p,s,n}^\Omega\lVert\ww\rVert_{\H^{s,p}(\Omega)},\quad\forall \ww\in\D_p(\AA_\mathcal{D}^{\sfrac{s}{2}} ).\nonumber
\end{align}
Usually, in concrete cases this kind of description can be achieved by showing that (in addition to sectoriality on $\X$) the operator $\AA_\mathcal{D}$ also has bounded imaginary powers (BIP). This means that for $\X_{\mathfrak{n},\sigma}(\Omega)=\overline{\D(\AA_\mathcal{D})}^{\lVert\cdot\rVert_{\X}}$ we have that
\begin{align}
    (\AA_{\mathcal{D}}^{i\tau})_{\tau\in\RR} \subset \mathcal{L}(\X_{\mathfrak{n},\sigma}(\Omega)).\label{eq:BIPstokesIntroSec4}
\end{align}
Consequently, the domains of fractional powers on $\X$ built a complex interpolation scale, with
\begin{align*}
    [\X_{\mathfrak{n},\sigma}(\Omega), \D(\AA_{\mathcal{D}})]_{\theta} = \D(\AA_{\mathcal{D}}^\theta),\quad \theta\in(0,1).
\end{align*}
In this case, it suffices to characterize the left hand-side in terms of standard function spaces. Note that under suitable growth conditions on the operator norms with respect to $\tau$ and certain assumptions on the space $\X$, \eqref{eq:BIPstokesIntroSec4} implies  $\L^q(\X)$--maximal regularity for $\AA_\mathcal{D}$ on $\X$, \textit{i.e.} the estimate \eqref{est1-2} holds, by the Dore--Venni Theorem, Theorem~\ref{thm:DoreVennithm}. The BIP property \eqref{eq:BIPstokesIntroSec4} can be reduced to the knowledge of the so-called bounded holomorphic (or $\mathbf{H}^\infty$-)functional calculus
\begin{align*}
    \lVert f(\AA_\mathcal{D})\rVert_{\X_{\mathfrak{n},\sigma}\rightarrow\X} \lesssim_{n,\Omega,\omega}^\X\lVert f\rVert_{\L^\infty(\Sigma_\omega)} \quad\text{ with }\quad f(\AA_\mathcal{D}):=\frac{1}{2i\pi}\int_{\Sigma_{\omega-\varepsilon}} f(z)(z\I-\AA_\mathcal{D})^{-1} \d z.
\end{align*}
Here $f\,:\,\Sigma_\omega\longrightarrow\CC$ is any bounded holomorphic functions and $\omega,\omega-\varepsilon\in(\pi-\mu, \frac{\pi}{2})$, for  $\varepsilon>0$ small enough. Again, we refer to \textbf{Section~\ref{sec:MaxRegIntro}} for more details on BIP, the $\mathbf{H}^\infty$-calculus and their connection to $\L^q$--maximal regularity.

\medbreak

We mention here that the interest in the $\mathbf{H}^\infty$-functional calculus goes way beyond the sole interest in $\L^q$--maximal regularity \eqref{est1-2}. Indeed, for $\X=\L^p$, according to the milestone result of Van Nerveen, Veraar and Weis \cite[Theorem~1.1]{VanNerveenVeraarWeis}, the bounded $\mathbf{H}^\infty$-functional calculus of $(\D(A),A)$ on $\X$ also allows for the so-called \textbf{stochastic} $\L^q$--maximal regularity. The latter means more precisely that one has the following estimate for $p\in[2,\infty),\, q\in(2,\infty)$,
\begin{align*}
    \mathbb{E}\bigg\lVert t\mapsto\int_{0}^{t}A^{\sfrac{1}{2}}e^{-(t-s)A}\mathbf{G}(s)\,\d W_{\H}(s)\bigg\rVert_{\L^q(\RR_+,\L^p(\Omega))}^q \lesssim_{p,q,\Omega}\mathbb{E}\lVert \mathbf{G} \rVert_{\L^q(\RR_+,\L^p(\Omega,\H))}^q.
\end{align*}
Here $\mathbb{E}$ is the expectation, $\H$ a Hilbert space and $W_{\H}$ is an $\H$-valued cylindrical Wiener process. Our current work does not emphasis on this consequence. However, this should be kept in mind as one of the central motivations for the investigation of the $\mathbf{H}^\infty$-functional calculus. In the case of the Stokes operator, this is of a paramount importance for the  stochastic analysis of incompressible viscous fluids, see for instance \cite{AgrestiVeraar2024} in the case of the multidimensional torus. For an exhaustive overview on Stochastic PDEs and the related maximal regularity theory, one could consult the recent survey by Agresti and Veraar \cite{AgrestiVeraarSurveySPDE2025} and the references therein. A purely practical reason for investigating the boundedness of $\mathbf{H}^\infty$-functional calculus rather than BIP is the ability of being stable under perturbation. This is the case for the property of the bounded  $\mathbf{H}^\infty$-functional calculus, while it usually does not apply for BIP. Hence we ask:

\begin{question}\label{quest:mainQ3} For which function spaces $\X$, possibly depending on $\Omega$, or for which regularity assumptions on $\Omega$, possibly depending on $\X$, does the Stokes--Dirichlet operator $\AA_{\mathcal{D}}$ admit either bounded imaginary powers (BIP) or an appropriate bounded $\mathbf{H}^\infty$-functional calculus on $\X(\Omega)$? 

\medskip

\noindent Furthermore, we wonder if one can describe the domain of the fractional powers $\D(\AA_{\mathcal{D}}^\alpha)$ in terms of standard function spaces. If so, we want to know for which values of $\alpha$, depending on $\X$ and $\Omega$, does such a characterization hold.
\end{question}

The BIPs  and the $\mathbf{H}^\infty$-functional calculus on $\L^p$-spaces, $1<p<\infty$, have been fairly investigated for the Stokes--Dirichlet operator. In the following we give a short historical review.

\begin{itemize}
    \item BIP for the Stokes--Dirichlet operator on $\L^p(\Omega)$ for all $1<p<\infty$, with the description \eqref{eq:IntroFracPowStokesDescrip} valid for all $s\in[0,2]$, was proved by Giga for bounded $\C^\infty$-domains. We refer to the seminal work \cite{Giga85} and to Giga and Sohr \cite{GigaSohr89} for the case of exterior doamins.
    \item The $\mathbf{H}^\infty$-functional calculus on $\L^p(\Omega)$ for all $1<p<\infty$ for bounded and exterior $\C^3$-domains with the description \eqref{eq:IntroFracPowStokesDescrip} for all $s\in[0,2]$ is due to Noll and Saal \cite{NollSaal2003}.
    \item The $\mathbf{H}^\infty$-functional calculus on $\L^p(\Omega)$ for all $1<p<\infty$ for bounded  $\C^{2,\alpha}$-domains, $\alpha>0$, with the description \eqref{eq:IntroFracPowStokesDescrip} for all $s\in[0,2]$ is due to Kalton, Kunstmann and Weis \cite{KaltonKunstmannWeis2006,KunstmannWeis2013Errata}.
    \item The $\mathbf{H}^\infty$-functional calculus for bounded Lipschitz domains on $\L^p(\Omega)$ for all $p$ such that $\frac{2n}{n+1}-\varepsilon_\Omega<p<\frac{2n}{n-1}+\varepsilon_\Omega$, for some $\varepsilon_\Omega>0$, is shown in the pioneering work of Kunstmann and Weis \cite{KunstmannWeis2017}. In this siuation  the description \eqref{eq:IntroFracPowStokesDescrip} is valid for all $s\in(-1+\sfrac{1}{p},1]$, thanks to the refinement by Tolksdorf \cite{Tolksdorf2018-1}. On $\L^2(\Omega)$, \eqref{eq:IntroFracPowStokesDescrip} remains valid for all $s\in(-\sfrac{1}{2},\sfrac{3}{2}]$, thanks to the result of Mitrea and Monniaux \cite{MitreaMonniaux2008}.
\end{itemize}

We conclude by claiming that similar properties for the resolvent problem, as well as the standard functional analytic properties, for the Stokes--Dirichlet operator were investigated in other domains. For the corresponding $\L^p$-theory in layers, see for instance the works by Abe and Shibata \cite{AbeShibata2003} and by Abels \cite{Abels2005PI,Abels2005PII}. In the framework of the \textit{relative} $\L^p$-theory, Farwig, Kozono, and Sohr did perform the analysis on general unbounded but uniform $\C^{1,1}$ domains for all $1<p<\infty$ \cite{FarwigKozonoSohr2005,FarwigKozonoSohr2007,FarwigKozonoSohr2008}. This has been refined very recently by Kunstmann and Tolksdorf for uniformly Lipschitz domains in \cite{KunstmannTolksdorf2025}, under the additional restriction $\frac{2n}{n+1}-\varepsilon_\Omega<p<\frac{2n}{n-1}+\varepsilon_\Omega$, for some $\varepsilon_\Omega>0$. For a similar range of indices, the standard $\L^p$-theory for exterior Lipschitz domains has been achieved by Tolksdorf and Watanabe in \cite{TolksdorfWatanabe}.

\subsection{Main results}
Our main results are concerned with estimates in the spirit of \eqref{est1} and \eqref{est1-2} in a wealth of function spaces.
In the following subsections we highlight the most crucial ones. Notably, the results for the unsteady Stokes problem hinge on the analysis of the Stokes resolvent problem \eqref{eq:IntroSystStokes}
Estimates for the Stokes resolvent problem are the content of \textbf{Chapters~\ref{Sec:StokesHalfSpace}}~and~\textbf{\ref{sec:stokessteady}}. They are turned into estimates for the unsteady Stokes problem in \textbf{Chapter~\ref{sec:MaxReg}}. Knowledge about the properties of the resolvent operator has also other interesting consequences. In particular,
a description of the Stokes--Dirichlet operator and its fractional powers, which induce more specific and sharp regularity results.

We warn the reader that the generic Stokes--Dirichlet operator (for instance on $\L^2(\Omega)$, which is the central (and easiest) case) on a domain $\Omega$ of $\RR^n$ is \textbf{not} defined as ``$\AA_\mathcal{D}=\PP_{\Omega}(-\Delta_\mathcal{D})$'' with the Leray projection $\PP_\Omega$ and the Dirichlet Laplace operator $\Delta_{\mathcal D}$.
        In general, for the Stokes--Dirichlet and Dirichlet Laplace operators defined through their respective sesquilinear forms, and domains are \textit{a priori} such that
    \begin{align*}
        \D_2(\Delta_\mathcal{D})\cap \L^2_{\mathfrak{n},\sigma}(\Omega) \subset \D_2(\AA_\mathcal{D}) \nsubseteq\D_2(\Delta_\mathcal{D}).
    \end{align*}
    Here, the subscripts $\mathfrak n$ and $\sigma$ refer to the zero traces in normal direction and solenoidability, respectively.
    The knowledge of $\D_2(\AA_\mathcal{D}) \subseteq\D_2(\Delta_\mathcal{D})$ is equivalent to know whether the pressure term $\nabla \mathfrak{p}$ itself is an element of $\L^2(\Omega,\CC^n)$, which may fail in general.
    
\subsubsection{Homogeneous Sobolev and Besov estimates on the half space}
The first main result is the estimate
\begin{align}\label{estintro:4.1}
   \lVert  (\partial_t\uu,\nabla^2\uu,\nabla\mathfrak{p})\rVert_{\L^q(\RR_+,\dot{\B}^{s}_{p,r}(\RR^n_+))} \less \lVert \ff \rVert_{\L^q(\RR_+,\dot{\B}^{s}_{p,r}(\RR_+^n))}.
\end{align}
for the Stokes problem \eqref{stokes} in the half space $\RR_+^n$ (similarly for Bessel potential spaces $\dot{\H}^{s,p}(\RR_+^n)$ provided $1<p<\infty$).
Here $p,q,r\in[1,\infty]$ and $s\in(-1+\sfrac{1}{p},\sfrac{1}{p})$, with some restrictions on $p,q$ and $r$.
The restrictions on $s$ could be relaxed by requiring compatibility conditions on the data in the spirit of \cite[Section~3]{Gaudin2023Hodge}. 

The full statement is given in \textbf{Theorem~\ref{thm:DaPratoGrisvardRn+}}. It is important to note that homogeneous endpoint Besov spaces $p=1,\infty$ are included which cannot be reached by interpolating the known results which mostly restrict the analysis to $1<p<\infty$. In particular, we cover H\"older spaces choosing $p=r=\infty$, $q\in(1,\infty]$. In the most general version also Sobolev and Besov regularity-in-time is possible.

Most of the effort in proving \eqref{estintro:4.1} is concerned with a counterpart for the resolvent problem \eqref{eq:IntroSystStokes}, \textit{i.e.},
\begin{align}\label{estintro41b}
    |\lambda|\lVert \uu\rVert_{\dot{\B}^{s}_{p,r}(\RR^n_+)}  + \lVert ( \nabla^2\uu,\nabla{\mathfrak{p}})\rVert_{\dot{\B}^{s}_{p,r}(\RR^n_+)} &\lesssim_{n,\mu} \lVert \ff\rVert_{\dot{\B}^{s}_{p,r}(\RR^n_+)}.
\end{align}
The full statement can be found in \textbf{Theorem~\ref{thm:MetaThmDirichletStokesRn+}}. The proof of the latter is based on the classical representation formula by Uka\"{i} \cite{Ukai1987}. It relates the solution to \eqref{eq:IntroSystStokes} to the resolvent problem of the Laplace operator and Riesz potentials.

Yet, \eqref{estintro41b} is far from sufficient to cover the analysis that allows to reach bounded domains by localisation arguments. When a localization procedure is applied to a solution $\uu$ of the Stokes resolvent system \eqref{eq:IntroSystStokes} on $\Omega$ with a locally flat boundary, for $\eta$ a smooth cut-off,  one has to study the system of the following form with non-zero prescribed divergence:
\begin{align*}
\lambda (\eta\uu)-\Delta(\eta\uu)+\nabla (\eta\mathfrak{p})=(\eta\ff) + \mathbf{F}(\eta,\nabla\uu,\mathfrak{p}),\quad\div(\eta\uu)=\nabla\eta\cdot \uu,\qquad (\eta \uu)_{|_{\partial\RR^n_+}}=0.
\end{align*}
Thus, we are involved in the analysis of the full non-divergence-free system
\begin{align*}
\lambda \uu-\Delta\uu+\nabla \mathfrak{p}=\ff ,\quad\div\uu=g,\qquad \uu_{|_{\partial\RR^n_+}}=0.
\end{align*}
For such a system, \textbf{Theorem~\ref{thm:MetaThmDirichletStokesRn+}} also provide a full estimate similar to \eqref{estintro41b}:
\begin{align}\label{estintro41b2}
    |\lambda|\lVert \uu\rVert_{\dot{\B}^{s}_{p,r}(\RR^n_+)}  + \lVert ( \nabla^2\uu,\nabla{\mathfrak{p}})\rVert_{\dot{\B}^{s}_{p,r}(\RR^n_+)} &\lesssim_{n,\mu} \lVert \ff\rVert_{\dot{\B}^{s}_{p,r}(\RR^n_+)} + |\lambda|\lVert g\rVert_{\dot{\B}^{s-1}_{p,r}(\RR^n_+)}  + \lVert g\rVert_{\dot{\B}^{s+1}_{p,r}(\RR^n_+)},
\end{align}
which again remains valid even when $p,r=1,\infty$, and also for homogeneous Bessel potential spaces $\dot{\H}^{s,p}(\RR_+^n)$. Actually, to reach such an estimate, one is even due to the analysis of the system with non-zero prescribed boundary value. The core strategy to reach such results is developed through the study of the Laplacian on the half space, the Hodge decomposition including the cases $p=1,\infty$, and a complete adapted function space theory. All of that  is established in prior \textbf{Chapters~\ref{Sec:Laplacians} and~\ref{Sec:FunctionSpaceTheoryDivFree}} and in the \textbf{Appendices~\ref{App:VectValFuncSpaces}, \ref{App:BesovSemigrpEst} and \ref{App:EquivNormBesov}}.

\subsubsection{\texorpdfstring{$\L^\infty$}{Loo}--regularity results on the half space}
One can certainly not expect a direct counterpart of \eqref{estintro:4.1} or \eqref{estintro41b} for the space $\L^\infty(\RR_+^n)$ in place of $\dot{\B}_{\infty,r}^s(\RR_+^n)$.
However, we can  improve the known results in the $\L^\infty$-case to
\begin{align}\label{estintro:4.2}
    |\lambda|\lVert \uu\rVert_{\L^\infty(\RR^n_+)} + |\lambda|^\frac{1}{2}\lVert \nabla\uu\rVert_{\L^\infty(\RR^n_+)} + \lVert ( \nabla^2\uu,\nabla{\mathfrak{p}})\rVert_{\dot{\B}^0_{\infty,\infty}(\RR^n_+)} &\lesssim_{n,\mu} \lVert \ff\rVert_{\L^\infty(\RR^n_+)},
\end{align}
provided $\ff\in\L_{\mathfrak{n},\sigma}^\infty(\RR^n_+)$, and with the additional pressure estimate
\begin{align*}
    |\lambda|^\frac{1}{2}\lVert \nabla \mathfrak{p}\rVert_{\dot{\B}^{-1}_{\infty,\infty}(\RR^n_+)}\lesssim_{n,\mu} \lVert \ff\rVert_{\L^\infty(\RR^n_+)},
\end{align*}
see \textbf{Theorem~\ref{thm:StokesDirRn+Linfty}}. This is completely new: {In \cite[Theorem~3.5]{DeshHieberPruss2001}} the estimate
\begin{align}\label{estintro:4.2singleest}
    |\lambda|\lVert \uu\rVert_{\L^\infty(\RR^n_+)} &\lesssim_{n,\mu} \lVert \ff\rVert_{\L^\infty(\RR^n_+)},
\end{align}
has been only proved. Here $\uu =\mathrm{K}_\lambda \ff$, where $\mathrm{K}_\lambda$ is an integral operator  of convolution-type yielding a representation formula valid on $\L^p$ for all $1<p<\infty$, and so that the operator $\lambda \mathrm{K}_\lambda$ remains uniformly bounded on $\L^\infty_{\mathfrak{n},\sigma}(\RR^n_+)$. The analysis performed did not clarify whether $\mathrm{K}_\lambda \ff$  actually solves the Stokes--Dirichlet problem in an appropriate sense. However, it is deduced in \cite[Theorem~4.4]{DeshHieberPruss2001} that there exists a (\textit{a priori} abstract) closed operator $(\D(A_{D}),A_{D})$ on $\L^\infty_{\mathfrak{n},\sigma}(\RR^n_+)$ such that $\lambda \mathrm{K}_\lambda = \lambda (\lambda \I+A_{D})^{-1}$, and this was, so far, the only accessible knowledge for the Stokes--Dirichlet problem on $\L^\infty_{\mathfrak{n},\sigma}(\RR^n_+)$. Since the formula $\lambda (\lambda \I+A_{D})^{-1}= \lambda \mathrm{K}_\lambda = \lambda (\lambda \I+\AA_{\mathcal{D}})^{-1}$ is valid on $\L^2_{\mathfrak{n},\sigma}\cap\L^\infty_{\mathfrak{n},\sigma}(\RR^n_+)$, we still write $(\D_\infty(\AA_\mathcal{D}),\AA_\mathcal{D}):=(\D(A_{D}),A_{D})$.

A fundamental consequence of the estimate \eqref{estintro:4.2} above is a description of the domain $\D_\infty(\AA_\mathcal{D})$ of the Stokes--Dirichlet operator $\AA_\mathcal{D}$ on $\L^\infty_{\mathfrak{n},\sigma}(\RR^n_+)$. We prove in \textbf{Theorem~\ref{thm:DomainStokesLInftyRn+}} that
    \begin{align}\label{eq:1809a}
        \D_\infty(\AA_\mathcal{D})=\{\, \uu\in\W^{1,\infty}_{\mathcal{D},\sigma}\cap\B^{2}_{\infty,\infty}(\RR^n_+)^n\,:\,\PP_{\RR^n_+}(-\Delta_\mathcal{D}\uu)\in\L^\infty_{\mathfrak{n},\sigma}(\RR^n_+)\,\},
    \end{align}
    where $\sigma$ refers to solenoidability and $\mathcal D$ to homogeneous Dirichlet boundary conditions.
  Furthermore, it holds for all $\uu\in\D_\infty(\AA_\mathcal{D})$,
    \begin{align}\label{eq:1809}
        \AA_\mathcal{D}\uu = \PP_{\RR^n_+}(-\Delta_\mathcal{D}\uu).
    \end{align}
    We recall that \eqref{eq:1809} \textbf{is not} the definition of the Stokes operator.
The description of the domain of the Stokes operator from \eqref{eq:1809a} seems completely new and does not follow directly from the results in \cite{DeshHieberPruss2001} and subsequent papers.
 It was even an open problem to describe the domain of the Stokes--Dirichlet operator on $\L^\infty_{\mathfrak{n},\sigma}(\RR^n_+)$ by a standard function space of order $2$. In this regard the expected conjecture is
\begin{align*}
    \D_\infty(\AA_\mathcal{D})=\{\, \uu\in\W^{1,\infty}_{\mathcal{D},\sigma}(\RR^n_+)\,:\,\nabla^2\uu\in\mathrm{BMO}(\RR^n_+,\CC^{n^3})\,\,\&\,\,\PP_{\RR^n_+}(-\Delta_\mathcal{D}\uu)\in\L^\infty_{\mathfrak{n},\sigma}(\RR^n_+)\,\}.
\end{align*}
We emphasize the basic inclusion $\{\, \uu\in\W^{1,\infty}(\RR^n_+,\CC^n)\,:\,\nabla^2\uu\in\mathrm{BMO}(\RR^n_+,\CC^{n^3})\}\subset \B^{2}_{\infty,\infty}(\RR^n_+,\CC^{n})$, such that it does not follow from our result, and still left the expected regularity result as an open question.

The result given by \eqref{eq:1809a} and \eqref{eq:1809} might seem surprising at first glance since the Leray projection is ill-defined on $\L^\infty(\RR^n_+,\CC^n)$. This is due to the membership relation
$$\text{``}\PP_{\RR^n_+}(-\Delta_\mathcal{D}\uu)\in\L^\infty_{\mathfrak{n},\sigma}(\RR^n_+)\text{''},$$
which holds provided the domain of the operator is unbounded.

The core strategy to reach such results, \eqref{estintro41b} and \eqref{eq:1809a}-\eqref{eq:1809} from  \textbf{Theorems~\ref{thm:StokesDirRn+Linfty}~\&~\ref{thm:DomainStokesLInftyRn+}}, lies in a uniqueness/Liouville-type result for weak solutions of the Stokes-resolvent problem due to Maekawa, Miura and Prange \cite{MaekawaMiuraPrange2020}, and a comparison between the Uka\"{i} and Desch--Hieber-Pr\"{u}ss representation formulas.  The Uka\"{i} formula combined with a non-trivial interpolation argument is what will allow us to make such a comparison and to obtain the gain of two derivatives in the endpoint Besov spaces.

\subsubsection{\texorpdfstring{$\L^p$}{Lp}--Sobolev and Besov regularity theory for domains of minimal regularity}
Understanding the problem in the half space is the first step towards the problem in bounded domains. This is due to the fact that the boundary of $\RR_+^n$ is flat and arbitrarily smooth.\footnote{The geometry of the half space brings, however, also problems as being unbounded. This requires the use of homogeneous function spaces -- an issue which is not relevant for bounded domains as already discussed above.}
Applying some localisation and perturbation arguments to flatten the boundary, one expects that the same result holds for bounded domains with smooth boundary. Studying rough boundaries and finding minimal assumptions on their regularity for well-posedness results is a different matter. In this regards, for $q,p,\kappa\in[1,\infty]$, $s\in(-1+\sfrac{1}{p},\sfrac{1}{p})$, we prove the estimate
\begin{align}\label{estintro:4.3}
   \lVert  (\partial_t\uu,\nabla^2\uu,\nabla\mathfrak{p})\rVert_{\L^q(\RR_+,{\B}^{s}_{p,\kappa}(\Omega))} \less \lVert \ff \rVert_{\L^q(\RR_+,{\B}^{s}_{p,\kappa}(\Omega))},
\end{align}
for the Stokes problem \eqref{stokes} in a bounded Lipschitz domain $\Omega$ whose boundary is of class $\mathcal M_{\W}^{s+2-1/p,p}$ with small norm
(similarly for Bessel potential spaces ${\H}^{s,p}(\Omega)$ provided $1<p<\infty$). The full statements are given in \textbf{Theorems~\ref{thm:bounded}~and~\ref{thm:UMDMaxRegBddDomain}}. If one is only interested in controlling
\begin{align*}
    \AA_\mathcal{D}\uu=-\Delta\uu+\nabla\mathfrak p \quad \text{ instead of }\quad (\nabla^2\uu,\nabla\mathfrak p),
\end{align*}
it is even sufficient to have a boundary of class $\mathcal M_{\W}^{1+\alpha,r}$ for some arbitrary $\alpha>0$ and arbitrary $r\geqslant 1$, with small norm. This latter condition for the boundary includes domains which might  not even be $\C^1$-domains, but that always have at least a Lipschitz boundary (with an appropriate small multiplier norm).

The estimate \eqref{estintro:4.3} (and its Sobolev counterpart) gives a comprehensive picture concerning the solvability of the non-stationary Stokes problem in domains of minimal regularity.
The only comparable result is given in \cite[Theorem 3.1]{Breit2025}, where estimates in $\L^p(\Omega)$, $1<p<\infty$.

\medbreak

Note that the latter multiplier condition, for the boundary $\partial\Omega$ to be of the class $\mathcal M_{\W}^{1+\alpha,r}$ with small norm  for some arbitrary $\alpha>0$ and arbitrary $r\geqslant 1$, is also sufficient in order to obtain boundedness of the $\mathbf{H}^\infty$-functional calculus of $\AA_\mathcal{D}$ on $\L^p_{\mathfrak{n},\sigma}(\Omega)$, for all $1<p<\infty$. Furthermore, a characterization of its fractional powers is available, that is, for which $s$, depending on $p$ and $\Omega$, one has
\begin{align}\label{eq:IntroEquivNormFracPowHsp}
    \lVert \uu \rVert_{\H^{s,p}(\Omega)}\sim_{p,s,\Omega} \lVert \AA_\mathcal{D}^{\sfrac{s}{2}}\uu \rVert_{\L^{p}(\Omega)},\quad \uu\in\D_{p}(\AA_\mathcal{D}^{\sfrac{s}{2}}),
\end{align}
see \textbf{Theorem~\ref{thm:FinalResultSobolev1}} for the full statement (see also \textbf{Meta-Theorem~\ref{thm:metaThmremoveC1}} to lower the $\C^1$-regularity assumption for the boundary into a (small) Lipschitz-boundary regularity assumption). In \textbf{Theorem~\ref{thm:FinalResultSobolevC1q}}, it is precised that for bounded $\C^{1,\alpha}$-domains, $\alpha\in(0,1]$, \eqref{eq:IntroEquivNormFracPowHsp} holds whenever
\begin{align*}
    -1+\sfrac{1}{p}<s<1+\alpha+ \sfrac{1}{p}.
\end{align*}

In \cite{Breit2025}, the approach is based on a modification of the classical approach from \cite{Solonnikov1973}: The idea is to first solve the problem in the flat geometry and to build the solution of the original problem by concatenating the (transformed) solutions.  This leads to various lower order error terms which can be controlled for small times. A global-in-time solution can be obtained by gluing local solutions together.

\medbreak

Here, our approach is instead based on an analysis of the resolvent problem \eqref{eq:IntroSystStokes} and thus closer to the elliptic theory from \cite[Theorem 3.1]{Breit2024FSI2D}, see \textbf{Theorem~\ref{thm:StokesResolvent}}. This approach is somewhat opposite to the one mentioned above. Here, one  reduces the problem in the bounded domain locally to that in the half space by a coordinate transform using the charts describing the boundary.
This leads to a Stokes--type system with coefficients and a non--zero prescribed divergence, that can be written as a perturbation of the original one (homogeneous, divergence-free).
This is reminiscent of \cite[Theorem 3.1]{Breit2024FSI2D} (which is in turn inspired by the classical analysis from \cite{MazyaShaposhnikova2009}).
However, there are various obstacles. For instance, the localisation operation (multiplication by smooth cut-off functions) is not continuous on homogeneous Sobolev and Besov spaces of too low regularity. This is related to the non-trivial divergence for the resolvent problem (after localisation) and does not appear in \cite{Breit2024FSI2D}. Additionally, we have to control lower order terms which are ``scaling'' with respect $\lambda$ in order to close the estimate. This can either be done by a very recent result from \cite{GengShen2024} in the case of bounded $\C^1$-domains or by our refined and independent analysis for the resolvent problem in Lipschitz domains with small multiplier norms given in Section~\ref{Sec:RemoveC1}, \textbf{Meta-Theorem~\ref{thm:metaThmremoveC1}}.

Having in mind the use of multiplier theory, the key strategy to obtain all the desired sharp results when the domain has low boundary regularity is \textbf{not} to perform the analysis and localization arguments in $\L^p$-spaces. Instead, we carry them out directly in negative regularity Sobolev and Besov spaces, such as $\H^{s,p}(\Omega)$ or $\B^{s}_{p,q}(\Omega)$ with $s<0$. It is precisely this point of view that will ultimately provide the sharp regularity theory in $\L^p(\Omega)$ for $1<p<\infty$, and even in the limiting cases $p=1$ and $p=\infty$; see the next subsections.

\medbreak

To put it in a nutshell, the mantra for our overall strategy is the following: ``To reach optimal elliptic regularity in a given function space, one should forget about using it directly as a ground space for the corresponding analysis''.

However, in the $\L^p$-case, $1<p<\infty$, it is obtained in \textbf{Theorem~\ref{thm:FinalResultSobolevC1q}} that
\begin{align*}
    \D_p(\AA_\mathcal{D}) = \W^{1,p}_{0,\sigma}(\Omega)\cap\W^{2,p}(\Omega)
\end{align*}
with equivalence of norms, whenever $\Omega$ is a bounded $\B^{2-\sfrac{1}{p}}_{\infty,p}$-domain. We recall that for any $\varepsilon>0$,
\begin{align*}
    \C^{1,1-\sfrac{1}{p}+\varepsilon}_{b}(\RR^{n-1}) \subsetneq\B^{2-\sfrac{1}{p}}_{\infty,p}(\RR^{n-1}) \subsetneq \C^{1,1-\sfrac{1}{p}}_{b}(\RR^{n-1}).
\end{align*}

\subsubsection{\texorpdfstring{$\L^\infty$}{Loo}-regularity theory on bounded \texorpdfstring{$\C^{1,\alpha}$}{C1a}-domains}
Similarly to the situation in Besov spaces discussed in the previous subsection, one may wonder what the minimal assumptions on the boundary are such that a version of \eqref{estintro:4.2} holds for the problem in a bounded domain and right hand-side $\ff\in\L^\infty_{\mathfrak{n},\sigma}(\Omega)$. In \textbf{Theorem~\ref{thm:FinalResultLinfty1}}, we prove the estimate
\begin{align*}
     (1+\lvert\lambda\rvert) \lVert \uu\rVert_{\L^\infty(\Omega)}+(1+\lvert\lambda\rvert)^\frac{1-\alpha}{2}\lVert ( \nabla \uu, \mathfrak{p} )\rVert_{\B^{\alpha}_{\infty,\infty}(\Omega)} \nonumber \lesssim_{n,\alpha,\mu}^{\Omega} \lVert \ff\rVert_{\L^\infty(\Omega)}
\end{align*}
for the Stokes--Dirichlet resolvent problem \eqref{eq:IntroSystStokes} in a bounded $\C^{1,\alpha}$-domain with small $(1,\alpha)$-H\"{o}lderian constant. Here $\alpha\in(0,1)$ is arbitrary. We are also allowed
to describe the domain of the Stokes--Dirichlet operator on $\L^\infty(\Omega)$ (under the same assumption on the domain). With notations similar to that in \eqref{eq:1809a}, we prove in \textbf{Theorem~\ref{thm:FinalResultLinfty2}} that 
    \begin{align*}
        \D_{\infty}(\AA_\mathcal{D}):=\{\,\uu\in\C^{1,\alpha}_{\mathcal{D},\sigma}(\overline{\Omega})\,:\,\PP_{\Omega}(-\Delta_\mathcal{D}\uu)\in\L^\infty_{\mathfrak{n},\sigma}(\Omega)\,\},
    \end{align*}
    as well as
    \begin{align*}
        \AA_\mathcal{D}\uu = \PP_{\Omega}(-\Delta_\mathcal{D}\uu),\quad \text{ for all }\uu\in\D_{\infty}(\AA_{\mathcal{D}}).
    \end{align*}
These results improve those by Abe and Giga \cite[Theorem~1.2]{AbeGiga2013}, and Abe, Giga and Hieber \cite[Theorem~1.1]{AbeGigaHieber2015}, where bounded $\C^3$-domains and bounded $\C^2$-domains were considered respectively.  Also, only a Lipschitz bound can be derived from the resolvent estimate of the former, even in the more recent work by Geng and Shen \cite[Theorems~1.1~\&~9.4]{GengShen2025} for bounded $\C^{1,\alpha}$-domains. Under the same regularity assumptions on the boundary, we obtain a precise description of the exact regularity of the solution that matches the regularity of the boundary.

As a consequence of our analysis of the resolvent problem on $\L^\infty$, in \textbf{Theorem~\ref{thm:MaxRegBesov(Linfty)}}, it is obtained
\begin{align*}
 \lVert  (\uu,\partial_t\uu,-\Delta\uu+\nabla\mathfrak{p})\rVert_{\B^\beta_{q,\kappa}(\RR_+,\L^\infty(\Omega))}\lesssim_{q,\beta,\Omega} \lVert \ff \rVert_{\B^\beta_{q,\kappa}(\RR_+,\L^\infty(\Omega))}
\end{align*}
for the non-stationary Stokes problem
in a bounded $\C^{1,\alpha}$-domain, with $\alpha\in(0,1]$ arbitrary. 

\subsubsection{\texorpdfstring{$\L^1$}{L1}-type regularity theory on bounded \texorpdfstring{$\C^{1,\alpha}$}{C1alpha}-domains}

As a consequence of the study on Besov spaces $\B^{s}_{1,q}$ and on $\L^\infty$, 
\textbf{Corollaries~\ref{cor:StokesL1-semigroup}} and~\textbf{\ref{cor:StokesWs1}} establish, for the first time, a coherent $\L^1$-type theory for the Stokes--Dirichlet semigroup on bounded domains with low boundary regularity, relying on the identification $\B^{s}_{1,1}=\W^{s,1}$ for all $s\notin\NN$.

In the existing literature—essentially restricted to the half-space (see~\cite{GigaMatsuiShimizu1999,ShibataShimizu2001}) or to smooth bounded domains (see~\cite{Maremonti2011})—some results can extend parts of the classical $\L^p$-framework down to the endpoint case $p=1$. 
For $\Omega$ being either $\RR^n_+$ or a sufficiently smooth bounded domain, one has for instance
\begin{align*}
    \| \nabla e^{-t\AA_{\mathcal{D}}}\uu \|_{\L^1(\Omega)} 
    \lesssim_{n,\Omega} t^{-1/2}\|\uu\|_{\L^1(\Omega)},
\end{align*}
even though neither standard duality arguments nor heat kernel bounds are available in this setting. 
Such an estimate is the best one can expect, since in general the uniform boundedness of the Stokes–-Dirichlet semigroup may fail on $\L^1_{\mathfrak{n},\sigma}$. See for instance~\cite[Theorem~5.1]{DeshHieberPruss2001} for the case of $\L^1_{\mathfrak{n},\sigma}(\RR^{n}_+)$.

\medbreak

\textbf{Corollary~\ref{cor:StokesL1-semigroup}} shows that the Stokes--Dirichlet semigroup $e^{-t\AA_{\mathcal{D}}}$ still enjoys parabolic $\L^1$-smoothing effects in bounded $\C^{1,\alpha}$-domains, $\alpha\in(0,1]$:
\begin{align*}
    \| e^{-t\AA_{\mathcal{D}}}\uu \|_{\W^{s,1}(\Omega)} 
    \lesssim_{s,n,\Omega} t^{-s/2}\|\uu\|_{\L^1(\Omega)}, 
    \qquad s\in(0,2+\alpha), \ t>0.
\end{align*}
In particular, one recovers the gradient decay estimate 
$\| \nabla e^{-t\AA_{\mathcal{D}}}\uu \|_{\L^1(\Omega)} 
\lesssim_{n,\Omega} t^{-1/2}\|\uu\|_{\L^1(\Omega)}$, 
mirroring the diffusive behavior of the heat semigroup. 
This provides the first quantitative $\L^1$-regularization result for the Stokes--Dirichlet system in bounded, nonsmooth domains. 
More generally, $\L^1$–$\L^p$ decay estimates for all $1<p\leqslant \infty$ are obtained in \textbf{Proposition~\ref{prop:Optimaldecay}}.

\medbreak

\textbf{Corollary~\ref{cor:StokesWs1}} complements this picture by providing a detailed functional-analytic characterization of $\AA_{\mathcal{D}}$ in the $\L^1$-type setting for bounded Lipschitz domains $\Omega$ belonging to the rougher class $\mathcal{M}^{1+\alpha,1}_{\W}(\epsilon)$\footnote{This class includes Lipschitz domains with small Lipschitz constant that are not of class $\C^1$.} for some $\epsilon>0$. 
It shows that $\AA_{\mathcal{D}}$ acts as an isomorphism between $\W^{s+2,1}(\Omega)$- and $\W^{s,1}$-divergence-free spaces for $s>0$, admits a bounded $\mathbf{H}^{\infty}$-functional calculus, and that its fractional domain of order $\sfrac{\beta}{2}$ coincides with the natural $\W^{s+\beta,1}$-spaces, $s+\beta\neq1,2$, (with normal or no--slip boundary conditions depending on whether $s+\beta>1$). 
In particular, the full real interpolation scale associated with $\AA_{\mathcal{D}}$ is now identified in $\L^1$.

\medbreak

Together, these two results provide a robust analytic foundation for the $\L^1$-type theory of the Stokes operator in rough bounded domains.

\subsubsection{The Leray projection for bounded domains of low regularity, including the endpoint spaces}

The results of \textbf{Proposition~\ref{prop:HodgeDiracbdd}}, 
\textbf{Theorem~\ref{thm:SharpHodgeDecompC1}} 
and \textbf{Proposition~\ref{prop:NeumannPbLipschitzMult}} 
extend the classical \emph{Hodge decomposition} to the full scale of Besov and Sobolev spaces 
on bounded Lipschitz domains with low regularity, 
and provide a unified operator-theoretic approach valid even at the endpoint cases $p=1$ and $p=\infty$.

\medbreak

\noindent
The classical $\L^p$-decomposition
\begin{align*}
    \L^p(\Omega,\CC^n)
= \L^p_{\mathfrak{n},\sigma}(\Omega)
{\oplus}\overline{\nabla \H^{1,p}(\Omega)}
\end{align*}
ensures that any $\uu\in\L^p(\Omega,\CC^n)$ can be written as $\uu=\vv+\nabla\mathfrak{q}$, with $\div\vv=0$ in $\Omega$ and $\vv\cdot\nu=0$ on $\partial\Omega$. The projection map onto $\L^p_{\mathfrak{n},\sigma}(\Omega)$ is denoted $\PP_{\Omega}$ and called the Leray projection, and one has
\begin{align*}
    \PP_{\Omega}\uu = \vv =\uu -\nabla\mathfrak{q}.
\end{align*}
It was proved by Simader and Sohr~\cite{SimaderSohr1992} that the decomposition holds for all $p\in(1,\infty)$ on bounded and exterior $\C^1$-domains, 
and by Fabes, Mendez, and Mitrea~\cite{FabesMendezMitrea1998} for bounded Lipschitz domains when $p$ lies near~$2$.  
Bogovski\u{\i}~\cite{Bogovskii1986} constructed counterexamples in unbounded settings.  
The key analytic tool is the Neumann problem for the Laplacian,
\begin{align}\label{eq:IntroNeumannPbHodgeProj}
-\Delta \mathfrak{q} = -\div\uu, \qquad\text{and}\qquad 
\nu\cdot[\nabla \mathfrak{q}-\uu]=0, \, \text{ on } \partial\Omega,
\end{align}
which determines the gradient component.

\medbreak

\noindent
As a particular case of \textbf{Theorem~\ref{thm:SharpHodgeDecompC1}}, we extend this decomposition to the full Sobolev–Besov scale:
\begin{align}\label{eq:IntroHodgeDecomp}
    \H^{s,p}(\Omega,\CC^n)
= \H^{s,p}_{\mathfrak{n},\sigma}(\Omega)
\oplus \overline{\nabla \H^{s+1,p}(\Omega,\CC)}, 
\qquad s\in(-1+\tfrac{1}{p},\tfrac{1}{p}),
\end{align}
and analogously for $\B^s_{p,q}$, including $p,q=1,\infty$.  
This completes the results of Mitrea and Monniaux~\cite{MitreaMonniaux2008} for bounded $\C^{1}$-domains, and of Fujiwara and Yamazaki~\cite{FujiwaraYamazaki2007} for $\C^{2,1}$-domains, the latter being the only prior work treating endpoint Besov spaces.

\medbreak

\noindent
Following M${}^\text{c}$Intosh and Monniaux~\cite{McintoshMonniaux2018}, we reformulate the problem in the setting of differential forms:
\begin{align}\label{eq:IntroHodgeDecompDiff}
    {\H}^{s,p}(\Omega,\Lambda^k)
= {\H}^{s,p}_{\mathfrak{n},\sigma}(\Omega,\Lambda^k)
\oplus \overline{\d\,\H^{s+1,p}(\Omega,\Lambda^{k-1})},
\end{align}
where the exact part $\d\boldsymbol{\omega}$ arises from solving the Hodge–Laplace problem
\begin{align*}
    -\Delta\boldsymbol{\omega} = \delta u,
    \qquad \text{and}\qquad
    \nu\iprod\boldsymbol{\omega}=0,\quad 
    \nu\iprod[\d\boldsymbol{\omega}-u]=0,\, \text{ on } \partial\Omega,
\end{align*}
where $-\Delta_{\mathcal{H}}=(\d+\d^\ast)^2$, $\d^\ast =\delta$ with a prescribed boundary condition, and ${\H}^{s,p}_{\mathfrak{n},\sigma}(\Omega,\Lambda^k)={\H}^{s,p}(\Omega,\Lambda^k)\cap\N(\d^\ast)$. For $k=1$, this reduces to the vector-field case presented above, with $\delta=-\div$ on $1$-forms/vector fields and $\d = \nabla$ on $0$-forms/scalar functions.

\medbreak

\noindent
The analysis focuses instead on the \textit{first order} system induced by the \textbf{Hodge–Dirac operator} $D_\mathfrak{n}=\d+\d^\ast$, whose resolvent problem
\begin{align*}
    \lambda \omega+i(\d + \delta)\omega =f, \qquad \nu\iprod\omega=0,
\end{align*}
yields sharp regularity across Besov and Sobolev scales.  
\textbf{Proposition~\ref{prop:HodgeDiracbdd}} provides these resolvent estimates in all admissible $\B^{s}_{p,q}$ spaces for bounded Lipschitz domains of class $\mathcal{M}^{1+\alpha,r}_\W(\epsilon)$, $\epsilon>0$ small.  
As a consequence, \textbf{Theorem~\ref{thm:SharpHodgeDecompC1}} establishes \eqref{eq:IntroHodgeDecompDiff} and the bounded Hodge–Leray projections on divergence-free fields with vanishing normal trace.  
For $k=1$, this corresponds to the solvability of the Neumann problem for the potential $\mathfrak{p}$, as in \textbf{Proposition~\ref{prop:NeumannPbLipschitzMult}}, yielding quantitative $\B^{s}_{p,q}$-estimates on $\nabla\mathfrak{p}$ and the validity of \eqref{eq:IntroHodgeDecomp} on any bounded Lipschitz domain of class $\mathcal{M}^{1+\alpha,1}_\W(\epsilon)$. This even includes domains whose boundary does not belong to the class $\C^1$.

\medbreak

\noindent
As a consequence, \textbf{Theorem~\ref{thm:SharpHodgeDecompC1}} 
yields the Hodge decomposition \eqref{eq:IntroHodgeDecompDiff} in the full Besov and Sobolev scales for a wide subclass of bounded Lipschitz domains, with bounded Hodge–Leray projection $\PP_\Omega$ onto divergence-free fields with vanishing normal component on the boundary.

In the case of vector fields ($k=1$), the decomposition translates into the \emph{solvability of the Neumann problem} \eqref{eq:IntroNeumannPbHodgeProj} for the scalar potential $\mathfrak{q}$, established in \textbf{Proposition~\ref{prop:NeumannPbLipschitzMult}}, which provides quantitative $\B^{s}_{p,q}$-estimates on $\nabla\mathfrak{q}$. In particular, \eqref{eq:IntroHodgeDecomp} holds for all Sobolev spaces $\H^{s,p}$, $1<p<\infty$, and all Besov spaces $\B^{s}_{p,q}$, $1\leqslant p,q\leqslant \infty$, on any bounded Lipschitz domain of class $\mathcal{M}^{1+\alpha,r}_\W(\epsilon)$, for sufficiently small $\epsilon>0$, which in particular need not be of class $\C^1$.

\medbreak

\noindent
These results offer the first complete theory of Hodge(--Helmholtz) decompositions in a wide class of low-regularity bounded domains, valid across the full $\L^p$–Besov–Sobolev scales and including endpoint function spaces for non-smooth and bounded geometries.

\subsection{Technical biases, technical issues and subsequent additional results}

We aim to highlight the technical motivations behind the present work. In order to address the problems from our chosen perspective, we must handle several technical issues and shortcomings that have not yet been fully clarified in the literature. Here we explain their impact in order to assess the novelty and relevance of our approach, and we highlight a few independent yet central results that may be of broader interest.

\subsubsection{The overall strategy: robustness and strengths of the established framework}\label{sec:IntroOverallStrat}

The strategy of proof presented here is very general and not at all restricted to the Stokes--Dirichlet problem. It can be straightforward applied to many other elliptic/parabolic problems subjects to various suitable boundary conditions. For instance, one could consider the Stokes problem subject to Neumann or other boundary conditions, with variable viscosity of minimal smoothness (such as: given in space of multipliers), as well as many other variant problems. Our strategy even simplifies for problems with no prescribed algebraic condition such as the divergence-free constraint. So that for those system one should be able to obtain again a similar full theory on all function spaces simultaneously including $\H^{s,p}$ and $\B^{s}_{p,q}$ as well as the endpoint spaces $\L^\infty$, $\L^1$, $\B^{s}_{\infty,q}$ and $\B^s_{1,q}$. 

\medbreak

In particular, as exhibited in the steps below --except for the specific space $\L^\infty(\Omega)$-- our main results obtained throughout the current document are \textbf{completely independent} of the current state of literature on the Stokes--Dirichlet problem on bounded smooth (\textit{e.g.}, between $\C^3$ and $\C^{1,1}$) and non-smooth (\textit{e.g.}, $\C^{1}$ or Lipschitz) domains. Moreover, such knowledge would not simplify by any means the proofs of our results, since we reach some exclusive regularity threshold.

\medbreak
In the following we describe the general approach in order to make it applicable for other problems.

We want to investigate the generic problem
\begin{align*}
    \lambda u+{\L}u =f, \quad\text{ and }\quad \mathcal{B}u=0,\text{ on }\partial\Omega,
\end{align*}
for an elliptic operator ${\L}$ of order $2$  (not necessarily local: it could include some algebraic condition such as the divergence-free constraint) with a set of prescribed boundary conditions given by the boundary operator $\mathcal{B}$. The arising operator when defined on $\RR^{n}_+$ is denoted by $\L_{\mathcal{B},+}$, the sole $\L_{\mathcal{B}}$ notation is for  the problem on the domain $\Omega$.

The main ingredients are the following:
\begin{enumerate}
    \item A complete function space theory that takes into account specificity of the elliptic operator such as interpolation theory of Sobolev and Besov spaces with boundary and algebraic condition (such as divergence-free or curl-free) on $\RR^n$, $\RR^n_+$ and domains. Appropriate density results for smooth functions might also be needed. This corresponds to \textbf{Chapter~\ref{Sec:FunctionSpaceTheoryDivFree}}.

    \item A suitable and consistent representation formula $\mathrm{K}^{{\L}_{\mathcal{B},+}}_\lambda$ for a corresponding half space operator $(\lambda \I+{\L}_{\mathcal{B},+})^{-1}$, valid on all $\dot{\H}^{s,p}(\RR^n_+)$, $1<p<\infty$, and $\dot{\B}^{s}_{p,q}(\RR^n_+)$, for \textit{some } $-1+\sfrac{1}{p}<s<\sfrac{1}{p}$, as well as --whenever its possible-- $\L^\infty(\RR^n_+)$, $\L^1(\RR^n_+)$. Although a representation formula is not necessary, a key point of our strategy is to figure out at least one $s_p\in(-1+\sfrac{1}{p},\sfrac{1}{p}]$, depending on the regularity of the coefficients, possibly arbitrarily close do $-1+\sfrac{1}{p}$, such that for all $$-1+\sfrac{1}{p}<s<s_p,$$
    for all $f\in\dot{\H}^{s,p}(\RR^n_+)$
    \begin{align*}
        \lVert\nabla^2(\lambda \I+{\L}_{\mathcal{B},+})^{-1} f\rVert_{\dot{\H}^{s,p}(\RR^n_+)} \lesssim_{p,s,n}^{\L,\mathcal{B}}\lVert
        f\rVert_{\dot{\H}^{s,p}(\RR^n_+)}, 
    \end{align*}
    and similarly for $\dot{\B}^{s}_{p,q}$, $p,q\in[1,\infty]$. To this end, a representation formula is the most common approach and will often allow to reach (almost) all integrability indeces.
    
    If one uses a representation formula approach, one should use the $\L^2$-theory and check that $\mathrm{K}^{{\L}_{\mathcal{B},+}}_\lambda=(\lambda \I+{\L}_{\mathcal{B},+})^{-1}$ on $\L^2(\RR^n_+)$ or the corresponding closed subspace. This is \textbf{Chapter~\ref{Sec:StokesHalfSpace}} with $s_p=\sfrac{1}{p}$.

    \item Depending on the specific structure of ${\L}$ and  $\mathcal{B}$ on $\RR^{n}_+$, one might have to consider the \textit{ad-hoc} boundary value problems
    \begin{align*}
    \lambda u+{\L}u =0 \quad\text{ and }\quad \mathcal{B}u=g,\text{ on }\partial\RR^{n}_+,
\end{align*}
and
    \begin{align*}
    \lambda u-\Delta u =0 \quad\text{ and }\quad \mathcal{B}u=g,\text{ on }\partial\RR^{n}_+.
\end{align*}
Whenever $-1+\sfrac{1}{p}<s<s_p$, for both one shall prove an estimate of the form
\begin{align*}
    |\lambda|\lVert u\rVert_{\dot{\H}^{s,p}(\RR^n_+)}+ \lVert \nabla^2 u\rVert_{\dot{\H}^{s,p}(\RR^n_+)} \lesssim_{p,s,n}^{\L,\mathcal{B}} |\lambda|^{\frac{k}{2}+\frac{s}{2}-\frac{1}{2p}}\lVert
        g\rVert_{{\L}^{p}(\RR^{n-1})}+\lVert
        g\rVert_{\dot{\B}^{s+k}_{p,p}(\RR^{n-1})} 
\end{align*}
for $k=1$ or $2$ depending on the structure of $\mathcal{B}$. Similarly for $\dot{\B}^{s}_{p,q}$, $p,q\in[1,\infty]$, instead of $\dot{\H}^{s,p}$.

    When $\mathcal{B}$ has constant coefficients this yields Neumann or Dirichlet-type boundary conditions, the latter ties to the behavior of the boundary semigroups $(e^{-x_n({\lambda\I-\Delta'})^\frac{1}{2}})_{x_n\geqslant 0}$ and $(e^{-x_n({-\Delta'})^\frac{1}{2}})_{x_n\geqslant 0}$  acting on $\dot{\B}^{s-\sfrac{1}{p}}_{p,q}(\RR^{n-1})$ with values in $\dot{\H}^{s,p}(\RR^{n}_+)$ or $\dot{\B}^{s}_{p,q}(\RR^{n}_+)$, $s\in\RR$, with $p,q=1,\infty$ for Besov spaces. This is achieved in  \textbf{Chapters~\ref{Sec:Laplacians} and~\ref{Sec:FunctionSpaceTheoryDivFree},~Section~\ref{subsec:HodgeRn+}} by means of the \textbf{Appendices~\ref{App:VectValFuncSpaces}, \ref{App:BesovSemigrpEst} and \ref{App:EquivNormBesov}}.

    If the coefficients in $\L$ and $\mathcal{B}$ are non-constant and in a suitable multiplier class one can then proceed by ``perturbation arguments'' following for instance \cite[Chapter~7]{KunstmannWeis2004} and \cite[Chapter~6]{PrussSimonett2016}.

    \item Consider $\Omega$ a bounded  domain whose boundary is in the multiplier class $\mathcal{M}_\W^{1+\alpha,r}$, $\alpha>0$, $r\in[1,\infty]$, with small norm. Thanks to the multiplier theory by Maz'ya and Shaposhnikova \cite{MazyaShaposhnikova2009} reintroduced in \textbf{Section~\ref{sec:SM}}, we use standard localisation schemes and perturbation arguments on a solution $u\in\H^{1,2}(\Omega)$ for a given $f\in\L^2(\Omega)\cap\H^{s,p}(\Omega)$, to prove first that
    \begin{align*}
        |\lambda|\lVert u \rVert_{{\H}^{s,p}(\Omega)}+\lVert\nabla^2u \rVert_{{\H}^{s,p}(\Omega)} \lesssim_{p,s,n,\Omega}^{\L,\mathcal{B}}\lVert
        f\rVert_{{\H}^{s,p}(\Omega)},
    \end{align*}
    provided $s$ satisfies $-\sfrac{1}{p}<s<\min\big(s_p,-1+\frac{1+\min(r,p)\alpha}{p}\big)$. Then, up to consider the corresponding adjoint problem $\L^\ast$ with the system of boundary conditions $\mathcal{B}^\ast$, one obtains then the same result and deduces sectoriality of $\L_\mathcal{B}$ on $\H^{s,p}(\Omega)$ for all $s\in(-1+\sfrac{1}{p},\sfrac{1}{p})$ by duality and interpolation. However, for $s$ close to $1/p$, \textit{a priori}, we only have
    \begin{align*}
        |\lambda|\lVert u \rVert_{{\H}^{s,p}(\Omega)} \lesssim_{p,s,n,\Omega}^{\L,\mathcal{B}}\lVert
        f\rVert_{{\H}^{s,p}(\Omega)}.
    \end{align*}
    The same can be deduced for Besov spaces $\B^{s}_{p,q}$, including the cases $p=1$ and $p=\infty$.
    
    This is achieved for the Stokes--Dirichlet operator in \textbf{Chapter~\ref{sec:stokessteady}, Section~\ref{Sec:ResolventEstGain2Deriv}}. We recall that in our case, we have $s_p=\sfrac{1}{p}$.

    \item For the bounded $\mathbf{H}^\infty$-functional calculus, one first starts to show the following resolvent estimates on $\RR^{n}_+$, for $s\in(-1+\sfrac{1}{p},s_p)$, $\mu\in(\frac{\pi}{2},\pi)$, for all $(\lambda_\ell)_{\ell\in\NN}\in\{\,z\in\CC^\ast\,:\, |\arg(z)|<\mu\,\}$, all $(f_\ell)_{\ell\in\NN}\subset\dot{\H}^{s,p}(\RR^n_+)$, setting $u_\ell = (\lambda_\ell \I+{\L}_{\mathcal{B},+})^{-1} f_\ell$
    \begin{align}\label{eq:RboundIntroRn+}
        \bigg\lVert \sum_{\ell=1}^N r_\ell\lambda_\ell u_\ell \bigg\rVert_{\L^p(0,1;\dot{\H}^{s,p}(\RR^n_+))}+\bigg\lVert \sum_{\ell=1}^N r_\ell\nabla^2u_\ell\bigg\rVert_{\L^p(0,1;\dot{\H}^{s,p}(\RR^n_+))} \lesssim_{p,s,n}^{\L,\mathcal{B}}\bigg\lVert
        \sum_{\ell=1}^N r_\ell f_\ell\bigg\rVert_{\L^p(0,1;\dot{\H}^{s,p}(\RR^n_+))}, 
    \end{align}
    with a bound uniform with respect to $N\in\NN$ and where the definition of $(r_\ell)_{\ell\in\NN}$ is given in \textbf{Definition~\ref{def:rbound}}. Such a bound \eqref{eq:RboundIntroRn+} is called an $\mathcal{R}$-bound, and depending on the structure of $\L$ and the boundary conditions $\mathcal{B}$, one may have to consider a variation including non-zero boundary conditions $(g_\ell)_{\ell\in\NN}$.

    For the Stokes--Dirichlet problem on $\dot{\H}^{s,p}(\RR^n_+)$, this is achieved in \textbf{Chapter~\ref{Sec:StokesHalfSpace}, Section~\ref{sec:RboundRn+Hsp}}. For the elliptic operator $\L$ with a system of boundary conditions $\mathcal{B}$, one should instead draw inspiration from \cite[Chapter~6]{bookDenkHieberPruss2003}, \cite[Chapter~6]{KunstmannWeis2004} or   \cite[Chapter~6,~Section~6.2]{PrussSimonett2016}, but these do not apply directly in to Stokes-type problems due to the divergence(-free) constraint.

    A localisation argument applied to \eqref{eq:RboundIntroRn+} shows that for the corresponding problem on the domain $\Omega$ with a boundary in the class $\mathcal{M}_\W^{1+\alpha,r}$, one reaches
    \begin{align}\label{eq:RboundIntroOmega}
        \bigg\lVert \sum_{\ell=1}^N r_\ell\lambda_\ell u_\ell \bigg\rVert_{\L^p(0,1;{\H}^{s,p}(\Omega))}+\bigg\lVert \sum_{\ell=1}^N r_\ell\nabla^2u_\ell\bigg\rVert_{\L^p(0,1;{\H}^{s,p}(\Omega))} \lesssim_{p,s,n}^{\L,\mathcal{B}}\bigg\lVert
        \sum_{\ell=1}^N r_\ell f_\ell\bigg\rVert_{\L^p(0,1;{\H}^{s,p}(\Omega))}, 
    \end{align}
    whenever $-\sfrac{1}{p}<s<\min\big(s_p,-1+\frac{1+\min(r,p)\alpha}{p}\big)$.

    \medbreak

    One reproduces the procedure for $\L^\ast$ with boundary values $\mathcal{B}^\ast$. Consequently, thanks to \eqref{eq:RboundIntroOmega}, one applies the Kunstmann--Weis comparison principle given in \cite{KunstmannWeis2017}, in order to extrapolate the bounded $\mathbf{H}^\infty$-functional calculus of $\L_\mathcal{B}$ from $\L^2(\Omega)$ to $\H^{s,p}(\Omega)$, $-\sfrac{1}{p}<s<\min\big(s_p,-1+\frac{1+\min(r,p)\alpha}{p}\big)$. Finally, complex interpolation will yield the same result on $\L^p(\Omega)$ with a characterization of the domain of fractional powers in terms of standard Sobolev spaces $\H^{\beta,p}$ whenever $\beta<\min\big(s_p+2,1+\frac{1+\min(r,p)\alpha}{p}\big)$.

    For the Stokes--Dirichlet operator, this final step including \eqref{eq:RboundIntroOmega} and the boundedness of the $\mathbf{H}^\infty$-functional calculus --with an explanation of the Kunstmann--Weis comparison principle--  can be found in \textbf{Chapter~\ref{sec:stokessteady}, Section~\ref{sec:ProofLPbddDomain}}.
\end{enumerate}

\paragraph{Where is the imoprovement with respect to the standard approach ?}Consider $\alpha>0$, $r\in[1,\infty]$, the condition $-1+\sfrac{1}{p}<s<\min\big(s_p,-1+\frac{1+\min(r,p)\alpha}{p}\big)$ and coefficients for $\L$ and $\mathcal{B}$ belonging to an appropriate multiplier class, with boundary $\partial\Omega$ in the multiplier class $\mathcal{M}_\W^{1+\alpha,r}$ with small norm. Performing the analysis with the localisation argument directly in negative Sobolev and Besov spaces $\H^{s,p}(\Omega)$ and $\B^{s}_{p,q}(\Omega)$ in combination with multiplier theory is what allowed us to gain a full derivative on the known results could. This refers either to the boundary regularity or to the regularity of the coefficients.

However, for $\Omega\subset \RR^n$ with purely Lipschitz or $\C^{1}$-boundary as above, (or if the coefficients for $\L$ and $\mathcal{B}$ are very rough), this strategy may not apply anymore ($s_p$ might not exist). However the analysis, allows for domains that are ``in between'' Lipschitz and $\C^{1}$ in some sense.

Obviously, there are more technical underlying details, such as:
\begin{itemize}[label=--]
    \item It is not clear how to put away from $0$ for the localised resolvent estimates (for $\lambda$'s). Indeed, elementary localisation arguments for the resolvent system imply the need to consider $\lambda\I + c_0\I+\L_\mathcal{B}$ for some $c_0>0$ large instead of $\lambda\I +\L_\mathcal{B}$.
    \item For the Stokes--Dirichlet problem, a more specific issue is related to the fact that one has to deal with the pressure function and a non-trivial prescribed divergence condition for the localised system.
    \item A more generic issue is the transition from homogeneous to inhomogeneous function spaces while performing the localisation procedure. This becomes critical when combined with the issue arising from the necessity to control the divergence part of the localised system.
\end{itemize}

\subsubsection{Technical issue 1: Homogeneous function spaces for the half space analysis}

The first step towards an analysis on domains for the Stokes--Dirichlet problem is a suitable theory on the half space.

To this end, the use of homogeneous Sobolev and Besov spaces $\dot{\H}^{s,p}$, $\dot{\W}^{s,p}$ and $\dot{\B}^{s}_{p,q}$, instead of the inhomogeneous and more standard ones ${\H}^{s,p}$, ${\W}^{s,p}$ and ${\B}^{s}_{p,q}$, is necessary  for a suitable theory or several reasons.

\begin{itemize}[label=--]
    \item As mentioned before  (see the discussion around \eqref{estintro41b2}), for the said suitable theory, one should consider the system with non-zero prescribed divergence, which on $\L^p$ translates as the estimate
    \begin{align}\label{eq:IntroNondivfreeResvEstRn+}
    |\lambda|\lVert \uu\rVert_{\L^p(\RR^n_+)}  + \lVert ( \nabla^2\uu,\nabla{\mathfrak{p}})\rVert_{\L^p(\RR^n_+)} &\lesssim_{n,\mu} \lVert \ff\rVert_{\L^p(\RR^n_+)} + |\lambda|\lVert g\rVert_{\dot{\W}^{-1,p}(\RR^n_+)}  + \lVert g\rVert_{\dot{\W}^{1,p}(\RR^n_+)}.
\end{align}
The measure of $\div \uu =g$ in homogeneous Sobolev spaces $\dot{\W}^{-1,p}_0(\RR^n_+)\cap \dot{\W}^{1,p}(\RR^n_+)$ rather than inhomogeneous one is necessary to preserve the ``regularity scaling'' with respect to $\lambda$. When it comes to fractional Sobolev spaces ${\H}^{s,p}$, $-1+\sfrac{1}{p}<s<\sfrac{1}{p}$, this scaling property is only preserved if all quantities are measured with homogeneous $\dot{\H}^{s,p}$--norms, that is one has an estimate of the form:
\begin{align}\label{eq:IntroNondivfreeResvEstHspRn+}
    |\lambda|\lVert \uu\rVert_{\dot{\H}^{s,p}(\RR^n_+)}  + \lVert ( \nabla^2\uu,\nabla{\mathfrak{p}})\rVert_{\dot{\H}^{s,p}(\RR^n_+)} &\lesssim_{n,\mu} \lVert \ff\rVert_{\dot{\H}^{s,p}(\RR^n_+)} + |\lambda|\lVert g\rVert_{\dot{\H}^{s-1,p}(\RR^n_+)}  + \lVert g\rVert_{\dot{\H}^{s+1,p}(\RR^n_+)}.
\end{align}
Otherwise, identities like $\L^p\cap\dot{\H}^{s,p}={\H}^{s,p}$ $s>0$, or $\L^p+\dot{\H}^{s,p}={\H}^{s,p}$, $s<0$, will produce some anisotropy with respect to the parameter $\lambda$ either on the lower or the higher regularity part of the norms that measure the size of $g$. The same goes for Besov spaces $\dot{\B}^{s}_{p,q}$. If $g=0$ and $1<p<\infty$, one can measure all quantitities with respect to inhomogeneous norms, removing the dots in \eqref{eq:IntroNondivfreeResvEstHspRn+}.

\item Since we want to reach an analysis in the endpoint Besov spaces over $\L^1$ and $\L^\infty$, one also needs to consider exclusively the homogeneous Besov spaces $\dot{\B}^{s}_{1,q}$, $\dot{\B}^{-s}_{\infty,q}$, $s\in(0,1)$, and not their inhomogeneous counterparts.

If one considers $\RR^n$, this restriction appears naturally if one wants to take advantage of the boundedness of the Leray projection
\begin{align*}
    \PP_{\RR^n}=\I+\nabla(-\Delta)^{-1}\div.
\end{align*}
It is known \textbf{not} to be bounded on ${\B}^{s}_{1,q}(\RR^n)$, ${\B}^{-s}_{\infty,q}(\RR^n)$, but it is on $\dot{\B}^{s}_{1,q}(\RR^n)$, $\dot{\B}^{-s}_{\infty,q}(\RR^n)$.

Since the corresponding operator $\PP_{\RR^n_+}$ on $\RR^n_+$ is  known to be the composition of an extension operator by reflection and $\PP_{\RR^n}$, similarly, one cannot expect its boundedness on ${\B}^{s}_{1,q}(\RR^n_+)$, ${\B}^{-s}_{\infty,q}(\RR^n_+)$, while it is bounded on $\dot{\B}^{s}_{1,q}(\RR^n_+)$ and $\dot{\B}^{-s}_{\infty,q}(\RR^n_+)$, $s\in(0,1)$.

\item When one considers the regularity theory on $\L^\infty(\RR^n_+)$ with divergence--free condition, we show the estimate \eqref{estintro:4.2} in \textbf{Theorem~\ref{thm:StokesDirRn+Linfty}} and it is also proved that the $\dot{\B}^{0}_{\infty,\infty}$-norm for the measurement of $\nabla^2 \uu$ and $\nabla\mathfrak{p}$ can not be replaced by an the inhomogeneous  ${\B}^{0}_{\infty,\infty}$-norm (unless one accepts a logarithmic-loss of regularity and scaling with respect to the resolvent parameter $\lambda$).
\end{itemize}

\subsubsection{Technical issue 2: Homogeneous function spaces and prescribed divergence}

We address here a minor but systematic flaw concerning the Stokes-type system, due to the divergence--free condition. We explain it for the resolvent problem, but the issue remains if one aims to perform localisation for the evolution problem.

\medbreak

The standard strategy to prove higher and $\L^p$-regularity of a solution $(\uu,\mathfrak{p})\in\H^{1,2}_{0,\sigma}(\Omega)\times\L^2(\Omega)$ to \eqref{eq:IntroSystStokes}, where  $\Omega$ has a locally flat boundary, is to consider $\eta\in\Ccinfty(\RR^n)$ sufficiently localised, so one can apply \eqref{eq:IntroNondivfreeResvEstRn+} to $(\eta \uu, \eta \mathfrak{p})$, and obtain
\begin{align*}
    |\lambda|\lVert \eta\uu\rVert_{\L^p(\RR^n_+)}  + \lVert ( \nabla^2(\eta\uu),\nabla({\eta\mathfrak{p}}))\rVert_{\L^p(\RR^n_+)} &\lesssim_{n,\mu} \lVert \eta\ff + \mathbf{F}(\eta,\uu,\nabla\uu,\mathfrak{p})\rVert_{\L^p(\RR^n_+)} 
    \\ &\quad\qquad+ |\lambda|\lVert \div(\eta\uu)\rVert_{\dot{\W}^{-1,p}(\RR^n_+)}  + \lVert \div(\eta\uu)\rVert_{\dot{\W}^{1,p}(\RR^n_+)}
\end{align*}
for a certain function $\mathbf F$.
To close the estimate, the main issue here is to estimate properly the part in the last line arising from the prescribed divergence due to the localisation argument.
The second term is easy to estimate with the identity $\div(\eta\uu) = \nabla \eta \cdot \uu $, thanks to $\div \uu =0$, so that
\begin{align*}
    \lVert \div(\eta\uu)\rVert_{\dot{\W}^{1,p}(\RR^n_+)} \leqslant \lVert \nabla \eta \cdot \uu\rVert_{{\W}^{1,p}(\RR^n_+)}=\lVert \nabla \eta \cdot \uu\rVert_{{\W}^{1,p}(\Omega)} \lesssim_{n} \lVert \nabla \eta \rVert_{{\W}^{1,\infty}(\RR^n)} \lVert \uu\rVert_{{\W}^{1,p}(\Omega)}.
\end{align*}
However, and this is our main point here, the estimate for the first term $\lVert\div(\eta\uu)\rVert_{\dot{\W}^{-1,p}(\RR^n_+)}$ for which we want to prove
\begin{align}\label{eq:IntroEstDivLocal}
    \lVert\div(\eta\uu)\rVert_{\dot{\W}^{-1,p}(\RR^n_+)}\lesssim_{\Omega,p,n}\lVert\uu\rVert_{{\W}^{-1,p}(\Omega)}.
\end{align}
This is more subtle and often overlooked in the literature. Indeed, the standard strategy might fail to apply here, which is an insidious systematic flaw. 

Indeed, writing $\div(\eta\uu)=\nabla \eta \cdot \uu$, one usually performs (often implicitly) one of the two following incomplete arguments, which are both similar in nature:
\begin{itemize}
    \item Claiming that homogeneous and inhomogeneous norms are equivalent for compactly support distributions and using Leibniz rule writing naively
    \begin{align*}
        \lVert\div(\eta\uu)\rVert_{\dot{\W}^{-1,p}(\RR^n_+)} \sim_{p,n,\Omega} \lVert\div(\eta\uu)\rVert_{{\W}^{-1,p}(\Omega)} \lesssim_{p,n,\Omega} \lVert \nabla \eta \rVert_{{\W}^{1,\infty}(\RR^n)} \lVert \uu \rVert_{{\W}^{-1,p}(\Omega)} ; \text{ or}
    \end{align*}
    \item First, claiming that $\dot{\W}^{-1,p}(\RR^{n}_+)$ is stable under multiplication by smooth cut-offs and, secondly, claiming again that inhomogeneous and homogeneous norms are equivalent for compactly supported distributions
    \begin{align*}
        \lVert\div(\eta\uu)\rVert_{\dot{\W}^{-1,p}(\RR^n_+)} = \lVert \nabla \eta \cdot \uu\rVert_{\dot{\W}^{-1,p}(\RR^n_+)} \lesssim_{p,n,\eta} \lVert  \uu\rVert_{\dot{\W}^{-1,p}(\RR^n_+)} \lesssim_{p,n,\eta,\Omega}  \lVert \uu \rVert_{{\W}^{-1,p}(\Omega)}.
    \end{align*}
\end{itemize}
These are valid strategies \textbf{only} when $-1>-n/p'$, \textit{i.e.} $p>\frac{n}{n-1}$,   and \textbf{always false otherwise}. For more details, we refer to the current introduction \textbf{Section~\ref{sec:intro2}}, \textbf{Lemma~\ref{lem:CompactEquiHomInhom}} and the subsequent remark. The issue remains for homogeneous Sobolev and Besov spaces $\dot{\H}^{s-1,p}$ and $\dot{\B}^{s-1}_{p,q}$ whenever  $s-1\leqslant -{n}/{p'}$. Fortunately, we shall see that this is not critical.

Recently, it has been circumvented by Geng and Shen \cite[Lemma~5.2]{GengShen2024} in a quite elementary way. They proved that \eqref{eq:IntroEstDivLocal} holds for all $p\in(1,\infty)$ with a suitable control on the continuity constant.

However, we need to deal here with the complete problem in our framework, incorporating a treatment of homogeneous Sobolev and Besov spaces of fractional orders in combination with a suitable change of coordinates. This is \textit{la raison d'être} of \textbf{Lemma~\ref{lem:preservingGengShen}} prior to the central result \textbf{Theorem~\ref{thm:StokesResolvent}}.

\subsubsection{Technical issue 3: Boundary value problems in endpoint Besov spaces on the half space}

For our analysis of the Stokes--Dirichlet problem, and more generally resolvent elliptic boundary value problem, we are due to study the problem
\begin{equation*}\tag{DS${}_\lambda$}
    \left\{ \begin{array}{rllr}
         \lambda \uu - \Delta \uu +\nabla \mathfrak{p} &= 0 \text{, }&&\text{ in } \RR^n_+\text{,}\\
        \div \uu &= g\text{, } &&\text{ in } \RR^n_+\text{,}\\
        \uu_{|_{\partial\RR^n_+}} &=0\text{, } &&\text{ on } \partial\RR^n_+\text{.}
    \end{array}\right.
\end{equation*}
Removing the divergence part in Step 3.2 from the proof of \textbf{Theorem~\ref{thm:MetaThmDirichletStokesRn+}}, one has to study
\begin{equation*}\tag{DS${}_\lambda$'}\label{eq:IntroEqDirStokesSystemtn+BVP}
    \left\{ \begin{array}{rllr}
         \lambda \ww - \Delta \ww +\nabla \pi &= 0 \text{, }&&\text{ in } \RR^n_+\text{,}\\
        \div \uu &= 0\text{, } &&\text{ in } \RR^n_+\text{,}\\
        \uu'_{|_{\partial\RR^n_+}} &= \hh'\text{, } &&\text{ on } \partial\RR^n_+\\
        {\uu_n} _{|_{\partial\RR^n_+}} &= 0\text{, } &&\text{ on } \partial\RR^n_+\text{.}
    \end{array}\right.
\end{equation*}
for some $\hh'$ depending linearly on $g$.

Following the strategy from \cite[Theorem~1.3,~Proof]{FarwigSohr1994} and \cite[Lemma~4.2.1,~Proof]{DanchinMucha2015} or more recently \cite[Eq~(3.10)]{GengShen2024}, one ends up with the following representation formula, on the tangential Fourier side
\begin{align}
    \F'[\ww'](\xi',x_n) &= 
    \left[ \frac{\sqrt{\lambda +|\xi'|^2}e^{-x_n\sqrt{\lambda +|\xi'|^2}}-|\xi'|e^{-x_n|\xi'|}}{\sqrt{\lambda +|\xi'|^2}-|\xi'|}\right] \frac{\xi'\,\prescript{t}{}{\xi'}}{|\xi'|^2} \F'[\hh'](\xi')\nonumber\\ &\qquad\qquad\qquad\qquad\qquad\qquad\qquad + \left[\mathbf{I}_{n-1}-\frac{\xi'\,\prescript{t}{}{\xi'}}{|\xi'|^2}\right] e^{-x_n\sqrt{\lambda +|\xi'|^2}}\F'[\hh'](\xi')\label{eq;IntroRepFormula1}\\
    \F'[\ww_n](\xi',x_n) &= \left[ \frac{e^{-x_n\sqrt{\lambda +|\xi'|^2}}-e^{-x_n|\xi'|}}{\sqrt{\lambda +|\xi'|^2}-|\xi'|}\right] i\xi'\cdot\F'[\hh'](\xi'),\label{eq;IntroRepFormula2}
\end{align}
and a similar expression for $\nabla\pi$. The standard strategy to obtain suitable $\L^p$-bounds  of such \textbf{parameter-dependent tangential Fourier multiplier} operators is the use of adapted H\"{o}rmander-Mikhlin-type results such as \cite[Lemma~2.4]{FarwigSohr1994}.

\medbreak

What seems to be a theory for such multipliers began to take shape in the work of Shibata and Shimizu \cite[Section~5]{ShibataShimizu2010}, Enomoto and Shibata \cite{EnomotoShibata2013}. It has been successfully applied for the linear theory of the Neumann--Stokes resolvent problem and the coupled Lame-hyperbolic resolvent system. But results and strategies of such type restrict the analysis to
\begin{align}\label{eq:IntroRescticLp}
    \L^p(\RR^n_+),\qquad 1<p<\infty,
\end{align}
see also the work of Shibata \cite{Shibata2015}, Denk and  Shibata \cite{DenkShibata2019} for other applications of such results. There are also known strategy to extend manually the analysis for the boundedness of such parameter-dependent tangential Fourier multiplier operators as given in \eqref{eq;IntroRepFormula1} and \eqref{eq;IntroRepFormula2} to Besov spaces $\B^{s}_{p,q}$, $p\in(1,\infty)$, see for instance the recent work by Kuo and Shibata \cite[Section~4.2]{KuoShibata2025}.

\medbreak

However, recall that we aim for results and a theory that encapsulate the endpoint homogeneous Besov spaces including the endpoint cases 
\begin{align*}
    \dot{\B}^{s}_{1,q}(\RR^{n}_+)\quad\text{ and }\quad\dot{\B}^{s}_{\infty,q}(\RR^{n}_+)
\end{align*}
for the Stokes--Dirichlet problem. Therefore, all of the strategies and results previously mentioned above fail to apply.

\medbreak

In order to reach our goals, we change the point of view: Instead of the Fourier multipliers representation as in \eqref{eq;IntroRepFormula1} and \eqref{eq;IntroRepFormula2}, we express $\ww$ in a purely operator theoretic way:
\begin{align*}
    \ww'(\cdot,x_n) &= \frac{1}{\lambda}[(\lambda\I-\Delta')^\frac{1}{2}+(-\Delta')^\frac{1}{2}]
    \Big[ (\lambda\I-\Delta')^\frac{1}{2}e^{-x_n(\lambda\I-\Delta')^\frac{1}{2}}-(-\Delta')^\frac{1}{2}e^{-x_n(-\Delta')^\frac{1}{2}}\Big][\I-\PP_{\RR^{n-1}}]\hh'\nonumber\\ &\qquad\qquad\qquad\qquad\qquad\qquad\qquad + \quad e^{-x_n(\lambda\I-\Delta')^\frac{1}{2}}\PP_{\RR^{n-1}}\hh'\\
    \ww_n(\cdot,x_n) &= \frac{1}{\lambda}[(\lambda\I-\Delta')^\frac{1}{2}+(-\Delta')^\frac{1}{2}]
    \Big[ e^{-x_n(\lambda\I-\Delta')^\frac{1}{2}}-e^{-x_n(-\Delta')^\frac{1}{2}}\Big] \div'\hh',
\end{align*}
and a similar expression for $\nabla \pi$. Due to the mapping properties of $(\lambda\I-\Delta')^\frac{1}{2}$ and $(-\Delta')^\frac{1}{2}$, one can essentially reduce the analysis to the boundedness of the semigroup $T_\lambda=(e^{-x_n(\lambda\I-\Delta')^\frac{1}{2}})_{x_n\geqslant 0}$, for $\lambda\in \{0\}\cup\{ \,e^{i\theta}\,,\,|\theta|\leqslant\mu\,\}$, $\mu\in(0,\pi)$, as a linear map
\begin{align}\label{eq:IntroPoissonSemigroupBesov}
    T_\lambda\,:\,\dot{\B}^{s-\sfrac{1}{p}}_{p,q}(\RR^{n-1})\longrightarrow \dot{\B}^{s}_{p,q}(\RR^{n}_+),\qquad s\in\RR,\,p,q\in[1,\infty].
\end{align}
To this end, we proceed as follows
\begin{itemize}
    \item In \textbf{Appendix~\ref{App:VectValFuncSpaces}}, we provide a quick construction of homogeneous vector-valued Besov spaces $\dot{\B}^{s}_{p,q}(\RR^n,\X)$ for a normed space $X$ and gather few of their basic properties for their restriction on $\RR^n_+$;
    \item In \textbf{Appendix~\ref{App:BesovSemigrpEst}}, we provide some abstract Besov-in-time (the role played by the $x_n$ variable) semigroup estimates for generic semigroups of injective linear operator $(\D(A),A)$ on a Banach space $\X$;
    \item In \textbf{Appendix~\ref{App:EquivNormBesov}}, we specify and apply the previous results for $A=(\lambda\I-\Delta')^\frac{1}{2}$  provided $\lambda\in \{0\}\cup\{ \,e^{i\theta}\,,\,|\theta|\leqslant\mu\,\}$, $\mu\in(0,\pi)$. We also derive other related equivalent norms for scalar-valued Besov spaces.
\end{itemize}
We do not claim by any means that our boundedness results are more general than the ones already available in the literature. However, they allow to reach endpoint spaces, and to treat directly the case of homogeneous function spaces of fractional order.  

\medbreak

From the results obtained in the aforementioned appendices, and especially the boundedness property \eqref{eq:IntroPoissonSemigroupBesov}, one can derive proper regularity estimates for boundary value resolvent problems in various appropriate forms by a combination of the standard trace theorem and elementary interpolation inequalities. For the simplest cases of applications, see for instance the proofs for the case of the resolvent problem for the Dirichlet  and the Neumann Laplacian  \textbf{Propositions~\ref{prop:DirResolventPbRn+} and~\ref{prop:NeumannResolventPbRn+}}.

\subsubsection{Solenoidal function spaces theory}

As mentioned in Point~\textit{(i)} from \textbf{Section~\ref{sec:IntroOverallStrat}}, the first step in our analysis consists in developing a precise theory of Sobolev and Besov spaces adapted to the resolvent problem under consideration. Since we work with a Stokes-type system, it is natural to introduce function spaces that incorporate the algebraic divergence--free constraint, namely
\medbreak
\resizebox{0.93\textwidth}{!}{$
\H^{s,p}_\sigma(\Omega):=\left\{ \uu\in\H^{s,p}(\Omega,\CC^n) \,:\, \div \uu =0 \right\}\text{, }{\H}^{s,p}_{0,\sigma}(\Omega):=\left\{ \uu\in\H^{s,p}(\RR^n,\CC^n) \,:\, \supp \uu \subset \overline{\Omega}\,\,\&\,\,\div \uu =0\right\}\text{,}$}
\medbreak
\resizebox{0.97\textwidth}{!}{$
\B^{s,\sigma}_{p,q}(\Omega):=\left\{ \uu\in\B^{s}_{p,q}(\Omega,\CC^n) \,:\, \div \uu =0 \right\}\text{, }\B^{s,\sigma}_{p,q,0}(\Omega):=\left\{ \uu\in\B^{s}_{p,q}(\RR^n,\CC^n) \,:\, \supp \uu \subset \overline{\Omega}\,\,\&\,\,\div \uu =0\right\}\text{,}
$}
\medbreak
\noindent and their homogeneous counterparts (with $\dot{\H}$ and $\dot{\B}$ replacing $\H$ and $\B$), whenever $\Omega$ is a special Lipschitz domain.  

\medbreak

For these spaces we collect known, refined, and in some cases new, properties. A global overview appears at the beginning of \textbf{Chapter~\ref{Sec:FunctionSpaceTheoryDivFree}}, together with \textbf{Sections~\ref{subsec:PoincBog} and~\ref{sec:ExtOpDivPreser}}. The refinements and novelty of our results concern mostly the treatment of endpoint or homogeneous function spaces on special Lipschitz domains.

\medbreak

Concerning Bogovski\u{\i}-type potential operators, \textbf{Theorems~\ref{thm:PotentialOpBddLipDom} and~\ref{thm:PotentialOpBddLipDomCcinfty}} address the case of bounded Lipschitz domains, whereas \textbf{Theorem~\ref{thm:PotentialOpSpeLipDom}} yields corresponding (but more involved) results for (homogeneous) function spaces on special Lipschitz domains.  

\medbreak

We also mention several extension operators on half-spaces that preserve the divergence--free property and provide corresponding boundedness results: \textbf{Theorem~\ref{thm:SteinsExtensionOpDivFree}} concerns special Lipschitz domains, and for a specific purpose, \textbf{Remark~\ref{rem:FlatExtOpDiv}} discusses the flat case $\RR^{n}_+$.  

\medbreak

Density results are established, for instance
\begin{align*}
\overline{\Ccinftydiv(\overline{\Omega})}^{\lVert \cdot \rVert_{\H^{s,p}(\Omega)}} = \H^{s,p}_\sigma(\Omega),
\qquad
\overline{\Ccinftydiv({\Omega})}^{\lVert \cdot \rVert_{\H^{s,p}(\RR^n)}} = \H^{s,p}_{0,\sigma}(\Omega),\qquad s\in\RR,
\end{align*}
and their analogues for (endpoint) Besov spaces are treated in \textbf{Theorem~\ref{thm:DivergenceFreeSpacesbddLipDensity}} for bounded Lipschitz domains and \textbf{Theorem~\ref{thm:DivergenceFreeSpacesSpeiclaLipDensity} }for homogeneous function spaces on special Lipschitz domains.

\medbreak

We furthermore establish interpolation identities such as
\begin{align*}
({\B}^{s_0,\sigma}_{p,q_0}(\Omega),{\B}^{s_1,\sigma}_{p,q_1}(\Omega))_{\theta,q}
=
({\H}^{s_0,p}_{\sigma}(\Omega),{\H}^{s_1,p}_{\sigma}(\Omega))_{\theta,q}
=
{\B}^{s,\sigma}_{p,q}(\Omega),\qquad s_0<s<s_1,\quad 1\leqslant p,q \leqslant \infty,
\end{align*}
together with their homogeneous analogues, proved in \textbf{Theorem~\ref{thm:InterpHomSpacesBddLip}} for bounded Lipschitz domains and in \textbf{Theorem~\ref{thm:InterpHomSpacesLip}} for special Lipschitz domains, with similar results for the spaces $\H^{s,p}_{0,\sigma}$ and $\B^{s,\sigma}_{p,q,0}$.

\medbreak

Finally, all these results can be combined with the identification
\begin{align*}
\B^{s,\sigma}_{p,q,0}(\Omega) \simeq  
\begin{cases}
\left\{ \uu\in\B^{s}_{p,q}(\Omega,\CC^n) \,:\, \div \uu =0 \ \text{and}\ \nu\cdot\uu_{{|_{\partial\Omega}}} = 0 \right\}, 
& \text{if } s\in(-1+\sfrac{1}{p},\sfrac{1}{p}),\\[0.3em]
\left\{ \uu\in\B^{s}_{p,q}(\Omega,\CC^n) \,:\, \div \uu =0 \ \text{and}\ \uu_{|_{\partial\Omega}} = 0 \right\}, 
& \text{if } s\in(\sfrac{1}{p},1+\sfrac{1}{p}),
\end{cases}
\end{align*}
obtained in \textbf{Proposition~\ref{prop:IdentifVanishingDivFree}}, with a similar results for other function spaces. This is also shown to remain valid in the framework of homogeneous function spaces on special Lipschitz domains.

All such results are mostly stated in the more general framework of differential forms, that also includes the corresponding analysis of $\curl$--free vector fields of arbitrary dimension. Such a point of view also allows to simplify the computations concerning some operators and algebraic relations, especially for extension and Bogovski\u{\i}-type operators. However, this may be at cost of clarity for the lay reader, non-initiated to such formalism.

\medbreak

A subsequent question, so that above definitions are meaningful is the validity of the integration by parts formula
\begin{align*}
    \big\langle  \uu\cdot\nu_{|_{\partial\Omega}}, \varphi_{|_{\partial\Omega}} \big\rangle_{\partial\Omega}=\big\langle  \uu, \nabla \varphi \big\rangle_{\Omega}+ \big\langle \div \uu, \varphi\big\rangle_{\Omega}
\end{align*}
when $\uu \in\dot{\H}^{s,p}(\Omega)^n$ is such that $\div \uu\in \dot{\H}^{s,p}(\Omega)$, $s\in(-1+\sfrac{1}{p},\sfrac{1}{p})$, $1<p<\infty$, and similarly if one considers the Besov space $\dot{\B}^{s}_{p,q}(\Omega)$. This is the main point of \textbf{Appendix~\ref{App:Traces}}, giving the appropriate results for inhomogeneous and homogeneous function spaces on the half space $\RR^n_+$, see \textbf{Theorem~\ref{thm:partialtracesDiffFormRn+}}, for inhomogeneous function spaces on Lipschitz domains as well as \textbf{Theorem~\ref{thm:partialtracesDiffFormLip}} and \textbf{Theorem~\ref{thm:partialtracesDiffFormHomSpeLip}} for homogeneous function spaces on special Lipschitz domains.

When $\Omega$ is a bounded Lipschitz domain, a standard duality argument (depending on the choice of function spaces), as explained in \textbf{Appendix~\ref{App:Traces}}, for \eqref{eq:VaritationnalFormulaDefPartialtraceDiv}, shows that there exists a linear functional $\mathfrak{n}_{\partial\Omega}(\uu)$ depending linearly on $\uu$ such that
\begin{align*}
    \big\langle  \mathfrak{n}_{\partial\Omega}(\uu), \varphi_{|_{\partial\Omega}} \big\rangle_{\partial\Omega}=\big\langle  \uu, \nabla \varphi \big\rangle_{\Omega}+ \big\langle \div \uu, \varphi\big\rangle_{\Omega}
\end{align*}
inducing a bounded linear map for $s\in(-1+\sfrac{1}{p},\sfrac{1}{p})$, $1<p<\infty$,
\begin{align*}
    \big\{ \vv\in{\H}^{s,p}(\Omega,\CC^n) \,:\, \div \vv \in\H^{s,p}(\Omega,\CC) \big\} &\longrightarrow \B^{s-{\sfrac{1}{p}}}_{p,p}(\partial\Omega)\\
    \uu \qquad& \longmapsto  \mathfrak{n}_{\partial\Omega}(\uu),
\end{align*}
and similarly for the Besov spaces ${\B}^{s}_{p,q}(\Omega)$, $s\in(-1+\sfrac{1}{p},\sfrac{1}{p})$, $1\leqslant p,q \leqslant \infty$.

The major issue not entirely addressed in the literature is the truthfulness of the distributional equality
\begin{align}\label{eq:IntroEqubastractTrace}
    \mathfrak{n}_{\partial\Omega}(\uu) = \nu \cdot \uu_{|_{\partial\Omega}},
\end{align}
which is \textbf{not} understood as a definition but as a consequence of the definition. Although the identity \eqref{eq:IntroEqubastractTrace} would not be surprising, it still remains to be checked in appropriate functional frameworks.

This is nearly unknown in the case of homogeneous and/or endpoint function spaces, even in the case of the flat half-space. For some results in this direction, see \cite{MitreaMitreaShaw2008} for the case of inhomogeneous function spaces on (bounded) Lipschitz domains,  \cite{FujiwaraYamazaki2007} for some partial results in the case of homogeneous function spaces on exterior domains, or \cite[Appendix]{Gaudin2023Hodge} for non-optimal and  partial results for homogeneous function spaces on the half-space.

Again, all results are stated in the general framework of differential forms such that they include the corresponding analysis of $\curl$--contrained vector fields of arbitrary dimension, with the adapted partial trace.

In order to obtain general results on partial traces in \textbf{Theorems~\ref{thm:partialtracesDiffFormRn+},~\ref{thm:partialtracesDiffFormLip}},~and~\textbf{\ref{thm:partialtracesDiffFormHomSpeLip}}, the scheme of proof is nonlinear. This is due to the fact that, throughout the whole document, we present groups of results that share common thematics exactly once. Furthermore, additional technicalities are involved when homogeneous function spaces on special Lipschitz domains (even the half space) are considered.

The path to provide the appropriate proofs is the following:
\begin{enumerate}
    \item For \textbf{Theorem~\ref{thm:partialtracesDiffFormRn+}} that concerns $\RR^{n}_+$. Follow the proof of Theorem~\ref{thm:partialtracesDiffFormRn+} for the half space, including homogeneous function spaces, up to the end of Step 4. This yields the functional $\uu\mapsto \mathfrak{n}_{\partial\RR^{n}_+}(\uu)$ with the appropriate boundedness properties, and dependency on some specific component of $\uu$;
    \item From the the exact component dependency of $\uu\mapsto \mathfrak{n}_{\partial\RR^{n}_+}(\uu)$, we can prove Proposition~\ref{prop:ExtOpDiffFormRn+}. It gives the existence of an appropriate extension operator
    \begin{align*}
    \E_{\mathcal{H}_\ast}\,:\,\big\{ \vv\in\dot{\H}^{s,p}(\RR^n_+)^n \,:\, \div \vv \in\dot{\H}^{s,p}(\RR^n_+) \big\} &\longrightarrow \big\{ \vv\in\dot{\H}^{s,p}(\RR^n)^n \,:\, \div \vv \in\dot{\H}^{s,p}(\RR^n) \big\},
\end{align*}
also valid for other relevant function spaces. In particular, smooth functions are dense. Moreover, such function spaces are an interpolation scale according to Corollary~\ref{cor:InterpolationDiffFormRn+};
    \item Resume the proof of \textbf{Theorem~\ref{thm:partialtracesDiffFormRn+}} from Step 5, and deduce \eqref{eq:IntroEqubastractTrace} for which one has $\nu=-\mathfrak{e}_n$. This is a consequence of the strong density of smooth functions arising naturally from  Proposition~\ref{prop:ExtOpDiffFormRn+}. The same applies to inhomogeneous function spaces. The case of Besov spaces ${\B}^{s}_{p,\infty}$ and $\dot{\B}^{s}_{p,\infty}$, follows by real interpolation thanks to Corollary~\ref{cor:InterpolationDiffFormRn+}. Remaining function spaces, such as $\L^\infty$, follow by natural embeddings;
    \item For \textbf{Theorem~\ref{thm:partialtracesDiffFormLip}} for inhomogeneous function spaces on Lipschitz domains when $s\geqslant 0$: one applies the standard duality argument. The case of function spaces of  order $s = 0$ such as $\L^p$, $1\leqslant p<\infty$, or $\C^{0}_{ub}$, the result follows by a density argument given in Corollary~\ref{cor:densityLpLipDiffForms} (which is a consequence of Proposition~\ref{prop:ExtOpDiffFormRn+}). By decreasing inclusion of Sobolev and Besov spaces when $s\geqslant 0$, \eqref{eq:IntroEqubastractTrace}  remains valid for non-negative regularity indices.
    
    \item For \textbf{Theorem~\ref{thm:partialtracesDiffFormLip}} for inhomogeneous function spaces on (bounded or Special) Lipschitz domains when $s< 0$:  by Proposition~\ref{prop:DiffFormDomainsclosedEtc..}, due to the boundedness of Bogoski\u{\i} operators given in Theorems~\ref{thm:PotentialOpBddLipDom} and \ref{thm:PotentialOpSpeLipDom}, one has that the family
    \begin{align*}
        \Big(\big\{ \vv\in{\B}^{s}_{p,q}(\Omega)^n \,:\, \div \vv \in{\B}^{s}_{p,q}(\Omega) \big\}\Big)_{-1+\sfrac{1}{p}<s<\sfrac{1}{p},\,1<p<\infty}
    \end{align*}
    is an interpolation scale. Therefore, an abstract density result for interpolation scales combined with case $s\geqslant 0$, yields \eqref{eq:IntroEqubastractTrace} when $s<0$;
    \item For \textbf{Theorem~\ref{thm:partialtracesDiffFormHomSpeLip}} which concerns the case of homogeneous function spaces on special Lipschitz domains: the case $s<0$ is a straightforward consequence of Theorem~\ref{thm:partialtracesDiffFormLip}, and the case $s>0$ is obtained by a variation of the usual duality argument, interpolation yielding the case $s=0$, thanks to Proposition~\ref{prop:DiffFormDomainsclosedEtc..}.
\end{enumerate}

\subsubsection{Technical aspects around non-complete function spaces}

When performing the analysis in our framework with homogeneous function spaces on the whole space or the half spaces, function spaces of high regularity are not complete. This has many consequences on how one should state results or prove ``standard properties''. The main consequence of the lack of completeness is the lack of ability to construct limits. Therefore, for each linear operator one should ensure that it is well-defined on the whole function space, without resorting to a restriction to a dense subspace. 

\medbreak

We illustrate this with two fundamental and simplified examples. For simplicity, we restrict ourselves to standard $\L^p$-based Sobolev spaces.
\begin{enumerate}
    \item  The Laplacian is a bounded linear operator
    \begin{align*}
        -\Delta\,:\,\dot{\W}^{2,p}(\RR^n)\longrightarrow \L^p(\RR^n),\quad 1<p<\infty,
    \end{align*}
    but it is only an isomorphism with a well-defined inverse
    \begin{align*}
        (-\Delta)^{-1}\,:\,\L^p(\RR^n)\longrightarrow \dot{\W}^{2,p}(\RR^n),\quad 1<p<n/2,
    \end{align*}
    due to the lack of completeness of $\dot{\W}^{2,p}(\RR^n)$ when $p\geqslant n/2$. However, if
    \begin{align*}
        -\Delta u = f\in \L^q(\RR^n)\cap \L^p(\RR^n),
    \end{align*}
    with $q\in(1,n/2)$ and $p\in(1,\infty)$, one still has $u=(-\Delta)^{-1}f\in \dot{\W}^{2,q}(\RR^n)\cap\dot{\W}^{2,p}(\RR^n)$. Thus, intersecting with a complete space makes the range complete, and we obtain the estimate
    \begin{align*}
        \lVert \nabla^2 u\rVert_{\L^p(\RR^n)}\lesssim_{p,n} \lVert f\rVert_{\L^p(\RR^n)},\quad f\in\L^q(\RR^n)\cap\L^p(\RR^n),\quad p,q>1,\quad q<n/2.
    \end{align*}
    The possible lack of completeness of $\dot{\W}^{2,p}(\RR^n)$ prevents us from extending the solution operator to all $f\in\L^p(\RR^n)$.

    \item We now consider the Dirichlet–Laplace problem on the flat half-space
    \begin{align*}
        -\Delta u = f,\quad \text{in }\RR^n_+,\quad \text{and}\quad u_{|_{\partial\RR^n_+}} = h\ \text{on }\partial\RR^n_+.
    \end{align*}
    By linearity, one may write $u=u_{\circ}+u_{\partial}$, where
    \begin{align*}
        -\Delta u_\circ = f,\quad \text{in }\RR^n_+,\quad \text{and}\quad {u_\circ}_{|_{\partial\RR^n_+}}=0\ \text{on }\partial\RR^n_+,
    \end{align*}
    and
    \begin{align*}
        -\Delta u_\partial =0,\quad \text{in }\RR^n_+,\quad \text{and}\quad {u_\partial}_{|_{\partial\RR^n_+}}=h\ \text{on }\partial\RR^n_+.
    \end{align*}
    By odd extension and solvability of the Laplacian on $\RR^n$, one has uniqueness for $u_\circ$ provided $f\in\L^p(\RR^n_+)$, and existence whenever $f\in\L^q(\RR^n_+)\cap\L^p(\RR^n_+)$ with $p,q>1$ and $q<n/2$.

    However, as soon as $h\in\dot{\W}^{2-\sfrac{1}{p},p}(\partial\RR^{n}_+)=\dot{\B}^{2-\sfrac{1}{p}}_{p,p}(\RR^{n-1})$, one has existence and uniqueness for $u_\partial$ for all $1<p<\infty$ without any further restriction, given by
    \begin{align*}
        u_\partial(\cdot,x_n) = e^{-x_n(-\Delta')^\frac{1}{2}}h\in \dot{\W}^{2,p}(\RR^n_+),
    \end{align*}
    thanks to Proposition~\ref{prop:PoissonSemigroup3}, where it is shown that $(e^{-x_n(-\Delta')^\frac{1}{2}})_{x_n\geqslant 0}$ takes values in the appropriate homogeneous Sobolev space without any completeness assumption, see Proposition~\ref{prop:DirPbRn+} for the complete solvability of the Dirichlet–Laplace problem on the half-space. 
\end{enumerate}

In the particular case of solving (linear) PDEs, the main issue is about existence rather than uniqueness. More generally, the question is whether, given a linear operator, we can make sense of its image in an appropriate sense \emph{beforehand}, that is, prior to establishing the relevant norm bounds (which, we recall, is unnecessary if the range lies in a complete normed vector space). For a linear operator $\T$ to be well-defined from a normed space $\X$ to $\Y(\RR^n_+)$, possibly a non-complete function space, we require
\begin{align*}
    \T(\X)\subset \S'_h(\RR^n_+).
\end{align*}
This is the main technical threshold, as it requires us to show that the elements can actually be computed in a meaningful way. Once this is achieved, it suffices to prove the norm bound from $\X$ to $\Y(\RR^n_+)$, and one may even restrict the proof of boundedness to a dense subspace, since the limit is already well-defined (i.e., there is no need to construct the limit). However, in several cases one needs to restrict to complete spaces, or to provide sufficient conditions ensuring that the range of the operator is well-defined. Hence the corresponding statement might become substantially more involved than what one would initially expect.

The statements or proofs of well-definedness and boundedness beyond a complete range rely either on
\begin{itemize}
    \item considering intersections in order to deal with a complete space, as suggested in Proposition~\ref{prop:DirPbRn+}, Proposition~\ref{prop:GenNeuPbRn+} and Theorem~\ref{thm:MetaThmSteadyDirichletStokesRn+} for the conditions of existence; 
    \item proving that the range or the limit makes sense prior to the estimate, see for instance Lemma~\ref{lem:ExtDirNeuRn+}, Proposition~\ref{prop:bddExtOpHomFreeDivSpeLip} (following a strategy established in \cite[Sections~3.3~\&~3.4]{Gaudin2023Lip}), Remark~\ref{rem:FlatExtOpDiv}, and Proposition~\ref{prop:PoissonSemigroup3}.
\end{itemize}
The most involved result in this regard is Theorem~\ref{thm:PotentialOpSpeLipDom}, together with the subsequent Remark~\ref{rem:convergencePotOpSpeLipNoncomplete}, concerning Bogovski\u{\i}-type operators on special Lipschitz domains, where both types of approaches are combined. It has strong consequences for the proofs of subsequent interpolation and density results such as Theorems~\ref{thm:DivergenceFreeSpacesSpeiclaLipDensity}~\&~\ref{thm:InterpHomSpacesLip}. For these statements the arguments become more intricate than initially expected when dealing with specific function spaces, such as endpoint inhomogeneous ones or homogeneous functions spaces beyond the regularity threshold.

\newpage
\section{Functional Setting}

\subsection{Notations}

\paragraph{General setting.}

Throughout this paper, the dimension is $n\geqslant 2$, and $\NN$ is the set of non-negative integers. For $a,b\in\RR$ with $a\leqslant b$, we write $\llb a,b\rrb:=[a,b]\cap\ZZ$.

\medbreak

For two real numbers $A,B\in\RR$, $A\lesssim_{a,b,c} B$ means that there exists a constant $C>0$ depending on ${a,b,c}$ such that $A\leqslant C B$. When both $A\lesssim_{a,b,c} B$ and $B \lesssim_{a,b,c} A$ hold, we simply write $A\sim_{a,b,c} B$. When there are many indices, we may write $A\lesssim_{a,b,c}^{d,e,f} B$ instead of $A\lesssim_{a,b,c,d,e,f} B$.

\medbreak

Let $\S(\RR^n,\CC)$ be the space of complex-valued Schwartz functions, and $\S'(\RR^n,\CC)$ its dual, called the space of tempered distributions. The Fourier transform on $\S'(\RR^n,\CC)$ is written $\F$, and it is pointwise defined for any $f\in\L^1(\RR^n,\CC)$ as
\begin{align*}
  \F f(\xi) :=\int_{\RR^n} f(x)\,e^{-ix\cdot\xi}\,\d x\text{, } \xi\in\RR^n\text{. }
\end{align*}
Additionally, for $p\in[1,\infty]$, the quantity $p'=\tfrac{p}{p-1}$ is its \textit{\textbf{H\"{o}lder conjugate}}.

\medbreak

For any $m\in\NN$, the map $\nabla^m\,:\,\S'(\RR^n,\CC)\longrightarrow \S'(\RR^n,\CC^{n^m})$ is defined as $\nabla^m u := (\partial^\alpha u)_{|\alpha|=m}$. The Laplace operator on $\RR^n$ is given as $\Delta= \partial_{x_1}^2+\partial_{x_2}^2+\ldots+\partial_{x_{n-1}}^2+\partial_{x_n}^2$.

\medbreak

 We introduce the operator $\nabla'$ which stands for the gradient on $\RR^{n-1}$ identified with the first $n-1$ variables of $\RR^n$, \textit{i.e.} $\nabla'=(\partial_{x_1}, \ldots, \partial_{x_{n-1}})$. Similarly, one defines $\Delta' = \partial_{x_1}^2+\partial_{x_2}^2+\ldots+\partial_{x_{n-1}}^2$. We denote by $(e^{-t(-\Delta')^\frac{1}{2}})_{t\geqslant0}$ the Poisson semigroup on $\RR^{n-1}$, replacing it either by $T_0$ or $(-\Delta_{\mathcal{D},\partial})^{-1}$ depending on context.

 \medbreak

If $\Omega$ is an open set in $\RR^n$, $\Ccinfty(\Omega,\CC)$ denotes the set of smooth, compactly supported functions in $\Omega$, and $\mathcal{D}'(\Omega,\CC)$ is its topological dual. For $p\in[1,\infty)$, $\L^p(\Omega,\CC)$  is the normed vector space of complex-valued (Lebesgue-) measurable functions whose $p$-th power is integrable with respect to the Lebesgue measure. $\S(\overline{\Omega},\CC)$ (\textit{resp.} $\Ccinfty(\overline{\Omega},\CC)$) stands for functions which are restrictions on $\Omega$ of elements of $\S(\RR^n,\CC)$ (\textit{resp.} $\Ccinfty(\RR^n,\CC)$). Unless the contrary is explicitly stated, we always identify $\L^p(\Omega,\CC)$ (resp. $\Ccinfty(\Omega,\CC)$) as the subspace of functions in $\L^p(\RR^n,\CC)$ (resp. $\Ccinfty(\RR^n,\CC)$) supported in $\overline{\Omega}$ (resp. $\Omega$) through the extension by $0$ outside $\Omega$. $\L^\infty(\Omega,\CC)$ stands for the space of essentially bounded (Lebesgue-) measurable functions.

For $\Omega'\subset\Omega$ being two open sets of $\RR^n$, the restriction homomorphism from $\Distrib(\Omega')$ to $\Distrib(\Omega)$ is denoted by the linear operator $\R_{\Omega'}$. When applied to an element $u\in\L^{1}_{\textrm{loc}}(\Omega)$, one might also write $u_{|_{\Omega'}}:=\R_{\Omega'}u\in\L^{1}_{\textrm{loc}}(\Omega')$.

\medbreak

For $s\in\RR$, $p\in[1,\infty)$, $\ell^p_s(\ZZ,\CC)$ denotes the normed vector space of $p$-summable sequences of complex numbers with respect to the counting measure $2^{ksp}\d k$;  $\ell^\infty_s(\ZZ,\CC)$ denotes the vector space of sequences $(x_k)_{k\in\ZZ}$ such that $(2^{ks}x_k)_{k\in\ZZ}$  is bounded.
More generally, when $\mathrm{\X}$ is a Banach space, for $p\in[1,\infty]$, one can also consider $\L^p(\Omega,\mu,\mathrm{\X})$ denotes the space of (Bochner-)measurable functions $u\,:\,\Omega\longrightarrow \mathrm{\X}$, such that $t\mapsto\lVert u(t)\rVert_\mathrm{\X} \in \L^p(\Omega,\mu,\RR)$. Similarly, one can consider $\ell^p_s(\ZZ,\mathrm{\X})$. We set $\C^0(\Omega,\mathrm{\X})$ to be the space of continuous functions on $\Omega\subset \RR^n$ with values in $\mathrm{\X}$. The subspace $\C^0_{ub}(\RR^n,\mathrm{\X})$  consists of bounded and uniformly continuous functions. $\C^0_0(\RR^n,\mathrm{\X})$ denotes the subspace of continuous functions that vanish at infinity. For $\mathscr{C}\in\{\C,\C_c,\C_{ub},\C_0\}$, we define $\mathscr{C}^0(\overline{\Omega}, \mathrm{\X})$ as the vector space of continuous functions on $\overline{\Omega}$ which are restrictions of elements that belongs to $\mathscr{C}^0(\RR^n, \mathrm{\X})$. The same goes for the function spaces of continuously differentiable functions up to the order $k\in\NN^\ast$: $\mathscr{C}^k(\overline{\Omega}, \mathrm{\X}):=\{\,u\in \mathscr{C}^{k-1}(\overline{\Omega}, \mathrm{\X})\,:\,\nabla u \in \mathscr{C}^{k-1}(\overline{\Omega}, \mathrm{\X})\,\}$. For $k\in\NN$, the space of continuously differentiable functions up to the order $k$, such that the derivatives of order $k$ are Lipschitz is denoted by $\mathscr{C}^{k,1}(\overline{\Omega}, \mathrm{\X}):=\{\,u\in \mathscr{C}^{k}(\overline{\Omega}, \mathrm{\X})\,:\,\nabla^{k+1} u \in \L^{\infty}({\Omega}, \mathrm{\X})\,\}$. We set $\mathscr{C}^{k,0}:=\mathscr{C}^k$. We keep similar notations with $\Omega$ instead of $\overline{\Omega}$, where we lose information concerning the behavior near the boundary. Finally, for $k\in\NN$, $\alpha\in[0,1]$, $\mathfrak{a}\in\{ub,0\}$, we define
\begin{align*}
    \C^{k,\alpha}_{\mathfrak{a},0}({\Omega}):=\{\,u\in \C^{k,\alpha}_{\mathfrak{a}}(\RR^n)\,:\,\supp u \subset \overline{\Omega}\,\},
\end{align*}
where, when $0<\alpha<1$, $\C^{k,\alpha}_{b}(\Omega)$ stands for smooth functions $u\in \C^{k}_{ub}(\overline{\Omega})$ such that
\begin{align*}
    \sup_{\substack{x,y\in\Omega\\x\neq y}} \frac{|\nabla^ku(x)-\nabla^ku(y)|}{|x-y|^\alpha} <\infty.
\end{align*}
We call the previous quantity the $(k,\alpha)$-\textbf{H\"{o}lderian constant} of $u$. If $k=0$, we reduce the name to $\alpha$-\textbf{H\"{o}lderian constant}.
Finally, we also introduce the little H\"{o}lder space $\mathcal{C}^{k,\alpha}(\Omega)$, which stands for the elements $u\in\C^{k,\alpha}_b(\Omega)$, such that
\begin{align*}
     \frac{|\nabla^ku(x)-\nabla^ku(y)|}{|x-y|^\alpha}\xrightarrow[{\substack{|x-y|\rightarrow 0\\x,y\in\Omega\\x\neq y}}]{} 0.
\end{align*}
When it is meaningful, for $I\in\{(0,T)\,[0,T),\,[0,T]\}$, $\C^{\mathrm{w}\ast}(I;\X)$ stands for functions $t\longmapsto u(t)$ with value in a Banach space $\X$, which are weakly-$\ast$ continuous  on $I$.

\medbreak

For $\Omega$ an open set of $\RR^n$, we say that $\Omega$ is a special Lipschitz domain, if there exists, up to a rotation, a globally Lipschitz function $\phi\,:\,\RR^{n-1}\longrightarrow\RR$, such that
\begin{align*}
    \Omega =\{\,(x',x_n)\in\RR^{n-1}\times\RR\,:\, x_n>\phi(x')\,\}\text{.}
\end{align*}
In other words, a special Lipschitz domain of $\RR^n$ is the epigraph of real valued Lipschitz function defined on $\RR^{n-1}$. Furthermore, if $\phi\in\C_c^{0,1}(\RR^{n-1})$, $\Omega$ is said to be a small special Lipschitz  (\textit{resp.} a small special $\C^1$) domain. In this case, we denote by
\begin{align*}
    \begin{array}{cccccccc}
     \Psi &:& \RR^{n-1}\times\RR &\longrightarrow& \RR^{n-1}\times\RR  \\
     & &(x',x_n) &\longmapsto & (x',x_n+\phi(x')) 
\end{array}
\end{align*}
the associated bi-Lipschitz map such that $\Psi(\RR^n_+)=\Omega$ and $\Psi(\partial\RR^n_+)=\partial\Omega$.

From now on, and until the end of this paper $\Omega$ stands for a domain of $\RR^n$, with at least, if not empty, Lipschitz boundary. Recall briefly that $\partial\RR^n = \emptyset$, and $\partial\RR^n_+ = \RR^{n-1}\times \{0\}$. We also recall that the outer normal unit at $\partial\RR^n_+$ is ${\nu} = -\mathfrak{e}_n$, where $(\mathfrak{e}_k)_{k\in\llb 1, n\rrb}$ is the canonical basis of $\RR^n$, identified with its dual basis denoted by $(\d x_{k})_{k\in\llb 1,n\rrb}$, where $\d x_{k}(\mathfrak{e}_j) = \mathbbm{1}_{\{k\}}(j)$, $(k,j)\in\llb 1,n\rrb^2$.

Given $\mathcal{L}$ an elliptic differential operator on a domain $\Omega$ of $\RR^n$, equipped with a system of boundary conditions $\mathcal{B}$, when it  is meaningful, we denote by $\mathcal{L}_{\mathcal{B},\partial}^{-1}$ the solution operator to
\begin{equation*}
    \left\{ \begin{array}{rllr}
         \mathcal{L}u &= 0 \text{, }&&\text{ in } \Omega\text{,}\\
        \mathcal{B} u_{|_{\partial\Omega}} &=g\text{, } &&\text{ on } \partial\Omega\text{.}
    \end{array}
    \right.
\end{equation*}

\medbreak

\subsection{Interpolation theory}

For more details, the interested reader is invited to consult \cite{BerghLofstrom1976,bookTriebel1978,bookLunardiInterpTheory} or \cite[Section~1.3]{EgertPhDThesis2015}.

Let $(\X,\left\lVert\cdot\right\rVert_\X)$ and $(\Y,\left\lVert\cdot\right\rVert_\Y)$ be two normed vector spaces. We write $\X\hookrightarrow \Y$ to indicate that $\X$ embeds continuously into $\Y$. We briefly recall the basics of interpolation theory. If there exists a Hausdorff topological vector space $Z$, such that $\X,\Y\subset Z$, then $\X\cap \Y$ and $\X+\Y$ are normed vector spaces with their canonical norms. One can define the $K$-functional of $z\in \X+\Y$, for any $t>0$ by
\begin{align*}
    K(t,z,\X,\Y) := \underset{\substack{(x,y)\in \X\times \Y,\\ z=x+y}}{\inf}\left({\left\lVert{x}\right\rVert_{\X}+t\left\lVert{y}\right\rVert_{\Y}}\right)\text{. }
\end{align*}
This enables us to construct, for any $\theta\in(0,1)$, $q\in[1,\infty]$, the real interpolation spaces between $\X$ and $\Y$ with indices $\theta,q$ as
\begin{align*}
    (\X,\Y)_{\theta,q} := \left\{\, x\in \X+\Y\,\Big{|}\,t\longmapsto t^{-\theta}K(t,x,\X,\Y)\in\L^q_\ast(\mathbb{R}_+)\,\right\}\text{, }
\end{align*}
where $\L^q_\ast(\mathbb{R}_+):=\L^q(0,\infty,\mathrm{d}t/t)$. When $q=\infty$, one can consider
\begin{align*}
    (\X,\Y)_{\theta} := \left\{\, x\in \X+\Y\,\Big{|}\,\lim_{t,t^{-1}\rightarrow0_+} t^{-\theta}K(t,x,\X,\Y)=0\,\right\}\text{, }
\end{align*}
endowed with the norm $\lVert \cdot \rVert_{(\X,\Y)_{\theta,\infty}}$.

If moreover we assume that $\X$ and $\Y$ are complex Banach spaces, one can consider $\mathrm{F}(\X,\Y)$ the set of all continuous functions $f:\overline{S}\longmapsto \X+\Y$, $S$ being the strip of complex numbers whose real part is between $0$ and $1$, with $f$ holomorphic in $S$, and such that
\begin{align*}
    t\longmapsto f(it)\in \mathrm{C}^0_b(\mathbb{R},\X) \quad\text{ and }\quad t\longmapsto f(1+it)\in \mathrm{C}^0_b(\mathbb{R},\Y)\text{. }
\end{align*}
We can endow the space $\mathrm{F}(\X,\Y)$ with the norm
\begin{align*}
    \left\lVert{f}\right\rVert_{\mathrm{F}(\X,\Y)}:=\max\left(\underset{t\in\mathbb{R}}{\mathrm{sup}} \left\lVert {f(it)}\right\rVert_{\X},\underset{t\in\mathbb{R}}{\mathrm{sup}} \left\lVert {f(1+it)}\right\rVert_{\Y}\right)\text{, }
\end{align*}
which makes $\mathrm{F}(\X,\Y)$ a Banach space since it is a closed subspace of $\mathrm{C}^0(\overline{S},\X+\Y)$.
Hence, for $\theta\in(0,1)$, the normed vector space given by
\begin{align*}
    [\X,\Y]_{\theta} &:= \left\{\,f(\theta)\,\big{|}\,f\in \mathrm{F}(\X,\Y)\,\right\} 
    \text{, }\\
    \left\lVert{x}\right\rVert_{[\X,\Y]_{\theta}} &:= \underset{\substack{f\in \mathrm{F}(\X,\Y),\\ f(\theta)=x}}{\inf} \left\lVert{f}\right\rVert_{\mathrm{F}(\X,\Y)}\text{, }
\end{align*}
is a Banach space called the complex interpolation space between $\X$ and $\Y$ associated with $\theta$. 

We present here a useful and powerful argument to compute concrete real and complex interpolation spaces when one can realize a couple of normed vector spaces as some subspace of another one in a proper way (explained below). This procedure is called \textbf{the retraction and co-rectraction argument}. See for instance \cite[Sections~2.1~\&~6.4]{BerghLofstrom1976}, \cite[Section~1.2]{bookTriebel1978} and \cite[Section~1.3.1]{EgertPhDThesis2015}.

\begin{theorem}[ {\cite[Theorem,~Section~1.2.4]{bookTriebel1978}} ]\label{thm:RetractionThm} Let $(\X_0,\X_1)$ and $(\Y_0,\Y_1)$ be two pairs of compatible couples, and let
\begin{align*}
    \mathfrak{R}\,:\,\Y_0+\Y_1 \longrightarrow \X_0+\X_1\quad\text{and}\quad \mathfrak{E}\,:\,\X_0+\X_1 \longrightarrow \Y_0+\Y_1
\end{align*}
be two bounded linear operators such that $\mathfrak{R}\mathfrak{E}=\I_{\X_0+\X_1}$, \textit{i.e. } $\mathfrak{E}$ is a right bounded inverse for $\mathfrak{R}$.

Then for all $\theta\in(0,1)$, $p\in[1,\infty]$, one has
\begin{align*}
    (\X_0,\X_1)_{\theta,p}=\mathfrak{R}(\Y_0,\Y_1)_{\theta,p}\text{.}
\end{align*}
The result still holds, replacing $(\cdot,\cdot)_{\theta,p}$ by $(\cdot,\cdot)_{\theta}$.

Furthermore, if $(\X_0,\X_1)$ and $(\Y_0,\Y_1)$ are two interpolation couples of complex Banach spaces, then the result still holds for $[\cdot,\cdot]_{\theta}$, instead of $(\cdot,\cdot)_{\theta,p}$.
\end{theorem}

In practice, the pair of retraction and co-retraction $(\mathfrak{E},\mathfrak{R})$ is often an extension operator and a restriction operator $(\E,\R)$ or and injection (canonical embedding) and projection operator $(\iota,\P)$.

\subsection{Function spaces on the whole space and on Lipschitz domains}\label{subsec:FuncSpaceRnLipdomains}

\medbreak

\subsubsection{On the whole space.}\label{subsec:SobolevBesovRn} First, wet set $\S'_h(\RR^n)$, to be the subspace of tempered distributions $u\in\S'(\RR^n)$, such that
\begin{align*}
   \forall \Theta\in\Ccinfty(\RR^n),\,\lVert \Theta(\lambda\mathfrak{D})u\rVert_{\L^\infty(\RR^n)}\xrightarrow[\lambda\rightarrow\infty]{} 0.
\end{align*}
Let  $\varphi\in \mathrm{C}_c^\infty(\RR^n)$, a radial, real-valued, non-negative function such that
\begin{itemize}[label={$\bullet$}]
    \item $\supp \varphi \subset \B_{4/3}(0)$;
    \item ${\varphi}_{|_{\B_{3/4}(0)}}=1$.
\end{itemize}
We define the following functions for any $j\in\ZZ$ and for all $\xi\in\RR^n$,
\begin{align*}
    \varphi_j(\xi):=\varphi(2^{-j}\xi)\text{, }\qquad \psi_j(\xi) := \varphi_{j}(\xi/2)-\varphi_{j}(\xi)\text{,}
\end{align*}
and we consider the following two families of operators associated with their Fourier multipliers:
\begin{itemize}[label={$\bullet$}]
    \item The \textit{\textbf{homogeneous}} family of Littlewood-Paley dyadic decomposition operators $(\dot{\Delta}_j)_{j\in\ZZ}$, where
    \begin{align*}
        \dot{\Delta}_j := \F^{-1}\psi_j\F = (\F^{-1}\psi_j)\ast \text{,}
    \end{align*}
    \item The \textit{\textbf{inhomogeneous}} family of Littlewood-Paley dyadic decomposition operators $({\Delta}_k)_{k\in\ZZ}$, where
    \begin{align*}
       {\Delta}_{-1} := \F^{-1}\varphi\F = (\F^{-1}\varphi)\ast\text{,}
    \end{align*}
    $\Delta_k:=\dot{\Delta}_k$ for any $k\geqslant 0$, and $\Delta_k:=0$ for any $k\leqslant-2$.
    \item The low frequency cut-off operators $(\dot{S}_k)_{k\in\ZZ}$, where
    \begin{align*}
       \dot{S}_k := \F^{-1}\varphi_k\F\text{.}
    \end{align*}
    Note that $\dot{S}_0 = \Delta_{-1}$.
\end{itemize}

For $(\dot{\Delta}_{j})_{j\in\ZZ}$, $({\Delta}_{k})_{k\in\NN}$, a given \textbf{Littlewood-Paley decomposition}, for all $p,q\in[1,\infty]$, $s\in\RR$, we set:
\begin{align*}
    \left\lVert {u} \right\rVert_{\H^{s,p}(\RR^n)} := \left\lVert {(\I-\Delta)^\frac{s}{2} u} \right\rVert_{\L^{p}(\RR^n)}\text{ and } \left\lVert {u} \right\rVert_{\dot{\H}^{s,p}(\RR^n)} := \Big\lVert \sum_{j\in\ZZ} (-\Delta)^\frac{s}{2}\dot{\Delta}_{j} u  \Big\rVert_{{\L}^{p}(\RR^n)}\text{, }
\end{align*}
where for $u\in\S'(\RR^n)$,
\begin{align*}
    (-\Delta)^\frac{s}{2}\dot{\Delta}_j u := \F^{-1}|\xi|^s\F\dot{\Delta}_j u\text{ and } (\I-\Delta)^\frac{s}{2}u = \F^{-1}(1+|\xi|^2)^\frac{s}{2}\F u.
\end{align*}
The standard Sobolev norms are given by
\begin{align*}
    \lVert {u} \rVert_{\dot{\W}^{k,p}(\RR^n)}:=\lVert\nabla^k u\rVert_{\L^p(\RR^n)}\text{ and } \lVert {u} \rVert_{{\W}^{k,p}(\RR^n)}:=\lVert u\rVert_{\L^p(\RR^n)} +\lVert\nabla^k u\rVert_{\L^p(\RR^n)}.
\end{align*}
Finally the Besov norms are
\begin{center}
     \resizebox{0.97\textwidth}{!}{$
        \left\lVert u \right\rVert_{\B^{s}_{p,q}(\RR^n)}:= \left\lVert(2^{ks}\left\lVert {\Delta}_k u \right\rVert_{\L^{p}(\RR^n)})_{k\in\NN}\right\rVert_{\ell^{q}(\NN)}\text{ and }
    \left\lVert u \right\rVert_{\dot{\B}^{s}_{p,q}(\RR^n)} := \left\lVert \big(2^{js}\lVert \dot{\Delta}_j u \rVert_{\L^{p}(\RR^n)}\big)_{j\in\ZZ}\right\rVert_{\ell^{q}(\ZZ)}\text{.}$}
\end{center}
From those given norms we can define the following normed vector spaces
\begin{itemize}[label={$\bullet$}]
    \item the inhomogeneous and homogeneous Sobolev (Bessel and Riesz potential) spaces,\\
    \resizebox{0.93\textwidth}{!}{$
        \H^{s,p}(\RR^n)=\left\{ u\in\S'(\RR^n) \,:\, \left\lVert {u} \right\rVert_{\H^{s,p}(\RR^n)}<\infty \right\}\text{, }\dot{\H}^{s,p}(\RR^n)=\left\{ u\in\S'_h(\RR^n) \,:\, \left\lVert {u} \right\rVert_{\dot{\H}^{s,p}(\RR^n)}<\infty \right\}\text{;}$}
    \item the standard inhomogeneous and homogeneous Sobolev spaces,\\
    \resizebox{0.93\textwidth}{!}{$
        \mathrm{\W}^{k,p}(\RR^n)=\left\{ u\in\S'(\RR^n) \,:\, \left\lVert {u} \right\rVert_{{\W}^{k,p}(\RR^n)}<\infty \right\}\text{, }\dot{\W}^{k,p}(\RR^n)=\left\{ u\in\S'_h(\RR^n) \,:\, \lVert {u} \rVert_{\dot{\W}^{k,p}(\RR^n)}<\infty \right\}\text{;}$}
    \item and the inhomogeneous and homogeneous Besov spaces,\\
    \resizebox{0.93\textwidth}{!}{$\B^{s}_{p,q}(\RR^n)=\left\{ u\in\S'(\RR^n) \,:\, \left\lVert {u} \right\rVert_{\B^{s}_{p,q}(\RR^n)}<\infty \right\}\text{, }\dot{\B}^{s}_{p,q}(\RR^n)=\left\{ u\in\S'_h(\RR^n) \,:\, \left\lVert {u} \right\rVert_{\dot{\B}^{s}_{p,q}(\RR^n)}<\infty \right\}\text{.}
    $}
\end{itemize}
It is common to set the Sobolev-Soblodeckij spaces to be for all $s\in\RR\setminus\ZZ$, $p\in[1,\infty]$,
\begin{align*}
    \W^{s,p}(\RR^n):=\B^{s}_{p,p}(\RR^n)\text{ and }\dot{\W}^{s,p}(\RR^n):=\dot{\B}^{s}_{p,p}(\RR^n).
\end{align*}

We also consider the following endpoint function spaces:
 for $p,q\in[1,\infty]$, $s\in\RR$:
\begin{align*}
    \dot{\mathcal{B}}^{s}_{p,\infty}(\RR^n):=&\Big\{u\in \dot{\B}^{s}_{p,\infty}(\RR^n)\,:\, \lim\limits_{|j|\rightarrow\infty} 2^{js} \lVert\dot{\Delta}_{j}u\rVert_{\L^p(\RR^n)}=0\Big\},\\
    \dot{\B}^{s,0}_{\infty,q}(\RR^n):=&\Big\{u\in \dot{\B}^{s}_{\infty,q}(\RR^n)\,:\, (\dot{\Delta}_ju)_{j\in\ZZ}\subset\mathrm{C}_0^0(\RR^n)\Big\},\\
    \dot{\mathcal{B}}^{s,0}_{\infty,\infty}(\RR^n):=& \dot{\mathcal{B}}^{s}_{\infty,\infty}(\RR^n)\cap\dot{\B}^{s,0}_{\infty,\infty}(\RR^n).
\end{align*}
as well as their inhomogeneous counterparts.  It is well known that, with equivalence of norms, $$\H^{m,p}(\RR^n)=\W^{m,p}(\RR^n) \quad\text{as well as}\quad\dot{\H}^{m,p}(\RR^n)=\dot{\W}^{m,p}(\RR^n)$$ for all $m\in\NN$, $p\in(1,\infty)$. And when $s>0$, $p\in[1,\infty)$, $\dot{\X}^{s,p}(\RR^n)\cap\L^p(\RR^n)=\X^{s,p}(\RR^n)$, with equivalence of norms, for $\X^{s,p}$ to be either $\B^{s}_{p,q}$ or $\H^{s,p}$, assuming $p>1$ for the latter. For Besov spaces with $p=\infty$, the equivalence of norms for intersections still holds, but the equality of sets requires to take the intersection with $\S'_h(\RR^n)$, and one also has with equivalence of norms
\begin{align*}
    \B^{k+\alpha}_{\infty,\infty}(\RR^{n})&=\C^{k,\alpha}_{b}(\RR^n),\qquad  &\BesSmo^{k+\alpha}_{\infty,\infty}(\RR^{n})= \mathcal{C}^{k,\alpha}_{b}(\RR^n),\\
    \B^{k+\alpha,0}_{\infty,\infty}(\RR^{n})&=\C^{k}_0\cap\C^{k,\alpha}_{b}(\RR^n),\qquad  &\BesSmo^{k+\alpha,0}_{\infty,\infty}(\RR^{n})= \mathcal{C}^{k,\alpha}_{0}(\RR^n),
\end{align*}
whenever $k\in\NN$, $\alpha\in(0,1)$.

We also set for any $k\in\NN$
\begin{align*}
    \L^\infty_h(\RR^n):=\S'_h(\RR^n)\cap\L^\infty(\RR^n)\quad\text{and}\quad \C^{k}_{ub,h}(\RR^n):=\S'_h(\RR^n)\cap\C^{k}_{ub}(\RR^n),
\end{align*}
which are respectively closed subspaces of $\L^\infty(\RR^n)$ and $\C^{k}_{ub}(\RR^n)$. We recall that for homogeneous function spaces over $\L^\infty$, one has the following set equalities for all $s>0$, all $k\in\mathbb{N}$ and all $q\in[1,\infty]$:
\begin{align*}
    \dot{\B}^{s}_{\infty,q}(\RR^n)={\B}^{s}_{\infty,q}(\RR^n)\cap\S'_h(\RR^n),\,\text{ and }\,\dot{\W}^{k,\infty}(\RR^n)={\W}^{k,\infty}(\RR^n)\cap\S'_h(\RR^n).
\end{align*}

\textbf{A major issue of homogeneous function spaces, which stands also for their major feature,} is that  to be meaningful in the study of non-linear partial differential equations and boundary value problems, they are no longer complete and \textbf{one cannot consider their completion.} See \cite[Introduction]{Gaudin2023Lip} for more details. However, we still have the next result

\begin{theorem}\label{thm:Compelteness} Let $p,q\in[1,\infty]$, $s\in\RR$. Then the spaces
\begin{align*}
    \dot{\B}^{s}_{p,q}(\RR^n)\text{, }\dot{\W}^{s,p}(\RR^n)\text{ and }\dot{\H}^{s,p}(\RR^n),
\end{align*}
are complete if and only if
\begin{align*}\tag{$\mathcal{C}_{s,p,q}$}\label{AssumptionCompletenessExponents}
    \left[ s<\frac{n}{p} \right]\text{ or }\left[q=1\text{ and } s\leqslant\frac{n}{p} \right]\text{, }
\end{align*}
assuming $p=q$ in the case of Sobolev spaces. Furthermore, the intersection of any space with a complete one ensures completeness of the resulting space.
\end{theorem}

The loss of completeness actually relies on the fact that function spaces of negative regularity lack test functions. For instance, is easy to check that
\begin{align*}
    \Ccinfty(\RR^n)\not\subset\dot{\H}^{-\frac{n}{2},2}(\RR^n),
\end{align*}
just taking one smooth compactly supported function with integral $1$. Therefore, any completion of it embeds in a space \textbf{larger} than $\mathcal{D}'(\RR^n)$. This issue relates to the loss of control over low frequencies, called infrared divergence. However, the intersection $\Ccinfty\cap\X^{s,p}(\RR^n)$ is still dense when expected (\textit{e.g.}, when $p,q<\infty$, one also has the density result for $\B^{s,0}_{\infty,q}$, $\BesSmo^{s}_{p,\infty}$ and $\BesSmo^{s,0}_{\infty,\infty}$. Also $\C^\infty\cap\B^{s}_{\infty,q}$ is dense in $\B^{s}_{\infty,q}$ and it goes similarly for their homogeneous counterpart). Otherwise for homogeneous function spaces the standard universal dense subspace of smooth functions is $\S_0(\RR^n)$, which stands for the space of Schwartz function such that their Fourier transform is compactly supported away from 0, that is
\begin{align*}
    \S_0(\RR^n):=\{\,u\in\S(\RR^n)\,:\,0\notin\supp \F u ,\text{ and }\supp \F u\text{ is compact}\,\}.
\end{align*}
See \cite[Proposition~2.8]{Gaudin2023Lip} for more details.

A good news is that, regardless of the completeness, for $s_0<s<s_1$, $p,q\in[1,\infty]$, we always have the real interpolation identities
\begin{align*}
    {\B}^{s}_{p,q}(\RR^n)=({\X}^{s_0,p}(\RR^n),{\Y}^{s_1,p}(\RR^n))_{\theta,q}\text{ and }\dot{\B}^{s}_{p,q}(\RR^n)=(\dot{\X}^{s_0,p}(\RR^n),\dot{\Y}^{s_1,p}(\RR^n))_{\theta,q}
\end{align*}
where $\X,\Y\in\{\H,\B_{\cdot,r}, r\in[1,\infty]\}$. See \cite[Theorem~2.12]{Gaudin2023Lip} for more details. Complex interpolation is also available but requires completeness to be defined.

One also has the following, nice duality relations:
\begin{itemize}
    \item $(\H^{-s,p'}(\RR^n))' = \H^{s,p}(\RR^n)$, $s\in\RR$, $1<p<\infty$;
    \item $(\B^{-s}_{p',q'}(\RR^n))' = \B^{s}_{p,q}(\RR^n)$,$s\in\RR$, $p,q\in(1,\infty]$;
    \item $(\BesSmo^{-s}_{p',\infty,0}(\RR^n))' = \B^{s}_{p,1}(\RR^n)$, $s\in\RR$, $p\in(1,\infty]$;
    \item $(\B^{-s,0}_{\infty,q',0}(\RR^n))' = \B^{s}_{1,q}(\RR^n)$, $s\in\RR$, $q\in(1,\infty]$;
    \item $(\BesSmo^{-s,0}_{\infty,\infty,0}(\RR^n))' = \B^{s}_{1,1}(\RR^n)$, $s\in\RR$.
\end{itemize}
and for homogeneous function spaces
\begin{itemize}
    \item $(\dot{\H}^{-s,p'}(\RR^n))' = \dot{\H}^{s,p}(\RR^n)$, $s\in\RR$, $1<p<\infty$, $s<\sfrac{n}{p}$;
    \item $(\dot\B^{-s}_{p',q'}(\RR^n))' = \dot\B^{s}_{p,q}(\RR^n)$, $s\in\RR$, $p,q\in(1,\infty]$, with \eqref{AssumptionCompletenessExponents};
    \item $(\dot\BesSmo^{-s}_{p',\infty}(\RR^n))' = \dot\B^{s}_{p,1}(\RR^n)$, $s\leqslant \sfrac{n}{p}$, $p\in(1,\infty]$;
    \item $(\dot\B^{-s,0}_{\infty,q'}(\RR^n))' = \dot\B^{s}_{1,q}(\RR^n)$, $s<n$, $q\in(1,\infty]$;
    \item $(\dot\BesSmo^{-s,0}_{\infty,\infty}(\RR^n))' = \dot\B^{s}_{1,1}(\RR^n)$, $s\leqslant n$.
\end{itemize}
Note that the conditions on each line in the duality identities for homogeneous function spaces correspond to \eqref{AssumptionCompletenessExponents}. However, if we go beyond the completeness threshold, the duality bracket is still relevant and sufficient to compute what should be the norm of a "dual element" that belongs to $\S'_h(\RR^n)$, see \textit{e.g.} \cite[Proposition~2.29]{bookBahouriCheminDanchin}.

\medbreak

\subsubsection{On Lipschitz domains.} For any $\X^{s,p}\in\{ \B^{s}_{p,q}, \dot{\B}^{s}_{p,q},  {\H}^{s,p}, \dot{\H}^{s,p},  {\W}^{s,p}, \dot{\W}^{s,p}\}$, we define
\begin{align*}
    \X^{s,p}(\Omega):= \X^{s,p}(\RR^n)_{|_{\Omega}}\text{, with the quotient norm }\lVert u \rVert_{\X^{s,p}(\Omega)}:= \inf\limits_{\substack{\Tilde{u}\in \X^{s,p}(\RR^n),\\ \tilde{u}_{|_{\Omega}}=u\, .}} \lVert \Tilde{u} \rVert_{\X^{s,p}(\RR^n)}.
\end{align*}
Since there is a natural restriction homomorphism, we do have $\S'(\RR^n)\hookrightarrow \Distrib(\Omega)$, we also have the inclusion,
\begin{align*}
    \X(\Omega) \subset \Distrib(\Omega)\text{.}
\end{align*}
We denote $\S'(\Omega)$ (resp. $\S'_h(\Omega)$), the image of $\S'(\RR^n)$ (resp. $\S'_h(\RR^n)$) in $\Distrib(\Omega)$ by the restriction homomorphism.

A direct consequence  is the density of $\S_0(\overline{\Omega})\subset\S(\overline{\Omega})$ in each of them, and the completeness and reflexivity when their counterpart on $\RR^n$ are complete and reflexive. We also define
\begin{align*}
    \X^{s,p}_0(\Omega):= \left\{\,u\in \X^{s,p}(\RR^n) \,\,:\,\, \supp u \subset \overline{\Omega} \right\}\text{, }
\end{align*}
with its natural induced norm $\lVert  u \rVert_{\X_0(\Omega)}:= \lVert  u \rVert_{\X^{s,p}(\RR^n)}$. We always have the canonical continuous injection,
\begin{align*}
    \X^{s,p}_0(\Omega)\hookrightarrow \X^{s,p}(\Omega) \text{. }
\end{align*}
We recall that in most cases, $p,q\in[1,\infty)$, one generally has $\X^{s,p}_0(\Omega)=\overline{\Ccinfty(\Omega)\cap \X^{s,p}(\RR^n)}^{\lVert\cdot\rVert_{\X^{s,p}(\RR^n)}}$, and still holds weakly-$\ast$ when $q=\infty$ or $p=\infty$. One still has strong density for the spaces $\BesSmo^{s}_{p,\infty,0}(\Omega)$, $\B^{s,0}_{\infty,q,0}(\Omega)$, and $\BesSmo^{s,0}_{\infty,\infty,0}(\Omega)$, $p,q<\infty$, $s\in\RR$, and it also holds for their homogeneous counterparts.

An explicit and nearly exhaustive treatment of  inhomogeneous function space on  bounded Lipschitz domains (and for which, still a wide part of it remains valid for exterior and special Lipschitz domains) has been achieved by Kalton, Mitrea and Mitrea \cite{KaltonMayborodaMitrea2007}. A treatment of the full Triebel-Lizorkin scale, including Hardy spaces, and other endpoint function spaces is also considered.

It includes, in particular the natural interpolation identities: Provided $(p_0,q_0),(p_1,q_1),(p,q)\in[1,\infty]^2$, $s_0\neq s_1$ and two real numbers, and $\theta\in(0,1)$ are such that
\begin{align*}
    \left(s,\frac{1}{p_\theta},\frac{1}{q_\theta}\right):= (1-\theta)\left(s_0,\frac{1}{p_0},\frac{1}{q_0}\right)+ \theta\left(s_1,\frac{1}{p_1},\frac{1}{q_1}\right)\text{,}
\end{align*}
one has
\begin{itemize}
    \item $[\H^{s_0,p_0}(\Omega),\H^{s_1,p_1}(\Omega)]_\theta=\H^{s,p_\theta}(\Omega)$, $1<p_0,p_1<\infty$;
    
    \item $(\H^{s_0,p}(\Omega),\H^{s_1,p}(\Omega))_{\theta,q}=(\B^{s_0}_{p,q_0}(\Omega),\B^{s_1}_{p,q_1}(\Omega))_{\theta,q} = \B^{s}_{p,q}(\Omega)$;

    \item $[\B^{s_0}_{p_0,q_0}(\Omega),\B^{s_1}_{p_1,q_1}(\Omega)]_{\theta} = \B^{s}_{p_\theta,q_\theta}(\Omega)$, $q_\theta<\infty$.
\end{itemize}
One can also consider $\B^{\cdot,0}_{\infty,q}$ instead of $\B^{\cdot}_{\infty,q}$, and  $(\B^{s}_{p,\infty},\,(\,,\,)_{\theta,\infty})$ by $(\BesSmo^{s}_{p,\infty},\,(\,,\,)_{\theta})$. Obviously similar results of interpolation can be considered with
\begin{center}
    $\H_0$ and $\B_{\cdot,\cdot,0}$ instead of $\H$, $\B$.
\end{center}

One also has the standard duality relations
\begin{itemize}
    \item $(\H^{-s,p'}_0(\Omega))' = \H^{s,p}(\Omega)$, and $(\H^{-s,p'}(\Omega))' = \H^{s,p}_0(\Omega)$, $s\in\RR$, $1<p<\infty$;
    \item $(\B^{-s}_{p',q',0}(\Omega))' = \B^{s}_{p,q}(\Omega)$, and $(\B^{-s}_{p',q',0}(\Omega))' = \B^{s}_{p,q,0}(\Omega)$, $s\in\RR$, $p,q\in(1,\infty]$;
    \item $(\BesSmo^{-s}_{p',\infty,0}(\Omega))' = \B^{s}_{p,1}(\Omega)$, and $(\BesSmo^{-s}_{p',\infty}(\Omega))' = \B^{s}_{p,1,0}(\Omega)$, $s\in\RR$, $p\in(1,\infty]$;
    \item $(\B^{-s,0}_{\infty,q',0}(\Omega))' = \B^{s}_{1,q}(\Omega)$, and $(\B^{-s,0}_{\infty,q'}(\Omega))' = \B^{s}_{1,q,0}(\Omega)$, $s\in\RR$, $q\in(1,\infty]$;
    \item $(\BesSmo^{-s,0}_{\infty,\infty,0}(\Omega))' = \B^{s}_{1,1}(\Omega)$, and $(\BesSmo^{-s,0}_{\infty,\infty}(\Omega))' = \B^{s}_{1,1,0}(\Omega)$, $s\in\RR$.
\end{itemize}
Note that these duality relations can be strengthened in a surprising way whenever $\Omega$ is bounded, see Theorem~\ref{thm:Reflexiveendpointspaces} below.

For interpolation as well as duality results, the corresponding ones are available for homogeneous Sobolev and Besov spaces on special Lipschitz domains in \cite[Section~3]{Gaudin2023}, \textit{i.e.}
\begin{center}
     with $\dot{\H}$ and $\dot{\B}$ that replace $\H$ and $\B$.
\end{center}

As in the case of the whole space $\Omega=\RR^n$,  concerning \textbf{homogeneous} function spaces, the completeness of any involved function spaces is required to define properly \textbf{complex interpolation}. 

\medbreak 

Concerning \textbf{duality},  similarly to what happens on $\Omega=\RR^n$, one requires the completeness of the output space (so just assuming \eqref{AssumptionCompletenessExponents} for space that can be computed as the dual of another one). Beyond the completeness threshold, both spaces are still in duality but, since the other one is not complete, it can obviously not realize all the bounded continuous linear functionals. 

\medbreak

\textbf{Traces.} When $\Omega=\RR^n_+$, and $s>\sfrac{1}{p}$, one can usually define the restriction on the boundary $\RR^{n-1}\times\{0\}$ of an element $u\in\X^{s,p}(\RR^n_+)$, $u(\cdot,0)=:u_{|_{\partial\RR^n_+}}$ as an element of $\X_\partial^{s-\frac{1}{p},p}(\partial\RR^{n}_+)=\X_\partial^{s-\sfrac{1}{p},p}(\RR^{n-1})$, called the trace of $u$ where
\begin{align*}
    \X_\partial^{s-\sfrac{1}{p},p} =\begin{cases}
        \B^{s-\sfrac{1}{p}}_{p,p}\text{,  if }\X=\H\text{;}\\
        \B^{s-\sfrac{1}{p}}_{p,q}\text{,  if }\X=\B_{\cdot,q}\text{.}\\
    \end{cases}
\end{align*}

When $\Omega$ is a Lipschitz domain, one can define the corresponding function spaces on the boundary $\partial\Omega$ according to the surface measure $\d \sigma_{\Omega}$ by localisation and then transformation. For $p,q\in[1,\infty]$, $0<s<1$, one can define $\L^p(\partial\Omega)$, $\W^{1,p}(\partial\Omega)$, by localisation and change of coordinates, then one can construct corresponding Besov spaces by interpolation:
\begin{align*}
    \B^{s}_{p,q}(\partial\Omega):=(\L^p(\partial\Omega),\W^{1,p}(\partial\Omega))_{s,q}.
\end{align*}

In this configuration, it can be proved that one has a well-defined trace operator
\begin{align}\label{eq:TraceOpLipschitz}
   {[\cdot]}_{|_{\partial\Omega}}\,:\,\X^{s,p}(\Omega)\longrightarrow\X_\partial^{s-\sfrac{1}{p},p}(\partial\Omega)\text{, } \sfrac{1}{p}<s<1+\sfrac{1}{p}\text{, }p,q\in[1,\infty],
\end{align}
and similarly for homogeneous function spaces on a special Lipschitz domain, see in this case \cite[Section~4]{Gaudin2023Lip}. An improvement in the case of homogeneous function spaces on the half-space $\RR^n_+$ is given in Appendix~\ref{App:Traces}. 

For Lipschitz boundaries, while the charts may lack of regularity, one can still define inhomogeneous function spaces of higher order on the boundary, but the naive definition mentioned above is now longer suitable. Following \cite[Chapters~V,~VI~\&~VII]{JonssonWallin1984}, one can still make sense of \eqref{eq:TraceOpLipschitz} for any $s\in\RR$, $p,q\in[1,\infty]$ such that
\begin{align*}
    k+\sfrac{1}{p}<s<k+1+\sfrac{1}{p}\text{, }k\in\NN,
\end{align*}
up to some appropriate changes of the definition. In this setting, the trace of higher regularity Sobolev or Besov functions still has improved regularity in some sense: the regularity is not confined by~$s-\frac{1}{p}<1$. However, we mention here that in general, for Lipschitz domains and other irregular domains, whenever $s>1+\sfrac{1}{p}$ one can no longer properly identify by local transformations
\begin{align}\label{eq:IntroIdentBesovBoundaryBesovRn-1}
    \B^{s-\sfrac{1}{p}}_{p,q}(\partial\Omega)\quad\quad\text{ with }\quad\quad\B^{s-\sfrac{1}{p}}_{p,q}(\RR^{n-1}).
\end{align}
See \cite[End of Section~2]{MitreaMitreaMonniaux2008} for a more recent exposition that fits better our current presentation. However, more regular domains allow to improve this identification: for instance, $\C^{1,\alpha}$-domains, with $\alpha\in(0,1]$, allow the identification for $s<1+\alpha+\sfrac{1}{p}$, etc..\footnote{This can also be seen as one of the hidden purposes of the Sobolev multiplier theory developped by Maz'ya and Shaposhnikova: for which abstract rough boundaries, depending on $p$ and $s$, such a local identification \eqref{eq:IntroIdentBesovBoundaryBesovRn-1} is possible.} 

Inhomogeneous Besov spaces of negative regularity $-1<s<0$ on the boundary to be the elements $u\in(\C^{0,1}_c(\partial\Omega))'$, such that
\begin{align*}
    \lVert u \rVert_{\B^{s}_{p,q}(\partial\Omega)}:=\sup_{\substack{\varphi \in\C_c^{0,1}(\partial\Omega),\\\lVert \varphi\rVert_{\B^{-s}_{p',q'}(\partial\Omega)}\leqslant1}} |\langle u,\varphi\rangle_{\partial\Omega}|<\infty
\end{align*}
where $\langle \cdot, \cdot\rangle_{\partial\Omega}$ extends the $\L^2(\partial\Omega)$-pairing (with respect to the surface measure). Besov spaces of negative regularity are necessary in order to define partial trace of vector fields of low regularity. See the the second part of Appendix~\ref{App:Traces}.

For $p,q\in[1,\infty]$ and $s>-1+{\sfrac{1}{p}}$, one can set
\begin{align*}
    &{\X}^{s,p}_{\mathcal{D}}(\Omega):= \begin{cases} {\X}^{s,p}(\Omega) &\text{, if } s < {\sfrac{1}{p}}\text{,}\\
    \left\{ \,u\in{\X}^{s,p}(\Omega)\, \,:\,\, u_{|_{\partial\Omega}}=0 \,\right\} &\text{, if }s > {\sfrac{1}{p}}\text{,}
    \end{cases}\\
    &{\X}^{s,p}_{\mathcal{N}}(\Omega):= \begin{cases} {\X}^{s,p}(\Omega) &\text{, if } s < 1+ {\sfrac{1}{p}}\text{,}\\
    \left\{ \,u\in {\X}^{s,p}(\Omega)\, \,:\,\, \partial_{\nu}u_{|_{\partial\Omega}}=0\,\right\} &\text{, if }s >1+{\sfrac{1}{p}}\text{,}
    \end{cases}
\end{align*}
As well as their homogeneous counterparts when $\Omega$ is a special Lipschitz domain.

In the case of a homogeneous Dirichlet boundary condition it is well known extending by $0$ can still provide a suitable distribution over the whole space $\RR^n$ preserving the regularity. This is a very well known result, however, see \cite[Proposition~4.21]{Gaudin2023Lip} for the case of homogeneous function spaces.

\begin{proposition}\label{prop:Ext0DirLip}Let $p\in(1,\infty)$, $s\in(\sfrac{1}{p},1+\sfrac{1}{p})$. Let $\Omega$ be a bounded or special Lipschitz domain. Then extension map by $0$ to the whole space yields a canonical isomorphism
\begin{align*}
    \H^{s,p}_{\mathcal{D}}(\Omega) \simeq {\H}^{s,p}_0(\Omega) \text{.}
\end{align*}
The result still holds replacing ${\H}^{s,p}$ by ${\B}^{s}_{p,q}$, ${\BesSmo}^{s}_{p,\infty}$, ${\B}^{s,0}_{\infty,q}$, ${\W}^{1,1}$, ${\C}^{0}_0$, $\C^{0}_{ub}$, $p,q\in[1,\infty]$.

\medbreak

Additionally, the whole result remains valid for homogeneous function spaces (assuming $1<p<\infty$ in the case of Bessel potential spaces):
\begin{align*}
    \dot{\H}^{s,p},\, \dot{\B}^{s}_{p,q},\,\dot{\BesSmo}^{s}_{p,\infty},\,\dot{\B}^{s,0}_{\infty,q},\,\C^{0}_{ub,h},\text{ and }\,\dot{\W}^{1,1}\,,\,p,q\in[1,\infty].
\end{align*}

As a consequence, when $p,q<\infty$ ($\C^{0}_{ub}$, $\C^{0}_{ub,h}$ being appart), $\Ccinfty(\Omega)$ is a strongly dense subspace of the nullspace of the trace operator, allowing $p=\infty$ and $\C^{0}_{ub}$ whenever $\Omega$ is bounded.
\end{proposition}

\newpage
\textbf{Some fundamental results.} The very well known Proposition~\ref{prop:Ext0DirLip} is mainly due to the \textbf{fundamental Sobolev multiplier theorem} concerning stability of function spaces by the multiplication by $\mathbbm{1}_\Omega$, initially due to Strichartz for $\Omega=\RR^n_+$ in the case inhomogeneous (Bessel potential) Sobolev spaces of fractional order over $\RR^n$ \cite[Chapter~II,~Corollary~3.7]{Strichartz1967}. It has been extended several times to rougher boundaries and further (inhomogeneous) functions spaces, see for instance the work by Runst and Sickel \cite[Chapter~4,~Section~4.6.3]{RunstSickel96}  and \cite[Proposition~4.8]{Sickel1999}. For homogeneous function spaces, one can check \cite[Appendix,~Lemma~12]{DanchinMucha2009} for $\Omega=\RR^n_+$ and homogeneous Besov spaces, \cite[Proposition~3.5]{Gaudin2023Lip}, for $\Omega$ a special Lipschitz domain in the case of homogeneous Sobolev and Besov spaces. We restate it here and provide further of its consequences. A vector-valued version for Besov spaces is available in Appendix, see Theorem~\ref{thm:fundamentalmultBesovSpacesX}.

\medbreak

\textbf{A very important consequence from the result below} is that when $\Omega$ is a bounded Lipschitz domain, it will make then no differences to consider homogeneous or inhomogeneous Sobolev and Besov spaces whenever $-1+\sfrac{1}{p}<s<\sfrac{1}{p}$.

\begin{proposition}\label{prop:FundamentalExtby0HomFuncSpaces}Let $p,q\in[1,\infty]$ and $s\in(-1+{\sfrac{1}{p}},{\sfrac{1}{p}})$. Let $\Omega$ be either a special, bounded or exterior Lipschitz domain. For $(\mathfrak{H},\mathfrak{B})\in\{\,(\H,\B),\,(\dot{\H},\dot{\B})\,\}$, one does have:
\begin{enumerate}
    \item Provided additionally that $1<p<\infty$, for all $u\in\mathfrak{H}^{s,p}(\RR^n)$, 
    \begin{align*}
        \lVert \car_{\Omega} u\rVert_{\mathfrak{H}^{s,p}(\RR^n)} \less_{p,s,n}^{\partial\Omega} \lVert u\rVert_{\mathfrak{H}^{s,p}(\RR^n)}. 
    \end{align*}
    \item For all $u\in\mathfrak{B}^{s}_{p,q}(\RR^n)$,
    \begin{align*}
        \lVert \car_{\Omega} u\rVert_{\mathfrak{B}^{s}_{p,q}(\RR^n)} \less_{p,s,n}^{\partial\Omega} \lVert u\rVert_{\mathfrak{B}^{s}_{p,q}(\RR^n)}. 
    \end{align*}
    Moreover, $\car_{\Omega}\BesSmo_{p,\infty}^{s}(\RR^n)\subset\BesSmo_{p,\infty}^{s}(\RR^n)$, and $\car_{\Omega}\B_{\infty,q}^{s,0}(\RR^n)\subset\B_{\infty,q}^{s,0}(\RR^n)$, and similarly for the corresponding homogeneous Besov spaces.
    \item Consequently, one has the following canonical identification with equivalence of norms
    \begin{align*}
        \mathfrak{H}^{s,p}(\Omega) &= \mathfrak{H}^{s,p}_0(\Omega)\text{, }1<p<\infty\text{,}\\
        \mathfrak{B}^{s}_{p,q}(\Omega)&=\mathfrak{B}^{s}_{p,q,0}(\Omega)\text{.}
    \end{align*}
    \item If moreover $\Omega$ is bounded, we also have
    \begin{align*}
        \H^{s,p}(\Omega) &= \dot{\H}^{s,p}(\Omega)\text{, }1<p<\infty\text{,}\\
        \B^{s}_{p,q}(\Omega)&=\dot{\B}^{s}_{p,q}(\Omega)\text{,}
    \end{align*}
    with an equivalence of norms whose constants depend on the finite measure of $\Omega$.
\end{enumerate}
In each case, the constant depends only on the Lipschitz constant of $\Omega$.
\end{proposition}

\begin{proof}Points \textit{(i)} and \textit{(ii)} for homogeneous function spaces over special Lipschitz domains $\Omega$ are given by \cite[Proposition~3.5]{Gaudin2023Lip}. The case of special Lipschitz domains for inhomogeneous function spaces can be deduced by similar arguments. For $\X^{s,p}\in\{ \H^{s,p},\,\dot{\H}^{s,p},\,\B^{s}_{p,q},\,\dot{\B}^{s}_{p,q}\}$, $\Omega=\{ (x',x_n)\in\RR^n\,:\,x_n>\phi(x')\,\}$ to be a special Lipschitz domain, one obtains the operator-norm bound
\begin{align*}
    \lVert \car_{\Omega} \cdot[\cdot]\rVert_{\X^{s,p}(\RR^n)\rightarrow\X^{s,p}(\RR^n)}&\lesssim_{p,s,n} (1+\lVert\nabla'\phi\rVert_{\L^\infty(\RR^{n-1})})^{|s|}, \, \text{ if } s\neq 0, \\
    \text{ and }\quad\lVert \car_{\Omega} \cdot[\cdot]\rVert_{\X^{0,p}(\RR^n)\rightarrow\X^{0,p}(\RR^n)}&\lesssim_{p,\varepsilon,n} (1+\lVert\nabla'\phi\rVert_{\L^\infty(\RR^{n-1})})^{\varepsilon},
\end{align*}
for some arbitrary $\varepsilon\in(0,1)$. If Points \textit{(i)} and \textit{(ii)} are satisfied then Point \textit{(iii)} is a straightforward consequence. Thus, we focus on proving Points \textit{(i)} and \textit{(ii)} for bounded and exterior Lipschitz domains and Point \textit{(iv)}.

\textbf{Step 1:} The points \textit{(i)} and \textit{(ii)} for  bounded and exterior Lipschitz domains. We do it for $\X^{s,p}\in\{ \H^{s,p},\,\dot{\H}^{s,p},\,\B^{s}_{p,q},\,\dot{\B}^{s}_{p,q}\}$.

If $\Omega$ is a bounded Lipschitz domain, then one can cover $\partial\Omega$ by a finite number of balls $(\B_{j})_{j\in\llb1,N\rrb}$, such that up to rotation and translation, there exists a global Lipschitz function $\phi_j\,:\,\RR^{n-1}\longrightarrow\RR$, such that
\begin{align*}
    \Omega\cap\B_j = \{ (x',x_n)\in\RR^n\,:\,x_n>\phi_j(x')\,\}\cap \B_j.
\end{align*}
For $j\in\llb1,n \rrb$, we define $\widetilde{\Omega}_j := \{ (x',x_n)\in\RR^n\,:\,x_n>\phi_j(x')\,\}$, and we set $\Omega_0:=\Omega\setminus\Big[\bigcup_{j=1}^N\widetilde{\Omega}_j\cap\B_j\Big]$, and we consider $(\varphi_j)_{j\in\llb0,N\rrb}\in\Ccinfty(\RR^n)$ a parition of unity subordinated to $(\Omega_0, (\B_j)_{j\in\llb1,N\rrb})$.

One can write
\begin{align*}
    \car_{\Omega} = \varphi_0 + \sum_{j=1}^N \varphi_j\car_{\widetilde{\Omega}_j}.
\end{align*}
Thus, for $s\neq0$, by Propositions~\ref{prop:EmbeddingsHomogeneousSobolevBesovinMultiplierSpaces}~\&~\ref{prop:EmbeddingsSobolevBesovinMultiplierSpaces} combined with Proposition~\ref{prop:basicMultipliers} and the Sobolev embedding ${\dot{\B}^{|s|}_{\frac{n}{|s|},1}(\RR^n)}\hookrightarrow\L^\infty(\RR^n)$, for $s\neq 0$, $u\in{\X}^{s,p}(\RR^n)$,
\begin{align*}
    \lVert \car_\Omega u\rVert_{{\X}^{s,p}(\RR^n)} &\leqslant \lVert \varphi_0 u\rVert_{{\X}^{s,p}(\RR^n)} + \sum_{j=1}^N \lVert \varphi_j\car_{\widetilde{\Omega}_j} u\rVert_{{\X}^{s,p}(\RR^n)}\\
    &\lesssim_{p,s,n} \Big(\lVert \varphi_0 \rVert_{\dot{\B}^{|s|}_{\frac{n}{|s|},1}(\RR^n)} + \sum_{j=1}^N \lVert \varphi_j \rVert_{\dot{\B}^{|s|}_{\frac{n}{|s|},1}(\RR^n)}\lVert \car_{\widetilde{\Omega}_j} \cdot[\cdot]\rVert_{\X^{s,p}(\RR^n)\rightarrow\X^{s,p}(\RR^n)}\Big) \lVert u\rVert_{{\X}^{s,p}(\RR^n)}.
\end{align*}
The case $s=0$, available when $1<p<\infty$, follows by interpolation.  

Note that due to the scale invariance of the $\dot{\B}^{|s|}_{\frac{n}{|s|},1}(\RR^n)$-norms, this means that the operator norm $\lVert \car_{\Omega} \cdot[\cdot]\rVert_{\X^{s,p}(\RR^n)\rightarrow\X^{s,p}(\RR^n)}$ has a bound that only depends on $p$, $s$, $n$, the Lipschitz constant of $\partial\Omega$ (of each $\partial\widetilde{\Omega}_j$), and the minimal number of charts to describe the boundary. It does not depend on the measure of $\Omega$.

If one writes $\car_{\Omega^{c}} = \mathbf{1} - \car_{\Omega}$, then we obtain the result for exterior Lipschitz domains.

\textbf{Step 2:} We prove Point \textit{(iv)}. The equality inhomogeneous and homogeneous function spaces within the range $s\in(-1+\sfrac{1}{p},\sfrac{1}{p})$ for bounded Lipschitz domains.

Let $0<s<1$ and $u\in\dot{\B}^{s}_{1,1}(\Omega)$, by Point \textit{(iii)}, one has $\mathcal{E}_0u \in \dot{\B}^{s}_{1,1}(\RR^n)$, where $\mathcal{E}_0u$ is the extension of $u$ to the whole $\RR^n$ by $0$, with
\begin{align*}
    \lVert\mathcal{E}_0u \rVert_{\dot{\B}^{s}_{1,1}(\RR^n)} \leqslant \lVert \car_{\Omega} \cdot[\cdot]\rVert_{\dot{\B}^{s}_{1,1}(\RR^n)\rightarrow\dot{\B}^{s}_{1,1}(\RR^n)} \lVert u \rVert_{\dot{\B}^{s}_{1,1}(\Omega)}.
\end{align*}
By the Sobolev embedding, and compact support of $\mathcal{E}_0u$ in $\overline{\Omega}$, one obtains
\begin{align*}
    \lVert \mathcal{E}_0 u\rVert_{\L^1(\RR^n)} \leqslant |\Omega|^{\frac{s}{n}} \lVert \mathcal{E}_0 u\rVert_{\L^\frac{n}{n-s}(\RR^n)} \lesssim_{s,n} |\Omega|^{\frac{s}{n}} \lVert \mathcal{E}_0 u\rVert_{\dot{\B}^{s}_{1,1}(\RR^n)}\lesssim_{s,n,\partial\Omega} |\Omega|^{\frac{s}{n}} \lVert  u\rVert_{\dot{\B}^{s}_{1,1}(\Omega)}.
\end{align*}
Thus, by the definition of function spaces by restriction, we have obtained
\begin{align*}
    \lVert u \rVert_{{\B}^{s}_{1,1}(\Omega)} \leqslant  \lVert \mathcal{E}_0 u\rVert_{{\B}^{s}_{1,1}(\RR^n)} \lesssim_{s,n} \lVert \mathcal{E}_0 u\rVert_{\L^1(\RR^n)} + \lVert \mathcal{E}_0 u\rVert_{\dot{\B}^{s}_{1,1}(\RR^n)} \lesssim_{s,n,\partial\Omega} (1+|\Omega|^{\frac{s}{n}}) \lVert  u\rVert_{\dot{\B}^{s}_{1,1}(\Omega)}.
\end{align*}
Thus, for all $0<s<1$, one has 
\begin{align*}
    {\B}^{s}_{1,1}(\Omega)=\dot{\B}^{s}_{1,1}(\Omega)
\end{align*}
with equivalence of norms. By real interpolation, (apply interpolation to the identity map) combined with Point \textit{(iii)}, one obtains 
\begin{align*}
    {\B}^{s}_{1,q,0}(\Omega)={\B}^{s}_{1,q}(\Omega)=\dot{\B}^{s}_{1,q}(\Omega)=\dot{\B}^{s}_{1,q,0}(\Omega)
\end{align*}
for all $q\in[1,\infty]$ and all $0<s<1$, with equivalence of norms.

Actually, the same arguments apply to ${\H}^{s,p}$, $\B^{s}_{p,q}$, whenever $p\in(1,\infty)$ and $0<s<\sfrac{1}{p}$. Duality yields the case $p\in(1,\infty]$, $-1+\sfrac{1}{p}<s<0$.  By interpolation again, the equality remains valid for $s=0$ holds.
\end{proof}

\begin{corollary}\label{cor:Reflexiveendpointspaces}  Let $\Omega$ be a bounded Lipschitz domain. For $0<s<1$, $q\in[1,\infty]$, one has
\begin{align*}
    \B^{-s}_{\infty,q}(\Omega) = \B^{-s,0}_{\infty,q}(\Omega)\text{ and }\BesSmo^{-s}_{\infty,\infty}(\Omega) = \BesSmo^{-s,0}_{\infty,\infty}(\Omega).
\end{align*}
In particular, $\B^{-s}_{\infty,q}(\Omega)$ and $\B^{s}_{1,q}(\Omega)$ are reflexive whenever $1<q<\infty$.
\end{corollary}

\begin{proof} This is a direct consequence of the previous proposition. Indeed, if $q<\infty$ let us consider the extension of $u\in \B^{-s}_{\infty,q}(\Omega)=\dot{\B}^{-s}_{\infty,q}(\Omega)$ by $0$ to the whole space $\tilde{u}$, and then a sequence $(\tilde{u}_{\ell})_{\ell\in\NN}\subset\C^\infty_{ub,h}(\RR^n)$ converging towards $\tilde{u}$ in $\dot{\B}^{-s}_{\infty,q}(\RR^n)$. By the previous proposition, again, $(\mathbbm{1}_{\Omega}\tilde{u}_{\ell})_{\ell\in\NN}$ converges towards $\tilde{u}$ in $\dot{\B}^{-s}_{\infty,q}(\RR^n)$. By construction, and since $\Omega$ is bounded, one has $(\mathbbm{1}_{\Omega}\tilde{u}_{\ell})_{\ell\in\NN}\subset\dot{\B}^{-s,0}_{\infty,q}(\RR^n)$, since $\dot{\Delta}_j\mathbbm{1}_{\Omega}\tilde{u}_{\ell}\in\C^\infty_0(\RR^n)$.

Since $\dot{\B}^{-s,0}_{\infty,q}(\RR^n)$ is closed in $\dot{\B}^{-s}_{\infty,q}(\RR^n)$, one obtains $\tilde{u}\in\dot{\B}^{-s,0}_{\infty,q}(\RR^n)$, here, this is because by passing to the limit one obtains $\dot{\Delta}_j\mathbbm{1}_{\Omega}\tilde{u}\in\C^0_0(\RR^n)$, which is the definition of belonging to $\dot{\B}^{-s,0}_{\infty,q}(\RR^n)$. By restriction, it holds that $\tilde{u}_{|_{\Omega}} =u\in \dot{\B}^{-s,0}_{\infty,q}(\Omega)={\B}^{-s,0}_{\infty,q}(\Omega)$.

By real interpolation the remaining identities for $\B^{-s}_{\infty,\infty}$ and $\BesSmo^{-s}_{\infty,\infty}$ also hold. This finishes the proof. 
\end{proof}

We now start to mention a first important theorem of this paper, which is is a direct consequence of previous Corollary~\ref{cor:Reflexiveendpointspaces}. This result will be of paramount importance to reach optimal results in endpoint function spaces.

\begin{theorem}\label{thm:Reflexiveendpointspaces} Let $\Omega$ be a bounded Lipschitz domain. For all  $s\in\RR$, $p\in[1,\infty]$, $q\in(1,\infty)$ the spaces $\B^{s}_{p,q}(\Omega)$ and $\B^{s}_{p,q,0}(\Omega)$ are reflexive. Each of them admits respectively $\C^\infty(\overline{\Omega})$ and $\Ccinfty(\Omega)$ as strongly dense subspaces whenever $q\in[1,\infty)$. The strong density result remains valid for the spaces $\BesSmo^{s}_{p,\infty}(\Omega)$ and $\BesSmo^{s}_{p,\infty,0}(\Omega)$.
\end{theorem}

\begin{remark}The knowledge of endpoint Besov spaces $\B^{s}_{1,q}(\Omega)$, $\B^{s}_{\infty,q}(\Omega)$, $1<q<\infty$ being reflexive may seem quite surprising at the first glance, even false. However,  in the case of the torus $\mathbb{T}^n$, see \cite[Definition~1.32~\&~Theorem~1.37] {bookTriebel2008}, \cite[Introduction,~p.283]{DilworthFreemanOdellSchlumprecht2011} and the references therein, it is known that we have an isomorphism
\begin{align*}
    \iota_s\,:\, \B^{s}_{p,q}(\mathbb{T}^n) \longrightarrow \ell^q(\oplus_{n=1}^{\infty}\ell^{p}_n(\llb 1,n\rrb))
\end{align*}
where $\ell^q(\oplus_{n=1}^{\infty}\ell^{p}_n(\llb 1,n\rrb))$ has explicit norm
\begin{align*}
    \big\lVert (\alpha_{j,n})_{1\leqslant j\leqslant n\leqslant \infty}\big\rVert_{\ell^q(\oplus_{n=1}^{\infty}\ell^{p}_n(\llb 1,n\rrb))} = \left(\sum_{n\geqslant 1}\Big( \sum_{j=1}^{n} |2^{jn}\alpha_{j,n}|^p\Big)^{\sfrac{q}{p}}\right)^{\sfrac{1}{q}}
\end{align*}
The latter space is known to be reflexive\footnote{Understand it as the direct $\ell^q$-sum  of many countable weight copies of $\RR^n$ equipped with the finite-dimensional $\ell^p$-norm.} for all $p\in[1,\infty]$, $q\in(1,\infty)$, while not being super-reflexive when $p=1,\infty$, and in particular not UMD (every UMD Banach is necessarily super-reflexive). This provides endpoint Besov spaces on compact domains/manifolds as simple and explicit examples of reflexive but non-UMD Banach spaces. Super-reflexive function spaces that are non-UMD were constructed but lead to really involved counter-examples, see for instance the one due to Pisier \cite{Pisier75SuperRefnonUMD}. This example of endpoint Besov spaces shows that such spaces are not artificial and  can emerge quite naturally in the world of Partial Differential Equations.
\end{remark}

\subsubsection{Additionnal knowledge: Compact support, Poincare inequalities, and some extension properties on the half-space}

The next Lemma is fundamental since it allows to play between homogeneous and inhomogeneous function spaces in several cases, when one deals with localisation procedures.

\begin{lemma}\label{lem:CompactEquiHomInhom}Let $s\in\RR$, $p,q\in[1,\infty]$, and let $u\in\S'_h(\RR^n)$ be compactly supported. If $s>-n/p'$, then it holds that
\begin{align*}
    u\in\H^{s,p}(\RR^n) &\iff u\in\dot{\H}^{s,p}(\RR^n)\text{, }1<p<\infty\text{;}\\
     u\in\B^{s}_{p,q}(\RR^n) &\iff u\in\dot{\B}^{s}_{p,q}(\RR^n)\text{, }1\leqslant p,q\leqslant\infty\text{;}
\end{align*}
with equivalence of norms, for $\X^{s,p}\in\{{\H}^{s,p},\,\B^{s}_{p,q}\}$,
\begin{align*}
    \lVert u \rVert_{\dot{\X}^{s,p}(\RR^n)} &\lesssim_{p,s,n} \lVert u \rVert_{{\X}^{s,p}(\RR^n)} \lesssim_{p,s,n} (1+|K|^\frac{s}{n})\lVert u \rVert_{\dot{\X}^{s,p}(\RR^n)},\qquad \text{ if }s>0;\\
    \lVert u \rVert_{{\X}^{s,p}(\RR^n)} &\lesssim_{p,s,n} \lVert u \rVert_{\dot{\X}^{s,p}(\RR^n)} \lesssim_{p,s,n} (1+|K|^\frac{-s}{n})\lVert u \rVert_{{\X}^{s,p}(\RR^n)},\qquad \text{ if }s<0;\\
     \lVert u \rVert_{\dot{\B}^{0}_{p,q}(\RR^n)} \lesssim_{p,n} (1+|K|^\frac{1}{n})&\lVert u \rVert_{{\B}^{0}_{p,q}(\RR^n)},\quad \text{ and }\quad \lVert u \rVert_{{\B}^{0}_{p,q}(\RR^n)} \lesssim_{p,n} (1+|K|^\frac{2}{n})\lVert u \rVert_{\dot{\B}^{0}_{p,q}(\RR^n)}.
\end{align*}
Furthermore, in the case of Besov spaces, the result still holds if $s=-\sfrac{n}{p'}$ and $q=\infty$ whenever $p\in(1,\infty]$.
\end{lemma}

\begin{remark}The result is sharp. Because, as mentioned in the introduction, $\Ccinfty(\RR^n)\subset{\H}^{-n/2,2}(\RR^n)$, while $\Ccinfty(\RR^n)\not\subset\dot{\H}^{-n/2,2}(\RR^n)$. Similarly, $\Ccinfty(\RR^n)\subset \B^0_{1,1}(\RR^n)$, while $\Ccinfty(\RR^n)\not\subset \dot{\B}^0_{1,1}(\RR^n)$. Things got similar for other function spaces $\dot{\B}^{s}_{p,q}$, $\dot{\H}^{s,p}$, $s<-\frac{n}{p'}$.
\end{remark}

\begin{proof}Let $u\in\S'_h(\RR^n)$ with $\supp u=:K$ compact. The equivalence between the membership of inhomogeneous and homogeneous function spaces follows directly from \cite[Lemma~2.7]{Gaudin2023Lip} and compactness of $K$. It remains to check the equivalence of norms. Thanks to translation invariance of Lebesgue norms, without loss of generality, we can assume $0\in K$. 

\textbf{Step 1:} Now, we show equivalence of norms for Sobolev spaces. Let $1<p<\infty$, assume first $s > 0$. One always has $\H^{s,p}(\RR^{n})\hookrightarrow\dot{\H}^{s,p}(\RR^{n})$. If $0<s<\frac{n}{p}$, by H\"{o}lder's inequality and Sobolev embeddings, one obtains
\begin{align}\label{eq:ProofEquivCompactSupport1}
    \lVert u\rVert_{\L^p(\RR^n)} \leqslant |K|^{\frac{s}{n}}\lVert u\rVert_{\L^{\frac{np}{n-ps}}(\RR^n)}\lesssim_{p,s,n} |K|^{\frac{s}{n}}\lVert u\rVert_{\dot{\H}^{s,p}(\RR^n)}.
\end{align}
Now, if $s\geqslant\frac{n}{p}$, choose $0<s_0<\frac{n}{p}$. By interpolation inequalities, it holds 
\begin{align}\label{eq:ProofEquivCompactSupport2}
    \lVert u\rVert_{\L^p(\RR^n)}\lesssim_{p,s_0,n} |K|^{\frac{s_0}{n}}\lVert u\rVert_{\dot{\H}^{s_0,p}(\RR^n)}\lesssim_{p,s_0,n} |K|^{\frac{s_0}{n}}\lVert u\rVert_{{\L}^{p}(\RR^n)}^{1-\frac{s_0}{s}}\lVert u\rVert_{\dot{\H}^{s,p}(\RR^n)}^{\frac{s_0}{s}}.
\end{align}
Thus, one has for some $C_{s_0,p,n}>0$, and for all $\epsilon>0$, 
\begin{align}\label{eq:ProofEquivCompactSupport3}
    \lVert u\rVert_{\L^p(\RR^n)}\leqslant C_{s_0,p,n} |K|^{\frac{s_0}{n}}\left( \epsilon\lVert u\rVert_{{\L}^{p}(\RR^n)} + \epsilon^{1-s/s_0} \lVert u\rVert_{\dot{\H}^{s,p}(\RR^n)}\right).
\end{align}
Thereby, choosing $\epsilon = \frac{1}{2 C_{s_0,p,n}|K|^{\frac{s_0}{n}}}$ it holds
\begin{align}\label{eq:ProofEquivCompactSupport4}
    \lVert u\rVert_{\L^p(\RR^n)}\lesssim_{p,s,s_0,n} |K|^{\frac{s}{n}}\lVert u\rVert_{\dot{\H}^{s,p}(\RR^n)}.
\end{align}
If $-\frac{n}{p'}<s<0$, one always has $\dot{\H}^{s,p}(\RR^n)\hookrightarrow\H^{s,p}(\RR^n)$. Thus for the reverse embedding, for $\varphi\in\Ccinfty(\RR^{n})$ with $\varphi = 1$ on $K$, $\supp \varphi \subset 2K$, thus $\varphi u =u$ in $\Distrib(\RR^n)$, and by the postponed Proposition~\ref{prop:EmbeddingsHomogeneousSobolevBesovinMultiplierSpaces} and duality, one has for $u=a+b\in\dot{\H}^{s,p}(\RR^{n})+\L^p(\RR^n)$,
\begin{align*}
    \lVert u\rVert_{\dot{\H}^{s,p}(\RR^n)} = \lVert \varphi u\rVert_{\dot{\H}^{s,p}(\RR^n)}&\lesssim_{p,s,n} \lVert \varphi a\rVert_{\dot{\H}^{s,p}(\RR^n)} + \lVert \varphi b\rVert_{\dot{\H}^{s,p}(\RR^n)}\\
    &\lesssim_{p,s,n} \lVert\varphi\rVert_{\dot{\B}^{-s}_{\frac{n}{-s},1}(\RR^n)}\lVert a\rVert_{\dot{\H}^{s,p}(\RR^n)} + \lVert \varphi b\rVert_{\L^{\frac{np}{n+ps}}(\RR^n)}\\
    &\lesssim_{p,s,n} \lVert\varphi\rVert_{\dot{\B}^{-s}_{\frac{n}{-s},1}(\RR^n)}\lVert a\rVert_{\dot{\H}^{s,p}(\RR^n)} + |2K|^\frac{-s}{n}\lVert \varphi b\rVert_{\L^{p}(\RR^n)}\\
    &\lesssim_{p,s,n} \lVert\varphi\rVert_{\dot{\B}^{-s}_{\frac{n}{-s},1}(\RR^n)}\lVert a\rVert_{\dot{\H}^{s,p}(\RR^n)} + |K|^\frac{-s}{n}\lVert\varphi\rVert_{\L^\infty(\RR^n)}\lVert  b\rVert_{\L^{p}(\RR^n)}\\
    &\lesssim_{p,s,n}  \lVert\varphi\rVert_{\dot{\B}^{-s}_{\frac{n}{-s},1}(\RR^n)}(1+|K|^{\frac{-s}{n}})(\lVert a\rVert_{\dot{\H}^{s,p}(\RR^n)} + \lVert b\rVert_{\L^{p}(\RR^n)}).
\end{align*}
Here, we did take advantage of the Sobolev embedding $\dot{\B}^{-s}_{\frac{n}{-s},1}(\RR^n)\hookrightarrow\L^\infty(\RR^n)$. Since $\H^{s,p}(\RR^{n})=\dot{\H}^{s,p}(\RR^{n})+\L^p(\RR^n)$ with equivalence of norms, the result follows by taking the infimum on all such pairs $(a,b)$ which yields
\begin{align*}
    \lVert u\rVert_{\dot{\H}^{s,p}(\RR^n)}&\lesssim_{p,s,n} (1+|K|^{\frac{-s}{n}})\lVert u\rVert_{{\H}^{s,p}(\RR^n)}.
\end{align*}

\textbf{Step 2:} Now, for Besov spaces, $p,q\in[1,\infty]$.

\textbf{Step 2.1:} We start with the case $s>0$, one always has ${\B}^{s}_{p,q}(\RR^n)\hookrightarrow\dot{\B}^{s}_{p,q}(\RR^n)$. If $p=1$, as before, by Sobolev embeddings for $s\in(0,n]$, $\dot{\B}^{s}_{1,1}(\RR^n)\hookrightarrow\L^{\frac{n}{n-s}}(\RR^n)$,
\begin{align*}
    \lVert u\rVert_{\L^1(\RR^n)} \leqslant |K|^{\frac{s}{n}}\lVert u\rVert_{\L^{\frac{n}{n-s}}(\RR^n)}\lesssim_{p,s,n} |K|^{\frac{s}{n}}\lVert u\rVert_{\dot{\B}^{s}_{1,1}(\RR^n)}.
\end{align*}
If $s>n$, or $s>0$ with $q\in[1,\infty]$, for $0<s_0<\min(s,n)$, by interpolation inequalities, following \eqref{eq:ProofEquivCompactSupport1}--\eqref{eq:ProofEquivCompactSupport4} from Step 1, one obtains as before
\begin{align*}
    \lVert u\rVert_{\L^1(\RR^n)}\lesssim_{s,s_0,n} |K|^\frac{s}{n}\lVert u\rVert_{\dot{\B}^{s}_{1,q}(\RR^n)}.
\end{align*}
Similar arguments apply for ${\B}^{s}_{p,q}$, $p\in(1,\infty)$, $q\in[1,\infty]$, $s>0$ (replacing the use of $\L^1$, $s\in(0,n]$ and ${\B}^{s}_{1,1}$ by  respectively $\L^p$, $s\in(0,\sfrac{n}{p}]$ and ${\B}^{s}_{p,1}$).

\textbf{Step 2.2:} If $-n/p'<s<0$ and $p\in(1,\infty]$, $q\in[1,\infty]$ since $\B^{s}_{p,q}(\RR^{n})=\dot{\B}^{s}_{p,q}(\RR^{n})+\L^{p}(\RR^{n})$, then for $u=a+b\in\dot{\B}^{s}_{p,q}(\RR^{n})+\L^{p}(\RR^{n})$, as before by the Sobolev embeddings $\L^{\frac{np}{n-sp},q}(\RR^n)\hookrightarrow \dot{\B}^{s}_{p,q}(\RR^n)$ and $\dot{\B}^{-s}_{\frac{n}{-s},1}(\RR^n)\hookrightarrow\L^\infty(\RR^n)$, and by interpolation inequalities, it holds that
\begin{align*}
    \lVert u\rVert_{\dot{\B}^{s}_{p,q}(\RR^{n})} &\leqslant \lVert \varphi a\rVert_{\dot{\B}^{s}_{p,q}(\RR^{n})} +\lVert \varphi b\rVert_{\dot{\B}^{s}_{p,q}(\RR^{n})}\\
    &\lesssim_{p,s,n} \lVert \varphi a\rVert_{\dot{\B}^{s}_{p,q}(\RR^{n})} +\lVert \varphi b\rVert_{\L^{\frac{np}{n-ps},q}(\RR^{n})}\\
    &\lesssim_{p,s,n}\lVert\varphi\rVert_{\dot{\B}^{-s}_{\frac{n}{-s},1}(\RR^n)}\lVert a\rVert_{\dot{\B}^{s}_{p,q}(\RR^{n})} +\lVert  \varphi b\rVert_{\L^{1}(\RR^{n})}^\frac{-sp'}{n}\lVert  \varphi b\rVert_{\L^{p}(\RR^{n})}^{1+\frac{sp'}{n}}\\
    &\lesssim_{p,s,n}\lVert\varphi\rVert_{\dot{\B}^{-s}_{\frac{n}{-s},1}(\RR^n)}(1+|K|^{\frac{-s}{n}})\left( \lVert a\rVert_{\dot{\B}^{s}_{p,q}(\RR^{n})} +\lVert  b\rVert_{\L^{p}(\RR^{n})}\right).
\end{align*}
 Above, $\L^{r,q}(\RR^n):=(\L^1(\RR^n),\L^\infty(\RR^n))_{1-\frac{1}{r},q}$ is the standard Lorentz space of parameters $r\in(1,\infty)$, $q\in[1,\infty]$. Taking the infimum over all such pairs $(a,b)$ yields 
\begin{align*}
    \lVert u\rVert_{\dot{\B}^{s}_{p,q}(\RR^{n})} &\lesssim_{s,n}(1+|K|^{\frac{-s}{n}}) \lVert u\rVert_{{\B}^{s}_{p,q}(\RR^{n})}.
\end{align*}

\textbf{Step 2.3:} For $s=0$ and $p=\infty$, if $u$ is compactly supported then so is $\nabla u$ and by the previously established estimates it holds that
\begin{align*}
    \lVert u\rVert_{\dot{\B}^{0}_{\infty,q}(\RR^{n})} \sim_{s,n} \lVert \nabla u\rVert_{\dot{\B}^{-1}_{\infty,q}(\RR^{n})} \lesssim_{s,n} (1+|K|^\frac{1}{n})\lVert \nabla u\rVert_{{\B}^{-1}_{\infty,q}(\RR^{n})} \lesssim_{s,n}  (1+|K|^\frac{1}{n}) \lVert u\rVert_{{\B}^{0}_{\infty,q}(\RR^{n})}.
\end{align*}
Similarly, for the reverse estimates
\begin{align*}
    \lVert u\rVert_{{\B}^{0}_{\infty,q}(\RR^{n})} &\sim_{
    \varepsilon,n} \lVert u\rVert_{{\B}^{-\sfrac{1}{2}}_{\infty,q}(\RR^{n})} + \lVert \nabla u\rVert_{{\B}^{-1}_{\infty,q}(\RR^{n})}\\
    &\lesssim_{\varepsilon,n} (1+|K|^\frac{1}{2n}) \lVert u\rVert_{\dot{\B}^{-\sfrac{1}{2}}_{\infty,q}(\RR^{n})} + (1+|K|^\frac{1}{n}) \lVert \nabla u\rVert_{\dot{\B}^{-1}_{\infty,q}(\RR^{n})}.
\end{align*}
One clearly has $\lVert \nabla u\rVert_{\dot{\B}^{-1}_{\infty,q}(\RR^{n})}\lesssim_{s,n}\lVert u\rVert_{\dot{\B}^{0}_{\infty,q}(\RR^{n})}$. Thus, we focus on estimating $\lVert u\rVert_{\dot{\B}^{-\sfrac{1}{2}}_{\infty,q}(\RR^{n})}$. By interpolation inequalities, provided $0<\varepsilon<\theta$, one obtains by the case $s<0$, for all $\epsilon>0$
\begin{align*}
    \lVert u\rVert_{\dot{\B}^{-\sfrac{1}{2}}_{\infty,q}(\RR^{n})}  &\lesssim_{n} \lVert u\rVert_{\dot{\B}^{-1}_{\infty,q}(\RR^{n})}^{\frac{1}{2}} \lVert u\rVert_{\dot{\B}^{0}_{\infty,q}(\RR^{n})}^{\frac{1}{2}} \\
    &\lesssim_{n} \Big((1+|K|^\frac{1}{n})\lVert u\rVert_{{\B}^{-1}_{\infty,q}(\RR^{n})}\Big)^{\frac{1}{2}} \lVert u\rVert_{\dot{\B}^{0}_{\infty,q}(\RR^{n})}^{\frac{1}{2}}\\
    &\leqslant C_{n} (1+|K|^\frac{1}{n})^\frac{1}{2}\left[ \epsilon \lVert u\rVert_{{\B}^{-1}_{\infty,q}(\RR^{n})}+ \frac{1}{\epsilon}  \lVert u\rVert_{\dot{\B}^{0}_{\infty,q}(\RR^{n})}\right].
\end{align*}
Recalling that ${\B}^{0}_{\infty,q}(\RR^{n}) \hookrightarrow {\B}^{-1}_{\infty,q}(\RR^{n})$. Choosing $\epsilon=\frac{1}{2C_{n} (1+|K|^\frac{1}{n})^\frac{1}{2} (1+|K|^\frac{1}{n})}$, it yields
\begin{align*}
    \lVert u\rVert_{{\B}^{0}_{\infty,q}(\RR^{n})} 
    &\leqslant \frac{1}{2}\lVert u\rVert_{{\B}^{0}_{\infty,q}(\RR^{n})} + 2C_{n}^2(1+|K|^\frac{1}{n})(1+|K|^\frac{1}{2n})^2 \lVert  u\rVert_{\dot{\B}^{0}_{\infty,q}(\RR^{n})}.
\end{align*}
Hence,
\begin{align*}
    \lVert u\rVert_{{\B}^{0}_{\infty,q}(\RR^{n})} 
    &\lesssim_{n}(1+|K|^\frac{2}{n}) \lVert  u\rVert_{\dot{\B}^{0}_{\infty,q}(\RR^{n})}.
\end{align*}

\textbf{Step 2.4 :} The case $p=\infty$, $s>0$. First, we assume $s\in(0,n)$.

By inrterpolation inequalities, the case $s<0$ from Step 2.2, and the embedding $\L^\infty(\RR^n)\hookrightarrow\B^{-s}_{\infty,q}(\RR^n)$, we have for all $\epsilon>0$
\begin{align*}
    \lVert u\rVert_{\L^\infty(\RR^n)} &\lesssim_{n}\lVert u\rVert_{\dot{\B}^{-s}_{\infty,q}(\RR^{n})}^{\frac{1}{2}} \lVert u\rVert_{\dot{\B}^{s}_{\infty,q}(\RR^{n})}^{\frac{1}{2}} \\
    &\lesssim_{s,n} \Big((1+|K|^\frac{s}{n})\lVert u\rVert_{{\B}^{-s}_{\infty,q}(\RR^{n})}\Big)^{\frac{1}{2}} \lVert u\rVert_{\dot{\B}^{s}_{\infty,q}(\RR^{n})}^{\frac{1}{2}}\\
    &\lesssim_{s,n} \Big((1+|K|^\frac{s}{n})\lVert u\rVert_{{\L}^{\infty}(\RR^{n})}\Big)^{\frac{1}{2}} \lVert u\rVert_{\dot{\B}^{s}_{\infty,q}(\RR^{n})}^{\frac{1}{2}}\\
    &\leqslant C_{s,n} (1+|K|^\frac{s}{n})^\frac{1}{2}\left[ \epsilon \lVert u\rVert_{\L^\infty(\RR^{n})}+ \frac{1}{\epsilon}  \lVert u\rVert_{\dot{\B}^{s}_{\infty,q}(\RR^{n})}\right].
\end{align*}
As in previous steps, one obtains
\begin{align*}
     \lVert u\rVert_{\L^\infty(\RR^n)} \lesssim_{s,n} (1+|K|^\frac{s}{n})\lVert u\rVert_{\dot{\B}^{s}_{\infty,q}(\RR^{n})}.
\end{align*}
Bootstrapping the interpolation inequality argument for $s\geqslant n$, choosing $s_0\in(0,n)$, one obtains for all $\epsilon>0$
\begin{align*}
     \lVert u\rVert_{\L^\infty(\RR^n)} &\lesssim_{s_0,n} (1+|K|^\frac{s_0}{n})\lVert u\rVert_{\dot{\B}^{s_0}_{\infty,q}(\RR^{n})}\\
     &\lesssim_{s_0,s,n} (1+|K|^\frac{s_0}{n}) \lVert u\rVert_{\L^\infty(\RR^n)}^{1-\frac{s_0}{s}}\lVert u\rVert_{\dot{\B}^{s}_{\infty,q}(\RR^{n})}^{\frac{s_0}{s}}\\
     &\leqslant C_{s_0,s,n}(1+|K|^\frac{s_0}{n}) \left[ \epsilon \lVert u\rVert_{\L^\infty(\RR^n)}+ \epsilon^{1-s/s_0}  \lVert u\rVert_{\dot{\B}^{s}_{\infty,q}(\RR^{n})}\right].
\end{align*}
Choosing $\epsilon= \frac{1}{2C_{s_0,s,n}(1+|K|^\frac{s_0}{n})}$, one obtains
\begin{align*}
     \lVert u\rVert_{\L^\infty(\RR^n)} \lesssim_{s_0,s,n}(1+|K|^\frac{s_0}{n})^\frac{s}{s_0} \lVert u\rVert_{\dot{\B}^{s}_{\infty,q}(\RR^{n})} \lesssim_{s_0,s,n}(1+|K|^\frac{s}{n})\lVert u\rVert_{\dot{\B}^{s}_{\infty,q}(\RR^{n})}.
\end{align*}

\textbf{Step 2.5:} The case $s=0$, $p\in(1,\infty)$. Since $K$ is compact, there exists a finite family of disjoint cubes $(Q_j)_{j\in\llb1, N\rrb}$, such that
\begin{align*}
    K \subset \bigcup_{j=1}^N Q_j=:\mathcal{Q}
\end{align*}
with $\frac{1}{2}|\mathcal{Q}|\leqslant |K|\leqslant |\mathcal{Q}|$.

By Proposition~\ref{prop:FundamentalExtby0HomFuncSpaces}, $\mathcal{Q}$ is a Lipschitz domain, and  one has $u=\mathbbm{1}_{\mathcal{Q}} u \in{\dot{\B}^{0}_{p,q,0}(\mathcal{Q})}={\dot{\B}^{0}_{p,q}(\mathcal{Q})}$ with 
\begin{align*}
    \lVert u\rVert_{\dot{\B}^{0}_{p,q}(\RR^{n})}\sim_{p,n} \lVert u\rVert_{\dot{\B}^{0}_{p,q}(\mathcal{Q})}  \lesssim_{p,\varepsilon,n} (1+|\mathcal{Q}|^\frac{\varepsilon}{n})\lVert u\rVert_{{\B}^{0}_{p,q}(\mathcal{Q})} \lesssim_{p,\varepsilon,n}  (1+|\mathcal{Q}|^\frac{\varepsilon}{n})\lVert u\rVert_{{\B}^{0}_{p,q}(\RR^n)},
\end{align*}
for some $0<\varepsilon<\sfrac{1}{p}$. The same argument yields
\begin{align*}
    \lVert u\rVert_{{\B}^{0}_{p,q}(\RR^{n})}\lesssim_{p,\varepsilon,n}  (1+|\mathcal{Q}|^\frac{\varepsilon}{n})\lVert u\rVert_{\dot{\B}^{0}_{p,q}(\RR^n)}.
\end{align*}
This concludes the argument since $|\mathcal{Q}|\leqslant 2|K|$.

Finally, the remaining case $s=-\sfrac{n}{p'}$ with $q=\infty$, $p\in(1,\infty]$ admits a proof to the one of Step 1 and Step 2.1, due to the Sobolev embedding $\L^1(\RR^n)\hookrightarrow\dot{\BesSmo}^{-\sfrac{n}{p'}}_{p,\infty}(\RR^n)$.
\end{proof}

We can derive the following Poincaré-type inequalities for compactly supported distributions that do not necessarily vanish on the boundary.

\begin{corollary}\label{cor:GeneralizedPoincaréHalfSpace}Let $p,q\in[1,\infty]$. Let $u\in\S'_h(\RR^n_+)$ with relatively compact support $K\subset \overline{\RR^n_+}$ (\textit{i.e.} $K\cap\partial\RR^{n}_+$ migh not be empty). If $-1+{\sfrac{1}{p}}<s<{\sfrac{1}{p}}$, then it holds that
\begin{align*}
    \nabla u\in\dot{\H}^{s,p}(\RR^n_+) &\implies u\in{\H}^{s+1,p}(\RR^n_+)\text{, }1<p<\infty\text{;}\\
     \nabla u\in\dot{\B}^{s}_{p,q}(\RR^n_+) &\implies u\in{\B}^{s+1}_{p,q}(\RR^n_+)\text{, }1\leqslant p,q\leqslant\infty\text{;}
\end{align*}
with an estimate
\begin{align*}
    \lVert u \rVert_{{\X}^{s,p}(\RR^n_+)} \sim_{p,s,n,K}\lVert u \rVert_{\dot{\X}^{s,p}(\RR^n_+)} \lesssim_{p,s,n,K}  \lVert \nabla u \rVert_{\dot{\X}^{s,p}(\RR^n_+)}\sim_{p,s,n,K} \lVert \nabla u \rVert_{{\X}^{s,p}(\RR^n_+)}.
\end{align*}
\end{corollary}

\begin{proof} Provided $\E:=\E_\mathcal{N}$ is the extension operator from Lemma~\ref{lem:ExtDirNeuRn+} below, it holds that $\E u$ has compact support in $\tilde{K} = K\cup\{(x',- x_n)\,:\,(x', x_n)\in K\}$, and choose a ball $\B$ such that $\tilde{K}\subset B$, by the definition of function spaces by restriction and Lemma~\ref{lem:CompactEquiHomInhom},
\begin{align*}
    \lVert u \rVert_{\dot{\X}^{s,p}(\RR^n_+)} \leqslant \lVert \E u \rVert_{\dot{\X}^{s,p}(\RR^n)}  \lesssim_{p,s,n,\B} \lVert \E u \rVert_{{\X}^{s,p}(\RR^n)}=\lVert \E u \rVert_{{\X}^{s,p}_0(\B)}.
\end{align*}
Since $\E u$ has compact support, so does $\nabla\E u$, and by Lemma~\ref{lem:CompactEquiHomInhom}, one has
\begin{align*}
    \lVert \nabla \E u \rVert_{{\X}^{s,p}(\RR^n)}=\lVert \nabla \E u \rVert_{{\X}^{s,p}_0(B)}\sim_{p,s,n,B} \lVert \nabla \E u \rVert_{\dot{\X}^{s,p}(\RR^n)} \lesssim_{p,s,n} \lVert \nabla  u \rVert_{\dot{\X}^{s,p}(\RR^n_+)} .
\end{align*}
Therefore, it suffices to prove that
\begin{align*}
    \lVert v \rVert_{{\X}^{s,p}_0(\B)} \lesssim_{p,s,n} \lVert \nabla v \rVert_{{\X}^{s,p}_0(\B)},\,\forall v\in{\X}^{s+1,p}_0(\B).
\end{align*}
Up to translation, without loss of generality, we can assume $B\subset\RR^n_+$ and $0\in\partial B$. In this case, for $A=\mathrm{diam}(\B)$, one has
\begin{align*}
    v = \mathbbm{1}_{\B}[\mathbbm{1}_{[-A,0]}\ast_{n} \partial_{x_n}v]
\end{align*}
where $\ast_{n}$ stands for the convolution operator with respect to the $n$-th variable. By Proposition~\ref{prop:FundamentalExtby0HomFuncSpaces}, it holds that
\begin{align*}
    \lVert v \rVert_{{\X}^{s,p}_0(\B)} \lesssim_{p,s,n,\B} \lVert \mathbbm{1}_{[-A,0]}\ast_{n} \partial_{x_n}v \rVert_{{\X}^{s,p}(\RR^n)}. 
\end{align*}
If ${\X}^{s,p}={\H}^{s,p}$, then it holds that $(\I-\Delta)^{\sfrac{s}{2}}[\mathbbm{1}_{[-A,0]}\ast_{n} \partial_{x_n}v] = \mathbbm{1}_{[-A,0]}\ast_{n} [(\I-\Delta)^{\sfrac{s}{2}}\partial_{x_n}v] $. Otherwise, if ${\X}^{s,p}={\B}^{s}_{p,q}$, then for all $j\geqslant -1$, $\Delta_j[\mathbbm{1}_{[-A,0]}\ast_{n} \partial_{x_n}v] = \mathbbm{1}_{[-A,0]}\ast_{n} [\Delta_j\partial_{x_n}v] $. Consequently, applying Young's inequality for the convolution with respect to  the last variable, one deduces
\begin{align*}
    \lVert \mathbbm{1}_{[-A,0]}\ast_{n} \partial_{x_n}v \rVert_{{\X}^{s,p}(\RR^n)} \leqslant \lVert \mathbbm{1}_{[-A,0]}\rVert_{{\L}^{1}(\RR)} \lVert \partial_{x_n}v \rVert_{{\X}^{s,p}(\RR^n)} = \mathrm{diam}(B)\lVert \partial_{x_n}v \rVert_{{\X}^{s,p}_0(B)}.
\end{align*}
Which yields the claimed result.
\end{proof}

A similar proof yields the following
\begin{lemma}[Poincaré inequality]\label{lem:Poincaré0}Let $p,q\in[1,\infty]$, $s\in(-1+\sfrac{1}{p},\sfrac{1}{p})$ and $\Omega$ be a bounded Lipschitz domain. For all $u\in\B^{s+1}_{p,q,0}(\Omega)$, one has
\begin{align*}
    \lVert u\rVert_{\B^{s}_{p,q}(\Omega)}\lesssim_{p,s,n,\Omega} \lVert\nabla u\rVert_{\B^{s}_{p,q}(\Omega)}.
\end{align*}
The result still holds with $\H^{s,p}$ instead of $\B^{s}_{p,q}$ whenever $1<p<\infty$.
\end{lemma}

We now provide Poincaré-Wirtinger type inequalities, in particular the Ne\v{c}as negative norm Theorem.
\begin{lemma}[Ne\v{c}as negative norm Theorem]\label{lem:NegativeNormThm}Let $p,q\in[1,\infty]$, $s\in(-1+\sfrac{1}{p},\infty)$. Let $\Omega$ be a bounded Lipschitz domain. For all $u\in\B^{s+1}_{p,q}(\Omega)$,  one has
\begin{align*}
    \lVert u-u_{\Omega}\rVert_{\B^{s}_{p,q}(\Omega)}\lesssim_{p,s,n,\Omega} \lVert\nabla u\rVert_{\B^{s-1}_{p,q}(\Omega)},
\end{align*}
where $u_{\Omega}= \frac{1}{|\Omega|}\langle u,\mathbbm{1}_{\Omega}\rangle_{\Omega}$.

The result still holds with $\H^{s,p}$ instead of $\B^{s}_{p,q}$ whenever $1<p<\infty$.
\end{lemma}

\begin{lemma}[Poincaré-Wirtinger inequality]\label{lem:Poincaré}Let $p,q\in[1,\infty]$, $s\in(-2+\sfrac{1}{p},\infty)$. Let $\Omega$ be a bounded Lipschitz domain. For all $u\in\B^{s+1}_{p,q}(\Omega)$,  one has
\begin{align*}
    \lVert u-u_{\Omega}\rVert_{\B^{s}_{p,q}(\Omega)}\lesssim_{p,s,n,\Omega} \lVert\nabla u\rVert_{\B^{s}_{p,q}(\Omega)},
\end{align*}
where $u_{\Omega}= \frac{1}{|\Omega|}\langle u,\mathbbm{1}_{\Omega}\rangle_{\Omega}$.

The result still holds with $\H^{s,p}$ instead of $\B^{s}_{p,q}$ whenever $1<p<\infty$.
\end{lemma}

\begin{remark} The restriction $s>-2+\sfrac{1}{p}$ (or $s>-1+{\sfrac{1}{p}}$ in Lemma~\ref{lem:NegativeNormThm}) is due to $\langle u,\mathbbm{1}_{\Omega}\rangle_{\Omega}$ to be meaningful. Indeed $\mathbbm{1}_{\Omega}\in(\X^{s+1,p}(\Omega))'=\X'^{-s-1,p'}_{0}(\Omega)$ requires $-s-1<\sfrac{1}{p'}$.
\end{remark}

\begin{proof}[of Lemmas~\ref{lem:NegativeNormThm}~\&~\ref{lem:Poincaré}]Without loss of generality, we can assume $u_\Omega=0$. Consider for $\alpha\in\RR$ the closed subspace
\begin{align*}
   \X^{\alpha,p}_{\mfree}(\Omega)=\{u\in\X^{\alpha,p}(\Omega)\,:\,\langle u,\mathbbm{1}_{\Omega}\rangle_{\Omega}=0\}.
\end{align*}
Let $u\in\B^{s}_{p,q,\mfree}(\Omega)$, we assume $s<1/p$, then one can write
\begin{align*}
    \lVert u\rVert_{\B^{s}_{p,q}(\Omega)} = \sup_{\substack{\varphi\in\B^{-s}_{p',q',0}(\Omega),\\ \lVert \varphi\rVert_{\B^{-s}_{p',q'}(\RR^n)}\leqslant 1}}\,|\langle u,\varphi\rangle_\Omega| = \sup_{\substack{\varphi\in\B^{-s}_{p',q',0}(\Omega),\\ \lVert \varphi\rVert_{\B^{-s}_{p',q'}(\RR^n)}\leqslant 1}}\,|\langle u,\varphi-\varphi_{\Omega}\rangle_\Omega|.
\end{align*}
By the use of Bogoski\u{i} operators, Theorem~\ref{thm:PotentialOpBddLipDomCcinfty},
\begin{align*}
    \lVert u\rVert_{\B^{s}_{p,q}(\Omega)} & \sim_{s,n} \sup_{\substack{\varphi\in\B^{-s}_{p',q',0}(\Omega),\\ \lVert \varphi\rVert_{\B^{-s}_{p',q'}(\RR^n)}\leqslant 1}}\,|\langle u,\varphi-\varphi_{\Omega}\rangle_\Omega|\\
    &\sim_{p,s,n,\Omega} \sup_{\substack{\boldsymbol{\varphi}\in\B^{-s+1}_{p',q',0}(\Omega)^n,\\ \lVert \boldsymbol{\varphi}\rVert_{\B^{-s+1}_{p',q'}(\RR^n)}\leqslant 1}}\,|\langle u,\div \boldsymbol{\varphi}\rangle_\Omega|\\
     &\sim_{p,s,n,\Omega}  \sup_{\substack{\boldsymbol{\varphi}\in\B^{-s+1}_{p',q',0}(\Omega)^n,\\ \lVert \boldsymbol{\varphi}\rVert_{\B^{-s+1}_{p',q'}(\RR^n)}\leqslant 1}}\,|\langle \nabla u,\boldsymbol{\varphi}\rangle_\Omega|\\
     &\lesssim_{p,s,n,\Omega}  \lVert \nabla u\rVert_{\B^{s-1}_{p,q}(\Omega)} .
\end{align*}
Now, for $1/p<s+1<1+1/p$, $u\in{\B^{s+1}_{p,q,\mfree}(\Omega)}$, 
\begin{align*}
    \lVert u\rVert_{\B^{s+1}_{p,q}(\Omega)}&\sim_{p,s,n,\Omega} \lVert u\rVert_{\B^{s}_{p,q}(\Omega)}+ \lVert \nabla u\rVert_{\B^{s}_{p,q}(\Omega)}\\
    &\lesssim_{p,s,n,\Omega} \lVert \nabla u\rVert_{\B^{s-1}_{p,q}(\Omega)}+ \lVert \nabla u\rVert_{\B^{s}_{p,q}(\Omega)}\\
    &\lesssim_{p,s,n,\Omega} \lVert \nabla u\rVert_{\B^{s}_{p,q}(\Omega)}.
\end{align*}
Therefore, by induction, one obtains the result for all $s\in(-1+\sfrac{1}{p},\infty)\setminus (\sfrac{1}{p}+\NN)$. It remains to prove the case $s=\frac{1}{p}$, the remaining cases will follow by induction.

For $u\in{\B^{\frac{1}{p}+1}_{p,q,\mfree}(\Omega)}$, for $\varepsilon>0$ sufficiently small fixed, for any $\kappa>0$,
\begin{align*}
    \lVert u\rVert_{\B^{\frac{1}{p}+1}_{p,q}(\Omega)}&\sim_{p,s,n,\Omega} \lVert u\rVert_{\B^{\frac{1}{p}}_{p,q}(\Omega)}+ \lVert \nabla u\rVert_{\B^{\frac{1}{p}}_{p,q}(\Omega)}\\
    &\lesssim_{p,s,n,\Omega} \lVert u\rVert_{\B^{\frac{1}{p}-\varepsilon}_{p,q}(\Omega)}^\frac{1}{2}\lVert u\rVert_{\B^{\frac{1}{p}+\varepsilon}_{p,q}(\Omega)}^\frac{1}{2}+ \lVert \nabla u\rVert_{\B^{\frac{1}{p}}_{p,q}(\Omega)}\\
    &\lesssim_{p,s,n,\Omega} \lVert u\rVert_{\B^{\frac{1}{p}-\varepsilon}_{p,q}(\Omega)}^\frac{1}{2}\lVert \nabla u\rVert_{\B^{\frac{1}{p}-1+\varepsilon}_{p,q}(\Omega)}^\frac{1}{2}+ \lVert \nabla u\rVert_{\B^{\frac{1}{p}}_{p,q}(\Omega)}\\
    &\lesssim_{p,s,n,\Omega} \kappa\lVert u\rVert_{\B^{\frac{1}{p}-\varepsilon}_{p,q}(\Omega)} +\frac{1}{\kappa}\lVert \nabla u\rVert_{\B^{\frac{1}{p}-1+\varepsilon}_{p,q}(\Omega)}+ \lVert \nabla u\rVert_{\B^{\frac{1}{p}}_{p,q}(\Omega)}\\&\lesssim_{p,s,n,\Omega} \kappa\lVert u\rVert_{\B^{\frac{1}{p}-\varepsilon}_{p,q}(\Omega)} +\big(\frac{1}{\kappa}+1\big) \lVert \nabla u\rVert_{\B^{\frac{1}{p}}_{p,q}(\Omega)}.
\end{align*}
Choosing $\kappa$ small enough yields the result. Other function spaces admit a similar proof.
\end{proof}

\paragraph{Some additional extension properties in the case of the flat half-space.}

Following \cite[Section~5]{Gaudin2022}, we introduce the following extension operators defined for any measurable function $u\,:\,\RR^n_+\longrightarrow\CC$, for almost every $x=(x',x_n)\in\RR^{n}$:
\begin{align*}
    \E_\mathcal{D}u (x',x_n) &:= \begin{cases}
  u(x',x_n)\,&\text{, if } (x',x_n)\in\RR^{n-1}\times\RR_+\text{,}\\    
  -u(x',-x_n)\,&\text{, if } (x',x_n)\in\RR^{n-1}\times\RR_-^*\text{;}
\end{cases}\\
    \E_\mathcal{N}u (x',x_n) &:= \begin{cases}
  u(x',x_n)\,&\text{, if } (x',x_n)\in\RR^{n-1}\times\RR_+\text{,}\\    
  u(x',-x_n)\,&\text{, if } (x',x_n)\in\RR^{n-1}\times\RR_-^* \text{;}
\end{cases}\\
    \mathcal{E}_0u (x',x_n) &:= \begin{cases}
  u(x',x_n)\,&\text{, if } (x',x_n)\in\RR^{n-1}\times\RR_+\text{,}\\    
  0\,&\text{, if } (x',x_n)\in\RR^{n-1}\times\RR_-^* \text{.}
\end{cases}
\end{align*}

The sharp boundedness of these extension operators in the case homogeneous functions spaces --and including the endpoint cases-- will be of a paramount importance for the study of Dirichlet-Stokes problem in Section~\ref{Sec:StokesHalfSpace}.

\begin{lemma}[ {\cite[Lemma~3.6]{Gaudin2023Lip}} ] \label{lem:ExtDirNeuLinftyRn+} Let $\mathcal{J\in\{\mathcal{D},\mathcal{N}\}}$. For all $u\in\L^\infty_h(\RR^n_+)$, one has $\E_{\mathcal{J}}u\in \L^\infty_h(\RR^n)$ with 
\begin{align*}
    \lVert \E_{\mathcal{J}}u \rVert_{\L^\infty(\RR^n)} =  \lVert u \rVert_{\L^\infty(\RR^n_+)}.
\end{align*}
Furthermore,
\begin{itemize}
    \item $u\in\L^\infty_h(\RR^n_+)$ if and only if  $ \mathcal{E}_0{u}\in\L^\infty_h(\RR^n)$;
    \item $u\in\C^{0}_{ub,h}(\overline{\RR^n_+})$ if and only if $\E_\mathcal{N}{u}\in\C^{0}_{ub,h}({\RR}^n)$;
    \item $u\in\C^{0}_{ub,h,\mathcal{D}}(\overline{\RR^n_+})$ if and only if $\mathcal{E}_0 u\in\C^{0}_{ub,h}({\RR}^n)$.
\end{itemize}
\end{lemma}

\begin{lemma}\label{lem:ExtDirNeuRn+}Let $p,q\in[1,\infty]$, $s\in(-1+{\sfrac{1}{p}},2+{\sfrac{1}{p}})$, such that $s\neq{\sfrac{1}{p}}, 1+{\sfrac{1}{p}}$, and $\mathcal{J}\in\{\mathcal{D},\,\mathcal{N}\}$. For all $u\in\dot{\B}^{s}_{p,q,\mathcal{J}}(\RR^n_+)$, one has $\E_{\mathcal{J}}u\in \dot{\B}^{s}_{p,q}(\RR^n)$ with an estimate
\begin{align*}
    \lVert \E_{\mathcal{J}}u \rVert_{\dot{\B}^{s}_{p,q}(\RR^n)} \less_{p,s,n} \lVert u \rVert_{\dot{\B}^{s}_{p,q}(\RR^n_+)}
\end{align*}
Furthermore,
\begin{itemize}
    \item When $p=\infty$, we can replace $\dot{\B}^{s}_{\infty,q}$ by $\dot{\B}^{s,0}_{\infty,q}$.
    \item When $q=\infty$, we can replace $\dot{\B}^{s}_{p,\infty}$ by $\dot{\BesSmo}^{s}_{p,\infty}$.
    \item The result still holds for if we replace $\dot{\B}^{s}_{p,q}$ by $\dot{\H}^{s,p}$, $1<p<\infty$. 
    \item The result still holds for their inhomogeneous counterparts.
\end{itemize}
\end{lemma}

\begin{proof}[of Lemma~\ref{lem:ExtDirNeuLinftyRn+}] \textbf{Step 0:} The case $s\in(-1+{\sfrac{1}{p}},{\sfrac{1}{p}})$ from Lemma~\ref{lem:ExtDirNeuRn+}  is a direct consequence of Proposition~\ref{prop:FundamentalExtby0HomFuncSpaces} for $\Omega=\RR^n_+$.

\textbf{Step 1:} If $s\in({\sfrac{1}{p}},1+{\sfrac{1}{p}})$ and $\mathcal{J}=\mathcal{D}$, for $u\in\dot{\B}^{s}_{p,q,\mathcal{D}}(\RR^n_+)$, $\tilde{u}$, the extension by $0$ of $u$ to the whole $\RR^n$, belongs to by $\dot{\B}^{s}_{p,q}(\RR^n)$ by \cite[Proposition~4.21]{Gaudin2023Lip}. Therefore, it holds that $\E_\mathcal{D}u = \tilde{u}-\tilde{u}(\cdot,-\cdot) \in\dot{\B}^{s}_{p,q}(\RR^n)$. The case of other function spaces is similar.

\textbf{Step 2:} Let $s\in({\sfrac{1}{p}},1+{\sfrac{1}{p}})$ and $\mathcal{J}=\mathcal{N}$.

\textbf{Step 2.1:} We start with the case $p=\infty$, $0 < s<1$. By Lemma~\ref{lem:ExtDirNeuRn+}, one has the boundedness of
\begin{align}\label{eq:ProofEnOpC0ubh}
    \E_\mathcal{N}\,:\,\C_{ub,h}^0(\overline{\RR^n_+})\longrightarrow\C_{ub,h}^0({\RR^n}),
\end{align}
For all $u\in\dot{\W}^{1,\infty}(\mathbb{R}^n_+)\subset\C_{ub,h}^0({\overline{\RR^n_+}})$, since it holds almost everywhere that
\begin{align*}
    \nabla' \E_\mathcal{N} u =  \E_\mathcal{N} [\nabla'u]\text{, and } \partial_{x_n} \E_\mathcal{N} u =  \E_\mathcal{D} [ \partial_{x_n}  u]\text{,}
\end{align*}
one deduces
\begin{align*}
    \lVert \nabla \E_\mathcal{N} u \rVert_{\L^\infty(\mathbb{R}^n)} \leqslant \lVert \E_\mathcal{N}[ \nabla' u] \rVert_{\L^\infty(\mathbb{R}^n)} + \lVert \E_\mathcal{N}[ \partial_{x_n}u] \rVert_{\L^\infty(\mathbb{R}^n)}\leqslant 2 \lVert \nabla u \rVert_{\L^\infty(\mathbb{R}^n_+)}.
\end{align*}
Thus, combined with \eqref{eq:ProofEnOpC0ubh}, $\E_\mathcal{N} u\in{\W}^{1,\infty}(\mathbb{R}^n)\cap\C_{ub,h}^0({{\RR^n}})=\dot{\W}^{1,\infty}(\RR^n)$. Hence, we also did obtain the boundedness  
\begin{align}\label{eq:ProofEnOpLiph}
    \E_\mathcal{N}\,:\,\dot{\W}^{1,\infty}(\mathbb{R}^n_+)\longrightarrow\dot{\W}^{1,\infty}({\RR^n}).
\end{align}
By real interpolation applied between \eqref{eq:ProofEnOpLiph} and Step 1, one obtains the bounded linear operator
\begin{align}\label{eq:ProofEnOpBesovinfty-1<s<1}
    \E_\mathcal{N}\,:\, \dot{\B}^{s}_{\infty,q}({\RR^n_+})\longrightarrow\dot{\B}^{s}_{\infty,q}({\RR^n}),\, -1<s<1,\, q\in[1,\infty].
\end{align}

\textbf{Step 2.2:} The case $p\in[1,\infty)$, ${\sfrac{1}{p}}<s<1+{\sfrac{1}{p}}$.
First if $s<{\sfrac{n}{p}}$, then $\dot{\B}^{s}_{p,q,\mathcal{N}}(\RR^n_+)=\dot{\B}^{s}_{p,q}(\RR^n_+)$ and $\dot{\B}^{s}_{p,q}(\RR^n)$ are complete, therefore it suffices to perform the estimate.

If $s \geqslant {\sfrac{n}{p}}$   for $u\in\dot{\B}^{s}_{p,q}(\RR^n_+)$, by Sobolev embeddings, one has $u\in\dot{\B}^{s-\sfrac{n}{p}}_{\infty,q}(\RR^n_+)$. Therefore $\E_\mathcal{N} u\in\dot{\B}^{s-\sfrac{n}{p}}_{\infty,\infty}(\RR^n)\subset\C_{ub,h}^0(\RR^n)\subset\S'_h(\RR^n)$ by the end of the previous Step 2.1 \eqref{eq:ProofEnOpBesovinfty-1<s<1}. Hence, it reduces again to only show the estimates to prove membership and boundedness.

One can show in $\S'(\RR^n)$, that
\begin{align*}
    \nabla' \E_\mathcal{N} u =  \E_\mathcal{N} [\nabla'u]\text{, and } \partial_{x_n} \E_\mathcal{N} u =  \E_\mathcal{D} [ \partial_{x_n}  u]\text{.}
\end{align*}
Consequently, we can estimate
\begin{align*}
    \lVert \E_\mathcal{N} u \rVert_{\dot{\B}^{s}_{p,q}(\RR^n)}\sim_{s,p,n}\lVert \nabla \E_\mathcal{N} u \rVert_{\dot{\B}^{s-1}_{p,q}(\RR^n)} = \lVert ( \E_\mathcal{N}[\nabla' u],\E_\mathcal{D} [ \partial_{x_n}  u] ) \rVert_{\dot{\B}^{s-1}_{p,q}(\RR^n)} &\lesssim_{s,p,n}  \lVert \nabla u \rVert_{\dot{\B}^{s-1}_{p,q}(\RR^n_+)}\\ &\lesssim_{s,p,n}  \lVert u \rVert_{\dot{\B}^{s}_{p,q}(\RR^n_+)}. 
\end{align*}

\textbf{Step 3:} The case $p\in[1,\infty]$, $1+{\sfrac{1}{p}}<s<2+{\sfrac{1}{p}}$.

\textbf{Step 3.1:} We prove $\E_{\mathcal{N}}\dot{\B}^{s}_{p,q,\mathcal{J}}(\RR^n_+)\subset \S'_h(\RR^n)$.

As in Step 2.2, if $s<{\sfrac{n}{p}}$, then $\dot{\B}^{s}_{p,q,\mathcal{J}}(\RR^n_+)$ and $\dot{\B}^{s}_{p,q}(\RR^n)$ are complete, otherwise if $s>{\sfrac{n}{p}}$, $\dot{\B}^{s}_{p,q,\mathcal{J}}(\RR^n_+)\hookrightarrow \dot{\B}^{s-\sfrac{n}{p}}_{\infty,q}(\RR^n_+)= \S'_h\cap{\B}^{s-\sfrac{n}{p}}_{\infty,q}(\RR^n_+)$, so that by the end of Step 2.1 for $\mathcal{J}=\mathcal{N}$, or Step 1 for $\mathcal{J}=\mathcal{D}$, $\E_\mathcal{J} u\in\S'_h\cap{\B}^{s-\sfrac{n}{p}}_{\infty,q}(\RR^n)$ for all $u\in\dot{\B}^{s}_{p,q,\mathcal{J}}(\RR^n_+)\subset{\B}^{\min(\alpha,s-\sfrac{n}{p})}_{\infty,q}(\RR^n_+)$, provided $0<\alpha<1$.
If $s={\sfrac{n}{p}}$ and $\mathcal{J}=\mathcal{N}$, then the same argument applies with $\dot{\B}^{0}_{\infty,q}$ thanks to \eqref{eq:ProofEnOpBesovinfty-1<s<1}.

If $s\neq {\sfrac{n}{p}}$ and $\mathcal{J}=\mathcal{D}$,  the same arguments apply and are carried by the Step 1 instead of Step 2.1.

If $s={\sfrac{n}{p}}$ and $\mathcal{J}=\mathcal{D}$, by Sobolev embeddings and  following \cite[Proof of Proposition~4.20]{Gaudin2023Lip}, for $r\in(\max(p,n-1),\infty)$ such that ${\sfrac{n}{r}}<1+{\sfrac{1}{r}}$, it holds that
\begin{align*}
    \dot{\B}^{{\sfrac{n}{p}}}_{p,q,\mathcal{D}}(\RR^n_+) \hookrightarrow \dot{\B}^{{\sfrac{n}{r}}}_{r,q,\mathcal{D}}(\RR^n_+).
\end{align*}
Thereby, by Step 1, $\E_\mathcal{D}\dot{\B}^{{\sfrac{n}{p}}}_{p,q,\mathcal{D}}(\RR^n_+)\subset \E_\mathcal{D}\dot{\B}^{{\sfrac{n}{r}}}_{r,q,\mathcal{D}}(\RR^n_+) \subset\dot{\B}^{{\sfrac{n}{r}}}_{r,q}(\RR^n) \subset \S'_h(\RR^n)$.

\textbf{Step 3.2} We prove the estimates.

In any case, one has $\E_{\mathcal{J}}\dot{\B}^{s}_{p,q,\mathcal{J}}(\RR^n_+)\subset \S'_h(\RR^n)$. Finally, it suffices to perform the estimates. First notice that, one has for all $\varphi\in\S(\RR^n)$,
\begin{align*}
    \langle -\Delta \E_\mathcal{D}u, \varphi \rangle_{\RR^n} &= \langle -\Delta u, [\varphi-\varphi(\cdot,-\cdot)] \rangle_{\RR^n_+} = \langle\E_\mathcal{D}[ -\Delta u], \varphi \rangle_{\RR^n}\text{, if } u\in\dot{\B}^{s}_{p,q,\mathcal{D}}(\RR^n_+)\text{, }\\
    \langle -\Delta \E_\mathcal{N}u, \varphi \rangle_{\RR^n} &= \langle -\Delta u, [\varphi+\varphi(\cdot,-\cdot)] \rangle_{\RR^n_+} = \langle\E_\mathcal{N}[ -\Delta u], \varphi \rangle_{\RR^n}\text{, if } u\in\dot{\B}^{s}_{p,q,\mathcal{N}}(\RR^n_+)\text{.}
\end{align*}
So that, provided $u\in\dot{\B}^{s+2}_{p,q,\mathcal{J}}(\RR^n_+)$ where $\mathcal{J}\in\{\mathcal{D},\mathcal{N}\}$, by Step 1, it holds
\begin{align*}
    \lVert \E_\mathcal{J}u\rVert_{\dot{\B}^{s+2}_{p,q}(\RR^n)}\sim_{p,s,n}\lVert \Delta \E_\mathcal{J}u\rVert_{\dot{\B}^{s}_{p,q}(\RR^n)} = \lVert \E_\mathcal{J} \Delta u\rVert_{\dot{\B}^{s}_{p,q}(\RR^n)} \lesssim_{p,s,n} \lVert \Delta u\rVert_{\dot{\B}^{s}_{p,q}(\RR^n_+)}\lesssim_{p,s,n} \lVert u\rVert_{\dot{\B}^{s+2}_{p,q}(\RR^n_+)},
\end{align*}
which ends the proof.
\end{proof}

\begin{lemma}\label{lem:ExtOpNegativeBesovSpaces}Let $p,q\in[1,\infty]$, $s\in(-2+{\sfrac{1}{p}},{\sfrac{1}{p}})$. One has the bounded maps 
\begin{align*}
    \E_\mathcal{N}\,:\,\dot{\B}^{s}_{p,q,0}(\RR^n_+) \longrightarrow \dot{\B}^{s}_{p,q}(\RR^n),\text{ and }\E_\mathcal{D}\,:\,\dot{\B}^{s}_{p,q}(\RR^n_+) \longrightarrow \dot{\B}^{s}_{p,q}(\RR^n).
\end{align*}
Furthermore,
\begin{itemize}
    \item When $q=\infty$, we can replace $\dot{\B}^{s}_{p,\infty}$ by $\dot{\BesSmo}^{s}_{p,\infty}$.
    \item The result still holds if we replace $\dot{\B}^{s}_{p,q}$ by either $\dot{\B}^{s,0}_{\infty,q}$, or $\dot{\H}^{s,p}$ with $1<p<\infty$.
    \item The result still holds for their inhomogeneous counterparts.
\end{itemize}
\end{lemma}

\begin{proof}\textbf{Step 1:} If $s\in(-1+{\sfrac{1}{p}},{\sfrac{1}{p}})$, then one has the canonical identification $\dot{\B}^{s}_{p,q,0}(\RR^n_+)=\dot{\B}^{s}_{p,q}(\RR^n_+)$ and the result follows by Lemma~\ref{lem:ExtDirNeuRn+}.

\textbf{Step 2:} If $s\in(-2+{\sfrac{1}{p}},-1+{\sfrac{1}{p}})$. If $\mathcal{J}=\mathcal{N}$, then for $u\in\dot{\B}^{s}_{p,q,0}(\RR^n_+)\subset\dot{\B}^{s}_{p,q}(\RR^n)$, one has $\E_\mathcal{N}u = u + u(\cdot,-\cdot)\in\dot{\B}^{s}_{p,q}(\RR^n)$, and
\begin{align*}
    \lVert \E_\mathcal{N}u\rVert_{\dot{\B}^{s}_{p,q}(\RR^n)}\leqslant 2 \lVert  u\rVert_{\dot{\B}^{s}_{p,q}(\RR^n)} = 2\lVert u\rVert_{\dot{\B}^{s}_{p,q,0}(\RR^n_+)}.
\end{align*}
If $u\in\dot{\B}^{s}_{p,q}(\RR^n_+)$, for $\varphi\in\Ccinfty(\RR^n)$, it can be checked that
\begin{align*}
    \big\langle \E_\mathcal{D}u, \varphi \big\rangle_{\RR^n} = \big\langle u,  [\varphi-\varphi(\cdot,-\cdot)]\big\rangle_{\RR^n_+} = \big\langle u, [\I-\E_\mathcal{N}^-]\varphi_{|_{\RR^n_+}}\big\rangle_{\RR^n_+}.
\end{align*}
where we defined $\E_\mathcal{N}^{-}\varphi:= [\E_\mathcal{N}\varphi](\cdot',-\cdot)$. Indeed, by construction, one has $[\I-\E_\mathcal{N}^-]\varphi_{|_{\RR^n_+}}\in\dot{\B}^{-s}_{p',q',\mathcal{D}}(\RR^n_+)$, here ${\sfrac{1}{p'}}<-s<{1+\sfrac{1}{p'}}$, with the estimate
\begin{align*}
    \lVert [\I-\E_\mathcal{N}^-]\varphi_{|_{\RR^n_+}}\rVert_{\dot{\B}^{-s}_{p',q'}(\RR^n_+)}\leqslant \lVert [\I-\E_\mathcal{N}^-]\varphi\rVert_{\dot{\B}^{-s}_{p',q'}(\RR^n)}\lesssim_{p',-s,n} \lVert \varphi\rVert_{\dot{\B}^{-s}_{p',q'}(\RR^n)}.
\end{align*}
Since one has the canonical identification $\dot{\B}^{-s}_{p',q',\mathcal{D}}(\RR^n_+)=\dot{\B}^{-s}_{p',q',0}(\RR^n_+)$, see \cite[Proposition~4.21]{Gaudin2023Lip}, one can extend $[\I-\E_\mathcal{N}^-]\varphi_{|_{\RR^n_+}}$ to the whole space by $0$, yielding an element $\widetilde{[\I-\E_\mathcal{N}^-]\varphi}_{|_{\RR^n_+}}\in\dot{\B}^{-s}_{p',q',0}(\RR^n_+)$. Thereby, by duality \cite[Proposition~2.15~\&~Theorem~3.37]{Gaudin2023Lip} (for Sobolev spaces one uses \cite[Propositions~2.16~\&~3.43]{Gaudin2023Lip} instead), one deduces consequently
\begin{align*}
    |\big\langle \E_\mathcal{D}u, \varphi \big\rangle_{\RR^n}|\lesssim_{p,s,n}  \lVert u\rVert_{\dot{\B}^{s}_{p,q}(\RR^n_+)} \Big\lVert \widetilde{[\I-\E_\mathcal{N}^-]\varphi}_{|_{\RR^n_+}}\Big\rVert_{\dot{\B}^{-s}_{p',q',0}(\RR^n_+)} \lesssim_{p,s,n}  \lVert u\rVert_{\dot{\B}^{s}_{p,q}(\RR^n_+)} \lVert \varphi\rVert_{\dot{\B}^{-s}_{p',q'}(\RR^n)},
\end{align*}
and then taking the supremum over all $\varphi$ whose $\dot{\B}^{-s}_{p',q'}$-norm is less than $1$,
\begin{align*}
    \lVert \E_\mathcal{D}u\rVert_{\dot{\B}^{s}_{p,q}(\RR^n)}  \lesssim_{p,s,n}  \lVert u\rVert_{\dot{\B}^{s}_{p,q}(\RR^n_+)}.
\end{align*}
Finally, the case $s=-1+{\sfrac{1}{p}}$ follows by interpolation.
\end{proof}

\subsection{Sectorial operators, holomorphic functional calculus and maximal regularity}\label{sec:MaxRegIntro}

For more details about Sectorial operators, bounded olomorphic functional calculus, and $\L^q$-maximal regularity theory, motivated by the study of the Navier-Stokes and related equations, as well as the study of the underlying Stokes operator, the reader is encouraged to consult \cite{Monniaux2009MaxReg,PrussSimonett2016,BechtelBuiKunstmann2024} and the references therein.

\subsubsection{Sectorial Operators and Functional calculus.}

We introduce the following subsets of the complex plane
\begin{align*}
    \Sigma_\mu &:=\{ \,z\in\mathbb{C}^\ast\,:\,\lvert\mathrm{arg}(z)\rvert<\mu\,\}\text{, if } \mu\in(0,\pi)\text{,}\\
    \text{and} \quad\quad \mathrm{S}_\mu&:=\Sigma_\mu\cup(-\Sigma_\mu)\quad\quad\quad\,\,\quad\quad\text{, provided } \mu\in(0,\frac{\pi}{2}),
\end{align*}
we also define $\Sigma_0:= (0,\infty)$, $\mathrm{S}_0:=\RR$ and later we are going to consider their respective closure $\overline{\Sigma}_\mu$ and $\overline{\mathrm{S}}_\mu$.

An operator $(\D(A),A)$ on a Banach space $\X$, over the field of complex numbers, is said to be $\omega$-\textit{\textbf{sectorial}} if for a fixed $\omega\in [0,\pi)$ both conditions are satisfied
\begin{enumerate}
    \item $\sigma(A)\subset \overline{\Sigma}_\omega $, where $\sigma(A)$ stands for the spectrum of $A$;
    \item For all $\mu\in(\omega,\pi)$, $\sup_{\lambda\in \mathbb{C}\setminus\overline{\Sigma}_\mu}\lVert \lambda(\lambda \I-A)^{-1}\rVert_{\X\rightarrow \X} < \infty$.
\end{enumerate}

Similarly, $(\D(B),B)$ is said to be $\omega$-\textit{\textbf{bisectorial}} on $\X$ if for a fixed $\omega\in [0,\frac{\pi}{2})$ both conditions are satisfied
\begin{enumerate}
    \item $\sigma(B)\subset \overline{\mathrm{S}}_\omega $;
    \item For all $\mu\in(\omega,\frac{\pi}{2})$, $\sup_{\lambda\in \mathbb{C}\setminus\overline{\mathrm{S}}_\mu}\lVert \lambda(\lambda \I-B)^{-1}\rVert_{\X\rightarrow \X} < \infty$.
\end{enumerate}
In what follows, we mainly focus on the main features and properties of sectorial operators exclusively, because they stem for the main feature of this work (the Stokes operator is sectorial). This, even if some bisectorial operator will be involved at some point in this work (the said Hodge-Dirac operator is a bisectorial operator). However, bisectorial operators enjoy the exact same, or at least very similar, features, and quite extensive exposure dedicated to them can be found in \cite[Chapter~3]{EgertPhDThesis2015}.

One of the biggest feature of sectorial operators is about to be able to define the Dunford Functional calculus, that is, to make sense of the quantity
\begin{align}\label{eq:DunfordintegralFuncCalc}
    f(A)x = \frac{1}{2 i \pi}\int_{\partial \Sigma_\theta} f(z)(z\I-A)^{-1}x\,\mathrm{d}z\text{,}
\end{align}
for all $x\in \X$, provided $\theta \in(\omega,\pi)$ and $f\in\L^1_\ast(\Sigma_{\theta})\cap\mathrm{\mathbf{H}}^\infty(\Sigma_\theta)$. Here, $\mathrm{\mathbf{H}}^\infty(\Sigma_\theta)$ stands for the sets of holomorphic functions over $\Sigma_\theta$ which are uniformly bounded. For the following class of decaying holomorphic functions
\begin{align*}
    \mathbf{H}^\infty_0(\Sigma_\theta):=\Big\{\,\psi\in\mathbf{H}^\infty(\Sigma_\theta)\,:\exists\varepsilon_0>0,\,\forall z \in\Sigma_\theta,\,|\psi(z)|\lesssim_{\theta,\psi,\varepsilon_0}\frac{|z|^{\varepsilon_0}}{(1+|z|)^{2\varepsilon_0}}\,\Big\},
\end{align*}
it is easy to check that \eqref{eq:DunfordintegralFuncCalc} makes sense for any $\omega$-sectorial operator $(\D(A),A)$ on $\X$ and for any $x\in\X$.

If $(\D(A),A)$ is $\omega$-sectorial with $\omega\in[0,\pi)$, for $\mu\in(\omega,\pi)$, one says that $A$
\begin{enumerate}
    \item  has  \textit{\textbf{bounded imaginary powers}} (BIP) of type $\theta_A\geqslant 0$, if for $f(z)=z^{is}$ plugged in \eqref{eq:DunfordintegralFuncCalc}, $x\mapsto A^{is} x$ yields a bounded linear operator for all $s\in\RR$, and
\begin{align*}
    \theta_A := \inf \left\{ \nu\geqslant 0\,: \, \sup_{s\in\RR} e^{-\nu|s|}\lVert A^{is}\rVert_{\X \rightarrow \X}<\infty\right\}\text{. }
\end{align*}
One always has $\theta_A\geqslant \omega$.
    \item admits a \textit{\textbf{bounded}} (or $\mathrm{\mathbf{H}}^\infty(\Sigma_\mu)$-)\textit{\textbf{holomorphic functional calculus}} on $X$ (of angle $\mu$), if for $\theta\in (\omega,\mu)$, there exists a constant $K_\theta$, such that for all $f\in\mathrm{\mathbf{H}}^\infty(\Sigma_\theta)$,
\begin{align*}
    \lVert f(A)\rVert_{\X\rightarrow \X}\leqslant K_\theta \lVert f \rVert_{\L^\infty}\text{.}
\end{align*}
\end{enumerate}

Of course bounded functional calculus of angle $\mu$ for $A$, implies bounded imaginary powers of angle $\theta_A\in[\omega, \mu]$. One of the interesting consequences of the BIP property is the following result
\begin{theorem}[ {\cite[Theorem~6.6.9]{bookHaase2006}} ]\label{thm:BIP}Let $(\D(A),A)$ be a sectorial operator on $\X$ such that it has BIP. Then for all $\theta\in(0,1)$, $\zeta:=(1-\theta)\alpha+\theta\beta$, it holds that
\begin{align*}
    [\D(A^\alpha),\D(A^\beta)]_{\theta} = \D(A^{\zeta})
\end{align*}
for all $0\leqslant \Re \alpha\leqslant \Re \beta$, with either $\Re \alpha>0$ or $\alpha =0$.
\end{theorem}

An important result that will be used in Section~\ref{sec:stokessteady}, is that any invertible sectorial operator of angle less than $\frac{\pi}{2}$ does have bounded $\mathbf{H}^{\infty}$- functional calculus of the same angle (and in particular it has BIP) on the space $(\X,\D(A))_{\theta,q}$, $\theta\in(0,1)$, $q\in[1,\infty]$,  see \cite[Corollary~6.4.3~\&~Theorem~6.6.8]{bookHaase2006}.

\subsubsection{Abstract parabolic maximal regularity and \texorpdfstring{$\mathcal{R}$}{R}-boundedness.}\label{sec:AsbtractMaxReg}

Now, we consider  a closed operator $(\D(A),A)$ on a Banach space $\X$. We recall, see \cite[Theorem~3.7.11]{ArendtBattyHieberNeubranker2011}, that the two following assertions are equivalent:\textit{
\begin{enumerate}
    \item $A$ is $\omega$-sectorial on $\X$, for some $\omega\in [0,\tfrac{\pi}{2})$;
    \item $-A$ generates a bounded holomorphic semigroup on $\X$, denoted by $(e^{-tA})_{t\geqslant0}$.
\end{enumerate}}

Thus, provided that $A$ is $\omega$-sectorial on $\X$ for some $\omega\in [0,\tfrac{\pi}{2})$, for $T\in(0,\infty]$, we look at the following abstract Cauchy problem, 
\begin{align}\tag{ACP}\label{ACP}
    \left\{\begin{array}{rl}
            \partial_t u(t) +Au(t)  =& f(t) \,\text{, } 0<t<T\\
            u(0) =& u_0\text{,}
    \end{array}
    \right.
\end{align}
where $f\in \L^1_{\mathrm{loc}}(0,T;\X)$, $u_0 \in \Y$, with $\Y$ being some normed vector space depending on $\X$ and $\D(A)$.

Provided, $\Y$ is included in $\X$, the mild solution $u\in\C^0([0,T),\X)$ to \eqref{ACP} is unique and given  by the Duhamel formula
\begin{align*}
    u(t)=e^{-tA}u_0 + \int_{0}^{t} e^{-(t-s)A}f(s) \d s\text{, }t\geqslant 0.
\end{align*}

A very central question concerns the truthfulness of the inequality
\begin{align*}
    \lVert (\partial_tu, A u )\rVert_{\L^q(0,T;\X)}\lesssim_{T} \lVert f\rVert_{\L^q(0,T;\X)} + \lVert u_0\rVert_{\Y}
\end{align*}
provided provided $f\in\L^q(0,T;\X)$ and $u_0\in \Y$, for some choice of $\Y$.

The problem being linear, for $u_0$ satisfying
\begin{align}\label{eq:InitialDataLqMaxReg}
    \int_{0}^{\infty} \lVert Ae^{-tA} u_0 \rVert_{\X}^q \d t<\infty
\end{align}
it follows that we can reduce to the case $u_0=0$, and this leads to the following definition:

\begin{definition}We say that an $\omega$-sectorial operator $(\D(A),A)$ on $\X$, where $\omega\in[0,{\sfrac{\pi}{2}})$, satisfies the $\L^q$-\textit{\textbf{maximal regularity property}} over $(0,T)$ on $\X$, if the Banach-valued singular integral operator
\begin{align*}
    A(\partial_t+A)^{-1}\,:\,f\longmapsto \left[t\mapsto A\int_{0}^{t} e^{-(t-s)A}f(s) \d s\right]
\end{align*}
is bounded on $\L^q(0,T;\X)$, where $q\in[1,\infty]$ is fixed.
\end{definition}
Actually, for $q\in(1,\infty)$, it can be shown that the $\L^q$-{{maximal regularity property}} does not depends on $q$: if it holds for one particular $q_0\in(1,\infty)$ then it holds for all $q\in(1,\infty)$. This is however not the case for the endpoint cases $q=1,\infty$. However, the case $q_0=1$ always implies the cases $q\in(1,\infty)$.

A nice sufficient condition in order for $(\D(A),A)$ to satisfy the $\L^q$-\textit{\textbf{maximal regularity property}} is the Dore-Venni Theorem. See \cite[Theorem~9.3.11~\&~Corollary~9.3.12]{bookHaase2006}, originally proved in a weaker version asking for invertibility of $A$ \cite{DoreVenni1987}.

\begin{theorem}[ Dore-Venni ]\label{thm:DoreVennithm} Let $\X$ be a Banach space that belongs to the UMD class. If $(\D(A),A)$ is a sectorial operator that admits BIP of angle $ \theta_A\in[0,{\sfrac{\pi}{2}})$ on $\X$, then it has the $\L^q$-maximal regularity property over $(0,T)$ on $\X$, for all $q\in(1,\infty)$, and all $T\in(0,\infty]$.
\end{theorem}

A {Banach space} $\X$ is said to be of the \textbf{UMD} class, if for one (or equivalently all) $q\in(1,\infty)$, the Hilbert transform
\begin{align*}
    \H f(x) = \frac{1}{\pi}\mathrm{p.v.} \int_\RR f(x-t)\frac{\d t}{t} 
\end{align*}
yields a bounded linear operator on $\L^q(\RR,\X)$. Here, $\mathrm{p.v.}$ means that the corresponding integral is understood in the sense of Cauchy principal value.

The acronym UMD stands for \textbf{Unconditional Martingale Differences}, this denomination being a consequence of the equivalence proved by Bourgain \cite[Lemma~2]{Bourgain83} and Burkholder \cite[Theorem~2]{BurkholderUMD-HT83}, see also \cite[Chapters~4~\&~5, Theorem~5.1.1]{HytonenNeervenVeraarWeisbookVolI2016}. For a more recent exposition about the link between $\L^q$-maximal regularity, $\mathcal{R}$-boundedness and holomorphic functional calculus, see \cite{HytonenNeervenVeraarWeisbookVolIII2023}  see e.g.\cite[Chapters~3~\&~4]{PrussSimonett2016}. 

A definitive characterization of $\L^q$-maximal regularity within the framework of UMD Banach spaces as been obtained by Weis, asserting that $\mathcal{R}$-boundedness of the resolvent of $(\D(A),A)$ on $\X$ is a necessary and sufficient condition. Recall that all Lebesgue, Sobolev and Besov spaces, $\L^p(\Omega)$, $\H^{s,p}(\Omega)$ and $\B^{s}_{p,q}(\Omega)$ as well as any of their closed subspaces, are all UMD spaces for any $1<p,q<\infty$, $s\in\RR$.

\begin{definition}\label{def:rbound} Let $\X$ and $\Y$ be Banach spaces.  We say that a family $\mathcal{T}$ of bounded linear operators from $\X$ to $\Y$ is $\boldsymbol{\mathcal{R}}$\textbf{-bounded}, if there exists a constant $C\geqslant 0$ such that for all $N\in\NN^\ast$, $(T_j)_{j\in\llb1,N\rrb}\subset\mathcal{T}$ and $(x_j)_{j\in\llb1,N\rrb}\subset\X$, the inequality
\begin{align}\label{def:Rboundedness}
    \Bigg\lVert \sum_{j=1}^N r_j(\cdot)T_jx_j\Bigg\rVert_{\L^q(0,1;\Y)} \leqslant C\Bigg\lVert \sum_{j=1}^N r_j(\cdot)x_j\Bigg\rVert_{\L^q(0,1;\X)}
\end{align}
holds for one (or, equivalently, all) $q\in[1,\infty)$. Here, for all $j\in\NN$, $r_j$ are defined by $r_j(t):=\sgn(\sin(2\pi^j t))$, $t\in[0,1]$, called the Rademacher functions. The optimal constant $C\geqslant 0$, for which \eqref{def:Rboundedness} holds for the family of operators $\mathcal{T}$ is denoted by
\begin{align*}
    \mathcal{R}(\mathcal{T}) = \mathcal{R}\{T, T\in\mathcal{T}\}. 
\end{align*}
\end{definition}

\begin{theorem}[ {\cite[Theorem~4.2]{Weis2001}} ]\label{thm:WeisMaxRegRbound}Let $\X$ be an UMD Banach space. Let $(\D(A),A)$ be an unbounded linear operator on $\X$. Then $(\D(A),A)$ has the $\L^q$-maximal regularity property for one (or equivalently all) $q\in(1,\infty)$ if and only if there exists
$\theta\in(0,\sfrac{\pi}{2})$ such that the family
\begin{align*}
    \{ \lambda(\lambda\I-A)^{-1},\, \lambda\in\Sigma_{\pi-\theta}\,\}
\end{align*}
is $\mathcal{R}$-bounded on $\X$. Such an operator $(\D(A),A)$ is called \textbf{$\boldsymbol{\mathcal{R}}$-sectorial} of angle $\theta$, or just $\theta$-{$\mathcal{R}$}-sectorial. If it holds true for any $\theta\in(0,\frac{\pi}{2})$, one says that is $(\D(A),A)$ $0$-{$\mathcal{R}$}-sectorial.
\end{theorem}

Now, with dealing with initial data in \eqref{eq:InitialDataLqMaxReg} as a motivation, we introduce two quantities for $v\in \X + \D(A)$,
\begin{align*}
    \lVert v\rVert_{\mathring{\mathcal{D}}_{A}(\theta,q)}:= \left( \int_{0}^{\infty} \lVert t^{1-\theta}Ae^{-tA} v \rVert_{\X}^q \frac{\mathrm{d}t}{t}\right)^\frac{1}{q}\text{, and }  \lVert v\rVert_{{\mathcal{D}}_{A}(\theta,q)}:= \lVert v \rVert_\X + \lVert v\rVert_{\mathring{\mathcal{D}}_{A}(\theta,q)}\text{,}
\end{align*}
where $\theta \in(0,1)$, $q\in[1,\infty]$, with the special case
\begin{align*}
    \lVert v\rVert_{\mathring{\mathcal{D}}_{A}(\theta,\infty)}:= \sup_{t>0}\,\lVert t^{1-\theta}Ae^{-tA} v \rVert_{\X}.
\end{align*}
This leads to the construction of the vector space
\begin{align*}
    {\mathcal{D}}_{A}(\theta,q):= \{ v\in \X\,:\, \lVert v\rVert_{\mathring{\mathcal{D}}_{A}(\theta,q)}<\infty \} \text{.}
\end{align*}
The vector space ${\mathcal{D}}_{A}(\theta,q)$ is known to be a Banach space under the norm $\lVert \cdot\rVert_{{\mathcal{D}}_{A}(\theta,q)}$ and, moreover, it satisfies the following equality with equivalence of norms
\begin{align}\label{eq:InterpolationXDomainAIntro}
    {\mathcal{D}}_{A}(\theta,q) = (\X,\D(A))_{\theta,q}\text{,}
\end{align}
see \cite[Theorem~6.2.9]{bookHaase2006}. If moreover $0\in \rho(A)$, it has been proved, \cite[Corollary~6.5.5]{bookHaase2006}, that $\lVert \cdot\rVert_{\mathring{\mathcal{D}}_{A}(\theta,q)}$ and $\lVert \cdot\rVert_{{\mathcal{D}}_{A}(\theta,q)}$ are two equivalent norms on ${\mathcal{D}}_{A}(\theta,q)$. So we restrict ourselves to the case of injective but not invertible operators.

\begin{assumption}\label{asmpt:homogeneousdomaindef} The operator $(\D(A),A)$ is injective on $\X$, and there exists a normed vector space $(\Y, \left\lVert \cdot \right\rVert_{\Y})$, such that $\D(A)\subset \Y$, and for all $x\in \D(A)$,
\begin{align}
        \left\lVert Ax\right\rVert_{\X} \sim_{\X,\Y,A} \left\lVert x \right\rVert_{\Y} \text{.}
\end{align}
\end{assumption}

The idea is to construct a homogeneous version of $A$ denoted $\mathring{A}$, defining first its domain
\begin{align*}
    \D(\mathring{A}):= \{\, y\in \Y\, :\, \exists (x_n)_{n\in\mathbb{N}}\subset \D(A),\, \left\lVert y-x_n\right\rVert_{\Y}\xrightarrow[n\rightarrow \infty]{}0   \,\}\text{.}
\end{align*}
So that, for all $y\in \D(\mathring{A})$, $\X$ being complete, it is meaning full to set
\begin{align*}
    \mathring{A}y:= \lim_{n\rightarrow \infty} Ax_n \text{.}
\end{align*}
Constructed this way, the operator $\mathring{A}$ is then injective on $\D(\mathring{A})$. We notice that $\D(\mathring{A})$, endowed with the norm $\lVert\mathring{A}\cdot\rVert_{\X}$, is a normed vector space, but not necessarily complete. We also need the existence of a Hausdorff topological vector space $\Z$, such that $\X,\Y\subset \Z$, and to consider the following assumption. 
\begin{assumption}\label{asmpt:homogeneousdomainintersect} The operator $(\D(A),A)$ and the normed vector space $\Y$ are such that
\begin{align}
        \X\cap\D(\mathring{A}) = \D({A}) \text{.}
\end{align}
\end{assumption}
As a consequence of all above assumptions, we can extend naturally, see \cite[Remark~2.7]{DanchinHieberMuchaTolk2020}, $(e^{-tA})_{t\geqslant0}$ to a $\mathrm{C}_0$-semigroup,
\begin{align*}
        e^{-tA}\,:\, \X+\D(\mathring{A}) \longrightarrow \X+\D(\mathring{A})\text{, } t\geqslant 0 \text{,}
\end{align*}
so that, one can fully make sense of the following vector space,
\begin{align*}
    \mathring{\mathcal{D}}_{A}(\theta,q):= \{ v\in \X+\D(\mathring{A})\,:\, \lVert v\rVert_{\mathring{\mathcal{D}}_{A}(\theta,q)}<\infty \} \text{.}
\end{align*}
Similarly to what happens for ${\mathcal{D}}_{A}(\theta,q)$ in \eqref{eq:InterpolationXDomainAIntro}, it has been proved in \cite[Proposition~2.12]{DanchinHieberMuchaTolk2020}, that the following equality holds with equivalence of norms,
\begin{align}
    \mathring{\mathcal{D}}_{A}(\theta,q) = (\X,\D(\mathring{A}))_{\theta,q}\text{.}
\end{align}
When $\X$ does not have the UMD property or if $A$ does not satisfy $\mathcal{R}$-boundedness of its resolvent on $\X$ being UMD, it turns out that $\mathring{\mathcal{D}}_{A}(\theta,q)$ is a very nice surrogate to $\X$, since one always has $\L^q$-maximal regularity on such vector space, \textbf{even when $q=1,\infty$ !} This result is due to Da Prato and Grisvard \cite[Theorems~4.7~\&~4.15]{DaPratoGrisvard1975}, it was then recently extended by Danchin, Hieber, Mucha and Tolksdorf \cite[Theorem~2.20]{DanchinHieberMuchaTolk2020}, where they allowed global-in-time estimates without requiring invertibility of $A$.

\begin{theorem}[ \textbf{Da Prato--Grisvard - Danchin \& al.} ]\label{thm:DaPratoGrisvard} Let $q\in[1,\infty]$, $\theta\in(0,1)$. Consider $(\D(A),A)$ an $\omega$-sectorial operator on a Banach space $\X$, with $\omega\in[0,\frac{\pi}{2})$, and such that it satisfies Assumptions \ref{asmpt:homogeneousdomaindef} and \ref{asmpt:homogeneousdomainintersect}.  For simplicity we set $\theta_q:=1+\theta-{\sfrac{1}{q}}$

Let $T\in(0,\infty]$. Then, provided $f\in\L^{q}(0,T;{\mathcal{D}}_{A}(\theta,q))$ and $u_0\in \mathring{\mathcal{D}}_{A}(\theta_q,q)$, the problem \eqref{ACP} admits a unique mild solution $u\in \C^{0}_{ub}([0,T]; \mathring{\mathcal{D}}_{A}(\theta_q,q))$, such that $\partial_t u, Au \in \L^{q}(0,T; \mathring{\mathcal{D}}_{A}(\theta,q))$, with the estimates
\begin{align*}
    \lVert u\rVert_{\L^{\infty}(0,T;\mathring{\mathcal{D}}_{A}(\theta_q,q))}\lesssim_{\theta,q,\X,A} \lVert (\partial_t u, Au) \rVert_{\L^q(0,T;\mathring{\mathcal{D}}_{A}(\theta,q))} \lesssim_{\theta,q,\X,A} \lVert f \rVert_{\L^q(0,T;\mathring{\mathcal{D}}_{A}(\theta,q))}+ \lVert u_0\rVert_{\mathring{\mathcal{D}}_{A}(\theta_q,q)}.
\end{align*}
In particular, one can say that $A$ has the $\L^q$-maximal regularity property on $\mathring{\mathcal{D}}_{A}(\theta,q)$. When $q=\infty$, one only has uniform boundedness in-time.

\medbreak

Furthermore, if either $T<\infty$ or $(\D(A),A)$ is invertible the Assumptions \ref{asmpt:homogeneousdomaindef} and \ref{asmpt:homogeneousdomainintersect} can be removed.
\end{theorem}

However, the possible lack of completeness of $\D(\mathring{A})$ implies that $\mathring{\mathcal{D}}_{A}(\theta,q)$ is not necessarily complete\footnote{And having applications in concrete cases in mind --such as solvability of IBVP-- one cannot consider such a completion.}.
This has consequences on how to consider the forcing term $f$ in \eqref{ACP}, choosing $f\in \L^q(0,T;{\mathcal{D}}_{A}(\theta,q))$ instead of $f\in \L^q(0,T;\mathring{\mathcal{D}}_{A}(\theta,q))$ to avoid definition issues, the latter choice being possible when $\mathring{\mathcal{D}}_{A}(\theta,q)$ is a Banach space.

For standard $\L^q$-maximal regularity on UMD Banach space the main result of the theory reads is as follows as a consequence of Theorem~\ref{thm:WeisMaxRegRbound}. We recall that the condition of $\mathcal{R}$-boundedness could be replaced by BIP.

\begin{theorem}\label{thm:LqMaxRegUMDHomogeneous}Let $\omega\in [0,\frac{\pi}{2})$, and $(\D(A),A)$ be an $\omega$-$\mathcal{R}$-sectorial operator on a UMD Banach space $\X$ such that it satisfies the Assumptions \ref{asmpt:homogeneousdomaindef} and \ref{asmpt:homogeneousdomainintersect}. Let $q\in(1,\infty)$. 

 Let $T\in(0,\infty]$. For all $f\in\L^q(0,T;\X)$, all $u_0\in \mathring{\mathcal{D}}_{A}(1-{{\sfrac{1}{q}}},q)$, the problem \eqref{ACP} admits a unique mild solution $u\in \mathrm{C}^0_{ub}([0,T];\mathring{\mathcal{D}}_{A}(1-{\sfrac{1}{q}},q))$ such that $\partial_t u$, $Au \in \L^q(0,T;\X)$ with the estimates
\begin{align}\label{eq:BoundLqMaxReghomogeneous}
    \lVert u \rVert_{\L^\infty(0,T;\mathring{\mathcal{D}}_{A}(1-{\sfrac{1}{q}},q))} \lesssim_{A,q} \lVert (\partial_t u, Au)\rVert_{\L^q(0,T;\X)} \lesssim_{A,q} \lVert f\rVert_{\L^q(0,T;\X)} + \lVert u_0\rVert_{\mathring{\mathcal{D}}_{A}(1-{\sfrac{1}{q}},q)}\text{.}
\end{align}
Furthermore, if either $u_0=0$, or $T<\infty$, or $(\D(A),A)$ is invertible, then the Assumptions \ref{asmpt:homogeneousdomaindef} and \ref{asmpt:homogeneousdomainintersect} can be removed.
\end{theorem}

Note that in \cite[Theorems~4.7~\&~4.15]{DaPratoGrisvard1975}, one can still use $\X$ -- without the UMD assumption -- as ground space for maximal regularity if one uses instead a Besov (Sobolev-Soblodeckij) space in-time, and therefore claiming about $\W^{\alpha,q}(\X)$--maximal regularity whenever $\alpha$ is not an integer.

When $\X$ is UMD, global-in-time (homogeneous) Sobolev (Bessel potential) $\dot{\H}^{\alpha,q}(\X)$--maximal regularity under the assumption of the Dore-Venni Theorem, Theorem~ \ref{thm:DoreVennithm}, has been obtained by the second author, see \cite[Theorem~4.7]{Gaudin2023}. In this case, the assumption of having a UMD space cannot be dropped since it is necessary to have a proper construction of (homogeneous) vector-valued Bessel potential spaces.\footnote{For instance, provided $p\in(1,\infty)$, one has $\W^{1,p}(\X)=\H^{1,p}(\X)$ if and only if $\X$ is UMD.}

\subsubsection{Specific notation for the domains of operators over families of function spaces}

\medbreak

Let $p,q\in[1,\infty]$, $s\in\RR$.

\medbreak

We introduce domains for an operator $A$ acting on (closed subspaces of) Sobolev or Besov spaces, denoting
\begin{itemize}
    \item $\D^{s}_{p}(A)$ (resp. $\dot{\D}^{s}_{p}(A)$)  its domain on $\H^{s,p}$ (resp. $\dot{\H}^{s,p}$), $1<p<\infty$;
    \item $\D^{s}_{p,q}(A)$ (resp. $\dot{\D}^{s}_{p,q}(A)$)  its domain on $\B^{s}_{p,q}$ (resp. $\dot{\B}^{s}_{p,q}$);
    \item $\D^{k}_{p}(A)$ (resp. $\dot{\D}^{k}_{p}(A)$)  its domain on $\W^{k,p}$ (resp. $\dot{\W}^{k,p}$);
    \item $\D_{p}(A) = \D^{0}_{p}(A) = \dot{\D}^{0}_{p}(A)$ its domain on $\L^p$.
\end{itemize}
Since $\B^{s}_{p,p}=\W^{s,p}$, $s\notin\ZZ$, $p\in[1,\infty]$, and $\H^{k,p}=\W^{k,p}$, when $k\in\ZZ$, $1<p<\infty$, the case of $\W^{s,p}$, $s$ is already covered by other cases, and there is no ambiguity in the definitions of notations.

Similarly, $\N^{s}_{p}(A)$, $\N^{s}_{p,q}(A)$ will stand for its nullspace on $\H^{s,p}$ and $\B^{s}_{p,q}$, and range spaces will be given respectively by $\R^{s}_{p}(A)$ and $\R^{s}_{p,q}(A)$. We replace $\N$ and $\R$ by $\dot{\N}$ and $\dot{\R}$ for their corresponding sets on homogeneous function spaces. When it is relevant, its resolvent set on $\X^{s,p}$ will be denoted by $\rho_{\X^{s,p}}(A)\subset \CC$.

If the operator $A$ has different realizations depending on various function spaces and on the considered open set, we may write its domain $\D(A,\Omega)$, and similarly for its nullspace $\N$ and range space $\R$. We omit the open set $\Omega$ if there is no possible confusion.

\medbreak

For an abstract space $\X^{s,p}\in\{\,\H^{s,p},\,\B^{s}_{p,q},\,\W^{s,p}, q\in[1,\infty]\}$, $k:=s$, if $s\in\NN$, we denote in a unified way the domain of $A$ on $\X^{s,p}$, by
\begin{align*}
    \D^{s,p}(A):=\begin{cases} \D^{s}_{p}(A),\,\, \, \, \qquad, \textrm{if} \, \X^{s,p}=\H^{s,p};\\
    \D^{s}_{p,q}(A),\,\qquad, \textrm{if}\, \X^{s,p}=\B^{s}_{p,q};\\
    \D^{k}_{p}(A),\,\, \, \, \qquad, \textrm{if} \, \X^{s,p}=\W^{k,p}.
    \end{cases}
\end{align*}
We proceed similarly with the notations $\R^{s,p}$ and $\N^{s,p}$. When one wants to notify that $\X$ is a homogeneous Sobolev or Besov space, one writes a similar definition for $\dot{\D}^{s,p}(A)$ to be the domain of $A$ on the homogeneous function space $\dot{\X}^{s,p}\in\{\,\dot{\H}^{s,p},\,\dot{\B}^{s}_{p,q},\,\dot{\W}^{s,p}, q\in[1,\infty]\}$.

We emphasize that one should \textbf{really distinguish} the domain of an operator on a homogeneous function space $\dot{\D}(A)$, the homogeneous domain of an operator $\D(\mathring{A})$, and the homogeneous domain of an operator on a homogeneous function space $\dot{\D}(\mathring{A})$. Here is an example.

\medbreak

We set $A=-\Delta$, on $\X={\H}^{1,p}(\RR^n)$ and  $\dot\X=\dot{\H}^{1,p}(\RR^n)$, $1<p<\infty$, then one has
\begin{itemize}
    \item ${\D}_{p}^1(A) =\{\,u\in{\H}^{1,p}(\RR^n)\,:\,\Delta u\in {\H}^{1,p}(\RR^n)\,\} = {\H}^{3,p}(\RR^n)$;
    \item $\dot{\D}_{p}^1(A) =\{\,u\in\dot{\H}^{1,p}(\RR^n)\,:\,\Delta u\in \dot{\H}^{1,p}(\RR^n)\,\} = \dot{\H}^{1,p}\cap\dot{\H}^{3,p}(\RR^n)$;
    \item ${\D}_{p}^1(\mathring{A})  = \overline{\H^{3,p}(\RR^n)}^{\lVert \Delta \cdot\rVert_{\H^{1,p}(\RR^n)}} = \overline{\H^{3,p}(\RR^n)}^{\lVert  \cdot\rVert_{\dot\H^{2,p}(\RR^n)}+ \lVert \cdot\rVert_{\dot\H^{3,p}(\RR^n)}} = \dot{\H}^{2,p}\cap\dot{\H}^{3,p}(\RR^n)$;
    \item $\dot{\D}_{p}^1(\mathring{A})  = \overline{\dot{\H}^{1,p}\cap\dot{\H}^{3,p}(\RR^n)}^{\lVert \Delta \cdot\rVert_{\dot{\H}^{1,p}(\RR^n)}} = \overline{\dot{\H}^{1,p}\cap\dot{\H}^{3,p}(\RR^n)}^{\lVert \cdot\rVert_{\dot{\H}^{3,p}(\RR^n)}}= \dot{\H}^{3,p}(\RR^n)$.
\end{itemize}

\subsection{Sobolev and Besov multipliers  and further rough domains}\label{sec:SM}

\subsubsection{The Sobolev and Besov multipliers}
 
In accordance with \cite[Chapter 14]{MazyaShaposhnikova2009} the multiplier norm for Sobolev (or Besov) spaces $\X^{s,p}\in\{\W^{s,p},\H^{s,p}\}$\footnote{when $p=\infty$ one can also consider $\mathcal{B}^{s}_{\infty,\infty}$, $\mathcal{B}^{s,0}_{\infty,\infty}$, $\B^{s,0}_{\infty,\infty}$ as well as their homogeneous counterparts.}  is given for $\varphi\in\W^{1,1}_{\text{loc}}(\Omega)$ by
\begin{align}\label{eq:SoMo}
\|\varphi\|_{\mathcal{M}^{s,p}_{\X}(\Omega)}:=\sup_{\substack{\vv\in\X^{s-1,p}({\Omega},\CC^n)\\\|\vv\|_{\X^{s-1,p}(\Omega)}=1}}\|\nabla\varphi\cdot\mathbf{v}\|_{\X^{s-1,p}(\Omega)},
\end{align}
where $p\in[1,\infty]$ and $s\in\RR$ (we assume $1<p<\infty$ if $\X=\H$).
The space $\mathcal{M}^{s,p}_{\X}(\Omega)$ of Sobolev-Besov multipliers is defined as those objects for which the $\mathcal{M}^{s,p}_{\X}(\Omega)$-norm is finite.

In the case of (non-diagonal) Besov spaces, we write $\mathcal{M}^{s,p,q}_{\B}(\Omega)$, for $\X^{s-1,p}= \B^{s-1}_{p,q}$, provided $p,q\in[1,\infty]$.

For $\varepsilon>0$, we denote by $\mathcal{M}^{s,p}_{\X}(\Omega,\varepsilon)$ the ball of radius $\varepsilon$ with respect to the $\mathcal{M}^{s,p}_{\X}(\Omega)$-norm.

We also introduce the \textbf{original} Sobolev-Besov multiplier norm for $\varphi\in\L^1_{\text{loc}}(\Omega)$:
\begin{align}\label{eq:SoMo'}
\|\varphi\|_{\mathcal{M}^{s,p}_{\X, \tt{or}}(\Omega)}:=\sup_{\substack{w\in\X^{s,p}(\Omega)\\\|{w}\|_{\X^{s,p}(\Omega)}=1}}\|\varphi \,w\|_{\X^{s,p}(\Omega)} <\infty,
\end{align}
up to appropriate modification in the case of non-diagonal Besov spaces.

Note that for all $s\in\RR$, all $p,q\in[1,\infty]$, one always has
\begin{align*}
    \mathcal{M}^{s,p}_{\X, \tt{or}}(\RR^n)\hookrightarrow\L^\infty(\RR^n).
\end{align*}

The quantity \eqref{eq:SoMo'} also serves as a customary definition of the Sobolev-Besov multiplier norm in the literature, but \eqref{eq:SoMo} is more suitable for our purposes, noting that
\begin{align*}
    \|\varphi\|_{\mathcal{M}^{s,p}_{\X}(\Omega)}=\| \nabla \varphi\|_{\mathcal{M}^{s-1,p}_{\X, \tt{or}}(\Omega)}.
\end{align*}

We start with the following proposition.
\begin{proposition}\label{prop:basicMultipliers}Let $s\in\RR$, $p\in[1,\infty]$. It holds that
\begin{enumerate}
    \item One has the continuous inclusion ${\mathcal{M}}^{{s},p}_{\X,\tt{or}}(\Omega)\subset{\mathcal{M}}^{-s,p'}_{\X',\tt{or}}(\Omega)$;
    \item If $s\geqslant -1$, and $\X$ is an inhomogeneous function space, for all $\phi\in\C^{0,1}_{b}(\overline{\Omega})$, such that $\phi\in{\mathcal{M}}^{s,p}_{\X}(\Omega)$, one has $\phi\in{\mathcal{M}}^{s,p}_{\X,\tt{or}}(\Omega)$.\label{pt:implicationLipMultipliers}
\end{enumerate}
\end{proposition}

\begin{proof} Point  \textit{(i)} is clear. We only check point \textit{(ii)}. The case $p=\infty$ is straightforward. Let us first assume that $s\in[0,1]$, $p\in[1,\infty]$, $u\mapsto u\phi$ is bounded on $\L^p(\Omega)$ and on $\W^{1,p}(\Omega)$, so by real and complex interpolation, it is bounded on $\W^{s,p}(\Omega)$ and on $\H^{s,p}(\Omega)$ for all $s\in[0,1]$, assuming $p>1$ for the latter.

Now, for $s\geqslant 1$, we proceed by induction. Indeed, by the Leibniz rule
\begin{align*}
    \lVert u\phi \rVert_{\X^{s,p}(\Omega)}&\lesssim_{s,p,n,\Omega} \lVert u\phi\rVert_{\X^{s-1,p}(\Omega)} + \lVert \phi (\nabla u) \rVert_{\X^{s-1,p}(\Omega)} + \lVert u (\nabla\phi)\rVert_{\X^{s-1,p}(\Omega)}  \\
    &\lesssim_{p,s,n,\Omega} \lVert \phi\rVert_{\mathcal{M}^{s-1,p}_{\X,\tt{or}}(\Omega)} \lVert u\rVert_{\X^{s,p}(\Omega)} + \lVert \phi\rVert_{\mathcal{M}^{s,p}_{\X}(\Omega)}\lVert u\rVert_{\X^{s-1,p}(\Omega)}\\
    &\lesssim_{p,s,n,\Omega} (\lVert \phi\rVert_{\mathcal{M}^{s-1,p}_{\X,\tt{or}}(\Omega)} + \lVert \phi\rVert_{\mathcal{M}^{s,p}_{\X}(\Omega)})\lVert u\rVert_{\X^{s,p}(\Omega)}.
\end{align*}
For the case $-1\leqslant s<0$, one can check that $\mathcal{M}^{-s,p'}_{\X,\tt{or}}(\Omega)$ also preserves $\X_{0}^{-s,p'}(\Omega)$, when $-s\neq {\sfrac{1}{p'}}$ (it preserves homogeneous boundary conditions). Therefore, by duality from the case $-s\in[0,1]$ and interpolation, the result remains valid for $s\geqslant -1$.
\end{proof}

Note that in our applications, we always assume that the functions in question are uniformly Lipschitz continuous, such that the assertion \ref{pt:implicationLipMultipliers} from Proposition~\ref{prop:basicMultipliers} above is always satisfied. 

\medbreak

Let us collect some useful properties of Sobolev multipliers.

We start with the following rule about the composition with Sobolev multipliers in a way similar to \cite[Lemma 9.4.1]{MazyaShaposhnikova2009}.
\begin{lemma}\label{lem:multiplierComposition}Let $p,q\in[1,\infty]$, $s\in\RR$. Let $\Omega$, $\Omega'$ be two Lipschitz domains of $\RR^n$ and $\phi\,:\,\Omega\longrightarrow\Omega'$ be a globally bi-Lipschitz map. Then it holds that the map
\begin{align*}
    u\longmapsto u\circ\phi
\end{align*}
is well-defined and bounded
\begin{enumerate}
    \item from $\W^{s,p}(\Omega')$ to $\W^{s,p}(\Omega)$, provided $s\in[-1,1]$, and similarly for $\H^{s,p}$, if additionally $1<p<\infty$;
    \item from $\B^{s}_{p,q}(\Omega')$ to $\B^{s}_{p,q}(\Omega)$, provided $s\in(-1,1)$.
\end{enumerate}
Furthermore,
\begin{itemize}
    \item  if $s\geqslant 1$ and additionally $\phi\in\mathcal{M}^{s,p}_{\X}(\Omega)$ then composition by $\phi$ maps $\X^{s,p}(\Omega')$ to $\X^{s,p}(\Omega)$, including the case of Besov spaces.
\end{itemize}
\end{lemma}

\begin{proof} If $s\in[-1,1]$, this is straightforward from the case $s=0,1$. Indeed, one can check that whenever $p\in[1,\infty]$, 
\begin{align*}
    \lVert u\circ \phi\rVert_{\L^p(\Omega)} &\leqslant \lVert \det(\nabla (\phi^{-1}))\rVert_{\L^\infty(\Omega')}^{\sfrac{1}{p}} \lVert u\rVert_{\L^p(\Omega')}.\\
    \lVert \nabla (u\circ \phi)\rVert_{\L^p(\Omega)} &\leqslant \lVert \det(\nabla (\phi^{-1}))\rVert_{\L^\infty(\Omega')}^{\sfrac{1}{p}}\lVert \nabla \phi\rVert_{\L^\infty(\Omega)}  \lVert \nabla u\rVert_{\L^p(\Omega')}.
\end{align*}
So that, by complex interpolation, and recalling that $\lVert \det(\nabla (\phi^{-1}))\rVert_{\L^\infty(\Omega')}= \lVert \det((\nabla \phi)^{-1})\rVert_{\L^\infty(\Omega)}$, provided $1<p<\infty$, $s\in[0,1]$, this yields the bound
\begin{align*}
   \lVert u\circ \phi\rVert_{\H^{s,p}(\Omega)} &\lesssim_{p,s,n,\Omega}   \lVert \det((\nabla \phi)^{-1})\rVert_{\L^\infty(\Omega)}^{\sfrac{1}{p}}\left[ 1+\lVert \nabla \phi\rVert_{\L^\infty(\Omega)}\right]^{s}  \lVert u\rVert_{\H^{s,p}(\Omega')}.
\end{align*}
It also preserves $\H^{s,p}_0$, ${\sfrac{1}{p}}<s\leqslant 1$. Therefore, by duality and complex interpolation for $s\in[-1,1]$
\begin{align*}
   \lVert u\circ \phi\rVert_{\H^{s,p}(\Omega)} &\lesssim_{p,s,n,\Omega}  \lVert \det((\nabla \phi)^{-1})\rVert_{\L^\infty(\Omega)}^{\sfrac{1}{p}}\left[ 1+\lVert \nabla \phi\rVert_{\L^\infty(\Omega)}\right]^{|s|}  \lVert u\rVert_{\H^{s,p}(\Omega')}.
\end{align*}
By real interpolation, one obtains for $0<|s|<1$, $p,q\in[1,\infty]$,
\begin{align*}
   \lVert u\circ \phi\rVert_{\B^{s}_{p,q}(\Omega)} &\lesssim_{p,s,n,\Omega}  \lVert \det((\nabla \phi)^{-1})\rVert_{\L^\infty(\Omega)}^{\sfrac{1}{p}}\left[ 1+\lVert \nabla \phi\rVert_{\L^\infty(\Omega)}\right]^{|s|}  \lVert u\rVert_{\B^{s}_{p,q}(\Omega')} .
\end{align*}
If $s=0$, again for some $0<\varepsilon<1$,  by interpolation between the cases $s=\pm\varepsilon$,
\begin{align*}
   \lVert u\circ \phi\rVert_{\B^{0}_{p,q}(\Omega)} &\lesssim_{p,\varepsilon,n,\Omega} \lVert \det((\nabla \phi)^{-1})\rVert_{\L^\infty(\Omega)}^{\sfrac{1}{p}}\left[ 1+\lVert \nabla \phi\rVert_{\L^\infty(\Omega)}\right]^{\varepsilon}  \lVert u\rVert_{\B^{0}_{p,q}(\Omega')} .
\end{align*}
The case $s\geqslant 1$ follows by induction. We just give the bound for $s\in[1,2)$ in the case of Besov spaces, $s\in(1,2]$ in the case of Bessel potential spaces or standard Sobolev spaces of integer order. It holds
\begin{align*}
   &\lVert u\circ \phi\rVert_{\X^{s,p}(\Omega)} \\&\lesssim_{p,s,n,\Omega} \lVert u\circ \phi\rVert_{\X^{s-1,p}(\Omega)} + \lVert \nabla ( u\circ \phi )\rVert_{\X^{s-1,p}(\Omega)}\\
    &\lesssim_{p,s,n,\Omega}  \lVert \det((\nabla \phi)^{-1})\rVert_{\L^\infty(\Omega)}^{\sfrac{1}{p}}\left[ 1+\lVert \phi\rVert_{\mathcal{M}^{s,p}_{\X}(\Omega)}\right]\left[ 1+\lVert \nabla \phi\rVert_{\L^\infty(\Omega)}\right]^{\max(\varepsilon,s-1)}   \lVert u\rVert_{\X^{s,p}(\Omega')}.
\end{align*}
This ends the proof. Note that the dependency of the implicit constant with respect to $\Omega$ (not $\Omega'$) only depends on the interpolation constants, which relies itself on the choice of the extension operator from $\Omega$ to $\RR^n$.
\end{proof}

We also have the following natural continuous inclusions of Sobolev multipliers:
\begin{proposition}\label{prop:embeddingmultipliers}Let $p,q\in[1,\infty]$, $p<\infty$, $0<s_0<s<s_1$. Then it holds for any combination $\X,\Y,\Z\in\{\H,\B_{\cdot,q},\W\}$, that
\begin{align*}
    \mathcal{M}_{\X,\tt{or}}^{s_1,p}(\RR^n)\hookrightarrow\mathcal{M}_{\Y,\tt{or}}^{s,p}(\RR^n)\hookrightarrow\mathcal{M}_{\Z,\tt{or}}^{s_0,p}(\RR^n),
\end{align*}
assuming additionally $p>1$ whenever Bessel potential spaces are involved.
\end{proposition}

\begin{proof} Assume $1<p<\infty$, first. Let $\varphi\in\mathcal{M}^{s,p}_{\X,\tt{or}}(\RR^n)$, then by \cite[Lemma~2.14]{Sickel1999}, one has for $\varepsilon>0$,
\begin{align*}
    \lVert \varphi\rVert_{\mathcal{M}^{s-\varepsilon,p}_{\H,\tt{or}}(\RR^n)} &\leqslant \lVert \varphi\rVert_{\mathcal{M}^{s,p}_{\X,\tt{or}}(\RR^n)}, \\
    \text{ and }\quad\lVert \varphi\rVert_{\mathcal{M}^{s-\varepsilon,p,p}_{\B,\tt{or}}(\RR^n)} &\leqslant \lVert \varphi\rVert_{\mathcal{M}^{s,p}_{\X,\tt{or}}(\RR^n)}, 
\end{align*}
so that by real interpolation, for all $0<s-\varepsilon<s_0<s$, $q\in[1,\infty]$,
\begin{align*}
    \lVert u\varphi\rVert_{\B^{s_0}_{p,q}(\RR^n)} \lesssim_{p,s,n,\varepsilon} \lVert \varphi\rVert_{\mathcal{M}^{s,p}_{\X,\tt{or}}(\RR^n)}\lVert u\rVert_{\B^{s_0}_{p,q}(\RR^n)},
\end{align*}
and we did obtain
\begin{align*}
    \mathcal{M}^{s,p}_{\H,\tt{or}}+\mathcal{M}^{s,p,p}_{\B,\tt{or}}(\RR^n)\hookrightarrow\mathcal{M}^{s_0,p,q}_{\B,\tt{or}}(\RR^n).
\end{align*}
By real interpolation, for all $\varphi\in\mathcal{M}^{s-\varepsilon,p}_{\H,\tt{or}}\cap\mathcal{M}^{s_1,p,q}_{\B,\tt{or}}(\RR^n)$, one has the estimate of operator norms
\begin{align*}
    \lVert \varphi\rVert_{\mathcal{M}^{s,p}_{\X,\tt{or}}(\RR^n)}\lesssim_{p,s,s_1,n,\varepsilon}\lVert \varphi\rVert_{\mathcal{M}^{s-\varepsilon,p}_{\X,\tt{or}}(\RR^n)}^{1-\theta}\lVert \varphi\rVert_{\mathcal{M}^{s_1,p,q}_{\B,\tt{or}}(\RR^n)}^{\theta}.
\end{align*}
We set $\varphi_N:=S_N\Psi$, for all $N\in\NN$ where $\Psi\in \mathcal{M}^{s_1,p,q}_{\B,\tt{or}}(\RR^n)\setminus\{0\}$. Since  $\mathcal{M}^{s_1,p,q}_{\B,\tt{or}}(\RR^n)\hookrightarrow\L^\infty(\RR^n)$, it turns out that $\varphi_N\in\C^{\infty}_{ub}(\RR^n)$ for all $N\in\NN$, so that $\varphi_N\in\mathcal{M}^{s,p}_{\H,\tt{or}}(\RR^n)\hookrightarrow\mathcal{M}^{s-\varepsilon,p}_{\X,\tt{or}}(\RR^n)$, and one obtains
\begin{align*}
    \lVert \varphi_N\rVert_{\mathcal{M}^{s,p}_{\X,\tt{or}}(\RR^n)}^\theta\lesssim_{p,s,s_1,n,\varepsilon}\lVert \varphi_N\rVert_{\mathcal{M}^{s_1,p,q}_{\B,\tt{or}}(\RR^n)}^{\theta}.
\end{align*}
Hence, by \cite[Lemma~2.1,~\textit{(ii)}]{LiSickelYangYuan2024}, it holds
\begin{align*}
    \lVert \varphi_N\rVert_{\mathcal{M}^{s,p}_{\X,\tt{or}}}\lesssim_{p,s,s_1,n,\varepsilon}\lVert \Psi\rVert_{\mathcal{M}^{s_1,p,q}_{\B,\tt{or}}(\RR^n)}
\end{align*}
Consequently, for all $u\in\H^{s,p}(\RR^n)$, since $(\varphi_N)_{N\in\NN}\subset\mathcal{M}^{s_1,p,q}_{\B,\tt{or}}(\RR^n)\hookrightarrow\L^\infty(\RR^n)$ converges in $\S'(\RR^n)$ to $\Psi$, it holds
\begin{align*}
    \lVert u\Psi\rVert_{\X^{s,p}(\RR^n)} \leqslant \liminf_{N\longrightarrow \infty} \, \lVert u\varphi_N\rVert_{\X^{s,p}(\RR^n)} \lesssim_{p,s,s_1,n} \lVert \Psi\rVert_{\mathcal{M}^{s_1,p,q}_{\B,\tt{or}}(\RR^n)} \lVert u\rVert_{\X^{s,p}(\RR^n)}.
\end{align*}
Finally, we did obtain
\begin{align*}
    \lVert \Psi\rVert_{\mathcal{M}^{s,p}_{\X,\tt{or}}(\RR^n)} \lesssim_{p,s,s_1,n}\lVert \Psi\rVert_{\mathcal{M}^{s_1,p,q}_{\B,\tt{or}}(\RR^n)}.
\end{align*}
The same arguments are valid in the case $p=1$, for the spaces $\W^{s,1}$ and $\B^{s}_{1,q}$.
\end{proof}

\medbreak

When it comes to concrete function spaces, in order to illustrate, we provide several embeddings.

\begin{proposition}\label{prop:EmbeddingsHomogeneousSobolevBesovinMultiplierSpaces} Let $p,q\in[1,\infty]$, $s > 0$ and $m\in\NN^\ast$.
\begin{enumerate}
    \item Provided $1<p<\infty$, $\dot{\H}^{s,\frac{m}{s}}\cap \L^\infty(\RR^m)\hookrightarrow\mathcal{M}^{s,p}_{\dot{\H},\tt{or}}(\RR^m)$ whenever $0<s< {\sfrac{m}{p}}$.
    \item $\dot{\B}^{s}_{\frac{m}{s},1}(\RR^m)+\dot{\B}^{s}_{\frac{m}{s},q}\cap \L^\infty(\RR^m) \hookrightarrow\mathcal{M}^{s,p,q}_{\dot{\B},\tt{or}}(\RR^m)$, $0< s < {\sfrac{m}{p}}$.
\end{enumerate}
\end{proposition}

\begin{proof}For point \textit{(i)} with homogeneous function spaces, this is a direct consequence from \cite[Theorem~1.2~(2)]{Li-KatoPonceEst2019} with $s=s_1\in(0,\sfrac{m}{p})$, $s_2=0$, and Sobolev embeddings ${\L}^{\frac{m}{s}}(\RR^{m})\hookrightarrow \dot{\H}^{-k,\frac{m}{s-k}}(\RR^{m})$, $k<s$, yielding the estimate
\begin{align*}
    \lVert\varphi u\rVert_{\dot{\H}^{s,p}(\RR^{m})} \lesssim_{p,s,m}\left( \lVert\varphi \rVert_{\L^\infty(\RR^{m})} + \lVert\varphi \rVert_{\dot{\H}^{s,\frac{m}{s}}(\RR^{m})}\right) \lVert u\rVert_{\dot{\H}^{s,p}(\RR^{m})}.
\end{align*}

One may use the embedding $\dot{\H}^{\frac{m}{q},q}\hookrightarrow \mathrm{BMO}$, $1<q<\infty$.

For point \textit{(ii)} with homogeneous function spaces, we are going to use the scale invariance property of the involved norms. By \cite[Section~4.9.2,~Corollary~1]{RunstSickel96}, one has for all $\varphi\in \B^{s}_{\frac{m}{s},q}\cap\L^{\infty}(\RR^{m})$, all $u\in\B^{s}_{p,q}(\RR^m)$,
\begin{align*}
    \lVert \varphi u \rVert_{\dot{\B}^{s}_{p,q}(\RR^m)} \lesssim_{p,s,m,q}\lVert \varphi u \rVert_{\B^{s}_{p,q}(\RR^m)}\lesssim_{p,s,m,q} \left(\lVert \varphi \rVert_{\L^\infty(\RR^m)}+\lVert \varphi\rVert_{\B^{s}_{\frac{m}{s},q}(\RR^m)}\right)\lVert u \rVert_{\B^{s}_{p,q}(\RR^m)}.
\end{align*}
Since $s>0$, for $\lambda>0$, plugging $u_\lambda:=u(\lambda\cdot)$ and $\varphi_\lambda:=\varphi(\lambda\cdot)$, by the equivalence of norm $\B^{s}_{r,q}=\L^r\cap\dot{\B}^{s}_{r,q}$, for $r\in\{p,\frac{m}{s}\}$, we deduce
\begin{align*}
    \lambda^{s-\frac{m}{p}}\lVert \varphi u \rVert_{\dot{\B}^{s}_{p,q}(\RR^m)}&\lesssim_{p,s,n,q} \left(\lVert \varphi  \rVert_{\L^\infty(\RR^m)}+ \lambda^{-\frac{m}{q}}\lVert \varphi  \rVert_{\L^q(\RR^m)} + \lVert \varphi\rVert_{\dot{\B}^{s}_{\frac{m}{s},q}(\RR^m)}\right)\\&\qquad\qquad\qquad\qquad\qquad\times \left( \lambda^{-\frac{m}{p}}\lVert u \rVert_{\L^p(\RR^m)} + \lambda^{s-\frac{m}{p}}\lVert u \rVert_{\dot{\B}^{s}_{p,q}(\RR^m)} \right).
\end{align*}
Dividing by $\lambda^{s-\frac{m}{p}}$, then taking the limit as $\lambda$ goes to infinity yields
\begin{align*}
    \lVert \varphi u \rVert_{\dot{\B}^{s}_{p,q}(\RR^m)}\lesssim_{p,s,n,q} \left(\lVert \varphi  \rVert_{\L^\infty(\RR^m)}+ \lVert \varphi\rVert_{\dot{\B}^{s}_{\frac{m}{s},q}(\RR^m)}\right)\lVert u \rVert_{\dot{\B}^{s}_{p,q}(\RR^m)}.
\end{align*}
Now, we relax the estimate to all $\varphi \in\dot{\B}^{s}_{\frac{m}{s},q}\cap\L^{\infty}(\RR^{m})$ and all $u\in\dot{\B}^{s}_{p,q}(\RR^m)$. Note that, since one has $\varphi\in\L^{\infty}(\RR^{m})$, and $u\in\dot{\B}^{s}_{p,q}(\RR^m)\hookrightarrow\L^{\frac{mp}{m-ps},q}(\RR^m)$, both belonging to $\S'_h(\RR^m)$, the following limit holds in $\S'(\RR^m)$
\begin{align*}
    u_N \varphi_N = \left(\sum_{|j|\leqslant N} \dot{\Delta}_j u\right) \left(\sum_{|j|\leqslant N} \dot{\Delta}_j \varphi\right) \xrightarrow[N\rightarrow\infty]{} u\varphi.
\end{align*}
Therefore,
\begin{align*}
    \lVert \varphi u \rVert_{\dot{\B}^{s}_{p,q}(\RR^m)} &\leqslant \liminf_{N\rightarrow\infty} \,\lVert \varphi_N\,  u_N \rVert_{\dot{\B}^{s}_{p,q}(\RR^m)}\\
    &\lesssim_{p,s,n,q} \liminf_{N\rightarrow\infty} \,\left(\lVert \varphi_N  \rVert_{\L^\infty(\RR^m)}+ \lVert \varphi_N\rVert_{\dot{\B}^{s}_{\frac{m}{s},q}(\RR^m)}\right)\lVert u_N \rVert_{\dot{\B}^{s}_{p,q}(\RR^m)}\\
    &\lesssim_{p,s,n,q} \left(\lVert \varphi  \rVert_{\L^\infty(\RR^m)}+ \lVert \varphi\rVert_{\dot{\B}^{s}_{\frac{m}{s},q}(\RR^m)}\right)\lVert u \rVert_{\dot{\B}^{s}_{p,q}(\RR^m)}.
\end{align*}
Finally, the embeddings $\dot{\B}^{s}_{\frac{m}{s},1}(\RR^m)\hookrightarrow \L^{\infty}(\RR^m)$, $\dot{\B}^{s}_{\frac{m}{s},q}(\RR^m)$ yields the result.
\end{proof}

\begin{proposition}\label{prop:EmbeddingsSobolevBesovinMultiplierSpaces} Let $p,q\in[1,\infty]$, $s > 0$ and $m\in\NN^\ast$.
\begin{enumerate}
    \item Provided $1<p<\infty$, $\mathcal{M}^{s,p}_{\H,\tt{or}}(\RR^m)$ contains continuously
    \begin{enumerate}
        \item $\dot{\H}^{s,\frac{m}{s}}\cap \L^\infty(\RR^m)+{\H}^{s,r}\cap \L^\infty(\RR^m)$, whenever ${\sfrac{m}{r}}\leqslant s< {\sfrac{m}{p}}$;
        \item $\H^{s,p}(\RR^m)$, whenever $s>{\sfrac{m}{p}}$;
    \end{enumerate}
    \item Provided $1<p\leqslant\infty$, $\mathcal{M}^{s,p}_{\W,\tt{or}}(\RR^m)$ contains continuously
    \begin{enumerate}\item $\dot{\B}^{s}_{\frac{m}{s},p}\cap \L^\infty(\RR^m)$, whenever $s<\frac{m}{p}$;
        \item ${\B}^{s}_{r,p}\cap \L^\infty(\RR^m)$, provided $s\notin\NN$, whenever $ {\sfrac{m}{r}}\leqslant s < {\sfrac{m}{p}}$; 
        
        or  $r\in(p,\infty]$ and $s={\sfrac{m}{p}}$;
        \item $\W^{s,p}(\RR^m)$ whenever $s>{\sfrac{m}{p}}$;
    \end{enumerate}
    \item $\mathcal{M}^{s,p,q}_{\B,\tt{or}}(\RR^m)$ contains continuously
    \begin{enumerate}
        \item $\dot{\B}^{s}_{\frac{m}{s},q}\cap \L^\infty(\RR^m)$, whenever $s<\frac{m}{p}$;
        \item ${\B}^{s}_{r,q}\cap \L^\infty(\RR^m)$,  whenever $ {\sfrac{m}{r}}\leqslant s < {\sfrac{m}{p}}$; 
        
        \item ${\B}^{s}_{r,1}(\RR^m)$, if $r\in[{\frac{m}{s}},\infty]$, $s=\frac{m}{p}$;
        \item $\B^{s}_{p,q}(\RR^m)$ whenever $s>{\sfrac{m}{p}}$;
    \end{enumerate}
    \item If $p=1$, $\mathcal{M}^{s,1}_{\W,\tt{or}}(\RR^m)$ contains continuously
    \begin{enumerate}
        \item $\dot{\B}^{s}_{\frac{m}{s},1}(\RR^m)+{\B}^{s}_{r,1}(\RR^m)$, provided ${\sfrac{m}{r}}\leqslant s<m$;
        \item $\W^{s,1}(\RR^m)$ whenever $s\geqslant m$;
    \end{enumerate}
\end{enumerate}
\end{proposition}

\begin{remark}\label{rem:MultipliersonHalfSpaces}By the definition of function spaces by restriction and applying Stein's extension operator, see \cite[Sections~3~\&~4]{Gaudin2023Lip} and \cite[Chapter~9,~Section~9.2]{MazyaShaposhnikova2009}, the result remains valid if one replaces $\RR^m$ by a special Lipschitz domain $\Omega$. Furthermore, when $s>m/p$, $q\leqslant p$, so when the function spaces are algebras, one has
\begin{align*}
    \mathcal{M}^{s,p}_{\X,\tt{or}}(\RR^m)=\X_{\textrm{unif}}^{s,p}(\RR^m),
\end{align*}
where $\X_{\textrm{unif}}$ is defined in \eqref{eq:introUnifSpaces}. For more details, we refer to \cite[Theorem~2.5]{Sickel1999} and \cite[Theorem~1.2]{NguyenSickel2018}.
\end{remark}

\begin{proof} For the point \textit{(i)-(a)}, \textit{(ii)-(a)} and \textit{(iii)-(a)}, when homogeneous function spaces are involved, this follows from Proposition~\ref{prop:EmbeddingsHomogeneousSobolevBesovinMultiplierSpaces}. See \cite[Theorem~3.3.1]{MazyaShaposhnikova2009} and \cite[Theorem~4.4.4]{MazyaShaposhnikova2009} for respectively the second part of \textit{(i)-(a)} and the point \textit{(ii)-(b)}.

The points \textit{(i)-(b)}, \textit{(ii)-(c)}, \textit{(iii)-(d)} and \textit{(iv)-(b)} follow from the well known fact that $\H^{s,p}(\RR^m)$, $\W^{s,p}(\RR^m)$ and $\B^{s}_{p,q}(\RR^m)$ are algebras in this precise case.

The point \textit{(iii)-(c)} follows from \cite[Section~4.6.1,~Theorem~2,~eq.~(20)]{RunstSickel96}.

Therefore, it remains to show the second part of \textit{(iv)-(a)}. We are going to take advantage of Proposition~\ref{prop:EmbeddingsHomogeneousSobolevBesovinMultiplierSpaces}.

Now, the case ${\B}^{s}_{r,1}$ where $r\geqslant \frac{m}{s}$ and $\dot{\B}^{s}_{\frac{m}{s},1}$ when $s\in\NN$. Note that, one has
\begin{align*}
    {\B}^{k}_{\infty,1}(\RR^{m})\cdot\W^{k,1}(\RR^{m})\hookrightarrow\W^{k,\infty}(\RR^{m})\cdot\W^{k,1}(\RR^{m})\hookrightarrow\W^{k,1}(\RR^{m})\text{, }k\in\NN,
\end{align*}
so that by bilinear real interpolation, one obtains
\begin{align}\label{proof:MultipilerWs1-1}
    \B^{s}_{\infty,1}(\RR^{m})\cdot\W^{s,1}(\RR^{m})\hookrightarrow\W^{s,1}(\RR^{m})\text{, }
\end{align}
for all $s\geqslant 0$. By Leibniz rule and Sobolev embeddings it can be checked that one has for $k\in\llb0,m\rrb$
\begin{align*}
    \B^{k}_{\frac{m}{k},1}\cdot\W^{k,1}(\RR^{m})\hookrightarrow\W^{k,1}(\RR^{m}).
\end{align*}
Note that the same dilation argument as presented in the proof of Proposition~\ref{prop:EmbeddingsHomogeneousSobolevBesovinMultiplierSpaces} yields
\begin{align*}
    \dot{\B}^{k}_{\frac{m}{k},1}(\RR^{m})\cdot\dot{\W}^{k,1}(\RR^{m})\hookrightarrow\dot{\W}^{k,1}(\RR^{m}).
\end{align*}
Therefore, using the embedding $\dot{\B}^{k}_{\frac{m}{k},1}\hookrightarrow\L^\infty$, the fact that $\L^\infty\cdot\L^1\hookrightarrow \L^1$, and the equivalence of norms $\L^1\cap\dot{\W}^{k,1}=\W^{k,1}$, we deduce
\begin{align}\label{proof:MultipilerWs1-2}
    \dot{\B}^{k}_{\frac{m}{k},1}(\RR^{m})\cdot{\W}^{k,1}(\RR^{m})\hookrightarrow{\W}^{k,1}(\RR^{m}).
\end{align}
So that by bilinear real interpolation applied top \eqref{proof:MultipilerWs1-2}, one obtains, even in integer cases, by what precedes,
\begin{align}\label{proof:MultipilerWs1-3}
    \dot{\B}^{s}_{\frac{m}{s},1}(\RR^{m})\cdot\W^{s,1}(\RR^{m})\hookrightarrow\W^{s,1}(\RR^{m})\text{.}
\end{align}

Up to replace $\dot{\B}^{s}_{\frac{m}{s},1}(\RR^{m})$ by ${\B}^{s}_{\frac{m}{s},1}(\RR^{m})$, complex interpolation between \eqref{proof:MultipilerWs1-1}  and \eqref{proof:MultipilerWs1-3} yields the case $r\geqslant \frac{m}{s}$.
\end{proof}

\subsubsection{Parametrisation of some rough domains, traces of multipliers}\label{sec:para}
We follow the presentation from \cite[Section~3]{Breit2024FSI2D}.
 Let $\Omega\subset\RR^n$ be a bounded open set.
We assume that $\partial{\Omega}$ can be covered by a finite
number of open sets $\mathcal U^1,\dots,\mathcal U^\ell$ for some $\ell\in\mathbb N$, such that
the following holds. For each $j\in\llb1,\ell\rrb$ there is a reference point
$y_j\in\RR^n$ and a local coordinate system $\{\mathfrak{e}^j_m,\,m\in\llb 1,n\rrb\}$ (which we assume
to be orthonormal and set $\mathcal{Q}_j=(\mathfrak{e}_1^j|\dots |\mathfrak{e}_n^j)\in\RR^{n\times n}$), a function
$\varphi_j:\RR^{n-1}\rightarrow\RR$
and $r_j>0$
with the following properties:
\begin{enumerate}[label={\bf (A\arabic{*})}]
\item\label{A1} There is $h_j>0$ such that
$$\mathcal{U}^j=\{x=\mathcal Q_jz+y_j\in\RR^n:\,z=(z',z_n)\in\RR^n,\,|z'|<r_j,\,
|z_n-\varphi_j(z')|<h_j\}.$$
\item\label{A2} For $x\in\mathcal U^j$ we have with $z=\prescript{t}{}{\mathcal {Q}_j}(x-y_j)$
\begin{itemize}
\item $x\in\partial{\Omega}$ if and only if $z_n=\varphi_j(z')$;
\item $x\in{\Omega}$ if and only if $0<z_n-\varphi_j(z')<h_j$;
\item $x\notin{\Omega}$ if and only if $0>z_n-\varphi_j(z')>-h_j$.
\end{itemize}
\item\label{A3} We have that
$$\partial{\Omega}\subset \bigcup_{j=1}^\ell\mathcal U^j.$$
\end{enumerate}
In other words, for any $x_0\in\partial{\Omega}$, there is a neighborhood $\mathcal{U}$ of $x_0$ and a function $\varphi:\RR^{n-1}\rightarrow\RR$ such that after translation and rotation\footnote{By translation via $y_j$ and rotation via $\mathcal Q_j$ we can assume that $x_0=0$ and that the outer normal at~$x_0$ is pointing in the negative $x_n$-direction.}
 \begin{align}\label{eq:3009}
 \mathcal{U} \cap {\Omega} = \mathcal{U} \cap \tilde{\Omega},\quad \tilde{\Omega} = \{(x',x_n)\in \RR^n \,:\, x' \in \RR^{n-1}, x_n > \varphi(x')\}.
 \end{align}
 The regularity of $\partial{\Omega}$ will be described by means of local coordinates as just described.
 \begin{definition}\label{def:besovboundary}
 Let ${\Omega}\subset\RR^n$ be a bounded domain, $\alpha>0$, $1\leqslant \rho,r\leqslant\infty$, and $\varepsilon>0$. If there exists $\ell\in\mathbb N$ and functions $\varphi_1,\dots,\varphi_\ell\in\C^{0}(\RR^{n-1})$   satisfying \ref{A1}--\ref{A3}, we say that $\partial{\Omega}$ belongs, either, to the class 
 \begin{enumerate}
     \item $\mathbf{B}^\alpha_{\rho,r}$ if for all $1\leqslant j\leqslant N $, $\varphi_j\in\B^\alpha_{\rho,r}(\RR^{n-1})$;
     \item $\mathbf{W}^{\alpha,\rho}$ if for all $1\leqslant j\leqslant N $, $\varphi_j\in\W^{\alpha,\rho}(\RR^{n-1})$;
     \item $\mathcal{M}^{\alpha,\rho}_{\X}$ if for all $1\leqslant j\leqslant N $, $\varphi_j\in\mathcal{M}^{\alpha,\rho}_{\X}(\RR^{n-1})$;
     \item $\mathcal{M}^{\alpha,\rho}_{\X}(\varepsilon)$ if for all $1\leqslant j\leqslant N $, $\varphi_j\in\mathcal{M}^{\alpha,\rho}_{\X}(\RR^{n-1},\varepsilon)$.
 \end{enumerate}
 \end{definition}
Clearly, a similar definition applies for a Lipschitz boundary (or even a $\C^1$ or a $\C^{1,\alpha}$-boundary with $\alpha\in(0,1]$) by requiring that $\varphi_1,\dots,\varphi_\ell\in \C^{0,1}_{b}(\RR^{n-1})$ (or $\varphi_1,\dots,\varphi_\ell\in \C^{1,\alpha}(\RR^{n-1})$). We say that the local Lipschitz constant of $\partial{\Omega}$, denoted by $L_{\Omega})$, is (smaller or) equal to some number $L>0$ provided the Lipschitz constants of $\varphi_1,\dots,\varphi_\ell$ are not exceeding $L$. Our main results depend on the assumption of a sufficiently small local Lipschitz constant.

While this seems rather restrictive at first glance, it appears quite natural when looking closer. Indeed, it holds, for instance, if the regularity
of $\partial{\Omega}$ is better than Lipschitz (such as $\C^{1}$ or $\C^{1,\alpha}$, $\alpha\in(0,1]$).

By means of the transformations $\mathcal Q_j$ introduced above, we can assume that the reference point $y^j$ in question is the origin and that $\nabla\varphi_j(0)=0$. So for a $\C^{1}$ bounded domain, by continuity of the derivative, the Lipschitz constant can be made small up to increase the number of charts. For $\C^{1,\alpha}$, this can be made more explicitly: choosing $r_j$ in \ref{A1} small enough (which can be achieved simply by allowing more sets in the cover $\mathcal U^1,\dots,\mathcal U^l$) we have
\begin{align*}
|\nabla\varphi_j(z')|=|\nabla\varphi_j(z')-\nabla\varphi_j(0)|\leqslant\,r_j^\alpha\lVert\nabla\varphi_j\rVert_{\dot{\C}^\alpha}\ll 1
\end{align*}
for all $z'$ with $|z'|\leqslant r_j$. Similarly to the case of $\C^1$-domains, up to increase the number of charts, for $\varphi_j\in\mathcal{C}^{1,\alpha}$, $0<\alpha<1$, (the little H\"{o}lder-Lipschitz space), the $\alpha$-H\"{o}lder-Lipschitz norm $\lVert\nabla\varphi_j\rVert_{\dot{\C}^\alpha}$ can be made arbitrarily small, since by definition, one has
\begin{align*}
    \lim_{|x'-y'|\rightarrow 0} \frac{|\nabla \varphi_j(x')-\nabla \varphi_j(y')|}{|x'-y'|^\alpha} =0.
\end{align*}
Note that,  for $0<\alpha<\beta<1$, the $\C^{1,\beta}$-domains are in particular $\mathcal{C}^{1,\alpha}$-domains.

In order to describe precisely the behaviour of functions defined in ${\Omega}$ close to the boundary, we need to extend the functions $\varphi_1,\dots,\varphi_\ell$  from \ref{A1}--\ref{A3} to the half space
$\RR^n_+ := \{z = (z',z_n)\,:\, z_n > 0\}$. Hence, we are confronted with the task of extending a function~$\phi\,:\, \RR^{n-1}\to \RR$ to a mapping $\Phi\,:\, \RR^n_+ \to \RR^n$ that maps the 0-neighborhood in~$\RR^n_+$ to the $x_0$-neighborhood in~${\Omega}$. The mapping $(z',0) \mapsto (z',\phi(z'))$ locally maps the boundary of~$\RR^n_+$ to the one of~$\partial {\Omega}$. We extend this mapping using the extension operator of Maz'ya and Shaposhnikova~\cite[Section 9.4.3]{MazyaShaposhnikova2009}. Let $\zeta \in \Ccinfty(\B_1(0'))$ with $\zeta \geqslant 0$ and $\int_{\RR^{n-1}} \zeta(x')\,\d x'=1$ and $\zeta$ to be radial\footnote{We also mention here that in \cite{MazyaShaposhnikova2009}, $\zeta$ also has some additional conditions on its moments.}. Let $\zeta_t(x') := t^{-(n-1)} \zeta(x'/t)$ denote the induced family of mollifiers. We define the extension operator 
\begin{equation}\label{eq:ExtOpMazya}
    (\mathcal{T}\phi)(z',z_n)=\int_{\RR^{n-1}} \zeta_{z_n}(z'-y')\phi(y')\,\d y',\quad (z',z_n) \in \RR^n_+,
\end{equation}
where~$\phi:\RR^{n-1}\to \RR$ is a Lipschitz function with Lipschitz constant~$K$. 

\textbf{A major feature of the Sobolev-Besov multiplier theory} developed by Maz'ya and Shaposhnikova is that multipliers are stable under taking traces. That is, if one has a sufficient amount of regularity to define the trace, say  on $\W^{s,p}(\RR^n_+)$, $s>\frac{1}{p}$, then it holds that one can define naturally $\varphi_{|_{\partial\RR^n_+}} \in \mathcal{M}^{s-{\sfrac{1}{p}},p}_{\W,\tt{or}}(\RR^{n-1})$, provided $\varphi \in \mathcal{M}^{s,p}_{\W,\tt{or}}(\RR^{n}_+)$, thanks to the ontoness of the trace operator from $\W^{s,p}(\RR^n_+)$ to $\W^{s-\sfrac{1}{p},p}(\partial\RR^n_+)$, $s\notin \NN + \sfrac{1}{p}$.

This is even better. More precisely, \cite[Theorems~8.7.1,~8.7.2~\&~8.8.1]{MazyaShaposhnikova2009} and \cite[Theorem~1]{Shaposhnikova1989} can be summarized has follows:
\begin{theorem}[ Maz'ya \& Shaposhnikova ]\label{thm:LiftTraceMultiplier}Let $p\in[1,\infty)$, $s > {\sfrac{1}{p}}$ and $\X^{s,p}\in\{\H^{s,p},\W^{s,p}\}$, assuming
\begin{enumerate}
    \item $p>1$ when $\X=\H$; or
    \item $s\notin\NN$ if $p=1$.
\end{enumerate}

Then the trace operator $[\cdot]_{|_{\partial\RR^n_+}}\,:\,\mathcal{M}^{s,p}_{\X,\tt{or}}(\RR^{n}_+)\longrightarrow\mathcal{M}^{s-\sfrac{1}{p},p,p}_{\B,\tt{or}}(\RR^{n-1})$ is well-defined, bounded and onto with right bounded inverse $\mathcal{T}$ given by \eqref{eq:ExtOpMazya}, \textit{i.e.} for all $\varphi\in\mathcal{M}^{s-\sfrac{1}{p},p,p}_{\B,\tt{or}}(\RR^{n-1})$,
\begin{align*}
    \lVert \mathcal{T}\varphi\rVert_{\mathcal{M}^{s,p}_{\X,\tt{or}}(\RR^{n}_+)}\lesssim_{p,s,n}\lVert \varphi\rVert_{\mathcal{M}^{s-{\sfrac{1}{p}},p,p}_{\B,\tt{or}}(\RR^{n-1})}.
\end{align*}
Furthermore, if $\varphi\,:\,\RR^{n-1}\longrightarrow\CC$ is such that $\nabla'\varphi\in \mathcal{M}^{s-{\sfrac{1}{p}},p,p}_{\B,\tt{or}}(\RR^{n-1})$, one also has $\nabla\mathcal{T}\varphi\in\mathcal{M}^{s,p}_{\X,\tt{or}}(\RR^{n}_+)$ with an estimate
\begin{align*}
    \lVert \nabla \mathcal{T}\varphi\rVert_{\mathcal{M}^{s,p}_{\X,\tt{or}}(\RR^{n}_+)}\lesssim_{p,s,n}\lVert \nabla'\varphi\rVert_{\mathcal{M}^{s-{\sfrac{1}{p}},p,p}_{\B,\tt{or}}(\RR^{n-1})}.
\end{align*}
\end{theorem}

\begin{remark}The case $p=\infty$ is not treated, notice point \textit{(v)} of Lemma~\ref{lem:SmoothinOPMaSh} below is an appropriate surrogate, since in this case those function spaces are algebras (with substitute $s>0$ to the condition $s>\frac{1}{p}$).
\end{remark}

It is shown in \cite[Lemma 9.4.5]{MazyaShaposhnikova2009} that (for sufficiently large~$N$, i.e., $N \geqslant c(\phi) K+1$) the mapping
\begin{align*}
  \alpha_{z'}(z_n) \mapsto N\,z_n+(\mathcal{T} \phi)(z',z_n)
\end{align*}
is for every $z' \in \RR^{n-1}$ one to one and the inverse is Lipschitz with its gradient
bounded by $(N-K)^{-1}$.
Now, we define the mapping~$\boldsymbol{\Phi}\,:\, \RR^n_+ \to \RR^n$ as a rescaled version of the latter one by setting
\begin{align}\label{eq:Phi}
  \boldsymbol{\Phi}(z',z_n)
  &:=
    \big(z',
    \alpha_{z'}(z_n)\big) = 
    \big(z',
    \,z_n + (\mathcal{T} \phi)(z',z_n/K)\big).
\end{align}

Thus, $\boldsymbol{\Phi}$ is one-to-one (for sufficiently large~$N=N(K)$) and we can define its inverse $\boldsymbol{\Phi}^{-1}$.
The Jacobi matrix of the mapping $\boldsymbol{\Phi}$ is of the form
\begin{align}\label{J}
   \nabla \mathbf{\Phi} = 
  \begin{pmatrix}
    \mathbf{I}_{n-1}&0
    \\
    \partial_{z'} (\mathcal{T}  \phi)& 1+ 1/N\partial_{z_n}\mathcal{T}  \phi
  \end{pmatrix}.
\end{align}
Since $|\partial_{z_n}\mathcal{T}  \phi| \leqslant K$, we have \begin{align}\label{eq:detJ}\frac{1}{2} < 1-K/N \leqslant |\det(\nabla \mathbf{\Phi})| \leqslant 1+K/N\leqslant 2\end{align}
using that $N$ is large compared to~$K$. Finally, when $s\geqslant 1$, we note the implications
\begin{align} \label{eq:SMPhiPsi}
\nabla\boldsymbol{\Phi}\in\mathcal{M}^{s-1,p}_{\X,\tt{or}}(\RR^n_+)\,\,&\implies \,\,(\nabla\boldsymbol{\Phi})^{-1}\in\mathcal{M}^{s-1,p}_{\X,\tt{or}}(\RR^n_+),\\
\boldsymbol{\Phi}\in\mathcal{M}^{s,p}_{\X}(\RR^n_+)\,\,&\implies \,\,\boldsymbol{\Phi}^{-1}\in\mathcal{M}^{s,p}_{\X}(\RR^n_+),
\end{align}
which hold, for instance, if $\boldsymbol{\Phi}$ is Lipschitz continuous, cf. \cite[Lemma 9.4.2]{MazyaShaposhnikova2009}. In fact, one can prove \eqref{eq:SMPhiPsi} with the help of \eqref{eq:detJ}.

Note that the use of Theorem~\ref{thm:LiftTraceMultiplier} when $1\leqslant p<\infty$ with an abstract multiplier norm is somewhat necessary when one deals with those rough boundaries. Indeed, in standard ``concrete function spaces'': from such point of view the extension operator $\mathcal{T}$ does not behave well on the boundary charts, since it induces a loss of localisation in an optimal way due to scaling invariance, except when $p=\infty$, as shown by Lemma~\ref{lem:SmoothinOPMaSh} below. 

It is followed by Example~\ref{ex:"Conrete"VS"Abstract"MultiplierTheory} which demonstrates on how relevant is using abstract multiplier theory instead of multiplier estimates in standard function spaces. 

\begin{lemma}\label{lem:SmoothinOPMaSh}Let $s \in\RR$, $p,q\in[1,\infty]$. Then it holds that $\mathcal{T}$ maps continuously
\begin{enumerate}
    \item $\dot{\B}^{s-{\sfrac{1}{p}}}_{p,p}(\RR^{n-1})$ to $\dot{\H}^{s,p}(\RR^{n}_+)$, provided $1<p<\infty$;
    \item $\dot{\B}^{k-{\sfrac{1}{p}}}_{p,p}(\RR^{n-1})$ to $\dot{\W}^{k,p}(\RR^{n}_+)$, provided $k\in\NN$, $p<\infty$;
    \item $\W^{k,\infty}(\RR^{n-1})$ to ${\W}^{k,\infty}(\RR^{n}_+)$, provided $k\in\NN$, with homogeneous estimates;
    \item $\dot{\B}^{s-{\sfrac{1}{p}}}_{p,q}(\RR^{n-1})$ to $\dot{\B}^{s}_{p,q}(\RR^{n}_+)$;
    \item $\B^{s}_{\infty,q}(\RR^{n-1})$ to $\B^{s}_{\infty,q}(\RR^{n}_+)$.
\end{enumerate}
Furthermore, a similar result holds for endpoint function spaces ${\B}^{\bullet,0}_{\infty,q}$, $\dot{\BesSmo}^{\bullet}_{p,\infty}$, ${\BesSmo}^{\bullet,0}_{\infty,\infty}$ and $\C^{k}_0$.
\end{lemma}

\begin{proof}Before we actually start the proof, we note that one has boundedness of $\mathcal{T}$ as a bounded linear operator from  $\L^{\infty}(\RR^{n-1})$ to $\L^{\infty}(\RR^{n}_+)$, but also from $\C^{0}_0(\RR^{n-1})$ to $\C^0_{0}(\overline{\RR^{n}_+})$.

\textbf{Step 1:} We prove the boundedness from $\dot{\B}^{-{\sfrac{1}{p}}}_{p,p}(\RR^{n-1})$ to $\L^{p}(\RR^n_+)$ provided $1\leqslant p<\infty$. We deal with the case $p=1$ first. First, using the identity $u\ast v = [(-\Delta')u]\ast [(-\Delta')^{-1}v]$, we can write
\begin{align*}
    \lVert \mathcal{T}\phi\rVert_{\L^{1}(\RR^n_+)}&= \int_{\RR^{n}_+} \left|\int_{\RR^{n-1}} [-\Delta'\zeta] (y') \frac{1}{z_n^{2}}[(-\Delta')^{-1}\phi]( y' z_n + z')\d y' \right| \d z' \d z_n\\
    &= \int_{\RR^{n}_+} \left|\int_{\RR^{n-1}} [-\Delta'\zeta] (y') \frac{1}{z_n^{2}}[\tilde{\phi}( y' z_n + z')-\tilde{\phi}( z')]\d y' \right| \d z' \d z_n.
\end{align*}
where $\tilde{\phi}:=[(-\Delta')^{-1}\phi]$ and noting that the compact support and smoothness of $\zeta$ implies $\int_{\RR^{n-1}} {\Delta'}\zeta =0$. Since $\zeta$ is an even function and compactly supported in $\B_{1}(0')$, it can be deduced
\begin{align*}
    \lVert \mathcal{T}\phi\rVert_{\L^{1}(\RR^n_+)} &=\frac{1}{2}\int_{\RR_+}\int_{\RR^{n-1}} \left|\int_{\RR^{n-1}} [-\Delta'\zeta] (y') \frac{1}{t^{2}}[\tilde{\phi}(z'+ y' t ) +\tilde{\phi}( z'-y' t ) - 2\tilde{\phi}( z')]\d y' \right| \d z' \d t\\
    &=\frac{1}{2}\int_{\RR_+}\int_{\RR^{n-1}} \left|\int_{\RR^{n-1}} \frac{1}{t^{n+1}}[-\Delta'\zeta] \left(\frac{y'}{t}\right)[\tilde{\phi}(z'+ y' ) +\tilde{\phi}( z'-y' ) - 2\tilde{\phi}( z')]\d y' \right| \d z' \d t\\
    &\leqslant \frac{1}{2}\int_{\RR^{n-1}} \int_{\RR^{n-1}} \int_{|y'|}^{\infty}\left|\frac{1}{t^{n+1}}[-\Delta'\zeta] \left(\frac{y'}{t}\right)\right|\left|\tilde{\phi}(z'+ y' ) +\tilde{\phi}( z'-y' ) - 2\tilde{\phi}( z') \right| \d t\, \d y' \d z' 
\end{align*}
where the triangle inequality as well as Minkowski's integral inequality have been applied to obtain the last line.
Therefore, one obtains by \cite[Theorem~2.37]{bookBahouriCheminDanchin}:
\begin{align*}
    \lVert \mathcal{T}\phi\rVert_{\L^{1}(\RR^n_+)} &\leqslant \frac{1}{2n}\lVert \Delta'\zeta\rVert_{\L^\infty(\RR^{n-1})}\int_{\RR^{n-1}} \int_{\RR^{n-1}} \frac{|\tilde{\phi}(z'+ y' ) +\tilde{\phi}( z'-y' ) - 2\tilde{\phi}( z') |}{|y'|^{n+1}}  \d z' \d y' \\
    &\lesssim_{\zeta,n} \lVert \tilde{\phi}\rVert_{\dot{\B}^{1}_{1,1}(\RR^{n-1})}  =\lVert (-\Delta')^{-1}{\phi}\rVert_{\dot{\B}^{1}_{1,1}(\RR^{n-1})}\\
    &\lesssim_{\zeta,n} \lVert {\phi}\rVert_{\dot{\B}^{-1}_{1,1}(\RR^{n-1})}. 
\end{align*}

Now, since one also has
\begin{align*}
     \lVert \mathcal{T}\phi\rVert_{\L^{\infty}(\RR^n_+)} \leqslant \lVert \phi\rVert_{\L^{\infty}(\RR^{n-1})} \leqslant \lVert \phi\rVert_{\dot{\B}^{0}_{\infty,1}(\RR^{n-1})}
\end{align*}
for all $\phi\in\dot{\B}^{0}_{\infty,1}(\RR^{n-1})$, it follows by the real interpolation identities
\begin{align*}
    (\dot{\B}^{-1}_{1,1},\dot{\B}^{0}_{\infty,1})_{1-{\sfrac{1}{p}},p}=\dot{\B}^{-{\sfrac{1}{p}}}_{p,p}\text{ and }(\L^1,\L^\infty)_{1-{\sfrac{1}{p}},p}=\L^{p}
\end{align*}
true for all $1<p<\infty$, that $\mathcal{T}$ maps $\dot{\B}^{-{\sfrac{1}{p}}}_{p,p}(\RR^{n-1})$ to ${\L}^{p}(\RR^{n}_+)$.

\textbf{Step 2:} We prove that $\mathcal{T}$ maps $\dot{\B}^{1-{\sfrac{1}{p}}}_{p,p}(\RR^{n-1})$ to $\dot{\W}^{1,p}(\RR^{n}_+)$, if $p<\infty$ and ${\W}^{1,\infty}(\RR^{n-1})$ to ${\W}^{1,\infty}(\RR^{n}_+)$ with homogeneous norms. Since $\nabla'$ commutes with the convolution, it amounts to bound the $\L^p$-norm of $\partial_{z_n} \mathcal{T}\phi$. Note that one has the following formula
\begin{align*}
    \partial_{z_n}\mathcal{T}\phi(z') = -(n-1)\int_{\RR^{n-1}} \frac{1}{z_n} \zeta_{z_n}(z'-y')\phi(y') \d y' -\sum_{k=1}^{n-1}\int_{\RR^{n-1}} \frac{y_k}{z_n^2}[\partial_{y_k}\zeta]_{z_n}\left({y'}\right) \phi(z'-y') \d y'.
\end{align*}
We focus on the second term, which can be written as
\begin{align*}
    -\sum_{k=1}^{n-1}\int_{\RR^{n-1}} \frac{y_k}{z_n^2}[\partial_{y_k}\zeta]_{z_n}\left({y'}\right) \phi(z'-y') \d y'  &= -\sum_{k=1}^{n-1}\int_{\RR^{n-1}} \frac{1}{z_n}[\partial_{y_k}(y_k \zeta_{z_n})- \zeta_{z_n}]\left({y'}\right) \phi(z'-y') \d y'\\
    &= \sum_{k=1}^{n-1}\Bigg[\int_{\RR^{n-1}} \frac{1}{z_n}[y_k \zeta_{z_n}](y') [\partial_{y_k}\phi](z'-y') \d y' \\
    &\qquad\qquad\qquad +\int_{\RR^{n-1}} \frac{1}{z_n}\zeta_{z_n} \left({y'}\right) \phi(z'-y') \d y'\Bigg].
\end{align*}
This yields the equality
\begin{align*}
    \partial_{z_n}\mathcal{T}\phi(z') = \sum_{k=1}^{n-1}\int_{\RR^{n-1}} [y_k \zeta]_{z_n}(y') [\partial_{y_k}\phi](z'-y') \d y'.
\end{align*}
For the space $\W^{1,\infty}$, the result is straightforward from Young's inequality for convolutions. Otherwise, one can apply the result from Step 1, for $p=1,\infty$, $n-1$ times -- $\partial_{z_n}\mathcal{T}$ behaves as the operator $\mathcal{T}$ with kernel by $y\mapsto y\zeta(y)$ replacing $\zeta$-- in order to obtain
\begin{align*}
    \lVert \partial_{z_n}\mathcal{T}\phi\rVert_{\L^p(\RR^{n}_+)} \lesssim_{\zeta,n,p} \lVert \nabla' \phi\rVert_{\dot{\B}^{-{\sfrac{1}{p}}}_{p,p}(\RR^{n-1})}.
\end{align*}
Higher order estimates follow by induction. This yields the cases $s=k\in\NN$.

Now for negative regularity estimates, consider the (pre-)dual operator
\begin{align*}
    \mathcal{T}^{\ast}v(x') = \int_{\RR^{n-1}}\int_{0}^{\infty} \zeta_{z_n}(y')\phi(x'-y',z_n) \d z_n \,\d y'.
\end{align*}
By the previous step, and duality, one obtains for all $v\in\Ccinfty(\RR^{n}_+)$, $p\in(1,\infty]$,
\begin{align*}
    \lVert \mathcal{T}^{\ast}v \rVert_{\dot{\B}^{{\sfrac{1}{p'}}}_{p,p}(\RR^{n-1})}\lesssim_{\zeta,p,n} \lVert v \rVert_{{\L}^{p}(\RR^{n}_+)},
\end{align*}
Obviously $\nabla'$ commutes with $\mathcal{T}^{\ast}$, which yields easily boundedness from $\dot{\W}^{k,p}_0(\RR^n_+)$ to $\dot{\B}^{k+{\sfrac{1}{p'}}}_{p,p}(\RR^{n-1})$. By Fubini-Tonelli, obtains boundedness of $\mathcal{T}^{\ast}$ from $\L^{1}(\RR^{n}_+)$ to $\L^{1}(\RR^{n-1})$, then boundedness from $\dot{\W}^{k,1}_{0}(\RR^{n}_+)$ to $\dot{\W}^{k,1}(\RR^{n-1})$.

Complex and real interpolation for $\mathcal{T}$ and $\mathcal{T}^{\ast}$ yields the full result. For real and complex interpolation when $p<\infty$, due to the lack of completeness one should use the trick introduced in \cite[Proposition~3.12~\&~Remark~3.13]{Gaudin2023Lip}.
\end{proof}

\begin{example}\label{ex:"Conrete"VS"Abstract"MultiplierTheory} The goal of this example is to illustrate the use of the Sobolev-Besov multiplier theory in contrast to "by hands" multiplier estimates with standard Sobolev-Besov norms.
We consider $\varphi \in \C^{0,\alpha}_{ub}(\RR^{n-1})$, $\alpha\in(0,1)$, we want to deal with multiplier properties of $\mathcal{T}\phi$ on $\W^{s,1}(\RR^n_+)=\B^{s}_{1,1}(\RR^n_+)$, $0<s<1$.
\begin{itemize}
    \item The naive approach by standard methods --Lemma~\ref{lem:SmoothinOPMaSh}-- : one can only obtain \textit{a priori} $\mathcal{T}\phi\in\C^{0,\alpha}_{ub}(\RR^{n}_+)=\B^{\alpha}_{\infty,\infty}(\RR^n_+)$, with an estimate
    \begin{align*}
        \lVert\mathcal{T}\phi \rVert_{\B^{\alpha}_{\infty,\infty}(\RR^n_+)}\lesssim_{n,\alpha,\zeta} \lVert\phi \rVert_{\B^{\alpha}_{\infty,\infty}(\RR^{n-1})}.
    \end{align*}
    Therefore, in order to expect
    \begin{align*}
         \lVert (\mathcal{T}\phi) u \rVert_{\W^{s,1}(\RR^n_+)} \lesssim_{s,n,\alpha} \lVert\mathcal{T}\phi \rVert_{\B^{\alpha}_{\infty,\infty}(\RR^n_+)}\lVert u \rVert_{\W^{s,1}(\RR^n_+)} \lesssim_{s,n,\alpha,\zeta} \lVert\phi \rVert_{\B^{\alpha}_{\infty,\infty}(\RR^{n-1})}\lVert u \rVert_{\W^{s,1}(\RR^n_+)}
    \end{align*}
    from these sole information and by applying Proposition~\ref{prop:EmbeddingsSobolevBesovinMultiplierSpaces} (thanks to Remark~\ref{rem:MultipliersonHalfSpaces}), this requires the condition
    \begin{align*}
        s\in(0,\alpha).
    \end{align*}
    \item The approach by Multiplier Theory --Theorem~\ref{thm:LiftTraceMultiplier}--: we start not considering $\mathcal{T}\varphi$, but only $\varphi$. For all $s\in(0,\alpha)$, one has by Proposition~\ref{prop:EmbeddingsSobolevBesovinMultiplierSpaces}
    \begin{align*}
        \varphi\in\mathcal{M}^{s,1}_{\W,\tt{or}}(\RR^{n-1}).
    \end{align*}
    So by Theorem~\ref{thm:LiftTraceMultiplier}, one has $\mathcal{T}\varphi\in\mathcal{M}^{s+1,1}_{\W}(\RR^{n}_+)$ with an estimate
    \begin{align*}
        \lVert\mathcal{T}\phi \rVert_{\mathcal{M}^{s+1,1}_{\W,\tt{or}}(\RR^{n}_+)}\lesssim_{s,n,\zeta} \lVert\phi \rVert_{\mathcal{M}^{s,1}_{\W,\tt{or}}(\RR^{n-1})}\lesssim_{s,\alpha,n,\zeta}\lVert\phi \rVert_{\B^{\alpha}_{\infty,\infty}(\RR^{n-1})}.
    \end{align*}
    valid for all $s\in(0,\alpha)$. However, by Proposition~\ref{prop:embeddingmultipliers}, one deduces even 
    \begin{align*}
        \mathcal{T}\varphi\in\mathcal{M}^{s,1}_{\W,\tt{or}}(\RR^{n}_+),\qquad \textrm{for all } s\in(0,1+\alpha).
    \end{align*}
\end{itemize}
More generally on $\X^{s,p}$, $\mathcal{T}\phi\in\mathcal{M}^{s,p}_{\X,\tt{or}}(\RR^{n}_+)$, $s\in(0,\alpha+\sfrac{1}{p})$, while being restricted by regularity at most $\alpha$ when using multiplier estimates in standard function spaces ! 
\end{example}

We conclude this section with a result that allows to deal "uniformly" with all function spaces within the regularity range $s\in(-1+\sfrac{1}{p},\sfrac{1}{p})$, assuming just arbitrary small regularity with respect to an arbitrary function space.

\begin{proposition}\label{prop:MultipliersintheLplikeRange} Let $r\in[1,\infty]$, $\alpha\in(0,1)$. Then for all $\varphi\in\mathcal{M}_{\W}^{1+\alpha,r}(\RR^{n-1})$, one has
\begin{align*}
    \nabla \mathcal{T}\varphi \in \mathcal{M}_{\X,\tt{or}}^{s,p}(\RR^{n}_+)
\end{align*}
for all $p,q\in[1,\infty]$, all $s\in\RR$ such that either
 \begin{enumerate}
    \item $p\in(r,\infty]$, $s\in(-1+\frac{1}{p},\frac{1+r\alpha}{p})$; or
    \item $p\in[1,r]$, $s\in(-1+\frac{1}{p},\alpha+\frac{1}{p})$; or
    \item $p=q=r$, $s=\alpha+\frac{1}{p}$; 
\end{enumerate}
assuming $1<p<\infty$ whenever $\X=\H$. One also has the corresponding estimate
\begin{align*}
    \lVert \nabla \mathcal{T}\varphi\rVert_{\mathcal{M}^{s,p}_{\X,\tt{or}}(\RR^n_+)}\lesssim_{p,s,n}^{r,\alpha}\lVert \varphi\rVert_{\mathcal{M}^{1+\alpha,r}_{\W}(\RR^{n-1})}.
\end{align*}
In particular it holds for all $s\in(-1+\sfrac{1}{p},\sfrac{1}{p})$.
\end{proposition}

\begin{proof} \textbf{Step 1:} First, we deal with the cases $r=1,\infty$. If $r=1$, for $\varphi\in\mathcal{M}_{\W}^{1+\alpha,1}(\RR^{n-1})$, by Theorem~\ref{thm:LiftTraceMultiplier}, one obtains
\begin{align*}
    \nabla \mathcal{T}\varphi \in \mathcal{M}_{\W,\tt{or}}^{\alpha+1,1}(\RR^{n}_+)
\end{align*}
with the estimate
\begin{align*}
    \lVert \nabla \mathcal{T}\varphi\rVert_{\mathcal{M}^{\alpha+1,1}_{\W,\tt{or}}(\RR^n_+)}\lesssim_{n}^{\alpha}\lVert \varphi\rVert_{\mathcal{M}^{1+\alpha,1}_{\W}(\RR^{n-1})}.
\end{align*}
Since $\nabla\mathcal{T}\varphi\in\L^\infty(\RR^{n}_+)$, with the estimate
\begin{align*}
    \lVert \nabla \mathcal{T}\varphi\rVert_{\mathcal{M}^{1}_{\L,\tt{or}}(\RR^n_+)} = \lVert \nabla \mathcal{T}\varphi\rVert_{\L^\infty(\RR^n_+)} \lesssim_{n}\lVert \nabla' \varphi\rVert_{\L^\infty(\RR^{n-1})} \lesssim_{n}^{\alpha}\lVert \varphi\rVert_{\mathcal{M}^{1+\alpha,1}_{\W}(\RR^n)},
\end{align*}
by real interpolation, $\nabla \mathcal{T}\varphi\in\mathcal{M}_{\B,\tt{or}}^{s,1,q}(\RR^{n}_+)$ for all $0<s<1$, all $q\in[1,\infty]$, with an estimate
\begin{align*}
    \lVert \nabla \mathcal{T}\varphi\rVert_{\mathcal{M}^{s,1,q}_{\B,\tt{or}}(\RR^n_+)}\lesssim_{n,s}^{\alpha} \lVert \nabla \mathcal{T}\varphi\rVert_{\mathcal{M}^{1}_{\L,\tt{or}}(\RR^n_+)}^{1-\theta}\lVert \nabla \mathcal{T}\varphi\rVert_{\mathcal{M}^{\alpha+1,1}_{\W,\tt{or}}(\RR^n_+)}^{\theta} \lesssim_{n,s}^{\alpha}\lVert \varphi\rVert_{\mathcal{M}^{1+\alpha,1}_{\W}(\RR^{n-1})}.
\end{align*}
By duality, see Proposition~\ref{prop:basicMultipliers}, one obtains all $q\in(1,\infty)$, all $s\in(-1,0)$,
\begin{align*}
    \lVert \nabla \mathcal{T}\varphi\rVert_{\mathcal{M}^{s,\infty,q}_{\B,\tt{or}}(\RR^n_+)} \leqslant \lVert \nabla \mathcal{T}\varphi\rVert_{\mathcal{M}^{-s,1,q'}_{\B,\tt{or}}(\RR^n_+)}\lesssim_{n,s}^{\alpha}\lVert \varphi\rVert_{\mathcal{M}^{1+\alpha,1}_{\W}(\RR^{n-1})}.
\end{align*}
By complex and real interpolation between the case $p=1,\infty$, then real interpolation to get $q=1,\infty$, one obtains the result for all $s\in(-1+\sfrac{1}{p},\sfrac{1}{p})$:
\begin{align*}
    \lVert \nabla \mathcal{T}\varphi\rVert_{\mathcal{M}^{s,p,q}_{\B,\tt{or}}(\RR^n_+)} \lesssim_{n,s}^{\alpha}\lVert \varphi\rVert_{\mathcal{M}^{1+\alpha,1}_{\W}(\RR^{n-1})}.
\end{align*}
The case of spaces $\H^{s,p}$, $s\in(-1+\sfrac{1}{p},\sfrac{1}{p})$, follows from Proposition~\ref{prop:embeddingmultipliers} applied to the case $s\geqslant 0$ and a duality argument.

If $r=\infty$, note that $\varphi\in\mathcal{M}_{\W}^{1+\alpha,\infty}(\RR^{n-1})$, means $\nabla'\varphi\in \mathcal{M}_{\W,\tt{or}}^{\alpha,\infty}(\RR^{n-1})=\C^{0,\alpha}_b(\RR^{n-1})=\B^{\alpha}_{\infty,\infty}(\RR^{n-1})$. Following the second bullet point of Example~\ref{ex:"Conrete"VS"Abstract"MultiplierTheory}, this reduces to the case $r=1$.

\textbf{Step 2:} The case $r\in(1,\infty)$. Similarly, for $\varphi\in\mathcal{M}_{\W}^{1+\alpha,r}(\RR^{n-1})$, by Theorem~\ref{thm:LiftTraceMultiplier}, one obtains
\begin{align*}
    \nabla \mathcal{T}\varphi \in \mathcal{M}_{\W,\tt{or}}^{\alpha+\sfrac{1}{r},r}(\RR^{n}_+)
\end{align*}
with the estimate
\begin{align*}
    \lVert \nabla \mathcal{T}\varphi\rVert_{\mathcal{M}^{\alpha+\sfrac{1}{r},r}_{\W,\tt{or}}(\RR^n_+)} +\lVert \nabla \mathcal{T}\varphi\rVert_{\mathcal{M}^{\alpha+\sfrac{1}{r},r}_{\H,\tt{or}}(\RR^n_+)}\lesssim_{n}^{\alpha}\lVert \varphi\rVert_{\mathcal{M}^{1+\alpha,r}_{\W}(\RR^{n-1})}.
\end{align*}
By Proposition~\ref{prop:embeddingmultipliers} and Lemma~\ref{lem:EmbeddingMultipliersIntoMWs1}, one also obtains for all $\beta\in(0,\alpha)$, for all $p\in[1,r]$,
\begin{align*}
    \mathcal{M}^{\alpha,r}_{\W,\tt{or}}(\RR^{n-1}) \hookrightarrow \mathcal{M}^{\beta,p}_{\W,\tt{or}}(\RR^{n-1}).
\end{align*}
Hence, all $\beta\in(0,\alpha)$, for all $p\in[1,r]$,
\begin{align}\label{eq:ProofMultiplierEstLpRange p < r}
    \lVert \nabla \mathcal{T}\varphi\rVert_{\mathcal{M}^{\beta+\sfrac{1}{p},p}_{\W,\tt{or}}(\RR^n_+)}\lesssim_{n}^{\alpha}\lVert \varphi\rVert_{\mathcal{M}^{1+\alpha,r}_{\W}(\RR^{n-1})}.
\end{align}
By \cite[Proposition~3.5.3~\&~Theorem~9.2.1]{MazyaShaposhnikova2009}, one has the embeddings
\begin{align*}
    \mathcal{M}^{\alpha+\sfrac{1}{r},r}_{\H,\tt{or}}(\RR^n_+) \hookrightarrow \mathcal{M}^{s,q}_{\H,\tt{or}}(\RR^n_+)
\end{align*}
for all $0\leqslant s \leqslant \alpha + \sfrac{1}{r}$ and $p\in(1,\infty)$, such that either
\begin{enumerate}
    \item $\alpha r+1 >n$, and $sq\leqslant q \left( \alpha -\frac{n-1}{r}\right)+n$; or
    \item $\alpha r+ 1 \leqslant n$, and $sq<\alpha r +1$.
\end{enumerate}
Note that for \textit{(i)}, one always has $q \left( \alpha -\frac{n-1}{r}\right)+n>1$ due to $\alpha r+1 >n $, so that in any case it holds, for any $0<s<\frac{r\alpha+1}{q}$, any $q\in(r,\infty)$, that
\begin{align}\label{eq:ProofMultiplierEstLpRange q >r}
    \lVert \nabla \mathcal{T}\varphi\rVert_{\mathcal{M}^{s,q}_{\H,\tt{or}}(\RR^n_+)} \lesssim_{q,s,r}^{\alpha,n} \lVert \nabla \mathcal{T}\varphi\rVert_{\mathcal{M}^{\alpha+\frac{1}{r},r}_{\H,\tt{or}}(\RR^n_+)}\lesssim_{q,s,r}^{\alpha,n} \lVert  \varphi\rVert_{\mathcal{M}^{1+\alpha,r}_{\W}(\RR^{n-1})}.
\end{align}
Combining \eqref{eq:ProofMultiplierEstLpRange p < r} and \eqref{eq:ProofMultiplierEstLpRange q >r}, by Proposition~\ref{prop:embeddingmultipliers}, it holds for all $p\in(1,\infty)$, all $s\in[0,\sfrac{1}{p})$ that:
\begin{align*}
    \lVert \nabla \mathcal{T}\varphi\rVert_{\mathcal{M}^{s,p}_{\H,\tt{or}}(\RR^n_+)} \lesssim_{p,s,r}^{\alpha,n}\lVert  \varphi\rVert_{\mathcal{M}^{1+\alpha,r}_{\W}(\RR^{n-1})}.
\end{align*}
By duality, Proposition~\ref{prop:basicMultipliers} yields the case $s\in(-1+\sfrac{1}{p},0)$. Real interpolation yields for all $p\in(1,\infty)$, all $s\in(-1+\sfrac{1}{p},\sfrac{1}{p})$, and all $q\in[1,\infty]$, 
\begin{align*}
    \lVert \nabla \mathcal{T}\varphi\rVert_{\mathcal{M}^{s,p,q}_{\B,\tt{or}}(\RR^n_+)} \lesssim_{p,s,r}^{\alpha,n}\lVert  \varphi\rVert_{\mathcal{M}^{1+\alpha,r}_{\W}(\RR^{n-1})}.
\end{align*}
It remains to deal with the case of Besov spaces with $p=\infty$, which follows by duality from the case $p=1$.
\end{proof}

\begin{remark}\label{rem:BetterMultipliers}Actually, we did prove a slightly stronger statement, provided $\varphi\in\mathcal{M}^{1+\alpha,r}_{\X}(\RR^{n-1})$, for an arbitrary $r\in[1,\infty]$, and arbitrary $\alpha>0$, then that is $\nabla\mathcal{T}\varphi\in\mathcal{M}^{\alpha+\frac{1}{r},r}_{\X,\tt{or}}(\RR^{n-1})$ also belongs to the class of multipliers of \textbf{all} the function spaces in the \textbf{blue and green areas} of the figure below:
\begin{figure}[H]
\centering
\begin{tikzpicture}[yscale=3,xscale=7]
  \draw[->] (-0.1,0) -- (1.1,0) node[right,yshift=-2mm] {$1/p$};
  \draw[->] (0,-1.20) -- (0,1.15) node[above] {$s$};
  \draw[domain=0:1,smooth,variable=\x,blue] plot ({\x},{-1+\x}) node[right,yshift=2mm] {$s=-1+1/p$};
  \draw[domain=0:1,smooth,variable=\x,blue] plot ({\x},{\x}) node[right]{$s=1/p$};
  \draw[domain=0.2:1,line width=0.6mm,variable=\x,teal] plot ({\x},{\x+2/14}) ;
  \draw[domain=0:0.2,line width=0.6mm,variable=\x,teal] plot ({\x},{12*\x/7}) ;
  \draw[domain=0.8:1,line width=0.6mm,variable=\x,teal] plot ({\x},{12*(\x-1)/7}) ;
  \draw[domain=0:0.8,line width=0.6mm,variable=\x,teal] plot ({\x},{\x-2/14-1}) ;
  
  \draw[domain=0:1,dashed,line width=0.4mm,variable=\x,teal] plot ({\x},{2/14}) ;

  \fill[blue!30,opacity=0.3]  (0,0)-- (1,1)  -- (1,0) -- (0,-1) -- cycle;
  \fill[teal!30,opacity=0.5]  (0,0)-- (1/5,12/35)   -- (1,1+2/14)  -- (1,1)  -- cycle;
  \fill[teal!30,opacity=0.5]  (0,-1)-- (1,0)   -- (0.8,-12/35)  -- (0,-1-2/14)  -- cycle;

  \draw[dashed] (1,1) -- (0,1)  node[left] {$s=1$};

  \draw[dashed] (1/5,12/35) -- (1/5,0)  node[below] {$\boldsymbol{{\L}^{r}}$};
  
  \node[circle,fill,inner sep=1.5pt,label=below:$\L^2$] at (0.5,0) {};
  \node[circle,fill,inner sep=1.5pt,label=above:${\mathrm{H}}^1$] at (0.5,1) {};
  \node[circle,fill,inner sep=1.5pt,label=above:{\color{teal}$\boldsymbol{{\X}^{\alpha+{\sfrac{1}{r}},r}}$},teal] at (1/5,12/35) {};
  \node[circle,fill,inner sep=1.5pt,label=above right:{$\boldsymbol{{\X}^{\alpha,r}}$}] at (1/5,2/14) {};
  \node[circle,fill,inner sep=1.5pt,label=below right:{$\boldsymbol{{\X}^{-\alpha-\sfrac{1}{r},r'}}$}] at (0.8,-12/35) {};
  \node[circle,fill,inner sep=1.5pt,label=above:{$\boldsymbol{{\B}^{\alpha+1}_{1,1}={\W}^{\alpha+1,1}}$}] at (1,1+2/14) {};
  \node[circle,fill,inner sep=1.5pt,label=below left:{$\boldsymbol{{\L}^\infty}$}] at (0,0) {};
  \node[circle,fill,inner sep=1.5pt,label=below right:{$\boldsymbol{{\L}^1}$}] at (1,0) {};
  \node[circle,fill,inner sep=1.5pt,label=left:{$\boldsymbol{{\B}^{-\alpha-1}_{\infty,\infty}}$}] at (0,-1-2/14) {};
\end{tikzpicture}
\caption{Function spaces for which $\nabla\mathcal{T}\varphi$ has membership in the corresponding multiplier class. Borderline cases of regularity $s=\alpha+\sfrac{1}{p}$, $p< r$, and $s=\frac{r\alpha+1}{p}$, $p> r$, might not be reached. When $s\geqslant 0$, it corresponds to \eqref{eq:ProofMultiplierEstLpRange p < r} and \eqref{eq:ProofMultiplierEstLpRange q >r}. For negative $s$ this is just the dual range.}
\end{figure}
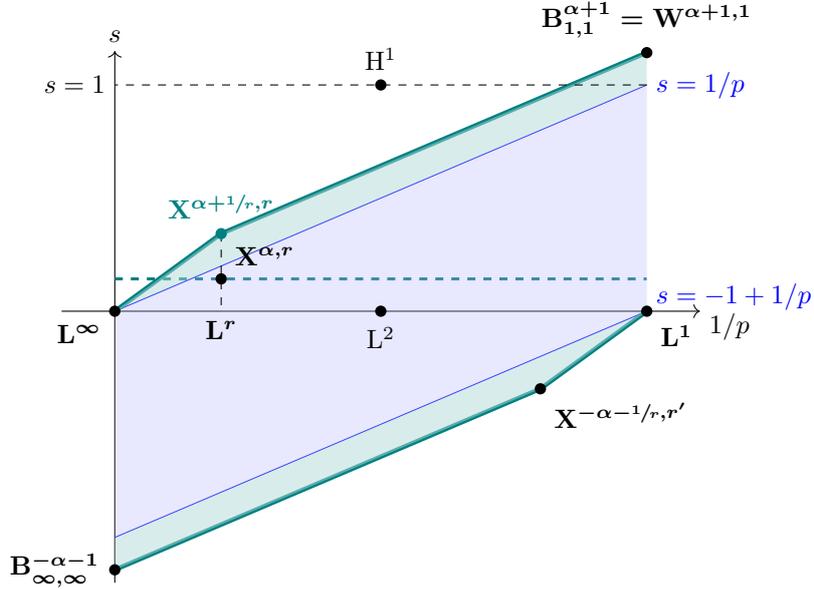

This suggests that given any domains with a boundary in a given multiplier class, one should be always allowed to ``perturb and propagate'' regularity beyond the sole class of function spaces associated to the chosen multiplier norm (\textit{i.e.} when dealing with a particular PDE, choosing one $r\in[1,\infty]$, so for a study \textit{a priori} on  a $\L^r$-type space, one could still obtain existence and uniqueness results in $\L^p$-type spaces for many $p\in[1,\infty]$). Mostly, interpolation techniques are then fully available.

It also tells us that minimal regularity assumptions that provide a margin to play with are given by the $\W^{s,1}$-multiplier norm for arbitrary small $s>0$.
\end{remark}

\begin{remark}\label{rem:MultiplierBoundarynotC1} We highlight now that, with such a theory, we are able to deal with domain with boundaries that might have a weaker smoothness requirement than $\C^1$ in some sense. Indeed, according to  \cite[Theorem~2.9]{Sickel1999} (we precise here $\W^{s,1}(\RR^n)=\mathrm{F}^{s}_{1,1}(\RR^n)$), having a domain $\Omega$ with local boundary charts $(\varphi_j)_{j\in\NN}\subset\mathcal{M}_{\W}^{1+\alpha,1}(\RR^{n-1})$ implies that $\varphi_j$ is globally Lipschitz and that the following Morrey-type condition is satisfied
\begin{align*}
    \sup_{\substack{\B_r\subset\RR^{n-1},\\ r\in(0,1]}} \frac{1}{r^{n-1-\alpha}} \int_{\B_r}\int_{\B_r} \frac{|\nabla'\varphi(x')-\nabla'\varphi(y')|}{|x'-y'|^{n-1+\alpha}}\,\d x' \,\d y'<\infty.
\end{align*}
Obviously, there exists continuously differentiable $\C^1$ functions such that this latter condition is not met.
\end{remark}

\subsubsection{Sufficient conditions for small multiplier norms}\label{Subsec:SmallMultNorm}

A key point in the strategy developed by Maz'ya and Shaposhnikova to perform analysis of elliptic PDEs in rough bounded domains over a prescribed function space $\X^{s,p}$ is about having small multiplier norms for the charts describing the boundary with respect to the trace space, the condition $\mathcal{M}^{s,p}_{\X_{\partial}}(\varepsilon)$ from Definition~\ref{def:besovboundary}, for small $\varepsilon>0$.

We provide here a statement giving sufficient conditions in terms of standard regularity up to increase the number of charts $(\varphi_{j})_{1\leqslant j\leqslant N}$. We did take advantage of Proposition~\ref{prop:embeddingmultipliers}, in order to transfer the result to Besov spaces $\B^{s}_{p,q}$, $p\neq q$, inducing the need of some arbitrarily small margin with respect to the regularity of the considered multiplier norm.

\begin{proposition}\label{prop:MspepsClassDomainsSufficientCdtn}Let $\alpha>0$, $\beta>\alpha$, $k\geqslant2$, $p,q\in[1,\infty]$. Let $\Omega$ be a bounded  Lipschitz domain. If
\begin{enumerate}
    \item $\Omega$ also belongs to the class $\mathbf{B}^{1+\alpha}_{r,p}$ or $\C^{\lfloor1+\beta \rfloor,\{\beta\}}$, then $\Omega$ is of class $\mathcal{M}^{1+\alpha,p}_{\W}(\varepsilon)$ for any $\varepsilon>0$, provided $r\in(\max(\frac{n-1}{\alpha},p),\infty]$;
    \item $\Omega$ also belongs to the class $\mathbf{B}^{1+\beta}_{r,q}$  or $\C^{\lfloor1+\beta \rfloor,\{\beta\}}$, then $\Omega$ is of class $\mathcal{M}^{1+\alpha,p,q}_{\B}(\varepsilon)$ for any $\varepsilon>0$, provided $r\in(\max(\frac{n-1}{\alpha},p),\infty]$;
    \item $\Omega$ also belongs to the class $\mathbf{B}^{k-{\sfrac{1}{p}}}_{\infty,p}$, then $\Omega$ is of class $\mathcal{M}^{k,p}_{\W}(\varepsilon)$ for any $\varepsilon>0$, whenever $1<p<\infty$;
    \item $\Omega$ also belongs to the class $\mathbf{W}^{1+\alpha,p}$, then $\Omega$ is of class $\mathcal{M}^{1+\alpha,p}_{\W}(\varepsilon)$ for any $\varepsilon>0$, whenever $1<p<\infty$, $\alpha>\sfrac{n-1}{p}$;
    \item $\Omega$ also belongs either to the class $\C^{\lfloor1+\beta \rfloor,\{\beta\}}$ or $\mathcal{C}^{\lfloor1+\alpha \rfloor,\{\alpha\}}$, then $\Omega$ is of class $\mathcal{M}^{1+\alpha,\infty,\infty}_{\B}(\varepsilon)$ for any $\varepsilon>0$.
\end{enumerate}
\end{proposition}

\begin{remark}\label{rem:ActuallyAlwayC1,a} Note that in Points \textit{(i)} and \textit{(ii)}, the condition on $r$ and $\alpha$, $r>\frac{n-1}{\alpha}$, is such that
\begin{align*}
     \mathbf{B}^{1+\alpha}_{r,p}\hookrightarrow \mathbf{B}^{1+\alpha - \frac{n-1}{r}}_{\infty,p} \hookrightarrow \C^{1,\alpha_0}.
\end{align*}
where $0<\alpha_0<\alpha-\frac{n-1}{r}$.
\end{remark}
We provide a short explanation for the content of Proposition~\ref{prop:MspepsClassDomainsSufficientCdtn}. Indeed, it is of special importance to us is describing the regularity of the boundary of a given domain by means of Sobolev multipliers. The multiplier classes gives minimal assumptions for the maximal (elliptic) regularity theory of PDEs \cite[Chapter 14]{MazyaShaposhnikova2009}. On the other hand, multiplier spaces are quite abstract. Hence, it is desirable to express by means of more customary function spaces. We will do this in the following
with a particular focus on those spaces appearing for the boundary regularity, provided a (bounded) domain $\Omega$ of $\RR^n$ with boundary $\partial\Omega$, the latter is described by local charts $\phi\,:\,\RR^{n-1}\longrightarrow\RR$.
As in \cite[Chapter 14]{MazyaShaposhnikova2009}, we are exhibiting here some of conditions that appears in Proposition~\ref{prop:MspepsClassDomainsSufficientCdtn} making a certain multiplier norm of $\phi$ small. In what follows below, we mention the following concordance of notations with respect to Proposition~\ref{prop:MspepsClassDomainsSufficientCdtn}
\begin{align*}
    s=1+\alpha,\quad \alpha>0.
\end{align*}
By \cite[Corollary 14.6.2]{MazyaShaposhnikova2009}, we have for $\phi\in \B^{s}_{\varrho,p}(\RR^{n-1})$ compactly supported and $\varepsilon\ll1$ that
\begin{align}\label{eq:MSa}
\|\phi\|_{\mathcal{M}_{\W}^{s,p}(\RR^{n-1})}\lesssim_{\|\phi\|_{\B^{s}_{\varrho,p}}}\varepsilon,
\end{align}
provided that $\|\nabla\phi\|_{L^{\infty}(\RR^{n-1})}\leqslant \varepsilon$, $s=l-{\sfrac{1}{p}}$ for some $l\in\NN$ and one of the following conditions holds:
\begin{itemize}
\item $p(s-1)<n-1$ and $\phi\in \B^{s}_{\varrho,p}(\RR^{n-1})$ with $\varrho\in\big[\frac{p(n-1)}{p(l-1)-1},\infty\big]$;
\item $p(s-1)=n-1$ and $\phi\in \B^{s}_{\varrho,p}(\RR^{n-1})$ with $\varrho\in(p,\infty]$.
\end{itemize}
Now, provided $s\geqslant 1$ ($s>1$ in the case of Besov spaces), by \cite[Corollary 4.3.8]{MazyaShaposhnikova2009} it holds
\begin{align}\label{eq:MSb}
\|\phi\|_{\mathcal{M}_{\W}^{s,p}(\RR^{m})}\sim_{p,s,m}
\|\nabla\phi\|_{\W^{s-1,p}(\RR^{m})} 
\end{align}
for $p(s-1)>m$. This is a simple consequence
of $\W^{s-1,p}(\RR^{m})$ being an algebra in this case. Thus, one obtains similarly
\begin{align}\label{eq:MSb'}
\|\phi\|_{\mathcal{M}_{\X}^{s,p}(\RR^{m})}\sim_{p,s,m}
\|\nabla\phi\|_{\X^{s-1,p}(\RR^{m})} 
\end{align}
for $p(s-1)>m$ and $\X\in\{\W,\,\B\}$.
If $\Omega\subset\RR^m$ is a Lipschitz domain with diameter $r$ and $ps>m$, we have
\begin{align}\label{eq:MSc}
\|\phi\|_{\mathcal{M}^{s,p}_{\X,\tt or}(\Omega)}\sim_{p,s,m}
r^{s-m/p}\|\phi\|_{\X^{s,p}(\Omega)}.
\end{align} 
This follows actually from a simple scaling argument since the left-hand side is scaling invariant (see \cite[Theorem~9.6.1,~\textit{(ii)}]{MazyaShaposhnikova2009} for the case $\X=\W$).
If $p(s-1)>n-1$,  it holds in accordance with \eqref{eq:MSc} that
\begin{align}\label{eq:0908b}
\partial\Omega\in \mathbf{X}^{s,p}\quad\Rightarrow\quad \partial\Omega \in \mathcal{M}_\X^{s,p}(\varepsilon)
\end{align}
for all $\varepsilon>0$ and $\X\in\{\W,\,\B\}$.

\medbreak

In order to take advantage of the theory of Elliptic Partial Differential Equations in smooth domains, in combination with the multiplier theory, the next result allows to consider regularized boundaries of a given Lipschitz domain, with boundary and domain constants that remain under control.

\medbreak

We mention explicitly here, and see for instance \cite[Lemma~2.1,~\textit{(ii)}]{LiSickelYangYuan2024}, that multiplier spaces for function spaces over $\RR^{n-1}$ are stable by convolution and have their own Young's inequality, that is
\begin{align}
    \lVert \eta\ast\varphi\rVert_{\mathcal{M}^{s,p}_{\X,\tt{or}}(\RR^{n-1})} \leqslant  \lVert \eta\rVert_{\L^1(\RR^{n-1})}\lVert \varphi\rVert_{\mathcal{M}^{s,p}_{\X,\tt{or}}(\RR^{n-1})}.\label{eq:mollifMultip}
\end{align}

\medbreak

The story of boundary regularization and approximation of Lipschitz domains by ones with smooth boundary in the interior seems to go back to Verchota, see \cite[Appendix]{VerchotaThesis}. We highlight here a very recent result that is more convenient with respect to our framework.

\medbreak

The boundary regularization of the result below being given by an explicit convolution in the proof as initially presented, we can then state it in a slightly more general version allowing Besov boundaries.   

\begin{theorem}[ {\cite[Theorem~1]{Antonini2024}} ]\label{thm:regularizeddomain}
    Let $\Omega\subset \RR^n$ be a bounded, Lipschitz domain with Lipschitz constant $L_\Omega$, and covering radius $R_\Omega$.
    \begin{enumerate}
        \item There exist sequences of bounded domains $\{\omega_m\},\{\Omega_m\}$, such that $\partial\omega_m\in \C^\infty,\,\partial\Omega_m\in \C^\infty$, and
    \begin{equation*}
        \omega_m\Subset\Omega\Subset \Omega_m\quad\text{for all $m\in \NN$.}
    \end{equation*}
    Their diameters satisfy
    \begin{equation}\label{diameters}
    \mathrm{diam}{(\Omega_m)}\leqslant c(n)\,\mathrm{diam}(\Omega)\,,\quad \mathrm{diam}(\omega_m)\leqslant c(n)\,\mathrm{diam}(\Omega)\,,
\end{equation}   
the following convergence property hold true  \begin{equation}\label{leb:distance}
        \lim_{m\to\infty} |\Omega_m\setminus\Omega|=0\,,\quad \lim_{m\to\infty}|\Omega\setminus \omega_m|=0\,,
    \end{equation} 
  the Hausdorff distances safisfy
\begin{equation}\label{hauss:dist}
        \mathrm{dist}_\H(\omega_m,\Omega)+\mathrm{dist}_\H(\Omega_m,\Omega)\leq\frac{12\,L_\Omega\sqrt{1+L_\Omega^2}}{m}\quad\text{for all $m\in \NN$,}
    \end{equation}
    and we may choose their Lipschitz characteristics $(L_{\Omega_m},R_{\Omega_m})$ and $(L_{\omega_m},R_{\omega_m})$ such that
    \begin{equation}\label{lip:car}
    \begin{split}
        L_{\Omega_m}\leqslant c(n)(1+L_\Omega^2)\,, & \quad   R_{\Omega_m}\geqslant R_\Omega/\big(c(n)(1+L_\Omega^2)\big) 
        \\
         L_{\omega_m}\leqslant c(n)(1+L_\Omega^2)\,, &\quad  R_{\omega_m}\geqslant R_\Omega/\big(c(n)(1+L_\Omega^2)\big)\,,\quad\text{for all $m\in \N$.}
        \end{split}
    \end{equation}
    Moreover, the smooth boundaries $\partial \omega_m, \partial \Omega_m$ are described with the help of the same co-ordinate 
    systems as $\partial \Omega$, i.e. there exists a finite number of local boundary charts $\{\phi^j\}_{j=1}^N,\{\psi^j_m\}_{j=1}^N$ and $\{\varphi^j_m\}_{j=1}^N$ which describe $\partial \Omega,\,\partial \Omega_m$ and $\partial \omega_m$ respectively, such that for each $j=1,\dots,N$ the functions $\psi^j_m,\varphi^j_m\in \C^\infty$ are defined on the same reference system as $\phi^j$, and 
    \begin{equation}\label{lchart:conv1}
        \psi^j_m\xrightarrow{m\to\infty} \phi^j\quad\text{and}\quad \varphi^j_m \xrightarrow{m\to\infty} \phi^j\quad\text{in $ \W^{1,p}(\B'_{R_\Omega-\varepsilon_0})$}\,,
    \end{equation}
    for all $p\in [1,\infty)$, for all $j=1,\dots,N$, and any fixed constant $\varepsilon_0\in (0,R_\Omega/2)$.
\item  If in addition $\partial \Omega\in \mathbf{B}^{s}_{p,q}$ (resp. $\boldsymbol{\BesSmo}^{s}_{p,\infty}$) for some $p\in[1,\infty]$, $q\in [1,\infty)$, $s\geqslant 1$, then
\begin{equation}\label{lchart:conv2}
    \psi^j_m\xrightarrow[m\to\infty]{} \phi^j\quad\text{and}\quad \varphi^j_m \xrightarrow[m\to\infty]{}\phi^j\quad\text{in $\B^{s}_{p,q}(\B'_{R_\Omega-\varepsilon_0})$ (resp. $\BesSmo^{s}_{p,\infty}(\B'_{R_\Omega-\varepsilon_0})$) }.
\end{equation}
A similar conclusion holds for $\partial \Omega\in\boldsymbol{\W}^{k,p}$, $\boldsymbol{\C}^{k}$, $\boldsymbol{\mathcal{C}}^{k,\alpha}$ and for $\boldsymbol{\C}^{k,\alpha}$, $p\in[1,\infty)$, $k\in\NN^\ast$, $\alpha\in(0,1]$, having only weak-$\ast$ convergence for the latter.
\end{enumerate}
\end{theorem}

\subsection{Differential forms in the context of fluid mechanics}\label{Subsec:DiffFormsIntro}

In the study of incompressible fluid dynamics, one may want to incorporate the "algebraic" condition of incompressibility for the fluid, namely
\begin{align*}
    \div u =0,
\end{align*}
in the function spaces we consider.

For the sake of completeness, and since it maybe of interest since it carries over additional geometric and algebraic properties, we will investigate instead the corresponding algebraic conditions at the level of differential forms
\begin{align*}
    \delta u = 0,
\end{align*}
which generalizes, and contains, the divergence free condition. The advantage of such generalized and abstract expression is that it can also carry over the behavior of curl-free functions and related function spaces in a single study.

Following \cite{McintoshMonniaux2018,Monniaux2021,Gaudin2023Hodge}, the complexified exterior algebra of $\RR^{n}$ is denoted by $\Lambda=\Lambda^{0} \oplus \Lambda^{1} \oplus \cdots \oplus \Lambda^{n}$ , where $\Lambda^k$ stands for all (up to complexification) alternating $k$-linear forms over $\RR^n$.

We also recall that for $k\in \llb 0,n\rrb$, $u \in \Lambda^k$ can be uniquely determined by $(u_{I})_{I\in\mathcal{I}^{n}_{k}}\in\CC^{\binom{n}{k}}$ such that
\begin{align*}
    u = \sum_{I\in\mathcal{I}^{n}_{k}} u_{I} \,\d x_{I} \text{,}
\end{align*}
where $\mathcal{I}^{n}_{k}=\{ (\ell_j)_{j\in\llb 1,k\rrb}\in\llb 1,n\rrb^k\,:\, \ell_j< \ell_{j+1}\}$, with $\lvert\mathcal{I}^{n}_{k}\rvert = {\binom{n}{k}}$, and $u_{I}$ and $\d x_{I}$ stands respectively for $u_{\ell_1 \ell_2\ldots \ell_k}$ and $\d x_{\ell_1}\wedge\d x_{\ell_2}\wedge\ldots\wedge\d x_{\ell_k}$ whenever $I=(\ell_j)_{j\in\llb 1,k\rrb}$, and the latter stands for the uniquely determined alternating $k$-linear form over $\RR^n$, such that
\begin{align*}
    \d x_{\ell_1}\wedge\d x_{\ell_2}\wedge\ldots\wedge\d x_{\ell_k}\left( \sum_{j=1}^n \alpha^{\ell_1}_{j}\mathfrak{e}_j,\, \sum_{j=1}^n \alpha^{\ell_1}_{j}\mathfrak{e}_j,\, \ldots,\, \sum_{j=1}^n \alpha^{\ell_k}_{j}\mathfrak{e}_j\right) = \prod_{m=1}^k \alpha^{\ell_m}_{\ell_m}.
\end{align*}
for all $(v_{\ell})_{\ell\in\llb1,k\rrb}\subset \RR^n$, with $v^\ell = \sum_{j=1}^n \alpha^{\ell}_{j}\mathfrak{e}_j$, $\alpha^{\ell}_{j}\in\RR$, $1\leqslant\ell\leqslant k$, $1 \leqslant j\leqslant n$.

One may also notice that such representation of $k$-differential forms with increasing index is possible due to symmetry properties (\textit{i.e.} $\d x_{\ell}\wedge\d x_{k} = - \d x_{k}\wedge\d x_{\ell}$ for all $k,\ell\in\llb 1,n\rrb$).

In particular, remark that $\Lambda^{0}\simeq \CC$, the space of complex scalars, and more generally $\Lambda^k\simeq \CC^{{\binom{n}{k}}}$, so that $\Lambda \simeq \CC^{2^n}$. We also set $\Lambda^{\ell}=\{0\}$ if $\ell<0$ or $\ell>n$.

On the exterior algebra $\Lambda$, the basic operations are
\begin{enumerate}[label=($\roman*$)]
    \item the exterior product $\wedge: \Lambda^{k} \times \Lambda^{\ell} \rightarrow \Lambda^{k+\ell}$,
    \item the interior product $\iprod: \Lambda^{k} \times \Lambda^{\ell} \rightarrow \Lambda^{\ell-k}$,
    \item the Hodge star operator $\star: \Lambda^{\ell} \rightarrow \Lambda^{n-\ell}$,
    \item the inner product $\langle\cdot, \cdot\rangle: \Lambda^{\ell} \times \Lambda^{\ell} \rightarrow \CC$.
\end{enumerate}
The operations detailed in \textit{(i)} and \textit{(ii)}, are naturally understood to be identically $0$ if one does have respectively $k+\ell\notin\llb 0,n\rrb$ and $\ell-k\notin\llb 0,n\rrb$.

If $a \in \Lambda^{1}, u \in \Lambda^{\ell}$ and $v \in \Lambda^{\ell+1}$, then
\begin{align*}
\langle a \wedge u, v\rangle=\langle u, a\iprod v\rangle.
\end{align*}
For more details, we refer to, \textit{e.g.}, \cite[Section~2]{AxelssonMcIntosh2004} and  \cite[Section~3]{CostabelMcIntosh2010}, noting that both papers contain some historical background (and being careful that ${\delta}$ has the opposite sign in \cite{AxelssonMcIntosh2004}). One may also consult \cite{DoCarmo1994} for an introduction from the euclidean setting point of view, and \cite[Section~1-3]{Jost2011} for basic and usual properties in the more general Riemannian setting\footnote{Notice that the Riemannian setting presented by Jost deals with compact manifold but a lot of computations remain true in their full generality, due to local behavior of each operation (Hodge-star operator, exterior and interior products etc.)}. 
In dimension $n=3$, this gives (see \cite{CostabelMcIntosh2010}) for a vector $a \in \RR^{3}$ identified with a $1$-form
\begin{itemize}
    \item   $u$ scalar, interpreted as 0 -form: $a \wedge u=u a$, $a\iprod u=0$;
    \item   $u$ scalar, interpreted as 3 -form: $a \wedge u=0$, $a\iprod u=u a$;
    \item   $u$ vector, interpreted as 1 -form: $a \wedge u=a \times u$, $a\iprod u=a \cdot u$;
    \item   $u$ vector, interpreted as 2 -form: $a \wedge u=a \cdot u$, $a\iprod u=-a \times u$.
\end{itemize}

From now and until the end of the present paper, if $p,q\in[1,\infty]$, $s\in\RR$, $k\in\llb 0,n \rrb$ and $\X^{s,p}\in\{\H^{s,p},\, \dot{\H}^{s,p},\, \B^{s}_{p,q},\, \dot{\B}^{s}_{p,q}\}$, then $\X(\Omega,\Lambda^k)$ stands for $k$-differential forms whose coefficients lie in $\X^{s,p}(\Omega,\CC)$, \textit{i.e.} $\X^{s,p}(\Omega,\Lambda^k)\simeq \X^{s,p}(\Omega,\CC^{\binom{n}{k}})$. One may also consider similarly $\X^{s,p}_0(\Omega,\Lambda^k)$.

We introduce the \textit{exterior derivative} $\d:=\nabla \wedge=\sum_{k=1}^{n} \partial_{x_k} \mathfrak{e}_{k} \wedge$ and the \textit{interior derivative} (or \textit{coderivative}) ${\delta}:=-\nabla\iprod=-\sum_{k=1}^{n} \partial_{x_k} \mathfrak{e}_{k}\iprod$ acting on differential forms on a domain $\Omega \subset \RR^{n}$, the operators $\d$ and ${\delta}$ are differential operators such that $\d^2 = \d \circ \d = 0$ and ${\delta}^2 = {\delta} \circ {\delta} = 0$, and each of them are bounded seen as linear operators $\X^{s,p}(\Omega,\Lambda)\longrightarrow \X^{s-1,p}(\Omega,\Lambda)$.

We allow us a slight abuse of notation: here we will not distinguish vectors of $\RR^n$, vector fields, and $1$-differential forms, and we recall the relation between $\d$ and ${\delta}$ via the Hodge star operator:
\begin{align*}
    \star {\delta} u=(-1)^{\ell} \d(\star u) \quad\text{ and }\quad \star \d u=(-1)^{\ell-1} {\delta}(\star u) \quad\text{ for an }\ell\text{-form }u.
\end{align*}

We recall the following integration by parts formula, for all $u,v\in \S({\overline{\Omega}},\Lambda)$,
\begin{align}
    \int_{\Omega} \langle \d u(x), v(x)\rangle \d x = \int_{\Omega} \langle u(x), {\delta} v(x)\rangle \d x + \int_{\partial\Omega} \langle u(x), \nu\iprod v(x)\rangle \d \sigma_x\text{,}\label{eq:IntbyParts1}\\
    \int_{\Omega} \langle {\delta} u(x), v(x)\rangle \d x = \int_{\Omega} \langle  u(x), \d v(x)\rangle \d x + \int_{\partial\Omega} \langle u(x), \nu\wedge v(x)\rangle \d \sigma_x\text{,}\label{eq:IntbyParts2}
\end{align}
which are true since we are in the Euclidean setting (otherwise the definition of $\delta$ has to be changed) and where $\nu$ is the  outer unit normal of $\Omega$ identified as a $1$-form. More generally, for all $T\in \mathcal{D}'(\Omega,\Lambda^k)\simeq \mathcal{D}'(\Omega,\CC^{\binom{n}{k}})$, $k\in\llb 0,n\rrb$, we define
    \begin{align*}
        \big\langle \d T, \phi \big\rangle_{\Omega} &:= \big\langle T, {\delta}\phi \big\rangle_{\Omega} \text{ for all }\phi \in \Ccinfty(\Omega,\Lambda^{k+1})\text{,}\\
        \big\langle {\delta}T, \psi \big\rangle_{\Omega} &:= \big\langle T,\d\psi \big\rangle_{\Omega} \text{ for all }\psi \in \Ccinfty(\Omega,\Lambda^{k-1})\text{.}
\end{align*}
\begin{definition}A distribution $u\in \Distrib(\Omega,\Lambda)$ is said to be
\begin{itemize}
    \item \textbf{closed} if $\d u =0$;
    \item \textbf{co-closed} if $\delta u =0$;
    \item \textbf{exact} if there exists $\omega\in\Distrib(\Omega,\Lambda)$ such that $u = \d \omega$;
    \item \textbf{co-exact} if there exists $\omega\in\Distrib(\Omega,\Lambda)$ such that $u = \delta \omega$.
\end{itemize}
\end{definition}

\medbreak

For arbitrary dimension $n$, the operator $\d$ restricted to its action on $\X^{s,p}(\Omega,\Lambda^1)$, with value in $\X^{s-1,p}(\Omega,\Lambda^2)$, and $\delta$ restricted to its action on $\X^{s,p}(\Omega,\Lambda^{n-1})$, with value in $\X^{s-1,p}(\Omega,\Lambda^{n-2})$, are fair consistent generalizations of the $\curl$ operator on $\RR^3$. Since in dimension $n$ higher than $4$, $n-1\neq 2$, we also have to distinguish their dual operators: the operator $\d$ restricted to its action on $\X^{s,p}(\Omega,\Lambda^{n-2})$ and the operator $\delta$ restricted to its action on $\X^{s,p}(\Omega,\Lambda^{2})$ which are fair consistent generalizations of the dual operator $\prescript{t}{}{\curl}$ (usually fully identified with the $\curl$ operator) on $\RR^3$.

Written differently, for $u\in\mathcal{D}'(\Omega,\Lambda)$, writing $ u =(u_0,\ldots, u_k,\ldots, u_n)\in\mathcal{D}'(\Omega,\oplus_{k=0}^n\Lambda^k)$,
\begin{enumerate}
    \item if $n=2$ :
    \begin{align*}
    \d u = (0,\, \nabla u_0,\, \curl u_1),\\
    \delta u = (-\div u_1,\, \nabla^\perp u_2,\, 0);
\end{align*}

    \item if $n=3$ :
    \begin{align*}
    \d u = (0,\, \nabla u_0,\, \curl u_1,\, -\div u_{2}),\\
    \delta u = (-\div u_1,\, {\curl}u_{2},\,  \nabla u_{3},\, 0);
\end{align*}

    \item if $n\geqslant 4$:
    \begin{align*}
    \d u = (0,\, \nabla u_0,\, \curl u_1,\, \ldots,\,  \prescript{t}{}{\curl}u_{n-2},\, -\div u_{n-1}),\\
    \delta u = (-\div u_1,\, \prescript{t}{}{\curl}u_{2},\, \ldots,\,  \curl u_{n-1},\,  \nabla u_{n},\, 0).
\end{align*}
\end{enumerate}

In particular, one may see those operators as unbounded ones and introduce their respective domains on $\L^p(\Omega, \Lambda^{k})$, $k\in\llb 0,n\rrb$, denoted by $\D_p(\d,\Lambda^k)$ and $\D_p({\delta},\Lambda^k)$ defined as
\begin{align*}
    \D_{p}(\d,\Lambda^k)&:=\left\{u \in \L^p(\Omega,\Lambda^k) \,\left|\, \d u \in \L^p(\Omega,\Lambda^{k+1})\right.\right\}\\\text { and } \D_{p}({\delta},\Lambda^k)&:=\left\{u \in \L^p(\Omega, \Lambda^k)\,\left|\, {\delta} u \in \L^p(\Omega, \Lambda^{k-1})\right.\right\},\\\text { and } \D_{p}(\underline{\d},\Lambda^k)&:=\left\{u \in \D_{p}(\d,\Lambda^k)\,\left|\, \nu\wedge u_{|_{\partial\Omega}} =0\right.\right\},\\\text { and } \D_{p}(\underline{\delta},\Lambda^k)&:=\left\{ u \in \D_{p}(\delta,\Lambda^k)\,\left|\, \nu\iprod u_{|_{\partial\Omega}} =0\right.\right\},
\end{align*}
as well as their ranges
\begin{align*}
    \R_{p}(\d,\Lambda^k)&:=\left\{v \in \L^p(\Omega,\Lambda^k) \,\left|\, v=\d u,\, u\in \D_{p}(\d,\Lambda^{k-1})\right.\right\}\\ \text { and } \R_{p}({\delta},\Lambda^k)&:=\left\{ v \in \L^p(\Omega,\Lambda^k) \,\left|\,v=\delta u,\, u\in \D_{p}(\d,\Lambda^{k+1})\right.\right\}.
\end{align*}
and similarly for $\underline{\d}$ and $\underline{\delta}$. The meaning of the boundary condition is understood in the distributional sense, see Appendix~\ref{App:Traces}.
Similar definitions are available for other inhomogeneous function spaces replacing $(\D_{p},\L^p)$ by either $(\D^{s}_p,\H^{s,p})$ or $(\D^{s}_{p,q},\B^{s}_{p,q})$ provided $p,q\in[1,\infty]$, $s\in(-1+{\sfrac{1}{p}},{\sfrac{1}{p}})$. When $\Omega$ is a special Lipschitz domain (like the half-space $\RR^n_+$) a similar definition is available for homogeneous function spaces replacing $(\D,\H,\B)$ by $(\dot{\D},\dot{\H},\dot{\B})$. The same goes for their respective range.

Before going further notice that at an early stage, the identification of the partial traces on the boundary  as $\nu\iprod u_{|_{\partial\Omega}}$ and $\nu\wedge u_{|_{\partial\Omega}}$ as in \eqref{eq:IntbyParts1} and \eqref{eq:IntbyParts2}, but for non-smooth differential forms is not entirely clear. See the first part of the proof of Theorem~\ref{thm:partialtracesDiffFormLip} for a proto-construction as abstract traces given by a distribution on the boundary. The full identification of partial traces as such is only available after Proposition~\ref{prop:DiffFormDomainsclosedEtc..} below is proved, which unlocks the second part of the proof of Theorem~\ref{thm:partialtracesDiffFormLip}. Before, having access to the definitive results, when it is relevant to make the distinction, we may write
\begin{align*}
    \mathfrak{n}_{\partial\Omega}(u), \qquad\text{ and }\qquad\mathfrak{t}_{\partial\Omega}(u) 
\end{align*}
for the distributions on the boundary that formally coincide respectively with  
\begin{align*}
    \nu\iprod u_{|_{\partial\Omega}}, \qquad\text{ and }\qquad\nu\wedge u_{|_{\partial\Omega}}.
\end{align*}

As a sequence of (densely defined but not necessarily closed) unbounded operators, we get:
\begin{center}
    \resizebox{\textwidth}{!}{
$\begin{array}{ccccccccccccccc}
     \d&:& \X^s(\Omega, \Lambda^0) & \overset{}{\longrightarrow} & \X^s(\Omega,\Lambda^1) &\overset{}{\longrightarrow} & \X^s(\Omega,\Lambda^2) & \overset{}{\longrightarrow} & \ldots &\overset{}{\longrightarrow} & \X^s(\Omega,\Lambda^{n-1})& \overset{}{\longrightarrow} & \X^s(\Omega,\Lambda^{n}) & {\longrightarrow} & 0\\
     0 & \longleftarrow & \X^s(\Omega, \Lambda^0) & \overset{}{\longleftarrow} & \X^s(\Omega,\Lambda^1) &\overset{}{\longleftarrow} & \X^s(\Omega,\Lambda^2) & \overset{}{\longleftarrow} & \ldots &\overset{}{\longleftarrow} & \X^s(\Omega,\Lambda^{n-1})& \overset{}{\longleftarrow} & \X^s(\Omega,\Lambda^{n}) &: &{\delta}\text{.}
\end{array}$}
\end{center}

We can use \eqref{eq:IntbyParts1} and \eqref{eq:IntbyParts2} to consider adjoints of $\d$ and ${\delta}$ in the sense of maximal adjoint operators in the Hilbert space $\L^2(\Omega,\Lambda)$, so that we will see later, that they have the following exact description of their domains
\begin{align*}
    &\D_{2}(\d^\ast,\Omega,\Lambda^k) = \{\, u\in \D_{2}({\delta},\Omega,\Lambda^k)\, \,:\, \nu\iprod u_{|_{\partial\Omega}} =0\, \} =\D_{2}(\underline{\delta},\Omega,\Lambda^k)\\
    \text{ and }\quad
    &\D_{2}({\delta}^\ast,\Omega,\Lambda^k) = \{\, u\in \D_{2}(\d,\Omega,\Lambda^k)\, \,:\, \nu\wedge u_{|_{\partial\Omega}} =0\, \}=\D_{2}(\underline{\d},\Omega,\Lambda^k)\text{.}
\end{align*}
One can also see these adjoint operators through the following $\L^2$-closures of unbounded operators,
\begin{align*}
    (\D_{2}(\d^\ast,\Omega,\Lambda^k),\d^\ast) = \overline{(\Ccinfty(\Omega,\Lambda^k),{\delta})}\quad\text{ and }\quad
    (\D_{2}({\delta}^\ast,\Omega,\Lambda^k),{\delta}^\ast) = \overline{(\Ccinfty(\Omega,\Lambda^k),\d)} \text{.}
\end{align*}
and similarly for other function spaces. See Proposition~\ref{prop:DiffFormDomainsclosedEtc..} below for the general statement.

\begin{proposition}\label{prop:Ext0PartialTraceVanish}Let $p,q\in[1,\infty]$, $s\in(-1+\sfrac{1}{p},\sfrac{1}{p})$. Let $\Omega$ be a bounded Lipschitz domain or special Lipschitz domain. The extension by $0$ from $\Omega$ to the whole space yields the canonical identification
\begin{align*}
    \{\,u\in\dot{\B}^{s}_{p,q}(\Omega,\Lambda)\,:\,\delta u\in\dot{\B}^{s}_{p,q}(\Omega,\Lambda)~\&~\nu\iprod u_{|_{\partial\Omega}}=0 \,\} \simeq \{\,u\in\dot{\B}^{s}_{p,q,0}(\Omega,\Lambda)\,:\,\delta u\in\dot{\B}^{s}_{p,q}(\RR^n,\Lambda)\,\}.
\end{align*}
The result still holds
\begin{itemize}
    \item for the Sobolev spaces $\dot{\H}^{s,p}$, assuming $1<p<\infty$;
    \item for the Besov spaces $\dot{\BesSmo}^{s}_{p,\infty}$, $\dot{\B}^{s,0}_{\infty,q}$;
    \item for the Lebesgue spaces $\L^1$ and $\L^\infty$;
    \item for the corresponding inhomogeneous function spaces.
\end{itemize}
A similar result holds for $(\d,\nu\wedge\cdot)$ instead of $(\delta,\nu\iprod\cdot)$.
\end{proposition}

\begin{remark}At this moment, even if one might acknowledge the existence of uniquely well defined partial trace $\mathfrak{n}_{\partial\Omega}(u)$ such that
\begin{align*}
    \langle \mathfrak{n}_{\partial\Omega}(u),\Psi_{|_{\partial\Omega}}\rangle_{\partial\Omega} = \langle u,\d\Psi\rangle_{\Omega}-\langle \delta u,\Psi\rangle_{\Omega},
\end{align*}
we insist that it is not clear at this stage that one can \textbf{truly} identify
\begin{align*}
    \mathfrak{n}_{\partial\Omega}(u) = \nu\iprod u_{|_{\partial\Omega}},
\end{align*}
and write $\nu\iprod u_{|_{\partial\Omega}}=0$. However, the notation $\nu\iprod u_{|_{\partial\Omega}}$ is already used here for better readability.
\end{remark}

\begin{proof}Let $u\in\dot{\B}^{s}_{p,q}(\Omega,\Lambda)$ such that $\delta u\in\dot{\B}^{s}_{p,q}(\Omega,\Lambda)$ and $\mathfrak{n}_{\partial\Omega}(u)=\nu\iprod u_{|_{\partial\Omega}}=0$. Since one has $(-1+\sfrac{1}{p},\sfrac{1}{p})$, this gives the canonical identification $\dot{\B}^{s}_{p,q}(\Omega,\Lambda)=\dot{\B}^{s}_{p,q,0}(\Omega,\Lambda)$, so that if one write $\widetilde{v}\in\dot{\B}^{s}_{p,q,0}(\Omega,\Lambda)$ the extension by $0$ to the whole space of $v\in\dot{\B}^{s}_{p,q}(\Omega,\Lambda)$, it holds for all $\varphi\in\Ccinfty(\RR^n,\Lambda)$,
\begin{align*}
    \big\langle \widetilde{u}, \d \varphi \big\rangle_{\RR^n} &= \big\langle {u}, \d \varphi \big\rangle_{\Omega}\\
    &=\big\langle \delta {u}, \varphi \big\rangle_{\Omega},\qquad\text{ (since }\quad\mathfrak{n}_{\partial\Omega}(u)=0)\\
    &=\big\langle \widetilde{\delta {u}}, \varphi \big\rangle_{\RR^n}.
\end{align*}
Therefore, it holds
\begin{align*}
    \delta\widetilde{u} = \widetilde{\delta {u}},\qquad\text{ in }\mathcal{D}'(\RR^n,\Lambda).
\end{align*}
By Proposition~\ref{prop:FundamentalExtby0HomFuncSpaces}, it gives the canonical embedding
\begin{align*}
    \{\,u\in\dot{\B}^{s}_{p,q}(\Omega,\Lambda)\,:\,\delta u\in\dot{\B}^{s}_{p,q}(\Omega,\Lambda)~\&~\nu\iprod u_{|_{\partial\Omega}}=0 \,\} \hookrightarrow \{\,u\in\dot{\B}^{s}_{p,q,0}(\Omega,\Lambda)\,:\,\delta u\in\dot{\B}^{s}_{p,q}(\RR^n,\Lambda)\,\}.
\end{align*}
The exact same computations, for $v\in \dot{\B}^{s}_{p,q,0}(\Omega,\Lambda)$ such that $\delta v\in\dot{\B}^{s}_{p,q}(\RR^n,\Lambda)$, for all $\varphi\in\Ccinfty(\RR^n,\Lambda)$
\begin{align*}
    \big\langle v_{|_{\Omega}}, \d \varphi \big\rangle_{\Omega} &= \big\langle v, \d \varphi \big\rangle_{\RR^n}\\
    &=\big\langle \delta v, \varphi \big\rangle_{\RR^n}\\
    &=\big\langle [\delta {v}]_{|_{\Omega}}, \varphi \big\rangle_{\Omega}.
\end{align*}
since $[\delta {v}]_{|_{\Omega}} = \delta v_{|_{\Omega}}$ in $\mathcal{D}'(\Omega,\Lambda)$, this yields the equality
\begin{align*}
    \big\langle v_{|_{\Omega}}, \d \varphi \big\rangle_{\Omega} &=\big\langle \delta {v}_{|_{\Omega}}, \varphi \big\rangle_{\Omega},\qquad\text{ for all } \varphi\in\Ccinfty(\overline{\Omega},\Lambda).
\end{align*}
This last line implies $\mathfrak{n}_{\partial\Omega}( v_{|_{\Omega}}) =0$, which gives the reverse embedding. The same proof applies to other function spaces.
\end{proof}

\begin{lemma}\label{lem:IntermediatedensityResultInhomogeneous}Let $\Omega$ be a bounded or special Lipschitz domain. Let $s\in\RR$, $p,q\in[1,\infty)$, with $m:=s$ if $s\in\NN$. For $\X^{s,p}\in\{\,\H^{s,p},\,\B^{s}_{p,q},\,\B^{s,0}_{\infty,q},\,\BesSmo^{s}_{p,\infty},\,\BesSmo^{s,0}_{\infty,\infty},\,\C^{m}_0,\,\W^{m,1}\}$, the space $\Ccinfty(\Omega,\Lambda)$ is a strongly dense subspace of
\begin{align*}
    \{\,u\in\X^{s,p}_0(\Omega,\Lambda)\,:\,\delta u\in\X^{s,p}(\RR^n,\Lambda)\,\}
\end{align*}
endowed with the graph norm. The result still holds with $\C^\infty_{ub,0}(\Omega,\Lambda)$ as a strong dense subspace if $\X^{s,p}=\B^{s}_{\infty,q}$ or $\BesSmo^{s}_{\infty,\infty}$. Finally, the result remains valid with $\d$ instead of $\delta$.
\end{lemma}

\begin{proof}By localisation and rotation, it suffices to prove it for a special Lipschitz domain $$\Omega=\{\,(x',x_n)\in\RR^{n-1}\times\RR\,:\,x_n>\phi(x')\}.$$ Let $u\in\X^{s,p}_0(\Omega)$ such that $\delta u \in\X^{s,p}_0(\Omega)$, by translation and convolution, for $t>0$ and $\varepsilon>0$ such that $\varepsilon<\frac{t}{2\big(1+\lVert\nabla'\phi\rVert_{\L^\infty(\RR^{n-1})}\big)}$, $u_{t}^{\varepsilon} := \varphi_{\varepsilon}\ast  u_t$ is smooth with support strictly contained in $\Omega$ and $u_{t}^{\varepsilon},\delta u_{t}^{\varepsilon}$ converges towards $u,\delta u$ in $\X^{s,p}(\RR^n)$, and then in $\X^{s,p}_0(\Omega)$. This gives the case $\B^{s}_{\infty,q}$ and $\BesSmo^{s}_{\infty,\infty}$ for $\C^\infty_{ub,0}$ functions. We consider $u^{R,\varepsilon}_t:=\Theta(\cdot/R)u^{\varepsilon}_t$, one has $\supp u^{R,\varepsilon}_t \subset \Omega\cap \overline{\B_{R}(0)}$, so that it is smooth with compact support. Finally,  taking the limit as $R$ goes to $\infty$, $u^{R,\varepsilon}_t$ and $\delta u^{R,\varepsilon}_t$ converge respectively to $u^{\varepsilon}_t$ and $\delta u^{\varepsilon}_t$. This gives the result.
\end{proof}

\begin{corollary}\label{cor:densityCcinfty} Let $p,q\in[1,\infty)$, $s\in(-1+\sfrac{1}{p},\sfrac{1}{p})$. Let $\Omega$ be a bounded or special Lipschitz domain. It holds that $\Ccinfty(\Omega,\Lambda)$ is a strong dense subspace of
\begin{align*}
    \D_{p,q}^{s}(\underline{\delta},\Omega,\Lambda)\text{, }\D_{\infty,q}^{s,0}(\underline{\delta},\Omega,\Lambda)\text{, }\mathcal{D}_{p,\infty}^{s}(\underline{\delta},\Omega,\Lambda)\text{, }\mathcal{D}_{\infty,\infty}^{s,0}(\underline{\delta},\Omega,\Lambda)\text{ and }\D_{p}^{s}(\underline{\delta},\Omega,\Lambda) 
\end{align*}
assuming $p>1$ or $[p=1$ and $s=0]$ for the latter. The result still holds if one replaces $\underline{\delta}$ by $\underline{\d}$.
\end{corollary}

\begin{proof} This is a direct consequence of Proposition~\ref{prop:Ext0PartialTraceVanish} and Lemma~\ref{lem:IntermediatedensityResultInhomogeneous}.
\end{proof}

\begin{remark}\label{rmk:DiffFormsanistropicEndBeginSpaces}The exterior and interior derivatives, as well as the Hodge-Dirac operators are \textit{a priori} unbounded operators only defined on the biggest space $\L^2(\Omega,\Lambda)$. However, throughout this paper we will use some semantic distortion referring to the interior or exterior derivatives as “unbounded operators on $\L^2(\Omega,\Lambda^k)$”, $k\in\llb0,n\rrb$, by their natural restriction to differential forms of degree $k$, such as
    \begin{align*}
        \d \,:\, \D_2(\d,\Lambda^k)\longrightarrow \L^2(\Omega,\Lambda^{k+1})\subset \L^2(\Omega,\Lambda),
    \end{align*}
    even if the range is not a subset of $\L^2(\Omega,\Lambda^k)$. This misuse will remain also for other function spaces that could replace $\L^2$. More generally, for $m\in\llb0,n\rrb$, for $0\leqslant k_0<k_1<\ldots<k_m\leqslant n$, we always use the canonical identification
    \begin{align*}
        \Ccinfty(\Omega,\Lambda^{k_0} \oplus\ldots\oplus \Lambda^{k_m}) \hookrightarrow \Ccinfty(\Omega,\Lambda) \text{ and }\mathcal{D}'(\Omega,\Lambda^{k_0} \oplus\ldots\oplus \Lambda^{k_m}) \hookrightarrow \mathcal{D}'(\Omega,\Lambda) 
    \end{align*}
    with identically zero coefficients on remaining indices.
\end{remark}

\begin{definition} \begin{enumerate}[label=($\roman*$)]
    \item The orthogonal projector defined on $\L^2(\Omega,\Lambda^k)$ onto $\N_2(\d^\ast,\Lambda^k)$ is denoted by $\mathbb{P}_\Omega$ and called the \textbf{generalized} \textbf{Hodge}\textbf{-}\textbf{Leray} (or \textbf{Leray}) \textbf{projector}.
    \item The (bounded) orthogonal projector defined on $\L^2(\Omega,\Lambda^k)$ onto $\N_2(\delta,\Lambda^k)$ is denoted by $\mathbb{Q}_\Omega$.
    \item For $p\in(1,\infty)$, $s\in\RR$ and $k\in\llb 0,n\rrb$, we say that $\H^{s,p}(\Omega,\Lambda^k)$ admits a \textbf{Hodge decomposition} if $(\D_p^s(\d,\Lambda^k),\d)$, $(\D_p^s(\delta,\Lambda^k),\delta)$ and their respective adjoints are closable and
    \begin{align*}
       \H^{s,p}(\Omega,\Lambda^k) &= \N_p^s(\mathfrak{d},\Lambda^k)\oplus \overline{\R_p^s(\mathfrak{d}^\ast,\Lambda^k)}\text{,}\\
        &=   \overline{\R_p^s(\mathfrak{d},\Lambda^k)} \oplus {\N_p^s(\mathfrak{d}^\ast,\Lambda^k)} \text{,}
\end{align*}
holds in the topological sense, provided $\mathfrak{d}\in\{\d,\delta\}$. We keep the same definition of the Hodge decomposition on other function spaces replacing $({\H}^{s,p},{\D}_p^s,{\R}_p^s,{\N}_p^s)$ by either $(\dot{\H}^{s,p},\dot{\D}_p^s,\dot{\R}_p^s,\dot{\N}_p^s)$, $({\B}^{s}_{p,q},{\D}_{p,q}^s,{\R}_{p,q}^s,{\N}_{p,q}^s)$, or even by $(\dot{\B}^{s}_{p,q},\dot{\D}_{p,q}^s,\dot{\R}_{p,q}^s,\dot{\N}_{p,q}^s)$, with $p\in[1,\infty]$, $q\in[1,\infty)$. When $q=\infty$, one considers weak-$\ast$ closure instead.
\end{enumerate}
\end{definition}

Following \cite[Sections~~2~\&~3]{McintoshMonniaux2018}, we state the following standard definitions and results.

\begin{definition} We set
\begin{itemize}
    \item The \textbf{Hodge-Dirac operator} on $\Omega$ with \textbf{tangential boundary conditions} is defined as
    \begin{align*}
        D_\mathfrak{t}:= \delta^\ast + \delta \text{.}
    \end{align*}
    Its square is the elliptic operator $-\Delta_{\mathcal{H}_\ast}:= D_\mathfrak{t}^2 = \delta^\ast\delta + \delta\delta^\ast$, is called the (negative) \textbf{Hodge Laplacian} with \textbf{relative boundary conditions} (also called \textbf{generalised Dirichlet boundary conditions})
    \begin{align*}
        \nu\wedge u _{|_{\partial \Omega}} = 0 \text{, and } \nu\wedge \delta u _{|_{\partial \Omega}} = 0\text{.}
    \end{align*}
    The restriction to scalar functions $u\,:\,\Omega\longrightarrow\Lambda^0$ gives $-\Delta_{\mathcal{H}_\ast}u = \delta\delta^\ast u = -\Delta_{\mathcal{D}}u$, where $-\Delta_{\mathcal{D}}$ is the Dirichlet Laplacian.
    \item The \textbf{Hodge-Dirac operator} on $\Omega$ with \textbf{normal boundary conditions} is defined as
    \begin{align*}
        D_\mathfrak{n}:= \d+\d^\ast \text{.}
    \end{align*}
    Its square is the elliptic operator $-\Delta_{\mathcal{H}}:=  \d\d^\ast + \d^\ast\d$, is called the (negative) \textbf{Hodge Laplacian} with \textbf{absolute boundary conditions} (also called \textbf{generalized Neumann boundary conditions})
    \begin{align*}
        \nu\iprod u _{|_{\partial \Omega}} = 0 \text{, and } \nu\iprod \d u _{|_{\partial \Omega}} = 0\text{.}
    \end{align*}
    The restriction to scalar functions $u\,:\,\Omega\longrightarrow\Lambda^0$ gives $-\Delta_{\mathcal{H}}u = \d^\ast\d u = -\Delta_{\mathcal{N}}u$,  where $-\Delta_{\mathcal{N}}$ is the Neumann Laplacian.
\end{itemize}
\end{definition}

\begin{proposition}\label{prop:HodgeDiracL2bdd}The Hodge-Dirac operator $(\D_2(D_\cdot,\Omega,\Lambda),D_\cdot)$  is a self-adjoint $0$-bisectorial operator on $\L^2(\Omega,\Lambda)$ so that it satisfies the following bound, for all $\theta\in (0,\frac{\pi}{2})$,
\begin{align}\label{eq:resolvEstDiracL2}
    \forall \lambda\in \mathrm{S}_\theta\text{, } \left\lVert \lambda(\lambda\I+i D_\cdot)^{-1}\right\rVert_{\L^2(\Omega)\rightarrow \L^2(\Omega)}\leqslant \frac{1}{\cos(\theta)} \text{.}
\end{align}
Moreover, for all $\mu\in(0,\frac{\pi}{2})$, it admits a bounded $\mathbf{H}^\infty(\mathrm{S}_\mu)$-functional calculus on $\L^2(\Omega,\Lambda)$ with bound $1$, \textit{i.e.} for all $f\in\mathbf{H}^\infty(\mathrm{S}_\mu)$, all $u\in\L^2(\Omega,\Lambda)$,
\begin{align}
    \left\lVert f(D_\cdot)u\right\rVert_{\L^2(\Omega)}\leqslant \left\lVert f\right\rVert_{\L^\infty(\mathrm{S}_\mu)}\lVert u\rVert_{\L^2(\Omega)}\text{.}
\end{align}
\end{proposition}

Actually, as a square of the Hodge-Dirac operator for $k\in\llb 0,n\rrb$, the Hodge Laplacian $(\D(\Delta_\mathcal{H},\Omega,\Lambda^k),-\Delta_\mathcal{H})$ can be realized on $\L^2(\Omega,\Lambda^k)$ by means of a densely defined, symmetric, accretive, continuous, closed, sesquilinear form on $\L^2(\Omega,\Lambda^k)$:
\begin{align}\label{eq:sesqlinformLaplacianHodge}
    \mathfrak{a}_\mathcal{H}\,:\, \D_2(\mathfrak{a}_\mathcal{H},\Lambda^k)^2\ni(u,v) \longmapsto \int_{\Omega}\langle \d u(x), {\d v(x)}\rangle \mathrm{~d}x + \int_{\Omega}\langle \delta u(x), {\delta v(x)}\rangle\mathrm{~d}x
\end{align}
with $\D_2(\mathfrak{a}_{\mathcal{H}},\Lambda^k)= \D_2(D_\cdot,\Omega,\Lambda^k)$. It is straightforward to see that the Hodge Laplacian is then a closed, densely defined, non-negative self-adjoint operator on $\L^2(\Omega,\Lambda^k)$. See \cite[Chapter~1]{bookOuhabaz2005} for more details about realization of operators via sesquilinear forms on a Hilbert space.

\begin{remark} We recall here that, if $\Omega\subset \RR^3$ is an open set with, say at least, Lipschitz boundary, one has formally for $u$ with value in $\Lambda^1\simeq \CC^3$ or $\Lambda^2\simeq \CC^3$,
\begin{align*}
    -\Delta_{\mathcal{H}}u = \curl \curl u - \nabla \div u
\end{align*}
with either one of the following couple of boundary conditions
\begin{align*}
    [ u\cdot\nu_{|_{\partial\Omega}} = 0 \text{, } \nu\times\curl u_{|_{\partial\Omega}} = 0]&\text{ or }[ u\times\nu_{|_{\partial\Omega}} = 0 \text{, } (\div u)\nu_{|_{\partial\Omega}} = 0] \text{.}
\end{align*}
\end{remark}

\begin{proposition}\label{prop:FullHodgedecompL2bdd}Let $\Omega$ be a bounded Lipschitz domain and $k\in\llb0,n\rrb$. 
\begin{enumerate}
    \item One has the direct orthogonal decomposition
\begin{align*}
    \L^2(\Omega,\Lambda^k)= \overline{\R_2(\d,\Lambda^k)}\oplus \overline{\R_2(\d^\ast,\Lambda^k)}\oplus (\N_2(\d,\Lambda^k)\cap \N_2(\d^\ast,\Lambda^k)),
\end{align*}
\item There exists a bounded linear operator $\mathcal{R}^{\sigma}\,:\,\L^2(\Omega,\Lambda^{k})\longrightarrow \R_{2}(\d,\Omega,\Lambda^{k+1})$, such that for for all $(v,w)\in \overline{\R_2(\d^\ast,\Lambda^k)} \times[  \overline{\R_2(\d,\Lambda^k)} \oplus (\N_2(\d,\Lambda^k)\cap \N_2(\d^\ast,\Lambda^k))]$
\begin{align*}
     \underline{\delta} \mathcal{R}^{\sigma} v =v, \quad  \text{ and }\quad \mathcal{R}^{\sigma} w=0.
\end{align*}
Additionally, $ \underline{\delta} \mathcal{R}^{\sigma}$ (\textit{resp.} $\mathcal{R}^{\sigma}\underline{\delta}$) extends uniquely as the unique orthogonal projection onto $\overline{\R_2(\d^\ast,\Lambda^k)}$ (\textit{resp.} onto $\overline{\R_2(\d,\Lambda^k)}$), with
\begin{align*}
    \mathcal{R}^{\sigma} \D_2(\underline{\delta},\Lambda) \subset  \D_2(\underline{\delta},\Lambda).
\end{align*}
\item Up to consider the extension from the dense subspaces $\D_2(\underline{\delta},\Lambda^k)=\D_2(\d^\ast,\Lambda^k)$ of $\L^{2}(\Omega,\Lambda^k)$, one has the following valid identity on $\L^{2}(\Omega,\Lambda^k)$
\begin{align*}
    \I=  \underline{\delta}\mathcal{R}^{\sigma} + \mathcal{R}^{\sigma} \underline{\delta} + \P_{\H},
\end{align*}
where $\P_{\H}$ is the orthogonal projection onto $\N_2(\d,\Lambda^k)\cap \N_2(\d^\ast,\Lambda^k)=\N_2(D_\mathfrak{n},\Lambda^k)$, the latter being a \textbf{finite dimensional subspace}.
\end{enumerate}
A similar result holds if one replaces $(\delta,\d,\sigma)$ by $(\d,\delta,\gamma)$ in which case one denotes by $\P_{\H_\ast}$ the counterparts of $\P_\H$.
\end{proposition}

\begin{remark}Notice that the uniqueness of orthogonal projectors onto closed subspace of a Hilbert space yields the equality $\PP_{\Omega} = \underline{\delta}\mathcal{R}^{\sigma}  + \P_{\H}$.
\end{remark}

\begin{lemma}\label{lem:approxHodgeLeray} Let $\Omega$ be a Lipschitz domain. For any $f\in\L^2(\Omega,\Lambda)$, the following limit holds strongly in $\L^2$:
\begin{align*}
    \PP_\Omega f = \P_\H f +  \lim_{t \rightarrow 0}\, \delta(it+D_\mathfrak{n})^{-1}\d (-it+D_\mathfrak{n})^{-1} [\I-\P_\H] f.
\end{align*}
\end{lemma}

\begin{proof}Up to consider $[\I-\P_\H] f$, we can assume $\P_\H f =0$. This means, we suppose $f\in\overline{\R_2(D_\mathfrak{n},\Omega,\Lambda)}$. Since $D_\mathfrak{n}$ is bi-sectorial by \cite[Proposition~3.2.2]{EgertPhDThesis2015}, $\d^2=0$ and $\delta^2=0$, one obtains
\begin{align*}
    \delta(it+D_\mathfrak{n})^{-1}\d (-it+D_\mathfrak{n})^{-1}f + \d(it+D_\mathfrak{n})^{-1}&\delta (-it+D_\mathfrak{n})^{-1}f \\
    &= D_\mathfrak{n}(it+D_\mathfrak{n})^{-1}D_\mathfrak{n} (-it+D_\mathfrak{n})^{-1} f \xrightarrow[t\rightarrow 0]{} f.
\end{align*}
However, one has
\begin{align*}
    \delta(it+D_\mathfrak{n})^{-1}\d (-it+D_\mathfrak{n})^{-1}f &= \PP_{\Omega}[ \delta(it+D_\mathfrak{n})^{-1}\d (-it+D_\mathfrak{n})^{-1} \d(it+D_\mathfrak{n})^{-1}\delta (-it+D_\mathfrak{n})^{-1} ]f.
\end{align*}
Therefore by continuity of $\PP_\Omega$,
\begin{align*}
    \delta(it+D_\mathfrak{n})^{-1}\d (-it+D_\mathfrak{n})^{-1}f \xrightarrow[t\rightarrow 0]{} \PP_{\Omega} f,
\end{align*}
which is the desired result.
\end{proof}

\begin{proposition}\label{prop:DiffFormDomainsclosedEtc..} Let $p,p_0,p_1,q,q_0,q_1\in[1,\infty]$ and $-1+{\sfrac{1}{p}}<s_0,s,s_1<{\sfrac{1}{p}}$, and $s=(1-\theta)s_0+\theta s_1$ and $\sfrac{1}{r} = (1-\theta)/r_0+\theta /r_1$, $r\in\{p,q\}$. Let $\Omega$ be either a bounded or special Lipschitz domain. The following hold
\begin{enumerate}
    \item The operators $\d$, $\delta$, $\underline{\d}$ and  $\underline{\delta}$ are closed, if defined as unbounded operators on respectively $\dot{\B}^{s}_{p,q}(\Omega,\Lambda)$, $\dot{\BesSmo}^{s}_{p,\infty}(\Omega,\Lambda)$, $\dot{\B}^{s,0}_{\infty,q}(\Omega,\Lambda)$, $\dot{\BesSmo}^{s,0}_{\infty,\infty}(\Omega,\Lambda)$,  $\dot{\BesSmo}^{s}_{p,\infty}(\Omega,\Lambda)$ and $\dot{\H}^{s,p}(\Omega,\Lambda)$  requiring additionally $1<p<\infty$ for the latter.
    
    If $q<\infty$, each operator is densely defined.

    \item The dual operators of $\d$, $\delta$, $\underline{\d}$ and  $\underline{\delta}$ defined as unbounded operators on $\dot{\X}^{s,p}(\Omega,\Lambda)$, are respectively $\underline{\delta}$, $\underline{\d}$, $\delta$ and $\d$ defined as unbounded operators on $\dot{\X}'{}^{-s,p'}(\Omega,\Lambda)$.

    \item For $\mathfrak{d}\in \{\d,\delta,\underline{\d},\underline{\delta}\}$, defined on $\dot{\X}^{s,p}(\Omega,\Lambda)$ with domain $\dot{\D}^{s,p}(\mathfrak{d},\Omega,\Lambda)$ one has the interpolation identities
    \begin{enumerate}
        \item $(\dot{\D}^{s_0,p}(\mathfrak{d},\Omega,\Lambda),\dot{\D}^{s_1,p}(\mathfrak{d},\Omega,\Lambda))_{\theta,q}= \dot{\D}^{s}_{p,q}(\mathfrak{d},\Omega,\Lambda)$, $s_0\neq s_1$;
        \item $[\dot{\D}^{s_0,p_0}(\mathfrak{d},\Omega,\Lambda),\dot{\D}^{s_1,p_1}(\mathfrak{d},\Omega,\Lambda)]_{\theta}= \dot{\D}^{s,p}(\mathfrak{d},\Omega,\Lambda)$, $q<\infty$.
    \end{enumerate}

    \item  If $\Omega$ is bounded, the operators $\d$, $\delta$, $\underline{\d}$ and  $\underline{\delta}$ are such that their range are closed on $\dot{\B}^{s}_{p,q}(\Omega,\Lambda)$, $\dot{\BesSmo}^{s}_{p,\infty}(\Omega,\Lambda)$, and ${\H}^{s,p}(\Omega,\Lambda)$ the latter requiring $1<p<\infty$. 
    
    Additionally for $\mathfrak{d}\in \{\d,\delta\}$, defined on ${\X}^{s,p}(\Omega,\Lambda)$ with range ${\R}^{s,p}(\mathfrak{d},\Omega,\Lambda)$ one has
    \begin{align*}
        \mathfrak{d}{\X}^{s+1,p}(\Omega,\Lambda) = {\R}^{s,p}(\mathfrak{d},\Omega,\Lambda),\qquad\text{ and }\qquad \underline{\mathfrak{d}}{\X}^{s+1,p}_\mathcal{D}(\Omega,\Lambda) = {\R}^{s,p}(\underline{\mathfrak{d}},\Omega,\Lambda).
    \end{align*}
\end{enumerate}
\end{proposition}

\begin{proof} The point \textit{(i)} is standard.

\textbf{Point \textit{(ii)}:} Let $1\leqslant p\leqslant \infty$, $s\in (-1+\sfrac{1}{p},\sfrac{1}{p})$. Let $u\in \dot{\D}_{p',\infty}^{-s}(\underline{\delta}^\ast,\Omega,\Lambda)$, for all $\varphi\in \dot{\D}_{p,1}^s(\underline{\delta},\Omega,\Lambda)$,
\begin{align*}
    \big\langle u, \underline{\delta} \varphi \big\rangle_{\Omega}=\big\langle \underline{\delta}^\ast u,  \varphi \big\rangle_{\Omega}.
\end{align*}
Taking $\varphi\in\Ccinfty(\Omega,\Lambda)$,
\begin{align*}
    \big\langle \underline{\delta}^\ast u,  \varphi \big\rangle_{\Omega}=\big\langle u, \underline{\delta} \varphi \big\rangle_{\Omega}=\big\langle u, {\delta} \varphi \big\rangle_{\Omega}=\big\langle \d u,  \varphi \big\rangle_{\Omega}.
\end{align*}
Therefore $\underline{\delta}^\ast u =\d u$ in $\mathcal{D}'(\Omega,\Lambda)$. Thus, one obtains $(\D_{p',\infty}^{-s}(\underline{\delta}^\ast,\Omega,\Lambda),\delta^\ast) \subset (\D_{p',\infty}^{-s}(\d,\Omega,\Lambda),\d)$. For the converse inclusion, let $u\in\dot{\D}_{p',\infty}^{-s}(\d,\Omega,\Lambda)$, for all $\varphi\in \dot{\D}_{p,1}^s(\underline{\delta},\Omega,\Lambda)$, since $\mathfrak{n}_{\partial\Omega}(u)=0$
\begin{align*}
    \big\langle u, \underline{\delta} \varphi \big\rangle_{\Omega}=\big\langle \d u,  \varphi \big\rangle_{\Omega}.
\end{align*}
Thus $\varphi\longmapsto \big\langle u, \underline{\delta} \varphi \big\rangle_{\Omega}$ extends as a bounded linear functional on $\dot{\B}^{s}_{p,1}(\Omega,\Lambda)$ since $\d u \in \dot{\B}^{-s}_{p',\infty}(\Omega,\Lambda)$. This yields the reverse inclusion  $(\D_{p',\infty}^{-s}(\d,\Omega,\Lambda),\d) \subset (\D_{p',\infty}^{-s}(\underline{\delta}^\ast,\Omega,\Lambda),\delta^\ast)$. The same goes for $\underline{\d}$, $\d$ and $\delta$, and for other function spaces.

\textbf{Points \textit{(iii)} \& \textit{(iv)}:} We consider the Bogovski\u{\i} and Poincaré operators from Theorems~\ref{thm:PotentialOpBddLipDom} and \ref{thm:PotentialOpSpeLipDom}. We do the proof for the domains of $\delta$ on homogeneous Sobolev and Besov spaces on special Lipschitz domains, the proof for other operators and on bounded Lipschitz domains is similar. We have the bounded linear maps
\begin{align*}
    \mathfrak{E}\,:\,\dot{\D}^{s,p}(\delta,\Omega,\Lambda) &\,\longrightarrow \dot{\X}^{s,p}(\Omega,\Lambda)\times\dot{\X}^{s,p}\cap\dot{\X}^{s+1,p}(\Omega,\Lambda)\\
    u &\,\longmapsto \qquad(\delta \T^{\sigma}u, \T^{\sigma}\delta u).
\end{align*}
and
\begin{align*}
    \mathfrak{R}\,:\,\dot{\X}^{s,p}(\Omega,\Lambda)\times\dot{\X}^{s,p}\cap\dot{\X}^{s+1,p}(\Omega,\Lambda) &\,\longrightarrow \dot{\D}^{s,p}(\delta,\Omega,\Lambda)\\
    (v,w)\qquad \qquad \qquad \qquad &\,\longmapsto \delta \T^{\sigma}v+w
\end{align*}
with
\begin{align*}
    \mathfrak{R}\mathfrak{E} = \I,\qquad \text{ on } \dot{\D}^{s,p}(\delta,\Omega,\Lambda).
\end{align*}
Since $(\dot{\X}^{s,p}(\Omega,\Lambda))_{s,p}$ and $(\dot{\X}^{s,p}\cap\dot{\X}^{s+1,p}(\Omega,\Lambda))_{s,p}$, $-1+1/p<s<1/p$, are interpolation scales, thanks to Proposition~\ref{prop:InterpHybridBesov} and Corollary~\ref{cor:InterpHybridBesovSpeLip}, this allows to transfer the result by a retraction and co-retraction argument. The case of the range of operators (on bounded domains) and of null spaces also admits a similar and simpler proof.

When $\Omega$ is bounded closedness of the range of $\d$, $\delta$, $\underline{\d}$ and $\underline{\delta}$ is provided by Theorem~\ref{thm:PotentialOpBddLipDom} giving bounded linear maps
\begin{align*}
    \mathcal{B}^{\gamma}&\,:\,\X^{s,p}(\Omega,\Lambda)\longrightarrow \X^{s+1,p}(\Omega,\Lambda);\\
    \mathcal{B}^{\sigma}&\,:\,\X^{s,p}(\Omega,\Lambda)\longrightarrow \X^{s+1,p}(\Omega,\Lambda);\\
    \mathcal{B}^{0,\gamma}&\,:\,\X^{s,p}_0(\Omega,\Lambda)\longrightarrow \X^{s+1,p}_0(\Omega,\Lambda);\\
    \mathcal{B}^{0,\sigma}&\,:\,\X^{s,p}_0(\Omega,\Lambda)\longrightarrow \X^{s+1,p}_0(\Omega,\Lambda):
\end{align*}
so that
\begin{align*}
    \d\mathcal{B}^{\gamma}u &=\d\mathcal{B}^{\gamma} \d v = v,\,\quad \forall u =\d v \in\R^{s,p}(\d,\Omega,\Lambda);\\
    \delta\mathcal{B}^{\sigma}u &=\delta\mathcal{B}^{\gamma} \delta v = v,\,\quad \forall u =\delta v \in\R^{s,p}(\delta,\Omega,\Lambda);\\
    \underline{\d}[\mathcal{B}^{0,\gamma}\tilde{u}]_{|_{\Omega}} &=[\d\mathcal{B}^{\gamma} {\d} \tilde{v}]_{|_{\Omega}} = v,\,\quad \forall u =\underline{\d} v \in\R^{s,p}(\underline{\d},\Omega,\Lambda);\\
    \underline{\delta}[\mathcal{B}^{0,\sigma}\tilde{u}]_{|_{\Omega}} &=[\delta \mathcal{B}^{\gamma} {\delta} \tilde{v}]_{|_{\Omega}} = v,\,\quad \forall u =\underline{\delta} v \in\R^{s,p}(\underline{\delta},\Omega,\Lambda),
\end{align*}
where the two last lines took twice advantages of Lemma~\ref{lem:IntermediatedensityResultInhomogeneous}. This shows the closedness of the range, due to the mapping properties as bounded operators.
\end{proof}

\paragraph{Further information: pullback, pushforward and Leibniz rules}

\begin{definition}
\label{def:changevar}  
Let $\Omega$ and $\Omega'$ be open sets of ${\RR}^n$ and
$\mathbf{\Psi}:\Omega\to\Omega'$ a bi-Lipschitz transformation. 
Denote by $\nabla{\Psi}(y)$ the Jacobian matrix of $\mathbf{\Psi}$ at a point $y\in \Omega$ 
and extend it, for almost every $y\in\Omega$, to an isomorphism $\nabla{\mathbf{\Psi}}(y):\Lambda^k\to\Lambda^k$ such that
$$
{\nabla{\Psi}(y)}(\d{x_{i_1}}\wedge\dots\wedge \d{x_{i_k}})=(\nabla{\Psi}(y)\d{x_{i_1}})\wedge\ldots
\wedge(\nabla\mathbf{\Psi}(y)\d{x_{i_k}}), \quad (i_\ell)_{\ell\in\llb 1,k\rrb}\subset\llb 1,n\rrb.
$$
The {\textbf{pullback}} of a differential form $u:\Omega'\to \Lambda$ is denoted by  
$\mathbf{\Psi}^*u:\Omega\to \Lambda$, the {\textbf{push forward}} of a differential form 
$f:\Omega\to\Lambda$ by $\mathbf{\Psi}_*f:\Omega'\to\Lambda$ and 
$\tilde{\mathbf{\Psi}}_*^{-1}u:=\det(\nabla \mathbf{\Psi})\mathbf{\Psi}_*^{-1}u:\Omega\to\Lambda$ are 
defined by
$$
(\mathbf{\Psi}^*u)(y):=\prescript{t}{}{[\nabla \mathbf{\Psi}(y)]}\bigl(u\circ\mathbf{\Psi}(y)\bigr)\quad\mbox{and}\quad
(\mathbf{\Psi}_*^{-1}u)(y):=[\nabla \mathbf{\Psi}(y)]^{-1}\bigl(u\circ\mathbf{\Psi}(y)\bigr),
\quad y\in \Omega',
$$
and where $\det(\nabla \mathbf{\Psi})(y)$ denotes the Jacobian determinant of $\mathbf{\Psi}$ at a
point $y\in \Omega$.
\end{definition}

\begin{remark}
\label{rem:bounds}
Note that for all $p\in[1,\infty]$, 
$\mathbf{\Psi}^*: \L^p(\Omega',\Lambda)\to \L^p(\Omega,\Lambda)$ and
$(\mathbf{\Psi}_*)^{-1}:\L^p(\Omega',\Lambda)\to \L^p(\Omega,\Lambda)$ are 
bounded with norms controlled by 
$\displaystyle{\lVert \nabla{\mathbf{\Psi}}\rVert_{\L^\infty(\Omega)}}$ and
$\displaystyle{\lVert (\nabla{\mathbf{\Psi}})^{-1}\rVert_{\L^\infty(\Omega)}}$, 
and hence by the Lipschitz constants of $\mathbf{\Psi}$ and $\mathbf{\Psi}^{-1}$.
\end{remark} 

\begin{remark}
\label{rem:drho=rhod}
For $\mathbf{\Psi}$ as in Definition~\ref{def:changevar} and a differential form
$u:\Omega'\to\Lambda$ the following commutation properties hold:
\begin{equation}
\label{eq:com}
\d(\mathbf{\Psi}^*u)=\mathbf{\Psi}^*(\d u)\quad\mbox{and}\quad
\delta(\tilde{\mathbf{\Psi}}_*^{-1}u)=\tilde{\mathbf{\Psi}}_*^{-1}(\delta u).
\end{equation}
In particular, if $u\in {\D}_p(\d,\Omega')$, then 
$\mathbf{\Psi}^*u\in{\D}_p(\d,\Omega)$ and if $u\in {\D}_p(\delta,\Omega')$, then 
$\tilde{\mathbf{\Psi}}_*^{-1}u\in{\D}_p(\delta,\Omega)$.

We also have the following homomorphism properties:
$$
\begin{array}{lcl}
\mathbf{\Psi}^*(u\wedge v)=\mathbf{\Psi}^*u\wedge\mathbf{\Psi}^*v, &\quad&
\mathbf{\Psi}_*^{-1}(u\wedge v)=\mathbf{\Psi}_*^{-1}u\wedge\mathbf{\Psi}_*^{-1}v,
\\[4pt]
\mathbf{\Psi}^*(u\iprod\, v)=\mathbf{\Psi}_*^{-1}u\iprod\,\mathbf{\Psi}^*v,
&&
\mathbf{\Psi}_*^{-1}(u\iprod\, v)=\mathbf{\Psi}^*u\iprod\,\mathbf{\Psi}_*^{-1}v.
\end{array}
$$

A very important particular case, if $\Omega=\{ (x',x_n)\in\RR^{n-1}\times\RR=\RR^n\,:\, x_n>\phi(x')\}$ is a special Lipschitz domain, so that $\mathbf{\Psi}=(x',x_n+\phi(x'))$ with outward unit normal $\mathbf{\nu}$, and $u\in\D_{p}(\d,\Omega)$, then one has on $\partial\Omega$
\begin{align*}
    \sqrt{1+|\nabla'\phi|^2} \mathbf{\Psi}_*^{-1}(\nu \iprod\, u) =  (-\mathfrak{e}_n)\iprod(\mathbf{\Psi}_*^{-1} u)
\end{align*}
since it can be checked that $\sqrt{1+|\nabla'\phi|^2}\mathbf{\Psi}^*(\nu)=-\mathfrak{e}_n$ almost everywhere on $\partial\Omega$. 
\end{remark}

\begin{remark}
\label{product-rule}
By the product rule for the exterior derivative and the interior derivative we 
have that for all bounded Lipschitz scalar-valued  functions $\eta$, for all 
$u\in {\D}^p(\d,\Omega)$ and $v\in{\D}^p({\delta},\Omega)$, 
then $\eta u\in {\D}^p(\d,\Omega)$, 
$\eta v\in{\D}^p({\delta},\Omega)$ with
\begin{equation}
\label{eq:product-rule}
\d(\eta u)=\eta\, \d u+\nabla\eta\wedge u\quad\mbox{and}\quad
\delta(\eta v)=\eta\,\delta v-\nabla\eta\iprod\,v .
\end{equation}
More generally, for $u$ a bounded Lipschitz $\ell$-form, for all 
$v\in {\D}^p(\d,\Omega)$, it holds
\begin{equation}
\label{eq:product-rule2}
\d(u\wedge v)=\d u\wedge v+(-1)^\ell u\wedge \d v,
\end{equation}
which gives also for all bounded Lipschitz scalar-valued  functions $\eta$, 
and for all $u\in {\D}^p(\d,\Omega)$ and all $v\in{\D}^p(\delta,\Omega)$:
\begin{equation}
\label{eq:product-rule3}
\d(\nabla\eta\wedge u)=-\nabla\eta\wedge \d u\text{ and }\delta(\nabla\eta\iprod v)=-\nabla\eta \iprod \delta v.
\end{equation}
\end{remark} 

When dealing with the natural bi-Lipschitz map arising from the description of the boundary, one is constrained to $\L^p$ spaces since the Jacobian of $\Psi$ only contains $\L^\infty$ coefficients, which destroys differential forms of higher and lower regularity, even $\C^{1}$ considered as such, does not seem to be enough especially to reach the endpoint spaces. However, it turns out that, a slight improvement of the regularity is enough. However, one has to consider a different global change of coordinates: see the next result.

\begin{lemma}\label{lem:SmoothingPullback} Let $\alpha>0$, $r\in[1,\infty]$, and let $\Omega$  be a special Lipschitz domain of class $\mathcal{M}^{1+\alpha,r}_\W$, \textit{i.e.} such that there exists $\phi\in\C^{0,1}(\RR^{{n-1}})$ satisfying additionally $\nabla'\phi\in\mathcal{M}^{\alpha,r}_{\W,\tt{or}}(\RR^{n-1})$ and such that, up to a rotation,
\begin{align*}
    \Omega=\{ (x',x_n)\in\RR^{n-1}\times\RR=\RR^n\,:\, x_n>\phi(x')\}.
\end{align*}

Then, the pullback and the (inverse) pushforward maps $\mathbf{\Phi}^{\ast}$ and $\mathbf{\Phi}_{\ast}^{-1}$ induced by $\phi$ according to \eqref{eq:Phi} are isomorphisms from ${\B}^{s}_{p,q}(\Omega,\Lambda^k)$ onto ${\B}^{s}_{p,q}(\RR^n_+,\Lambda^k)$ whenever either
\begin{enumerate}
    \item $p\in(r,\infty]$, $s\in(-1+\frac{1}{p},\frac{1+r\alpha}{p})$; or
    \item $p\in[1,r]$, $s\in(-1+\frac{1}{p},\alpha+\frac{1}{p})$; or
    \item $p=q=r$, $s=\alpha+\sfrac{1}{p}$.
\end{enumerate}
Furthermore, it holds that
\begin{itemize}
    \item one has the point-wise equality
    \begin{align*}
    -\mathfrak{e}_n = {\sqrt{1+|\nabla'\phi|^2}}\mathbf{\Phi}^{\ast}_{|_{\partial\RR^n_+}}(\nu);
\end{align*}
    \item the results still holds if we replace $\B^{s}_{p,q}$, by either $\BesSmo^{s}_{p,\infty}$, $\B^{s,0}_{\infty,q}$, $\BesSmo^{s,0}_{\infty,\infty}$ or $\H^{s,p}$ requiring for  the latter $1<p=q<\infty$.
    \item when $p=\infty$, for $\nabla'\phi \in \B^{\alpha}_{\infty,q}(\RR^{n-1})$, then the result still holds for the function spaces $\B^{s}_{\infty,q}$, $-1-\alpha < s\leqslant \alpha$. The same goes if we replace $\B^{\bullet}_{\infty,q}$ by either $\B^{\bullet,0}_{\infty,q}$, $\BesSmo^{\bullet, 0}_{\infty,\infty}$.
\end{itemize}
\end{lemma}

\begin{proof} We recall that the formula for the pullback and the (inverse) pushforward are given by
\begin{align*}
    \boldsymbol{\Phi}^\ast u = \prescript{t}{}{(\nabla \boldsymbol{\Phi})} [ u\circ \boldsymbol{\Phi}] \text{ and } \boldsymbol{\Phi}_\ast^{-1} u = {(\nabla \boldsymbol{\Phi})^{-1}} [ u\circ \boldsymbol{\Phi}].
\end{align*}
The result is then a direct consequence of Proposition~\ref{prop:MultipliersintheLplikeRange} and Lemma~\ref{lem:multiplierComposition}.
\end{proof}

The next straightforward corollary is fundamental for our later purposes.

\begin{corollary}\label{cor:SmoothPullbackdomainDiffForms}With the same notations and assumptions established in Lemma~\ref{lem:SmoothingPullback}, one has isomorphisms
\begin{align*}
    \begin{array}{rllr}
    \mathbf{\Phi}^{\ast}&\,:\, \D^{s}_{p,q}(\d, \Omega,\Lambda )&\longmapsto\D^{s}_{p,q}(\d, \RR^n_+,\Lambda ),\\
    \tilde{\mathbf{\Phi}}_{\ast}^{-1}&\,:\, \D^{s}_{p,q}(\delta, \Omega,\Lambda )&\longmapsto\D^{s}_{p,q}(\delta, \RR^n_+,\Lambda ).
    \end{array}
\end{align*}
Furthermore, the result still holds
\begin{itemize}
    \item if we replace $\D^{s}_{p,q}$ by either $\mathcal{D}^{s}_{p,\infty}$, $\D^{s,0}_{\infty,q}$, $\mathcal{D}^{s,0}_{\infty,\infty}$  or even $\D^{s}_{p}$ requiring $1<p<\infty$ for the latter;
    \item  if we replace $\d$ and $\delta$ by respectively $\underline{\d}$ and $\underline{\delta}$, whenever $-1+{\sfrac{1}{p}}<s<{\sfrac{1}{p}}$.
\end{itemize}
\end{corollary}

\begin{remark}\label{rem:SobBesovDiffFormPullback}Recall that when $k=0$, one has
\begin{align*}
    \D^{s}_{p,q}(\d, \Omega,\Lambda^{0} ) = \B^{s+1}_{p,q}(\Omega,\CC)
\end{align*}
and similarly for other inhomogeneous function spaces. 
\end{remark}

\newpage
\section{Prelude: Further information on the Dirichlet and Neumann Laplacian}\label{Sec:Laplacians}

This very short preliminary section is about gathering further information on the regularity and the well-posedness for the Dirichlet and Neumann problem and their resolvent version. First on  the half-space $\RR^n_+$, then on bounded domains. In this specific section, all the function spaces are over $\CC$.

\begin{definition}The (negative) Dirichlet Laplacian $-\Delta_\mathcal{D}$ on a Lipschitz domain $\Omega$, is the $\L^2$-realization of the bounded linear operator $-\Delta\,:\,\H^{1,2}_{0}(\Omega)\longrightarrow\H^{-1,2}(\Omega)$ denoted by $(\D_2(\Delta_{\mathcal{D}}),-\Delta_{\mathcal{D}})$, and induced by the sesquilinear form
\begin{align*}
    \mathfrak{a}_{\mathcal{D}}\,:\,\D(\mathfrak{a}_\mathcal{D})\times \D(\mathfrak{a}_\mathcal{D}) &\longrightarrow \CC\\
    (u ,v)\qquad&\longmapsto \langle \nabla  u, \nabla  v\rangle_{\Omega},
\end{align*}
with form domain $\D(\mathfrak{a}_\mathcal{D}):=\H^{1,2}_{0}(\Omega)$.
\end{definition}

\begin{definition}The (negative) Neumann Laplacian $-\Delta_\mathcal{N}$ on a Lipschitz domain $\Omega$, is the $\L^2$-realization of the bounded linear operator $-\Delta\,:\,\H^{1,2}(\Omega)\longrightarrow\H^{-1,2}_0(\Omega)$ denoted by $(\D_2(\Delta_{\mathcal{N}}),-\Delta_{\mathcal{N}})$, and induced by the sesquilinear form
\begin{align*}
    \mathfrak{a}_{\mathcal{N}}\,:\,\D(\mathfrak{a}_\mathcal{N})\times \D(\mathfrak{a}_\mathcal{N}) &\longrightarrow \CC\\
    (u ,v)\qquad&\longmapsto \langle \nabla  u, \nabla  v\rangle_{\Omega},
\end{align*}
with form domain $\D(\mathfrak{a}_\mathcal{N}):=\H^{1,2}(\Omega)$. 
\end{definition}

So that it is easy to see, and well-known, that both the Dirichlet and Neumann Laplacians, are closed, densely defined, non-negative self-adjoint operators on $\mathrm{L}^2(\Omega,\mathbb{C})$, see \cite[Chapter~1,~Section~1.2]{bookOuhabaz2005}.

The questions highlighted and answered in this section are usual
\begin{enumerate}
    \item regularity properties of solutions on $\L^p$ and other function spaces, including endpoint ones;
    \item regularity properties of the semigroup, boundedness of the Holomorphic functional calculus;
    \item identification of the domain of its fractional powers in various context when it is relevant.
\end{enumerate}
Although most of the questions addressed here have already been studied extensively in prior work, certain subtleties—well known to experts—remain sparsely documented by explicit references. In this note, we gather several of these refinements and provide explicit proofs or, at least a clear sketch of the central arguments.

\subsection{The full Dirichlet problem for the Laplacian on the half-space}

Since the solution on the flat half-space is given by odd reflection with respect to the boundary $\partial\RR^n_+$, the proofs below mostly focus on proving uniqueness and how to reach the case of endpoint function spaces and how to deal with non-zero boundary values. For preliminary results in this direction, where the endpoint cases and non-zero boundary values are omitted, see \cite[Section~5]{Gaudin2022}.

\begin{proposition}\label{prop:DirResolventPbRn+}Let $p,q\in[1,\infty]$, $-1+{\sfrac{1}{p}}<s<{\sfrac{1}{p}}$. Let $\mu\in[0,\pi)$, consider $f\in\dot{\B}^{s}_{p,q}(\RR^n_+)$ and assume one of the following conditions is satisfied
\begin{enumerate}
    \item $h\in\dot{\B}^{s-{\sfrac{1}{p}}}_{p,q}\cap \dot{\B}^{s+2-{\sfrac{1}{p}}}_{p,q}(\partial\RR^n_+)$; or
    \item $h\in{\L}^{p}\cap \dot{\B}^{s+2-{\sfrac{1}{p}}}_{p,q}(\partial\RR^n_+)$; or
    \item $h:=H_{|_{\partial\RR^n_+}}$ with $H\in\dot{\B}^{s}_{p,q}\cap \dot{\B}^{s+2}_{p,q}(\RR^n_+)$.
\end{enumerate}
Then, for all $\lambda\in\Sigma_\mu$, the Dirichlet resolvent problem
\begin{equation*}\tag{$\mathcal{DL}_{\lambda}$}\label{ResolvDirLap}
    \left\{ \begin{array}{rllr}
         \lambda u - \Delta u &= f \text{, }&&\text{ in } \RR^n_+\text{,}\\
        u_{|_{\partial\RR^n_+}} &=h\text{, }&&\text{ on } \partial\RR^n_+\text{.}
    \end{array}
    \right.
\end{equation*}
admits a unique solution $u\in\dot{\B}^{s}_{p,q}\cap \dot{\B}^{s+2}_{p,q}(\RR^n_+)$ which obeys the estimate
\begin{align*}
    |\lambda|\lVert  u\rVert_{\dot{\B}^{s}_{p,q}(\RR^n_+)}+\lVert \nabla^2 u\rVert_{\dot{\B}^{s}_{p,q}(\RR^n_+)}&\lesssim_{p,n,s}^{\mu} \lVert  f\rVert_{\dot{\B}^{s}_{p,q}(\RR^n_+)}\\ &\qquad\qquad + \min \begin{cases}
  |\lambda|\lVert  H\rVert_{\dot{\B}^{s}_{p,q}(\RR^n_+)}+\lVert  H\rVert_{\dot{\B}^{s+2}_{p,q}(\RR^n_+)}\\    
  |\lambda|\lVert  h\rVert_{\dot{\B}^{s-{\sfrac{1}{p}}}_{p,q}(\partial\RR^n_+)}+\lVert  h\rVert_{\dot{\B}^{s+2-{\sfrac{1}{p}}}_{p,q}(\partial\RR^n_+)}\\    
  |\lambda|^{1-(\sfrac{1}{2p}-\sfrac{s}{2})}\lVert  h\rVert_{\L^p(\partial\RR^n_+)}+\lVert  h\rVert_{\dot{\B}^{s+2-{\sfrac{1}{p}}}_{p,q}(\partial\RR^n_+)}
  \end{cases}.
\end{align*}
If additionally, for some $r,q\in[1,\infty]$, $\alpha\in(-1+\sfrac{1}{r},2+\sfrac{1}{r})$, $\alpha\neq\sfrac{1}{r},1+\sfrac{1}{r}$, one also has $f\in\dot{\B}^{\alpha}_{r,\kappa,\mathcal{D}}(\RR^n_+)$ and $h=0$, then one  obtains 
$u\in\dot{\B}^{\alpha}_{r,\kappa}\cap \dot{\B}^{\alpha+2}_{r,\kappa}(\RR^n_+)$ with the estimate
\begin{align*}
    |\lambda|\lVert  u\rVert_{\dot{\B}^{\alpha}_{r,\kappa}(\RR^n_+)}+\lVert \nabla^2 u\rVert_{\dot{\B}^{\alpha}_{r,\kappa}(\RR^n_+)}\lesssim_{r,n,\alpha}^{\mu} \lVert  f\rVert_{\dot{\B}^{\alpha}_{r,\kappa}(\RR^n_+)}.
\end{align*}
Furthermore
\begin{itemize}
    \item a similar result holds for the remaining endpoint Besov spaces $\dot{\BesSmo}^{s}_{p,\infty}$, $\dot{\B}^{s,0}_{\infty,q}$, and  $\dot{\BesSmo}^{s,0}_{\infty,\infty}$;
    
    \item In the case $p=q\in(1,\infty)$,  the result still holds replacing the spaces $(\dot{\B}^{\bullet}_{p,q}(\RR^n_+),\dot{\B}^{\bullet}_{p,q}(\partial\RR^n_+))$ by $(\dot{\H}^{\bullet,p}(\RR^n_+),\dot{\B}^{\bullet}_{p,p}(\partial\RR^n_+))$, and similarly whenever $r=\kappa\in(1,\infty)$;

    \item For any $\theta\in(0,\pi)$, the Dirichlet Laplacian has a $\mathbf{H}^\infty(\Sigma_\theta)$-functional calculus on any previously mentioned function spaces;
    
    \item  If either $h=0$ or $\lambda=1$, everything remains valid for the corresponding inhomogeneous function spaces whenever $1<p<\infty$.
\end{itemize}
\end{proposition}

\begin{proof} The case $p\in(1,\infty)$ when $h=0$ has been already treated in \cite[Section~5,~Proposition~5.3]{Gaudin2022}. We treat the remaining cases. Let $\lambda\in\Sigma_\mu$.

\textbf{Step 1:} The case $h=0$.

\textbf{Step 1.1:} Uniqueness. We deal with the case $p=\infty$, the case $p=1$ is similar.

Let $u\in\dot{\B}^{s}_{\infty,q}\cap \dot{\B}^{s+2}_{\infty,q}(\RR^n_+)$ be such that it satisfies \eqref{ResolvDirLap} with $f=0$, $h=0$. By interpolation inequalities, recalling that $-1<s<0$, it holds that $u\in {\B}^{s+2}_{\infty,q,\mathcal{D}}(\RR^n_+)\subset \C^{1,1+s}_{ub}(\overline{\RR^{n}_+})$, with $u_{|_{\partial\RR^n_+}}=0$. We consider $U:=\E_{\mathcal{D}}u$ which belongs to $\dot{\B}^{s}_{\infty,q}\cap \dot{\B}^{s+2}_{\infty,q}(\RR^n)$ by Lemma \ref{lem:ExtDirNeuRn+}, and by the trace theorem $U(\cdot,0) = U_{|_{\partial\RR^n_+}}= U_{|_{\partial\RR^n_-}}=0$. Let  $\varphi\in\Ccinfty(\RR^n)$, then it is usual to figure out
\begin{align*}
    \langle \lambda U - \Delta U, \varphi \rangle_{\RR^n} &= \int_{\RR^n} \Big[ \lambda U(x) \varphi(x) + \nabla U(x) \cdot \nabla \varphi(x)\Big] \,\d x\\
    &= \int_{\RR^n_+} \lambda u(x) [\varphi(x)-\varphi(x',-x_n)] \,\d x\\
    & \qquad+ \int_{\RR^n_+}  \Big[ \nabla' u(x) \cdot \nabla' [\varphi(x)-\varphi(x',-x_n)] + \partial_{x_n}u(x) \partial_{x_n}[\varphi(x)-\varphi(x',-x_n)]\Big]\,\d x\\
    &= \langle \lambda u - \Delta u, [\varphi - \varphi(\cdot,-\cdot)] \rangle_{\RR^n_+} =0,
\end{align*}
since $\mathbbm{1}_{\RR^n_+}[\varphi-\varphi(\cdot,-\cdot)] = \mathrm{P}_0 \varphi := [\I-\E_{\mathcal{N}}^{-}]\varphi$ is an element of $\dot{\B}^{-s}_{1,r,0}(\RR^{n}_+)\cap\dot{\B}^{-s+1}_{1,r,0}(\RR^{n}_+)$ for all $r\in[1,\infty]$, recalling that $\dot{\B}^{\alpha}_{\infty,q}(\RR^{n}_+)= (\dot{\B}^{-\alpha}_{1,q',0}(\RR^{n}_+))'$, and $\dot{\B}^{\alpha}_{\infty,1}(\RR^{n}_+)= (\dot{\BesSmo}^{-\alpha}_{1,\infty,0}(\RR^{n}_+))$, when $\alpha<0$, $q\in(1,\infty],$. Thus $\lambda U - \Delta U=0$ in $\mathcal{S}'(\RR^n)$, hence $U=0$, and then $u=0$.

\textbf{Step 1.2:} The existence is well known and standard, given by $u:=[(\lambda\I-\Delta)^{-1}\E_{\mathcal{D}}f]_{|_{\RR^n_+}}$ as in, \textit{e.g.}, the proof of \cite[Section~5,~Proposition~5.3]{Gaudin2022}.

\textbf{Step 2:} Existence for $f=0$, $h\neq0$, $p\in[1,\infty]$. Up to a rescaling argument, and by uniqueness, one may replace $(u,h)$ by $\big(|\lambda|u(\cdot/|\lambda|^{1/2}), |\lambda|h(\cdot/|\lambda|^{1/2}))$ and assume $|\lambda|=1$. It suffices to set $u:=e^{-x_n(\lambda\I-\Delta')^{\sfrac{1}{2}}}h$, and by Proposition~\ref{prop:DumpedPoissonSemigroup2}, one obtains
\begin{align*}
    \lVert u\rVert_{\dot{\B}^{s}_{p,q}(\RR^n_+)} + \lVert \nabla^2 u\rVert_{\dot{\B}^{s}_{p,q}(\RR^n_+)} &\lesssim_{p,s,n}^{\mu} \lVert h\rVert_{{\B}^{s-{\sfrac{1}{p}}}_{p,q}(\partial\RR^{n}_+)}  +\lVert h\rVert_{{\B}^{s+2-{\sfrac{1}{p}}}_{p,q}(\partial\RR^{n}_+)}\\ &\lesssim_{p,s,n}^{\mu} \lVert h\rVert_{{\B}^{s+2-{\sfrac{1}{p}}}_{p,q}(\partial\RR^{n}_+)}\\
    &\lesssim_{p,s,n}^{\mu} \lVert h\rVert_{{\L}^{p}(\partial\RR^{n}_+)} + \lVert h\rVert_{\dot{\B}^{s+2-{\sfrac{1}{p}}}_{p,q}(\partial\RR^{n}_+)}.
\end{align*}
From this point, if $h\in\dot{\B}^{s-{\sfrac{1}{p}}}_{p,q}\cap\dot{\B}^{s+2-{\sfrac{1}{p}}}_{p,q}(\partial\RR^{n}_+)$, from interpolation inequalities with $\theta={\sfrac{1}{2}}({\sfrac{1}{p}}-s)$, one deduces
\begin{align*}
    \lVert h\rVert_{{\L}^{p}(\partial\RR^{n}_+)} + \lVert h\rVert_{\dot{\B}^{s+2-{\sfrac{1}{p}}}_{p,q}(\partial\RR^{n}_+)} &\lesssim_{p,s,n} \lVert h\rVert_{\dot{\B}^{s-{\sfrac{1}{p}}}_{p,q}(\partial\RR^{n}_+)}^{1-\theta}\lVert h\rVert_{\dot{\B}^{s+2-{\sfrac{1}{p}}}_{p,q}(\partial\RR^{n}_+)}^\theta+\lVert h\rVert_{\dot{\B}^{s+2-{\sfrac{1}{p}}}_{p,q}(\partial\RR^{n}_+)} \\&\lesssim_{p,s,n} \lVert h\rVert_{\dot{\B}^{s-{\sfrac{1}{p}}}_{p,q}(\partial\RR^{n}_+)}+\lVert h\rVert_{\dot{\B}^{s+2-{\sfrac{1}{p}}}_{p,q}(\partial\RR^{n}_+)}.
\end{align*}
If instead $h=H_{|_{\partial\RR^n_+}}$, provided $H\in\dot{\B}^{s}_{p,q}\cap\dot{\B}^{s+2}_{p,q}(\RR^{n}_+)$, by the trace theorem and then again by interpolation inequalities,
\begin{align*}
    \lVert h\rVert_{{\L}^{p}(\partial\RR^{n}_+)} + \lVert h\rVert_{\dot{\B}^{s+2-{\sfrac{1}{p}}}_{p,q}(\partial\RR^{n}_+)} &\lesssim_{p,s,n} \lVert H\rVert_{\dot{\B}^{{\sfrac{1}{p}}}_{p,1}(\RR^{n}_+)}+\lVert H\rVert_{\dot{\B}^{s+2}_{p,q}(\partial\RR^{n}_+)}\\&\lesssim_{p,s,n} \lVert H\rVert_{\dot{\B}^{s}_{p,q}(\RR^{n}_+)}^{1-\theta}\lVert H\rVert_{\dot{\B}^{s+2}_{p,q}(\RR^{n}_+)}^\theta+\lVert H\rVert_{\dot{\B}^{s+2}_{p,q}(\partial\RR^{n}_+)}\\ &\lesssim_{p,s,n} \lVert H\rVert_{\dot{\B}^{s}_{p,q}(\partial\RR^{n}_+)}+\lVert H\rVert_{\dot{\B}^{s+2}_{p,q}(\partial\RR^{n}_+)}.
\end{align*}
This ends the proof.
\end{proof}

\begin{proposition}\label{prop:DirPbRn+}Let $p,q\in[1,\infty]$, $-2+{\sfrac{1}{p}}<s<{\sfrac{1}{p}}$. Let $f\in\dot{\B}^{s}_{p,q}(\RR^n_+)$ and $h\in\dot{\B}^{s+2-{\sfrac{1}{p}}}_{p,q}(\partial\RR^n_+)$. Then, the following Dirichlet problem
\begin{equation*}\tag{$\mathcal{DL}_{0}$}\label{DirLap}
    \left\{ \begin{array}{rllr}
         - \Delta u &= f \text{, }&&\text{ in } \RR^n_+\text{,}\\
        u_{|_{\partial\RR^n_+}} &=h\text{, }&&\text{ on } \partial\RR^n_+\text{.}
    \end{array}
    \right.
\end{equation*}
if it admits a solution $u\in\dot{\B}^{s+2}_{p,q}(\RR^n_+)$, it is unique. Moreover, it obeys the estimates,
\begin{align*}
    \lVert u\rVert_{\L^\infty_{x_n}(\RR_+,\dot{\B}^{s+2-{\sfrac{1}{p}}}_{p,q}(\RR^{n-1}))}\lesssim_{p,n,s}\lVert \nabla^2 u\rVert_{\dot{\B}^{s}_{p,q}(\RR^n_+)}\lesssim_{p,n,s} \lVert  f\rVert_{\dot{\B}^{s}_{p,q}(\RR^n_+)} + \lVert  h\rVert_{\dot{\B}^{s+2-{\sfrac{1}{p}}}_{p,q}(\partial\RR^n_+)}.
\end{align*}
Furthermore,
\begin{itemize}
    \item  if $f\in\dot{\B}^{\alpha}_{r,\kappa}(\RR^n_+)$, for some $r,\kappa\in[1,\infty]$, $\alpha\in(-2+\sfrac{1}{r},\sfrac{1}{r})$ satisfying $(\mathcal{C}_{\alpha+2,r,\kappa})$, then the solution always exists;
    \item a similar result holds for the endpoint Besov spaces $\dot{\BesSmo}^{s}_{p,\infty}$, $\dot{\B}^{s,0}_{\infty,q}$, and  $\dot{\BesSmo}^{s,0}_{\infty,\infty}$, and the solution is strongly continuous with respect to $x_n$;
    \item in the case $p=q\in(1,\infty)$, the result still holds replacing the spaces $(\dot{\B}^{s+2}_{p,q}(\RR^n_+),\dot{\B}^{s+2-{\sfrac{1}{p}}}_{p,q}(\partial\RR^n_+))$ by $(\dot{\H}^{s+2,p}(\RR^n_+),\dot{\B}^{s+2-{\sfrac{1}{p}}}_{p,p}(\partial\RR^n_+))$.
\end{itemize}
\end{proposition}

\begin{proof} By the trace theorem and Proposition~\ref{prop:PoissonSemigroup3}, up to consider $u-(-\Delta_{\mathcal{D,\partial}})^{-1}h$, we can assume $h=0$. Assume $u\in\dot{\B}^{\alpha}_{r,\kappa}(\mathbb{R}^n_+)$, with $\alpha>1/r$ is a solution to \eqref{DirLap} with $h=0$ and $f\in\dot{\B}^{s}_{p,q}(\mathbb{R}^n_+)$. We want to show that, necessarily, $u\in \dot{\B}^{s+2}_{p,q}(\mathbb{R}^n_+)$ with the corresponding estimate. We assume ${\sfrac{1}{p}}<s+2<2+{\sfrac{1}{p}}$.

As for Lemma~\ref{lem:ExtDirNeuRn+}, by the same computations in its proof, it holds that $\E_\mathcal{D}u\in \dot{\B}^{s+2}_{p,q}(\mathbb{R}^n)$, since we have for any $\varphi\in\Ccinfty(\RR^n)$,
\begin{align*}
    \langle -\Delta \E_\mathcal{D}u, \varphi \rangle_{\RR^n} = \langle f, [\varphi-\varphi(\cdot,-\cdot)] \rangle_{\RR^n_+} = \langle \E_\mathcal{D}f, \varphi \rangle_{\RR^n},
\end{align*}
obtained since $[\varphi-\varphi(\cdot,-\cdot)]\in\dot{\B}^{-s}_{p',r,0}(\mathbb{R}^n_+)$, for all $r\in[1,\infty]$. By duality, and the boundedness of $\E_\mathcal{D}$, we obtain the result.
\end{proof}

\subsection{The Dirichlet problem on bounded domains}

We start with an important and known theorem for bounded Lipschitz domains and which will be of a certain importance for latter parts of the present work. This is a particular instance/consequence of a famous result of Fabes, Mendez and Mitrea  \cite[Theorem~10.1]{FabesMendezMitrea1998}.

\begin{theorem}\label{thm:DirLapC1Domains} Let  $\Omega$ is be a bounded Lipschitz domain. There exists $\varepsilon_{{}_\Omega}\in(0,1]$ depending on $\Omega$, such that for all $p\in(1,\infty)$, it holds that
\begin{enumerate}
    \item $-\Delta_\mathcal{D}\,:\,\H^{s+1,p}_{0}(\Omega)\longrightarrow\H^{s-1,p}(\Omega)$ is an isomorphism  for all $s\in(-1+\sfrac{1}{p},\sfrac{1}{p})$ that satisfies either
    \begin{itemize}
        \item $s\in(\sfrac{3}{p}-2-\varepsilon_{{}_\Omega},\sfrac{1}{p})$, if $p\in(1,\frac{2}{1+\varepsilon_{{}_\Omega}}]$; or
        \item $s\in(-1+\sfrac{1}{p},\sfrac{1}{p})$, if $p\in[\frac{2}{1+\varepsilon_{{}_\Omega}},\frac{2}{1-\varepsilon_{\Omega}}]$; or
        \item $s\in(-1+\sfrac{1}{p},\sfrac{3}{p}-1+\varepsilon_{{}_\Omega})$, if $p\in[\frac{2}{1-\varepsilon_{{}_\Omega}},\infty)$.
    \end{itemize}
    \item $(\D_p(-\Delta_{\mathcal{D}}),-\Delta_{\mathcal{D}})$ is an invertible $0$-sectorial operator on $\L^p(\Omega)$ which admits a bounded $\mathbf{H}^{\infty}(\Sigma_\theta)$-functional calculus for all $\theta\in(0,\pi)$;
    \item For all $\alpha \in(-1+\sfrac{1}{p},\sfrac{1}{p})$, one has 
    \begin{align*}
        \D_p((-\Delta_\mathcal{D})^{\sfrac{\alpha}{2}}) = \H^{\alpha,p}(\Omega).
\end{align*}
with equivalence of norms.
\end{enumerate}
By real interpolation similar results are available for Besov spaces $\B^{\bullet}_{p,q}$ and $\BesSmo^{\bullet}_{p,\infty}$, $p\in(1,\infty)$, $q\in[1,\infty]$.
\end{theorem}

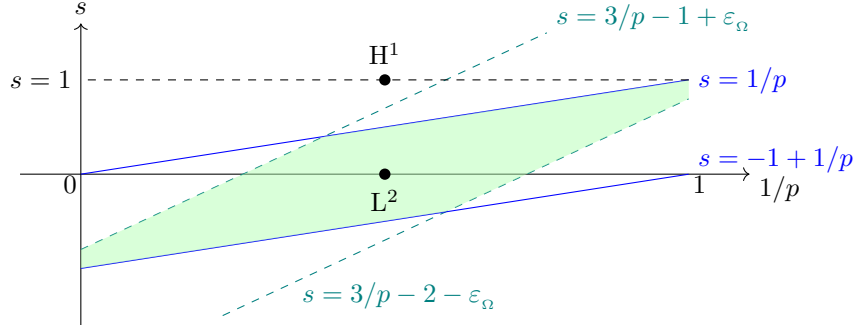
\begin{figure}[H]
\centering
\begin{tikzpicture}[yscale=1.25,xscale=8]

\def\littleparam{0.2}

  \draw[->] (-0.1,0) -- (1.1,0) node[right,yshift=-2mm] {$1/p$};
  \draw[->] (0,-1.6) -- (0,1.6) node[above] {$s$};

  \draw[domain=0:(5/6-\littleparam/3)),dashed,variable=\x,teal] plot (\x,3*\x-1+\littleparam) node[right,yshift=2mm] {$s=3/p-1+\varepsilon_{{}_\Omega}$};
  \draw[domain=(1/6+\littleparam/3):1,dashed,variable=\x,teal] plot (\x,3*\x-2-\littleparam);
  \node[anchor=north east] at (0.7,-1.05) {\color{teal}$s=3/p-2-\varepsilon_{{}_\Omega}$};
  
  \draw[domain=0:1,smooth,variable=\x,blue] plot ({\x},{-1+\x}) node[right,yshift=2mm] {$s=-1+1/p$};
  \draw[domain=0:1,smooth,variable=\x,blue] plot ({\x},{\x}) node[right]{$s=1/p$};
 \fill[green!30,opacity=0.5] (0,-1+\littleparam) -- plot[domain=0:{(1-\littleparam)/2}](\x,3*\x-1+\littleparam)-- plot[domain={(1-\littleparam)/2}:1](\x,\x)  -- (1,1-\littleparam)-- ({(1+\littleparam)/2},{-(1-\littleparam)/2})  -- (0,-1)-- cycle;

  \draw[dashed] (1,1) -- (0,1)  node[left] {$s=1$};
  \node[circle,fill,inner sep=1.5pt,label=below:$\mathrm{L}^2$] at (0.5,0) {};
  \node[circle,fill,inner sep=1.5pt,label=above:${\mathrm{H}}^1$] at (0.5,1) {};
  \draw[circle,fill,inner sep=1pt] (0,0) node[below left] {$0$};
  \draw[circle,fill,inner sep=1pt] (1,0) node[below right] {$1$};
\end{tikzpicture}
\caption{Representation of $(s,\sfrac{1}{p})$ for the isomorphism property in point \textit{(i)} of above Theorem~\ref{thm:DirLapC1Domains}.}
\label{Fig:RegularityDirLapLipDomains}
\end{figure}

\begin{proof} The really though part is Point \textit{(i)} that comes from \cite[Theorem~10.1]{FabesMendezMitrea1998}.
\medbreak
\textbf{For Point \textit{(ii)}:} as a positive self-adjoint operator $\L^2(\Omega)$, $(\D_2(-\Delta_{\mathcal{D}}),-\Delta_{\mathcal{D}})$ is an invertible $0$-sectorial operator on $\L^2(\Omega)$ which admits a bounded $\mathbf{H}^{\infty}(\Sigma_\theta)$-functional calculus for all $\theta\in(0,\pi)$. One can check that \cite[Corollary~2.18]{bookOuhabaz2005} applies, so that the generated semigroup $(e^{t\Delta_{\mathcal{D}}})_{t\geqslant 0}$ extends naturally  has a uniformly bounded holomorphic one on $\L^p(\Omega)$ for all $p\in[1,\infty]$. Also, the bound
\begin{align}\label{eq:LpLqDecay}
    \lVert e^{t\Delta_{\mathcal{D}}} u\rVert_{\L^q(\Omega)} \lesssim_{p,q,n,\Omega} \frac{1}{t^{\frac{n}{2}(\frac{1}{p}-\frac{1}{q})}}\lVert  u\rVert_{\L^p(\Omega)}, 
\end{align}
holds for all $u\in\L^p(\Omega)$ and all $t>0$, provided $1\leqslant p\leqslant q\leqslant \infty$.
Now, notice that as in the proof of \cite[Theorem~8]{Davies1995}, one can deduce there exists $c>0$, for all measurable sets $E,F\subset\Omega$, all $u\in\L^2(\Omega)$, all $t>0$
\begin{align}\label{eq:L2offdiag}
    \lVert \mathbbm{1}_{F} e^{t\Delta_{\mathcal{D}}}(\mathbbm{1}_{E} u)\rVert_{\L^2(\Omega)} \lesssim_{n,\Omega} e^{-c \frac{\d(E,F)^2}{t}}\lVert  u\rVert_{\L^2(E)}.
\end{align}
Having \eqref{eq:LpLqDecay} and \eqref{eq:L2offdiag} in hand, one can apply \cite[Theorem~2.2]{BlunckKuntsmann2003} to deduce boundedness of the $\mathbf{H}^\infty(\Sigma_\theta)$-functional calculus on $\L^p(\Omega)$, for all $\theta\in(0,\pi)$, all $p\in(1,\infty)$.

\medbreak

\textbf{For Point \textit{(iii)}:} From the boundedness of the $\mathbf{H}^\infty(\Sigma_\theta)$-functional calculus on $\L^p(\Omega)$, in particular $(\D_p(-\Delta_{\mathcal{D}}),-\Delta_{\mathcal{D}})$ has BIP on $\L^p(\Omega)$ for all $p\in (1,\infty)$. On $\L^2(\Omega)$, this yields directly for all $\alpha\in(0,1)$, $\alpha\neq\sfrac{1}{2}$,
\begin{align*}
        \H^{\alpha,2}_0(\Omega)=[\L^2(\Omega),\H^{1,2}_0(\Omega)]_\alpha=[\L^2(\Omega),\D_2((-\Delta_\mathcal{D})^{\sfrac{1}{2}})]_\alpha = \D_2((-\Delta_\mathcal{D})^{\sfrac{\alpha}{2}}).
\end{align*}
By duality, for all $s\in[0,1]$, $s\neq\sfrac{1}{2}$, $(-\Delta_\mathcal{D})^{\sfrac{s}{2}}\,:\,\L^2(\Omega)\longrightarrow\H^{-s}(\Omega)$ is also an isomorphism. One has 
\begin{enumerate}
    \item For $\Re (z)=0$, $(-\Delta_\mathcal{D})^{\sfrac{z}{2}}$ act as an isomorphism on $\L^p(\Omega)$, for all $p\in(1,\infty)$,
    \item For $\Re (z)=1$, $(-\Delta_\mathcal{D})^{\sfrac{z}{2}}$ act is an isomorphism mapping $\H^{1,2}_0(\Omega)$ to $\L^2(\Omega)$;
\end{enumerate}
By complex interpolation, and duality, for all $\alpha\in(0,1)$, $2\alpha\neq\sfrac{1}{p}+\NN$, such that $|\frac{1}{p}-\frac{1}{2}|+2|\alpha|<2$, we obtain
\begin{align*}
    \D_p((-\Delta_\mathcal{D})^{\sfrac{\alpha}{2}}) = \H^{\alpha,p}_0(\Omega).
\end{align*}

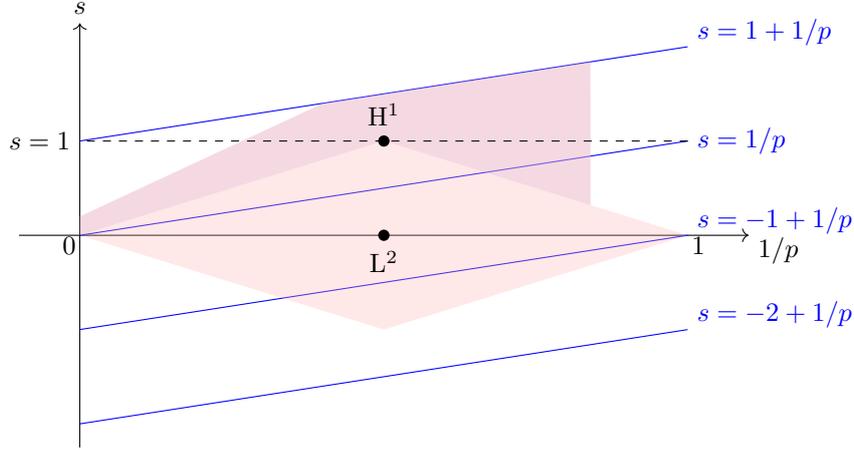
\begin{figure}[H]
\centering
\begin{tikzpicture}[yscale=1.25,xscale=8]
  \draw[->] (-0.1,0) -- (1.1,0) node[right,yshift=-2mm] {$1/p$};
  \draw[->] (0,-2.25) -- (0,2.25) node[above] {$s$};

 \draw[domain=0:1,smooth,variable=\x,blue] plot ({\x},{-2+\x}) node[right,yshift=2mm] {$s=-2+1/p$};
  \draw[domain=0:1,smooth,variable=\x,blue] plot ({\x},{-1+\x}) node[right,yshift=2mm] {$s=-1+1/p$};
  \draw[domain=0:1,smooth,variable=\x,blue] plot ({\x},{\x}) node[right]{$s=1/p$};
  \draw[domain=0:1,smooth,variable=\x,blue] plot ({\x},{1+\x}) node[right,yshift=2mm] {$s=1+1/p$};
   \fill[red!30,opacity=0.3] (0,0) -- plot[domain=0.5:1] (\x,2-2*\x) -- (1/2,0) -- cycle;
  \fill[red!30,opacity=0.3] (0,0) -- plot[domain=0:0.5] (\x,-2*\x)  -- (1,0) -- cycle;

    \def\littleparam{0.2}

   \fill[purple!30,opacity=0.5] (0,0)--(0,\littleparam) -- plot[domain=0:{(1-\littleparam)/2}](\x,3*\x+\littleparam)-- plot[domain={(1-\littleparam)/2}:{(4/5+\littleparam/5)}](\x,\x+1)  -- ({4/5+\littleparam/5},{2/5-2*\littleparam/5})-- (1/2,1) -- cycle;

  \draw[dashed] (1,1) -- (0,1)  node[left] {$s=1$};
  \node[circle,fill,inner sep=1.5pt,label=below:$\mathrm{L}^2$] at (0.5,0) {};
  \node[circle,fill,inner sep=1.5pt,label=above:${\mathrm{H}}^1$] at (0.5,1) {};
  \draw[circle,fill,inner sep=1pt] (0,0) node[below left] {$0$};
  \draw[circle,fill,inner sep=1pt] (1,0) node[below right] {$1$};
\end{tikzpicture}
\caption{Representation of $(s,\sfrac{1}{p})$, for Bessel potential spaces that coincides with the domain of an appropriate fractional power of Dirichlet Laplacian, the parallelogram such that $p\in(1,\infty)$ and $|\frac{1}{p}-\frac{1}{2}|+2|s|<2$. It also includes the first iteration below (the first bullet point).}
\end{figure}
One always has the embedding for $\alpha\in(0,1)$,
\begin{align}\label{eq:embeddingFracDomainDirLap}
        \H^{2\alpha,p}_0(\Omega)=[\L^p(\Omega),\H^{2,p}_0(\Omega)]_\alpha\hookrightarrow[\L^p(\Omega),\D_p(-\Delta_\mathcal{D})]_\alpha = \D_p((-\Delta_\mathcal{D})^{\alpha}).
\end{align}
For the reverse embedding, let $p\in(1,\infty)$, $s\in(-1+{\sfrac{1}{p}},{\sfrac{1}{p}})$, as in Point \textit{(i)}, see Figure~\ref{Fig:RegularityDirLapLipDomains}, let $\FF\in\Ccinfty(\Omega,\CC^{n})$. One has $\div (\FF) \in\H^{s-1,p}(\Omega)$ with
\begin{align*}
    \lVert \div (\FF)\rVert_{\H^{s-1,p}(\Omega)} \leqslant \lVert \FF\rVert_{\H^{s,p}(\Omega)}.
\end{align*}
By part \textit{(i)}, solving
\begin{equation*}
    \left\{ \begin{array}{rllr}
         - \Delta u &= \div (\FF) \text{, }&&\text{ in } \Omega\text{,}\\
        u_{|_{\partial\Omega}} &=0\text{, }&&\text{ on } \partial\Omega\text{.}
    \end{array}
    \right.
\end{equation*}
yields
\begin{align*}
    \lVert \nabla(-\Delta_\mathcal{D})^{-1}\div (\FF)\rVert_{\H^{s,p}(\Omega)} =\lVert \nabla u\rVert_{\H^{s,p}(\Omega)} \lesssim_{p,s,n,\Omega} \lVert \div (\FF)\rVert_{\H^{s-1,p}(\Omega)} \lesssim_{p,s,n,\Omega} \lVert \FF\rVert_{\H^{s,p}(\Omega)}.
\end{align*}
Therefore, writing $(-\Delta_\mathcal{D})^{-\sfrac{1}{2}}\div (\FF) = (-\Delta_\mathcal{D})^{\sfrac{1}{2}}(-\Delta_\mathcal{D})^{-1}\div (\FF) $, thanks to \eqref{eq:embeddingFracDomainDirLap} and Poincaré's inequality Lemma~\ref{lem:Poincaré0}, one obtains
\begin{align*}
    \lVert (-\Delta_\mathcal{D})^{-\sfrac{1}{2}}\div (\FF)\rVert_{\H^{s,p}(\Omega)} &=  \lVert (-\Delta_\mathcal{D})^{\sfrac{1}{2}}(-\Delta_\mathcal{D})^{-1}\div (\FF)\rVert_{\H^{s,p}(\Omega)}\\
    &\lesssim_{p,s,n,\Omega} \lVert \nabla(-\Delta_\mathcal{D})^{-1}\div (\FF)\rVert_{\H^{s,p}(\Omega)}\\
    &\lesssim_{p,s,n,\Omega} \lVert \FF\rVert_{\H^{s,p}(\Omega)}.
\end{align*}
Hence, by duality, one obtains for all $f\in\H^{s,p}(\Omega)$
\begin{align*}
    \lVert \nabla(-\Delta_\mathcal{D})^{-\sfrac{1}{2}} f\rVert_{\H^{s,p}(\Omega)}\lesssim_{p,s,n,\Omega} \lVert f\rVert_{\H^{s,p}(\Omega)}.
\end{align*}
This provides for all $p\in(1,\infty)$, $s\in(-1+\sfrac{1}{p},\sfrac{1}{p})$ as in Figure~\ref{Fig:RegularityDirLapLipDomains}, that $(-\Delta_\mathcal{D})^{\sfrac{1}{2}}\,:\,\H^{s+1,p}_0(\Omega)\longrightarrow\H^{s,p}(\Omega)$ is an isomorphism. Therefore, one deduces
\begin{itemize}
    \item  that we have for all $p\in(1,\infty)$, all $s\in(-1+\sfrac{1}{p},1+\sfrac{1}{p})$, $s\notin\ZZ+{\sfrac{1}{p}}$, $s<\sfrac{1}{p}$, $|\frac{1}{p}-\frac{1}{2}|+2|s|<2$ and requiring $s+1$ belongs to Figure~\ref{Fig:RegularityDirLapLipDomains} for the second identity
\begin{align*}
    \D_p((-\Delta_\mathcal{D})^{\sfrac{s}{2}}) = \H^{s,p}_0(\Omega)\text{ and }\D_p((-\Delta_\mathcal{D})^{\frac{1+s}{2}}) = \H^{s+1,p}_0(\Omega).
\end{align*}
which also imposes $p>\frac{5}{4+\varepsilon_{{}_\Omega}}$;
\item and for $p$ close to $1$, $s\in(-\varepsilon_{{}_\Omega}-1+\sfrac{1}{p},-1+\sfrac{1}{p})$, one has an isomorphism $$(-\Delta_\mathcal{D})^{\sfrac{1}{2}}\,:\,\H^{s+1,p}_0(\Omega)\longrightarrow\H^{s,p}(\Omega).$$
\end{itemize}
Applying complex interpolation between both bullet points, we reach the following isomorphism for all $r\in(1,2]$:
\begin{align*}
    (-\Delta_\mathcal{D})^{\sfrac{1}{2}}\,:\,\H^{1,r}_0(\Omega)\longrightarrow\L^r(\Omega).
\end{align*}
Due to the BIP property implied by the bounded holomorphic functional calculus on $\L^r$, one obtains the isomorphism
\begin{align*}
    (-\Delta_\mathcal{D})^{\sfrac{s}{2}}\,:\,\H^{s,r}(\Omega)\longrightarrow\L^r(\Omega).
\end{align*}
for all $0\leqslant s<\sfrac{1}{r}$. Hence, the result holds for all $p\in(1,\infty)$, $s\in[0,\sfrac{1}{p})$. By duality, one deduces the case of negative regularity indices.
\end{proof}

\subsection{The Neumann resolvent problem on the half-space}

As consequence of the study of the Dirichlet Laplacian on the half-space, one also obtains its counterpart for the Neumann, which will still be of use but of of much lesser importance, anecdotal.

\begin{proposition}\label{prop:NeumannResolventPbRn+}Let $p,q\in[1,\infty]$, $-1+{\sfrac{1}{p}}<s<{\sfrac{1}{p}}$. Let $\mu\in[0,\pi)$, consider $f\in\dot{\B}^{s}_{p,q}(\RR^n_+)$ and assume one of the following conditions is satisfied
\begin{enumerate}
    \item $h\in\dot{\B}^{s-1-{\sfrac{1}{p}}}_{p,q}\cap \dot{\B}^{s+1-{\sfrac{1}{p}}}_{p,q}(\partial\RR^n_+)$;
    \item $h\in{\L}^{p}\cap \dot{\B}^{s+1-{\sfrac{1}{p}}}_{p,q}(\partial\RR^n_+)$; or
    \item $h:=H_{|_{\partial\RR^n_+}}$ with $H\in\dot{\B}^{s-1}_{p,q}\cap \dot{\B}^{s+1}_{p,q}(\RR^n_+)$.
\end{enumerate}

Then, for all $\lambda\in\Sigma_\mu$, the Neumann resolvent problem
\begin{equation*}\tag{$\mathcal{NL}_{\lambda}$}\label{ResolvNeuLap}
    \left\{ \begin{array}{rllr}
         \lambda u - \Delta u &= f \text{, }&&\text{ in } \RR^n_+\text{,}\\
        \partial_{x_n}u_{|_{\partial\RR^n_+}} &=h\text{, }&&\text{ on } \partial\RR^n_+\text{.}
    \end{array}
    \right.
\end{equation*}
admits a unique solution $u\in\dot{\B}^{s}_{p,q}\cap \dot{\B}^{s+2}_{p,q}(\RR^n_+)$ which obeys the estimate
\begin{align*}
    |\lambda|\lVert  u\rVert_{\dot{\B}^{s}_{p,q}(\RR^n_+)}+\lVert \nabla^2 u\rVert_{\dot{\B}^{s}_{p,q}(\RR^n_+)} &\lesssim_{p,n,s}^{\mu} \lVert  f\rVert_{\dot{\B}^{s}_{p,q}(\RR^n_+)}\\
    &\qquad\qquad+ \min \begin{cases}
  |\lambda|\lVert  H\rVert_{\dot{\B}^{s-1}_{p,q}(\RR^n_+)}+\lVert  H\rVert_{\dot{\B}^{s+1}_{p,q}(\RR^n_+)}\\    
  |\lambda|\lVert  h\rVert_{\dot{\B}^{s-1-{\sfrac{1}{p}}}_{p,q}(\partial\RR^n_+)}+\lVert h\rVert_{\dot{\B}^{s+1-{\sfrac{1}{p}}}_{p,q}(\partial\RR^n_+)}\\    
  |\lambda|^{\sfrac{1}{2}-(\sfrac{1}{2p}-\sfrac{s}{2})}\lVert  h\rVert_{\L^p(\partial\RR^n_+)}+\lVert  h\rVert_{\dot{\B}^{s+1-{\sfrac{1}{p}}}_{p,q}(\partial\RR^n_+)}
  \end{cases}.
\end{align*}
If additionally, for some $r,q\in[1,\infty]$, $\alpha\in[\sfrac{1}{r},1+\sfrac{1}{r})$, one has $f\in\dot{\B}^{\alpha}_{r,\kappa}(\RR^n_+)$ then  one also obtains 
$u\in\dot{\B}^{\alpha}_{r,\kappa}\cap \dot{\B}^{\alpha+2}_{r,\kappa}(\RR^n_+)$ with the estimate
\begin{align*}
    |\lambda|\lVert  u\rVert_{\dot{\B}^{\alpha}_{r,\kappa}(\RR^n_+)}+\lVert \nabla^2 u\rVert_{\dot{\B}^{\alpha}_{r,\kappa}(\RR^n_+)}&\lesssim_{r,n,\alpha}^{\mu} \lVert  f\rVert_{\dot{\B}^{\alpha}_{r,\kappa}(\RR^n_+)}\\
    &\qquad\qquad+ \min \begin{cases}
  |\lambda|\lVert  H\rVert_{\dot{\B}^{\alpha-1}_{r,\kappa}(\RR^n_+)}+\lVert  H\rVert_{\dot{\B}^{\alpha+1}_{r,\kappa}(\RR^n_+)}\\    
  |\lambda|\lVert  h\rVert_{\dot{\B}^{\alpha-1-{\sfrac{1}{r}}}_{r,\kappa}(\partial\RR^n_+)}+\lVert  h\rVert_{\dot{\B}^{\alpha+1-{\sfrac{1}{r}}}_{r,\kappa}(\partial\RR^n_+)}\\    
  |\lambda|^{{\sfrac{1}{2}}+\sfrac{\alpha}{2}-\sfrac{1}{2r}}\lVert    h\rVert_{\L^r(\partial\RR^n_+)}+\lVert h\rVert_{\dot{\B}^{\alpha+1-{\sfrac{1}{r}}}_{r,\kappa}(\partial\RR^n_+)}
  \end{cases}.
\end{align*}
Furthermore
\begin{itemize}
    \item If $h=0$, the result remains valid for $-2+\sfrac{1}{p}<s\leqslant-1+\sfrac{1}{p}$: provided $f\in\dot{\B}^{s}_{p,q,0}(\RR^n_+)$ one has a unique solution $u$ lying in $\dot{\B}^{s}_{p,q,0}\cap\dot{\B}^{s+2}_{p,q}(\RR^n_+)$;
    
    \item a similar result holds for the remaining endpoint Besov spaces $\dot{\BesSmo}^{s}_{p,\infty}$, $\dot{\B}^{s,0}_{\infty,q}$, and  $\dot{\BesSmo}^{s,0}_{\infty,\infty}$;
    
    \item In the case $p=q\in(1,\infty)$,  the result still holds replacing the spaces $(\dot{\B}^{\bullet}_{p,q}(\RR^n_+),\dot{\B}^{\bullet}_{p,q}(\partial\RR^n_+))$ by $(\dot{\H}^{\bullet,p}(\RR^n_+),\dot{\B}^{\bullet}_{p,p}(\partial\RR^n_+))$, and similarly whenever $r=\kappa\in(1,\infty)$;

    \item For any $\theta\in(0,\pi)$, the Neumann Laplacian has a $\mathbf{H}^\infty(\Sigma_\theta)$-functional calculus on any previously mentioned function spaces;

    \item  If either $h=0$ or $\lambda=1$, everything remains valid for the corresponding inhomogeneous function spaces whenever $1<p<\infty$.
\end{itemize}
\end{proposition}

\begin{proof} When $h=0$, $-1+\sfrac{1}{p}<s<\sfrac{1}{p}$, $p,q\in[1,\infty]$, the solution is given by even reflexion, so one can prove uniqueness and existence as in the proof of Proposition~\ref{prop:DirResolventPbRn+} up to appropriate changes, taking $u:=[(\lambda\I-\Delta)^{-1}\E_{\mathcal{N}} f]_{|_{\RR^n_+}}$, so the result holds by Lemma~\ref{lem:ExtDirNeuRn+}. Still for $h=0$, estimates of higher order follow from applying separately $\nabla'$ which commutes, then noting that $\partial_{x_n}u$ solves a Dirichlet problem so that Proposition~\ref{prop:DirResolventPbRn+} applies.

We focus on existence for $f=0$, $h\neq0$. For now, we assume $-1+\sfrac{1}{p}<s<\sfrac{1}{p}$. Up to a scaling argument, we can assume $|\lambda|=1$. Then, the solution should be given by 
\begin{align*}
    u(\cdot,x_n):=(\lambda \I-\Delta')^{-\sfrac{1}{2}}e^{-x_n(\lambda \I-\Delta')^{\sfrac{1}{2}}}h=e^{-x_n(\lambda \I-\Delta')^{\sfrac{1}{2}}}\Tilde{h}\text{,}\quad x_n\geqslant 0.
\end{align*}
Therefore by Proposition~\ref{prop:DumpedPoissonSemigroup2},
\begin{align*}
    \lVert u\rVert_{\dot{\B}^{s}_{p,q}(\RR^n_+)} + \lVert \nabla^2 u\rVert_{\dot{\B}^{s}_{p,q}(\RR^n_+)} &\lesssim_{p,s,n}^{\mu} \lVert \Tilde{h}\rVert_{{\B}^{s-{\sfrac{1}{p}}}_{p,q}(\partial\RR^{n}_+)}  +\lVert \Tilde{h}\rVert_{{\B}^{s+2-{\sfrac{1}{p}}}_{p,q}(\partial\RR^{n}_+)}\\ &\lesssim_{p,s,n}^{\mu} \lVert \Tilde{h}\rVert_{{\B}^{s+2-{\sfrac{1}{p}}}_{p,q}(\partial\RR^{n}_+)}\\&\lesssim_{p,s,n}^{\mu} \lVert h\rVert_{{\B}^{s+1-{\sfrac{1}{p}}}_{p,q}(\partial\RR^{n}_+)}\\
    &\lesssim_{p,s,n}^{\mu} \lVert h\rVert_{{\L}^{p}(\partial\RR^{n}_+)} + \lVert h\rVert_{\dot{\B}^{s+1-{\sfrac{1}{p}}}_{p,q}(\partial\RR^{n}_+)}.
\end{align*}
As in the proof of Proposition~\ref{prop:DirResolventPbRn+}, interpolation inequalities yield the result. Now for the case of higher regularity estimates, since $\partial_{x_n}u =-e^{-x_n(\lambda \I-\Delta')^{\sfrac{1}{2}}}h$ and $\nabla'u = \nabla'(\lambda \I-\Delta')^{-\sfrac{1}{2}}e^{-x_n(\lambda \I-\Delta')^{\sfrac{1}{2}}}h$
\begin{align*}
    \lVert \partial_{x_n} u\rVert_{\dot{\B}^{s}_{p,q}(\RR^n_+)} + \lVert \partial_{x_n}\nabla^2 u\rVert_{\dot{\B}^{s}_{p,q}(\RR^n_+)} &\lesssim_{p,s,n}^{\mu} \lVert h\rVert_{{\B}^{s-{\sfrac{1}{p}}}_{p,q}(\partial\RR^{n}_+)}  +\lVert h\rVert_{{\B}^{s+2-{\sfrac{1}{p}}}_{p,q}(\partial\RR^{n}_+)}\\ &\lesssim_{p,s,n}^{\mu} \lVert h\rVert_{{\B}^{s+2-{\sfrac{1}{p}}}_{p,q}(\partial\RR^{n}_+)}.
\end{align*}
and
\begin{align*}
    \lVert \nabla' u\rVert_{\dot{\B}^{s}_{p,q}(\RR^n_+)} + \lVert \nabla'\nabla^2 u\rVert_{\dot{\B}^{s}_{p,q}(\RR^n_+)} &\lesssim_{p,s,n}^{\mu} \lVert \nabla' h\rVert_{{\B}^{s+1-{\sfrac{1}{p}}}_{p,q}(\partial\RR^{n}_+)}.
\end{align*}
Then one can conclude similarly by interpolation type inequalities as in Proposition~\ref{prop:DirResolventPbRn+}.
\end{proof}

We state the following proposition, which is an easy consequence, and for which the proof is left to the reader. It will be useful later on.
\begin{proposition}\label{prop:IsomNeuMannHomSobSpaRn+} Let $\mu\in[0,\pi)$. For all $p\in(1,\infty)$, all $s\in(-1+\sfrac{1}{p},\sfrac{1}{p})$, all $\lambda\in\Sigma_\mu$, it holds that
\begin{align*}
    (-\Delta_\mathcal{N})^{-\frac{1}{2}}(\lambda\I-\Delta_\mathcal{N})\,:\,\dot{\H}^{s-1,p}_0\cap\dot{\H}^{s+1,p}(\RR^n_+)\longrightarrow \dot{\H}^{s,p}(\RR^n_+)
\end{align*}
is an isomorphism of Banach spaces.
\end{proposition}

\newpage
\section{Generalized divergence-free and curl-free Sobolev and Besov spaces on Lipschitz domains}\label{Sec:FunctionSpaceTheoryDivFree}

For $\X\in\{\Ccinfty,\,\S,\,\S_0,\,\L^{p},\,\H^{s,p},\,\W^{m,p},\,\B^{s}_{p,q},\,\BesSmo^{s}_{p,\infty},\,\B^{s,0}_{\infty,q},\,{\C}^{m}_{ub},\,{\C}^{m}_{0}\}$, $s\in\RR$,  $p,q\in[1,\infty]$, with $m:=s$ if $s\in\NN$, we define
\begin{align*}
    &\X_{\sigma}(\Omega,\Lambda):= \{\, u\in\X(\Omega,\Lambda)\,:\, \delta u = 0 \text{ in } \mathcal{D}'(\Omega,\Lambda) \,\},\\
    \text{ and } &\X_{\gamma}(\Omega,\Lambda):= \{\, u\in\X(\Omega,\Lambda)\,:\, \d u = 0 \text{ in } \mathcal{D}'(\Omega,\Lambda) \,\},
\end{align*}
and similarly, whenever $\X\notin \{\C^\infty_c,\,\S,\,\S_0\}$, we can define for $\eta\in\{\sigma,\gamma\}$,
\begin{align*}
    \X_{0,\eta}(\Omega,\Lambda):=\X_0(\Omega,\Lambda)\cap\X_{\eta}(\mathbb{R}^n,\Lambda).
\end{align*}
When it is applicable, we also define
\begin{align*}
    \X_{\mathfrak{n}}(\Omega,\Lambda)&:= \{\, u\in\X(\Omega,\Lambda)\,:\, \nu\iprod u_{|_{\partial\Omega}} = 0 \,\},\\
    \text{ and }\,\, \X_{\mathfrak{t}}(\Omega,\Lambda)&:= \{\, u\in\X(\Omega,\Lambda)\,:\, \nu\wedge u_{|_{\partial\Omega}} = 0 \,\},\\
    \text{ and } \X_{\mathcal{D}}(\Omega,\Lambda)&:= \{\, u\in\X(\Omega,\Lambda)\,:\,u_{|_{\partial\Omega}} = 0 \,\},
\end{align*}
where for the latter the trace is considered component-wise.

For $k\in\llb 0,n\rrb$, one has the isometric identification through the Hodge-star operator $\star$
\begin{align*}
    \X_{\sigma}(\Omega,\Lambda^{k})=\,\star \,\X_{\gamma}(\Omega,\Lambda^{n-k}),\quad \text{ and }\quad\X_{\mathfrak{n}}(\Omega,\Lambda^k)=\,\star \,\X_{\mathfrak{t}}(\Omega,\Lambda^{n-k}).
\end{align*}
The same is achieved for the homogeneous function spaces $$\dot{\X}\in \{\dot{\H}^{s,p},\,\dot{\B}^{s}_{p,q},\,\dot{\BesSmo}^{s}_{p,\infty},\,\dot{\B}^{s,0}_{\infty,q},\,\dot{\W}^{m,p},\,\dot{\C}^{m}_{ub},,\,\dot{\C}^{m}_{ub,h},\,\dot{\C}^{m}_{0}\}.$$

\paragraph{The goals and motivations of the current chapter.}

\begin{itemize}
    \item The \textbf{Section~\ref{Sec:DivFreeSpacesRn}} below concerns the study of the function spaces only the case of the whole space $\Omega=\RR^n$. The main questions it aims to address are the same as the ones in the two bullet points below for \textbf{Sections~\ref{Sec:BogovExtOpDivFreeLipDomains:InterpDensity}} and \textbf{\ref{Sec:HodgeLapHodgeDirHodgeLeray}}, but in a simpler and somewhat quite known setting. The results established on $\RR^n$ are independent, some of them are new with such a generality, and will serve as preparatory results for the subsequent sections, mostly Section~\ref{Sec:BogovExtOpDivFreeLipDomains:InterpDensity}.
    
    \item In \textbf{Section~\ref{Sec:BogovExtOpDivFreeLipDomains:InterpDensity}}, we are going to give the \textbf{fundamental results} of \textbf{density} and \textbf{interpolation}. Those are well-known for vector fields, and were proven in several specific cases. In the case of vector fields $\Lambda^{1}\simeq \CC^n$, results like
\begin{align*}
    \L^p_{\mathfrak{n},\sigma}(\Omega,\CC^n)= \overline{\Ccinftydiv(\Omega,\CC^n)}^{\lVert \cdot \rVert_{\L^p(\Omega)}}, \quad\text{ and }\quad
    \W^{1,p}_{0,\sigma}(\Omega,\CC^n)= \overline{\Ccinftydiv(\Omega,\CC^n)}^{\lVert \cdot \rVert_{\W^{1,p}(\RR^n)}},
\end{align*}
for $p\in(1,\infty)$, $\Omega$ to be $\RR^n$, a special Lipschitz domain or a bounded domain are quite well-known but not that much exposed in the literature, see for instance \cite[Section~III.2,~Theorem~III.2.3]{bookGaldi2011} and \cite[Section~2.2,~Lemma~2.2.3]{SohrBook2001}. However, much less is known in the case of other function spaces of fractional order $\H^{s,p}_{0,\sigma}(\Omega,\CC^n)$ and $\W^{s,p}_{0,\sigma}(\Omega,\CC^n)$, $s\in\RR$, $p\in(1,\infty)$, see \cite[Proposition~2.10]{MitreaMonniaux2008} for the case of Bessel potential spaces, $s>-1+\sfrac{1}{p}$. But as far as the knowledge of the authors is concerned, it seems that almost nothing is known in the case of endpoint function spaces $p=1,\infty$, such as $\L^1_{\mathfrak{n},\sigma}$, $\W^{s,1}_{0,\sigma}$, $\C^{m}_{0,0,\sigma}$, $\BesSmo^{s,0,\sigma}_{\infty,\infty,0}=\mathcal{C}^{m,\alpha}_{0,0,\sigma}$, $s=m+\alpha$, $m\in\NN$, $\alpha\in(0,1)$, even when $\Omega=\RR^n$.

\medbreak

Also only few results are known for interpolation of divergence free function spaces with
\begin{align*}
    [\H^{s_0,p_0}_{0,\sigma}(\Omega,\CC^n),\H^{s_1,p_1}_{0,\sigma}(\Omega,\CC^n)]_{\theta} = \H^{s,p}_{0,\sigma}(\Omega,\CC^n), 
\end{align*}
$p_0,p_1\in(1,\infty), \, s_j>-1+\sfrac{1}{p_j},\, j\in\{0,1\}$, $s=(1-\theta)s_0+\theta s_1$, $1/p=(1-\theta)/p_0+\theta/p_1$, $\theta\in(0,1)$, see for instance \cite[Theorem~2.12]{MitreaMonniaux2008}.

\medbreak

Obviously, even less is known in the case of homogeneous function spaces over $\RR^n$, $\Omega=\RR^n_+$ the half-space or $\Omega$ a special Lipschitz domain: the only known results to the best of the authors' knowledge are for the whole and the flat half-space, in the case of real interpolation\cite[Propositions~3.24~\&~3.25]{DanchinHieberMuchaTolk2020}, and some additional results can be deduced from the work \cite{Gaudin2023Hodge}, .  As much as possible, we will try to provide the corresponding interpolation and density results for special and bounded Lipschitz domains in the case of differential forms of arbitrary degree, encompassing then comparable results for (higher-dimensional) curl-free vector fields.

\item In \textbf{Section~\ref{Sec:HodgeLapHodgeDirHodgeLeray}}, provided $\Omega$ is a bounded --at least Lipschitz-- domain, the central question is about to know whether the \textbf{Hodge-Leray projector} 
\[
\PP_\Omega \,:\, \L^p(\Omega,\CC^n) \longrightarrow \L^p_{\mathfrak{n},\sigma}(\Omega,\CC^n)
\]  
remains bounded for $p\neq 2$, but also beyond the solely class of $\L^p$-spaces, including possibly Sobolev spaces $\H^{s,p}$ and Besov spaces $\B^{s}_{p,q}$, possibly allowing $p=1,\infty$ for the latter.

\medbreak

For bounded Lipschitz domains, Fabes, Mendez and Mitrea \cite{FabesMendezMitrea1998} showed boundedness only for $p$ near $2$, while Simader and Sohr \cite{SimaderSohr1992} proved it for all $p\in(1,\infty)$ in $\C^1$-bounded and exterior domains. On less regular geometries the situation deteriorates: Bogovski\u{\i} \cite{Bogovskii1986} gave counterexamples on unbounded sets, and Tolksdorf \cite[Theorem~5.1.10]{TolksdorfPhDThesis2017} obtained boundedness only for restricted $p$-intervals on special Lipschitz domains.  

\medbreak

The extension to fractional Sobolev and Besov scales has been quite well studied. Mitrea--Monniaux \cite[Proposition~2.16]{MitreaMonniaux2008} established the decomposition in $\H^{s,p}$ for $\C^1$-bounded domains provided $s\in(-1+\sfrac{1}{p},\sfrac{1}{p})$, and near the Hilbert scale for bounded Lipschitz domains. In Besov spaces, the only general statement including $\B^{s}_{1,q}$ and $\B^{-s}_{\infty,q}$, $0<s<1$, $q\in[1,\infty]$ is due to Fujiwara and Yamazaki \cite[Theorem~3.1]{FujiwaraYamazaki2007} on $\C^{2,1}$-domains.

\medbreak

Our objective is therefore to establish the boundedness of $\PP_\Omega$ and the corresponding decompositions for a class of \textbf{sufficiently rough} domains in Sobolev and Besov spaces, including the \textbf{fractional and endpoint regimes}. Finally, in order to cover not only vector fields but also higher-degree objects, we will also formulate our results in the general framework of the \textbf{Hodge decomposition} for \textbf{differential forms}. This naturally connects to sharp elliptic regularity for boundary value problems, in particular the Neumann/Hodge Laplacian and the Hodge-Dirac operators.
\end{itemize}

\subsection{The case of the whole space.}\label{Sec:DivFreeSpacesRn}

As prototype results, we need to investigate first density and interpolation results in the case of the whole space $\RR^n$. 

\begin{theorem}\label{thm:DivergenceFreeSpacesRnDensity}Let $p,q\in[1,\infty)$, $s\in\RR$ and $m\in\NN$.
\begin{enumerate}
    \item The space $\S_{0,\sigma}(\RR^n,\Lambda)$ is a dense subspace of 
    \begin{enumerate}
        \item the inhomogeneous functions spaces:  ${\B}^{s,\sigma}_{p,q}(\RR^n,\Lambda)$, ${\BesSmo}^{s,\sigma}_{p,\infty}(\RR^n,\Lambda)$, ${\B}^{s,0,\sigma}_{\infty,q}(\RR^n,\Lambda)$,
        
        \noindent ${\BesSmo}^{s,0,\sigma}_{\infty,\infty}(\RR^n,\Lambda)$, ${\H}^{s,p}_{\sigma}(\RR^n,\Lambda)$, ${\W}^{m,1}_{\sigma}(\RR^n,\Lambda)$, and ${\C}^m_{0,\sigma}(\RR^n,\Lambda)$.
        \item the homogeneous function spaces: $\dot{\B}^{s,\sigma}_{p,q}(\RR^n,\Lambda)$, $\dot{\BesSmo}^{s,\sigma}_{p,\infty}(\RR^n,\Lambda)$, $\dot{\B}^{s,0,\sigma}_{\infty,q}(\RR^n,\Lambda)$,

        \noindent$\dot{\BesSmo}^{s,0,\sigma}_{\infty,\infty}(\RR^n,\Lambda)$, $\dot{\H}^{s,p}_{\sigma}(\RR^n,\Lambda)$, $\dot{\W}^{m,1}_{\sigma}(\RR^n,\Lambda)$, and $\dot{\C}^m_{0,\sigma}(\RR^n,\Lambda)$.
    \end{enumerate}
    
    \item The space $\Ccinftydiv(\RR^n,\Lambda)$ is a dense subspace of
    \begin{enumerate}
        \item the inhomogeneous functions spaces:  ${\B}^{s,\sigma}_{p,q}(\RR^n,\Lambda)$, ${\BesSmo}^{s,\sigma}_{p,\infty}(\RR^n,\Lambda)$, ${\B}^{s,0,\sigma}_{\infty,q}(\RR^n,\Lambda)$,
        
        \noindent${\BesSmo}^{s,0,\sigma}_{\infty,\infty}(\RR^n,\Lambda)$, ${\H}^{s,p}_{\sigma}(\RR^n,\Lambda)$, ${\W}^{m,1}_{\sigma}(\RR^n,\Lambda)$, and ${\C}^m_{0,\sigma}(\RR^n,\Lambda)$;
        \item the homogeneous Besov spaces: $\dot{\B}^{s,\sigma}_{p,q}(\RR^n,\Lambda)$, $s>-n/p'$,  $\dot{\B}^{s,0,\sigma}_{\infty,q}(\RR^n,\Lambda)$, $s>-n$, $\dot{\BesSmo}^{s,0,\sigma}_{\infty,\infty}(\RR^n,\Lambda)$, $s\geqslant-n$, $\dot{\BesSmo}^{s,\sigma}_{p,\infty}(\RR^n,\Lambda)$, $s\geqslant-n/p'$;
        \item the homogeneous Sobolev spaces: $\dot{\H}^{s,p}_{\sigma}(\RR^n,\Lambda)$, $s>-n/p'$, $\dot{\W}^{m,1}_{\sigma}(\RR^n,\Lambda)$ and 
        
        \noindent $\dot{\C}^m_{0,\sigma}(\RR^n,\Lambda)$.
    \end{enumerate}
    \item The space $\Ccinftydiv\cap\dot{\B}^{s}_{p,q}(\RR^n,\Lambda)$ is a dense subspace of $\dot{\B}^{s,\sigma}_{p,q}(\RR^n,\Lambda)$. 
    
    The result remains true if we replace $\dot{\B}^{s}_{p,q}$, by either $\dot{\B}^{s,0}_{\infty,q}$, $\dot{\BesSmo}^{s}_{p,\infty}$, $\dot{\BesSmo}^{s,0}_{\infty,\infty}$ or $\dot{\H}^{s,p}$.

    \item The space $\C_{ub,\sigma}^\infty\cap\dot{\B}^{s}_{\infty,q}(\RR^n,\Lambda)$ is a dense subspace of $\dot{\B}^{s,\sigma}_{\infty,q}(\RR^n,\Lambda)$. 
    
    The result remains true if we replace $\dot{\B}^{s}_{\infty,q}$ by $\dot{\BesSmo}^{s}_{\infty,\infty}$, but also their inhomogeneous counterparts.
\end{enumerate}
Everything still holds, replacing $\sigma$ by $\gamma$.
\end{theorem}

\begin{remark}In each case except the last one, the approximation procedure is universal.
\end{remark}

Before we prove Theorem \ref{thm:DivergenceFreeSpacesRnDensity}, we have to collect several intermediate technical results and to discuss the mapping properties of the (generalized) Leray projector $$\PP_{\RR^n}=\I-\d(-\Delta)^{-1}\delta = \I+ R\wedge( R\iprod \cdot ) = \delta(-\Delta)^{-1}\d =- R \iprod ( R\wedge \cdot ).$$
It can be written in the specific case of vector fields as
\begin{align*}
    \PP_{\RR^n}=\I+\nabla(-\Delta)^{-1}\div = \I+R\langle R,\cdot\rangle = \prescript{t}{}{\curl}(-\Delta)^{-1}\curl,
\end{align*}
where $R= \prescript{t}{}{(R_1,R_2,\ldots,R_n)}$, with $R_k=\partial_{x_k}(-\Delta)^{-1/2}$ the Riesz transform with respect to the $k$-th variable. While $R=\nabla (-\Delta)^{-1/2}$, we set similarly
\begin{align*}
    S:=\nabla'(-\Delta')^{-1/2} = \prescript{t}{}{(S_1, S_2, \ldots, S_{n-1})}.
\end{align*}

We notice that
\begin{align}
    \PP_{\RR^n} \S_{0}(\RR^n,\Lambda)= \S_{0,\sigma}(\RR^n,\Lambda).
\end{align}

In what follows, we set $\X^{s,p}(\RR^n)$ to be any of the following normed spaces
\begin{itemize}
    \item ${\H}^{s,p}(\RR^n)$, $\dot{\H}^{s,p}(\RR^n)$, ${\B}^{s}_{p,q}(\RR^n)$, for $1 <p< \infty$,  $q\in[1, \infty]$, $s\in\RR$; or
    \item $\dot{\B}^{s}_{p,q}(\RR^n)$, for $1\leqslant p,q \leqslant \infty$, $s\in\RR$; or
    \item ${\mathcal{H}}^{1}(\RR^n)$, for $p=1$, $s=0$.
\end{itemize}
Let $\Y^{\alpha,r}$ be any of the above spaces, where we have replaced the indices $(s,p,q)$ and the symbols $\{{\H}^{s,p},\,{\B}^{s}_{p,q},\,{\mathcal{H}}^{1}\}$ by $(\alpha,r,\kappa)$ and $\{{\H}^{\alpha,r},\,{\B}^{\alpha}_{r,\kappa},\,{\mathcal{H}}^{1}\}$.

\begin{lemma}\label{lem:FourierLocalLerayProj}Let $u\in \S'(\RR^n,\Lambda)$, such that $(\dot{\Delta}_ju)_{j\in\ZZ}\subset\L^p(\RR^n,\Lambda)$, for some $p\in[1,\infty]$. For $N\in\NN$, we set
\begin{align*}
    u_{N} = \sum_{|j|\leqslant N} \dot{\Delta}_j u.
\end{align*}

For all $N\in\NN$, $s\in\RR$, $q\in[1,\infty]$, one has $\PP u_N\in \dot{\B}^{s,\sigma}_{p,q}(\RR^n,\Lambda)$, with the uniform estimate
\begin{align*}
    \lVert \PP_{\RR^n} u_N \rVert_{\dot{\B}^{s}_{p,q}(\RR^n)} \less_{p,s,n} \lVert u_N \rVert_{\dot{\B}^{s}_{p,q}(\RR^n)}.
\end{align*}
\end{lemma}

\begin{theorem}\label{thm:RieszTransfRn} Let $p,q,r,\kappa\in[1,\infty]$, and $s,\alpha\in\RR$. Let $\mathcal{T}\in\{\,R_1,\,R_2,\dots,\,R_{n},\,S_1,\,S_2,\ldots ,\,S_{n-1}\,\}$.

Then, if $\X^{s,p}(\RR^n)$ is complete, for all $u\in\X^{s,p}(\RR^n)\cap\Y^{\alpha,r}(\RR^n)$, one has
\begin{align*}
    \lVert \mathcal{T} u \rVert_{\X^{s,p}(\RR^n)}\less_{p,s,n} \lVert u \rVert_{\X^{s,p}(\RR^n)} \quad\text{and}\quad\lVert \mathcal{T} u \rVert_{\Y^{\alpha,r}(\RR^n)}\less_{r,\alpha,n} \lVert u \rVert_{\Y^{\alpha,r}(\RR^n)}.
\end{align*}
Moreover, we also have $\mathcal{T} \dot{\mathfrak{B}} \subset \dot{\mathfrak{B}}$, for $\dot{\mathfrak{B}}\in\{ \dot{\BesSmo}_{\cdot,\infty}^{\cdot},\,\dot{\B}_{\infty,q}^{\cdot,0} \}$ whenever it is well-defined. 
\end{theorem}

\begin{proof}[of Lemma~\ref{lem:FourierLocalLerayProj}~\&~Theorem~\ref{thm:RieszTransfRn}] The proof is standard and for Besov spaces, including endpoints cases, one uses \cite[Lemmas~2.1~\&~2.2]{bookBahouriCheminDanchin}. For Sobolev and Hardy spaces see \cite[Chapter~2,~Theorem~1,~Chapter~3,~Section~1,~\&~Chapter~VII,~Section~3.2,~Corollary~1]{Stein1970}.
\end{proof}

\begin{theorem}\label{thm:BddlerayProjRn}Let $p,q,r,\kappa\in[1,\infty]$, and $s,\alpha\in\RR$.

Provided $\X^{s,p}(\RR^n,\Lambda)$ is complete, for all $u\in\X^{s,p}(\RR^n,\Lambda)\cap\Y^{\alpha,r}(\RR^n,\Lambda)$, one has
\begin{align*}
    \PP_{\RR^n} u\in\X^{s,p}_{\sigma}\cap\Y^{\alpha,r}_{\sigma}(\RR^n,\Lambda),
\end{align*}
with the estimates
\begin{align*}
    \lVert \PP_{\RR^n} u \rVert_{\X^{s,p}(\RR^n)}\less_{p,s,n} \lVert u \rVert_{\X^{s,p}(\RR^n)} \quad\text{and}\quad\lVert \PP_{\RR^n} u \rVert_{\Y^{\alpha,r}(\RR^n)}\less_{r,\alpha,n} \lVert u \rVert_{\Y^{\alpha,r}(\RR^n)}.
\end{align*}
Moreover, we also have $\PP_{\RR^n} \dot{\mathfrak{B}} \subset \dot{\mathfrak{B}}$, for $\dot{\mathfrak{B}}\in\{ \dot{\BesSmo}_{\cdot,\infty}^{\cdot},\,\dot{\B}_{\infty,q}^{\cdot,0} \}$ whenever it is well-defined.
\end{theorem}

\begin{lemma}\label{lem:0meanL1divfree}For all $u\in\L^1_\sigma(\RR^n,\Lambda^\ell)$ a $\ell$-form with $1\leqslant \ell\leqslant n-1$, one has component-wise
\begin{align*}
    \int_{\RR^n} u(x)\,\d x =0.
\end{align*}
\end{lemma}

\begin{proof}For $\e>0$, one has $u_\e:=e^{\e \Delta}u \in\L^1_\sigma(\RR^n)\cap\C^\infty_{0,\sigma}(\RR^n)\subset \L^p_\sigma(\RR^n)$, $1<p<\infty$, with
\begin{align*}
    \int_{\RR^n} u(x)\,\d x = \int_{\RR^n} u_\e(x)\,\d x.
\end{align*}
Because $u_\e\in\L^p_\sigma(\RR^n)$, provided $1<p<\infty$, we can write $$u_\e = \PP u_\e= \delta [e^{\frac{\e}{2} \Delta} (-\Delta)^{-1}\d e^{\frac{\e}{2} \Delta}u]=: \delta v_\e.$$  By construction: by Sobolev embeddings and regularizing properties of the convolution, one can check that $v_\e\in\C^\infty_{0}(\RR^n)$, 
 with $u_\e=\delta v_\e\in\L^1(\RR^n)$, and the equality holds everywhere. Since each component of $u_\e$ is a finite sum of derivatives of order one of components of $v_\e$, by Fubini and the fundamental theorem of calculus, it holds that
\begin{align*}
    \int_{\RR^n} u(x)\,\d x = \int_{\RR^n} u_\e(x)\,\d x = \int_{\RR^n} \delta v_\e(x)\,\d x =0.
\end{align*}
This ends the proof.
\end{proof}

\begin{proof} [of Theorem  \ref{thm:DivergenceFreeSpacesRnDensity}] The proof is partially inspired from \cite[Prop.~2.27]{bookBahouriCheminDanchin} and \cite[Prop.~2.8]{Gaudin2023Lip}. One has $\Ccinftydiv(\RR^n,\Lambda^0)=\Ccinfty(\RR^n,\CC)$, $\S_{0,\sigma}(\RR^n,\Lambda^0)=\S_{0}(\RR^n,\CC)$ and $\S_{0,\sigma}(\RR^n,\Lambda^n)=\Ccinftydiv(\RR^n,\Lambda^n)=\{0\}$. Therefore, we can focus on the case of $\ell$-forms with $1\leqslant \ell\leqslant n-1$. 

\textbf{Step 1:} Point \textit{(i)}. One can check, as in \cite[Prop.~2.8]{Gaudin2023Lip}, that $u_N= \sum\limits_{|j|\leqslant N}\dot{\Delta}_j u\in\C^\infty_{0,\sigma}(\RR^n,\Lambda)$ for $N\in\NN$, satisfies
\begin{align*}
    \lVert u -u_N\rVert_{\B^{s}_{p,q}(\RR^n)}\xrightarrow[N\rightarrow\infty]{}0,
\end{align*}
given $u\in\B^{s,\sigma}_{p,q}(\RR^n,\Lambda)$ and similarly for any other function spaces presented in \textit{(i)-(a)} and \textit{(i)-(b)}.

\textbf{Step 1.1:} However, we have to comment the case of inhomogeneous function spaces when $p=1$.

For $u\in\W^{k,1}_{\sigma}(\RR^n,\Lambda)$, one has in particular $u\in\L^1_{\sigma}(\RR^n,\Lambda)$ which, by Lemma \ref{lem:0meanL1divfree}, implies that $u$ has a $0$ integral over the whole space. By \cite[Prop.~2.8]{Gaudin2023Lip} and \cite[Lem.~E.5.2]{bookHaase2006}, one has
\begin{align*}
    \lVert u -u_N\rVert_{\L^1(\RR^n)}+\lVert \nabla^k u -\nabla^k u_N\rVert_{\L^1(\RR^n)}\xrightarrow[N\rightarrow\infty]{} 0.
\end{align*}
For $u\in\B^{s,\sigma}_{1,q}(\RR^n)$, one has $\Delta_{-1}u=\dot{S}_0u\in \L^1_\sigma(\RR^n)$. Consequently, one can write
\begin{align*}
    u_N = (\dot{S}_0u)_{N} + ([\I-\dot{S}_0]u)_N,
\end{align*}
with, as before, $\left((\dot{S}_0u)_{N}\right)_{N\in\NN}$ converging to $\dot{S}_0u$ in $\L^1_\sigma(\RR^n)$, by Lemma \ref{lem:0meanL1divfree} and \cite[Lem.~E.5.2]{bookHaase2006}, and $\left(([\I-\dot{S}_0]u)_N\right)_{N\in\NN}$ converging to $[\I-\dot{S}_0]u$ in $\dot{\B}^{s,\sigma}_{1,q}(\RR^n)$ by \cite[Prop.~2.8]{Gaudin2023Lip} (or equivalently in ${\B}^{s,\sigma}_{1,q}(\RR^n,\Lambda)$ by \cite[Thm.~6.3.2]{BerghLofstrom1976}). Therefore, for $N\geqslant 2$,
\begin{align*}
    \lVert u -u_N\rVert_{\B^{s}_{1,q}(\RR^n)} &\lesssim_{s,n} \lVert\dot{S}_0u- (\dot{S}_0u)_N\rVert_{\L^1(\RR^n)} + \left(\sum_{j\geqslant 0} 2^{jsq}\lVert \Delta_{j} u 
 - \Delta_{j}u_N\rVert_{\L^1(\RR^n)}^q\right)^\frac{1}{q}\\
 & \lesssim_{s,n}\lVert\dot{S}_0u- (\dot{S}_0u)_N\rVert_{\L^1(\RR^n)} + \left(\sum_{j\geqslant 1} 2^{jsq}\lVert \Delta_{j}  [\I-\dot{S}_0][u-u_N]\rVert_{\L^1(\RR^n)}^q\right)^\frac{1}{q}\\
 &\xrightarrow[N\rightarrow\infty]{} 0.
\end{align*}

\textbf{Step 1.2:} For $\e>0$ fixed, let $N\in \NN$ be such that
\begin{align*}
    \lVert u -u_N\rVert_{\B^{s}_{p,q}(\RR^n)} <\e.
\end{align*}
Following the proof of \cite[Prop.~2.8]{Gaudin2023Lip}, for $M\geqslant N+1$, $R>0$, provided $\Theta\in \C_c^\infty(\RR^n)$, real valued, supported in $\B_2(0)$ and  be such that $\Theta_{|_{\B_1(0)}}=1$, and $\Theta_R:=\Theta(\cdot /R)$, we introduce
\begin{align*}
    u_{N,M}^R:= \PP_{\RR^n}(\dot{S}_{M}-\dot{S}_{-M})[\Theta_{R}u_{N}]\in\S_{0,\sigma}(\RR^n)\text{.}
\end{align*}
Since $u_N\in\C^\infty_{0,\sigma}\cap\W^{k,p}_{\sigma}(\RR^n,\Lambda)$, for all $k\in\NN$, we do have $\PP u_N=u_N$ by \cite[Lem.~2.2]{bookBahouriCheminDanchin}, and since $\dot{\Delta}_{k} u_{N} = 0$, $k\leqslant -M-1$, we have $\dot{S}_{-M}u_{N}=0$, and $\dot{S}_{M}u_{N}=u_{N}$, and then
\begin{align*}
    u_{N,M}^R - u_{N} = \PP_{\RR^n}(\dot{S}_{M}-\dot{S}_{-M})[(\Theta_{R}-1)u_{N}]\text{.}
\end{align*}
By the triangle inequality,
\begin{align*}
    \lVert u -u_{N,M}^R\rVert_{\B^{s}_{p,q}(\RR^n)} < \e + \lVert u_N -u_{N,M}^R\rVert_{\B^{s}_{p,q}(\RR^n)}.
\end{align*}
Now, one applies \cite[Lem.~2.1]{bookBahouriCheminDanchin}, Lemma \ref{lem:FourierLocalLerayProj}, then  \cite[Lem.~2.1]{bookBahouriCheminDanchin} again, so provided $m=\max(0,\lceil s\rceil)+1$, one obtains
\begin{align*}
    \lVert u_N -u_{N,M}^R\rVert_{\B^{s}_{p,q}(\RR^n)} &\sim_{M,s} \lVert u_N -u_{N,M}^R\rVert_{\dot{\B}^{s}_{p,q}(\RR^n)} = \lVert \mathbb{P}(\dot{S}_{M}-\dot{S}_{-M})[(\Theta_{R}-1)u_{N}]\rVert_{\dot{\B}^{s}_{p,q}(\RR^n)}\\
    &\less_{M,s} \lVert (\dot{S}_{M}-\dot{S}_{-M})[(\Theta_{R}-1)u_{N}]\rVert_{\dot{\B}^{s}_{p,q}(\RR^n)}\\
    &\less_{M,s} \lVert (\dot{S}_{M}-\dot{S}_{-M})[(\Theta_{R}-1)u_{N}]\rVert_{{\B}^{s}_{p,q}(\RR^n)}\\
     &\less_{M,s} \lVert (\dot{S}_{M}-\dot{S}_{-M})[(\Theta_{R}-1)u_{N}]\rVert_{{\W}^{m,p}(\RR^n)}\\
     &\less_{M,s}\lVert [(\Theta_{R}-1)u_{N}]\rVert_{{\W}^{m,p}(\RR^n)}.
\end{align*}
The Dominated Convergence Theorem  yields the result as $R$ tends to infinity, since for $R$ sufficiently large,
\begin{align*}
    \lVert u -u_{N,M}^R\rVert_{\B^{s}_{p,q}(\RR^n)} < 2\e.
\end{align*}
The proof for other function spaces, inhomogeneous or homogeneous, is similar.

\textbf{Step 2:} We address the point \textit{(ii)}. For $\e>0$ fixed, let $N,M,R>0$ to be such that
\begin{align*}
    \lVert u -u_{N,M}^R\rVert_{\B^{s}_{p,q}(\RR^n)} < \e.
\end{align*}
For $\eta>0$, we set
\begin{align*}
    u_{N,M}^{R,\eta} = \delta [\Theta(\cdot/\eta)(-\Delta)^{-1}\d]u_{N,M}^{R} \in \Ccinftydiv(\RR^n),
\end{align*}
with $(-\Delta)^{-1}\d \, u_{N,M}^{R}\in\S_0(\RR^n)$. Since $\PP_{\RR^n} u_{N,M}^{R} = u_{N,M}^{R}$, one has
\begin{align}
    u_{N,M}^{R,\eta}-u_{N,M}^{R} = \delta [(\Theta(\cdot/\eta)-1)(-\Delta)^{-1}\d]u_{N,M}^{R}.\label{eq:ApproxDivFreeFactorDeriv}
\end{align}
One can conclude following the arguments in \cite[Prop.~2.8,~Proof,~Step~3]{Gaudin2023Lip} in each case so that for $\eta>0$, large enough, one has
\begin{align*}
    \lVert u -u_{N,M}^{R,\eta}\rVert_{\B^{s}_{p,q}(\RR^n)} &\leqslant \lVert u -u_{N,M}^{R}\rVert_{\B^{s}_{p,q}(\RR^n)}+\lVert u_{N,M}^{R}\ -u_{N,M}^{R,\eta}\rVert_{\B^{s}_{p,q}(\RR^n)}\\
    &< \e +\lVert u_{N,M}^{R}\ -u_{N,M}^{R,\eta}\rVert_{\B^{s}_{p,q}(\RR^n)} <2\e,
\end{align*}
and it goes similarly for other function spaces.

\textbf{Step 3:} The proof of \textit{(iii)} follows from a slight modification of the arguments present in Step 2, according to \cite[Prop.~2.8,~Proof,~Step~4]{Gaudin2023Lip}. Due to \textit{(ii)}, we can assume without loss of generality that $s\leqslant -n/p'$.

Let $\mathfrak{m}\in\NN$ such that $s+2\mathfrak{m}>0$. Since $u_{N,M}^R\in\S_0(\RR^n,\Lambda)$, we do have $(-\Delta)^{-\mathfrak{m}-1}\d\, u_{N,M}^R\in\S_0(\RR^n,\Lambda)$. Therefore, for any $\eta>0$, we introduce
  \begin{align*}
      \tilde{u}_{N,M}^{R,\eta} := \delta(-\Delta)^{\mathfrak{m}} \Theta(\cdot/\eta) (-\Delta)^{-\mathfrak{m}-1}\d\, u_{N,M}^R\in\Ccinftydiv(\RR^n)\cap \dot{\B}^{s}_{p,q}(\RR^n).
  \end{align*}
From this, we can conclude as in Step 2:
\begin{align*}
    \lVert \tilde{u}_{N,M}^{R,\eta} - u_{N,M}^R\rVert_{\dot{\B}^{s}_{p,q}(\RR^n)} \lesssim_{p,\mathfrak{m},n} \lVert [\Theta(\cdot/\eta)-1 ](-\Delta)^{-\mathfrak{m}-1}\d \,u_{N,M}^R\rVert_{\dot{\B}^{s+2\mathfrak{m}}_{p,q}(\RR^n)}  \xrightarrow[\eta\longrightarrow\infty]{}0\text{. }
\end{align*}
This finishes the proof.
\end{proof}

We improve the interpolation result \cite[Proposition~3.24]{DanchinHieberMuchaTolk2020}, removing the completeness assumption, and considering the endpoints. The next result also improves known interpolation results of endpoints such as the ones by Ogawa and Shimizu in \cite{OgawaShimizu2010}. Wee also include some inhomogeneous endpoint function spaces and many others.

We recall here that we deal with $\S'_h$-constructions of the homogeneous Besov spaces which are not complete if \eqref{AssumptionCompletenessExponents} is not satisfied. Hence, standard  proofs and argument for interpolation of standard functions spaces may not apply as such.

To improve the readability of the next result, we omit the mention of the exterior algebra $\Lambda$, writing only $\X_{\sigma}(\RR^n)$ instead of $\X_{\sigma}(\RR^n,\Lambda)$ for differential forms $u$, such that, both, the coefficients satisfy $\big((u_{I})_{I\in\mathcal{I}^{k}_{n}}\big)_{0\leqslant k\leqslant n}\subset \X(\RR^n,\CC)$, and $\delta u =0$.

\begin{theorem}\label{thm:InterpHomSpacesRn}Let $1\leqslant p_0,p_1,p,q,q_0,q_1\leqslant \infty$, $s_0,s_1\in\RR$, such that $s_0\neq s_1$, and for $\theta\in(0,1)$, let
\begin{align*}
    \left(s,\frac{1}{p_\theta},\frac{1}{q_\theta}\right):= (1-\theta)\left(s_0,\frac{1}{p_0},\frac{1}{q_0}\right)+ \theta\left(s_1,\frac{1}{p_1},\frac{1}{q_1}\right)\text{. }
\end{align*}
If $s_0,s_1\in\NN$, we write $m_0:=s_0$ and $m_1:=s_1$.
We have the following interpolation identities with equivalence of norms
\begin{enumerate}
    \item $(\dot{\B}^{s_0,\sigma}_{p,q_0}(\RR^n),\dot{\B}^{s_1,\sigma}_{p,q_1}(\RR^n))_{\theta,q}=(\dot{\H}^{s_0,p}_{\sigma}(\RR^n),\dot{\H}^{s_1,p}_{\sigma}(\RR^n))_{\theta,q}=\dot{\B}^{s,\sigma}_{p,q}(\RR^n)$\label{eq:realInterpHomBspqRn};
    \item $(\dot{\B}^{s_0,0,\sigma}_{\infty,q_0}(\RR^n),\dot{\B}^{s_1,0,\sigma}_{\infty,q_1}(\RR^n))_{\theta,q}=(\dot{\C}^{m_0}_{0,\sigma}(\RR^n),\dot{\C}^{m_1}_{0,\sigma}(\RR^n))_{\theta,q}=\dot{\B}^{s,0,\sigma}_{\infty,q}(\RR^n)$\label{eq:realInterpHomBsinftyqRn};
    \item $(\dot{\W}^{m_0,1}_{\sigma}(\RR^n),\dot{\W}^{m_1,1}_{\sigma}(\RR^n))_{\theta,q}=\dot{\B}^{s,\sigma}_{1,q}(\RR^n)$\label{eq:realInterpHomBspqRnL1};
    \item $(\dot{\C}^{m_0}_{ub,h,\sigma}(\RR^n),\dot{\C}^{m_1}_{ub,h,\sigma}(\RR^n))_{\theta,q}=(\dot{\W}^{m_0,\infty}_{\sigma}(\RR^n),\dot{\W}^{m_1,\infty}_{\sigma}(\RR^n))_{\theta,q}=\dot{\B}^{s,\sigma}_{\infty,q}(\RR^n)$\label{eq:realInterpHomBspqRnLinfty};
    \item $[\dot{\H}^{s_0,p_0}_{\sigma}(\RR^n),\dot{\H}^{s_1,p_1}_{\sigma}(\RR^n)]_{\theta} = \dot{\H}^{s,p_\theta}_{\sigma}(\RR^n)$, if $1<p_0,p_1<\infty$, and \hyperref[AssumptionCompletenessExponents]{$(\mathcal{C}_{s_j,p_j})$} is true for $j\in\{0,1\}$\label{eq:complexInterpHomSobRn};
    \item $[\dot{\B}^{s_0,\sigma}_{p_0,q_0}(\RR^n),\dot{\B}^{s_1,\sigma}_{p_1,q_1}(\RR^n)]_{\theta} = \dot{\B}^{s,\sigma}_{p_\theta,q_\theta}(\RR^n)$, if \hyperref[AssumptionCompletenessExponents]{$(\mathcal{C}_{s_j,p_j,q_j})$} is satisfied for $j\in\{0,1\}$, $q_\theta<\infty$.\label{eq:complexInterpHomBspqRn}
\end{enumerate}
Moreover, 
\begin{itemize}
    \item The identity \ref{eq:realInterpHomBspqRn} still holds if we replace $\dot{\B}^{s_j,\sigma}_{p,q_j}(\RR^n)$ by $\dot{\mathcal{B}}^{s_j,\sigma}_{p,\infty}(\RR^n)$, $j\in\{0,1\}$;
    \item When we replace $(\cdot,\cdot)_{\theta,q}$ by $(\cdot,\cdot)_{\theta}$, all the real interpolation identities still hold with $\dot{\mathcal{B}}^{s,\sigma}_{p,\infty}$ as an output space ($\dot{\mathcal{B}}^{s,0,\sigma}_{\infty,\infty}$ in \ref{eq:realInterpHomBsinftyqRn});
    \item The identity \ref{eq:complexInterpHomBspqRn} remains valid for $q_\theta=\infty$ with an output space $\dot{\BesSmo}^{s_\theta,\sigma}_{p_{\theta},\infty}(\RR^n)$.
    \item Every interpolation result \textbf{remains true for inhomogeneous function spaces}, \textit{i.e.}, replacing
    \begin{align*}
        \{\dot{\W},\,\dot{\H},\,\dot{\B},\,\dot{\BesSmo},\,\dot{\C}\}\text{ by }\{\W,\,\H,\,\B,\,{\BesSmo},\,\C\}.
    \end{align*}
    In this case, one removes the conditions \hyperref[AssumptionCompletenessExponents]{$(\mathcal{C}_{s_j,p_j})$}  and \hyperref[AssumptionCompletenessExponents]{$(\mathcal{C}_{s_j,p_j,q_j})$}  on the exponents that appear respectively in \ref{eq:complexInterpHomSobRn} and \ref{eq:complexInterpHomBspqRn}. However, considering \ref{eq:complexInterpHomBspqRn} for inhomogeneous function spaces,  one needs the additional result
    \begin{align}
        [\L^{p_0}_{\sigma}(\RR^n),\L^{p_1}_{\sigma}(\RR^n)]_{\theta} = \L^{p_\theta}_{\sigma}(\RR^n),\quad 1\leqslant p_0,p_1\leqslant\infty.\label{eq:interpFreeDivLpRn}
    \end{align}
    which we will prove below.
    \item Everything still holds with $\gamma$ replacing $\sigma$.
\end{itemize}
\end{theorem}

\begin{proof} \textbf{Step 1:} We just give the scheme of proof, inspired from the proof of \cite[Theorem~2.12]{Gaudin2023Lip} to prove \ref{eq:realInterpHomBspqRn}-\ref{eq:complexInterpHomSobRn} in the case of homogeneous function spaces\footnote{One cannot use in all cases a direct retraction and co-retraction argument by means of the Leray projection $\PP_{\RR^n}$ to reach the whole result. Indeed, for the homogeneous function spaces, since the spaces of high regularity are not complete, a retraction co-retraction argument by the Leray projection and \cite[Thm.~2.12]{Gaudin2023Lip} cannot be reached, since $\PP_{\RR^n}$ is not well-defined on the whole scale of homogeneous Sobolev and Besov spaces (more precisely: only (weak-$\ast$)-densely defined and (weak-$\ast$)-densely bounded when $s>n/p$, and it cannot be extended by density due the lack of completeness). In the case of inhomogeneous function spaces, the Leray projection does not act boundedly on the whole scale of inhomogeneous Sobolev and Besov spaces because of the endpoint cases $p=1,\infty$.}. Thanks to \cite[Theorem~2.12, Proof, Steps~2~\&~3.2]{Gaudin2023Lip}, it is not difficult to see that the homogeneous Littlewood-Paley decomposition $(\dot{\Delta}_j)_{j\in\ZZ}$, and the reconstruction operator $\Tilde{\Sigma}^{\bullet}$, formally defined as,
\begin{align*}
    \Tilde{\Sigma}^{\bullet}( (w_j)_{j\in\ZZ}) := \sum_{j\in\ZZ} \dot{\Delta}_j[w_{j-1}+w_{j}+w_{j+1}]
\end{align*}
are such that they restrict as bounded operators 
\begin{align}
    \bullet \, &\tilde{\Sigma}^{\bullet}\,:\,\ell^q_s(\ZZ,\L^p_\sigma(\RR^n))\longrightarrow \dot{\B}^{s,\sigma}_{p,q}(\RR^n)\text{, $p,q\in[1,\infty]$, $s\in\RR$ such that \eqref{AssumptionCompletenessExponents},}\nonumber\\
    \bullet \, &(\dot{\Delta}_j)_{j\in\ZZ}\,:\,\dot{\B}^{s,\sigma}_{p,q}(\RR^n)\longrightarrow \ell^q_s(\ZZ,\L^p_\sigma(\RR^n))\text{,  $p,q\in[1,\infty]$, $s\in\RR$,}\label{eq:liftOperatorLinftyRealInterpRn}\\
    \bullet \, &(\dot{\Delta}_j\tilde{\Sigma}^{\bullet})_{j\in\ZZ}\,:\,\ell^q_s(\ZZ,\L^p_\sigma(\RR^n))\longrightarrow \ell^q_s(\ZZ,\L^p_\sigma(\RR^n))\text{, $p,q\in[1,\infty]$, $s\in\RR$.}\nonumber
\end{align}
Therefore, one can reproduce \cite[Theorem~2.12, Proof, Steps~2.2~\&~4]{Gaudin2023Lip} in order to obtain \ref{eq:realInterpHomBspqRn} in the case of homogeneous Besov spaces. The identity for Besov spaces \ref{eq:realInterpHomBsinftyqRn} admits the same proof up to appropriate changes. The identity \ref{eq:realInterpHomBspqRn} in the case of homogeneous Sobolev spaces follows from itself in the specific case of homogeneous Besov spaces, and the well-known embedding
\begin{align*}
    \dot{\B}^{s_j,\sigma}_{p,1}(\RR^n) \hookrightarrow\dot{\H}^{s_j,p}_{\sigma}(\RR^n) \hookrightarrow \dot{\B}^{s_j,\sigma}_{p,\infty}(\RR^n), \text{ provided }s_j\in\RR,\,\,j\in\{0,1\},\,p\in[1,\infty].
\end{align*}
The same argument applies to obtain \ref{eq:realInterpHomBsinftyqRn} in the case of the spaces of continuous and smooth functions, as well as \ref{eq:realInterpHomBspqRnL1} and \ref{eq:realInterpHomBspqRnLinfty}.

\textbf{Step 2:} We prove \ref{eq:complexInterpHomSobRn} and \ref{eq:complexInterpHomBspqRn} in the case of homogeneous function spaces. By Theorem \ref{thm:BddlerayProjRn}, since $1<p_0,p_1<\infty$ and $s_j<n/p_j$, for $j\in\{0,1\}$, it turns out that the couple $\{\dot{\H}^{s_0,p_0}_{\sigma}(\RR^n),\dot{\H}^{s_1,p_1}_{\sigma}(\RR^n)\}$ is a retract of $\{\dot{\H}^{s_0,p_0}(\RR^n),\dot{\H}^{s_1,p_1}(\RR^n)\}$ by means of the Leray projection $\PP_{\RR^n}$.

For the same reason, for all $p_j,q_j\in[1,\infty]$, $s_j\in\RR$, such that \hyperref[AssumptionCompletenessExponents]{$(\mathcal{C}_{s_j,p_j,q_j})$} is satisfied, $j\in\{0,1\}$, the couple $\{\dot{\B}^{s_0,\sigma}_{p_0,q_0}(\RR^n),\dot{\B}^{s_1,\sigma}_{p_1,q_1}(\RR^n)\}$ is a retract of $\{\dot{\B}^{s_0}_{p_0,q_0}(\RR^n),\dot{\B}^{s_1}_{p_1,q_1}(\RR^n)\}$ by means of the Leray projection $\PP_{\RR^n}$.

In both cases, \ref{eq:complexInterpHomSobRn} and \ref{eq:complexInterpHomBspqRn} are then consequences of \cite[Thm.~2.12]{Gaudin2023Lip}.

\textbf{Step 3:} We give a scheme of proof to prove \ref{eq:realInterpHomBspqRn}-\ref{eq:realInterpHomBspqRnLinfty} in the case of inhomogeneous function spaces. We just indicate how to prove \ref{eq:realInterpHomBspqRn} in this specific case, the remaining ones admit arguments similar to those presented in Step 1.

The inhomogeneous Littlewood-Paley decomposition $({\Delta}_j)_{j\geqslant -1}$, and the reconstruction operator $\Tilde{\Sigma}$, formally defined as,
\begin{align*}
    \Tilde{\Sigma}( w_{-1},(w_j)_{j\geqslant 0}) :=  \Delta_{-1}[w_{-1}+w_0]+\sum_{j\geqslant 0} {\Delta}_j[w_{j-1}+w_{j}+w_{j+1}]
\end{align*}
are such that, for all $p,q\in[1,\infty]$, $s\in\RR$, one has the bounded operators 
\begin{align*}
    [\Delta_{-1},({\Delta}_j)_{j\in\NN}]\,:\,{\B}^{s,\sigma}_{p,q}(\RR^n)\longrightarrow \L^p_\sigma(\RR^n)\times\ell^q_s(\NN,\L^p_\sigma(\RR^n)),\\
    \tilde{\Sigma}\,:\,\L^p_\sigma(\RR^n)\times\ell^q_s(\NN,\L^p_\sigma(\RR^n))\longrightarrow {\B}^{s,\sigma}_{p,q}(\RR^n),
\end{align*}
such that $\tilde{\Sigma}[\Delta_{-1},({\Delta}_j)_{j\in\NN}] =\mathrm{I}$ on $\S'(\RR^n)$. Hence, the couple $\{{\B}^{s_0,\sigma}_{p,q_0}(\RR^n),{\B}^{s_1,\sigma}_{p,q_1}(\RR^n)\}$ is a retract of the couple $\{\L^p_\sigma(\RR^n)\times\ell^{q_0}_{s_0}(\NN,\L^p_\sigma(\RR^n)),\L^p_\sigma(\RR^n)\times\ell^{q_1}_{s_1}(\NN,\L^p_\sigma(\RR^n))\}$ through the Littlewood-Paley decomposition. Thus, the result follows from \cite[Thm.~5.6.1]{BerghLofstrom1976}.

The remaining identities \ref{eq:realInterpHomBspqRn}-\ref{eq:realInterpHomBspqRnLinfty}, as in Step 1, admit either a similar proof or can be deduced from embeddings such as
\begin{align*}
    {\B}^{k_j,\sigma}_{p,1}(\RR^n) \hookrightarrow{\W}^{k_j,p}_{\sigma}(\RR^n) \hookrightarrow {\B}^{k_j,\sigma}_{p,\infty}(\RR^n), \text{ provided }k_j\in\NN,\,\,j\in\{0,1\},\,p\in[1,\infty].
\end{align*}

\textbf{Step 4:} We provide the arguments to obtain \ref{eq:complexInterpHomSobRn} and \ref{eq:complexInterpHomBspqRn} in the case of homogeneous function spaces. First, \ref{eq:complexInterpHomSobRn} can be obtained by the exact same argument in Step 2.

Second, to obtain \ref{eq:complexInterpHomBspqRn}, one has to investigate the truthfulness of \eqref{eq:interpFreeDivLpRn} first. The task is about to include the cases $p_j=1,\infty$ for $j=0$ or $j=1$, otherwise one could proceed as in Step 2 for the identity \ref{eq:complexInterpHomSobRn}.

\textbf{Step 4.1:} We prove the identity \eqref{eq:interpFreeDivLpRn}. Let $1<q<\infty$, we recall that one has the standard complex interpolation identity
\begin{align*}
    [\mathcal{H}^{1}(\RR^n),\L^q(\RR^n)]_{\theta} = \L^{\frac{q}{1+\theta(q-1)}}(\RR^n),\, \forall\theta\in(0,1),
\end{align*}
where $\mathcal{H}^{1}(\RR^n)$ is the standard Hardy space. But, by Theorem \ref{thm:BddlerayProjRn}, the couple $\{\mathcal{H}^{1}_\sigma(\RR^n),\L^q_\sigma(\RR^n)\}$ is a rectract of the couple $\{\mathcal{H}^{1}(\RR^n),\L^q(\RR^n)\}$, therefore
\begin{align*}
    [\mathcal{H}^{1}_\sigma(\RR^n),\L^q_\sigma(\RR^n)]_{\theta} = \L^{\frac{q}{1+\theta(q-1)}}_\sigma(\RR^n),\, \forall\theta\in(0,1).
\end{align*}
Since $\mathcal{H}^{1}_\sigma(\RR^n)\hookrightarrow \L^1_\sigma(\RR^n)$, it holds that, for $\theta\in(0,1)$,
\begin{align*}
    \L^{\frac{q}{1+\theta(q-1)}}_\sigma(\RR^n) = [\mathcal{H}^{1}_\sigma(\RR^n),\L^q_\sigma(\RR^n)]_{\theta} &\hookrightarrow [\L^1_\sigma(\RR^n),\L^q_\sigma(\RR^n)]_{\theta} \\&\hookrightarrow [\L^1(\RR^n),\L^q(\RR^n)]_{\theta}=\L^{\frac{q}{1+\theta(q-1)}}(\RR^n).
\end{align*}
Thus, we have obtained
\begin{align*}
    \L^{\frac{q}{1+\theta(q-1)}}_\sigma(\RR^n)\hookrightarrow[\L^1_\sigma(\RR^n),\L^q_\sigma(\RR^n)]_{\theta}\subset (\L^1_\sigma(\RR^n)+\L^q_\sigma(\RR^n))\cap \L^{\frac{q}{1+\theta(q-1)}}(\RR^n) \hookrightarrow \L^{\frac{q}{1+\theta(q-1)}}_\sigma(\RR^n),
\end{align*}
so that, we did prove for all $\theta\in(0,1)$, all $q\in(1,\infty)$,
\begin{align}\label{eq:ProofComplInterpRnDivFreeL1}
    \L^{\frac{q}{1+\theta(q-1)}}_\sigma(\RR^n) = [\L^1_\sigma(\RR^n),\L^q_\sigma(\RR^n)]_{\theta}
\end{align}

Due to the fact $\L^{1}_{\sigma}(\RR^n)=\overline{\L^1\cap\L^q_\sigma(\RR^n)}^{\lVert\cdot\rVert_{\L^1(\RR^n)}}$, by \cite[Theorem~4.2.2]{BerghLofstrom1976}, it holds
\begin{align*}
    \L^{\frac{q}{1+\theta(q-1)}}_\sigma(\RR^n) = [\L^1_\sigma(\RR^n),\L^q_\sigma(\RR^n)]_{\theta} = [\L^1(\RR^n),\L^q_\sigma(\RR^n)]_{\theta}.
\end{align*}
By duality applied to \eqref{eq:ProofComplInterpRnDivFreeL1}, thanks to \cite[Corollary~4.5.2]{BerghLofstrom1976}, it holds
\begin{align*}
    \L^{\frac{q'}{1-\theta}}_\sigma(\RR^n) =  [\L^{q'}_\sigma(\RR^n),\L^\infty(\RR^n)]_{\theta}.
\end{align*}
Since it is clear that $\overline{\L^\infty\cap\L^{q'}_\sigma(\RR^n)}^{\lVert\cdot\rVert_{\L^\infty(\RR^n)}} \hookrightarrow\L^{\infty}_{\sigma}(\RR^n)$, by \cite[Theorem~4.2.2]{BerghLofstrom1976} again, it holds
\begin{align*}
    \L^{\frac{q'}{1-\theta}}_\sigma(\RR^n) =  [\L^{q'}_\sigma(\RR^n),\L^\infty(\RR^n)]_{\theta} &\hookrightarrow[\L^{q'}_\sigma(\RR^n),\L^\infty_\sigma(\RR^n)]_{\theta}\\
    &\hookrightarrow (\L^{q'}_\sigma(\RR^n)+\L^\infty_\sigma(\RR^n))\cap[\L^{q'}(\RR^n),\L^\infty(\RR^n)]_{\theta}\\
    &\hookrightarrow (\L^{q'}_\sigma(\RR^n)+\L^\infty_\sigma(\RR^n))\cap\L^{\frac{q'}{1-\theta}}(\RR^n)\\
    &\hookrightarrow \L^{\frac{q'}{1-\theta}}_\sigma(\RR^n).
\end{align*}
So that, up to replace $q'$ by $q$, we did obtain for all $q\in(1,\infty)$, all $\theta\in(0,1)$:
\begin{align}\label{eq:ProofComplInterpRnDivFreeLinfty}
    \L^{\frac{q}{1-\theta}}_\sigma(\RR^n) = [\L^q_\sigma(\RR^n),\L^\infty_\sigma(\RR^n)]_{\theta}.
\end{align}
Thanks to  \eqref{eq:ProofComplInterpRnDivFreeL1} and \eqref{eq:ProofComplInterpRnDivFreeLinfty}, one can apply Wolff's reiteration theorem \cite[Theorem~2]{Wolff1982} exactly as for the proof of \cite[Theorem~3]{Wolff1982}, in order to obtain for all $\theta\in(0,1)$:
\begin{align*}
    \L^{\frac{1}{1-\theta}}_\sigma(\RR^n) = [\L^1_\sigma(\RR^n),\L^\infty_\sigma(\RR^n)]_{\theta}.
\end{align*}
Thereby, the standard reiteration Theorem \cite[Theorem~4.6.1]{BerghLofstrom1976} yields \eqref{eq:interpFreeDivLpRn} in its entirety.

\textbf{Step 4.2:} The identity \ref{eq:complexInterpHomBspqRn} in the case of inhomogeneous function spaces follows from the retraction argument given in Step 3, the previous Step 4.1 and \cite[Theorem~5.6.3]{BerghLofstrom1976}.
\end{proof}

\subsection{The study of function spaces on Lipschitz domains}\label{Sec:BogovExtOpDivFreeLipDomains:InterpDensity}

We want to carry over density and interpolation results  for function spaces from the whole space to Lipschitz domains. To do so, we take advantage of several operators such as a divergence preserving extension operator, or Bogovski{\u{\i}}-like potential operators.

\subsubsection{Potential operators \textit{à la} Poincaré and Bogovski{\u{\i}} in Lipschitz domains}\label{subsec:PoincBog}

\paragraph{Overall picture.} The resolution of divergence-type equations is a cornerstone in the analysis of incompressible fluid mechanics and the theory of elliptic PDEs. Given a bounded domain $\Omega \subset \mathbb{R}^n$ and a function $f \in \L^p(\Omega)$
we are looking for a solution vector field $\uu \in \W^{1,p}_0(\Omega, \CC^n)$ satisfying
\begin{align}\label{eq:DivEqBogov}
    \div\, \uu =f,\qquad\text{ in } \Omega,
\end{align}
which is possible only when $f$ has vanishing mean value
\begin{align*}
    \int_\Omega f(x) \, \d x = 0.
\end{align*}

More precisely, two classes of explicit constructions of solution operators have emerged as fundamental tools:
\begin{itemize}
    \item \textbf{Bogovski{\u{\i}} operators}, which provide an explicit, bounded right-inverse to the divergence on suitable domains;

    Let $\Omega \subset \mathbb{R}^n$ be a bounded domain that is star-shaped with respect to a point $x_0 \in \Omega$, assuming without loss of generality $x_0 = 0$. It has been proved by Bogovski{\u{\i}} in \cite{Bogovskii1979}, that for any $1 < p < \infty$, there exists a linear operator 
\begin{align*}
    \mathcal{B}\,:\, \L^p(\Omega) \longrightarrow \W_0^{1,p}(\Omega, \CC^n),
\end{align*}
such that:
\begin{align*}
    \div \mathcal{B}(g) = g, \quad \text{for all } g \in\L^p(\Omega) \text{ such that }\int_{\Omega} g(x) \,\d x =0,
\end{align*}
solving continuously the equation with the bound
\begin{align*}
    \lVert \mathcal{B}(g)\rVert_{\W^{1,p}(\Omega)} \lesssim_{p,n,\Omega} \lVert g\rVert_{\L^p(\Omega)},\quad \text{ for all } g \in\L^p(\Omega).
\end{align*}
Note also that such a solution to \eqref{eq:DivEqBogov} is not unique.

When $\Omega$ is star-shaped with respect to $0 \in \Omega$, assuming without loss of generality $\overline{\B_1(0)}\subset \Omega$, then a Bogovski{\u{\i}} operator may be defined for $f\in\Ccinfty(\Omega)$ explicitly as:
\begin{align}\label{eq:BogovskiiOpDiv}
    \mathcal{B}(f)(x) := \int_{\B_1(0)}\int_{1}^{\infty} \eta(y)(x-y) t^{n}f(y+t(x-y)) \frac{\d t}{t} \d y
\end{align}
where $\eta \in \Ccinfty(\B_1(0))$ is a cutoff function satisfying $\int_{\B_1(0)} \eta(\tau)\,\d\tau = 1$. With such a definition, one obtains
\begin{align*}
    \div \mathcal{B}(f) = f - \bigg(\int_{\Omega} f(x)\d x\bigg)\eta, \quad \text{for all }f \in\L^p(\Omega),
\end{align*}
up to consider the extension of $f$ to the whole space by $0$ in \eqref{eq:BogovskiiOpDiv}.

The case of general bounded Lipschitz domains $\Omega$, can be obtained by a star-shaped covering of the domain.

One can also prove boundedness on $\W^{m,p}_0(\Omega)$, $m\geqslant 2$, $1<p<\infty$, and one has  boundedness on the negative and fractional Sobolev scales $\W^{s,p}_0(\Omega)$ and $\H^{s,p}_0(\Omega)$, $s>-2+\sfrac{1}{p}$, for such kind of operators by Gei{\ss}ert, Heck and Hieber \cite{GesseirtHeckHieber}.

One can also consider similar problems such as
\begin{align}\label{eq:DivEqBogovGen1}
\nabla u  =\ff,\qquad\text{ in } \Omega,\quad\ff\in\L^p(\Omega,\CC^n),
\end{align}
with unknown $u\in\W^{1,p}_0(\Omega)$, as well as
\begin{align}\label{eq:DivEqBogovGen2}
\curl \uu  =\ff,\qquad\text{ in } \Omega,\quad\ff\in\L^p(\Omega,\CC^3),
\end{align}
with unknown $\uu\in \W^{1,p}_0(\Omega,\CC^3)$. Note that both have been also considered by Bogovski{\u{\i}} in \cite{Bogovskii1980} for \eqref{eq:DivEqBogovGen1}, and Borchers and Sohr dealt with \eqref{eq:DivEqBogovGen2} in \cite{BorchersSohr1990}, with similar solution operators. Each requires additionnal necessary, but not sufficient, compatibility conditions such as $\curl \ff =0$ and $\nu\times \ff_{|_{\partial\Omega}}=0$ for \eqref{eq:DivEqBogovGen1}, and $\div \ff =0$ with $\nu\cdot \ff_{|_{\partial\Omega}}=0$ for \eqref{eq:DivEqBogovGen2} in order to provide the existence of a solution. This problem rather belongs to the realm of Poincaré-type results.
    
\item \textbf{Poincaré-type potential operators}, which pseudo-invert the exterior derivative $\d$ on closed differential forms, or co-closed differential forms if one considers the interior derivative $\delta$ instead. 

Indeed, it turns out that \eqref{eq:DivEqBogov}, \eqref{eq:DivEqBogovGen1} and \eqref{eq:DivEqBogovGen2} can be written in an unified way, for $k\in\llb 0,n\rrb$, $1<p<\infty$, either in exterior derivative form
\begin{align}\label{eq:DivEqBogovGen1DiffForm}
\d u  &= f,\qquad\text{ in } \Omega,\quad f\in\L^p(\Omega,\Lambda^k),\\
u_{|_{\partial\Omega}}  &= 0,\qquad\text{ on } \partial\Omega,\nonumber
\end{align}
with unknown $u\in\W^{1,p}_0(\Omega,\Lambda^{k-1})$, as well as in interior derivative form
\begin{align}\label{eq:DivEqBogovGen2DiffForm}
\delta u  &= f,\qquad\text{ in } \Omega,\quad f\in\L^p(\Omega,\Lambda^k),\\
 u_{|_{\partial\Omega}}  &= 0,\qquad\text{ on } \partial\Omega,\nonumber
\end{align}
with unknown $u\in\W^{1,p}_0(\Omega,\Lambda^{k+1})$, and both \eqref{eq:DivEqBogovGen1DiffForm} and \eqref{eq:DivEqBogovGen2DiffForm} are equivalent through the Hodge-star duality. Note that a necessary condition for the solvability of both equations is the membership of $f$ respectively to $\R_p(\underline{\d},\Omega,\Lambda^k)$ or $\R_p(\underline{\delta},\Omega,\Lambda^k)$, which means that, by definition,  $f$ is respectively an exact or a co-exact $k$-form.

A first study of potential operators at the differential form level has been achieved by Mitrea, Mitrea and Monniaux in \cite[Section~3]{MitreaMitreaMonniaux2008} including a treatment in the case of Triebel-Lizorkin and Besov spaces $\mathrm{F}^{s}_{p,q}$, $\B^{s}_{p,q}$, $1<p,q<\infty$, $s\in(-1+\sfrac{1}{p},\infty)$ for star-shaped bounded subdomains of a compact manifold with Lipschitz boundary.

Shortly after, Costabel and M${}^\text{c}$Intosh, in \cite{CostabelMcIntosh2010}, did provide a similar construction for arbitrary bounded Lipschitz domains $\Omega$ of the Euclidean space and they proved that both the Poincaré and the Bogovski{\u{\i}} operators are given by pseudo-differential operators of order $-1$. This allowed them to state the result for the whole Besov and Triebel-Lizorkin scales, $\B^{s}_{p,q}$ and $\mathrm{F}^{s}_{p,q}$, $0<p,q\leqslant\infty$, $s\in\RR$, with $p,q<\infty$ in the latter case. However, in the case of arbitrary bounded Lipschitz domains, not necessarily star-shaped, the potential operators are only given up to a compact perturbation in the sense that
\begin{align}\label{eq:DivEqBogovGen1DiffFormkompact}
\d \T u  &= u - \L u,\qquad\text{ in } \Omega,\quad u\in\R_p(\underline{\d},\Omega,\Lambda^k)
\end{align}
where $\T\,:\,\L^p(\Omega,\Lambda^k)\longrightarrow\W^{1,p}_0(\Omega,\Lambda^{k-1})$ is the "almost" potential operator and the "perturbation" $\L\,:\,\L^p(\Omega,\Lambda^k)\longrightarrow \C^{\infty}_{0,0}(\Omega,\Lambda)\subset \L^p(\Omega,\Lambda^k)$ has finite dimensional range (and is in particular compact).
\end{itemize}

Similar considerations can be made on unbounded domains such as the half-space $\Omega=\RR^n_+$, or $\Omega$ being a special Lipschitz domain. However, in this case, due to scaling issues, provided $p\in(1,\infty)$, if one considers a problem like
\begin{align*}
\delta u  &= f,\qquad\text{ in } \Omega,\quad f\in\overline{\R_p(\underline{\delta},\Omega,\Lambda^k)}^{\lVert \cdot\rVert_{\L^p(\Omega)}},\\
 u_{|_{\partial\Omega}}  &= 0,\qquad\text{ on } \partial\Omega,\nonumber
\end{align*}
it holds that necessarily the solution lies in a homogeneous function space, \textit{i.e.} $u\in \dot{\W}^{1,p}_0(\Omega,\Lambda^{k-1})$, more precisely
\begin{align*}
    u =\T f,\quad\text{ with } \T\,:\,\dot{\H}^{s,p}_0(\Omega,\Lambda^k)\longrightarrow\dot{\H}^{s+1,p}_0(\Omega,\Lambda^{k-1}),\, s\in\RR.
\end{align*}

In this case, but also for other scales of function spaces, this problem has been solved by Costabel, M${}^\text{c}$Intosh and Taggart in \cite{CostabelTaggartMcIntosh2013}, including even homogeneous Triebel-Lizorkin and Besov spaces, $\dot{\mathrm{F}}^{s}_{p,q}$ and $\dot{\B}^{s}_{p,q}$, $s\in\RR$, $p,q\in[1,\infty]$.

However, all the constructions highlighted in this short review do not fit exactly our purpose, due to several technical obstructions. For instance, and this is not exhaustive:
\begin{itemize}
    \item One may want $\L$ in \eqref{eq:DivEqBogovGen1DiffFormkompact} to vanish on $\R_{p}(\underline{\d},\Omega,\Lambda)$, so that \eqref{eq:DivEqBogovGen1DiffFormkompact} becomes $\d \T u =u$, and then $\T$ becomes a "true" potential operator;
    \item The construction of potential operators for special Lipschitz domains involving homogeneous function spaces in \cite{CostabelTaggartMcIntosh2013} is achieved given a construction of homogeneous function spaces modulo polynomials. This is not a suitable framework for non-linear and boundary value problems, as alluded to in the introduction.
\end{itemize}

The goal of this subsection is to provide a modified and revised version of the presented potential operators in order to reach endpoint function spaces and homogeneous function spaces on special Lipschitz domains with the desired properties.

\paragraph{On bounded Lipschitz domains.}

The next Theorem provides Poincaré and Bogovski{\u{\i}} type operators yielding bounded projections onto the range and the null space of $\delta$ and $\underline{\delta}$. Later, its main purpose will be to carry over the interpolation properties.

\begin{theorem}\label{thm:PotentialOpBddLipDom}Let $\Omega$ be a bounded Lipschitz domain. There exists linear operators $\mathcal{B}^\sigma$, $\K^{\sigma}$ and $\mathcal{B}^{0,\sigma}$, $\K^{0,\sigma}$ such that for any $p,q\in[1,\infty]$, $s\in\mathbb{R}$, for ${\X}^{s,p}\in \{\,{\H}^{s,p},\, {\B}^s_{p,q},\,{\BesSmo}^{s}_{p,\infty}\,\}$, assuming $1<p<\infty$ if $\X=\H$, one has 
\begin{enumerate}[label=($\roman*$)]
    \item The mapping properties, provided $k\in\llb 0,n\rrb$,
    \begin{align*}
        &\mathcal{B}^{\sigma}\,:\,\X^{s,p}(\Omega,\Lambda^k)\longrightarrow \X^{s+1,p}(\Omega,\Lambda^{k+1})\\
        \text{ and }\quad&\K^{\sigma}\,:\,\X^{s,p}(\Omega,\Lambda^k)\longrightarrow  \C^{\infty}(\overline{\Omega},\Lambda^{k})\subset \X^{s+1,p}(\Omega,\Lambda^k);\\
        &\mathcal{B}^{0,\sigma}\,:\,\X^{s,p}_0(\Omega,\Lambda^k)\longrightarrow \X^{s+1,p}_0(\Omega,\Lambda^{k+1})\\
        \text{ and }\quad&\K^{0,\sigma}\,:\,\X^{s,p}_0(\Omega,\Lambda^k)\longrightarrow  \C^{\infty}_{0,0}({\Omega},\Lambda^{k})\subset \X^{s+1,p}_0(\Omega,\Lambda^k);
    \end{align*}
    where $\K^{\sigma}$ and $\K^{0,\sigma}$ have finite-dimensional range with dimension depending on $k,n$ and $\Omega$, but independent of $s\in\RR$, $p,q\in[1,\infty]$;
    \item The identities
    \begin{align*}
        v &= \delta \mathcal{B}^{\sigma} v + \mathcal{B}^{\sigma}\delta v +\K^{\sigma}v&&\text{ with }\delta\K^{\sigma} v=0,\quad \text{ for all } v \in\X^{s,p}(\Omega,\Lambda^k),\\
        u &= \delta \mathcal{B}^{0,\sigma} u + \mathcal{B}^{0,\sigma}\delta u +\K^{0,\sigma}u&&\text{ with }\delta\K^{0,\sigma} u=0,\quad \text{ for all } u \in\X^{s,p}_0(\Omega,\Lambda^k),
    \end{align*}
    are valid;
    \item Additionally, $\K^{\sigma} \delta =0$ on $\X^{s,p}(\Omega,\Lambda)$ and $\K^{0,\sigma} \delta  =0$ on $\X^{s,p}_0(\Omega,\Lambda)$.
\end{enumerate}
Furthermore,
\begin{itemize}
    \item One has $\mathcal{B}^{\sigma}(\C^{\infty}(\overline{\Omega},\Lambda))\subset \C^{\infty}(\overline{\Omega},\Lambda)$ and $\mathcal{B}^{0,\sigma}(\C^{\infty}_{0,0}({\Omega},\Lambda))\subset \C^{\infty}_{0,0}({\Omega},\Lambda)$;
    \item A similar result holds for $(\d,\gamma,\Lambda^{k-1})$ replacing $(\delta,\sigma,\Lambda^{k+1})$. 
\end{itemize}
\end{theorem}

\begin{proof}We deal only with the case of the pair of operators $(\mathcal{B}^{0,\sigma},\K^{0,\sigma})$ on $\X^{s,p}_0(\Omega,\Lambda)$. The other case admits a similar and shorter proof.

\textbf{Step 1:} Construction.  We follow the strategy exhibited in \cite[Section~4,~p.~1725]{McintoshMonniaux2018}. By \cite[Theorem~4.6,~Remark~4.12]{CostabelMcIntosh2010}, up to consider the conjugation by the Hodge-Star operator, there exists linear operators $\T,\L$ with the mapping properties
\begin{align}
        &\T\,:\,\X^{s,p}_0(\Omega,\Lambda^k)\longrightarrow \X^{s+1,p}_0(\Omega,\Lambda^{k+1})\\
        \text{ and }\quad&\L\,:\,\X^{s,p}_0(\Omega,\Lambda^k)\longrightarrow  \C^{\infty}_{0,0}({\Omega},\Lambda^{k})\subset \X^{s+1,p}_0(\Omega,\Lambda^k),
\end{align}
the identity
\begin{align}\label{eq:ProofBogovUnchangedOp}
    \I = \delta \T +\T\delta + \L, \quad\text{ on } \X^{s,p}_0(\Omega,\Lambda),
\end{align}
and where $\L$ has finite-dimensional range with dimension depending on $k,n$ and $\Omega$, but independent of $s\in\RR$, $p,q\in[1,\infty]$ by \cite[Theorem~4.9,~Remark~4.10~\&~Theorem~1.1]{CostabelMcIntosh2010}. One has the commutation relation
\begin{align}\label{eq:ProofBogovCommutationRelationCompact}
        \delta \L = \L\delta, \quad\text{ on } \X^{s,p}_0(\Omega,\Lambda),
\end{align}
with the additionnal bounded mapping property, provided that $\mathcal{E}_0$ is the extension operator from $\Omega$ to the whole space by $0$,
\begin{align}\label{eq:ProofBogovCompactOpBddL2Xsp}
    \L\mathcal{E}_0\,:\,\L^{2}(\Omega,\Lambda) \longrightarrow \X^{s,p}_0(\Omega,\Lambda).
\end{align}
due to the fact that $\L$ is given by a pseudo-differential operator of order $-\infty$, see \cite[Theorem~4.6]{CostabelMcIntosh2010}. Therefore, following \cite[Section~4,~p.~1725]{McintoshMonniaux2018}, with the notations introduced in Proposition~\ref{prop:FullHodgedecompL2bdd}, we set
\begin{align*}
    \mathcal{B}^{0,\sigma} := \T + \L \T + \L \mathcal{E}_0 \mathcal{R}^{\sigma}  \R_{\Omega}\L, \quad\text{ and }\quad \K^{0,\sigma}:= \L \mathcal{E}_0 \P_{\H}\R_{\Omega}\L
\end{align*}
where  $\R_{\Omega}$ is the restriction operator from the whole space to $\Omega$. By the finite dimensional range of $\L$, equivalence of all norms in finite dimensions, and the boundedness properties of $\T$ and $\L$ \eqref{eq:ProofBogovUnchangedOp} and \eqref{eq:ProofBogovCompactOpBddL2Xsp}, and of $\mathcal{R}^{\sigma}$ given in Proposition~\ref{prop:FullHodgedecompL2bdd}, it holds that
\begin{align*}
    \lVert \mathcal{B}^{0,\sigma} u \rVert_{\X^{s+1,p}_0(\Omega)} \lesssim_{p,s,n,\Omega}\lVert  u \rVert_{\X^{s,p}_0(\Omega)} \text{ and } \lVert \K^{0,\sigma} u \rVert_{\X^{s+m,p}_0(\Omega)} \lesssim_{p,s,n,\Omega}^m\lVert  u \rVert_{\X^{s,p}_0(\Omega)}
\end{align*}
for all $m\in\NN$, all $u\in\X^{s,p}_0(\Omega,\Lambda)$.

\textbf{Step 2:} We check all the identities. Up to apply real interpolation, one can assume first that $q<\infty$, so that we can argue by strong density of $\Ccinfty(\Omega,\Lambda)$ in $\X^{s,p}_0(\Omega,\Lambda)$, since $\Omega$ is bounded. For $u\in \Ccinfty(\Omega,\Lambda)$, thanks to \eqref{eq:ProofBogovUnchangedOp} and \eqref{eq:ProofBogovCommutationRelationCompact}, it holds
\begin{align*}
    \delta \mathcal{B}^{0,\sigma} u + \mathcal{B}^{0,\sigma}\delta u + \K^{0,\sigma} u &= \delta  \T u + \L \delta \T u + \L \mathcal{E}_0 \underline{\delta}\mathcal{R}^{\sigma}  \R_{\Omega}\L u \\ &\qquad +  \T\delta u + \L \T{\delta} u + \L \mathcal{E}_0 \mathcal{R}^{\sigma}\underline{\delta}  \R_{\Omega}\L u + \L \mathcal{E}_0 \P_{\H}\R_{\Omega}\L u\\
    &= [\I-\L] u + \L[\I-\L] u + \L \mathcal{E}_0 [\underline{\delta}\mathcal{R}^{\sigma}+\mathcal{R}^{\sigma}\underline{\delta} +\P_{\H}] \R_{\Omega}\L u\\
    &= u - \L^2 u + \L \mathcal{E}_0 \R_{\Omega}\L u\\
    &= u.
\end{align*}
where we did apply Proposition~\ref{prop:FullHodgedecompL2bdd} and $\mathcal{E}_0 \R_{\Omega}=\I$ on $\C^{\infty}_{0,0}(\Omega,\Lambda)$. We also have, thanks to \eqref{eq:ProofBogovCommutationRelationCompact},
\begin{align*}
    \delta \K^{0,\sigma} u = \L \mathcal{E}_0 \delta \P_{\H}\R_{\Omega}\L u =0,
\end{align*}
as well as
\begin{align*}
    \K^{0,\sigma} \delta u = \L \mathcal{E}_0 [\P_{\H}\underline{\delta} ]\R_{\Omega}\L u =0,
\end{align*}
 which holds due to $\P_\H \underline{\delta} =0$ on $\D_2(\underline{\delta},\Omega,\Lambda)$ by construction,  according to Propositions~\ref{prop:Ext0PartialTraceVanish} and \ref{prop:FullHodgedecompL2bdd}.
\end{proof}

To provide optimal density results in the case of vector fields $\Lambda^1\simeq\CC^n$, we need a different version that preserves smooth compactly supported functions.

\begin{theorem}\label{thm:PotentialOpBddLipDomCcinfty}Let $\Omega$ be a bounded Lipschitz domain. There exists linear operators $\mathcal{B}^{0,\sigma}_c$ such that for any $p,q\in[1,\infty]$, $s>-2+1/p$, for ${\X}^{s,p}\in \{\,{\H}^{s,p},\, {\B}^s_{p,q},\,{\BesSmo}^{s}_{p,\infty}\,\}$, assuming $1<p<\infty$ if $\X=\H$, one has 
\begin{enumerate}[label=($\roman*$)]
    \item The mapping property,
    \begin{align*}
        \mathcal{B}^{0,\sigma}_c\,:\,\X^{s,p}_0(\Omega,\CC)\longrightarrow \X^{s+1,p}_0(\Omega,\CC^n);
    \end{align*}
    \item The identities
    \begin{align*}
      \div \mathcal{B}^{0,\sigma}_c [\div \vv ] = \div  \vv, \text{ for all } \vv \in\X^{s+1,p}_0(\Omega,\CC^n),
    \end{align*}
    and
    \begin{align*}
      \div \mathcal{B}^{0,\sigma}_c [\ww ] = \ww, \text{ for all } \ww \in\X^{s+1,p}_0(\Omega,\CC^n)\text{ such that } \langle\ww,\mathbbm{1}_{\Omega}\rangle_{\Omega}=0;
    \end{align*}
    \item Additionally, $\mathcal{B}^{0,\sigma}_c(\Ccinfty(\Omega,\CC))\subset\Ccinfty(\Omega,\CC^n)$.
\end{enumerate}
\end{theorem}

\begin{proof} It suffices to apply the strategy from the proof of \cite[Lemmas~2.3~\&~2.4,~ Theorem~2.5]{GesseirtHeckHieber} (see also \cite[Proposition~2.5]{MitreaMonniaux2008} for a similar argument) to \cite[Corollaries~3.3~\&~3.4, Remark~3.5]{CostabelMcIntosh2010}. The latter allows to deal with the case of Besov spaces including $p=1,\infty$.
\end{proof}

\paragraph{On special Lipschitz domains.}
The following theorem is a revised version of\cite[Section~3,~Theorem~3.3, Corollary~3.4,~Section~4,~Proposition~4.1,~Corollary~4.2]{CostabelTaggartMcIntosh2013}, in a slightly different framework than the one introduced in the given reference, but remains true replacing $\S'(\mathbb{R}^n)/\mathcal{P}(\RR^n)$ by $\S'(\RR^n)$ then using an appropriate realization of homogeneous Sobolev and Besov spaces over $\S'_h(\RR^n)$. \textbf{Due to the lack of completeness for function spaces of high regularity, our statement is a bit more involved.} While  the next result is not as functional as its counterparts on bounded Lipschitz domains, it will be sufficient for our purpose.

\begin{theorem}\label{thm:PotentialOpSpeLipDom}Let $\Omega$ be a special Lipschitz domain and $p,q\in[1,\infty]$, $s\in\mathbb{R}$. Consider $\dot{\X}^{s,p}\in \{\,\dot{\H}^{s,p},\, \dot{\B}^s_{p,q},\,\dot{\B}^{s}_{\infty,q},\,\dot{\BesSmo}^{s}_{p,\infty}\,\}$, assuming $1<p<\infty$ if $\X=\H$. There exists linear operators $\T^\sigma$, $(\T^\sigma_{a,b})_{0<a<b}$ respectively defined on $\S_0(\mathbb{R}^n,\Lambda)$ and $\S'_h(\RR^n,\Lambda)$ such that,
\begin{enumerate}[label=($\roman*$)]
    \item $\T^\sigma_{a,b}(\S'_h(\mathbb{R}^n,\Lambda^{k}))\subset \S'_h(\mathbb{R}^n,\Lambda^{k+1})$ and $\T^\sigma(\S_0(\mathbb{R}^n,\Lambda^{k}))\subset \S_0(\mathbb{R}^n,\Lambda^{k+1})$, for all $0<a<b$;
    \item For all $u\in\S'_h(\RR^n,\Lambda)$, it holds in $\S'(\RR^n,\Lambda)$
    \begin{align*}
        \delta\T^{\sigma}_{a,b}u+\T^{\sigma}_{a,b} \delta u \xrightarrow[\substack{a\longrightarrow 0_+\\b\longrightarrow\infty}]{} u.
    \end{align*}
    Furthermore, if $\delta u =0$ in $\Omega$, then it holds in $\mathcal{D}'(\Omega,\Lambda)$,
    \begin{align*}
        \delta\T^{\sigma}_{a,b}u \xrightarrow[\substack{a\longrightarrow0_+\\b\longrightarrow\infty}]{} u.
    \end{align*}
    \item For all $0<a<b$, there are induced by restriction  well-defined uniformly bounded linear operators $\T^\sigma_{a,b}:\,\dot{\X}^{s,p}(\Omega,\Lambda) \longrightarrow \dot{\X}^{s+1,p}(\Omega,\Lambda)$, and $\delta\T^\sigma_{a,b},\T^\sigma_{a,b}\delta:\,\dot{\X}^{s,p}(\Omega,\Lambda) \longrightarrow \dot{\X}^{s,p}(\Omega,\Lambda)$ \textit{ i.e. } for all $u\in \dot{\X}^{s,p}(\Omega,\Lambda)$,
    \begin{align*}
        \left\lVert \T^\sigma_{a,b} u \right\rVert_{\dot{\X}^{s+1,p}(\Omega)} &\lesssim_{n,s,p,q,\Omega} \left\lVert u \right\rVert_{\dot{\X}^{s,p}(\Omega)},\\
        \left\lVert \delta \T^\sigma_{a,b} u \right\rVert_{\dot{\X}^{s,p}(\Omega)}+\left\lVert \T^\sigma_{a,b}\delta u \right\rVert_{\dot{\X}^{s,p}(\Omega)} &\lesssim_{n,s,p,q,\Omega} \left\lVert u \right\rVert_{\dot{\X}^{s,p}(\Omega)}.
    \end{align*}
    Furthermore, the convergence properties of point \textit{(ii)} still hold in the sense of the present involved $\X^{\bullet,p}$-norms above, whenever $q<\infty$. 
    \item If \hyperref[AssumptionCompletenessExponents]{$(\mathcal{C}_{s+1,p,q})$} , then $\T^\sigma\,:\,\dot{\X}^{s,p}(\Omega,\Lambda) \longrightarrow \dot{\X}^{s+1,p}(\Omega,\Lambda)$ is the limit operator of $(\T^\sigma_{a,b})_{0<a<b}$ in the strong sense. If instead \eqref{AssumptionCompletenessExponents}, one considers the strong limits $\delta\T^\sigma,\T^\sigma\delta:\,\dot{\X}^{s,p}(\Omega,\Lambda) \longrightarrow \dot{\X}^{s,p}(\Omega,\Lambda)$ of $(\delta\T^\sigma_{a,b},\T^\sigma_{a,b}\delta)_{0<a<b}$, and also on their inhomogeneous counterparts\footnote{On inhomogeneous function spaces one can only consider $\delta\T^\sigma,\T^\sigma\delta$, and not $\T^\sigma$ itself due to homogeneity/scaling considerations.} $\X^{s,p}(\Omega,\Lambda)$ whenever $1<p<\infty$, $1\leqslant q\leqslant \infty$, $s\in\RR$ without any further assumption.
    \item If $($\ref{AssumptionCompletenessExponents}$)$ is satisfied, we have the following decomposition
    \begin{align*}
        \dot{\X}^{s,p}(\Omega,\Lambda^k) &= \delta \T^{\sigma}\dot{\X}^{s,p}(\Omega,\Lambda^{k})\oplus \T^{\sigma}\delta\dot{\X}^{s,p}(\Omega,\Lambda^{k}).
    \end{align*}
    with bounded projections $\delta \T^{\sigma}$ and $\T^{\sigma}\delta $. This translates  into
    \begin{align*}
        \dot{\X}^{s,p}(\Omega,\Lambda^k) &= \overline{\delta \dot{\X}^{s+1,p}(\Omega,\Lambda^{k+1})}\oplus \T^{\sigma}\delta\dot{\X}^{s,p}(\Omega,\Lambda^{k}),
    \end{align*}
    where the closure is taken weakly-$\ast$ whenever $q=\infty$. If additionally $1<p<\infty$, a similar result holds for inhomogeneous function spaces $\X^{s,p}(\Omega,\Lambda)$, for all $s\in\RR$, $q\in[1,\infty]$.
\end{enumerate}
Furthermore,
\begin{itemize}
    \item The whole result still holds for $\overline{\Omega}$-supported differential forms in $\RR^n$, replacing $\dot{\X}^{s,p}$ by $\dot{\X}^{s,p}_0$, with corresponding operators $\T^{0,\sigma}$ and $(\T^{0,\sigma}_{a,b})_{0<a<b}$;
    \item If $u\in\S'_h(\RR^n,\Lambda)$, is such that $\supp u \subset \overline{\Omega^c}$, then $\supp (\T^{\sigma}_{a,b} u)\subset \overline{\Omega^c}$ for all $0<a<b$. In particular, when it is well-defined, this property is preserved for the corresponding limit operator $\T^{\sigma}$;
    \item If $u\in\S'_h(\RR^n,\Lambda)$ is such that $\supp u \subset \overline{\Omega}$, then $\supp (\T^{0,\sigma}_{a,b} u)\subset \overline{\Omega}$ for all $0<a<b$. In particular, when it is well-defined, this property is preserved for the corresponding limit operator $\T^{0,\sigma}$;
    \item For all $0<a<b$, one also has $\T^{\sigma}_{a,b} (\Ccinfty(\RR^n,\Lambda)),\T^{0,\sigma}_{a,b} (\Ccinfty(\RR^n,\Lambda))\subset \Ccinfty(\RR^n,\Lambda)$.
    \item Everything still holds for $(\d,\gamma,\Lambda^{k-1})$ replacing $(\delta,\sigma,\Lambda^{k+1})$ and moreover $\T^{\gamma}_{a,b}=(\T^{0,\sigma}_{a,b})^\ast$ and $\T^{0,\gamma}_{a,b}=(\T^{\sigma}_{a,b})^\ast$.
\end{itemize}
\end{theorem}

\begin{remark}\label{rem:convergencePotOpSpeLipNoncomplete} We provide two remarks:
\begin{itemize}
\item  For all $u\in\S'_h(\RR^n,\Lambda)$ such that  $\delta u=0$, one has
\begin{align*}
    \delta \T^{\sigma}_{a,b} u = \delta \T^{\sigma}_{a,b} u + \T^{\sigma}_{a,b} \delta  u = u - \theta_{a}\ast u +\theta_b\ast u
\end{align*}
where $(\theta_t)_{t>0}$ is a suitable approximation of the identity. Consequently, on $\dot{\X}^{m,p}_{\sigma} \in\{\,\dot{\W}^{m,p}_{\sigma},\,\dot{\C}^{m}_{ub,h,\sigma},\,\dot{\C}^{m}_{0,\sigma}\}$, $p\in[1,\infty)$, one obtains trivially the following uniform bound for all $u\in \dot{\X}^{m,p}_{\sigma}(\RR^n,\Lambda)$
\begin{align*}
    \lVert \delta\T^{\sigma}_{a,b} u \rVert_{\dot{\X}^{m,p}(\RR^n)} \leqslant 3 \lVert  u \rVert_{\dot{\X}^{m,p}(\RR^n)}.
\end{align*}
Furthermore, as stated in the Theorem, one has convergence of $(\delta\T_{a,b}^{\sigma}u)_{0<a<b}$ towards $u\in\dot{\X}^{m,p}_{\sigma}(\RR^n,\Lambda)$. \textbf{However, this tells nothing about the boundedness of} $\delta\T_{a,b}^{\sigma}$, $\delta\T^{\sigma}$ \textbf{and} $\T_{a,b}^{\sigma}\delta$, $\T^{\sigma}\delta$,\textbf{ on  } $\dot{\X}^{m,p}(\RR^n)$ \textbf{itself} (\textit{i.e.} when no divergence-free/co-closed condition is assumed).

    \item For $\dot{\X}^{s,p}\in \{\,\dot{\H}^{s,p},\, \dot{\B}^s_{p,q},\,\dot{\B}^{s}_{\infty,q},\,\dot{\BesSmo}^{s}_{p,\infty}\,\}$, one can also consider $\T^{\sigma}$ as bounded linear operator from $\dot{\X}^{s,p}\cap\dot{\X}^{\alpha,r}(\Omega,\Lambda)$ to $\dot{\X}^{s+1,p}\cap\dot{\X}^{\alpha+1,r}(\Omega,\Lambda)$, as long as $\dot{\X}^{s+1,p}(\RR^n)$ is a complete space. In this case, one still has the corresponding decoupled estimates arising from point \textit{(iii)}, with
\begin{align*}
    \lVert\T^{\sigma} u \rVert_{\dot{\X}^{\alpha+1,r}(\Omega)} \lesssim_{\alpha,r,\partial\Omega}  \lVert u \rVert_{\dot{\X}^{\alpha,r}(\Omega)} \text{, for all }u\in\dot{\X}^{s,p}\cap\dot{\X}^{\alpha,r}(\Omega,\Lambda).
\end{align*}
And similarly when $\dot{\X}^{s,p}(\RR^n)$ is complete, with 
\begin{align*}
    \lVert \delta\T^{\sigma} u \rVert_{\dot{\X}^{\alpha,r}(\Omega)} + \lVert \T^{\sigma}\delta u \rVert_{\dot{\X}^{\alpha,r}(\Omega)} \lesssim_{\alpha,r,\partial\Omega}  \lVert u \rVert_{\dot{\X}^{\alpha,r}(\Omega)} \text{, for all }u\in\dot{\X}^{s,p}\cap\dot{\X}^{\alpha,r}(\Omega,\Lambda).
\end{align*}
In these cases, one still has convergence of $(\delta\T_{a,b}^{\sigma}u)_{0<a<b}$ and $(\T_{a,b}^{\sigma}\delta u)_{0<a<b}$ towards $\delta\T^{\sigma}u$ and $\T^{\sigma}\delta u$ respectively, with respect to the involved norms.
\end{itemize}

The same considerations apply to $\T^{0,\sigma}$ and $\dot{\X}_0$ instead of $\T^{\sigma}$ and $\dot{\X}$.
\end{remark}

\begin{proof} We give the argument for the convergence exhibited in point \textit{(ii)} using \cite[Propositions~3.1~\&~3.3]{CostabelTaggartMcIntosh2013}, the rest can be deduced from the analysis performed in \cite[Sections~3,~4~\&~7]{CostabelTaggartMcIntosh2013}. For all $u\in\S'_h(\RR^n,\Lambda)$,  one writes for all $0<a<b$, $\Gamma^{a,b}= \delta\T^{\sigma}_{a,b}+\T^{\sigma}_{a,b} \delta $ and
\begin{align*}
    u -\delta\T^{\sigma}_{a,b}u+\T^{\sigma}_{a,b} \delta u  &= [\I-\Gamma^{a,b}]\left(\sum_{j<0} \dot{\Delta}_j u\right) + [\I-\Gamma^{a,b}]\left(\sum_{j\geqslant 0} \dot{\Delta}_j u\right). 
\end{align*}
Note that $0<a<b<\infty$ ensures that $\T^{\sigma}_{a,b}u$, $\delta\T^{\sigma}_{a,b}u,\T^{\sigma}_{a,b} \delta u,$ and $\Gamma^{a,b}u\in\S'_h(\RR^n,\Lambda)$. It suffices then to check the convergence for the block of low frequencies, considering
\begin{align*}
    [\I-\Gamma^{a,b}]\left(\sum_{j<0} \dot{\Delta}_j u\right). 
\end{align*}
According to notations introduced in \cite[Propositions~3.1~\&~3.3]{CostabelTaggartMcIntosh2013}, for $u_{-}:=\sum_{j<0} \dot{\Delta}_j u$, $0<a<b$, one can write
\begin{align*}
    [\I-\Gamma^{a,b}]u_- = u_- - \theta_{a}\ast u_- +\theta_b\ast u_{-}, 
\end{align*}
where $(\theta_t)_{t>0}$ is a suitable approximation of the identity. Since $ u\in\S'_h(\RR^n,\Lambda)$, it turns out that $u_{-}\in\C^{0}_{ub,h}(\RR^n,\Lambda)$, so that
\begin{align*}
    \lVert u_- - \theta_{a}\ast u_-\rVert_{\L^\infty(\RR^n)} \xrightarrow[a\rightarrow 0_+]{} 0.
\end{align*}
Therefore, it remains to show that
\begin{align}\label{eq:proofBogovSpeLipVanishingLowFreq}
    \theta_{b}\ast u_- \xrightarrow[b\rightarrow \infty]{} 0\text{, in } \S'(\RR^n,\Lambda).
\end{align}
Since $u\in\S'_h(\RR^n,\Lambda)$, we recall that for $u_{-}^{N}:=\sum\limits_{-N\leqslant j<0} \dot{\Delta}_j u$, $N\in\NN^\ast$, one has by definition
\begin{align*}
    \lVert u_- - u_-^{N}\rVert_{\L^\infty(\RR^n)} \xrightarrow[N\rightarrow \infty]{} 0.
\end{align*}
Consequently, let $\varepsilon>0$ and fix $N\in\NN^\ast$ such that $\lVert u_- - u_-^{N}\rVert_{\L^\infty(\RR^n)}<\varepsilon$. Hence, for $\varphi\in\S(\RR^n,\Lambda)$, we write
\begin{align*}
     \langle\theta_{b}\ast u_-,\varphi\rangle =  \langle\theta_{b}\ast[ u_- - u_{-}^N],\varphi\rangle + \langle\theta_{b}\ast u_{-}^N,\varphi\rangle.
\end{align*}
However, by the properties of the convolution for tempered distributions, one also has
\begin{align*}
    \langle\theta_{b}\ast u_{-}^N,\varphi\rangle = \langle u,\theta_{b}\ast\varphi_{-}^N\rangle \xrightarrow[b\rightarrow \infty]{} 0,
\end{align*}
due to the fact that $\varphi_{-}^N \in \S_0(\RR^n,\Lambda)$ by construction. Now, by Young's Inequality for the convolution, we obtain
\begin{align*}
     \limsup_{b\rightarrow\infty}\,|\langle\theta_{b}\ast u_-,\varphi\rangle| &\leqslant    \limsup_{b\rightarrow\infty} \,|\langle\theta_{b}\ast[ u_- - u_{-}^N],\varphi\rangle|\\
     &\leqslant \lVert \theta_1 \rVert_{\L^1(\RR^n)}\lVert u_- - u_{-}^N \rVert_{\L^\infty(\RR^n)}\lVert \varphi \rVert_{\L^1(\RR^n)}\\
     &< \varepsilon \lVert \theta_1 \rVert_{\L^1(\RR^n)}\lVert \varphi \rVert_{\L^1(\RR^n)},
\end{align*}
which holds true for all $\varepsilon>0$, yielding \eqref{eq:proofBogovSpeLipVanishingLowFreq}.
\end{proof}

\subsubsection{A divergence preserving extension operator.}\label{sec:ExtOpDivPreser}

In \cite[Chap.~VI,~Sec.~3]{Stein1970}, Stein introduced an extension operator $\Ext$ from a special Lipschitz domain to the whole space, such that it has homogeneous estimates of integer order on any $\L^p$-spaces:
\begin{align*}
    \lVert \nabla^k (\Ext u)\rVert_{\L^p(\RR^n)}\lesssim_{p,n,k,\partial\Omega}\lVert \nabla^k u\rVert_{\L^p(\Omega)},\quad \forall u\in\W^{m,p}(\RR^n),\, k\in\llb 0,m\rrb,\, p\in[1,\infty],\, m\in\NN.
\end{align*}
This extension operator did play an important role in the regularity theory of partial differential equations on bounded Lipschitz domains, see for instance \cite{JerisonKenig1995}. The homogeneous estimates shown for Stein's extension operator allowed the second author to establish a suitable theory of homogeneous Sobolev and Besov spaces on special Lipschitz domains in \cite{Gaudin2023Lip}.

A bit more than a decade ago, Hiptmair, Li and Zou, in  \cite[Thm.~3.5]{HiptmairLiZou2012}, did improve Stein's construction, giving an operator $\Ext_\gamma$ which satisfies the following identity on differential forms with smooth coefficients
\begin{align*}
    \d \, \Ext_\gamma = \Ext_\gamma \, \d.
\end{align*}
In particular, they deduced the existence of an extension operator that commutes with the $\div$ and the $\curl$ operators on smooth functions. As a direct consequence, such an extension operator preserves the divergence free and the curl free properties.

On $k$-forms, one can consider instead
\begin{align*}
    \ExtDiv_{|_{\Lambda^k}} := (-1)^{k(n-k)}\star {\Ext_\gamma}_{|_{\Lambda^{n-k}}} \star
\end{align*}
and then introduce
\begin{align*}
    \ExtDiv := \bigoplus_{k=0}^{n} \ExtDiv_{|_{\Lambda^k}}.
\end{align*}
By construction, $\ExtDiv$ is an extension operator, and we do have $\delta \,\ExtDiv = \ExtDiv \,\delta$. One can then state the next result, where, as for Theorems~\ref{thm:InterpHomSpacesRn}~\&~\ref{thm:DivergenceFreeSpacesSpeiclaLipDensity}, we omit the mention of the exterior algebra $\Lambda$ in order to alleviate notations.

\begin{theorem}[ {\cite[Chapter~VI,~Section~3,~Theorem~5']{Stein1970}, \cite[Thm.~3.5]{HiptmairLiZou2012}} ]\label{thm:SteinsExtensionOpDivFree} For $\Omega$ a special Lipschitz domain, there exists a well-defined linear operator
\begin{align*}
    \ExtDiv\,:\,\L^1(\Omega)+\L^\infty(\Omega) \longrightarrow \L^1(\RR^n)+\L^\infty(\RR^n)
\end{align*}
such that $\delta \,\ExtDiv = \ExtDiv \,\delta$ on differential forms with smooth coefficients. For all $u\in\L^1(\Omega)+\L^\infty(\Omega)$,
\begin{align*}
    [{\ExtDiv u}]_{|_{\Omega}} = u\text{,}
\end{align*}
and for all $p\in[1,\infty]$, $m\in\NN$, if $u\in\mathrm{W}^{m,p}(\Omega)$, we have $\ExtDiv u\in \mathrm{W}^{m,p}(\RR^n)$, with the estimates
\begin{equation}\label{eq:SteinExtOpEstInteger}
    \lVert \nabla^k (\ExtDiv u)\rVert_{\L^p(\RR^n)}\lesssim_{p,n,k,\partial\Omega}\lVert \nabla^k u\rVert_{\L^p(\Omega)},\quad k\in\llb 0,m\rrb.
\end{equation}
If additionally $\delta u$ in $\mathrm{W}^{m,p}(\Omega)$, then $\delta \ExtDiv u= \ExtDiv \delta u\in \mathrm{W}^{m,p}(\RR^n)$, so that $ \ExtDiv u \in\mathrm{W}^{m,p}_{\sigma}(\RR^n)$ whenever $u \in\mathrm{W}^{m,p}_{\sigma}(\Omega)$.

We also have $\ExtDiv[\C^{m}_0(\overline{\Omega})] \subset \C^{m}_0(\RR^n)$ and $\ExtDiv[\C^{m}_{ub}(\overline{\Omega})] \subset \C^{m}_{ub}(\RR^n)$, for all $m\in\NN$.
\end{theorem}

 Before continuing, we recall that the issue for an extension operator $\mathcal{E}$ to be well-defined on $\L^\infty$-based homogeneous function spaces like $\dot{\B}^{s}_{\infty,q}$, say $s>0$, while
 \begin{align*}
    \mathcal{E}(\dot{\B}^{s}_{\infty,q}(\Omega))\subset \mathcal{E}({\B}^{s}_{\infty,q}(\Omega))\subset{\B}^{s}_{\infty,q}(\RR^n),
 \end{align*}
 is mostly about the knowledge of
 \begin{align*}
     \mathcal{E}(\dot{\B}^{s}_{\infty,q}(\Omega))\subset\mathcal{S}'_h(\RR^n),
 \end{align*}
 in order to fit the definition of function spaces by restriction. This issue  does not happen on $\C^{0}_0$ and $\L^p$-based homogeneous function spaces $1\leqslant p<\infty$. See \cite[Section~3]{Gaudin2023Lip} for more details.

The divergence-free/co-closed preserving extension operator from \cite{HiptmairLiZou2012} and Stein's extension operator do enjoy the same structure. Therefore, by complex interpolation and real interpolation for inhomogeneous function spaces, and by reproducing the analysis done in \cite[Lem.~3.10~\textit{(iii)}, Prop.~3.12~\&~3.34]{Gaudin2023Lip} for the homogeneous function spaces, we are able to obtain the next proposition.

\begin{proposition}\label{prop:bddExtOpHomFreeDivSpeLip} Let $s\geqslant 0$, $p,q\in[1,\infty]$, $p<\infty$ with $m:=s$, if $s\in\NN$. Consider the divergence-preserving extension operator $\ExtDiv$ introduced in Theorem \ref{thm:SteinsExtensionOpDivFree}. Let $\Ext^{c}_{\sigma}$ be the corresponding extension operator from $\overline{\Omega}^c$ to $\RR^n$, and introduce the projection operator
\begin{align*}
    \cp_{0,\sigma} := \I - \Ext^{c}_{\sigma}[ \R_{\overline{\Omega}^c}\cdot],
\end{align*}
satisfying $\delta\, \cp_{0,\sigma}= \cp_{0,\sigma}\,\delta$ on differential forms.

For $\dot{\X}^{s,p}\in\{\dot{\W}^{m,1},\,\dot{\H}^{s,p},\,\dot{\B}^{s}_{p,q},\,\dot{\B}^{s,0}_{\infty,q},\,\dot{\C}^{m}_{0}\}$, assuming either
\begin{itemize}
    \item $s>0$ if $\X=\B$; or
    \item $p>1$ if $\X=\H$;
\end{itemize}
then the extension operator $\ExtDiv$  maps boundedly  $\dot{\X}^{s,p}_{\sigma}(\Omega)$ to $\dot{\X}^{s,p}_{\sigma}(\RR^n)$, as well as the projection operator $\cp_{0,\sigma}$ maps boundedly $\dot{\X}^{s,p}_{\sigma}(\RR^n)$ to $\dot{\X}^{s,p}_{0,\sigma}(\Omega)$.

Moreover, one has $\ExtDiv \dot{\BesSmo}^{s,\sigma}_{p,\infty}(\Omega) \subset \dot{\BesSmo}^{s,\sigma}_{p,\infty}(\RR^n)$ (resp. $\cp_{0,\sigma} \dot{\BesSmo}^{s,\sigma}_{p,\infty}(\RR^n) \subset \dot{\BesSmo}^{s,\sigma}_{p,\infty,0}(\Omega)$) and similarly with $\dot{\BesSmo}^{s,0,\sigma}_{\infty,\infty}$. The result still holds for their inhomogeneous counterparts, including ${\B}^{s,\sigma}_{\infty,q}$, $\mathcal{B}^{s,\sigma}_{\infty,\infty}$ and ${\C}^{m}_{ub,\sigma}$.
\end{proposition}

\begin{remark}\label{rem:NotS'hBsinftyqSteinExtOp}To continue the discussion at the beginning of this subsection: In Proposition~\ref{prop:bddExtOpHomFreeDivSpeLip}, concerning the fonction spaces $\dot{\B}^{s}_{\infty,q}$, $s>0$, $q\in[1,\infty]$, the estimate
\begin{align*}
    \lVert \mathcal{E}_\sigma u\rVert_{\dot{\B}^{s}_{\infty,q}(\RR^n)}\lesssim_{p,n,\partial\Omega}\lVert u\rVert_{\dot{\B}^{s}_{\infty,q}(\Omega)}
\end{align*}
is valid, but it is not clear that $\mathcal{E}_\sigma u\in\S'_h(\RR^n,\Lambda)$.
\end{remark}

\paragraph{Concerning the Half-space.}

\begin{remark}\label{rem:FlatExtOpDiv}  When $\Omega=\RR^n_+$, following \cite[Lemma~3.15]{DanchinHieberMuchaTolk2020}, one can reach a similar result for the spaces $\dot{\B}^{s}_{\infty,q}$, $s>-1$, $q\in[1,\infty]$. For all $p,q\in[1,\infty]$, $s>-1+\sfrac{1}{p}$, for $m\in\NN$ fixed, such that $s<m+\sfrac{1}{p}$, there are extension and projection operators 
\begin{align*}
    \E_{\sigma}^m\,&:\,\dot{\B}^{s}_{p,q}({\RR^n_+},\CC^n) \longrightarrow\dot{\B}^{s}_{p,q}(\RR^n,\CC^n),\\
    \P_{0,\sigma}^m\,&:\,\dot{\B}^{s}_{p,q}({\RR^n},\CC^n) \longrightarrow\dot{\B}^{s}_{p,q,0}(\RR^n_+,\CC^n).
\end{align*}
and in particular each preserves free-divergence vector fields as bounded linear operators
\begin{align*}
    \E_{\sigma}^m\,&:\,\dot{\B}^{s,\sigma}_{p,q}({\RR^n_+},\CC^n) \longrightarrow\dot{\B}^{s,\sigma}_{p,q}(\RR^n,\CC^n),\\
    \P_{0,\sigma}^m\,&:\,\dot{\B}^{s,\sigma}_{p,q}({\RR^n},\CC^n) \longrightarrow\dot{\B}^{s,\sigma}_{p,q,0}(\RR^n_+,\CC^n).
\end{align*}
Where as in \cite[Lemma~3.15,~eq.~(3.7)]{DanchinHieberMuchaTolk2020},
\begin{align*}
    \E_{\sigma}^m\,:\,\L^{1}_{\text{loc}}(\overline{\RR^n_+},\CC^n) \longrightarrow\L^{1}_{\text{loc}}(\RR^n,\CC^n),
\end{align*}
is defined for all $\ff\in\L^{1}_{\text{loc}}(\overline{\RR^n_+},\CC^n)$, $\E_{\sigma}^m\ff(x):=\ff(x)$  for almost every $x\in\overline{\RR^n_+}$, and for almost every $x=(x',x_n)\in\RR^n_-$,
\begin{align*}
    \E_{\sigma}^m\ff(x',x_n) := \bigg( \sum_{j=0}^m \alpha_j\ff'\Big(x',-\frac{x_n}{j+1}\Big),\, \sum_{j=0}^m \beta_j\ff_n\Big(x',-\frac{x_n}{j+1}\Big)\bigg).
\end{align*}
Here, $(\alpha_j)_{0\leqslant j\leqslant m},(\beta_j)_{0\leqslant j \leqslant m}\subset\CC$ are such that for all $\ell\in\llb 0,m-1\rrb$, 
\begin{align*}
    \sum_{j=0}^m \alpha_j\Big(-\frac{1}{j+1}\Big)^\ell =1,\quad \sum_{j=0}^m \beta_j\Big(-\frac{1}{j+1}\Big)^\ell =1,\quad\text{ and }\quad (j+1)\alpha_j=-\beta_j\text{ for all }j\in\llb0,m\rrb,
\end{align*}
which are given by inverting a Vandermonde matrix. The condition $(j+1)\alpha_j=-\beta_j$, $j\in\llb0,m\rrb$, ensures that for all $\ff\in\L^{1}_{\text{loc}}(\overline{\RR^n_+},\CC^n)$,
\begin{align*}
    \div \E^{m}_{\sigma} \ff =  \E^{m}[ \div \ff].
\end{align*}
Where, for all $g\in\L^{1}_{\text{loc}}(\overline{\RR^n_+},\CC)$, $\E^m g(x):=g(x)$  for almost every $x\in\overline{\RR^n_+}$, and for almost every $x=(x',x_n)\in\RR^n_-$,
\begin{align*}
    \E^m g(x',x_n) := \sum_{j=0}^m \alpha_j g\Big(x',-\frac{x_n}{j+1}\Big).
\end{align*}
With this particular choice of $(\alpha_j)_{0\leqslant j\leqslant m}\subset\CC$ yielding extension and projection operators that satisfy
\begin{align*}
    \E^m\,&:\,\dot{\B}^{s}_{p,q}({\RR^n_+},\CC) \longrightarrow\dot{\B}^{s}_{p,q}(\RR^n,\CC),\\
    \P_{0}^m\,&:\,\dot{\B}^{s}_{p,q}({\RR^n},\CC) \longrightarrow\dot{\B}^{s}_{p,q,0}(\RR^n_+,\CC).
\end{align*}
Similar considerations apply for $\dot{\W}^{\ell,p}$, $\dot{\C}^\ell_{ub,h}$, $\ell\in\llb0,m\rrb$, $p\in[1,\infty]$, $\dot{\H}^{s,p}$, $1<p<\infty$, $s>-1+1/p$. Note that one has to check that $\E^m_{\sigma} \dot{\B}^{s}_{p,q}(\RR^n_+,\CC^n) \subset \S'_h(\RR^n,\CC^n)$, $s\geqslant0$, $p,q\in[1,\infty]$, which can be proved as in the proof of Lemma~\ref{lem:ExtDirNeuRn+}. The details are left to the reader. We also mention that with additional --tedious but not hard-- work one can build a differential form counterpart of above extension and projection operators.
\end{remark}

We give an additional result due to the behavior on constants --indeed $\P_{0,\sigma}^{m}(\CC^n) =\{0\}$--  which will be of use later on. The proof is left to the reader.
\begin{lemma}\label{lem:FlatProjOpDivModConst} Provided $m\in\NN$, $\P_{0,\sigma}^{m}$ induces a bounded linear quotient and co-restricted map
\begin{align*}
    \P_{0,\sigma}^{m}&\,:\,\left.{\L^{\infty}(\RR^n,\CC^n)}\middle/{\CC^n}\right. \longrightarrow \{\,\vv\in\L^\infty(\RR^n,\CC^n)\,:\,\supp\vv \subset \overline{\RR^n_+}\,\}.
\end{align*}
Furthermore, it is weakly-$\ast$ continuous and preserves divergence--free vector fields, and, additionally, its (pre-)dual operator satisfies the boundedness property
\begin{align*}
    (\P_{0,\sigma}^{m})^\ast&\,:\,{\L^{1}(\RR^n,\CC^n)} \longrightarrow {\L^{1}_\mfree(\RR^n,\CC^n)}=\bigg\{\,\vv\in\L^1(\RR^n,\CC^n)\,:\,\int_{\RR^n}\vv=\mathbf{0} \,\bigg\}.
\end{align*}
\end{lemma}

\subsubsection{Density}

As for Theorem~\ref{thm:InterpHomSpacesRn}, to alleviate notations and improve readability, we omit mentioning the exterior algebra $\Lambda$ in the upcoming next results.

\begin{theorem}\label{thm:DivergenceFreeSpacesbddLipDensity}Let $p,q\in[1,\infty)$, $s\in\RR$, $m\in\NN$ and $\Omega$ be a bounded Lipschitz domain of $\RR^n$. For $\CC^n$-valued function spaces ($\Lambda^1\simeq \CC^n$), it holds that:
\begin{enumerate}
    \item the space $\C_{\sigma}^\infty(\overline{\Omega})$ is a strongly dense subspace  of the function spaces:  ${\B}^{s,\sigma}_{p,q}(\Omega)$, ${\BesSmo}^{s,\sigma}_{p,\infty}(\Omega)$, ${\B}^{s,\sigma}_{\infty,q}(\Omega)$, ${\BesSmo}^{s,\sigma}_{\infty,\infty}(\Omega)$, ${\H}^{s,p}_{\sigma}(\Omega)$ with $1<p<\infty$, ${\W}^{m,1}_{\sigma}(\Omega)$ and ${\C}^m_{\sigma}(\overline{\Omega})$.

    \item the space $\Ccinftydiv(\Omega)$ is a strongly dense subspace of the function spaces:  ${\B}^{s,\sigma}_{p,q,0}(\Omega)$, ${\BesSmo}^{s,\sigma}_{p,\infty,0}(\Omega)$, ${\B}^{s,\sigma}_{\infty,q,0}(\Omega)$, ${\BesSmo}^{s,\sigma}_{\infty,\infty,0}(\Omega)$, $\W^{m,1}_{0,\sigma}(\Omega)$, ${\C}^m_{0,0,\sigma}(\Omega)$ and in ${\H}^{s,p}_{0,\sigma}(\Omega)$, assuming $1<p<\infty$ for the latter.
\end{enumerate}
Furthermore, for differential forms of any degree ($\Lambda=\Lambda^0\oplus\Lambda^1\oplus\ldots\oplus\Lambda^{n-1}\oplus\Lambda^{n}$) instead of vector fields only ($\Lambda^1\simeq\CC^n$), the result remains valid with the exact same proof if we replace $\Ccinftydiv(\Omega,\Lambda^1)$ by $\C^{\infty}_{0,0,\sigma}(\Omega,\Lambda)$\footnote{So the price to pay is the loss of compact support strictly included $\Omega$.}.
\end{theorem}

\begin{proof}[ of Theorem~\ref{thm:DivergenceFreeSpacesbddLipDensity} ] We deal with the case of $\overline{\Omega}$-supported divergence free vector fields, the other case admits a similar and simpler proof. In order to prove the result for bounded domains, especially in endpoint function spaces, we are going to make use of Lemma~\ref{lem:IntermediatedensityResultInhomogeneous}. 

\textbf{Step 1:} The case of ${\C}^m_{0,0,\sigma}(\Omega)$, $m\in\NN$. Let $\uu\in {\C}^m_{0,0,\sigma}(\Omega)$, by Lemma~\ref{lem:IntermediatedensityResultInhomogeneous} there exists $(\uu_k)_{k\in\NN}\subset\Ccinfty(\Omega,\CC^n)$, such that
\begin{align*}
    \lVert \uu- \uu_k \rVert_{\W^{m,\infty}(\RR^n)} + \lVert \div \uu_k \rVert_{\W^{m,\infty}(\RR^n)} \xrightarrow[k\longrightarrow\infty]{} 0
\end{align*}
Therefore, by Theorem~\ref{thm:PotentialOpBddLipDomCcinfty}, one has $\uu_k^\sigma:= \uu_k - \mathcal{B}^{0,\sigma}_c(\div \uu_k)\in\Ccinftydiv(\Omega)$. By Sobolev embeddings for $p>n$, it holds that
\begin{align*}
    \lVert \uu- \uu_k^\sigma \rVert_{\W^{m,\infty}(\RR^n)} &\leqslant \lVert \uu- \uu_k \rVert_{\W^{m,\infty}(\RR^n)} + \lVert \mathcal{B}^{0,\sigma}_c \div \uu_k \rVert_{\W^{m,\infty}(\RR^n)}\\
    &\lesssim_{p,m,n} \lVert \uu- \uu_k \rVert_{\W^{m,\infty}(\RR^n)} + \lVert \mathcal{B}^{0,\sigma}_c \div \uu_k \rVert_{\W^{m+1,p}(\RR^n)}\\
    &\lesssim_{p,m,n,\Omega} \lVert \uu- \uu_k \rVert_{\W^{m,\infty}(\RR^n)} + \lVert \div \uu_k \rVert_{\W^{m,p}(\RR^n)}\\
    &\lesssim_{p,m,n,\Omega} \lVert \uu- \uu_k \rVert_{\W^{m,\infty}(\RR^n)} + \lVert \div \uu_k \rVert_{\W^{m,\infty}(\RR^n)} \xrightarrow[k\longrightarrow\infty]{} 0.
\end{align*}
where the last estimate is obtained by the finiteness of the measure of $\Omega$. Therefore the result holds.

\textbf{Step 2:} The case of  Besov and Bessel potential spaces. In order to fix the ideas assume $\uu\in\B^{s,\sigma}_{p,q,0}(\Omega)$, the proof is similar for the other function spaces. We proceed as before, but since $s\in\RR$, we are specifically taking advantage of Theorem~\ref{thm:PotentialOpBddLipDom}, one only has $\uu_k^\sigma:= \uu_k - \mathcal{B}^{0,\sigma}(\div \uu_k) \in\C^\infty_{0,0,\sigma}(\Omega)$ with
\begin{align*}
    \lVert \uu- \uu_k^\sigma \rVert_{\B^{s}_{p,q}(\RR^n)}  &\leqslant \lVert \uu- \uu_k \rVert_{\B^{s}_{p,q}(\RR^n)} + \lVert \mathcal{B}^{0,\sigma} \div \uu_k \rVert_{\B^{s}_{p,q}(\RR^n)}\\
    &\lesssim_{p,s,n,\Omega} \lVert \uu- \uu_k \rVert_{\B^{s}_{p,q}(\RR^n)} + \lVert \div \uu_k \rVert_{\B^{s}_{p,q}(\RR^n)}\xrightarrow[k\longrightarrow\infty]{} 0.
\end{align*}
Now, for $m:=\max(0,\lceil s\rceil)+1$, one has $\uu_k^\sigma\in \C^\infty_{0,0,\sigma}(\Omega)\subset \C^m_{0,0,\sigma}(\Omega)$, so that by Step 1, there exists $(\uu_{k,\ell}^\sigma)_{\ell\in\NN}\subset \Ccinftydiv(\Omega)$ converging towards $\uu_k^\sigma$ in $\C^m_{0}(\RR^n,\CC^n)$, so that by the finiteness of the measure of $\Omega$,
\begin{align*}
     \lVert \uu- \uu_{k,\ell}^\sigma \rVert_{\B^{s}_{p,q}(\RR^n)} &\leqslant  \lVert \uu- \uu_k^\sigma \rVert_{\B^{s}_{p,q}(\RR^n)} +  \lVert \uu_{k}^\sigma- \uu_{k,\ell}^\sigma \rVert_{\B^{s}_{p,q}(\RR^n)} \\
      &\lesssim_{p,s,n}  \lVert \uu- \uu_k^\sigma \rVert_{\B^{s}_{p,q}(\RR^n)} +  \lVert \uu_{k}^\sigma- \uu_{k,\ell}^\sigma \rVert_{\W^{m,p}(\RR^n)} \\
      &\lesssim_{p,s,n,\Omega}  \lVert \uu- \uu_k^\sigma \rVert_{\B^{s}_{p,q}(\RR^n)} +  \lVert \uu_{k}^\sigma- \uu_{k,\ell}^\sigma \rVert_{\W^{m,\infty}(\RR^n)},
\end{align*}
choosing successively $k$ large enough, then $\ell$ large enough, yields the result.

\textbf{Step 3:} The case of the Sobolev space $\W^{m,1}_{0,\sigma}(\Omega)$. It starts similarly, but we take advantage of interpolation inequalities to deduce the following
\begin{align*}
    \lVert \uu- \uu_k^\sigma \rVert_{\W^{m,1}(\RR^n)}  &\leqslant \lVert \uu- \uu_k \rVert_{\W^{m,1}(\RR^n)} + \lVert \mathcal{B}^{0,\sigma} \div \uu_k \rVert_{\W^{m,1}(\RR^n)}\\
    &\lesssim_{p,s,n} \lVert \uu- \uu_k \rVert_{\W^{m,1}(\RR^n)} + \lVert \mathcal{B}^{0,\sigma} \div \uu_k \rVert_{\B^{m-1}_{1,\infty}(\RR^n)}^\frac{1}{2}\lVert \mathcal{B}^{0,\sigma} \div \uu_k \rVert_{\B^{m+1}_{1,\infty}(\RR^n)}^\frac{1}{2}\\
     &\lesssim_{p,s,n,\Omega} \lVert \uu- \uu_k \rVert_{\W^{m,1}(\RR^n)} + \lVert \mathcal{B}^{0,\sigma} \div \uu_k \rVert_{\B^{m}_{1,\infty}(\RR^n)}^\frac{1}{2} \lVert \div \uu_k \rVert_{\B^{m}_{1,\infty}(\RR^n)}^\frac{1}{2}\\
     &\lesssim_{p,s,n,\Omega} \lVert \uu- \uu_k \rVert_{\W^{m,1}(\RR^n)} + \lVert \uu_k \rVert_{\B^{m}_{1,\infty}(\RR^n)}^\frac{1}{2}\lVert \div \uu_k \rVert_{\B^{m}_{1,\infty}(\RR^n)}^\frac{1}{2}\\
     &\lesssim_{p,s,n,\Omega} \lVert \uu- \uu_k \rVert_{\W^{m,1}(\RR^n)} + \lVert \uu_k \rVert_{\W^{m,1}(\RR^n)}^\frac{1}{2}\lVert \div \uu_k \rVert_{\W^{m,1}(\RR^n)}^\frac{1}{2}\xrightarrow[k\longrightarrow\infty]{} 0.
\end{align*}
Finally, one can proceed as in Step 2 approximating $\uu_k^\sigma\in \C^\infty_{0,0,\sigma}(\Omega)\subset \C^m_{0,0,\sigma}(\Omega)$ by $(\uu_{k,\ell}^\sigma)_{\ell\in\NN}\subset \Ccinftydiv(\Omega)$ in $\C^m_{0,0,\sigma}(\Omega)$.
\end{proof}

\begin{theorem}\label{thm:DivergenceFreeSpacesSpeiclaLipDensity}Let $p,q\in[1,\infty)$, $s\in\RR$, $m\in\NN$ and $\Omega$ be a special Lipschitz domain of $\RR^n$.
It holds that:
\begin{enumerate}
    \item the space $\Ccinftydiv(\overline{\Omega})$ is a dense subspace  
    \begin{enumerate}
        \item of the homogeneous Besov spaces: $\dot{\B}^{s,\sigma}_{p,q}(\Omega)$, $s>-n/p'$,  $\dot{\B}^{s,0,\sigma}_{\infty,q}(\Omega)$, $s>-n$, $\dot{\BesSmo}^{s,0,\sigma}_{\infty,\infty}(\Omega)$, $s\geqslant-n$, $\dot{\BesSmo}^{s,\sigma}_{p,\infty}(\Omega)$, $s\geqslant-n/p'$;
        \item of the homogeneous Sobolev spaces: $\dot{\H}^{s,p}_{\sigma}(\Omega)$, $s>-n/p'$, $1<p<\infty$, $\dot{\W}^{m,1}_{\sigma}(\Omega)$ and $\dot{\C}^m_{0,\sigma}(\overline{\Omega})$;
        \item for all $s\in\RR$, $\Ccinftydiv(\overline{\Omega})\cap\dot{\B}^{s}_{p,q}(\Omega)$ is strongly dense in $\dot{\B}^{s,\sigma}_{p,q}(\Omega)$.

        The result remains true if we replace $\dot{\B}^{s}_{p,q}$, by either $\dot{\B}^{s,0}_{\infty,q}$, $\dot{\BesSmo}^{s}_{p,\infty}$, $\dot{\BesSmo}^{s,0}_{\infty,\infty}$ or $\dot{\H}^{s,p}$.
    \end{enumerate}
    
    \item the space $\Ccinftydiv(\Omega)$ is a dense subspace 
    \begin{enumerate}
        \item of the homogeneous Besov spaces: $\dot{\B}^{s,\sigma}_{p,q,0}(\Omega)$, $s>-n/p'$,  $\dot{\B}^{s,0,\sigma}_{\infty,q,0}(\Omega)$, $s>-n$, $\dot{\BesSmo}^{s,0,\sigma}_{\infty,\infty,0}(\Omega)$, $s\geqslant-n$, $\dot{\BesSmo}^{s,\sigma}_{p,\infty,0}(\Omega)$, $s\geqslant-n/p'$;
        \item  of the homogeneous Sobolev spaces: $\dot{\H}^{s,p}_{0,\sigma}(\Omega)$, $s>-n/p'$, $1<p<\infty$, $\dot{\W}^{m,1}_{0,\sigma}(\Omega)$ and $\dot{\C}^m_{0,0,\sigma}(\Omega)$;

        \item for all $s\in\RR$, $\Ccinftydiv({\Omega})\cap\dot{\B}^{s}_{p,q}(\RR^n)$ is strongly dense in $\dot{\B}^{s,\sigma}_{p,q,0}(\Omega)$.

        The result remains true if we replace $\dot{\B}^{s}_{p,q}$, by either $\dot{\B}^{s,0}_{\infty,q}$, $\dot{\BesSmo}^{s}_{p,\infty}$, $\dot{\BesSmo}^{s,0}_{\infty,\infty}$ or $\dot{\H}^{s,p}$.
    \end{enumerate}
    \item Additionally,
    \begin{enumerate}
        \item $\C^\infty_{ub,h,0,\sigma}({\Omega})$ is a strongly dense subspace of $\C^m_{ub,h,0,\sigma}({\Omega})$;
        \item $\C^\infty_{ub,h,\sigma}(\overline{\Omega})$ is a strongly dense subspace of $\C^m_{ub,h,\sigma}(\overline{\Omega})$;
        \item $\C^\infty_{ub,0,\sigma}({\Omega})\cap\dot{\B}^{s}_{\infty,q}(\RR^n)$ is a strongly dense subspace of $\dot{\B}^{s,\sigma}_{\infty,q,0}({\Omega})$;
        \item $\C^\infty_{ub,\sigma}(\overline{\Omega})\cap\dot{\B}^{s}_{\infty,q}(\Omega)$ is a strongly dense subspace of $\dot{\B}^{s,\sigma}_{\infty,q}(\Omega)$;
    \end{enumerate}
    and similarly with $\dot{\BesSmo}^{s}_{\infty,\infty}$ replacing $\dot{\B}^{s}_{\infty,q}$.
\end{enumerate}
Furthermore, similar results hold
\begin{itemize}
    \item  for the corresponding  inhomogeneous function spaces;
    \item  with $\gamma$ replacing $\sigma$.
\end{itemize}
\end{theorem}

\begin{proof} \textbf{Step 1:} We prove Point \textit{(i)}.

 \textbf{Step 1.1:}  Point \textit{(i)} for the case of homogeneous Besov and Bessel potential spaces with $s>0$, $\dot{\W}^{m,1}_{\sigma}$, and $\C^{m}_{0,\sigma}(\overline{\Omega})$, $m\in\NN$. This is a direct consequence from Theorem~\ref{prop:bddExtOpHomFreeDivSpeLip}. We focus on the case $\dot{\W}^{m,1}_{\sigma}(\Omega)$, the remaining cases admit a similar proof.
For say $u\in\dot{\W}^{m,1}_{\sigma}(\Omega)$, then it holds that $\mathcal{E}_{\sigma}u\in\dot{\W}^{m,1}_{\sigma}(\RR^n)$. Consequently, it suffices to approximate $\mathcal{E}_{\sigma}u$ by $(U_\ell)_{\ell\in\NN}\subset\Ccinftydiv(\RR^n)$, so that $(u_\ell)_{\ell\in\NN} =({U_\ell}_{|_{\Omega}})_{\ell\in\NN} \subset\Ccinftydiv(\overline{\Omega})$ converges towards $u$.

The same goes for other function spaces.

\textbf{Step 1.2:}  Point \textit{(i)} with then case of homogeneous Besov and Sobolev spaces with $s\leqslant 0$.

Let $u\in\dot{\B}^{s,\sigma}_{p,q}(\Omega)$, by Theorem~\ref{thm:PotentialOpSpeLipDom}, since $\delta u =0$ one has
\begin{align*}
    \lVert u - \delta \T^{\sigma}_{a,b}u\rVert_{\dot{\B}^{s}_{p,q}(\Omega)}\xrightarrow[\substack{a\rightarrow 0_+\\b\rightarrow\infty}]{} 0.
\end{align*}
Let $(u_\ell)_{\ell\in\NN}\subset\Ccinfty(\overline{\Omega})\cap\dot{\B}^{s}_{p,q}(\Omega)$ converging to $u$ in $\dot{\B}^{s}_{p,q}(\Omega)$. For $0<a<b$, $\ell\in\NN$, one obtains by Theorem~\ref{thm:PotentialOpSpeLipDom}
\begin{align*}
    \lVert u - \delta \T^{\sigma}_{a,b}u_\ell\rVert_{\dot{\B}^{s}_{p,q}(\Omega)} &\leqslant \lVert u - \delta \T^{\sigma}_{a,b}u\rVert_{\dot{\B}^{s}_{p,q}(\Omega)} + \lVert \delta \T^{\sigma}_{a,b} u - \delta \T^{\sigma}_{a,b}u_\ell\rVert_{\dot{\B}^{s}_{p,q}(\Omega)} \\
    &\lesssim_{p,s,n,\partial\Omega} \lVert u - \delta \T^{\sigma}_{a,b}u\rVert_{\dot{\B}^{s}_{p,q}(\Omega)} + \lVert u-u_\ell\rVert_{\dot{\B}^{s}_{p,q}(\Omega)}.
\end{align*}
Let $\varepsilon>0$, and consider $b>a>0$ such that $\lVert u - \delta \T^{\sigma}_{a,b}u\rVert_{\dot{\B}^{s}_{p,q}(\Omega)}<\varepsilon$, and consider $\ell\in\NN$ such that $\lVert u-u_\ell\rVert_{\dot{\B}^{s}_{p,q}(\Omega)}<\varepsilon$, it holds that 
\begin{align*}
    \lVert u - \delta \T^{\sigma}_{a,b}u_\ell\rVert_{\dot{\B}^{s}_{p,q}(\Omega)} \lesssim_{p,s,n,\partial\Omega} \varepsilon
\end{align*}
which yields the result since $\delta \T^{\sigma}_{a,b}u_\ell\in\Ccinftydiv(\overline{\Omega})\cap\dot{\B}^{s}_{p,q}(\Omega)$ by Theorem~\ref{thm:PotentialOpSpeLipDom}. Note that it also yields \textit{(i)-(a)} and \textit{(i)-(b)} due to the fact that in such cases one has the inclusion $\Ccinfty(\overline{\Omega})\subset\dot{\B}^{s}_{p,q}(\Omega)$.

\textbf{Step 2:} We prove point \textit{(ii)}.

\textbf{Step 2.1:} We deal with \textit{(ii)-(c)} with the case of homogeneous Besov and Bessel potential spaces assuming \eqref{AssumptionCompletenessExponents}. In order to fix the ideas we deal with the case $\dot{\B}^{s}_{p,q}$, the other cases admit a similar proof. As in Step 1.2,  for $u\in \dot{\B}^{s,\sigma}_{p,q,0}(\Omega)\subset\dot{\B}^{s,\sigma}_{p,q}(\RR^n)$,
\begin{align*}
    \lVert u - \delta \T^{0,\sigma}_{a,b}u\rVert_{\dot{\B}^{s}_{p,q}(\RR^n)}\xrightarrow[\substack{a\rightarrow 0_+\\b\rightarrow\infty}]{} 0.
\end{align*}
Let $(u_\ell)_{\ell\in\NN}\subset\Ccinfty({\Omega})\cap\dot{\B}^{s}_{p,q}(\RR^n)$ converging toward $u$ in $\dot{\B}^{s}_{p,q,0}(\Omega)\subset \dot{\B}^{s}_{p,q}(\RR^n)$. Let $\varepsilon>0$, and consider $b>a>0$ such that $\lVert u - \delta \T^{0,\sigma}_{a,b}u\rVert_{\dot{\B}^{s}_{p,q}(\RR^n)}<\varepsilon$, and consider $\ell\in\NN$ such that $\lVert u-u_\ell\rVert_{\dot{\B}^{s}_{p,q}(\RR^n)}<\varepsilon$, it holds that 
\begin{align*}
    \lVert u - \delta \T^{0,\sigma}_{a,b}u_\ell\rVert_{\dot{\B}^{s}_{p,q}(\RR^n)} \lesssim_{p,s,n,\partial\Omega} \varepsilon.
\end{align*}
The issue here concerns the compact support to be included in $\Omega$ for the elements that constitute  the approximating sequence $( \delta \T^{0,\sigma}_{a,b}u_\ell)_{a,b,\ell}$.

So, while $v_\ell^{a,b}:=\delta \T^{\sigma}_{a,b}u_\ell\in \Ccinftydiv\cap\dot{\B}^{s}_{p,q}(\RR^n)$ with compact support by Theorem~\ref{thm:PotentialOpSpeLipDom}, it holds its support is only included in $\overline{\Omega}$. Since $\Omega$  is a special Lipschitz domain, by strong continuity of translations for in the considered function spaces, one can introduce for all $t>0$, all $x\in\RR^n$, $v_{\ell,t}^{a,b}(x):=v_\ell^{a,b}(x',x_n-t)$ so that by construction $v_{\ell,t}^{a,b}\in\Ccinftydiv(\Omega)\cap\dot{\B}^{s}_{p,q}(\RR^n)$, with, for $t>0$ small enough,
\begin{align*}
    \lVert u - v_{\ell,t}^{a,b}\rVert_{\dot{\B}^{s}_{p,q}(\RR^n)}\lesssim_{p,s,n,\partial\Omega}\varepsilon.
\end{align*}
\textbf{Step 2.2:} We deal with homogeneous Besov and Bessel potential spaces in the case of Points \textit{(ii)-(a)} and \textit{(ii)-(b)} assuming $s> 0$ and $p<\infty$. Note the overlap with the case \eqref{AssumptionCompletenessExponents} already treated in the previous Step 2.1.

In order to fix the ideas, we deall with the case of the homogeneous Besov space $\dot{\B}^{s}_{p,q}$, $q<\infty$, the other cases admit a similar proof. 

Let $u\in\dot{\B}^{s,\sigma}_{p,q,0}(\Omega)$, by Theorem~\ref{prop:bddExtOpHomFreeDivSpeLip}, one has $\cp_{0,\sigma} u = u$, so that for $u_N:=  \cp_{0,\sigma} \big(\sum_{|j|\leqslant N} \dot{\Delta}_j u\big)\in\C^{\infty}_{0,0,\sigma}(\Omega)\cap \B^{s+m}_{p,q}(\RR^n)$ for all $m\in\NN$ and one obtains
\begin{align*}
    \lVert u-u_N \rVert_{\dot{\B}^{s}_{p,q}(\RR^n)} \xrightarrow[N\rightarrow\infty]{} 0.
\end{align*}
In particular, one has $\C^{\infty}_{0,0,\sigma}(\Omega)\cap \dot{\B}^{s_0}_{p,1}(\RR^n)\cap \dot{\B}^{s}_{p,q}(\RR^n)$, where $0<s_0\leqslant \sfrac{n}{p}$ and $s_0<s$, so that by Theorem~\ref{thm:PotentialOpSpeLipDom} and its following remark, Remark~\ref{rem:convergencePotOpSpeLipNoncomplete},
\begin{align*}
    \lVert u_N - \delta \T^{0,\sigma}_{a,b} u_N \rVert_{\dot{\B}^{s}_{p,q}(\RR^n)} \xrightarrow[\substack{a\rightarrow 0_+\\b\rightarrow\infty}]{} 0.
\end{align*}
For $\varepsilon>0$, we consider $b$ and $N$ large enough, and $a$ small enough such that
\begin{align*}
    \lVert u-\delta \T^{0,\sigma}_{a,b} u_N\rVert_{\dot{\B}^{s}_{p,q}(\RR^n)} < \varepsilon .
\end{align*}
Now, consider $(u_{N,\ell})_{\ell\in\NN}\subset \Ccinfty(\Omega)$, such that it converges towards $u_N$ in $\dot{\B}^{s}_{p,q,0}(\Omega)\hookrightarrow\dot{\B}^{s}_{p,q}(\RR^n)$. By Theorem~\ref{thm:PotentialOpSpeLipDom} one has
\begin{align*}
    \lVert u-\delta \T^{0,\sigma}_{a,b} u_{N,\ell}\rVert_{\dot{\B}^{s}_{p,q}(\RR^n)} &\leqslant \lVert u-\delta \T^{0,\sigma}_{a,b} u_{N}\rVert_{\dot{\B}^{s}_{p,q}(\RR^n)} + \lVert \delta \T^{0,\sigma}_{a,b} u_{N}-\delta \T^{0,\sigma}_{a,b} u_{N,\ell}\rVert_{\dot{\B}^{s}_{p,q}(\RR^n)}\\
    &\lesssim_{p,s,n,\partial\Omega} \varepsilon + \lVert  u_{N}- u_{N,\ell}\rVert_{\dot{\B}^{s}_{p,q}(\RR^n)},
\end{align*}
so that for $\ell$ large enough
\begin{align*}
    \lVert u-\delta \T^{0,\sigma}_{a,b} u_{N,\ell}\rVert_{\dot{\B}^{s}_{p,q}(\RR^n)} \lesssim_{p,s,n,\partial\Omega} \varepsilon.
\end{align*}
As in the previous Step 2.1, one can conclude by strong continuity of translations since, by Theorem~\ref{thm:PotentialOpSpeLipDom}, one obtains $\delta \T^{0,\sigma}_{a,b} u_{N,\ell}\in\Ccinftydiv\cap\dot{\B}^{s}_{p,q}(\RR^n)$ with compact support in $\overline{\Omega}$.

\textbf{Step 2.3:} We deal with the remaining function spaces: $\dot{\W}^{m,1}_{0,\sigma}(\Omega)$, $\dot{\C}^{m}_{0,0,\sigma}(\Omega)$, $\dot{\B}^{s,0,\sigma}_{\infty,q,0}(\Omega)$ and $\dot{\BesSmo}^{s,0,\sigma}_{\infty,\infty,0}(\Omega)$, $m\in\NN$, $s\geqslant 0$, $q\in[1,\infty)$.

We take inspiration from the previous Step 2.2, we consider $0<a<b$ and $N\in\NN$. Let $(u_{\ell})_{\ell\in\NN}\subset \Ccinftydiv(\RR^n)$ converging to $u$. For $\varepsilon>0$, provided $a$ is small enough and $b$ and $\ell$ are large enough, by Theorem~\ref{prop:bddExtOpHomFreeDivSpeLip} and Remark~\ref{rem:convergencePotOpSpeLipNoncomplete}, it holds
\begin{align*}
    \lVert u - \delta\T^{0,\sigma}_{a,b}u_\ell\rVert_{\dot{\W}^{m,1}(\RR^n)}<\varepsilon.
\end{align*}
Note that $\T^{0,\sigma}_{a,b}u_\ell\in\Ccinfty(\RR^n)$, therefore for $\Theta\in\Ccinfty(\B_{1}(0))$, with $\Theta =1$ on $\B_{\frac{1}{2}}(0)$, for $R>0$, we look at
\begin{align*}
    \lVert u - \delta(\Theta(\cdot/R)\cp_{0,\sigma} [\T^{0,\sigma}_{a,b}u_\ell])\rVert_{\dot{\W}^{m,1}(\RR^n)} &\leqslant \lVert u - \Theta(\cdot/R)\cp_{0,\sigma}\delta [\T^{0,\sigma}_{a,b}u_\ell])\rVert_{\dot{\W}^{m,1}(\RR^n)} \\ &\qquad+  \frac{1}{R}\lVert [\nabla \Theta](\cdot/R)\iprod \cp_{0,\sigma}[\T^{0,\sigma}_{a,b}u_\ell])\rVert_{\dot{\W}^{m,1}(\RR^n)}\\
    &\leqslant \lVert u - \cp_{0,\sigma}\delta [\T^{0,\sigma}_{a,b}u_\ell]\rVert_{\dot{\W}^{m,1}(\RR^n)}\\ &\qquad + \lVert \cp_{0,\sigma}\delta [\T^{0,\sigma}_{a,b}u_\ell]- \Theta(\cdot/R)\cp_{0,\sigma}\delta [\T^{0,\sigma}_{a,b}u_\ell] \rVert_{\dot{\W}^{m,1}(\RR^n)}\\ &\qquad +  \frac{1}{R}\lVert [\nabla \Theta](\cdot/R)\iprod \cp_{0,\sigma}[\T^{0,\sigma}_{a,b}u_\ell]\rVert_{\dot{\W}^{m,1}(\RR^n)}\\
    &\lesssim_{m,n,\partial\Omega} \varepsilon + \lVert \cp_{0,\sigma}\delta [\T^{0,\sigma}_{a,b}u_\ell]- \Theta(\cdot/R)\cp_{0,\sigma}\delta [\T^{0,\sigma}_{a,b}u_\ell] \rVert_{{\W}^{m,1}(\RR^n)}\\ &\qquad +  \frac{1}{R}\lVert [\nabla \Theta](\cdot/R)\iprod \cp_{0,\sigma}[\T^{0,\sigma}_{a,b}u_\ell]\rVert_{{\W}^{m,1}(\RR^n)}.
\end{align*}
Hence, for $R$ large enough, recalling the fact that $\T^{0,\sigma}_{a,b}u_\ell\in\Ccinfty(\RR^n)$, the dominated convergence theorem yields
\begin{align*}
    \lVert u - \delta(\Theta(\cdot/R)\cp_{0,\sigma} [\T^{0,\sigma}_{a,b}u_\ell])\rVert_{\dot{\W}^{m,1}(\RR^n)}\lesssim_{m,n,\partial\Omega} \varepsilon.
\end{align*}
It holds that $v_{\ell}^{a,b,R}:=\delta(\Theta(\cdot/R)\cp_{0,\sigma} [\T^{0,\sigma}_{a,b}u_\ell])$ belongs to $\Ccinftydiv\cap\dot{\W}^{m,1}(\RR^n)$, with compact support included in $\overline{\Omega}$. Again, as in the end of previous steps, strong continuity of translations yields the result. The proof for $\C^m_{0,0,\sigma}(\Omega)=\dot{\C}^m_{0,0,\sigma}(\Omega)$ is similar.

Now, let $u\in\dot{\B}^{s,0,\sigma}_{\infty,q,0}(\Omega)\hookrightarrow\dot{\B}^{s,0,\sigma}_{\infty,q}(\RR^n)$, note that for all $j\in\ZZ$, $\dot{\Delta}_j u \in\C^{\infty}_{0,\sigma}(\RR^n)$, and for all $b>a>0$, $\delta\T^{0,\sigma}_{a,b} u \in\S'_h(\RR^n)$ and $\dot{\Delta}_j\delta\T^{0,\sigma}_{a,b} u=\delta\T^{0,\sigma}_{a,b}\dot{\Delta}_j u \in\C^{\infty}_{0,\sigma}(\RR^n)$, and it holds
\begin{align*}
    \lVert u- \delta\T^{0,\sigma}_{a,b} u \rVert_{\dot{\B}^{s}_{\infty,q}(\RR^n)} &= \Big(\sum_{j\in\ZZ} 2^{jsq} \lVert \dot{\Delta}_ju- \delta\T^{0,\sigma}_{a,b} \dot{\Delta}_ju \rVert_{\L^\infty(\RR^n)}^q \Big)^{\sfrac{1}{q}} \\
    &\leqslant 2 \Big(\sum_{j\in\ZZ} 2^{jsq} \lVert \dot{\Delta}_ju \rVert_{\L^\infty(\RR^n)}^q \Big)^{\sfrac{1}{q}} = 2 \lVert u\rVert_{\dot{\B}^{s}_{\infty,q}(\RR^n)}.
\end{align*}
where we did apply Remark~\ref{rem:convergencePotOpSpeLipNoncomplete} to obtain the bounds
\begin{align*}
    \lVert \dot{\Delta}_ju- \delta\T^{0,\sigma}_{a,b} \dot{\Delta}_ju \rVert_{\L^\infty(\RR^n)} \leqslant 2 \lVert \dot{\Delta}_ju \rVert_{\L^\infty(\RR^n)},\\
     \lVert \delta\T^{0,\sigma}_{a,b} \dot{\Delta}_ju \rVert_{\L^\infty(\RR^n)} \leqslant 3 \lVert \dot{\Delta}_ju \rVert_{\L^\infty(\RR^n)}.
\end{align*}
Since $\dot{\Delta}_ju\in\C^{\infty}_{0,\sigma}(\RR^n)$,  Remark~\ref{rem:convergencePotOpSpeLipNoncomplete} yields
\begin{align*}
    \lVert \dot{\Delta}_ju- \delta\T^{0,\sigma}_{a,b} \dot{\Delta}_ju \rVert_{\L^\infty(\RR^n)} \xrightarrow[\substack{a\rightarrow 0_+\\b\rightarrow\infty}]{} 0.
\end{align*}
By the dominated convergence theorem, it holds that
\begin{align*}
    \lVert u- \delta\T^{0,\sigma}_{a,b} u \rVert_{\dot{\B}^{s}_{\infty,q}(\RR^n)} \xrightarrow[\substack{a\rightarrow 0_+\\b\rightarrow\infty}]{} 0.
\end{align*}
Now, assume $(u_{\ell})_{\ell\in\NN}\subset \Ccinfty(\Omega)$ converging to $u$ in $\dot{\B}^{s,0}_{\infty,q,0}(\Omega)\hookrightarrow\dot{\B}^{s}_{\infty,q}(\RR^n)$
\begin{align*}
    \lVert u- \delta\T^{0,\sigma}_{a,b} u_{\ell} \rVert_{\dot{\B}^{s}_{\infty,q}(\RR^n)}\leqslant \lVert u- \delta\T^{0,\sigma}_{a,b} u \rVert_{\dot{\B}^{s}_{\infty,q}(\RR^n)}+ 2\lVert  u - u_{\ell} \rVert_{\dot{\B}^{s}_{\infty,q}(\RR^n)}.
\end{align*}
So that the result follows ending this steps as the previous ones: taking $a$ small enough, $b$ and $\ell$ large enough, and concluding as in Step 2.1 by strong continuity of translations.

All the remaining cases, including co-closed differential forms that does not vanish on the boundary, and the cases of inhomogeneous function spaces, can be proved by similar arguments.
\end{proof}

\subsubsection{Interpolation}

This first interpolation result only concerns bounded Lipschitz domains.

\begin{theorem}\label{thm:InterpHomSpacesBddLip}Let $1\leqslant p_0,p_1,p,q,q_0,q_1\leqslant \infty$, $s_0,s_1\in\RR$, such that $s_0\neq s_1$, and for $\theta\in(0,1)$. Consider $\Omega$ to be a bounded Lipschitz domain, and let
\begin{align*}
    \left(s,\frac{1}{p_\theta},\frac{1}{q_\theta}\right):= (1-\theta)\left(s_0,\frac{1}{p_0},\frac{1}{q_0}\right)+ \theta\left(s_1,\frac{1}{p_1},\frac{1}{q_1}\right)\text{. }
\end{align*}
If $s_0,s_1\in\NN$, we write $m_0:=s_0$ and $m_1:=s_1$.
We have the following interpolation identities with equivalence of norms
\begin{enumerate}
    \item $({\B}^{s_0,\sigma}_{p,q_0}(\Omega),{\B}^{s_1,\sigma}_{p,q_1}(\Omega))_{\theta,q}=({\H}^{s_0,p}_{\sigma}(\Omega),{\H}^{s_1,p}_{\sigma}(\Omega))_{\theta,q}=({\W}^{m_0,p}_{\sigma}(\Omega),{\W}^{m_1,p}_{\sigma}(\Omega))_{\theta,q}={\B}^{s,\sigma}_{p,q}(\Omega)$\label{eq:realInterpHomBspqBddLip};
    \item $({\C}^{m_0}_{\sigma}(\overline{\Omega}),{\C}^{m_1}_{\sigma}(\overline{\Omega}))_{\theta,q}={\B}^{s,\sigma}_{\infty,q}(\Omega)$\label{eq:realInterpHomBsinftyqBddLip};
    \item $[{\H}^{s_0,p_0}_{\sigma}(\Omega),{\H}^{s_1,p_1}_{\sigma}(\Omega)]_{\theta} = {\H}^{s,p_\theta}_{\sigma}(\Omega)$, if $1<p_0,p_1<\infty$\label{eq:complexInterpHomSobBddLip};
    \item $[{\B}^{s_0,\sigma}_{p_0,q_0}(\Omega),{\B}^{s_1,\sigma}_{p_1,q_1}(\Omega)]_{\theta} = {\B}^{s,\sigma}_{p_\theta,q_\theta}(\Omega)$,  $q_\theta<\infty$.\label{eq:complexInterpHomBspqBddLip}
\end{enumerate}
and similarly with
\begin{align*}
    \{{\W}_0,\,{\H}_0,\,{\B}_{\cdot,\cdot,0},\,{\C}_{0,0}\}\quad \text{ instead of  }\quad\{{\W},\,{\H},\,{\B},\,{\C}\}.
\end{align*}

Moreover, 
\begin{itemize}
    \item The identity \ref{eq:realInterpHomBspqBddLip} still holds if we replace ${\B}^{s_j,\sigma}_{p,q_j}(\Omega)$ by ${\mathcal{B}}^{s_j,\sigma}_{p,\infty}(\Omega)$, $j\in\{0,1\}$;
    \item When we replace $(\cdot,\cdot)_{\theta,q}$ by $(\cdot,\cdot)_{\theta}$, all the real interpolation identities still hold with ${\mathcal{B}}^{s,\sigma}_{p,\infty}$ as an output space  ;
    \item The identity \ref{eq:complexInterpHomBspqBddLip} remains valid for $q_\theta=\infty$ with an output space ${\BesSmo}^{s_\theta,\sigma}_{p_{\theta},\infty}(\Omega)$.
\end{itemize}
\end{theorem}

\begin{proof} We  deal with \ref{eq:realInterpHomBspqBddLip} and \ref{eq:realInterpHomBsinftyqBddLip} in the case of $\overline{\Omega}$-supported differential forms. By Theorem~\ref{thm:PotentialOpBddLipDom}, we introduce the linear operator
\begin{align*}
    \P_{0,\sigma}:= \I-\mathcal{B}^{0,\sigma}\delta.
\end{align*}
For all $p,q\in[1,\infty]$, $s\in\RR$, it is a bounded linear projection
\begin{align*}
    \P_{0,\sigma}\,:\,\B^{s}_{p,q,0}(\Omega)\longrightarrow\B^{s,\sigma}_{p,q,0}(\Omega).
\end{align*}
Therefore, by retraction and co-retraction, it holds
\begin{align*}
    ({\B}^{s_0,\sigma}_{p,q_0,0}(\Omega),{\B}^{s_1,\sigma}_{p,q_1,0}(\Omega))_{\theta,q} = \P_{0,\sigma}({\B}^{s_0}_{p,q_0,0}(\Omega),{\B}^{s_1}_{p,q_1,0}(\Omega))_{\theta,q} = \P_{0,\sigma}\B^{s}_{p,q,0}(\Omega) = \B^{s,\sigma}_{p,q,0}(\Omega).
\end{align*}
Finally, the natural embeddings for $j\in\{0,1\}$,
\begin{align*}
    \B^{s_j,\sigma}_{p,1,0}(\Omega)\hookrightarrow \H^{s_j,p}_{0,\sigma}(\Omega)\hookrightarrow \B^{s_j,\sigma}_{p,\infty,0}(\Omega),
\end{align*}
and
\begin{align*}
    \B^{m_j,\sigma}_{p,1,0}(\Omega)&\hookrightarrow \W^{m_j,p}_{0,\sigma}(\Omega)\hookrightarrow \B^{m_j,\sigma}_{p,\infty,0}(\Omega),\\
    \B^{m_j,\sigma}_{\infty,1,0}(\Omega)&\hookrightarrow \C^{m_j}_{0,0,\sigma}(\Omega)\hookrightarrow \B^{m_j,\sigma}_{\infty,\infty,0}(\Omega),
\end{align*}
yield  \ref{eq:realInterpHomBspqBddLip} and \ref{eq:realInterpHomBsinftyqBddLip}. The case of complex interpolation holds by the same retraction and co-retraction argument.

The proof still applies verbatim for co-closed differential forms that does not vanish on the boundary.
\end{proof}

This second result only concerns special Lipschitz domains.

\begin{theorem}\label{thm:InterpHomSpacesLip}Let $1\leqslant p_0,p_1,p,q,q_0,q_1\leqslant \infty$, $s_0,s_1\in\RR$, such that $s_0\neq s_1$, and for $\theta\in(0,1)$. Let $\Omega$  be a special Lipschitz domain, and
\begin{align*}
    \left(s,\frac{1}{p_\theta},\frac{1}{q_\theta}\right):= (1-\theta)\left(s_0,\frac{1}{p_0},\frac{1}{q_0}\right)+ \theta\left(s_1,\frac{1}{p_1},\frac{1}{q_1}\right)\text{. }
\end{align*}
If $s_0,s_1\in\NN$, we write $m_0:=s_0$ and $m_1:=s_1$.
We have the following interpolation identities with equivalence of norms
\begin{enumerate}
    \item $(\dot{\B}^{s_0,\sigma}_{p,q_0}(\Omega),\dot{\B}^{s_1,\sigma}_{p,q_1}(\Omega))_{\theta,q}=(\dot{\H}^{s_0,p}_{\sigma}(\Omega),\dot{\H}^{s_1,p}_{\sigma}(\Omega))_{\theta,q}=\dot{\B}^{s,\sigma}_{p,q}(\Omega)$\label{eq:realInterpHomBspqLip}, assuming $s_0,s_1<1$ if $p=\infty$;
    \item $(\dot{\B}^{s_0,0,\sigma}_{\infty,q_0}(\Omega),\dot{\B}^{s_1,0,\sigma}_{\infty,q_1}(\Omega))_{\theta,q}=(\dot{\C}^{m_0}_{0,\sigma}(\overline{\Omega}),\dot{\C}^{m_1}_{0,\sigma}(\overline{\Omega}))_{\theta,q}=\dot{\B}^{s,0,\sigma}_{\infty,q}(\Omega)$\label{eq:realInterpHomBsinftyqLip};
    \item $(\dot{\W}^{m_0,1}_{\sigma}(\Omega),\dot{\W}^{m_1,1}_{\sigma}(\Omega))_{\theta,q}=\dot{\B}^{s,\sigma}_{1,q}(\Omega)$\label{eq:realInterpHomBspqLipL1};
    \item $[\dot{\H}^{s_0,p_0}_{\sigma}(\Omega),\dot{\H}^{s_1,p_1}_{\sigma}(\Omega)]_{\theta} = \dot{\H}^{s,p_\theta}_{\sigma}(\Omega)$, if $1<p_0,p_1<\infty$, and \hyperref[AssumptionCompletenessExponents]{$(\mathcal{C}_{s_j,p_j})$} is true for $j\in\{0,1\}$\label{eq:complexInterpHomSobLip};
    \item $[\dot{\B}^{s_0,\sigma}_{p_0,q_0}(\Omega),\dot{\B}^{s_1,\sigma}_{p_1,q_1}(\Omega)]_{\theta} = \dot{\B}^{s,\sigma}_{p_\theta,q_\theta}(\Omega)$, if  \hyperref[AssumptionCompletenessExponents]{$(\mathcal{C}_{s_j,p_j,q_j})$} is satisfied for $j\in\{0,1\}$, $q_\theta<\infty$.\label{eq:complexInterpHomBspqLip}
\end{enumerate}
and similarly with
\begin{align*}
    \{\dot{\W}_0,\,\dot{\H}_0,\,\dot{\B}_{\cdot,\cdot,0},\,\dot{\C}_{0,0}\}\quad \text{ instead of  }\quad\{\dot{\W},\,\dot{\H},\,\dot{\B},\,\dot{\C}\}.
\end{align*}

Furthermore, 
\begin{itemize}
    \item The identity \ref{eq:realInterpHomBspqLip} still holds if we replace $\dot{\B}^{s_j,\sigma}_{p,q_j}(\Omega)$ by $\dot{\mathcal{B}}^{s_j,\sigma}_{p,\infty}(\Omega)$, $j\in\{0,1\}$;
    \item When one considers $\Omega=\RR^n_+$, one can remove the conditions in \ref{eq:realInterpHomBspqLip} and it also holds that
    \begin{align}\label{eq:realInterpHomBspqLipLinfty}
        (\dot{\C}^{m_0}_{ub,h,\sigma}(\overline{\RR^n_+}),\dot{\C}^{m_1}_{ub,h,\sigma}(\overline{\RR^n_+}))_{\theta,q}=(\dot{\W}^{m_0,\infty}_{\sigma}({\RR^n_+}),\dot{\W}^{m_1,\infty}_{\sigma}({\RR^n_+}))_{\theta,q}=\dot{\B}^{s,\sigma}_{\infty,q}({\RR^n_+}),\\
        (\dot{\C}^{m_0}_{ub,h,0,\sigma}({\RR^n_+}),\dot{\C}^{m_1}_{ub,h,0,\sigma}({\RR^n_+}))_{\theta,q}=(\dot{\W}^{m_0,\infty}_{0,\sigma}({\RR^n_+}),\dot{\W}^{m_1,\infty}_{0,\sigma}({\RR^n_+}))_{\theta,q}=\dot{\B}^{s,\sigma}_{\infty,q,0}({\RR^n_+});\nonumber
    \end{align}
    \item When we replace $(\cdot,\cdot)_{\theta,q}$ by $(\cdot,\cdot)_{\theta}$, all the real interpolation identities still hold with $\dot{\mathcal{B}}^{s,\sigma}_{p,\infty}$ as an output space ($\dot{\mathcal{B}}^{s,0,\sigma}_{\infty,\infty}$ in \ref{eq:realInterpHomBsinftyqLip});
    \item The identity \ref{eq:complexInterpHomBspqLip} remains valid for $q_\theta=\infty$ with an output space $\dot{\BesSmo}^{s_\theta,\sigma}_{p_{\theta},\infty}(\Omega)$.
    \item Every interpolation result and remark above \textbf{remains true} for \textbf{inhomogeneous spaces} \textit{i.e.}, replacing
    \begin{align*}
        \{\dot{\W},\,\dot{\H},\,\dot{\B},\,\dot{\BesSmo},\,\dot{\C}\}\text{ by }\{\W,\,\H,\,\B,\,{\BesSmo},\,\C\},
    \end{align*}
    but requiring additionally that, either, $1<p,p_0,p_1<\infty$ \textbf{or} $s,s_0,s_1>0$.  In this case, one also removes the completeness conditions \hyperref[AssumptionCompletenessExponents]{$(\mathcal{C}_{s_j,p_j})$} and \hyperref[AssumptionCompletenessExponents]{$(\mathcal{C}_{s_j,p_j,q_j})$} in \ref{eq:complexInterpHomSobLip} and \ref{eq:complexInterpHomBspqLip}. For vector fields on $\Omega=\RR^n_+$, the results remain valid provided $1\leqslant p,p_0,p_1\leqslant\infty$ \textbf{and} $s_j>-1+1/p_j$, $j\in\{ 0,1\}$.
\end{itemize}
\end{theorem}

\begin{proof}We mostly focus on $\overline{\Omega}$-supported co-closed differential forms in $\RR^n$: similar argument may be applied or reproduced for the other cases. As usual the main issue is the (loss of) completeness of homogeneous function spaces of higher order.

\textbf{Step 1:} The case of homogeneous function spaces. Real interpolation of Besov spaces.

By construction it can be checked that under the assumptions of \ref{eq:realInterpHomBspqLip}
\begin{align*}
    (\dot{\B}^{s_0,\sigma}_{p,q_0,0}(\Omega),\dot{\B}^{s_1,\sigma}_{p,q_1,0}(\Omega))_{\theta,q} \hookrightarrow (\dot{\B}^{s_0}_{p,q_0,0}(\Omega),\dot{\B}^{s_1}_{p,q_1,0}(\Omega))_{\theta,q} = \dot{\B}^{s}_{p,q,0}(\Omega).
\end{align*}
But since $(\dot{\B}^{s_0,\sigma}_{p,q_0,0}(\Omega),\dot{\B}^{s_1,\sigma}_{p,q_1,0}(\Omega))_{\theta,q}\subset \dot{\B}^{s_0,\sigma}_{p,q_0,0}(\Omega)+\dot{\B}^{s_1,\sigma}_{p,q_1,0}(\Omega)$, we obtain
\begin{align*}
    (\dot{\B}^{s_0,\sigma}_{p,q_0,0}(\Omega),\dot{\B}^{s_1,\sigma}_{p,q_1,0}(\Omega))_{\theta,q} \hookrightarrow  \dot{\B}^{s,\sigma}_{p,q,0}(\Omega).
\end{align*}

\textbf{Step 1.1:} For the reverse embedding, assume first that $s_0$ and $q_0$ are such that \hyperref[AssumptionCompletenessExponents]{$(\mathcal{C}_{s_0,p,q_0})$} is satisfied. Let $u\in\dot{\B}^{s,\sigma}_{p,q,0}(\Omega)\subset\dot{\B}^{s}_{p,q,0}(\Omega)\subset \dot{\B}^{s_0}_{p,q_0,0}(\Omega)+\dot{\B}^{s_1}_{p,q_1,0}(\Omega)$, and consider $(\mathfrak{a},\mathfrak{b})\in\dot{\B}^{s_0}_{p,q_0,0}(\Omega)\times\dot{\B}^{s_1}_{p,q_1,0}(\Omega)$ such that $u=\mathfrak{a}+\mathfrak{b}$. Since $\delta u =0$, by Theorem~\ref{thm:PotentialOpSpeLipDom}, for all $0<a<b<\infty$, we obtain
\begin{align*}
    \delta \T^{0,\sigma}_{a,b} u +\T^{0,\sigma}_{a,b} \delta u = \delta \T^{0,\sigma}_{a,b} u  = \delta \T^{0,\sigma}_{a,b}  \mathfrak{a}+ \delta \T^{0,\sigma}_{a,b} \mathfrak{b}
\end{align*}
By Remark~\ref{rem:convergencePotOpSpeLipNoncomplete}, the fact $\delta u=0$, and since \hyperref[AssumptionCompletenessExponents]{$(\mathcal{C}_{s_0,p,q_0})$}, by Theorem~\ref{thm:PotentialOpSpeLipDom} it holds that in $\S'_h(\RR^n)$
\begin{align*}
    \delta \T^{0,\sigma}_{a,b} u = \delta \T^{0,\sigma}_{a,b} u +\T^{0,\sigma}_{a,b} \delta u &\xrightarrow[\substack{a \longrightarrow 0_+\\ b \longrightarrow \infty}]{} u,\\
    \delta \T^{0,\sigma}_{a,b}  \mathfrak{a} &\xrightarrow[\substack{a \longrightarrow 0_+\\ b \longrightarrow \infty}]{} \delta \T^{0,\sigma}  \mathfrak{a}.
\end{align*}
Therefore, $(\delta \T^{0,\sigma}_{a,b}  \mathfrak{b})_{0<a<b<\infty} = (\delta \T^{0,\sigma}_{a,b}  [u-\mathfrak{a}])_{0<a<b<\infty}$ admits a limit in $\mathcal{S}'_h(\mathbb{R}^n)$, and we can write
\begin{align*}
    u = \delta \T^{0,\sigma}  \mathfrak{a}+ \delta \T^{0,\sigma} \mathfrak{b}
\end{align*}
with the estimates for $t>0$
\begin{align*}
    K(t,u,\dot{\B}^{s_0,\sigma}_{p,q_0,0}(\Omega),\dot{\B}^{s_1,\sigma}_{p,q_1,0}(\Omega)) &\leqslant \lVert \delta \T^{0,\sigma}  \mathfrak{a}\rVert_{\dot{\B}^{s_0}_{p,q_0}(\RR^n)} + t \lVert \delta \T^{0,\sigma}  \mathfrak{b}\rVert_{\dot{\B}^{s_1}_{p,q_1}(\RR^n)}\\
    &\lesssim_{p,s,n,\Omega}  \lVert  \mathfrak{a}\rVert_{\dot{\B}^{s_0}_{p,q_0}(\RR^n)} + t \lVert  \mathfrak{b}\rVert_{\dot{\B}^{s_1}_{p,q_1}(\RR^n)}.
\end{align*}
Taking the infimum on all such pairs $(\mathfrak{a},\mathfrak{b})$, one obtains
\begin{align*}
    K(t,u,\dot{\B}^{s_0,\sigma}_{p,q_0,0}(\Omega),\dot{\B}^{s_1,\sigma}_{p,q_1,0}(\Omega)) &\lesssim_{p,s,n,\Omega}  K(t,u,\dot{\B}^{s_0}_{p,q_0,0}(\Omega),\dot{\B}^{s_1}_{p,q_1,0}(\Omega)).
\end{align*}
Therefore, we deduce the reverse embedding by the inequalities
\begin{align*}
    \lVert u \rVert_{(\dot{\B}^{s_0,\sigma}_{p,q_0,0}(\Omega),\dot{\B}^{s_1,\sigma}_{p,q_1,0}(\Omega))_{\theta,q}} \lesssim_{p,s,n,\Omega} \lVert u \rVert_{(\dot{\B}^{s_0}_{p,q_0,0}(\Omega),\dot{\B}^{s_1}_{p,q_1,0}(\Omega))_{\theta,q}} \sim_{p,s,n,\Omega} \lVert u \rVert_{\dot{\B}^{s}_{p,q}(\RR^n)}.
\end{align*}

\textbf{Step 1.2:} Now, we remove the condition \hyperref[AssumptionCompletenessExponents]{$(\mathcal{C}_{s_0,p,q_0})$}. Consider arbitrary $s_0<s_1$, and by the previous step, let $\alpha<s_0$ , $r\in[1,\infty]$, $0<\eta<1$, such that \hyperref[AssumptionCompletenessExponents]{$(\mathcal{C}_{\alpha,p,r})$} is satisfied and
\begin{align*}
    (\dot{\B}^{\alpha,\sigma}_{p,r,0}(\Omega),\dot{\B}^{s_1,\sigma}_{p,q_1,0}(\Omega))_{\eta,q_0} = \dot{\B}^{s_0,\sigma}_{p,q_0,0}(\Omega).
\end{align*}
By the reiteration result \cite[Lemma~C.1]{Gaudin2023Lip}, and the previous step, it holds
\begin{align*}
    (\dot{\B}^{s_0,\sigma}_{p,q_0,0}(\Omega),\dot{\B}^{s_1,\sigma}_{p,q_1,0}(\Omega))_{\theta,q} &= \Big((\dot{\B}^{\alpha,\sigma}_{p,r,0}(\Omega),\dot{\B}^{s_1,\sigma}_{p,q_1,0}(\Omega))_{\eta,q_0},\dot{\B}^{s_1,\sigma}_{p,q_1,0}(\Omega)\Big)_{\theta,q}\\
    &=(\dot{\B}^{\alpha,\sigma}_{p,r,0}(\Omega),\dot{\B}^{s_1,\sigma}_{p,q_1,0}(\Omega))_{(1-\theta)\eta+\theta,q} = \dot{\B}^{s,\sigma}_{p,q,0}(\Omega).
\end{align*}

\textbf{Step 1.3:} We remove the condition $s_0,s_1<1$ in \ref{eq:realInterpHomBspqLip} provided $p=\infty$, whenever $\Omega=\RR^n_+$. This obstruction came from the special Lipschitz domain-case in \cite[Theorem~3.46]{Gaudin2023Lip}. Since in the same aforementioned statement, this restriction can be removed when $\Omega=\RR^n_+$, it holds that previous Step 1.1 and Step 1.2 still apply.

\textbf{Step 2:} Complex interpolation of homogeneous function spaces.

In the complex interpolation statements all function spaces are complete, and  Theorem~\ref{thm:PotentialOpSpeLipDom} provides a bounded projection operator
\begin{align*}
    \P^{0,\sigma}:=\delta \T^{0,\sigma}&\,:\,\dot{\H}^{s,p}_0(\Omega) \longrightarrow\dot{\H}^{s,p}_{0,\sigma}(\Omega),\qquad 1<p<\infty \\
    &\,:\,\dot{\B}^{s}_{p,q,0}(\Omega) \longrightarrow\dot{\B}^{s,\sigma}_{p,q,0}(\Omega).
\end{align*}
The result follows from a retraction and co-retraction as in the proof of Theorem~\ref{thm:InterpHomSpacesBddLip}:
\begin{align*}
    [\dot{\H}^{s_0,p_0}_{0,\sigma}(\Omega),\dot{\H}^{s_1,p_0}_{0,\sigma}(\Omega)]_{\theta} &= [\P^{0,\sigma}\dot{\H}^{s_0,p_0}_{0}(\Omega),\P^{0,\sigma}\dot{\H}^{s_1,p_1}_{0}(\Omega)]_{\theta}\\
    &= \P^{0,\sigma}[\dot{\H}^{s_0,p_0}_{0}(\Omega),\dot{\H}^{s_1,p_1}_{0}(\Omega)]_{\theta}\\
    &= \P^{0,\sigma}\dot{\H}^{s,p}_{0}(\Omega)\\
    &= \dot{\H}^{s,p}_{0,\sigma}(\Omega).
\end{align*}
The same goes for homogeneous Besov spaces.

\textbf{Step 3:} Complex and real interpolation of inhomogeneous function spaces.

\textbf{Step 3.1:} The case $1<p<\infty$. By Remark~\ref{rem:convergencePotOpSpeLipNoncomplete} and  Theorem~\ref{thm:PotentialOpSpeLipDom}, since $\H^{\beta,p}(\RR^n) = \dot{\H}^{\beta,p}(\RR^n)\cap\L^p(\RR^n)$ and $\H^{-\beta,p}(\RR^n) = \dot{\H}^{-\beta,p}(\RR^n)+\L^p(\RR^n)$, $\beta\geqslant 0$, for all $s\in\RR$, we obtain bounded linear maps
\begin{align*}
    \delta \T^{0,\sigma},\, \T^{0,\sigma}\delta  \,:\,\H^{s,p}_0(\Omega) \longrightarrow\H^{s,p}_{0}(\Omega) 
\end{align*}
and by real interpolation for $q\in[1,\infty]$,
\begin{align*}
    \delta \T^{0,\sigma},\, \T^{0,\sigma}\delta   \,:\,\B^{s}_{p,q,0}(\Omega) \longrightarrow\B^{s}_{p,q,0}(\Omega).
\end{align*}
with $\I =\delta \T^{0,\sigma} +\T^{0,\sigma}\delta$. In particular, this yields bounded linear projections
\begin{align*}
    \P^{0,\sigma}:=\delta \T^{0,\sigma}&\,:\,\H^{s,p}_{0}(\Omega) \longrightarrow\H^{s,p}_{0,\sigma}(\Omega), \\
    &\,:\,\B^{s}_{p,q,0}(\Omega) \longrightarrow\B^{s,\sigma}_{p,q,0}(\Omega).
\end{align*}
Then, the result follows from a retraction and co-retraction as in the proof of Theorem~\ref{thm:InterpHomSpacesBddLip}.

\textbf{Step 3.2:} The case  $s,s_0,s_1>0$. By Theorem~\ref{thm:SteinsExtensionOpDivFree}, for all $m_j\in\NN$, $p\in[1,\infty)$, $j\in\NN$,
\begin{align*}
    \mathcal{E}_{\sigma}\,&:\,\W^{m_j,p}(\Omega,\Lambda) \longrightarrow\W^{m_j,p}(\RR^n,\Lambda),\\
    \mathcal{P}_{0,\sigma}\,&:\,\W^{m_j,1}(\RR^n,\Lambda) \longrightarrow\W^{m_j,1}_0(\Omega,\Lambda)\\
    \mathcal{E}_{\sigma}\,&:\,\C^{m_j}_{ub}(\overline{\Omega},\Lambda) \longrightarrow\C^{m_j}_{ub}(\RR^n,\Lambda),\\
    \mathcal{P}_{0,\sigma}\,&:\,\C^{m_j}_{ub}(\RR^n,\Lambda) \longrightarrow\C^{m_j}_{ub,0}(\Omega,\Lambda)
\end{align*}
and both preserve co-closed forms. Therefore, by real interpolation for all $s>0$, all $p,q\in[1,\infty]$, and since $\mathcal{E}_{\sigma}(\L^p_{\sigma}(\Omega))\subset \L^p_{\sigma}(\RR^n)$ and $\mathcal{P}_{0,\sigma}(\L^p_{\sigma}(\RR^n))\subset \L^p_{\sigma}(\RR^n)$, one obtains bounded linear extension and projection operators
\begin{align*}
    \mathcal{E}_{\sigma}\,&:\,\B^{s,\sigma}_{p,q}(\Omega) \longrightarrow\B^{s,\sigma}_{p,q}(\RR^n),\\
    \mathcal{P}_{0,\sigma}\,&:\,\B^{s,\sigma}_{p,q}(\RR^n) \longrightarrow\B^{s,\sigma}_{p,q,0}(\Omega).
\end{align*}
The same applies to Bessel potential spaces by complex interpolation. 
The results of complex and real interpolation then follow from a retraction and co-retraction argument as in Step 2.

\textbf{Step 3.3:} The case of inhomgoeneous function spaces of vector fields, provided $\Omega=\RR^n_+$, whenever $1\leqslant p,p_0,p_1\leqslant\infty$  $s_j>-1+1/p_j$, $j\in\{ 0,1\}$. The proof is actually verbatim the one of Step 3.2 using the extension and projection operators $\E^{m}_{\sigma}$ and $\P_{0,\sigma}^m$ from Remark~ \ref{rem:FlatExtOpDiv}, provided $m$ is large enough (such that $s_0,s_1<m$ is sufficient).

\textbf{Step 4:} We deal with real interpolation of any other function spaces. It has been proved before that \ref{eq:realInterpHomBspqLip} and \ref{eq:realInterpHomBsinftyqLip} hold when Besov spaces are considered. The remaining cases are obtained through the natural embeddings for all $p\in[1,\infty]$, $m\in\NN$, $s\in\RR$,
\begin{align*}
    \dot{\B}^{s,\sigma}_{p,1,0}(\Omega)\hookrightarrow&\, \dot{\H}^{s,p}_{0,\sigma}(\Omega) \hookrightarrow \dot{\B}^{s,\sigma}_{p,\infty,0}(\Omega),\\
    \dot{\B}^{m,\sigma}_{p,1,0}(\Omega)\hookrightarrow&\, \dot{\W}^{m,p}_{0,\sigma}(\Omega) \hookrightarrow \dot{\B}^{m,\sigma}_{p,\infty,0}(\Omega),\\
    \dot{\B}^{m,0,\sigma}_{\infty,1,0}(\Omega)\hookrightarrow&\, \dot{\C}^{m}_{0,0,\sigma}(\Omega) \hookrightarrow \dot{\B}^{m,0,\sigma}_{\infty,\infty,0}(\Omega),\\
    \dot{\B}^{m,\sigma}_{\infty,1,0}(\Omega)\hookrightarrow&\, \dot{\C}^{m}_{ub,h,0,\sigma}(\Omega) \hookrightarrow \dot{\B}^{m,\sigma}_{\infty,\infty,0}(\Omega),
\end{align*}
and similarly for other function spaces.
\end{proof}

\subsubsection{Additional important results}

\begin{lemma}\label{lem:weakStartdensityCcinftydiv}Let $p\in[1,\infty]$, $s\in\RR$. Let $\Omega$ be either a bounded or special Lipschitz domain. Then the space $\Ccinftydiv(\Omega)$ is weakly-$\ast$ dense 
\begin{enumerate}
    \item in the function space ${\B}^{s,\sigma}_{p,\infty,0}(\Omega)$;
    \item  and if $\Omega$ is a special Lipschitz domain, in the homogeneous Besov spaces $\dot{\B}^{s,\sigma}_{p,\infty,0}(\Omega)$, provided $-n/p'\leqslant s<n/p$.
\end{enumerate}
\end{lemma}

\begin{proposition}\label{prop:IdentifVanishingDivFree}Let $p,q\in[1,\infty]$, $s\in(-1+\sfrac{1}{p},\sfrac{1}{p})$. Let $\Omega$ be a bounded Lipschitz domain or a special Lipschitz domain.  The extension by $0$ to the whole space provides the canonical identifications
\begin{enumerate}
    \item if $-1+\sfrac{1}{p}<s<\sfrac{1}{p}$, $\B^{s,\sigma}_{p,q,0}(\Omega) \simeq \B^{s,\sigma}_{p,q,\mathfrak{n}}(\Omega)$;
    \item if $\sfrac{1}{p}<s<1+\sfrac{1}{p}$, $\B^{s,\sigma}_{p,q,0}(\Omega) \simeq \B^{s,\sigma}_{p,q,\mathcal{D}}(\Omega)$.
\end{enumerate}
The result still holds
\begin{itemize}
    \item  for the Sobolev spaces $\H^{s,p}$, assuming $1<p<\infty$;
    \item  for the Besov spaces $\BesSmo^{s}_{p,\infty}$;
    \item  when $\Omega$ is a special Lipschitz domain, for homogeneous Sobolev and Besov spaces $\dot{\B}^{s}_{p,q}$, $\dot{\BesSmo}^{s}_{p,\infty}$, $\dot{\B}^{s,0}_{\infty,q}$, $\dot{\BesSmo}^{s,0}_{\infty,\infty}$, and $\dot{\H}^{s,p}$, assuming $1<p<\infty$ for the latter;
    \item in a similar way with
    \begin{align*}
        \{ u\in\L^1_{\sigma}(\RR^n)\,:\, \supp u \subset\overline{\Omega}\} &\simeq \L^1_{\mathfrak{n},\sigma}(\Omega)\\
        \dot{\W}^{1,1}_{0,\sigma}(\Omega)=\{ u\in\dot{\W}^{1,1}_{\sigma}(\RR^n)\,:\, \supp u \subset\overline{\Omega}\} &\simeq \dot{\W}^{1,1}_{\mathcal{D},\sigma}(\Omega)\\
        \C^{0}_{0,0,\sigma}(\Omega)=\{ u\in\C^{0}_{0,\sigma}(\RR^n)\,:\, \supp u\subset\overline{\Omega}\}  &\simeq \{ u\in\C^{0}_{0,\sigma}(\overline{\Omega})\,:\,  u_{|_{\partial\Omega}} =0\};\\
        \C^{0}_{ub,0,\sigma}(\Omega)=\{ u\in\C^{0}_{ub,\sigma}(\RR^n)\,:\, \supp u \subset\overline{\Omega}\} &\simeq \{ u\in\C^{0}_{ub,\sigma}(\overline{\Omega})\,:\,  u_{|_{\partial\Omega}} =0\};\\
        \C^{0}_{ub,h,0,\sigma}(\Omega)=\{ u\in\C^{0}_{ub,h,\sigma}(\RR^n)\,:\, \supp u \subset\overline{\Omega}\} &\simeq \{ u\in\C^{0}_{ub,h,\sigma}(\overline{\Omega})\,:\,  u_{|_{\partial\Omega}} =0\};\\
        \{ u\in\L^\infty_{\sigma}(\RR^n)\,:\, \supp u \subset\overline{\Omega}\} &\simeq \L^\infty_{\mathfrak{n},\sigma}(\Omega);\\
        \{ u\in\L^\infty_{h,\sigma}(\RR^n)\,:\, \supp u \subset\overline{\Omega}\} &\simeq \L^\infty_{h,\mathfrak{n},\sigma}(\Omega).
    \end{align*}
\end{itemize}
\end{proposition}

\begin{proof}Note that the cases $\L^1$, $\L^\infty$, $\dot{\B}^s_{p,q}(\Omega)$ and $\dot{\H}^{s,p}(\Omega)$, $-1+1/p<s<1/p$ are particular cases of Proposition~\ref{prop:Ext0PartialTraceVanish} where $\delta u=0$. In the case $1/p<s<1+1/p$, or if $\dot{\W}^{1,1}$ and $\C^0$-type spaces are considered, the result follows from Proposition~\ref{prop:Ext0DirLip}. Indeed, as in the proof of Proposition~\ref{prop:Ext0PartialTraceVanish}, one obtains $\delta \widetilde{u} = \widetilde{\delta u}$ in $\mathcal{D}'(\RR^n,\Lambda)$, where, here, $v\mapsto \tilde{v}$ stands for the extension from $\Omega$ to $\RR^n$ by $0$.
\end{proof}

This yields the following interpolation result:

\begin{proposition}\label{prop:InterpDivFreeC1} Let $p,p_0,p_1,q,q_0,q_1\in[1,\infty]$, $-1+{\sfrac{1}{p_j}}<s_j<2+{\sfrac{1}{p_j}}$, $s_j\neq 1/p_j,1+1/p_j$, $j\in\{0,1\}$ and $\Omega$ be a bounded Lipschitz domain. For all $\theta\in(0,1)$, such that $s:=s_0(1-\theta)+\theta s_1<1+\sfrac{1}{p}$, the following interpolation identities hold true
\begin{enumerate}
    \item If $s_0<s<s_1$, $p=p_0=p_1$,
    \begin{align*}
        ({\H}^{s_0,p}_{{\mathcal{D}},\sigma}(\Omega),{\H}^{s_1,p}_{{\mathcal{D}},\sigma}(\Omega))_{\theta,q} =({\B}^{s_0,\sigma}_{p,q_0,{\mathcal{D}}}(\Omega),{\B}^{s_1,\sigma}_{p,q_1,{\mathcal{D}}}(\Omega))_{\theta,q} = {\B}^{s,\sigma}_{p,q,\mathcal{D}}(\Omega);
    \end{align*}
    \item If $1<p_0,p_1<\infty$, and $\sfrac{1}{p}=\sfrac{1-\theta}{p_0}+\sfrac{\theta}{p_1}$, 
    \begin{align*}
        [{\H}^{s_0,p_0}_{{\mathcal{D}},\sigma}(\Omega),{\H}^{s_1,p_1}_{{\mathcal{D}},\sigma}(\Omega)]_{\theta}={\H}^{s,p}_{{\mathcal{D}},\sigma}(\Omega);
    \end{align*}
    \item Provided $\sfrac{1}{r} = \sfrac{1-\theta}{r_0}+\sfrac{\theta}{r_1}$, $r\in\{p,q\}$, $q<\infty$,
    \begin{align*}
        [{\B}^{s_0,\sigma}_{p_0,q_0,{\mathcal{D}}}(\Omega),{\B}^{s_1,\sigma}_{p_1,q_1,{\mathcal{D}}}(\Omega)]_{\theta}={\B}^{s,\sigma}_{p,q,{\mathcal{D}}}(\Omega).
    \end{align*}
\end{enumerate}
Furthermore,
\begin{itemize}
    \item the same real interpolation result \textit{(i)} still holds replacing $(\cdot,\cdot)_{\theta,\infty}$ and ${\B}^{s,\sigma}_{p,\infty,\mathcal{D}}$ by respectively $(\cdot,\cdot)_{\theta}$ and ${\BesSmo}^{s,\sigma}_{p,\infty,\mathcal{D}}$;
    \item the complex interpolation result \textit{(iii)} holds when $q=\infty$ up to replace ${\B}^{s,\sigma}_{p,\infty,\mathcal{D}}$ by ${\BesSmo}^{s,\sigma}_{p,\infty,\mathcal{D}}$.
\end{itemize}
\end{proposition}

\begin{proof}It is straightforward to notice that
\begin{align*}
    ({\B}^{s_0,\sigma}_{p,q_0,{\mathcal{D}}}(\Omega),{\B}^{s_1,\sigma}_{p,q_1,{\mathcal{D}}}(\Omega))_{\theta,q} \hookrightarrow {\B}^{s}_{p,q}(\Omega).
\end{align*}
But then since ${\B}^{s_1,\sigma}_{p,q_1,{\mathcal{D}}}(\Omega)={\B}^{s_0,\sigma}_{p,q_0,{\mathcal{D}}}(\Omega)\cap{\B}^{s_1,\sigma}_{p,q_1,{\mathcal{D}}}(\Omega)$ is dense in $({\B}^{s_0,\sigma}_{p,q_0,{\mathcal{D}}}(\Omega),{\B}^{s_1,\sigma}_{p,q_1,{\mathcal{D}}}(\Omega))_{\theta,q}$ whenever $q<\infty$, see \cite[Theorem~3.4.2]{BerghLofstrom1976}, one preserves the divergence-free/co-closed condition as well as the targeted(Dirichlet or no-slip) boundary  condition. 

For the reverse embedding, note that one always has with equivalence of norms
\begin{align*}
    {\B}^{\tilde{s},\sigma}_{p,r,0}(\Omega)\hookrightarrow {\B}^{\tilde{s},\sigma}_{p,r,\mathcal{D}}(\Omega)
\end{align*}
with canonical identification whenever $\tilde{s}\in(-1+\sfrac{1}{r},1+\sfrac{1}{r})$, $\tilde{s}\neq\sfrac{1}{r}$ by Proposition~\ref{prop:IdentifVanishingDivFree}. So that the reverse embedding holds, thanks to the condition $s<1+\sfrac{1}{p}$,
\begin{align*}
    {\B}^{s,\sigma}_{p,q,{\mathcal{D}}}(\Omega) = {\B}^{s,\sigma}_{p,q,0}(\Omega) = ({\B}^{s_0,\sigma}_{p,q_0,0}(\Omega),{\B}^{s_1,\sigma}_{p,q_1,0}(\Omega))_{\theta,q} \hookrightarrow ({\B}^{s_0,\sigma}_{p,q_0,{\mathcal{D}}}(\Omega),{\B}^{s_1,\sigma}_{p,q_1,{\mathcal{D}}}(\Omega))_{\theta,q}.
\end{align*}

The case $q=\infty$ follows from the reiteration theorem.

The case of real interpolation of Bessel potential spaces follows from the sandwich embedding
\begin{align*}
    {\B}^{\tilde{s},\sigma}_{r,1,\mathcal{D}}(\Omega)\hookrightarrow {\H}^{\tilde{s},r}_{\mathcal{D},\sigma}(\Omega)\hookrightarrow {\B}^{\tilde{s},\sigma}_{r,\infty,\mathcal{D}}(\Omega)
\end{align*}
valid for all $ \tilde{s}\in(-1+\sfrac{1}{r},2+\sfrac{1}{r})$, $\tilde{s}\neq\sfrac{1}{r},1+\sfrac{1}{r}$, $r\in[1,\infty]$.

For complex interpolation, the same arguments are still valid.
\end{proof}

The case of homogeneous function spaces admits a similar proof.

\begin{proposition}\label{prop:HomInterpDivFreeC1} Let $p,p_0,p_1,q,q_0,q_1\in[1,\infty]$, $-1+{\sfrac{1}{p_j}}<s_j<2+{\sfrac{1}{p_j}}$, $s_j\neq 1/p_j,1+1/p_j$, $j\in\{0,1\}$ and $\Omega$ be a special Lipschitz domain. For all $\theta\in(0,1)$, such that $s:=s_0(1-\theta)+\theta s_1<1+\sfrac{1}{p}$, the following interpolation identities hold true
\begin{enumerate}
    \item If $s_0<s<s_1$, $p=p_0=p_1$, assuming additionally $s_1<1$ if $p=\infty$,\label{eq:HomInterpDivFreeSpeLipBC}
    \begin{align*}
        (\dot{\H}^{s_0,p}_{{\mathcal{D}},\sigma}(\Omega),\dot{\H}^{s_1,p}_{{\mathcal{D}},\sigma}(\Omega))_{\theta,q} =(\dot{\B}^{s_0,\sigma}_{p,q_0,{\mathcal{D}}}(\Omega),\dot{\B}^{s_1,\sigma}_{p,q_1,{\mathcal{D}}}(\Omega))_{\theta,q} = \dot{\B}^{s,\sigma}_{p,q,\mathcal{D}}(\Omega);
    \end{align*}
    \item If $1<p_0,p_1<\infty$, and $\sfrac{1}{p}=\sfrac{1-\theta}{p_0}+\sfrac{\theta}{p_1}$, provided \hyperref[AssumptionCompletenessExponents]{$(\mathcal{C}_{s_j,p_j})$}, $j\in\{0,1\}$,
    \begin{align*}
        [\dot{\H}^{s_0,p_0}_{{\mathcal{D}},\sigma}(\Omega),\dot{\H}^{s_1,p_1}_{{\mathcal{D}},\sigma}(\Omega)]_{\theta}=\dot{\H}^{s,p}_{{\mathcal{D}},\sigma}(\Omega);
    \end{align*}
    \item Provided $\sfrac{1}{r} = \sfrac{1-\theta}{r_0}+\sfrac{\theta}{r_1}$, $r\in\{p,q\}$, $q<\infty$, and \hyperref[AssumptionCompletenessExponents]{$(\mathcal{C}_{s_j,p_j,q_j})$}, $j\in\{0,1\}$,\label{eq:HomCompInterpDivFreeSpeLipBC}
    \begin{align*}
        [\dot{\B}^{s_0,\sigma}_{p_0,q_0,{\mathcal{D}}}(\Omega),\dot{\B}^{s_1,\sigma}_{p_1,q_1,{\mathcal{D}}}(\Omega)]_{\theta}=\dot{\B}^{s,\sigma}_{p,q,{\mathcal{D}}}(\Omega).
    \end{align*}
\end{enumerate}
Furthermore,
\begin{itemize}
    \item the same real interpolation result \ref{eq:HomInterpDivFreeSpeLipBC} still holds replacing $(\cdot,\cdot)_{\theta,\infty}$ and $\dot{\B}^{s,\sigma}_{p,\infty,\mathcal{D}}$ by respectively $(\cdot,\cdot)_{\theta}$ and $\dot{\BesSmo}^{s,\sigma}_{p,\infty,\mathcal{D}}$;
    \item the complex interpolation result \ref{eq:HomCompInterpDivFreeSpeLipBC} holds when $q=\infty$ up to replace $\dot{\B}^{s,\sigma}_{p,\infty,\mathcal{D}}$ by $\dot{\BesSmo}^{s,\sigma}_{p,\infty,\mathcal{D}}$;
    \item if $\Omega=\RR^n_+$, then one can remove the condition "$s_1<1$ if $p=\infty$" in \ref{eq:HomInterpDivFreeSpeLipBC}.
\end{itemize}
\end{proposition}

\subsection{The Hodge Laplacian, the Hodge-Dirac operator and boundedness of the Hodge-Leray projection as consequences}\label{Sec:HodgeLapHodgeDirHodgeLeray}

This section is dedicated to the analysis of the Hodge Laplacian and the Hodge Dirac operator motivated by two main points:
\begin{itemize}
    \item The analysis of the Hodge Laplacian on the half-space, in the continuation of \cite{Gaudin2023Hodge}, will alow us later on to provide higher order estimates of the Stokes--Dirichlet resolvent problem in Section~\ref{Sec:StokesHalfSpace}. One will even consider endpoint function spaces such as $\L^\infty$ or $\B^{s}_{\infty,q}$.
    \item The analysis of the Hodge-Leray projection which, as exhibited below, ties to the careful analysis of the regularity properties of the Hodge Laplacian. Since the Hodge Laplacian $-\Delta_\mathcal{H}$ is the square of the Hodge-Dirac operator $D_\mathfrak{n}=\d+\d^\ast$, when on bounded domains we will reduce our analysis to the  resolvent problem of the corresponding first order system, following an idea of M${}^\text{c}$Intosh and Monniaux \cite{McintoshMonniaux2018}.
\end{itemize}

The analysis of incompressible flows, in particular the Navier--Stokes equations, relies on the \textbf{Helmholtz decomposition} 
\[
\uu = \vv+\nabla \mathfrak{q}, \qquad \div \vv=0 \; (\text{with possibly } \vv\cdot\nu_{|_{\partial\Omega}}=0),
\] 
which provides the "incompressible part" $v$ and the potential part $\nabla \mathfrak{q}$ of $u$.   In $\L^2$, one obtains the topological orthogonal decomposition
\[
\L^2(\Omega,\mathbb{C}^n)
= \L^2_{\mathfrak{n},\sigma}(\Omega)\overset{\perp}{\oplus}\overline{\nabla {\H}^{1,2}(\Omega)},
\]
via the Helmholtz--Leray\footnote{or Hodge--Leray, or Hodge--Helmholtz, or whatever combination. However, the convention seem to be the following: the name of Hodge is involved if one is more involved in the framework of Differential Geometry/De Rham Cohomology. The name of Helmholtz is concerned when one considers such a decomposition mostly at the level of vector fields. The name of Leray obviously comes in mind whenever the analysis incompressible fluid flows is involved.} projector $\PP_\Omega$, initially only defined on $\L^2(\Omega,\CC^n)$ as an orthogonal projection, see \textit{e.g.} \cite[Chapter~2]{SohrBook2001}.  

The extension to the \textbf{$\L^p$-theory}, $p\in(1,\infty)$, $p\neq 2$, is more delicate and deeply linked to Harmonic Analysis and sharp elliptic regularity. That is to know whether one has
\begin{align*}
    \L^p(\Omega,\mathbb{C}^n)
= \L^p_{\mathfrak{n},\sigma}(\Omega)\overset{\perp}{\oplus}\overline{\nabla {\H}^{1,p}(\Omega)}.
\end{align*}
For bounded Lipschitz domains, Fabes--Mendez--Mitrea \cite[Theorem~12.2]{FabesMendezMitrea1998} proved the decomposition for $p\in(3/2-\varepsilon,3+\varepsilon)$. Simader--Sohr \cite[Theorem~1.4]{SimaderSohr1992} obtained it for bounded and exterior $\C^1$-domains for all $p\in(1,\infty)$. However, Bogovski\u{\i} \cite{Bogovskii1986} produced counterexamples in unbounded domains. When special Lipschitz domains were concerned Tolksdorf \cite[Theorem~5.1.10]{TolksdorfPhDThesis2017} obtained the admissible range $p\in\big(\frac{2n}{n+1}-\varepsilon,\frac{2n}{n-1}+\varepsilon\big)$. More general situations are treated in \cite{FarwigKozonoSohr2005,FarwigKozonoSohr2007}.  

A key tool here is the \textbf{Neumann problem} for the Laplacian. On bounded domains, the solvability of
\[
-\Delta \mathfrak{q} = \div \uu \quad \text{in }\Omega, 
\qquad \partial_\nu \mathfrak{q} = \uu\cdot \nu \quad \text{on }\partial\Omega,
\]
provides the gradient component in the decomposition. The well-posedness of this problem (see \cite{SimaderSohr1992,FabesMendezMitrea1998}) underlies uniqueness of $v$ and $q$, up to constants, and shows how the decomposition is tied to elliptic theory.

A natural question is the extension to the Sobolev and Besov scales:
\[
{\H}^{s,p}(\Omega,\CC^n) = {\H}^{s,p}_{\mathfrak{n},\sigma}(\Omega)\oplus 
\overline{\nabla {\H}^{s+1,p}(\Omega,\CC)}, \qquad s\in(-1+\sfrac{1}{p},\sfrac{1}{p}),\; p\in(1,\infty).
\]
The result {\cite[Proposition~2.16,~Remark]{MitreaMonniaux2008}} by Mitrea and Monniaux implies such a decomposition for all indices mentioned above when $\Omega$ is a bounded $\C^1$-domain. It also holds for function spaces that are close the family $(\H^{s,2})_{|s|<\frac{1}{2}}$ when $\Omega$ is a bounded Lipschitz domain. Concerning a result that includes fractional endpoint function spaces, such as Besov spaces $\B^{s}_{1,q}$, $\B^{-s}_{\infty,q}$, $0<s<1$, $q\in[1,\infty]$, the only known result to the best of authors' knowledge was achieved in \cite[Theorem~4.2]{FujiwaraYamazaki2007}. However, this work by Fujiwara and Yamazaki only concerned vector fields and $\C^{2,1}$- bounded (and exterior) domains.

\medbreak

Following \cite[Introduction]{Gaudin2023Hodge}, for \textbf{differential forms} one could instead consider the general \textbf{Hodge decomposition} :
\[
\dot{\H}^{s,p}(\Omega,\Lambda^k)= 
\dot{\H}^{s,p}_{\mathfrak{n},\sigma}(\Omega,\Lambda^k)
\oplus \overline{\d\,{\H}^{s+1,p}(\Omega,\Lambda^{k-1})},
\]
looking for the analogous decomposition
\[
u = v+\d \boldsymbol{\omega}, \quad \delta v=0, \; \text{ and }\; \nu\iprod v_{|_{\partial\Omega}}=0.
\] 

Here, the "gradient"-like component $\d \boldsymbol{\omega}$, which is actually called the exact part,  arises by solving the \textbf{Hodge-Laplace elliptic problem}:
\[
-\Delta_{\mathcal{H}} \boldsymbol{\omega} = \delta u \quad \text{in }\Omega, 
\qquad \nu\iprod \boldsymbol{\omega} = 0, \;\; \nu\iprod [\d\boldsymbol{\omega}-u] = 0 \quad \text{on }\partial\Omega,
\]
for a $(k-1)$-differential form $\boldsymbol{\omega}$.  Then $\d\boldsymbol{\omega}$ provides the exact part of $u$, and $u-\d\boldsymbol{\omega}$ lies in the divergence-free/co-closed subspace. Note that for $k=1$, $\nu\iprod \boldsymbol{\omega} = 0$ simplifies trivially as $0 = 0$. This is the natural analogue of the Neumann problem in the setting of differential forms, and   the \textbf{Hodge Laplacian}
\[
-\Delta_{\mathcal{H}}=(\d+\d^\ast)^2
\]
respects the decomposition according to \cite{GeissertHeckTrunk2013,MonniauxShen2018} in the case of vector fields.

Classical results are due to Schwarz on smooth manifolds \cite{Schwartz1995} and more recently to M${}^\text{c}$Intosh  and  Monniaux \cite{McintoshMonniaux2018} for Lipschitz domains of $\RR^n$ on $\L^p$, within a range of exponent $p$ near $2$. Very recently Monniaux did obtain the Hodge decomposition for differential forms of arbitrary degrees is proved for bounded $\C^1$-domains on $\L^p$ spaces for any $1<p<\infty$, \cite[Theorem~2.2]{Monniaux2025} .

\medbreak

A important idea of  M${}^\text{c}$Intosh and Monniaux \cite{McintoshMonniaux2018} was about exhibiting that this analysis essentially reduces to the analysis of the Hodge-Dirac operator on the appropriate scale of function spaces, instead of the direct analysis of its square the Hodge Laplacian.

\medbreak

\subsubsection{A preliminary study on the flat half-space involving the Hodge Laplacian}\label{subsec:HodgeRn+}

When on the half-space, we can define for $\X\in \{\H^{s,p},\,\B^{s}_{p,q},\,\dot{\H}^{s,p},\,\dot{\B}^{s}_{p,q},\,\W^{m,p},\,\dot{\W}^{m,p}\}$, $k\in\llb 0,n\rrb$,
\begin{align}\label{eq:HodgeSpaces}
    \X_\mathcal{H}(\RR^n_+,\Lambda^k) &:=\{\,u\in\X(\RR^n_+,\Lambda^k)\,:\,\forall I'\in\mathcal{I}^{n-1}_{k-1},\,u_{I',n}\in\X_0(\RR^n_+,\CC)\,\},\\
    \X_{\mathcal{H}_\ast}(\RR^n_+,\Lambda^k) &:=\{\,u\in\X(\RR^n_+,\Lambda^k)\,:\,\forall I\in\mathcal{I}^{n-1}_{k},\,u_{I}\in\X_0(\RR^n_+,\CC)\,\},\nonumber
\end{align}
and the smooth compactly supported counterpart,
\begin{align}\label{eq:HodgeSpacesCinfty}
    \C^\infty_{\tilde{c},\mathcal{H}}(\overline{\RR^n_+},\Lambda^k) &:=\{\,u\in\Ccinfty(\overline{\RR^n_+},\Lambda^k)\,:\,\forall I'\in\mathcal{I}^{n-1}_{k-1},\,u_{I',n}\in\Ccinfty(\RR^n_+,\CC)\,\}\\
    \C^\infty_{\tilde{c},\mathcal{H}_\ast}(\overline{\RR^n_+},\Lambda^k) &:=\{\,u\in\Ccinfty(\overline{\RR^n_+},\Lambda^k)\,:\,\forall I\in\mathcal{I}^{n-1}_{k},\,u_{I}\in\Ccinfty(\RR^n_+,\CC)\,\}.\nonumber
\end{align}

For a measurable function $u\,:\,\mathbb{R}^n_+\longrightarrow \Lambda^k$, provided $k\in\llb 0,n\rrb$, $I\in\mathcal{I}^{k}_{n}$, we define
\begin{align}\label{eq:HodgeReflexExtOp}
    (\mathrm{E}_{\mathcal{H}}u)_{I}:= \begin{cases}
  \mathrm{E}_{\mathcal{D}}u_{I}\,&\text{, if } n\in I\text{,}\\    
  \mathrm{E}_{\mathcal{N}}u_{I}\,&\text{, if } n\notin I\text{;} 
    \end{cases}\, \text{ and }\, (\mathrm{E}_{\mathcal{H}_{\ast}}u)_{I}:= \begin{cases}
  \mathrm{E}_{\mathcal{N}}u_{I}\,&\text{, if } n\in I\text{,}\\    
  \mathrm{E}_{\mathcal{D}}u_{I}\,&\text{, if } n\notin I\text{;} 
    \end{cases}
\end{align}

\begin{proposition}\label{prop:ExtHodgeRn+} Let $p,q\in[1,\infty]$, $-2+{\sfrac{1}{p}}<s<1+{\sfrac{1}{p}}$, such that $s-\sfrac{1}{p}\notin\ZZ$ and $k\in\llb 0,n\rrb$.  One has a bounded map
\begin{align*}
    \E_\mathcal{H}\,:&\,\dot{\B}^{s}_{p,q,\mathcal{H}}(\RR^n_+,\Lambda^k) \longrightarrow \dot{\B}^{s}_{p,q}(\RR^n,\Lambda^k)\text{, }-1+{\sfrac{1}{p}}<s<1+{\sfrac{1}{p}},\\
    \,:&\,\dot{\B}^{s}_{p,q,\mathcal{H}_\ast}(\RR^n_+,\Lambda^k) \longrightarrow \dot{\B}^{s}_{p,q}(\RR^n,\Lambda^k)\text{, }-2+{\sfrac{1}{p}}<s<{\sfrac{1}{p}},
\end{align*}
Furthermore,
\begin{itemize}
    \item When $q=\infty$, we can replace $\dot{\B}^{s}_{p,\infty}$ by $\dot{\BesSmo}^{s}_{p,\infty}$.
    \item The result still holds if we replace $\dot{\B}^{s}_{p,q}$ by either $\dot{\B}^{s,0}_{\infty,q}$, $\dot{\W}^{1,1}$ or $\dot{\H}^{s,p}$ with $1<p<\infty$ for the latter. 
    \item The same result holds exchanging the roles of $\mathcal{H}$ and $\mathcal{H}_\ast$.
    \item The result still holds for their inhomogeneous counterparts.
\end{itemize}
\end{proposition}

\begin{proof} This is a direct consequence from Lemmas~\ref{lem:ExtDirNeuRn+}~\&~\ref{lem:ExtOpNegativeBesovSpaces}.
\end{proof}

\begin{proposition}\label{prop:HodgeFuncSpacesRn+} Let $p,q\in[1,\infty]$, $-2+{\sfrac{1}{p}}<s<1+{\sfrac{1}{p}}$, such that \eqref{AssumptionCompletenessExponents} is satisfied  and $k\in\llb 0,n\rrb$.  The following duality identities hold:
\begin{align*}
\begin{array}{rlrlll}
    (\dot{\H}^{-s,p'}_{\mathcal{H}_{\ast}}(\RR^n_+,\Lambda^k))'&=\dot{\H}^{s,p}_{\mathcal{H}}(\RR^n_+,\Lambda^k)\, \text{,}&&\text{where }1<p<\infty\text{;}\\
    (\dot{\B}^{-s}_{p',q',\mathcal{H}_{\ast}}(\RR^n_+,\Lambda^k))'&=\dot{\B}^{s}_{p,q,\mathcal{H}}(\RR^n_+,\Lambda^k)\, \text{,}&&\text{where } p,q>1\text{;}\\
    (\dot{\mathcal{B}}^{-s}_{p',\infty,\mathcal{H}_{\ast}}(\RR^n_+,\Lambda^k))'&=\dot{\B}^{s}_{p,1,\mathcal{H}}(\RR^n_+,\Lambda^k)\, \text{,}&&\text{where }p>1,q=1\text{;}\\
    (\dot{\B}^{-s,0}_{\infty,q'\mathcal{H}_{\ast}}(\RR^n_+,\Lambda^k))'&=\dot{\B}^{s}_{1,q,\mathcal{H}}(\RR^n_+,\Lambda^k)\, \text{,}&&\text{where }p=1,q>1\text{;}\\
    (\dot{\mathcal{B}}^{-s,0}_{\infty,\infty,\mathcal{H}_{\ast}}(\RR^n_+,\Lambda^k))'&=\dot{\B}^{s}_{1,1,\mathcal{H}}(\RR^n_+,\Lambda^k)\, \text{,}&&\text{where }p=q=1\text{.}
\end{array}
\end{align*}
Furthermore, if $p,q<\infty$, $\C_{\tilde{c},\mathcal{H}}^\infty(\overline{\RR^n_+},\Lambda^k)$ (up to taking the intersection) is a strongly dense subspace of $\dot{\H}^{s,p}_{\mathcal{H}}(\RR^n_+,\Lambda^k)$ ($1<p<\infty$), $\dot{\B}^{s}_{p,q,\mathcal{H}}(\RR^n_+,\Lambda^k)$, $\dot{\mathcal{B}}^{s}_{p,\infty,\mathcal{H}}(\RR^n_+,\Lambda^k)$, $\dot{\B}^{s,0}_{\infty,q,\mathcal{H}}(\RR^n_+,\Lambda^k)$ and in $\dot{\mathcal{B}}^{s,0}_{\infty,\infty,\mathcal{H}}(\RR^n_+,\Lambda^k)$.

Moreover,
\begin{itemize}
    \item the statement remains true for inhomogeneous function spaces;
    \item one can exchange the roles of $\mathcal{H}$ and $\mathcal{H}_\ast$.
\end{itemize}
\end{proposition}

\begin{proof} This follows from the definition of function spaces \eqref{eq:HodgeSpaces} and \eqref{eq:HodgeSpacesCinfty}.
\end{proof}

\begin{proposition}\label{prop:GenNeuPbRn+}Let $p,q\in[1,\infty]$, $-1+{\sfrac{1}{p}}<s<{\sfrac{1}{p}}$, $k\in\llb 0,n\rrb$. Let $f\in\dot{\B}^{s-1}_{p,q,\mathcal{H}_\ast}(\RR^n_+,\Lambda^k)$, $F\in\dot{\B}^{s}_{p,q}(\RR^n_+,\Lambda^{k+1})$, $g\in\dot{\B}^{s+1-1/p}_{p,q}(\partial\RR^n_+,\Lambda^{k-1})$  and $h\in\dot{\B}^{s-1/p}_{p,q}(\partial\RR^n_+,\Lambda^k)$ such that $\mathfrak{e}_n\iprod (g,h) =0$. 

Then, the following generalized Neumann problem
\begin{equation*}\tag{$\mathcal{GNL}_{0}$}\label{NeuLap}
    \left\{ \begin{array}{rllr}
         - \Delta \mathfrak{p} &= f+ \delta F \text{, }&&\text{ in } \RR^n_+\text{,}\\
        (-\mathfrak{e}_n)\iprod \mathfrak{p}_{|_{\partial\RR^n_+}} &= g\text{, }&&\text{ on } \partial\RR^n_+\text{,}\\
        (-\mathfrak{e}_n)\iprod \left(\d \mathfrak{p}-F\right)_{|_{\partial\RR^n_+}} &= h\text{, }&&\text{ on } \partial\RR^n_+\text{.}
    \end{array}
    \right.
\end{equation*}
if it admits a solution $\mathfrak{p}\in\dot{\B}^{s+1}_{p,q}(\RR^n_+,\Lambda^k)$, then it is unique. Moreover, it obeys the estimates,
\begin{align*}
  \lVert \nabla \mathfrak{p}\rVert_{\dot{\B}^{s}_{p,q}(\RR^n_+)}\lesssim_{p,n,s} \lVert  F\rVert_{\dot{\B}^{s}_{p,q}(\RR^n_+)}+ \lVert  f\rVert_{\dot{\B}^{s-1}_{p,q}(\RR^n_+)} + \lVert  g\rVert_{\dot{\B}^{s+1-1/p}_{p,q}(\partial\RR^n_+)}+ \lVert  h\rVert_{\dot{\B}^{s-1/p}_{p,q}(\partial\RR^n_+)}.
\end{align*}
Furthermore,
\begin{itemize}
    \item if $\delta f=0$, $\delta' g =0$ and $\delta'h =0$ then $\delta \mathfrak{p} = 0$;

    \item if $f\in\dot{\B}^{\alpha_0-1}_{r_0,\kappa_0,\mathcal{H}_\ast}(\RR^n_+,\Lambda^k)$, $F\in\dot{\B}^{\alpha_1}_{r_1,\kappa_1}(\RR^n_+,\Lambda^{k+1})$, $g\in\dot{\B}^{\alpha_2+1-\sfrac{1}{r_2}}_{r_2,\kappa_2}(\partial\RR^n_+,\Lambda^{k-1})$  and $h\in\dot{\B}^{\alpha_3-\sfrac{1}{r_3}}_{r_3,\kappa_3}(\partial\RR^n_+,\Lambda^k)$ for some $r_j,\kappa_j\in[1,\infty]$, $\alpha_j\in(-1+\sfrac{1}{r_j},\sfrac{1}{r_j})$ satisfying \hyperref[AssumptionCompletenessExponents]{$(\mathcal{C}_{\alpha_j+1,r_j,\kappa_j})$}  for all $j\in\llb 0,3\rrb$, then the solution always exists;

    \item a similar result holds for the endpoint Besov spaces $\dot{\BesSmo}^{s}_{p,\infty}$, $\dot{\B}^{s,0}_{\infty,q}$, and  $\dot{\BesSmo}^{s,0}_{\infty,\infty}$;

    \item in the case $p=q\in(1,\infty)$, the result still holds replacing the spaces $(\dot{\B}^{\bullet}_{p,q}(\RR^n_+),\dot{\B}^{\bullet}_{p,q}(\partial\RR^n_+))$ by $(\dot{\H}^{\bullet,p}(\RR^n_+),\dot{\B}^{\bullet}_{p,p}(\partial\RR^n_+))$.
\end{itemize}
\end{proposition}

\begin{remark} Due to low regularity issues, solving the problem \eqref{NeuLap}, especially the meaning of traces, might be unclear. We say that $\mathfrak{p}$ solves \eqref{NeuLap} if it satisfies $(-\mathfrak{e}_n)\iprod\mathfrak{p}_{|_{\partial\RR^n_+}}=g$ in the strong trace sense, and
\begin{align*}
    \big\langle \d\mathfrak{p}-F,\, \d\Psi\big\rangle_{\RR^n_+} + \big\langle \delta\mathfrak{p},\, \delta \Psi\big\rangle_{\RR^n_+} = \big\langle f,\, \Psi\big\rangle_{\RR^n_+} +\big\langle h,\, \Psi_{|_{\partial\RR^n_+}}\big\rangle_{\partial\RR^n_+}\text{, }\forall \Psi\in\C^\infty_{\tilde{c},\mathcal{H}}(\overline{\RR^n_+},\Lambda^k) .
\end{align*}
Furthermore, notice that if $\delta f=0$, then the meaning of the trace equality $(-\mathfrak{e}_n)\iprod \left(\d \mathfrak{p}-F\right)_{|_{\partial\RR^n_+}} = h$ is strong by Theorem~\ref{thm:partialtracesDiffFormRn+}.

The choice of $\C^\infty_{\tilde{c},\mathcal{H}}(\overline{\RR^n_+},\Lambda^k)$ for the space of test functions is a consequence of both Proposition~\ref{prop:HodgeFuncSpacesRn+} and the trace theorem Theorem~\ref{thm:partialtracesDiffFormRn+}.
\end{remark}

\begin{proof} \textbf{Step 1:} Uniqueness. Let $\mathfrak{p}\in\dot{\B}^{s+1}_{p,\infty}(\RR^n_+,\Lambda^k)$ such that
\begin{equation*}
    \left\{ \begin{array}{rllr}
         - \Delta \mathfrak{p} &= 0 \text{, }&&\text{ in } \RR^n_+\text{,}\\
        (-\mathfrak{e}_n)\iprod \mathfrak{p}_{|_{\partial\RR^n_+}} &= 0\text{, }&&\text{ on } \partial\RR^n_+\text{,}\\
        (-\mathfrak{e}_n)\iprod \d \mathfrak{p}_{|_{\partial\RR^n_+}} &= 0\text{, }&&\text{ on } \partial\RR^n_+\text{.}
    \end{array}
    \right.
\end{equation*}
Notice that the trace condition $(-\mathfrak{e}_n)\iprod \d \mathfrak{p}_{|_{\partial\RR^n_+}}=0$ is well-defined since $\div(\nabla\mathfrak{p})=0 \in\dot{\B}^{s}_{p,q}(\RR^n_+,\Lambda^{k})$ (the divergence and the gradient are applied component-wise). Hence, by the trace theorem, one obtains $\partial_{x_j}\mathfrak{p}\in\C^{\mathrm{w}\ast}_{ub,x_n}(\overline{\RR_+},\dot{\B}^{s-1/p}_{p,\infty}(\RR^{n-1},\Lambda^k))$ for all $j\in\llb1,n\rrb$. Thus,
\begin{align}\label{eq:traceMeaning}
    (-\mathfrak{e}_n)\iprod \d \mathfrak{p}_{|_{\partial\RR^n_+}} = (-\mathfrak{e}_n)\iprod\left(\sum_{j=1}^n\sum_{I\in\mathcal{I}_{k}^{n}} \partial_{x_j}\mathfrak{p}_{I}(\cdot,0)\, \d{x_j}\wedge\d{x_I}\right)\in\dot{\B}^{s-1/p}_{p,\infty}(\RR^{n-1},\Lambda^k).
\end{align}

Notice that the first condition $(-\mathfrak{e}_n)\iprod \mathfrak{p}_{|_{\partial\RR^n_+}} = 0$ becomes
\begin{align*}
    \sum_{I'\in\mathcal{I}_{k-1}^{n-1}} \mathfrak{p}_{I',n}(\cdot,0)\,\d{x_{I'}} =0.
\end{align*}
Thus $-\Delta \mathfrak{p}_{I',n} = 0$ with $\mathfrak{p}_{I',n}{}_{|_{\partial\RR^n_+}} =0$ for all $I'\in\mathcal{I}_{k-1}^{n-1}$. By injectivity of the Dirichlet Laplacian, Proposition~\ref{prop:DirPbRn+}, we obtain $\mathfrak{p}_{I',n}=0$. Plugging it in \eqref{eq:traceMeaning}, we obtain
\begin{align*}
    \sum_{I\in\mathcal{I}_{k}^{n-1}} \partial_{x_n}\mathfrak{p}_{I}(\cdot,0)\, \d{x_I}=0.
\end{align*}
Thus it reduces to the injectivity of the Neumann Laplacian. Let $\mathbf{p}\in\dot{\B}^{s+1}_{p,\infty}(\RR^n_+,\CC)$ be such that $-\Delta \mathbf{p} =0$. We show that $\Delta \E_{\mathcal{N}}\mathbf{p}=0$ in $\mathcal{D}'(\RR^n)$. Let $\varphi\in\Ccinfty(\RR^n)$, one has
\begin{align*}
     \big\langle -\Delta \E_\mathcal{N}\mathbf{p},\, \varphi \big\rangle_{\RR^n} = \big\langle  \E_\mathcal{N}\mathbf{p},\, -\Delta \varphi \big\rangle_{\RR^n} &= \big\langle  \mathbf{p},\, -\Delta \varphi \big\rangle_{\RR^n_+} + \big\langle  \mathbf{p},\, -\Delta \varphi(\cdot,-\cdot) \big\rangle_{\RR^n_+}\\
     &= \big\langle  \mathbf{p}, -\Delta [\varphi+\varphi(\cdot,-\cdot)] \big\rangle_{\RR^n_+}\\
     &= \big\langle  \nabla \mathbf{p},\, \nabla [\varphi+\varphi(\cdot,-\cdot)] \big\rangle_{\RR^n_+}\\
     &= \big\langle  -\Delta \mathbf{p},\, [\varphi+\varphi(\cdot,-\cdot)] \big\rangle_{\RR^n_+}=0.
\end{align*}
Hence, $\E_\mathcal{N}\mathbf{p}$ is harmonic on $\RR^n$, thus it is a constant by Liouville's Theorem. Since $\E_\mathcal{N}\mathbf{p}\in\S'_h(\RR^n)$, it ensures that $\E_\mathcal{N}\mathbf{p} =0$ and then $\mathbf{p}=0$.

\textbf{Step 2:} Existence. We can decompose the problem in a natural way by linearity.

\textbf{Step 2.1:} First notice that for $f,F=0$, $\mathfrak{p}$ should be given by
\begin{align*}
    \begin{array}{rlrlrl}
    \mathfrak{p}_{I',n}(\cdot,x_n) &:= (-1)^ke^{-x_n(-\Delta')^{\sfrac{1}{2}}}g_{I'}\text{,}&&x_n\geqslant 0\text{,} &&I'\in\mathcal{I}^{n-1}_{k-1},\\
    \mathfrak{p}_{I}(\cdot,x_n) &:= -(-\Delta')^{-\sfrac{1}{2}}e^{-x_n(-\Delta')^{\sfrac{1}{2}}}[h_{I}-(\d'g)_{I}]\text{,}&&x_n\geqslant 0\text{,}&& I\in\mathcal{I}^{n-1}_{k}.
    \end{array}
\end{align*}
Written differently:
\begin{align*}
    \mathfrak{p}(\cdot,x_n) = (-1)^k [e^{-x_n(-\Delta')^{\sfrac{1}{2}}}g] \wedge \d x_n - (-\Delta')^{-\sfrac{1}{2}}[e^{-x_n(-\Delta')^{\sfrac{1}{2}}}(h-\d'g)]\text{, }x_n\geqslant 0\text{.}
\end{align*}

One has $\Delta \mathfrak{p}=0$ and then the partial traces are well-defined in the strong sense  by Proposition~\ref{prop:PoissonSemigroup3} and the trace theorems, Theorems~\ref{thm:tracesRn+}~\&~\ref{thm:partialtracesDiffFormRn+}.

This procedure has to be modified a little bit if $\dot{\B}^{s+1}_{p,q}(\RR^n_+,\Lambda^k)$ is not a complete space. Assume first that $q<\infty$. For $N\in\NN$, we set
\begin{align*}
    \mathfrak{p}_N(\cdot,x_n) = (-1)^k [e^{-x_n(-\Delta')^{\sfrac{1}{2}}}g_N] \wedge \d x_n - (-\Delta')^{-\sfrac{1}{2}}[e^{-x_n(-\Delta')^{\sfrac{1}{2}}}(h_N-\d'g_N)]\text{, }x_n\geqslant 0\text{.}
\end{align*}
where we defined
\begin{align*}
    & g_N := \sum_{|j|\leqslant N} \Delta_j'g \in \bigcap_{\substack{\zeta\in\RR\\ 1\leqslant \kappa \leqslant \infty }}\dot{\B}^{\zeta}_{p,\kappa}(\partial\RR^n_+,\Lambda^{k-1})\text{,}\\
    \text{ and }& h_N := \sum_{|j|\leqslant N} \Delta_j'h \in  \bigcap_{\substack{\zeta\in\RR\\ 1\leqslant \kappa \leqslant \infty }}\dot{\B}^{\zeta}_{p,\kappa}(\partial\RR^n_+,\Lambda^{k})\text{.}
\end{align*}

For $N,N'\in\NN$, by Proposition~\ref{prop:PoissonSemigroup3} and one has
\begin{align*}
    \lVert \d \mathfrak{p}_N\rVert_{\dot{\B}^{s}_{p,q}(\RR^n_+)} &\leqslant\lVert \nabla \mathfrak{p}_N\rVert_{\dot{\B}^{s}_{p,q}(\RR^n_+)} \lesssim_{n,s} \lVert  g_N\rVert_{\dot{\B}^{s+1-\sfrac
    {1}{p}}_{p,q}(\partial\RR^n_+)} + \lVert  h_N\rVert_{\dot{\B}^{s-\sfrac
    {1}{p}}_{p,q}(\partial\RR^n_+)} ,\\
    \lVert \d \mathfrak{p}_N - \d \mathfrak{p}_{N'} \rVert_{\dot{\B}^{s}_{p,q}(\RR^n_+)}&\leqslant\lVert \nabla \mathfrak{p}_N - \nabla \mathfrak{p}_{N'}\rVert_{\dot{\B}^{s}_{p,q}(\RR^n_+)}  \lesssim_{n,s} \lVert  g_N-g_{N'}\rVert_{\dot{\B}^{s+1-\sfrac
    {1}{p}}_{p,q}(\partial\RR^n_+)}\\
    &\qquad\qquad\qquad\qquad\qquad\qquad\qquad\qquad\qquad+ \lVert  h_N-h_{N'}\rVert_{\dot{\B}^{s-\sfrac
    {1}{p}}_{p,q}(\partial\RR^n_+)}.
\end{align*}

Therefore, $(g_N)_{N\in\NN}$ and $( h_N)_{N\in\NN}$ being Cauchy sequences, an then so are $( \nabla {\mathfrak{p}}_N)_{N\in\NN}$ and $( \d {\mathfrak{p}}_N)_{N\in\NN}$. Consequently, there exists ${\mathfrak{q}}\in\mathcal{D}'(\RR^n_+,\Lambda^{k-1})$ uniquely determined up to a constant $k-1$-form such that $( \nabla \mathfrak{p}_N)_{N\in\NN}$ and $( \d \mathfrak{p}_N)_{j\in\NN}$ converge, respectively, to $\nabla \mathfrak{q}$ and $\d \mathfrak{q}$. If the equivalence class of $\mathfrak{q}$ contains an element in $\S'_h(\RR^n_+,\Lambda^{k-1})$ up to a constant $k-1$-form, for $\mathfrak{p}$ to be such a representative, by approximation by $(\mathfrak{p}_N)_{N\in\NN}$ and continuity of traces, it must solve the corresponding generalized Neumann problem.

It is straightforward to check that $\delta \mathfrak{p} =0$ whenever $\delta'(g,h)=0$.

\textbf{Step 2.2:} Now, if $f=0$, $g=0$ and $h=0$, for $F\in \dot{\B}^{s}_{p,q}(\RR^n_+,\Lambda^{k+1})$, one has $\E_\mathcal{H} F \in \dot{\B}^{s}_{p,q}(\RR^n,\Lambda^{k+1})$, and $\tilde{\mathfrak{p}}:=(-\Delta)^{-1}\delta (\E_\mathcal{H} F) \in\dot{\B}^{s+1}_{p,q}(\RR^n,\Lambda^{k})$, with $\delta \tilde{\mathfrak{p}} =0$, and $(-\Delta)\tilde{\mathfrak{p}} = \delta(\E_\mathcal{H} F)$ in $\S'(\RR^n,\Lambda^k)\subset\mathcal{D}'(\RR^n,\Lambda^k)$.

In particular, 
\begin{align*}
    \big\langle \d\tilde{\mathfrak{p}}-(\E_\mathcal{H} F),\, \d\Psi\big\rangle_{\RR^n} =\big\langle \d\tilde{\mathfrak{p}}-(\E_\mathcal{H} F),\, \d\Psi\big\rangle_{\RR^n} + \big\langle \delta\tilde{\mathfrak{p}},\, \delta \Psi\big\rangle_{\RR^n} = 0\text{, }\forall \Psi\in\Ccinfty({\RR^n},\Lambda^k).
\end{align*}

We can set $\mathfrak{p}:= \tilde{\mathfrak{p}}_{|_{\RR^n_+}}$, so, by restriction, we still have $\delta {\mathfrak{p}} =0$, and $(-\Delta){\mathfrak{p}} = \delta F$ in $\mathcal{D}'(\RR^n_+,\Lambda^k)$. It suffices to check that the boundary conditions are satisfied. Since $\delta(-\Delta)^{-1}=(-\Delta)^{-1}\delta$ on $\RR^n$, one can also write $\tilde{\mathfrak{p}} = (-\Delta)^{-\sfrac{1}{2}}\delta(-\Delta)^{-\sfrac{1}{2}} \E_\mathcal{H} F$. From the $\L^r$-theory \cite[Section~3]{Gaudin2023Hodge}, if we assume temporarily $F\in\L^r\cap \dot{\B}^{s}_{p,q}(\RR^n_+,\Lambda^{k+1})$ for $r\in(1,n)$, one has $\tilde{\mathfrak{p}} = \E_\mathcal{H}[ (-\Delta_\mathcal{H})^{-\sfrac{1}{2}}\d^\ast(-\Delta_\mathcal{H})^{-\sfrac{1}{2}} F]$, and consequently $\mathfrak{p} = (-\Delta_\mathcal{H})^{-\sfrac{1}{2}}\d^\ast(-\Delta_\mathcal{H})^{-\sfrac{1}{2}} F \in \dot{\H}^{1,r}\cap \dot{\B}^{s+1}_{p,q}(\RR^n_+,\Lambda^{k+1})$, with $(-\mathfrak{e}_n)\iprod \mathfrak{p}_{|_{\partial\RR^n_+}} = 0$.

Since $\delta (\d \mathfrak{p} - F)=0$, the partial trace $(-\mathfrak{e}_n)\iprod (\d\mathfrak{p}-F)_{|_{\partial\RR^n_+}}$ is well-defined by Theorem~\ref{thm:partialtracesDiffFormRn+}, and we do have for all $\Psi\in\C^\infty_{\tilde{c},\mathcal{H}}(\overline{\RR^n_+},\Lambda^k)$
\begin{align}\label{eq:ProofTraceEqualityHidgeLap}
    \big\langle \d\mathfrak{p}-F,\, \d\Psi\big\rangle_{\RR^n_+} + \big\langle \delta\mathfrak{p},\, \delta \Psi\big\rangle_{\RR^n_+} &= \big\langle \d (-\Delta_\mathcal{H})^{-\sfrac{1}{2}}\d^\ast(-\Delta_\mathcal{H})^{-\sfrac{1}{2}} F-F,\, \d\Psi\big\rangle_{\RR^n_+}\nonumber\\
    &=\big\langle \d^\ast(-\Delta_\mathcal{H})^{-\sfrac{1}{2}}\d (-\Delta_\mathcal{H})^{-\sfrac{1}{2}} F,\, \d\Psi\big\rangle_{\RR^n_+} = 0.
\end{align}
Thus $(-\mathfrak{e}_n)\iprod (\d\mathfrak{p}-F)_{|_{\partial\RR^n_+}}=0$. If $p,q<\infty$ and $(\mathcal{C}_{s+1,p,q})$ is satisfied, and we obtain the estimate
\begin{align*}
    \lVert \d \mathfrak{p}\rVert_{\dot{\B}^{s}_{p,q}(\RR^n_+)}\leqslant\lVert \nabla \mathfrak{p}\rVert_{\dot{\B}^{s}_{p,q}(\RR^n_+)} \leqslant \lVert \nabla(-\Delta)^{-\sfrac{1}{2}}\delta(-\Delta)^{-\sfrac{1}{2}} \E _\mathcal{H} F\rVert_{\dot{\B}^{s}_{p,q}(\RR^n)}\lesssim_{p,n,s} \lVert  F\rVert_{\dot{\B}^{s}_{p,q}(\RR^n_+)}.
\end{align*}
Thus, the equality \eqref{eq:ProofTraceEqualityHidgeLap}, the trace properties, as well as the estimate remain valid for all $F\in\dot{\B}^{s}_{p,q}(\RR^n_+,\Lambda^{k+1})$, since $\L^r\cap \dot{\B}^{s}_{p,q}(\RR^n_+,\Lambda^{k+1})$ is strongly dense in $\dot{\B}^{s}_{p,q}(\RR^n_+,\Lambda^{k+1})$. The case $q=\infty$ follows by real interpolation. If $p=\infty$, we do have $\L^r\cap\dot{\B}^{s}_{\infty,q}(\RR^n_+) = \L^r\cap\dot{\B}^{s,0}_{\infty,q}(\RR^n_+)$, and $\L^r\cap\dot{\B}^{s}_{\infty,\infty}(\RR^n_+) = \L^r\cap\dot{\BesSmo}^{s,0}_{\infty,\infty}(\RR^n_+)$. This yields the case of function spaces $\dot{\B}^{s,0}_{\infty,q}$, $q\in[1,\infty]$, and $\dot{\BesSmo}^{s,0}_{\infty,\infty}$. Now, for the function spaces $\dot{\B}^{s}_{\infty,q}$, $q\in[1,\infty)$ and $\dot{\BesSmo}^{s}_{\infty,\infty}$, provided $-1<s<0$, we consider for $j\in\NN$
\begin{align*}
    \tilde{F}_j := [\dot{S}_j-\dot{S}_{-j}]\E_\mathcal{H} F \in \bigcap_{\substack{\zeta\in\RR\\ 1\leqslant \kappa \leqslant \infty }}\dot{\B}^{\zeta}_{\infty,\kappa}(\RR^n,\Lambda^{k+1}),
\end{align*}
so that $(\tilde{F}_j)_{j\in\NN}$ converges to $\E_\mathcal{H} F$ in $\dot{\B}^{s}_{\infty,q}(\RR^n,\Lambda^{k+1})$ and then, by restriction $({F}_j)_{j\in\NN}:=(\tilde{F}_j{}_{|_{\RR^n_+}})_{j\in\NN}$ converges to  $F$ in $\dot{\B}^{s}_{\infty,q}(\RR^n_+,\Lambda^{k+1})$ (and similarly for $\dot{\BesSmo}^{s}_{\infty,\infty}$ up to appropriate changes). 

One sets for all $j\in\NN$, $\tilde{\mathfrak{p}}_j:=(-\Delta)^{-\sfrac{1}{2}}\delta(-\Delta)^{-\sfrac{1}{2}} \tilde{F}_j$, and ${\mathfrak{p}}_j:=\tilde{\mathfrak{p}}_j{}_{|_{\RR^n_+}}$, which are smooth and belong to $\dot{\B}^{\zeta}_{\infty,\kappa}$ for all $\zeta\in\RR$, all $\kappa\in[1,\infty]$. One has
\begin{align*}
    (-\mathfrak{e}_n)\iprod \tilde{\mathfrak{p}}_j &= (-\mathfrak{e}_n)\iprod\sum_{\ell=1}^{n}\sum_{I\in\mathcal{I}^{n}_{k+1}} \partial_{x_\ell}  (-\Delta)^{-1}[\dot{S}_j-\dot{S}_{-j}]\E_\mathcal{H} F_I [(-\mathfrak{e}_\ell)\iprod \d x_I]\\
    &= \sum_{\ell=1}^{n-1}\sum_{I'\in\mathcal{I}^{n-1}_{k}} (-1)^{k}\partial_{x_\ell}  (-\Delta)^{-1}[\dot{S}_j-\dot{S}_{-j}](\E_\mathcal{H} F)_{I',n} [(-\mathfrak{e}_\ell)\iprod \d x_{I'}]\\
    &= \sum_{\ell=1}^{n-1}\sum_{I'\in\mathcal{I}^{n-1}_{k}} (-1)^{k}  (-\Delta)^{-1}[\dot{S}_j-\dot{S}_{-j}]\E_\mathcal{D} [\partial_{x_\ell} F_{I',n}] [(-\mathfrak{e}_\ell)\iprod \d x_{I'}].
\end{align*}
Thus, it is straightforward to obtain $(-\mathfrak{e}_n)\iprod \tilde{\mathfrak{p}}_j(\cdot,0) = (-\mathfrak{e}_n)\iprod {\mathfrak{p}}_j{}_{|_{\partial\RR^n_+}} = 0$. In a similar way, writing
\begin{align*}
    (-\mathfrak{e}_n)\iprod (\d \tilde{\mathfrak{p}}_j - \Tilde{F}_j) =  (-\mathfrak{e}_n)\iprod \delta(-\Delta)^{-\sfrac{1}{2}}\d(-\Delta)^{-\sfrac{1}{2}}\Tilde{F}_j =  (-\mathfrak{e}_n)\iprod (-\Delta)^{-1} \delta\d \Tilde{F}_j.
\end{align*}
As before, it can be checked that $(-\mathfrak{e}_n)\iprod (\d \tilde{\mathfrak{p}}_j - \Tilde{F}_j)(\cdot,0) = (-\mathfrak{e}_n)\iprod (\d {\mathfrak{p}}_j - {F}_j)_{|_{\partial\RR^n_+}}=0$. For $j,j'\in\NN$, we have the estimates
\begin{align*}
    \lVert \d \mathfrak{p}_j\rVert_{\dot{\B}^{s}_{\infty,q}(\RR^n_+)} &\leqslant\lVert \nabla \mathfrak{p}_j\rVert_{\dot{\B}^{s}_{\infty,q}(\RR^n_+)} \lesssim_{n,s} \lVert  F_j\rVert_{\dot{\B}^{s}_{\infty,q}(\RR^n_+)},\\
    \lVert \d \mathfrak{p}_j - \d \mathfrak{p}_{j'} \rVert_{\dot{\B}^{s}_{\infty,q}(\RR^n_+)}&\leqslant\lVert \nabla \mathfrak{p}_j - \nabla \mathfrak{p}_{j'}\rVert_{\dot{\B}^{s}_{\infty,q}(\RR^n_+)}  \lesssim_{n,s} \lVert  F_j-F_{j'}\rVert_{\dot{\B}^{s}_{\infty,q}(\RR^n_+)}.
\end{align*}
Therefore, $(\tilde{F}_j)_{j\in\NN}$ and $( F_j)_{j\in\NN}$ being Cauchy sequences, so are $( \nabla \tilde{\mathfrak{p}}_j)_{j\in\NN}$ and $( \d \tilde{\mathfrak{p}}_j)_{j\in\NN}$, and then $( \nabla \mathfrak{p}_j)_{j\in\NN}$ and $( \d \mathfrak{p}_j)_{j\in\NN}$. Consequently, there exists $\tilde{\mathfrak{q}}\in\S'(\RR^n,\Lambda^{k-1})$ uniquely determined up to a constant $k-1$-form such that $( \nabla \mathfrak{p}_j)_{j\in\NN}$ and $( \d \mathfrak{p}_j)_{j\in\NN}$ converge, respectively, to $\nabla \mathfrak{q}$ and $\d \mathfrak{q}$, where $\mathfrak{q}:=\tilde{\mathfrak{q}}_{|_{\RR^n_+}}$. If the equivalence class of $\mathfrak{q}$ contains an element in $\S'_h(\RR^n_+,\Lambda^{k-1})$ up to a constant $k-1$-form, for $\mathfrak{p}$ to be such a representative, by approximation by $(\mathfrak{p}_j)_{j\in\NN}$ and continuity of traces, it must solve the generalized Neumann problem.

By linearity, for the case $f\in\dot{\B}^{s-1}_{p,q,\mathcal{H}_{\ast}}(\RR^n_+,\Lambda^k)$, and $F=0$, $G=0$, $h=0$, $g=0$, it suffices to reproduce everything with $\mathfrak{p}=[(-\Delta)^{-1}\E_\mathcal{H}f]_{|_{\RR^n_+}}$, thanks to Proposition~\ref{prop:ExtHodgeRn+}.

Finally, the case $q=\infty$, for the space $\dot{\B}^{s}_{\infty,\infty}$, follows by real interpolation.
\end{proof}
The next Theorem is then an immediate corollary.

\begin{theorem}\label{thm:HodgeDecompRn+1}Let $p,q\in[1,\infty]$, $s\in(-1+\sfrac{1}{p},\sfrac{1}{p})$ and $k\in\llb 0,n\rrb$. For all $u\in\dot{\B}^{s}_{p,q}(\RR^n_+,\Lambda^k)$, there exists $\mathfrak{p}\in\L^p_{\textit{loc}}(\RR^n_+,\Lambda^{k-1})$, unique up to a constant $k-1$-form, such that $\d\mathfrak{p}\in\dot{\B}^{s}_{p,q}(\RR^n_+,\Lambda^k)$, $\nabla\mathfrak{p}\in\dot{\B}^{s}_{p,q}(\RR^n_+,\Lambda^{k-1})^n$, and satisfying
\begin{enumerate}
    \item $\delta(u-\d\mathfrak{p}) = 0$, in $\mathcal{D}'(\RR^n_+,\Lambda^{k-1})$;
     \item $[\mathfrak{e}_n\iprod (u-\d\mathfrak{p})]{}_{|_{\partial\RR^n_+}} = 0$;
\end{enumerate}
with the estimate
\begin{align*}
    \lVert u-\d\mathfrak{p} \rVert_{\dot{\B}^{s}_{p,q}(\RR^n_+)} + \lVert\d\mathfrak{p} \rVert_{\dot{\B}^{s}_{p,q}(\RR^n_+)} \lesssim_{p,s,n} \lVert u\rVert_{\dot{\B}^{s}_{p,q}(\RR^n_+)}.
\end{align*}
Furthermore,
\begin{itemize}
    \item one has the identity
    \begin{align*}
        \E_\mathcal{H}[\d \mathfrak{p}]  = \d(-\Delta)^{-\sfrac{1}{2}}\delta (-\Delta)^{-\sfrac{1}{2}}\E_\mathcal{H}u   = [\I-\delta (-\Delta)^{-\sfrac{1}{2}}\d(-\Delta)^{-\sfrac{1}{2}}]\E_\mathcal{H}u ;
    \end{align*}
    \item a similar result holds for the Besov spaces $\dot{\BesSmo}^{s}_{p,\infty}$, $\dot{\B}^{s,0}_{\infty,q}$, and  $\dot{\BesSmo}^{s,0}_{\infty,\infty}$;

    \item in the case $p\in(1,\infty)$, the result still holds replacing the space $\dot{\B}^{s}_{p,q}$ by either ${\B}^{s}_{p,q}$, $\dot{\H}^{s,p}$, or ${\H}^{s,p}$.
\end{itemize}
\end{theorem}

Now, we provide here  resolvent estimates for the Hodge Laplacian on the flat half-space, including new end-point function spaces and higher order estimates.

\begin{proposition}\label{prop:HodgeResolventPbRn+}Let $p,q\in[1,\infty]$, $-1+{\sfrac{1}{p}}<s<{\sfrac{1}{p}}$. Let $\mu\in[0,\pi)$, $f\in\dot{\B}^{s}_{p,q}(\RR^n_+,\Lambda)$. The resolvent problem
\begin{equation*}\tag{$\mathcal{HL}_{\lambda}$}\label{ResolvHodgeLap}
    \left\{ \begin{array}{rllr}
         \lambda u - \Delta u &= f \text{, }&&\text{ in } \RR^n_+\text{,}\\
        (-\mathfrak{e}_n)\iprod u_{|_{\partial\RR^n_+}} &=0\text{, }&&\text{ on } \partial\RR^n_+\text{,}\\
        (-\mathfrak{e}_n)\iprod \d u_{|_{\partial\RR^n_+}} &=0\text{, }&&\text{ on } \partial\RR^n_+\text{,}\\
    \end{array}
    \right.
\end{equation*}
admits a unique solution $u\in\dot{\B}^{s}_{p,q}\cap \dot{\B}^{s+2}_{p,q}(\RR^n_+,\Lambda)$ which obeys the estimate
\begin{align*}
    |\lambda|\lVert  u\rVert_{\dot{\B}^{s}_{p,q}(\RR^n_+)}+\lVert \nabla^2 u\rVert_{\dot{\B}^{s}_{p,q}(\RR^n_+)} \lesssim_{p,n,s}^{\mu} \lVert  f\rVert_{\dot{\B}^{s}_{p,q}(\RR^n_+)}
\end{align*}
If additionally, for some $r,q\in[1,\infty]$, $\alpha\in(\sfrac{1}{r},1+\sfrac{1}{r})$, one has $f\in\dot{\B}^{\alpha}_{r,\kappa,\mathcal{H}}(\RR^n_+,\Lambda)$, then one also obtains 
$u\in\dot{\B}^{\alpha}_{r,\kappa}\cap \dot{\B}^{\alpha+2}_{r,\kappa}(\RR^n_+,\Lambda)$ with the estimate
\begin{align*}
    |\lambda|\lVert  u\rVert_{\dot{\B}^{\alpha}_{r,\kappa}(\RR^n_+)}+\lVert \nabla^2 u\rVert_{\dot{\B}^{\alpha}_{r,\kappa}(\RR^n_+)}\lesssim_{r,n,\alpha}^{\mu} \lVert  f\rVert_{\dot{\B}^{\alpha}_{r,\kappa}(\RR^n_+)}.
\end{align*}
Furthermore
\begin{itemize}
    \item The result remains valid for $-2+\sfrac{1}{p}<s<-1+\sfrac{1}{p}$: provided $f\in\dot{\B}^{s}_{p,q,\mathcal{H}_\ast}(\RR^n_+,\Lambda)$ one has a unique solution $u$ lying in $\dot{\B}^{s}_{p,q,\mathcal{H}_\ast}\cap\dot{\B}^{s+2}_{p,q,\mathcal{H}}(\RR^n_+,\Lambda)$;
    
    \item a similar result holds for the remaining endpoint Besov spaces $\dot{\BesSmo}^{s}_{p,\infty}$, $\dot{\B}^{s,0}_{\infty,q}$, and  $\dot{\BesSmo}^{s,0}_{\infty,\infty}$;
    
    \item In the case $p=q\in(1,\infty)$,  the result still holds replacing the spaces $\dot{\B}^{\bullet}_{p,q}(\RR^n_+)$ by $\dot{\H}^{\bullet,p}(\RR^n_+)$, and similarly whenever $r=\kappa\in(1,\infty)$;

    \item For any $\theta\in(0,\pi)$, the Hodge Laplacian has a $\mathbf{H}^\infty(\Sigma_\theta)$-functional calculus on any previously mentioned function spaces;

    \item  Everything remains valid for the corresponding inhomogeneous function spaces whenever $1<p<\infty$.
\end{itemize}
\end{proposition}

\begin{proof} Thanks to Propositions~\ref{prop:ExtHodgeRn+} and \ref{prop:ExtOpDiffFormRn+}, one can follow and improve the former results of the second author \cite[Lemma~2.26~\&~Proposition~2.28]{Gaudin2023Hodge}, allowing to include this broader range of exponents ($p=1,\infty$, as well as higher and lower regularity exponent on $-2+\sfrac{1}{p}<s<-1+\sfrac{1}{p}$ and $\sfrac{1}{p}<s<1+\sfrac{1}{p}$), by still setting
\begin{align*}
    u :=[(\lambda \I-\Delta)^{-1}\E_{\mathcal{H}}f]_{|_{\RR^n_+}},
\end{align*}
with the identity
\begin{align*}
    \E_{\mathcal{H}} u = (\lambda \I-\Delta)^{-1}\E_{\mathcal{H}}f.
\end{align*}
The details are left to the reader. Due to the structure of $\E_\mathcal{H}$, one can also take advantage of Propositions~\ref{prop:DirResolventPbRn+}~and~\ref{prop:NeumannResolventPbRn+}.
\end{proof}

\begin{proposition}\label{prop:HodgeResolventPbRn+Linfty}Let $\mu\in[0,\pi)$, $f\in\L^\infty(\RR^n_+,\Lambda)$. The resolvent problem
\begin{equation*}\tag{$\mathcal{HL}_{\lambda}$}\label{ResolvHodgeLapEndpointinfty}
    \left\{ \begin{array}{rllr}
         \lambda u - \Delta u &= f \text{, }&&\text{ in } \RR^n_+\text{,}\\
        (-\mathfrak{e}_n)\iprod u_{|_{\partial\RR^n_+}} &=0\text{, }&&\text{ on } \partial\RR^n_+\text{,}\\
        (-\mathfrak{e}_n)\iprod \d u_{|_{\partial\RR^n_+}} &=0\text{, }&&\text{ on } \partial\RR^n_+\text{,}\\
    \end{array}
    \right.
\end{equation*}
admits a unique solution $u\in\B^{2}_{\infty,\infty}(\RR^n_+,\Lambda)$ which obeys the estimate
\begin{align*}
    |\lambda|\lVert  u\rVert_{\L^\infty(\RR^n_+)}+|\lambda|^\frac{1}{2}\lVert \nabla u\rVert_{\L^\infty(\RR^n_+)}+\lVert  \Delta u\rVert_{\L^\infty(\RR^n_+)}+ \lVert \nabla^2 u\rVert_{\dot{\B}^0_{\infty,\infty}(\RR^n_+)} \lesssim_{n}^{\mu} \lVert  f\rVert_{\L^\infty(\RR^n_+)},
\end{align*}
and for any $c>0$
\begin{align*}
    \left(\frac{1}{1+ |\ln(\min(c,|\lambda|))|}\right)\lVert \nabla^2 u\rVert_{{\B}^0_{\infty,\infty}(\RR^n_+)} \lesssim_{n,c}^{\mu} \lVert  f\rVert_{\L^\infty(\RR^n_+)}.
\end{align*}
If additionally $f\in\W^{1,\infty}_{\mathcal{H}}(\RR^n_+,\Lambda)$, one also obtains $u\in\B^{3}_{\infty,\infty}(\RR^n_+,\Lambda)$ with the estimate
\begin{align*}
    |\lambda|\lVert \nabla u\rVert_{\L^\infty(\RR^n_+)}+|\lambda|^\frac{1}{2}\lVert \nabla^2 u\rVert_{\L^\infty(\RR^n_+)}+\lVert  \nabla \Delta u\rVert_{\L^\infty(\RR^n_+)}+ \lVert \nabla^3 u\rVert_{\dot{\B}^0_{\infty,\infty}(\RR^n_+)} \lesssim_{n}^{\mu} \lVert \nabla f\rVert_{\L^\infty(\RR^n_+)}
\end{align*}
and
\begin{align*}
    \left(\frac{1}{1+ |\ln(\min(c,|\lambda|))|}\right)\lVert \nabla^3 u\rVert_{{\B}^0_{\infty,\infty}(\RR^n_+)} \lesssim_{n,c}^{\mu} \lVert \nabla  f\rVert_{\L^\infty(\RR^n_+)}.
\end{align*}
\end{proposition}

\begin{remark}Note that the same proof as the one exhibited below also yields a more optimal standard/expected result, that is, provided $f\in\L^\infty(\RR^n_+,\Lambda)$, one has $\nabla^2 u \in \mathrm{BMO}(\RR^n_+,\Lambda\otimes\CC^{n^2})$. However, this will not be needed in the current work.
\end{remark}

\begin{proof} \textbf{Step 1:} We start with uniqueness for the weak resolvent problem when $f\in\L^\infty(\RR^n_+,\Lambda)$.

We prove uniqueness of weak solutions $u\in\D_{\infty}(\d,\RR^n_+)\cap\D_{\infty}(\underline{\delta},\RR^n_+)$ to \eqref{ResolvHodgeLapEndpointinfty} with $f=0$, \textit{i.e.}
\begin{align}\label{eq:PROOFWeakSolHodgeLinfty}
    \lambda \langle u,\varphi\rangle_{\RR^n_+} + \langle \d u, \d \varphi\rangle_{\RR^n_+} + \langle \delta u,\delta\varphi\rangle_{\RR^n_+} = 0,\qquad \forall \varphi\in\C_{\tilde{c},\mathcal{H}}^\infty(\overline{\RR^n_+},\Lambda).
\end{align} 
implies $u=0$. By density, one obtains a weak solution if and only if
\begin{align*}
    \lambda \langle u,\varphi\rangle_{\RR^n_+} + \langle \d u, \d \varphi\rangle_{\RR^n_+} + \langle \delta u,\delta\varphi\rangle_{\RR^n_+} = 0,\qquad \forall \varphi\in\W_{\mathcal{H}}^{1,1}({\RR^n_+},\Lambda).
\end{align*} 
So let $u\in \D_{\infty}(\d,\RR^n_+)\cap\D_{\infty}(\underline{\delta},\RR^n_+)$ be such a weak solution. For each $\Psi\in\Ccinfty(\RR^n_+,\Lambda)$, we set
\begin{align*}
    \varphi := [(\lambda \I-\Delta)^{-1}\E_{\mathcal{H}}\Psi]_{|_{\RR^n_+}}
\end{align*}
Therefore, $\varphi\in\H^{2,2}(\RR^n_+,\Lambda)$, with for all $k\in\llb0,n\rrb$,
\begin{itemize}
    \item $\varphi_{I',n}{}_{|_{\partial\RR^n_+}} =0,\qquad \forall I'\in\mathcal{I}^{n-1}_{k-1}$,
    \item $\partial_{x_n}\varphi_{I}{}_{|_{\partial\RR^n_+}} =0,\qquad \forall I\in\mathcal{I}^{n-1}_{k}$,
\end{itemize}
so that
\begin{align*}
    (-\mathfrak{e}_n)\iprod\varphi_{|_{\partial\RR^n_+}} =0,\qquad \text{ and }\qquad(-\mathfrak{e}_n)\iprod \d\varphi_{|_{\partial\RR^n_+}} =0,
\end{align*}
and $\varphi$ satisfies the equation
\begin{align*}
    \lambda\varphi - \Delta \varphi = \Psi,\qquad \text{ in }\L^2(\RR^n_+,\Lambda).
\end{align*}
Furthermore due to the definition of function spaces by restriction one also obtains $\varphi\in\W^{1,1}(\RR^n_+,\Lambda)$, with an estimate
\begin{align*}
    |\lambda|\lVert  \varphi\rVert_{\L^1(\RR^n_+)}+&|\lambda|^\frac{1}{2}\lVert  \nabla \varphi\rVert_{\L^1(\RR^n_+)} \\
    &\leqslant |\lambda|\lVert  (\lambda \I-\Delta)^{-1}\E_{\mathcal{H}}\Psi\rVert_{\L^1(\RR^n)}+|\lambda|^\frac{1}{2}\lVert  \nabla (\lambda \I-\Delta)^{-1}\E_{\mathcal{H}}\Psi\rVert_{\L^1(\RR^n)}\\
    &\lesssim_{\mu,n} \lVert  \E_{\mathcal{H}}\Psi\rVert_{\L^1(\RR^n)}\\
    &\lesssim_{\mu,n} \lVert  \Psi\rVert_{\L^1(\RR^n_+)}.
\end{align*}
Additionally, $\varphi\in\W^{2,1}(\RR^n_+,\Lambda)$, since
\begin{align*}
    |\lambda|^\frac{1}{2}\lVert  \nabla^2 \varphi\rVert_{\L^1(\RR^n_+)} &\leqslant |\lambda|^\frac{1}{2}\lVert  \nabla^2 (\lambda \I-\Delta)^{-1}\E_{\mathcal{H}}\Psi\rVert_{\L^1(\RR^n)}\\
    &\leqslant |\lambda|^\frac{1}{2}\lVert  \nabla (\lambda \I-\Delta)^{-1}\E_{\mathcal{H}}\nabla'\Psi\rVert_{\L^1(\RR^n)} +|\lambda|^\frac{1}{2}\lVert  \nabla (\lambda \I-\Delta)^{-1}\E_{\mathcal{H}_\ast}\partial_{x_n}\Psi\rVert_{\L^1(\RR^n)} \\
    &\lesssim_{\mu,n} \lVert \E_{\mathcal{H}}\nabla'\Psi\rVert_{\L^1(\RR^n)} +\lVert  \E_{\mathcal{H}_\ast}\partial_{x_n}\Psi\rVert_{\L^1(\RR^n)} \\
    &\lesssim_{\mu,n} \lVert \nabla\Psi\rVert_{\L^1(\RR^n_+)}.
\end{align*}
In particular, due to $\varphi\in\W^{2,1}(\RR^n_+,\Lambda)$ with the appropriate boundary conditions, we did obtain $\varphi\in\D_{1}(\d,\RR^n_+)\cap\D_{1}(\underline{\delta},\RR^n_+)$ with $\d \varphi\in\D_{1}(\underline{\delta},\RR^n_+)$ and $\underline{\delta} \varphi\in\D_{1}(\d,\RR^n_+)$. This legitimates the following integration by parts
\begin{align*}
    \langle u, \Psi\rangle_{\RR^n_+} =  \langle u, \lambda\varphi - \Delta \varphi\rangle_{\RR^n_+} = \lambda\langle u, \varphi\rangle_{\RR^n_+} + \langle \d u, \d\varphi\rangle_{\RR^n_+}+ \langle \delta u, \delta\varphi\rangle_{\RR^n_+} =0.
\end{align*}
Thus, $\langle u, \Psi\rangle_{\RR^n_+}=0$ for all $\Psi\in\Ccinfty(\RR^n_+,\Lambda)$, which implies $u=0$.

\textbf{Step 2:} We prove existence. So for $f\in\L^\infty(\RR^n_+,\Lambda)$, the solution is given by
\begin{align*}
    u:=[(\lambda \I-\Delta)^{-1}\E_{\mathcal{H}}f]_{|_{\RR^n_+}}.
\end{align*}
The estimates admit a proof similar to the one proving that the test function is $\W^{1,1}$, except one has to prove that $\nabla^2 u\in\dot{\B}^{0}_{\infty,\infty}(\RR^n_+,\Lambda\otimes\CC^{n^2})$.

Note that, $(\lambda \I-\Delta)^{-1}\E_{\mathcal{H}}f$ belongs to $\L^\infty(\RR^n,\Lambda)$. Therefore $\nabla^2 (\lambda \I-\Delta)^{-1}\E_{\mathcal{H}}f\in\dot{\B}^{-2}_{\infty,\infty}(\RR^n,\Lambda\otimes\CC^{n^2})\subset\S'_h(\RR^n,\Lambda\otimes\CC^{n^2})$, and
\begin{align*}
    \lVert \nabla^2 (\lambda \I-\Delta)^{-1}\E_{\mathcal{H}}f \rVert_{\dot{\B}^{0}_{\infty,\infty}(\RR^n)} \lesssim_{\mu,n} \lVert \E_{\mathcal{H}}f \rVert_{{\L}^{\infty}(\RR^n)}\lesssim_{\mu,n} \lVert f \rVert_{{\L}^{\infty}(\RR^n_+)}
\end{align*}
Since $\nabla^2 (\lambda \I-\Delta)^{-1}\E_{\mathcal{H}}f\in\S'_h(\RR^n,\Lambda\otimes\CC^{n^2})$ and $\lVert \nabla^2 (\lambda \I-\Delta)^{-1}\E_{\mathcal{H}}f \rVert_{\dot{\B}^{0}_{\infty,\infty}(\RR^n)}<\infty$, it holds that 
\begin{align*}
    \nabla^2 (\lambda \I-\Delta)^{-1}\E_{\mathcal{H}}f \in \dot{\B}^{0}_{\infty,\infty}(\RR^n,\Lambda\otimes\CC^{n^2}).
\end{align*}
Consequently by the definition of function spaces by restriction
\begin{align*}
    \lVert \nabla^2 u\rVert_{\dot{\B}^{0}_{\infty,\infty}(\RR^n_+)}\leqslant \lVert \nabla^2 (\lambda \I-\Delta)^{-1}\E_{\mathcal{H}}f\rVert_{\dot{\B}^{0}_{\infty,\infty}(\RR^n)}\lesssim_{\mu,n} \lVert f \rVert_{{\L}^{\infty}(\RR^n_+)}.
\end{align*}

\textbf{Step 3:} We prove the logarithmic loss-resolvent estimate when the solution is measured with an inhomogeneous Besov norm for the second order derivatives. This is achieved by a standard but non trivial dilation argument.

First, provided $f\in\L^\infty(\RR^n_+,\Lambda)$, note that $u$ is a solution to the problem \eqref{ResolvHodgeLapEndpointinfty} if and only if $u_{\lambda}= |\lambda|u(\cdot/|\lambda|^\frac{1}{2})$ is a solution to \eqref{ResolvHodgeLapEndpointinfty} with forcing term $f_\lambda=(\cdot/|\lambda|^\frac{1}{2})$, so
\begin{align*}
    \lVert \nabla^2 u_\lambda \rVert_{{\B}^{0}_{\infty,\infty}(\RR^n_+)} &\lesssim_n  \lVert  u_\lambda \rVert_{{\B}^{2}_{\infty,\infty}(\RR^n_+)} \\ 
    &\lesssim_{n} \lVert u_\lambda \rVert_{\L^\infty(\RR^n_+)} + \lVert \nabla^2 u_\lambda\rVert_{\dot{\B}^{0}_{\infty,\infty}(\RR^n_+)}\\
    &\lesssim_{n,\mu} \lVert f_\lambda \rVert_{\L^\infty(\RR^n_+)} = \lVert f\rVert_{\L^\infty(\RR^n_+)}.
\end{align*}
However, one has $\nabla^2 u_\lambda = (\nabla^2 u)(\cdot/|\lambda|^\frac{1}{2})$,  if $|\lambda|<1$, one obtains $1/|\lambda|^\frac{1}{2}>1$, and  by \cite[Theorem~3.11]{Vybiral2008}:
\begin{align*}
    \lVert \nabla^2 u \rVert_{{\B}^{0}_{\infty,\infty}(\RR^n_+)} &= \lVert (\nabla^2 u_\lambda)_{\lambda^{-1}}\rVert_{{\B}^{0}_{\infty,\infty}(\RR^n_+)}\\
    &\lesssim_{n} (1+|\ln(|\lambda|)|) \lVert \nabla^2 u_\lambda \rVert_{{\B}^{0}_{\infty,\infty}(\RR^n_+)}\\
    &\lesssim_{n,\mu} (1+|\ln(|\lambda|)|) \lVert f\rVert_{\L^\infty(\RR^n_+)}.
\end{align*}
Now if $|\lambda|\geqslant 1$, one obtains
\begin{align*}
    \lVert \nabla^2 u \rVert_{{\B}^{0}_{\infty,\infty}(\RR^n_+)} &\lesssim_{n} \lVert  u \rVert_{{\B}^{2}_{\infty,\infty}(\RR^n_+)}\\
    &\sim_{n} \lVert u\rVert_{{\L}^{\infty}(\RR^n_+)} + \lVert  \nabla^2 u \rVert_{\dot{\B}^{0}_{\infty,\infty}(\RR^n_+)}\\
    &\lesssim_{n,\mu} \left(\frac{1}{|\lambda|}+ 1\right) \lVert  f \rVert_{{\L}^\infty(\RR^n_+)}\\
    &\lesssim_{n,\mu} 2 \lVert  f \rVert_{{\L}^\infty(\RR^n_+)}.
\end{align*}
Summarizing both cases, we did obtain
\begin{align*}
    \left(\frac{1}{1+ |\ln(\min(1,|\lambda|)|)}\right)\lVert \nabla^2 u\rVert_{{\B}^0_{\infty,\infty}(\RR^n_+)} \lesssim_{n}^{\mu} \lVert  f\rVert_{\L^\infty(\RR^n_+)}.
\end{align*}

Concerning higher order estimates when $f\in\W^{1,\infty}_\mathcal{H}(\RR^n_+,\Lambda)$, just use the component-wise formulas  $\nabla'\E_{\mathcal{H}}f = \E_{\mathcal{H}}\nabla'f$ and $\partial_{x_n}\E_{\mathcal{H}}f = \E_{\mathcal{H}_\ast}\partial_{x_n}f$.
\end{proof}

\begin{corollary}\label{cor:HodgeResolventPbRn+LinftyBesov}Let $s\in(0,1)$ and $q\in[1,\infty]$. Let $\mu\in[0,\pi)$, $f\in\B^{s}_{\infty,q,\mathcal{H}}(\RR^n_+,\Lambda)$. For any $\lambda\in\Sigma_\mu$, the resolvent problem
\begin{equation*}\tag{$\mathcal{HL}_{\lambda}$}\label{ResolvHodgeLapEndpointinftyBesov}
    \left\{ \begin{array}{rllr}
         \lambda u - \Delta u &= f \text{, }&&\text{ in } \RR^n_+\text{,}\\
        (-\mathfrak{e}_n)\iprod u_{|_{\partial\RR^n_+}} &=0\text{, }&&\text{ on } \partial\RR^n_+\text{,}\\
        (-\mathfrak{e}_n)\iprod \d u_{|_{\partial\RR^n_+}} &=0\text{, }&&\text{ on } \partial\RR^n_+\text{,}\\
    \end{array}
    \right.
\end{equation*}
admits a unique solution $u\in\B^{s+2}_{\infty,q}(\RR^n_+,\Lambda)$ which obeys the estimate
\begin{align*}
    |\lambda|\lVert  u\rVert_{\B^{s}_{\infty,q}(\RR^n_+)}+\lVert  \Delta u\rVert_{\B^{s}_{\infty,q}(\RR^n_+)}+ \lVert \nabla^2 u\rVert_{\dot{\B}^0_{\infty,\infty}(\RR^n_+)} + \lVert \nabla^2 u\rVert_{\dot{\B}^s_{\infty,q}(\RR^n_+)} \lesssim_{n}^{\mu} \lVert  f\rVert_{\B^{s}_{\infty,q}(\RR^n_+)}
\end{align*}
and for any $c>0$
\begin{align*}
    \left(\frac{1}{1+ |\ln(\min(c,|\lambda|))}\right)\lVert \nabla^2 u\rVert_{{\B}^s_{\infty,q}(\RR^n_+)} \lesssim_{n}^{\mu} \lVert  f\rVert_{\B^{s}_{\infty,q}(\RR^n_+)}.
\end{align*}
\end{corollary}

\begin{proof}Just apply real interpolation to Proposition~\ref{prop:HodgeResolventPbRn+Linfty}.
\end{proof}

\subsubsection{The Leray projection on the half-space and bounded domains. Duality.}

\medbreak

 The next result is just a reformulation of Theorem~\ref{thm:HodgeDecompRn+1} for the Hodge decompositions of homogeneous function spaces on the half-space.
\begin{theorem}\label{thm:HodgeDecompRn+2}Let $p,q\in[1,\infty]$, $s\in(-1+\sfrac{1}{p},\sfrac{1}{p})$ and $k\in\llb 0,n\rrb$. The following assertions hold:
\begin{enumerate}
    \item The generalized Hodge--Leray projection $\PP_{\RR^n_+}$ is well-defined and bounded seen as an operator on $\dot{\B}^{s}_{p,q}(\RR^n_+,\Lambda^k)$, and we do have the identity
    \begin{align*}
        \E_\mathcal{H}\PP_{\RR^n_+} = \PP_{\RR^n}\E_\mathcal{H} ;
    \end{align*}
    \item The following Hodge decomposition holds:
    \begin{align*}
        \dot{\B}^{s}_{p,q}(\RR^n_+,\Lambda^k) = \dot{\B}^{s,\sigma}_{p,q,\mathfrak{n}}(\RR^n_+,\Lambda^k)\oplus\dot{\B}^{s,\gamma}_{p,q}(\RR^n_+,\Lambda^k);
    \end{align*}
\end{enumerate}
Furthermore,
\begin{itemize}
    \item a similar result holds for the Besov spaces $\dot{\BesSmo}^{s}_{p,\infty}$, $\dot{\B}^{s,0}_{\infty,q}$, and  $\dot{\BesSmo}^{s,0}_{\infty,\infty}$;

    \item in the case $p\in(1,\infty)$, the result still holds replacing the space $\dot{\B}^{s}_{p,q}$ by either ${\B}^{s}_{p,q}$, $\dot{\H}^{s,p}$, or ${\H}^{s,p}$;

    \item a similar result holds if one replaces $(\sigma,\gamma,\mathfrak{n},\PP,\mathcal{H})$ by $(\gamma,\sigma,\mathfrak{t},\QQ,\mathcal{H}_\ast)$.
\end{itemize}
\end{theorem}

We want to perform a counterpart of the Theorem above for $\C^{1,\alpha}$-bounded  and similar but rougher domains and which also includes the endpoint $p=1,\infty$ when it comes to Besov spaces.

\medbreak

As briefly exposed before, the behavior of the Hodge decomposition mainly ties to the behavior of the Hodge-Laplace operator. Due to several technical obstructions, we rather prefer to concentrate on the behavior of the Hodge-Dirac operator and propose here a first order approach to the Hodge decomposition. This idea is inpired by  the work of M${}^{\text{c}}$Intosh and Monniaux \cite{McintoshMonniaux2018}.

\medbreak

We start with a result on the half-space, aiming to reach the corresponding one on bounded domains by localisation techniques and Sobolev-Besov multiplier theory.
\begin{proposition}\label{prop:HodgeDiracRn+} Let $\mu\in(0,\frac{\pi}{2})$. Let $p,q\in[1,\infty]$, $s\in(-1+\sfrac{1}{p},\sfrac{1}{p})$. For all $\lambda\in \mathrm{S}_\mu$, all $f\in\dot{\B}^{s}_{p,q}(\RR^n_+,\Lambda)$, the resolvent problem
\begin{equation*}\tag{$\mathcal{HD}_{\lambda}$}\label{ResolvHodgeDirac}
    \left\{ \begin{array}{rllr}
         \lambda u + i(\d u + \delta u) &= f \text{, }&&\text{ in } \RR^n_+\text{,}\\
        (-\mathfrak{e}_n)\iprod u_{|_{\partial\RR^n_+}} &=0\text{, }&&\text{ on } \partial\RR^n_+\text{,}
    \end{array}
    \right.
\end{equation*}
admits a unique solution $u\in\dot{\B}^{s}_{p,q}\cap\dot{\B}^{s+1}_{p,q}(\RR^n_+,\Lambda)$, with an estimate
\begin{align*}
    |\lambda|\lVert u\rVert_{\dot{\B}^{s}_{p,q}(\RR^n_+)} + \lVert \nabla u\rVert_{\dot{\B}^{s}_{p,q}(\RR^n_+)} \lesssim_{p,s,n}^\mu \lVert f\rVert_{\dot{\B}^{s}_{p,q}(\RR^n_+)}.
\end{align*}
Furthermore
\begin{itemize}
    \item a similar result holds for the remaining endpoint Besov spaces $\dot{\BesSmo}^{s}_{p,\infty}$, $\dot{\B}^{s,0}_{\infty,q}$, ,$\dot{\BesSmo}^{s,0}_{\infty,\infty}$ as well as for Sobolev spaces $\dot{\H}^{s,p}$, assuming $1<p<\infty$ for the latter;

    \item a similar result holds with the boundary condition $(-\mathfrak{e}_n)\wedge u_{|_{\partial\RR^n_+}}$ instead of $(-\mathfrak{e}_n)\iprod u_{|_{\partial\RR^n_+}}$.
\end{itemize}
\end{proposition}

\begin{proof}Following the proof of Proposition~\ref{prop:HodgeResolventPbRn+}, since for $D=\d+\delta$, one has
\begin{align*}
    D^2 = (\d + \delta )^2 = -\Delta,
\end{align*}
we set
\begin{align*}
    u:= [(\lambda - i D)(\lambda^2-\Delta)^{-1}\E_\mathcal{H}f]_{|_{\RR^n_+}}.
\end{align*}
Therefore by the definition of function spaces by restriction and Proposition~\ref{prop:ExtHodgeRn+}
\begin{align*}
    \lVert u\rVert_{\dot{\B}^{s}_{p,q}(\RR^n_+)} &\leqslant \lVert (\lambda - i D)(\lambda^2-\Delta)^{-1}\E_\mathcal{H}f\rVert_{\dot{\B}^{s}_{p,q}(\RR^n)}\\
    &\leqslant \lVert \lambda(\lambda^2-\Delta)^{-1}\E_\mathcal{H}f\rVert_{\dot{\B}^{s}_{p,q}(\RR^n)} + \lVert \nabla(\lambda^2-\Delta)^{-1}\E_\mathcal{H}f\rVert_{\dot{\B}^{s}_{p,q}(\RR^n)}\\
    &\lesssim_{p,s,n}^\mu \frac{1}{|\lambda|}\lVert \E_\mathcal{H}f\rVert_{\dot{\B}^{s}_{p,q}(\RR^n)}\\
    &\lesssim_{p,s,n}^\mu \frac{1}{|\lambda|}\lVert f\rVert_{\dot{\B}^{s}_{p,q}(\RR^n_+)}.
\end{align*}
Similarly
\begin{align*}
    \lVert \nabla u\rVert_{\dot{\B}^{s}_{p,q}(\RR^n_+)} &\leqslant \lVert \nabla (\lambda - i D)(\lambda^2-\Delta)^{-1}\E_\mathcal{H}f\rVert_{\dot{\B}^{s}_{p,q}(\RR^n)}\\
    &\leqslant \lVert \nabla \lambda(\lambda^2-\Delta)^{-1}\E_\mathcal{H}f\rVert_{\dot{\B}^{s}_{p,q}(\RR^n)} + \lVert \nabla^2(\lambda^2-\Delta)^{-1}\E_\mathcal{H}f\rVert_{\dot{\B}^{s}_{p,q}(\RR^n)}\\
    &\lesssim_{p,s,n}^\mu\lVert f\rVert_{\dot{\B}^{s}_{p,q}(\RR^n_+)}.
\end{align*}
The same goes for Sobolev spaces. Uniqueness can be shown by an argument similar to the one in the proof of Proposition~\ref{prop:HodgeResolventPbRn+Linfty}, in Step 1, due to self-adjointness of the problem.
\end{proof}

\begin{proposition}\label{prop:HodgeDiracbdd} Let $\alpha\in(0,1)$, $r\in[1,\infty]$, and let $\Omega$ be a bounded Lipschitz domain of class $\mathcal{M}^{1+\alpha,r}_\W(\epsilon)$ for some sufficiently small $\epsilon>0$. Let $\mu\in(0,\frac{\pi}{2})$. Let $p,q\in[1,\infty]$, $s\in(-1+\sfrac{1}{p},\sfrac{1}{p})$. For all $\lambda\in \mathrm{S}_\mu$, all $f\in{\B}^{s}_{p,q}(\Omega,\Lambda)$, the resolvent problem
\begin{equation*}\tag{$\mathcal{HD}_{\lambda}$}\label{ResolvHodgeDiracDomain}
    \left\{ \begin{array}{rllr}
         \lambda u + i(\d u + \delta u) &= f \text{, }&&\text{ in } \Omega\text{,}\\
        \nu\iprod u_{|_{\partial\Omega}} &=0\text{, }&&\text{ on } \partial\Omega\text{,}
    \end{array}
    \right.
\end{equation*}
admits a unique solution $u\in{\D}^{s}_{p,q}(\d)\cap{\D}^{s}_{p,q}(\delta,\Omega,\Lambda)$, with an estimate
\begin{align*}
    |\lambda|\lVert u\rVert_{{\B}^{s}_{p,q}(\Omega)} + \lVert (\d u,\delta u)\rVert_{{\B}^{s}_{p,q}(\Omega)} \lesssim_{p,s,n}^{\mu,\Omega} \lVert f\rVert_{{\B}^{s}_{p,q}(\Omega)}.
\end{align*}
If additionally, either
\begin{enumerate}
    \item $p\in(r,\infty]$, $s\in(-1+\frac{1}{p},-1+\frac{1+r\alpha}{p})$; or
    \item $p\in[1,r]$, $s\in(-1+\frac{1}{p},-1+\alpha+\frac{1}{p})$; or
    \item $p=q=r$, $s=-1+\alpha+\frac{1}{p}$;
\end{enumerate}
then $u\in\B^{s+1}_{p,q}(\Omega)$ and
\begin{align*}
    \lVert  u\rVert_{{\B}^{s+1}_{p,q}(\Omega)} \lesssim_{p,s,n}^\Omega \lVert( u,\,\d u,\,\delta u)\rVert_{{\B}^{s}_{p,q}(\Omega)}.
\end{align*}
Furthermore
\begin{itemize}
    \item a similar result holds for the remaining endpoint Besov spaces ${\BesSmo}^{s}_{p,\infty}$, as well as for Sobolev spaces ${\H}^{s,p}$, assuming $1<p<\infty$ for the latter;

    \item a similar result holds with the tangential boundary condition $\nu\wedge u_{|_{\partial\Omega}}$ instead of the normal one $\nu\iprod u_{|_{\partial\Omega}}$.
\end{itemize}
\end{proposition}

\begin{proof} We deal with the case $r=1$, the proof for $r\in(1,\infty)$ is similar otherwise, and the case $r=\infty$ is simpler. We proceed by standard localisation techniques and decompose the proof in several parts: first, the case of the spaces $\H^{s,2}$, $-\frac{1}{2}<s<\frac{1}{2}$, then the case of the Besov and Sobolev spaces $\B^{s}_{p,q}(\Omega)$, $\H^{s,p}(\Omega)$, $1<p<\infty$, $q\in[1,\infty]$, $s\in(-1+\sfrac{1}{p},\sfrac{1}{p})$, and finally the cases $p=1,\infty$ for Besov spaces. The procedure will be iterative and relying on the fact that at each step the extrapolation of the Hodge-Dirac operator has a nullspace independent of the considered parameters, which we shall prove.

By assumption, there is $\ell\in\mathbb N$ and functions $\varphi_1,\dots,\varphi_\ell\in\mathcal{M}^{1+\alpha,1}_{\W}(\epsilon)$ satisfying \ref{A1}--\ref{A3} for a sufficiently small $\epsilon>0$.

We consider an open set $\mathcal U^0\Subset{\Omega}$ such that ${\Omega}\subset \cup_{j=0}^\ell \mathcal U^j$. Finally, we consider a decomposition of unity $(\eta_j)_{j\in\llb0,\ell\rrb}$ with respect to the covering
$\mathcal U^0,\dots,\mathcal U^\ell$ of ${\Omega}$. 
For $j\in\llb 1,\ell\rrb$, we consider the extension $\boldsymbol{\Phi}_j$ of $\varphi_j$ given by \eqref{eq:Phi} with inverse $\boldsymbol{\Phi}_j^{-1}$. So that $\nabla\boldsymbol{\Phi}_j\in\mathcal{M}^{\beta,p}_{\X,\tt{or}}(\epsilon)$  for all $-1+\frac{1}{p}<\beta<\frac{1+\alpha}{p}$, all $p\in[1,\infty)$, thanks to Remark~\ref{rem:BetterMultipliers}. That is $\nabla' \varphi_\ell$ has small multiplier norm.

Let $f\in\Ccinfty(\Omega,\Lambda)$ and consider the unique solution $u\in \D_2(\d,\Omega,\Lambda)\cap\D_2(\d^\ast,\Omega,\Lambda)$ provided by Proposition~\ref{prop:HodgeDiracL2bdd}.
\begin{equation}\label{eq:HodgeDirOr}
    \left\{ \begin{array}{rllr}
         \lambda u-i(\d+\delta)u &= f \text{, }&&\text{ in } \Omega\text{,}\\
        \nu\iprod{u}{}_{|_{\partial\Omega}} &=0\text{, } &&\text{ on } \partial\Omega\text{.}
    \end{array}
    \right.
\end{equation}

\textbf{Step 1:} \textbf{Regularization of the problem and Preparatory setting.} This preparatory steps aims to help of proving that solution to a family of regularized problems will converge to the desired solution and to rewrite the system as a suitable perturbation of its counterpart on the half-space. By Theorem~\ref{thm:regularizeddomain}, consider a sequence of approximating domains $(\Omega^{\varepsilon})_{\varepsilon>0}$ associated with the regularized charts $((\varphi_{j}^{\varepsilon})_{j\in\llb1,\ell\rrb})_{\varepsilon>0}$ obtained by convolution according to Theorem~\ref{thm:regularizeddomain}, and the associated diffeomorphism. Consider $u^{\varepsilon}$ be the solution to
\begin{equation}\label{eq:HodgeDirOrregularized}
    \left\{ \begin{array}{rllr}
         \lambda u^{\varepsilon}+i (\d+\delta)u^{\varepsilon} &= f \text{, }&&\text{ in } \Omega_\varepsilon\text{,}\\
        \nu_{\varepsilon}\iprod{u^{\varepsilon}}{}_{|_{\partial\Omega_\varepsilon}} &=0\text{, } &&\text{ on } \partial\Omega_\varepsilon\text{.}
    \end{array}
    \right.
\end{equation}

By Theorem~\ref{thm:regularizeddomain}, and since the characteristic constants of the domain are preserved through the approximation procedure for the boundary Theorem~\ref{thm:regularizeddomain}, up to consider the restriction of $u^{\varepsilon}$ to $\Omega$, one obtains $u^{\varepsilon}\in\L^2(\Omega,\Lambda)$, and according to Proposition~\ref{prop:HodgeDiracL2bdd}
\begin{align*}
    |\lambda|\lVert u^{\varepsilon} \rVert_{\L^{2}(\Omega)} + \lVert (\d u^{\varepsilon}, \delta u^\varepsilon )\rVert_{\L^{2}(\Omega)} \leqslant \frac{2}{\cos(\mu)}\lVert f \rVert_{\L^{2}(\Omega)},
\end{align*}
uniformly with respect to $\varepsilon>0$. By Banach-Alaoglu, up to extract many finitely subsequences, $(u^{\varepsilon})_\varepsilon$ converges weakly to some limit $v\in\L^{2}(\Omega,\Lambda)$  and similarly with $(\d u^{\varepsilon}, \delta u^\varepsilon )$ is converging weakly towards some couple $(v_\sigma,v_\gamma)\in\L^{2}(\Omega,\Lambda)^2$ and by the compact embedding $\L^2\hookrightarrow\H^{s,2}$, $-\frac{1}{2}<s<0$, it converges strongly in $\H^{s,2}(\Omega,\Lambda)$,  and by closedness of $\d$ and $\delta$ in $\H^{s,2}(\Omega,\Lambda)$, one obtains
\begin{align*}
    (u^{\varepsilon}, \d u^{\varepsilon}, \delta u^\varepsilon )\xrightarrow[\varepsilon\rightarrow 0]{} (v,v_\sigma,v_\gamma)= (v, \d v, \delta v). 
\end{align*}
Hence, $(u^{\varepsilon})_\varepsilon$ converges weakly towards $v\in\D_2(\d,\Omega,\Lambda)\cap\D_2(\delta,\Omega,\Lambda)$. It holds 
\begin{align*}
    \lambda v + i(\d + \delta) v =f,\quad\text{ in }\quad\mathcal{D}'(\Omega,\Lambda).
\end{align*}

By \cite[Theorem~11.2,~eq.~(11.13)]{MitreaMitreaTaylor2001} and the formula $\nu\wedge(\nu\iprod\cdot)+\nu\iprod(\nu\wedge\cdot) = |\nu|^2\I=\I$, it holds that for any $\varepsilon$, $u^{\varepsilon}_{|_{\partial\Omega_{\varepsilon}}}\in\L^2(\partial\Omega_\varepsilon,\Lambda)$, with a uniform estimate, thanks to Theorem~\ref{thm:regularizeddomain},
\begin{align*}
    \lVert u^{\varepsilon}_{|_{\partial\Omega_{\varepsilon}}}\rVert_{\L^2(\partial\Omega_\varepsilon)}\lesssim_{\Omega} \lVert (u^{\varepsilon},\,\d u^\varepsilon,\,\delta u^\varepsilon)\rVert_{\L^2(\Omega_\varepsilon)}.
\end{align*}
Up to localisation and to consider $(\eta_j u_\varepsilon) \circ \boldsymbol{\Phi}_{j}^\varepsilon$, $j\in\llb 1,\ell\rrb$, due to Theorem~\ref{thm:regularizeddomain} and the definition of $\boldsymbol{\Phi}_{j}^\varepsilon$, one obtains the weak convergence of $\Big([(\eta_j u_\varepsilon) \circ \boldsymbol{\Phi}_{j}^\varepsilon]_{|_{\partial\RR^n_+}}\Big)_\varepsilon$ towards $[(\eta_j v) \circ \boldsymbol{\Phi}_{j}]_{|_{\partial\RR^n_+}}$ in $\L^2(\partial\RR^n_+,\Lambda)$, and consequently
\begin{align*}
    \nu\iprod v_{|_{\partial\Omega}} =0.
\end{align*}
Therefore, along with $v\in\D_2(\d,\Omega,\Lambda)\cap\D_2(\delta,\Omega,\Lambda)$ and  $\lambda v + i(\d + \delta) v =f$, we deduce  that necessarily $u=v$.
Thus, we did obtain
\begin{align*}
    (u^{\varepsilon}, \d u^{\varepsilon}, \delta u^\varepsilon )\xrightarrow[\varepsilon\rightarrow 0]{} (u,\d u,\delta u)
\end{align*}
weakly in $\L^2$ and strongly in $\H^{s,2}$, for all $-\frac{1}{2}<s<0$.

 Let us fix $j\in\llb1,\ell\rrb$ and assume, without loss of generality, that the reference point $y_j=0$ and that the outer normal at~$0$ is pointing in the negative $x_n$-direction (this saves us some notation regarding the translation and rotation of the coordinate system).
We multiply $u^\varepsilon$ by $\eta_j$ and take the (inverse) push forward in order to obtain for $\tilde{u}_j^{\varepsilon}:=(\tilde{\boldsymbol{\Phi}}_{\ast j}^{\varepsilon})^{-1}u_j^{\varepsilon}=(\tilde{\boldsymbol{\Phi}}_{\ast j}^{\varepsilon})^{-1}(\eta_j u^{\varepsilon})$, and $f_j=\eta_jf$, and $\tilde{f}_j^\varepsilon:=(\tilde{\boldsymbol{\Phi}}_{\ast j}^{\varepsilon})^{-1} f_j$ the equations
\begin{equation}\label{eq:HodgeDir2}
    \left\{ \begin{array}{rllr}
         \lambda\tilde{u}_j^{\varepsilon}+ i((\tilde{\boldsymbol{\Phi}}_{\ast j}^{\varepsilon})^{-1}\d {u}_j^{\varepsilon}+ \delta \tilde{u}_j^\varepsilon)  &= i\,[(\tilde{\boldsymbol{\Phi}}_{\ast j}^{\varepsilon})^{-1}[\d+\delta,\eta_j]u^{\varepsilon}]+\tilde{f}_j^\varepsilon \text{, }&&\text{ in } \RR^n_+\text{,}\\
        (-\mathfrak{e}_n)\iprod{\tilde{u}_j^{\varepsilon}}{}_{|_{\partial\RR^n_+}} &=0\text{, } &&\text{ on } \partial\RR^n_+\text{.}
    \end{array}
    \right.
\end{equation}
with the commutators $[\d+\delta,\eta_j]= -\nabla\eta_j\wedge(\cdot) + \nabla \eta_j\iprod(\cdot)$, see \eqref{eq:product-rule}. We write $\tilde{g}_j^\varepsilon = [(\tilde{\boldsymbol{\Phi}}_{j,\ast}^{\varepsilon})^{-1}[\d+\delta,\eta_j]u^{\varepsilon}]$, and we also write
\begin{align*}
    \tilde{\boldsymbol{\Phi}}_{\ast j}^{-1}\d {u}_j = \tilde{\boldsymbol{\Phi}}_{\ast j}^{-1}({\boldsymbol{\Phi}}_{j}^\ast)^{-1}\d {\boldsymbol{\Phi}}_{j}^{\ast}(\tilde{\boldsymbol{\Phi}}_{\ast j}) \tilde{u}_j
\end{align*}
We write $\mathbf{A}_j:=\tilde{\boldsymbol{\Phi}}_{\ast j}^{-1} (\boldsymbol{\Phi}_{j}^{\ast})^{-1}=\det(\nabla\boldsymbol{\Phi}_j)(\nabla\boldsymbol{\Phi}_j)^{-1}\cdot\prescript{t}{}{(\nabla\boldsymbol{\Phi}_j)}^{-1}$ so that it also holds $\mathbf{A}_j^{-1}={\boldsymbol{\Phi}}_{j}^{\ast}(\tilde{\boldsymbol{\Phi}}_{\ast j})=\det(\nabla\boldsymbol{\Phi}_j)^{-1}\cdot\prescript{t}{}{(\nabla\boldsymbol{\Phi}_j)}\cdot\nabla\boldsymbol{\Phi}_j^{-1}$ where the action of matrices is given at the differential form level according to Definition~\ref{def:changevar}. In this case, putting $\tilde{g}_j^\varepsilon:= i(\tilde{\boldsymbol{\Phi}}_{\ast j}^{\varepsilon})^{-1}[-\nabla\eta_j\wedge p  + \nabla \eta_j\iprod p]$, \eqref{eq:HodgeDir2} becomes
\begin{equation}\label{eq:HodgeDir3}
    \left\{ \begin{array}{rllr}
         \lambda\tilde{u}_j^{\varepsilon}+ i(\d \tilde{u}_j^\varepsilon+ \delta \tilde{u}_j^\varepsilon)  &= i[\d - \mathbf{A}_j^\varepsilon \d {(\mathbf{A}_j^\varepsilon)}^{-1}]\tilde{u}_j^\varepsilon+ \tilde{f}_j^\varepsilon + \tilde{g}_j^\varepsilon\text{, }&&\text{ in } \RR^n_+\text{,}\\
        (-\mathfrak{e}_n)\iprod{\tilde{u}_j^{\varepsilon}}{}_{|_{\partial\RR^n_+}} &=0\text{, } &&\text{ on } \partial\RR^n_+\text{.}
    \end{array}
    \right.
\end{equation}

\textbf{Step 2:}\textbf{ We prove regularity of the problem on $\H^{s,2}$}, first dealing partially with the case $-\frac{1}{2}<s<-\frac{1}{2} + \frac{\alpha}{2}$, then reaching any $-\frac{1}{2}<s<\frac{1}{2}$. 

\textbf{Step 2.1:} \textbf{We focus on proving the solution $u$ is sufficiently regular when $-\frac{1}{2}<s<-\frac{1}{2} + \frac{\alpha}{2}$}, and provide a partial resolvent estimate. Note that since $f$ is smooth, by construction it holds that $u^{\varepsilon}$ is smooth in ${\Omega}_{\varepsilon}$ up to the boundary, and applying Proposition~\ref{prop:HodgeDiracRn+} to \eqref{eq:HodgeDir3}, it holds
\begin{align}\label{eq:Step2-ProofResolvHodgeDirHs2}
    |\lambda|\lVert\tilde{u}_j^{\varepsilon} \rVert_{\dot{\H}^{s,2}(\RR^n_+)} + \lVert \nabla\tilde{u}_j^{\varepsilon} \rVert_{\dot{\H}^{s,2}(\RR^n_+)} \lesssim_{s,n}^\mu \big\lVert [\d - \mathbf{A}_j^\varepsilon \d {(\mathbf{A}_j^\varepsilon)}^{-1}]\tilde{u}_j^\varepsilon \big\rVert_{\dot{\H}^{s,2}(\RR^n_+)}   + \lVert \tilde{f}_j^\varepsilon\rVert_{\dot{\H}^{s,2}(\RR^n_+)}  + \lVert \tilde{g}_j^\varepsilon\rVert_{\dot{\H}^{s,2}(\RR^n_+)}.
\end{align}
By Proposition~\ref{prop:FundamentalExtby0HomFuncSpaces}, Lemma~\ref{lem:CompactEquiHomInhom}  and Lemma~\ref{lem:SmoothingPullback}:
\begin{align*}
\lVert \tilde{f}_j^\varepsilon\rVert_{\dot{\H}^{s,2}(\RR^n_+)}  + \lVert \tilde{g}_j^\varepsilon\rVert_{\dot{\H}^{s,2}(\RR^n_+)} \lesssim_{p,s,n}^\Omega \lVert {f}\rVert_{{\H}^{s,2}(\Omega)}  + \lVert u^\varepsilon\rVert_{{\H}^{s,2}(\Omega_\varepsilon)}, 
\end{align*}
where, in order to apply Lemma~\ref{lem:SmoothingPullback}, we recall that we did assume 
\begin{align}\label{eq:ProofHodgeDircSmallReguAssumpL2}
    -\frac{1}{2}<s<-\frac{1}{2} + \frac{\alpha}{2}.
\end{align}
We write
\begin{align*}
    \d - \mathbf{A}_j \d \prescript{t}{}{\mathbf{A}_j^{-1}} =\d -   \mathbf{A}_j\d  +   \mathbf{A}_j\d - \mathbf{A}_j \d \prescript{t}{}{(\mathbf{A}_j)^{-1}} = [\I-{\mathbf{A}_j}]  \d + \mathbf{A}_j\d [\I- \prescript{t}{}{\mathbf{A}_j^{-1}}].
\end{align*}
Therefore, thanks to compact support Lemma~\ref{lem:CompactEquiHomInhom}, we can bound
\begin{align*}
    \big\lVert [\d - \mathbf{A}_j^\varepsilon \d {(\mathbf{A}_j^\varepsilon)}^{-1}]\tilde{u}_j^\varepsilon \big\rVert_{\dot{\H}^{s,2}(\RR^n_+)} &\lesssim_{s,n,\Omega} \lVert \I-{\mathbf{A}_j^\varepsilon}\rVert_{\mathcal{M}_{\H}^{s,2}(\RR^n_+)}^n\big\lVert\d{(\mathbf{A}_j^\varepsilon)}^{-1}\tilde{u}_j^\varepsilon \big\rVert_{{\H}^{s,2}(\RR^n_+)} \\
    &\qquad\qquad+ \lVert {\mathbf{A}_j^\varepsilon}\rVert_{\mathcal{M}_{\H}^{s,2}(\RR^n_+)}^n\big\lVert\d [\I- ({\mathbf{A}_j^{\varepsilon})^{-1}}]\tilde{u}_j^\varepsilon \big\rVert_{{\H}^{s,2}(\RR^n_+)}\\
    &\lesssim_{s,n,\Omega} \Big(\lVert \I-{\mathbf{A}_j^\varepsilon}\rVert_{\mathcal{M}_{\H}^{s,2}(\RR^n_+)}^n \lVert {(\mathbf{A}_j^\varepsilon)}^{-1}\rVert_{\mathcal{M}_{\H}^{s+1,2}(\RR^n_+)}^n\\
    &\qquad\qquad+ \lVert {\mathbf{A}_j^\varepsilon}\rVert_{\mathcal{M}_{\H}^{s,2}(\RR^n_+)}^n\lVert [\I- ({\mathbf{A}^{\varepsilon}_j)^{-1}}]\rVert_{\mathcal{M}_{\H}^{s,2}(\RR^n_+)}^n\Big)\big\lVert\tilde{u}_j^\varepsilon \big\rVert_{{\H}^{s+1,2}(\RR^n_+)}.
\end{align*}
Now, we write
\begin{align*}
    \I- \mathbf{A}_j&= \I-\det(\nabla \boldsymbol{\Phi}_j)(\nabla  \boldsymbol{\Phi}_j)^{-1}\cdot\prescript{t}{}{(\nabla  \boldsymbol{\Phi}_j)}^{-1}\\&= \I - \det(\nabla \boldsymbol{\Phi}_j)(\nabla  \boldsymbol{\Phi}_j)^{-1} +\det(\nabla \boldsymbol{\Phi}_j)(\nabla  \boldsymbol{\Phi}_j)^{-1}[\I-\prescript{t}{}{(\nabla  \boldsymbol{\Phi}_j)}^{-1}],
\end{align*}
so that, thanks to Proposition~\ref{prop:MultipliersintheLplikeRange}, it holds
\begin{align*}
    \lVert \I- \mathbf{A}_j^\varepsilon \rVert_{\mathcal{M}^{s,2}_{\H,\tt{or}}(\RR^n_+)} &\leqslant \lVert \I - \det(\nabla \boldsymbol{\Phi}_j^\varepsilon)(\nabla  \boldsymbol{\Phi}_j^\varepsilon)^{-1} \rVert_{\mathcal{M}^{s,2}_{\H,\tt{or}}(\RR^n_+)} \\ &\qquad+ \lVert\det(\nabla \boldsymbol{\Phi}_j^\varepsilon)(\nabla  \boldsymbol{\Phi}_j^\varepsilon)^{-1} \rVert_{\mathcal{M}^{s,2}_{\H,\tt{or}}(\RR^n_+)}\lVert\I-\prescript{t}{}{(\nabla  \boldsymbol{\Phi}_j^\varepsilon)}^{-1} \rVert_{\mathcal{M}^{s,2}_{\H,\tt{or}}(\RR^n_+)}\\
    &\lesssim_n \big( 1+ \lVert \nabla\mathcal{T}\varphi_j^\varepsilon \rVert_{\mathcal{M}^{s,2}_{\H,\tt{or}}(\RR^n_+)}\big)\lVert \nabla\mathcal{T}\varphi_j^\varepsilon \rVert_{\mathcal{M}^{s,2}_{\H,\tt{or}}(\RR^n_+)} \\
    &\lesssim_{s,n}^\Omega\lVert \varphi_j^\varepsilon \rVert_{\mathcal{M}^{1+\alpha,r}_{\W}(\RR^{n-1})}\\
    &\lesssim_{s,n}^\Omega\lVert \varphi_j \rVert_{\mathcal{M}^{1+\alpha,r}_{\W}(\RR^{n-1})} <\epsilon.
\end{align*}
The last estimate follows from the fact that the regularization of $(\varphi_j^\varepsilon)_{\varepsilon}$ is a regularization by convolution of $\varphi_j$, so it preserves uniform bounds with respect to $\varepsilon$. Similarly, still due to Proposition~\ref{prop:MultipliersintheLplikeRange} and \eqref{eq:ProofHodgeDircSmallReguAssumpL2}, and taking advantage of the fact $\epsilon\lesssim 1$,
\begin{align*}
    \lVert (\mathbf{A}_j^\varepsilon)^{-1} \rVert_{\mathcal{M}^{s+1,2}_{\H,\tt{or}}(\RR^n_+)} &\lesssim_{n,\Omega} (1+\lVert \nabla\mathcal{T}\varphi_j^\varepsilon \rVert_{\mathcal{M}^{s+1,2}_{\H,\tt{or}}(\RR^n_+)})^2\\
    &\lesssim_{n,\Omega} (1+\lVert \varphi_j \rVert_{\mathcal{M}^{1+\alpha,r}_{\W}(\RR^{n-1})})^2\lesssim_{n,\Omega} 1.
\end{align*}
One can employ the same arguments in order to bound
\begin{align*}
    \lVert \I- (\mathbf{A}_j^\varepsilon)^{-1} \rVert_{\mathcal{M}^{s+1,2}_{\H,\tt{or}}(\RR^n_+)} \lesssim_{s,n} \epsilon,\quad \text{ and }\quad \lVert \mathbf{A}_j^\varepsilon \rVert_{\mathcal{M}^{s,2}_{\H,\tt{or}}(\RR^n_+)}\lesssim_{n,\Omega} 1.
\end{align*}
Consequently, \eqref{eq:Step2-ProofResolvHodgeDirHs2} becomes
\begin{align*}
    |\lambda|\lVert\tilde{u}_j^{\varepsilon} \rVert_{\dot{\H}^{s,2}(\RR^n_+)} + \lVert \nabla\tilde{u}_j^{\varepsilon} \rVert_{\dot{\H}^{s,2}(\RR^n_+)} \lesssim_{s,n,\Omega}^\mu \epsilon\lVert \tilde{u}_j^{\varepsilon} \rVert_{{\H}^{s+1,2}(\RR^n_+)}   + \lVert u^\varepsilon\rVert_{{\H}^{s,2}(\Omega_\varepsilon)}  + \lVert f\rVert_{{\H}^{s,2}(\Omega)},
\end{align*}
which by compact support also implies
\begin{align*}
    |\lambda|\lVert\tilde{u}_j^{\varepsilon} \rVert_{\dot{\H}^{s,2}(\RR^n_+)} + \lVert \nabla\tilde{u}_j^{\varepsilon} \rVert_{\dot{\H}^{s,2}(\RR^n_+)} &\lesssim_{s,n,\Omega}^\mu \epsilon^n \left(\lVert \tilde{u}_j^{\varepsilon} \rVert_{\dot{\H}^{s,2}(\RR^n_+)}+ \lVert \nabla \tilde{u}_j^{\varepsilon} \rVert_{\dot{\H}^{s,2}(\RR^n_+)}\right)   + \lVert u\rVert_{{\H}^{s,2}(\Omega)}  + \lVert f\rVert_{{\H}^{s,2}(\Omega)}\\
     &\lesssim_{s,n,\Omega}^\mu \epsilon^n \left(\lVert u\rVert_{{\H}^{s,2}(\Omega)}+ \lVert \nabla \tilde{u}_j^{\varepsilon} \rVert_{\dot{\H}^{s,2}(\RR^n_+)}\right)   + \lVert u\rVert_{{\H}^{s,2}(\Omega)}  + \lVert f\rVert_{{\H}^{s,2}(\Omega)}\\
     &\lesssim_{s,n,\Omega}^\mu \epsilon^n \lVert \nabla \tilde{u}_j^{\varepsilon} \rVert_{\dot{\H}^{s,2}(\RR^n_+)} + (1+\epsilon^n )\lVert u^\varepsilon\rVert_{{\H}^{s,2}(\Omega)}  + \lVert f\rVert_{{\H}^{s,2}(\Omega)}\\
     &\lesssim_{s,n,\Omega}^\mu \epsilon^n \lVert \nabla \tilde{u}_j^{\varepsilon} \rVert_{\dot{\H}^{s,2}(\RR^n_+)} + \lVert u^\varepsilon\rVert_{{\H}^{s,2}(\Omega_\varepsilon)}  + \lVert f\rVert_{{\H}^{s,2}(\Omega)}.
\end{align*}
So for $\epsilon$ chosen small enough, one obtains
\begin{align*}
    |\lambda|\lVert u^{\varepsilon} \rVert_{{\H}^{s,2}(\Omega_\varepsilon)} + \lVert \nabla{u}^{\varepsilon} \rVert_{{\H}^{s,2}(\Omega_\varepsilon)} &\lesssim_{s,n,\Omega}^\mu  \lVert u^\varepsilon\rVert_{{\H}^{s,2}(\Omega_\varepsilon)}  + \lVert f\rVert_{{\H}^{s,2}(\Omega)}.
\end{align*}
Hence $(u_\varepsilon)_\varepsilon$ is uniformly bounded in ${\H}^{s+1,2}(\Omega,\Lambda)$\footnote{recall that $f\in\Ccinfty(\Omega,\Lambda)\subset\Ccinfty(\Omega_\varepsilon,\Lambda)$, so that $(u_\varepsilon)_\varepsilon$ is uniformly bounded in $\L^2(\Omega,\Lambda)$} which implies that its limit $u$ also belongs to ${\H}^{s+1,2}(\Omega,\Lambda)$. Consequently, one can reperform the whole argument without any regularization procedure, in order to obtain independently of $\varepsilon>0$,
\begin{align}\label{eq:Step2-ProofResolvHodgeDirHs2Close}
    |\lambda|\lVert u \rVert_{{\H}^{s,2}(\Omega)} + \lVert \nabla{u} \rVert_{{\H}^{s,2}(\Omega)} \lesssim_{s,n,\Omega}^\mu  \lVert u\rVert_{{\H}^{s,2}(\Omega)}  + \lVert f\rVert_{{\H}^{s,2}(\Omega)}.
\end{align}
\textbf{Step 2.2:} \textbf{We prove full bisectoriality for all $-\sfrac{1}{2}<s<\sfrac{1}{2}$, and that the null space of the extrapolated operator is independent of $s$.} As a consequence of \eqref{eq:Step2-ProofResolvHodgeDirHs2Close} provided $R_0>0$ is a large enough but fixed number, depending on $s$, $n$, $\mu$, and $\Omega$,  one obtains the following bounds,
\begin{align}\label{eq:ProofHodge-DircAlmostResolventbound}
    \sup_{\lambda\in\mathrm{S}_\mu\setminus\B_{R_0}(0)}\,\lVert \lambda (\lambda\I+iD_{\mathfrak{n}})^{-1} \rVert_{{\H}^{s,2}(\Omega)\rightarrow{\H}^{s,2}(\Omega)}\lesssim_{s,n,R_0}^{\mu,\Omega} 1,
\end{align}
and
\begin{align*}
      \sup_{\lambda\in\mathrm{S}_\mu\setminus\B_{R_0}(0)}\,\lVert \nabla  (\lambda\I+iD_{\mathfrak{n}})^{-1} \rVert_{{\H}^{s,2}(\Omega)\rightarrow{\H}^{s,2}(\Omega)} \lesssim_{s,n,R_0}^{\mu,\Omega} 1,
\end{align*}
 where we recall that \eqref{eq:ProofHodgeDircSmallReguAssumpL2} is assumed. This implies that
 \begin{align*}
     (\D_2^s({D_\mathfrak{n}},\Omega,\Lambda),D_\mathfrak{n})=(\D_2^s(\d)\cap\D_2^s(\d^\ast,\Omega,\Lambda), \d + \d^\ast)
 \end{align*}
 is a closed operator, and that furthermore,
 \begin{align}
     \D_2^s({D_\mathfrak{n}},\Omega,\Lambda) \hookrightarrow\H^{s+1,2}_{\mathfrak{n}}(\Omega,\Lambda) , \quad -\sfrac{1}{2}<s<\sfrac{(\alpha-1)}{2}.
 \end{align}
 By duality and complex interpolation, due to self-adjointness, \eqref{eq:ProofHodge-DircAlmostResolventbound} remains valid for all $s\in(-\frac{1}{2},\frac{1}{2})$. So that one obtains a net of closed operators consistent operators
 \begin{align*}
     (\D_2^{s'}({D_\mathfrak{n}},\Omega,\Lambda),D_\mathfrak{n})_{|_{\H^{s,2}(\Omega,\Lambda)}}=(\D_2^s({D_\mathfrak{n}},\Omega,\Lambda),D_\mathfrak{n})\hookrightarrow \H^{\frac{1}{2},2}(\Omega,\Lambda), \quad\text{ for }\quad -\sfrac{1}{2}<s'<s<\sfrac{1}{2}
 \end{align*}
 Thus by \cite[Appendix~A,~Proposition~A.2.8~\&~Corollary~A.2.9]{bookHaase2006} and Proposition~\ref{prop:HodgeDiracL2bdd},
 \begin{align}\label{eq:ProofHodgeDiracRsolvSetHs2}
     \mathrm{S}_\mu\subset \rho_{\L^2}({D_\mathfrak{n}})\subset\rho_{\H^{s,2}}({D_\mathfrak{n}})\quad\text{ for }\quad 0<s<\frac{1}{2},
 \end{align}
and by duality for all $s\in(-\sfrac{1}{2},\sfrac{1}{2})$. One still has to prove the uniform resolvent bound.

 Now, if we let $v\in\N_2^s({D_\mathfrak{n}},\Omega,\Lambda)$, $-\frac{1}{2}<s<\frac{1}{2}$ then it holds that $v\in \H^{(\frac{1}{2}+\frac{\alpha}{2})^-,2}(\Omega,\Lambda)\hookrightarrow \L^2(\Omega,\Lambda)$, with
 \begin{align*}
     \d v + \d^\ast v =0
 \end{align*}
 This implies $v\in \N_2(\d)\cap\N_2(\d^\ast,\Omega,\Lambda)=\N_2({D_\mathfrak{n}},\Omega,\Lambda)$, consequently
 \begin{align*}
     \N_2^s({D_\mathfrak{n}},\Omega,\Lambda) = \N_2({D_\mathfrak{n}},\Omega,\Lambda).
 \end{align*}
 Recall that $\N_2({D_\mathfrak{n}},\Omega,\Lambda)$ is finite dimensional so every well-defined norms are equivalent, this yields that the projector $\P_\H$ from Proposition~\ref{prop:FullHodgedecompL2bdd}, is also bounded as a linear operator 
 \begin{align*}
     \P_\H\,:\,\H^{s,2}(\Omega,\Lambda)\longrightarrow\N_2^s({D_\mathfrak{n}},\Omega,\Lambda).
 \end{align*}
Now, we assume again $-\sfrac{1}{2}<s<\sfrac{(\alpha-1)}{2}$. By a standard compactness argument by contradiction, thanks to the compact embedding $\H^{s+1}(\Omega,\Lambda)\hookrightarrow\H^{s}(\Omega,\Lambda)$, it holds
\begin{align*}
    \lVert v-\P_\H v \rVert_{\H^{s,2}(\Omega)} \lesssim_{s,n}^\Omega \lVert D_\mathfrak{n}v \rVert_{\H^{s,2}(\Omega)},  \quad\text{ for all }v \in\D_2^s({D_\mathfrak{n}},\Omega,\Lambda), \quad -\sfrac{1}{2}<s<\sfrac{(\alpha-1)}{2}.
\end{align*}
This also implies the closedness of the range. Now, for $u$ the solution of \eqref{eq:HodgeDirOr} for $f\in\Ccinfty(\Omega,\Lambda)$ and $\lambda\in\mathrm{S}_\mu$ fixed,
\begin{align*}
    \lVert u-\P_\H u \rVert_{\H^{s,2}(\Omega)} \lesssim_{s,n}^\Omega \lVert D_\mathfrak{n} [u-\P_\H u ] \rVert_{\H^{s,2}(\Omega)}&=\lVert D_\mathfrak{n} u \rVert_{\H^{s,2}(\Omega)}\\
    &= \big\lVert [\I-\P_\H ]D_\mathfrak{n} u \big\rVert_{\H^{s,2}(\Omega)}\\
    &= \big\lVert \lambda [\I-\P_\H ] u + [\I-\P_\H ] f \big\rVert_{\H^{s,2}(\Omega)}\\
    &\lesssim_{s,n,\Omega} |\lambda|\big\lVert [\I-\P_\H ] u\big\rVert_{\H^{s,2}(\Omega)}+\lVert  f \big\rVert_{\H^{s,2}(\Omega)},
\end{align*}
so that for $r_0>0$ fixed but small enough, for $\lambda\in\mathrm{S}_\mu\cap\B_{r_0}(0)$:
\begin{align*}
    \left(1-\frac{|\lambda|}{2r_0}\right)\lVert u-\P_\H u \rVert_{\H^{s,2}(\Omega)} \lesssim_{s,n,\Omega}^{r_0,\mu}\lVert  f \rVert_{\H^{s,2}(\Omega)}.
\end{align*}
Consequently, due to the identity $\lambda \P_\H u =\P_\H f$, one obtains
\begin{align}\label{eq:ProofHodge-DircAlmostResolventboundHs2small}
    |\lambda|\lVert u\rVert_{\H^{s,2}(\Omega)} &\leqslant |\lambda|\lVert  [\I-\P_\H] u\rVert_{\H^{s,2}(\Omega)} + \lVert \P_\H f\rVert_{\H^{s,2}(\Omega)}\nonumber\\&\lesssim_{s,n,} |\lambda|\lVert  u-\P_\H u\rVert_{\H^{s,2}(\Omega)} + \lVert  f\rVert_{\H^{s,2}(\Omega)}\nonumber\\
    &\lesssim_{s,n,\Omega}^{r_0,\mu} \frac{2r_0|\lambda|}{2r_0-|\lambda|}\lVert f\rVert_{\H^{s,2}(\Omega)} + \lVert  f\rVert_{\H^{s,2}(\Omega)}\nonumber\\
    &\lesssim_{s,n,\Omega}^{r_0,\mu} (2r_0+1) \lVert  f\rVert_{\H^{s,2}(\Omega)}.
\end{align}
From there, combining \eqref{eq:ProofHodge-DircAlmostResolventbound} and \eqref{eq:ProofHodge-DircAlmostResolventboundHs2small}
\begin{align}\label{eq:ProofHodgeDirResolventBoundexceptAnnuli}
    \sup_{\substack{\lambda\in\mathrm{S}_\mu,\\ \lambda\notin\overline{\B_{R_0}(0)}\setminus\B_{r_0}(0)}} \,\lVert \lambda (\lambda\I+iD_{\mathfrak{n}})^{-1} \rVert_{{\H}^{s,2}(\Omega)\rightarrow{\H}^{s,2}(\Omega)}\lesssim_{s,n,R_0,r_0}^{\mu,\Omega} 1, \quad -\sfrac{1}{2}<s<\sfrac{(\alpha-1)}{2}.
\end{align}
Concerning the annuli part $\mathcal{C}_{R_0,r_0}^{\mu}:=\mathrm{S}_\mu\cap [\overline{\B_{R_0}(0)}\setminus\B_{r_0}(0)]$, by \eqref{eq:ProofHodgeDiracRsolvSetHs2}, one already has
\begin{align*}
    \mathcal{C}_{R_0,r_0}^{\mu}\subset\rho_{\H^{s,2}}(D_\mathfrak{n}).
\end{align*}
so that, the number $\mu\in(0,\frac{\pi}{2})$ being arbitrary, by holomorphy of
\begin{align*}
    \lambda\longmapsto\lambda(\lambda\I+iD_\mathfrak{n})^{-1}
\end{align*}
in a relatively compact neighborhood of $\mathcal{C}_{R_0,r_0}^{\mu}$ in $\mathrm{S}_{\mu'}$ for $\mu'\in(\mu,\frac{\pi}{2})$, one obtains uniform boundedness of the resolvent on $\mathcal{C}_{R_0,r_0}^{\mu}$, and this, in combination with \eqref{eq:ProofHodgeDirResolventBoundexceptAnnuli}, yields the full resolvent bound
\begin{align}
    \sup_{\substack{\lambda\in\mathrm{S}_\mu}} \,\lVert \lambda (\lambda\I+iD_{\mathfrak{n}})^{-1} \rVert_{{\H}^{s,2}(\Omega)\rightarrow{\H}^{s,2}(\Omega)}\lesssim_{s,n,R_0,r_0}^{\mu,\Omega} 1, \quad -\sfrac{1}{2}<s<\sfrac{(-1+\alpha)}{2},
\end{align}
whenever $s\in(-\frac{1}{2},\frac{-1+\alpha}{2})$, complex interpolation with $\L^2$ and duality yields the result on $\H^{s,2}$ for all $s\in(-\frac{1}{2},\frac{1}{2})$.

\textbf{Step 3:} \textbf{We deal with the iterative scheme allowing to reach $\H^{s,p}$-spaces for all $1<p<\infty$.} Note that the whole argument will be similar to Step 2, thanks to Proposition~\ref{prop:MultipliersintheLplikeRange} and Remark~\ref{rem:BetterMultipliers}, so we just achieve the first iteration. We set 
\begin{align*}
    p_0:= \frac{n-1-\alpha}{n-1}2, \quad\text{ and }\quad p_{k+1}:=\max\Big(\frac{n-1-\alpha}{n-1}p_k,1\Big), \quad k\in\NN.
\end{align*}
Note that $p_k\in[1,2)$ for all $k\in\NN$ and also that $1$ is reached in a finite number of step. For the reader's convenience, we mention that $p_0\in(1,2)$ is the minimal number larger than 1 such that the following embedding holds
\begin{align*}
    \H^{-1+\frac{\alpha+1}{p_0},p_0}(\Omega)\hookrightarrow \H^{-\frac{1}{2},2}(\Omega).
\end{align*}
Similarly, $p_k$ is the minimal real number larger than $1$ such that 
\begin{align*}
    \H^{-1+\frac{1+\alpha}{p_{k}},p_{k}}(\Omega)\hookrightarrow \H^{-1+\frac{1}{p_{k-1}},p_{k-1}}(\Omega).
\end{align*}

For $p\in(p_0,p_0')$, $s\in(-1+\sfrac{1}{p},-1+\sfrac{(1+\alpha)}{p}$, one can reproduce Step 2.1 and re-obtain for $\lambda\in\mathrm{S}_\mu$, $\mu\in(0,\frac{\pi}{2})$, and $f\in\Ccinfty(\Omega,\Lambda)$, that the solution $u$ is such that $u\in\H^{s+1,p}(\Omega,\Lambda)$ and
\begin{align}\label{eq:Step2-ProofResolvHodgeDirHspClose}
    |\lambda|\lVert u \rVert_{{\H}^{s,p}(\Omega)} + \lVert \nabla{u} \rVert_{{\H}^{s,p}(\Omega)} \lesssim_{s,n,\Omega}^\mu  \lVert u\rVert_{{\H}^{s,p}(\Omega)}  + \lVert f\rVert_{{\H}^{s,p}(\Omega)}.
\end{align}
As in Step 2.1, this implies that
 \begin{align*}
     (\D_p^s({D_\mathfrak{n}},\Omega,\Lambda),D_\mathfrak{n})=(\D_p^s(\d)\cap\D_p^s(\d^\ast,\Omega,\Lambda), \d + \d^\ast)
 \end{align*}
 is a closed operator, with
 \begin{align}
     \D_p^s({D_\mathfrak{n}},\Omega,\Lambda) \hookrightarrow\H^{s+1,p}_{\mathfrak{n}}(\Omega,\Lambda) , \quad -1+\sfrac{1}{p}<s<-\sfrac{1}{p}+\sfrac{(1+\alpha)}{p},
 \end{align}
 with a resolvent estimate --provided  $R_0>0$ is depending on $s$, $p$, $n$, $\Omega$ and $\mu$--
 \begin{align}\label{eq:ProofHodge-DircAlmostResolventboundHsp}
    \sup_{\substack{\lambda\in\mathrm{S}_\mu\setminus\B_{R_0}(0)}}\lVert \lambda (\lambda\I+iD_{\mathfrak{n}})^{-1} \rVert_{{\H}^{s,p}(\Omega)\rightarrow{\H}^{s,p}(\Omega)} + \lVert \nabla  (\lambda\I+iD_{\mathfrak{n}})^{-1} \rVert_{{\H}^{s,p}(\Omega)\rightarrow{\H}^{s,p}(\Omega)} \lesssim_{p,s,n,R_0}^{\mu,\Omega} 1,
\end{align}
 
 By duality and complex interpolation, due to self-adjointness on $\L^2$, \eqref{eq:ProofHodge-DircAlmostResolventboundHsp} remains valid for all $s\in(-1+\frac{1}{p},\frac{1}{p})$ all $p\in(p_0,p_0')$. So that, one obtains a net of closed operators consistent operators
 \begin{align}\label{eq:RestrictOperatorHsp}
     (\D_p^{s'}({D_\mathfrak{n}},\Omega,\Lambda),D_\mathfrak{n})_{|_{\H^{s,p}(\Omega,\Lambda)}}=(\D_p^s({D_\mathfrak{n}},\Omega,\Lambda),D_\mathfrak{n}), \quad\text{ for }\quad -1+\sfrac{1}{p}<s'<s<\sfrac{1}{p}.
 \end{align}
 When $p>2$, by Sobolev embeddings one can always find $\tilde{s}\in(-{\sfrac{1}{2}},{\sfrac{1}{2}})$ depending on $s$ below, such that
  \begin{align*}
     (\D_2^{\tilde{s}}({D_\mathfrak{n}},\Omega,\Lambda),D_\mathfrak{n})_{|_{\H^{s,p}(\Omega,\Lambda)}}=(\D_p^s({D_\mathfrak{n}},\Omega,\Lambda),D_\mathfrak{n}), \quad\text{ for }\quad -1+\sfrac{1}{p}<s<-1+\sfrac{1+\alpha}{p},
 \end{align*}
and thanks to the fact $2<p<p_0'$
 \begin{align*}
    (\D_2^{\tilde{s}}({D_\mathfrak{n}},\Omega,\Lambda),D_\mathfrak{n}) \hookrightarrow \H^{\frac{1}{2},2}(\Omega,\Lambda) \hookrightarrow \H^{s,p}(\Omega,\Lambda), \quad\text{ for }\quad -1+\sfrac{1}{p}<s<-1+\sfrac{1+\alpha}{p}.
 \end{align*}
And in a similar but simpler way for $p_0<p<2$
  \begin{align*}
     (\D_2^{s}({D_\mathfrak{n}},\Omega,\Lambda),D_\mathfrak{n})_{|_{\H^{s,p}(\Omega,\Lambda)}}=(\D_p^s({D_\mathfrak{n}},\Omega,\Lambda),D_\mathfrak{n}), \quad\text{ for }\quad -1+\sfrac{1}{p}<s<-1+\sfrac{1+\alpha}{2}
 \end{align*}
So that one can apply again \cite[Appendix~A,~Proposition~A.2.8~\&~Corollary~A.2.9]{bookHaase2006} and Proposition~\ref{prop:HodgeDiracL2bdd}, yielding
 \begin{align}\label{eq:ProofHodgeDiracRsolvSetHsp}
     \mathrm{S}_\mu  \subset \rho_{\H^{s,p}}({D_\mathfrak{n}})\quad\text{ for all }\quad p\in(p_0,p_0'),\quad -1+\sfrac{1}{p}<s<-1+\sfrac{(1+\alpha)}{p},
 \end{align}
and by duality it also holds whenever $s\in(-\sfrac{\alpha}{p'}+\sfrac{1}{p},\sfrac{1}{p})$. Now, \eqref{eq:RestrictOperatorHsp} implies that \eqref{eq:ProofHodgeDiracRsolvSetHsp} holds for all $p\in(2,p_0')$ and all $s\in(-1+\sfrac{1}{p},\sfrac{1}{p})$, and by duality for all $p\in(p_0,p_0')$. It remains to control the behavior of the resolvent with respect to $\lambda\in \mathrm{S}_\mu$.

\medbreak

Now, assuming $ -1+\sfrac{1}{p}<s<\sfrac{1}{p}$, $p\in(p_0,p_0')$, by previous Step 2 we improve the description of $\N_2(D_\mathfrak{n},\Omega,\Lambda)$, so that by Sobolev embeddings on has
\begin{align*}
    \N_2(D_\mathfrak{n},\Omega,\Lambda) =\{ \, u \in \H^{(\frac{1+\alpha}{2})^-,2}(\Omega,\Lambda)\,:\,\d u = \delta u= 0,\, \nu\iprod u_{|_{\partial\Omega}}=0 \} \hookrightarrow \N_{p}^s(D_\mathfrak{n},\Omega,\Lambda).
\end{align*}
Now, let $v\in\N_{p}^s(D_\mathfrak{n},\Omega,\Lambda)$, then one has either $v\in\D_{p}^s(D_\mathfrak{n},\Omega,\Lambda)\subset\H^{s+1,p}(\Omega,\Lambda)$ if $ -1+\sfrac{1}{p}<s<-1+\sfrac{(1+\alpha)}{p}$, or $v\in\D_{p}^s(D_\mathfrak{n},\Omega,\Lambda)\subset\H^{(\frac{1+\alpha}{p})^-,p}(\Omega,\Lambda)$, if $s\in [-1+\sfrac{(1+\alpha)}{p},\sfrac{1}{p})$. In both cases one obtains $v\in \L^2(\Omega,\Lambda)$, $\nu\iprod v_{|_{\partial\Omega}}=0$ and
\begin{align*}
    \d v + \delta v=0.
\end{align*}
The latter implies $v\in \N_2(D_\mathfrak{n},\Omega,\Lambda)$. Consequently, we deduce
\begin{align*}
    \N_2(D_\mathfrak{n},\Omega,\Lambda) =\N_{p}^s(D_\mathfrak{n},\Omega,\Lambda)\quad\text{ for all }\quad p\in(p_0,p_0'),\quad -1+\sfrac{1}{p}<s<\sfrac{1}{p}.
\end{align*}
We recall that $\N_2(D_\mathfrak{n},\Omega,\Lambda)$ is finite dimensional, so all finite norms are equivalent. Therefore for all $w\in\H^{s,p}(\Omega,\Lambda)$,
\begin{align*}
    \lVert \P_\H w\rVert_{\H^{s,p}(\Omega)} \lesssim_{n,\Omega}\lVert  w\rVert_{\H^{s,p}(\Omega)}.
\end{align*}

Now, the exact same compactness argument as in Step 2.2 yields
\begin{align*}
    \lVert v-\P_\H v \rVert_{\H^{s,p}(\Omega)} \lesssim_{s,n}^\Omega \lVert D_\mathfrak{n}v \rVert_{\H^{s,p}(\Omega)},  \quad\text{ for all }v \in\D_p^s({D_\mathfrak{n}},\Omega,\Lambda), \quad -1+\sfrac{1}{p}<s<-1+\sfrac{(1+\alpha)}{p}.
\end{align*}

Therefore, one can conclude with the exact same arguments as in the end of Step 2.2, which yields the full bisectoriality estimate for all $p\in(p_0,p_0')$ and all $s\in(-1+\sfrac{1}{p},\sfrac{1}{p})$, by duality and complex interpolation. Iterating the argument we did provide at the next step $k\geqslant 1$, for $p\in(p_k,p_k')$, one should replace the use of the family $(\H^{s,2})_{|s|<\frac{1}{2}}$ by $(\H^{\tau,q})_{\tau\in(-1+\sfrac{1}{q},\sfrac{1}{q}),q\in(p_{k-1},p_{k-1}')}$. Thus, we obtain the $\H^{s,p}$-result for all $1<p<\infty$. The corresponding case for Besov spaces follows from real interpolation.

\textbf{Step 4:} \textbf{We deal with the case of Besov spaces $p=1,\infty$.} We concentrate on the spaces $\B^{s}_{1,q}$, $1< q<\infty$, $0<s<1$. Thanks to Sobolev embeddings and the case $1<p<\infty$, whenever $s\in(0,\alpha)$, one can reproduce Step 2.1, --recall that $\B^{s}_{1,q}$ is reflexive, see Theorem~\ref{thm:Reflexiveendpointspaces}-- and obtain that $u\in\B^{s+1}_{1,q}(\Omega,\Lambda)$, and
\begin{align}\label{eq:Step2-ProofResolvHodgeDirBs1qClose}
    |\lambda|\lVert u \rVert_{{\B}^{s}_{1,q}(\Omega)} + \lVert \nabla{u} \rVert_{{\B}^{s}_{1,q}} \lesssim_{s,n,\Omega}^\mu  \lVert u\rVert_{{\B}^{s}_{1,q}(\Omega)}  + \lVert f\rVert_{{\B}^{s}_{1,q}(\Omega)}.
\end{align}
So the following bound holds for any $\lambda\in\mathrm{S}_\mu\setminus\B_{R_0}(0)$:
\begin{align}\label{eq:ProofHodge-DircAlmostResolventboundBs1q}
    \lVert \lambda (\lambda\I+iD_{\mathfrak{n}})^{-1} \rVert_{{\B}^{s}_{1,q}(\Omega)\rightarrow{\B}^{s}_{1,q}(\Omega)} + \lVert \nabla  (\lambda\I+iD_{\mathfrak{n}})^{-1} \rVert_{{\B}^{s}_{1,q}(\Omega)\rightarrow{\B}^{s}_{1,q}(\Omega)} \lesssim_{s,n,R_0}^{\mu,\Omega} 1,
\end{align}
Note with $s<\alpha$, that one has
\begin{align*}
    \D_{1,q}^{s}(D_\mathfrak{n},\Omega,\Lambda)=\B^{s+1}_{1,q,\mathfrak{n}}(\Omega,\Lambda).
\end{align*}
Since $s>0$, one obtains directly from  the previous steps and Sobolev embeddings, and a suitable $r\in(1,\infty)$,
\begin{align*}
    \N_{2}(D_\mathfrak{n},\Omega,\Lambda) = \N_{2,q}^s(D_\mathfrak{n},\Omega,\Lambda) \hookrightarrow\N_{1,q}^s(D_\mathfrak{n},\Omega,\Lambda)\hookrightarrow \N_{r}(D_\mathfrak{n},\Omega,\Lambda)=\N_{2}(D_\mathfrak{n},\Omega,\Lambda).
\end{align*}
and from the case $p$ large, by Sobolev embeddings, one also obtains
\begin{align*}
    \N_{2}(D_\mathfrak{n},\Omega,\Lambda)=\N_{p}(D_\mathfrak{n},\Omega,\Lambda) \hookrightarrow \N_{\infty,q'}^\tau(D_\mathfrak{n},\Omega,\Lambda), \quad\text{for all } -1<\tau<0.
\end{align*}
Hence,  $\P_\H$ is well-defined and bounded on ${\B}^{s}_{1,q}(\Omega,\Lambda)$.

Thus, from the same compactness argument from Step 2.2 (again recall that Banach-Alaoglu is applicable due to Theorem~\ref{thm:Reflexiveendpointspaces}) one obtains the following resolvent bound for all $s\in(0,\alpha)$
\begin{align*}
    \sup_{\lambda\in \mathrm{S}_\mu}\, \lVert \lambda (\lambda\I+iD_{\mathfrak{n}})^{-1} \rVert_{{\B}^{s}_{1,q}(\Omega)}\lesssim_{s,n,R_0}^{\mu,\Omega} 1,
\end{align*}
It remains to prove the result when $s\in[\alpha,1)$. We consider, $s_0\in(0,\alpha)$ to be fixed, for $\theta\in(0,1)$ such that $s:=s_0+1-\theta<1$, one has
\begin{align*}
({\B}^{s_0}_{1,q}(\Omega,\Lambda),\D_{1,q}^{s_0}(D_\mathfrak{n},\Omega,\Lambda) )_{\theta,q} &= ({\B}^{s_0}_{1,q}(\Omega,\Lambda),{\B}^{s_0+1}_{1,q,\mathfrak{n}}(\Omega,\Lambda) )_{\theta,q}\\
&\hookrightarrow  ({\B}^{s_0}_{1,q}(\Omega,\Lambda),{\B}^{s_0+1}_{1,q}(\Omega,\Lambda) )_{\theta,q} = {\B}^{s}_{1,q}(\Omega,\Lambda),
\end{align*}
and since $0<s<1$,
\begin{align*}
    {\B}^{s}_{1,q}(\Omega,\Lambda) = {\B}^{s}_{1,q,0}(\Omega,\Lambda) = ({\B}^{s_0}_{1,q,0}(\Omega,\Lambda),{\B}^{s_0+1}_{1,q,0}(\Omega,\Lambda) )_{\theta,q} &\hookrightarrow({\B}^{s_0}_{1,q}(\Omega,\Lambda),{\B}^{s_0+1}_{1,q,\mathfrak{n}}(\Omega,\Lambda) )_{\theta,q}\\
    &= ({\B}^{s_0}_{1,q}(\Omega,\Lambda),{\B}^{s_0+1}_{1,q,\mathfrak{n}}(\Omega,\Lambda) )_{\theta,q}.
\end{align*}
Thus, we did obtain
\begin{align*}
    ({\B}^{s_0}_{1,q}(\Omega,\Lambda),\D_{1,q}^{s_0}(D_\mathfrak{n},\Omega,\Lambda) )_{\theta,q}= {\B}^{s}_{1,q}(\Omega,\Lambda),
\end{align*}
and consequently being bisectorial on ${\B}^{s_0}_{1,q}(\Omega,\Lambda)$ and $\D_{1,q}^{s_0}(D_\mathfrak{n},\Omega,\Lambda)$ by the beginning of the current step, one obtains the resolvent bound for $s\in[\alpha,1)$ by interpolation. By duality, one obtains the case $p=\infty$, and by real interpolation the cases $q=1,\infty$. If $r=\infty$, one can also reproduce the localisation procedure for the estimates to self improve regularity of solution when $s\in(-1,-1+\alpha)$. We mention that if $r=\infty$ one could achieve directly the case of the endpoint spaces $\B^{s}_{\infty,q}$ through the localisation techniques instead of a duality argument. The proof is now finished.
\end{proof}

Through the proof of Proposition~\ref{prop:HodgeDiracbdd}, we did derive the following fundamental result as a corollary.

\begin{corollary}\label{cor:HarmonicformsProj}Let $\alpha\in(0,1)$, $r\in[1,\infty]$, and let $\Omega$ be a bounded  Lipschitz  domain of class $\mathcal{M}^{1+\alpha,r}_\W(\epsilon)$ for some sufficiently small $\epsilon>0$. For all $p,q\in[1,\infty]$, $s\in(-1+\sfrac{1}{p},\sfrac{1}{p})$, the projection onto Harmonic forms
\begin{align*}
    \P_\H\,:\,\B^{s}_{p,q}(\Omega,\Lambda)\longrightarrow \N_2(\d)\cap\N_2(\d^\ast,\Omega,\Lambda) \subset \B^{s}_{p,q}(\Omega,\Lambda),
\end{align*}
is well-defined and bounded, and its finite-dimensional range is independent of $p,q,s$, and similarly with $\H^{s,p}(\Omega,\Lambda)$, $1<p<\infty$.

\medbreak

\noindent Furthermore, for all $q\in[1,\infty]$, all $p\in[1,\infty]$ and $s$ such that either
 \begin{enumerate}
    \item $p\in(r,\infty]$, $s\in(-1+\frac{1}{p},-1+\frac{1+r\alpha}{p})$; or
    \item $p\in[1,r]$, $s\in(-1+\frac{1}{p},-1+\alpha+\frac{1}{p})$; or
    \item $p=q=r$, $s=-1+\alpha+\frac{1}{p}$;
\end{enumerate}
one has
\begin{align*}
    \N_2(\d)\cap\N_2(\d^\ast,\Omega,\Lambda) \subset \B^{s+1}_{p,q}(\Omega,\Lambda).
\end{align*}
The same result holds if one replaces $(\d,\d^\ast)$ by $(\delta^\ast,\delta)$.
\end{corollary}

The novelty of the next result is the fact that it allows $p=1,\infty$, under the assumption on the domain to be in a multiplier class, it particular the result applies for bounded $\C^{1,\alpha}$-domains, $\alpha>0$, taking $r=\infty$. This seems to be  up to now the only result that includes endpoint function spaces for differential forms, even in the smooth setting.

\begin{theorem}\label{thm:SharpHodgeDecompC1}Let $\alpha\in(0,1)$, $r\in[1,\infty]$, and let $\Omega$ be a either $\RR^n_+$ or a bounded  Lipschitz  domain of class $\mathcal{M}^{1+\alpha,r}_\W(\epsilon)$ for some sufficiently small $\epsilon>0$.

\medbreak

\noindent For all $p,q\in[1,\infty]$, $s\in(-1+\sfrac{1}{p},\sfrac{1}{p})$, the following Hodge decompositions hold:
    \begin{align*}
        \dot{\B}^{s}_{p,q}(\Omega,\Lambda) &= \dot{\B}^{s,\sigma}_{p,q,\mathfrak{n}}(\Omega,\Lambda)\oplus\, \overline{\d {\B}^{s+1}_{p,q}(\Omega,\Lambda)};\\
        &= \dot{\B}^{s,\sigma}_{p,q}(\Omega,\Lambda)\oplus\, \overline{\d {\B}^{s+1}_{p,q,0}(\Omega,\Lambda)},
    \end{align*}
with respective bounded Hodge-Leray projections
\begin{align*}
            \PP_\Omega&\,:\,\dot{\B}^{s}_{p,q}(\Omega,\Lambda) \longrightarrow \dot{\B}^{s,\sigma}_{p,q,\mathfrak{n}}(\Omega,\Lambda),\\
\text{ and }\qquad\QQ_\Omega&\,:\,\dot{\B}^{s}_{p,q}(\Omega,\Lambda) \longrightarrow \dot{\B}^{s,\sigma}_{p,q}(\Omega,\Lambda).
    \end{align*}
Furthermore,
\begin{itemize}
    \item a similar result holds for the Besov spaces $\dot{\BesSmo}^{s}_{p,\infty}$, $\dot{\B}^{s,0}_{\infty,q}$, $\dot{\BesSmo}^{s,0}_{\infty,\infty}$ and Sobolev spaces $\dot{\H}^{s,p}$, assuming $1<p<\infty$ for the latter;

    \item the closure is taken weakly-$\ast$ whenever $q=\infty$.

    \item a similar result holds if one replaces $(\sigma,\mathfrak{n},\d)$ by $(\gamma,\mathfrak{t},\delta)$.
\end{itemize}
\end{theorem}

\begin{proof} We only check the case of bounded  Lipschitz  domain of class $\mathcal{M}^{1+\alpha,r}_\W(\epsilon)$ for some small $\epsilon>0$. By Lemma~\ref{lem:approxHodgeLeray}, for any $f\in\Ccinfty(\Omega,\Lambda)$, one has
\begin{align*}
    \PP_\Omega f &= \lim_{t\rightarrow 0}\,\d^\ast(it+ D_\mathfrak{n})^{-1}\d(-it+ D_\mathfrak{n})^{-1}[f-\P_\H f] + \P_\H f.
\end{align*}
Therefore, by the Fatou property for function spaces, Proposition~\ref{prop:HodgeDiracbdd}, and Corollary~\ref{cor:HarmonicformsProj}, it holds 
\begin{align*}
    \lVert \PP_\Omega f \rVert_{\B^{s}_{p,q}(\Omega)}&\leqslant \liminf_{t\rightarrow 0} \,\lVert \d^\ast(it+ D_\mathfrak{n})^{-1}\d(-it+ D_\mathfrak{n})^{-1}[f-\P_\H f] + \P_\H f\rVert_{\B^{s}_{p,q}(\Omega)}\\
    &\lesssim_{p,n,s,\Omega}  \lVert f-\P_\H f\rVert_{\B^{s}_{p,q}(\Omega)}+\lVert  \P_\H f\rVert_{\B^{s}_{p,q}(\Omega)}\\
    &\lesssim_{p,n,s,\Omega} \lVert  f\rVert_{\B^{s}_{p,q}(\Omega)}.
\end{align*}
assuming first $1<q<\infty$, one can extend the operator by density, and the same argument applies for the Sobolev spaces $\H^{s,p}$. The case $q=1,\infty$ holds by real interpolation. 
\end{proof}

In the case of vector fields, we obtain the following fundamental consequence.
\begin{proposition}\label{prop:NeumannPbLipschitzMult}Let $\alpha\in(0,1)$, $r\in[1,\infty]$, and let $\Omega$ be a bounded  Lipschitz  domain of class $\mathcal{M}^{1+\alpha,r}_\W(\epsilon)$ for some sufficiently small $\epsilon>0$.

\medbreak

Let $p,q\in[1,\infty]$, $s\in(-1+\sfrac{1}{p},\sfrac{1}{p})$. For all $\uu\in\B^{s}_{p,q}(\Omega,\CC^n)$,  the following Neumann problem
\begin{equation*}\tag{$\mathcal{NL}_{\Omega}$}\label{NeuLapDomain}
    \left\{ \begin{array}{rllr}
         - \Delta \mathfrak{p} &=-\div \uu \text{, }&&\text{ in } \Omega\text{,}\\ 
        \nu\cdot \left(\nabla \mathfrak{p}-\uu\right)_{|_{\partial\Omega}} &= 0\text{, }&&\text{ on } \partial\Omega\text{.}
    \end{array}
    \right.
\end{equation*}
 admits a unique solution $\mathfrak{p}\in{\B}^{s+1}_{p,q,\mfree}(\Omega,\CC)$, with the estimate
\begin{align*}
  \lVert  \nabla\mathfrak{p}\rVert_{{\B}^{s}_{p,q}(\Omega)}\lesssim_{p,n,s,\Omega} \lVert \uu\rVert_{{\B}^{s}_{p,q}(\Omega)}.
\end{align*}
Furthermore, a similar result holds for the endpoint Besov spaces ${\BesSmo}^{\bullet}_{p,\infty}$, and Bessel potential spaces ${\H}^{\bullet,p}(\Omega)$, assuming $1<p<\infty$ for the latter.
\end{proposition}

The next two results are straightforward consequences of Theorem~\ref{thm:SharpHodgeDecompC1}.

\begin{proposition}\label{prop:DualityDivFreeRn+} Let $\alpha\in(0,1)$, $r\in[1,\infty]$, and let $\Omega$ be a either $\RR^n_+$ or a bounded Lipschitz domain of class $\mathcal{M}^{1+\alpha,r}_\W(\epsilon)$ for some sufficiently small $\epsilon>0$. Let $1\leqslant p,q \leqslant \infty$, $s\in(-1+\sfrac{1}{p},\sfrac{1}{p})$. Then, the following duality identities hold:
\begin{align*}
\begin{array}{rlrlrlll}
    (\dot{\H}^{-s,p'}_{\mathfrak{n},\sigma}(\Omega))'&=\dot{\H}^{s,p}_{\mathfrak{n},\sigma}(\Omega)\, \text{,}&&(\dot{\H}^{-s,p'}_{\sigma}(\Omega))'&=\dot{\H}^{s,p}_{\sigma}(\Omega)\, \text{,}&&\text{where }1<p<\infty\text{;}\\
    (\dot{\B}^{-s,\sigma}_{p',q',\mathfrak{n}}(\Omega))'&=\dot{\B}^{s,\sigma}_{p,q,\mathfrak{n}}(\Omega)\, \text{,}&&(\dot{\B}^{-s,\sigma}_{p',q'}(\Omega))'&=\dot{\B}^{s,\sigma}_{p,q}(\Omega)\, \text{,}&& \text{where } p,q>1\text{;}\\
    (\dot{\mathcal{B}}^{-s,\sigma}_{p',\infty,\mathfrak{n}}(\Omega))'&=\dot{\B}^{s,\sigma}_{p,1,\mathfrak{n}}(\Omega)\, \text{,}&&(\dot{\mathcal{B}}^{-s,\sigma}_{p',\infty}(\Omega))'&=\dot{\B}^{s,\sigma}_{p,1}(\Omega)\, \text{,}&&\text{where }p>1,q=1\text{;}\\
    (\dot{\B}^{-s,0,\sigma}_{\infty,q',\mathfrak{n}}(\Omega))'&=\dot{\B}^{s,\sigma}_{1,q,\mathfrak{n}}(\Omega)\, \text{,}&&(\dot{\B}^{-s,0,\sigma}_{\infty,q'}(\Omega))'&=\dot{\B}^{s,\sigma}_{1,q}(\Omega)\, \text{,}&&\text{where }p=1,q>1\text{;}\\
    (\dot{\mathcal{B}}^{-s,0,\sigma}_{\infty,\infty,\mathfrak{n}}(\Omega))'&=\dot{\B}^{s,\sigma}_{1,1,\mathfrak{n}}(\Omega)\, \text{,}&&(\dot{\mathcal{B}}^{-s,0,\sigma}_{\infty,\infty}(\Omega))'&=\dot{\B}^{s,\sigma}_{1,1}(\Omega)\, \text{,}&&\text{where }p=q=1\text{.}
\end{array}
\end{align*}
Moreover, when $\Omega=\RR^n_+$ and $1<p<\infty$, the statement remains true for inhomogeneous function spaces. The result still holds if one replaces $(\sigma,\mathfrak{n})$ by $(\gamma,\mathfrak{t})$.
\end{proposition}

\begin{corollary}\label{cor:Weak*densityDivFreeRn+} Let $\alpha\in(0,1)$, $r\in[1,\infty]$, and $\Omega$ be either $\RR^n_+$ or a bounded Lipschitz domain of class $\mathcal{M}^{1+\alpha,r}_\W(\epsilon)$ for some sufficiently small $\epsilon>0$. Considering the case of vector fields only, $\Lambda^1\simeq\CC^n$, the space $\Ccinftydiv(\Omega)$ is weakly-$\ast$ dense in 
\begin{enumerate}
    \item $\dot{\B}^{s,\sigma}_{p,\infty,\mathfrak{n}}(\Omega)$,  $p\in[1,\infty]$, $s\in(-1+\sfrac{1}{p},\sfrac{1}{p})$;
    \item $\dot{\B}^{s,\sigma}_{\infty,q,\mathfrak{n}}(\Omega)$, $q\in[1,\infty]$, $s\in(-1,0)$.
\end{enumerate}
\end{corollary}

\newpage
\section{The Stokes--Dirichlet problem on the half-space}\label{Sec:StokesHalfSpace}

From now on, and until the end of the present work, all function spaces are either scalar-valued or vector fields-valued, \textit{i.e.} over $\Lambda^0\simeq\CC$ and $\Lambda^1\simeq\CC^n$, unless the contrary is explicitly stated. 

\begin{notation}
In this section, for $p,q,r,\kappa\in[1,\infty]$, $s,\uptau\in\RR$, we set $\dot{\X}^{s,p}$ to be any of the following normed spaces
\begin{itemize}
    \item $\dot{\H}^{s,p}$, for $1 <p< \infty$, $s\in\RR$; or
    \item $\dot{\B}^{s}_{p,q}$, for $1\leqslant p,q \leqslant \infty$, $s\in\RR$.
\end{itemize}
Similarly, we may consider $\dot{\Y}^{\uptau,r}$ be any of the above spaces, where we have replaced the indices $(s,p,q)$ and the symbols $\{\dot{\H}^{s,p},\,\dot{\B}^{s}_{p,q}\}$ by $(\uptau,r,\kappa)$ and $\{\dot{\H}^{\uptau,r},\,\dot{\B}^{\uptau}_{r,\kappa}\}$.
\end{notation}
Additionally,  when one writes $\X_{\cdot,\sigma}(\Omega)$ it is always implicit that the space is $\CC^n$-valued, and when $\Omega$ is a special or bounded Lipschitz domain, we set the new notation
\begin{align*}
   {\X}^{s,p}_{\mathcal{D},\sigma}(\Omega):= \begin{cases}{\X}^{s,p}_{\mathfrak{n},\sigma}(\Omega)=\left\{ \,\uu\in{\X}^{s,p}_\sigma(\Omega)\, \,:\,\, \uu\cdot\nu_{|_{\partial\Omega}}=0 \,\right\}  &\text{, if } -1+{\sfrac{1}{p}}<s < {\sfrac{1}{p}}\text{,}\\
    \left\{ \,\uu\in{\X}^{s,p}_\sigma(\Omega)\, \,:\,\, \uu_{|_{\partial\Omega}}=0 \,\right\} &\text{, if }s > {\sfrac{1}{p}}\text{,}
    \end{cases}
\end{align*}
so that by Proposition~\ref{prop:IdentifVanishingDivFree}, one always has the canonical identification
\begin{align*}
   {\X}^{s,p}_{\mathcal{D},\sigma}(\Omega)={\X}^{s,p}_{0,\sigma}(\Omega),\qquad \text{ whenever }-1+{\sfrac{1}{p}}<s < 1+{\sfrac{1}{p}},\quad s\neq{\sfrac{1}{p}}.
\end{align*}
When $\Omega$ is a special Lipschitz domain, we keep similar notations for homogeneous function spaces with the same consequences.

\medbreak

Before going to the heart of the matter for our study of the Stokes--Dirichlet problem on several domains-- in the present section $\RR^n_+$, and later on bounded domains with some minimal regularity assumptions--, we present it in the general $\L^2$-setting over Lipschitz domains.

\subsection{Presentation of the Stokes--Dirichlet problem on Lipschitz domains.}\label{subsec:GeneralStokesDirPbLipDomains}

Let $\Omega\subset\RR^n$, $n\geqslant 2$, be a Lipschitz domain. 
 We introduce the Stokes--Dirichlet operator drawing inspiration from \cite[Section~4]{MitreaMonniaux2008}, \cite[Section~2]{MonniauxShen2018}. On bounded Lipschitz domains, we endow $\H^{1,2}_{0}(\Omega)$ and $\H^{1,2}_{0,\sigma}(\Omega)$ with the norm $\lVert \cdot \rVert_{\dot{\H}^{1,2}}$. When $\Omega$ is a special Lipschitz domain, we endow respectively ${\H}^{1,2}_{0}(\Omega)$, ${\H}^{1,2}_{0,\sigma}(\Omega)$,  and $\dot{\H}^{1,2}_{0}(\Omega)$, $\dot{\H}^{1,2}_{0,\sigma}(\Omega)$ with the norms $\lVert \cdot \rVert_{{\H}^{1,2}}$ and $\lVert \cdot \rVert_{\dot{\H}^{1,2}}$. We set
\begin{enumerate}
    \item $\iota\,:\,\L^2_{\mathfrak{n},\sigma}(\Omega)\longrightarrow \L^2(\Omega,\CC^n)$ to be the canonical embedding;
    \item $\PP_{\Omega}\,:\,\L^2(\Omega,\CC^n)\longrightarrow\L^2_{\mathfrak{n},\sigma}(\Omega)$ to be the Leray projection; note that $\iota'=\PP_{\Omega}$ with $\PP_{\Omega}\iota=\mathbf{I}_{\L^2_{\mathfrak{n},\sigma}}$;
    \item $\iota_0\,:\,\H^{1,2}_{0,\sigma}(\Omega)\longrightarrow \H^{1,2}_{0}(\Omega,\CC^n)$ to be the canonical embedding, with the property $\iota_{|_{\H^{1,2}_{0,\sigma}}}=\iota_0$;
    \item $\PP_{1}\,:\,\H^{-1,2}(\Omega,\CC^n)\longrightarrow(\H^{1,2}_{0,\sigma}(\Omega))'$, the projection given by $\PP_1:=\iota_0'$.
\end{enumerate}
Note that when $\Omega$ is a special Lipschitz domain, one can check that $\iota_0$ coincides with the canonical embedding $\iota_0\,:\,\dot{\H}^{1,2}_{0,\sigma}(\Omega)\longrightarrow \dot{\H}^{1,2}_{0}(\Omega,\CC^n)$, so that $\mathbb{P}_1$ restricts and co-restricts as a bounded linear projection $\PP_1\,:\,\dot{\H}^{-1,2}(\Omega,\CC^n)\longrightarrow (\dot{\H}^{1,2}_{0,\sigma}(\Omega))'$.

We warn the reader that in general, with our current choice of notation, the Hilbert space $(\H^{1,2}_{0,\sigma}(\Omega))'$ may not be identifiable with $\H^{-1,2}_{\sigma}(\Omega)$. In general, the elements $(\H^{1,2}_{0,\sigma}(\Omega))'$ may not even be \textit{a priori} identifiable with distributions on $\Omega$, unless the domain has more regularity.

\begin{definition}The Stokes--Dirichlet operator is the unbounded $\L^2_{\mathfrak{n},\sigma}$-realization of the bounded operator $\AA_{\mathcal{D}}\,:\,\H^{1,2}_{0,\sigma}(\Omega)\longrightarrow (\H^{1,2}_{0,\sigma}(\Omega))'$, denoted by $(\D_2(\AA_{\mathcal{D}}),\AA_{\mathcal{D}})$, and induced by the sesquilinear form
\begin{align*}
    \mathfrak{a}_{\mathcal{D},\sigma}\,:\,\D(\mathfrak{a}_{\mathcal{D},\sigma})\times \D(\mathfrak{a}_{\mathcal{D},\sigma}) &\longrightarrow \CC\\
    (\uu ,\vv)\quad\qquad&\longmapsto \langle \nabla \iota \uu, \nabla \iota \vv\rangle_{\Omega},
\end{align*}
with form domain $\D(\mathfrak{a}_{\mathcal{D},\sigma}):=\H^{1,2}_{0,\sigma}(\Omega)$.
\end{definition}

\begin{proposition}[ {\cite[Theorem~4.7]{MitreaMonniaux2008},  \cite[Prop.~1]{MonniauxShen2018}} ]\label{Prop:L2StokesDir} The $\L^2_{\mathfrak{n},\sigma}$-realization of the Stokes--Dirichlet operator -- with the description:
\begin{align*}
        \D_2(\AA_{\mathcal{D}}) &= \{\, \uu\in \H^{1,2}_{0,\sigma}(\Omega) \,:\, \PP_1(-\Delta_\mathcal{D})\iota_0 \uu\in \L^2_{\mathfrak{n},\sigma}(\Omega) \,\},\\
        \AA_{\mathcal{D}}&=\PP_1(-\Delta_\mathcal{D})\iota_0 \uu, \quad \uu\in\D_2(\AA_{\mathcal{D}}),
\end{align*}
    where $\Delta_\mathcal{D}\,:\,\H^{1,2}_{0}(\Omega)\longrightarrow \H^{-1,2}(\Omega)$ is the (weak) Dirichlet Laplacian -- enjoys the following properties:
\begin{enumerate}
    \item It is a $0$-sectorial self-adjoint operator on $\L^2_{\mathfrak{n},\sigma}(\Omega)$. More precisely, for $\mu\in[0,\pi)$, for all $\lambda\in\Sigma_\mu$ and for all $\ff\in\L^2_{\mathfrak{n},\sigma}(\Omega)$, there exists a unique $\uu\in\D_2(\AA_{\mathcal{D}})\subset\H^{1,2}_{0,\sigma}(\Omega)$, such that
    \begin{align*}
        \lambda \uu +\AA_\mathcal{D}\uu =\ff\quad \text{ in } \L^2_{\mathfrak{n},\sigma}(\Omega)
    \end{align*}
    with the estimate
    \begin{align*}
        |\lambda|\lVert \uu \rVert_{\L^2(\Omega)} + |\lambda|^{\sfrac{1}{2}}\lVert \nabla \uu \rVert_{\L^2(\Omega)} +   \lVert \AA_\mathcal{D} \uu \rVert_{\L^2(\Omega)} \less_{\mu} \lVert \ff\rVert_{\L^2(\Omega)}.
    \end{align*}
    In particular, $\AA_\mathcal{D}$ generates a uniformly bounded holomorphic $\C_0$-semigroup on $\L^2_{\mathfrak{n},\sigma}(\Omega)$. Furthermore, $\AA_\mathcal{D}$ admits a bounded $\mathrm{\mathbf{H}}^\infty(\Sigma_\theta)$-functional calculus  for any $\theta\in(0,\pi)$.
    \item For all $ \uu \in\D_2(\AA_{\mathcal{D}})$, there exists $\mathfrak{p}\in\L^2_{\mathrm{loc}}(\overline{\Omega})$, such that
    \begin{align*}
        \iota \AA_\mathcal{D} \uu = - \Delta \iota_0 \uu + \nabla \mathfrak{p}\quad\text{ in } \H^{-1,2}(\Omega,\CC^n), 
    \end{align*}
    and the domain has the following description
    \begin{align*}
        \D_2(\AA_{\mathcal{D}})=\{\, \uu\in \H^{1,2}_{0,\sigma}(\Omega) \,:\, \exists \mathfrak{p} \in\L^2_{\mathrm{loc}}(\overline{\Omega}),\, -\Delta \iota_0 \uu + \nabla \mathfrak{p} \in \L^2_{\mathfrak{n},\sigma}(\Omega) \,\}.
    \end{align*}
    \item The operator is injective, and, when $\Omega$ is bounded, the operator is invertible. Its square root has domain $\D_2(\AA_\mathcal{D}^{\sfrac{1}{2}})=\H^{1,2}_{0,\sigma}(\Omega)$, with 
    \begin{align*}
        \lVert \AA_\mathcal{D}^{\sfrac{1}{2}} \uu \rVert_{\L^2(\Omega)}= \lVert \nabla \uu \rVert_{\L^2(\Omega)}, \quad \uu\in\D_2(\AA_\mathcal{D}^{\sfrac{1}{2}}).
    \end{align*}
     This induces an isomorphism $\AA_\mathcal{D}^{\sfrac{1}{2}}\,:\,\H^{1,2}_{0,\sigma}(\Omega)\longrightarrow\L^2_{\mathfrak{n},\sigma}(\Omega)$. ($\AA_\mathcal{D}^{\sfrac{1}{2}}\,:\,\dot{\H}^{1,2}_{0,\sigma}(\Omega)\longrightarrow\L^2_{\mathfrak{n},\sigma}(\Omega)$, when $\Omega$ is a special Lipschitz domain and $n\geqslant3$).
\end{enumerate}
\end{proposition}

In particular, according to Proposition~\ref{Prop:L2StokesDir} above, we have a suitable standard description on the Hilbert space setting $(\H^{1,2}_{0,\sigma}(\Omega),\L^2_{\mathfrak{n},\sigma}(\Omega), (\H^{1,2}_{0,\sigma}(\Omega))')$ to solve uniquely the following partial differential equation, provided $\mu\in(0,\pi)$, for any $\lambda\in\Sigma_\mu$ and any $\ff \in\L^2_{\mathfrak{n},\sigma}(\Omega)$,
\begin{equation*}\tag{DS${}_\lambda$}\label{eq:DirStokesSystemBddLip}
    \left\{ \begin{array}{rllr}
         \lambda \uu - \Delta \uu +\nabla \mathfrak{p} &= \ff \text{, }&&\text{ in } \Omega\text{,}\\
        \div \uu &= 0\text{, } &&\text{ in } \Omega\text{,}\\
        \uu_{|_{\partial\Omega}} &=0\text{, } &&\text{ on } \partial\Omega\text{.}
    \end{array}
    \right.
\end{equation*}
called the Stokes--Dirichlet (resolvent) problem.

\begin{remark}\begin{itemize}
     \item The Stokes--Dirichlet operator \textbf{\textit{is not}} defined as "$\AA_\mathcal{D}=\PP_{\Omega}(-\Delta_\mathcal{D})$":
        
        In general, for the Stokes--Dirichlet and Dirichlet Laplacian defined through their respective sesquilinear forms, their domains are such that
    \begin{align*}
        \D_2(\Delta_\mathcal{D})\cap \L^2_{\mathfrak{n},\sigma}(\Omega) \subset \D_2(\AA_\mathcal{D}) \nsubseteq\D_2(\Delta_\mathcal{D}).
    \end{align*}
    The knowledge of $\D_2(\AA_\mathcal{D}) \subseteq\D_2(\Delta_\mathcal{D})$ is equivalent to know whether the pressure term $\nabla \mathfrak{p}$ itself is an element of $\L^2(\Omega,\CC^n)$, which may fail in general. This is a matter of regularity properties of the Stokes--Dirichlet operator and the Dirichlet Laplacian, which ties to the geometric and regularity properties of the domain $\Omega$.
    \item In the case of unbounded domains, such as a special Lipschitz domain, and more specifically the flat half-space, the Stokes--Dirichlet operator may not be invertible on $\L^2_{\mathfrak{n},\sigma}(\Omega,\CC^n)$ while being injective. Thus, the square-root property holds with homogeneous norms, but the operator is invertible only at the “homogeneous level”. An additional major issue occurs when $n=2$, since, in general, $\dot{\H}^{1,2}(\Omega)$ is not a Banach space.
\end{itemize}
\end{remark}

In this particular setting, and this is standard \textit{e.g.} \cite[Chapter~2,~Section~2.6]{SohrBook2001}, one can even consider data $\ff$ with lower regularity in \eqref{eq:DirStokesSystemBddLip}. Although this is very well known and standard, we provide here a short proof since it will simplify our presentation of a quite tedious strategy in the late Section~\ref{Sec:RemoveC1}.

\begin{lemma}\label{lem:H-1EstBddLip} Let $\Omega$ be a bounded Lipschitz domain. Let $\mu\in(0,\pi)$, for any $\lambda\in\Sigma_\mu\cup\{0\}$, for any $\mathbf{F}\in\L^2(\Omega,\CC^{n^2})$, one has the resolvent estimate
\begin{align}\label{eq:EstA^-1PdivL²}
    (1+|\lambda|)^\frac{1}{2}\lVert (\lambda\I+\AA_\mathcal{D})^{-1}\PP_\Omega\div(\mathbf{F}) \rVert_{\L^2(\Omega)}\lesssim_{\mu,\Omega} \lVert \mathbf{F}\rVert_{\L^2(\Omega)}.
\end{align}
In particular, for any $\ff\in\H^{-1,2}(\Omega,\CC^n)$, the problem \eqref{eq:DirStokesSystemBddLip} admits a unique solution $(\uu,\mathfrak{p})\in\H^{1,2}_{0,\sigma}(\Omega)\times \L^2_\mfree(\Omega,\CC)$, and one has the estimate
\begin{align}\label{eq:H-1EstBddLip}
    (1+|\lambda|)^\frac{1}{2}\lVert \uu \rVert_{\L^2(\Omega)}\lesssim_{\mu,\Omega} \lVert \ff\rVert_{\H^{-1,2}(\Omega)}.
\end{align}
\end{lemma}

\begin{proof}Let $\FF\in\Ccinfty(\overline{\Omega},\CC^{n^2})$, let $\boldsymbol{\varphi}\in\Ccinfty(\Omega,\CC^n)$, by self-adjointness of both the Hodge-Leray projection and of the resolvent operator $(\lambda\I+\AA_\mathcal{D})^{-1}$,
\begin{align*}
    \langle  (\lambda\I+\AA_\mathcal{D})^{-1}\PP_\Omega\div(\mathbf{F}),\,\boldsymbol{\varphi}\rangle_\Omega &= \langle  \PP_{\Omega}(\lambda\I+\AA_\mathcal{D})^{-1}\PP_\Omega\div(\mathbf{F}),\,\boldsymbol{\varphi}\rangle_\Omega\\
     &= \langle  \div(\mathbf{F}),\,\PP_\Omega(\lambda\I+\AA_\mathcal{D})^{-1}\PP_{\Omega}\boldsymbol{\varphi}\rangle_\Omega\\
     &= \langle  \div(\mathbf{F}),\,(\lambda\I+\AA_\mathcal{D})^{-1}\PP_{\Omega}\boldsymbol{\varphi}\rangle_\Omega\\
     &= -\langle  \mathbf{F},\,\nabla(\lambda\I+\AA_\mathcal{D})^{-1}\PP_{\Omega}\boldsymbol{\varphi}\rangle_\Omega.
\end{align*}
Hence, by Cauchy-Schwarz, and Point \textit{(i)} of Proposition~\ref{Prop:L2StokesDir},
\begin{align*}
    |\langle  (\lambda\I+\AA_\mathcal{D})^{-1}\PP_\Omega\div(\mathbf{F}),\,\boldsymbol{\varphi}\rangle_\Omega|&\leqslant \lVert \mathbf{F}\rVert_{\L^2(\Omega)}\lVert \nabla(\lambda\I+\AA_\mathcal{D})^{-1}\PP_{\Omega}\boldsymbol{\varphi}\rVert_{\L^2(\Omega)}\\
    &\lesssim_{\mu,\Omega} \frac{1}{(1+|\lambda|)^\frac{1}{2}}\lVert \mathbf{F}\rVert_{\L^2(\Omega)}\lVert \boldsymbol{\varphi}\rVert_{\L^2(\Omega)}.
\end{align*}
Taking the supremum over $\boldsymbol{\varphi}$ with $\L^2$-norm $1$, and extending above inequality to the whole $\L^2(\Omega,\CC^n)$ by density, \eqref{eq:EstA^-1PdivL²} holds.

\medbreak

\noindent Now, we recall that $\H^{-1,2}(\Omega,\CC^n) = \L^2(\Omega,\CC^n)+ \div(\L^2(\Omega,\CC^{n^2}))$. Let $\ff\in\H^{-1,2}(\Omega,\CC^n)$ and $(\mathbf{a},\mathbf{b})\in\L^2(\Omega,\CC^n)\times\L^2(\Omega,\CC^{n^2})$ such that $\ff=\mathbf{a}+\div(\mathbf{b})$, it holds
\begin{align*}
     \lVert (\lambda\I+\AA_\mathcal{D})^{-1}\PP_\Omega[\mathbf{a}+\div(\mathbf{b})] \rVert_{\L^2(\Omega)} &\leqslant  \lVert (\lambda\I+\AA_\mathcal{D})^{-1}\PP_\Omega\mathbf{a} \rVert_{\L^2(\Omega)} + \lVert (\lambda\I+\AA_\mathcal{D})^{-1}\PP_\Omega\div(\mathbf{b}) \rVert_{\L^2(\Omega)} \\
     &\lesssim_{\mu,\Omega}\frac{1}{(1+|\lambda|)}\lVert\mathbf{a} \rVert_{\L^2(\Omega)} + \frac{1}{(1+|\lambda|)^\frac{1}{2}}\lVert\mathbf{b} \rVert_{\L^2(\Omega)}\\
     &\lesssim_{\mu,\Omega}\frac{1}{(1+|\lambda|)^\frac{1}{2}}\Big(\lVert\mathbf{a} \rVert_{\L^2(\Omega)} + \lVert\mathbf{b} \rVert_{\L^2(\Omega)} \Big).
\end{align*}
Taking the infimum on all such pairs  $(\mathbf{a},\mathbf{b})$ yields \eqref{eq:H-1EstBddLip}.
\end{proof}

\begin{theorem}[ {\cite[Theorem~5.1]{MitreaMonniaux2008}} ]\label{thm:StokesDirichletHsFracPower} Assume $\Omega$ is a bounded Lipschitz domain. For $-\sfrac{1}{2}<s<\sfrac{3}{2}$, $s\neq\sfrac{1}{2}$, we have an isomorphism
\begin{align*}
    \AA_{\mathcal{D}}^{\sfrac{s}{2}}\,:\,\dot{\H}^{s,2}_{0,\sigma}(\Omega)\longrightarrow\L^2_{\mathfrak{n},\sigma}(\Omega).
\end{align*}
In particular, one can realize $\AA_{\mathcal{D}}$ as a densely defined, closed, injective, $0$-sectorial operator on  $\dot{\H}^{s,2}_{0,\sigma}(\Omega)$.

\medbreak

\noindent Furthermore, the result remains true
\begin{itemize}
    \item if $\Omega=\RR^n_+$ and $-\sfrac{1}{2}<s<1$;
    \item if $\Omega$ is special Lipschitz domain and $0\leqslant s<1$,
\end{itemize}
the case $s=1$ being allowed for both if $n\geqslant 3$.
\end{theorem}

\begin{proof} When $\Omega$ is a bounded domain, the result is claimed and proved in \cite[Theorem~5.1]{MitreaMonniaux2008}.

The only remaining case is for $\Omega$ to be a special Lipschitz domain. We recall that, in this case the main issue of this result is the possible lack of completeness of $\dot{\H}^{1,2}_{0,\sigma}(\Omega)$ when $n=2$, and that one cannot take an abstract completion to preserve meaningfulness of various quantities.

\textbf{Step 1:} We show the equivalence of norms between $\lVert\cdot\rVert_{\dot{\H}^{s,2}(\Omega)}$ and $\lVert \AA^{\sfrac{s}{2}}_{\mathcal{D}}\cdot\rVert_{\L^2(\Omega)}$ By the sequilinear form, the definition of the square root in the Hilbertian case, we obtain a continuous mapping 
\begin{align*}
    \lVert\AA_{\mathcal{D}}^{\sfrac{1}{2}} \uu\rVert_{\L^2(\Omega)}^2=\lVert\nabla \uu\rVert_{\L^2(\Omega)}^2,\qquad \AA_{\mathcal{D}}^{\sfrac{1}{2}} \,:\, \dot{\H}^{1,2}_{0,\sigma}(\Omega)\longrightarrow\L^2_{\mathfrak{n},\sigma}(\Omega)
\end{align*}
since $\H^{1,2}_{0,\sigma}$ is dense in $\dot{\H}^{1,2}_{0,\sigma}$. The continuous inverse can only be constructed up to a dense subspace, the operator being injective but $\dot{\H}^{1,2}_{0,\sigma}$ may not be complete. For $s\in[0,1]$, we define the homogeneous domains of fractional powers of the Stokes--Dirichlet operator by
\begin{align*}
    \D_2( \mathring{\AA}^{\sfrac{s}{2}}_{\mathcal{D}}) := \{ \uu \in \L^{2}_{\mathfrak{n},\sigma}(\Omega)+\dot{\H}^{1,2}_{0,\sigma}(\Omega) \,:\, \lVert\uu\rVert_{\eus{D}_{{\AA}_\mathcal{D}}(s,2)}<\infty\}
\end{align*}
where
\begin{align*}
    \lVert \uu \rVert_{\eus{D}_{\AA_{\mathcal{D}}}(s,2)}^2 :&= \int_{0}^{\infty} \lVert t^{2-s}\AA_{\mathcal{D}}e^{-t\AA_{\mathcal{D}}^{\sfrac{1}{2}}} \uu \rVert_{\L^2(\Omega)}^2 \frac{\d t}{t} \\ 
    &= \int_{0}^{\infty} \lVert t^{2}\AA_{\mathcal{D}}e^{-t\AA_{\mathcal{D}}^{\sfrac{1}{2}}}\uu \rVert_{\L^2(\Omega)}^2 \frac{\d t}{{t}^{1+2s}}.
\end{align*}

Thanks to M${}^{\text{c}}$Intosh's theorem \cite[Theorem~8]{McIntosh1986}, we have for all $s\in[0,1]$, $\uu\in \D_2(\AA_{\mathcal{D}}^{\sfrac{s}{2}})$,
\begin{align}\label{eq:proofFracHomdomStokeDir}
    \lVert \uu  \rVert_{\eus{D}_{\AA_{\mathcal{D}}}(s,2)}\sim_{s} \lVert \AA_{\mathcal{D}}^{\sfrac{s}{2}} \uu \rVert_{\L^2(\Omega)}
\end{align}
and we recall that since $\AA_{\mathcal{D}}$ is injective and $\L^2_{\mathfrak{n},\sigma}(\Omega)$ is a Hilbert space, $\D_2(\AA_{\mathcal{D}}^{\sfrac{1}{2}})\cap\R_2(\AA_{\mathcal{D}}^{\sfrac{1}{2}})$ is a dense subspace of $\D_2( \AA_{\mathcal{D}}^{\sfrac{s}{2}})$. And by \cite[Theorem~8]{McIntosh1986}, by the $\mathbf{H}^{\infty}$-functional calculus, one also has the equivalence of norms,
\begin{align*}
    \lVert \uu  \rVert_{\eus{D}_{\AA_{\mathcal{D}}}(s,2)}^2\sim_{s} \int_{0}^{\infty} \lVert t^{1-s}\AA_{\mathcal{D}}^{\sfrac{1}{2}}e^{-t\AA_{\mathcal{D}}^{\sfrac{1}{2}}}\uu \rVert_{\L^2(\Omega)}^2 \frac{\d t}{t}
\end{align*}
for all $\uu\in \L^{2}_{\mathfrak{n},\sigma}(\Omega)+\dot{\H}^{1,2}_{0,\sigma}(\Omega)$. Therefore, by \cite[Proposition~2.12]{DanchinHieberMuchaTolk2020}, Proposition~\ref{prop:IdentifVanishingDivFree} and Theorem~\ref{thm:InterpHomSpacesLip} one obtains that for all $0<s<1$, $s\neq\frac{1}{2}$,
\begin{align}\label{eq:proofInterpHomdomStokeDir}
    \eus{D}_{\AA_{\mathcal{D}}}(s,2) = (\L^{2}_{\mathfrak{n},\sigma}(\Omega),\D_2( \mathring{\AA}^{\sfrac{1}{2}}_{\mathcal{D}}))_{s,2} = (\L^{2}_{\mathfrak{n},\sigma}(\Omega),\dot{\H}^{1,2}_{0,\sigma}(\Omega))_{s,2} = \dot{\H}^{s,2}_{0,\sigma}(\Omega)
\end{align}
with equivalence of norms. Combining \eqref{eq:proofInterpHomdomStokeDir} and \eqref{eq:proofFracHomdomStokeDir}, one obtains,
\begin{align*}
    \lVert \uu \rVert_{\dot{\H}^{s,2}(\Omega)}\sim_{s,n}\lVert \AA^{\sfrac{s}{2}}_{\mathcal{D}}\uu\rVert_{\L^2(\Omega)}
\end{align*}
for all $\uu\in\dot{\H}^{s,2}_{0,\sigma}(\Omega)$, provided $s\neq \frac{1}{2}$.

Finally, we show that for all $\alpha\in[0,1]$, $\alpha\neq1/2$, the canonical embedding
\begin{align*}
    \iota \,:\, \dot{\H}^{\alpha,2}_{0,\sigma}(\Omega)\longrightarrow\eus{D}_{\AA_\mathcal{D}}(\alpha,2).
\end{align*}
is actually an isomorphism as soon as $\alpha<1$. Hence, we need to show that $\AA_{\mathcal{D}}^{\sfrac{\alpha}{2}}$ is invertible.

We introduce the densely defined (holomorphic) family of operators, provided $\Re(z)\in[0,1]$,
\begin{align*}
    \AA^{-\sfrac{z}{2}}_{\mathcal{D}}\,:\, \mathrm{D}_2(\AA^{\sfrac{1}{2}}_{\mathcal{D}})\cap\mathrm{R}_2(\AA^{\sfrac{1}{2}}_{\mathcal{D}})\longrightarrow \L^2_{\mathfrak{n},\sigma}(\Omega)+\dot{\H}^{1,2}_{0,\sigma}(\Omega).
\end{align*}
By the mapping property $\AA^{\sfrac{1}{2}}_{\mathcal{D}}\,:\,\dot{\H}^{1,2}_{0,\sigma}(\Omega)\longrightarrow\L^2_{\mathfrak{n},\sigma}(\Omega)$, $\text{M}^{\text{c}}$Intosh's Theorem, and the fact that $\AA_{\mathcal{D}}$ has BIP on $\L^2_{\mathfrak{n},\sigma}(\Omega)$, we obtain 
\begin{align*}
    \lVert \AA_{\mathcal{D}}^{-\frac{1+i\tau}{2}}\uu  \rVert_{\dot{\H}^{1,2}_0(\Omega)}^2 \sim  \int_{0}^{\infty}\lVert t\AA_{\mathcal{D}}^{\frac{1}{2}-\frac{i\tau}{2}}e^{-t\AA_{\mathcal{D}}^{\frac{1}{2}}}\uu\rVert_{\L^2(\Omega)}^2 \frac{\d t}{t} &\lesssim_{\tau} \lVert \uu  \rVert_{\L^2(\Omega)}^2,\\
    \lVert \AA_{\mathcal{D}}^{-\frac{i\tau}{2}}\uu \rVert_{{\L}^{2}(\Omega)}^2 \sim   \int_{0}^{\infty}\lVert t^{2}\AA_{\mathcal{D}}^{1-\frac{i\tau}{2}}e^{-t\AA_{\mathcal{D}}}\uu \rVert_{\L^2(\Omega)}^2 \frac{\d t}{t} &\lesssim_{\tau} \lVert \uu  \rVert_{\L^2(\Omega)}^2.
\end{align*}
For $\Re(z)\in[0,1]$, we introduce the map
\begin{align*}
    T_{z} := (-\Delta)^{\sfrac{z}{2}}\mathcal{E}_0\AA^{-\sfrac{z}{2}}_{\mathcal{D}}
\end{align*}
where $\mathcal{E}_0$ is the extension operator from $\Omega$ to the whole space by $0$. For $k=0,1$, $\tau\in\RR$, the linear operator $T_{k+i\tau}$ is well-defined and bounded from $\L^2_{\mathfrak{n,\sigma}}(\Omega)$ to $\L^2_{\sigma}(\RR^n)$ by Proposition~\ref{prop:IdentifVanishingDivFree}.
Therefore, by Stein's complex interpolation, and by density, we have a well-defined bounded operator
\begin{align*}
    \AA_{\mathcal{D}}^{-\frac{\alpha}{2}}\,:\,\L^{2}_{\mathfrak{n},\sigma}(\Omega)\longrightarrow \dot{\H}^{\alpha,2}_{0,\sigma}(\Omega)
\end{align*}
for all $\alpha\in[0,1)$, $\alpha\neq 1/2$, since $\dot{\H}^{\alpha,2}_{0,\sigma}(\Omega)$ is complete. This completes the proof of Theorem~\ref{thm:StokesDirichletHsFracPower}, since the case $s\in(-\sfrac{1}{2},0)$ for $\Omega=\RR^n_+$ holds by duality.
\end{proof}

\subsection{New  estimates for the Stokes--Dirichlet resolvent problem on the half-space}

During the last few decades for more than fifty years, people did analyze a lot the linear Stokes equations on the half-space $\RR^n_+$. Among all these studies, people were able to produce many representation formula, thanks to Fourier analysis. Here, we are going to take advantage of the notable contribution by Uka\"{i} \cite{Ukai1987}.

Before the next definition, we recall the notations for Riesz Transforms, that is
\begin{align*}
    R=\nabla(-\Delta)^{-\frac{1}{2}},\quad R'=\nabla'(-\Delta)^{-\frac{1}{2}} \quad\text{ and }\quad S=\nabla'(-\Delta')^{-\frac{1}{2}}.
\end{align*}

\begin{definition}\label{def:UkaiOP} We define the Uka\"{i} operators for the velocity field of the Stokes--Dirichlet problem on the flat half-space according to \cite[Eq. (1.2)--(1.5),~(1.10)~\&~Thm~1.1]{Ukai1987}:
\begin{align*}
    \mathrm{E}_{0,\sigma}^{\mathcal{U}} &:= \begin{bmatrix}\mathbf{I}_{n-1}& {S}\\ \prescript{t}{}{(-S)} & \I\end{bmatrix} [\mathrm{E}_{\mathcal{D}}\otimes\mathbf{I}_n]\\ \text{ and }\qquad \Gamma_{\sigma}^{\mathcal{U}}&:= \begin{bmatrix}\mathbf{I}_{n-1}& (-S)\\ 0 & \I\end{bmatrix}\begin{bmatrix}\mathbf{I}_{n-1}& 0\\ 0 & (R'\cdot S)(R'\cdot S+ R_n) \end{bmatrix}\begin{bmatrix}\mathbf{I}_{n-1}& 0\\ 0 & \car_{\RR^n_+} \end{bmatrix},
\end{align*}
The operators for the pressure are given by
\begin{align*}
    \Xi_{\gamma}^{\mathcal{U}}:=(-\Delta_{\mathcal{D},\partial})^{-1}\mathcal{T}_{\partial\RR^n_+}\partial_{x_n}\text{ and } \Pi_{\gamma}^{\mathcal{U}}:=\mathrm{E}_{\mathcal{D}}\prescript{t}{}{(-S+\I\otimes\mathfrak{e}_n)}.
\end{align*}
\end{definition}
Recall here that $\mathcal{T}_{\partial\RR^n_+}$ stands for the trace operator.

\medbreak

The main purpose of Uka\"{i} to introduce such operators, was to give an explicit formula for the Stokes--Dirichlet semigroup on $\L^2_{\mathfrak{n},\sigma}(\RR^n_+)$, which can be summarized as
\begin{align*}
    e^{-t\AA_\mathcal{D}}f = [\Gamma_{\sigma}^{\mathcal{U}}\,e^{t\Delta} \,\mathrm{E}_{0,\sigma}^{\mathcal{U}}f]_{\RR^n_+},\qquad f\in\L^2_{\mathfrak{n},\sigma}(\RR^n_+).
\end{align*}
Thanks to the analysis performed in \cite{Ukai1987}-- $\Gamma_{\sigma}^{\mathcal{U}}$ and $\mathrm{E}_{0,\sigma}^{\mathcal{U}}$ being made of (partial) Riesz transform and rough cut-off-- this allowed Uka\"{i} to extend easily, and then recover, in a simpler way, the boundedness and regularity properties of
\begin{align*}
    (e^{-t\AA_\mathcal{D}})_{t\geqslant 0}\qquad \text{on }\L^p_{\mathfrak{n},\sigma}(\RR^n_+)\text{, }1<p<\infty.
\end{align*}

Our point of view here is about to say that, with this strategy in mind, we can go even further and reach directly, without any major other tools than suitable function space theory,  all (homogeneous) Sobolev spaces of fractional orders $\dot{\H}^{s,p}$, $1<p<\infty$, including also Besov spaces $\dot{\B}^{s}_{p,q}$, even the endpoint ones $p,q=1,\infty$, also encompassing $\dot{\W}^{s,1}=\dot{\B}^{s}_{1,1}$ and $\dot{\C}^{s}_{ub}=\dot{\B}^{s}_{\infty,\infty}$, $s\notin\NN$. Instead of dealing with the semigroup, we deal with the resolvent operator, which is equivalent, but allows more flexibility when it comes to localisation techniques.

The reachability of fractional Sobolev and Besov spaces is possible thanks to the next lemma.

\begin{lemma}\label{lem:UkaiOperators}Let $p,q\in[1,\infty]$ and $-2+{\sfrac{1}{p}}<s<{\sfrac{1}{p}}$. The following Ukaï operators are well-defined and bounded
\begin{align*}
    \mathrm{E}_{0,\sigma}^{\mathcal{U}}&\,:\, \dot{\X}^{s,p}(\RR^n_+,\CC^n) \longrightarrow \dot{\X}^{s,p} (\RR^n,\CC^n),\\
     \Gamma_{\sigma}^{\mathcal{U}}&\,:\, \dot{\X}^{s,p}(\RR^n, \CC^n ) \longrightarrow \dot{\X}^{s,p} (\RR^n, \CC^n ),\\
     \nabla\Xi_{\gamma}^{\mathcal{U}}&\,:\, \dot{\X}^{s+3,p}(\RR^n_+,\CC) \longrightarrow \dot{\X}^{s+1,p} (\RR^n_+,\CC),\\
     \Pi_{\gamma}^{\mathcal{U}}&\,:\,\dot{\X}^{s,p} (\RR^n_+,\CC^n)\longrightarrow \dot{\X}^{s,p} (\RR^n_+,\CC^n).
\end{align*}

Furthermore, for all $\uu\in\dot{\X}^{s,p}\cap\dot{\X}^{s+1,p}_{\mathcal{D},\sigma}(\RR^n_+,\CC^n)$, one has $\mathrm{E}_{0,\sigma}^{\mathcal{U}}\uu$, $\Gamma_{\sigma}^{\mathcal{U}}\mathrm{E}_{0,\sigma}^{\mathcal{U}}\uu\in\dot{\X}^{s,p}\cap\dot{\X}^{s+1,p}(\RR^n,\CC^n)$, with the estimates
\begin{align*}
    \lVert  \mathrm{E}_{0,\sigma}^{\mathcal{U}}\uu\rVert_{\dot{\X}^{s+k,p}(\RR^n)} +  \lVert\Gamma_{\sigma}^{\mathcal{U}}\mathrm{E}_{0,\sigma}^{\mathcal{U}}\uu\rVert_{\dot{\X}^{s+k,p}(\RR^n)} \lesssim_{p,s,n}^{k} \lVert  \uu\rVert_{\dot{\X}^{s+k,p}(\RR^n_+)}\text{, } k\in\{0,1\}\text{.}
\end{align*}
with $[\Gamma_{\sigma}^{\mathcal{U}}\mathrm{E}_{0,\sigma}^{\mathcal{U}}\uu]_{|_{\RR^n_+}}=\uu$.
\end{lemma}

\begin{proof} The boundedness of the operators $\mathrm{E}_{0,\sigma}^{\mathcal{U}}$, $\Gamma_{\sigma}^{\mathcal{U}}$ and  $\Pi_{\gamma}^{\mathcal{U}}$ when $-2+{\sfrac{1}{p}}<s<{\sfrac{1}{p}}$ is a consequence of Lemmas~\ref{lem:ExtDirNeuRn+} and~\ref{lem:ExtOpNegativeBesovSpaces} and Theorem~\ref{thm:RieszTransfRn}. The fact that $\uu\mapsto[\Gamma_{\sigma}^{\mathcal{U}}\mathrm{E}_{0,\sigma}^{\mathcal{U}}\uu]_{|_{\RR^n_+}}$ preserves divergence-free functions with zero (partial) trace follows from the analysis performed by Uka\"{i} in \cite{Ukai1987}. Higher order estimates of order $1$ are also given by Lemmas~\ref{lem:ExtDirNeuRn+}. 
\end{proof}

We now prove one of the main theorems. We reprove well-definedness of the formula starting from the $\L^2$ setting and deduce all the other results from that. We think it is worth mentioning the previous recent iteration for the analysis of the resolvent problem by Watanabe \cite{Watanabe2025}, but only the divergence-free and reflexive case was considered, that is, their analysis was restricted to Besov spaces $\dot{\B}^{s,\sigma}_{p,1,\mathfrak{n}}(\RR^n_+)$, $1<p<\infty$. We also mention the recent partial analysis by Geng and Shen in the non-divergence free case where only $\L^q$-spaces are concerned \cite{GengShen2024}.  Dealing directly with the evolution problem on $\dot{\B}^{s}_{p,1}(\RR^n_+,\CC^n)$, $1<p<\infty$, we also mention the work of Danchin and Mucha \cite{DanchinMucha2009,DanchinMucha2015}. Here we provide a full result with a prescribed non-zero divergence, and non-zero boundary data.

\begin{theorem}\label{thm:MetaThmDirichletStokesRn+}Let $p,q\in[1,\infty]$, $s\in (-1+\sfrac{1}{p},\sfrac{1}{p})$. Let $\mu\in[0,\pi)$ and $\lambda\in\Sigma_\mu$.
\begin{enumerate}[label=($\roman*$)]
    \item Provided $\ff\in {\dot{\X}}^{s,p}(\RR^n_+,\CC^n)$, $g\in[\dot{\X}^{s-1,p}_0\cap\dot{\X}^{s+1,p}](\RR^n_+,\CC)$, we assume either $\hh\in \dot{\X}_{\partial}^{s-1/p,p}\cap\dot{\X}_{\partial}^{s+2-1/p,p}(\partial\RR^n_+,\CC^{n})$,  or $\hh:=\mathbf{H}_{|_{\partial\RR^n_+}}$ for $\mathbf{H}\in[\dot{\X}^{s,p}\cap\dot{\X}^{s+2,p}](\RR^n_+,\CC^n)$, then the Stokes--Dirichlet resolvent problem
    \begin{equation*}\tag{DS${}_\lambda$}\label{eq:DirStokesSystemtn+}
    \left\{ \begin{array}{rllr}
         \lambda \uu - \Delta \uu +\nabla \mathfrak{p} &= \ff \text{, }&&\text{ in } \RR^n_+\text{,}\\
        \div \uu &= g\text{, } &&\text{ in } \RR^n_+\text{,}\\
        \uu_{|_{\partial\RR^n_+}} &=\hh\text{, } &&\text{ on } \partial\RR^n_+\text{.}
    \end{array}
    \right.
\end{equation*}
         admits a unique solution $(\uu,\nabla \mathfrak{p})\in {\dot{\X}}^{s,p}\cap{\dot{\X}}^{s+2,p}(\RR^n_+,\CC^n) \times {\dot{\X}}^{s,p}(\RR^n_+,\CC^n)$ with  the estimates,
        \begin{align}\label{eq:FullResolvEstRn+}
            \lvert\lambda\rvert\lVert  \uu\rVert_{{\dot{\X}}^{s,p}(\RR^n_+)}+&\lVert \nabla^2 \uu\rVert_{{\dot{\X}}^{s,p}(\RR^n_+)} + \lVert \nabla \mathfrak{p}\rVert_{{\dot{\X}}^{s,p}(\RR^n_+)} \nonumber\\&\lesssim_{p,n,s,\mu} \lVert \ff\rVert_{{\dot{\X}}^{s,p}(\RR^n_+)} + \lvert\lambda\rvert\lVert  g\rVert_{{\dot{\X}}^{s-1,p}(\RR^n_+)}+\lVert \nabla g\rVert_{{\dot{\X}}^{s,p}(\RR^n_+)}\\
            &\qquad + \min \begin{cases}
  \lvert\lambda\rvert\lVert  \mathbf{H}\rVert_{{\dot{\X}}^{s,p}(\RR^n_+)}+\lVert \mathbf{H}\rVert_{{\dot{\X}}^{s+2,p}(\RR^n_+)}\\    
 \lvert\lambda\rvert\lVert  \hh\rVert_{{\dot{\X}}_\partial^{s-1/p,p}(\partial\RR^n_+)}+\lVert \hh\rVert_{{\dot{\X}}_\partial^{s+2-1/p,p}(\partial\RR^n_+)} 
\end{cases}.\nonumber
        \end{align}
        
        \item  The operator $(\dot{\D}^{s,p}(\AA_\mathcal{D},\RR^n_+),\AA_{\mathcal{D}})$ is an injective $0$-sectorial operator on $\dot{\X}^{s,p}_{\mathfrak{n},\sigma}(\RR^n_+)$ with domain
        \begin{align*}
            \dot{\D}^{s,p}(\AA_\mathcal{D},\RR^n_+) = \dot{\X}^{s,p}_{\mathfrak{n},\sigma}\cap\dot{\X}^{s+2,p}_{\mathcal{D},\sigma}(\RR^n_+) = \dot{\X}^{s,p}_{\mathfrak{n},\sigma}\cap\dot{\X}^{s+1,p}_{0}(\RR^n_+)\cap\dot{\X}^{s+2,p}(\RR^n_+)
        \end{align*}
        and expression
        \begin{align*}
            \AA_\mathcal{D}\uu = \PP_{\RR^n_+}(-\Delta_{\mathcal{D}}\uu),\qquad \forall \uu \in \dot{\D}^{s,p}(\AA_\mathcal{D},\RR^n_+).
        \end{align*}
        
        Moreover, the following resolvent identity holds true on $\dot{\X}^{s,p}_{\mathfrak{n},\sigma}(\RR^n_+)$,
        \begin{align}\label{eq:ResolventEqualityUkai}
             (\lambda \I + \AA_{\mathcal{D}})^{-1} = \R_{\RR^n_+}\Gamma_{\sigma}^{\mathcal{U}}(\lambda \I - \Delta)^{-1}\E_{0,\sigma}^{\mathcal{U}}\text{.}
        \end{align}
\end{enumerate}
\end{theorem}

\begin{proof} \textbf{Step 1:} We start with $\textit{(i)}$, assuming first $\ff\in\dot{\X}^{s,p}_{\mathcal{D},\sigma}(\RR^n_+)$, $-1+{\sfrac{1}{p}}<s<{\sfrac{1}{p}}$, and $g=0$ and $\hh=0$. Indeed, by linearity of the problem, continuity of the Leray projection, and up to correct the pressure term, one can always replace $\ff$ by $\PP_{\RR^n_+}\ff$ and then assume $\div \ff =0$ and $\nu\cdot \ff_{|_{\partial\RR^n_+}} =0$.

\medbreak

\textbf{Step 1.1:}  A short investigation of the case $p=2$, $s=0$. 

By Proposition \ref{Prop:L2StokesDir}, there exists a unique solution $(\uu,\nabla\mathfrak{p})\in\H^{1,2}_{0,\sigma}(\RR^n_+)\times\H^{-1,2}(\RR^n_+,\CC^n)$. 

However, the analysis performed in \cite{Ukai1987} shows that necessarily, by uniqueness,
\begin{align*}
    \uu = [\Gamma_{\sigma}^{\mathcal{U}}(\lambda \I - \Delta)^{-1}\mathrm{E}_{0,\sigma}^{\mathcal{U}} \ff]_{|_{\RR^n_+}} \text{ and } \nabla \mathfrak{p}= \nabla\Xi_{\gamma}^{\mathcal{U}}(\lambda \I - \Delta)^{-1}\Pi_{\gamma}^{\mathcal{U}} \ff.
\end{align*}
Having, these formulas in mind, we show that $\uu\in\H^{2,2}(\RR^n_+,\CC^n)$.

Notice that, the operators $\Gamma_{\sigma}^{\mathcal{U}}$, $(\lambda \I - \Delta)^{-1}$, $\mathrm{E}_{0,\sigma}^{\mathcal{U}}$ are bounded on $\L^2$, in a way that
\begin{align*}
    \lVert \uu\rVert_{\L^2(\RR^n_+)} \leqslant \lVert \Gamma_{\sigma}^{\mathcal{U}}(\lambda \I - \Delta)^{-1}\mathrm{E}_{0,\sigma}^{\mathcal{U}} \ff \rVert_{\L^2(\RR^n)}&\lesssim_{n} \lVert 
 (\lambda \I - \Delta)^{-1}\mathrm{E}_{0,\sigma}^{\mathcal{U}} \ff \rVert_{\L^2(\RR^n)}\\
 &\less_{\mu,n} \frac{1}{|\lambda|}\lVert 
\mathrm{E}_{0,\sigma}^{\mathcal{U}} \ff \rVert_{\L^2(\RR^n)} \less_{\mu,n} \frac{1}{|\lambda|}\lVert 
\ff \rVert_{\L^2(\RR^n_+)}.
\end{align*}

To prove $\uu\in\H^{2,2}(\RR^n_+,\CC^n)$, since $ \uu = [\Gamma_{\sigma}^{\mathcal{U}}(\lambda \I - \Delta)^{-1}\mathrm{E}_{0,\sigma}^{\mathcal{U}} \ff]_{|_{\RR^n_+}}$, it suffices to check that
\begin{align*}
    \Delta\Gamma_{\sigma}^{\mathcal{U}}(\lambda \I - \Delta)^{-1}\mathrm{E}_{0,\sigma}^{\mathcal{U}} \ff \in\L^2(\RR^n,\CC^n).
\end{align*}
We proceed by density, assuming $\ff\in\Ccinftydiv(\RR^n_+)$ for simplicity. Since $\Delta'$ commutes with $\Gamma_{\sigma}^{\mathcal{U}}$, it reduces to show that
\begin{align*}
    \partial_{x_n}^2[(R'\cdot S)&(R'\cdot S+R_n)]\mathbbm{1}_{\RR^n_+}(\lambda \I - \Delta)^{-1}[\mathrm{E}_{0,\sigma}^{\mathcal{U}} \ff]_n \\ &=  \partial_{x_n}^2[(R'\cdot S)(R'\cdot S+R_n)]\mathbbm{1}_{\RR^n_+}(\lambda \I - \Delta)^{-1}[ \prescript{t}{}{(-S)}\cdot\E_\mathcal{D}\ff'+\E_\mathcal{D}\ff_n] \in\L^2(\RR^n,\CC).
\end{align*}
Since $\div(\ff)=0$,
\begin{align}
    \partial_{x_n}^2&[(R'\cdot S)(R'\cdot S+R_n)]\mathbbm{1}_{\RR^n_+}(\lambda \I - \Delta)^{-1}[ \prescript{t}{}{(-S)}\cdot \E_\mathcal{D}\ff'+ \E_\mathcal{D}\ff_n] \label{eq:ProofCalculationsUkaiOPStokesEstRn+}\\
    &=\partial_{x_n}^2\left[\frac{(-\Delta')^\frac{1}{2}}{(-\Delta)^\frac{1}{2}}\Bigg(\frac{(-\Delta')^\frac{1}{2}}{(-\Delta)^\frac{1}{2}}+\frac{-\partial_{x_n}}{(-\Delta)^\frac{1}{2}}\Bigg)\right]\mathbbm{1}_{\RR^n_+}(\lambda \I - \Delta)^{-1}\left[ \frac{\partial_{x_n}}{(-\Delta')^\frac{1}{2}} \E_\mathcal{N}\ff_n+ \E_\mathcal{D}\ff_n\right]\nonumber\\
    &= \partial_{x_n}^2\left[\frac{1}{(-\Delta)^\frac{1}{2}}\Bigg(\frac{(-\Delta')^\frac{1}{2}}{(-\Delta)^\frac{1}{2}}+\frac{-\partial_{x_n}}{(-\Delta)^\frac{1}{2}}\Bigg)\right]\mathbbm{1}_{\RR^n_+}\partial_{x_n}(\lambda \I - \Delta)^{-1}\E_\mathcal{N}\ff_n\nonumber\\ &\qquad  \partial_{x_n}^2\left[\frac{(-\Delta')^\frac{1}{2}}{(-\Delta)^\frac{1}{2}}\Bigg(\frac{(-\Delta')^\frac{1}{2}}{(-\Delta)^\frac{1}{2}}+\frac{-\partial_{x_n}}{(-\Delta)^\frac{1}{2}}\Bigg)\right]\mathbbm{1}_{\RR^n_+}(\lambda \I - \Delta)^{-1}\E_\mathcal{D}\ff_n\nonumber\\
    &= \left[\frac{\partial_{x_n}}{(-\Delta)^\frac{1}{2}}\Bigg(\frac{(-\Delta')^\frac{1}{2}}{(-\Delta)^\frac{1}{2}}+\frac{-\partial_{x_n}}{(-\Delta)^\frac{1}{2}}\Bigg)\right] \partial_{x_n}\mathbbm{1}_{\RR^n_+}\partial_{x_n}(\lambda \I - \Delta)^{-1}\E_\mathcal{N}\ff_n\nonumber\\ &\qquad  +\left[\frac{(-\Delta')^\frac{1}{2}}{(-\Delta)^\frac{1}{2}}\Bigg(\frac{(-\Delta')^\frac{1}{2}}{(-\Delta)^\frac{1}{2}}+\frac{-\partial_{x_n}}{(-\Delta)^\frac{1}{2}}\Bigg)\right]\partial_{x_n}^2\mathbbm{1}_{\RR^n_+}(\lambda \I - \Delta)^{-1}\E_\mathcal{D}\ff_n\nonumber\nonumber\\
    &= [R_n(R'\cdot S+R_n)] \partial_{x_n}\mathbbm{1}_{\RR^n_+}\partial_{x_n}(\lambda \I - \Delta)^{-1}\E_\mathcal{N}\ff_n\nonumber\\ &\qquad  +[(R'\cdot S)(R'\cdot S+R_n)]\partial_{x_n}^2\mathbbm{1}_{\RR^n_+}(\lambda \I - \Delta)^{-1}\E_\mathcal{D}\ff_n\nonumber\\
    &=[R_n(R'\cdot S+R_n)] \mathbbm{1}_{\RR^n_+}\partial_{x_n}^2(\lambda \I - \Delta)^{-1}\E_\mathcal{N}\ff_n\nonumber\\ &\qquad  +[(R'\cdot S)(R'\cdot S+R_n)]\mathbbm{1}_{\RR^n_+}\partial_{x_n}^2(\lambda \I - \Delta)^{-1}\E_\mathcal{D}\ff_n.\label{eq:ProofCalculationsUkaiOPStokesEstRn+end}
\end{align}
Thus,
\begin{align*}
    \lVert \nabla^2 \uu\rVert_{\L^2(\RR^n_+)} &\leqslant \lVert \nabla^2\Gamma_{\sigma}^{\mathcal{U}}(\lambda \I - \Delta)^{-1}\mathrm{E}_{0,\sigma}^{\mathcal{U}} \ff \rVert_{\L^2(\RR^n)}\\
    &\leqslant \lVert 
 \Delta\Gamma_{\sigma}^{\mathcal{U}} (\lambda \I - \Delta)^{-1}\mathrm{E}_{0,\sigma}^{\mathcal{U}} \ff \rVert_{\L^2(\RR^n)}\\
 &\lesssim_{n} \lVert  \Delta'\Gamma_{\sigma}^{\mathcal{U}} 
 (\lambda \I - \Delta)^{-1}\mathrm{E}_{0,\sigma}^{\mathcal{U}} \ff \rVert_{\L^2(\RR^n)} \\
 &\qquad\qquad+ \lVert \partial_{x_n}^2 \Gamma_{\sigma}^{\mathcal{U}} 
 (\lambda \I - \Delta)^{-1}\mathrm{E}_{0,\sigma}^{\mathcal{U}} \ff \rVert_{\L^2(\RR^n)} \\
 &\lesssim_{n} \lVert  \Gamma_{\sigma}^{\mathcal{U}} 
 \Delta'(\lambda \I - \Delta)^{-1}\mathrm{E}_{0,\sigma}^{\mathcal{U}} \ff \rVert_{\L^2(\RR^n)} \\
 &\qquad\qquad+ \lVert \partial_{x_n}^2 
 (\lambda \I - \Delta)^{-1}(\mathrm{E}_{0,\sigma}^{\mathcal{U}} \ff)'\rVert_{\L^2(\RR^n)} \\
 &\qquad\qquad+ \lVert  \partial_{x_n}^2 S [(R'\cdot S)(R'\cdot S+R_n)]\mathbbm{1}_{\RR^n_+}(\lambda \I - \Delta)^{-1}[\mathrm{E}_{0,\sigma}^{\mathcal{U}} \ff]_n\rVert_{\L^2(\RR^n)} \\
 &\qquad\qquad+ \lVert  \partial_{x_n}^2 [(R'\cdot S)(R'\cdot S+R_n)]\mathbbm{1}_{\RR^n_+}(\lambda \I - \Delta)^{-1}[\mathrm{E}_{0,\sigma}^{\mathcal{U}} \ff]_n\rVert_{\L^2(\RR^n)} \\
 &\lesssim_{\mu,n} \lVert 
(\mathrm{E}_{0,\sigma}^{\mathcal{U}} \ff)' \rVert_{\L^2(\RR^n)} +\lVert 
\mathrm{E}_{\mathcal{D}} \ff_n \rVert_{\L^2(\RR^n)} +\lVert 
\mathrm{E}_{\mathcal{N}} \ff_n \rVert_{\L^2(\RR^n)} \\
&\lesssim_{\mu,n} \lVert 
\ff \rVert_{\L^2(\RR^n_+)}.
\end{align*}
The transition from the antepenultimate estimate to the penultimate one did take advantage of preceding calculations \eqref{eq:ProofCalculationsUkaiOPStokesEstRn+}--\eqref{eq:ProofCalculationsUkaiOPStokesEstRn+end}.

In particular, \eqref{eq:ResolventEqualityUkai} holds true on $\L^2_{\mathfrak{n},\sigma}(\RR^n_+)$. The estimate on the pressure follows from the equality $\nabla\mathfrak{p}=-\lambda \uu +\Delta \uu+\ff$, and the fact that $\Delta \uu+\ff\in\L^2(\RR^n_+,\CC^n)$, since $\uu\in\H^{2,2}(\RR^{n}_+,\CC^n)$.

\medbreak

\textbf{Step 1.2:} Extrapolation to other function spaces. We perform the extrapolation procedure for the spaces 
\begin{equation}\label{eq:StokesDirRn+ProofListSmoothSpaces}
    \dot{\B}^{s,\sigma}_{p,q,\mathfrak{n}}\text{, } \dot{\BesSmo}^{s,\sigma}_{p,\infty,\mathfrak{n}}\text{, } 1\leqslant p <\infty\text{, } \dot{\B}^{s,0,\sigma}_{\infty,q,\mathfrak{n}}\text{, } \dot{\BesSmo}^{s,0,\sigma}_{\infty,\infty,\mathfrak{n}}\text{, } 1\leqslant q<\infty\text{, and } \dot{\H}^{s,p}_{\mathfrak{n},\sigma}\text{, } 1<p<\infty.
\end{equation}
To fix the ideas, we perform the case $\dot{\B}^{s,\sigma}_{p,q,\mathfrak{n}}$, the other ones admit the same proof.  Let $\ff\in\dot{\B}^{s,\sigma}_{p,q,\mathfrak{n}}(\RR^{n}_+)$, then by Lemma \ref{lem:UkaiOperators},
\begin{align*}
    \lVert  \Gamma_{\sigma}^{\mathcal{U}} \lambda(\lambda \I - \Delta)^{-1}\mathrm{E}_{0,\sigma}^{\mathcal{U}} \ff\rVert_{\dot{\B}^{s}_{p,q}(\RR^{n})} \less_{p,s,n,\mu} \lVert \ff \rVert_{\dot{\B}^{s}_{p,q}(\RR^{n}_+)}.
\end{align*}
Now, up to choose $\ff\in\L^2_{\mathfrak{n},\sigma}(\RR^n_+)\cap\dot{\B}^{s,\sigma}_{p,q,\mathfrak{n}}(\RR^{n}_+)$, by --and as in-- the previous Step 1.1, using the definition of function spaces by restriction and Lemma~\ref{lem:UkaiOperators}:
\begin{align*}
    \lVert \nabla^2 \uu\rVert_{\dot{\B}^{s}_{p,q}(\RR^{n}_+)} \leqslant \lVert \nabla^2\Gamma_{\sigma}^{\mathcal{U}}(\lambda \I - \Delta)^{-1}\mathrm{E}_{0,\sigma}^{\mathcal{U}} \ff \rVert_{\dot{\B}^{s}_{p,q}(\RR^{n})}&\less_{p,s,n} \lVert \Delta 
 \Gamma_{\sigma}^{\mathcal{U}} (\lambda \I - \Delta)^{-1}\mathrm{E}_{0,\sigma}^{\mathcal{U}} \ff \rVert_{\dot{\B}^{s}_{p,q}(\RR^{n})}\\
 &\less_{p,s,n,\mu} \lVert 
\ff \rVert_{\dot{\B}^{s}_{p,q}(\RR^{n}_+)}.
\end{align*}

Therefore, we obtain that $\uu\in\H^{2,2}_{\mathcal{D},\sigma}(\RR^n_+)\cap [\dot{\B}^{s,\sigma}_{p,q,\mathfrak{n}}\cap \dot{\B}^{s+2,\sigma}_{p,q,\mathcal{D}}](\RR^n_+)$, thus
the equality $\nabla\mathfrak{p}=-\lambda \uu +\Delta\uu +\ff$, which holds in $\L^2(\RR^n_+,\CC^n)$, implies $\nabla \mathfrak{p}\in \dot{\B}^{s}_{p,q}(\RR^n_+,\CC^n)$, and we recover the estimate
\begin{align*}
    \lVert \nabla \mathfrak{p}\rVert_{\dot{\B}^{s}_{p,q}(\RR^{n}_+)} \less_{p,s,n,\mu} \lVert 
\ff \rVert_{\dot{\B}^{s}_{p,q}(\RR^{n}_+)}.
\end{align*}

We recall that $\L^2_{\mathfrak{n},\sigma}(\RR^n_+)\cap\dot{\B}^{s,\sigma}_{p,q,\mathfrak{n}}(\RR^{n}_+)$ is dense in $\dot{\B}^{s,\sigma}_{p,q,\mathfrak{n}}(\RR^{n}_+)$, since it contains $\Ccinftydiv(\RR^n_+)$, see Theorem~\ref{thm:DivergenceFreeSpacesSpeiclaLipDensity}, therefore the resolvent estimate
\begin{align*}
      |\lambda|\lVert \uu\rVert_{\dot{\B}^{s}_{p,q}(\RR^{n}_+)}+\lVert \nabla^2 \uu\rVert_{\dot{\B}^{s}_{p,q}(\RR^{n}_+)}+\lVert \nabla \mathfrak{p}\rVert_{\dot{\B}^{s}_{p,q}(\RR^{n}_+)} \less_{p,s,n,\mu} \lVert 
\ff \rVert_{\dot{\B}^{s}_{p,q}(\RR^{n}_+)}.
\end{align*}
extends by density to all $\ff\in\dot{\B}^{s,\sigma}_{p,q,\mathfrak{n}}(\RR^{n}_+)$.

Note that during the limiting procedure, approximating $\ff\in\dot{\B}^{s,\sigma}_{p,q,\mathfrak{n}}(\RR^{n}_+)$ by elements in $\L^2_{\mathfrak{n},\sigma}(\RR^n_+)\cap\dot{\B}^{s,\sigma}_{p,q,\mathfrak{n}}(\RR^{n}_+)$, we preserve the divergence free and the Dirichlet boundary conditions for $\uu$.

Finally, for $\ff$ that belongs to one of the spaces from \eqref{eq:StokesDirRn+ProofListSmoothSpaces},  $(\lambda\I +\AA_\mathcal{D})^{-1}\ff$ is a solution to \eqref{eq:DirStokesSystemtn+} (with $g=0$ and $\hh=0$).

The boundedness of the extrapolated resolvent in the case of Besov spaces $\dot{\B}^{s,\sigma}_{p,\infty,\mathfrak{n}}(\RR^n_+)$ is a consequence of real interpolation. The case of Besov spaces $\dot{\B}^{s,\sigma}_{\infty,q,\mathfrak{n}}$, $q\in[1,\infty]$ is deferred to the end of Step~2.2.

\medbreak

\textbf{Step 2:} For $-1+\sfrac{1}{p}<s<\sfrac{1}{p}$, we prove that the  resolvent extended on other function spaces actually solves \eqref{eq:DirStokesSystemtn+} with uniqueness. We also prove the injectivity of the extrapolated operator $\AA_\mathcal{D}$ itself.

\medbreak

\textbf{Step 2.1:} We check the uniqueness for solutions to \eqref{eq:DirStokesSystemtn+} in the class $[\dot{\B}^{s,\sigma}_{p,q,\mathfrak{n}}\cap \dot{\B}^{s+2,\sigma}_{p,q,\mathcal{D}}](\RR^n_+)$ (here $1\leqslant p,q \leqslant \infty$). Let $(\vv,\nabla\mathfrak{q})\in[\dot{\B}^{s,\sigma}_{p,\infty,\mathfrak{n}}\cap \dot{\B}^{s+2,\sigma}_{p,\infty,\mathcal{D}}](\RR^n_+)\times\dot{\B}^{s}_{p,\infty}(\RR^n_+,\CC^n)$ be a solution of
\begin{equation}\label{eq:ProofInjectiveStokesRn+}
    \left\{ \begin{array}{rllr}
         \lambda \vv - \Delta \vv +\nabla \mathfrak{q} &= 0 \text{, }&&\text{ in } \RR^n_+\text{,}\\
        \div \vv &= 0\text{, }&&\text{ in } \RR^n_+\text{,}\\
        \vv_{|_{\partial\RR^n_+}} &=0\text{, } &&\text{ on } \partial\RR^n_+\text{.}
    \end{array}
    \right.
\end{equation}
If $p=2$, $q\in[1,\infty]$, $s\in(-{\sfrac{1}{2}},{\sfrac{1}{2}})$ this follows by Theorem~\ref{thm:StokesDirichletHsFracPower}, duality, then real interpolation. Now, we assume $1\leqslant p < 2 $. First, if $n=2$, by interpolation inequality and Sobolev embeddings, one has $$\vv\in [\dot{\B}^{s}_{p,\infty,\mathfrak{n}}\cap\dot{\B}^{s+2,\sigma}_{p,\infty,\mathcal{D}}](\RR^n_+)\hookrightarrow \dot{\H}^{\sfrac{n}{p}-\sfrac{n}{2},p}_{\mathcal{D},\sigma}\cap\dot{\H}^{\sfrac{n}{p}-\sfrac{n}{2}+1,p}_{\mathcal{D},\sigma}(\RR^n_+)\hookrightarrow {\H}^{1,2}_{0,\sigma}(\RR^n_+).$$
Thus, from the $\L^2$-theory, we obtain $\vv=0$. If $n=2$, it remains to check the case $p=1$, $0<s<1$, which follows from the Sobolev embbeding $\dot{\B}^{s}_{1,\infty}\hookrightarrow \dot{\B}^{0}_{\frac{2}{2-s},\infty}$.

Now, if $n\geqslant 3$, notice that $p\in[1,2)\subseteq[1,n-1)$, implies that $\dot{\B}^{s+1}_{p,\infty}(\RR^n)$ is a complete space, so that the pressure itself (not only its gradient) $\mathfrak{q}$ is well-defined and uniquely determined, since it is necessarily given by $\mathfrak{q}(\cdot,x_n):=(-\Delta')^{-\sfrac{1}{2}}e^{-x_n(-\Delta')^{\sfrac{1}{2}}}\mathcal{T}_{\partial\RR^n_+}[\Delta \uu_n]$ by Proposition~\ref{prop:GenNeuPbRn+}. Hence, by the definition of $\E_{0,\sigma}^\mathcal{U}$ and the fact that $S=\nabla'(-\Delta')^{-\sfrac{1}{2}}$ it can be checked that
\begin{align*}
    \E_{0,\sigma}^\mathcal{U}[\nabla \mathfrak{q}] =0.
\end{align*}
Therefore, applying $\E_{0,\sigma}^\mathcal{U}$ to the first equation of \eqref{eq:ProofInjectiveStokesRn+}, it becomes
\begin{align*}
    \lambda [\E_{0,\sigma}^\mathcal{U} \uu] - [\E_{0,\sigma}^\mathcal{U} \Delta \uu] = 0.
\end{align*}
Since $\uu\in \dot{\B}^{s}_{p,\infty}\cap \dot{\B}^{s+2}_{p,\infty}(\RR^n_+,\CC^n)$ with $\uu_{|_{\partial\RR^n_+}}=0$, by Lemma~\ref{lem:ExtDirNeuRn+} and its proof, we obtain $\E_\mathcal{D}\uu\in \dot{\B}^{s}_{p,\infty}\cap \dot{\B}^{s+2}_{p,\infty}(\RR^n,\CC^n)$,  with the identity
\begin{align*}
    \Delta \E_\mathcal{D}\uu = \E_\mathcal{D}\Delta\uu \text{, in } \S'(\RR^n,\CC^n).
\end{align*}
Thus, by the definition of $\E_{0,\sigma}^\mathcal{U}$ again, we deduce that 
\begin{align*}
    \lambda [\E_{0,\sigma}^\mathcal{U} \uu] -  \Delta [\E_{0,\sigma}^\mathcal{U}  \uu] = 0\text{, in } \S'(\RR^n,\CC^n).
\end{align*}
So that, $\E_{0,\sigma}^\mathcal{U} \uu=0$, and since $[\Gamma_{\sigma}^\mathcal{U}\E_{0,\sigma}^\mathcal{U}\uu]_{|_{\RR^n_+}}=\uu$, see Lemma~\ref{lem:UkaiOperators}, it follows that $\uu=0$.

It follows that we have uniqueness for the problem \eqref{eq:DirStokesSystemtn+} in $\dot{\B}^{s}_{p,\infty}\cap \dot{\B}^{s+2}_{p,\infty}(\RR^n_+,\CC^n)$ for all $p\in[1,2]$, and then in $\dot{\H}^{s,p}\cap \dot{\H}^{s+2,p}(\RR^n_+,\CC^n)$  and in  $\dot{\B}^{s}_{p,q}\cap \dot{\B}^{s+2}_{p,q}(\RR^n_+,\CC^n)$ for all $p\in[1,2]$, $q\in[1,\infty]$.

Hence, for $n\geqslant 2$, it follows that we do have existence and uniqueness in $\dot{\B}^{s,\sigma}_{p,q,\mathfrak{n}}\cap \dot{\B}^{s+2,\sigma}_{p,q,\mathcal{D}}(\RR^n_+,\CC^n)$ for all $p\in[1,2]$, $q\in[1,\infty]$, $s\in(-1+{\sfrac{1}{p}},{\sfrac{1}{p}})$ provided $\ff\in\dot{\B}^{s,\sigma}_{p,q,\mathfrak{n}}(\RR^n_+)$, $g=0$, $\hh=0$. If $p\in(2,\infty]$, $q\in[1,\infty]$, since the Stokes--Dirichlet resolvent problem is self-adjoint on $\L^2_{\mathfrak{n},\sigma}(\RR^n_+)$, and that we have existence for the (pre-)dual problem by Step 1.2 and the current Step 2.1, this imposes uniqueness by duality, see Proposition~\ref{prop:DualityDivFreeRn+}.

Consequently, we did extrapolate $\AA_\mathcal{D}$ as a densely defined and closed $0$-sectorial operator on the spaces
\begin{itemize}
    \item $\dot{\B}^{s,\sigma}_{p,q,\mathfrak{n}}$, $\dot{\BesSmo}^{s,\sigma}_{p,\infty,\mathfrak{n}}$, $1\leqslant p <\infty$, $1\leqslant q< \infty$  and $\dot{\H}^{s,p}_{\mathfrak{n},\sigma}$, $1<p<\infty$;
    \item $\dot{\B}^{s,\sigma}_{\infty,q,\mathfrak{n}}$, $\dot{\B}^{s,0,\sigma}_{\infty,q,\mathfrak{n}}$, $\dot{\BesSmo}^{s,0,\sigma}_{\infty,\infty,\mathfrak{n}}$, $1\leqslant q<\infty$,
    \item $\dot{\B}^{s,\sigma}_{p,\infty,\mathfrak{n}}$, $1\leqslant p<\infty$, $\dot{\B}^{s,0,\sigma}_{\infty,\infty,\mathfrak{n}}$, but only weakly-$\ast$ densely defined.
\end{itemize}
As a direct consequence if $\dot{\X}^{s,p}$ is any of the spaces above, then the domain of $\AA_\mathcal{D}$ in $\dot{\X}^{s,p}_{\mathfrak{n},\sigma}$ is its naive domain given by
\begin{align*}
    &\dot{\D}^{s,p}(\AA_{\mathcal{D}},\RR^n_+) = \{\,\uu\in[\dot{\X}^{s,p}\cap\dot{\X}^{s+2,p}](\RR^n_+)\,:\,\div \uu =0 \text{, and }\, \uu_{|_{\partial\RR^n_+}}=0 \,\}, \\
    &\AA_{\mathcal{D}}\uu=\PP_{\RR^n_+}(-\Delta_\mathcal{D}\uu)\text{, }\uu\in \dot{\D}^{s,p}(\AA_{\mathcal{D}},\RR^n_+).
\end{align*}

Indeed, if $\uu\in\dot{\D}^{s,p}(\AA_{\mathcal{D}},\RR^n_+)\subset[\dot{\X}^{s,p}\cap\dot{\X}^{s+2,p}_{\mathcal{D},\sigma}](\RR^n_+)$, let $\nabla\mathfrak{p}\in\dot{\X}^{s,p}(\RR^n_+,\CC^n)$ such that
\begin{align*}
    \AA_\mathcal{D}u =-\Delta \uu + \nabla \mathfrak{p}=:\ff\in\dot{\X}^{s,p}_{\mathfrak{n},\sigma}(\RR^n_+).
\end{align*}
By the boundedness of the Leray projection $\PP_{\RR^n_+}$ on $\dot{\X}^{s,p}(\RR^n_+,\CC^n)$, one deduces
\begin{align*}
    \ff = \PP_{\RR^n_+}\ff = \PP_{\RR^n_+}(-\Delta \uu + \nabla \mathfrak{p}) = \PP_{\RR^n_+}(-\Delta \uu).
\end{align*}

It remains to check existence with the appropriate meaning in the case of the spaces $\dot{\B}^{s,\sigma}_{\infty,q,\mathfrak{n}}(\RR^n_+)$ and the accompanying resolvent estimate. The description of the operator will have the same proof. 

\medbreak

\textbf{Step 2.2:} We show the extrapolated operator $\AA_\mathcal{D}$ is injective. Uniqueness in the following spaces
\begin{itemize}
    \item $\dot{\B}^{s,\sigma}_{p,q,\mathfrak{n}}$, $\dot{\BesSmo}^{s,\sigma}_{p,\infty,\mathfrak{n}}$, $1\leqslant p,q<\infty$, and $\dot{\H}^{s,p}_{\mathfrak{n},\sigma}$, $1<p<\infty$;
    \item $\dot{\B}^{s,0,\sigma}_{\infty,q,\mathfrak{n}}$, $\dot{\BesSmo}^{s,0,\sigma}_{\infty,\infty,\mathfrak{n}}$, $1\leqslant q<\infty$.
\end{itemize}
will follow from the uniqueness in the case of the Besov spaces $\dot{\B}^{s,\sigma}_{p,\infty,\mathfrak{n}}$, $1\leqslant p\leqslant \infty$. We consider a solution  $(\vv,\nabla\mathfrak{q})\in[\dot{\B}^{s,\sigma}_{p,\infty,\mathfrak{n}}\cap\dot{\B}^{s+2,\sigma}_{p,\infty,\mathcal{D}}](\RR^n_+)\times\dot{\B}^{s}_{p,\infty}(\RR^n_+,\CC^n)$ such that 
\begin{equation}\label{eq:SelfSystemStokesInjectivity}
     \left\{ \begin{array}{rllr}
           - \Delta \vv +\nabla \mathfrak{q} &= 0 \text{, }&&\text{ in } \RR^n_+\text{,}\\
        \div \vv &= 0\text{, }&&\text{ in } \RR^n_+\text{,}\\
        \uu_{|_{\partial\RR^n_+}} &=0\text{, } &&\text{ on } \partial\RR^n_+\text{.}
    \end{array}
    \right.
\end{equation}
By taking the divergence of the first equation, for $\mathfrak{q}$ determined up to a constant, it should satisfies
\begin{equation*}
    \left\{ \begin{array}{rllr}
         - \Delta \mathfrak{q} &= 0 \text{, }&&\text{ in } \RR^n_+\text{,}\\
        \partial_{x_n}\mathfrak{q}_{|_{\partial\RR^n_+}} &= \Delta \vv_n{}_{|_{\partial\RR^n_+}}\text{, }  &&\text{ on } \partial\RR^n_+\text{.}
    \end{array}
    \right.
\end{equation*}
Consequently, one has $\nabla\mathfrak{q}(\cdot',x_n)=-S e^{-x_n(-\Delta')^\frac{1}{2}}[\Delta \vv_n{}_{|_{\partial\RR^n_+}}]+e^{-x_n(-\Delta')^\frac{1}{2}}[\Delta \vv_n{}_{|_{\partial\RR^n_+}}]\mathfrak{e}_n$, so that applying the Uka\"{i} operator $\mathrm{E}_{0,\sigma}^{\mathcal{U}}$ to $\nabla \mathfrak{q}$, one obtains
\begin{align*}
    \mathrm{E}_{0,\sigma}^{\mathcal{U}}\nabla \mathfrak{q} =\mathbf{0}.
\end{align*}
Therefore, applying $\mathrm{E}_{0,\sigma}^{\mathcal{U}}$ to the first equation of \eqref{eq:SelfSystemStokesInjectivity}, we deduce
\begin{align*}
   \mathrm{E}_{0,\sigma}^{\mathcal{U}}\Delta \vv =\mathbf{0}.
\end{align*}
On the other hand, due to $\uu_{|_{\partial\RR^n_+}}=0$, it holds in $\S'(\RR^n,\CC^n)$:
\begin{align*}
    \mathrm{E}_{0,\sigma}^{\mathcal{U}}(-\Delta \vv) = \begin{bmatrix}\mathbf{I}_{n-1}& {S}\\ \prescript{t}{}{(-S)} & \I\end{bmatrix}\Delta \E_{\mathcal{D}}\vv = \Delta \begin{bmatrix}\mathbf{I}_{n-1}& {S}\\ \prescript{t}{}{(-S)} & \I\end{bmatrix}\E_{\mathcal{D}}\vv =-\Delta \mathrm{E}_{0,\sigma}^{\mathcal{U}}\vv
\end{align*}
So that, $\E_{0,\sigma}^\mathcal{U} \vv=0$ due to $\E_{0,\sigma}^\mathcal{U} \vv\in\dot{\B}^{s}_{p,\infty}(\RR^n)$ by  Lemma~\ref{lem:UkaiOperators}. Since $[\Gamma_{\sigma}^\mathcal{U}\E_{0,\sigma}^\mathcal{U}\vv]_{|_{\RR^n_+}}=\vv$, see Lemma~\ref{lem:UkaiOperators}, it follows that $\vv=0$.

\medbreak

\textbf{Step 2.3:} It remains to deal with solving \eqref{eq:DirStokesSystemtn+} in the spaces $\dot{\B}^{s,\sigma}_{\infty,q,\mathfrak{n}}$, $q\in[1,\infty]$, $s\in(-1,0)$, and to obtain the resolvent estimate. Actually, it relies on the duality result Proposition \ref{prop:DualityDivFreeRn+}. The dense range of $\AA_\mathcal{D}$ as an operator on $\dot{\B}^{s,\sigma}_{1,q',\mathfrak{n}}(\RR^n_+)$, $1\leqslant q'<\infty$, and on $\dot{\BesSmo}^{s,\sigma}_{1,\infty,\mathfrak{n}}(\RR^n_+)$, implies uniqueness of such a solution $\uu\in[\dot{\B}^{s,\sigma}_{\infty,q,\mathfrak{n}}\cap\dot{\B}^{s+2,\sigma}_{\infty,q,\mathcal{D}}](\RR^n_+)$  to \eqref{eq:DirStokesSystemtn+} (with $g=0$ and $\hh=0$). By Lemma \ref{lem:UkaiOperators}, for all $\ff\in\dot{\B}^{s,\sigma}_{\infty,q,\mathfrak{n}}(\RR^n_+)$:
\begin{align*}
    |\lambda|\lVert  \Gamma_{\sigma}^{\mathcal{U}} (\lambda \I - \Delta)^{-1}\mathrm{E}_{0,\sigma}^{\mathcal{U}} \ff\rVert_{\dot{\B}^{s}_{\infty,q}(\RR^{n})} \less_{p,s,n,\mu} \lVert \ff \rVert_{\dot{\B}^{s}_{\infty,q}(\RR^{n}_+)},\\
    \lVert\Delta\Gamma_{\sigma}^{\mathcal{U}} (\lambda \I - \Delta)^{-1}\mathrm{E}_{0,\sigma}^{\mathcal{U}} \ff\rVert_{\dot{\B}^{s}_{\infty,q}(\RR^{n})} \less_{p,s,n,\mu} \lVert \ff \rVert_{\dot{\B}^{s}_{\infty,q}(\RR^{n}_+)}.
\end{align*}
So we have to check that $ \Gamma_{\sigma}^{\mathcal{U}} (\lambda \I - \Delta)^{-1}\mathrm{E}_{0,\sigma}^{\mathcal{U}} \ff$ solves \eqref{eq:DirStokesSystemtn+} in $\RR^n_+$ (with $g=0$ 
 and $\hh=0$). Notice that for $\mathcal{T}\in \{\,\Gamma_{\sigma}^{\mathcal{U}},\, (\lambda \I - \Delta)^{-1},\,\mathrm{E}_{0,\sigma}^{\mathcal{U}}\,\}$, one can show on $\dot{\B}^{s}_{\infty,q}$ that
\begin{align*}
    \mathcal{T} =\begin{cases}
                \mathcal{T}_{|_{\dot{\B}^{s,0}_{\infty,q'}}}^{\ast\ast}\text{, if }1<q<\infty;\\
                \mathcal{T}_{|_{\dot{\BesSmo}^{s,0}_{\infty,\infty}}}^{\ast\ast}\text{, if }q=\infty;\\
                \big( (\mathcal{T}_{|_{\dot{\B}^{s,0}_{\infty,1}}}^{\ast})_{|_{\dot{\BesSmo}^{s}_{1,\infty}}}\big)^\ast\text{, if }q=1.
                \end{cases}
\end{align*}
Thus, the map $\ff\mapsto \Gamma_{\sigma}^{\mathcal{U}} (\lambda \I - \Delta)^{-1}\mathrm{E}_{0,\sigma}^{\mathcal{U}} \ff$ is weak-$\ast$ continuous on $\dot{\B}^{s}_{\infty,q}$. By Corollary \ref{cor:Weak*densityDivFreeRn+}, $\Ccinftydiv(\RR^n_+)$ is weakly-$\ast$ dense in $\dot{\B}^{s,\sigma}_{\infty,q,\mathfrak{n}}(\RR^{n}_+)$. From this point, one can check for $(\ff_k)_{k\in\NN}\subset\Ccinftydiv(\RR^n_+)$ converging weakly-$\ast$ to $\ff$ in $\dot{\B}^{s,\sigma}_{\infty,q,\mathfrak{n}}(\RR^{n}_+)$, that $\uu_k:= [\Gamma_{\sigma}^{\mathcal{U}} (\lambda \I - \Delta)^{-1}\mathrm{E}_{0,\sigma}^{\mathcal{U}} \ff_k]_{|_{\RR^n_+}}$ and $\nabla\mathfrak{p}_k=-[\I-\PP_{\RR^n_+}][(-\Delta)\uu_k-\ff_k]$ converge weakly-$\ast$ to a solution $(\uu,\nabla\mathfrak{p})$ of \eqref{eq:DirStokesSystemtn+} with data $\ff$ (and with $g=0$ 
 and $\hh=0$). By uniqueness of the limit and by weak-$\ast$ continuity of involved operators, it turns out that $\uu=[\Gamma_{\sigma}^{\mathcal{U}} (\lambda \I - \Delta)^{-1}\mathrm{E}_{0,\sigma}^{\mathcal{U}} \ff]_{|_{\RR^n_+}}\in[\dot{\B}^{s,\sigma}_{\infty,q,\mathfrak{n}}\cap\dot{\B}^{s+2,\sigma}_{\infty,q,\mathcal{D}}](\RR^n_+)$ is such a solution. 

Thus, the formula
\begin{align*}
    \uu =[\Gamma_{\sigma}^{\mathcal{U}} (\lambda \I - \Delta)^{-1}\mathrm{E}_{0,\sigma}^{\mathcal{U}} \ff]_{|_{\RR^n_+}}
\end{align*}
remains valid on $\dot{\B}^{s,\sigma}_{\infty,q,\mathfrak{n}}(\RR^n_+)$ to provide a (the) solution.

\textbf{Step 3:} Now, still for $p,q\in[1,\infty]$, $s\in(-1+\sfrac{1}{p},\sfrac{1}{p})$, since we already have uniqueness by linearity and the previous steps, we only have to construct a solution to \eqref{eq:DirStokesSystemtn+} for $g\neq0$ and $\hh\neq0$. By boundedness of the Leray projection, then by linearity we can assume $\ff=0$. Up to consider ${|\lambda|}\uu(\cdot/|\lambda|^{\sfrac{1}{2}})$, $|\lambda|^{\sfrac{1}{2}}\mathfrak{ p}(\cdot/|\lambda|^{\sfrac{1}{2}})$, $|\lambda|^{\sfrac{1}{2}} g(\cdot/|\lambda|^{\sfrac{1}{2}})$ and ${|\lambda|}\hh(\cdot/|\lambda|^{\sfrac{1}{2}})$, instead of $\uu$, $\mathfrak{p}$, $g$ and $\hh$ respectively, we can assume, without loss of generality $|\lambda|=1$.

We split the system in two parts, writing $g=0+g$ and $\hh = \hh' +\hh_n\mathfrak{e}_n$, $\uu=\vv+\ww$, and $\mathfrak{p}=\mathfrak{q}+\pi$.

\textbf{Step 3.1:} We consider the system
\begin{equation*}
    \left\{ \begin{array}{rllr}
         \lambda \vv - \Delta \vv +\nabla \mathfrak{q} &= 0 \text{, }&&\text{ in } \RR^n_+\text{,}\\
        \div \vv &= 0\text{, }&&\text{ in } \RR^n_+\text{,}\\
        \vv_{|_{\partial\RR^n_+}} &=(\hh',0)\text{, } &&\text{ on } \partial\RR^n_+\text{.}
    \end{array}
    \right.
\end{equation*}

The ansatz for $\vv=(\vv',\vv_n)$ and $\nabla \mathfrak{q}$ are given, for $x_n\geqslant0$, by
\begin{align*}
    \vv'(\cdot,x_n)&:=\frac{1}{\lambda} [(\lambda\I-\Delta')^{{\sfrac{1}{2}}}+(-\Delta')^{{\sfrac{1}{2}}}]((\lambda\I-\Delta')^{{\sfrac{1}{2}}}e^{-x_n(\lambda\I-\Delta')^{{\sfrac{1}{2}}}}-(-\Delta')^{{\sfrac{1}{2}}}e^{-x_n(-\Delta')^{{\sfrac{1}{2}}}})(\I-\PP')\hh'\\ &\qquad \qquad + e^{-x_n(\lambda\I-\Delta')^{{\sfrac{1}{2}}}}\PP'\hh';\\
   \vv_n (\cdot,x_n)&:= \frac{1}{\lambda}[(\lambda\I-\Delta')^{{\sfrac{1}{2}}}+(-\Delta')^{{\sfrac{1}{2}}}](e^{-x_n(\lambda\I-\Delta')^{{\sfrac{1}{2}}}}-e^{-x_n(-\Delta')^{{\sfrac{1}{2}}}})\div' \hh';\\
   \nabla \mathfrak{q} (\cdot,x_n) &:= \nabla(-\Delta')^{-{\sfrac{1}{2}}}[(\lambda\I-\Delta')^{{\sfrac{1}{2}}}+(-\Delta')^{{\sfrac{1}{2}}}]e^{-x_n(-\Delta')^{{\sfrac{1}{2}}}}\div' \hh'.
\end{align*}
This can be derived from \cite[Lemma~4.2.1,~Proof]{DanchinMucha2015}\footnote{in this reference, one should somehow follow the proof replacing $i\xi_0$ by $\lambda$}. One may also consult \cite[Theorem~1.3,~Proof]{FarwigSohr1994}. For simplification, we write here $\PP':=\PP_{\RR^{n-1}}$.

By Propositions~\ref{prop:DumpedPoissonSemigroup2} and \ref{prop:PoissonSemigroup3}, for $\hh'\in\dot{\B}^{s-1/p}_{p,p}\cap\dot{\B}^{s+2-1/p}_{p,p}(\partial\RR^{n}_+,\CC^{n-1})$, recalling that $|\lambda|=1$,
\begin{align}
    \lVert \vv \rVert_{\dot{\B}^{s}_{p,q}(\RR^n_+)} &\lesssim_{p,s,n}^{\mu} \lVert [(\lambda\I-\Delta')^{{\sfrac{1}{2}}}+(-\Delta')^{{\sfrac{1}{2}}}](\I-\PP')\hh' \rVert_{{\B}^{s+1-1/p}_{p,q}(\RR^{n-1})} + \lVert \PP' \hh' \rVert_{{\B}^{s-1/p}_{p,q}(\RR^{n-1})}\nonumber\\ &\qquad +\lVert [(\lambda\I-\Delta')^{{\sfrac{1}{2}}}+(-\Delta')^{{\sfrac{1}{2}}}] (\I-\PP')\hh' \rVert_{\dot{\B}^{s+1-1/p}_{p,q}(\RR^{n-1})} \nonumber\\ &\qquad +\lVert [(\lambda\I-\Delta')^{{\sfrac{1}{2}}}+(-\Delta')^{{\sfrac{1}{2}}}]\div' \hh' \rVert_{\dot{\B}^{s-1/p}_{p,q}(\RR^{n-1})} \nonumber\\
    &\lesssim_{p,s,n}^{\mu} \lVert (\I-\PP')\hh' \rVert_{{\B}^{s+2-1/p}_{p,q}(\RR^{n-1})} + \lVert \PP' \hh' \rVert_{{\B}^{s+2-1/p}_{p,q}(\RR^{n-1})}+\lVert \hh' \rVert_{{\B}^{s+2-1/p}_{p,q}(\RR^{n-1})}\nonumber\\
    &\lesssim_{p,s,n}^{\mu} \lVert (\hh',(\I-\PP')\hh',\PP'\hh')\rVert_{\L^p(\RR^{n-1})}+ \lVert (\hh',(\I-\PP')\hh',\PP'\hh') \rVert_{\dot{\B}^{s+2-1/p}_{p,q}(\RR^{n-1})}\label{eq:EstimateTraceProofNonDivFreeRn+}
\end{align}
By interpolation inequalities,
\begin{align*}
    \lVert \vv \rVert_{\dot{\B}^{s}_{p,q}(\RR^n_+)} &\lesssim_{p,s,n}^{\mu} \lVert (\hh',(\I-\PP')\hh',\PP'\hh')\rVert_{\dot{\B}^{s-1/p}_{p,q}(\RR^{n-1})}^{1-\theta}\lVert (\hh',(\I-\PP')\hh',\PP'\hh') \rVert_{\dot{\B}^{s+2-1/p}_{p,q}(\RR^{n-1})}^{\theta}\nonumber\\ &\qquad\qquad\qquad\qquad\qquad\qquad\qquad\qquad\qquad\qquad+ \lVert (\hh',(\I-\PP')\hh',\PP'\hh') \rVert_{\dot{\B}^{s+2-1/p}_{p,q}(\RR^{n-1})}\nonumber\\
    &\lesssim_{p,s,n}^{\mu} \lVert \hh' \rVert_{\dot{\B}^{s-1/p}_{p,p}(\RR^{n-1})}+ \lVert \hh' \rVert_{\dot{\B}^{s+2-1/p}_{p,q}(\RR^{n-1})}\nonumber.
\end{align*}

Now, if $\hh':=\mathbf{H'}_{|_{\partial\RR^n_+}}$, for some $\mathbf{H'}\in\dot{\B}^{s}_{p,q}(\RR^n_+,\CC^{n-1})\cap\dot{\B}^{s+2}_{p,q}(\RR^n_+,\CC^{n-1})$, we start from \eqref{eq:EstimateTraceProofNonDivFreeRn+}, so that by the trace theorem : 
\begin{align*}
    \lVert \vv \rVert_{\dot{\B}^{s}_{p,q}(\RR^n_+)} &\lesssim_{p,s,n}^{\mu} \lVert (\mathbf{H'},(\I-\PP')\mathbf{H'},\PP'\mathbf{H'})\rVert_{\dot{\B}^{1/p}_{p,1}(\RR^{n}_+)}+ \lVert (\mathbf{H'},(\I-\PP')\mathbf{H'},\PP'\mathbf{H'}) \rVert_{\dot{\B}^{s+2}_{p,q}(\RR^{n}_+)}\\
    &\lesssim_{p,s,n}^{\mu} \lVert \mathbf{H'}\rVert_{\dot{\B}^{s}_{p,q}(\RR^{n}_+)}^{1-\theta} \lVert \mathbf{H'}\rVert_{\dot{\B}^{s+2}_{p,q}(\RR^{n}_+)}^{\theta}+ \lVert \mathbf{H'} \rVert_{\dot{\B}^{s+2}_{p,q}(\RR^{n}_+)}\\
    &\lesssim_{p,s,n}^{\mu} \lVert \mathbf{H'}\rVert_{\dot{\B}^{s}_{p,q}(\RR^{n}_+)} + \lVert \mathbf{H'} \rVert_{\dot{\B}^{s+2}_{p,q}(\RR^{n}_+)}.
\end{align*}

For the pressure term, one has, similarly:
\begin{align*}
     \lVert \nabla \mathfrak{q} \rVert_{\dot{\B}^{s}_{p,q}(\RR^n_+)} &\lesssim_{p,s,n} \lVert [(\lambda\I-\Delta')^{{\sfrac{1}{2}}}+(-\Delta')^{{\sfrac{1}{2}}}] \div \hh' \rVert_{\dot{\B}^{s-1/p}_{p,q}(\RR^{n-1})} \\ &\lesssim_{p,s,n}^{\mu} \lVert \hh' \rVert_{\dot{\B}^{s+1-1/p}_{p,q}(\RR^{n-1})} + \lVert \hh' \rVert_{\dot{\B}^{s+2-1/p}_{p,q}(\RR^{n-1})}.
\end{align*}
Now, the second order derivatives of $\vv$, by Proposition \ref{prop:DirPbRn+}:
\begin{align*}
     \lVert \nabla^2 \vv \rVert_{\dot{\B}^{s}_{p,q}(\RR^n_+)} &\lesssim_{p,s,n}  \lVert \Delta \vv \rVert_{\dot{\B}^{s}_{p,q}(\RR^n_+)} + \lVert \hh'\rVert_{\dot{\B}^{s+2-1/p}_{p,q}(\RR^{n-1})}\\
     &\lesssim_{p,s,n}  \lVert \vv \rVert_{\dot{\B}^{s}_{p,q}(\RR^n_+)} + \lVert \nabla \mathfrak{q} \rVert_{\dot{\B}^{s}_{p,q}(\RR^n_+)} + \lVert \hh'\rVert_{\dot{\B}^{s+2-1/p}_{p,q}(\RR^{n-1})},
\end{align*}
where we did use  $\Delta \vv = \lambda \vv + \nabla\mathfrak{q}$ to obtain the last estimate, and recall that we assumed $|\lambda|=1$.

\medbreak

\noindent Summing up everything, we have obtained
\begin{align*}
    \lVert  \vv\rVert_{\dot{\B}^{s}_{p,q}(\RR^n_+)}+&\lVert \nabla^2 \vv\rVert_{\dot{\B}^{s}_{p,q}(\RR^n_+)} + \lVert \nabla \mathfrak{q}\rVert_{\dot{\B}^{s}_{p,q}(\RR^n_+)} \lesssim_{p,n,s,\mu}  \min \begin{cases}
  \lVert  \mathbf{H}'\rVert_{\dot{\B}^{s}_{p,q}(\RR^n_+)}+\lVert \mathbf{H}'\rVert_{\dot{\B}^{s}_{p,q}(\RR^n_+)}\\    
  \lVert  \hh'\rVert_{\dot{\B}^{s-1/p}_{p,q}(\partial\RR^n_+)}+\lVert \hh'\rVert_{\dot{\B}^{s+2-1/p}_{p,q}(\partial\RR^n_+)}
\end{cases}.
\end{align*}

The case of homogeneous Sobolev space $\dot{\H}^{s,p}(\RR^n_+)$, $p\in(1,\infty)$, admits a similar proof.

\textbf{Step 3.2:} We consider the system
\begin{equation}\label{eq:ProofRn+StokesSystemnondivfreeonly}
    \left\{ \begin{array}{rllr}
         \lambda \ww - \Delta \ww +\nabla \pi &= 0 \text{, }&&\text{ in } \RR^n_+\text{,}\\
        \div \ww &= g\text{, }&&\text{ in } \RR^n_+\text{,}\\
        \ww_{|_{\partial\RR^n_+}} &= \hh_n\mathfrak{e}_n\text{, } &&\text{ on } \partial\RR^n_+\text{.}
    \end{array}
    \right.
\end{equation}

Up to consider again $|\lambda|=1$, if we set $\ww_1:=\ww-\ww_0$, with $\ww_0(\cdot,x_n):= e^{-x_n(\lambda\I-\Delta')^{\sfrac{1}{2}}}\hh_n \mathfrak{e}_n$, it reduces to the system
\begin{equation}\label{eq:ProofRn+StokesSystemnondivLift}
    \left\{ \begin{array}{rllr}
         \lambda \ww_1 - \Delta \ww_1 +\nabla \pi &= 0 \text{, }&&\text{ in } \RR^n_+\text{,}\\
        \div \ww_1 &= \Tilde{g}\text{, }&&\text{ in } \RR^n_+\text{,}\\
        {\ww_1}_{|_{\partial\RR^n_+}} &= 0\text{, } &&\text{ on } \partial\RR^n_+\text{,}
    \end{array}
    \right.
\end{equation}
where $\Tilde{g}=g-\partial_{x_n}\ww_0\in\dot{\B}^{s-1}_{p,q,0}\cap\dot{\B}^{s+1}_{p,q}(\RR^n_+)$, and the estimate provided by Proposition~\ref{prop:DirResolventPbRn+}:
\begin{align}\label{eq:EstStokesDirNonDivFree0}
    \lVert  \ww_0 \rVert_{\dot{\B}^{s}_{p,q}(\RR^n_+)} + \lVert\nabla^2 \ww_0\rVert_{\dot{\B}^{s}_{p,q}(\RR^n_+)} \lesssim_{p,s,n}^{\mu} \min \begin{cases}
  \lVert  \mathbf{H}_n\rVert_{\dot{\B}^{s}_{p,q}(\RR^n_+)}+\lVert \nabla^2 \mathbf{H}_n\rVert_{\dot{\B}^{s}_{p,q}(\RR^n_+)}\\    
 \lVert  \hh_n\rVert_{\dot{\B}^{s-1/p}_{p,q}(\partial\RR^n_+)}+\lVert \nabla'^2 \hh_n\rVert_{\dot{\B}^{s-1/p}_{p,q}(\partial\RR^n_+)}
  \end{cases}.
\end{align}
Following the beginning of the proof of \cite[Theorem~3.1]{GengShen2024}, we plug a new ansatz to solve \eqref{eq:ProofRn+StokesSystemnondivLift}. We consider $\tilde{\ww}=[(-\Delta)^{-{\sfrac{1}{2}}}\nabla(-\Delta)^{-{\sfrac{1}{2}}}\E_{\mathcal{N}}\tilde{g}]_{|_{\RR^n_+}}$\footnote{written that way to avoid issues of definition due to the lack of completeness for function spaces with high regularities.} and pressure term $\nabla \tilde{\pi}= \nabla\tilde{g}-\lambda\tilde{\ww}$, which by definition of function spaces by restriction and Lemma~\ref{lem:ExtOpNegativeBesovSpaces}, satisfy the estimate
\begin{align}\label{eq:EstStokesDirNonDivFree1}
    \lVert \tilde{\ww}\rVert_{\dot{\B}^{s}_{p,q}(\RR^n_+)}+\lVert \nabla^2\tilde{\ww}\rVert_{\dot{\B}^{s}_{p,q}(\RR^n_+)} +\lVert \nabla\tilde{\pi}\rVert_{\dot{\B}^{s}_{p,q}(\RR^n_+)}  \lesssim_{p,s,n}  \lVert \tilde{g}\rVert_{\dot{\B}^{s-1}_{p,q}(\RR^n_+)} + \lVert  \tilde{g}\rVert_{\dot{\B}^{s+1}_{p,q}(\RR^n_+)}.
\end{align}
Moreover, one can check that ${ \tilde{\ww}_n }{}_{|_{\partial\RR^n_+}}=0$, thus for $\ww_2:=\ww_1-\tilde{\ww}$ and $\nabla \pi' = \nabla\pi-\nabla \tilde{\pi}$, the system reduces to
\begin{equation*}
    \left\{ \begin{array}{rllr}
         \lambda \ww_2 - \Delta \ww_2 +\nabla \pi' &= 0 \text{, }&&\text{ in } \RR^n_+\text{,}\\
        \div \ww_2 &= 0\text{, }&&\text{ in } \RR^n_+\text{,}\\
        {\ww_2}_{|_{\partial\RR^n_+}} &= -({\tilde{\ww}'}_{|_{\partial\RR^n_+}},0)\text{, } &&\text{ on } \partial\RR^n_+\text{,}
    \end{array}
    \right.
\end{equation*}
But this amounts to the previous Step 3.1, which states that there exists a unique solution $(\ww_2,\nabla\pi')$, which satisfies the estimate
\begin{align}\label{eq:EstStokesDirNonDivFree2}
    \lVert {\ww}_2\rVert_{\dot{\B}^{s}_{p,q}(\RR^n_+)}+\lVert \nabla^2{\ww}_2\rVert_{\dot{\B}^{s}_{p,q}(\RR^n_+)} +\lVert \nabla{\pi'}\rVert_{\dot{\B}^{s}_{p,q}(\RR^n_+)}  \lesssim_{p,s,n}^{\mu}  \lVert \tilde{\ww}'\rVert_{\dot{\B}^{s}_{p,q}(\RR^n_+)} + \lVert  \tilde{\ww}'\rVert_{\dot{\B}^{s+2}_{p,q}(\RR^n_+)}
\end{align}
Since $\ww= \ww_0+\ww_1 = \ww_0+\tilde{\ww}+\ww_2$ and $\nabla\pi=\nabla\tilde{\pi}+\nabla\pi'$, we did construct a (the unique) solution $(\ww,\nabla \pi)$ to \eqref{eq:ProofRn+StokesSystemnondivfreeonly}, which, thanks to \eqref{eq:EstStokesDirNonDivFree0}, \eqref{eq:EstStokesDirNonDivFree1} and \eqref{eq:EstStokesDirNonDivFree2}, satisfies the estimate
\begin{align*}
    \lVert {\ww}\rVert_{\dot{\B}^{s}_{p,q}(\RR^n_+)}+\lVert &\nabla^2{\ww}\rVert_{\dot{\B}^{s}_{p,q}(\RR^n_+)} +\lVert \nabla{\pi}\rVert_{\dot{\B}^{s}_{p,q}(\RR^n_+)} \\ &\lesssim_{p,s,n}^{\mu}  \lVert g\rVert_{\dot{\B}^{s-1}_{p,q}(\RR^n_+)} + \lVert g\rVert_{\dot{\B}^{s+1}_{p,q}(\RR^n_+)}\\ &\quad\quad\quad + \min \begin{cases}
  \lVert  \mathbf{H}_n\rVert_{\dot{\B}^{s}_{p,q}(\RR^n_+)}+\lVert \nabla^2 \mathbf{H}_n\rVert_{\dot{\B}^{s}_{p,q}(\RR^n_+)}\\    
 \lVert  \hh_n\rVert_{\dot{\B}^{s-1/p}_{p,q}(\partial\RR^n_+)}+\lVert \nabla'^2 \hh_n\rVert_{\dot{\B}^{s-1/p}_{p,q}(\partial\RR^n_+)}.
  \end{cases}
\end{align*}
This ends the Step 3, and the proof for the resolvent estimate.
\end{proof}

\begin{proposition}\label{prop:MetaThmDirichletStokesRn+HigherReg}Let $p,q,r,\kappa\in[1,\infty]$, $s\in (-1+\sfrac{1}{p},\sfrac{1}{p})$, $\uptau\in(-1+\sfrac{1}{r},1+\sfrac{1}{r})$, such that $\uptau\neq\sfrac{1}{r}$. Let $\ff\in {\dot{\X}}^{s,p}_{\mathfrak{n},\sigma}\cap\dot{\Y}^{\uptau,r}_{\mathcal{D},\sigma}(\RR^n_+)$, then the solution $(\uu,\nabla \mathfrak{p})$ given by Theorem~\ref{thm:MetaThmDirichletStokesRn+} to the Stokes--Dirichlet resolvent problem 
    \begin{equation*}\tag{DS${}_\lambda$}\label{eq:DirStokesSystemtn+Hom}
    \left\{ \begin{array}{rllr}
         \lambda \uu - \Delta \uu +\nabla \mathfrak{p} &= \ff \text{, }&&\text{ in } \RR^n_+\text{,}\\
        \div \uu &= 0\text{, } &&\text{ in } \RR^n_+\text{,}\\
        \uu_{|_{\partial\RR^n_+}} &=0\text{, } &&\text{ on } \partial\RR^n_+\text{.}
    \end{array}
    \right.
\end{equation*}
        also satisfies $(\uu,\nabla \mathfrak{p})\in \dot{\Y}^{\uptau,r}\cap\dot{\Y}^{\uptau+2,r}(\RR^n_+,\CC^n) \times \dot{\Y}^{\uptau,r}(\RR^n_+,\CC^n)$ with  the estimates
        \begin{align*}
            \lvert\lambda\rvert\lVert  \uu\rVert_{\dot{\Y}^{\uptau,r}(\RR^n_+)}+\lVert \nabla^2 \uu\rVert_{\dot{\Y}^{\uptau,r}(\RR^n_+)} + \lVert \nabla \mathfrak{p}\rVert_{\dot{\Y}^{\uptau,r}(\RR^n_+)} \lesssim_{r,n,\uptau,\mu} \lVert \ff\rVert_{{\dot{\Y}}^{\uptau,r}(\RR^n_+)}.
        \end{align*}
\end{proposition}

\begin{proof} One only deal with the case  $\uptau\in(\sfrac{1}{r},1+\sfrac{1}{r})$. For $p,q\in[1,\infty]$, $s\in(-1+\sfrac{1}{p},\sfrac{1}{p})$, we assume first $p=r$, $q=\kappa$ et $\uptau =s+1$. Let $\ff\in\dot{\X}^{s,p}_{\mathfrak{n},\sigma}\cap\dot{\X}^{s+1,p}_{\mathcal{D},\sigma}(\RR^n_+)$, and consider $(\uu,\nabla \mathfrak{p})\in\dot{\X}^{s,p}_{\mathfrak{n},\sigma}\cap\dot{\X}^{s+2,p}_{\mathcal{D},\sigma}(\RR^n_+)\times\dot{\X}^{s,p}(\RR^n_+,\CC^n)$ be the unique solution to \eqref{eq:DirStokesSystemtn+Hom}.

Note that in particular, writing for $h>0$
\begin{align*}
    \partial_{j}^{h} w := \frac{w(\cdot+h\mathfrak{e}_j)-w}{h},
\end{align*}
one obtains
\begin{equation*}
    \left\{ \begin{array}{rllr}
         \lambda  \partial_{j}^{h}\uu  - \Delta  \partial_{j}^{h}\uu +\nabla  \partial_{j}^{h}\mathfrak{p} &=  \partial_{j}^{h}\ff \text{, }&&\text{ in } \RR^n_+\text{,}\\
        \div  (\partial_{j}^{h}\uu) &= 0\text{, } &&\text{ in } \RR^n_+\text{,}\\
        \partial_{j}^{h} \uu_{|_{\partial\RR^n_+}} &=0\text{, } &&\text{ on } \partial\RR^n_+\text{.}
    \end{array}
    \right.
\end{equation*}
Since $(\partial_{j}^{h}\uu)_{h>0}$ converges towards $\partial_{x_j}\uu$ strongly in $\dot{\X}^{s,p}({\RR^n_+})$, then the convergence holds in the distributional sense, so that for $\D'^{h} := (\partial_{j}^{h})_{1\leqslant j \leqslant n-1}$:
\begin{align*}
    |\lambda|\lVert \nabla' \uu\rVert_{\dot{\X}^{s,p}({\RR^n_+})} &+  \lVert \nabla' \nabla^2 \uu\rVert_{\dot{\X}^{s,p}({\RR^n_+})} \\
    &\leqslant \liminf_{h\longrightarrow 0_+} \, |\lambda|\lVert \D'^{h}\uu\rVert_{\dot{\X}^{s,p}({\RR^n_+})} + \lVert \nabla^2 \D'^{h}\uu\rVert_{\dot{\X}^{s,p}({\RR^n_+})}\\ &\lesssim_{p,s,n,\mu} \liminf_{h\longrightarrow 0_+} \,\lVert \D'^{h}\ff\rVert_{\dot{\X}^{s,p}({\RR^n_+})} \\
    &\lesssim_{p,s,n,\mu}\, \lVert \nabla' \ff\rVert_{\dot{\X}^{s,p}({\RR^n_+})}.
\end{align*}
To bound $\partial_{x_n}\uu_n$, we use the divergence free condition:
\begin{align*}
    |\lambda|\lVert \partial_{x_n}\uu_n\rVert_{\dot{\X}^{s,p}({\RR^n_+})} &+ \lVert \partial_{x_n} \nabla^2 \uu_n\rVert_{\dot{\X}^{s,p}({\RR^n_+})}\\
    &=  |\lambda|\lVert \div' \uu'\rVert_{{\L}^{\infty}(\RR^n_+)}+ \lVert \nabla^2\div' \uu'\rVert_{\dot{\X}^{s,p}({\RR^n_+})}\\
    &\lesssim_{n}  |\lambda|\lVert \nabla' \uu'\rVert_{\dot{\X}^{s,p}({\RR^n_+})}  +\lVert \nabla'\nabla^2 \uu'\rVert_{\dot{\X}^{s,p}({\RR^n_+})}\\
    &\lesssim_{s,n,\mu} \lVert \nabla \ff\rVert_{\dot{\X}^{s,p}({\RR^n_+})}.
\end{align*}
It remains to obtain the estimates on $\partial_{x_n}\uu'=(\partial_{x_n}\uu_j)_{1\leqslant j\leqslant n-1}$. One considers $\omega := \curl \uu(=\d \uu)$, since $\ff\in\dot{\X}^{s+1,p}_{\mathcal{D},\sigma}({\RR^n})$, for $\tilde{\ff}$ the extension of ${\ff}$ to the whole $\RR^n$ by $0$, one has $\tilde{\ff} \in \dot{\X}^{s+1,p}_{0,\sigma}({\RR^n})$ and
\begin{align*}
    \curl \tilde{\ff} = \widetilde{\curl \ff} \in  \dot{\X}^{s,p}_{0}({\RR^n})
\end{align*}
with $\supp \curl \tilde{\ff}\subset \overline{\RR^n_+}$, so that
\begin{align*}
    (-\mathfrak{e}_n)\iprod \curl {\ff}_{|_{\partial\RR^n_+}} = 0
\end{align*}
and similarly $(-\mathfrak{e}_n)\iprod \omega_{|_{\partial\RR^n_+}} = 0$. Since $\d \omega = \d^2 \uu =0$, one also obtains for free $(-\mathfrak{e}_n)\iprod \d \omega_{|_{\partial\RR^n_+}} = 0$, so that
\begin{equation*}
    \left\{ \begin{array}{rllr}
         \lambda  \omega  - \Delta  \omega &=  \curl \ff \text{, }&&\text{ in } \RR^n_+\text{,}\\
         (-\mathfrak{e}_n)\iprod \omega_{|_{\partial\RR^n_+}}  &=0\text{, } &&\text{ on } \partial\RR^n_+,\\
         (-\mathfrak{e}_n)\iprod\d \omega_{|_{\partial\RR^n_+}}  &=0\text{, } &&\text{ on } \partial\RR^n_+\text{.}
    \end{array}
    \right.
\end{equation*}
Now, applying Proposition~\ref{prop:HodgeResolventPbRn+Linfty}, we deduce from $(\omega)_{kn}= \partial_{x_n}\uu_k -\partial_{x_k}\uu_n$ and writing $ \partial_{x_n}\uu' = \sum_{k=1}^{n-1}(\omega)_{kn} \mathfrak{e}_k + \nabla'\uu_n$
\begin{align*}
    |\lambda|\lVert \partial_{x_n}\uu' \rVert_{\dot{\X}^{s,p}({\RR^n_+})} &+  \lVert \nabla^2\partial_{x_n}\uu' \rVert_{\dot{\X}^{s,p}({\RR^n_+})}\\
    &\leqslant |\lambda|\lVert \omega \rVert_{\dot{\X}^{s,p}({\RR^n_+})} + \lVert \nabla^2 \omega \rVert_{\dot{\X}^{s,p}({\RR^n_+})}\\ & \qquad+ |\lambda|\lVert \nabla'\uu_n \rVert_{\dot{\X}^{s,p}({\RR^n_+})} + \lVert \nabla^2\nabla'\uu_n  \rVert_{\dot{\X}^{s,p}({\RR^n_+})}\\
    &\lesssim_{\mu,p,s,n} \lVert \curl \ff \rVert_{\dot{\X}^{s,p}({\RR^n_+})} +\lVert \nabla'\ff\rVert_{\dot{\X}^{s,p}({\RR^n_+})}\\
    &\lesssim_{\mu,p,s,n} \lVert \nabla \ff \rVert_{\dot{\X}^{s,p}({\RR^n_+})}.
\end{align*}
By linearity, one obtains the remaining estimate on $\nabla\mathfrak{p}$.

We now, remove the assumption $\dot{\X}^{s+1,p}_{\mathcal{D},\sigma}=\dot{\Y}^{\uptau,r}_{\mathcal{D},\sigma}$. Now, provided $r\in[1,\infty]$, $\uptau\in(\sfrac{1}{r},1+\sfrac{1}{r})$ and $\kappa$ are now arbitrary, if $q<\infty$ or $\kappa<\infty$, one has $\dot{\X}^{s,p}_{\mathfrak{n},\sigma}\cap\dot{\Y}^{\uptau-1,r}_{\mathcal{D},\sigma}\cap\dot{\X}^{s+1,p}_{\mathcal{D},\sigma}\cap\dot{\Y}^{\uptau,r}_{\mathcal{D},\sigma}(\RR^n_+)$ as a strongly dense subspace of $\dot{\X}^{s,p}_{\mathfrak{n},\sigma}\cap\dot{\Y}^{\uptau,r}_{\mathcal{D},\sigma}(\RR^n_+)$, which yields the result. If $q=\kappa=\infty$, one consider the sequence $(\ff_N)_{N\in\NN}$  defined for all $N\in\NN$ by
\begin{align*}
    \ff_N := \Big[\Gamma_{\sigma}^\mathcal{U}\Big[ \sum_{|j|\leqslant N} \dot{\Delta}_j \E_{0,\sigma}^{\mathcal{U}}\ff\Big]\Big]_{|_{\RR^n_+}},
\end{align*}
which converges weakly-$\ast$ towards $\ff$ in $\dot{\X}^{s,p}_{\mathfrak{n},\sigma}\cap\dot{\Y}^{\uptau,r}_{\mathcal{D},\sigma}(\RR^n_+)$, thanks to the arguments in the proof of Theorem~\ref{thm:MetaThmDirichletStokesRn+}, Step 2.3, and which remains uniformly bounded by the corresponding norms of $\ff$ thanks to Lemma~\ref{lem:UkaiOperators}. One can than check by the same arguments that the corresponding family of solutions $(\uu_N)_{N\in\NN}$ converges weakly-$\ast$ towards $\uu$ the unique solution in $\dot{\X}^{s,p}_{\mathfrak{n},\sigma}\cap\dot{\X}^{s+2,p}_{\mathfrak{n},\sigma}$. By the Fatou property for the intersection of homogeneous Besov spaces \cite[Theorem~2.25]{bookBahouriCheminDanchin} --which remains valid on $\RR^n_+$-- (consider the intersection $\dot{\X}^{s,p}\cap\dot{\Y}^{\uptau,r}(\RR^n_+)$ as a ground space which is complete due to $s<\sfrac{1}{p}$), one also obtains $\uu\in\dot{\Y}^{\uptau,r}_{\mathcal{D},\sigma}\cap\dot{\Y}^{\uptau+2,r}(\RR^n_+)$ and the corresponding estimate.
\end{proof}

\begin{proposition}\label{prop:MetaPropDirichletStokesRn+2}Let $p,q,r,\kappa\in[1,\infty]$, $s\in (-1+\sfrac{1}{p},\sfrac{1}{p})$, $\uptau\in(-1+\sfrac{1}{r},1+\sfrac{1}{r})$, such that $\uptau\neq\sfrac{1}{r}$.
\begin{enumerate}[label=($\roman*$)]
        \item  For any $\theta\in(0,\pi)$, the Stokes--Dirichlet operator $\AA_{\mathcal{D}}$  admits a bounded $\mathbf{H}^\infty(\Sigma_\theta)$-functional calculus  on the function spaces: ${\dot{\X}}^{s,p}_{\mathfrak{n},\sigma}(\RR^n_+)$, ${\dot{\X}}^{s,p}_{\mathfrak{n},\sigma}(\RR^n_+)\cap {\dot{\Y}}^{\uptau,r}_{\mathcal{D},\sigma}(\RR^n_+)$.

         \item  The operator $\AA_{\mathcal{D}}$ on $\dot{\X}^{s,p}_{\mathfrak{n},\sigma}(\RR^n_+)$ has the following square-root property
         \begin{align*}
             \dot{\D}^{s,p}(\AA_\mathcal{D}^{\sfrac{1}{2}},\RR^n_+) = \dot{\X}^{s,p}_{\mathfrak{n},\sigma}\cap\dot{\X}^{s+1,p}_{\mathcal{D},\sigma}(\RR^n_+)
         \end{align*}
         and an estimate
         \begin{align*}
             \lVert \AA_\mathcal{D}^{\sfrac{1}{2}} \uu\rVert_{\dot{\X}^{s,p}(\RR^n_+)}\sim_{s,p,n} \lVert \nabla \uu\rVert_{\dot{\X}^{s,p}(\RR^n_+)},\qquad \forall \uu \in \dot{\D}^{s,p}(\AA_\mathcal{D}^{\sfrac{1}{2}},\RR^n_+).
         \end{align*}
         One also has a bounded injective map
         \begin{align*}
             \AA_\mathcal{D}^{\sfrac{1}{2}}\,:\,\dot{\X}^{s+1,p}_{\mathcal{D},\sigma}(\RR^n_+)\longrightarrow\dot{\X}^{s,p}_{\mathfrak{n},\sigma}(\RR^n_+).
         \end{align*}
\end{enumerate}
\end{proposition}

\begin{proof} \textbf{Step 1:} The bounded holomorphic functional calculus whenever $p\in[1,\infty]$, $-1+\sfrac{1}{p}<s<\sfrac{1}{p}$.

We start recalling that the (negative) Laplacian admits bounded holomorphic functional calculus of angle $0$ on $\dot{\X}^{s,p}(\RR^n)$. More precisely, for all $\mu\in[0,\pi)$, all $\zeta\in(\mu,\pi)$, all $f\in\mathbf{H}^{\infty}(\Sigma_\zeta)$, and all $u\in\dot{\X}^{s,p}(\RR^n)$,
\begin{align*}
    \lVert f(-\Delta) u\rVert_{\dot{\X}^{s,p}(\RR^n)}\lesssim_{p,s,n}^{\zeta} \lVert f\rVert_{\L^{\infty}(\Sigma_\zeta)} \lVert u\rVert_{\dot{\X}^{s,p}(\RR^n)}.
\end{align*}
Now, thanks to the resolvent formula \eqref{eq:ResolventEqualityUkai} proven to be valid on $\dot{\X}^{s,p}_{\mathcal{D},\sigma}(\RR^n_+)$ in Step 2, by linearity, for all $f\in\mathbf{H}^{\infty}_0(\Sigma_\zeta)$, we can make sense of
\begin{align*}
    f(\AA_\mathcal{D}) &= \frac{1}{2\pi i}\int_{\partial\Sigma_\zeta} f(z)(z\I-\AA_\mathcal{D})^{-1}\,\d z\\
    &= \frac{1}{2\pi i}\int_{\partial\Sigma_\zeta} f(z)\R_{\RR^n_+}\Gamma_{\sigma}^{\mathcal{U}}(z\I+\Delta)^{-1}\mathrm{E}_{0,\sigma}^{\mathcal{U}}\,\d z\\
    &= \R_{\RR^n_+}\Gamma_{\sigma}^{\mathcal{U}}\left(\frac{1}{2\pi i}\int_{\partial\Sigma_\zeta} f(z)(z\I+\Delta)^{-1}\,\d z\right) \mathrm{E}_{0,\sigma}^{\mathcal{U}}\\
    &=\R_{\RR^n_+}\Gamma_{\sigma}^{\mathcal{U}}f(-\Delta)\mathrm{E}_{0,\sigma}^{\mathcal{U}}.
\end{align*}
So by the definition of function spaces by restriction and by Lemma~\ref{lem:UkaiOperators}, for all $\uu\in\dot{\X}^{s,p}_{\mathfrak{n},\sigma}(\RR^n_+)$, it holds
\begin{align*}
    \lVert f(\AA_\mathcal{D}) \uu\rVert_{\dot{\X}^{s,p}(\RR^n_+)}\lesssim_{p,s,n} \lVert f(-\Delta)\rVert_{\dot{\X}^{s,p}(\RR^n)\rightarrow\dot{\X}^{s,p}(\RR^n)} \lVert \mathrm{E}_{0,\sigma}^{\mathcal{U}} \uu\rVert_{\dot{\X}^{s,p}(\RR^n)}  \lesssim_{p,s,n}^{\zeta}  \lVert f\rVert_{\L^{\infty}(\Sigma_\zeta)} \lVert \uu\rVert_{\dot{\X}^{s,p}(\RR^n_+)}.
\end{align*}

\textbf{Step 2:} Functional Calculus of higher regularity by the squareroot property. We start noticing the particular formula initially deduced by Ukaï \eqref{eq:ResolventEqualityUkai}, but now valid on $\dot{\X}^{s,p}_{\mathfrak{n},\sigma}(\RR^n_+)$:
\begin{align*}
    e^{-t\AA_\mathcal{D}} = [\Gamma_{\sigma}^{\mathcal{U}}e^{t\Delta}\mathrm{E}_{0,\sigma}^{\mathcal{U}}\cdot]_{|_{\RR^n_+}}.
\end{align*}
The unbounded operator $\AA_\mathcal{D}$ being $0$-sectorial on $\dot{\X}^{s,p}_{\mathfrak{n},\sigma}(\RR^n_+)$. Hence, the inverse square-root of $\AA_\mathcal{D}$ is given by the formula
\begin{align*}
    \AA_\mathcal{D}^{-\sfrac{1}{2}} = \frac{1}{\sqrt{\pi}} \int_{0}^{\infty} e^{-t\AA_\mathcal{D}} \frac{\d\,  t}{\sqrt{t}},
\end{align*}
valid on $\R^{s,p}(\AA_\mathcal{D}^{\sfrac{1}{2}})$, the range of $\AA_\mathcal{D}^{\sfrac{1}{2}}$ on $\dot{\X}^{s,p}_{\mathcal{D},\sigma}(\RR^n_+)$, by M${}^{\text{c}}$Intosh's approximation Theorem \cite[Theorem~5.2.6]{bookHaase2006}. We prove that for all $\ff\in\dot{\X}^{s,p}_{\mathcal{D},\sigma}(\RR^n_+)$, $-1+\sfrac{1}{p}<s<\sfrac{1}{p}$, $q<\infty$, that
\begin{align*}
    \lVert \nabla \AA_\mathcal{D}^{-\sfrac{1}{2}} \ff \rVert_{\dot{\X}^{s,p}(\RR^n_+)}\lesssim_{p,s,n} \lVert \ff \rVert_{\dot{\X}^{s,p}(\RR^n_+)}.
\end{align*}

Since $q$ is finite, one approximates $\ff$ by $(\ff_N)_{N\in\NN}$, defined for all $N\in\NN$ by
\begin{align*}
    \ff_N := \Big[\Gamma_{\sigma}^\mathcal{U}\Big[ \sum_{|j|\leqslant N} \dot{\Delta}_j \E_{0,\sigma}^{\mathcal{U}}\ff\Big]\Big]_{|_{\RR^n_+}},
\end{align*}
which converges to $\ff$ by Lemma~\ref{lem:UkaiOperators} and \cite[Proposition~2.8]{Gaudin2023Lip}. By the functional calculus identity \eqref{eq:ResolventEqualityUkai} and the finiteness of the sum, one obtains $\ff_N \in{\D}^{s,p}(\AA_\mathcal{D})\cap{\R}^{s,p}(\AA_\mathcal{D})$. 
By the definition of function spaces by restriction and Lemma~\ref{lem:UkaiOperators}, it holds
\begin{align*}
    \lVert\nabla\AA_\mathcal{D}^{-\sfrac{1}{2}} \ff_N\rVert_{\dot{\X}^{s,p}(\RR^n_+)} &\leqslant \liminf_{\varepsilon\rightarrow0} \,\lVert\nabla(\varepsilon\I+\AA_\mathcal{D})^{-\sfrac{1}{2}} \ff_N\rVert_{\dot{\X}^{s,p}(\RR^n_+)}\\
    &\leqslant \liminf_{\varepsilon\rightarrow0} \,\lVert\nabla \Gamma^{\mathcal{U}}_{\sigma}(\varepsilon\I-\Delta)^{-\sfrac{1}{2}}\E_{0,\sigma}^\mathcal{U} \ff_N\rVert_{\dot{\X}^{s,p}(\RR^n)}\\
    &\lesssim_{p,s,n}  \liminf_{\varepsilon\rightarrow0} \,\lVert(\varepsilon\I-\Delta)^{-\sfrac{1}{2}}\E_{0,\sigma}^\mathcal{U} \ff_N\rVert_{\dot{\X}^{s+1,p}(\RR^n)}\\
    &\lesssim_{p,s,n}  \lVert\E_{0,\sigma}^\mathcal{U} \ff_N\rVert_{\dot{\X}^{s,p}(\RR^n)}\\
    &\lesssim_{p,s,n}  \lVert\ff_N\rVert_{\dot{\X}^{s,p}(\RR^n_+)}.
\end{align*}
From this point, one can extend the map by density. The case $q=\infty$ follows by real interpolation.
Now, we prove that, including the case $q=\infty$,
\begin{align*}
    \lVert \AA_\mathcal{D}^{\sfrac{1}{2}} \ff \rVert_{\dot{\X}^{s,p}(\RR^n_+)}\lesssim_{p,s,n} \lVert \nabla \ff \rVert_{\dot{\X}^{s,p}(\RR^n_+)}.
\end{align*}
But one can proceed similarly, writing again up to a restriction on $\RR^n_+$
\begin{align*}
    \AA_\mathcal{D}^{\sfrac{1}{2}} \ff_N&= \frac{1}{\sqrt{\pi}} \int_{0}^{\infty} \AA_\mathcal{D} e^{-t\AA_\mathcal{D}}\ff_N \frac{\d\,  t}{\sqrt{t}}\\
     &= \frac{1}{\sqrt{\pi}} \int_{0}^{\infty} \Big[\Gamma_{\sigma}^\mathcal{U} (-\Delta) e^{t\Delta}\E_{0,\sigma}^{\mathcal{U}}\ff_N\Big] \frac{\d\,  t}{\sqrt{t}}
     \\
     &= \frac{1}{\sqrt{\pi}} \int_{0}^{\infty} \Big[\Gamma_{\sigma}^\mathcal{U} (-\Delta)^{\sfrac{1}{2}} e^{t\Delta}  (-\Delta)^{\sfrac{1}{2}} \E_{0,\sigma}^{\mathcal{U}}\ff_N\Big] \frac{\d\,  t}{\sqrt{t}}\\
     &= \Big[\Gamma_{\sigma}^\mathcal{U}(-\Delta)^{\sfrac{1}{2}} \E_{0,\sigma}^{\mathcal{U}}\ff_N\Big].
\end{align*}
By the definition of function spaces by restriction and Lemma~\ref{lem:UkaiOperators},
it holds
\begin{align*}
    \lVert \AA_\mathcal{D}^{\sfrac{1}{2}} \ff \rVert_{\dot{\X}^{s,p}(\RR^n_+)}\leqslant \liminf_{N\rightarrow\infty}\,\lVert \AA_\mathcal{D}^{\sfrac{1}{2}} \ff_{N} \rVert_{\dot{\X}^{s,p}(\RR^n_+)} &\leqslant \liminf_{N\rightarrow\infty}\,\lVert \Gamma_{\sigma}^\mathcal{U}(-\Delta)^{\sfrac{1}{2}} \E_{0,\sigma}^{\mathcal{U}}\ff_N \rVert_{\dot{\X}^{s,p}(\RR^n)}\\
    &\lesssim_{p,s,n} \liminf_{N\rightarrow\infty}\,\lVert (-\Delta)^{\sfrac{1}{2}} \E_{0,\sigma}^{\mathcal{U}}\ff_N  \rVert_{\dot{\X}^{s,p}(\RR^n)}\\
    &\lesssim_{p,s,n} \lVert \nabla \ff \rVert_{\dot{\X}^{s,p}(\RR^n_+)}.
\end{align*}
Then, we obtain the case of the spaces $\dot{\Y}^{\uptau,r}(\RR^n_+)$, $\uptau\in(-1+\sfrac{1}{r},1+\sfrac{1}{r})$, $s\neq\sfrac{1}{r}$, $r,\kappa\in[1,\infty]$ in \textit{(i)}, since we have obtained
\begin{align*}
    \lVert \uu\rVert_{\dot{\X}^{s+1,p}(\RR^n_+)} \sim_{p,s,n} \lVert \AA_{\mathcal{D}}^{\sfrac{1}{2}} \uu\rVert_{\dot{\X}^{s,p}(\RR^n_+)},
\end{align*}
for all $\uu\in\dot{\X}^{s+1}_{\mathcal{D},\sigma}(\RR^n_+)$. This finishes the proof.
\end{proof}

\subsection{Improvement of known results: the steady problem, in middle and low regularity and a consequence}

Now, we deal with the steady-state problem.

\begin{theorem}\label{thm:MetaThmSteadyDirichletStokesRn+}Let $p,q\in[1,\infty]$, $s\in (-2+{\sfrac{1}{p}},{\sfrac{1}{p}})$, $s\neq-1+\sfrac{1}{p}$. Let $\ff\in {\dot{\X}}^{s,p}(\RR^n_+,\CC^n)$, $g\in\dot{\X}^{s+1,p}(\RR^n_+,\CC)$ and $\hh\in \dot{\X}_{\partial}^{s+2-{\sfrac{1}{p}},p}(\partial\RR^n_+,\CC^{n})$ then the non-homogeneous Stokes--Dirichlet steady-state problem
    \begin{equation*}\tag{DS${}_0$}\label{eq:SteadyDirStokesSystemtn+}
    \left\{ \begin{array}{rllr}
         - \Delta \uu +\nabla \mathfrak{p} &= \ff \text{, }&&\text{ in } \RR^n_+\text{,}\\
        \div \uu &= g\text{, } &&\text{ in } \RR^n_+\text{,}\\
        \uu_{|_{\partial\RR^n_+}} &=\hh\text{, } &&\text{ on } \partial\RR^n_+\text{.}
    \end{array}
    \right.
\end{equation*}
         admits at most one solution $(\uu,\nabla \mathfrak{p})\in {\dot{\X}}^{s+2,p}(\RR^n_+,\CC^n) \times {\dot{\X}}^{s,p}(\RR^n_+,\CC^n)$, and it necessarily satisfies the estimates,
        \begin{align}\label{eq:FullSteadyEstRn+}
            \lVert \nabla^2 \uu\rVert_{{\dot{\X}}^{s,p}(\RR^n_+)} + \lVert \nabla \mathfrak{p}\rVert_{{\dot{\X}}^{s,p}(\RR^n_+)} \nonumber\lesssim_{p,n,s} \lVert \ff\rVert_{{\dot{\X}}^{s,p}(\RR^n_+)} +\lVert \nabla g\rVert_{{\dot{\X}}^{s,p}(\RR^n_+)}  +\lVert \hh\rVert_{\dot{\X}_\partial^{s+2-1/p,p}(\partial\RR^n_+)}.
        \end{align}
Furthermore,
\begin{itemize}
    \item if $\ff\in\dot{\B}^{\alpha_0}_{r_0,\kappa_0}(\RR^n_+,\CC^n)$, $g\in\dot{\B}^{\alpha_1+1}_{r_1,\kappa_1}(\RR^n_+,\CC)$  and $\hh\in\dot{\B}^{\alpha_2+2-\sfrac{1}{r_2}}_{r_2,\kappa_2}(\partial\RR^n_+,\CC^n)$ for some $r_j,\kappa_j\in[1,\infty]$, $\alpha_j\in(-2+\sfrac{1}{r_j},\sfrac{1}{r_j})$ satisfying \hyperref[AssumptionCompletenessExponents]{$(\mathcal{C}_{\alpha_j+2,r_j,\kappa_j})$} for all $j\in\llb 0,2\rrb$, then the solution always exists;

    \item a similar result holds for the endpoint Besov spaces $\dot{\BesSmo}^{s}_{p,\infty}$, $\dot{\B}^{s,0}_{\infty,q}$, and  $\dot{\BesSmo}^{s,0}_{\infty,\infty}$;

    \item in the case $p=q\in(1,\infty)$, the result still holds replacing the spaces $(\dot{\B}^{\bullet}_{p,q}(\RR^n_+),\dot{\B}^{\bullet}_{p,q}(\partial\RR^n_+))$ by $(\dot{\H}^{\bullet,p}(\RR^n_+),\dot{\B}^{\bullet}_{p,p}(\partial\RR^n_+))$.
\end{itemize}
\end{theorem}

\begin{proof} \textbf{Step 1:} Uniqueness when $s>-1+{\sfrac{1}{p}}$. Let $(\vv,\nabla \mathfrak{q})\in\dot{\B}^{s+2}_{p,\infty}(\RR^n_+,\CC^n)\times \dot{\B}^{s}_{p,\infty}(\RR^n_+,\CC^n)$ be a solution of
\begin{equation*}\tag{DS${}_00$}\label{eq:SteadyDirStokesSystemtn+0}
    \left\{ \begin{array}{rllr}
         - \Delta \vv +\nabla \mathfrak{q} &= \mathbf{0} \text{, }&&\text{ in } \RR^n_+\text{,}\\
        \div \vv &= 0\text{, } &&\text{ in } \RR^n_+\text{,}\\
        \vv_{|_{\partial\RR^n_+}} &=\mathbf{0}\text{, } &&\text{ on } \partial\RR^n_+\text{.}
    \end{array}
    \right.
\end{equation*}
By taking the divergence of the first equation, for $\mathfrak{q}$ determined up to a constant, it should satisfies
\begin{equation*}
    \left\{ \begin{array}{rllr}
         - \Delta \mathfrak{q} &= 0 \text{, }&&\text{ in } \RR^n_+\text{,}\\
        \partial_{x_n}\mathfrak{q}_{|_{\partial\RR^n_+}} &= \Delta \vv_n{}_{|_{\partial\RR^n_+}}\text{, }  &&\text{ on } \partial\RR^n_+\text{.}
    \end{array}
    \right.
\end{equation*}
Consequently, one has $\nabla\mathfrak{q}(\cdot',x_n)=-S e^{-x_n(-\Delta')^\frac{1}{2}}[\Delta \vv_n{}_{|_{\partial\RR^n_+}}]+e^{-x_n(-\Delta')^\frac{1}{2}}[\Delta \vv_n{}_{|_{\partial\RR^n_+}}]\mathfrak{e}_n$, so that applying the Uka\"{i} operator $\mathrm{E}_{0,\sigma}^{\mathcal{U}}$ to $\nabla \mathfrak{q}$, one obtains
\begin{align*}
    \mathrm{E}_{0,\sigma}^{\mathcal{U}}\nabla \mathfrak{q} =\mathbf{0}.
\end{align*}
Therefore, applying $(-\Delta')^\frac{1}{2}\mathrm{E}_{0,\sigma}^{\mathcal{U}}$ to the first equation of \eqref{eq:SteadyDirStokesSystemtn+0}, we deduce
\begin{align*}
    (-\Delta')^\frac{1}{2}\mathrm{E}_{0,\sigma}^{\mathcal{U}}\Delta \vv =\mathbf{0}.
\end{align*}
On the other hand, due to $\vv_{|_{\partial\RR^n_+}}=0$, it holds in $\S'(\RR^n)$:
\begin{align*}
    (-\Delta')^\frac{1}{2}\mathrm{E}_{0,\sigma}^{\mathcal{U}}(-\Delta \vv) = \begin{bmatrix}(-\Delta')^\frac{1}{2}\mathbf{I}_{n-1}& {\nabla'}\\ {-\div'} & (-\Delta')^\frac{1}{2}\end{bmatrix}\Delta \E_{\mathcal{D}}\vv = \begin{pmatrix}-\Delta\E_{\mathcal{D}}[(-\Delta')^\frac{1}{2}\vv'+\nabla'\vv_n]\\ \Delta\E_{\mathcal{D}}[-\div'\vv'+(-\Delta')^\frac{1}{2}\vv_n]\end{pmatrix}
\end{align*}
Since each element belongs to $\S'_h(\RR^n)$ as a byproduct of our assumption on the solution, it holds that
\begin{align*}
    \E_{\mathcal{D}}[-\div'\vv'+(-\Delta')^\frac{1}{2}\vv_n]=0.
\end{align*}
Due to the divergence free condition on $\vv$, one obtains
\begin{align*}
    \partial_{x_n}\vv_n+(-\Delta')^\frac{1}{2}\vv_n=0
\end{align*}
so, applying $\partial_{x_n}-(-\Delta')^\frac{1}{2}$, we deduce that
\begin{align*}
    \Delta \vv_n =0.
\end{align*}
By uniqueness for the Dirichlet Laplacian, this yields $\vv_n =0$. Thus, getting back on the pressure, we obtain $$\nabla\mathfrak{q}(\cdot',x_n)=-S e^{-x_n(-\Delta')^\frac{1}{2}}[\Delta \vv_n{}_{|_{\partial\RR^n_+}}]+e^{-x_n(-\Delta')^\frac{1}{2}}[\Delta \vv_n{}_{|_{\partial\RR^n_+}}]\mathfrak{e}_n =\mathbf{0}.$$
Thereby, $\vv'$ is a solution to the homogeneous Dirichlet Laplacian with a $\mathbf{0}'$ forcing term. This implies $\vv'=0$.

\textbf{Step 2:} Existence by a solution formula for $-2+{\sfrac{1}{p}}<s<{\sfrac{1}{p}}$. We give a representation formula for the solution. In order to fix the ideas, we take $\ff\in\dot{\B}^{s}_{p,q}(\RR^n_+,\CC^n)$, $q<\infty$, $g\in\dot{\B}^{s+1}_{p,q}(\RR^n_+,\CC)$ and $\hh\in\dot{\B}^{s+2-\sfrac{1}{p}}_{p,q}(\RR^n_+,\CC^n)$. Up to proceed by regularization of the data $\ff$, $g$ and $\hh$ as in the proof of Proposition~\ref{prop:GenNeuPbRn+}, by defining for $N\in\NN$,
\begin{itemize}
    \item $\ff_N := \left[\sum_{|j|\leqslant N} \dot{\Delta}_j\E_{\mathcal{D}} \ff\right]_{|_{\RR^n_+}}$;
    \item $g_N := \left[\sum_{|j|\leqslant N} \dot{\Delta}_j\E_{\mathcal{N}} g\right]_{|_{\RR^n_+}}$;
    \item $\hh_N := \sum_{|j|\leqslant N} \dot{\Delta}'_j\hh$.
\end{itemize}
which are all well-defined by Lemma~\ref{lem:ExtDirNeuRn+} and converging, respectively, towards $\ff$, $g$ and $\hh$ by \cite[Proposition~2.8]{Gaudin2023Lip}, we can omit the subscript $N$.

We look for a solution $\uu=\uu_1+\uu_2$, \eqref{eq:SteadyDirStokesSystemtn+} where $\uu_2 = (-\Delta_\mathcal{D})^{-1}\ff + (-\Delta_{\mathcal{D},\partial})^{-1}\hh$ with $ \uu_2\in\dot{\B}^{s+2}_{p,q}(\RR^n_+,\CC^n)$ by Proposition~\ref{prop:DirPbRn+}, it amounts to solve
\begin{equation*}\tag{DS${}^{(1)}_0$}\label{eq:ProofSteadyDirStokesSystemRn+1}
    \left\{ \begin{array}{rllr}
         - \Delta \uu_1 +\nabla \mathfrak{p} &= 0 \text{, }&&\text{ in } \RR^n_+\text{,}\\
        \div \uu_1 &= \tilde{g} \text{, } &&\text{ in } \RR^n_+\text{,}\\
        \uu_1{}_{|_{\partial\RR^n_+}} &=0\text{, } &&\text{ on } \partial\RR^n_+\text{.}
    \end{array}
    \right.
\end{equation*}
where $\tilde{g} := g - \div(\uu_2)$. Now, following the idea in Step 3.2 from the proof of Theorem~\ref{thm:MetaThmSteadyDirichletStokesRn+}, we write $\nabla\mathfrak{p}= \nabla\pi+\nabla\tilde{\pi}$, $\uu_1=\ww + \tilde{\ww}$, where $\tilde{\ww}=-[\nabla(-\Delta)^{-1}\E_{\mathcal{N}}\tilde{g}]_{|_{\RR^n_+}}$ and $\nabla\tilde{\pi} = -[\nabla\E_{\mathcal{N}}\tilde{g}]_{|_{\RR^n_+}}$, and one has necessarily $\tilde{\ww}_n{}_{|_{\partial\RR^n_+}}=0$. This leaves obviously to solve the following system
\begin{equation*}
    \left\{ \begin{array}{rllr}
         - \Delta \ww +\nabla \pi &= 0 \text{, }&&\text{ in } \RR^n_+\text{,}\\
        \div \ww &= 0 \text{, } &&\text{ in } \RR^n_+\text{,}\\
        \ww{}_{|_{\partial\RR^n_+}} &=(-\tilde{\ww}'{}_{|_{\partial\RR^n_+}} ,0)\text{, } &&\text{ on } \partial\RR^n_+\text{.}
    \end{array}
    \right.
\end{equation*}
The ansatz for $\ww=(\ww',\ww_n)$ and $\nabla \pi$ are given, for $x_n\geqslant0$, and $\tilde{\hh}'= -\tilde{\ww}'{}_{|_{\partial\RR^n_+}}$ by
\begin{align*}
    \ww'(\cdot,x_n)&:= [e^{-x_n(-\Delta')^{{\sfrac{1}{2}}}}-x_n(-\Delta')^{{\sfrac{1}{2}}}e^{-x_n(-\Delta')^{{\sfrac{1}{2}}}}](\I-\PP')\tilde{\hh}' + e^{-x_n(-\Delta')^{{\sfrac{1}{2}}}}\PP'\tilde{\hh}';\\
   \ww_n (\cdot,x_n)&:= -x_n e^{-x_n(-\Delta')^{{\sfrac{1}{2}}}}\div' \tilde{\hh}';\\
   \nabla \pi (\cdot,x_n) &:= 2\nabla e^{-x_n(-\Delta')^{{\sfrac{1}{2}}}}\div' \tilde{\hh}'.
\end{align*}
One can obtain the estimates applying Proposition~\ref{prop:PoissonSemigroup3}.

\textbf{Step 3:} Uniqueness if $s\in(-2+{\sfrac{1}{p}},-1+{\sfrac{1}{p}})$. We proceed by a duality argument. Let $(\vv,\nabla\mathfrak{q})\in\dot{\B}^{s+2}_{p,\infty}(\RR^n_+,\CC^n)\times \dot{\B}^{s}_{p,\infty}(\RR^n_+,\CC^n)$ be such that it satisfies \eqref{eq:SteadyDirStokesSystemtn+0}. In particular, $\vv\in\dot{\B}^{s+2,\sigma}_{p,\infty,\mathcal{D}}(\RR^n_+)$, and for all $\varphi\in\Ccinfty({\RR^n_+},\CC^n)\subset \dot{\B}^{-s}_{p',1}(\RR^n_+,\CC^n)$,
\begin{align*}
    \langle -\Delta \vv + \nabla \mathfrak{q},\, \varphi\rangle_{\RR^{n}_+} = 0
\end{align*}
In particular if $\varphi\in\Ccinftydiv(\RR^n_+)$, one obtains
\begin{align*}
    \langle  \vv ,\, -\Delta\varphi\rangle_{\RR^{n}_+} = 0
\end{align*}
Since $\vv\in \dot{\B}^{s+2,\sigma}_{p,\infty,\mathcal{D}}(\RR^n_+)$ with ${\sfrac{1}{p}}<s+2<1+{\sfrac{1}{p}}$, one can approximate $\vv$ (weakly-$\ast$) by elements from $\Ccinftydiv(\RR^n_+)$, (note that $-s\in(\sfrac{1}{p'},1+\sfrac{1}{p'})$), in order to obtain for all $\varphi\in\dot{\B}^{-s,\sigma}_{p',1,\mathcal{D}}(\RR^n_+)$ and all $\nabla \Pi \in \dot{\B}^{-2-s}_{p',1}(\RR^n_+,\CC^n)$,
\begin{align*}
     \langle  \vv ,\, -\Delta\varphi+\nabla \Pi\rangle_{\RR^{n}_+} = 0
\end{align*}
By the previous Step 2, and since we are dealing with a Besov space of third index $1$, the range of the map $(\varphi,\nabla\Pi)\mapsto -\Delta \varphi+\nabla \Pi$ is strongly dense in $\dot{\B}^{-2-s}_{p',1}(\RR^n_+,\CC^n)$, thus it holds that $\vv =0$, and then $\nabla\mathfrak{q}=0$, yielding the uniqueness.
\end{proof}

\begin{corollary}\label{cor:HomogeneousDomainStokesRn+}Let $p,q\in[1,\infty]$, $q<\infty$, $-1+{\sfrac{1}{p}}<s<{\sfrac{1}{p}}$. Then the unbounded operator $(\dot{\D}^{s}_{p,q}(\AA_\mathcal{D}),{\AA}_\mathcal{D})$ on $\dot{\B}^{s,\sigma}_{p,q,\mathfrak{n}}(\RR^n_+,\CC^n)$ satisfies Assumptions \ref{asmpt:homogeneousdomaindef} and \ref{asmpt:homogeneousdomainintersect} and the homogeneous domain of $\AA_\mathcal{D}$ satisfies
\begin{align*}
    \dot{\D}^{s}_{p,q}(\mathring{\AA}_\mathcal{D}) = \dot{\B}^{s+2,\sigma}_{p,q,\mathcal{D}}(\RR^n_+,\CC^n),
\end{align*}
with equivalence of norms.
In particular, one has
\begin{align*}
    \dot{\B}^{s+2,\sigma}_{p,q,\mathcal{D}}(\RR^n_+,\CC^n) = \overline{{\B}^{s+2,\sigma}_{p,q,\mathcal{D}}(\RR^n_+,\CC^n)}^{\lVert\cdot\rVert_{\dot{\B}^{s+2}_{p,q}(\RR^n_+)}}.
\end{align*}
Furthermore,
\begin{itemize}
    \item The result still holds for $(\dot{\mathcal{D}}^{s}_{p,\infty},\dot{\BesSmo}^{s,\sigma}_{p,\infty,\mathfrak{n}},\dot{\BesSmo}^{s+2,\sigma}_{p,\infty,\mathcal{D}})$ and $(\dot{\D}^{s,0}_{\infty,q},\dot{\B}^{s,0,\sigma}_{\infty,q,\mathfrak{n}},\dot{\B}^{s+2,0,\sigma}_{\infty,q,\mathcal{D}})$.
    \item The result still holds for Sobolev spaces replacing $(\dot{\D}^{s}_{p,q},\dot{\B}^{s,\sigma}_{p,q,\mathfrak{n}},\dot{\B}^{s+2,\sigma}_{p,q,\mathcal{D}})$ by $(\dot{\D}^{s}_p,\dot{\H}^{s,p}_{\mathfrak{n},\sigma},\dot{\H}^{s+2,p}_{\mathcal{D},\sigma})$, provided $1<p<\infty$.
\end{itemize}
\end{corollary}

\begin{remark} From the proof below, in the case of $(\dot{\D}^{s}_{p,\infty},\dot{\B}^{s,\sigma}_{p,\infty,\mathfrak{n}},\dot{\B}^{s+2,\sigma}_{p,\infty,\mathcal{D}})$,  one only obtains
    \begin{align*}
        \dot{\BesSmo}^{s+2,\sigma}_{p,\infty,\mathcal{D}}(\RR^n_+,\CC^n)\subseteq \dot{\D}^{s}_{p,\infty}(\mathring{\AA}_\mathcal{D}) \subseteq \dot{\B}^{s+2,\sigma}_{p,\infty,\mathcal{D}}(\RR^n_+,\CC^n),
    \end{align*}
    and similarly for $(\dot{\D}^{s,0}_{\infty,\infty},\dot{\B}^{s,0,\sigma}_{\infty,\infty,\mathfrak{n}},\dot{\B}^{s+2,0,\sigma}_{\infty,\infty,\mathcal{D}})$. The issue being that the domain is only weakly-$\ast$ dense for the considered ground space.
\end{remark}

\begin{proof} By Theorem~\ref{thm:MetaThmSteadyDirichletStokesRn+}, for all $\uu\in\dot{\D}^{s}_{p,q}({\AA}_\mathcal{D})$, one has
\begin{align*}
    \lVert \AA_\mathcal{D}\uu\rVert_{\dot{\B}^{s}_{p,q}(\RR^n_+)} \sim_{p,s,n}\lVert \nabla ^2 \uu\rVert_{\dot{\B}^{s}_{p,q}(\RR^n_+)} \sim_{p,s,n}\lVert \uu\rVert_{\dot{\B}^{s+2}_{p,q}(\RR^n_+)}.
\end{align*}
Thus, to check the Assumption \ref{asmpt:homogeneousdomaindef}, we choose ${\Y}=\dot{\B}^{s+2}_{p,q}(\RR^n_+,\CC^n)$.

The inclusion $\dot{\D}^{s}_{p,q}(\mathring{\AA}_\mathcal{D})\subset \dot{\B}^{s+2,\sigma}_{p,q,\mathcal{D}}(\RR^n_+)$ holds by construction. Indeed, if $(\uu_\ell)_{\ell\in\NN}\subset\dot{\D}^{s}_{p,q}({\AA}_\mathcal{D})$ converges towards an element $\uu\in\dot{\B}^{s+2}_{p,q}(\RR^n_+,\CC^n)$, then it holds by strong continuity of the trace operator, Theorem~\ref{thm:tracesRn+}, that $\uu_{|_{\partial\RR^n_+}}=0$, since $\uu_\ell {}_{|_{\partial\RR^n_+}}=0$, for all $\ell\in\NN$. It remains to check the divergence free condition for $\uu$, which holds trivially since $0=\div \uu_\ell$ converges to $\div \uu$ in $\dot{\B}^{s+1}_{p,q}(\RR^n_+,\CC)$, and thus $\div \uu =0$.

Hence, we just prove the reverse inclusion $\dot{\B}^{s+2,\sigma}_{p,q,\mathcal{D}}(\RR^n_+,\CC^n)\subset \dot{\D}^{s}_{p,q}(\mathring{\AA}_\mathcal{D})$. Let $\uu\in \dot{\B}^{s+2,\sigma}_{p,q,\mathcal{D}}(\RR^n_+,\CC^n)$, we set $\ff:= \PP_{\RR^n_+}(-\Delta_{\mathcal{D}} \uu) = -\Delta \uu + \nabla \mathfrak{p} \in  \dot{\B}^{s,\sigma}_{p,q,\mathfrak{n}}(\RR^n_+,\CC^n)$. By Theorem~\ref{thm:MetaThmDirichletStokesRn+}, we set for all $\varepsilon>0$, $\uu_\varepsilon:=(\varepsilon\I+\AA_\mathcal{D})^{-1}\ff\in \dot{\D}^{s}_{p,q}({\AA}_\mathcal{D})$. So that by Theorem~\ref{thm:MetaThmSteadyDirichletStokesRn+}, it holds
\begin{align*}
    \lVert \nabla^2(\uu-\uu_\varepsilon)\rVert_{\dot{\B}^{s}_{p,q}(\RR^n_+)} \lesssim_{p,s,n} \lVert \ff - \AA_\mathcal{D}(\varepsilon\I+\AA_\mathcal{D})^{-1}\ff\rVert_{\dot{\B}^{s}_{p,q}(\RR^n_+)}.
\end{align*}
If $q<\infty$, then by Lemma~\ref{lem:UkaiOperators}, and the Functional Calculus identity from Theorem~\ref{thm:MetaThmDirichletStokesRn+}, one obtains
\begin{align*}
    \lVert \ff - \AA_\mathcal{D}(\varepsilon\I+\AA_\mathcal{D})^{-1}\ff\rVert_{\dot{\B}^{s}_{p,q}(\RR^n_+)} \lesssim_{p,s,n} \lVert \E_{0,\sigma}^\mathcal{U}\ff - (-\Delta)(\varepsilon\I-\Delta)^{-1}\E_{0,\sigma}^\mathcal{U}\ff\rVert_{\dot{\B}^{s}_{p,q}(\RR^n)} \xrightarrow[\varepsilon\rightarrow0]{} 0.
\end{align*}
This yields the reverse inclusion. Therefore, Theorem~\ref{thm:MetaThmDirichletStokesRn+} giving $\dot{\D}^{s}_{p,q}({\AA}_\mathcal{D})= \dot{\B}^{s,\sigma}_{p,q,\mathfrak{n}}\cap\dot{\B}^{s+2,\sigma}_{p,q,\mathcal{D}}(\RR^n_+,\CC^n)$, it implies that the equality $\dot{\D}^{s}_{p,q}({\AA}_\mathcal{D})=\dot{\B}^{s,\sigma}_{p,q,\mathfrak{n}}(\RR^n_+)\cap\dot{\D}^{s}_{p,q}(\mathring{\AA}_\mathcal{D})$ holds trivially. Additionally, if one approximates $\ff$ by $(\ff_N)_{N\in\NN}$, defined for all $N\in\NN$ by
\begin{align*}
    \ff_N := \Big[\Gamma_{\sigma}^\mathcal{U}\Big[ \sum_{|j|\leqslant N} \dot{\Delta}_j \E_{0,\sigma}^{\mathcal{U}}\ff\Big]\Big]_{|_{\RR^n_+}},
\end{align*}
it gives us the claimed density result by Lemma~\ref{lem:UkaiOperators} and \cite[Proposition~2.8]{Gaudin2023Lip}.
\end{proof}

\subsection{The \texorpdfstring{$\L^\infty$}{Loo}-theory: a characterization of the domain and higher-order estimates.}

We can  improve the results in the $\L^\infty$-case. However, a major issue of the Stokes operator on $\L^\infty(\RR^n_+)$ (and $\C^{0}_{ub}(\overline{\RR^n_+})$) is that both the Stokes problem, and also its resolvent counterpart, does have a nontrivial nullspace.

Indeed, for any $\mathbf{a}'\in\CC^{n-1}$, any $\lambda\in\Sigma_{\pi}\cup\{0\}$, any $\kappa\in\CC$, we set
\begin{align}\label{eq:nontrivialSolStokesRn+}
    \ww_{\mathbf{a}'}(x',x_n):= \frac{1}{\lambda}(1-e^{-x_n\sqrt{\lambda}})\begin{bmatrix}
           \mathbf{a}'\\
           0
         \end{bmatrix}\,\quad\text{ and }\quad\mathfrak{q}_{\mathbf{a}'}(x',x_n):=-\mathbf{a}'\cdot x' + \kappa.
\end{align}

The couple $(\ww_{\mathbf{a}'},\nabla \mathfrak{q}_{\mathbf{a}'})\in\W^{2,\infty}_{\mathcal{D},\sigma}(\RR^n_+)\times\L^\infty(\RR^n_+,\CC^n)$ is a non-trivial couple that satisfies the system
\begin{equation*}
    \left\{ \begin{array}{rllr}
         \lambda \ww_{\mathbf{a}'} - \Delta \ww_{\mathbf{a}'} +\nabla \mathfrak{q}_{\mathbf{a}'} &= 0 \text{, }&&\text{ in } \RR^n_+\text{,}\\
        \div \ww_{\mathbf{a}'} &= 0\text{, } &&\text{ in } \RR^n_+\text{,}\\
        \ww_{\mathbf{a}'} {}_{|_{\partial\RR^n_+}} &=0\text{, } &&\text{ on } \partial\RR^n_+\text{.}
    \end{array}
    \right.
\end{equation*}

This means that on $\L^\infty_{\mathfrak{n},\sigma}(\RR^n_+)$, the Stokes operator given with domain
\begin{align*}
    \big\{ \,\uu\in \L^{\infty}_{\mathfrak{n},\sigma}(\RR^n_+)\,:\, \exists \nabla\mathfrak{p}\in\B^{-2}_{\infty,\infty}(\RR^n_+,\CC^n),\, -\Delta \uu +\nabla \mathfrak{p} \in\L^{\infty}_{\mathfrak{n},\sigma}(\RR^n_+)\,\big\}
\end{align*}
must be at least a multivalued operator, but not only, this also implies that we have a multivalued resolvent operator whenever the later exists. The graph arising from the resolvent problem is in fact a dual-coverage operator graph. It is quite astonishing: while the Stokes--Dirichlet operator on function spaces of bounded functions has been  studied for several decades, that this fact has not been that much highlighted in the literature until recently. See for instance \cite[Introduction,~Sec.~1B]{MaekawaMiuraPrange2020} for a brief discussion.

More precisely, we say that $(\uu,\mathfrak{p})\in\W^{1,1}_{\text{loc}}(\overline{\RR^n_+})\times \L^1_{\text{loc}}({\RR^n_+})$ is a weak solution to the Stokes--Dirichlet resolvent problem  with forcing term $\ff\in\L^1_{\mathfrak{n},\sigma,\text{loc}}(\overline{\RR^n_+})$ if it satisfies $\uu_{|_{\partial\RR^n_+}}=0$,
\begin{align}\label{eq:weakStokesResRn+1}
   \lambda  \langle\uu,\,  \boldsymbol{\varphi}\rangle_{\RR^n_+} + \langle\nabla\uu,\, \nabla\boldsymbol{\varphi}\rangle_{\RR^n_+}-\langle \mathfrak{p},\, \div \boldsymbol{\varphi}\rangle_{\RR^n_+} &=\langle \ff,\,\boldsymbol{\varphi}\rangle_{\RR^n_+},\qquad \forall \boldsymbol{\varphi} \in\Ccinfty({\RR^n_+},\CC^n),\,
\end{align}
and
\begin{align}\label{eq:weakStokesResRn+2}
    \langle\uu,\, \nabla{\phi}\rangle_{\RR^n_+} &= 0,\qquad \forall {\phi} \in\Ccinfty(\overline{\RR^n_+}).
\end{align}
Maekawa, Miura and Prange did obtain the following Liouville-type/uniqueness-type result:
\begin{theorem}[ {\cite[Theorem~4]{MaekawaMiuraPrange2020}} ]\label{thm:MaekawaMiruaPrangeUniqueness} Let $\mu\in[0,\pi)$ and let $\lambda\in\Sigma_\mu$. Let $(\uu,\mathfrak{p})\in\W^{1,\infty}({\RR^n_+},\CC^n)\times \L^1_{\mathrm{loc}}({\RR^n_+})$ be a solution to the weak Stokes--Dirichlet resolvent problem \eqref{eq:weakStokesResRn+1}-\eqref{eq:weakStokesResRn+2} with $\ff =0$. Then it holds that
\begin{align*}
    (\uu,\mathfrak{p})\in\{ (\ww_{\mathbf{a}'},\mathfrak{q}_{\mathbf{a}'}),\,\mathbf{a}'\in\CC^{n-1}\},
\end{align*}
where the expression of $\ww_{\mathbf{a}'}$ and $\mathfrak{q}_{\mathbf{a}'}$, provided $\mathbf{a}'\in\CC^{n-1}$, are given in \eqref{eq:nontrivialSolStokesRn+}.
\end{theorem}
Note that our definition of weak solutions is a bit different than theirs, but it can be checked that their proof still applies up to some minor modifications.

This suggests that the function space $\L^{\infty}_{\mathfrak{n},\sigma}(\RR^n_+)$ itself is not really a suitable ground space for the study of the Stokes problem considering bounded functions on the half-space. In order to reach properly a mono-valued  version of the Stokes--Dirichlet operator and its resolvent operator, a proposition is to replace
\begin{align*}
    \L^{\infty}_{\mathfrak{n},\sigma}(\RR^n_+)\quad\text{ by }\quad \L^{\infty}_{h,\mathfrak{n},\sigma}(\RR^n_+)
\end{align*}
where we recall that, by Lemma~\ref{lem:ExtDirNeuLinftyRn+},  $\L^{\infty}_{h}(\RR^n)=\L^\infty(\RR^n)\cap\S'_h(\RR^n)$ is such that $$\L^\infty_h(\RR^n)\cap[\CC\mathbbm{1}_{\RR^n_+}\oplus\CC\mathbbm{1}_{\RR^n_-}]=\{0\},$$
so that by restriction by Lemma~\ref{lem:ExtDirNeuLinftyRn+} again, $\CC^n\cap\L^{\infty}_{h,\mathfrak{n},\sigma}(\RR^n_+)=\{0\}.$

In fact, we don't need to  dive into such considerations to provide uniqueness, even if using $\L^{\infty}_{h,\mathfrak{n},\sigma}(\RR^n_+)$ as a ground space would provide uniqueness for free, according to Theorem~\ref{thm:MaekawaMiruaPrangeUniqueness}.

To fit our purpose, and to follow our previous idea, we aim to put the ``bad low frequency behavior'' away from the pressure term, which should be enough to claim uniqueness. We introduce the appropriate (complete) function  spaces for our analysis
\begin{itemize}
    \item provided $s\in\RR$, $q\in[1,\infty]$, the ``decaying'' endpoint inhomogeneous Besov spaces
    \begin{align*}
        \tilde{\B}^{s}_{\infty,q}(\RR^n_+) := \B^{s}_{\infty,q}(\RR^n_+)\cap\S'_h(\RR^n_+)
    \end{align*}
    endowed with the norm $\lVert\cdot\rVert_{\B^{s}_{\infty,q}(\RR^n_+)}$. Note that one has the algebraic equality
    \begin{align*}
        \dot{\B}^{s}_{\infty,q}(\RR^n_+) = \tilde{\B}^{s}_{\infty,q}(\RR^n_+),\quad\text{whenever }s>0,\, q\in[1,\infty];
    \end{align*}
    Again, for all $s\in\RR$, $\tilde{\B}^{s}_{\infty,q}(\RR^n_+,\CC^n)\cap\CC^n =\{0\}$.
    \item the hybrid Besov space
\begin{align*}
    \hat{\B}^{1}_{\infty,\infty}(\RR^n_+,\CC) = \left 
\{\,\mathfrak{q}\in\mathcal{D}'(\RR^n_+)\,:\,\nabla \mathfrak{q}\in\dot{\B}^{-1}_{\infty,\infty}\cap{\B}^{0}_{\infty,\infty}(\RR^n_+,\CC^n)\,\}\middle/ \CC\right.
\end{align*}
endowed with the norm
\begin{align*}
    \mathfrak{q}\longmapsto \lVert \nabla \mathfrak{q}\rVert_{\dot{\B}^{-1}_{\infty,\infty}(\RR^n_+)}+\lVert \nabla \mathfrak{q}\rVert_{{\B}^{0}_{\infty,\infty}(\RR^n_+)}.
\end{align*}

This later function space does not contain non-zero Lipschitz continuous functions with at most linear growth, and one has
\begin{align*}
    \nabla\hat{\B}^{1}_{\infty,\infty}(\RR^n_+,\CC)\subset \tilde{\B}^{0}_{\infty,\infty}(\RR^n_+,\CC^n).
\end{align*}

Additionally, it can be checked that it can be interpreted by the duality relation
\begin{align*}
    \hat{\B}^{1}_{\infty,\infty}(\RR^n_+,\CC)=(\mathring{\B}^{-1}_{1,1,0}(\RR^{n}_+,\CC))'
\end{align*}
where
\begin{align*}
    \mathring{\B}^{-1}_{1,1,0}(\RR^{n}_+,\CC)&:=-\div(\dot{\B}^{1}_{1,1,0}(\RR^n_+,\CC^n)+{\B}^{0}_{1,1,0}(\RR^n_+,\CC^n)).
\end{align*}
This yields that $\hat{\B}^{1}_{\infty,\infty}(\RR^n_+,\CC)$ is weakly-$\ast$ complete. Note also that, by construction, $\div(\Ccinfty(\RR^{n}_+,\CC^n))$ is strongly dense in $\mathring{\B}^{-1}_{1,1,0}(\RR^{n}_+)$.
\end{itemize}

The overall strategy of the next proof is guided by the following simple idea:
\begin{center}
    to (real) interpolate $\dot{\B}^{s}_{\infty,q}$ and  $\dot{\B}^{s+1}_{\infty,q}$-estimates from Theorem~\ref{thm:MetaThmDirichletStokesRn+} and Proposition~\ref{prop:MetaThmDirichletStokesRn+HigherReg} for $-1<s<0$, and use the (false) ``embedding''
    \begin{align}\label{}
        \L^{\infty}_{\mathfrak{n},\sigma}(\RR^n_+)\hookrightarrow\dot{\B}^{0,\sigma}_{\infty,\infty,0}(\RR^n_+),
    \end{align}
    to recover ${\B}^{2}_{\infty,\infty}$-regularity of "the" solution.
\end{center}
However, there are several technical obstructions:
\begin{enumerate}
    \item Since $\L^{\infty}(\RR^n)\not\subset\dot{\B}^{0}_{\infty,\infty}(\RR^n)$, due to $\L^{\infty}(\RR^n)\not\subset\S'_h(\RR^n)$, the embedding above cannot be reached.
    \item Even if above embedding was possible, the real interpolation of the $\dot{\B}^{s}_{\infty,q}$ and  $\dot{\B}^{s+1}_{\infty,q}$-estimates from Proposition~\ref{prop:MetaThmDirichletStokesRn+HigherReg} is flawed, since the  $\dot{\B}^{s+1}_{\infty,q}$-estimates are not true on the whole $\dot{\B}^{s+1}_{\infty,q}$ but rather only for elements in spaces like $\dot{\B}^{s}_{\infty,q}\cap\dot{\B}^{s+1}_{\infty,q}$.
    \item Additionally, the already known representation formulas giving existence to the resolvent problem by Uka\"{i} \cite{Ukai1987} and Desch, Hieber and Pr\"{u}ss \cite{DeshHieberPruss2001} are \textit{a priori} not the same, even the one by Uka\"{i} encounters issues for the boundedness on $\L^\infty$. So, we will need a priori to compare the two type of solutions in an appropriate way.
\end{enumerate}

The representation formula of Desch, Hieber and Pr\"{u}ss \cite[Section~4,~Theorems~4.3~\&~4.4]{DeshHieberPruss2001} yields an estimate
\begin{align*}
    |\lambda|\lVert \uu\rVert_{\L^\infty(\RR^n_+)}  \lesssim_{n,\mu} \lVert \ff\rVert_{\L^\infty(\RR^n_+)},
\end{align*}
provided $\ff\in\L^\infty_{\mathfrak{n},\sigma}(\RR^n_+)$, $\lambda\in\Sigma_{\mu}$, $\mu\in[0,\pi)$, where, following \cite[Section~2~\&~Section~3,~p.120--121]{DeshHieberPruss2001}, $\uu$ is given by a solution operator
\begin{align}\label{eq:DeshHieberPrussFormula}
    \uu = (\lambda\I-\Delta_\mathcal{D})^{-1}\ff + \mathrm{T}_\lambda \ff,
\end{align}
where for any $p\in(1,\infty]$, any $\gg\in\L^p(\RR^n_+,\CC^n)$,
\begin{align*}
    \mathrm{T}_\lambda \gg := \mathrm{T}_\lambda' \gg' + (\mathrm{T}_{\lambda,n} \gg')\mathfrak{e}_n ,
\end{align*}
with for almost every $(x',x_n)\in\RR^n_+$,
\begin{align*}
    \mathrm{T}_\lambda' \gg'(x',x_n) &:= \int_{0}^{\infty}\int_{\RR^{n-1}} \mathfrak{r}'_\lambda(x'-y',x_n,y_n)\gg'(y',y_n)\,\d y' \d y_n\\
    \mathrm{T}_{\lambda,n} \gg'(x',x_n) &:= \int_{0}^{\infty}\int_{\RR^{n-1}} \mathfrak{r}_{\lambda,n}(x'-y',x_n,y_n) \cdot \gg'(y',y_n)\,\d y' \d y_n.
\end{align*}
The integral kernel being given by
\begin{align*}
    \mathfrak{r}'_\lambda (x',x_n,y_n)&:=\frac{1}{(2\pi)^{n-1}}\int_{\RR^{n-1}}\frac{e^{-x_n|\xi'|}-e^{-x_n\sqrt{\lambda+|\xi'|^2}}}{\sqrt{\lambda+|\xi'|^2}-|\xi'|}\frac{1}{\sqrt{\lambda+|\xi'|^2}|\xi'|}\xi'\,\prescript{t}{}{\xi'}\,e^{-y_n\sqrt{\lambda+|\xi'|^2}} e^{i \,x'\cdot\xi'} \d\xi',\\
    \mathfrak{r}_{\lambda,n}(x',x_n,y_n) &:= \frac{1}{(2\pi)^{n-1}}\int_{\RR^{n-1}}\frac{e^{-x_n|\xi'|}-e^{-x_n\sqrt{\lambda+|\xi'|^2}}}{\sqrt{\lambda+|\xi'|^2}-|\xi'|}\frac{i\xi'}{\sqrt{\lambda+|\xi'|^2}}e^{-y_n\sqrt{\lambda+|\xi'|^2}} e^{i \,x'\cdot\xi'} \d\xi'.
\end{align*}

Desch, Hieber and Pr\"{u}ss did prove by kernel estimates that, on $\L^\infty_{\mathfrak{n},\sigma}(\RR^n_+)$, the holomorphic family of operators arising from {\eqref{eq:DeshHieberPrussFormula}} can be realized by a unique unbounded closed operator on $\L^\infty_{\mathfrak{n},\sigma}(\RR^n_+)$, they denote $(\D(A_{\L^\infty_\sigma}),A_{\L^\infty_\sigma})$, such that on the whole  $\L^\infty_{\mathfrak{n},\sigma}(\RR^n_+)$:
\begin{align*}
    (\lambda\I-\Delta_\mathcal{D})^{-1} + \mathrm{T}_\lambda = (\lambda\I+A_{\L^\infty_\sigma})^{-1}.
\end{align*}
For consistency, since $(\lambda\I-\Delta_\mathcal{D})^{-1} + \mathrm{T}_\lambda = (\lambda\I+\AA_\mathcal{D})^{-1}$ on $\L^{2}_{\mathfrak{n},\sigma}\cap\L^{\infty}_{\mathfrak{n},\sigma}(\RR^{n}_+)$, we still denote here the arising operator by
\begin{align*}
    (\D_\infty(\AA_\mathcal{D}),\AA_\mathcal{D})=(\D(A_{\L^\infty_\sigma}),A_{\L^\infty_\sigma}).
\end{align*}
This analysis has been pursued by Maekawa, Miura and Prange \cite{MaekawaMiuraPrange2020} for the analysis of the Stokes--Dirichlet problem in the spaces $\L^p_{\textrm{uloc}}(\RR^n_+)$, $1\leqslant p\leqslant \infty$, with refined point-wise estimates on the kernel. More precisely, they did give explicit Poisson-type kernel bounds on $\mathfrak{r}_\lambda$ that we will not use here.

Our next result improves known regularity estimates, and prove that the domain of the Stokes operator lies in a standard space of order $2$, with a suitable description of its domain as a direct consequence, see Theorem~\ref{thm:DomainStokesLInftyRn+} below. We investigate existence and uniqueness for the resolvent problem building on the work and results by  Desch, Hieber and Pr\"{u}ss \cite{DeshHieberPruss2001} and Maekawa, Miura and Prange \cite{MaekawaMiuraPrange2020}, proving also that in this case the pressure is log-Lipschitz continuous, see \cite[Definition~2.106~\&~Proposition~2.107]{bookBahouriCheminDanchin} for a definition and a suitable sufficient condition ($\B^{0}_{\infty,\infty}$-membership of the gradient).

\begin{theorem}\label{thm:StokesDirRn+Linfty} Let $\mu\in[0,\pi)$. For all $\ff\in\L^\infty_{\mathfrak{n},\sigma}(\RR^n_+)$, for all $\lambda\in\Sigma_{\mu}$, there exists a unique solution $(\uu,\mathfrak{p})\in{\B}^{2,\sigma}_{\infty,\infty,\mathcal{D}}(\RR^n_+)\times\hat{\B}^{1}_{\infty,\infty}(\RR^n_+,\CC)$ to
\begin{equation*}\tag{DS${}_\lambda$}\label{eq:DirStokesSystemtn+Infty}
    \left\{ \begin{array}{rllr}
         \lambda \uu - \Delta \uu +\nabla \mathfrak{p} &= \ff \text{, }&&\text{ in } \RR^n_+\text{,}\\
        \div \uu &= 0\text{, } &&\text{ in } \RR^n_+\text{,}\\
        \uu_{|_{\partial\RR^n_+}} &=0\text{, } &&\text{ on } \partial\RR^n_+\text{.}
    \end{array}
    \right.
\end{equation*}
satisfying the estimates
\begin{align*}
    |\lambda|\lVert \uu\rVert_{\L^\infty(\RR^n_+)} + |\lambda|^\frac{1}{2}\lVert \nabla\uu\rVert_{\L^\infty(\RR^n_+)} + \lVert ( \nabla^2\uu,\nabla{\mathfrak{p}})\rVert_{\dot{\B}^0_{\infty,\infty}(\RR^n_+)} &\lesssim_{n,\mu} \lVert \ff\rVert_{\L^\infty(\RR^n_+)},
\end{align*}
and with the additional pressure estimate
\begin{align*}
    |\lambda|^\frac{1}{2}\lVert \nabla \mathfrak{p}\rVert_{\dot{\B}^{-1}_{\infty,\infty}(\RR^n_+)}\lesssim_{n,\mu} \lVert \ff\rVert_{\L^\infty(\RR^n_+)},
\end{align*}
as well as for any $c>0$
\begin{align*}
    \left(\frac{1}{1+|\ln(\min(c,|\lambda|))|}\right)\lVert (\nabla^2\uu,\nabla \mathfrak{p})\rVert_{{\B}^0_{\infty,\infty}(\RR^n_+)} \lesssim_{n,\mu,c} \lVert \ff\rVert_{\L^\infty(\RR^n_+)}.
\end{align*}
Furthermore, 
\begin{itemize}
    \item if $\ff\in\C^{0}_{ub,\mathcal{D},\sigma}({\RR^{n}_+})$, one has $\uu\in\BesSmo^{2}_{\infty,\infty}(\RR^n_+)$;
    \item if $\ff\in\L^{\infty}_{h,\mathfrak{n},\sigma}({\RR^{n}_+})$, one has $\uu\in\tilde{\B}^{2}_{\infty,\infty}(\RR^n_+)$;
    \item if $\ff\in\C^{0}_{ub,h,\mathcal{D},\sigma}({\RR^{n}_+})$, one has $\uu\in\tilde{\BesSmo}^{2}_{\infty,\infty}(\RR^n_+)$;
    \item if $\ff\in\C^{0}_{0,0,\sigma}({\RR^{n}_+})$, one has $\uu\in\BesSmo^{2,0}_{\infty,\infty}(\RR^n_+)$.
\end{itemize}
\end{theorem}

Before proving Theorem~\ref{thm:StokesDirRn+Linfty}, we exhibit the following fundamental consequence:

\begin{theorem}\label{thm:DomainStokesLInftyRn+} The Stokes--Dirichlet operator $(\D_\infty(\AA_\mathcal{D}),\AA_\mathcal{D})$ on $\L^\infty_{\mathfrak{n},\sigma}(\RR^n_+)$ is such that its domain satisfies the following description
    \begin{align*}
        \D_\infty(\AA_\mathcal{D})=\{\, \uu\in\B^{2,\sigma}_{\infty,\infty,\mathcal{D}}(\RR^n_+)\,:\,\PP_{\RR^n_+}(-\Delta_\mathcal{D}\uu)\in\L^\infty_{\mathfrak{n},\sigma}(\RR^n_+)\,\},
    \end{align*}
    and for all $\uu\in\D_\infty(\AA_\mathcal{D})$,
    \begin{align*}
        \AA_\mathcal{D}\uu = \PP_{\RR^n_+}(-\Delta_\mathcal{D}\uu).
    \end{align*}
\end{theorem}
Before, we start the proofs, we provide several comments.
\begin{itemize}
    \item This result might seems surprising, since at the first glance the Hodge-Leray projection is ill-defined on $\L^\infty(\RR^n_+,\CC^n)$. It is because, the membership relation
    $$\text{``}\PP_{\RR^n_+}(-\Delta_\mathcal{D}\uu)\in\L^\infty_{\mathfrak{n},\sigma}(\RR^n_+)\text{''},$$
    provided the knowledge of the domain of the operator, is an unbounded one.

    \item To the best of the authors' knowledge, over the last two decades, this seems to be the  only attempt  in the literature for an explicit description of the domain of the Stokes--Dirichlet operator on $\L^\infty(\RR^n_+)$.

    \item It was already an open question to know whether it would be possible to describe the domain of the Stokes--Dirichlet operator on $\L^\infty_{\mathfrak{n},\sigma}(\RR^n_+)$ by a standard  function space of order $2$. The reasonable and expected conjecture by the community is the following
\end{itemize}
\begin{conjecture} The domain of the Stokes--Dirichlet operator on $\L^\infty_{\mathfrak{n},\sigma}(\RR^n_+)$ would be given by
\begin{align*}
    \D_\infty(\AA_\mathcal{D})=\{\, \uu\in\W^{1,\infty}_{\mathcal{D},\sigma}(\RR^n_+)\,:\,\nabla^2\uu\in\mathrm{BMO}(\RR^n_+,\CC^{n^3})\,\,\&\,\,\PP_{\RR^n_+}(-\Delta_\mathcal{D}\uu)\in\L^\infty_{\mathfrak{n},\sigma}(\RR^n_+)\,\}.
\end{align*}
\end{conjecture}
One notices that $\{\, \uu\in\W^{1,\infty}(\RR^n_+,\CC^n)\,:\,\nabla^2\uu\in\mathrm{BMO}(\RR^n_+,\CC^{n^3})\}\subset \B^{2}_{\infty,\infty}(\RR^n_+,\CC^{n})$.

\medbreak

While above conjecture is fairly well expected to be true, it would be way more surprising to know if it holds
\begin{align*}
    \D_\infty(\AA_\mathcal{D})=\{\,\uu\in \D_{\infty}(\Delta_\mathcal{D})\cap\L^\infty_{\mathfrak{n},\sigma}(\RR^n_+)\,:\,\PP_{\RR^n_+}(-\Delta_\mathcal{D}\uu)\in\L^\infty_{\mathfrak{n},\sigma}(\RR^n_+)\,\},
\end{align*}
which would imply by linearity that the pressure term satisfies $\nabla\mathfrak{p}\in\L^\infty(\RR^n_+,\CC^n)$. This is certainly something that one would not expect.

\medbreak

\noindent Now, we start the proofs of Theorems~\ref{thm:StokesDirRn+Linfty} and \ref{thm:DomainStokesLInftyRn+}. 

\begin{proof}[of Theorem~\ref{thm:DomainStokesLInftyRn+}] By Theorem~\ref{thm:StokesDirRn+Linfty} above, one has the following temporary description for the domain of the Stokes--Dirichlet operator
\begin{align*}
    \D_\infty(\AA_\mathcal{D})=\{\, \vv\in\B^{2,\sigma}_{\infty,\infty,\mathcal{D}}(\RR^n_+)\,:\,\exists \nabla\mathfrak{q} \in\tilde{\B}^{0}_{\infty,\infty}(\RR^n_+,\CC^n),\,-\Delta \vv + \nabla \mathfrak{q}\in\L^\infty_{\mathfrak{n},\sigma}(\RR^n_+)\,\}.
\end{align*}
Note that the condition "$\nabla\mathfrak{q} \in\tilde{\B}^{0}_{\infty,\infty}(\RR^n_+)$" is due to the uniqueness from Step 0 in the proof of Theorem~\ref{thm:StokesDirRn+Linfty} below, which requires only this better condition (asking for weaker localisation), uniqueness implying that $\mathfrak{q}\in\hat{\B}^{1}_{\infty,\infty}(\RR^n_+)$ which is a weaker condition (giving stronger localisation).

\textbf{Step 1:} The left to right embedding. Let $\uu\in \D_\infty(\AA_\mathcal{D})\subset \B^{2,\sigma}_{\infty,\infty,\mathcal{D}}(\RR^n_+)$, and let $\mathfrak{p} \in\hat{\B}^{1}_{\infty,\infty}(\RR^n_+)$, such that
\begin{align*}
    \ff:=-\Delta \uu + \nabla \mathfrak{p} \in\L^\infty_{\mathfrak{n},\sigma}(\RR^n_+).
\end{align*}
One has $\uu\in \B^{s+2,\sigma}_{\infty,1,\mathcal{D}}(\RR^n_+)$, all for $-1<s<0$, and then
\begin{align*}
    \Delta \uu \in \dot{\B}^{s,\sigma}_{\infty,1}(\RR^n_+).
\end{align*}
On the other hand, since  $\mathfrak{p} \in\hat{\B}^{1}_{\infty,\infty}(\RR^n_+)$, it holds $\nabla\mathfrak{p} \in\dot{\B}^{0}_{\infty,\infty}\cap\dot{\B}^{-1}_{\infty,\infty}(\RR^n_+,\CC^n)$,  one also has by interpolation identities
\begin{align*}
    \nabla\mathfrak{p} \in\dot{\B}^{s}_{\infty,1}(\RR^n_+,\CC^n).
\end{align*}
Therefore,
\begin{align*}
    -\Delta \uu + \nabla \mathfrak{p}  = \ff\in\dot{\B}^{s}_{\infty,1}\cap\L^\infty_{\mathfrak{n},\sigma}(\RR^n_+,\CC^n)\subset \dot{\B}^{s,\sigma}_{\infty,1,\mathfrak{n}}(\RR^n_+).
\end{align*}
Recalling, that here $-1<s<0$, it is legitimate to apply the Hodge-Leray projection and deduce
\begin{align*}
    \PP_{\RR^n_+}(-\Delta \uu) = \PP_{\RR^n_+}(-\Delta \uu+ \nabla \mathfrak{p}) = \PP_{\RR^n_+}\ff =\ff \in\L^\infty_{\mathfrak{n},\sigma}(\RR^n_+).
\end{align*}
\textbf{Step 2:} For the reverse embedding, let $\uu\in\B^{2,\sigma}_{\infty,\infty,\mathcal{D}}(\RR^n_+)$ be such that
\begin{align*}
    \PP_{\RR^n_+}(-\Delta \uu)\in\L^\infty_{\mathfrak{n},\sigma}(\RR^n_+).
\end{align*}
So by the same argument as in Step 1, there exists $\nabla \mathfrak{p}\in\dot{\B}^{s}_{\infty,1}(\RR^n_+,\CC^n)$, for any $-1<s<0$, such that
\begin{align*}
    -\Delta \uu + \nabla \mathfrak{p} \in\dot{\B}^{s,\sigma}_{\infty,1}\cap\L^\infty_{\mathfrak{n},\sigma}(\RR^n_+).
\end{align*}
Since $\uu\in \B^{2,\sigma}_{\infty,\infty,\mathcal{D}}(\RR^n_+)$, by linearity it holds
\begin{align*}
     \nabla \mathfrak{p} \in \B^{0}_{\infty,\infty}(\RR^n_+,\CC^n).
\end{align*}
Therefore,  $\nabla \mathfrak{p} \in \B^{0}_{\infty,\infty}\cap\dot{\B}^{s,\sigma}_{\infty,1}(\RR^n_+,\CC^n)\subset \tilde{\B}^{0}_{\infty,\infty}(\RR^n_+,\CC^n)$. This ends the argument for  the reverse inclusion, and thereby the proof.
\end{proof}

\begin{proof}[of Theorem~\ref{thm:StokesDirRn+Linfty}]Before we actually start the proof: the overall spirit of the proof is to prove uniqueness then by comparing the solution given by Ukaï's representation formula to the one by the different representation formula of Desch, Hieber and Pr\"{u}ss \cite{DeshHieberPruss2001}, up to weakly-$\ast$ dense subset of $\L^\infty_{\mathfrak{n},\sigma}(\RR^n_+)$.

Furthermore, by a scaling argument, and without loss of generality, as in the proof of Theorem~\ref{thm:MetaThmDirichletStokesRn+}, we can assume $|\lambda|=1$. See the end of the proof for the non-trivial and adapted dilation argument.

\textbf{Step 0:} We prove uniqueness. Let $(\uu,\nabla\mathfrak{p})\in{\B}^{2,\sigma}_{\infty,\infty,\mathcal{D}}(\RR^n_+)\times\tilde{\B}^{0}_{\infty,\infty}(\RR^n_+,\CC^n)$ be a solution to \eqref{eq:DirStokesSystemtn+Infty} with $\ff=0$. Since such a solution is in particular a weak solution to \eqref{eq:weakStokesResRn+1}-\eqref{eq:weakStokesResRn+2} with $\ff =0$, by Theorem~\ref{thm:MaekawaMiruaPrangeUniqueness}, it holds that, for some $\mathbf{a}'\in\CC^{n-1}$,
\begin{align*}
    (\uu,\mathfrak{p})=(\ww_{\mathbf{a}'},\mathfrak{q}_{\mathbf{a}'})
\end{align*}
where the expressions of $\ww_{\mathbf{a}'}$ and $\mathfrak{q}_{\mathbf{a}'}$ are given in \eqref{eq:nontrivialSolStokesRn+}. Since $\nabla'\mathfrak{p}=\nabla'\mathfrak{q}_{\mathbf{a}'}\in\tilde{\B}^{0}_{\infty,\infty}(\RR^{n}_+,\CC^{n-1})$, this implies that $\mathbf{a}'=\mathbf{0}'$. Therefore, $\nabla\mathfrak{p}=0$ and $\uu$ is a solution to the Dirichlet Laplacian resolvent problem, which implies $\uu=0$.

\textbf{Step 1:} We provide a preparatory step, proving first a Lipschitz bound for $\uu$ and a low frequency control on the gradient of the pressure $\nabla\mathfrak{p}$. 

In this current Step 1, we aim to prove first two facts for the solution
\begin{itemize}
    \item The Lipschitz bound
    \begin{align*}
        |\lambda|^\frac{1}{2}\lVert \nabla \uu\rVert_{\L^\infty(\RR^n_+)}  \lesssim_{n,\mu} \lVert \ff\rVert_{\L^\infty(\RR^n_+)},
    \end{align*}
    which stands for the goal of Step 1.1;
    \item The low frequency control on the pressure term $\nabla \mathfrak{p} = -(\lambda\I-\Delta)\T_\lambda \ff'$:
    \begin{align*}
        |\lambda|^\frac{1}{2}\lVert \nabla \mathfrak{p}\rVert_{\dot{\B}^{-1}_{\infty,\infty}(\RR^n_+)}  \lesssim_{n,\mu} \lVert \ff\rVert_{\L^\infty(\RR^n_+)},
    \end{align*}
    which is the goal of Step 1.2.
    \end{itemize}
This will prove in particular that $(\uu,\nabla\mathfrak{p})$ given by  the Desch, Hieber and Pr\"{u}ss formula solves the resolvent problem in $\B^{-1}_{\infty,\infty}(\RR^n_+,\CC^n)$.

\textbf{Step 1.1:} For the Lipschitz bound we proceed by standard kernel estimates. By \cite[Proposition~3.1~\&~Remark~3.2]{DeshHieberPruss2001}, and following the proof of \cite[Proposition~3.4]{DeshHieberPruss2001}, one obtains --recall that here $|\lambda|=1$--, for some $\mathfrak{c}>0$, all $x_n,y_n>0$,
\begin{align*}
    \int_{\RR^{n-1}} |\nabla' \mathfrak{r}_\lambda(x',x_n,y_n)|\d x' \lesssim_{n,\mu} e^{-\mathfrak{c} y_n}\frac{x_n}{1+x_n}\frac{1}{(x_n+y_n)^2},
\end{align*}
so that
\begin{align*}
    \int_{0}^\infty\int_{\RR^{n-1}} |\nabla' \mathfrak{r}_\lambda(x',x_n,y_n)|\d x' \d y_n&\lesssim_{n,\mu} \int_{0}^\infty e^{-\mathfrak{c} y_n}\frac{x_n}{1+x_n}\frac{1}{(x_n+y_n)^2}\d y_n \\
    &\lesssim_{n,\mu} \int_{0}^\infty e^{-\mathfrak{c} x_n z}\frac{1}{1+x_n}\frac{1}{(1+z)^2}\d z\\
    &\lesssim_{n,\mu} \frac{1}{1+x_n}.
\end{align*}
Therefore by \cite[Proposition~3.3]{DeshHieberPruss2001} $\nabla'\mathrm{T}_\lambda$ is bounded as an operator from $\L^\infty(\RR^n_+,\CC^n)$ to $\L^\infty(\RR^n_+,\CC^{n\times(n-1)})$. Consequently, by Proposition~\ref{prop:HodgeResolventPbRn+Linfty}\footnote{the Hodge Laplacian on $n$-forms gives the Dirichlet Laplacian for scalar functions.} and the identity \eqref{eq:DeshHieberPrussFormula}, due to the divergence-free condition, it holds that
\begin{align*}
    \lVert \nabla' \uu\rVert_{\L^\infty(\RR^n_+)} &\leqslant \lVert \nabla' (\lambda\I-\Delta_\mathcal{D})^{-1}\ff\rVert_{\L^\infty(\RR^n_+)} + \lVert \nabla' \mathrm{T}_\lambda\ff\rVert_{\L^\infty(\RR^n_+)}\\
    &\lesssim_{n,\mu} \lVert \ff\rVert_{\L^\infty(\RR^n_+)}.
\end{align*}
Now, we want to obtain the estimates on $\partial_{x_n}\uu'=(\partial_{x_n}\uu_j)_{1\leqslant j\leqslant n-1}$.
For this, we study the linear operator $\partial_{x_n} \mathrm{T}_\lambda'$. It turns out that for almost every $x'\in\RR^{n-1}$, $x_n,y_n>0$, one can write
\begin{align*}
    \partial_{x_n}&\mathfrak{r}'_\lambda (x',x_n,y_n)\\
    &=\frac{1}{(2\pi)^{n-1}}\int_{\RR^{n-1}}\frac{|\xi'|e^{-x_n|\xi'|}}{\sqrt{\lambda+|\xi'|^2}}\frac{\xi'\,\prescript{t}{}{\xi'}}{|\xi'|^2}\,e^{-y_n\sqrt{\lambda+|\xi'|^2}} e^{i \,x'\cdot\xi'} \d\xi'\\
    &\quad+\frac{1}{(2\pi)^{n-1}}\int_{\RR^{n-1}} x_n |\xi'|e^{-x_n|\xi'|}\phi(x_n(\sqrt{\lambda+|\xi'|^2}-|\xi'|))\frac{\xi'\,\prescript{t}{}{\xi'}}{|\xi'|^2}\,e^{-y_n\sqrt{\lambda+|\xi'|^2}} e^{i \,x'\cdot\xi'} \d\xi',
\end{align*}
where $\phi(z)=\frac{1-e^{-z}}{z}$, $z\in\CC^{\ast}$. We follow \cite[Section~3,~p.120--122,~eq.~(3.5)]{DeshHieberPruss2001}: for $x'\in\RR^{n-1}$ fixed, considering an orthogonal matrix $\mathrm{Q}$ with real coefficients such that $\mathrm{Q} x' = |x'|\mathfrak{e}_1$, then we perform a change of variables
\begin{align*}
    \RR^{n-1}\setminus\RR \mathrm{Q}\mathfrak{e}_1\,&\longrightarrow\, \RR\times\SS_{n-2}\times(0,\infty)\\
    \xi'\,&\longmapsto\, \Big(\mathrm{Q}^{-1}\xi'\cdot\mathfrak{e}_1,\, \frac{\mathrm{Q}^{-1}\xi'-(\mathrm{Q}^{-1}\xi'\cdot\mathfrak{e}_1)\mathfrak{e}_1}{|\mathrm{Q}^{-1}\xi'-(\mathrm{Q}^{-1}\xi'\cdot\mathfrak{e}_1)\mathfrak{e}_1\,|}, |\mathrm{Q}^{-1}\xi'-(\mathrm{Q}^{-1}\xi'\cdot\mathfrak{e}_1)\mathfrak{e}_1| \Big)=(\eta,\, \Theta,\, r),
\end{align*}
yielding 
\begin{align*}
    &\partial_{x_n}\mathfrak{r}'_\lambda (x',x_n,y_n)\\
    &=\frac{1}{(2\pi)^{n-1}}\int_{0}^\infty r^{n-3}\int_{\SS_{n-2}}\int_{\RR} e^{i |x'|\eta}\frac{\sqrt{r^2+\eta^2}e^{-x_n\sqrt{r^2+\eta^2}}}{\sqrt{\lambda+r^2+\eta^2}}\\ &\quad\quad\quad\quad\quad\quad\quad\quad\quad\quad\quad\quad\quad\quad\quad\quad\quad\quad\times\frac{(\eta \mathfrak{e}_1+r\Theta)\,\prescript{t}{}{(\eta \mathfrak{e}_1+r\Theta)}}{r^2+\eta^2}\,e^{-y_n\sqrt{\lambda+r^2+\eta^2}} \,\d\eta \,\d \Theta\, \d r\\
    &\quad+\frac{1}{(2\pi)^{n-1}}\int_{0}^\infty r^{n-3}\int_{\SS_{n-2}}\int_{\RR}e^{i |x'|\eta}x_n \sqrt{r^2+\eta^2}e^{-x_n\sqrt{r^2+\eta^2}}\phi(x_n(\sqrt{\lambda+r^2+\eta^2}-\sqrt{r^2+\eta^2}))\\&\quad\quad\quad\quad\quad\quad\quad\quad\quad\quad\quad\quad\quad\quad\quad\quad\quad\quad\times\frac{(\eta \mathfrak{e}_1+r\Theta)\,\prescript{t}{}{(\eta \mathfrak{e}_1+r\Theta)}}{r^2+\eta^2}\,e^{-y_n\sqrt{\lambda+r^2+\eta^2}} \,\d\eta \,\d \Theta\, \d r.
\end{align*}
Provided $\varepsilon_0>0$ is small enough, for $\varepsilon\in(0,\varepsilon_0)$, we perform a change in integration contour replacing $\eta$ by $s+i\varepsilon(r+|s|)$, which, by holomorphy of the integrands with respect to $\eta$, does not change the value of the integral, and proceding as \cite[Section~3,~p.122]{DeshHieberPruss2001}, it yields
\begin{align*}
    |\partial_{x_n}&\mathfrak{r}'_\lambda (x',x_n,y_n)|\\
    &\lesssim_{\mu,\varepsilon_0,n} \int_{0}^\infty\int_{0}^\infty r^{n-3} \frac{r+s}{1+r+s}e^{-(\varepsilon |x'|+c(x_n+y_n))(r+s)}e^{-cy_n}\d s\, \d r\\
    & \qquad\qquad+\int_{0}^\infty\int_{0}^\infty r^{n-3} \frac{(r+s)^2}{(1+r+s)(1+r+s+x_n)}e^{-(\varepsilon |x'|+\frac{c}{2}(x_n+y_n))(r+s)}e^{-cy_n}\d s\, \d r\\
     &\lesssim_{\mu,\varepsilon_0,n}  e^{-cy_n}\Big(\int_{0}^\infty r^{n-2}\int_{r}^\infty  \frac{s^2}{1+s}e^{-s(\varepsilon |x'|+c(x_n+y_n))}\frac{\d s}{s}\, \frac{\d r}{r}\\
    & \qquad\qquad+\int_{0}^\infty r^{n-2}\int_{r}^\infty  \frac{s^3}{(1+s)(1+s+x_n)}e^{-s(\varepsilon |x'|+\frac{c}{2}(x_n+y_n))}\frac{\d s}{s}\, \frac{\d r}{r}\Big)\\
    &\lesssim_{\mu,\varepsilon_0,n}  e^{-cy_n}\Big(\int_{0}^\infty \frac{r^{n-1}}{1+r}e^{-r(\varepsilon |x'|+c(x_n+y_n))}{\d r}\\
    & \qquad\qquad+\int_{0}^\infty  \frac{r^{n}}{(1+r)(1+r+x_n)}e^{-r(\varepsilon |x'|+\frac{c}{2}(x_n+y_n))}{\d r}\Big)\\
    &\lesssim_{\mu,\varepsilon_0,n}  e^{-cy_n}\int_{0}^\infty  \frac{r^{n-1}}{(1+r)}e^{-r(\varepsilon |x'|+\frac{c}{2}(x_n+y_n))}{\d r}
\end{align*}
This yields
\begin{align*}
    \int_{0}^{\infty}\int_{\RR^{n-1}}&|\partial_{x_n}\mathfrak{r}'_\lambda (x',x_n,y_n)| \,\d x'\d y_n\\
    &\lesssim_{\mu,\varepsilon_0,n}  \int_{0}^{\infty}\int_{\RR^{n-1}} e^{-cy_n}\int_{0}^\infty  \frac{r^{n-1}}{(1+r)}e^{-r(\varepsilon |x'|+\frac{c}{2} (x_n+y_n))}\,{\d r}\,{\d x'}\,{\d y_n}\\
    &\lesssim_{\mu,\varepsilon_0,n}  \int_{0}^{\infty}\int_{0}^\infty \int_{0}^\infty e^{-cy_n} \frac{r^{n-1}\rho^{n-2}}{(1+r)}e^{-r(\varepsilon \rho+\frac{c}{2} (x_n+y_n))}\,{\d r}\,{\d \rho}\,{\d y_n}\\
    &\lesssim_{\mu,\varepsilon_0,n}  \int_{0}^{\infty}\int_{0}^\infty \int_{0}^\infty  e^{-cy_n}\frac{\tau}{(1+r)}e^{-\varepsilon\tau}e^{-r\frac{c}{2} (x_n+y_n)}\,{\d \tau}\,{\d r}\,{\d y_n}\\
    &\lesssim_{\mu,\varepsilon_0,\varepsilon,n}  \int_{0}^{\infty}\int_{0}^\infty  e^{-cy_n}\frac{1}{(1+r)}e^{-r\frac{c}{2} (x_n+y_n)}\,{\d y_n}\,{\d r}\\
    &\lesssim_{\mu,\varepsilon_0,\varepsilon,n}  \int_{0}^\infty  \frac{1}{(1+r)^2}e^{-r\frac{c}{2} x_n}\,{\d r},
\end{align*}
thus
\begin{align*}
     \sup_{x_n>0}\int_{0}^{\infty}\int_{\RR^{n-1}}|\partial_{x_n}\mathfrak{r}'_\lambda (x',x_n,y_n)|\, \d x'\d y_n \lesssim_{\mu,\varepsilon_0,\varepsilon,n}  \int_{0}^\infty  \frac{1}{(1+r)^2}{\d r} <\infty.
\end{align*}
Therefore, we can apply \cite[Lemma~3.3]{DeshHieberPruss2001}, giving the boundedness of $\partial_{x_n}\mathrm{T}_\lambda'$ as a linear operator from $\L^\infty(\RR^n_+,\CC^{n-1})$ to $\L^\infty(\RR^n_+,\CC^{n-1})$, and consequently, due to the divergence-free condition on $\uu$:
\begin{align*}
    \lVert \nabla \uu\rVert_{\L^\infty(\RR^{n}_+)}&\lesssim_{n} \lVert \nabla' \uu\rVert_{\L^\infty(\RR^{n}_+)} + \lVert \partial_{x_n} \uu_n\rVert_{\L^\infty(\RR^{n}_+)}+ \lVert \partial_{x_n} \uu'\rVert_{\L^\infty(\RR^{n}_+)}\\
    &\lesssim_{n} \lVert \nabla' \uu\rVert_{\L^\infty(\RR^{n}_+)} + \lVert \div' \uu'\rVert_{\L^\infty(\RR^{n}_+)}+ \lVert \partial_{x_n} \uu'\rVert_{\L^\infty(\RR^{n}_+)}\\
    &\lesssim_{n} 2\lVert \nabla' \uu\rVert_{\L^\infty(\RR^{n}_+)} + \lVert \partial_{x_n} (\lambda\I-\Delta_{\mathcal{D}})^{-1}\ff'\rVert_{\L^\infty(\RR^{n}_+)} +\lVert \partial_{x_n} \mathrm{T}_\lambda'\ff'\rVert_{\L^\infty(\RR^{n}_+)}\\
    &\lesssim_{\mu,n} \lVert \ff\rVert_{\L^\infty(\RR^{n}_+)}.
\end{align*}

\textbf{Step 1.2:} For the low frequency control on the pressure, we change the point of view, writing
\begin{align*}
    \nabla'\mathfrak{p}(\cdot,x_n) &= -(\lambda\I-\Delta)\T'_\lambda\ff'(\cdot,x_n)\\
    &= \nabla'\left(\frac{\div'}{(\lambda\I-\Delta')^{\frac{1}{2}}}+\frac{\div'}{(-\Delta')^{\frac{1}{2}}}\right)e^{-x_n(-\Delta')^\frac{1}{2}}\int_{0}^\infty e^{-y_n(\lambda\I-\Delta')^\frac{1}{2}} \ff'(\cdot,y_n)\d y_n,
\end{align*}
so that
\begin{align*}
    \lVert\nabla'\mathfrak{p}\rVert_{\dot{\B}^{-1}_{\infty,\infty}(\RR^n_+)} &=  \bigg\lVert \nabla'\left(\frac{\div'}{(\lambda\I-\Delta')^{\frac{1}{2}}}+\frac{\div'}{(-\Delta')^{\frac{1}{2}}}\right)e^{-x_n(-\Delta')^\frac{1}{2}}\int_{0}^\infty e^{-y_n(\lambda\I-\Delta')^\frac{1}{2}} \ff'(\cdot,y_n)\d y_n\bigg\rVert_{\dot{\B}^{-1}_{\infty,\infty}(\RR^n_+)}\\
    &\lesssim_{n} \bigg\lVert \nabla'\left(\frac{\div'}{(\lambda\I-\Delta')^{\frac{1}{2}}}\right)e^{-x_n(-\Delta')^\frac{1}{2}}\int_{0}^\infty e^{-y_n(\lambda\I-\Delta')^\frac{1}{2}} \ff'(\cdot,y_n)\d y_n\bigg\rVert_{\dot{\B}^{-1}_{\infty,\infty}(\RR^n_+)}\\
    &\qquad\qquad+ \bigg\lVert \nabla'\left(\frac{\div'}{(-\Delta')^{\frac{1}{2}}}\right)e^{-x_n(-\Delta')^\frac{1}{2}}\int_{0}^\infty e^{-y_n(\lambda\I-\Delta')^\frac{1}{2}} \ff'(\cdot,y_n)\d y_n\bigg\rVert_{\dot{\B}^{-1}_{\infty,\infty}(\RR^n_+)}\\
    &\lesssim_{n} \bigg\lVert \left(\frac{\div'}{(\lambda\I-\Delta')^{\frac{1}{2}}}\right)e^{-x_n(-\Delta')^\frac{1}{2}}\int_{0}^\infty e^{-y_n(\lambda\I-\Delta')^\frac{1}{2}} \ff'(\cdot,y_n)\d y_n\bigg\rVert_{{\B}^{0}_{\infty,\infty}(\RR^n_+)}\\
    &\qquad\qquad+ \bigg\lVert \nabla'e^{-x_n(-\Delta')^\frac{1}{2}}\int_{0}^\infty e^{-y_n(\lambda\I-\Delta')^\frac{1}{2}} \ff'(\cdot,y_n)\d y_n\bigg\rVert_{\dot{\B}^{-1}_{\infty,\infty}(\RR^n_+)}\\
    &\lesssim_{n} \bigg\lVert \left(\frac{\div'}{(\lambda\I-\Delta')^{\frac{1}{2}}}\right)e^{-x_n(-\Delta')^\frac{1}{2}}\int_{0}^\infty e^{-y_n(\lambda\I-\Delta')^\frac{1}{2}} \ff'(\cdot,y_n)\d y_n\bigg\rVert_{{\B}^{0}_{\infty,\infty}(\RR^n_+)}\\
    &\qquad\qquad+ \bigg\lVert e^{-x_n(-\Delta')^\frac{1}{2}}\int_{0}^\infty e^{-y_n(\lambda\I-\Delta')^\frac{1}{2}} \ff'(\cdot,y_n)\d y_n\bigg\rVert_{{\B}^{0}_{\infty,\infty}(\RR^n_+)}\\
    &\lesssim_{n} \bigg\lVert \left(\frac{\div'}{(\lambda\I-\Delta')^{\frac{1}{2}}}\right)\int_{0}^\infty e^{-y_n(\lambda\I-\Delta')^\frac{1}{2}} \ff'(\cdot,y_n)\d y_n\bigg\rVert_{{\B}^{0}_{\infty,\infty}(\RR^{n-1})}\\
    &\qquad\qquad+ \bigg\lVert \int_{0}^\infty e^{-y_n(\lambda\I-\Delta')^\frac{1}{2}} \ff'(\cdot,y_n)\d y_n\bigg\rVert_{{\B}^{0}_{\infty,\infty}(\RR^{n-1})}\\
    &\lesssim_{n} \bigg\lVert \int_{0}^\infty e^{-y_n(\lambda\I-\Delta')^\frac{1}{2}} \ff'(\cdot,y_n)\d y_n\bigg\rVert_{{\L}^{\infty}(\RR^{n-1})}\\
    &\lesssim_{n,\mu} \frac{1}{|\lambda|^\frac{1}{2}}\lVert \ff'\rVert_{{\L}^{\infty}(\RR^{n}_+)}
\end{align*}
noting that to obtain the fourth inequality from the third one, we did apply Corollary~\ref{cor:PoissonSmigrpInhomSpaces}.

\medbreak

\noindent Now, for $\partial_{x_n}\mathfrak{p}$, one has the representation
\begin{align*}
    \partial_{x_n}\mathfrak{p}(\cdot,x_n) &= -(\lambda\I-\Delta)\T_{\lambda,n}\ff'(\cdot,x_n)\\
    &= -\div'e^{-x_n(-\Delta')^\frac{1}{2}}\left(\I+\frac{(-\Delta')^{\frac{1}{2}}}{(\lambda\I-\Delta')^{\frac{1}{2}}}\right)\int_{0}^\infty e^{-y_n(\lambda\I-\Delta')^\frac{1}{2}} \ff'(\cdot,y_n)\d y_n,
\end{align*}
so that
\begin{align*}
    \lVert\partial_{x_n}\mathfrak{p}\rVert_{\dot{\B}^{-1}_{\infty,\infty}(\RR^n_+)} &=  \bigg\lVert \div'e^{-x_n(-\Delta')^\frac{1}{2}}\left(\I+\frac{(-\Delta')^{\frac{1}{2}}}{(\lambda\I-\Delta')^{\frac{1}{2}}}\right)\int_{0}^\infty e^{-y_n(\lambda\I-\Delta')^\frac{1}{2}} \ff'(\cdot,y_n)\d y_n\bigg\rVert_{\dot{\B}^{-1}_{\infty,\infty}(\RR^n_+)}\\
    &\lesssim_{n}  \bigg\lVert e^{-x_n(-\Delta')^\frac{1}{2}}\left(\I+\frac{(-\Delta')^{\frac{1}{2}}}{(\lambda\I-\Delta')^{\frac{1}{2}}}\right)\int_{0}^\infty e^{-y_n(\lambda\I-\Delta')^\frac{1}{2}} \ff'(\cdot,y_n)\d y_n\bigg\rVert_{{\B}^{0}_{\infty,\infty}(\RR^n_+)}\\
    &\lesssim_{n}  \bigg\lVert \left(\I+\frac{(-\Delta')^{\frac{1}{2}}}{(\lambda\I-\Delta')^{\frac{1}{2}}}\right)\int_{0}^\infty e^{-y_n(\lambda\I-\Delta')^\frac{1}{2}} \ff'(\cdot,y_n)\d y_n\bigg\rVert_{{\B}^{0}_{\infty,\infty}(\RR^{n-1})}\\
    &\lesssim_{n}  \bigg\lVert \int_{0}^\infty e^{-y_n(\lambda\I-\Delta')^\frac{1}{2}} \ff'(\cdot,y_n)\d y_n\bigg\rVert_{{\L}^{\infty}(\RR^{n-1})}\\
    &\lesssim_{n,\mu} \frac{1}{|\lambda|^\frac{1}{2}}\lVert \ff'\rVert_{{\L}^{\infty}(\RR^{n}_+)}.
\end{align*}
This finishes this substep.

\textbf{Step 2:} We prove full regularity bounds of order 2 and the appropriate regularity bound for the pressure.

Therefore, we want to show first that $\nabla^2 \uu\in{\B}^0_{\infty,\infty}(\RR^n_+)$, $\Delta \uu\in\L^\infty(\RR^n_+)$ with the estimates
\begin{align*}
    \left(\frac{1}{1+|\ln(\min(1,|\lambda|))|}\right)\lVert \nabla^2\uu\rVert_{{\B}^0_{\infty,\infty}(\RR^n_+)} \lesssim_{n,\mu} \lVert \ff\rVert_{\L^\infty(\RR^n_+)},
\end{align*}
and 
\begin{align*}
    \lVert \nabla^2\uu\rVert_{\dot{\B}^0_{\infty,\infty}(\RR^n_+)} \lesssim_{n,\mu} \lVert \ff\rVert_{\L^\infty(\RR^n_+)}.
\end{align*}
Again, we recall that we assumed temporarily $|\lambda|=1$. We will remove this assumption in Step 3.

\medbreak

\noindent By Theorem~\ref{thm:MetaThmDirichletStokesRn+} and Proposition~\ref{prop:MetaThmDirichletStokesRn+HigherReg}, for $-1<s<0$, one has
\begin{align*}
    \lVert(\lambda\I+\AA_\mathcal{D})^{-1}\ff\rVert_{\dot{\B}^{s}_{\infty,1}(\RR^n_+)}+\lVert \nabla^2(\lambda\I+\AA_\mathcal{D})^{-1}\ff\rVert_{\dot{\B}^{s}_{\infty,1}(\RR^n_+)} &\lesssim_{s,n,\mu}\lVert \ff\rVert_{\dot{\B}^{s}_{\infty,1}(\RR^n_+)},\\ &\qquad\qquad\qquad \ff\in\dot{\B}^{s,\sigma}_{\infty,1,\mathfrak{n}}(\RR^n_+),\\
    \lVert(\lambda\I+\AA_\mathcal{D})^{-1}\ff\rVert_{\dot{\B}^{s+1}_{\infty,1}(\RR^n_+)}+\lVert \nabla^2(\lambda\I+\AA_\mathcal{D})^{-1}\ff\rVert_{\dot{\B}^{s+1}_{\infty,1}(\RR^n_+)} &\lesssim_{s,n,\mu}\lVert \ff\rVert_{\dot{\B}^{s+1}_{\infty,1}(\RR^n_+)},\\ &\qquad\qquad\qquad \ff\in\dot{\B}^{s,\sigma}_{\infty,1,\mathfrak{n}}\cap\dot{\B}^{s+1,\sigma}_{\infty,1,\mathcal{D}}(\RR^n_+).
\end{align*}
By Remark~\ref{rem:FlatExtOpDiv}, for $m\geqslant 3$ fixed, we consider a lifted and bounded  (extended by density) extension of $(\lambda\I+\AA_\mathcal{D})^{-1}$:
\begin{align*}
    \tilde{A}_\lambda^{\sigma}\,:\, \ell_{s+\mathfrak{k}}^{1}(\ZZ,\L^\infty_\sigma(\RR^n)) &\longrightarrow \ell_{s+\mathfrak{k}}^{1}(\ZZ,\L^\infty_\sigma(\RR^n))\cap\ell_{s+2+\mathfrak{k}}^{1}(\ZZ,\L^\infty_\sigma(\RR^n))\\
    (\ff_j)_{j\in\ZZ}\qquad& \longmapsto \left[\dot{\Delta}_{k}\E^{m}_{\sigma}(\lambda\I+\AA_\mathcal{D})^{-1}\R_{\RR^n_+}\P_{0,\sigma}^m\sum_{j\in\ZZ}\dot{\Delta}_j[\ff_{j-1}+\ff_{j}+\ff_{j+1}]\right]_{k\in\ZZ},
\end{align*}
provided $\mathfrak{k}=0,1$, and $\mathcal{E}_0$ is the extension from to the whole space by $0$. Note that for $\ff_j=\dot{\Delta}_j\mathcal{E}_0\ff$, $\ff\in\dot{\B}^{s,\sigma}_{\infty,1,\mathfrak{n}}(\RR^n_+)$, one has
\begin{align*}
    \sum_{k\in\ZZ} &[\dot{\Delta}_{k-1}+\dot{\Delta}_{k}+\dot{\Delta}_{k+1}]\tilde{A}_\lambda^{\sigma,k}(\dot{\Delta}_j\mathcal{E}_0\ff)_{j\in\ZZ} \\
    &= \sum_{k\in\ZZ} \dot{\Delta}_{k}\E^{m}_{\sigma}(\lambda\I+\AA_\mathcal{D})^{-1}\R_{\RR^n_+}\P_{0,\sigma}^m\sum_{j\in\ZZ}\dot{\Delta}_j[[\dot{\Delta}_{j-1}+\dot{\Delta}_{j}+\dot{\Delta}_{j+1}]\mathcal{E}_0\ff]\\
    &= \sum_{k\in\ZZ} \dot{\Delta}_{k}\E^{m}_{\sigma}(\lambda\I+\AA_\mathcal{D})^{-1}\R_{\RR^n_+}\P_{0,\sigma}^m\sum_{j\in\ZZ}\dot{\Delta}_j\mathcal{E}_0\ff\\
    &= \sum_{k\in\ZZ} \dot{\Delta}_{k}\E^{m}_{\sigma}(\lambda\I+\AA_\mathcal{D})^{-1}\R_{\RR^n_+}\P_{0,\sigma}^m\mathcal{E}_0\ff\\
    &= \sum_{k\in\ZZ} \dot{\Delta}_{k}\E^{m}_{\sigma}(\lambda\I+\AA_\mathcal{D})^{-1}\ff\\
    &= \E^{m}_{\sigma}(\lambda\I+\AA_\mathcal{D})^{-1}\ff.
\end{align*}
By real interpolation between the cases $\mathfrak{k}=0$ and $\mathfrak{k}=1$, one obtains the bounded linear map
\begin{align*}
    \tilde{A}_\lambda^{\sigma}\,:\, \ell^{\infty}(\ZZ,\L^\infty_\sigma(\RR^n)) &\longrightarrow \ell_{2}^{\infty}(\ZZ,\L^\infty_\sigma(\RR^n)),
\end{align*}
with bound
\begin{align*}
    \lVert [\tilde{A}_\lambda^{\sigma}(\ff_j)_{j\in\ZZ}]\rVert_{\ell_{2}^{\infty}(\ZZ,\L^\infty(\RR^n))} \lesssim_{n,\mu} \lVert(\ff_j)_{j\in\ZZ}\rVert_{\ell^{\infty}(\ZZ,\L^\infty(\RR^n))}.
\end{align*}
Therefore, for $\ff\in\L^\infty_{\mathfrak{n},\sigma}(\RR^n_+)$, one sets $\ff_j:=\dot{\Delta}_j\mathcal{E}_0\ff\in\L^\infty_{\sigma}(\RR^n)$ and the estimates
\begin{align*}
    \lVert [\tilde{A}_\lambda^{\sigma}(\dot{\Delta}_j\mathcal{E}_0\ff)_{j\in\ZZ}]\rVert_{\ell_{2}^{\infty}(\ZZ,\L^\infty(\RR^n))} &\lesssim_{n,\mu} \lVert(\dot{\Delta}_j\mathcal{E}_0\ff)_{j\in\ZZ}\rVert_{\ell^{\infty}(\ZZ,\L^\infty(\RR^n))}\\ &\lesssim_{n,\mu} \lVert\mathcal{E}_0\ff\rVert_{\L^\infty(\RR^n)}\\ &\lesssim_{n,\mu} \lVert\ff\rVert_{\L^\infty(\RR^n_+)} .
\end{align*}
We set $\uu:=(\lambda\I-\Delta_\mathcal{D})^{-1}\ff+\T_{\lambda}\ff$, and $\nabla\mathfrak{p}=-(\lambda\I-\Delta) \T_{\lambda}\ff\in \dot{\B}^{-1}_{\infty,\infty}(\RR^n_+)$. For for all $N\in\NN$, thanks to Theorem~\ref{thm:MetaThmDirichletStokesRn+},  Remark~\ref{rem:FlatExtOpDiv}, and Step 1, one can set
\begin{align*}
    \ff_N&:= \R_{\RR^n_+}\left[\P_{0,\sigma}^m\sum_{j\geqslant - N} \dot{\Delta}_j\mathcal{E}_0\ff\right]\in\L^\infty_{\mathfrak{n},\sigma}(\RR^n_+)\cap\dot{\B}^{s,\sigma}_{\infty,1,\mathfrak{n}}(\RR^n_+)\subset\L^\infty_{h,\mathfrak{n},\sigma}(\RR^n_+),\\
    \uu_N&:=(\lambda\I+\AA_\mathcal{D})^{-1}\ff_N\in\dot{\B}^{s,\sigma}_{\infty,1,\mathfrak{n}}(\RR^n_+)\cap\dot{\B}^{s+2,\sigma}_{\infty,1,\mathcal{D}}(\RR^n_+)\subset\tilde{\B}^{s+2,\sigma}_{\infty,1,\mathcal{D}}(\RR^n_+)\subset\W^{1,\infty}_{\mathcal{D},\sigma}(\RR^n_+)\cap\L^\infty_{h,\mathfrak{n},\sigma}(\RR^n_+),\\
    \vv_N&:=(\lambda\I-\Delta_\mathcal{D})^{-1}\ff_N+\T_{\lambda}\ff_N\in\W^{1,\infty}_{\mathcal{D},\sigma}(\RR^n_+),
\end{align*}
where $s$ could be any $-1<s<0$. Thanks to Lemma~\ref{lem:FlatProjOpDivModConst}, it can be checked that $(\ff_N)_{N\in\NN}$ converges weakly-$\ast$ to $\ff$ in $\L^\infty(\RR^n_+)=(\L^1(\RR^n_+))'$, and one has the uniform bound
\begin{align*}
    \lVert \ff_N\rVert_{\L^\infty(\RR^n_+)}\lesssim_{m,n} \lVert \ff\rVert_{\L^\infty(\RR^n_+)}.
\end{align*}

We also consider the respective pressure terms $\nabla\mathfrak{p}_N\in\dot{\B}^{s}_{\infty,1}(\RR^n_+,\CC^n)$, for any $-1<s<0$ given by Theorem~\ref{thm:MetaThmDirichletStokesRn+}, and $\nabla\mathfrak{q}_N = -(\lambda\I-\Delta)\T_\lambda \ff_N\in\dot{\B}^{-1}_{\infty,\infty}(\RR^{n}_+,\CC^n)$, so that both couples $(\uu_N,\mathfrak{p}_N)$ and $(\vv_N,\mathfrak{q}_N)$ are solving \eqref{eq:DirStokesSystemtn+Infty} with forcing term $\ff_N$. In particular, both are weak solutions to \eqref{eq:weakStokesResRn+1}-\eqref{eq:weakStokesResRn+2} with forcing term $\ff_N$.
So that, for all $N\in\NN$, $\E_{\sigma}^m\uu_N\in{\B}^{2}_{\infty,\infty}(\RR^n,\CC^n)$, yielding $\uu_N\in{\B}^{2}_{\infty,\infty}(\RR^n_+,\CC^n)$ and this holds with the uniform bound
\begin{align*}
    \lVert \nabla^2\E^{m}_{\sigma} \uu_N\rVert_{\dot{\B}^{0}_{\infty,\infty}(\RR^n)}\sim_{n} \lVert [\tilde{A}_\lambda^{\sigma}(\dot{\Delta}_j\mathcal{E}_0\ff_N)_{j\in\ZZ}]\rVert_{\ell_{2}^{\infty}(\ZZ,\L^\infty(\RR^n))} &\lesssim_{n,\mu} \lVert\ff_N\rVert_{\L^\infty(\RR^n_+)}\\ &\lesssim_{n,\mu} \lVert\ff\rVert_{\L^\infty(\RR^n_+)}.
\end{align*}
Consequently, in the distributional sense (even in $\B^{-1}_{\infty,\infty}(\RR^n_+,\CC^n)$), it holds
\begin{align*}
    \lambda (\uu_N-\vv_N) - \Delta (\uu_N-\vv_N) + \nabla(\mathfrak{p}_N-\mathfrak{q}_N) = 0.
\end{align*}
Since $\uu_N-\vv_N\in\W^{1,\infty}_{\mathcal{D}}(\RR^n_+)$, one obtains $\Delta (\uu_N-\vv_N)\in\dot{\W}^{-1,\infty}(\RR^n_+,\CC^n)$, and one deduces
\begin{align*}
    \nabla(\mathfrak{p}_N-\mathfrak{q}_N) \in \W^{1,\infty}_{\mathcal{D}}\cap\L^\infty_h(\RR^n_+,\CC^n) + \dot{\W}^{-1,\infty}(\RR^n_+,\CC^n).
\end{align*}
Therefore, ${\mathfrak{p}_N-\mathfrak{q}_N}$ admits a representative in $\L^1_{\mathrm{loc}}(\RR^n_+)$, and so does $\mathfrak{q}_N$.

Hence, we can apply Theorem~\ref{thm:MaekawaMiruaPrangeUniqueness}, so that there exists $\mathbf{a}'_N\in\CC^{n-1}$ such that $(\uu_N,\mathfrak{p}_N)=(\vv_N,\mathfrak{q}_N) + (\ww_{\mathbf{a}'_N},\mathfrak{q}_{\mathbf{a}'_N})$, which implies $\vv_N\in\B^{2}_{\infty,\infty}(\RR^n_+,\CC^n)$, and by linearity $\vv_N,\ff_N\in\L^\infty(\RR^n_+,\CC^n)\subset \B^{0}_{\infty,\infty}(\RR^n_+)$ and $\Delta\vv_N\in\B^{0}_{\infty,\infty}(\RR^n_+,\CC^n)$ are implying 
\begin{align*}
    \nabla \mathfrak{q}_N= -\lambda \vv_N +\Delta \vv_N +\ff_N \in\B^{0}_{\infty,\infty}(\RR^n_+).
\end{align*}
Thus, $\nabla'\mathfrak{q}_N\in \dot{\B}^{-1}_{\infty,\infty}\cap\B^{0}_{\infty,\infty}(\RR^n_+,\CC^{n-1})\subset \tilde{\B}^{0}_{\infty,\infty}(\RR^n_+,\CC^{n-1})$, and similarly $\nabla'\mathfrak{p}_N\in \dot{\B}^{s}_{\infty,1}\cap\B^{0}_{\infty,\infty}(\RR^n_+,\CC^{n-1})\subset \tilde{\B}^{0}_{\infty,\infty}(\RR^n_+,\CC^{n-1})$  so that by linearity
\begin{align*}
    \mathbf{a}'_N =\nabla'\mathfrak{p}_N-\nabla'\mathfrak{q}_N\in\tilde{\B}^{0}_{\infty,\infty}(\RR^n_+,\CC^{n-1})
\end{align*}
which implies $\mathbf{a}'_N=\mathbf{0}'$ and consequently for all $N\in\NN$
\begin{align*}
    (\uu_{N},\nabla\mathfrak{p}_N)=(\vv_{N},\nabla\mathfrak{q}_N).
\end{align*}

Therefore, we deduce a bound uniform with respect to $N\in\NN$ given by
\begin{align}
    \lVert\uu_N\rVert_{{\B}^{2}_{\infty,\infty}(\RR^n_+)} &\leqslant \lVert \E^m_\sigma\uu_N\rVert_{{\B}^{2}_{\infty,\infty}(\RR^n)}\nonumber\\
    &\lesssim_{n} \lVert \E^{m}_{\sigma} \uu_N\rVert_{\L^\infty(\RR^n)}+\lVert \nabla^2\E^{m}_{\sigma} \uu_N\rVert_{\dot{\B}^{0}_{\infty,\infty}(\RR^n)}\nonumber\\
    &\lesssim_{\mu,n}\lVert  \uu_N\rVert_{\L^\infty(\RR^n_+)} + \lVert\ff_N\rVert_{\L^\infty(\RR^n_+)} \label{eq:LinftyStokesProofUniformBound}\\ &\lesssim_{\mu,n}\lVert  \vv_N\rVert_{\L^\infty(\RR^n_+)} + \lVert\ff_N\rVert_{\L^\infty(\RR^n_+)} \nonumber\\ &\lesssim_{\mu,n}\lVert\ff_N\rVert_{\L^\infty(\RR^n_+)}\nonumber\\
    &\lesssim_{\mu,n}\lVert\ff\rVert_{\L^\infty(\RR^n_+)}.\nonumber
\end{align}
One also obtains by linearity
\begin{align*}
    \lVert\nabla \mathfrak{p}_N\rVert_{{\B}^{0}_{\infty,\infty}(\RR^n_+)} \lesssim_{\mu,n}\lVert\ff\rVert_{\L^\infty(\RR^n_+)}.
\end{align*}
Consequently, there exists a subsequence $(\uu_{N_k}, \mathfrak{p}_{N_k})_{k\in\NN}$ that converges weakly-$\ast$ towards some $(\tilde{\uu},\tilde{\mathfrak{p}})$ in $\B^{2}_{\infty,\infty}(\RR^n_+,\CC^n)\times\hat{\B}^{1}_{\infty,\infty}(\RR^n_+)=(\B^{-2}_{1,1,0}(\RR^n_+)\times\mathring{\B}^{-1}_{1,1,0}(\RR^n_+))'$. For all $\boldsymbol{\varphi}\in\Ccinfty(\RR^n_+,\CC^n)$, it can be checked that one has
\begin{align*}
    &\lambda \int_{\RR^n_+}  \uu_{N_k}(x) \cdot \boldsymbol{\varphi}(x)\, \d x + \int_{\RR^n_+} \nabla \uu_{N_k}(x) : \nabla \boldsymbol{\varphi}(x)\, \d x - \int_{\RR^n_+} \mathfrak{p}_{N_k}(x) \div \boldsymbol{\varphi}(x)\, \d x  =  \int_{\RR^n_+} \ff_{N_k}(x) \cdot \boldsymbol{\varphi}(x)\, \d x\\
    &\xrightarrow[k\longrightarrow\infty]{} \lambda \int_{\RR^n_+}  \tilde{\uu}(x) \cdot \boldsymbol{\varphi}(x)\, \d x + \int_{\RR^n_+} \nabla \tilde{\uu}(x) : \nabla \boldsymbol{\varphi}(x)\, \d x - \int_{\RR^n_+} \tilde{\mathfrak{p}}(x) \div \boldsymbol{\varphi}(x)\, \d x =  \int_{\RR^n_+} \ff(x) \cdot \boldsymbol{\varphi}(x)\, \d x.
\end{align*}
Consequently, by Theorem~\ref{thm:MaekawaMiruaPrangeUniqueness},  this implies that  $(\uu-\tilde{\uu},\mathfrak{p}-\tilde{\mathfrak{p}}) = (\ww_\mathbf{a'},\mathfrak{q}_\mathbf{a'})$ for some $\mathbf{a'}\in\CC^{n-1}$.
Since $\hat{\B}^{1}_{\infty,\infty}(\RR^n_+)$ is weak-$\ast$ complete,
one obtains in particular
\begin{align*}
   \mathbf{a}'= \nabla'\mathfrak{p}-\nabla'\tilde{\mathfrak{p}}\in\dot{\B}^{-1}_{\infty,\infty}(\RR^n_+,\CC^{n-1}),
\end{align*}
this implying $\mathbf{a}'=\mathbf{0}'$, we deduce ultimately
\begin{align*}
    (\uu,\nabla\mathfrak{p})= (\tilde{\uu},\nabla \tilde{\mathfrak{p}}).
\end{align*}
In particular, this is the solution given by the Desch, Hieber and Pr\"uss representation formula. We conclude with $\uu\in\B^{2}_{\infty,\infty}(\RR^{n}_+,\CC^n)$ and $\nabla\mathfrak{p}\in\tilde{\B}^{0}_{\infty,\infty}(\RR^{n}_+,\CC^{n})$, also having the estimates
\begin{align*}
    \lVert \uu\rVert_{{\B}^{2}_{\infty,\infty}(\RR^n_+)} \leqslant \liminf_{k\longrightarrow\infty} \lVert \uu_{N_k}\rVert_{{\B}^{2}_{\infty,\infty}(\RR^n_+)}\lesssim_{n,\mu} \lVert\ff\rVert_{\L^\infty(\RR^n_+)},\\
    \lVert \nabla\mathfrak{p}\rVert_{\dot{\B}^{0}_{\infty,\infty}(\RR^n_+)} \leqslant \lVert \nabla\mathfrak{p}\rVert_{{\B}^{0}_{\infty,\infty}(\RR^n_+)} \leqslant \liminf_{k\longrightarrow\infty} \lVert \nabla\mathfrak{p}_{N_k}\rVert_{{\B}^{0}_{\infty,\infty}(\RR^n_+)}\lesssim_{n,\mu} \lVert\ff\rVert_{\L^\infty(\RR^n_+)}.
\end{align*}
Note that $\uu\in {\B}^{2}_{\infty,\infty}(\RR^n_+,\CC^n)\subset\L^\infty(\RR^n_+,\CC^n)$ implies that $\nabla^2 \E^{m}_{\sigma} \uu \in \dot{\B}^{-2}_{\infty,\infty}(\RR^n,\CC^{n^3})\subset \S'_h(\RR^n,\CC^{n^3})$, therefore $\nabla^2 \E^{m}_{\sigma} \uu\in\dot{\B}^{0}_{\infty,\infty}(\RR^n,\CC^{n^3})$, and by definition of function spaces by restriction
\begin{align}\label{eq:ProofLinftyLogHomEst}
    \lVert \uu\rVert_{{\L}^{\infty}(\RR^n_+)}+\lVert \nabla\uu\rVert_{{\L}^{\infty}(\RR^n_+)}+\lVert \nabla^2\uu\rVert_{\dot{\B}^{0}_{\infty,\infty}(\RR^n_+)}  \lesssim_{n}  \lVert \uu\rVert_{{\B}^{2}_{\infty,\infty}(\RR^n_+)}  \lesssim_{n,\mu} \lVert\ff\rVert_{\L^\infty(\RR^n_+)}.
\end{align}
Otherwise, we consider the estimate
\begin{align}\label{eq:ProofLinftyLogHomEst2}
    \lVert (\nabla^2 \uu,\nabla\mathfrak{p})\rVert_{{\B}^{0}_{\infty,\infty}(\RR^n_+)} \lesssim_{\mu,n} \lVert \uu\rVert_{{\B}^{2}_{\infty,\infty}(\RR^n_+)} + \lVert\ff\rVert_{\L^\infty(\RR^n_+)}\lesssim_{\mu,n} \lVert\ff\rVert_{\L^\infty(\RR^n_+)}.
\end{align}
\textbf{Step 3:} Now, we get rid of the assumption $|\lambda|=1$ by a dilation argument. Starting from \eqref{eq:ProofLinftyLogHomEst} the dilation procedure is the standard one, we omit it. However, the dilation procedure applied to \eqref{eq:ProofLinftyLogHomEst2} is different. To end the proof, we are distinguishing  two cases $|\lambda|<1$ and $|\lambda|\geqslant 1$. The main issue here concerns the fact that we deal with an inhomogeneous Besov space $\B^{0}_{\infty,\infty}$ which is not dilation-invariant.
We assume $\lambda\in\Sigma_\mu$, note that for $\tilde{\uu}_{\lambda} :={|\lambda|}\uu(\cdot/|\lambda|^\frac{1}{2}) =: {|\lambda|}\uu_{|\lambda|^\frac{1}{2}}$ and $\tilde{\ff}_{\lambda} :=\ff_{|\lambda|^\frac{1}{2}}$, one has 
\begin{align*}
    \nabla^2 \left({|\lambda|}\uu_{|\lambda|^\frac{1}{2}}\right) = (\nabla^2\uu)_{|\lambda|^\frac{1}{2}}
\end{align*}
so that if $|\lambda|<1$, one obtains $1/|\lambda|^\frac{1}{2}>1$, and  by \cite[Theorem~3.11]{Vybiral2008}:
\begin{align*}
    \lVert\nabla^2\uu \rVert_{\B^{0}_{\infty,\infty}(\RR^n_+)} &= \lVert ((\nabla^2\uu)_{|\lambda|^\frac{1}{2}})_{|\lambda|^{-\frac{1}{2}}}\rVert_{\B^{0}_{\infty,\infty}(\RR^n_+)}\\
    & \lesssim_{n} \big({1+|\ln(|\lambda|)|}\big)\,\lVert (\nabla^2\uu)_{|\lambda|^\frac{1}{2}}\rVert_{\B^{0}_{\infty,\infty}(\RR^n_+)}= \big({1+|\ln(|\lambda|)|}\big)\,\lVert \nabla^2\tilde{\uu}_{\lambda}\rVert_{\B^{0}_{\infty,\infty}(\RR^n_+)}\\
    &\lesssim_{n,\mu} \big({1+|\ln(|\lambda|)|}\big)\,\lVert\tilde{\ff}_{\lambda}\rVert_{\L^\infty(\RR^n_+)}= \big({1+|\ln(|\lambda|)|}\big)\,\lVert\ff\rVert_{\L^\infty(\RR^n_+)}.
\end{align*}
Now, if $|\lambda|\geqslant 1$, one can reperform the whole proof replacing \eqref{eq:LinftyStokesProofUniformBound} by the estimate
\begin{align*}
    \lVert\uu\rVert_{{\B}^{2}_{\infty,\infty}(\RR^n_+)} 
    &\lesssim_{\mu,n}\left(\frac{1}{|\lambda|}+1\right)\lVert\ff\rVert_{\L^\infty(\RR^n_+)}\\
    &\lesssim_{\mu,n} 2\lVert\ff\rVert_{\L^\infty(\RR^n_+)}.
\end{align*}
Consequently,
\begin{align*}
     \lVert \nabla^2 \uu\rVert_{{\B}^{0}_{\infty,\infty}(\RR^n_+)}\leqslant \lVert \uu\rVert_{{\B}^{2}_{\infty,\infty}(\RR^n_+)} \lesssim_{n,\mu} \lVert\ff\rVert_{\L^\infty(\RR^n_+)}.
\end{align*}
The facts that $\uu,\ff\in\L^\infty(\RR^n_+)\subset \B^{0}_{\infty,\infty}(\RR^n_+)$ and $\Delta\uu\in\B^{0}_{\infty,\infty}(\RR^n_+)$ are implying 
\begin{align*}
    \nabla \mathfrak{p}= -\lambda \uu +\Delta \uu +\ff \in\B^{0}_{\infty,\infty}(\RR^n_+).
\end{align*}
Summing up everything, we get
\begin{align*}
    \left(\frac{1}{1+|\ln(\min(1,|\lambda|))|}\right)\lVert (\nabla^2 \uu,  \nabla\mathfrak{p})\rVert_{{\B}^{0}_{\infty,\infty}(\RR^n_+)} \lesssim_{n,\mu} \lVert\ff\rVert_{\L^\infty(\RR^n_+)}.
\end{align*}
\textbf{Step 4:} We just check that $\uu\in\tilde{\B}^2_{\infty,\infty}(\RR^n_+,\CC^n)$ whenever $\ff\in\L^\infty_{h,\mathfrak{n},\sigma}(\RR^n_+)$. In order to see it, we write
\begin{align*}
    \uu = \frac{1}{\lambda} \ff + \frac{1}{\lambda}\left(\Delta \uu-\nabla\mathfrak{p}\right)\in \L^\infty_{h,\mathfrak{n},\sigma}(\RR^n_+) + \tilde{\B}^0_{\infty,\infty}(\RR^n_+,\CC^n)\subset \tilde{\B}^0_{\infty,\infty}(\RR^n_+,\CC^n)
\end{align*}
Thus, 
\begin{align*}
    \uu \in \tilde{\B}^0_{\infty,\infty}\cap {\B}^2_{\infty,\infty}(\RR^n_+,\CC^n)\subset\tilde{\B}^2_{\infty,\infty}(\RR^n_+,\CC^n).
\end{align*}
This finishes the proof.
\end{proof}

We can recover additionnal regularity properties  in the case of inhomogeneous end-point Besov spaces. Even if it is interesting in itself, this will be more relevant and useful later on for endpoint interpolation theory.

\begin{proposition}\label{prop:DampedStokesDirRn+Lip} Let $\mu\in[0,\pi)$. For all $\lambda\in\Sigma_{\mu}$, for all $\ff\in\W^{1,\infty}_{\mathcal{D},\sigma}(\RR^n_+)$ the system
\begin{equation*}\tag{DDS${}_\lambda$}\label{eq:DampedDirStokesSystemtn+InftyLip}
    \left\{ \begin{array}{rllr}
         \lambda \uu  - \Delta \uu +\nabla \mathfrak{p} &= \ff \text{, }&&\text{ in } \RR^n_+\text{,}\\
        \div \uu &= 0\text{, } &&\text{ in } \RR^n_+\text{,}\\
        \uu_{|_{\partial\RR^n_+}} &=0\text{, } &&\text{ on } \partial\RR^n_+\text{.}
    \end{array}
    \right.
\end{equation*}
admits a unique solution $(\uu,\nabla\mathfrak{p})\in{\B}^{3,\sigma}_{\infty,\infty,\mathcal{D}}(\RR^n_+)\times\tilde{\B}^{1}_{\infty,\infty}(\RR^n_+,\CC^n)$ satisfying the estimates
\begin{align*}
    |\lambda|\lVert \nabla \uu\rVert_{\L^\infty(\RR^n_+)} +|\lambda|^\frac{1}{2}\lVert \nabla^2 \uu\rVert_{\L^\infty(\RR^n_+)}+ \lVert (\nabla^3 \uu,\nabla^{2}\mathfrak{p})\rVert_{\dot{\B}^{0}_{\infty,\infty}(\RR^n_+)} \lesssim_{n,\mu} \lVert \nabla \ff\rVert_{\L^\infty(\RR^n_+)}
\end{align*}
and for any $c>0$
\begin{align*}
    \left(\frac{1}{1+|\ln(\min(c,|\lambda|))|}\right)\lVert (\nabla^3\uu,\nabla^2 \mathfrak{p})\rVert_{{\B}^0_{\infty,\infty}(\RR^n_+)} \lesssim_{n,\mu,c} \lVert \nabla\ff\rVert_{\L^\infty(\RR^n_+)}.
\end{align*}
\end{proposition}

\begin{proposition}\label{prop:DampedStokesDirRn+BesovInfty} Let $q\in[1,\infty]$, $s\in(-1,2)$, $s\neq0$. Let $\mu\in[0,\pi)$. For all $\ff\in{\B}^{s,\sigma}_{\infty,q,\mathcal{D}}(\RR^n_+)$, for all $\lambda\in\Sigma_{\mu}$, there exists a unique solution $(\uu,\nabla\mathfrak{p})\in{\B}^{s+2,\sigma}_{\infty,q,\mathcal{D}}(\RR^n_+)\times\tilde{\B}^{s}_{\infty,q}(\RR^n_+,\CC^n)$ to \eqref{eq:DampedDirStokesSystemtn+InftyLip} satisfying the estimates
\begin{align*}
    |\lambda|\lVert \uu\rVert_{{\B}^{s}_{\infty,q}(\RR^n_+)} +|\lambda|^\frac{1}{2}\lVert \nabla\uu\rVert_{{\B}^{s}_{\infty,q}(\RR^n_+)}  \lesssim_{s,n,\mu} \lVert \ff\rVert_{{\B}^{s}_{\infty,q}(\RR^n_+)},
\end{align*}
and for any $c>0$
\begin{align*}
    \left(\frac{1}{1+|\ln(\min(c,|\lambda|))|}\right)\lVert (\nabla^2 \uu,\nabla\mathfrak{p})\rVert_{{\B}^{s}_{\infty,q}(\RR^n_+)} \lesssim_{s,n,\mu,c} \lVert \ff\rVert_{{\B}^{s}_{\infty,q}(\RR^n_+)}.
\end{align*}
Moreover, when $s>0$, one also has
\begin{align*}
    \lVert (\nabla^2 \uu,\nabla\mathfrak{p})\rVert_{\dot{\B}^{0}_{\infty,\infty}(\RR^n_+)}+ \lVert( \nabla^2 \uu,\nabla\mathfrak{p}) \rVert_{\dot{\B}^{s}_{\infty,q}(\RR^n_+)}\lesssim_{s,n,\mu} \lVert \ff\rVert_{{\B}^{s}_{\infty,q}(\RR^n_+)}.
\end{align*}

If additionally $\ff\in\tilde{\B}^{s,\sigma}_{\infty,q,\mathcal{D}}(\RR^n_+)$ (\textit{resp.} $\ff\in\tilde{\BesSmo}^{s,\sigma}_{\infty,\infty,\mathcal{D}}, {\B}^{s,0,\sigma}_{\infty,q,\mathcal{D}}({\RR^{n}_+})$) one has $\uu\in\tilde{\B}^{s+2,\sigma}_{\infty,q,\mathcal{D}}(\RR^n_+)$  (\textit{resp.} $\uu\in\tilde{\BesSmo}^{s+2,\sigma}_{\infty,\infty}(\RR^n_+),\B^{s+2,0,\sigma}_{\infty,q}(\RR^n_+)$).
\end{proposition}

\begin{proof}[of Propositions~\ref{prop:DampedStokesDirRn+Lip}~\&~\ref{prop:DampedStokesDirRn+BesovInfty}]\textbf{Step 1:} We prove the Lipschitz estimate of Proposition~\ref{prop:DampedStokesDirRn+Lip}.  First, we assume $\ff\in\W^{1,\infty}_{\mathcal{D},\sigma}(\RR^n_+)$, by Theorem~\ref{thm:StokesDirRn+Linfty}, one has existence and uniqueness of a solution $\uu\in\B^{2}_{\infty,\infty}(\RR^n_+,\CC^n)$, with estimates
\begin{align*}
    |\lambda|\lVert \uu\rVert_{{\L}^{\infty}(\RR^n_+)} +|\lambda|^\frac
    {1}{2}\lVert \nabla \uu\rVert_{{\L}^{\infty}(\RR^n_+)} + \lVert \nabla^2 \uu\rVert_{\dot{\B}^{0}_{\infty,\infty}(\RR^n_+)} \lesssim_{s,n,\mu} \lVert \ff\rVert_{{\L}^{\infty}(\RR^n_+)}.
\end{align*}
and
\begin{align*}
    \left(\frac{1}{1+|\ln(\min(c,|\lambda|))|}\right)\lVert (\nabla^2 \uu,\nabla\mathfrak{p})\rVert_{{\B}^{0}_{\infty,\infty}(\RR^n_+)} \lesssim_{s,n,\mu,c} \lVert \ff\rVert_{{\L}^{\infty}(\RR^n_+)}.
\end{align*}

Note that in particular, one has $\uu\in\C^{1,\beta}_{ub,\mathcal{D},\sigma}(\overline{\RR^n_+})$, for all $0<\beta<1$. For $j\in\llb1,n-1\rrb$, writing for $h>0$
\begin{align*}
    \partial_{j}^{h} w := \frac{w(\cdot+h\mathfrak{e}_j)-w}{h},
\end{align*}
one obtains
\begin{equation*}
    \left\{ \begin{array}{rllr}
         \lambda  \partial_{j}^{h}\uu  - \Delta  \partial_{j}^{h}\uu +\nabla  \partial_{j}^{h}\mathfrak{p} &=  \partial_{j}^{h}\ff \text{, }&&\text{ in } \RR^n_+\text{,}\\
        \div  (\partial_{j}^{h}\uu) &= 0\text{, } &&\text{ in } \RR^n_+\text{,}\\
        \partial_{j}^{h} \uu_{|_{\partial\RR^n_+}} &=0\text{, } &&\text{ on } \partial\RR^n_+\text{.}
    \end{array}
    \right.
\end{equation*}
Since $(\partial_{j}^{h}\uu)_{h>0}$ converges towards $\partial_{x_j}\uu$ in $\C^{0,\beta}_{ub}(\overline{\RR^n_+})$, then the convergence holds in the distributional sense, so that for $\D'^{h} := (\partial_{j}^{h})_{1\leqslant j \leqslant n-1}$:
\begin{align*}
    |\lambda|\lVert \nabla' \uu\rVert_{{\L}^{\infty}(\RR^n_+)} &+ |\lambda|^\frac{1}{2}\lVert \nabla'\nabla \uu\rVert_{{\L}^{\infty}(\RR^n_+)} + \lVert \nabla' \nabla^2 \uu\rVert_{\dot{\B}^{0}_{\infty,\infty}(\RR^n_+)} \\
    &\leqslant \liminf_{h\longrightarrow 0_+} \, |\lambda|\lVert \D'^{h}\uu\rVert_{{\L}^{\infty}(\RR^n_+)} + |\lambda|^\frac{1}{2}\lVert \nabla\D'^{h} \uu\rVert_{{\L}^{\infty}(\RR^n_+)} + \lVert \nabla^2 \D'^{h}\uu\rVert_{\dot{\B}^{0}_{\infty,\infty}(\RR^n_+)}\\ &\lesssim_{s,n,\mu} \liminf_{h\longrightarrow 0_+} \,\lVert \D'^{h}\ff\rVert_{{\L}^{\infty}(\RR^n_+)} \\
    &\lesssim_{s,n,\mu}\, \lVert \nabla' \ff\rVert_{{\L}^{\infty}(\RR^n_+)}.
\end{align*}
Similarly, we deduce
\begin{align*}
    \left(\frac{1}{1+|\ln(\min(c,|\lambda|))|}\right)\lVert \nabla' \nabla^2 \uu\rVert_{{\B}^{0}_{\infty,\infty}(\RR^n_+)}
    &\leqslant \liminf_{h\longrightarrow 0_+} \, \left(\frac{1}{1+|\ln(\min(c,|\lambda|))|}\right)\lVert \nabla^2 \D'^{h}\uu\rVert_{{\B}^{0}_{\infty,\infty}(\RR^n_+)}\\ &\lesssim_{s,n,\mu} \liminf_{h\longrightarrow 0_+} \,\lVert \D'^{h}\ff\rVert_{{\L}^{\infty}(\RR^n_+)} \\
    &\lesssim_{s,n,\mu}\, \lVert \nabla' \ff\rVert_{{\L}^{\infty}(\RR^n_+)}.
\end{align*}
By linearity, one also obtains
\begin{align*}
    \lVert \nabla' \nabla \mathfrak{p}\rVert_{{\B}^{0}_{\infty,\infty}(\RR^n_+)} \lesssim_{s,n,\mu} \lVert \nabla' \ff\rVert_{{\L}^{\infty}(\RR^n_+)}.
\end{align*}
To bound $\partial_{x_n}\uu_n$, we use the divergence free condition:
\begin{align*}
    |\lambda|\lVert \partial_{x_n}\uu_n\rVert_{{\L}^{\infty}(\RR^n_+)} &+ |\lambda|^\frac{1}{2}\lVert \partial_{x_n}\nabla'\uu_n\rVert_{{\L}^{\infty}(\RR^n_+)}+ \lVert \partial_{x_n} \nabla^2 \uu_n\rVert_{\dot{\B}^{0}_{\infty,\infty}(\RR^n_+)}\\
    &=  |\lambda|\lVert \div' \uu'\rVert_{{\L}^{\infty}(\RR^n_+)} +|\lambda|^\frac{1}{2}\lVert \nabla \div' \uu'\rVert_{{\L}^{\infty}(\RR^n_+)} + \lVert \nabla^2\div' \uu'\rVert_{\dot{\B}^{0}_{\infty,\infty}(\RR^n_+)}\\
    &\lesssim_{n}  |\lambda|\lVert \nabla' \uu'\rVert_{{\L}^{\infty}(\RR^n_+)} +|\lambda|^\frac{1}{2}\lVert \nabla \nabla' \uu'\rVert_{{\L}^{\infty}(\RR^n_+)} +\lVert \nabla'\nabla^2 \uu'\rVert_{\dot{\B}^{0}_{\infty,\infty}(\RR^n_+)}\\
    &\lesssim_{s,n,\mu} \lVert \nabla \ff\rVert_{{\L}^{\infty}(\RR^n_+)},
\end{align*}
as well as
\begin{align*}
    \left(\frac{1}{1+|\ln(\min(c,|\lambda|))|}\right)\lVert \partial_{x_n} \nabla^2 \uu_n\rVert_{{\B}^{0}_{\infty,\infty}(\RR^n_+)}\lesssim_{s,n,\mu} \lVert \nabla \ff\rVert_{{\L}^{\infty}(\RR^n_+)},
\end{align*}
By linearity, this also yields
\begin{align*}
    \left(\frac{1}{1+|\ln(\min(c,|\lambda|))|}\right)&\lVert \partial_{x_n}^2\mathfrak{p}\rVert_{{\B}^{0}_{\infty,\infty}(\RR^n_+)}\\
    &\lesssim_{n,c}  |\lambda|\lVert \partial_{x_n}\uu_n\rVert_{{\L}^{\infty}(\RR^n_+)}+ \lVert \partial_{x_n} \ff_n\rVert_{{\L}^{\infty}(\RR^n_+)}\\
    &\qquad\qquad+ \left(\frac{1}{1+|\ln(\min(c,|\lambda|))|}\right)\lVert \partial_{x_n} \Delta \uu_n\rVert_{{\B}^{0}_{\infty,\infty}(\RR^n_+)}\\
    &\lesssim_{s,n,\mu,c} \lVert \nabla \ff\rVert_{{\L}^{\infty}(\RR^n_+)}.
\end{align*}
Therefore, we have obtained
\begin{align*}
    \left(\frac{1}{1+|\ln(\min(c,|\lambda|))|}\right)\lVert \nabla^2\mathfrak{p}\rVert_{{\B}^{0}_{\infty,\infty}(\RR^n_+)} \lesssim_{s,n,\mu,c} \lVert \nabla \ff\rVert_{{\L}^{\infty}(\RR^n_+)}.
\end{align*}
It remains to obtain the estimates on $\partial_{x_n}\uu'=(\partial_{x_n}\uu_j)_{1\leqslant j\leqslant n-1}$. One considers $\omega := \curl \uu(=\d \uu)$, since $\ff\in\W^{1,\infty}_{\mathcal{D},\sigma}(\RR^n_+)$, for $\tilde{\ff}$ the extension of ${\ff}$ to the whole $\RR^n$ by $0$, one has $\tilde{\ff} \in \W^{1,\infty}_{\sigma}(\RR^n)$ and
\begin{align*}
    \curl \tilde{\ff} = \widetilde{\curl \ff} \in  \L^{\infty}(\RR^n)
\end{align*}
with $\supp \curl \tilde{\ff}\subset \overline{\RR^n_+}$, so that
\begin{align*}
    (-\mathfrak{e}_n)\iprod \curl {\ff}_{|_{\partial\RR^n_+}} = 0
\end{align*}
and similarly $(-\mathfrak{e}_n)\iprod \omega_{|_{\partial\RR^n_+}} = 0$. Since $\d \omega = \d^2 \uu =0$, one also obtains for free $(-\mathfrak{e}_n)\iprod \d \omega_{|_{\partial\RR^n_+}} = 0$, so that
\begin{equation*}
    \left\{ \begin{array}{rllr}
         \lambda  \omega  - \Delta  \omega &=  \curl \ff \text{, }&&\text{ in } \RR^n_+\text{,}\\
         (-\mathfrak{e}_n)\iprod \omega_{|_{\partial\RR^n_+}}  &=0\text{, } &&\text{ on } \partial\RR^n_+,\\
         (-\mathfrak{e}_n)\iprod\d \omega_{|_{\partial\RR^n_+}}  &=0\text{, } &&\text{ on } \partial\RR^n_+\text{.}
    \end{array}
    \right.
\end{equation*}
Now, applying Proposition~\ref{prop:HodgeResolventPbRn+Linfty}, we deduce from $(\omega)_{kn}= \partial_{x_n}\uu_k -\partial_{x_k}\uu_n$ and writing $ \partial_{x_n}\uu' = \sum_{k=1}^{n-1}(\omega)_{kn} \mathfrak{e}_k + \nabla'\uu_n$
\begin{align*}
    |\lambda|\lVert \partial_{x_n}\uu' \rVert_{\L^\infty(\RR^n_+)} &+ |\lambda|^\frac{1}{2}\lVert \nabla\partial_{x_n}\uu' \rVert_{\L^\infty(\RR^n_+)} +  \lVert \nabla^2\partial_{x_n}\uu' \rVert_{\dot{\B}^0_{\infty,\infty}(\RR^n_+)}\\
    &\leqslant |\lambda|\lVert \omega \rVert_{\L^\infty(\RR^n_+)} + |\lambda|^\frac{1}{2}\lVert \nabla \omega \rVert_{\L^\infty(\RR^n_+)}+ \lVert \nabla^2 \omega \rVert_{\dot{\B}^0_{\infty,\infty}(\RR^n_+)}\\ & \qquad+ |\lambda|\lVert \nabla'\uu_n \rVert_{\L^\infty(\RR^n_+)} + |\lambda|^\frac{1}{2}\lVert \nabla\nabla'\uu_n\rVert_{\L^\infty(\RR^n_+)}+ \lVert \nabla^2\nabla'\uu_n  \rVert_{\dot{\B}^0_{\infty,\infty}(\RR^n_+)}\\
    &\lesssim_{\mu,n} \lVert \curl \ff \rVert_{\L^\infty(\RR^n_+)} +\lVert \nabla'\ff\rVert_{\L^\infty(\RR^n_+)}\\
    &\lesssim_{\mu,n} \lVert \nabla \ff \rVert_{\L^\infty(\RR^n_+)}.
\end{align*}
Hence, the result claimed in Proposition~\ref{prop:DampedStokesDirRn+Lip}.

\textbf{Step 2:} Now, we prove Proposition~\ref{prop:DampedStokesDirRn+BesovInfty} when $s\in(0,2)$, $s\neq 1$. By Theorem~\ref{thm:StokesDirRn+Linfty} and Step 1, provided first that $s\in(0,1)$ and $q\in[1,\infty]$, by real interpolation --see Theorem~\ref{thm:InterpHomSpacesLip}-- one obtains for any $\lambda\in\Sigma_\mu$, any $\ff\in{\B}^{s,\sigma}_{\infty,q,\mathcal{D}}(\RR^n_+)$
\begin{align*}
    |\lambda|\lVert \uu\rVert_{{\B}^{s}_{\infty,q}(\RR^n_+)} + \lVert \nabla^2 \uu\rVert_{\dot{\B}^{0}_{\infty,\infty}(\RR^n_+)} + \lVert \nabla^2 \uu\rVert_{\dot{\B}^{s}_{\infty,q}(\RR^n_+)} \lesssim_{s,n,\mu} \lVert \ff\rVert_{{\B}^{s}_{\infty,q}(\RR^n_+)},
\end{align*}
and for any $c>0$
\begin{align*}
    \left(\frac{1}{1+|\ln(\min(c,|\lambda|))|}\right)\lVert (\nabla^2 \uu,\nabla\mathfrak{p})\rVert_{{\B}^{s}_{\infty,q}(\RR^n_+)} \lesssim_{s,n,\mu,c} \lVert \ff\rVert_{{\B}^{s}_{\infty,q}(\RR^n_+)}.
\end{align*}
From this point, the estimates when $s\in(1,2)$ are deduced by reproducing \textit{verbatim} Step 1 by means of Corollary~\ref{cor:HodgeResolventPbRn+LinftyBesov} instead of Proposition~\ref{prop:HodgeResolventPbRn+Linfty}.

\textbf{Step 3:} In order to prove Proposition~\ref{prop:DampedStokesDirRn+BesovInfty} when $s\in(-1,0)$. One takes advantage of the following fact which is standard to prove: that for any $-1<s<0$, it holds
\begin{align*}
   {\B}^{s,\sigma}_{\infty,q,\mathfrak{n}}(\RR^n_+) = \L^\infty_{\mathfrak{n,\sigma}}(\RR^n_+) + \dot{\B}^{s,\sigma}_{\infty,q,\mathfrak{n}}(\RR^n_+)
\end{align*}
with equivalence of norms.

Due to real interpolation, we assume without loss of generality $q<\infty$, and consider $\ff\in\C^{\infty}_{ub,0,\sigma}(\RR^n_+)$, the latter space being strongly dense in ${\B}^{s,\sigma}_{\infty,q,\mathfrak{n}}(\RR^n_+)$, by Theorem~\ref{thm:DivergenceFreeSpacesSpeiclaLipDensity} and Proposition~\ref{prop:IdentifVanishingDivFree}. Let $\ff=\mathbf{a}+\mathbf{b}\in\L^\infty_{\mathfrak{n,\sigma}}(\RR^n_+) + \dot{\B}^{s,\sigma}_{\infty,q,\mathfrak{n}}(\RR^n_+)={\B}^{s,\sigma}_{\infty,q,\mathfrak{n}}(\RR^n_+)$, by Theorems~\ref{thm:MetaThmDirichletStokesRn+} and \ref{thm:StokesDirRn+Linfty}, writing $\uu = \uu_{\mathbf{a}}+\uu_{\mathbf{b}}$, it holds
\begin{align*}
    |\lambda|\lVert \uu \rVert_{{\B}^{s}_{\infty,q}(\RR^n_+)} &\lesssim_{n,s} |\lambda|\lVert \uu_\mathbf{a} \rVert_{\L^\infty(\RR^n_+)} + |\lambda|\lVert \uu_\mathbf{b} \rVert_{\dot{\B}^{s}_{\infty,q}(\RR^n_+)}\\
    &\lesssim_{s,n,\mu} \lVert \mathbf{a} \rVert_{\L^\infty(\RR^n_+)} +\lVert \mathbf{b} \rVert_{\dot{\B}^{s}_{\infty,q}(\RR^n_+)}.
\end{align*}
Taking the infimum on all such pairs $(\mathbf{a},\mathbf{b})$, one obtains
\begin{align*}
    |\lambda|\lVert \uu \rVert_{{\B}^{s}_{\infty,q}(\RR^n_+)}
    &\lesssim_{s,n,\mu}\lVert \ff \rVert_{{\B}^{s}_{\infty,q}(\RR^n_+)}.
\end{align*}
So the resolvent operator also extends by density on the whole ${\B}^{s,\sigma}_{\infty,q,\mathfrak{n}}(\RR^n_+)$.

For higher order regularity estimates, we start by recalling that since $-1<s<0$, we also have ${\B}^{s}_{\infty,q}(\RR^n_+) = {\B}^{0}_{\infty,\infty}(\RR^n_+) + \dot{\B}^{s}_{\infty,q}(\RR^n_+)$. Therefore, by Theorems~\ref{thm:MetaThmDirichletStokesRn+} and \ref{thm:StokesDirRn+Linfty} it holds
\begin{align*}
    \lVert \nabla^2  \uu \rVert_{{\B}^{s}_{\infty,q}(\RR^n_+)} &\lesssim_{n,s} \lVert \nabla^2\uu_\mathbf{a} \rVert_{\B^0_{\infty,\infty}(\RR^n_+)} + \lVert \nabla^2\uu_\mathbf{b} \rVert_{\dot{\B}^{s}_{\infty,q}(\RR^n_+)}\\
    &\lesssim_{s,c,n,\mu} (1+|\ln(\min(c,|\lambda|)|)\lVert \mathbf{a} \rVert_{\L^\infty(\RR^n_+)} +\lVert \mathbf{b} \rVert_{\dot{\B}^{s}_{\infty,q}(\RR^n_+)}\\
    &\lesssim_{s,c,n,\mu} (1+|\ln(\min(c,|\lambda|)|)\left[\lVert \mathbf{a} \rVert_{\L^\infty(\RR^n_+)} +\lVert \mathbf{b} \rVert_{\dot{\B}^{s}_{\infty,q}(\RR^n_+)}\right].
\end{align*}
One can now take the infimum on all such pairs $(\mathbf{a},\mathbf{b})$, and obtains, thanks to linearity,
\begin{align*}
    \left(\frac{1}{(1+|\ln(\min(c,|\lambda|)|)}\right)\lVert ( \nabla^2 \uu, \nabla \mathfrak{p} )\rVert_{{\B}^{s}_{\infty,q}(\RR^n_+)}
    &\lesssim_{s,n}\lVert \ff \rVert_{{\B}^{s}_{\infty,q}(\RR^n_+)}.
\end{align*}

It remains to check that $\uu\in\tilde
{\B}^{s+2}_{\infty,q}(\RR^n_+,\CC^n)$ whenever $\ff\in\tilde
{\B}^{s,\sigma}_{\infty,q,\mathcal{D}}(\RR^n_+)$. If $s>0$, $\ff\in\tilde
{\B}^{s,\sigma}_{\infty,q,\mathcal{D}}(\RR^n_+)\subset \C^{0}_{ub,h,0,\sigma}(\RR^n_+)$, by Theorem~\ref{thm:StokesDirRn+Linfty}, one obtains
\begin{align*}
    \uu \in {\B}^{s+2}_{\infty,q}\cap \dot{\B}^{2}_{\infty,\infty}(\RR^n_+,\CC^n) = \tilde{\B}^{s+2}_{\infty,q}(\RR^n_+,\CC^n),
\end{align*}
which yields the result. If $s<0$ and $q<\infty$, one argues by density from the case $s>0$, since $\tilde{\B}^{s+2}_{\infty,q}$ is a closed subspace of ${\B}^{s+2}_{\infty,q}$. The case $q=\infty$ follows from an interpolation argument.
\end{proof}

\begin{proposition}\label{prop:BesovInftySqrtDirStokesRn+} For all $q\in[1,\infty]$, all $-1<s<1$, $s\neq 0$, it holds that 
\begin{align*}
    (\I+\AA_\mathcal{D})^{\sfrac{1}{2}} \,:\,\B^{s+1,\sigma}_{\infty,q,\mathcal{D}}(\RR^n_+)\longrightarrow \B^{s,\sigma}_{\infty,q,\mathcal{D}}(\RR^n_+) 
\end{align*}
is an isomorphism.

Furthermore,
\begin{itemize}
    \item a similar result holds for $\B^{\bullet,0,\sigma}_{\infty,q,\mathcal{D}}(\RR^n_+)$, $\BesSmo^{\bullet,\sigma}_{\infty,\infty,\mathcal{D}}(\RR^n_+)$ and $\BesSmo^{\bullet,0,\sigma}_{\infty,\infty,\mathcal{D}}(\RR^n_+)$;
    \item additionally, one has the mapping property
    \begin{align*}
        (\I+\AA_\mathcal{D})^{\sfrac{1}{2}} \tilde{\B}^{s+1,\sigma}_{\infty,q,\mathcal{D}}(\RR^n_+) = \tilde{\B}^{s,\sigma}_{\infty,q,\mathcal{D}}(\RR^n_+) .
    \end{align*}
\end{itemize}
\end{proposition}

\begin{proof} By Proposition~\ref{prop:DampedStokesDirRn+BesovInfty}, for all $q\in[1,\infty]$, all $-1<s<1$, $s\neq 0$ $\I+\AA_\mathcal{D}$ is an invertible $0$-sectorial operator on $\B^{s,\sigma}_{\infty,q,\mathcal{D}}(\RR^n_+)$, whose domain is a closed subspace of
\begin{align*}
    \B^{s+2,\sigma}_{\infty,q,\mathcal{D}}(\RR^n_+).
\end{align*}

By Theorem~\ref{thm:InterpHomSpacesLip}, due to $-1<s<1$, $s\neq 0$, whenever $s+2\theta<1$, one has 
\begin{align*}
    \B^{s+2\theta,\sigma}_{\infty,q,\mathcal{D}}(\RR^n_+) =\B^{s+2\theta,\sigma}_{\infty,q,0}(\RR^n_+) &= (\B^{s,\sigma}_{\infty,q,0}(\RR^n_+), \B^{s+2,\sigma}_{\infty,q,0}(\RR^n_+))_{\theta,q}\\ 
    &\hookrightarrow (\B^{s,\sigma}_{\infty,q,\mathcal{D}}(\RR^n_+), \D^{s}_{\infty,q}(\I+\AA_\mathcal{D}))_{\theta,q} \\ 
    &\hookrightarrow (\B^{s,\sigma}_{\infty,q,\mathcal{D}}(\RR^n_+),  \B^{s+2,\sigma}_{\infty,q,\mathcal{D}}(\RR^n_+))_{\theta,q} \\ 
    &\hookrightarrow \B^{s+2\theta,\sigma}_{\infty,q,\mathcal{D}}(\RR^n_+),
\end{align*}
the second embedding being obtained through Proposition~\ref{prop:DampedStokesDirRn+BesovInfty}. It yields the following equality with equivalence of norms provided, $-1<s<1$, $s\neq 0$,  $s+2\theta<1$:
\begin{align*}
    \B^{s+2\theta,\sigma}_{\infty,q,\mathcal{D}}(\RR^n_+)=(\B^{s,\sigma}_{\infty,q,\mathcal{D}}(\RR^n_+), \D^{s}_{\infty,q}(\I+\AA_\mathcal{D}))_{\theta,q}.
\end{align*}
By \cite[Proposition~6.4.1,~a),~\&~Corollary~6.5.5]{bookHaase2006}, for any $\theta\in(\frac{1}{2},1)$, $q\in[1,\infty]$, one has an isomorphism
\begin{align*}
   (\I+\AA_\mathcal{D})^\frac{1}{2}\,:\,(\B^{s,\sigma}_{\infty,q,\mathcal{D}}(\RR^n_+), \D^{s}_{\infty,q}(\I+\AA_\mathcal{D}))_{\theta,q} \longrightarrow (\B^{s,\sigma}_{\infty,q,\mathcal{D}}(\RR^n_+), \D^{s}_{\infty,q}(\I+\AA_\mathcal{D}))_{\theta-\frac{1}{2},q},
\end{align*}
which combined with the preceding interpolation identity yields the result.

Finally, the fact that it maps $\tilde{\B}^{s+1,\sigma}_{\infty,q,\mathcal{D}}(\RR^n_+)$ to $\tilde{\B}^{s,\sigma}_{\infty,q,\mathcal{D}}(\RR^n_+)$, follows from the fact that for all $\lambda\in\Sigma_\mu\cup\{0\}$, according to Proposition~\ref{prop:DampedStokesDirRn+BesovInfty}, the resolvent preserves the $\tilde{\B}^{\bullet,\sigma}_{\infty,q}$-spaces:
\begin{align*}
    (\lambda\I+\I+\AA_\mathcal{D})^{-1}\tilde{\B}^{s,\sigma}_{\infty,q,\mathcal{D}}(\RR^n_+) \subset \tilde{\B}^{s+2,\sigma}_{\infty,q,\mathcal{D}}(\RR^n_+).
\end{align*}
Thereby, so do the corresponding fractional powers.
\end{proof}

\subsection{Miscellaneous: A \texorpdfstring{$\mathcal{R}-$}{R-}bound on homogeneous Sobolev spaces}\label{sec:RboundRn+Hsp}

With an interest for the bounded $\mathbf{H}^\infty$-functional calculus and the $\L^q$-maximal regularity property, from the representation formula exposed in the proof of Theorem~\ref{thm:MetaThmDirichletStokesRn+}, we derive the following $\mathcal{R}$-bound. It will be of use to extrapolate the boundedness of the $\mathbf{H}^{\infty}$-functional calculus of the Stokes operator in bounded rougher domains in Section~\ref{sec:ProofLPbddDomain}, following \cite{KunstmannWeis2017,GabelTolksdorf2022}.

\begin{proposition}\label{prop:RboundHspRn+}Let $p\in(1,\infty)$, $s\in(-1+\sfrac{1}{p},\sfrac{1}{p})$. Let $\mu\in(0,\pi)$, and $(\lambda_\ell)_{\ell\in\NN}\subset\Sigma_{\mu}$. For all $(\ff_{\ell})_{\ell\in\NN}\subset\dot{\H}^{s,p}(\RR^n_+,\CC^n)$, $(g_{\ell})_{\ell\in\NN}\subset\dot{\H}^{s-1,p}_0\cap\dot{\H}^{s+1,p}(\RR^n_+)$, consider $(\uu_\ell,\nabla{\mathfrak{p}}_{\ell})_{\ell\in\NN}$ the corresponding family of solutions provided by Theorem~\ref{thm:MetaThmDirichletStokesRn+}. Then for all $N\in\NN$, one has
\begin{align*}
    \bigg\lVert \sum_{\ell=1}^{N} r_\ell|\lambda_\ell| \uu_\ell\bigg\rVert_{\L^p(0,1;\dot{\H}^{s,p}(\RR^n_+))} + &\bigg\lVert \sum_{\ell=1}^{N} r_\ell (\nabla^2 \uu_\ell,\nabla\mathfrak{p}_\ell)\bigg\rVert_{\L^p(0,1;\dot{\H}^{s,p}(\RR^n_+))}\\ \qquad & \lesssim_{p,s,n,\mu}\bigg\lVert \sum_{\ell=1}^{N} r_\ell\ff_\ell\bigg\rVert_{\L^p(0,1;\dot{\H}^{s,p}(\RR^n_+))} \\ & + \bigg\lVert \sum_{\ell=1}^{N} r_\ell |\lambda_\ell| g_\ell\bigg\rVert_{\L^p(0,1;\dot{\H}^{s-1,p}(\RR^n_+))} + \bigg\lVert \sum_{\ell=1}^{N} r_\ell \nabla g_\ell\bigg\rVert_{\L^p(0,1;\dot{\H}^{s,p}(\RR^n_+))}.
\end{align*}
\end{proposition}

Before providing the proof, we recall several known facts for $\mathcal{R}$-boundedness of families of operators, and their direct consequences. The first proposition is about fundamental knowledge the second is about direct consequences. See \cite[Chapter~2]{KunstmannWeis2004}, \cite[Chapter~4,~Section~4.1]{PrussSimonett2016} and  \cite[Chapter~8]{HytonenNeervenVeraarWeisbookVolII2018} for more details.

\begin{proposition}\label{prop:RboundednessFudnamentalProp}Let $\X$, $\Y$ and $\Z$ be three Banach spaces. It holds that
\begin{enumerate}
    \item Provided $T\,:\,\X\longrightarrow\Y$ is linear and bounded then the singleton family $\{T\}$ is $\mathcal{R}$-bounded with the estimate
    \begin{align*}
        \mathcal{R}\{T\} \leqslant \lVert T\rVert_{\X\rightarrow\Y};
    \end{align*}
    \item Provided $\mathcal{T}$ and $\mathcal{S}$ are both $\mathcal{R}$-bounded from $\X$ to $\Y$, then it holds that the families
    \begin{align*}
        \mathcal{S}\cup\mathcal{T},\,\quad\text{ and }\quad \mathcal{S}+\mathcal{T}:=\{\,S+T,\,(S,T)\in\mathcal{T}\times\mathcal{S}\}
    \end{align*}
    are both $\mathcal{R}$-bounded from $\X$ to $\Y$ with the $\mathcal{R}$-bounds
    \begin{align*}
        \mathcal{R}(\mathcal{S}\cup\mathcal{T})\leqslant\mathcal{R}(\mathcal{S}+\mathcal{T})\leqslant \mathcal{R}(\mathcal{S})+\mathcal{R}(\mathcal{T});
    \end{align*}
    \item Provided $\mathcal{T}$ and $\mathcal{S}$ are respectively $\mathcal{R}$-bounded from $\X$ to $\Y$ and from $\Y$ to $\Z$, then it holds that the family
    \begin{align*}
        \mathcal{S}\mathcal{T}:=\{\,ST,\,(S,T)\in\mathcal{T}\times\mathcal{S}\}
    \end{align*}
    is $\mathcal{R}$-bounded from $\X$ to $\Z$ with the $\mathcal{R}$-bound
    \begin{align*}
        \mathcal{R}(\mathcal{S}\mathcal{T})\leqslant \mathcal{R}(\mathcal{S})\mathcal{R}(\mathcal{T});
    \end{align*}
    \item Let $p\in(1,\infty)$, $\mu\in[0,\pi)$, $N\geqslant 1$, it holds that the following families are $\mathcal{R}$-bounded on $\L^p(\RR^n,\CC^N)$:
    \begin{itemize}
        \item $\{\lambda(\lambda\I-\Delta)^{-1},\, \lambda\in\Sigma_\mu\}$;
        \item $\{\partial_{x_k}\partial_{x_j}(\lambda\I-\Delta)^{-1},\, \lambda\in\Sigma_\mu\}$, for any $k,j\in\llb1,n\rrb$;
        \item $\{\partial_{x_k}(\lambda\I-\Delta)^{-\frac{1}{2}},\, \lambda\in\Sigma_\mu\}$, for any $k\in\llb1,n\rrb$;
        \item $\{\lambda^\frac{1}{2}(\lambda\I-\Delta)^{-\frac{1}{2}},\, \lambda\in\Sigma_\mu\}$;
        \item $\{(\lambda\I-\Delta')^\frac{1}{2}(\lambda\I-\Delta)^{-\frac{1}{2}},\, \lambda\in\Sigma_\mu\}$;
        \item $\{(-\Delta')^\frac{1}{2}(\lambda\I-\Delta)^{-\frac{1}{2}},\, \lambda\in\Sigma_\mu\}$, for any $k,j\in\llb1,n\rrb$;
    \end{itemize}
and everything remains valid by replacing arbitrarily some $\lambda$-expressions by $|\lambda|$.
\end{enumerate}
\end{proposition}

\begin{corollary}\label{cor:RboundsTwistedNeumanResolvHspRn+} Let $p\in(1,\infty)$, $-1+\sfrac{1}{p}<s<\sfrac{1}{p}$ and $N\in\NN^\ast$. Then it holds that, for all $\mu\in[0,\pi)$, the following sets of operators are $\mathcal{R}$-bounded on $\dot{\H}^{s,p}(\RR^n_+,\CC^N)$:
\begin{itemize}
    \item $\{ \,\lambda(|\lambda| \I-\Delta_\mathcal{N})^{-1},\,\lambda\in\Sigma_\mu\}$,
    \item  $\{ \,(-\Delta')^\frac{1}{2}[(\lambda\I-\Delta')^\frac{1}{2}+(-\Delta')^\frac{1}{2}](|\lambda| \I-\Delta_\mathcal{N})^{-1},\,\lambda\in\Sigma_\mu\}$;
    \item $\{ \,(\lambda\I-\Delta')^\frac{1}{2}[(\lambda\I-\Delta')^\frac{1}{2}+(-\Delta')^\frac{1}{2}](|\lambda| \I-\Delta_\mathcal{N})^{-1},\,\lambda\in\Sigma_\mu\}$;
    \item $\{ \,\partial_{x_n}[(\lambda\I-\Delta')^\frac{1}{2}+(-\Delta')^\frac{1}{2}](|\lambda| \I-\Delta_\mathcal{N})^{-1},\,\lambda\in\Sigma_\mu\}$.
\end{itemize}
\end{corollary}

\begin{proof} Recall that the resolvent for the Neumann Laplacian can be written as
\begin{align*}
    (\lambda\I-\Delta_\mathcal{N})^{-1}=\R_{\RR^n_+}(\lambda\I-\Delta)^{-1}\E_{\mathcal{N}},
\end{align*}
where $\E_\mathcal{N}$ is bounded from $\dot{\H}^{s,p}(\RR^n_+,\CC)$ to $\dot{\H}^{s,p}(\RR^n,\CC)$. Therefore, this reduces to the study of same families of operators on $\RR^n$ with the Laplacian on the whole space instead of the Neumann Laplacian, asking then for $\mathcal{R}$-boundedness on $\dot{\H}^{s,p}(\RR^n,\CC)$. Due to the isomorphism $(-\Delta)^\frac{s}{2}\,:\,\dot{\H}^{s,p}(\RR^n,\CC)\longrightarrow\L^p(\RR^n,\CC)$, this reduces to study of the same families in the $\L^p(\RR^n,\CC)$-case, for which such properties are known to hold due to Proposition~\ref{prop:RboundednessFudnamentalProp}.
\end{proof}

\begin{proof}[of Proposition~\ref{prop:RboundHspRn+}] \textbf{Step 1:} We start by analyzing the problem to reduce the problem to the $\mathcal{R}$-boundedness of simpler families of operators. If $(g_{\ell})_{\ell\in\NN}=0$ the result follows from the functional calculus identity from Theorem~\ref{thm:MetaThmDirichletStokesRn+}
\begin{align*}
    \uu_\ell = (\lambda_\ell\I+\AA_\mathcal{D})^{-1}\PP_{\RR^{n}_+}\ff_\ell = \R_{\RR^n_+}\Gamma^\mathcal{U}(\lambda_\ell\I-\Delta)^{-1}\E_{\sigma}^{\mathcal{U}}\PP_{\RR^{n}_+}\ff_\ell
\end{align*}
and the $\mathcal{R}$-boundedness of the Laplacian on $\RR^n$, and the fact $\R_{\RR^n_+},\Gamma^\mathcal{U},\E_{\sigma}^{\mathcal{U}}$ and $\PP_{\RR^{n}_+}$ are singleton operators, see  Proposition~\ref{prop:RboundednessFudnamentalProp}.

\medbreak

\noindent Hence, by linearity and without loss of generality we can assume $\ff_\ell=0$.  For $(\lambda,\uu,\mathfrak{p},g)\in\{(\lambda_\ell,\uu_\ell,\mathfrak{p}_\ell,g_\ell),\, \ell\in\NN\}$, setting $\ww=(-\Delta_{\mathcal{N}})^{-1}g$, following Steps 3.1 and 3.2 in the proof of Theorem~\ref{thm:MetaThmDirichletStokesRn+}, one can write
\begin{align*}
    &\uu'(x',x_n)\\
    &\,\,:=\frac{1}{\lambda} [(\lambda\I-\Delta')^{{\sfrac{1}{2}}}+(-\Delta')^{{\sfrac{1}{2}}}]((\lambda\I-\Delta')^{{\sfrac{1}{2}}}e^{-x_n(\lambda\I-\Delta')^{{\sfrac{1}{2}}}}-(-\Delta')^{{\sfrac{1}{2}}}e^{-x_n(-\Delta')^{{\sfrac{1}{2}}}})[\nabla'\ww](x',0)\\ &\qquad\qquad -\nabla'\ww(x',x_n)\\
   &\uu_n (x',x_n)\\
   &\,\,:= -\frac{1}{\lambda}[(\lambda\I-\Delta')^{{\sfrac{1}{2}}}+(-\Delta')^{{\sfrac{1}{2}}}](e^{-x_n(\lambda\I-\Delta')^{{\sfrac{1}{2}}}}-e^{-x_n(-\Delta')^{{\sfrac{1}{2}}}})[(-\Delta')\ww](x',0)\\&\qquad\qquad-\partial_{x_n}\ww(x',x_n) ;\\
   &\mathfrak{p} (x',x_n) := -(-\Delta')^{{\sfrac{1}{2}}}[(\lambda\I-\Delta')^{{\sfrac{1}{2}}}+(-\Delta')^{{\sfrac{1}{2}}}]e^{-x_n(-\Delta')^{{\sfrac{1}{2}}}}\ww(x',0)+\lambda\ww(x',x_n).
\end{align*}
By linearity it suffices to prove the $\mathcal{R}$-bound on $(\lambda_\ell\uu_\ell,\nabla^2\uu_\ell)_{\ell\in\NN}$. Furthermore, we take a focus on $\uu'$, $\uu_n$ admits a similar treatment. Using partial Fourier transform, one may write
\begin{align*}
    \mathcal{F}'\uu'(\xi',x_n) &= -\int_{0}^\infty \partial_{y_n}m_\lambda(\xi',x_n+y_n)i\xi'\mathcal{F}'\ww(\xi',y_n) \d y_n\\
    &\qquad-\int_{0}^{\infty} m_\lambda(\xi',x_n+y_n)i\xi'\mathcal{F}'\partial_{x_n}\ww(\xi',y_n)\d y_n-i\xi'\mathcal{F}'\partial_{x_n}\ww(\xi',x_n)
\end{align*}
where
\begin{align*}
    m_\lambda(\xi',z) := \frac{1}{\lambda} [(\lambda\I+|\xi'|^2)^\frac{1}{2} + |\xi'|][(\lambda\I+|\xi'|^2)^\frac{1}{2}e^{-z(\lambda\I+|\xi'|^2)^\frac{1}{2}}- |\xi'|e^{-z|\xi'|}]
\end{align*}
and similarly for $\uu_n$. Using this trick, setting $A_\lambda=(\lambda\I-\Delta')^\frac{1}{2}$, one can rewrite
\begin{align*}
    &\uu'(\cdot',x_n)= \nabla' \ww(\cdot',x_n) \\
    &+ \frac{1}{\lambda} [A_\lambda+A_0] \Bigg(\int_{-\infty}^{x_n} A_\lambda^2 e^{-(x_n-y_n)A_\lambda}\nabla'\ww(\cdot',-y_n) \d y_n -\int_{-\infty}^{x_n} A_0^2 e^{-(x_n-y_n)A_0}\nabla'\ww(\cdot',-y_n) \d y_n\\
    &\qquad-\int_{0}^{x_n} A_\lambda^2 e^{-(x_n-y_n)A_\lambda}\nabla'\ww(\cdot',-y_n) \d y_n +\int_{0}^{x_n} A_0^2 e^{-(x_n-y_n)A_0}\nabla'\ww(\cdot',-y_n) \d y_n\\
    & \qquad +\int_{-\infty}^{x_n} A_\lambda e^{-(x_n-y_n)A_\lambda} \partial_{x_n}\nabla'\ww(\cdot',-y_n) \d y_n -\int_{-\infty}^{x_n} A_0 e^{-(x_n-y_n)A_0} \partial_{x_n}\nabla'\ww(\cdot',-y_n) \d y_n\\
    &\qquad-\int_{0}^{x_n} A_\lambda e^{-(x_n-y_n)A_\lambda}\partial_{x_n}\nabla'\ww(\cdot',-y_n) \d y_n +\int_{0}^{x_n} A_0 e^{-(x_n-y_n)A_0}\partial_{x_n}\nabla'\ww(\cdot',-y_n) \d y_n\Bigg)
\end{align*}
We continue to rewrite differently each term, writing first
\begin{align*}
    \nabla \ww = \nabla(-\Delta_\mathcal{N})^{-1}g &=\nabla(|\lambda| \I-\Delta_\mathcal{N})^{-1}(\lambda \I-\Delta_\mathcal{N})(-\Delta_\mathcal{N})^{-1}g\\
    &=   \nabla(-\Delta_\mathcal{N})^{-\frac{1}{2}}(|\lambda| \I-\Delta_\mathcal{N})^{-1}(-\Delta_\mathcal{N})^{-\frac{1}{2}}(|\lambda| \I-\Delta_\mathcal{N})g\\
    &=   [A_\lambda+A_0]^{-1}\nabla(-\Delta_\mathcal{N})^{-\frac{1}{2}}[A_\lambda+A_0](|\lambda| \I-\Delta_\mathcal{N})^{-1}(-\Delta_\mathcal{N})^{-\frac{1}{2}}(|\lambda| \I-\Delta_\mathcal{N})g.
\end{align*}
Note that by Proposition~\ref{prop:IsomNeuMannHomSobSpaRn+} $(-\Delta_\mathcal{N})^{-\frac{1}{2}}(|\lambda| \I-\Delta_\mathcal{N})$ is an isomorphism from $\dot{\H}^{s-1,p}_0\cap\dot{\H}^{s+1,p}(\RR^n_+)$ to $\dot{\H}^{s,p}(\RR^n_+)$. Now, we introduce some notations:
\begin{itemize}
    \item $R_{\mathcal{N}}:=\nabla(-\Delta_\mathcal{N})^{-\frac{1}{2}}$;
    \item $R_{\mathcal{N}}^{\zeta}(\lambda) := A_\zeta[A_\lambda+A_0](|\lambda| \I-\Delta_\mathcal{N})^{-1}$, $\zeta\in\{0,\lambda\}$;
    \item $R_{\mathcal{N}}^{n}(\lambda) := \partial_{x_n}[A_\lambda+A_0](|\lambda| \I-\Delta_\mathcal{N})^{-1}$;
    \item $g_\lambda := (-\Delta_\mathcal{N})^{-\frac{1}{2}}(|\lambda| \I-\Delta_\mathcal{N})g\in\dot{\H}^{s,p}(\RR^n_+)$;
    \item $\mathcal{M}_{\zeta,\mathfrak{a}} f(x_n):= \int_{\mathfrak{a}}^{x_n} A_\zeta e^{-(x_n-y_n)A_\zeta}g(-y_n) \d y_n$, provided $x_n\geqslant 0$, $\zeta\in\{0,\lambda\}$, $\mathfrak{a}\in\{-\infty,0\}$.
\end{itemize}
Hence, provided $\mathfrak{a}\in\{-\infty,0\}$, $\zeta\in\{0,\lambda\}$, one has
\begin{align*}
    &[A_\lambda+A_0]\int_{\mathfrak{a}}^{x_n} A_\zeta^2 e^{-(x_n-y_n)A_\zeta}\nabla'\ww(\cdot',-y_n) \d y_n\\
    &= \int_{\mathfrak{a}}^{x_n} A_\zeta e^{-(x_n-y_n)A_\lambda}[\nabla'(-\Delta_\mathcal{N})^{-\frac{1}{2}}A_\zeta[A_\lambda+A_0](|\lambda| \I-\Delta_\mathcal{N})^{-1}g_\lambda](\cdot',-y_n) \d y_n \\
    &= \int_{\mathfrak{a}}^{x_n} A_\lambda e^{-(x_n-y_n)A_\lambda}[R_{\mathcal{N}}'R_{\mathcal{N}}^{\zeta}(\lambda)g_\lambda](\cdot',-y_n) \d y_n \\
    &= \mathcal{M}_{\zeta,\mathfrak{a}}R_{\mathcal{N}}'R_{\mathcal{N}}^{\zeta}(\lambda)g_\lambda.
\end{align*}
Thus, we write
\begin{align}\label{eq:ProofrBOUNDrewrited}
    \uu' = -\frac{1}{\lambda}R_\mathcal{N}'\lambda(|\lambda| \I-\Delta_\mathcal{N})^{-1}g_\lambda + \frac{1}{\lambda}\sum_{\substack{\zeta\in\{0,\lambda\},\\\mathfrak{a}\in\{-\infty,0\}}} \pm\mathcal{M}_{\zeta,\mathfrak{a}}R_{\mathcal{N}}'R_{\mathcal{N}}^{\zeta}(\lambda)g_\lambda \mp \mathcal{M}_{\zeta,\mathfrak{a}}R_{\mathcal{N}}^{n}(\lambda)R_{\mathcal{N}}'g_\lambda.
\end{align}
Since by Corollary~\ref{cor:RboundsTwistedNeumanResolvHspRn+}, the families $\{ \lambda(|\lambda| \I-\Delta_\mathcal{N})^{-1},\,\lambda\in\Sigma_\mu\}$, $\{ R_{\mathcal{N}}^{0}(\lambda),\,\lambda\in\Sigma_\mu\}$, $\{ R_{\mathcal{N}}^{\lambda}(\lambda),\,\lambda\in\Sigma_\mu\}$ and \{$R_{\mathcal{N}}^{n}(\lambda),\,\lambda\in\Sigma_\mu\}$ are known to be $\mathcal{R}$-bounded families of operators on $\dot{\H}^{s,p}(\RR^n_+)$, and since $\{\mathcal{M}_{0,0}\}$ and $\{\mathcal{M}_{-\infty,0}\}$ are singleton operators, the problem reduce to the knowledge of the $\mathcal{R}$-boundedness of the two families
\begin{align*}
    \{ \mathcal{M}_{\lambda,\mathfrak{a}},\,\lambda\in\Sigma_\mu\},\qquad \mathfrak{a}\in\{-\infty,0\}.
\end{align*}
\textbf{Step 2:} We did reduce the problem to the following
\begin{itemize}
    \item Proving the boundedness of $\mathcal{M}_{0,0}$ and $\mathcal{M}_{0,-\infty}$ on $\dot{\H}^{s,p}(\RR^n_+)$;
    \item Proving the $\mathcal{R}$-boundedness of $\{ \mathcal{M}_{\lambda,\mathfrak{a}},\,\lambda\in\Sigma_\mu\},$ for $ \mathfrak{a}\in\{-\infty,0\}$ on $\dot{\H}^{s,p}(\RR^n_+)$.
\end{itemize}
\textbf{Step 2.1:} We prove boundedness of $\mathcal{M}_{0,0}$ and $\mathcal{M}_{0,-\infty}$ on $\dot{\H}^{s,p}(\RR^n_+)$, recalling that for all (say for simplicity) $f\in\Ccinfty(\RR^n)$, one writes for all $a\in(-\infty,0]$ (not to be mistaken with $\mathfrak{a}\in\{-\infty,0\}$ !)
\begin{align*}
    \mathcal{M}_{0,a}f(x_n) = \int_{a}^{x_n} A_0e^{-(x_n-y_n)A_0}f(\cdot',-y_n)\d y_n.
\end{align*}

Recall that $A_0=(-\Delta')^\frac{1}{2}$ has bounded $\mathbf{H}^\infty(\Sigma_\theta)$-functional calculus on $\dot{\H}^{s,p}(\RR^{n-1})$, $-1+\sfrac{1}{p}<s<\sfrac{1}{p}$, for all $\theta\in(0,\pi)$. Therefore, it has BIP on it, and by global-in-time homogeneous Sobolev maximal regularity \cite[Theorem~4.7]{Gaudin2023}, for all $-1+\sfrac{1}{p}<\beta,s<\sfrac{1}{p}$, for all $a\leqslant 0$,
\begin{align*}
    \lVert \mathcal{M}_{0,a}f \lVert_{\dot{\H}^{\beta,p}(\RR,\dot{\H}^{s,p}(\RR^{n-1}))} &\lesssim_{\beta,n,p} \lVert \mathcal{M}_{0,a}f \lVert_{\dot{\H}^{\beta,p}([a,\infty),\dot{\H}^{s,p}(\RR^{n-1}))}\\
    &\lesssim_{\beta,n,p} \lVert f(\cdot',-\cdot) \lVert_{\dot{\H}^{\beta,p}([a,\infty),\dot{\H}^{s,p}(\RR^{n-1}))}\\
    &\lesssim_{\beta,n,p} \lVert f \lVert_{\dot{\H}^{\beta,p}(\RR,\dot{\H}^{s,p}(\RR^{n-1}))}.
\end{align*}
Before we continue, note that thanks to \cite[Proposition~3.6~(ii)]{Gaudin2023}, we did extend $f$ to the whole line by $0$, preserving the $\dot{\H}^{\beta,p}(\RR,\dot{\H}^{s,p}(\RR^{n-1}))$-regularity, but also $\mathcal{M}_{0,a}f$ from $[a,\infty)$ to the whole line by $0$, and the bound of the extension does not depend on $a\in(-\infty,0]$, due to the invariance by translation.

Now, if $0\leqslant s< \sfrac{1}{p}$,
\begin{align*}
    \lVert \mathcal{M}_{0,a}f \lVert_{\dot{\H}^{s,p}(\RR^{n})} &\lesssim_{s,n,p} \lVert \mathcal{M}_{a,0}f \lVert_{\dot{\H}^{s,p}(\RR,\L^p(\RR^{n-1}))}+ \lVert \mathcal{M}_{0,a}f \lVert_{\L^p(\RR,\dot{\H}^{s,p}(\RR^{n-1}))}\\
    &\lesssim_{s,n,p} \lVert f \lVert_{\dot{\H}^{s,p}(\RR,\L^p(\RR^{n-1}))}+ \lVert f \lVert_{\L^p(\RR,\dot{\H}^{s,p}(\RR^{n-1}))}\\
    &\lesssim_{s,n,p} \lVert f \lVert_{\dot{\H}^{s,p}(\RR^{n})}.
\end{align*}
Now, if $0<s<1-\sfrac{1}{p}$, we write $f=g+h$ for $h_0\in\dot{\H}^{-s,p}(\RR,\L^p(\RR^{n-1}))$, $h_1\in\L^p(\RR,\dot{\H}^{-s,p}(\RR^{n-1}))$, it holds
\begin{align*}
    \lVert \mathcal{M}_{0,a}f \lVert_{\dot{\H}^{-s,p}(\RR^{n})} &\leqslant \lVert \mathcal{M}_{0,a}h_0 \lVert_{\dot{\H}^{-s,p}(\RR^{n})}+\lVert \mathcal{M}_{0,a}h_1 \lVert_{\dot{\H}^{-s,p}(\RR^{n})}\\
    &\lesssim_{s,n,p}\lVert \mathcal{M}_{0,a}h_0 \lVert_{\dot{\H}^{-s,p}(\RR,\L^p(\RR^{n-1}))}+ \lVert \mathcal{M}_{a,0}h_1 \lVert_{\L^p(\RR,\dot{\H}^{-s,p}(\RR^{n-1}))}\Big)\\
    &\lesssim_{s,n,p} \lVert h_0 \lVert_{\dot{\H}^{-s,p}(\RR,\L^p(\RR^{n-1}))}+\lVert h_1 \lVert_{\L^p(\RR,\dot{\H}^{-s,p}(\RR^{n-1}))},
\end{align*}
taking the infimum on all pairs $(h_0,h_1)$ yields
\begin{align*}
    \lVert \mathcal{M}_{0,a}f \lVert_{\dot{\H}^{-s,p}(\RR^{n})}\lesssim_{s,n,p} \lVert f \lVert_{\dot{\H}^{-s,p}(\RR^{n})},\qquad -s\in(-1+\sfrac{1}{p},0).
\end{align*}
Therefore for any $s\in(-1+\sfrac{1}{p},\sfrac{1}{p})$, any $f\in\Ccinfty(\RR^n)$, since all the previous bounds were independent of $a\in(-\infty,0]$, it holds
\begin{align*}
    \lVert \mathcal{M}_{0,-\infty}f \lVert_{\dot{\H}^{s,p}(\RR^{n})}\leqslant \liminf_{a\longrightarrow -\infty}\,\lVert \mathcal{M}_{0,a}f \lVert_{\dot{\H}^{s,p}(\RR^{n})}\lesssim_{s,n,p} \lVert f \lVert_{\dot{\H}^{s,p}(\RR^{n})}.
\end{align*}
One can conclude by density of $\Ccinfty(\RR^n)$ in $\dot{\H}^{s,p}(\RR^n)$.

\textbf{Step 2.2:} We prove the $\mathcal{R}$-boundedness of the family $\{ \mathcal{M}_{\lambda,-\infty},\,\lambda\in\Sigma_\mu\}$ on $\dot{\H}^{s,p}(\RR^n)$. Since $\mathcal{M}_{\lambda,-\infty}$ is a convolution operator, one has $(-\Delta)^\frac{s}{2}\mathcal{M}_{\lambda,-\infty}=\mathcal{M}_{\lambda,-\infty}(-\Delta)^\frac{s}{2}$, so it suffices to prove it on $\L^p(\RR^n)$. We write thanks to the Fubini theorem, $\L^p(\RR^n)=\L^p(\RR,\L^p(\RR^{n-1}))$, we take the Fourier transform with respect to $x_n$  and we obtain for all $f\in\S(\RR,\L^p(\RR^{n-1}))$,
\begin{align*}
    \mathcal{F}_n[\mathcal{M}_{\lambda,-\infty}f](x',\tau) = (\lambda \I-\Delta')^{\frac{1}{2}}(i\tau+(\lambda \I-\Delta')^{\frac{1}{2}})^{-1}\mathcal{F}_nf(x',\tau).
\end{align*}
We want to check that the families
\begin{align*}
    &\{ \,(\lambda \I-\Delta')^{\frac{1}{2}}(i\tau+(\lambda \I-\Delta')^{\frac{1}{2}})^{-1},\,\lambda\in\Sigma_\mu\},\\
    \text{ and }\quad&\{ \,i\tau(\lambda \I-\Delta')^{\frac{1}{2}}(i\tau+(\lambda \I-\Delta')^{\frac{1}{2}})^{-2},\,\lambda\in\Sigma_\mu\},
\end{align*}
are both $\mathcal{R}$-bounded on $\L^p(\RR^{n-1})$ with a $\mathcal{R}$-bound uniform with respect to $\tau\in\RR$, in which case the result would follow immediately from \cite[Statement~5.2,~Fourier~Multipliers,~(b)]{KunstmannWeis2004}. Applying Fourier Transform with respect to $x'\in\RR^{n-1}$, it is then sufficient to prove a bound uniform with respect to $\tau\in\RR$ and $\lambda\in\Sigma_\mu$ for the families
\begin{align*}
    &\{ \,(\lambda \I+|\xi'|^2)^{\frac{1}{2}}(i\tau+(\lambda \I+|\xi'|^2)^{\frac{1}{2}})^{-1},\,\lambda\in\Sigma_\mu\},\\
    \text{ and }\quad&\{ \,i\tau(\lambda \I+|\xi'|^2)^{\frac{1}{2}}(i\tau+(\lambda \I+|\xi'|^2)^{\frac{1}{2}})^{-2},\,\lambda\in\Sigma_\mu\},
\end{align*}
and their higher order partial derivatives with respect to $\xi'$ up to the order $n$, looking for inequalities
\begin{align*}
    &|\xi'|^k\nabla^{k}_{\xi'}\left[\frac{(\lambda \I+|\xi'|^2)^{\frac{1}{2}}}{i\tau+(\lambda \I+|\xi'|^2)^{\frac{1}{2}}}\right] \lesssim_{\mu,n} 1,\qquad k\in\llb 1,n\rrb;\\
    \text{ and }\qquad& |\xi'|^k\nabla^{k}_{\xi'}\left[\frac{i\tau(\lambda \I+|\xi'|^2)^{\frac{1}{2}}}{(i\tau+(\lambda \I+|\xi'|^2)^{\frac{1}{2}})^2}\right] \lesssim_{\mu,n} 1,\qquad k\in\llb 1,n\rrb.
\end{align*}
This is true and can be checked by induction. Thus, by \cite[Statement~5.2,~Fourier~Multipliers,~(a)]{KunstmannWeis2004}, one has $\mathcal{R}$-bounded families on $\L^p(\RR^{n-1})$:
\begin{align*}
    &\sup_{\tau \in\RR} \,\mathcal{R}\left\{ \,(\lambda \I+|\xi'|^2)^{\frac{1}{2}}(i\tau+(\lambda \I+|\xi'|^2)^{\frac{1}{2}})^{-1},\,\lambda\in\Sigma_\mu\right\} \lesssim_{\mu,n} 1,\\
    \text{ and }\qquad&\sup_{\tau \in\RR} \,\mathcal{R}\left\{ \,i\tau(\lambda \I+|\xi'|^2)^{\frac{1}{2}}(i\tau+(\lambda \I+|\xi'|^2)^{\frac{1}{2}})^{-2},\,\lambda\in\Sigma_\mu\right\} \lesssim_{\mu,n} 1.
\end{align*}
Now applying \cite[Statement~5.2,~Fourier~Multipliers,~(a)]{KunstmannWeis2004}, on $\L^p(\RR,\L^p(\RR^{n-1}))=\L^p(\RR^n)$, we obtain the $\mathcal{R}$-bounded family
\begin{align*}
    \mathcal{R}\left\{ \,\mathcal{M}_{\lambda,-\infty},\,\lambda\in\Sigma_\mu\right\} \lesssim_{p,\mu,n} 1.
\end{align*}

\textbf{Step 3:} The $\mathcal{R}$-bounds. Let $N\in\NN$, and let $(\lambda_\ell,\uu_\ell,\mathfrak{p}_\ell,g_\ell)_{\ell\in\NN}$. We aim to prove the $\mathcal{R}$-bounds for $(\lambda_\ell \uu_\ell)_{\ell\in\NN}$ and $(\nabla^2\uu_\ell)_{\ell\in\NN}$. Before we start, for any $\ell\in\NN$, here, one has $g_{\lambda_\ell} = (-\Delta_{\mathcal{N}})^{-1}(|\lambda_\ell|\I-\Delta)g_\ell$.

\textbf{Step 3.1:} We provide the $\mathcal{R}$-bound for $(\uu'_\ell)_{\ell\in\NN}$, the one for $(\uu_{n,\ell})_{\ell\in\NN}$ can be obtained similarly. Thanks to \eqref{eq:ProofrBOUNDrewrited},
\begin{align*}
    &\bigg\lVert \sum_{\ell=1}^{N} r_\ell|\lambda_\ell|\uu_\ell'\bigg\rVert_{\L^p(0,1;\dot{\H}^{s,p}(\RR^n_+))}\\
    &\leqslant \bigg\lVert \sum_{\ell=1}^{N} r_\ell\frac{|\lambda_\ell|}{\lambda_\ell}R_\mathcal{N}'\lambda_\ell(|\lambda_\ell| \I-\Delta_\mathcal{N})^{-1}g_\lambda \bigg\rVert_{\L^p(0,1;\dot{\H}^{s,p}(\RR^n_+))}\\
    &\qquad+ \bigg\lVert \sum_{\ell=1}^N r_\ell\frac{|\lambda_\ell|}{\lambda_\ell}\sum_{\substack{\zeta\in\{0,\lambda_\ell\},\\\mathfrak{a}\in\{-\infty,0\}}} \Big(\pm\mathcal{M}_{\zeta,\mathfrak{a}}R_{\mathcal{N}}'R_{\mathcal{N}}^{\zeta}(\lambda_\ell)g_{\lambda_\ell} \mp \mathcal{M}_{\zeta,\mathfrak{a}}R_{\mathcal{N}}^{n}({\lambda_\ell})R_{\mathcal{N}}'g_{\lambda_\ell}\Big) \bigg\rVert_{\L^p(0,1;\dot{\H}^{s,p}(\RR^n_+))}\\
    &\leqslant \bigg\lVert \sum_{\ell=1}^{N} r_\ell R_\mathcal{N}'|\lambda_\ell|(|\lambda_\ell| \I-\Delta_\mathcal{N})^{-1}g_\lambda \bigg\rVert_{\L^p(0,1;\dot{\H}^{s,p}(\RR^n_+))}  \\
    &\qquad+ \sum_{\mathfrak{a}\in\{0,-\infty\}}\bigg(\bigg\lVert \sum_{\ell=1}^{N} r_\ell \mathcal{M}_{\lambda_\ell,\mathfrak{a}}R_{\mathcal{N}}'R_{\mathcal{N}}^{\lambda_\ell}(\lambda_\ell)g_{\lambda_\ell} \bigg\rVert_{\L^p(0,1;\dot{\H}^{s,p}(\RR^n_+))}\\
    &\qquad\qquad\qquad+\bigg\lVert \sum_{\ell=1}^{N} r_\ell \mathcal{M}_{\lambda_\ell,\mathfrak{a}}R_{\mathcal{N}}^{n}({\lambda_\ell})R_{\mathcal{N}}'g_{\lambda_\ell} \bigg\rVert_{\L^p(0,1;\dot{\H}^{s,p}(\RR^n_+))}\\
    &\qquad\qquad\qquad+\bigg\lVert \sum_{\ell=1}^{N} r_\ell \mathcal{M}_{0,\mathfrak{a}}R_{\mathcal{N}}'R_{\mathcal{N}}^{0}(\lambda_\ell)g_{\lambda_\ell} \bigg\rVert_{\L^p(0,1;\dot{\H}^{s,p}(\RR^n_+))}\\
    &\qquad\qquad\qquad+\bigg\lVert \sum_{\ell=1}^{N} r_\ell \mathcal{M}_{0,\mathfrak{a}}R_{\mathcal{N}}^{n}({\lambda_\ell})R_{\mathcal{N}}'g_{\lambda_\ell} \bigg\rVert_{\L^p(0,1;\dot{\H}^{s,p}(\RR^n_+))}\Bigg)\\
    &\lesssim_{p,s,n,\mu} \bigg\lVert \sum_{\ell=1}^{N} r_\ell g_{\lambda_\ell} \bigg\rVert_{\L^p(0,1;\dot{\H}^{s,p}(\RR^n_+))}\\
    &\lesssim_{p,s,n,\mu} \lVert (-\Delta_\mathcal{N})^{-\frac{1}{2}}\rVert_{\dot{\H}^{s-1,p}_0(\RR^n_+)\rightarrow \dot{\H}^{s,p}(\RR^n_+)}\bigg\lVert \sum_{\ell=1}^{N} r_\ell (|\lambda_\ell|\I-\Delta)g_\ell \bigg\rVert_{\L^p(0,1;\dot{\H}^{s-1,p}(\RR^n_+))}\\
    &\lesssim_{p,s,n,\mu} \bigg\lVert \sum_{\ell=1}^{N} r_\ell |\lambda_\ell|g_\ell \bigg\rVert_{\L^p(0,1;\dot{\H}^{s-1,p}(\RR^n_+))}+\bigg\lVert \sum_{\ell=1}^{N} r_\ell \nabla g_\ell \bigg\rVert_{\L^p(0,1;\dot{\H}^{s,p}(\RR^n_+))}.
\end{align*}
The implicit constant in the third inequality above is given by
\begin{align*}
    &\lVert R_\mathcal{N}'\rVert_{\dot{\H}^{s,p}(\RR^n_+)\rightarrow \dot{\H}^{s,p}(\RR^n_+)} \Bigg( \mathcal{R}\{ |\lambda|(|\lambda| \I-\Delta_\mathcal{N})^{-1},\,\lambda\in\Sigma_\mu\} \\
    &+ \sum_{\mathfrak{a}\in\{0,-\infty\}}\bigg( \mathcal{R}\{ \mathcal{M}_{\lambda,\mathfrak{a}},\,\lambda\in\Sigma_\mu\}\big(\mathcal{R}\{ R_{\mathcal{N}}^{\lambda}(\lambda),\,\lambda\in\Sigma_\mu\} +\mathcal{R}\{ R_{\mathcal{N}}^{n}(\lambda),\,\lambda\in\Sigma_\mu\}\big)\\ 
    &\qquad+\lVert \mathcal{M}_{0,\mathfrak{a}}\rVert_{\dot{\H}^{s,p}(\RR^n_+)\rightarrow\dot{\H}^{s,p}(\RR^n_+)}\big(\mathcal{R}\{ R_{\mathcal{N}}^{0}(\lambda),\,\lambda\in\Sigma_\mu\} +\mathcal{R}\{ R_{\mathcal{N}}^{n}(\lambda),\,\lambda\in\Sigma_\mu\}\big)\bigg)\Bigg)
\end{align*}
\textbf{Step 3.2:} We provide the argument for the $\mathcal{R}$-bound for $(\nabla^2\uu'_\ell)_{\ell\in\NN}$, again the proof for $(\nabla^2\uu_{n,\ell})_{\ell\in\NN}$ will be similar. Note that for any $\ell\in\NN$, that $\uu_\ell \in \dot{\D}^{s}_{p}(-\Delta_\mathcal{D})\subset\dot{\H}^{s,p}\cap\dot{\H}^{s+2,p}_\mathcal{D}(\RR^n_+)$, therefore, it holds
\begin{align*}
    \bigg\lVert \sum_{\ell=1}^{N} r_\ell\nabla^2\uu_\ell'\bigg\rVert_{\L^p(0,1;\dot{\H}^{s,p}(\RR^n_+))}\sim_{p,s,n} \bigg\lVert \sum_{\ell=1}^{N} r_\ell(-\Delta)\uu_\ell'\bigg\rVert_{\L^p(0,1;\dot{\H}^{s,p}(\RR^n_+))}.
\end{align*}
In this case, it can be checked that
\begin{align*}
     -\Delta \uu'_\ell(\cdot',x_n)= -[A_{\lambda_\ell}+A_0] A_{\lambda_\ell} e^{-x_n A_{\lambda_\ell}}[\nabla'(-\Delta_\mathcal{N})^{-1} g_\ell](\cdot',0) - [\nabla'(-\Delta_\mathcal{N})^{-1} g_\ell](\cdot',x_n).
\end{align*}
The whole proof for the $\mathcal{R}$-bound of $(-\Delta \uu'_\ell)_{\ell\in\NN}$ is then similar to, and even shorter than, the one of $( \uu'_\ell)_{\ell\in\NN}$, thus left to the reader.
\end{proof}

\newpage

\section{The Stokes--Dirichlet problem in bounded rough domains}\label{sec:stokessteady}

In this section, we consider the steady (resolvent) Stokes system
\begin{equation*}\tag{DS${_\lambda}$}\label{eq:Stokes}
    \left\{ \begin{array}{rllr}
         \lambda \uu - \Delta \uu +\nabla \mathfrak{p} &= \ff \text{, }&&\text{ in } \Omega\text{,}\\
        \div \, \uu &= 0\text{, } &&\text{ in } \Omega\text{,}\\
        \uu_{|_{\partial\Omega}} &=0\text{, } &&\text{ on } \partial\Omega\text{.}
    \end{array}
    \right.
\end{equation*}
in a domain $\Omega\subset\RR^n$ with unit normal $\nu$. The first main result given in following Theorem~\ref{thm:StokesResolvent} is a maximal elliptic regularity estimate for the solution in terms of the right-hand side and the boundary datum under nearly minimal assumption on the regularity of $\Omega$. This will allow us to derive the full picture for the regularity properties of the Stokes--Dirichlet operators and the underlying PDEs in the following sections. The corresponding multiplier spaces are introduced in Section \ref{sec:SM}.

The overall spirit is the following:
\begin{itemize}
    \item Provided a boundary of low regularity in the class $\mathcal{M}^{1+\alpha,r}_{\W}(\epsilon)$, $\epsilon>0$, the idea is to perform the analysis on a low regularity Sobolev (or Besov) space $\H^{s,p}$, $s>-1+{\sfrac{1}{p}}$, $s$ sufficiently close to $-1+{\sfrac{1}{p}}$, for which one can ensure a full gain of two derivatives for the solution of \eqref{eq:Stokes}. This is just the content of Section~\ref{Sec:ResolventEstGain2Deriv} with Theorem~\ref{thm:StokesResolvent}.
    \item Thanks to this exact gain of two derivatives on these low regularity Sobolev spaces, $\H^{s,p}$, $s>-1+{\sfrac{1}{p}}$, $s$ sufficiently close to $-1+{\sfrac{1}{p}}$, one will obtain in this case the following description for the domain of the Stokes--Dirichlet operator arising from \eqref{eq:Stokes}
    \begin{align}\label{eq:descritptionDomainStokesIntroSecStokes}
        \D_{p}^s(\AA_\mathcal{D}) =  \H^{s+2,p}_{\mathcal{D},\sigma}(\Omega) = \H^{s+1,p}_{0,\sigma}\cap\H^{s+2,p}(\Omega),\quad\text{ on } \H^{s,p}_{\mathfrak{n},\sigma}(\Omega).
    \end{align}
    \item With this exact description of the domain, one can easily (but tediously) extrapolate the boundedness of the $\mathbf{H}^{\infty}$-functional calculus from the scale $(\H^{s,2})_{|s|<\frac{1}{2}}$ to $\H^{s,p}$, $s>-1+{\sfrac{1}{p}}$, $p\in(1,\infty)$, $s$ sufficiently close to $-1+{\sfrac{1}{p}}$, resulting in particular in the BIP property. This particular step is achieved by applying the Kunstmann-Weis comparison principle following \cite{KunstmannWeis2017,GabelTolksdorf2022}. This is the first part of Theorem~\ref{thm:FinalResultSobolev1} in Section~\ref{sec:ProofLPbddDomain}.
    \item The BIP property yielding
    \begin{align*}
        [\H^{s,p}_{\mathfrak{n},\sigma}(\Omega),\D_{p}^s(\AA_\mathcal{D})]_{\alpha} = \D_{p}^s(\AA_\mathcal{D}^\alpha),\quad \text{ for all } \alpha\in(0,1),
    \end{align*}
    one obtains from Proposition~\ref{prop:InterpDivFreeC1} and \eqref{eq:descritptionDomainStokesIntroSecStokes}, the following identification, naturally accompanied by the underlying isomorphism (we recall that $s$ is negative),
    \begin{align*}
        \L^{p}_{\mathfrak{n},\sigma}(\Omega) = [\H^{s,p}_{\mathfrak{n},\sigma}(\Omega),\H^{s+2,p}_{\mathcal{D},\sigma}(\Omega)]_{\theta} = [\H^{s,p}_{\mathfrak{n},\sigma}(\Omega),\D_{p}^s(\AA_\mathcal{D})]_{\theta} =  \D_{p}^s(\AA_{\mathcal{D}}^\theta),\quad\text{ for some } \theta\in(0,1).
    \end{align*}
    Therefore, one obtains directly the extension of the (part of the) Dirichlet--Stokes operator in $\L^p_{\mathfrak{n},\sigma}(\Omega)$, and then exact solvability of \eqref{eq:Stokes} for $\ff\in\L^p_{\mathfrak{n},\sigma}(\Omega)$. In this case, our approach will lead to a precise description of the optimal gain of regularity for $\uu$ solving  \eqref{eq:Stokes} depending on the parameters $r\in[1,\infty]$ and $\alpha\in(0,1)$ that describe the regularity of the boundary. This is the second part of Theorem~\ref{thm:FinalResultSobolev1} in Section~\ref{sec:ProofLPbddDomain}.
    \item A similar but simpler analysis is performed on Besov spaces including endpoints $p=1,\infty$ in Theorem~\ref{thm:FinalResultBesov2} of Section~\ref{sec:BesovResolvResultsProof}.
    
    \item We will then precise our results and their consequences within the scale of $\L^p$-spaces $1\leqslant p\leqslant\infty$, for bounded $\C^{1,\alpha}$-domains, $\alpha\in(0,1]$, and provide a wide list of standard but sharpened consequences such has $\L^p$-$\L^q$ decay estimates for the semigroup, including $p,q=1,\infty$. This is mainly achieved in Theorems~\ref{thm:FinalResultSobolevC1q}~and~\ref{thm:FinalResultLinfty2} and Proposition~\ref{prop:Optimaldecay}. 

    \item Ultimately, in Section~\ref{Sec:RemoveC1}, Meta-Theorem~\ref{thm:metaThmremoveC1},  we will remove the extra-assumption of $\C^1$-boundary initially assumed in Theorem~\ref{thm:StokesResolvent} and its subsequent results.
\end{itemize}

\subsection{A first regularity result for the resolvent problem}\label{Sec:ResolventEstGain2Deriv}

Before presenting the next result, which we consider as an important step in this work,we state here a lemma that allows to deal with estimates for compactly supported perturbations of the divergence free condition in negative homogeneous Sobolev spaces, with a gain of one derivative.

\begin{lemma}\label{lem:preservingGengShen} Let $p,q\in[1,\infty]$, $s\in(-1+{\sfrac{1}{p}},{\sfrac{1}{p}})$. Let $\Omega$ be a bounded Lipschitz domain, and $\eta\in\C_c^{\infty}(\RR^n)$. Assume there exists a special Lipschitz domain ${\Omega}'$, such that $\Omega\subset\Omega'$, and one has, up to a rotation, the following description
\begin{align*}
    \Omega'=\{\,(x',x_n)\in\RR^n\,:\,\phi(x')>x_n\,\},
\end{align*}
for some $\phi\in\C^{0,1}_{ub}(\RR^{n-1})$. If $s\leqslant 0$ we assume additionally that for some $\alpha>0$, $r\in[1,\infty]$,
\begin{itemize}
    \item $\phi\in\mathcal{M}^{1+\alpha,r}_{\W}(\RR^{n-1})$.
\end{itemize}

\medbreak

We consider the pullback and the (inverse) pushforward maps $\mathbf{\Phi}^{\ast}$ and $\mathbf{\Phi}_{\ast}^{-1}$ induced by  $\phi$ according to \eqref{eq:Phi}.

\medbreak

Then, for all $\uu\in\B^{s+1,\sigma}_{p,q,\mathcal{D}}(\Omega)$, one has $\div [\tilde{\boldsymbol{\Phi}}^{-1}_\ast (\eta \uu)]\in\dot{\B}^{s-1}_{p,q,0}(\RR^n_+)$, with the estimate
\begin{align}\label{eq:PushforwardEstLocalisedDiv}
    \lVert\div [\tilde{\boldsymbol{\Phi}}^{-1}_\ast (\eta \uu)]\rVert_{\dot{\B}^{s-1}_{p,q}(\RR^n)} \lesssim_{p,s,n,\Omega} \lVert\nabla \eta\rVert_{\W^{1,\infty}(\RR^n)}\Big(1+\lVert \phi\rVert_{\mathcal{M}^{1+\alpha,r}_{\W}(\RR^{n-1})}\Big)^2 \lVert\uu\rVert_{{\B}^{s-1}_{p,q}(\Omega)}, 
\end{align}
and the result still holds replacing $\B^{\bullet}_{p,q}$ by $\H^{\bullet,p}$ whenever $1<p<\infty$.

Furthermore, if $s>0$ ($s\geqslant 0$ in the case of Bessel potential spaces), in \eqref{eq:PushforwardEstLocalisedDiv}, one can replace
\begin{align*}
    \lVert \phi\rVert_{\mathcal{M}^{1+\alpha,r}_{\W}(\RR^{n-1})}\qquad\text{ by }\qquad\lVert \nabla'\phi\rVert_{\L^\infty(\RR^{n-1})}.
\end{align*}

\end{lemma}

\begin{proof}One writes $\tilde{\Omega}:=\boldsymbol{\Phi}^{-1}(\Omega)$, and recall that one has $\RR^{n}_+=\boldsymbol{\Phi}^{-1}(\Omega')$. For $h\in\S_0(\RR^n,\CC)$, we have
\begin{align*}
    \langle \div [\tilde{\boldsymbol{\Phi}}^{-1}_\ast (\eta \uu)], h\rangle_{\RR^n} = \langle \div [\tilde{\boldsymbol{\Phi}}^{-1}_\ast (\eta \uu)],  h-(h)_{\tilde{\Omega}}\rangle_{\RR^n}.
\end{align*}
Since $\div \uu =0$, it holds that $\div [\tilde{\boldsymbol{\Phi}}^{-1}_\ast  \uu] =\tilde{\boldsymbol{\Phi}}^{-1}_\ast \div   \uu =0$, hence $\div [\tilde{\boldsymbol{\Phi}}^{-1}_\ast (\eta \uu)] = \det(\nabla {\boldsymbol{\Phi}}) [\nabla\eta \cdot\uu]\circ \boldsymbol{\Phi}$. We also mention that $\tilde{\boldsymbol{\Phi}}^{-1}_\ast (\eta \uu)$ has support contained in $\overline{\RR^{n}_+}$, so does its divergence.
One can write, thanks to a change of variable,
\begin{align*}
    \langle \div [\tilde{\boldsymbol{\Phi}}^{-1}_\ast (\eta \uu)], h\rangle_{\RR^n}  = \langle \uu, \nabla\eta [(h-(h)_{\tilde{\Omega}})\circ\boldsymbol{\Phi}^{-1}]\rangle_{\RR^n}.
\end{align*}
Finally by duality, since $\uu\in \B^{s+1,\sigma}_{p,q,\mathcal{D}}(\Omega)=\B^{s+1,\sigma}_{p,q,0}(\Omega)\hookrightarrow\B^{s-1}_{p,q,0}(\Omega,\CC^n)$, it holds
\begin{align*}
    |\langle \div [\tilde{\boldsymbol{\Phi}}^{-1}_\ast (\eta \uu)], h\rangle_{\RR^n}| &\leqslant\lVert \uu\rVert_{\B^{s-1}_{p,q}(\Omega)}\lVert \nabla\eta [(h-(h)_{\tilde{\Omega}})\circ\boldsymbol{\Phi}^{-1}]\rVert_{\B^{-s+1}_{p',q'}(\Omega)}\\
    &\lesssim_{p,s,n}\lVert \uu\rVert_{\B^{s-1}_{p,q}(\Omega)} \lVert \nabla \eta\rVert_{\W^{1,\infty}(\RR^n)}\lVert [(h-(h)_{\tilde{\Omega}})\circ\boldsymbol{\Phi}^{-1}]\rVert_{\B^{-s+1}_{p',q'}(\Omega)}.
\end{align*}
From this point, one just has to estimate the term $\lVert [(h-(h)_\Omega)\circ\boldsymbol{\Phi}^{-1}]\rVert_{\B^{-s+1}_{p',q'}(\Omega)}$. By Lemma~\ref{lem:multiplierComposition}, if $s >0$ (or $s\geqslant 0$ in the case of Bessel potential spaces)
\begin{align*}
    \lVert [(h-(h)_{\tilde{\Omega}})\circ\boldsymbol{\Phi}^{-1}]\rVert_{\B^{-s+1}_{p',q'}(\Omega)} &\lesssim_{p,s,n,\Omega} \lVert \det(\nabla\boldsymbol{\Phi})\rVert_{\L^\infty(\Omega')}^{\sfrac{1}{p}}\left[ 1+\lVert \nabla \boldsymbol{\Phi}^{-1}\rVert_{\L^\infty(\Omega')}\right] \lVert h-(h)_{\tilde{\Omega}}\rVert_{\B^{-s+1}_{p',q'}(\tilde{\Omega})}\\
    &\lesssim_{p,s,n,\Omega} 2^{1/p}\left[ 1+\lVert \nabla \boldsymbol{\Phi}\rVert_{\L^\infty(\RR^n_+)}\right] \lVert h-(h)_{\tilde{\Omega}}\rVert_{\B^{-s+1}_{p',q'}(\tilde{\Omega})}.\\
    &\lesssim_{p,s,n,\Omega} \left[ 1+\lVert \nabla'\phi\rVert_{\L^\infty(\RR^{n-1})}\right] \lVert h-(h)_{\tilde{\Omega}}\rVert_{\B^{-s+1}_{p',q'}(\tilde{\Omega})}.
\end{align*}
and if $s\leqslant 0$ ($s<0$)
\begin{align*}
    \lVert [(h-(h)_{\tilde{\Omega}})\circ\boldsymbol{\Phi}^{-1}]\rVert_{\B^{-s+1}_{p',q'}(\Omega)} &\lesssim_{p,s,n,\Omega} \lVert \det(\nabla\boldsymbol{\Phi})\rVert_{\L^\infty(\Omega')}^{\sfrac{1}{p}}\left[ 1+\lVert \nabla \boldsymbol{\Phi}^{-1}\rVert_{\L^\infty(\Omega')}\right]\\ &\qquad\qquad\qquad\qquad\qquad\qquad\times\left[ 1+\lVert \boldsymbol{\Phi}\rVert_{\mathcal{M}^{-s+1,p',q'}_{\B}(\RR^n_+)}\right]  \lVert h-(h)_{\tilde{\Omega}}\rVert_{\B^{-s+1}_{p',q'}(\tilde{\Omega})}\\
    &\lesssim_{p,s,n,\Omega} \left[1+\lVert \boldsymbol{\Phi}\rVert_{\mathcal{M}^{-s+1,p',q'}_{\B}(\RR^n_+)}\right]^2  \lVert h-(h)_{\tilde{\Omega}}\rVert_{\B^{-s+1}_{p',q'}(\tilde{\Omega})}.
\end{align*}
On the first hand, by Lemma~\ref{lem:NegativeNormThm}, Proposition~\ref{prop:FundamentalExtby0HomFuncSpaces} and the definition of function spaces by restriction, one deduces
\begin{align*}
    \lVert h-(h)_{\tilde{\Omega}}\rVert_{\B^{-s+1}_{p',q'}(\tilde{\Omega})} &\lesssim_{p,s,n,\phi,\Omega} \lVert \nabla h\rVert_{\B^{-s}_{p',q'}(\tilde{\Omega})}\\
    &\sim_{p,s,n,\phi,\Omega} \lVert \nabla h\rVert_{\dot{\B}^{-s}_{p',q'}(\tilde{\Omega})} \\
    &\lesssim_{p,s,n,\phi,\Omega} \lVert \nabla h\rVert_{\dot{\B}^{-s}_{p',q'}(\RR^n)} \\
    &\lesssim_{p,s,n,\phi,\Omega} \lVert h\rVert_{\dot{\B}^{-s+1}_{p',q'}(\RR^n)}. 
\end{align*}
On the other hand,
\begin{align*}
    \lVert \boldsymbol{\Phi}\rVert_{\mathcal{M}^{-s+1,p',q'}_{\B}(\RR^n_+)} = \lVert \nabla \mathcal{T}\phi\rVert_{\mathcal{M}^{-s,p',q'}_{\B}(\RR^n_+)}
\end{align*}
Thus, if $s< 0$, by Proposition~\ref{prop:MultipliersintheLplikeRange},  one obtains
\begin{align*}
    \lVert \boldsymbol{\Phi}\rVert_{\mathcal{M}^{-s+1,p',q'}_{\B}(\RR^n_+)} = \lVert \nabla \mathcal{T}\phi\rVert_{\mathcal{M}^{-s,p',q'}_{\B}(\RR^n_+)}  &\lesssim_{p,s,n,{\alpha}} \lVert  \phi\rVert_{\mathcal{M}^{1+{\alpha},r}_{\W}(\RR^{n-1})}.
\end{align*}

Summarizing up everything, and taking the supremum over $h\in\S_0(\RR^n,\CC)\subset\dot{\B}^{-s+1}_{p',q'}(\RR^n,\CC)$, we did obtain for all $s< 0$, 
\begin{align*}
    \lVert\div [\tilde{\boldsymbol{\Phi}}^{-1}_\ast (\eta \uu)]\rVert_{\dot{\B}^{s-1}_{p,q}(\RR^n,\CC)} &\lesssim_{p,s,n}^{\Omega}\lVert \nabla \eta\rVert_{\W^{1,\infty}(\RR^n)}\left[1+\lVert  \phi\rVert_{\mathcal{M}^{1+{\alpha},r}_{\W}(\RR^{n-1})}\right]^2 \lVert \uu\rVert_{\B^{s-1}_{p,q}(\Omega)} .
\end{align*}
If $s>0$, ($s\geqslant0$ in the case of Bessel potential spaces), one has similarly
\begin{align*}
    \lVert\div [\tilde{\boldsymbol{\Phi}}^{-1}_\ast (\eta \uu)]\rVert_{\dot{\B}^{s-1}_{p,q}(\RR^n,\CC)} &\lesssim_{p,s,n}^{\Omega,\zeta}\lVert \nabla \eta\rVert_{\W^{1,\infty}(\RR^n)}\left[1+\lVert  \phi\rVert_{\mathcal{M}^{1+{\alpha},r}_{\W}(\RR^{n-1})}\right] \lVert \uu\rVert_{\B^{s-1}_{p,q}(\Omega)} .
\end{align*}
It remains to deal with the case $s=0$ for Besov spaces $1<p<\infty$. Since $p$, $\eta$ and $\phi$ are fixed and independent of $s$, by real interpolation between $s=\pm\varepsilon$, one obtains
\begin{align*}
    \lVert\div [\tilde{\boldsymbol{\Phi}}^{-1}_\ast (\eta \uu)]\rVert_{\dot{\B}^{-1}_{p,q}(\RR^n,\CC)} &\lesssim_{p,s,n}^{\Omega}\lVert \nabla \eta\rVert_{\W^{1,\infty}(\RR^n)}\left[1+\lVert \phi\rVert_{\mathcal{M}^{1+\alpha,r}_{\W}(\RR^{n-1})}\right]^{2(1-\theta)}\\
    &\qquad\qquad\qquad\qquad\qquad\times\left[1+\lVert \nabla'\phi\rVert_{\L^\infty(\RR^{n-1})}\right]^{\theta} \lVert \uu\rVert_{\B^{-1}_{p,q}(\Omega)}\\
    &\lesssim_{p,s,n}^{\Omega}\lVert \nabla \eta\rVert_{\W^{1,\infty}(\RR^n)}\left[1+\lVert  \phi\rVert_{\mathcal{M}^{1+{\alpha},r}_{\W}(\RR^{n-1})}\right]^{2-\theta} \lVert \uu\rVert_{\B^{-1}_{p,q}(\Omega)}.
\end{align*}
This finishes the proof.
\end{proof}

\newpage
\begin{theorem}[Regularity of the Stokes resolvent problem]\label{thm:StokesResolvent}  Let $\alpha\in(0,1)$, $r\in[1,\infty]$, and let $\Omega$ be a bounded $\C^1$-domain of class $\mathcal{M}^{1+\alpha,r}_\W(\epsilon)$ for some sufficiently small $\epsilon>0$.

\medbreak

\noindent Let $p,q\in[1,\infty]$ $s\in(-1+{\sfrac{1}{p}},\sfrac{1}{p})$, such that either
\begin{enumerate}
    \item $p\in(r,\infty]$, $s\in(-1+\frac{1}{p},-1+\frac{1+r\alpha}{p})$; or
    \item $p\in[1,r]$, $s\in(-1+\frac{1}{p},-1+\alpha+\frac{1}{p})$; or
    \item $p=q=r$, $s=-1+\alpha+\frac{1}{p}$.
\end{enumerate}

Then, for $\mu\in[0,\pi)$, for all $\lambda\in\Sigma_{\mu}\cup\{0\}$, for all $f\in\B^{s,\sigma}_{p,q,\mathfrak{n}}(\Omega)$, the resolvent problem \eqref{eq:Stokes} admits a unique solution $(\uu,\mathfrak{p})\in\B^{s+2}_{p,q}(\Omega,\CC^n)\times\B^{s+1}_{p,q,\mfree}(\Omega)$ satisfying the estimate
\begin{align}\label{eq:thm:ResolvEstbbdDomain}
     (1+\lvert\lambda\rvert) \lVert \uu\rVert_{\B^{s}_{p,q}(\Omega)}+\lVert \nabla^2 \uu\rVert_{\B^{s}_{p,q}(\Omega)} + \lVert \nabla \mathfrak{p}\rVert_{\B^{s}_{p,q}(\Omega)} \nonumber \lesssim_{p,n,s,\mu}^{\Omega} \lVert \ff\rVert_{\B^{s}_{p,q}(\Omega)}.
\end{align}
Furthermore, the result still holds if one replaces $\B^{s}_{p,q}$ by either $\mathcal{B}^{s}_{p,\infty}$ or $\H^{s,p}$, assuming $1<p<\infty$ for the latter.
\end{theorem}

\begin{remark}In comparison to other results presented before, we did assume $\C^1$-boundary instead of Lipschitz boundary. This is in order to simplify the proof by using a result by Geng and Shen \cite[Theorem~5.1]{GengShen2024} especially the crucial estimate
\begin{align*}
(1+|\lambda|)^\frac{1}{2}\lVert \uu \rVert_{\L^\kappa(\Omega)}\lesssim_{\kappa,n,\Omega,\mu}\lVert \ff \rVert_{\W^{-1,\kappa}(\Omega)},\quad \kappa\in(1,\infty),
\end{align*}
for a bounded $\C^1$-domain $\Omega$. This inequality just being the $\L^\kappa$-counterparts of the one obtained in Lemma~\ref{lem:H-1EstBddLip}.

We make such an additional assumption for the convenience of the reader to draw the overall spirit of the localisation argument, then to focus on how to obtain the uniform estimates in $\lambda$ \eqref{eq:thm:ResolvEstbbdDomain}, as well as the other subsequent properties in next sections. The proof for Lipschitz domains is more intricate, but uses similar arguments. See Section~\ref{Sec:RemoveC1}, Theorem~\ref{thm:metaThmremoveC1}, where we draw a scheme of proof.
\end{remark}

\begin{proof} We deal with the case $r=1$. For the case $r=\infty$, one replace everywhere the quantity $\sfrac{(1+\alpha)}{p}$ by $\alpha+\sfrac{1}{p}$, and  the condition $s\in(-1+\sfrac{1}{p},-1+\sfrac{(1+\alpha)}{p})$ by $s\in(-1+\sfrac{1}{p},-1+\alpha+\sfrac{1}{p})$, allowing to perform the case $p=\infty$ with the exact same proof, since the interval is then non-empty. The case $r\in(1,\infty)$ is an in-between case which admits a similar proof, see Proposition~\ref{prop:MultipliersintheLplikeRange} and Remark~\ref{rem:BetterMultipliers}.

\medbreak

Note that one not only needs to show the estimates, we also need to show $\uu\in\B^{s+2}_{p,q}(\Omega)$ prior to perform them. For simplicity, here, we assume $1<q<\infty$, the endpoint cases being reachable through real interpolation. Therefore, in this setting, we exclusively deal with reflexive spaces, see Theorem~\ref{thm:Reflexiveendpointspaces}. We assume temporarily $p>1$ and $s< \min(-1+\frac{1+\alpha}{p},0)$, or that we are dealing with the case of $\L^p$ spaces with $s=0$. The case of other function spaces with $s\geqslant0$, and of Besov spaces with $p=1$ is deferred to Step 4 of the proof.

\medbreak

\textbf{Step 1: Preliminary setting.}  
By assumption, there is $\ell\in\mathbb N$ and functions $\varphi_1,\dots,\varphi_\ell\in\mathcal{M}^{1+\alpha,1}_{\W}(\epsilon)$ satisfying \ref{A1}--\ref{A3} for some sufficiently small $\epsilon>0$. By Proposition~\ref{prop:MultipliersintheLplikeRange}, we have $\nabla \mathcal{T}\varphi_j\in\mathcal{M}^{s,p}_{\X,\tt{or}}(\RR^{n-1})$ for all $-1+\sfrac{1}{p}<s<\sfrac{(1+\alpha)}{p}$, all $p\in[1,\infty]$, with an estimate
\begin{align}\label{eq:proofStokesDirbbdDomainSmallMultNorm}
    \lVert\nabla \mathcal{T}\varphi_j \rVert_{\mathcal{M}^{s,p}_{\X,\tt{or}}(\RR^{n}_+)}\lesssim_{p,s,n} \lVert\nabla' \varphi_j \rVert_{\mathcal{M}^{\alpha,1}_{\W,\tt{or}}(\RR^{n-1})} =  \lVert \varphi_j \rVert_{\mathcal{M}^{1+\alpha,1}_{\W}(\RR^{n-1})} < \epsilon.
\end{align}
That is $\nabla' \varphi_j$ has small multiplier norm and so does $\nabla \mathcal{T}\varphi_j$, for $j\in\llb 1,\ell\rrb$.

We consider an open set $\mathcal U^0\Subset{\Omega}$ and a family of open sets $(\mathcal{U}_j)_{j\in\llb 1,\ell\rrb}$ such that ${\Omega}\subset \cup_{j=0}^\ell \mathcal U^j$. Naturally, we consider a related partition of unity $(\eta_j)_{j\in\llb0,\ell\rrb}$ with respect to the covering
$\mathcal U^0,\dots,\mathcal U^\ell$ of ${\Omega}$. 
For $j\in\llb 1,\ell\rrb$, we consider the extension $\boldsymbol{\Phi}_j$ of $\varphi_j$ given by \eqref{eq:Phi} with inverse $\boldsymbol{\Phi}_j^{-1}$. Let $\ff\in\Ccinftydiv(\Omega)$ and consider the unique solution $\uu\in \cap_{r>1}\W^{1,r}_{0,\sigma}(\Omega)$ thanks to \cite[Theorem~5.1]{GengShen2024}.

\begin{equation*}\label{eq:StokesOr}
    \left\{ \begin{array}{rllr}
         \lambda\uu-\Delta \uu+\nabla\mathfrak{p} &= \ff \text{, }&&\text{ in } \Omega\text{,}\\
        \div\, \uu &= 0\text{, } &&\text{ in } \Omega\text{,}\\
        {\uu}{}_{|_{\partial\Omega}} &=0\text{, } &&\text{ on } \partial\Omega\text{.}
    \end{array}
    \right.
\end{equation*}

By Theorem~\ref{thm:regularizeddomain}, we consider a sequence of inner approximating domains $(\Omega^{\varepsilon})_{\varepsilon>0}$ associated with the regularised charts $((\varphi_{j}^{\varepsilon})_{j\in\llb1,\ell\rrb})_{\varepsilon>0}$ obtained by convolution according to Theorem~\ref{thm:regularizeddomain}, and the associated (lifted) diffeomorphisms $((\boldsymbol{\Phi}_{j}^{\varepsilon})_{j\in\llb1,\ell\rrb})_{\varepsilon>0}$. Let $\varepsilon_0>0$ such that for all $0<\varepsilon<\varepsilon_0$, it holds $\supp \ff \subset \Omega_{\varepsilon}$, and consider the solution $(\uu^{\varepsilon},\mathfrak{p}^{\varepsilon})$  to
\begin{equation}\label{eq:StokesOrregularized}
    \left\{ \begin{array}{rllr}
         \lambda\uu^{\varepsilon}-\Delta \uu^{\varepsilon}+\nabla\mathfrak{p}^{\varepsilon} &= \ff \text{, }&&\text{ in } \Omega_\varepsilon\text{,}\\
        \div\, \uu^{\varepsilon} &= 0\text{, } &&\text{ in } \Omega_\varepsilon\text{,}\\
        {\uu^{\varepsilon}}{}_{|_{\partial\Omega_\varepsilon}} &=0\text{, } &&\text{ on } \partial\Omega_\varepsilon\text{.}
    \end{array}
    \right.
\end{equation}

By \cite[Theorem~5.1]{GengShen2024}, and since the characteristic constants of the domain are preserved through the approximation procedure for the boundary by Theorem~\ref{thm:regularizeddomain}, up to consider the extension of $\uu^{\varepsilon}$ by $0$ to $\Omega$, one obtains $\uu^{\varepsilon}\in\W^{1,r}_{0,\sigma}(\Omega)$ for all $r\in(1,\infty)$, with
\begin{align*}
    \lVert \uu^{\varepsilon} \rVert_{\W^{1,r}(\Omega)} \lesssim_{r,n,\Omega} \lVert \ff \rVert_{\L^{r}(\Omega)},
\end{align*}
uniformly with respect to $\varepsilon>0$. By Banach-Alaoglu, up to extract finitely many subsequences, $\uu^{\varepsilon}$ converges weakly to some limit $\vv\in\W^{1,r}_{0,\sigma}(\Omega)$ (for any $r>1$) and by compact embedding, strongly in $\B^{s+1,\sigma}_{p,q,0}(\Omega),\H^{s+1,p}_{0,\sigma}(\Omega)$, recall that $s+1<1$\footnote{this is the reason of the temporary restrictions $s<0$ and $p>1$: to be able to obtain convergence, --or, at least, boundedness-- in function spaces with at least one more derivative than the ground space.}. For $\boldsymbol{\varphi}\in\Ccinftydiv(\Omega)$ can check that if $\Omega_\varepsilon$ contains the support of $\boldsymbol{\varphi}$,
\begin{align*}
    \lambda\int_{\Omega} \uu^\varepsilon\cdot\boldsymbol{\varphi} + &\int_{\Omega} \nabla \uu^\varepsilon \cdot \nabla \boldsymbol{\varphi} =\int_{\Omega} \ff\cdot\boldsymbol{\varphi}\\
    &\xrightarrow[\varepsilon\longrightarrow 0]{} \lambda\int_{\Omega} \vv\cdot\boldsymbol{\varphi} + \int_{\Omega} \nabla \vv \cdot \nabla \boldsymbol{\varphi} =\int_{\Omega} \ff\cdot\boldsymbol{\varphi}.
\end{align*}
The limit $\vv\in\W^{1,2}_{0,\sigma}(\Omega)$ satisfying the weak Stokes--Dirichlet equation for all $\boldsymbol{\varphi}\in\Ccinftydiv(\Omega)$, by uniqueness in $\H^{1,2}_{0,\sigma}(\Omega)$, one obtains $\vv=\uu$. And similarly $\nabla\mathfrak{p_\varepsilon}$ converges towards $\nabla\mathfrak{p}$ in $\B^{s-1}_{p,q}(\Omega)$.

Note that since $\ff$ is smooth and $\Omega_\varepsilon$ has smooth boundary, by construction it holds that $\uu^{\varepsilon}$ and $\mathfrak{p}^\varepsilon$ are smooth in $\overline{\Omega}_{\varepsilon}$, up to the boundary.
 
 Let us fix $j\in\llb1,\ell\rrb$ and assume, without loss of generality, that the reference point $y_j=0$ and that the outer normal at~$0$ is pointing in the negative $x_n$-direction (this saves us some notation regarding the translation and rotation of the coordinate system).
We multiply $\uu^\varepsilon$ by $\eta_j$ and obtain for $\uu_j^{\varepsilon}:=\eta_j\uu^{\varepsilon}$, $\mathfrak{p}_j:=\eta_j\mathfrak{p}^{\varepsilon}$ and $\ff_j:=\eta_j\ff$ the equations
\begin{equation}\label{eq:Stokes2}
    \left\{ \begin{array}{rllr}
         \lambda\uu_j^{\varepsilon}-\Delta \uu_j^{\varepsilon}+\nabla\mathfrak{p}_j^{\varepsilon} &= [-\Delta,\eta_j]\uu^{\varepsilon}+[\nabla,\eta_j]\mathfrak{p}^{\varepsilon}+\ff_j \text{, }&&\text{ in } \Omega_{\varepsilon}\text{,}\\
        \div\, \uu_j^{\varepsilon} &= \div(\eta_j \uu^{\varepsilon})\text{, } &&\text{ in } \Omega_{\varepsilon}\text{,}\\
        {\uu_j^{\varepsilon}}{}_{|_{\partial\Omega_{\varepsilon}}} &=0\text{, } &&\text{ on } \partial\Omega_{\varepsilon}\text{.}
    \end{array}
    \right.
\end{equation}
with the commutators $[\Delta,\eta_j]=\Delta\eta_j+2\nabla\eta_j\cdot\nabla$ and $[\nabla,\eta_j]=\nabla\eta_j$.
Now, we set $\tilde{\uu}_j^{\varepsilon}:=\uu_j^{\varepsilon}\circ\boldsymbol{\Phi}_j^{\varepsilon}$, $\tilde{\mathfrak{p}}_j^{\varepsilon}:=\mathfrak{p}_j^{\varepsilon}\circ\boldsymbol{\Phi}_j^{\varepsilon}$, $\mathbf{g}_j:=\det (\nabla\boldsymbol{\Phi}_j^{\varepsilon})([-\Delta,\eta_j]\uu+[\nabla,\eta_j]\mathfrak{p}^{\varepsilon}+\ff_j)\circ\boldsymbol{\Phi}_j^{\varepsilon}$, and
\begin{align*}
    h_j^\varepsilon:=\det(\nabla\boldsymbol{\Phi}_j^{\varepsilon})(\div(\eta_j\uu^{\varepsilon}))\circ\boldsymbol{\Phi}_j^{\varepsilon} = (\tilde{\boldsymbol{\Phi}}_{j,\ast}^{\varepsilon})^{-1}\div(\eta_j\uu^{\varepsilon}) = \div[ (\tilde{\boldsymbol{\Phi}}_{j,\ast}^{\varepsilon})^{-1}(\eta_j\uu^{\varepsilon})].
\end{align*}
We obtain the equations
\begin{equation}\label{eq:Stokes3}
    \left\{ \begin{array}{rllr}
        \lambda\det(\nabla\boldsymbol{\Phi}_j^{\varepsilon})\tilde{\uu}_j^{\varepsilon}-\div \big(\mathbf{A}_j^{\varepsilon}\nabla\tilde{\uu}_j)+\div(\mathbf{B}_j^\varepsilon\tilde{\mathfrak{p}}_j^\varepsilon) &= \gg_j^\varepsilon \text{, }&&\text{ in } \RR^n_+\text{,}\\
        \div\, \tilde{\uu}_j^{\varepsilon} + \div([{\mathbf{B}_j^{\varepsilon}}-\mathbf{I}_{n}] \tilde{\uu}_j^{\varepsilon})  &= h_j^{\varepsilon}\text{, } &&\text{ in } \RR^n_+\text{,}\\
        {\tilde{\uu}_j^{\varepsilon}}{}_{|_{\partial\RR^n_+}} &=0\text{, } &&\text{ on } \partial\RR^n_+\text{.}
    \end{array}
    \right.
\end{equation}
where $\mathbf{A}_j:=\tilde{\boldsymbol{\Phi}}_{\ast j}^{-1} (\boldsymbol{\Phi}_{j}^{\ast})^{-1}=\det(\nabla\boldsymbol{\Phi}_j)(\nabla\boldsymbol{\Phi}_j)^{-1}(\prescript{t}{}{\nabla\boldsymbol{\Phi}_j})$ and $\mathbf{B}_j:=\det(\nabla\boldsymbol{\Phi}_j)(\nabla\boldsymbol{\Phi}_j)^{-1}$
 (note that, here, $\tilde{\boldsymbol{\Phi}}_{\ast j}^{-1}$ and $(\boldsymbol{\Phi}_{j}^{\ast})^{-1}$ are applied component-wise identifying $\nabla\mathfrak{p}=\div(\mathfrak{p}\mathbf{I}_{n})$\footnote{\textit{i.e.} considering each scalar equation $\lambda {v}_k - \div(\nabla v_k - \mathfrak{q} \mathfrak{e}_k)=\ldots$, $k\in\llb1,n\rrb$. However, for the divergence constraint consider the (inverse) pushforward $\tilde{\boldsymbol{\Phi}}^{-1}_\ast$ acting as such on $\div(\eta\vv)=...$.}).
 This can be rewritten as
\begin{equation}\label{eq:Stokes5}
    \left\{ \begin{array}{rllr}
        \lambda \tilde{\uu}_j^{\varepsilon}-\Delta\tilde{\uu}_j^{\varepsilon} + \nabla \tilde{\mathfrak{p}}_j^{\varepsilon} -\div \big([\mathbf{A}_j^{\varepsilon}-\mathbf{I}_{n}]\nabla\tilde{\uu}_j^{\varepsilon})+\div([{\mathbf{B}_j^{\varepsilon} - \mathbf{I}_{n}}]\tilde{\mathfrak{p}}_j^{\varepsilon}) &= \gg_j^{\varepsilon} + \lambda[1-\det(\nabla\boldsymbol{\Phi}_j)]\tilde{\uu}_j^{\varepsilon} \text{,}\\
        \div\, \tilde{\uu}_j^{\varepsilon} + \div([{\mathbf{B}_j^{\varepsilon}}-\mathbf{I}_{n}] \tilde{\uu}_j^{\varepsilon})  &= h_j^{\varepsilon}\text{,}\\
        {\tilde{\uu}_j^{\varepsilon}}{}_{|_{\partial\RR^n_+}} &=0\text{.}
    \end{array}
    \right.
\end{equation}
Setting
\begin{align*}
\mathcal{S}^{j} (\tilde{\uu},\tilde{\mathfrak{p}})&=\mathcal{S}^{j} _0(\tilde{\uu})+\mathcal{S}^{j}_1 (\tilde{\uu})+ \mathcal{S}^{j} _2(\tilde{\mathfrak{p}}),\\
\mathcal{S}^{j}_0 (\tilde{\uu})&=\lambda (1-\det(\nabla\boldsymbol{\Phi}_j)) \tilde{\uu},\\
\mathcal{S}^{j}_1  (\tilde{\uu})&=-\div \big((\mathbf{I}_{n}-\mathbf{A}_j)\nabla\tilde{\uu}),\\
\mathcal{S}^{j}_2 (\tilde{\mathfrak{p}})&=\div ((\mathbf{B}_j-\mathbf{I}_n)\tilde{\mathfrak{p}}),\\
\mathfrak{s}^{j}(\tilde{\uu})&=\div([\mathbf{I}_{n}-{}{\mathbf{B}_j}]\tilde{\uu})=[\mathbf{I}_{n}-{}{\mathbf{B}_j}]:\nabla \tilde{\uu},\\
\end{align*} we can finally write
\eqref{eq:Stokes5} as
\begin{align}\label{eq:Stokes6}
\lambda\tilde{\uu}_j^{\varepsilon}-\Delta\tilde{\uu}_j^{\varepsilon}+\nabla\tilde{\mathfrak{p}}_j^{\varepsilon}&=\mathcal S^{j,\varepsilon} (\tilde{\uu}_j^{\varepsilon},\tilde{\mathfrak{p}}_j^{\varepsilon})+\mathbf{g}_j^{\varepsilon},
\quad\div \tilde{\uu}_j^{\varepsilon}=\mathfrak s^{j,\varepsilon}(\tilde{\uu}_j^{\varepsilon})+h_j^{\varepsilon},\quad
{\tilde{\uu}_j^{\varepsilon}}{}_{|_{\partial\RR^n_+}}=0.
\end{align}
The estimate for the half-space (see Theorem \ref{thm:MetaThmDirichletStokesRn+}) yields
\begin{align}
|\lambda|\|\tilde{\uu}_j^{\varepsilon}\|_{\dot{\B}^{s}_{p,q}(\RR^n_+)}+\| \nabla^2 \tilde{\uu}_j^{\varepsilon}\|_{\dot{\B}^{s}_{p,q}(\RR^n_+)}+\|\tilde \nabla {\mathfrak{p}}_j^{\varepsilon}\|_{\dot{\B}^{s}_{p,q}(\RR^n_+)}&\lesssim_{p,s,n,\mu} \|\mathcal S^{j,\varepsilon} (\tilde{\uu}_j^{\varepsilon},\tilde{\mathfrak{p}}_j^{\varepsilon})+\mathbf{g}_j^{\varepsilon}\|_{\dot{\B}^{s}_{p,q}(\RR^n_+)}\label{eq:0201c}\\ 
 &\qquad\qquad+ |\lambda   |\|\mathfrak{s}^{j,\varepsilon}(\tilde{\uu}_j^{\varepsilon})+h_j^{\varepsilon}\|_{\dot{\B}^{s-1}_{p,q}(\RR^n_+)}\nonumber\\ &\qquad\qquad\qquad +\|\mathfrak s^{j,\varepsilon}(\tilde{\uu}_j^{\varepsilon})+h_j^{\varepsilon}\|_{\dot{\B}^{s+1}_{p,q}(\RR^n_+)}\nonumber.
\end{align}

\textbf{Step 2: Uniform estimates and proof of regularity.} To prove ${\B}^{s+2}_{p,q}$--regularity, we need to show the boundedness of the ${\B}^{s+2}_{p,q}$--norm of $\uu^{\varepsilon}$ by its lower order ${\B}^{s+1}_{p,q}$, and ${\B}^{s}_{p,q}$-norms uniformly with respect to $\varepsilon$. Therefore, our task consists in estimating the right-hand side uniformly with respect to $\varepsilon$.

In order to estimate $\|\mathcal S^{j,\varepsilon}(\tilde{\uu}^\varepsilon_j,\tilde{\mathfrak{p}}^\varepsilon_j)\|_{\dot{\B}^{s+1}_{p,q}(\RR^n_+)}$,  $|\lambda|\|\mathfrak s^{j,\varepsilon}(\tilde{\uu}^\varepsilon)\|_{\dot{\B}^{s-1}_{p,q}(\RR^n_+)}$ and $\|\mathfrak s^{j,\varepsilon}(\tilde{\uu}^\varepsilon)\|_{\dot{\B}^{s+1}_{p,q}(\RR^n_+)}$ we
use the Besov multiplier norm introduced in \eqref{eq:SoMo}.  Hence we obtain
by \eqref{J}, and the definitions of $\mathbf{A}_j$ and $\boldsymbol{\Phi}_j$ (by compact support in $\mathcal{U}_j$, one can apply Lemma~\ref{lem:CompactEquiHomInhom} and Proposition~\ref{prop:FundamentalExtby0HomFuncSpaces}) it holds
\begin{align*}
&\|\mathcal S_1^{j,\varepsilon}(\tilde{\uu}^{\varepsilon}_j)\|_{\dot{\B}^{s}_{p,q}(\RR^n_+)}\lesssim_{p,s,n,\Omega} \|\mathcal S_1^{j,\varepsilon}(\tilde{\uu}^{\varepsilon}_j)\|_{{\B}^{s}_{p,q}(\RR^n_+)}\\
&\quad\lesssim_{p,s,n,\Omega}  \|\mathbf I_{n}-\mathbf{A}_j^\varepsilon\|_{\mathcal{M}_{\B,\tt or}^{s+1,p,q}(\RR^n_+)}\|\nabla\tilde{\uu}^\varepsilon_j\|_{\B^{s+1}_{p,q}(\RR^n_+)}\\
&\quad\lesssim_{p,s,n,\Omega}  \|\mathbf I_{n}-\mathbf{A}_j^\varepsilon\|_{\mathcal{M}_{\B,\tt or}^{s+1,p,q}(\RR^n_+)}\|{\uu}^\varepsilon_j\|_{\B^{s+2}_{p,q}(\Omega_\varepsilon)}.
\end{align*}
As in the proof of Proposition~\ref{prop:HodgeDiracbdd}, Step 2.1, writing 
\begin{align*}
    \I- \mathbf{A}_j&= \I-\det(\nabla \boldsymbol{\Phi}_j)(\nabla  \boldsymbol{\Phi}_j)^{-1}\cdot\prescript{t}{}{(\nabla  \boldsymbol{\Phi}_j)}\\
    &= \I - \det(\nabla \boldsymbol{\Phi}_j)(\nabla  \boldsymbol{\Phi}_j)^{-1} +\det(\nabla \boldsymbol{\Phi}_j)(\nabla  \boldsymbol{\Phi}_j)^{-1}[\I-\prescript{t}{}{(\nabla  \boldsymbol{\Phi}_j)}]\\
    &= [\I-(\nabla  \boldsymbol{\Phi}_j)^{-1}] +[1- \det(\nabla \boldsymbol{\Phi}_j)](\nabla  \boldsymbol{\Phi}_j)^{-1} +\det(\nabla \boldsymbol{\Phi}_j)(\nabla  \boldsymbol{\Phi}_j)^{-1}[\I-\prescript{t}{}{(\nabla  \boldsymbol{\Phi}_j)}],
\end{align*}
so that for $\epsilon>0$ small enough (not to be mistaken with $\varepsilon$ ! see the beginning of the proof for notations)  since $s+1<\frac{1+\alpha}{p}$, thanks to \eqref{eq:proofStokesDirbbdDomainSmallMultNorm}, it holds
\begin{align*}
    \lVert \I- \mathbf{A}_j^\varepsilon \rVert_{\mathcal{M}^{s+1,p,q}_{\B,\tt{or}}(\RR^n_+)} &\leqslant \big\lVert \I-(\nabla  \boldsymbol{\Phi}_j^\varepsilon)^{-1}\big\rVert_{\mathcal{M}^{s+1,p,q}_{\B,\tt{or}}(\RR^n_+)}+\big\lVert (1 - \det(\nabla \boldsymbol{\Phi}_j^\varepsilon))(\nabla  \boldsymbol{\Phi}_j^\varepsilon)^{-1} \big\rVert_{\mathcal{M}^{s+1,p,q}_{\B,\tt{or}}(\RR^n_+)} \\ &\qquad+ \big\lVert\det(\nabla \boldsymbol{\Phi}_j^\varepsilon)(\nabla  \boldsymbol{\Phi}_j^\varepsilon)^{-1} \big\rVert_{\mathcal{M}^{s+1,p,q}_{\B,\tt{or}}(\RR^n_+)}\big\lVert\I-\prescript{t}{}{(\nabla  \boldsymbol{\Phi}_j^\varepsilon)} \big\rVert_{\mathcal{M}^{s+1,p,q}_{\B,\tt{or}}(\RR^n_+)}\\
    &\lesssim_{n}  \lVert \nabla\mathcal{T}\varphi_j^\varepsilon \rVert_{\mathcal{M}^{s+1,p,q}_{\B,\tt{or}}(\RR^n_+)} + \lVert\partial_{x_n}\mathcal{T}\varphi_j^\varepsilon\rVert_{\mathcal{M}^{s+1,p,q}_{\B,\tt{or}}(\RR^n_+)}\left( 1+\lVert\nabla\mathcal{T}\varphi_j^\varepsilon \rVert_{\mathcal{M}^{s+1,p,q}_{\B,\tt{or}}(\RR^n_+)}\right)\\
    &\qquad+ \big(1+\big\lVert  \nabla\mathcal{T}\varphi_j^\varepsilon\big\rVert_{\mathcal{M}^{s+1,p,q}_{\B,\tt{or}}(\RR^n_+)}\big) \big\lVert\nabla \mathcal{T}\varphi_j^\varepsilon \big\rVert_{\mathcal{M}^{s+1,p,q}_{\B,\tt{or}}(\RR^n_+)}\\
    &\lesssim_{n,\Omega}  \lVert  \varphi_j^\varepsilon\rVert_{\mathcal{M}^{1+\alpha,1}_{\W}(\RR^{n-1})} + \lVert\varphi_j^\varepsilon\rVert_{\mathcal{M}^{1+\alpha,1}_{\W}(\RR^{n-1})}\left( 1+\lVert\varphi_j^\varepsilon \rVert_{\mathcal{M}^{1+\alpha,1}_{\W}(\RR^{n-1})}\right)\\
    &\qquad+ \big(1+\lVert  \varphi_j^\varepsilon\rVert_{\mathcal{M}^{1+\alpha,1}_{\W}(\RR^{n-1})}\big) \lVert\varphi_j^\varepsilon \rVert_{\mathcal{M}^{1+\alpha,1}_{\W}(\RR^{n-1})}\\
    &\lesssim_{n,\Omega}  \lVert  \varphi_j\rVert_{\mathcal{M}^{1+\alpha,1}_{\W}(\RR^{n-1})} + \lVert\varphi_j\rVert_{\mathcal{M}^{1+\alpha,1}_{\W}(\RR^{n-1})}\left( 1+\lVert\varphi_j\rVert_{\mathcal{M}^{1+\alpha,1}_{\W}(\RR^{n-1})}\right)\\
    &\qquad+ \big(1+\lVert  \varphi_j\rVert_{\mathcal{M}^{1+\alpha,1}_{\W}(\RR^{n-1})}\big) \lVert\varphi_j \rVert_{\mathcal{M}^{1+\alpha,1}_{\W}(\RR^{n-1})}\\
    &\lesssim_{\Omega} \epsilon + 2\epsilon(1+\epsilon).
\end{align*}
Here, we used that $\varphi_j^{\varepsilon}$ is a mollified version of $\varphi_j$ according to Theorem~\ref{thm:regularizeddomain} and Young's inequality for convolutions for multipliers \eqref{eq:mollifMultip}, so that the bound is uniform with respect to $\varepsilon$. So, we finally have
\begin{align*}
\|\mathcal S_1^{j,\varepsilon}(\tilde{\uu}_j^{\varepsilon})\|_{\B^{s-2}_{p,q}(\RR^n_+)}&\lesssim_{p,s,n,\Omega} \epsilon \|{\uu}^{\varepsilon}_j\|_{\B^{s+2}_{p,q}(\Omega_\varepsilon)}
\end{align*}
and, similarly, 
\begin{align*}
\|\mathcal S_0^{j,\varepsilon}(\tilde{\uu}_j^\varepsilon)\|_{\B^{s}_{p,q}(\RR^n_+)}&\lesssim_{p,s,n,\Omega}  |\lambda|\|1-\det(\nabla\boldsymbol{\Phi}_j^\varepsilon)\|_{\mathcal{M}^{s,p,q}_{\B,\tt or}(\RR^n_+)}\|{\uu}^\varepsilon_j\|_{\B^{s}_{p,q}(\Omega_\varepsilon)}\\
&\lesssim_{p,s,n,\Omega}    |\lambda|\|\partial_{x_n}\mathcal{T}\varphi^{\varepsilon}_j\|_{\mathcal{M}^{s,p,q}_{\B,\tt or}(\RR^n_+)}\|{\uu}^\varepsilon_j\|_{\B^{s}_{p,q}(\Omega_\varepsilon)}\\
&\lesssim_{p,s,n,\Omega}   |\lambda|\|\varphi^{\varepsilon}_j\|_{\mathcal{M}^{1+\alpha,1}_{\W}(\RR^{n-1})}\|{\uu}^\varepsilon_j\|_{\B^{s}_{p,q}(\Omega_\varepsilon)}\\
&\lesssim_{p,s,n,\Omega}   \epsilon|\lambda|\|{\uu}^\varepsilon_j\|_{\B^{s}_{p,q}(\Omega_\varepsilon)}.
\end{align*}
Now, since $s+1<\frac{1+\alpha}{p}$, one obtains in the exact same way
\begin{align*}
\|\mathcal S_2^{j,\varepsilon}(\tilde{\mathfrak{p}}^{\varepsilon})_j\|_{\dot{\B}^{s}_{p,q}(\RR^n_+)}&\lesssim_{p,s,n,\Omega}  \|\mathbf{B}_j^\varepsilon-\mathbf{I}_{n}\|_{\mathcal{M}_{\B,\tt or}^{s+1,p,q}(\RR^n_+)}\|\tilde{\mathfrak{p}}_{j}^\varepsilon\|_{\B^{s+1}_{p,q}(\RR^n_+)}\\
&\lesssim_{p,s,n,\Omega} \lVert \partial_{x_n}\mathcal{T}\varphi_j^\varepsilon \rVert_{\mathcal{M}^{s+1,p,q}_{\B,\tt{or}}(\RR^n_+)}\|{\mathfrak{p}}^\varepsilon\|_{\B^{s+1}_{p,q}(\Omega_{\varepsilon})}\\
&\lesssim_{p,s,n,\Omega} \epsilon\|{\mathfrak{p}}^\varepsilon\|_{\B^{s+1}_{p,q}(\Omega_{\varepsilon})},
\end{align*}
as well as, 
\begin{align*}
\|\mathfrak s^{j,\varepsilon}(\tilde{\uu}^\varepsilon_j)\|_{\dot{\B}^{s+1}_{p,q}(\RR^n_+)}
&\lesssim_{p,s,n,\Omega}\|\mathbf I_{n}-\mathbf{B}_j^\varepsilon\|_{\mathcal M^{s+1,p,q}_{\B,\tt or}(\RR^n_+)}\|\nabla\tilde{\uu}_{j}^\varepsilon\|_{\B^{s+1}_{p,q}(\RR^n_+)}\\
&\lesssim_{p,s,n,\Omega} \|\nabla\mathcal{T}\varphi_j^{\varepsilon}\|_{\mathcal M^{s+1,p,q}_{\B,\tt or}(\RR^{n}_+)}\|{\uu}^{\varepsilon}_j\|_{\B^{s+2}_{p,q}(\Omega_\varepsilon)}\\
&\lesssim_{p,s,n,\Omega} \epsilon\|{\uu}^{\varepsilon}_j\|_{\B^{s+2}_{p,q}(\Omega_\varepsilon)}.
\end{align*}
thanks to the equality $s^{j,\varepsilon}(\tilde{\uu}^\varepsilon_j)= [\mathbf{I}_n-\mathbf{B}_j^\varepsilon]:\nabla \tilde{\uu}_{j}^\varepsilon$. Now, taking advantage of the equality $s^{j,\varepsilon}(\tilde{\uu}^\varepsilon_j)=\div([\mathbf{I}_n-\mathbf{B}_j^\varepsilon]\tilde{\uu}_{j}^\varepsilon)$, we can estimate
\begin{align*}
|\lambda|\|\mathfrak s^{j,\varepsilon}(\tilde{\uu}^\varepsilon_j)\|_{\dot{\B}^{s-1}_{p,q}(\RR^n_+)}
&\lesssim_{p,s,n,\Omega} |\lambda|\|{\mathbf{B}}_j^{\varepsilon}-\mathbf I_{n}\|_{\mathcal M^{s,p,q}_{\B,\tt or}(\RR^n_+)}\|\tilde{\uu}_{j}^\varepsilon\|_{\B^{s}_{p,q}(\RR^n_+)}\\
&\lesssim_{p,s,n,\Omega} |\lambda|\lVert \partial_{x_n}\mathcal{T}\varphi_j^\varepsilon \rVert_{\mathcal{M}^{s,p,q}_{\B,\tt{or}}(\RR^n_+)}\|{\uu}^{\varepsilon}_j\|_{\B^{s+2}_{p,q}(\Omega_\varepsilon)}\\
&\lesssim_{p,s,n,\Omega} \epsilon|\lambda|\|{\uu}^{\varepsilon}_j\|_{\B^{s}_{p,q}(\Omega_\varepsilon)}.
\end{align*}
We conclude that
 \begin{align}\label{eq:0301}
 \|\mathcal S^{j,\varepsilon}(\tilde{\uu}_j^\varepsilon,\tilde{\mathfrak{p}}_j^\varepsilon)\|_{\dot{\B}^{s}_{p,q}(\RR^n_+)}+&|\lambda|\|\mathfrak s^{j,\varepsilon}(\tilde{\uu}_j^{\varepsilon})\|_{\dot{\B}^{s-1}_{p,q}(\RR^n_+)} + \|\mathfrak s^{j,\varepsilon}(\tilde{\uu}_j^{\varepsilon})\|_{\dot{\B}^{s+1}_{p,q}(\RR^n_+)}\\
 &\leqslant C_{p,s,n,\Omega} \cdot \epsilon \big(|\lambda|\|{\uu}^{\varepsilon}_j\|_{\B^{s}_{p,q}(\Omega_\varepsilon)}+\|{\uu}^\varepsilon\|_{\B^{s+2}_{p,q}(\Omega_{\varepsilon})}+\|{\mathfrak{p}}^\varepsilon_j\|_{\B^{s+1}_{p,q}(\Omega_{\varepsilon})}\big)\nonumber
 \end{align}
 for sufficiently small $\epsilon>0$. On the other hand, we have
 \begin{align*}
 \|\mathbf{g}_j^\varepsilon\|_{\dot{\B}^{s}_{p,q}(\RR^n_+)}&\lesssim_{p,s,n,\Omega} \| [\Delta \eta_j \uu^\varepsilon]\circ\boldsymbol{\Phi}_{j}^\varepsilon\|_{\dot{\B}^{s}_{p,q}(\RR^n_+)}+ \| [\nabla \eta_j \nabla\uu^\varepsilon]\circ\boldsymbol{\Phi}_{j}^\varepsilon\|_{\dot{\B}^{s}_{p,q}(\RR^n_+)}\\ & \qquad\qquad+ \| [\nabla{\eta_j} \mathfrak{p}^{\varepsilon}]\circ\boldsymbol{\Phi}_j^\varepsilon\|_{\dot{\B}^{s}_{p,q}(\RR^n_+)}+ \|\ff_j\circ\boldsymbol{\Phi}_j^{\varepsilon}\|_{\B^{s}_{p,q}(\RR^n_+)}\\
& \lesssim_{p,s,n,\Omega} \|\uu^{\varepsilon}\|_{\B^{s}_{p,q}(\Omega_\varepsilon)}+ \|\mathfrak{p}^{\varepsilon}\|_{\B^{s+1}_{p,q}(\Omega_\varepsilon)}+\|\ff\|_{\B^{s}_{p,q}(\Omega_\varepsilon)},
 \end{align*}
 where the hidden constant depends on $\det(\nabla\boldsymbol{\Phi}_j^{\varepsilon})$ and  $\|\boldsymbol{\Phi}_j^{\varepsilon}\|_{ \mathcal{M}^{s+1,p,q}_{\X}(\RR^n_+)}$ being uniformly controlled by the corresponding one for $\varphi_j$.
 
Similarly, concerning $h_j^\varepsilon$, by keeping the expression $h_j^\varepsilon = \div((\tilde{\boldsymbol{\Phi}}^{\varepsilon}_{j,\ast})^{-1}(\eta_j\uu^\varepsilon))$, one applies Lemma~\ref{lem:preservingGengShen} in combination with Lemma~\ref{lem:SmoothingPullback}, to obtain
\begin{align*}
|\lambda| \|h_j^\varepsilon\|_{\dot{\B}^{s-1}_{p,q}(\RR^n_+)}&\lesssim_{p,s,n,\Omega} |\lambda| \| \uu^{\varepsilon}\|_{{\B}^{s-1}_{p,q}(\Omega_\varepsilon)}.
\end{align*}
Now, writing $h_j^\varepsilon = (\tilde{\boldsymbol{\Phi}}^{\varepsilon}_{j,\ast})^{-1}\div(\eta_j\uu^\varepsilon)=(\tilde{\boldsymbol{\Phi}}^{\varepsilon}_{j,\ast})^{-1}[\nabla \eta_j \cdot \uu^{\varepsilon}]$, since one has compact support in $(\boldsymbol{\Phi}^{\varepsilon}_j)^{-1}({\Omega}_\varepsilon)$ by Lemmas~\ref{lem:CompactEquiHomInhom}, \ref{lem:SmoothingPullback} and \ref{lem:multiplierComposition},  
\begin{align*}
     \|h_j^\varepsilon\|_{\dot{\B}^{s+1}_{p,q}(\RR^n_+)}\lesssim_{p,s,n}^{\Omega} \|h_j^\varepsilon\|_{{\B}^{s+1}_{p,q}(\RR^n_+)} \lesssim_{p,s,n}^{\Omega} \| \div(\eta_j \uu^\varepsilon)\|_{{\B}^{s+1}_{p,q}(\mathcal{U}_j)}\lesssim_{p,s,n}^{\Omega} \|\uu^\varepsilon\|_{{\B}^{s+1}_{p,q}(\Omega_\varepsilon)}.
\end{align*}
Plugging this and \eqref{eq:0301} into \eqref{eq:0201c} shows for all $j\in\llb1,\ell\rrb$
 \begin{align}\label{almost0}
|\lambda|\|{\uu}^{\varepsilon}_j\|_{\B^{s}_{p,q}(\Omega_\varepsilon)}+\|{\uu}_j^{\varepsilon}\|_{\B^{s+2}_{p,q}(\Omega_\varepsilon)}+\|{\mathfrak{p}}_j^\varepsilon\|_{\B^{s+1}_{p,q}(\Omega_{\varepsilon})}&\lesssim_{p,s,n}^{\Omega,\mu} |\lambda| \|\uu^\varepsilon\|_{\B^{s-1}_{p,q}(\Omega_\varepsilon)} +
 \|\uu^\varepsilon\|_{\B^{s+1}_{p,q}(\Omega_\varepsilon)}\\& \qquad\qquad+ \|\mathfrak{p}^{\varepsilon}\|_{\B^{s}_{p,q}(\Omega_{\varepsilon})}+\|\ff\|_{\B^{s}_{p,q}(\Omega_\varepsilon)},\nonumber
\end{align}
provided $\epsilon$ is sufficiently small. Clearly, the same estimate (even without the first two terms on the right-hand side) holds for $j=0$ by local regularity theory for the Stokes system (the case of the whole space).

Summing up everything we obtain
 \begin{align}\label{almost1}
|\lambda|\|{\uu}^{\varepsilon}\|_{\B^{s}_{p,q}(\Omega_\varepsilon)}+\|{\uu}^{\varepsilon}\|_{\B^{s+2}_{p,q}(\Omega_\varepsilon)}+\|{\mathfrak{p}}^\varepsilon\|_{\B^{s+1}_{p,q}(\Omega_{\varepsilon})}&\lesssim_{p,s,n}^{\Omega,\mu} |\lambda| \|\uu^\varepsilon\|_{\B^{s-1}_{p,q}(\Omega_\varepsilon)} +
 \|\uu^\varepsilon\|_{\B^{s+1}_{p,q}(\Omega_\varepsilon)}\\& \qquad\qquad+ \|\mathfrak{p}^{\varepsilon}\|_{\B^{s}_{p,q}(\Omega_{\varepsilon})}+\|\ff\|_{\B^{s}_{p,q}(\Omega_\varepsilon)},\nonumber
\end{align}

Now, choose a universal extension operator $\mathcal{E}_{\Omega_{\varepsilon}}$ from $\Omega_{\varepsilon}$ to $\Omega$, for instance the one from \cite{Rychkov1999}, where the boundedness constants only depends on the measure and the diameter of the domain as well as its Lipschitz constant. One obtains the following uniform bound thanks the work performed in Step 1
\begin{align*}
   \limsup_{\varepsilon\rightarrow 0} \,\lVert \mathcal{E}_{\Omega_{\varepsilon}}\uu^{\varepsilon}\rVert_{\B^{s+2}_{p,q}(\Omega)} &\lesssim_{p,s,n,\Omega}  \limsup_{\varepsilon\rightarrow 0} \, \lVert \uu\rVert_{\B^{s+2}_{p,q}(\Omega_\varepsilon)}\\ &\lesssim_{p,s,n,\Omega,\mu} |\lambda| \|\uu\|_{\B^{s-1}_{p,q}(\Omega)} +
 \|\uu\|_{\B^{s}_{p,q}(\Omega)}+ \|\mathfrak{p}\|_{\B^{s}_{p,q}(\Omega)}+\|\ff\|_{\B^{s}_{p,q}(\Omega)} <\infty,
\end{align*}
where the right-hand side is obtained through the convergence  given by the extension by $0$ to the whole domain $\Omega$ (as exhibited in Step 1). The result holds for all $-1+\sfrac{1}{p}<s<\min(-1+\sfrac{(1+\alpha)}{p},0)$. Therefore, by boundedness and compact embeddings $(\mathcal{E}_{\Omega_{\varepsilon}}\uu^{\varepsilon})_{\varepsilon}$ admits a limit in $\B^{s+2}_{p,q}(\Omega)$ (only weak in $\W^{2,p}$ in the case of $\L^p$ spaces). Necessarily, this limit coincides with $\uu$ in $\Omega_\varepsilon$ for any $\varepsilon$, thus in $\Omega$. Consequently, the solution is such that $\uu\in\B^{s+2}_{p,q}(\Omega)$ (if $s=0$, on $\L^p$, one obtains regularity $\W^{2,p}$). 

It follows that the solution has the desired smoothness and one can reperform the whole procedure without any changes, just removing all the $\varepsilon$ and dealing with the non-smoothed charts instead. The constants being uniform with respect to $\Omega$ and $\ff \in\Ccinftydiv(\Omega)$ being an arbitrary function, we obtain 

\begin{align*}
(1+|\lambda|)\|{\uu}\|_{\B^{s}_{p,q}(\Omega)}+\|{\uu}\|_{\B^{s+2}_{p,q}(\Omega)}+\|{\mathfrak{p}}\|_{\B^{s+1}_{p,q}(\Omega)}&\lesssim_{p,s,n}^{\Omega,\mu}  |\lambda| \|\uu\|_{\B^{s-1}_{p,q}(\Omega)} +
 \|\uu\|_{\B^{s+1}_{p,q}(\Omega)}\\& \qquad\qquad+ \|\mathfrak{p}\|_{\B^{s}_{p,q}(\Omega)}+\|\ff\|_{\B^{s}_{p,q}(\Omega)}\\
 &\lesssim_{p,s,n}^{\Omega,\mu}  
 \|\uu\|_{\B^{s+1}_{p,q}(\Omega)}+ \|\mathfrak{p}\|_{\B^{s}_{p,q}(\Omega)}+\|\ff\|_{\B^{s}_{p,q}(\Omega)},\nonumber
\end{align*}
the last inequality coming from the fact that $\lambda \uu = \Delta \uu -\nabla \mathfrak{p} +\ff$.
Exactly the same goes for the spaces $\H^{s,p}$, $p\in(1,\infty)$, for which the case $s= 0$ is allowed whenever $1< p\leqslant 1+\alpha$.

 \textbf{Step 3: The final resolvent estimate.}

 We start by dealing with the pressure term on the right hand side. By Ne\v{c}as' negative norm theorem, Lemma~\ref{lem:NegativeNormThm}, one obtains
\begin{align*}
(1+|\lambda|)\|{\uu}\|_{\B^{s}_{p,q}(\Omega)}+\|{\uu}\|_{\B^{s+2}_{p,q}(\Omega)}+\|{\mathfrak{p}}\|_{\B^{s+1}_{p,q}(\Omega)}&\lesssim_{p,s,n}^{\Omega,\mu} \|{\uu}\|_{\B^{s+1}_{p,q}(\Omega)}+\|\nabla\mathfrak{p}\|_{\B^{s-1}_{p,q}(\Omega)}+\|\ff\|_{\B^{s}_{p,q}(\Omega)}\\
&\lesssim_{p,s,n,\theta}^{\Omega,\mu} \|{\uu}\|_{\B^{s+1}_{p,q}(\Omega)}+ \|\nabla\mathfrak{p}\|_{\B^{s-{\vartheta}}_{p,q}(\Omega)}+\|\ff\|_{\B^{s}_{p,q}(\Omega)},\nonumber
\end{align*}
provided $\vartheta>0$ is such that $s-{\vartheta}>-1+\sfrac{1}{p}$.

Applying the divergence operator to  the equation $\lambda \uu - \Delta \uu +\nabla\mathfrak{p} = \ff$ yields
\begin{align*}
    -\Delta \mathfrak{p} =0,\quad\text{in }\Omega,\quad\text{ and }  \quad \partial_\nu[\mathfrak{p}-\Delta \uu] =0,\quad\text{on }\partial\Omega.
\end{align*}
These equations hold and are meaningful due  to the fact that $\Delta\uu\in\B^{s}_{p,q}(\Omega)$ is divergence free. By Proposition~\ref{prop:NeumannPbLipschitzMult}, it holds
\begin{align*}
   \|\nabla\mathfrak{p}\|_{\B^{s-{\vartheta}}_{p,q}(\Omega)}\lesssim_{p,s,\vartheta}^{n,\Omega} \|\Delta \uu\|_{\B^{s-{\vartheta}}_{p,q}(\Omega)}\lesssim_{p,s,\vartheta}^{n,\Omega} \|\uu\|_{\B^{s+2-{\vartheta}}_{p,q}(\Omega)}.
\end{align*}
Hence, it is deduced that
\begin{align*}
(1+|\lambda|)\|{\uu}\|_{\B^{s}_{p,q}(\Omega)}+\|{\uu}\|_{\B^{s+2}_{p,q}(\Omega)}+\|{\mathfrak{p}}\|_{\B^{s+1}_{p,q}(\Omega)}&\lesssim_{p,s,n,\theta}^{\Omega,\mu} \|\uu\|_{\B^{s+2-{\vartheta}}_{p,q}(\Omega)}+\|\ff\|_{\B^{s}_{p,q}(\Omega)}.
\end{align*}

Now, we fix $s_0<s$ and $p_0\in(1,\infty)$ such that the following embeddings are true (if $p\in(1,\infty)$ one can choose  $p_0=p$):
\begin{align}
\L^{p_0}({\Omega})&\hookrightarrow \B^{s_0}_{p,q}({\Omega}),\label{eq:emb1}\\
\B^{s}_{p,q}(\Omega)& \hookrightarrow \W^{-1,p_0}(\Omega). \label{eq:emb2}
\end{align}
For the first embedding, it is sufficient to choose $s_0$ sufficiently small (it is always true if $s_0\leqslant -n$). For the second,
we choose $p_0$ 'close' to $p$ (or sufficiently large if $p=\infty$) and recall the assumption $-1+\sfrac{1}{p}<s<\min(0,-1+\sfrac{(1+\alpha)}{p})$.
Now, there is $\theta\in(0,1)$ such that
 \begin{align*}
\|\uu\|_{\B^{s+2-{\vartheta}}_{p,q}(\Omega)}&\lesssim_{p,s,s_0}^{n,\Omega} \|\uu\|_{\B^{s+2}_{p,q}(\Omega)}^{\theta}\|\uu\|_{\B^{s_0}_{p,q}(\Omega)}^{1-\theta}\lesssim_{p,s,s_0}^{n,\Omega} \|\uu\|_{\B^{s+2}_{p,q}(\Omega)}^{\theta}\|\uu\|_{\L^{p_0}(\Omega)}^{1-\theta} \\&\lesssim_{p,s,s_0}^{n,\Omega}\|\uu\|_{\B^{s+2}_{p,q}(\Omega)}^{\theta}\|\ff\|_{\W^{-1,p_0}(\Omega)}^{1-\theta}\lesssim_{p,s,s_0}^{n,\Omega}\|\uu\|_{\B^{s+2}_{p,q}(\Omega)}^{\theta}\|\ff\|_{\B^{s}_{p,q}(\Omega)}^{1-\theta}
 \end{align*}
 using the estimate 
 \begin{align}\label{eq:shen}
 (1+|\lambda|)^{\sfrac{1}{2}}\|\uu\|_{\L^{p_0}(\Omega)}\lesssim_{p_0,s,n,\Omega} \|\ff\|_{\W^{-1,p_0}(\Omega)}
 \end{align}
 from \cite[Theorem~5.1]{GengShen2024}, obtained for the resolvent problem in bounded $\C^1$-domains in the penultimate and the embedding \eqref{eq:emb2} in the ultimate step.
 Hence, by interpolation inequalities, we obtain
  \begin{align*}
 \|\uu\|_{\B^{s+1}_{p,q}(\Omega)}&\leqslant\kappa\|\uu\|_{\B^{s}_{p,q}(\Omega)}+C_\kappa\|\ff\|_{\B^{s}_{p,q}(\Omega)}
 \end{align*}
for any $\kappa>0$, and consequently
\begin{align*}
(1+|\lambda|)\|{\uu}\|_{\B^{s}_{p,q}(\Omega)}+\|{\uu}\|_{\B^{s+2}_{p,q}(\Omega)}+\|{\mathfrak{p}}\|_{\B^{s+1}_{p,q}(\Omega)}&\lesssim_{p,s,n}^{\Omega,\mu} 
 \kappa\|\uu\|_{\B^{s+2}_{p,q}(\Omega)}+(1+C_\kappa)\|\ff\|_{\B^{s}_{p,q}(\Omega)},\nonumber
\end{align*}
Choosing $\kappa$ sufficiently small yields
\begin{align*}
(1+|\lambda|)\|{\uu}\|_{\B^{s}_{p,q}(\Omega)}+\|{\uu}\|_{\B^{s+2}_{p,q}(\Omega)}+\|{\mathfrak{p}}\|_{\B^{s+1}_{p,q}(\Omega)}&\lesssim_{p,s,n}^{\Omega,\mu} \|\ff\|_{\B^{s}_{p,q}(\Omega)}.
\end{align*}
Obviously, the same goes for the spaces $\H^{s,p}$, $s\leqslant 0$.

\medbreak

\textbf{Step 4: Reaching the remaining cases $\H^{s,p}$, $0<s<-1+\sfrac{(1+\alpha)}{p}$, $1<p<1+\alpha$ and $\B^{s}_{p,q}$, $0\leqslant s<-1+\sfrac{(1+\alpha)}{p}$, $1\leqslant p<1+\alpha$ ($s>0$, if $p=1$).} The idea here is to bootstrap the strategy from Steps 1-3. By Step $3$, choosing $\ff\in\Ccinftydiv(\Omega)$ ensures that the solution lies in $\W^{2,\tilde{p}}(\Omega,\CC^n)$ for all $\tilde{p}\in(1,1+\alpha)$, in particular $\uu\in\B^{s+1}_{\tilde{p},q}(\Omega,\CC^n)$. The same approximation procedure and estimates as in Steps 1-3 yields that $\uu$ necessarily lies in $\B^{s+2}_{p,q}(\Omega,\CC^n)$. Consequently, we did obtain that, for $\ff\in\Ccinftydiv(\Omega)$, $\uu$ its associated energy solution in $\W^{1,2}_{0,\sigma}(\Omega)$ also belongs to any  Sobolev and Besov spaces and that, for all $p,q\in[1,\infty]$, $-1+\sfrac{1}{p}<s<-1+\sfrac{(1+\alpha)}{p}$,
\begin{align*}
(1+|\lambda|)\|{\uu}\|_{\B^{s}_{p,q}(\Omega)}+\|{\uu}\|_{\B^{s+2}_{p,q}(\Omega)}+\|{\mathfrak{p}}\|_{\B^{s+1}_{p,q}(\Omega)}&\lesssim_{p,s,n}^{\Omega,\mu}  \|\ff\|_{\B^{s}_{p,q}(\Omega)},\nonumber
\end{align*}
and whenever $1<p<\infty$:
\begin{align*}
(1+|\lambda|)\|{\uu}\|_{\H^{s,p}(\Omega)}+\|{\uu}\|_{\H^{s+2,p}(\Omega)}+\|{\mathfrak{p}}\|_{\H^{s+1,p}(\Omega)}&\lesssim_{p,s,n}^{\Omega,\mu}  \|\ff\|_{\H^{s,p}(\Omega)}.
\end{align*}

Therefore, we did extrapolate the $\L^2$-resolvent of $(\D_2(\AA_\mathcal{D}),\AA_\mathcal{D})$ to  $\B^{s}_{p,q}(\Omega)$ and $\H^{s,p}(\Omega)$, preserving consistency, and we obtain the domain inclusion of the extrapolated operator
\begin{align*}
    \D^{s}_p(\AA_\mathcal{D})\subset\H^{s+2,p}_{\mathcal{D},\sigma}(\Omega).
\end{align*}
Similarly for Besov spaces, real interpolation giving boundedness of the resolvent in the case $q=1,\infty$.

\textbf{Step 5: Uniqueness}

For $\lambda\in\Sigma_\mu\cup\{0\}$. If $p\geqslant 2$, since the domain is bounded one always have an embedding $\H^{s+2,p}_{\mathcal{D},\sigma}(\Omega),\B^{s+2,\sigma}_{p,q,\mathcal{D}}(\Omega) \hookrightarrow\H^{1,2}_{0,\sigma}(\Omega)$, yielding uniqueness. If $p\leqslant 2$, the same argument applies using iteratively Sobolev embeddings instead.
\end{proof}

\begin{remark}\label{rem:AboutTheProof}A direct consequence of the end of last of the previous proof is that $\D^{s}_p(\AA_\mathcal{D})=\H^{s+2,p}_{\mathcal{D},\sigma}(\Omega)$, so that for all $\uu\in \D^{s}_p(\AA_\mathcal{D})$, one has necessarily $$\AA_\mathcal{D}\uu=\PP_{\Omega}(-\Delta_{\mathcal{D}}\uu).$$ Note that this equality does \textbf{not} arise from the definition but is rather a consequence of the regularity properties given by the theorem above. There, we did  extrapolate the resolvent operator from the $\L^2$-case, and our theorem gives bound for a typical energy solutions given a datum in a universally dense subspace. At this stage, even if one considers $\C^{1,\alpha}$-domains, if $\alpha\leqslant\frac{1}{2}$, the equality $\AA_\mathcal{D}=\PP_{\Omega}(-\Delta_{\mathcal{D}})$ is \textit{a priori} not accessible on $\L^2$.
\end{remark}

\newpage

\subsection{The definitive results  for the elliptic theory}\label{sec:InterpdomainsHinftyCalc}

For simplicity, we define the following set of conditions for $p,q,r\in[1,\infty]$, $s\in\RR$:
\begin{align*}\tag{$\prescript{\alpha}{r}{\ast}_{p,q}^s$}\label{eq:ConditionRegularity}
    \begin{cases}
    p\in(r,\infty]\text{, }&s<-1+\frac{1+r\alpha}{p}\,\,\,;\text{ or }\\
    p\in[1,r]\text{, }&s <-1+\alpha+\frac{1}{p};\text{ or}\\
     p=q=r\text{, }&s=-1+\alpha+\frac{1}{p}.
           \end{cases}
\end{align*}
This condition corresponds to the highlighted area at the bottom of the figure below.

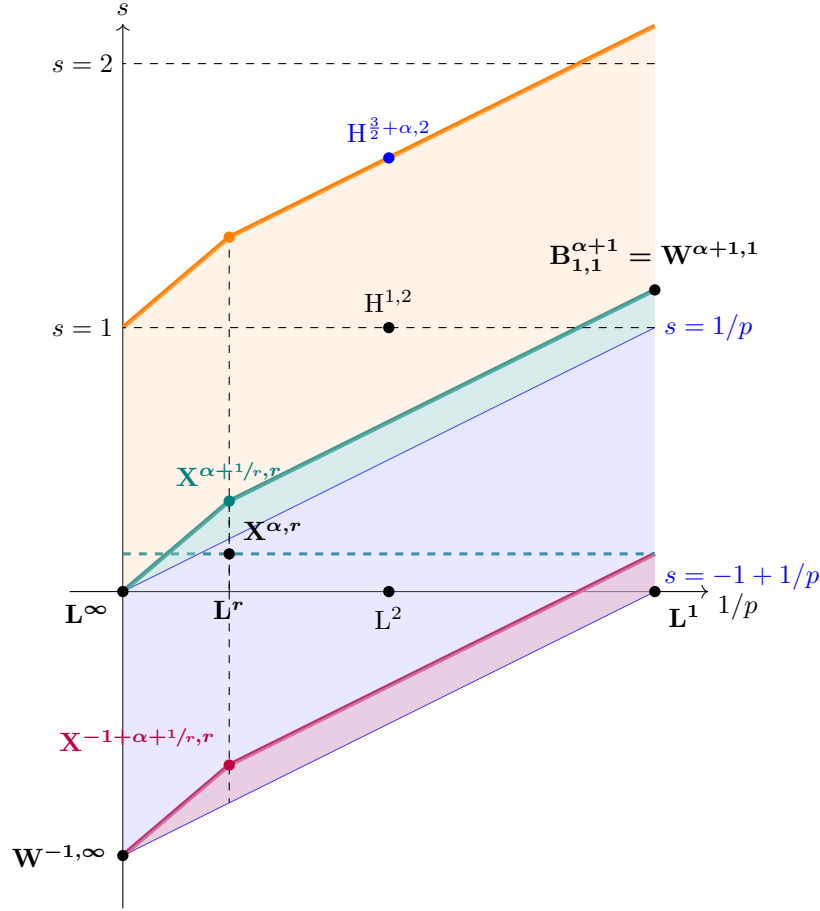
\begin{figure}[H]
\centering
\begin{tikzpicture}[yscale=3.5,xscale=7]
  \draw[->] (-0.1,0) -- (1.1,0) node[right,yshift=-2mm] {$1/p$};
  \draw[->] (0,-1.20) -- (0,2.15) node[above] {$s$};
  \draw[domain=0:1,smooth,variable=\x,blue] plot ({\x},{-1+\x}) node[right,yshift=2mm] {$s=-1+1/p$};
  \draw[domain=0:1,smooth,variable=\x,blue] plot ({\x},{\x}) node[right]{$s=1/p$};
  \draw[domain=0.2:1,line width=0.6mm,variable=\x,teal] plot ({\x},{\x+2/14}) ;
  \draw[domain=0:0.2,line width=0.6mm,variable=\x,teal] plot ({\x},{12*\x/7}) ;
  \draw[domain=0.2:1,line width=0.6mm,variable=\x,purple] plot ({\x},{\x+2/14-1}) ;
  \draw[domain=0:0.2,line width=0.6mm,variable=\x,purple] plot ({\x},{12*\x/7-1}) ;
  \draw[domain=0.2:1,line width=0.6mm,variable=\x,orange] plot ({\x},{\x+2/14+1}) ;
  \draw[domain=0:0.2,line width=0.6mm,variable=\x,orange] plot ({\x},{12*\x/7+1}) ;
  
  \draw[domain=0:1,dashed,line width=0.4mm,variable=\x,teal] plot ({\x},{2/14}) ;

  \fill[blue!30,opacity=0.3]  (0,0)-- (1,1)  -- (1,0) -- (0,-1) -- cycle;
  \fill[teal!30,opacity=0.5]  (0,0)-- (1/5,12/35)   -- (1,1+2/14)  -- (1,1)  -- cycle;
  \fill[purple!30,opacity=0.5]  (0,-1)-- (1/5,12/35-1)   -- (1,1+2/14-1)  -- (1,0)  -- cycle;
  \fill[orange!30,opacity=0.3]  (0,0)-- (0,1)  -- (1/5,12/35+1) -- (1,2+2/14) --  (1,1+2/14) -- (1/5,12/35) -- cycle;

  \draw[dashed] (1,1) -- (0,1)  node[left] {$s=1$};
  \draw[dashed] (1,2) -- (0,2)  node[left] {$s=2$};

  \draw[dashed] (1/5,12/35) -- (1/5,0)  node[below] {$\boldsymbol{{\L}^{r}}$};
  \draw[dashed] (1/5,0) -- (1/5,12/35-1-2/14) ;
  \draw[dashed] (1/5,0) -- (1/5,12/35+1) ;
  
  \node[circle,fill,inner sep=1.5pt,label=below:$\L^2$] at (0.5,0) {};
  \node[circle,fill,inner sep=1.5pt,label=above:${\H}^{1,2}$] at (0.5,1) {};
  \node[circle,fill,inner sep=1.5pt,label=above:{\color{blue}${\H}^{\frac{3}{2}+\alpha,2}$},blue] at (0.5,0.5+2/14+1) {};
  \node[circle,fill,inner sep=1.5pt,orange] at (1/5,12/35+1) {};
  \node[circle,fill,inner sep=1.5pt,label=above:{\color{teal}$\boldsymbol{{\X}^{\alpha+{\sfrac{1}{r}},r}}$},teal] at (1/5,12/35) {};
  \node[circle,fill,inner sep=1.5pt,label=above right:{$\boldsymbol{{\X}^{\alpha,r}}$}] at (1/5,2/14) {};
  \node[circle,fill,inner sep=1.5pt,label=above left:{\color{purple}$\boldsymbol{{\X}^{-1+\alpha+\sfrac{1}{r},r}}$},purple] at (1/5,12/35-1) {};
  \node[circle,fill,inner sep=1.5pt,label=above:{$\boldsymbol{{\B}^{\alpha+1}_{1,1}={\W}^{\alpha+1,1}}$}] at (1,1+2/14) {};
  \node[circle,fill,inner sep=1.5pt,label=below left:{$\boldsymbol{{\L}^\infty}$}] at (0,0) {};
  \node[circle,fill,inner sep=1.5pt,label=below right:{$\boldsymbol{{\L}^1}$}] at (1,0) {};
  \node[circle,fill,inner sep=1.5pt,label=left:{$\boldsymbol{{\W}^{-1,\infty}}$}] at (0,-1) {};
\end{tikzpicture}
\caption{Stabilitised function spaces for a given boundary $\partial\Omega$ in the multiplier class $\mathcal{M}^{1+\alpha,r}_\W$ and associated maximal gain of Sobolev regularity for the solutions of the Stokes--Dirichlet (resolvent) system.}
\end{figure}
The upper threshold is the optimal regularity one can reach on the scale for a fixed Lebesgue index $p$,  provided a boundary in the multiplier class $\mathcal{M}^{1+\alpha,r}_\W$. Here it is illustrated when one wants to solve the system \eqref{eq:Stokes} with say $\ff$ that lies in $\L^2(\Omega)$(here $r$ fixed is such that $r>2$). If instead $\ff\in\X^{s,p}(\Omega)$ for $s,p,q$ satisfying \eqref{eq:ConditionRegularity}, \textit{i.e.} in the bottom area, one always gain an exact amount of two full derivatives in the usual (Sobolev) sense, according to Theorem~\ref{thm:StokesResolvent}. The results we aim to prove  are just the rigorous presentation of such statements and even beyond.

\subsubsection{The results on \texorpdfstring{$\L^p$}{Lp} and Sobolev spaces}\label{sec:ProofLPbddDomain}

The next Theorem is a significant improvement over the current state of the art, providing a gain of between $1$ and $2$ derivatives with respect to the regularity of the boundary\footnote{In fact, even slightly more in some sense, since we will remove the $\C^1$-assumption in Section~\ref{Sec:RemoveC1}}. For the boundedness of the $\mathbf{H}^{\infty}$-functional calculus on $\L^p$-spaces in the case of smooth bounded domains (which implies BIP, the characterization of the domain of fractional powers due to smoothness, and $\L^q$-maximal regularity), we refer to the work of Noll and Saal \cite[Theorem~3]{NollSaal2003} for $\C^3$-domains, and to Kalton, Kunstmann, and Weis \cite[Theorem~9.17]{KaltonKunstmannWeis2006} together with \cite[Theorem~2.1]{KunstmannWeis2013Errata} for $\C^{2,\alpha}$-domains with $0<\alpha<1$. For bounded Lipschitz domains, the boundedness of the $\mathbf{H}^{\infty}$-functional calculus holds on $\L^p$ only for 
\begin{align*}
     p \in \Bigl(\tfrac{2n}{n+1}-\varepsilon,\,\tfrac{2n}{n-1}+\varepsilon\Bigr),
\end{align*}
where $\varepsilon>0$ depends on $\Omega$, see the work of Kunstmann and Weis \cite[Theorem~16]{KunstmannWeis2017}. However, in the work of Kunstmann and Weis, the identification of the domain of fractional powers is restricted to a particular range of regularity for $s$, with a subsequent improvement by Tolksdorf to the range $s\in(-1+\sfrac{1}{p},1]$, see \cite[Theorem~1.1]{Tolksdorf2018-1}. The precise behavior could it be in the case of $\mathbf{H}^\infty$-calculus as well as the characterization of the domain fractional powers, for boundary regularity that lies between Lipschitz and $\C^{2,\alpha}$, $\alpha>0$, was still not solved yet.

For the sole property of bounded imaginary powers in the case of bounded $\C^\infty$-domains---and consequently the identification of the domain of fractional powers---it goes back originally to Giga \cite[Theorems~1~\&~3]{Giga85} by a Pseudo-Differential Operators approach.

\begin{theorem}\label{thm:FinalResultSobolev1}   Let $\alpha\in(0,1)$, $r\in[1,\infty]$, and let $\Omega$ be a bounded $\C^1$-domain of class $\mathcal{M}^{1+\alpha,r}_\W(\epsilon)$ for some sufficiently small $\epsilon>0$.

\medbreak

Let $p\in(1,\infty)$. Then it holds that
\begin{enumerate}
    \item $\AA_{\mathcal{D}}\,:\,\H^{s+2,p}_{\mathcal{D},\sigma}(\Omega)\longrightarrow\H^{s,p}_{\mathfrak{n},\sigma}(\Omega)$ is an isomorphism, whenever $s>-1+\frac{1}{p}$ satisfies \hyperref[eq:ConditionRegularity]{$(\prescript{\alpha}{r}{\ast}_{p,p}^s)$}.
    
    \item The operator $(\D_p(\AA_\mathcal{D}),\AA_\mathcal{D})$ is densely defined, invertible and $0$-sectorial on $\L^p_{\mathfrak{n},\sigma}(\Omega)$, and it admits a bounded $\mathbf{H}^{\infty}(\Sigma_\theta)$-functional calculus for all $\theta\in(0,\pi)$;
    \item For all $s\in(-1+\sfrac{1}{p},\sfrac{1}{p})$, one has
    \begin{align*}
        \D_p^s(\AA_\mathcal{D})=\begin{cases}
           \H^{s+2,p}_{\mathcal{D},\sigma}(\Omega)= \H^{s+1,p}_{0,\sigma}\cap\H^{s+2,p}(\Omega),&\text { if  \hyperref[eq:ConditionRegularity]{$(\prescript{\alpha}{r}{\ast}_{p,p}^s)$} is satisfied};\\
            \{\, \uu\in\H^{1+\frac{1}{p}+\frac{\alpha}{2p}}_{\mathcal{D},\sigma}(\Omega)\,:\,\PP_{\Omega}(-\Delta_\mathcal{D}\uu)\in\H^{s,p}_{\mathfrak{n},\sigma}(\Omega)\,\},&\text { otherwise.}
        \end{cases}
    \end{align*}
    So that, in any case, for all $\uu\in\D_p^s(\AA_\mathcal{D})$,
    \begin{align*}
        \AA_\mathcal{D}\uu = \PP_{\Omega}(-\Delta_\mathcal{D}\uu);
    \end{align*}
    \item For all $s\in(-1+\sfrac{1}{p},2+\sfrac{1}{p})$, $s\neq\sfrac{1}{p},1+\sfrac{1}{p}$, such that \hyperref[eq:ConditionRegularity]{$(\prescript{\alpha}{r}{\ast}_{p,p}^{s-2})$} is satisfied, one has with equivalence of norms
    \begin{align*}
        \D_p(\AA_\mathcal{D}^{\sfrac{s}{2}}) =  \H^{s,p}_{\mathcal{D},\sigma}(\Omega).
  \end{align*}
\end{enumerate}
\end{theorem}

\begin{theorem}\label{thm:FinalResultSobolev2} Let $\alpha\in(0,1)$, $r\in[1,\infty]$, and let $\Omega$ be a bounded $\C^1$-domain of class $\mathcal{M}^{1+\alpha,r}_\W(\epsilon)$ for some sufficiently small $\epsilon>0$. Let $p\in(1,\infty)$ and $s\in(-1+{\sfrac{1}{p}},\sfrac{1}{p})$. For $\mu\in[0,\pi)$, for all $\lambda\in\Sigma_{\mu}\cup\{0\}$, for all $f\in\H^{s,p}_{\mathfrak{n},\sigma}(\Omega)$, 
\begin{enumerate}
    \item If $s$ is such that \hyperref[eq:ConditionRegularity]{$(\prescript{\alpha}{r}{\ast}_{p,p}^{s})$}, then the resolvent problem \eqref{eq:Stokes} admits a unique solution $$(\uu,\mathfrak{p})\in\H^{s+2,p}(\Omega,\CC^n)\times\H^{s+1,p}_{\mfree}(\Omega),$$ which also satisfies the following estimate
\begin{align*}
     (1+\lvert\lambda\rvert) \lVert \uu\rVert_{\H^{s,p}(\Omega)}+\lVert ( \nabla^2 \uu, \nabla \mathfrak{p})\rVert_{\H^{s,p}(\Omega)} \nonumber \lesssim_{p,n,s,\mu}^{\Omega} \lVert \ff\rVert_{\H^{s,p}(\Omega)};
\end{align*}
    \item Otherwise the resolvent problem \eqref{eq:Stokes} admits a unique solution $(\uu,\mathfrak{p})$ satisfying the membership 
    $$(\uu,\mathfrak{p})\in\H^{1+\beta,p}(\Omega,\CC^n)\times\H^{\beta,p}_\mfree(\Omega),$$
    for all $\beta\in[0,1+s]$ such that \hyperref[eq:ConditionRegularity]{$(\prescript{\alpha}{r}{\ast}_{p,p}^{\beta-1})$} and satisfies the following estimate for all such $\beta$:
\begin{align*}
     (1+\lvert\lambda\rvert) \lVert \uu\rVert_{\H^{s,p}(\Omega)}+(1+\lvert\lambda\rvert)^\frac{s-(\beta-1)}{2}\lVert (\nabla \uu,\, \mathfrak{p})\rVert_{\H^{\beta,p}(\Omega)}  \nonumber \lesssim_{p,n,s,\mu}^{\Omega,\beta} \lVert \ff\rVert_{\H^{s,p}(\Omega)}.
\end{align*}
\end{enumerate}
\end{theorem}

We mainly aim to prove Theorem~\ref{thm:FinalResultSobolev1}, the other results will be then more or less straightforward consequences. As for the proof of Theorem~\ref{thm:StokesResolvent}, we just deal with the case $r=1$, the cases $r\in(1,\infty]$ admitting either a similar or easier proof. We recall the overall scheme of the proof is the following
\begin{itemize}
    \item We prove first the boundedness of the $\mathbf{H}^\infty(\Sigma_\theta)$-functional calculus for $(\D_p^s(\AA_\mathcal{D}),\AA_\mathcal{D})$ on $\H^{s,p}_{\mathfrak{n},\sigma}(\Omega)$, $-1+\sfrac{1}{p}<s<-1+\sfrac{(1+\alpha)}{p}$, $p\geqslant 2$.
    \item Thanks to BIP inherited from the bounded holomorphic functional calculus, we compute explicitly the domain of fractional powers, first when $2<p<\infty$, and secondly for $1<p<2$ by duality.
    \item Thanks to characterisation of the domain of fractional powers, we provide the explicit domain and explicit general formula for the extrapolated Stokes operator on $\H^{s,p}_{\mathfrak{n},\sigma}$, even when $s\geqslant -1+\sfrac{(1+\alpha)}{p}$.
\end{itemize}

For $\mathcal{R}$-boundedness and the properties and definition we are using, we refer to Section~\ref{sec:AsbtractMaxReg}, Definition~\ref{def:rbound}, Section~\ref{sec:RboundRn+Hsp} and Proposition~\ref{prop:RboundednessFudnamentalProp}.

\begin{proof}[of Theorems~\ref{thm:FinalResultSobolev1}~\&~Theorems~\ref{thm:FinalResultSobolev2}]Note that Point~\textit{(i)} from Theorem~\ref{thm:FinalResultSobolev1} is just Theorem~\ref{thm:StokesResolvent}. We recall that as for the proof of Theorem~\ref{thm:StokesResolvent}, we just deal with the case $r=1$. 

\textbf{Step 1:} \textbf{The boundedness of the $\mathbf{H}^\infty(\Sigma_\theta)$-functional calculus for $(\D_p^s(\AA_\mathcal{D}),\AA_\mathcal{D})$ on $\H^{s,p}_{\mathfrak{n},\sigma}(\Omega)$, $-1+\sfrac{1}{p}<s<-1+\sfrac{(1+\alpha)}{p}$, $p\geqslant 2$.}

The true goals and substeps for the boundedness of the $\mathbf{H}^\infty(\Sigma_\theta)$-functional calculus are the following
\begin{itemize}
    \item We prove $\mathcal{R}$-boundedness of the resolvent of the operator $\AA_\mathcal{D}$ on $\H^{s,p}_{\mathfrak{n,\sigma}}(\Omega)$, $-1+\sfrac{1}{p}<s<-1+\sfrac{(1+\alpha)}{p}$.
    \item We use the $\mathcal{R}$-boundedness of the resolvent for the Stokes operator and the $\mathbf{H}^\infty(\Sigma_\theta)$-functional calculus of the Dirichlet Laplacian Theorem~\ref{thm:DirLapC1Domains}, to extrapolate Dirichlet--Stokes operator's  $\mathbf{H}^\infty(\Sigma_\theta)$-functional calculus from $\H^{s,2}_{\mathfrak{n},\sigma}(\Omega)$, $s\in(-\sfrac{1}{2},\sfrac{1}{2})$, to $\H^{s,p}_{\mathfrak{n,\sigma}}(\Omega)$, $-1+\sfrac{1}{p}<s<-1+\sfrac{(1+\alpha)}{p}$, $p\in(2,\infty)$.
\end{itemize}

\medbreak

\textbf{Step 1.1:} \textbf{We prove the $\mathcal{R}$-boundedness of the resolvent on $\H^{s,p}_{\mathfrak{n},\sigma}(\Omega)$,  $-1+\sfrac{1}{p}<s<-1+\sfrac{(1+\alpha)}{p}$, $1<p<\infty$.} Let $\mu\in(0,\pi)$. Since the resolvent set of an unbounded operator is open and since $0\in\rho(\AA_\mathcal{D})$ by  Theorem~\ref{thm:StokesResolvent}, it holds that
\begin{align*}
    \lambda\longmapsto \lambda(\lambda\I+\AA_\mathcal{D})^{-1}
\end{align*}
is then holomorphic in a neighborhood of $\overline{\Sigma_\mu}$. Consequently, by \cite[Example~2.19]{KunstmannWeis2004}, it holds that
\begin{align*}
    \{ \lambda(\lambda\I+\AA_\mathcal{D})^{-1}\,,\, \lambda\in\Sigma_\mu\cap\overline{\B_R(0)}\cup\{0\} \}
\end{align*}
is $\mathcal{R}$-bounded on $\H^{s,p}_{\mathfrak{n,\sigma}}(\Omega)$ for any $R>0$, and that it suffices then to show the $\mathcal{R}$-boundedness of the family
\begin{align*}
    \{ \lambda(\lambda\I+\AA_\mathcal{D})^{-1}\,,\, \lambda\in\Sigma_\mu,\,|\lambda|>R_0 \}
\end{align*}
on $\H^{s,p}_{\mathfrak{n,\sigma}}(\Omega)$, for some $R_0>0$ to be chosen later on.

\medbreak

Let $(\ff_\ell)_{\ell\in\NN}\subset\H^{s,p}_{\mathfrak{n},\sigma}(\Omega)$, $(\lambda_\ell)_{\ell\in\NN}\subset\Sigma_\mu\cap\overline{\B_{R_0}(0)}^c$, and consider the corresponding family  of solutions $(\uu_\ell,\mathfrak{p}_\ell)_{\ell\in\NN}\subset \H^{s+2,p}_{\mathcal{D},\sigma}(\Omega)\times\H^{s+1,p}_\mfree(\Omega)$ given by Theorem~\ref{thm:StokesResolvent}. Let $\N\in\NN^\ast$, and with the notations introduced in the proof of  Theorem~\ref{thm:StokesResolvent}, $L\in\NN$, and $(\mathcal{U}_j)_{j\in\llb0,L\rrb}$ the covering of $\Omega$ with patch $(\eta_j)_{j\in\llb 0,L\rrb}$, with local (lifted) diffeomorphism $(\boldsymbol{\Phi}_j)_{j\in\llb 0,L\rrb}$, for the localised couples $(\tilde{\uu}_{j,\ell},\tilde{\mathfrak{p}}_{j,\ell}) := [\eta_j(\uu_\ell,\mathfrak{p}_\ell)] \circ \boldsymbol{\Phi}_j$, and corresponding localised data $\ff_{j,\ell}:=(\eta_j\ff_\ell)$, by Proposition~\ref{prop:RboundHspRn+}, writing
\begin{align*}
    \gg_{j,\ell} := \det(\nabla\boldsymbol{\Phi}_j)[\Delta \eta_j \uu_\ell+\nabla \eta_j \cdot \nabla\uu_\ell+ \nabla{\eta_j} \mathfrak{p}_\ell+\ff_{j,\ell}]\circ \boldsymbol{\Phi}_j
\end{align*}
it holds,
\begin{align*}
    \bigg\lVert &\sum_{\ell=1}^{N} r_\ell\lambda_\ell \tilde{\uu}_{j,\ell}\bigg\rVert_{\L^p(0,1;\dot{\H}^{s,p}(\RR^n_+))} + \bigg\lVert \sum_{\ell=1}^{N} r_\ell (\nabla^2 \tilde{\uu}_{j,\ell},\nabla\tilde{\mathfrak{p}}_{j,\ell})\bigg\rVert_{\L^p(0,1;\dot{\H}^{s,p}(\RR^n_+))}\\ \qquad & \lesssim_{p,s,n,\mu}\bigg\lVert \sum_{\ell=1}^{N} r_\ell \lambda_\ell (1-\det(\nabla\Phi_j))\tilde{\uu}_{j,\ell}\bigg\rVert_{\L^p(0,1;\dot{\H}^{s,p}(\RR^n_+))} \\ 
    &\qquad\qquad +\bigg\lVert \sum_{\ell=1}^{N} r_\ell \div([\mathbf{I}_n-\mathbf{A}_j]\tilde{\uu}_{j,\ell})\bigg\rVert_{\L^p(0,1;\dot{\H}^{s,p}(\RR^n_+))}  \\   
    &\qquad\qquad +\bigg\lVert \sum_{\ell=1}^{N} r_\ell \div([\mathbf{I}_n-\mathbf{B}_j]\tilde{\mathfrak{p}}_{j,\ell})\bigg\rVert_{\L^p(0,1;\dot{\H}^{s,p}(\RR^n_+))}  \\
    &\qquad\qquad +\bigg\lVert \sum_{\ell=1}^{N} r_\ell \gg_{j,\ell}\bigg\rVert_{\L^p(0,1;\dot{\H}^{s,p}(\RR^n_+))}  \\  
    &\qquad\qquad + \bigg\lVert \sum_{\ell=1}^{N} r_\ell \lambda_\ell\left( \div([\mathbf{I}_n-\mathbf{B}_j]\tilde{\uu}_{j,\ell}) + \tilde{\Phi}_{j,\ast}^{-1}(\div\eta_j\uu_\ell)\right)\bigg\rVert_{\L^p(0,1;\dot{\H}^{s-1,p}(\RR^n_+))}\\  
    &\qquad\qquad  + \bigg\lVert \sum_{\ell=1}^{N} r_\ell \nabla \left[\div([\mathbf{I}_n-\mathbf{B}_j]\tilde{\uu}_{j,\ell}) + \tilde{\Phi}_{j,\ast}^{-1}(\div\eta_j\uu_\ell)\right]\bigg\rVert_{\L^p(0,1;\dot{\H}^{s,p}(\RR^n_+))}.
\end{align*}
Since the operators arising from the localisation and transformations are singletons and independent of $\ell$, thanks to Proposition~\ref{prop:RboundednessFudnamentalProp}, one can follow the arguments in the proof of Theorem~\ref{thm:StokesResolvent}, in particular Step 2. This yields the similar $\mathcal{R}$-estimate
\begin{align*}
    &\bigg\lVert \sum_{\ell=1}^{N} r_\ell\lambda_\ell \tilde{\uu}_{j,\ell}\bigg\rVert_{\L^p(0,1;\dot{\H}^{s,p}(\RR^n_+))} + \bigg\lVert \sum_{\ell=1}^{N} r_\ell (\nabla^2 \tilde{\uu}_{j,\ell},\nabla\tilde{\mathfrak{p}}_{j,\ell})\bigg\rVert_{\L^p(0,1;\dot{\H}^{s,p}(\RR^n_+))}\\
     \qquad & \lesssim_{p,s,n}^{\mu,\Omega} \epsilon\left[\bigg\lVert \sum_{\ell=1}^{N} r_\ell \lambda_\ell{\uu}_{\ell}\bigg\rVert_{\L^p(0,1;{\H}^{s,p}(\Omega))} + \bigg\lVert \sum_{\ell=1}^{N} r_\ell {\uu}_{\ell}\bigg\rVert_{\L^p(0,1;{\H}^{s+2,p}(\Omega))}+\bigg\lVert \sum_{\ell=1}^{N} r_\ell \mathfrak{p}_{\ell}\bigg\rVert_{\L^p(0,1;{\H}^{s+1,p}(\Omega))} \right] \\ 
    &\qquad\quad +\left[\bigg\lVert \sum_{\ell=1}^{N} r_\ell \lambda_\ell{\uu}_{\ell}\bigg\rVert_{\L^p(0,1;{\H}^{s-1,p}(\Omega))} + \bigg\lVert \sum_{\ell=1}^{N} r_\ell {\uu}_{\ell}\bigg\rVert_{\L^p(0,1;{\H}^{s,p}(\Omega))}+\bigg\lVert \sum_{\ell=1}^{N} r_\ell \mathfrak{p}_{\ell}\bigg\rVert_{\L^p(0,1;{\H}^{s,p}(\Omega))} \right] \\  
    &\qquad\quad +\bigg\lVert \sum_{\ell=1}^{N} r_\ell {\ff}_{\ell}\bigg\rVert_{\L^p(0,1;{\H}^{s,p}(\Omega))}.
\end{align*}
Thus, by the triangle inequality and since $\epsilon>0$ can be chosen small enough, it reduces to the estimate
\begin{align*}
    \bigg\lVert &\sum_{\ell=1}^{N} r_\ell\lambda_\ell {\uu}_{\ell}\bigg\rVert_{\L^p(0,1;{\H}^{s,p}(\Omega))} + \bigg\lVert \sum_{\ell=1}^{N} r_\ell (\nabla^2 {\uu}_{\ell},\nabla{\mathfrak{p}}_{\ell})\bigg\rVert_{\L^p(0,1;{\H}^{s,p}(\Omega))}\\
     \qquad & \lesssim_{p,s,n}^{\mu,\Omega} \left[\bigg\lVert \sum_{\ell=1}^{N} r_\ell \lambda_\ell{\uu}_{\ell}\bigg\rVert_{\L^p(0,1;{\H}^{s-1,p}(\Omega))} + \bigg\lVert \sum_{\ell=1}^{N} r_\ell {\uu}_{\ell}\bigg\rVert_{\L^p(0,1;{\H}^{s,p}(\Omega))}+\bigg\lVert \sum_{\ell=1}^{N} r_\ell \mathfrak{p}_{\ell}\bigg\rVert_{\L^p(0,1;{\H}^{s,p}(\Omega))} \right] \\  
    &\qquad\qquad +\bigg\lVert \sum_{\ell=1}^{N} r_\ell {\ff}_{\ell}\bigg\rVert_{\L^p(0,1;{\H}^{s,p}(\Omega))}.
\end{align*}
Since $\lambda_\ell\uu_\ell=  \Delta\uu_\ell - \nabla\mathfrak{p}_\ell+\ff_\ell$, by linearity it holds,
\begin{align}\label{eq:ProofRboundC1aInterm}
    \bigg\lVert &\sum_{\ell=1}^{N} r_\ell\lambda_\ell {\uu}_{\ell}\bigg\rVert_{\L^p(0,1;{\H}^{s,p}(\Omega))} + \bigg\lVert \sum_{\ell=1}^{N} r_\ell (\nabla^2 {\uu}_{\ell},\nabla{\mathfrak{p}}_{\ell})\bigg\rVert_{\L^p(0,1;{\H}^{s,p}(\Omega))}\\
     \qquad & \lesssim_{p,s,n}^{\mu,\Omega} \bigg\lVert \sum_{\ell=1}^{N} r_\ell {\uu}_{\ell}\bigg\rVert_{\L^p(0,1;{\H}^{s,p}(\Omega))}+\bigg\lVert \sum_{\ell=1}^{N} r_\ell \mathfrak{p}_{\ell}\bigg\rVert_{\L^p(0,1;{\H}^{s,p}(\Omega))} +\bigg\lVert \sum_{\ell=1}^{N} r_\ell {\ff}_{\ell}\bigg\rVert_{\L^p(0,1;{\H}^{s,p}(\Omega))}.\nonumber
\end{align}
By Nec\v{a}s' negative norm theorem Lemma~\ref{lem:NegativeNormThm} applied to $\sum_{\ell=1}^{N} r_\ell \mathfrak{p}_{\ell}$, one obtains successively
\begin{align*}
    \bigg\lVert \sum_{\ell=1}^{N} r_\ell \mathfrak{p}_{\ell}\bigg\rVert_{\L^p(0,1;{\H}^{s,p}(\Omega))} &\lesssim_{p,s,n}^{\Omega} \bigg\lVert \sum_{\ell=1}^{N} r_\ell \nabla \mathfrak{p}_{\ell}\bigg\rVert_{\L^p(0,1;{\H}^{s-1,p}(\Omega))}\\
    &\lesssim_{p,s,n}^{\varepsilon,\Omega} \bigg\lVert \sum_{\ell=1}^{N} r_\ell \nabla \mathfrak{p}_{\ell}\bigg\rVert_{\L^p(0,1;{\H}^{s-\varepsilon,p}(\Omega))}.
\end{align*}
provided a given sufficiently small $\varepsilon>0$, is such that $s-\varepsilon>-1+\sfrac{1}{p}$. For any $\ell$, one has
\begin{equation*}\label{eq:PoissonPbPressureProof}
    \left\{ \begin{array}{rllr}
         -\div (\nabla \mathfrak{p}_\ell) &= 0 \text{, }&&\text{ in } \Omega\text{,}\\
        \nabla\mathfrak{p}_\ell \cdot\nu_{|_{\partial\Omega}} &= \Delta \uu_\ell\cdot\nu _{|_{\partial\Omega}}\text{, } &&\text{ on } \partial\Omega\text{.}
    \end{array}
    \right.
\end{equation*}
By Proposition~\ref{prop:NeumannPbLipschitzMult}, one obtains $\nabla \mathfrak{p}_{\ell} =  [\I-\PP_\Omega](\Delta\uu_\ell)$, and since the Hodge-Leray projection is a singleton of a bounded operator on $\H^{s,p}(\Omega,\CC^n)$ by Theorem~\ref{thm:SharpHodgeDecompC1}, it holds
\begin{align*}
    \bigg\lVert \sum_{\ell=1}^{N} r_\ell \nabla \mathfrak{p}_{\ell}\bigg\rVert_{\L^p(0,1;{\H}^{s-\varepsilon,p}(\Omega))} &\lesssim_{p,s,\varepsilon,n,\Omega}  \bigg\lVert \sum_{\ell=1}^{N} r_\ell \Delta \uu_\ell\bigg\rVert_{\L^p(0,1;{\H}^{s-\varepsilon,p}(\Omega))}.
\end{align*}
due to Proposition~\ref{prop:RboundednessFudnamentalProp}. Now, by interpolation inequality, it yields for all $\kappa>0$,
\begin{align}\label{eq:ProofRboundC1aPressureEst}
     \bigg\lVert \sum_{\ell=1}^{N} r_\ell \mathfrak{p}_{\ell}\bigg\rVert_{\L^p(0,1;{\H}^{s,p}(\Omega))} \leqslant \kappa  \bigg\lVert \sum_{\ell=1}^{N} r_\ell \Delta \uu_\ell\bigg\rVert_{\L^p(0,1;{\H}^{s,p}(\Omega))} + C_\kappa \bigg\lVert \sum_{\ell=1}^{N} r_\ell \uu_\ell\bigg\rVert_{\L^p(0,1;{\H}^{s-1,p}(\Omega))}.
\end{align}
Plugging \eqref{eq:ProofRboundC1aPressureEst} in \eqref{eq:ProofRboundC1aInterm}, by the embedding $\H^{s,p}(\Omega)\hookrightarrow\H^{s-1,p}(\Omega)$, and choosing $\kappa$ small enough, one reaches
\begin{align*}
    \bigg\lVert &\sum_{\ell=1}^{N} r_\ell\lambda_\ell {\uu}_{\ell}\bigg\rVert_{\L^p(0,1;{\H}^{s,p}(\Omega))} + \bigg\lVert \sum_{\ell=1}^{N} r_\ell (\nabla^2 {\uu}_{\ell},\nabla{\mathfrak{p}}_{\ell})\bigg\rVert_{\L^p(0,1;{\H}^{s,p}(\Omega))}\\
     \qquad & \leqslant C_{p,s,n,\mu,\Omega} \bigg( \bigg\lVert \sum_{\ell=1}^{N} r_\ell {\uu}_{\ell}\bigg\rVert_{\L^p(0,1;{\H}^{s,p}(\Omega))}+\bigg\lVert \sum_{\ell=1}^{N} r_\ell {\ff}_{\ell}\bigg\rVert_{\L^p(0,1;{\H}^{s,p}(\Omega))}\bigg)
\end{align*}
for some constant $C_{p,s,n,\mu,\varepsilon,\Omega}>0$. If one chooses $R_0> 4C_{p,s,n,\mu,\varepsilon,\Omega} +2$, it reduces to 
\begin{align*}
    \bigg\lVert &\sum_{\ell=1}^{N} r_\ell\lambda_\ell {\uu}_{\ell}\bigg\rVert_{\L^p(0,1;{\H}^{s,p}(\Omega))} + \bigg\lVert \sum_{\ell=1}^{N} r_\ell (\nabla^2 {\uu}_{\ell},\nabla{\mathfrak{p}}_{\ell})\bigg\rVert_{\L^p(0,1;{\H}^{s,p}(\Omega))} \lesssim_{p,s,\varepsilon}^{n,\mu,\Omega} \bigg\lVert \sum_{\ell=1}^{N} r_\ell {\ff}_{\ell}\bigg\rVert_{\L^p(0,1;{\H}^{s,p}(\Omega))},
\end{align*}
which means in particular that the family
\begin{align*}
    \{ \lambda(\lambda\I+\AA_\mathcal{D})^{-1}\,,\, \lambda\in\Sigma_\mu,\,|\lambda|>R_0 \}
\end{align*}
is $\mathcal{R}$-bounded on $\H^{s,p}_{\mathfrak{n},\sigma}(\Omega)$. 

\textbf{Additional point:} Provided $p\in(1,\infty)$, note that the self-adjointness of $\AA_\mathcal{D}$ on $\L^{2}_{\mathfrak{n},\sigma}(\Omega)$ and duality are allowing to conclude that for all $\mu\in[0,{\pi})$
\begin{align}
    \{ (1+|\lambda|)(\lambda\I+\AA_\mathcal{D})^{-1}\,,\, \lambda\in\Sigma_\mu\cup\{0\} \}\label{eq:DualRboundStokesDir}
\end{align}
extends uniquely by density as a $\mathcal{R}$-bounded family on $\H^{s,p'}_{\mathfrak{n},\sigma}(\Omega)$ , all $s\in(-\sfrac{\alpha}{p}+\sfrac{1}{p'},\sfrac{1}{p'})$.

\medbreak

\textbf{Step 1.2:} We prove extrapolation of the $\mathbf{H}^{\infty}$-functional calculus by the Kunstmann-Weis \textbf{\textit{comparison principle}} \cite[Theorem~9,~Propositions~10~\&~11]{KunstmannWeis2017}. Although the strategy is already established, we present a full version of the argument for convenience. We follow the presentation from \cite[Section~3,~p.250--256]{GabelTolksdorf2022}, where the comparison principle is presented as follows.
\begin{theorem}[ {{Kunstmann \& Weis}, \cite[Theorem~9]{KunstmannWeis2017}} ]\label{Thm:Kunstmann-WeisHinfty}
Let $\X$ and $\Y$ be Banach spaces. Let $\mathfrak{R}\,:\,\Y \longrightarrow \X$ and $\mathfrak{E} \,:\,\X \longrightarrow \Y$ be bounded linear operators satisfying $ \mathfrak{R}\mathfrak{E} = \I$ on $\X$. 

\medbreak

Let $(\D(B),B)$ have a bounded $\mathbf{H}^{\infty} (\Sigma_{\tau})$-functional calculus in $\Y$ for some $\tau \in (0,\pi)$, and let $(\D(A),A)$ be $\mathcal{R}$-sectorial in $\X$. 
Assume that there are functions $\varphi , \psi \in \mathbf{H}_0^{\infty} (\Sigma_{\mu}) \setminus \{ 0 \}$ where $\mu \in (0 , \tau)$ and $C_1, C_2 > 0$ are such that, for some $\beta > 0$ and all $\ell\in \ZZ$,
\begin{align}
 \sup_{1 \leqslant \tau , t \leqslant 2} \mathcal{R} \big\{ \varphi (\tau 2^{j + \ell} A) \mathfrak{R} \psi(t 2^j B) : j \in \ZZ \big\} &\leqslant C_1 2^{- \beta \lvert \ell \rvert} \label{eq:cond}, \\
 \text{and}\qquad\sup_{1 \leqslant \tau , t \leqslant 2} \mathcal{R} \big\{ \varphi (\tau 2^{j + \ell} A)^{\prime} \mathfrak{E}^{\prime} \psi (t 2^j B)^{\prime} : j \in \ZZ \big\} &\leqslant C_2 2^{- \beta \lvert \ell \rvert}. \label{eq:DualCond}
\end{align}
Then $A$ has a bounded $\mathbf{H}^{\infty} (\Sigma_{\mu})$-functional calculus on $\X$.
\end{theorem}

Kunstmann and Weis provided also tools for establishing~\eqref{eq:cond} and~\eqref{eq:DualCond} as given below, summarized in this form following \cite[Section~3,~p.250--256]{GabelTolksdorf2022}. 

\begin{proposition}[ {\cite[Propositions~10~\&~11]{KunstmannWeis2017}} ]\label{Prop:Fractionalcomparison}
In the setting of Theorem~\ref{Thm:Kunstmann-WeisHinfty}, suppose that there exist $\alpha_0 > 0$ and $C > 0$ such that,
for $\alpha = \pm \alpha_0$, we have
\begin{align}\label{Eq: Retraction Condition}
 \mathfrak{R} \, \D(B^{\alpha}) \subset \D(A^{\alpha}), \qquad \| A^{\alpha}  \mathfrak{R} y \|_X \leqslant C \| B^{\alpha} y \|_Y \quad \text{for all} \quad y \in \D(B^{\alpha}),
\end{align}
and
\begin{align}\label{Eq: Injection Condition}
 \mathfrak{E} \, \D(A^{\alpha}) \subset \D(B^{\alpha}), \qquad \| B^{\alpha} \mathfrak{E} x \|_Y \leqslant C \| A^{\alpha} x \|_X \quad \text{for all} \quad x \in \D(A^{\alpha}).
\end{align}
Then Conditions~\eqref{eq:cond} and~\eqref{eq:DualCond} hold for the choice $C_1 = C_2 = C$, $\beta = \alpha_0$, and $\varphi(\lambda) = \psi(\lambda) = \lambda^{2\alpha_0}(1 + \lambda)^{-4 \alpha_0}$.
\end{proposition}

The idea of our argument will be close to the strategy exhibited in \cite[Section~3]{GabelTolksdorf2022}, but instead of extrapolating the $\mathbf{H}^{\infty}$-functional calculus on the pure $\L^p$-scale, we will do it for negative regularity Bessel potential spaces.

We want to apply Theorem~\ref{Thm:Kunstmann-WeisHinfty} to $\X=\H^{s,p}_{\mathfrak{n},\sigma}(\Omega)$, $\Y=\H^{s,p}(\Omega,\CC^n)$, $-1+\sfrac{1}{p}<s<-1+\sfrac{(1+\alpha)}{p}$,
\begin{align*}
    \mathfrak{E}:= \iota\,:\,\H^{s,p}_{\mathfrak{n},\sigma}(\Omega) \longrightarrow \H^{s,p}(\Omega,\CC^n),\qquad\text{ and }\qquad\mathfrak{R}:= \PP_\Omega\,:\, \H^{s,p}(\Omega,\CC^n)\longrightarrow \H^{s,p}_{\mathfrak{n},\sigma}(\Omega),
\end{align*}
and operators
\begin{align*}
    (\D(A),A) = (\D_{p}^s(\AA_\mathcal{D}),\AA_\mathcal{D})\qquad\text{ and }\qquad (\D(B),B) = (\D_{p}^s(-\Delta_\mathcal{D}),-\Delta_\mathcal{D}).
\end{align*}
Note that, by Proposition~\ref{prop:DualityDivFreeRn+} and according to Section~\ref{subsec:GeneralStokesDirPbLipDomains}, one has dual operators
\begin{align*}
    \mathfrak{R}'= (\PP_\Omega)' = \iota&\,:\,\H^{-s,p'}_{\mathfrak{n},\sigma}(\Omega) \longrightarrow \H^{-s,p'}(\Omega,\CC^n),\\
    \text{ and }\qquad\mathfrak{E}'= \iota' =\PP_\Omega&\,:\, \H^{-s,p'}(\Omega,\CC^n)\longrightarrow \H^{-s,p'}_{\mathfrak{n},\sigma}(\Omega),
\end{align*}
and, by self-adjointness, unbounded dual operators
\begin{align*}
    (\D(A'),A') = (\D_{p'}^{-s}(\AA_\mathcal{D}),\AA_\mathcal{D})\qquad\text{ and }\qquad (\D(B),B) = (\D_{p'}^{-s}(-\Delta_\mathcal{D}),-\Delta_\mathcal{D}).
\end{align*}
acting respectively on $\H^{-s,p'}_{\mathfrak{n},\sigma}(\Omega)$, and $\H^{-s,p'}(\Omega,\CC^n)$. This can be checked by universal density of $\Ccinftydiv(\Omega)$ and $\Ccinfty(\Omega,\CC^n)$ in both scale of spaces, extending each resolvent this way. The case of the Stokes--Dirichlet is a particular instance of \eqref{eq:DualRboundStokesDir}. The good behavior of the Dirichlet Laplacian is due to Theorem~\ref{thm:DirLapC1Domains}.

\medbreak

We note that for $\varphi , \psi \in \mathbf{H}_0^{\infty} (\Sigma_{\mu}) \setminus \{ 0 \}$, provided $\mu\in(0,\pi)$, we have
\begin{itemize}
    \item $\{\varphi(t 2^{j+\ell}\AA_\mathcal{D}),\,t>0,\,j,\ell\in\ZZ\}$ is $\mathcal{R}$-bounded in  $\H^{s,p}_{\mathfrak{n},\sigma}(\Omega)$ and $\H^{-s,p'}_{\mathfrak{n},\sigma}(\Omega)$ by $\mathcal{R}$-sectoriality: use Step 1.1 and apply \cite[Lemma~3.3]{KaltonKunstmannWeis2006};
    \item $\{\psi(-t 2^{j+\ell}\Delta_\mathcal{D}),\,t>0,\,j,\ell\in\ZZ\}$ is $\mathcal{R}$-bounded in  $\H^{s,p}(\Omega,\CC^n)$ and $\H^{-s,p'}(\Omega,\CC^n)$ by Theorems~\ref{thm:DirLapC1Domains},~\ref{thm:DoreVennithm},~\ref{thm:WeisMaxRegRbound} (\textit{i.e.}, BIP implies Maximal regularity, implying $\mathcal{R}$-boundedness of the resolvent) and \cite[Lemma~3.3]{KaltonKunstmannWeis2006}.
\end{itemize}
Consequently, by Theorem~\ref{thm:SharpHodgeDecompC1} and Proposition~\ref{prop:RboundednessFudnamentalProp}, we obtain the following uniform $\mathcal{R}$-boundedness 
\begin{align}\label{eq:ProofBddRboundHsp}
    \sup_{\ell\in\ZZ} \sup_{1 \leqslant \tau , t \leqslant 2} \mathcal{R} \big\{ \varphi (\tau 2^{j + \ell} \AA_\mathcal{D}) \PP_\Omega \psi(-t 2^j\Delta_\mathcal{D} ) : j \in \ZZ \big\} \lesssim_{\mu,\varphi,\psi} 1,
\end{align}
as a family of bounded linear operators from $\H^{s,p}(\Omega,\CC^n)$ to $\H^{s,p}_{\mathfrak{n},\sigma}(\Omega)$, and from $\H^{-s,p'}(\Omega,\CC^n)$ to $\H^{-s,p'}_{\mathfrak{n},\sigma}(\Omega)$, for any $s\in(-1+\sfrac{1}{p},-1+\sfrac{(1+\alpha)}{p})$.

In \eqref{eq:ProofBddRboundHsp}, we only have a uniform bound, while we want some "off-diagonal" decay in order to apply Theorem~\ref{Thm:Kunstmann-WeisHinfty}. The idea, is to put $\H^{s,p}$ and $\H^{-s,p'}$ as interior point of an interpolation procedure between the couple of spaces $(\H^{s_0,q},\H^{\tilde{s},2})$ and $(\H^{-\tilde{s},2},\H^{-s_0,q'})$, having uniform boundedness for one and  exponential decay for the other, it should provide the desired decay for the interior spaces $\H^{s,p}$ and $\H^{-s,p'}$.

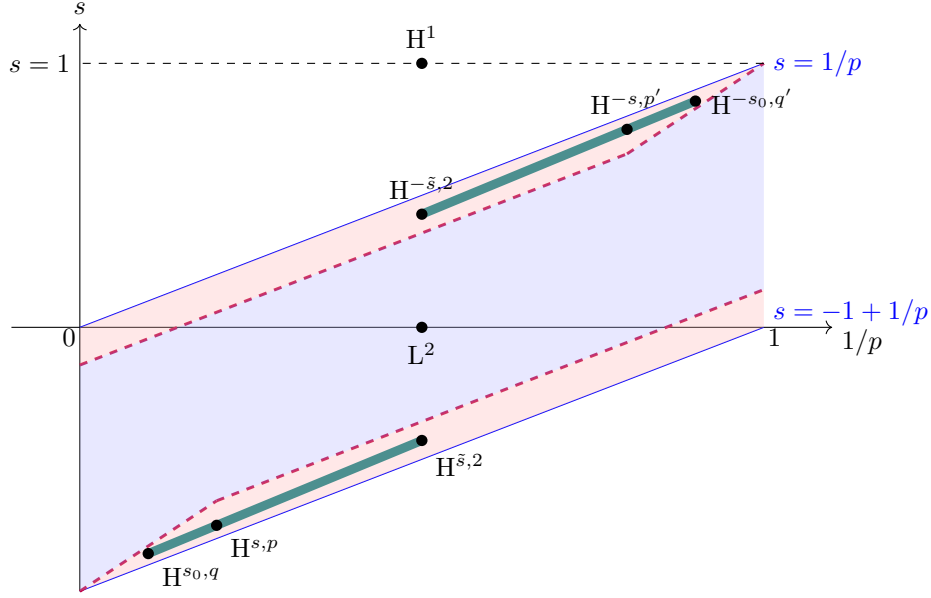
\begin{figure}[H]
\centering
\begin{tikzpicture}[yscale=3.5,xscale=9]
  \draw[->] (-0.1,0) -- (1.1,0) node[right,yshift=-2mm] {$1/p$};
  \draw[->] (0,-1.0) -- (0,1.15) node[above] {$s$};
  \draw[domain=0:1,smooth,variable=\x,blue] plot ({\x},{-1+\x}) node[right,yshift=2mm] {$s=-1+1/p$};
  \draw[domain=0:1,smooth,variable=\x,blue] plot ({\x},{\x}) node[right]{$s=1/p$};

  \draw[domain=0.2:1,dashed,line width=0.4mm,variable=\x,purple] plot ({\x},{\x+2/14-1}) ;
  \draw[domain=0:0.2,dashed,line width=0.4mm,variable=\x,purple] plot ({\x},{12*\x/7-1}) ;\draw[domain=0:0.8,dashed,line width=0.4mm,variable=\x,purple] plot ({\x},{\x-2/14}); 
  \draw[domain=0.8:1,dashed,line width=0.4mm,variable=\x,purple] plot ({\x},{12*\x/7-10/14});
 \draw[domain=0.5:0.90,line width=1.2mm,variable=\x,teal] plot ({\x},{15*\x/14-3/28}) ;
 \draw[domain=0.10:0.5,line width=1.2mm,variable=\x,teal] plot ({\x},{15*\x/14-27/28}) ;

\fill[red!30,opacity=0.3] (0,-1) -- (1/5,12/35-1) -- (1,2/14) -- (1,0) -- cycle;
\fill[red!30,opacity=0.3] (0,0) -- (0,-2/14) -- (0.8,1-12/35)  -- (1,1) -- cycle;
  \fill[blue!30,opacity=0.3]  (0,-1)-- (0,-2/14)  -- (0.8,1-12/35) -- (1,1) --  (1,2/14) -- (1/5,12/35-1) -- cycle;
  \draw[dashed] (1,1) -- (0,1)  node[left] {$s=1$};
  \node[circle,fill,inner sep=1.5pt,label=below:$\L^2$] at (0.5,0) {};
  \node[circle,fill,inner sep=1.5pt,label=above:${\H}^1$] at (0.5,1) {};
  \node[circle,fill,inner sep=1.5pt,label=above:${\H}^{-\tilde{s},2}$] at (0.5,6/14) {};
  \node[circle,fill,inner sep=1.5pt,label=below right:${\H}^{\tilde{s},2}$] at (0.5,-6/14) {};
  \node[circle,fill,inner sep=1.5pt,label=above:${\H}^{-s,p'}$] at (0.8,10.5/14) {};
  \node[circle,fill,inner sep=1.5pt,label=below right:${\H}^{s,p}$] at (0.2,-10.5/14) {};
  \node[circle,fill,inner sep=1.5pt,label=right:${\H}^{-s_0,q'}$] at (0.90,15*0.9/14-3/28) {};
  \node[circle,fill,inner sep=1.5pt,label=below right:${\H}^{s_0,q}$] at (0.10,15*0.1/14-27/28) {} ;
  \draw[circle,fill,inner sep=1pt] (0,0) node[below left] {$0$};
  \draw[circle,fill,inner sep=1pt] (1,0) node[below right] {$1$};
\end{tikzpicture}
\caption{Representation of the interpolation procedure to extrapolate the ``off-diagonal'' decay, provided a given fixed $r\in[1,\infty]$. (Recall that in the current proof $r=1$ is assumed.)}
\end{figure}
Therefore, one has to check the off-diagonal decay on a part of the scale $(\H^{s,2})_{|s|<\frac{1}{2}}$, the idea is to apply Proposition~\ref{Prop:Fractionalcomparison}.

We let $\tilde{s}\in(-\sfrac{1}{2},-1+\sfrac{(1+\alpha)}{2})$, by Theorems~\ref{thm:DirLapC1Domains}~\&~\ref{thm:StokesDirichletHsFracPower}, since the domain has at least $\C^{1}$-boundary (for this particular argument a Lipschitz boundary is sufficient), one has for $\varepsilon_0>0$, such that $\tilde{s}\pm\varepsilon_0\in(-\sfrac{1}{2},-1+\sfrac{(1+\alpha)}{2})$, one has
\begin{itemize}
    \item $\D_{2}^{\tilde{s}}((-\Delta_\mathcal{D})^{\pm\sfrac{\varepsilon_0}{2}}) = \H^{\tilde{s}\pm\varepsilon_0,2}(\Omega,\CC^n)$ on $\H^{\tilde{s},2}(\Omega,\CC^n)$; and
    \item $\D_{2}^{\tilde{s}}(\AA_\mathcal{D}^{\pm\sfrac{\varepsilon_0}{2}}) = \H^{\tilde{s}\pm\varepsilon_0,2}_{\mathfrak{n},\sigma}(\Omega)$ on $\H^{\tilde{s},2}_{\mathfrak{n},\sigma}(\Omega)$;
\end{itemize}
with the underlying isomorphisms. Thus, by Theorem~\ref{thm:SharpHodgeDecompC1}, for $\vv\in \D_{2}^{\tilde{s}}((-\Delta_\mathcal{D})^{\pm\sfrac{\varepsilon_0}{2}}) = \H^{\tilde{s}\pm\varepsilon_0,2}(\Omega,\CC^n)$, one has $\PP_\Omega\vv\in\H^{\tilde{s}\pm\varepsilon_0,2}_{\mathfrak{n},\sigma}(\Omega)$, and by the isomorphisms properties it holds
\begin{align*}
    \lVert \AA_\mathcal{D}^{\pm\sfrac{\varepsilon_0}{2}} \PP_\Omega\vv\rVert_{\H^{\tilde{s},2}(\Omega)}\sim_{\varepsilon_0,\Omega,\tilde{s}}\lVert  \PP_\Omega\vv\rVert_{\H^{\tilde{s}\pm\varepsilon_0,2}(\Omega)}&\lesssim_{\varepsilon_0,\Omega,\tilde{s}}\lVert  \vv\rVert_{\H^{\tilde{s}\pm\varepsilon_0,2}(\Omega)}\\
    &\lesssim_{\varepsilon_0,\Omega,\tilde{s}}\lVert (-\Delta_\mathcal{D})^{\pm\sfrac{\varepsilon_0}{2}}  \vv\rVert_{\H^{\tilde{s},2}(\Omega)}.
\end{align*}
Similarly, due to the natural embedding $\iota\,:\,\H^{\tilde{s}\pm\varepsilon_0,2}_{\mathfrak{n},\sigma}(\Omega)\longrightarrow\H^{\tilde{s}\pm\varepsilon_0,2}(\Omega,\CC^n)$, for all $\ww\in \D_{2}^{\tilde{s}}(\AA_\mathcal{D}^{\pm\sfrac{\varepsilon_0}{2}}) = \H^{\tilde{s}\pm\varepsilon_0,2}_{\mathfrak{n,\sigma}}(\Omega)$,
\begin{align*}
    \lVert (-\Delta_\mathcal{D})^{\pm\sfrac{\varepsilon_0}{2}}  \ww\rVert_{\H^{\tilde{s},2}(\Omega)}\sim_{\varepsilon_0,\Omega,\tilde{s}}\lVert  \ww\rVert_{\H^{\tilde{s}\pm\varepsilon_0,2}(\Omega)}&\lesssim_{\varepsilon_0,\Omega,\tilde{s}}\lVert \AA_\mathcal{D}^{\pm\sfrac{\varepsilon_0}{2}}  \ww\rVert_{\H^{\tilde{s},2}(\Omega)}.
\end{align*}
So, since \eqref{Eq: Retraction Condition} and \eqref{Eq: Injection Condition} are satisfied on $\X=\H^{\tilde{s},2}_{\mathfrak{n},\sigma}(\Omega)$, $\Y=\H^{\tilde{s},2}(\Omega,\CC^n)$, with $\alpha_0=\varepsilon_0$, Proposition~\ref{Prop:Fractionalcomparison} yielding
\begin{align}\label{eq:ProofBddRboundHs2}
     \sup_{1 \leqslant \tau , t \leqslant 2} \mathcal{R} \big\{ \varphi (\tau 2^{j + \ell} \AA_\mathcal{D}) \PP_\Omega \psi(-t 2^j\Delta_\mathcal{D} ) : j \in \ZZ \big\} \lesssim_{\mu,\varepsilon_0,\Omega} 2^{-\varepsilon_0|\ell|},
\end{align}
as a families of bounded operators from $\H^{\pm\tilde{s},2}(\Omega,\CC^n)$ to $\H^{\pm\tilde{s},2}_{\mathfrak{n},\sigma}(\Omega)$. Choosing $q\in(2,\infty)$ arbitrarily large and $s_0\in(-1+\sfrac{1}{q},-1+\sfrac{(1+\alpha)}{q})$, complex interpolation between \eqref{eq:ProofBddRboundHs2} on $\H^{\tilde{s},2},\H^{-\tilde{s},2}$ with \eqref{eq:ProofBddRboundHsp} on $\H^{{s}_0,q},\H^{-s_0,q'}$, yields the $\mathcal{R}$-bound
\begin{align}\label{eq:ProofBddRboundinterp}
     \sup_{1 \leqslant \tau , t \leqslant 2} \mathcal{R} \big\{ \varphi (\tau 2^{j + \ell} \AA_\mathcal{D}) \PP_\Omega \psi(-t 2^j\Delta_\mathcal{D} ) : j \in \ZZ \big\} \lesssim_{\mu,\varepsilon_0,\Omega} 2^{-\theta\varepsilon_0|\ell|},\qquad\text{for all}\qquad \ell\in\ZZ,
\end{align}
on $\H^{{s},p},\H^{-s,p'}$. Here, $\theta\in(0,1)$ is such that
\begin{align*}
    [\H^{{s}_0,q}(\Omega),\H^{\tilde{s},2}(\Omega)]_\theta = \H^{s,p}(\Omega).
\end{align*}
Finally, Theorem~\ref{Thm:Kunstmann-WeisHinfty} applies, and since the choice of parameters is arbitrary the boundedness of the $\mathbf{H}^\infty(\Sigma_\mu)$-functional calculus on $\H^{s,p}_{\mathfrak{n},\sigma}(\Omega)$ holds for all $p\in[2,\infty)$, and $s\in(-1+\sfrac{1}{p},-1+\sfrac{(1+\alpha)}{p})$, all $\mu\in(0,\pi)$.

\medbreak

\textbf{Step 2:} \textbf{We compute fractional powers  of the Stokes--Dirichlet operator} on $\H^{s,p}_{\mathfrak{n,\sigma}}(\Omega)$, $-1+\sfrac{1}{p}<s<-1+\sfrac{(1+\alpha)}{p}$, when $p\in[2,\infty)$, and deduce the result for all $1<p<\infty$.

We first assume $p\in[2,\infty)$, $s\in(-1+\sfrac{1}{p},-1+\sfrac{(1+\alpha)}{p})$, yielding the boundedness of the $\mathbf{H}^\infty(\Sigma_\mu)$-functional calculus for $\mu\in(0,{\pi})$ for $(\D^{s}_p(\AA_\mathcal{D}),\AA_\mathcal{D})$. In particular, it has BIP, so that for all $\theta\in(0,1)$
\begin{align*}
    \D^{s}_p(\AA_\mathcal{D}^\theta) = [\H^{s,p}_{\mathfrak{n},\sigma}(\Omega),\D^{s}_p(\AA_\mathcal{D})]_\theta.
\end{align*}
Taking into account the first point  of Remark~\ref{rem:AboutTheProof}, thanks to Proposition~\ref{prop:InterpDivFreeC1}, if $s+2\theta<1+\sfrac{1}{p}$, $s+2\theta\neq\sfrac{1}{p}$,
\begin{align*}
    \D^{s}_p(\AA_\mathcal{D}^\theta)=[\H^{s,p}_{\mathfrak{n},\sigma}(\Omega),\D^{s}_p(\AA_\mathcal{D})]_\theta = [\H^{s,p}_{\mathfrak{n},\sigma}(\Omega),\H^{s+2,p}_{\mathcal{D},\sigma}(\Omega)]_\theta = \H^{s+2\theta,p}_{\mathcal{D},\sigma}(\Omega).
\end{align*}
In particular, for all $s\in(-1+\sfrac{1}{p},1+\sfrac{(1+\alpha)}{p})$,
\begin{align}\label{eq:StokesSqrtBddDomains}
    \AA_\mathcal{D}^{\sfrac{1}{2}}&\,:\,\H^{s+2,p}_{\mathcal{D},\sigma}(\Omega)\longrightarrow\H^{s+1,p}_{\mathcal{D},\sigma}(\Omega),\\
    &\,:\,\H^{s+1,p}_{\mathcal{D},\sigma}(\Omega)\longrightarrow\H^{s,p}_{\mathfrak{n},\sigma}(\Omega)\nonumber
\end{align}
are isomophisms. Now, if $s+2\theta\in(1+\sfrac{1}{p},1+\sfrac{(1+\alpha)}{p})$, by invertibility of $\AA_\mathcal{D}$,
\begin{align*}
    \D^{s}_p(\AA_\mathcal{D}^{\theta})=\AA_\mathcal{D}^{-\sfrac{1}{2}}\D^{s}_p(\AA_\mathcal{D}^{\theta-\sfrac{1}{2}})= \AA_\mathcal{D}^{-\sfrac{1}{2}}\H^{s+2\theta-1,p}_{\mathcal{D},\sigma}(\Omega) = \H^{s+2\theta,p}_{\mathcal{D},\sigma}(\Omega),
\end{align*}
the last equality being supported by \eqref{eq:StokesSqrtBddDomains} and the fact that $s+2\theta-1\in(\sfrac{1}{p},\sfrac{(1+\alpha)}{p})$.

As a direct consequence, one can extrapolate $(\D_p(\AA_\mathcal{D}),\AA_\mathcal{D})$ on $\L^p_{\mathfrak{n},\sigma}(\Omega) = \AA_\mathcal{D}^{\sfrac{s}{2}}\H^{s,p}_{\mathfrak{n},\sigma}(\Omega)$ as a $0$-sectorial invertible operator with bounded $\mathbf{H}^\infty(\Sigma_\mu)$-functional calculus for $\mu\in(0,{\pi})$, and by consistency of the extrapolated operator, we have
\begin{align}\label{eq:FracPowersStokes}
    \D_p(\AA_\mathcal{D}^{\sfrac{\tau}{2}}) = \H^{\tau,p}_{\mathcal{D},\sigma}(\Omega),\qquad \text{for all}\quad\tau\in(-1+\sfrac{1}{p},1+\sfrac{(1+\alpha)}{p}),\,\,\text{all }\, p\in[2,\infty).
\end{align}

Since $(\D_2(\AA_\mathcal{D}),\AA_\mathcal{D})$ is self-adjoint, by density, and duality one obtains isomorphisms
\begin{align*}
    \AA_\mathcal{D}^{\sfrac{\tau}{2}}\,:\,\H^{\tau,p}_{\mathcal{D},\sigma}(\Omega)\longrightarrow\L^p_{\mathfrak{n},\sigma}(\Omega),
\end{align*}
for all $\tau \in(-1+\sfrac{1}{p},\sfrac{1}{p})$, all $p\in(1,2)$. One also obtains boundedness of the $\mathbf{H}^\infty(\Sigma_\mu)$-functional calculus for $\mu\in(0,\pi)$ on $\L^p_{\mathfrak{n},\sigma}(\Omega)$, $p\in(1,2)$. Starting again from the case $s\in(-1+\sfrac{1}{p},-1+\sfrac{(1+\alpha)}{p})$, as in the beginning of the present step, one can reproduce the arguments with Bounded Imaginary Powers including \eqref{eq:StokesSqrtBddDomains}, in order to deduce the validity of \eqref{eq:FracPowersStokes} for all $p\in(1,2)$ and all $\tau\in(-1+\sfrac{1}{p},1+\sfrac{(1+\alpha)}{p})$.

\textbf{Step 3:} \textbf{Description of the domain $\D^{s}_p(\AA_\mathcal{D})$ and expression of the operator $\AA_\mathcal{D}$.} By Theorem~\ref{thm:StokesResolvent}, the description as $\H^{s+1,p}_{0,\sigma}\cap\H^{s+2,p}(\Omega)$ is a direct consequence whenever $p\in(1,\infty)$ and $s>-1+\sfrac{1}{p}$ are satisfying \hyperref[eq:ConditionRegularity]{$(\prescript{\alpha}{r}{\ast}_{p,p}^{s})$}. In particular, the Hodge-Leray projection is well-defined as acting on  $\Delta \uu,\nabla\mathfrak{p}\in\H^{s,p}(\Omega,\CC^n)$ provided $\uu\in\D_{p}^s(\AA_\mathcal{D})$, setting
\begin{align*}
    \AA_\mathcal{D}\uu =-\Delta \uu+\nabla\mathfrak{p}=:\ff
\end{align*}
one obtains
\begin{align*}
    \AA_\mathcal{D}\uu =\ff = \PP_\Omega \ff = \PP_\Omega [-\Delta \uu+\nabla\mathfrak{p}] = \PP_\Omega [-\Delta \uu].
\end{align*}
This is legitimate since $\PP_\Omega$ acts properly on the elements $-\Delta \uu,\nabla\mathfrak{p}\in\H^{s,p}(\Omega,\CC^n)$, $s>-1+\sfrac{1}{p}$ by Theorem~\ref{thm:SharpHodgeDecompC1}.

The description in the remaining cases. Recall that previously  \hyperref[eq:ConditionRegularity]{$(\prescript{\alpha}{r}{\ast}_{p,p}^{s})$} was assumed. Now we just deal with the case $s=1+\sfrac{(1+\alpha)}{p}$. The remaining cases for larger $s$ admit a similar argument (and if $r\neq p$ then this corresponds to the case $p>r$, $s=1+\sfrac{(1+r\alpha)}{p}$, otherwise, one deals with $s=1+\alpha+\sfrac{1}{p}$, $p<r$). By the previous considerations on fractional powers, one has for all $s'<s$
\begin{align*}
    (\D^{s'}_p(\AA_\mathcal{D}),\AA_\mathcal{D})_{|_{\H^{s,p}_{\mathfrak{n},\sigma}(\Omega)}} = (\D^{s}_p(\AA_\mathcal{D}),\AA_\mathcal{D}),
\end{align*}
\textit{i.e.}
\begin{align*}
    \D^{s}_p(\AA_\mathcal{D})=\big\{\,\uu\in\D^{s'}_p(\AA_\mathcal{D})\,:\,\AA_\mathcal{D}\uu\in\H^{s,p}_{\mathfrak{n},\sigma}(\Omega)\,\big\}
\end{align*}
where the membership is an unbounded membership relation. However, due to the previous result on the characterisation of the domain of fractional powers
\begin{align*}
    \D^{s'}_p(\AA_\mathcal{D}) = \H^{s'+1,p}_{0,\sigma}\cap\H^{s'+2,p}(\Omega) \hookrightarrow \H^{1+\frac{1}{p}+\frac{\alpha}{2p},p}_{\mathcal{D},\sigma}(\Omega)=\D^{-1+\frac{1}{p}+\frac{\alpha}{2p}}_p(\AA_\mathcal{D}).
\end{align*}
This latter embedding is independent of the choice of $r$, and remains true even if $r>1$. Hence, it holds
\begin{align*}
    \D^{s}_p(\AA_\mathcal{D})\hookrightarrow\big\{\,\uu\in\H^{1+\frac{1}{p}+\frac{\alpha}{2p},p}_{\mathcal{D},\sigma}(\Omega)\,:\,\AA_\mathcal{D}\uu=\PP_{\Omega}(-\Delta_{\mathcal{D}} \uu)\in\H^{s,p}_{\mathfrak{n},\sigma}(\Omega)\,\big\}
\end{align*}
where $\PP_{\Omega}(-\Delta_{\mathcal{D}} \uu)$ makes sense (boundedly) in $\H^{-1+\frac{1}{p}+\frac{\alpha}{2p},p}(\Omega)$, and the membership relation to $\H^{s,p}_{\mathfrak{n},\sigma}(\Omega)$ is still an unbounded one in the latter description. Let us prove the reverse embedding. Let $\uu\in \H^{1+\frac{1}{p}+\frac{\alpha}{2p},p}_{\mathcal{D},\sigma}(\Omega)$, $\PP_\Omega(-\Delta \uu) \in \H^{-1+\frac{1}{p}+\frac{\alpha}{2p},p}_{\mathfrak{n},\sigma}\cap\H^{s,p}_{\mathfrak{n},\sigma}(\Omega)$, then, in particular, one obtains
\begin{align*}
    \ff:=\PP_\Omega(-\Delta \uu) =-\Delta\uu + \nabla \mathfrak{p} \in\H^{s,p}_{\mathfrak{n},\sigma}(\Omega),
\end{align*}
by the previously proved elliptic regularity results and uniqueness. Since $\H^{s,p}_{\mathfrak{n},\sigma}(\Omega)\subset \H^{s',p}_{\mathfrak{n},\sigma}(\Omega)$, for all $s'<s$, one obtains $\uu\in\H^{s'+1,p}_{0,\sigma}\cap\H^{s'+2,p}(\Omega)$. This yields the result for the description of the domain, and for the identity $\PP_\Omega(-\Delta)=\AA_\mathcal{D}$.
\end{proof}

\subsubsection{The results on Besov, H\"{o}lder and \texorpdfstring{$\L^1$}{L1}-type Sobolev spaces}\label{sec:BesovResolvResultsProof}

Now we present the full result for Besov spaces $\B^{s}_{p,q}$, for which we recall that it contains specific subclasses that are the Sobolev-Soblodeckij spaces $\W^{s,p}=\B^{s}_{p,p}$, $s\notin\ZZ$, $p\in[1,\infty)$, and the Lipschitz-H\"{o}lder spaces $\C^{k,\beta}=\B^{k+\beta}_{\infty,\infty}$, $k\in\NN$, $\beta\in(0,1)$.

We provide a detailed proof primarily for Theorem~\ref{thm:FinalResultBesov2}, and we will briefly explain the particular case \(p = r = \infty\), which yields Theorem~\ref{thm:FinalResultBesov1}. The remaining results follow as straightforward consequences of Theorem~\ref{thm:FinalResultBesov2}. When \(r < \infty\) and \(p = \infty\), note that the condition \hyperref[eq:ConditionRegularity]{$(\prescript{\alpha}{r}{\ast}_{p,q}^s)$} becomes empty whenever $s\in(-1+\sfrac{1}{p},\sfrac{1}{p})$.

As the reader will see, the proof simplifies significantly compared to its counterpart on $\L^p$ when it comes to identify the domain of fractional powers. This is due to the fact that the boundedness of the $\mathbf{H}^\infty$-functional calculus is immediate, thanks to an identification of the real interpolation scale associated with the operator (by Dore's Theorem \cite[Theorem~6.1.3]{bookHaase2006}), and the stability of Besov spaces under real interpolation, see the end of Remark~\ref{rem:FuncCalcSimplAbsFC} for additional insights. Therefore, the following results are essentially straightforward consequences of Theorem~\ref{thm:StokesResolvent}. Nevertheless, for the sake of completeness, the arguments to compute the domain of fractional powers are still provided.

\begin{theorem}\label{thm:FinalResultBesov0}Let $\alpha\in(0,1)$, $r\in[1,\infty]$, and let $\Omega$ be a bounded $\C^1$-domain of class $\mathcal{M}^{1+\alpha,r}_\W(\epsilon)$ for some sufficiently small $\epsilon>0$.

\medbreak

\noindent Let $p,q\in[1,\infty]$. Let $s\in(-1+{\sfrac{1}{p}},\sfrac{1}{p})$. For $\mu\in[0,\pi)$, for all $\lambda\in\Sigma_{\mu}\cup\{0\}$, for all $f\in\B^{s,\sigma}_{p,q,\mathfrak{n}}(\Omega)$, 
\begin{itemize}
    \item If \eqref{eq:ConditionRegularity} is satisfied: the resolvent problem \eqref{eq:Stokes} admits a unique solution $$(\uu,\mathfrak{p})\in\B^{s+2}_{p,q}(\Omega,\CC^n)\times\B^{s+1}_{p,q,\mfree}(\Omega),$$ which also satisfies the following estimate
\begin{align*}
     (1+\lvert\lambda\rvert) \lVert \uu\rVert_{\B^{s}_{p,q}(\Omega)}+\lVert ( \nabla^2 \uu, \nabla \mathfrak{p})\rVert_{\B^{s}_{p,q}(\Omega)} \nonumber \lesssim_{p,n,s,\mu}^{\Omega} \lVert \ff\rVert_{\B^{s}_{p,q}(\Omega)};
\end{align*}
    \item Otherwise the resolvent problem \eqref{eq:Stokes} admits a unique solution $(\uu,\mathfrak{p})$ satisfying the membership 
    $$(\uu,\mathfrak{p})\in\B^{1+\beta}_{p,q}(\Omega,\CC^n)\times\B^{\beta}_{p,q,\mfree}(\Omega),$$
    for all $\beta\in[0,1+s]$ such that \hyperref[eq:ConditionRegularity]{$(\prescript{\alpha}{r}{\ast}_{p,q}^{\beta-1})$} and that satisfies the following estimate for all such $\beta$:
\begin{align*}
     (1+\lvert\lambda\rvert) \lVert \uu\rVert_{\B^{s}_{p,q}(\Omega)}+(1+\lvert\lambda\rvert)^\frac{s+1-\beta}{2}\lVert \nabla \uu\rVert_{\B^{\beta}_{p,q}(\Omega)} + \lVert \mathfrak{p}\rVert_{\B^{\beta}_{p,q}(\Omega)} \nonumber \lesssim_{p,n,s,\mu}^{\Omega,\beta} \lVert \ff\rVert_{\B^{s}_{p,q}(\Omega)}.
\end{align*}
\end{itemize}
\end{theorem}

From the previous theorem, we highlight the following particular case:

\begin{theorem}\label{thm:FinalResultBesov1}Let $\alpha\in(0,1)$, and consider $\Omega$ to be a bounded $\mathcal{C}^{1,\alpha}$-domain or bounded $\C^{1,\alpha}$-domain with small $(1,\alpha)$-H\"{o}lderian constant.

For $\mu\in[0,\pi)$, for $\beta\in(0,\alpha]$, for all $\lambda\in\Sigma_{\mu}\cup\{0\}$,  for all $f\in\B^{-1+\beta,\sigma}_{\infty,\infty,\mathfrak{n}}(\Omega)$, the resolvent problem \eqref{eq:Stokes} admits a unique solution $(\uu,\mathfrak{p})\in\C^{1,\beta}(\overline{\Omega},\CC^n)\times\C^{0,\beta}_{\mfree}(\overline{\Omega})$ satisfying the estimate
\begin{align*}
     (1+\lvert\lambda\rvert) \lVert \uu\rVert_{\B^{-1+\beta}_{\infty,\infty}(\Omega)}+\lVert (\nabla \uu, \mathfrak{p} )\rVert_{\B^{\beta}_{\infty,\infty}(\Omega)} \nonumber \lesssim_{n,\alpha,\mu}^{\Omega} \lVert \ff\rVert_{\B^{-1+\beta}_{\infty,\infty}(\Omega)}.
\end{align*}
If additionally $\ff\in\BesSmo^{-1+\beta,\sigma}_{\infty,\infty,\mathfrak{n}}(\Omega)$, then it holds that $\uu\in \mathcal{C}^{1,\beta}(\overline{\Omega},\CC^n)$, whenever $\beta<\alpha$, or  if $\beta=\alpha$ and $\Omega$ has a $\mathcal{C}^{1,\alpha}$-boundary.
\end{theorem}

\begin{theorem}\label{thm:FinalResultBesov2} Let $\alpha\in(0,1)$, $r\in[1,\infty]$, and let $\Omega$ be a bounded $\C^1$-domain of class $\mathcal{M}^{1+\alpha,r}_\W(\epsilon)$ for some sufficiently small $\epsilon>0$.

\medbreak

Let $p\in[1,\infty]$. Then it holds that
\begin{enumerate}
    \item $\AA_{\mathcal{D}}\,:\,\B^{s+2,\sigma}_{p,q,\mathcal{D}}(\Omega)\longrightarrow\B^{s,\sigma}_{p,q,\mathfrak{n}}(\Omega)$, is an isomorphism for all $q\in[1,\infty]$, $s>-1+\sfrac{1}{p}$ whenever it is satisfied \hyperref[eq:ConditionRegularity]{$(\prescript{\alpha}{r}{\ast}_{p,q}^s)$};
    \item For all $q\in[1,\infty]$, $-1+\sfrac{1}{p}<s<\sfrac{1}{p}$, $(\D_{p,q}^{s}(\AA_\mathcal{D}),\AA_\mathcal{D})$ is an invertible $0$-sectorial on $\B^{s,\sigma}_{p,q,\mathfrak{n}}(\Omega)$ which admits bounded $\mathbf{H}^{\infty}(\Sigma_\theta)$-functional calculus for all $\theta\in(0,\pi)$;
    \item For all $s\in(-1+\sfrac{1}{p},\sfrac{1}{p})$, all $\beta\in(0,2)$, such that  \hyperref[eq:ConditionRegularity]{$(\prescript{\alpha}{r}{\ast}_{p,q}^{s+\beta-2})$}, $s+\beta\neq\sfrac{1}{p}$, one has with equivalence of norms, for all $q\in[1,\infty]$,
    \begin{align*}
        \D^s_{p,q}(\AA_\mathcal{D}^{\sfrac{\beta}{2}}) =  \B^{s+\beta,\sigma}_{p,q,\mathcal{D}}(\Omega).
  \end{align*}
  \item If $\Omega$ is a bounded domain with $\mathcal{C}^{1,\alpha}$ boundary, or $\C^{1,\alpha}$ boundary with small $(1,\alpha)$-H\"{o}lderian,  constant then it holds
  \begin{align*}
      \D_{\infty,\infty}^{\alpha-1}(\AA_\mathcal{D}) =  \B^{1+\alpha,\sigma}_{\infty,\infty,\mathcal{D}}(\Omega) = \C^{1,\alpha}_{\mathcal{D},\sigma}(\overline{\Omega}).
  \end{align*}
\end{enumerate}
If $q<\infty$, then the Stokes--Dirichlet operator is densely defined. The whole result still holds for the Besov spaces $\BesSmo^{s}_{p,\infty}$ up to appropriate modifications.
\end{theorem}

For the next result --being an immediate consequence by setting $p=q=1$--, one can assume without loss of generality that $r=1$, since $\mathcal{M}^{1+\alpha,1}_{\W}$ is a lower regularity condition for the boundary than any $\mathcal{M}^{1+\alpha,r}_{\W}$-condition, $r\in(1,\infty]$.

\begin{corollary}\label{cor:StokesWs1}Let $\alpha\in(0,1)$, and consider $\Omega$ to be bounded $\C^1$-domain of class $\mathcal{M}^{1+\alpha,1}_{\W}(\epsilon)$ for some sufficiently small $\epsilon>0$. 

\medbreak

Then it holds that
\begin{enumerate}
    \item $\AA_{\mathcal{D}}\,:\,\W^{1,1}_{0,\sigma}(\Omega)\cap\W^{s+2,1}(\Omega)\longrightarrow\W^{s,1}_{\mathfrak{n},\sigma}(\Omega)$, is an isomorphism for all $s\in(0,\alpha)$;
    \item For all $s\in(0,1)$, $(\D_{1,1}^{s}(\AA_\mathcal{D}),\AA_\mathcal{D})$ is an invertible $0$-sectorial operator on $\W^{s,1}_{\mathfrak{n},\sigma}(\Omega)$ which admits bounded $\mathbf{H}^{\infty}(\Sigma_\theta)$-functional calculus for all $\theta\in(0,\pi)$;
    \item For all $s\in(0,1)$, all $\beta\in(0,2)$, such that $s+\beta\leqslant 2+\alpha$, $s+\beta\neq 1$, one has with equivalence of norms
    \begin{align*}
        \D^s_{1,1}(\AA_\mathcal{D}^{\sfrac{\beta}{2}}) =  \begin{cases}
            \W^{s+\beta,1}_{\mathfrak{n},\sigma}(\Omega),\quad &\,\textrm{if }s+\beta<1;\\
            \W^{1,1}_{0,\sigma}(\Omega)\cap\W^{s+\beta,1}(\Omega),\quad &\,\textrm{if }s+\beta>1.
        \end{cases}
  \end{align*}
\end{enumerate}
\end{corollary}

\begin{proof} \textbf{Step 1:} the case $r<\infty$. By Theorem~\ref{thm:StokesResolvent}, when $p<\infty$, the Stokes operator is an invertible $0$-sectorial operator on $\B^{s,\sigma}_{p,q,\mathfrak{n}}(\Omega)$ whenever $s\in(-1+\sfrac{1}{p},\sfrac{1}{p})$, $p$ and $q$ are satisfying \hyperref[eq:ConditionRegularity]{$(\prescript{\alpha}{r}{\ast}_{p,q}^s)$}, with domain
\begin{align*}
    \D^{s}_{p,q}(\AA_\mathcal{D}) = \B^{s+2,\sigma}_{p,q,\mathcal{D}}(\Omega).
\end{align*}
Consequently, by the Dore Theorem \cite[Theorem~6.1.3]{bookHaase2006}, and real interpolation Proposition~\ref{prop:InterpDivFreeC1}, one obtains that $\AA_\mathcal{D}$ remains an invertible $0$-sectorial operator with  bounded $\mathbf{H}^{\infty}(\Sigma_\theta)$-functional calculus, for all $\theta\in(0,\pi)$, on the Besov spaces
\begin{align*}
    \B^{s,\sigma}_{p,q,\mathfrak{n}}(\Omega), \quad \text{ for all }s\in(-1+\sfrac{1}{p},\sfrac{1}{p}),\, p\in[1,\infty),q\in[1,\infty].
\end{align*}
The case $p=\infty$ is obtained by duality and self-adjointness.
Moreover, by the bounded holomorphic functional calculus,  \cite[Proposition~6.4.1,~a),~\&~Corollary~6.5.5]{bookHaase2006} and Proposition~\ref{prop:InterpDivFreeC1}, for all $s\in(-1+\sfrac{1}{p},\sfrac{1}{p})$ such that the condition \hyperref[eq:ConditionRegularity]{$(\prescript{\alpha}{r}{\ast}_{p,q}^s)$}  is satisfied, the obtained isomorphism 
\begin{align*}
    \AA_\mathcal{D}^\frac{\beta}{2}\,:\,(\B^{s,\sigma}_{p,q,\mathfrak{n}}(\Omega),\D^{s}_{p,q}(\AA_\mathcal{D}))_{\theta,q}\longrightarrow\,(\B^{s,\sigma}_{p,q,\mathfrak{n}}(\Omega),\D^{s}_{p,q}(\AA_\mathcal{D}))_{\theta-\frac{\beta}{2},q}
\end{align*}
translates into an isomorphism
\begin{align*}
    \AA_\mathcal{D}^\frac{\beta}{2}\,:\,\B^{s+2\theta,\sigma}_{p,q,\mathcal{D}}(\Omega)\longrightarrow\,\B^{s+2\theta-\beta,\sigma}_{p,q,\mathcal{D}}(\Omega)
\end{align*}
whenever $s+2\theta<1+\sfrac{1}{p}$. The choices of $\beta\in(0,2)$ and $s\in(-1+\sfrac{1}{p},\sfrac{1}{p}$) such that the condition \hyperref[eq:ConditionRegularity]{$(\prescript{\alpha}{r}{\ast}_{p,q}^s)$} is satisfied are arbitrary. Therefore, Point \textit{(iii)} from Theorem~\ref{thm:FinalResultBesov2} holds for all $s\in(-1+\sfrac{1}{p},\sfrac{1}{p})$, $\beta\in (0,2)$ such that $s+\beta<1+\sfrac{1}{p}$.

\medbreak

\noindent Now, let $s\in(-1+\sfrac{1}{p},\sfrac{1}{p})$ and $\beta\in (0,2)$, such that $s+\beta>1+\sfrac{1}{p}$ also satisfies \hyperref[eq:ConditionRegularity]{$(\prescript{\alpha}{r}{\ast}_{p,q}^{s+\beta-2})$}. One has $\sfrac{1}{p}<s+\beta-1<1+\sfrac{1}{p}$, and by previous considerations, one also has an isomorphism
\begin{align*}
    \AA_\mathcal{D}^{-\frac{1}{2}}\,:\,\B^{s+\beta-1,\sigma}_{p,q,\mathcal{D}}(\Omega)\longrightarrow\,\B^{s+\beta,\sigma}_{p,q,\mathcal{D}}(\Omega),
\end{align*}
Therefore, writing $\AA_\mathcal{D}^{\frac{\beta}{2}} = \AA_\mathcal{D}^{\frac{\beta+1}{2}}\AA_\mathcal{D}^{-\frac{1}{2}}$, yields Point \textit{(iii)} of Theorem~\ref{thm:FinalResultBesov2} in any cases.

\textbf{Step 2:} The case $r=\infty$. The proof remain identical, but then since $\mathcal{M}^{1+\alpha,\infty}_{\W}=\B^{1+\alpha}_{\infty,\infty}=\C^{1,\alpha}$, the smallness assumption on the multiplier norm $\mathcal{M}^{1+\alpha,\infty}_{\W}(\epsilon)$ for $\epsilon>0$ simply becomes a smallness assumption on the $(1,\alpha)$-H\"{o}lderian constant of the charts describing the locally the boundary. That is, if $\varphi\,:\,\RR^{n-1}\longrightarrow\RR$ describes locally the boundary of $\Omega$, it can be chosen so that it satisfies
\begin{align*}
    \lVert\nabla' \varphi\rVert_{\L^\infty}+\sup_{\substack{x',y'\in\RR^{n-1}\\x'\neq y'}}\frac{|\nabla'\varphi(x')-\nabla'\varphi(y')|}{|x'-y'|^\alpha}<\epsilon.
\end{align*}
Note that if the boundary lies in the little H\"{o}lder class $\mathcal{C}^{1,\alpha}$, the choices for $\varphi$ can be made in such way that this condition is always satisfied for any fixed, but arbitrarily small, $\epsilon>0$. See Sections~\ref{sec:para}~and~\ref{Subsec:SmallMultNorm}, as well as Proposition~\ref{prop:MspepsClassDomainsSufficientCdtn} for more details.

\noindent Note also that the condition \hyperref[eq:ConditionRegularity]{$(\prescript{\alpha}{\infty}{\ast}_{\infty,\infty}^s)$} just translates into
\begin{align*}
    s \leqslant -1+\alpha.
\end{align*}
Consequently, all arguments in the previous step remain valid, and the result follows.
\end{proof}

\subsection{Other instances of the main results}

Here, we specify the result in the cases that could be of specific interests. Especially, a study of the precise regularity properties of the Stokes operator in bounded $\C^{1,\alpha}$-domains, which seems still be an open question up to now.  As a by-product of our analysis, we also derive several consequences for the pure $\L^1$-theory of the Stokes--Dirichlet operator in such domains.

\medbreak

Furthermore, we also provide the arguments that allow to remove the regularity condition of a $\C^1$-boundary. From Section~\ref{Sec:RemoveC1}, one deduces that the result for the $\L^\infty$-theory below, Theorem~\ref{thm:FinalResultLinfty2}~and~\ref{thm:FinalResultLinfty1}, are the only results that depend on Geng and Shen’s landmark work on the Dirichlet--Stokes operator \cite{GengShen2024,GengShen2025}, and more generally of the existing literature about the Stokes operator on bounded domains.

\subsubsection{The particular case of the \texorpdfstring{$\L^p$}{Lp}-theory, \texorpdfstring{$1\leqslant p\leqslant \infty$}{1<p<oo}, for bounded \texorpdfstring{$\C^{1,\alpha}$}{C1a}-domains.}\label{sec:LpTheoryC1a}

We precise our previous theorems and result in the particular case of bounded $\C^{1,\alpha}$-domains. Note that by Proposition~\ref{prop:MspepsClassDomainsSufficientCdtn}, any $\C^{1,\alpha}$-bounded domain is a domain of class $\mathcal{M}^{1+\beta,p}_{\W}(\epsilon)$ for any $\beta\in(0,\alpha)$, and any $\epsilon>0$. We would like to bring to the attention of the reader that as classes of domains, one has for any $\varepsilon>0$
\begin{align*}
    \C^{1,1-\sfrac{1}{p}+\varepsilon}\subsetneq\B^{2-\sfrac{1}{p}}_{\infty,p}\subsetneq\C^{1,1-\sfrac{1}{p}}
\end{align*}
while having $\mathbf{B}^{2-\sfrac{1}{p}}_{\infty,p}\subset\mathcal{M}^{2-\sfrac{1}{p},p}_{\W}(\epsilon)$, for any $\epsilon>0$, see Proposition~\ref{prop:MspepsClassDomainsSufficientCdtn}.

\noindent For bounded $\C^{1,\alpha}$-domains, $\alpha\in(0,1]$, having $s>-1+\sfrac{1}{p}$ to satisfy the condition \hyperref[eq:ConditionRegularity]{$(\prescript{\beta}{p}{\ast}_{p,p}^{s})$} for any $\beta\in(0,\alpha)$ is equivalent to 
\begin{align*}
    s<-1+\alpha+\sfrac{1}{p}.
\end{align*}

\noindent For a domain in the class  $\mathbf{B}^{2-\sfrac{1}{p}}_{\infty,p}\subset\mathcal{M}^{2-\sfrac{1}{p},p}_{\W}(\epsilon)$, \hyperref[eq:ConditionRegularity]{$\big(\prescript{1-\frac{1}{p}}{p}{\ast}_{p,p}^{s}\big)$} becomes
\begin{align*}
    s\leqslant 0.
\end{align*}

\noindent Consequently, Theorem~\ref{thm:FinalResultSobolev1} for the $\L^p$-theory can be restated and simplified as follows

\begin{theorem}\label{thm:FinalResultSobolevC1q}   Let $\alpha\in(0,1]$ and let $\Omega$ be a bounded $\C^{1,\alpha}$-domain. For any $p\in(1,\infty)$, it holds that
\begin{enumerate}
    \item $\AA_{\mathcal{D}}\,:\,\H^{s+2,p}_{\mathcal{D},\sigma}(\Omega)\longrightarrow\H^{s,p}_{\mathfrak{n},\sigma}(\Omega)$, is an isomorphism for all $s\in(-1+\sfrac{1}{p},-1+\alpha+\sfrac{1}{p})$.
    
    \item The operator $(\D_p(\AA_\mathcal{D}),\AA_\mathcal{D})$ is a densely defined, invertible and $0$-sectorial on $\L^p_{\mathfrak{n},\sigma}(\Omega)$, and it admits bounded $\mathbf{H}^{\infty}(\Sigma_\theta)$-functional calculus for all $\theta\in(0,\pi)$;
    
    \item For all $s\in(-1+\sfrac{1}{p},1+\alpha+\sfrac{1}{p})$, $s\neq\sfrac{1}{p},1+\sfrac{1}{p}$, one has with equivalence of norms
    \begin{align*}
        \D_p(\AA_\mathcal{D}^{\sfrac{s}{2}}) =  \H^{s,p}_{\mathcal{D},\sigma}(\Omega).
  \end{align*}

    \item \label{StatementW2pReg} If additionally $\Omega$ is a bounded domain of class $\mathbf{B}^{2-\sfrac{1}{p}}_{\infty,p}$, then it holds
    \begin{align*}
        \D_p(\AA_\mathcal{D}) = \W^{1,p}_{0,\sigma}(\Omega)\cap\W^{2,p}(\Omega).
    \end{align*}
\end{enumerate}
\end{theorem}
Concerning Point \ref{StatementW2pReg} only, if $p\in(1,\infty)$ is arbitrary,  $\B^{2-\sfrac{1}{p}}_{\infty,p}$-regularity condition for the boundary is known to be sharp in the case of the $\W^{2,p}$--regularity of the Dirichlet Laplacian on $\L^p$ spaces, see \cite[Chapter~14,~Section~14.6,~p.513--514]{MazyaShaposhnikova2009}. The same conclusion for the $\W^{2,p}$--regularity on $\L^p$ for the Stokes--Dirichlet operator, can also be derived from a  very recent result in the 3-dimensional case by the first author \cite[Theorem~3.1]{Breit2025}. If $p>n$,  Point \ref{StatementW2pReg} still holds for $\W^{2-\sfrac{1}{p},p}$-regularity of the boundary.

\medbreak

Before this, the best known result for $\W^{2,p}$-regularity on $\L^p$ goes back to Solonnikov \cite{Solonnikov1977}, for which it was required to have a $\C^{1,1}$-domain, which yields for us a gain of $\frac{1}{p}$-derivative for the boundary regularity, and even slightly more if one considers a boundary in the Besov or Sobolev class whenever one has a high Lebesgue index, $p>n$, see Proposition~\ref{prop:MspepsClassDomainsSufficientCdtn}.

\medbreak

Now, our previous analysis in Section~\ref{sec:BesovResolvResultsProof}, Step 2 in the proof, in combination with a recent work of Geng and Shen \cite[Theorems~1.1~\&~9.4]{GengShen2025} allows to obtain the sharpest result known up to now concerning the $\L^\infty$-theory. We recall here for convenience that any bounded  $\C^{1,\alpha}$-domain is also a $\mathcal{C}^{1,\beta}$-domain for any $\beta\in(0,\alpha)$.

\begin{theorem}\label{thm:FinalResultLinfty2} Let $\alpha\in(0,1)$, and consider $\Omega$ to be a $\mathcal{C}^{1,\alpha}$-bounded domain or $\C^{1,\alpha}$-bounded domain with small $(1,\alpha)$-H\"{o}lderian constant. Then it holds that
\begin{enumerate}
    \item $(\D_{\infty}(\AA_\mathcal{D}),\AA_\mathcal{D})$ is an invertible $0$-sectorial operator on $\L^\infty_{\mathfrak{n},\sigma}(\Omega)$ and restricts as such on $\C^{0}_{0,0,\sigma}({\Omega})$;
    \item For all $\beta\in[0,\alpha]$, it holds
    \begin{align*}
      \D_{\infty}(\AA_\mathcal{D}^\frac{1+\beta}{2}) \hookrightarrow  \B^{1+\beta,\sigma}_{\infty,\infty,\mathcal{D}}(\Omega) = \C^{1,\beta}_{\mathcal{D},\sigma}(\overline{\Omega}).
  \end{align*}
    In particular, for all $\uu\in\D_{\infty}(\AA_\mathcal{D})$, one has $\uu\in\C^{1,\beta}(\overline{\Omega})$ and an estimate 
    \begin{align*}
        (1+|\lambda|)^\frac{1-\beta}{2}\lVert \uu \rVert_{\B^{1,\beta}_{\infty,\infty}(\Omega)} \lesssim_{s,n,\Omega}^{\alpha,\beta,\mu} \lVert \lambda\uu +\AA_{\mathcal{D}}\uu\rVert_{\L^\infty(\Omega)}
    \end{align*}
    valid for  all $\lambda\in\Sigma_{\mu}\cup\{0\}$, $\mu\in[0,\pi)$ being fixed.
    \item One has the following domain description and expression for Dirichlet--Stokes operator:
    \begin{align*}
        \D_{\infty}(\AA_\mathcal{D}):=\{\,\uu\in\C^{1,\alpha}_{\mathcal{D},\sigma}(\overline{\Omega})\,:\,\PP_{\Omega}(-\Delta_\mathcal{D}\uu)\in\L^\infty_{\mathfrak{n},\sigma}(\Omega)\,\},
    \end{align*}
    and
    \begin{align*}
        \AA_\mathcal{D}\uu = \PP_{\Omega}(-\Delta_\mathcal{D}\uu),\quad \text{ for all }\uu\in\D_{\infty}(\AA_{\mathcal{D}}).
    \end{align*}
\end{enumerate}
Furthermore, the Stokes--Dirichlet operator is densely defined on  $\C^{0}_{0,0,\sigma}({\Omega})$.
\end{theorem}

In terms of solvability of \eqref{eq:Stokes}, one also obtains the next result.

\begin{theorem}\label{thm:FinalResultLinfty1}Let $\alpha\in(0,1)$, and consider $\Omega$ to be a $\mathcal{C}^{1,\alpha}$-bounded domain or $\C^{1,\alpha}$-bounded domain with small $(1,\alpha)$-H\"{o}lderian constant.

For $\mu\in[0,\pi)$, for all $\lambda\in\Sigma_{\mu}\cup\{0\}$, for all $f\in\L^\infty_{\mathfrak{n},\sigma}(\Omega)$, the resolvent problem \eqref{eq:Stokes} admits a unique solution $(\uu,\mathfrak{p})\in\C^{1,\alpha}(\overline{\Omega},\CC^n)\times\C^{0,\alpha}_{\mfree}(\overline{\Omega})$ that also satisfies the estimate
\begin{align*}
     (1+\lvert\lambda\rvert) \lVert \uu\rVert_{\L^\infty(\Omega)}+(1+\lvert\lambda\rvert)^\frac{1-\alpha}{2}\lVert ( \nabla \uu, \mathfrak{p} )\rVert_{\B^{\alpha}_{\infty,\infty}(\Omega)} \nonumber \lesssim_{n,\alpha,\mu}^{\Omega} \lVert \ff\rVert_{\L^\infty(\Omega)}.
\end{align*}
If additionally $\ff\in\C^{0}_{0,0,\sigma}(\Omega)$, then it holds that $\uu\in \mathcal{C}^{1,\alpha}(\overline{\Omega},\CC^n)$ whenever $\Omega$ is of class $\mathcal{C}^{1,\alpha}$.
\end{theorem}

The theorems above, Theorem~\ref{thm:FinalResultLinfty2}~\&~\ref{thm:FinalResultLinfty1}, improve the known results by Abe and Giga \cite[Theorem~1.2]{AbeGiga2013}, and Abe, Giga and Hieber \cite[Theorem~1.1]{AbeGigaHieber2015}, where bounded $\C^3$-domains and bounded $\C^2$-domains were considered respectively. Only Lipschitz bound can be derived from the resolvent estimate of the former. Note that our result does not rely on the notion of admissible domains as introduced by Abe and Giga. Finally, we require less regularity for the boundary than for their results and we obtain a precise description of the exact regularity of the solution that matches the regularity of the boundary.

\begin{proof}[of Theorems~\ref{thm:FinalResultLinfty2}~\&~\ref{thm:FinalResultLinfty1}] We start from Theorem~\ref{thm:FinalResultBesov1}. For $s\in(-1,-1+\alpha]$, by real interpolation, $\AA_\mathcal{D}$ remains an invertible $0$-sectorial operator on the space
\begin{align*}
    (\BesSmo^{s,\sigma}_{\infty,\infty,\mathfrak{n}}(\Omega),\mathcal{D}^{s}_{\infty,\infty}(\AA_\mathcal{D}))_{-\frac{s}{2},\infty}.
\end{align*}
So, we have in particular the following resolvent estimate, by setting $\uu:=(\lambda\I+\AA_\mathcal{D})^{-1}\ff$ for $\ff\in(\BesSmo^{s,\sigma}_{\infty,\infty,\mathfrak{n}}(\Omega),\mathcal{D}^{s}_{\infty,\infty}(\AA_\mathcal{D}))_{-\frac{s}{2},\infty}$, for all $\lambda\in\Sigma_\mu\cup\{0\}$, $\mu\in(0,\pi)$:
\begin{align*}
   &(1+|\lambda|)\lVert \uu \rVert_{(\BesSmo^{s,\sigma}_{\infty,\infty,\mathfrak{n}}(\Omega),\mathcal{D}^{s}_{\infty,\infty}(\AA_\mathcal{D}))_{-\frac{s}{2},\infty}}\\ &\qquad+  (1+|\lambda|)^\frac{1-\alpha}{2}\lVert \AA_\mathcal{D}^{\frac{1+\alpha}{2}} \uu \rVert_{(\BesSmo^{s,\sigma}_{\infty,\infty,\mathfrak{n}}(\Omega),\mathcal{D}^{s}_{\infty,\infty}(\AA_\mathcal{D}))_{-\frac{s}{2},\infty}} \lesssim_{\mu,\Omega} \lVert \ff \rVert_{(\BesSmo^{s,\sigma}_{\infty,\infty,\mathfrak{n}}(\Omega),\mathcal{D}^{s}_{\infty,\infty}(\AA_\mathcal{D}))_{-\frac{s}{2},\infty}}.
\end{align*}

From this point, due to Propositions~\ref{prop:IdentifVanishingDivFree}~and~\ref{prop:InterpDivFreeC1}, Theorems~\ref{thm:InterpHomSpacesBddLip}~and~\ref{thm:FinalResultBesov1} as well as\cite[Proposition~6.4.1,~a),~\&~Corollary~6.5.5]{bookHaase2006}, one can canonically identify with equivalence of norms
\begin{align*}
    (\BesSmo^{s,\sigma}_{\infty,\infty,\mathfrak{n}}(\Omega),\mathcal{D}^{s}_{\infty,\infty}(\AA_\mathcal{D}))_{-\frac{s}{2},\infty} = \B^{0,\sigma}_{\infty,\infty,0}(\Omega).
\end{align*} 
Hence if $\ff\in\L^\infty_{\mathfrak{n},\sigma}(\Omega)$, by Proposition~\ref{prop:IdentifVanishingDivFree}, the extension by $0$ to the whole $\RR^n$, implies that we have a canonical embedding
\begin{align*}
    \L^\infty_{\mathfrak{n},\sigma}(\Omega) \hookrightarrow\B^{0,\sigma}_{\infty,\infty,0}(\Omega),
\end{align*}
yielding
\begin{align*}
    (1+|\lambda|)^\frac{1-\alpha}{2}\lVert \AA_\mathcal{D}^{\frac{1+\alpha}{2}} \uu \rVert_{(\BesSmo^{s,\sigma}_{\infty,\infty,\mathfrak{n}}(\Omega),\mathcal{D}^{s}_{\infty,\infty}(\AA_\mathcal{D}))_{-\frac{s}{2},\infty}} \lesssim_{\mu,\Omega} \lVert \ff \rVert_{(\BesSmo^{s,\sigma}_{\infty,\infty,\mathfrak{n}}(\Omega),\mathcal{D}^{s}_{\infty,\infty}(\AA_\mathcal{D}))_{-\frac{s}{2},\infty}}\lesssim_{\Omega,\mu} \lVert \ff \rVert_{\L^\infty(\Omega)}.
\end{align*}
By \cite[Proposition~6.4.1,~a),~\&~Corollary~6.5.5]{bookHaase2006}, Theorem~\ref{thm:FinalResultBesov1} and Proposition~\ref{prop:InterpDivFreeC1}, one obtains
\begin{align*}
    \lVert \AA_\mathcal{D}^{\frac{1+\alpha}{2}} \uu \rVert_{(\BesSmo^{s,\sigma}_{\infty,\infty,\mathfrak{n}}(\Omega),\mathcal{D}^{s}_{\infty,\infty}(\AA_\mathcal{D}))_{-\frac{s}{2},\infty}}&\sim_{\varepsilon,\Omega} \lVert \AA_\mathcal{D}^{\frac{1+\alpha+\varepsilon}{2}} \uu \rVert_{(\BesSmo^{s,\sigma}_{\infty,\infty,\mathfrak{n}}(\Omega),\mathcal{D}^{s}_{\infty,\infty}(\AA_\mathcal{D}))_{-\frac{s}{2}-\varepsilon,\infty}}\\
    &\sim_{\varepsilon,\Omega,\alpha} \lVert \AA_\mathcal{D}^{\frac{1+\alpha+\varepsilon}{2}} \uu \rVert_{\B^{-\varepsilon}_{\infty,\infty}(\Omega)}\\
    &\sim_{\varepsilon,\Omega,\alpha} \lVert \uu \rVert_{\B^{1+\alpha}_{\infty,\infty}(\Omega)}.
\end{align*}
This yields the estimate
\begin{align*}
   (1+|\lambda|)^\frac{1-\alpha}{2} \lVert \uu \rVert_{\B^{1+\alpha}_{\infty,\infty}(\Omega)}\lesssim_{\Omega,\mu} \lVert \ff \rVert_{\L^\infty(\Omega)}.
\end{align*}
Since $\ff\in\L^\infty_{\mathfrak{n},\sigma}(\Omega)\hookrightarrow\B^{-1+\alpha,\sigma}_{\infty,\infty,\mathfrak{n}}(\Omega)$, by Theorem~\ref{thm:FinalResultBesov1}, one has
\begin{align*}
    \lVert \uu \rVert_{\B^{1+\alpha}_{\infty,\infty}(\Omega)}\sim_{\alpha,n,\Omega} \lVert (\nabla^2 \uu, \nabla \mathfrak{p}) \rVert_{\B^{-1+\alpha}_{\infty,\infty}(\Omega)}\sim_{\alpha,n,\Omega} \lVert (\nabla \uu, \mathfrak{p}) \rVert_{\B^{\alpha}_{\infty,\infty}(\Omega)},
\end{align*}
so we deduce
\begin{align*}
   (1+|\lambda|)^\frac{1-\alpha}{2} \lVert (\nabla \uu,\mathfrak{p} )\rVert_{\B^{\alpha}_{\infty,\infty}(\Omega)}\lesssim_{\Omega,\mu} \lVert \ff \rVert_{\L^\infty(\Omega)}.
\end{align*}
The estimate and $0$-sectoriality of $\AA_\mathcal{D}$ on $\L^\infty_{\mathfrak{n},\sigma}(\Omega)$ itself is a recent result of Geng and Shen  \cite[Theorems~1.1~\&~9.4]{GengShen2025} providing exactly the estimate
\begin{align*}
     (1+|\lambda|) \lVert \uu \rVert_{\L^\infty(\Omega)}\lesssim_{\Omega,\mu} \lVert \ff \rVert_{\L^\infty(\Omega)}.
\end{align*}
We end the proof discussing about the extra assertion concerning the smoothness of $\uu\in \mathcal{C}^{1,\alpha}(\overline{\Omega},\CC^n)$ whenever $\ff\in\C^{0}_{0,0,\sigma}(\Omega)$ and $\Omega$ is of class $\mathcal{C}^{1,\alpha}$. This can be seen from the localisation argument for a given a smooth compactly supported function $\ff\in\Ccinftydiv(\Omega)$ in the proof of Theorem~\ref{thm:StokesResolvent}, the function space $\mathcal{C}^{1,\alpha}$ being stable by composition and multiplication.
\end{proof}

Note that when $\lambda=0$, the proof is still independent of the knowledge of Geng and Shen's work, yielding a particular result which reads as follows
\begin{proposition}\label{prop:LinftySteady}Let $\alpha\in(0,1)$, and consider $\Omega$ to be a $\mathcal{C}^{1,\alpha}$-bounded domain or a $\C^{1,\alpha}$-bounded domain with small $(1,\alpha)$-H\"{o}lderian constant.

\medbreak

\noindent For all $\ff\in\L^\infty_{\mathfrak{n},\sigma}(\Omega)$, the steady Stokes system
\begin{equation*}\tag{DS}\label{eq:StokesLinfty}
    \left\{ \begin{array}{rllr}
         - \Delta \uu +\nabla \mathfrak{p} &= \ff \text{, }&&\text{ in } \Omega\text{,}\\
        \div \, \uu &= 0\text{, } &&\text{ in } \Omega\text{,}\\
        \uu_{|_{\partial\Omega}} &=0\text{, } &&\text{ on } \partial\Omega\text{.}
    \end{array}
    \right.
\end{equation*}
admits a unique solution $(\uu,\mathfrak{p})\in\C^{1,\alpha}(\overline{\Omega},\CC^n)\times\C^{0,\alpha}_\mfree(\overline{\Omega})$, with the estimate
\begin{align*}
    \lVert \uu \rVert_{\B^{1+\alpha}_{\infty,\infty}(\Omega)}+ \lVert \mathfrak{p} \rVert_{\B^{\alpha}_{\infty,\infty}(\Omega)} \lesssim_{\Omega} \lVert \ff\rVert_{\L^\infty(\Omega)}.
\end{align*}
\end{proposition}

 Note that, Cacciopoli-type inequalities and standard Moser's iteration argument \cite[Part~II,~Theorem~1.3,~b)]{GiaquintaModica1982} allows to remove the smallness assumption on the $(1,\alpha)$-H\"{o}lderian constant in Proposition~\ref{prop:LinftySteady}, so that the result was already known. The same argument still applies for the regularity of the resolvent problem  Theorem~\ref{thm:FinalResultLinfty1}. However, one then loses the regularity "scaling" with respect to resolvent parameter $\lambda$ (the $(1+|\lambda|)^\frac{1-\alpha}{2}$-factor in the estimate disappears). 

An interpolation argument as shown in the proof above yields the next Lemma.
\begin{lemma}\label{lem:lossLinftyDecayEstC1a} Let $\alpha\in(0,1]$, and $\Omega$ be a bounded $\C^{1,\alpha}$-domain. Let $\beta\in(0,1)$ and let $q\in[1,\infty]$. For all $\uu\in\B^{-\beta,\sigma}_{\infty,q,\mathfrak{n}}(\Omega)$, one has
\begin{align*}
    \lVert e^{-t\AA_\mathcal{D}}\uu\rVert_{\L^\infty(\Omega)} &\lesssim_{p,n,\Omega} \frac{1}{t^{\frac{\beta}{2}}}\lVert\uu\rVert_{\B^{-\beta}_{\infty,q}(\Omega)}
\end{align*}
\end{lemma}

\begin{proof}Consider $\beta\in(0,1)$, then for any $\varepsilon\in(0,1)$, since the semigroup generated by $\AA_\mathcal{D}$ is analytic on $\B^{-\varepsilon,\sigma}_{\infty,\infty,\mathfrak{n}}(\Omega)$, by the isomorphism property given by the characterization of fractional powers in Theorem~\ref{thm:FinalResultBesov2}, so it is on $\B^{\varepsilon,\sigma}_{\infty,\infty,\mathcal{D}}(\Omega)$. If we assume $\varepsilon<\beta<1-\varepsilon$, one has both inequalities as well as their underlying mapping properties:
\begin{align*}
    \lVert e^{-\AA_\mathcal{D}}\uu\rVert_{\B^{\varepsilon}_{\infty,\infty}(\Omega)}\lesssim_{\varepsilon,\beta,\Omega}\lVert \AA_\mathcal{D}^{\sfrac{\beta}{2}} e^{-\AA_\mathcal{D}}\uu\rVert_{\B^{\varepsilon-\beta}_{\infty,\infty}(\Omega)}\lesssim_{\varepsilon,\Omega}\frac{1}{t^\frac{\beta}{2}}\lVert \uu\rVert_{\B^{\varepsilon-\beta}_{\infty,\infty}(\Omega)}
\end{align*}
and
\begin{align*}
    \lVert e^{-\AA_\mathcal{D}}\uu\rVert_{\B^{-\varepsilon}_{\infty,\infty}(\Omega)}\lesssim_{\varepsilon,\beta,\Omega}\lVert \AA_\mathcal{D}^{\sfrac{\beta}{2}} e^{-\AA_\mathcal{D}}\uu\rVert_{\B^{-\varepsilon-\beta}_{\infty,\infty}(\Omega)}\lesssim_{\varepsilon,\Omega}\frac{1}{t^\frac{\beta}{2}}\lVert \uu\rVert_{\B^{-\varepsilon-\beta}_{\infty,\infty}(\Omega)}.
\end{align*}
By real interpolation, one obtains
\begin{align}
    \lVert e^{-\AA_\mathcal{D}}\uu\rVert_{(\B^{-\varepsilon,\sigma}_{\infty,\infty,\mathfrak{n}}(\Omega),\B^{\varepsilon,\sigma}_{\infty,\infty,\mathcal{D}}(\Omega))_{\frac{1}{2},1}}\lesssim_{\varepsilon,\beta,\Omega}\frac{1}{t^\frac{\beta}{2}}\lVert \uu\rVert_{(\B^{-\beta-\varepsilon,\sigma}_{\infty,\infty,\mathfrak{n}}(\Omega),\B^{-\beta+\varepsilon,\sigma}_{\infty,\infty,\mathfrak{n}}(\Omega))_{\frac{1}{2},1}}.\label{eq:lossLinftyDecayEstC1a}
\end{align}
However, one has
\begin{align*}
    (\B^{-\varepsilon,\sigma}_{\infty,\infty,\mathfrak{n}}(\Omega),\B^{\varepsilon,\sigma}_{\infty,\infty,\mathcal{D}}(\Omega))_{\frac{1}{2},1} \hookrightarrow (\B^{-\varepsilon}_{\infty,\infty}(\Omega,\CC^n),\B^{\varepsilon}_{\infty,\infty}(\Omega,\CC^n))_{\frac{1}{2},1} =\B^{0}_{\infty,1}(\Omega,\CC^n)\hookrightarrow\L^\infty(\Omega,\CC^n),
\end{align*}
and by Theorem~\ref{prop:InterpDivFreeC1}
\begin{align*}
    (\B^{-\beta-\varepsilon,\sigma}_{\infty,\infty,\mathfrak{n}}(\Omega),\B^{-\beta+\varepsilon,\sigma}_{\infty,\infty,\mathfrak{n}}(\Omega))_{\frac{1}{2},1} = \B^{-\beta,\sigma}_{\infty,1,\mathfrak{n}}(\Omega).
\end{align*}
Therefore, \eqref{eq:lossLinftyDecayEstC1a} becomes
\begin{align*}
    \lVert e^{-t\AA_\mathcal{D}}\uu\rVert_{\L^\infty(\Omega)} &\lesssim_{\varepsilon,\beta,\Omega} \frac{1}{t^{\frac{\beta}{2}}}\lVert\uu\rVert_{\B^{-\beta}_{\infty,1}(\Omega)}.
\end{align*}
Since $\varepsilon>0$ is arbitrary, the estimate above holds for all $\beta\in(0,1)$. Applying real interpolation again yields the case $\B^{-\beta}_{\infty,q}(\Omega)$ , $q\in[1,\infty]$, $\beta\in(0,1)$.
\end{proof}

We provide here a consequence for the $\L^1$-theory. To the best of the authors' knowledge, this is the first result for bounded domains whose boundary is not of high regularity. Even then, only several occurrences did appear in the literature. See \cite[Theorem~0.1]{GigaMatsuiShimizu1999} and \cite[Theorem~4.1]{ShibataShimizu2001} for $\L^1$-type estimates in the case of the half-space, \cite[Theorem~3.1]{Maremonti2011} for smooth domains.

\begin{corollary}\label{cor:StokesL1-semigroup} Let $\alpha\in(0,1]$, and let $\Omega$ be a bounded $\C^{1,\alpha}$-domain. For any $\uu\in\L^1_{\mathfrak{n},\sigma}(\Omega)$, all $s\in(0,2+\alpha)$, one has
\begin{align*}
    \lVert e^{-t\AA_\mathcal{D}}\uu\rVert_{\W^{s,1}(\Omega)} \lesssim_{n,s,\Omega} \frac{1}{t^\frac{s}{2}} \lVert \uu\rVert_{\L^{1}(\Omega)},\quad \text{ for all }t>0. 
\end{align*}
In particular,
\begin{align*}
    \lVert \nabla e^{-t\AA_\mathcal{D}}\uu\rVert_{\L^1(\Omega)} &\lesssim_{n,\Omega} \frac{1}{t^\frac{1}{2}} \lVert \uu\rVert_{\L^{1}(\Omega)},\quad \text{ for all }t>0.\\ 
    \lVert \nabla^2 e^{-t\AA_\mathcal{D}}\uu\rVert_{\L^1(\Omega)} &\lesssim_{n,\Omega} \frac{1}{t} \lVert \uu\rVert_{\L^{1}(\Omega)},
\end{align*}
\end{corollary}

\begin{proof} Let $s\in(0,1)$. By Lemma~\ref{lem:lossLinftyDecayEstC1a}, for any $\boldsymbol{\varphi}\in\Ccinfty(\Omega,\CC^n)$, one has
\begin{align}\label{eq:DecayEst}
    \lVert e^{-t\AA_\mathcal{D}}\PP_{\Omega}\boldsymbol{\varphi}\rVert_{\L^\infty(\Omega)} &\lesssim_{s,n,\Omega} \frac{1}{t^{\frac{s}{2}}}\lVert\boldsymbol{\varphi}\rVert_{\B^{-s}_{\infty,\infty}(\Omega)}.
\end{align}
For $\uu\in\Ccinftydiv(\Omega)$, we have the duality identity
\begin{align*}
    \langle e^{-t\AA_\mathcal{D}} \uu,\, \boldsymbol{\varphi}\rangle_{\Omega} = \langle  \uu,\, e^{-t\AA_\mathcal{D}}\PP_{\Omega}\boldsymbol{\varphi}\rangle_{\Omega}.
\end{align*}
Therefore by H\"{o}lder's inequality and \eqref{eq:DecayEst}
\begin{align*}
    |\langle e^{-t\AA_\mathcal{D}} \uu,\, \boldsymbol{\varphi}\rangle_{\Omega}|\leqslant \lVert \uu\rVert_{\L^1(\Omega)}\lVert e^{-t\AA_\mathcal{D}}\PP_{\Omega}\boldsymbol{\varphi} \rVert_{\L^\infty(\Omega)} \lesssim_{s,n,\Omega} \frac{1}{t^\frac{s}{2}}\lVert \uu\rVert_{\L^1(\Omega)}\lVert \boldsymbol{\varphi} \rVert_{\B^{-s}_{\infty,\infty}(\Omega)}.
\end{align*}
Taking the supremum over all $\boldsymbol{\varphi}$ with $\B^{-s}_{\infty,\infty}$-norm $1$, yields for all $\uu\in\Ccinftydiv(\Omega)$,
\begin{align}\label{eq:proofWs1Decay1}
    \lVert  e^{-t\AA_\mathcal{D}}\uu\rVert_{\W^{s,1}(\Omega)}=\lVert \uu\rVert_{\B^{s}_{1,1}(\Omega)}\lesssim_{s,n,\Omega} \frac{1}{t^\frac{s}{2}}\lVert \uu\rVert_{\L^1(\Omega)}.
\end{align}
Now, consider $\beta>0$ such that $s+\beta\in(1,2+\alpha)$, $s+\beta\neq 1$. By Corollary~\ref{cor:StokesWs1}, it holds
\begin{align}\label{eq:proofWs1Decay2}
    \lVert  e^{-t\AA_\mathcal{D}} \uu\rVert_{\W^{s+\beta,1}(\Omega)}\sim_{\beta,\Omega} \lVert \AA_\mathcal{D}^{\sfrac{\beta}{2}} e^{-t\AA_\mathcal{D}} \uu \rVert_{\W^{s,1}(\Omega)} \lesssim_{s,n,\Omega} \frac{1}{t^\frac{s+\beta}{2}}\lVert \uu\rVert_{\L^1(\Omega)}.
\end{align}
By real interpolation between \eqref{eq:proofWs1Decay1} and \eqref{eq:proofWs1Decay2}, for $k=1,2$ one reaches the last case with
\begin{align*}
    \lVert \nabla^k e^{-t\AA_\mathcal{D}}\uu\rVert_{\L^1(\Omega)} \lesssim_{n} \lVert e^{-t\AA_\mathcal{D}}\uu\rVert_{\B^k_{1,1}(\Omega)} \lesssim_{s,n,\Omega} \frac{1}{t^\frac{k}{2}}\lVert \uu\rVert_{\L^1(\Omega)}.
\end{align*}
This ends the proof.
\end{proof}

As another direct consequence of Lemma~\ref{lem:lossLinftyDecayEstC1a}, we obtain the following important $\L^p$-$\L^q$ decay estimates for the semigroup which includes the end-point spaces $\L^1$ and $\L^\infty$. We mention that fractional versions are possible.

\begin{proposition}\label{prop:Optimaldecay} Let $\alpha\in(0,1]$, and let $\Omega$ be a bounded $\C^{1,\alpha}$-domain. Let $1\leqslant p < q\leqslant \infty$. It holds that
\begin{enumerate}
    \item  for all $\uu\in\L^p_{\mathfrak{n},\sigma}(\Omega)$, one has for all $t>0$
\begin{align*}
    \lVert e^{-t\AA_\mathcal{D}}\uu\rVert_{\L^q(\Omega)} &\lesssim_{p,q,n,\Omega} \frac{1}{t^{\frac{n}{2}\left(\frac{1}{p}-\frac{1}{q}\right)}}\lVert\uu\rVert_{\L^p(\Omega)},\\
    \quad \lVert \nabla e^{-t\AA_\mathcal{D}}\uu\rVert_{\L^q(\Omega)} &\lesssim_{p,q,n,\Omega} \frac{1}{t^{\frac{n}{2}\left(\frac{1}{p}-\frac{1}{q}\right)+\frac{1}{2}}}\lVert\uu\rVert_{\L^p(\Omega)},\\
    \text{ and }\quad \lVert \AA_\mathcal{D} e^{-t\AA_\mathcal{D}}\uu\rVert_{\L^q(\Omega)} &\lesssim_{p,q,n,\Omega} \frac{1}{t^{\frac{n}{2}\left(\frac{1}{p}-\frac{1}{q}\right)+1}}\lVert\uu\rVert_{\L^p(\Omega)};
\end{align*}
\item for all $\mathbf{F}\in\L^p(\Omega,\CC^{n^2})$, one has for all $t>0$,
\begin{align*}
    \lVert  e^{-t\AA_\mathcal{D}}\PP_{\Omega}\div(\mathbf{F})\rVert_{\L^q(\Omega)} &\lesssim_{p,q,n,\Omega} \frac{1}{t^{\frac{n}{2}\left(\frac{1}{p}-\frac{1}{q}\right)+\frac{1}{2}}}\lVert \mathbf{F}\rVert_{\L^p(\Omega)}.
\end{align*}
\end{enumerate}
Furthermore, the two last estimates in \textit{(i)} remain valid whenever $p=q$.
\end{proposition}

\begin{proof}\textbf{Step 1:} The only main obstacle is about reaching the endpoint cases. First note that for $r>n$, thanks to Lemma~\ref{lem:lossLinftyDecayEstC1a} and Sobolev' embeddings
\begin{align}
    \lVert e^{-t\AA_\mathcal{D}}\uu\rVert_{\L^\infty(\Omega)} 
    &\lesssim_{\Omega,r} \frac{1}{t^{\frac{n}{2r}}}\lVert e^{-\frac{t}{2}\AA_\mathcal{D}}\uu\rVert_{\B^{-\frac{n}{r}}_{\infty,\infty}(\Omega)}\nonumber\\
    &\lesssim_{\Omega,r} \frac{1}{t^{\frac{n}{2r}}}\lVert e^{-\frac{t}{2}\AA_\mathcal{D}}\uu\rVert_{\L^{r}(\Omega)}\nonumber\\
    &\lesssim_{\Omega,r} \frac{1}{t^{\frac{n}{2r}}}\lVert\uu\rVert_{\L^{r}(\Omega)}.\label{eq:ProofEstLinftyLr}
\end{align}
Up to replace  $\uu$, by $\PP_\Omega \uu$, by duality, considering temporarily $\boldsymbol{\varphi}\in\Ccinftydiv(\Omega)$, for $\uu\in\L^r(\Omega,\CC^n)$ one obtains
\begin{align*}
    \langle e^{-t\AA_\mathcal{D}}\boldsymbol{\varphi},\uu\rangle_{\Omega} &= \langle \PP_{\Omega} e^{-t\AA_\mathcal{D}}\boldsymbol{\varphi},\,\uu\rangle_{\Omega}\\
    &= \langle \boldsymbol{\varphi},\, e^{-t\AA_\mathcal{D}}\PP_{\Omega}\uu\rangle_{\Omega},
\end{align*}
so that by H\"{o}lder's inequality it holds
\begin{align*}
   | \langle e^{-t\AA_\mathcal{D}}\boldsymbol{\varphi},\uu\rangle_{\Omega}|&\leqslant \lVert \boldsymbol{\varphi}\rVert_{\L^1(\Omega)}\lVert e^{-t\AA_\mathcal{D}}\PP_{\Omega}\uu\rVert_{\L^\infty(\Omega)}\\
   &\lesssim_{r,n,\Omega} \frac{1}{t^{\frac{n}{2r}}}\lVert \boldsymbol{\varphi}\rVert_{\L^1(\Omega)}\lVert \uu\rVert_{\L^r(\Omega)}.
\end{align*}
Taking the supremum over all $\uu$ with $\L^r$-norm $1$, we did obtain for all $1<r'<\frac{n}{n-1}$, all $\boldsymbol{\varphi}\in\L^{1}_{\mathfrak{n},\sigma}(\Omega,\CC^n)$,
\begin{align}\label{eq:ProofEstLinftyL1Lr'}
    \lVert e^{-t\AA_\mathcal{D}}\boldsymbol{\varphi}\rVert_{\L^{r'}(\Omega)} \lesssim_{\Omega,r}  \frac{1}{t^{\frac{n}{2}-\frac{n}{2r'}}}\lVert\boldsymbol{\varphi}\rVert_{\L^{1}(\Omega)}.
\end{align}

\textbf{Step 2:} The case $1\leqslant p <q\leqslant \infty$. By Theorem~\ref{thm:FinalResultSobolevC1q}, for any $r\in(1,\infty)$,
\begin{align*}
    \lVert \nabla \vv\rVert_{\L^{r}(\Omega)}\sim_{r,\Omega} \lVert\AA_\mathcal{D}^{\sfrac{1}{2}}\vv\rVert_{\L^{r}(\Omega)}.
\end{align*}
Therefore, for $1<p<n$, $\frac{1}{r}=\frac{1}{p}-\frac{1}{n}$ by Sobolev embeddings and analyticity of the Stokes--Dirichlet semigroup\footnote{One can check the end of the proof for a more elementary argument.}
\begin{align}\label{eq:ProofEstLinftyLpLp*}
    \lVert  e^{-t\AA_\mathcal{D}}\boldsymbol{\varphi}\rVert_{\L^{r}(\Omega)} \lesssim_{r,n}\lVert \nabla e^{-t\AA_\mathcal{D}}\boldsymbol{\varphi}\rVert_{\L^{p}(\Omega)} \lesssim_{\Omega,r}  \frac{1}{t^\frac{1}{2}}\lVert \boldsymbol{\varphi}\rVert_{\L^{p}(\Omega)}.
\end{align}
Hence for $p\in(\frac{n}{2},n)$, it holds $r>n$, and combining \eqref{eq:ProofEstLinftyLr} with \eqref{eq:ProofEstLinftyLpLp*}, for any $p\in(\frac{n}{2},n)\cup(n,\infty)$, it holds
\begin{align*}
    \lVert e^{-t\AA_\mathcal{D}}\boldsymbol{\varphi}\rVert_{\L^\infty(\Omega)} \lesssim_{p,\Omega} \frac{1}{t^{\frac{n}{2p}}}\lVert \boldsymbol{\varphi}\rVert_{\L^{p}(\Omega)},
\end{align*}
complex interpolation allowing the case $p=n$. One can iterate the argument, so that it holds for any $p\in(\frac{n}{k},\infty)$, for the $k^\text{th}$ iteration, $k\in\llb 2,n\rrb$. Therefore the $\L^p$-$\L^\infty$ estimates holds for all $p\in(1,\infty)$. By applying  complex interpolation with the $\L^p$-$\L^p$ estimate for $1<p<\infty$ (uniform boundedness of the semigroup), one obtains the $\L^p$-$\L^q$ estimate for all $1<p<q\leqslant \infty$. Now, taking $1<p<\frac{n}{n-1}$, and $q\in(p,\infty]$, combining the $\L^p$-$\L^q$ estimate with \eqref{eq:ProofEstLinftyL1Lr'} (here $r'=p$) one obtains therefore the $\L^1$-$\L^q$ estimate for any $1<q\leqslant \infty$. All the remaining estimates are straightforward consequences, except the gradient estimate
\begin{align*}
    \lVert \nabla e^{-t\AA_\mathcal{D}}\uu\rVert_{\L^\infty(\Omega)} \lesssim_{\Omega,p} \frac{1}{t^{\frac{n}{2p}+\frac{1}{2}}}\lVert \uu\rVert_{\L^p(\Omega)}.
\end{align*}
Due to the semigroup property $e^{-t\AA_\mathcal{D}}=e^{-\frac{t}{2}\AA_\mathcal{D}}e^{-\frac{t}{2}\AA_\mathcal{D}}$, and the already proven $\L^p$-$\L^q$ estimate, it is sufficient to check
\begin{align}\label{eq:LinftyLinftygradientDecay}
    \lVert \nabla e^{-t\AA_\mathcal{D}}\uu\rVert_{\L^\infty(\Omega)} \lesssim_{\Omega} \frac{1}{t^\frac{1}{2}}\lVert \uu\rVert_{\L^\infty(\Omega)}.
\end{align}
This is however a direct consequence of interpolation inequalities applied to Theorem~\ref{thm:FinalResultLinfty1}, yielding
\begin{align*}
    (1+|\lambda|)^\frac{1}{2}\lVert \nabla(\lambda\I+\AA_\mathcal{D})^{-1} \ff\rVert_{\L^\infty(\Omega)}\lesssim_{\mu,\Omega}\lVert \ff\rVert_{\L^\infty(\Omega)},\quad \ff\in\L^\infty_{\mathfrak{n},\sigma}(\Omega),\quad \lambda\in\Sigma_\mu\cup\{0\},\quad \mu\in(0,\pi).
\end{align*}
The latter can be tuned into \eqref{eq:LinftyLinftygradientDecay} by plugging it into the Cauchy integral formula that represents the semigroup through integration of its resolvent, see \cite[Proposition~3.7]{TolksdorfWatanabe}. We finish the proof claiming that the estimates in Point \textit{(ii)} are obtained by duality  and the same arguments one can find in the first half of the current step. We conclude: by duality arguments one obtains the $\L^p-\L^q$,  $\L^1-\L^q$ estimates for $\mathbf{F}\longmapsto e^{-t\AA_\mathcal{D}}\PP_{\Omega}\div(\mathbf{F})$, whenever $1<p\leqslant q<\infty$. The remaining cases of $\L^p-\L^\infty$ and $\L^1-\L^\infty$-estimates are obtained by writing $ e^{-t\AA_\mathcal{D}}\PP_{\Omega}\div(\mathbf{F})= e^{-\frac{t}{2}\AA_\mathcal{D}}e^{-\frac{t}{2}\AA_\mathcal{D}}\PP_{\Omega}\div(\mathbf{F})$.
\end{proof}

\subsubsection{Lowering the regularity: removing the \texorpdfstring{$\C^1$}{C1}-regularity assumption}\label{Sec:RemoveC1}

We state the following Meta-Theorem.

\begin{theorem}\label{thm:metaThmremoveC1}In \textbf{Theorem~\ref{thm:StokesResolvent}} and all its subsequent results, one can replace the assumption
\begin{center}
    ``Let $\alpha\in(0,1)$, $r\in[1,\infty]$, and let $\Omega$ be a bounded \textbf{$\C^1$}-domain of class $\mathcal{M}^{1+\alpha,r}_\W(\epsilon)$ for some sufficiently small $\epsilon>0$.''
\end{center}
 by 
\begin{center}
    ``Let $\alpha\in(0,1)$, $r\in[1,\infty]$, and let $\Omega$ be a bounded \textbf{Lipschitz} domain of class $\mathcal{M}^{1+\alpha,r}_\W(\epsilon)$ for some sufficiently small $\epsilon>0$.''
\end{center}
\end{theorem}

\noindent Note that this statement is not artificial due to Remark~\ref{rem:MultiplierBoundarynotC1}.

For better readability, so the reader's convenience, we provide a detailed scheme of proof.

\begin{proof}Note the that the main crucial point in the proof of \eqref{eq:H-1EstBddLip} asking for a $\C^1$-domain, was the estimate \eqref{eq:shen}
\begin{align}\label{eq:GengShen}
    (1+|\lambda|)^\frac{1}{2}\lVert \uu \rVert_{\L^{q}(\Omega)}\lesssim_{\Omega,q,\mu}\lVert\ff\rVert_{\W^{-1,q}(\Omega)},\,\quad \ff\in\W^{-1,q}(\Omega,\CC^n),
\end{align}
in combination with the embeddings \eqref{eq:emb1} and \eqref{eq:emb2},
\begin{align*}
    \H^{s,p}(\Omega)&\hookrightarrow\W^{-1,q}(\Omega),\\
   \text{ and } \quad\quad\quad\quad \L^{q}(\Omega)&\hookrightarrow\H^{s_0,p}(\Omega),
\end{align*}
for a well chosen $q\in(1,\infty)$ and $s_0< 0$,  depending on $s$ and $p$ ($s_0=-(n+1)$ is sufficient).

\noindent Note that according to the proof of Lemma~\ref{lem:H-1EstBddLip}, to reach \eqref{eq:GengShen}, it is sufficient to know for all $\ff\in\Ccinfty(\Omega,\CC^n)$
\begin{align*}
(1+|\lambda|)\lVert \uu \rVert_{\L^{q}(\Omega)}&\lesssim_{\Omega,q,\mu}\lVert\ff\rVert_{\L^{q}(\Omega)}\\
    \text{ and }\quad(1+|\lambda|)^\frac{1}{2}\lVert \nabla \uu \rVert_{\L^{q'}(\Omega)}&\lesssim_{\Omega,\mu,q}\lVert\ff\rVert_{\L^{q'}(\Omega)}.
\end{align*}
Finally, to obtain the latter gradient estimate it is sufficient to know the way stronger result
\begin{align*}
    \D_{q'}(\AA_\mathcal{D}^{\sfrac{1}{2}})=\H^{1,q'}_{0,\sigma}(\Omega)\quad\text{ and }\quad\lVert\AA_\mathcal{D}^{\sfrac{1}{2}}\vv\rVert_{\L^{q'}(\Omega)}\sim_{\Omega,q'}\lVert \nabla \vv\rVert_{\L^{q'}(\Omega)}.
\end{align*}

Our strategy is as follows: with an idea similar to the one in the proof Proposition~\ref{prop:HodgeDiracbdd} we proceed by a bootstrap. We proceed for the first main step as follows (as in the proof of Theorem~\ref{thm:StokesResolvent}, we assume without loss of generality $r=1$):
\begin{enumerate}
    \item Choose $p_1:=\frac{2(n-1)}{n}$, so that for all $p\in(p_1,p_1')$, all $s\in(-1+\sfrac{1}{p},-1+\sfrac{(1+\alpha)}{p})$, one has an embeddings
    \begin{align*}
    \H^{s,p}(\Omega)&\hookrightarrow\H^{-1,2}(\Omega),\\
   \text{ and } \quad\quad\quad\quad \L^{2}(\Omega)&\hookrightarrow\H^{-(n+1),p}(\Omega).
\end{align*}
\item Consider $\ff\in\Ccinftydiv(\Omega)$, and proceed as in Step 1 and Step 2, in order to obtains for some given $\vartheta\in(0,1)$, and for all $p\in(p_1,p_1')$ 
\begin{align*}
    (1+|\lambda|)\|{\uu}\|_{\H^{s,p}(\Omega)}+\|{\uu}\|_{\H^{s+2,p}(\Omega)}+\|{\mathfrak{p}}\|_{\H^{s+1,p}(\Omega)}&\lesssim_{p,s,n,\theta}^{\Omega,\mu} \|\uu\|_{\H^{s+2-\vartheta,p}(\Omega)}+\|\ff\|_{\B^{s}_{p,q}(\Omega)}.
\end{align*}
\item By interpolation inequalities and Lemma~\ref{lem:H-1EstBddLip}
 \begin{align*}
\|\uu\|_{\H^{s+2-{\vartheta},p}(\Omega)}&\lesssim_{p,s}^{n,\Omega} \|\uu\|_{\H^{s+2,p}(\Omega)}^{\theta}\|\uu\|_{\H^{-(n+1),p}(\Omega)}^{1-\theta}\lesssim_{p,s,}^{n,\Omega} \|\uu\|_{\H^{s+2,p}(\Omega)}^{\theta}\|\uu\|_{\L^{2}(\Omega)}^{1-\theta} \\&\lesssim_{p,s,\mu}^{n,\Omega}\|\uu\|_{\H^{s+2,p}(\Omega)}^{\theta}\|\ff\|_{\H^{-1,2}(\Omega)}^{1-\theta}\lesssim_{p,s,\mu}^{n,\Omega}\|\uu\|_{\H^{s+2,p}(\Omega)}^{\theta}\|\ff\|_{\H^{s,p}(\Omega)}^{1-\theta}.
 \end{align*}
 \item As in Step 3 in the proof of Theorem~\ref{thm:StokesResolvent}, by interpolation inequalities, one obtains for all $\ff\in\Ccinftydiv(\Omega)$, all $p\in(p_1,p_1')$, $s\in(-1+\sfrac{1}{p},-1+\sfrac{(1+\alpha)}{p})$,
\begin{align*}
    (1+|\lambda|)\|{\uu}\|_{\H^{s,p}(\Omega)}+\|{\uu}\|_{\H^{s+2,p}(\Omega)}+\|{\mathfrak{p}}\|_{\H^{s+1,p}(\Omega)}&\lesssim_{p,s,n,\theta}^{\Omega,\mu} \|\ff\|_{\H^{s,p}(\Omega)}.
\end{align*}
 Therefore, for all $p\in(p_1,p_1')$, $s\in(-1+\sfrac{1}{p},-1+\sfrac{(1+\alpha)}{p})$, one obtains that one can promote the Dirichlet--Stokes operators as an injective $0$-sectorial operator $(\D_p^s(\AA_\mathcal{D}),\AA_\mathcal{D})$ on $\H^{s,p}_{\mathfrak{n},\sigma}(\Omega)$, with domain
 \begin{align*}
     \D_p^s(\AA_\mathcal{D})=\H^{s+1}_{0,\sigma}(\Omega)\cap\H^{s+2,p}(\Omega).
 \end{align*}
 \item Now, one can reproduce Step 1 of the proof in Section~\ref{sec:ProofLPbddDomain}: the extrapolation method of the $\mathbf{H}^\infty$-functional calculus by the Kunstmann-Weis comparison principle (for the interpolation of the $\mathcal{R}$-bound take a $q_0\in(p_1,2)$ arbitrarily close to $p_1$). Therefore, one obtains for all $p\in(p_1,p_1')$, $s\in(-1+\sfrac{1}{p},-1+\sfrac{(1+\alpha)}{p})$, that  $(\D_p^s(\AA_\mathcal{D}),\AA_\mathcal{D})$ has a bounded $\mathbf{H}^\infty$-functional calculus on $\H^{s,p}_{\mathfrak{n},\sigma}(\Omega)$.

 \medbreak
 
 In particular it has BIP, and from Proposition~\ref{prop:InterpDivFreeC1}, on obtains that for all $p\in(p_1,p_1')$, a $0$-sectorial injective operator $(\D_p(\AA_\mathcal{D}),\AA_\mathcal{D})$ on $\L^p_{\mathfrak{n},\sigma}(\Omega)$ with the squareroot property
 \begin{align*}
     \D_{p}(\AA_\mathcal{D}^{\sfrac{1}{2}})=\H^{1,p}_{0,\sigma}(\Omega)\quad\text{ and }\quad\lVert\AA_\mathcal{D}^{\sfrac{1}{2}}\vv\rVert_{\L^{p}(\Omega)}\sim_{\Omega,p}\lVert \nabla \vv\rVert_{\L^{p}(\Omega)}.
 \end{align*}
 \item The previous step \textit{(v)} implies in particular that for all $p\in(p_1,p_1')$, all $\ff\in\L^p_{\mathfrak{n},\sigma}(\Omega)$
 \begin{align*}
     \lVert&\nabla (\lambda\I+\AA_\mathcal{D})^{-1}\ff \rVert_{\L^{p}(\Omega)}\sim_{p,s,n}^{\Omega} \lVert \AA_\mathcal{D}^{\sfrac{1}{2}} (\lambda\I+\AA_\mathcal{D})^{-1}\ff \rVert_{\L^{p}(\Omega)}\lesssim_{p,s,n}^{\Omega,\mu} (1+|\lambda|)^{-\frac{1}{2}}\lVert\ff\rVert_{\L^{p}(\Omega)}.
 \end{align*}
 So that with the exact same proof of Lemma~\ref{lem:H-1EstBddLip}, it can be turned into the estimate
 \begin{align*}
     (1+|\lambda|)^\frac{1}{2}\lVert \uu \rVert_{\L^{p}(\Omega)}\lesssim_{p,s,n}^{\Omega,\mu} \lVert\ff\rVert_{\W^{-1,p}(\Omega)},\,\quad \ff\in\W^{-1,p}(\Omega,\CC^n).
 \end{align*}
 valid for all $\ff\in\W^{-1,p}(\Omega,\CC^n)$, $p\in(p_1,p_1')$.
\end{enumerate}

This ends the first Step of the bootstrap procedure. Now, one sets for all integer $k\geqslant1$, $p_{k+1}:=\max\big(\frac{n-1}{n}p_k,1\big)$, so that one has the sharp Sobolev embedding
\begin{align*}
    \H^{-1+\frac{1}{p_{k+1}},p_{k+1}}(\Omega)&\hookrightarrow\H^{-1,p_k}(\Omega).
\end{align*}

By the exact last step of the first iteration, one can reproduce the whole procedure and iterate for all $p\in(p_2,p_2')$, replacing $(\L^2,\H^{-1,2})$ by $(\L^q,\H^{-1,q})$ for some $q\in(p_1,p_1')$.

\medbreak

Similarly, for all $k\geqslant 2$, for all $p\in(p_k,p_k')$, choosing then an appropriate $(\L^q,\H^{-1,q})$ for some $q\in(p_{k-1},p_{k-1}')$. This yields the result for Sobolev spaces $\H^{s,p}(\Omega)$ for all $p\in(1,\infty)$, all $s\in(-1+\sfrac{1}{p},\sfrac{1}{p})$ yielding also directly the characterization for the domains of fractional powers. 

\medbreak

Finally, real interpolation yields the cases of Besov spaces $\B^{s}_{p,q}(\Omega)$ for all $p\in(1,\infty)$, all $q\in[1,\infty]$, all $s\in(-1+\sfrac{1}{p},\sfrac{1}{p})$. The case of Besov spaces $p=1,\infty$ can be then achieved separately, reproducing the proof of Theorem~\ref{thm:StokesResolvent} for their sole purpose (assuming implicitly $r=\infty$ when $p=\infty$). 
\end{proof}

\begin{remark}\label{rem:FuncCalcSimplAbsFC}We did choose to present the proof above which also proves iteratively the extrapolation of the bounded holomorphic calculus on $\L^p$-spaces by the Kunstmann-Weis argument for only one reason. That is: it only uses and involves Sobolev and $\L^p$-spaces. No use of Besov spaces was needed except for their own purpose. This is designed for the reader which is only interested by the $\L^p$-theory, $1<p<\infty$.

\medbreak

\noindent It goes without saying that the argument here can be simplified by removing the Kunstmann-Weis extrapolation method, if we replace the embeddings
\begin{align*}
    (1+|\lambda|)^\frac{1}{2}\lVert \uu \rVert_{\L^{q}(\Omega)}\lesssim_{\Omega,r,\mu}\lVert\ff\rVert_{\W^{-1,q}(\Omega)},\,\quad \ff\in\W^{-1,q}(\Omega,\CC^n).
\end{align*}
in combination with embeddings
\begin{align*}
    \H^{s,p}(\Omega)&\hookrightarrow\W^{-1,q}(\Omega),\\
   \text{ and } \quad\quad\quad\quad \L^{q}(\Omega)&\hookrightarrow\H^{s_0,p}(\Omega),
\end{align*}
by their following Besov counterparts
\begin{align}\label{eq:GengShenBesov}
    (1+|\lambda|)^\frac{1}{2}\lVert \uu \rVert_{\B^{0}_{q,\kappa}(\Omega)}\lesssim_{\Omega,q,\mu}\lVert\ff\rVert_{\B^{-1}_{q,\kappa}(\Omega)},\,\quad \ff\in\B^{-1}_{q,\kappa}(\Omega,\CC^n),
\end{align}
as well as the embeddings
\begin{align*}
    \H^{s,p}(\Omega)&\hookrightarrow\B^{-1}_{q,\kappa}(\Omega),\\
   \text{ and } \quad\quad\quad\quad \B^{0}_{q,\kappa}(\Omega)&\hookrightarrow\H^{-(n+1),p}(\Omega).
\end{align*}
This is due to the fact that Dore's Theorem, \cite[Theorem~6.1.3]{bookHaase2006} and results such as \cite[Proposition~6.4.1,~a), \& Corollary~6.5.5]{bookHaase2006}, are making the bounded holomorphic functional calculus --and then the identifications of fractional powers-- for the Stokes--Dirichlet operator automatic on the scale of Besov spaces, thanks to  Proposition~\ref{prop:InterpDivFreeC1}. This implies automatically \eqref{eq:GengShenBesov}.

As a real interpolation scale, Besov spaces  behaves extremely well with semigroup theory and sectorial operators. If one restrict themselves to the study in the sole setting of Besov spaces, then it should be possible to reach any result above and perform all the proofs in a ``quite elementary'' (but quite long) way, without invoking explicitly any abstract argument concerning the (bounded) holomorphic functional calculus or the BIP property on such spaces. (However Dore's Theorem and the so called absolute functional calculus remain to be the culprits, just the consequences of them for the study of some "concrete operators" can be derived without such general statements in mind. The notion of Absolute Functional Calculus is due to Kalton and Kucherenko \cite{KaltonKucherenko2010}.) 
\end{remark}

Finally, this iterative procedure shows that our strategy is completely unrelated to the specific Stokes--Dirichlet (resolvent) problem and that it can be adapted to many other kind of (non-local) elliptic problems, as long as one has a suitable representation formula, or at least estimates, for a corresponding problem on $\RR^n_+$. One can even deal with a sufficiently wide class of non-self-adjoint problems, say, for a given second order elliptic operator $\mathcal{L}$ with prescribed sufficiently natural boundary values\footnote{\textit{i.e.} in the Lopatinskii-Shapiro sense. For more details about these types of boundary conditions, see, for instance, \cite{bookDenkHieberPruss2003,Guidetti1991Ell,Guidetti1991Interp,KunstmannWeis2004,PrussSimonett2016}.}, if one considers simultaneously the elliptic problems arising from $\mathcal{L}$ and $\mathcal{L}^\ast$.

\newpage

\section{Maximal regularity results for the Stokes--Dirichlet problem}\label{sec:MaxReg}

This particular and last section aims to present  optimal solvability results for the Stokes--Dirichlet evolution problem on bounded domains $\Omega$ with minimal regularity assumption and on $\RR^n_+$. It mostly exhibits nearly straightforward consequences of our established  abstract-analytic and regularity theory of the Stokes--Dirichlet (resolvent) problem and its underlying operator.

\medbreak

Indeed, thanks to the standard $\L^q$-maximal regularity results Theorems~\ref{thm:LqMaxRegUMDHomogeneous} and \ref{thm:DaPratoGrisvard}, and Theorems~\ref{thm:FinalResultSobolev1} and \ref{thm:FinalResultBesov1}, one can consider the Stokes--Dirichlet operator $\AA_\mathcal{D}$ on a wide family of function spaces $\{ \X^{s,p}(\Omega),\,{s,p}\}$ for the solvability of the linear evolution problem. In particular, for any $T\in(0,\infty]$, one can solve uniquely the abstract Cauchy problem
\begin{align}\tag{ACP}\label{ACPStokes}
    \left\{\begin{array}{rl}
            \partial_t  \uu(t) +\AA_\mathcal{D}\uu(t)  =& \ff(t) \,\text{, } 0<t<T,\\
            u(0) =& u_0\text{,}
    \end{array}
    \right.
\end{align}
with estimates
\begin{align*}
    \lVert \uu \rVert_{\L^{\infty}(\mathcal{D}_{\AA_\mathcal{D}}(1+\theta-{\sfrac{1}{q},q}))}\less_{\theta,q,\Omega} \lVert  (\partial_t\uu,\AA_\mathcal{D}\uu)\rVert_{\L^q(\mathcal{D}_{\AA_\mathcal{D}}(\theta,q))} \less_{\theta,q,\Omega} \lVert \ff \rVert_{\L^q(\mathcal{D}_{\AA_\mathcal{D}}(\theta,q))} + \lVert \uu_0 \rVert_{\mathcal{D}_{\AA_\mathcal{D}}(1+\theta-{\sfrac{1}{q},q})}.
\end{align*}
We recall that $\mathcal{D}_{\AA_\mathcal{D}}(\theta,q)=(\X_{\mathfrak{n},\sigma},\D(\AA_\mathcal{D}))_{\theta,q}$, $\theta\in(0,1)$, $q\in[1,\infty]$, and 
\begin{align*}
    \lVert \uu \rVert_{\L^{\infty}(\mathcal{D}_{\AA_\mathcal{D}}(1-{\sfrac{1}{q},q}))}\less_{q,\Omega} \lVert  (\partial_t\uu,\AA_\mathcal{D}\uu)\rVert_{\L^q(\X)} \less_{q,\Omega} \lVert \ff \rVert_{\L^q(\X)} + \lVert \uu_0 \rVert_{\mathcal{D}_{\AA_\mathcal{D}}(1-{\sfrac{1}{q},q})},
\end{align*}
whenever $\X_{\mathfrak{n},\sigma}$ is UMD and $1<q<\infty$.

It appears clearly that to have a full description of the solvability of the problem in terms of sharp regularity of solutions, one has to describe the space of initial data when possible. Therefore, one does need to compute the real interpolation spaces
\begin{align*}
        \mathcal{D}_{\AA_\mathcal{D}}(\theta,q)=(\X_{\mathfrak{n},\sigma},\D(\AA_\mathcal{D}))_{\theta,q},
\end{align*}
for all $\theta\in(0,1)$, all $q\in[1,\infty]$, and a suitable choice of $\X_{\mathfrak{n},\sigma}$.

\textbf{The goals of this chapter are the following:}
\begin{itemize}
    \item As briefly outlined above, the first section of this chapter, \textbf{Section~\ref{Sec:InterpTheoryDomainStokes}}, is devoted to the interpolation theory for the domain of the Stokes--Dirichlet operator. 

    \medbreak
    
    In the case of the flat half-space $\Omega=\RR^n_+$, the situation turns out to be more delicate. In order to apply the framework of Theorem~\ref{thm:DaPratoGrisvard} and obtain global-in-time estimates, one is naturally led to consider interpolation of (abstract) homogeneous function spaces with boundary conditions, for which very few results are available in the literature. We therefore provide a complete description of the interpolation spaces associated with the Stokes operator on $\RR^n_+$, including the endpoint cases $p,q=1,\infty$.

    \item The second section, \textbf{Section~\ref{sec:StandardMaxRegStokes}}, is concerned with solvability results for the Stokes--Dirichlet evolution problem with sharp regularity description of solutions. We begin with the well-known case of the half-space and then proceed to bounded domains under minimal regularity assumptions.

    \item Finally, the $\L^\infty$-setting for the Stokes--Dirichlet operator does not fall within any of the above frameworks. In \textbf{Section~\ref{Sec:MAxRegforLinfty}}, we highlight a lesser-known variant of the Da Prato--Grisvard theory which nevertheless ensures global-in-time maximal regularity with $\L^\infty$ as a ground space for the space variable.
\end{itemize}

\subsection{Interpolation theory for the domain of the Dirichlet-Stokes operator on \texorpdfstring{$\RR^n_+$}{the half-space} and bounded rough domains}\label{Sec:InterpTheoryDomainStokes}

We start with the case of bounded domains with minimal regularity assumption.

\begin{proposition}\label{prop:InterpDivFreeC1a} Let $\alpha\in(0,1)$, $r\in[1,\infty]$, and let $\Omega$ be a bounded Lipschitz domain of class $\mathcal{M}^{1+\alpha,r}_\W(\epsilon)$ for some sufficiently small $\epsilon>0$.

Let $p,q,q_0,q_1\in[1,\infty]$, $-1+{\sfrac{1}{p}}<s_j<2+{\sfrac{1}{p}}$, $s_j\neq \sfrac{1}{p}$. For all $\theta\in(0,1)$, such that $s:=s_0(1-\theta)+\theta s_1$ satisfies \hyperref[eq:ConditionRegularity]{$(\prescript{\alpha}{r}{\ast}_{p,q}^{s-2})$}, $s\neq\sfrac{1}{p}$, the following interpolation identity holds true
\begin{align}\label{eq:realInterpDivFreebddRoughDomain}
        \big({\H}^{s_0,p}_{{\mathcal{D}},\sigma}(\Omega),{\H}^{s_1,p}_{{\mathcal{D}},\sigma}(\Omega)\big)_{\theta,q} =\big({\B}^{s_0,\sigma}_{p,q_0,{\mathcal{D}}}(\Omega),{\B}^{s_1,\sigma}_{p,q_1,{\mathcal{D}}}(\Omega)\big)_{\theta,q} = {\B}^{s,\sigma}_{p,q,\mathcal{D}}(\Omega);
\end{align}
Still under the assumption \hyperref[eq:ConditionRegularity]{$(\prescript{\alpha}{r}{\ast}_{p,q}^{s-2})$},  for $s_1=s_0+2$, $s_0\in(-1+\sfrac{1}{p},\sfrac{1}{p})$, for $s\neq1/p$, one also has the identity
\begin{align}\label{eq:realInterpDomainStokesbddDomain}
        \big({\H}^{s_0,p}_{{\mathfrak{n}},\sigma}(\Omega),\D^{s_0}_{p}(\AA_\mathcal{D})\big)_{\theta,q} =\big({\B}^{s_0,\sigma}_{p,q_0,{\mathcal{D}}}(\Omega),\D^{s_0}_{p,q_0}(\AA_\mathcal{D})\big)_{\theta,q} = {\B}^{s,\sigma}_{p,q,\mathcal{D}}(\Omega).
\end{align}
Furthermore,  
\begin{itemize}
    \item the result still holds if one replaces $(\cdot,\cdot)_{\theta,\infty}$ and ${\B}^{s,\sigma}_{p,\infty,\mathcal{D}}$ by respectively $(\cdot,\cdot)_{\theta}$ and ${\BesSmo}^{s,\sigma}_{p,\infty,\mathcal{D}}$;
    
    \item in the particular case  of $\Omega$ being a bounded $\C^{1,\alpha}$-domain with small constant, one recalls here that the condition \hyperref[eq:ConditionRegularity]{$(\prescript{\alpha}{r}{\ast}_{p,q}^{s-2})$} becomes $s<1+\alpha+\sfrac{1}{p}$, and $s\leqslant 1+\alpha$ whenever $p=q=\infty$.
\end{itemize}
\end{proposition}

\begin{remark}Carefully note that while the two interpolation identities yield the same output space, both identities are different. Indeed, for $s_0\in(-1+\sfrac{1}{p},\sfrac{1}{p})$ that does \textbf{not} satisfy \hyperref[eq:ConditionRegularity]{$(\prescript{\alpha}{r}{\ast}_{p,q}^{s_0})$}, it might happen that
\begin{align*}
{\H}^{s_0+2,p}_{\mathcal{D},\sigma}(\Omega)\subsetneq &\D^{s_0}_{p}(\AA_\mathcal{D}) \not\subset{\H}^{s_0+2,p}_{\mathcal{D},\sigma}(\Omega),\\
   \text{ and }\quad {\B}^{s_0+2,\sigma}_{p,q_0,\mathcal{D}}(\Omega)\subsetneq &\D^{s_0}_{p,q_0}(\AA_\mathcal{D}) \not\subset{\B}^{s_0+2,\sigma}_{p,q_0,\mathcal{D}}(\Omega).
\end{align*}
\end{remark}

\begin{proof} The identity \eqref{eq:realInterpDomainStokesbddDomain} has been proved during the proof Theorem~\ref{thm:FinalResultBesov2} in the case of Besov spaces. The case of Sobolev spaces follows from the natural embeddings
\begin{align*}
    {\B}^{s_0,\sigma}_{p,1,\mathfrak{n}}(\Omega) \hookrightarrow &{\H}^{s_0,p}_{{\mathfrak{n}},\sigma}(\Omega) \hookrightarrow  {\B}^{s_0,\sigma}_{p,\infty,\mathfrak{n}}(\Omega),\\
    \text{ and }\quad \D^{s_0}_{p,1}(\AA_\mathcal{D})\hookrightarrow&\D^{s_0}_{p}(\AA_\mathcal{D})\hookrightarrow\D^{s_0}_{p,\infty}(\AA_\mathcal{D}).
\end{align*}
From this point, \eqref{eq:realInterpDivFreebddRoughDomain} follows from the reiteration theorem \cite[Theorem~3.5.3]{BerghLofstrom1976} applied to \eqref{eq:realInterpDomainStokesbddDomain} and Proposition~\ref{prop:InterpDivFreeC1}.
\end{proof}

The next result on the half space $\RR^n_+$ was previously known only for indices $s \in \bigl(-1+\sfrac{1}{p},\,\sfrac{1}{p}\bigr)$, $1<p<\infty$, $q\in[1,\infty]$, mostly due to the boundedness of the Hodge--Leray projection (Theorem~\ref{thm:HodgeDecompRn+2}). See for instance \cite{Gaudin2023Hodge}. One can also deduce parts of the identity \eqref{eq:1stInterpolSolenoidalDir} below from the uniqueness of the Cauchy problem and the proof of \cite[Theorem~2.9]{Watanabe2025}, but only under the assumption $s \in \bigl(-1+\sfrac{1}{p},\,\sfrac{1}{p}\bigr)$ (first assuming $q=1$, then applying a reiteration theorem). Note that for its (fully) inhomogeneous counterpart, the following proposition is simple to deduce and well known. However, interpolation of homogeneous function spaces incorporating boundary conditions, even on the flat half-space, is non-trivial, but most importantly poorly investigated since it is a fairly recent topic.

\begin{proposition}\label{prop:InterpHomDivFreeRn+} Let $p,q,q_0,q_1\in[1,\infty]$, and $s_0,s_1\in(-1+{\sfrac{1}{p}},2+{\sfrac{1}{p}})$, $s_0\neq s_1$. For all $\theta\in(0,1)$,  one has
\begin{align}
    (\dot{\H}^{s_0,p}_{{\mathcal{D}},\sigma}(\RR^n_+),\dot{\H}^{s_1,p}_{{\mathcal{D}},\sigma}(\RR^n_+))_{\theta,q} =(\dot{\B}^{s_0,\sigma}_{p,q_0,{\mathcal{D}}}(\RR^n_+),\dot{\B}^{s_1,\sigma}_{p,q_1,{\mathcal{D}}}(\RR^n_+))_{\theta,q} &= \dot{\B}^{s,\sigma}_{p,q,\mathcal{D}}(\RR^n_+)\text{,}\label{eq:1stInterpolSolenoidalDir}
\end{align}
with equivalence of norms, whenever $s,s_0,s_1\neq 1/p$.

\medbreak

\noindent In particular, under the assumption $s_1=s_0+2$, $s_0\in(-1+\sfrac{1}{p},\sfrac{1}{p})$, it holds
\begin{align*}
        (\dot{\H}^{s_0,p}_{{\mathfrak{n}},\sigma}(\RR^n_+),\dot{\D}^{s_0}_{p}(\mathring{\AA}_\mathcal{D}))_{\theta,q} =(\dot{\B}^{s_0,\sigma}_{p,q_0,{\mathcal{D}}}(\RR^n_+),\dot{\D}^{s_0}_{p,q_0}(\mathring{\AA}_\mathcal{D}))_{\theta,q} = \dot{\B}^{s,\sigma}_{p,q,\mathcal{D}}(\RR^n_+).
\end{align*}
Furthermore,
\begin{itemize}
    \item the same result holds replacing $(\cdot,\cdot)_{\theta,\infty}$ and $\dot{\B}^{s,\sigma}_{p,\infty,\mathcal{D}}$ by respectively $(\cdot,\cdot)_{\theta}$ and $\dot{\BesSmo}^{s,\sigma}_{p,\infty,\mathcal{D}}$.
    \item the same result holds replacing $\dot{\B}^{\bullet,\sigma}_{\infty,q,\mathcal{D}}$ and $\dot{\BesSmo}^{\bullet,\sigma}_{\infty,\infty,\mathcal{D}}$ by respectively $\dot{\B}^{\bullet,0,\sigma}_{\infty,q,\mathcal{D}}$ and $\dot{\BesSmo}^{\bullet,0,\sigma}_{\infty,\infty,\mathcal{D}}$;
\end{itemize}
\end{proposition}

\begin{proof} We first notice that the second interpolation identity is a straightforward consequence of Corollary~\ref{cor:HomogeneousDomainStokesRn+} and the first identity \eqref{eq:1stInterpolSolenoidalDir}. So it reduces to prove the first interpolation identity \eqref{eq:1stInterpolSolenoidalDir}.

\medbreak

\noindent The result is already known whenever $-1+1/p<s_0<s<1+1/p$ thanks to Proposition~\ref{prop:HomInterpDivFreeC1}. By Theorem~\ref{thm:InterpHomSpacesBddLip} and the fact that
\begin{align}\label{eq:ProofBesovDivFreeRn+Embedding}
    (\dot{\B}^{s_0,\sigma}_{p,q_0,{\mathcal{D}}}(\RR^n_+),\dot{\B}^{s_1,\sigma}_{p,q_1,{\mathcal{D}}}(\RR^n_+))_{\theta,q} \hookrightarrow \dot{\B}^{s,\sigma}_{p,q,\mathcal{D}}(\RR^n_+),
\end{align}
one only needs to deal with the reverse embedding in the case $1+1/p<s<s_1<2+1/p$.

\textbf{Step 1:} The case $p<\infty$. Due to
\begin{align*}
    (\dot{\B}^{s_0,\sigma}_{p,1,{\mathcal{D}}}(\RR^n_+),\dot{\B}^{s_1,\sigma}_{p,1,{\mathcal{D}}}(\RR^n_+))_{\theta,q}\hookrightarrow(\dot{\B}^{s_0,\sigma}_{p,q_0,{\mathcal{D}}}(\RR^n_+),\dot{\B}^{s_1,\sigma}_{p,q_1,{\mathcal{D}}}(\RR^n_+))_{\theta,q},
\end{align*}
it suffices to prove
\begin{align*}
   \dot{\B}^{s,\sigma}_{p,q,\mathcal{D}}(\RR^n_+) \hookrightarrow(\dot{\B}^{s_0,\sigma}_{p,1,{\mathcal{D}}}(\RR^n_+),\dot{\B}^{s_1,\sigma}_{p,1,{\mathcal{D}}}(\RR^n_+))_{\theta,q}.
\end{align*}
Hence, without loss od generality we can assume $q_0=q_1=1$.

\textbf{Step 1.1:} Assume temporarily \hyperref[AssumptionCompletenessExponents]{ $(\mathcal{C}_{s_0,p,1})$}, $s_0>1/p$ (\textit{i.e.} $1/p<s_0\leqslant n/p$).

Let $\uu\in\dot{\B}^{s,\sigma}_{p,q,\mathcal{D}}(\RR^n_+)$.
We aim to control the $K$-functional of $\uu$, that is, to prove that one has for any $t>0$
\begin{align*}
    K(t,\uu,\dot{\B}^{s_0,\sigma}_{p,1,\mathcal{D}}(\RR^n_+),\dot{\B}^{s_1,\sigma}_{p,1,\mathcal{D}}(\RR^n_+)) \lesssim_{p,s_0,s_1,n} K(t,\uu,\dot{\B}^{s_0,\sigma}_{p,1}(\RR^n_+),\dot{\B}^{s_1,\sigma}_{p,1}(\RR^n_+)).
\end{align*}

\medbreak
 
\noindent By Theorem~\ref{thm:InterpHomSpacesLip}, $\uu\in\dot{\B}^{s,\sigma}_{p,q}(\RR^n_+) \subset  \dot{\B}^{s_0,\sigma}_{p,1}(\RR^n_+) + \dot{\B}^{s_1,\sigma}_{p,1}(\RR^n_+)$. Let $(\mathfrak{a},\mathfrak{b})\in \dot{\B}^{s_0,\sigma}_{p,1}(\RR^n_+)\times\dot{\B}^{s_1,\sigma}_{p,1}(\RR^n_+)$, such that 
\begin{align*}
    \uu=\mathfrak{a}+\mathfrak{b}.
\end{align*}
 The goal is to find $\mathfrak{a}_{\mathcal{D}},\mathfrak{b}_{\mathcal{D}}$ such that $\mathfrak{a}_{\mathcal{D}}{}_{|_{\partial\RR^n_+}}=\mathfrak{b}_{\mathcal{D}}{}_{|_{\partial\RR^n_+}} =0$, $\div \mathfrak{a}_{\mathcal{D}}= \div\mathfrak{b}_{\mathcal{D}} =0$, and
\begin{align*}
    \uu=\mathfrak{a}_{\mathcal{D}}+\mathfrak{b}_{\mathcal{D}}.
\end{align*}
with the estimates
\begin{align*}
    \lVert \mathfrak{a}_{\mathcal{D}}\rVert_{\dot{\B}^{s_0}_{p,1}(\RR^n_+)} \lesssim_{p,s_0,n}\lVert \mathfrak{a}\rVert_{\dot{\B}^{s_0}_{p,1}(\RR^n_+)},\, \text{ and }\lVert \mathfrak{b}_{\mathcal{D}}\rVert_{\dot{\B}^{s_1}_{p,1}(\RR^n_+)} \lesssim_{p,s_1,n}\lVert \mathfrak{b}\rVert_{\dot{\B}^{s_1}_{p,1}(\RR^n_+)}.
\end{align*}
First, since $\mathfrak{a}\in \dot{\B}^{s_0,\sigma}_{p,1}(\RR^n_+)$ with $1/p<s_0\leqslant n/p$, by Theorem~\ref{thm:MetaThmSteadyDirichletStokesRn+}, there exists a unique $(\ww,\mathfrak{q})\in \dot{\B}^{s_0,\sigma}_{p,1,\mathcal{D}}(\RR^n_+)\times \dot{\B}^{s_0-1}_{p,1}(\RR^n_+,\CC^n)$ such that
\begin{equation*}
    \left\{ \begin{array}{rllr}
         - \Delta\ww +\nabla \mathfrak{q} &= 0 \text{, }&&\text{ in } \RR^n_+\text{,}\\
        \div \ww &= 0\text{, } &&\text{ in } \RR^n_+\text{,}\\
        \ww{}_{|_{\partial\RR^n_+}} &=\mathfrak{a}_{|_{\partial\RR^n_+}} = -\mathfrak{b}_{|_{\partial\RR^n_+}}\text{, } &&\text{ on } \partial\RR^n_+\text{.}
    \end{array}
    \right.
\end{equation*}
with  the estimates
\begin{align*}
            \lVert \ww\rVert_{\dot{\B}^{s_0}_{p,1}(\RR^n_+)} + \lVert \mathfrak{q}\rVert_{\dot{\B}^{s_0-1}_{p,1}(\RR^n_+)} \nonumber\lesssim_{p,n,s_0} \lVert  \mathfrak{a}_{|_{\partial\RR^n_+}}\rVert_{\dot{\B}^{s_0-1/p}_{p,1}(\partial\RR^n_+)}\lesssim_{p,n,s_0} \lVert  \mathfrak{a}\rVert_{\dot{\B}^{s_0}_{p,1}(\RR^n_+)}.
\end{align*}
Note that since $\uu_{|_{\partial\RR^n_+}}=0$, it implies $\mathfrak{a}_{|_{\partial\RR^n_+}}=-\mathfrak{b}_{|_{\partial\RR^n_+}}\in\dot{\B}^{s_0-1/p}_{p,1}\cap\dot{\B}^{s_1-1/p}_{p,1}(\partial\RR^n_+)$, which, by Theorem~\ref{thm:MetaThmSteadyDirichletStokesRn+} and the trace theorem, also implies
\begin{align*}
            \lVert \ww\rVert_{\dot{\B}^{s_1}_{p,1}(\RR^n_+)} + \lVert \mathfrak{q}\rVert_{\dot{\B}^{s_1-1}_{p,1}(\RR^n_+)} \nonumber\lesssim_{p,n,s_0} \lVert  \mathfrak{b}_{|_{\partial\RR^n_+}}\rVert_{\dot{\B}^{s_1-1/p}_{p,1}(\partial\RR^n_+)}\lesssim_{p,n,s_1} \lVert  \mathfrak{b}\rVert_{\dot{\B}^{s_1}_{p,1}(\RR^n_+)}.
\end{align*}
Therefore, one may write
\begin{align*}
    \uu =\mathfrak{a}+\mathfrak{b} &= (\mathfrak{a} - \ww) + (\mathfrak{b}+\ww)\\
    &=: \mathfrak{a}_{\mathcal{D}} + \mathfrak{b}_{\mathcal{D}}.
\end{align*}
By construction, $(\mathfrak{a}_{\mathcal{D}},\mathfrak{b}_{\mathcal{D}})\in\dot{\B}^{s_0,\sigma}_{p,1,\mathcal{D}}(\RR^n_+)\times\dot{\B}^{s_1,\sigma}_{p,1,\mathcal{D}}(\RR^n_+)$, with the estimates
\begin{align*}
            &\lVert \mathfrak{a}_{\mathcal{D}}\rVert_{\dot{\B}^{s_0}_{p,1}(\RR^n_+)}=\lVert \mathfrak{a}- \ww\rVert_{\dot{\B}^{s_0}_{p,1}(\RR^n_+)} \lesssim_{p,n,s_0} \lVert  \mathfrak{a}\rVert_{\dot{\B}^{s_0}_{p,1}(\RR^n_+)}\\
            \,\text{ and }&\lVert \mathfrak{b}_{\mathcal{D}}\rVert_{\dot{\B}^{s_1}_{p,1}(\RR^n_+)}=\lVert \mathfrak{b}+ \ww\rVert_{\dot{\B}^{s_1}_{p,1}(\RR^n_+)} \lesssim_{p,n,s_1} \lVert  \mathfrak{b}\rVert_{\dot{\B}^{s_1}_{p,1}(\RR^n_+)}.
\end{align*}
Everything above yields $\uu\in\dot{\B}^{s_0,\sigma}_{p,1,\mathcal{D}}(\RR^n_+)+\dot{\B}^{s_1,\sigma}_{p,1,\mathcal{D}}(\RR^n_+)$, and for $t>0$, by the definition of the $K$-functional
\begin{align*}
    K(t,\uu,\dot{\B}^{s_0,\sigma}_{p,1,\mathcal{D}}(\RR^n_+),\dot{\B}^{s_1,\sigma}_{p,1,\mathcal{D}}(\RR^n_+)) &\leqslant \lVert \mathfrak{a}_{\mathcal{D}}\rVert_{\dot{\B}^{s_0}_{p,1}(\RR^n_+)} +t \lVert \mathfrak{b}_{\mathcal{D}}\rVert_{\dot{\B}^{s_1}_{p,1}(\RR^n_+)}\\
    &\lesssim_{p,s_0,s_1,n} \lVert \mathfrak{a}\rVert_{\dot{\B}^{s_0}_{p,1}(\RR^n_+)} +t \lVert \mathfrak{b}\rVert_{\dot{\B}^{s_1}_{p,1}(\RR^n_+)}.
\end{align*}
Taking the infimum on all such pairs $(\mathfrak{a},\mathfrak{b})$, one deduces for all $t>0$
\begin{align*}
    K(t,\uu,\dot{\B}^{s_0,\sigma}_{p,1,\mathcal{D}}(\RR^n_+),\dot{\B}^{s_1,\sigma}_{p,1,\mathcal{D}}(\RR^n_+)) \lesssim_{p,s_0,s_1,n} K(t,\uu,\dot{\B}^{s_0,\sigma}_{p,1}(\RR^n_+),\dot{\B}^{s_1,\sigma}_{p,1}(\RR^n_+)).
\end{align*}
Now, multiplying by $t^{-\theta}$, and then taking the $\L^q_{\ast}$-norm on both sides, by Theorem~\ref{thm:InterpHomSpacesLip}, one obtains
\begin{align*}
    \lVert \uu \rVert_{(\dot{\B}^{s_0,\sigma}_{p,1,\mathcal{D}}(\RR^n_+),\dot{\B}^{s_1,\sigma}_{p,1,\mathcal{D}}(\RR^n_+))_{\theta,q}}\lesssim_{p,s_0,s_1,n}  \lVert \uu \rVert_{(\dot{\B}^{s_0,\sigma}_{p,1}(\RR^n_+),\dot{\B}^{s_1,\sigma}_{p,1}(\RR^n_+))_{\theta,q}}\sim_{p,s_0,s_1,n} \lVert \uu \rVert_{\dot{\B}^{s}_{p,q}(\RR^n_+)}.
\end{align*}
This gives then the reverse embedding and consequently the result if $1/p <s_0 \leqslant n/p$, $p<\infty$.

\textbf{Step 1.2:} Now the case $s_0>1/p$ without assuming \hyperref[AssumptionCompletenessExponents]{ $(\mathcal{C}_{s_0,p,1})$}. Let $1/p<\alpha<s_0$, such that $\alpha\leqslant n/p$, by the previous Step 1.1
\begin{align*}
    (\dot{\B}^{s_0,\sigma}_{p,1,\mathcal{D}}(\RR^n_+),\dot{\B}^{s_1,\sigma}_{p,1,\mathcal{D}}(\RR^n_+))_{\theta,q} &= \Big((\dot{\B}^{\alpha,\sigma}_{p,1,\mathcal{D}}(\RR^n_+),\dot{\B}^{s_1,\sigma}_{p,1,\mathcal{D}}(\RR^n_+))_{\frac{s_0-\alpha}{s_1-\alpha},1},\dot{\B}^{s_1,\sigma}_{p,1,\mathcal{D}}(\RR^n_+)\Big)_{\theta,q}\\
    &=(\dot{\B}^{\alpha,\sigma}_{p,1,\mathcal{D}}(\RR^n_+),\dot{\B}^{s_1,\sigma}_{p,1,\mathcal{D}}(\RR^n_+))_{\frac{s-\alpha}{s_1-\alpha},q}\\
    &= \dot{\B}^{s,\sigma}_{p,q,\mathcal{D}}(\RR^n_+),
\end{align*}
where we did apply \cite[Lemma~C.1]{Gaudin2023Lip}. This yields then the fulls result whenever $p<\infty$, $1/p<s_0,s,s_1<2+1/p$.

\textbf{Step 1.3:} Now the case $-1+1/p<s_0<1/p$. Let $1/p<s<s_1<2+1/p$, note that if $1/p<\alpha<\min(s,1+1/p)$, we already know
\begin{align*}
    \dot{\B}^{\alpha,\sigma}_{p,q,\mathcal{D}}(\RR^n_+)\hookrightarrow \dot{\B}^{\alpha,\sigma}_{p,q,0}(\RR^n_+) &= (\dot{\B}^{s_0,\sigma}_{p,1,0}(\RR^n_+),\dot{\B}^{s_1,\sigma}_{p,1,0}(\RR^n_+))_{\frac{\alpha-s_0}{s_1-s_0},q}\\
    &\hookrightarrow(\dot{\B}^{s_0,\sigma}_{p,1,\mathcal{D}}(\RR^n_+),\dot{\B}^{s_1,\sigma}_{p,1,\mathcal{D}}(\RR^n_+))_{\frac{\alpha-s_0}{s_1-s_0},q}\\
    &\hookrightarrow \dot{\B}^{\alpha,\sigma}_{p,q,\mathcal{D}}(\RR^n_+).
\end{align*}
By the extremal reiteration property \cite[Lemma~C.1]{Gaudin2023Lip} and Step 1.2, we obtain
\begin{align*}
    \dot{\B}^{s,\sigma}_{p,q,\mathcal{D}}(\RR^n_+)&=(\dot{\B}^{\alpha,\sigma}_{p,q,\mathcal{D}}(\RR^n_+),\dot{\B}^{s_1,\sigma}_{p,1,\mathcal{D}}(\RR^n_+))_{\frac{s-\alpha}{s_1-\alpha},q} \\
    &= \Big((\dot{\B}^{s_0,\sigma}_{p,1,\mathcal{D}}(\RR^n_+),\dot{\B}^{s_1,\sigma}_{p,1,\mathcal{D}}(\RR^n_+))_{\frac{\alpha-s_0}{s_1-s_0},q},\dot{\B}^{s_1,\sigma}_{p,1,\mathcal{D}}(\RR^n_+)\Big)_{\frac{s-\alpha}{s_1-\alpha},q}\\
    &= (\dot{\B}^{s_0,\sigma}_{p,1,\mathcal{D}}(\RR^n_+),\dot{\B}^{s_1,\sigma}_{p,1,\mathcal{D}}(\RR^n_+))_{\frac{s-s_0}{s_1-s_0},q}.
\end{align*}
This ends the case $p<\infty$.

\textbf{Step 2:} The case $p=\infty$. Again, by Theorem~\ref{thm:InterpHomSpacesBddLip} and Proposition~\ref{prop:IdentifVanishingDivFree}, the result is known whenever $-1<s_0<s<1$  and $1<s<s_1<2$.

\textbf{Step 2.1:} The case $0<s_0<s<s_1<2$. Prior to the desired result, we want to prove that, for inhomogeneous function spaces:
\begin{align}\label{eq:proofInterpInftyDivFreecaseInhom}
    (\C^0_{ub,h,0,\sigma}(\RR^n_+),{\B}^{2,\sigma}_{\infty,\infty,\mathcal{D}}(\RR^n_+))_{\theta,q} = \tilde{\B}^{2\theta,\sigma}_{\infty,q,\mathcal{D}}(\RR^n_+)
\end{align}
with equivalence of norms, where, provided $s>0$, and $1\leqslant q\leqslant\infty$, we recall that $\tilde{\B}^{s}_{\infty,q}(\RR^n_+):={\B}^{s}_{\infty,q}\cap\S'_h(\RR^n_+)$ equipped with the inhomogeneous ${\B}^{s}_{\infty,q}$-norm. Recall that, as sets, one has the equality $\tilde{\B}^{s}_{\infty,q}(\RR^n_+)=\dot{\B}^{s}_{\infty,q}(\RR^n_+)$, $s>0$, $1\leqslant q\leqslant\infty$.

We notice that one trivially obtains an embedding
\begin{align*}
    (\C^0_{ub,h,0,\sigma}(\RR^n_+),{\B}^{2,\sigma}_{\infty,\infty,\mathcal{D}}(\RR^n_+))_{\theta,q} \hookrightarrow {\B}^{2\theta,\sigma}_{\infty,q,\mathcal{D}}(\RR^n_+).
\end{align*}
However, by \cite[Theorem~3.4.2]{BerghLofstrom1976}, if $q<\infty$, $\tilde{\B}^{2,\sigma}_{\infty,\infty,\mathcal{D}}(\RR^n_+)=\C^0_{ub,h,0,\sigma}(\RR^n_+)\cap{\B}^{2,\sigma}_{\infty,\infty,\mathcal{D}}(\RR^n_+)$ is dense in $(\C^0_{ub,h,0,\sigma}(\RR^n_+),{\B}^{2,\sigma}_{\infty,\infty,\mathcal{D}}(\RR^n_+))_{\theta,q}$, and the closure of  $\tilde{\B}^{2}_{\infty,\infty}(\RR^n_+)$ in ${\B}^{2\theta}_{\infty,q}(\RR^n_+)$ is $\tilde{\B}^{2\theta}_{\infty,q}(\RR^n_+)$.
Both together yields
\begin{align*}
    (\C^0_{ub,h,0,\sigma}(\RR^n_+),{\B}^{2,\sigma}_{\infty,\infty,\mathcal{D}}(\RR^n_+))_{\theta,q} \hookrightarrow {\B}^{2\theta,\sigma}_{\infty,q,\mathcal{D}}\cap\tilde{\B}^{2\theta}_{\infty,q} (\RR^n_+) =\tilde{\B}^{2\theta,\sigma}_{\infty,q,\mathcal{D}}(\RR^n_+).
\end{align*}
The reiteration theorem yields inclusion in $\S'_h(\RR^n_+)$ when $q=\infty$, and therefore the corresponding embedding. 

For the reverse embedding for $2\theta<1$, one has
\begin{align*}
    \tilde{\B}^{2\theta,\sigma}_{\infty,q,\mathcal{D}}(\RR^n_+) = \tilde{\B}^{2\theta,\sigma}_{\infty,q,0}(\RR^n_+) &= (\C^0_{ub,h,0,\sigma}(\RR^n_+),\tilde{\B}^{2,\sigma}_{\infty,\infty,0}(\RR^n_+))_{\theta,q}\\ &\hookrightarrow (\C^0_{ub,h,0,\sigma}(\RR^n_+),{\B}^{2,\sigma}_{\infty,\infty,\mathcal{D}}(\RR^n_+))_{\theta,q},
\end{align*}
giving \eqref{eq:proofInterpInftyDivFreecaseInhom}, for $2\theta<1$, $1\leqslant q\leqslant \infty$.

For the reverse embedding for $2\theta>1$, according to Theorem~\ref{thm:StokesDirRn+Linfty}, we can introduce the unbounded invertible $0$-sectorial operator on $\C^0_{ub,h,0,\sigma}(\RR^n_+)$:
\begin{align*}
    \D(\widetilde{\AA}_\mathcal{D})&:= \D_\infty({\AA_\mathcal{D}})\cap\C^0_{ub,h,0,\sigma}(\RR^n_+)\hookrightarrow \tilde{\BesSmo}^{2,\sigma}_{\infty,\infty,\mathcal{D}}(\RR^n_+),\,\\
    \widetilde{\AA}_\mathcal{D}\vv&:= (\I +{\AA_\mathcal{D}}) \vv = \vv - \Delta \vv + \nabla \mathfrak{q},\qquad\qquad\qquad \vv\in\D(\widetilde{\AA}_\mathcal{D}).
\end{align*}
By previous considerations,
\begin{align*}
    (\C^0_{ub,h,0,\sigma}(\RR^n_+),\D(\widetilde{\AA}_\mathcal{D}))_{\theta-\frac{1}{2},q}=\tilde{\B}^{2\theta-1,\sigma}_{\infty,q,\mathcal{D}}(\RR^n_+).
\end{align*}
By Proposition~\ref{prop:BesovInftySqrtDirStokesRn+}, one has the isomorphism
\begin{align*}
    \widetilde{\AA}_\mathcal{D}^{-\frac{1}{2}}\tilde{\B}^{2\theta-1,\sigma}_{\infty,q,\mathcal{D}}(\RR^n_+) = \tilde{\B}^{2\theta,\sigma}_{\infty,q,\mathcal{D}}(\RR^n_+),
\end{align*}
and the isomorphism
\begin{align*}
    \widetilde{\AA}_\mathcal{D}^{\frac{1}{2}}(\C^0_{ub,h,0,\sigma}(\RR^n_+),\D(\widetilde{\AA}_\mathcal{D}))_{\theta,q} = (\C^0_{ub,h,0,\sigma}(\RR^n_+),\D(\widetilde{\AA}_\mathcal{D}))_{\theta-\frac{1}{2},q},
\end{align*}
the latter being given by \cite[Proposition~6.4.1,~a),~\&~Corollary~6.5.5]{bookHaase2006}. Therefore, we deduce
\begin{align*}
    (\C^0_{ub,h,0,\sigma}(\RR^n_+),\D(\widetilde{\AA}_\mathcal{D}))_{\theta,q} &=  \widetilde{\AA}_\mathcal{D}^{-\frac{1}{2}}(\C^0_{ub,h,0,\sigma}(\RR^n_+),\D(\widetilde{\AA}_\mathcal{D}))_{\theta-\frac{1}{2},q} \\
    &= \widetilde{\AA}_\mathcal{D}^{-\frac{1}{2}}\tilde{\B}^{2\theta-1,\sigma}_{\infty,q,\mathcal{D}}(\RR^n_+)\\
    &=  \tilde{\B}^{2\theta,\sigma}_{\infty,q,\mathcal{D}}(\RR^n_+).
\end{align*}
By the reiteration Theorem \cite[Theorem~3.5.4]{BerghLofstrom1976}, for any $0<s_0<s<s_1<2$, we obtain
\begin{align*}
    (\tilde{\B}^{s_0,\sigma}_{\infty,q_0,\mathcal{D}}(\RR^n_+),&\tilde{\B}^{s_1,\sigma}_{\infty,q_1,\mathcal{D}}(\RR^n_+))_{\frac{s-s_0}{s_1-s_0},q}\\
    &= \big((\C^0_{ub,h,0,\sigma}(\RR^n_+),\D(\widetilde{\AA}_\mathcal{D}))_{\frac{s_0}{2},q},(\C^0_{ub,h,0,\sigma}(\RR^n_+),\D(\widetilde{\AA}_\mathcal{D}))_{\frac{s_1}{2},q}\big)_{\frac{s-s_0}{s_1-s_0},q}\\
    &= (\C^0_{ub,h,0,\sigma}(\RR^n_+),\D(\widetilde{\AA}_\mathcal{D}))_{\frac{s}{2},q}\\
    &= \tilde{\B}^{s}_{\infty,q,\mathcal{D}}(\RR^n_+).
\end{align*}
Therefore, combining it with \eqref{eq:ProofBesovDivFreeRn+Embedding}, one obtains for $0<s_0<s<s_1<2$
\begin{align*}
    \tilde{\B}^{s,\sigma}_{\infty,q,\mathcal{D}}(\RR^n_+) &= (\tilde{\B}^{s_0,\sigma}_{\infty,q_0,\mathcal{D}}(\RR^n_+),\tilde{\B}^{s_1,\sigma}_{\infty,q_1,\mathcal{D}}(\RR^n_+))_{\theta,q}\\
    &\hookrightarrow (\dot{\B}^{s_0,\sigma}_{\infty,q_0,\mathcal{D}}(\RR^n_+),\dot{\B}^{s_1,\sigma}_{\infty,q_1,\mathcal{D}}(\RR^n_+))_{\theta,q}\\
    &\hookrightarrow\dot{\B}^{s,\sigma}_{\infty,q,\mathcal{D}}(\RR^n_+).
\end{align*}
Recall that algebraically $\tilde{\B}^{s,\sigma}_{\infty,q,\mathcal{D}}(\RR^n_+)=\dot{\B}^{s,\sigma}_{\infty,q,\mathcal{D}}(\RR^n_+)$, since $s>0$, so that one has for any $\uu\in \dot{\B}^{s}_{\infty,q,\mathcal{D}}(\RR^n_+)$,
\begin{align*}
    \lVert \uu \rVert_{\dot{\B}^{s}_{\infty,q}(\RR^n_+)} \lesssim_{s,n} \lVert \uu \rVert_{(\dot{\B}^{s_0,\sigma}_{\infty,q_0,\mathcal{D}}(\RR^n_+),\dot{\B}^{s_1,\sigma}_{\infty,q_1,\mathcal{D}}(\RR^n_+))_{\theta,q}} &\lesssim_{s,n} \lVert \uu \rVert_{{\B}^{s}_{\infty,q}(\RR^n_+)}\\ &\lesssim_{s,n} \lVert \uu \rVert_{{\L}^{\infty}(\RR^n_+)}+  \lVert \uu \rVert_{\dot{\B}^{s}_{\infty,q}(\RR^n_+)}.
\end{align*}
For $\lambda>0$, we consider the dilation $\uu_\lambda = \uu(\lambda\cdot)$, so that
\begin{align}\label{eq:ProofrealInterpLInftydivFreeRn+}
    \lambda^{s}\lVert \uu \rVert_{\dot{\B}^{s}_{\infty,q}(\RR^n_+)} \lesssim_{s,n} \lVert \uu_\lambda \rVert_{(\dot{\B}^{s_0,\sigma}_{\infty,q_0,\mathcal{D}}(\RR^n_+),\dot{\B}^{s_1,\sigma}_{\infty,q_1,\mathcal{D}}(\RR^n_+))_{\theta,q}} &\lesssim_{s,n} \lVert \uu \rVert_{{\B}^{s}_{\infty,q}(\RR^n_+)}\\ &\lesssim_{s,n} \lVert \uu \rVert_{{\L}^{\infty}(\RR^n_+)}+  \lambda^s\lVert \uu \rVert_{\dot{\B}^{s}_{\infty,q}(\RR^n_+)}\nonumber.
\end{align}
Due to the dilation invariance of $\dot{\B}^{s_0,\sigma}_{\infty,q_0,\mathcal{D}}(\RR^n_+)$ and $\dot{\B}^{s_1,\sigma}_{\infty,q_0,\mathcal{D}}(\RR^n_+)$, and homogeneity of the $K$-functional, it can be checked that
\begin{align*}
    \lVert \uu_\lambda \rVert_{(\dot{\B}^{s_0,\sigma}_{\infty,q_0,\mathcal{D}}(\RR^n_+),\dot{\B}^{s_1,\sigma}_{\infty,q_1,\mathcal{D}}(\RR^n_+))_{\theta,q}} \sim_{s_0,s_1,n,s} \lambda^s\lVert \uu \rVert_{(\dot{\B}^{s_0,\sigma}_{\infty,q_0,\mathcal{D}}(\RR^n_+),\dot{\B}^{s_1,\sigma}_{\infty,q_1,\mathcal{D}}(\RR^n_+))_{\theta,q}}.
\end{align*}
Thus, diving by $\lambda^s$, \eqref{eq:ProofrealInterpLInftydivFreeRn+} becomes
\begin{align*}
    \lVert \uu \rVert_{\dot{\B}^{s}_{\infty,q}(\RR^n_+)} \lesssim_{s,s_0,s_1,n} \lVert \uu \rVert_{(\dot{\B}^{s_0,\sigma}_{\infty,q_0,\mathcal{D}}(\RR^n_+),\dot{\B}^{s_1,\sigma}_{\infty,q_1,\mathcal{D}}(\RR^n_+))_{\theta,q}} &\lesssim_{s,s_0,s_1,n}  \frac{1}{\lambda^s}\lVert \uu \rVert_{{\L}^{\infty}(\RR^n_+)}+  \lVert \uu \rVert_{\dot{\B}^{s}_{\infty,q}(\RR^n_+)}.
\end{align*}
Taking the limit as $\lambda$ goes to infinity, it yields
\begin{align*}
    \lVert \uu \rVert_{\dot{\B}^{s}_{\infty,q}(\RR^n_+)} \lesssim_{s,s_0,s_1,n} \lVert \uu \rVert_{(\dot{\B}^{s_0,\sigma}_{\infty,q_0,\mathcal{D}}(\RR^n_+),\dot{\B}^{s_1,\sigma}_{\infty,q_1,\mathcal{D}}(\RR^n_+))_{\theta,q}} &\lesssim_{s,s_0,s_1,n}  \lVert \uu \rVert_{\dot{\B}^{s}_{\infty,q}(\RR^n_+)},
\end{align*}
so that algebraically and with equivalence of norms, it holds for all $0<s_0<s<s_1<2$, and all $q_0,q_1,q\in[1,\infty]$,
\begin{align*}
    (\dot{\B}^{s_0,\sigma}_{\infty,q_0,\mathcal{D}}(\RR^n_+),\dot{\B}^{s_1,\sigma}_{\infty,q_1,\mathcal{D}}(\RR^n_+))_{\theta,q} = \dot{\B}^{s,\sigma}_{\infty,q,\mathcal{D}}(\RR^n_+).
\end{align*}
The same arguments yield algebraically and with equivalence of norms, for any $0<s<s_1<2$, and any $q_1,q\in[1,\infty]$,
\begin{align*}
    (\C^0_{ub,h,0,\sigma}(\RR^n_+),\dot{\B}^{s_1,\sigma}_{\infty,q_1,\mathcal{D}}(\RR^n_+))_{\theta,q} = \dot{\B}^{s,\sigma}_{\infty,q,\mathcal{D}}(\RR^n_+).
\end{align*}
\textbf{Step 2.2:} Now, it remains to deal with the case $-1<s_0<0<s<s_1<2$. But this amounts to reproduce Step 1.2 just setting $p=\infty$, and taking advantage of Step 2.1 instead of Step 1.1, which finishes the proof.
\end{proof}

\subsection{Global-in-time maximal regularity results in the standard settings}\label{sec:StandardMaxRegStokes}

In this whole section, in order to alleviate the notation, we do not explicitly write the dependency of the domain $\Omega$ and the time interval in the estimates, since the domain is fixed, and one can always be reduce our result to the well-posedness globally-in-time .

We start with global-in-time maximal regularity on $\RR^n_+$. The result for UMD Banach spaces reads as follows:

\begin{theorem}\label{thm:UMDMaxRegRn+} Let $p,q\in(1,\infty)$, $s\in(-1+\sfrac{1}{p},\sfrac{1}{p})$ and $\beta\in(-1+\sfrac{1}{q},\sfrac{1}{q})$ such that $s+2+2\beta-\sfrac{2}{q}\neq\sfrac{1}{p}$. For simplicity, we set $\beta_q:=1+\beta-{\sfrac{1}{q}}$.

\medbreak

 Let $T\in(0,\infty]$. For all $\ff\in\dot{\H}^{\beta,q}(0,T;\dot{\H}^{s,p}(\RR^n_+,\CC^n))$, all $\uu_0\in\dot{\B}^{s+2\beta_q,\sigma}_{p,q,\mathcal{D}}(\RR^n_+)$, the problem
\begin{equation*}\tag{DS}
    \left\{ \begin{array}{rllr}
         \partial_t \uu - \Delta \uu +\nabla \mathfrak{p} &= \ff \text{, }&&\text{ in } (0,T)\times{\RR}^n_+\text{,}\\
        \div \uu &= 0\text{, } &&\text{ in } (0,T)\times{\RR}^n_+\text{,}\\
        \uu_{|_{\partial{\RR}^n_+}} &=\mathbf{0}\text{, } &&\text{ on } (0,T)\times\partial{\RR}^n_+\text{,}\\
        \uu(0) &=\uu_0\text{, } &&\text{ in } \dot{\B}^{s+2\beta_q}_{p,q}(\RR^n_+,\CC^n)\text{,}
    \end{array}
    \right.
\end{equation*}
admits a unique (mild) solution $(\uu,\nabla \mathfrak{p})\in \C_{ub}^0([0,T];\dot{\B}^{s+2\beta_q}_{p,q}(\RR^n_+,\CC^n))\times \dot{\H}^{\beta,q}(0,T;\dot{\H}^{s,p}({\RR}^n_+,\CC^n))$, such that $\partial_t \uu ,\nabla^2\uu, \in\dot{\H}^{\beta,q}(0,T;\dot{\H}^{s,p}({\RR}^n_+,\CC^n))$, and satisfying estimates
\begin{align*}
    \lVert \uu \rVert_{\L^{\infty}(\dot{\B}^{s+2\beta_q}_{p,q})}\less_{p,s,n}^{\beta,q} &\lVert  (\partial_t\uu,\nabla^2\uu,\nabla\mathfrak{p})\rVert_{\dot{\H}^{\beta,q}(\dot{\H}^{s,p})}\less_{p,s,n}^{\beta,q} \lVert \ff \rVert_{\dot{\H}^{\beta,q}(\dot{\H}^{s,p})} + \lVert \uu_0 \rVert_{\dot{\B}^{s+2\beta_q}_{p,q}}.
\end{align*}
The whole result still holds if one replaces $\dot{\H}^{s,p}$ by $\dot{\B}^{s}_{p,r}$, $1<p,r<\infty$.
\end{theorem}

Note that obviously $\dot{\H}^{0,q}=\L^q$. 

\begin{proof} Up to apply the Leray projection and to modify the pressure term, one can consider $\ff\in\dot{\H}^{\beta,q}(0,T;\dot{\H}^{s,p}_{\mathfrak{n},\sigma}(\RR^n_+))$, so that the problem reduces to the analysis of the Stokes--Dirichlet operator on $\dot{\H}^{s,p}_{\mathfrak{n},\sigma}(\RR^n_+)$. To obtain Homogeneous Sobolev global-in-time maximal regularity, we did apply \cite[Theorem~4.7]{Gaudin2023}, recalling that both $\dot{\H}^{s,p}_{\mathfrak{n},\sigma}(\RR^n_+)$ and $\dot{\B}^{s,\sigma}_{p,r,\mathfrak{n}}(\RR^n_+)$ are UMD Banach spaces, $1<p,r<\infty$, $s\in(-1+\sfrac{1}{p},\sfrac{1}{p})$, and the Stokes--Dirichlet operator has bounded $\mathbf{H}^\infty$-functional calculus by Proposition~\ref{prop:MetaPropDirichletStokesRn+2}. The trace space for initial data being given by Corollary~\ref{cor:HomogeneousDomainStokesRn+} and Proposition~\ref{prop:InterpHomDivFreeRn+}.
\end{proof}

\medbreak

Now, we provide the Da Prato--Grisvard maximal regularity result on $\RR^n_+$, by applying Theorem~\ref{thm:DaPratoGrisvard}.

\begin{theorem}\label{thm:DaPratoGrisvardRn+} Let $p,q,r\in[1,\infty]$, and $s\in(-1+\sfrac{1}{p},\sfrac{1}{p})$ such that $s+2-\sfrac{2}{q}\neq\sfrac{1}{p}$. We assume either
\begin{itemize}
    \item $r,q\in(1,\infty)$; or
    \item $r=1$, $q\in[1,\infty)$; or
    \item $r=\infty$, $q\in(1,\infty]$.
\end{itemize}

\medbreak

 Let $T\in(0,\infty]$. For all $\ff\in\L^q(0,T;\dot{\B}^{s}_{p,r}(\RR^n_+,\CC^n))$, all $\uu_0\in\dot{\B}^{s+2-\sfrac{2}{q},\sigma}_{p,q,\mathcal{D}}(\RR^n_+)$, the problem
\begin{equation*}\tag{DS}
    \left\{ \begin{array}{rllr}
         \partial_t \uu - \Delta \uu +\nabla \mathfrak{p} &= \ff \text{, }&&\text{ in } (0,T)\times\RR^n_+\text{,}\\
        \div \uu &= 0\text{, } &&\text{ in } (0,T)\times\RR^n_+\text{,}\\
        \uu_{|_{\partial\Omega}} &=\mathbf{0}\text{, } &&\text{ on } (0,T)\times\partial\RR^n_+\text{,}\\
        \uu(0) &=\uu_0\text{, } &&\text{ in } \dot{\B}^{s+2-\sfrac{2}{q}}_{p,q}(\RR^n_+,\CC^n)\text{,}
    \end{array}
    \right.
\end{equation*}
admits a unique (mild) solution $(\uu,\nabla\mathfrak{p})\in \C^{0}_{ub}([0,T];\dot{\B}^{s+2-\sfrac{2}{q}}_{p,q}(\RR^n_+,\CC^n))\times\L^q(0,T;\dot{\B}^{s}_{p,r}({\Omega},\CC^n))$, such that $\partial_t \uu ,\nabla^2\uu, \nabla \mathfrak{p}\in\L^q(0,T;\dot{\B}^{s}_{p,r}(\RR^n_+,\CC^n))$, with estimates
\begin{align*}
    \lVert \uu \rVert_{\L^{\infty}(\dot{\B}^{s+2-\sfrac{2}{q}}_{p,q})}\less_{p,s,n}^{q} \lVert  (\partial_t\uu,\nabla^2\uu,\nabla\mathfrak{p})\rVert_{\L^q(\dot{\B}^{s}_{p,r})} \less_{p,s,n}^{q} \lVert \ff \rVert_{\L^q(\dot{\B}^{s}_{p,r})} + \lVert \uu_0 \rVert_{\dot{\B}^{s+2-\sfrac{2}{q}}_{p,q}}.
\end{align*}
When $q=\infty$, one only has weak-$\ast$ continuity in time.
\end{theorem}

\begin{proof} When $p,r\in(1,\infty)$, the result was already obtained in Theorem~\ref{thm:UMDMaxRegRn+}. If $q=r\in[1,\infty]$, this follows from the Da Prato--Grisvard Theorem~\ref{thm:DaPratoGrisvard}, and the characterization of the homogeneous domains of the Stokes--Dirichlet operator and the interpolation spaces Corollary~\ref{cor:HomogeneousDomainStokesRn+} and Proposition~\ref{prop:InterpHomDivFreeRn+}.

We mention the fact that $\L^1$-in-time maximal regularity implies $\L^q$-maximal regularity for all $q\in[1,\infty)$ by \cite[Theorem~2.13]{BechtelBuiKunstmann2024}. The same result gives similarly that $\L^\infty$-maximal regularity implies $\L^q$-maximal regularity for all $q\in(1,\infty]$.
\end{proof}

Now, we proceed to the statement on bounded domains with minimal regularity. First the UMD case. Before stating our next result, and according to Proposition~\ref{prop:InterpDivFreeC1a} we set the following for all $p,q,r\in[1,\infty]$, $s\in(-1+\sfrac{1}{p},\sfrac{1}{p})$, all $\theta\in(0,1)$
\begin{align*}
    {\mathcal{D}}^{p}_{\AA_\mathcal{D}}\Big(\frac{s}{2}+\theta,\,q\Big):=(\B^{s,\sigma}_{p,r,\mathfrak{n}}(\Omega),\D^{s}_{p,r}(\AA_\mathcal{D}))_{\theta,q}.
\end{align*}
By  Proposition~\ref{prop:InterpDivFreeC1a}, ${\mathcal{D}}^{p}_{\AA_\mathcal{D}}\Big(\frac{s}{2}+\theta,\,q\Big)$ is well-defined and independent of $r\in[1,\infty]$.

The next theorem, Theorem~\ref{thm:UMDMaxRegBddDomain}, is obtained through the application of \cite[Theorem~4.7]{Gaudin2023}, where the invertibility of the operator allows to consider inhomogeneous function spaces. For the intermediate regularity (mixed derivative) estimates see for instance the work by Pr\"{u}ss \cite{Pruss2002}, and Pr\"{u}ss and Simonett \cite[Corollary~4.5.10]{PrussSimonett2016}.

Otherwise the proofs for the next two results are exactly the same as the ones of  Theorems~\ref{thm:UMDMaxRegRn+}~and~\ref{thm:DaPratoGrisvardRn+}, taking advantage of Theorem~\ref{thm:FinalResultSobolev1} and Proposition~\ref{prop:InterpDivFreeC1a}.

\begin{theorem}\label{thm:UMDMaxRegBddDomain}Let $\alpha\in(0,1)$, $r\in[1,\infty]$, and let $\Omega$ be a bounded Lipschitz domain of class $\mathcal{M}^{1+\alpha,r}_\W(\epsilon)$ for some sufficiently small $\epsilon>0$.

\medbreak

Let $p,q\in(1,\infty)$, $s\in(-1+\sfrac{1}{p},\sfrac{1}{p})$ and $\beta\in(-1+\sfrac{1}{q},\sfrac{1}{q})$ such that $s+2+2\beta-\sfrac{2}{q}\neq \sfrac{1}{p}$. For simplicity, set $\beta_q:=1+\beta-\sfrac{1}{q}$.

\medbreak

 Let $T\in(0,\infty]$. For all $\ff\in\H^{\beta,q}(0,T;\H^{s,p}({\Omega},\CC^n))$, all $\uu_0\in{\mathcal{D}}^{s,p}_{\AA_\mathcal{D}}(\sfrac{s}{2}+\beta_q,q)$, the problem
\begin{equation*}\tag{DS}
    \left\{ \begin{array}{rllr}
         \partial_t \uu - \Delta \uu +\nabla \mathfrak{p} &= \ff \text{, }&&\text{ in } (0,T)\times\Omega\text{,}\\
        \div \uu &= 0\text{, } &&\text{ in } (0,T)\times\Omega\text{,}\\
        \uu_{|_{\partial\Omega}} &=\mathbf{0}\text{, } &&\text{ on } (0,T)\times\partial\Omega\text{,}\\
        \uu(0) &=\uu_0\text{, } &&\text{ in } {\mathcal{D}}^{s,p}_{\AA_\mathcal{D}}(\beta_q,q)\text{,}
    \end{array}
    \right.
\end{equation*}
admits a unique (mild) solution $(\uu,\mathfrak{p})$ such that $\uu\in \C_{ub}^0([0,T];{\mathcal{D}}^{p}_{\AA_\mathcal{D}}(\sfrac{s}{2}+\beta_q,q))$, and $\partial_t \uu$, $-\Delta \uu+\nabla\mathfrak{p}\in\H^{\beta,q}(0,T;\H^{s,p}({\Omega},\CC^n))$, with estimates
\begin{align*}
    \lVert \uu \rVert_{\L^{\infty}(\mathcal{D}_{\AA_\mathcal{D}}^p(\sfrac{s}{2}+\beta_q,q))}\less_{p,s,n}^{\Omega,\beta,q} \lVert  (\uu,\partial_t\uu,-\Delta\uu+\nabla\mathfrak{p})&\rVert_{\H^{\beta,q}(\H^{s,p})}\\
    &\less_{p,s,n}^{\Omega,\beta,q} \lVert \ff \rVert_{\H^{\beta,q}(\H^{s,p})} + \lVert \uu_0 \rVert_{{\mathcal{D}}^{p}_{\AA_\mathcal{D}}(\sfrac{s}{2}+\beta_q,q)}.
\end{align*}
Additionally, 
\begin{itemize}
    \item for any $\uptau \in[0,1]$ such that \hyperref[eq:ConditionRegularity]{$(\prescript{\alpha}{r}{\ast}_{p,q}^{s-2+2\uptau})$}, one has $\uu\in \H^{1+\beta-\uptau,q}(0,T;\H^{s+2\uptau,p}(\Omega,\CC^n))$ and the estimate
    \begin{align*}
        \lVert  \uu \rVert_{\H^{1+\beta-\uptau,q}(\H^{s+2\uptau,p})}\less_{p,s,n}^{\Omega,\beta,q} \lVert \ff \rVert_{\H^{\beta,q}(\H^{s,p})} + \lVert \uu_0 \rVert_{{\mathcal{D}}^{p}_{\AA_\mathcal{D}}(\sfrac{s}{2}+\beta_q,q)}.
    \end{align*}
    \item for any $\beta_{o}\in[0,\beta_q]$ such that \hyperref[eq:ConditionRegularity]{$(\prescript{\alpha}{r}{\ast}_{p,q}^{s+2\beta_0-2})$}, one has $\uu\in\C_{ub}^0([0,T];\B^{s+2\beta_o}_{p,q}(\Omega,\CC^n))$, with the estimate
    \begin{align*}
        \lVert \uu \rVert_{\L^{\infty}(\B^{s+2\beta_o}_{p,q})}\less_{p,s,n}^{\Omega,\beta,q}\lVert \uu \rVert_{\L^{\infty}(\mathcal{D}^p_{\AA_\mathcal{D}}(\sfrac{s}{2}+\beta_q,q))},
    \end{align*}
    \item a similar result holds, if one replaces $\H^{s,p}$ by $\B^{s}_{p,\kappa}$, $1<p,\kappa<\infty$.
\end{itemize}
\end{theorem}

Now, we consider its Da Prato--Grisvard counterpart.

\begin{theorem}\label{thm:bounded}Let $\alpha\in(0,1)$, $r\in[1,\infty]$, and let $\Omega$ be a bounded Lipschitz domain of class $\mathcal{M}^{1+\alpha,r}_\W(\epsilon)$ for some sufficiently small $\epsilon>0$.

\medbreak

Let $p,q,\kappa\in[1,\infty]$, $s\in(-1+\sfrac{1}{p},\sfrac{1}{p})$. We assume that $s+2-\sfrac{2}{q}\neq \sfrac{1}{p}$,  and that either
\begin{itemize}
    \item $\kappa,q\in(1,\infty)$; or
    \item $\kappa=1$, $q\in[1,\infty)$; or
    \item $\kappa=\infty$, $q\in(1,\infty]$.
\end{itemize}
For simplicity, we set $\tilde{s}_q:=\frac{s}{2}+1-\frac{1}{q}$.

\medbreak

 Let $T\in(0,\infty]$. For all $\ff\in\L^{q}(0,T;\B^{s}_{p,\kappa}({\Omega},\CC^n))$, all $\uu_0\in{\mathcal{D}}^{p}_{\AA_\mathcal{D}}(\tilde{s}_q,q)$, the problem
\begin{equation*}\tag{DS}
    \left\{ \begin{array}{rllr}
         \partial_t \uu - \Delta \uu +\nabla \mathfrak{p} &= \ff \text{, }&&\text{ in } (0,T)\times\Omega\text{,}\\
        \div \uu &= 0\text{, } &&\text{ in } (0,T)\times\Omega\text{,}\\
        \uu_{|_{\partial\Omega}} &=\mathbf{0}\text{, } &&\text{ on } (0,T)\times\partial\Omega\text{,}\\
        \uu(0) &=\uu_0\text{, } &&\text{ in }  {\mathcal{D}}^{p}_{\AA_\mathcal{D}}(\tilde{s}_q,q)\text{,}
    \end{array}
    \right.
\end{equation*}
admits a unique (mild) solution $(\uu,\mathfrak{p})$ such that $\uu\in \C_{ub}^0([0,T];{\mathcal{D}}^{p}_{\AA_\mathcal{D}}(\tilde{s}_q,q))$, and $\partial_t \uu$, $-\Delta \uu+\nabla\mathfrak{p}\in\L^q(0,T;\B^{s}_{p,\kappa}({\Omega},\CC^n))$, with the estimates
\begin{align*}
    \lVert \uu \rVert_{\L^{\infty}({\mathcal{D}}^{p}_{\AA_\mathcal{D}}(\tilde{s}_q,\,q))}\less_{p,s,n}^{\Omega,\theta,q} \lVert  (\uu,\partial_t\uu,-\Delta\uu+\nabla\mathfrak{p})\rVert_{\L^{q}(\B^{s}_{p,\kappa})}\less_{p,s,n}^{\Omega,\theta,q} \lVert \ff \rVert_{\L^{q}(\B^{s}_{p,\kappa})} + \lVert \uu_0 \rVert_{{\mathcal{D}}^{p}_{\AA_\mathcal{D}}(\tilde{s}_q,q)}.
\end{align*}
When $q=\infty$, one only has weak-$\ast$ continuity.

Additionally, provided \hyperref[eq:ConditionRegularity]{$(\prescript{\alpha}{r}{\ast}_{p,\kappa}^{s-2})$} is satisfied, one has $\nabla^2\uu,\,\nabla \mathfrak{p}\in \L^{q}(0,T;\B^{s}_{p,\kappa}(\Omega,\CC^n))$ and the estimate
    \begin{align*}
        \lVert  -\Delta \uu + \nabla \mathfrak{p} \rVert_{\L^{q}(\B^{s}_{p,\kappa})}\sim_{p,s,n}^{\Omega,\beta,q}  \lVert  (\uu,\,\nabla^2\uu,\, \nabla \mathfrak{p}) \rVert_{\L^{q}(\B^{s}_{p,\kappa})}.
    \end{align*}
\end{theorem}

\subsection{Global-in-time maximal regularity results: a possible way to use the \texorpdfstring{$\L^\infty$}{Linfty}-theory}\label{Sec:MAxRegforLinfty}

While establishing an $\L^\infty$-theory for the semigroup of Stokes--Dirichlet (or any other elliptic) operator is interesting in itself, it turns out to be inapplicable in practice in the standard setting built for evolution equations. Indeed, $\L^\infty$ is neither UMD (it is not reflexive, hence not UMD), and also not the interior point of some non-trivial real interpolation scale. Therefore, standard $\L^q$-maximal regularity and results such as Theorem~\ref{thm:DaPratoGrisvard} are out of reach.

However, we provide here a not so much well-known result (consequence) of the Da Prato--Grisvard Theory, that allows to deal with an arbitrary $\omega$-sectorial operator, $\omega\in[0,\frac{\pi}{2})$, $(\D(A),A)$ on an arbitrary Banach space $\X$. However, the price to pay, is that one has to leave the $\L^q$-scale to measure regularity in-time, and have then to deal with Besov (Sobolev-Soblodeckij) in-time regularity, however it will allow to have $q=1,\infty$. We recall the coincidental (it is actually a definition) identity $\B^{\beta}_{q,q}(0,T;\X)=\W^{\beta,q}(0,T;\X)$ for any $\beta\in(0,\infty)\setminus\NN$, any $q\in[1,\infty]$. In particular, we can recover H\"{o}lder-in-time maximal regularity results !

\begin{theorem}[Da Prato--Grisvard]\label{thm:DaPratoGrisvard2}Let $\omega\in[0,\frac{\pi}{2})$, and $(\D(A),A)$ be an invertible $\omega$-sectorial operator on a Banach space $\X$. Let $q,\kappa\in[1,\infty]$, $\beta\in(-1+\sfrac{1}{q},\sfrac{1}{q})$. For simplicity, we set $\beta_q:=1+\beta-\sfrac{1}{q}\in(0,1)$. 

Let $T\in(0,\infty]$. For any $f\in\B^{\beta}_{q,\kappa}(0,T;\X)$ and any $u_0\in\mathcal{D}_{A}(\beta_q,\kappa)$ the problem
\begin{align}\tag{ACP}\label{ACPBesov}
    \left\{\begin{array}{rl}
            \partial_t u(t) +Au(t)  =& f(t) \,\text{, } 0<t<T,\\
            u(0) =& u_0\text{,}
    \end{array}
    \right.
\end{align}
admits a unique (mild) solution $u\in\C^{0}_{ub}([0,T];\mathcal{D}_{A}(\beta_q,\kappa))$, such that $u,\partial_t u$ and $Au\in\B^{\beta}_{q,\kappa}(0,T;\X)$, with the estimates
\begin{align*}
    \lVert u\rVert_{\L^{\infty}(0,T;\mathcal{D}_{A}(\beta_q,\kappa))}\lesssim_{\beta,q,A} \lVert  (u,\,\partial_tu,\,A u)\rVert_{\B^\beta_{q,\kappa}(0,T;\X)} \lesssim_{\beta,q,A} \lVert \ff \rVert_{\B^\beta_{q,\kappa}(0,T;\X)} + \lVert \uu_0 \rVert_{\mathcal{D}_{A}(\beta_q,\kappa)}.
\end{align*}
Furthermore, 
\begin{itemize}
    \item when $\kappa=\infty$, one only has $u\in\L^\infty([0,T];\mathcal{D}_{A}(\beta_q,\infty))$ with the strong continuity property
    \begin{align*}
        u\in\C^0_{ub}([0,T];\mathcal{D}_{A}(\vartheta,1)),\quad\text{ for all }\vartheta\in(0,\beta_q);
    \end{align*}
    \item if $T<\infty$, one can remove the invertibility assumption on $(\D(A),A)$, yielding  implicit constants with an exponential growth with respect to $T$.
\end{itemize}
\end{theorem}

\begin{proof} We deal with the case $T=\infty$, otherwise the proof is similar. By assumption $(\D(A),A)$ is an invertible $\omega$-sectorial operator on a Banach space $\X$. Furthermore, thanks to an analysis similar to \cite[Chapter~3,~Section~3.2]{PrussSimonett2016}, one can deduce that the operator $(\D_{q,1}^{\beta_0}(\partial_t,\RR_+,\X),\partial_t)$ on $\B^{\beta_0}_{q,1}(\RR_+,\X)$, with domain
\begin{align*}
    \D_{q,1}^{\beta_0}(\partial_t,\RR_+,\X) :=\B^{\beta_0+1}_{q,1,0}(\RR_+,\X) =\{ \,u\in\B^{\beta_0+1}_{q,1}(\RR_+,\X) \,:\, u(0)=0\, \}
\end{align*}
is a densely defined closed injective $\frac{\pi}{2}$-sectorial operator for any $\beta_0\in(-1+\sfrac{1}{q},\sfrac{1}{q})$.

Since one has $\omega+\frac{\pi}{2}<\pi$, and $(\D(A),A)$ is invertible on $\X$, by \cite[Theorem~3.11]{DaPratoGrisvard1975}, it holds that, for any $\kappa\in[1,\infty]$, any $\theta\in(0,1)$, such that $\beta:=\beta_0+\theta<{\sfrac{1}{q}}$, the operator $A(\partial_{t}+A)^{-1}$ is bounded on
\begin{align*}
    (\B^{\beta_0}_{q,1}(\RR_+;\X),\D_{q,1}^{\beta_0}(\partial_t,\RR_+;\X))_{\theta,\kappa} =(\B^{\beta_0}_{q,1}(\RR_+;\X),\D_{q,1}^{\beta_0}(\partial_t,\RR_+;\X))_{\theta,\kappa}=\B^{\beta}_{q,\kappa}(\RR_+,\X).
\end{align*}
In particular, one has for any $f\in\B^{\beta}_{q,\kappa}(\RR_+;\X)$, $u_0=0$, a unique mild solution $v\in\C^{0}_{ub,0}({\RR_+};\X)$ to \eqref{ACPBesov}, with an estimate
\begin{align*}
    \lVert  (v,\,\partial_tv,\,A v)\rVert_{\B^\beta_{q,\kappa}(\RR_+;\X)} \lesssim_{\beta,q,A} \lVert f \rVert_{\B^\beta_{q,\kappa}(\RR_+;\X)}. 
\end{align*}
We consider $e^{-tA}u_0$, for $u_0\neq 0$, so we can set $u(t):=v(t)+e^{-tA}u_0$, for $t\geqslant0$, so that by Lemma~\ref{lem:BesovintimeSemigroup}, and one obtains the estimate
\begin{align*}
    \lVert  (u,\,\partial_tu,\,A u)\rVert_{\B^\beta_{q,\kappa}(\RR_+;\X)} \lesssim_{\beta,q,A} \lVert f \rVert_{\B^\beta_{q,\kappa}(\RR_+;\X)} + \lVert u_0\rVert_{\mathcal{D}_A(\beta_q,\kappa)}.
\end{align*}

Now, we set $\kappa=1$. By \cite[Theorem~3.1]{MeyriesVeraar2014} (whose proof remains valid for $q=1,\infty$ when no weight is involved), one has the chain of inequalities
\begin{align*}
    \lVert u \rVert_{\L^\infty(\RR_+;\mathcal{D}_A(\beta_q,1))}\lesssim_{\X,A,q} \lVert u \rVert_{\B^{\sfrac{1}{q}}_{q,1}(\RR_+;\mathcal{D}_A(\beta_q,\kappa))} &\lesssim_{\X,A,q}\lVert u \rVert_{\B^{\beta+1}_{q,1}(\RR_+;\X)}^{1-\beta_q} \lVert A u \rVert_{\B^{\beta}_{q,1}(\RR_+;\X)}^{\beta_q}\\
    &\lesssim_{\beta,\X,A,q}\lVert (u,\,\partial_t u,\,  A u) \rVert_{\B^{\beta}_{q,1} (\RR_+;\X)}\\
    &\lesssim_{\beta,\X,A,q} \lVert f \rVert_{\B^\beta_{q,1}(\RR_+;\X)} + \lVert u_0\rVert_{\mathcal{D}_A(\beta_q,1)}.
\end{align*}
Hence, for an arbitrary fixed $t\geqslant 0$, by real interpolation of the linear map
$(f,u_0)\,\longmapsto\,u(t)$, one obtains for $\kappa\in[1,\infty]$
\begin{align*}
    \lVert u \rVert_{\L^\infty(\RR_+;\mathcal{D}_A(\beta_q,\kappa))}\lesssim_{\beta,\X,A,q} 
     \lVert f \rVert_{\B^\beta_{q,\kappa}(\RR_+;\X)} + \lVert u_0\rVert_{\mathcal{D}_A(\beta_q,\kappa)}.
\end{align*}
Therefore, the closed graph theorem applies and one obtains
\begin{align*}
    \lVert u\rVert_{\L^{\infty}(0,T;\mathcal{D}_{A}(\beta_q,\kappa))}\lesssim_{\beta,\X,q,A} \lVert  (u,\,\partial_tu,\,A u)\rVert_{\B^\beta_{q,\kappa}(0,T;\X)}.
\end{align*}
This was the last required estimate.
\end{proof}

It would be of great interest to obtain the homogeneous counterpart of Theorem~\ref{thm:DaPratoGrisvard2} in a vein similar to Theorem~\ref{thm:DaPratoGrisvard} allowing then the same estimates for operators that are injective but not invertible such has the Stokes operator on $\L^\infty_{\mathfrak{n},\sigma}(\RR^n_+)$.

This is the result for bounded domains:

\begin{theorem}\label{thm:MaxRegBesov(Linfty)} Let $q,\kappa\in[1,\infty]$, $\alpha\in(0,1)$ and $\beta\in(-1+\sfrac{1}{q},\sfrac{1}{q})$ such that $2+2\beta-\sfrac{2}{q}\neq0,1$. Let $\Omega$ be either  a bounded $\mathcal{C}^{1,\alpha}$-domain, or bounded $\C^{1,\alpha}$ domain with small $(1,\alpha)$-H\"{o}lderian constant.  For simplicity, we set $\beta_q:=1+\beta-\sfrac{1}{q}$.

\medbreak

 Let $T\in(0,\infty]$. For all $\ff\in\B^{\beta}_{q,\kappa}(0,T;\L^\infty_{\mathfrak{n},\sigma}({\Omega}))$, all $\uu_0\in{\mathcal{D}}_{\AA_\mathcal{D}}(\beta_q,\kappa)$, the problem
\begin{equation*}\tag{DS}
    \left\{ \begin{array}{rllr}
         \partial_t \uu - \Delta \uu +\nabla \mathfrak{p} &= \ff \text{, }&&\text{ in } (0,T)\times\Omega\text{,}\\
        \div \uu &= 0\text{, } &&\text{ in } (0,T)\times\Omega\text{,}\\
        \uu_{|_{\partial\Omega}} &=\mathbf{0}\text{, } &&\text{ on } (0,T)\times\partial\Omega\text{,}\\
        \uu(0) &=\uu_0\text{, } &&\text{ in } {\mathcal{D}}_{\AA_\mathcal{D}}(\beta_q,\kappa)\text{,}
    \end{array}
    \right.
\end{equation*}
admits a unique (mild) solution $(\uu,\mathfrak{p})\in \C_{ub}^0([0,T];\mathcal{D}_{\AA_\mathcal{D}}(\beta_q,\kappa))\times\B^\beta_{q,\kappa}(0,T;\C^{1,\alpha}_\mfree(\overline{\Omega}))$, such that both $\partial_t \uu$ and $-\Delta \uu+\nabla\mathfrak{p}\in\B^\beta_{q,\kappa}(0,T;\L^\infty({\Omega},\CC^n))$, with estimates
\begin{align*}
    \lVert \uu \rVert_{\L^{\infty}(\mathcal{D}_{\AA_\mathcal{D}}(\beta_q,\kappa))}\less_{\beta,q,n}^{\Omega} \lVert  (\uu,\partial_t\uu,-\Delta\uu+\nabla\mathfrak{p})\rVert_{\B^\beta_{q,\kappa}(\L^\infty)}\less_{\beta,q,n}^{\Omega} \lVert \ff \rVert_{\B^\beta_{q,\kappa}(\L^\infty)} + \lVert \uu_0 \rVert_{\mathcal{D}_{\AA_\mathcal{D}}(\beta_q,\kappa)}.
\end{align*}
Additionally, if either $2\beta_q< 1+ \alpha$, or $2\beta_q=1+\alpha$ with $\kappa=\infty$, one has 
\begin{align*}
    \mathcal{D}_{\AA_\mathcal{D}}(\beta_q,\kappa)=\B^{2\beta_q,\sigma}_{\infty,\kappa,\mathcal{D}}(\Omega).
\end{align*}
Except that, when $\kappa=\infty$, one only has weak-$\ast$ continuity in-time.
\end{theorem}

We mention that one can even apply Theorem~\ref{thm:DaPratoGrisvard2} for the Dirichlet-Stokes operator with $\X=\B^{s}_{p,\tilde{r}}$. This yields $\B^{\beta}_{q,\kappa}(\B^{s}_{p,\tilde{r}})$-maximal regularity for the Stokes--Dirichlet problem on a given bounded Lipschitz domain that belongs that the class $\mathcal{M}^{1+\alpha,r}_\W(\epsilon)$, provided $p,r,\tilde{r}\in[1,\infty]$, $s\in(-1+\sfrac{1}{p},\sfrac{1}{p})$. However, we shall not inflict this atrocity on the reader.



\appendix

\section{An embedding for spaces of multipliers.}

We provide here an additional result of an embedding for spaces of multipliers which does not fit the spirit of the presentation in Section~\ref{sec:SM}.

\begin{lemma}\label{lem:EmbeddingMultipliersIntoMWs1} For all $0<s<\tau<1$, all $p\in(1,\infty]$, one has the continuous embedding
\begin{align*}
    \mathcal{M}_{\W,\tt{or}}^{\tau,p}(\RR^n)\hookrightarrow \mathcal{M}_{\W,\tt{or}}^{s,1}(\RR^n).
\end{align*}
The result still holds with $ \mathcal{M}_{\H,\tt{or}}^{\tau,p}(\RR^n)$ whenever $p<\infty$.
\end{lemma}

\begin{proof} \textbf{Step 1:} We deal with the case $1<p<n$. We recall that in our specific case, multiplier norms verify the following estimates
\begin{align*}
    \lVert\varphi\rVert_{ \mathcal{M}_{\W,\tt{or}}^{s,1}(\RR^n)} &\sim_{s,n} \lVert \varphi\rVert_{\L^\infty(\RR^n)}+\sup_{\substack{\B_r\subset\RR^n,\\ r\in(0,1]}} \frac{1}{r^{n-s}} \int_{\B_r}\int_{\B_r} \frac{|\varphi(x)-\varphi(y)|}{|x-y|^{n+s}} \,\d x  \, \d y,\qquad \text{ and }\\
    \lVert\varphi\rVert_{ \mathcal{M}_{\W,\tt{or}}^{\tau,p}(\RR^n)} &\gtrsim_{p,\tau,n} \lVert \varphi\rVert_{\L^\infty(\RR^n)}+\sup_{\substack{\B_r\subset\RR^n,\\ r\in(0,1]}} \left(\frac{1}{r^{n-p\tau}} \int_{\B_r}\int_{\RR^n} \frac{|\varphi(x+h)-\varphi(x)|^p}{|h|^{n+p\tau}} \,\d x \,  \d h\right)^\frac{1}{p},
\end{align*}
according to \cite[Proposition~4.3.1]{MazyaShaposhnikova2009} and \cite[Theorem~2.9]{Sickel1999} ($\W^{s,1}(\RR^n)=\B_{1,1}^{s}(\RR^n)=\mathrm{F}^{s}_{1,1}(\RR^n)$).

Therefore, writing $\tau =s+\varepsilon$, for $\varepsilon>0$, it is sufficient to show that
\begin{align}
    r^{s-n} \int_{\B_r}\int_{\B_r} \frac{|\varphi(x)-\varphi(y)|}{|x-y|^{n+s}}&\, \d x  \, \d y \label{eq:ProofEmbeddingWs1-0} \\
    &\lesssim_{p,s,n}^\varepsilon r^{s+\varepsilon-\frac{n}{p}} \left(\int_{\B_r}\int_{\RR^n} \frac{|\varphi(x+h)-\varphi(x)|^p}{|h|^{n+p(s+\varepsilon)}}\, \d x \, \d h\right)^\frac{1}{p} +  \lVert \varphi\rVert_{\L^\infty(\RR^n)}.\nonumber
\end{align}
By change of variable, and since $x, x+h\in \B_r$ implies $|h|\leqslant 2r$, it holds
\begin{align}\label{eq:ProofEmbeddingWs1-0.5}
    \int_{\B_r}\int_{\B_r} \frac{|\varphi(x)-\varphi(y)|}{|x-y|^{n+s}} \,\d x \,\d y \leqslant \int_{\B_{2r}}\int_{\B_r} \frac{|\varphi(x+h)-\varphi(h)|}{|h|^{n+s}} \,\d x \,\d h. 
\end{align}
Applying Fubini, we consider and split
\begin{align*}
    \int_{\B_{r}}\int_{\B_{2r}} &\frac{|\varphi(x+h)-\varphi(h)|}{|h|^{n+s}}  \d h  \, \d x\\
    &= \int_{\B_{r}}\int_{\B_r} \frac{|\varphi(x+h)-\varphi(h)|}{|h|^{n+s}}  \d h \, \d x + \int_{\B_{r}}\int_{\B_{2r}\setminus\B_r} \frac{|\varphi(x+h)-\varphi(h)|}{|h|^{n+s}}  \d h \,  \d x,
\end{align*}
so that by H\"{o}lder inequality on the first integral
\begin{align*}
    \int_{\B_r} \frac{|\varphi(x+h)-\varphi(h)|}{|h|^{n+s}}  \d h &\leqslant \left(\int_{\B_r} \frac{|\varphi(x+h)-\varphi(h)|^p}{|h|^{n+p(s+\varepsilon)}}  \d h\right)^\frac{1}{p}\left(\int_{\B_r} \frac{|h|^{\varepsilon p'}}{|h|^{n}}  \d h\right)^\frac{1}{p'}\\
    &\lesssim_{p,s,n}^{\varepsilon} r^{\varepsilon}\left(\int_{\B_r} \frac{|\varphi(x+h)-\varphi(h)|^p}{|h|^{n+p(s+\varepsilon)}}  \d h\right)^\frac{1}{p}.
\end{align*}
Hence, integrating again on $\B_r$, applying again H\"{o}lder inequality:
\begin{align}
    \int_{\B_{r}}\int_{\B_r} \frac{|\varphi(x+h)-\varphi(h)|}{|h|^{n+s}}\, \d x\, \d h &\lesssim_{p,s,n}^{\varepsilon} r^{\varepsilon} \int_{\B_{r}}\left(\int_{\B_r} \frac{|\varphi(x+h)-\varphi(h)|^p}{|h|^{n+p(s+\varepsilon)}} \, \d h\right)^\frac{1}{p} \,\d x\nonumber\\
    &\lesssim_{p,s,n}^{\varepsilon} r^{\varepsilon+\frac{n}{p'}} \left(\int_{\B_{r}}\int_{\B_r} \frac{|\varphi(x+h)-\varphi(h)|^p}{|h|^{n+p(s+\varepsilon)}} \, \d h  \,  \d x\right)^\frac{1}{p}.\label{eq:ProofEmbeddingWs1-1}
\end{align}
Now, for the second integral
\begin{align}
    \int_{\B_r}\int_{\B_{2r}\setminus\B_r} \frac{|\varphi(x+h)-\varphi(h)|}{|h|^{n+s}}  \d h \,\d x &\leqslant 2 \lVert \varphi\rVert_{\L^\infty(\RR^n)} |\B_r| \int_{\B_{2r}\setminus\B_{r}}\frac{1}{|h|^{n+s}} \, \d h\nonumber\\
    &\leqslant 2 |\B_1| \left(\int_{\B_{2}\setminus\B_{1}}\frac{1}{|h|^{n+s}} \, \d h\right) \lVert \varphi\rVert_{\L^\infty(\RR^n)} r^{n-s}\nonumber\\
    &\lesssim_n r^{n-s}\lVert \varphi\rVert_{\L^\infty(\RR^n)}.\label{eq:ProofEmbeddingWs1-2}
\end{align}
Replacing \eqref{eq:ProofEmbeddingWs1-1} and \eqref{eq:ProofEmbeddingWs1-2} in \eqref{eq:ProofEmbeddingWs1-0.5} yields \eqref{eq:ProofEmbeddingWs1-0}.

\textbf{Step 2:} The case $p\in[n,\infty]$. If $p<\infty$, by  Proposition~\ref{prop:embeddingmultipliers}, \cite[Proposition~3.5.3]{MazyaShaposhnikova2009} and Proposition~\ref{prop:embeddingmultipliers} again, choosing $q\in(1,n)$, it holds
\begin{align*}
    \mathcal{M}_{\W,\tt{or}}^{s+\varepsilon,p}(\RR^n)\hookrightarrow \mathcal{M}_{\H,\tt{or}}^{s+\frac{2}{3}\varepsilon,p}(\RR^n) \hookrightarrow \mathcal{M}_{\H,\tt{or}}^{s+\frac{2}{3}\varepsilon,q}(\RR^n)\hookrightarrow \mathcal{M}_{\W,\tt{or}}^{s+\frac{1}{3}\varepsilon,q}(\RR^n)\hookrightarrow\mathcal{M}_{\W,\tt{or}}^{s,1}(\RR^n)
\end{align*}
where the last embedding is a consequence of Step 1. The case $p=\infty$ holds trivially since $\mathcal{M}_{\W,\tt{or}}^{s+\varepsilon,\infty}(\RR^n)=\B^{s+\varepsilon}_{\infty,\infty}(\RR^n)\hookrightarrow\mathcal{M}_{\W,\tt{or}}^{s,1}(\RR^n)$.
\end{proof}

\section{On interpolation of intersections.}

Following the exact same ideas in \cite[Appendix~A]{Gaudin2022} and the proof of \cite[Theorem~2.12]{Gaudin2023Lip}, we obtain the next result.

\begin{proposition}\label{prop:InterpHybridBesov}Let $q_0,q_1,r_0,r_1,p,\Tilde{q}\in[1,\infty]$, $s_0,s_1,\alpha_0,\alpha_1\in\RR$ such that $s_j\leqslant \alpha_j$, $j\in\{0,1\}$. For $\theta\in(0,1)$, we set for $\beta\in\{s,\alpha\}$, $\beta:= (1-\theta)\beta_0+\theta \beta_1$, $1/q:= (1-\theta)/q_0+\theta /q_1$.

\begin{enumerate}
    \item We assume $s_0\neq s_1$ and $\alpha_0\neq \alpha_1$, one has, in this case,
    \begin{align*}
        \Big( \dot{\H}^{s_0,p}\cap\dot{\H}^{\alpha_0,p}(\RR^{n}), \dot{\H}^{s_1,p}\cap\dot{\H}^{\alpha_1,p}(\RR^{n}) \Big)_{\theta,\Tilde{q}}&=\Big( \dot{\B}_{p,q_0}^{s_0}\cap\dot{\B}_{p,r_0}^{\alpha_0}(\RR^{n}), \dot{\B}_{p,q_1}^{s_1}\cap\dot{\B}_{p,r_1}^{\alpha_1}(\RR^{n}) \Big)_{\theta,\Tilde{q}}\\ &= \dot{\B}_{p,\Tilde{q}}^{s}\cap\dot{\B}_{p,\Tilde{q}}^{\alpha}(\RR^{n}).
    \end{align*}
    A similar result holds for the corresponding intersections of spaces $\dot{\B}^{\bullet,0}_{\infty,\bullet}$, $\dot{\BesSmo}^{\bullet}_{\bullet,\infty}$ and $\dot{\BesSmo}^{\bullet,0}_{\infty,\infty}$.
    \item We assume \hyperref[AssumptionCompletenessExponents]{$(\mathcal{C}_{s_j,p,q_j})$} , for $j\in\{0,1\}$, $1/r:= (1-\theta)/r_0+\theta /r_1$, and $(q_0,q_1),(r_0,r_1)\neq(\infty,\infty)$. In this case, one has
    \begin{align*}
        \Big[ \dot{\B}_{p,q_0}^{s_0}\cap\dot{\B}_{p,r_0}^{\alpha_0}(\RR^{n}), \dot{\B}_{p,q_1}^{s_1}\cap\dot{\B}_{p,r_1}^{\alpha_1}(\RR^{n}) \Big]_{\theta} = \dot{\B}_{p,q}^{s}\cap\dot{\B}_{p,r}^{\alpha}(\RR^{n}).
    \end{align*}
    \item We assume \hyperref[AssumptionCompletenessExponents]{$(\mathcal{C}_{s_j,p})$}, for $j\in\{0,1\}$, and $1<p<\infty$. In this case, one has
    \begin{align*}
        \Big[ \dot{\H}^{s_0,p}\cap\dot{\H}^{\alpha_0,p}(\RR^{n}), \dot{\H}^{s_1,p}\cap\dot{\H}^{\alpha_1,p}(\RR^{n}) \Big]_{\theta} =\dot{\H}^{s,p}\cap\dot{\H}^{\alpha,p}(\RR^{n}).
    \end{align*}
\end{enumerate}
\end{proposition}

We aim to transfer the result above to homogeneous function spaces on special Lipschitz domains, whenever \textit{at least} $s_0,\alpha_0,s_1,\alpha_1>-1+\sfrac{1}{p}$.

We start with a Lemma that enhance the boundedness of Stein's extension operator, denoted here by $\tilde{\mathcal{E}}$ obtained in \cite[Chapter~VI,~Section~3,~Theorem~5']{Stein1970}, and improved in \cite[Propositions~3.12~\&~3.34]{Gaudin2023Lip}.

\begin{lemma}\label{lem:SteinExtOpNegative} Let $\Omega$ be a special Lipschitz domain, and consider $\tilde{\mathcal{E}}$ to be Stein's extension operator from \cite[Chapter~VI,~Section~3,~Theorem~5']{Stein1970}.  Then for all $p,q\in[1,\infty]$, $s\in(-1,\infty)$, $\tilde{\mathcal{E}}$ maps continuously
\begin{itemize}
    \item $\dot{\B}^{s}_{p,q}(\Omega)$ to $\dot{\B}^{s}_{p,q}(\RR^n)$;
    \item $\dot{\H}^{s,p}(\Omega)$ to $\dot{\H}^{s,p}(\RR^n)$, if $1<p<\infty$.
\end{itemize}
\end{lemma}

\begin{proof} Note that due to \cite[Propositions~3.12~\&~3.34]{Gaudin2023Lip}, and due to interpolation, we only need to prove it for $-1<s<0$ and $q<\infty$. The key point here is to notice that the divergence preserving extension operator (here restricted to vector fields $\Lambda^1\simeq\CC^n$) given by Hiptmair, Li and Zou $\mathcal{E}_\sigma$ --exposed in Section~\ref{sec:ExtOpDivPreser}, Proposition~\ref{prop:bddExtOpHomFreeDivSpeLip}-- is an actual instance of Stein's extension operator $\tilde{\mathcal{E}}$, and it can be checked that in fact for all smooth vector fields $\uu$
\begin{align*}
    \div\mathcal{E}_{\sigma}\uu = \tilde{\mathcal{E}}[\div \uu].
\end{align*}
This is a direct consequence of the definition by Hodge-Star duality and pullback, as given for the operators in Section~\ref{sec:ExtOpDivPreser}, \cite{HiptmairLiZou2012} and \cite[Chapter~VI,~Section~3]{Stein1970}, see also \cite[Section~3]{Gaudin2023Lip}. The details are left to the reader.

In order to fix the ideas, we assume $q<\infty$ and we let $(U_k)_{k\in\NN}\subset\C^{\infty}_{ub,h}\cap\dot{\B}^{s}_{p,q}\cap\dot{\B}^{s+1}_{p,q}(\RR^n)$, such that $(U_k)_{k\in\NN}$ converges towards $U\in\dot{\B}^{s}_{p,q}(\RR^n)$ and arbitrary extension of $u\in\dot{\B}^{s}_{p,q}({\Omega})$. We set $(u_k)_{k\in\NN}:=(U_k{}_{|_\Omega})_{k\in\NN}$, so $(u_k)_{k\in\NN}$ converges to $u$. By Theorem~\ref{thm:PotentialOpSpeLipDom}, one can write
\begin{align*}
    u =[\div \T^{\sigma}] u,\,\text{ and }\,  u_k =\div \T^{\sigma} u_k.
\end{align*}
It holds that
\begin{align*}
    \lVert \tilde{\mathcal{E}}u_k\rVert_{\dot{\B}^{s}_{p,q}(\RR^n)}=\lVert \tilde{\mathcal{E}} \div \T^{\sigma} u_k\rVert_{\dot{\B}^{s}_{p,q}(\RR^n)} &\leqslant \liminf_{\substack{a\rightarrow 0_+\\ b\rightarrow\infty}} \lVert \div {\mathcal{E}}_{\sigma} \T^{\sigma}_{a,b} u_k\rVert_{\dot{\B}^{s}_{p,q}(\RR^n)}\\
    &\lesssim_{p,s,n} \liminf_{\substack{a\rightarrow 0_+\\ b\rightarrow\infty}}\lVert {\mathcal{E}}_{\sigma}\T^{\sigma}_{a,b}u_k\rVert_{\dot{\B}^{s+1}_{p,q}(\RR^n)}\\
    &\lesssim_{p,s,n,\partial\Omega} \lVert \T^{\sigma}u_k\rVert_{\dot{\B}^{s+1}_{p,q}(\Omega)}\\
    &\lesssim_{p,s,n,\partial\Omega} \lVert u_k\rVert_{\dot{\B}^{s}_{p,q}(\Omega)}.
\end{align*}
Where, in order to obtain the estimates above we did apply Theorem~\ref{thm:PotentialOpSpeLipDom} and Remark~\ref{rem:NotS'hBsinftyqSteinExtOp}.
Therefore, by extension by density, one obtains for all $u\in\dot{\B}^{s}_{p,q}(\Omega)$
\begin{align*}
    \lVert \tilde{\mathcal{E}}u\rVert_{\dot{\B}^{s}_{p,q}(\RR^n)} &\lesssim_{p,s,n} \lVert u\rVert_{\dot{\B}^{s}_{p,q}(\Omega)}.
\end{align*}
This yields the result, since the proof for the homogeneous Sobolev spaces is similar, and the case of Besov spaces.
\end{proof}

A standard retraction and co-retraction argument allowed by Lemma~\ref{lem:SteinExtOpNegative} allows the following counterparts of 

\begin{corollary}\label{cor:InterpHybridBesovSpeLip} Proposition~\ref{prop:InterpHybridBesov} still holds for $\Omega$ to be a special Lipschitz domain instead of $\RR^n$ if one assumes additionally that $s_0,s_1,\alpha_0,\alpha_1>-1$.
\end{corollary}

\section{Traces and partial traces including endpoint spaces}\label{App:Traces}

\paragraph{About standard strong traces on the flat half-space.} Following the strategy exhibited in \cite[Section~4]{Gaudin2023Lip}, we can remove to ask for the completeness of function spaces in order to define the trace on the boundary $\partial\RR^n_+$.

\begin{theorem}\label{thm:tracesRn+}Let $p,q\in[1,\infty]$, $q<\infty$, $s>{\sfrac{1}{p}}$. The following assertions hold true
\begin{enumerate}
    \item for $1<p<\infty$,
    \begin{align*}
        \dot{\H}^{s,p}(\RR^n_+)\hookrightarrow \C^{0}_{0,x_n}(\overline{\RR_+},\dot{\B}^{s-\sfrac{1}{p}}_{p,p}(\RR^{n-1}));
    \end{align*}
    \item for $p<\infty$,
    \begin{align*}
        \dot{\B}^{s}_{p,q}(\RR^n_+) &\hookrightarrow \C^{0}_{0,x_n}(\overline{\RR_+},\dot{\B}^{s-\sfrac{1}{p}}_{p,q}(\RR^{n-1}));\\
        \dot{\BesSmo}^{s}_{p,\infty}(\RR^n_+)&\hookrightarrow \C^{0}_{0,x_n}(\overline{\RR_+},\dot{\BesSmo}^{s-\sfrac{1}{p}}_{p,\infty}(\RR^{n-1}));\\
        \dot{\B}^{s}_{p,\infty}(\RR^n_+)&\hookrightarrow \C^{\mathrm{w}\ast}_{0,x_n}(\overline{\RR_+},\dot{\B}^{s-\sfrac{1}{p}}_{p,\infty}(\RR^{n-1})).
    \end{align*}
    \item for $p=\infty$,
    \begin{align*}
        {\B}^{s}_{\infty,q}(\RR^n_+) &\hookrightarrow \C^{0}_{ub,x_n}(\overline{\RR_+},{\B}^{s}_{\infty,q}(\RR^{n-1}));\\
        {\B}^{s}_{\infty,\infty}(\RR^n_+) &\hookrightarrow \C^{\mathrm{w}\ast}_{ub,x_n}(\overline{\RR_+},{\B}^{s}_{\infty,\infty}(\RR^{n-1}));
    \end{align*}
    with the homogeneous estimates
    \begin{align*}
        \sup_{x_n\geqslant 0} \lVert u(\cdot,x_n)\rVert_{\dot{\B}^{s}_{\infty,r}(\RR^{n-1})} \lesssim_{s,n} \lVert u\rVert_{\dot{\B}^{s}_{\infty,r}(\RR^{n}_+)} \text{, }u\in\B^{s}_{\infty,r}(\RR^n_+)\text{, }r\in[1,\infty].
    \end{align*}
    Additionally, a similar result holds if we replace $({\B}^{s}_{\infty,q},\C^{0}_{ub})$ by $({\B}^{s,0}_{\infty,q},\C^0_0)$, $({\BesSmo}^{s}_{\infty,\infty},\C^0_{ub})$ or $({\BesSmo}^{s,0}_{\infty,\infty},\C^0_0)$, and ${\B}^{s}_{\infty,\infty}$ by ${\B}^{s,0}_{\infty,\infty}$;
    \item for $p<\infty$,
    \begin{align*}
        \dot{\B}^{\sfrac{1}{p}}_{p,1}(\RR^n_+) &\hookrightarrow \C^{0}_{0,x_n}(\overline{\RR_+},{\L}^{p}(\RR^{n-1}));\\
        \dot{\B}^{0}_{\infty,1}(\RR^n_+)&\hookrightarrow \C^{0}_{ub,x_n}(\overline{\RR_+},\C^{0}_{ub}(\RR^{n-1}));\\
        \dot{\B}^{0,0}_{\infty,1}(\RR^n_+)&\hookrightarrow \C^{0}_{0,x_n}(\overline{\RR_+},\C^{0}_{0}(\RR^{n-1})).
    \end{align*}
\end{enumerate}
Moreover, 
\begin{itemize}
    \item the induced trace operator $u\mapsto u_{|_{\partial\RR^n_+}}=u(\cdot,0)$ admits a right bounded inverse in the cases \textit{(i)} and \textit{(ii)};
    \item in any of above cases, if additionally $u\in\L^\infty(\RR^{n}_+)$, one also has $u(\cdot,0)\in\L^\infty(\RR^{n-1})$ with the estimate
    \begin{align*}
        \lVert u(\cdot,0)\rVert_{\L^{\infty}(\RR^{n-1})}\leqslant \lVert u\rVert_{\L^{\infty}(\RR^{n}_+)}.
    \end{align*}
\end{itemize}
\end{theorem}

\begin{remark} The right bounded inverse of the trace operator in the cases \textit{(i)} and \textit{(ii)} is given by the operator $T_0=(-\Delta_{\mathcal{D},\partial})^{-1}$ from Proposition~\ref{prop:PoissonSemigroup3} below.

\textit{La raison d'être} of point \textit{(iii)}'s formulation is due to the following facts, we remind to the reader here.
\begin{itemize}
    \item We have, for $s>0$ and $q\in[1,\infty]$, the algebraic equalities $\dot{\B}^{s}_{\infty,q}(\RR^n_+) = {\B}^{s}_{\infty,q}(\RR^n_+)\cap \S'_h(\RR^n_+)\subset {\B}^{s}_{\infty,q}(\RR^n_+)$, and ${\B}^{s,0}_{\infty,q}(\RR^n_+)=\dot{\B}^{s,0}_{\infty,q}(\RR^n_+)\subset\mathcal{S}'_h(\RR^n_+)$, hence the provided embeddings and the estimate are more general.
    \item An important fact is that, in general, the trace of an element from the space ${\B}^{s}_{\infty,q}(\RR^n_+)\cap \S'_h(\RR^n_+)$ may not lie in $\S'_h(\RR^{n-1})$. Consider for instance
    \begin{align*}
        (x',x_n)\mapsto e^{-x_n (\I-\Delta')^{\sfrac{1}{2}} } \mathbf{1} (x') = e^{-x_n} \in \W^{1,\infty}(\RR^n_+)\cap\S'_h(\RR^n_+),
    \end{align*}
    recalling that we have $\mathbf{1}\notin \S'_h(\RR^{n-1})$. This is a consequence of Proposition~\ref{prop:DumpedPoissonSemigroup2} below.
    \item However, for $p\in[1,\infty)$ and $s>{\sfrac{1}{p}}$, we do have the inclusions
    \begin{align*}
        \dot{\B}^{s}_{p,q}(\RR^n_+) \subset {\B}^{s}_{p,q}(\RR^n_+) + \C^{\lceil s \rceil}_0(\overline{\RR^n_+}) \subset \C^{0}_{0,x_n}(\overline{\RR_+},{\B}^{s-\sfrac{1}{p}}_{p,q}(\RR^{n-1})) + \C^{0}_{0,x_n}(\overline{\RR_+},\C^{\lceil s \rceil}_0(\RR^{n-1}))
    \end{align*}
    replacing $\C^{0}_{0,x_n}$ by $\C^{\mathrm{w}\ast}_{0,x_n}$ in the last inclusion, whenever $q=\infty$. This implies that for all $u\in\dot{\B}^{s}_{p,q}(\RR^n_+)$, one has $u(\cdot,x_n)\in\B^{s-\sfrac{1}{p}}_{p,q}(\RR^{n-1})+\C^{\lceil s \rceil}_0(\RR^{n-1})\subset\mathcal{S}'_h(\RR^{n-1})$, for all $x_n\geqslant0$.
    \item The extension operator $T_\lambda$ from the boundary $\partial\RR^n_+$ to $\RR^n_+$, given in Proposition~\ref{prop:DumpedPoissonSemigroup2} below, is also such that for all $f\in \dot{\B}^{s}_{\infty,q}(\partial\RR^n_+)\subset \S'_h(\RR^{n-1})$, one has $T_\lambda f\in\dot{\B}^{s}_{\infty,q}(\partial\RR^n_+) \subset \S'_h(\RR^n_+)$, whenever $q\in[1,\infty]$ with the corresponding estimate.
\end{itemize}
\end{remark}

\begin{corollary}\label{cor:TraceIntersectHomSpacesRn+}Let $p,q\in[1,\infty)$, and $s\in(-1+\sfrac{1}{p},\sfrac{1}{p})$.
\begin{enumerate}
    \item $\dot{\H}^{s,p}\cap\dot{\H}^{s+1,p}(\RR^n_+)\hookrightarrow \C^{0}_{0,x_n}(\overline{\RR_+},{\B}^{s+1-\sfrac{1}{p}}_{p,p}(\RR^{n-1}))$, provided $p>1$;
    \item $\dot{\B}^{s}_{p,q}\cap \dot{\B}^{s+1}_{p,q}(\RR^n_+) \hookrightarrow \C^{0}_{0,x_n}(\overline{\RR_+},{\B}^{s+1-\sfrac{1}{p}}_{p,q}(\RR^{n-1}));$
    \item $\dot{\BesSmo}^{s}_{p,\infty}\cap\dot{\BesSmo}^{s+1}_{p,\infty}(\RR^n_+)\hookrightarrow \C^{0}_{0,x_n}(\overline{\RR_+},{\BesSmo}^{s+1-\sfrac{1}{p}}_{p,\infty}(\RR^{n-1}))$;
    \item $\dot{\B}^{s}_{\infty,q}\cap \dot{\B}^{s+1}_{\infty,q}(\RR^n_+) \hookrightarrow \C^{0}_{ub,x_n}(\overline{\RR_+},{\B}^{s+1-\sfrac{1}{p}}_{\infty,q}(\RR^{n-1}))$, $-1<s<0$;
    \item  $\dot{\B}^{s}_{p,\infty}\cap \dot{\B}^{s+1}_{p,\infty}(\RR^n_+)\hookrightarrow \C^{\mathrm{w}\ast}_{0,x_n}(\overline{\RR_+},{\B}^{s+1-\sfrac{1}{p}}_{p,\infty}(\RR^{n-1})).$
\end{enumerate}
Moreover, the induced trace operator $u\mapsto u_{|_{\partial\RR^n_+}}=u(\cdot,0)$ is onto in each of the cases above, and admits a right bounded inverse that can be given by the linear operator $T_\lambda$ as in Proposition~\ref{prop:DumpedPoissonSemigroup2} below.

Furthermore, the result still holds if we replace ${\B}^{\cdot}_{p,q}$ by either ${\B}^{\cdot,0}_{\infty,q}$ or ${\BesSmo}^{\cdot,0}_{\infty,\infty}$, and $({\B}^{\cdot}_{p,\infty},\C^{\mathrm{w}\ast}_{0,x_n})$ by either $({\B}^{\cdot,0}_{\infty,\infty},\C^{\mathrm{w}\ast}_{0,x_n})$ or $({\B}^{\cdot}_{\infty,\infty},\C^{\mathrm{w}\ast}_{ub,x_n})$.
\end{corollary}

\paragraph{Overall strategy to define and identify partial traces on the boundary.} When dealing with vector fields of low regularity, say $\uu\in\X^{s,p}(\Omega,\CC^n)$, say for $-1+{\sfrac{1}{p}}<s<{\sfrac{1}{p}}$ and such that $\div u\in\X^{s,p}(\Omega,\CC)$, it is of standard interest to figure out how to define the normal trace of $\uu$, ${\nu\cdot \uu}_{|_{\partial\Omega}}$ on the boundary. The very standard and well-known procedure starting from the variational/duality functional
    \begin{align}\label{eq:VaritationnalFormulaDefPartialtraceDiv}
        (\uu,\varphi)\longmapsto\langle  \uu, \nabla \varphi\rangle_{\Omega}+ \langle \div \uu, \varphi\rangle_{\Omega}
    \end{align}
reads as follows:
\begin{enumerate}
    \item One proves that $\Ccinfty(\overline{\Omega})$ is dense in $\X'^{-s,p'}\cap\X'^{-s+1,p'}(\Omega,\CC^n)$ with surjective trace in $\tilde{\X}^{-s+1-{\sfrac{1}{p'}},p'}_{\partial}(\partial\Omega)$ (or in nice closed subspaces, still denoted the same)
    \item One proves that $\Ccinfty({\Omega})$ is dense in the subspaces constituted of the elements $\varphi$ of $\X'^{-s,p'}\cap\X'^{-s+1,p'}(\Omega,\CC^n)$ such that $\varphi_{|_{\partial\Omega}}=0$.
    \item Then, $\uu$ being fixed, one deduces that the formula \eqref{eq:VaritationnalFormulaDefPartialtraceDiv} induces a distribution on the boundary that acts only depending on the value of $\varphi_{|_{\partial\Omega}}$.
    \item "Conclude" from point \textit{(i)} and \textit{(iii)} in combination with the Open mapping Theorem there exists a unique distribution on the boundary $\mathfrak{n}(\uu)\in\tilde{\X}'^{s-1/p,p}_{\partial}(\partial\Omega)$ (the dual space of $\tilde{\X}^{-s+1-{\sfrac{1}{p'}},p'}_{\partial}(\partial\Omega)$) such that
    \begin{align*}
        \big\langle  \mathfrak{n}_{\partial\Omega}(\uu), \varphi_{|_{\partial\Omega}} \big\rangle_{\partial\Omega}=\big\langle  \uu, \nabla \varphi \big\rangle_{\Omega}+ \big\langle \div \uu, \varphi\big\rangle_{\Omega}
    \end{align*}
    for all $\varphi\in\Ccinfty(\overline{\Omega})$.
\end{enumerate}
While this procedure is very well-known and standard, it is usually additionally claimed that one can identify $\mathfrak{n}_{\partial\Omega}(\uu)$ with $\nu\cdot \uu_{|_{\partial\Omega}}$, while this is true on the closure of smooth functions in the $\X^{s,p}$-domain of the divergence operator, one still has to prove this is true for the whole $\X^{s,p}$-domain of the divergence operator\footnote{This fact obviously even more neglected when it comes to homogeneous function space, even the second author didn't explicitly removed this ambiguity in his previous work \cite[Appendix]{Gaudin2023Hodge}}. When $\X^{s,p}$ is nice enough, this amounts to prove that
\begin{align*}
    \D^{s,p}(\div,\Omega) = \overline{\Ccinfty(\overline{\Omega},\CC^n)}^{\lVert\cdot\rVert_{\X^{s,p}(\Omega)}+\lVert \div (\cdot)\rVert_{\X^{s,p}(\Omega)}}.
\end{align*}
But his might require additional unattainable knowledge, except in several particular cases. The point will be to take advantage of these particular cases.

Therefore, one will have to stop temporarily the proof, to perform the following arguments
\begin{itemize}
    \item Use the variationnal/duality formula and the new abstract trace information on the boundary \eqref{eq:VaritationnalFormulaDefPartialtraceDiv} to deduce density of smooth function in $\D^{s,p}(\div,\RR^n_+)$ first then $\D^{s,p}(\div,\Omega)$ by transformation and localization.
    \item Once well-defined for standard nice enough function spaces, in the case of endpoint spaces for which smooth functions are not dense, one may perform interpolation arguments. Hence, one may need to show additionally that the family of function spaces $(\D^{s,p}(\div,\Omega))_{s,p}$ is an interpolation scale. Other methods, such has embedding-type results, might be of use to deduce the distributional equality $\mathfrak{n}_{\partial\Omega}(\uu)=\nu\cdot \uu_{|_{\partial\Omega}}$. Sometimes, one might use a combination of both.
\end{itemize}
Eventually, the same considerations are still meaningful for the $\curl$ operator and its associated partial trace $\nu\times(\cdot)_{|_{\partial\Omega}}$, as well as their higher-dimensional counterparts represented by differential forms.

\paragraph{The case of the flat half-space.} Hence, we start by investigating the case of the flat half-space.

\begin{theorem}\label{thm:partialtracesDiffFormRn+} Let $p,q\in[1,\infty]$, and $s\in(-1+\sfrac{1}{p},\sfrac{1}{p})$ and $k\geqslant 1$.
\begin{enumerate}
    \item For $u\in{\L}^{1}(\RR^n_+,\Lambda^k)$ such that $\delta u \in {\L}^{1}(\RR^n_+,\Lambda^{k-1})$, one can define the partial trace
    \begin{align*}
    (-\mathfrak{e}_n)\iprod u_{|_{\partial\RR^n_+}} \in \W^{-1,1}(\partial\RR^n_+,\Lambda^{k-1}),
\end{align*}
with the estimate
\begin{align*}
    \lVert  (-\mathfrak{e}_n)\iprod u_{|_{\partial\RR^n_+}} \rVert_{ \W^{-1,1}(\partial\RR^n_+)} \lesssim_{n} \lVert  (u,\,\delta u) \rVert_{{\L}^{1}(\RR^n_+)}.
\end{align*}

\item For $u\in{\L}^{\infty}(\RR^n_+,\Lambda^k)$ such that $\delta u \in {\L}^{\infty}(\RR^n_+,\Lambda^{k-1})$, one can define the partial trace
\begin{align*}
    (-\mathfrak{e}_n)\iprod u_{|_{\partial\RR^n_+}} \in \L^{\infty}(\partial\RR^n_+,\Lambda^{k-1}),
\end{align*}
with the estimate
\begin{align*}
    \lVert  (-\mathfrak{e}_n)\iprod u_{|_{\partial\RR^n_+}} \rVert_{{\L}^{\infty}(\partial\RR^n_+)} \lesssim_{n} \lVert  (u,\,\delta u) \rVert_{{\L}^{\infty}(\RR^n_+)}.
\end{align*}
    \item For $u\in\dot{\B}^{s}_{p,q}(\RR^n_+,\Lambda^k)$ such that $\delta u\in \dot{\B}^{s}_{p,q}(\RR^n_+,\Lambda^{k-1})$, one can define the partial trace
\begin{align*}
    (-\mathfrak{e}_n)\iprod u_{|_{\partial\RR^n_+}} \in {\B}^{s-\sfrac{1}{p}}_{p,q}(\partial\RR^n_+,\Lambda^{k-1}),
\end{align*}
with the estimate
\begin{align*}
    \lVert  (-\mathfrak{e}_n)\iprod u_{|_{\partial\RR^n_+}} \rVert_{{\B}^{s-\sfrac{1}{p}}_{p,q}(\partial\RR^n_+)} \lesssim_{s,n} \lVert  (u,\,\delta u) \rVert_{\dot{\B}^{s}_{p,q}(\RR^n_+)}.
\end{align*}
If $p=q\in(1,\infty)$, the result still holds with $\dot{\H}^{s,p}$ instead of $\dot{\B}^{s}_{p,p}$.
\end{enumerate}
Furthermore, 
\begin{itemize}
    \item if $\delta u =0$ one can replace the trace spaces ${\B}^{s-\sfrac{1}{p}}_{p,q}(\partial\RR^n_+)$ and $\W^{-1,1}(\partial\RR^n_+)$, by their homogeneous counterpart $\dot{\B}^{s-\sfrac{1}{p}}_{p,q}(\partial\RR^n_+)$ and $\dot{\W}^{-1,1}(\partial\RR^n_+)$;

    \item a similar result holds with $((-\mathfrak{e}_n)\wedge\cdot,\, \d,\, \Lambda^{k+1})$ that replaces $((-\mathfrak{e}_n)\iprod\cdot,\, \delta,\, \Lambda^{k-1})$;
    
    \item up to appropriate changes, the whole result remains valid for inhomogeneous function spaces.
\end{itemize}
\end{theorem}

\begin{remark}Actually the proof below shows that the map
\begin{align*}
\left\{ \begin{array}{cllc}
          \D^{s}_{p,q}(\delta,\RR^n_+,\Lambda)& \longrightarrow && \C^{0}_{ub,x_n}(\overline{\RR_+},\B^{s-\frac{1}{p}}_{p,q}(\RR^{n-1},\Lambda))\\
    u&\longmapsto && [x_n\mapsto (-\mathfrak{e}_n)\iprod u (\cdot,x_n)].
    \end{array}
    \right.
\end{align*}
is well-defined and continuous provided $p\in[1,\infty]$, $q\in[1,\infty)$. For the other function spaces the result still hold up to appropriate changes. For instance, in the case of $\L^\infty$ or endpoint Besov spaces $\B^{s}_{p,\infty}$, one only ask for weak-$\ast$ continuity with respect to the last variable. In particular, the meaning of partial traces on the boundary is \textbf{stronger} than one might initially expect \textbf{!}
\end{remark}

\begin{proof} \textbf{Step 1:} The case \textit{(i)}. Let $u\in{\L}^{1}(\RR^n_+,\Lambda^k)$ such that $\delta u \in {\L}^{1}(\RR^n_+,\Lambda^{k-1})$. We can define for all $\psi\in\mathcal{S}(\partial\mathbb{R}^{n}_+,\Lambda^{k-1})$, and $\Psi\in \mathrm{W}^{1,\infty}(\mathbb{R}^n_+,\Lambda^{k-1})$ such that $\Psi_{|_{\partial\mathbb{R}^n_+}}=\psi$, the following functional,
\begin{align*}
    \kappa_{u}(\Psi):=  \int_{\mathbb{R}^{n}_+}  u(x) \cdot \mathrm{d} {\Psi}(x)\mathrm{~d}x - \int_{\mathbb{R}^{n}_+}  \mathbf{\delta} u(x)\cdot  {\Psi}(x)\mathrm{~d}x\text{.}
\end{align*}
First, the map $(u,\Psi)\mapsto\kappa_{u}(\Psi)$ is well-defined and bilinear on ${\mathrm{D}}_{1}( \mathbf{\delta},\mathbb{R}^{n}_+,\Lambda^k)\times\mathrm{W}^{1,\infty}(\mathbb{R}^n_+,\Lambda^{k-1})$, and only depends on the boundary value $\psi$ of $\Psi$. It is straightforward from H\"{o}lder's inequality that,
\begin{align*}
    \lvert \kappa_{u}(\Psi) \rvert &\leqslant\lVert u \rVert_{{\L}^{1}(\mathbb{R}^n_+)}\lVert \mathrm{d}\Psi \rVert_{{\L}^{\infty}(\mathbb{R}^n_+)} + \lVert {\delta} u \rVert_{{\L}^{1}(\mathbb{R}^n_+)}\lVert\Psi \rVert_{{\L}^{\infty}(\mathbb{R}^n_+)}\\
    &\leqslant\lVert (u,\,\delta u) \rVert_{\L^1(\mathbb{R}^n_+)} \lVert \Psi \rVert_{{\W}^{1,\infty}(\mathbb{R}^n_+)}\text{.}
\end{align*}
In particular restricts as a linear functional on $\C^{1}_{0}(\overline{\RR^n_+})$, the  action of $\kappa_u$, on $\Psi$, $u$ being fixed, only depends on the boundary value $\psi$ of $\Psi$ since $\kappa_u$ cancels out on ${\C}^{1}_{0,0}(\mathbb{R}^n_+)=\overline{\Ccinfty(\RR^n_+)}^{\lVert \cdot \rVert_{{\W}^{1,\infty}(\mathbb{R}^n_+)}}$. Since the trace operator from ${\W}^{1,\infty}(\mathbb{R}^n_+)$ to ${\W}^{1,\infty}(\partial\RR^n_+)$ is onto, and so is its restriction from $\C^{1}_{0}(\overline{\RR^n_+})$ to $\C^{1}_{0}(\partial\RR^n_+)$. Therefore, by the open mapping theorem, $\kappa_u$ induces and restricts as a linear functional $\tilde{\kappa}_{u} \in (\C^{1}_{0}(\partial\RR^n_+))'$, with $\kappa_u(\Psi) = \tilde{\kappa}_{u} (\psi)$ and one has an estimate
\begin{align*}
    \lvert \tilde{\kappa}_{u}(\psi) \rvert  &\lesssim_{n}\lVert (u,\,\delta u) \rVert_{\L^1(\mathbb{R}^n_+)} \lVert \psi \rVert_{\W^{1,\infty}(\partial\RR^n_+)}\text{.}
\end{align*}
By duality, there exists a unique element $\mathfrak{n}(u)\in \mathcal{M}^{-1}(\partial\RR^n_+,\Lambda^{k-1})$, depending linearly on $u$ such that
\begin{align}\label{eq:IBPPartialTraceProof}
    \langle \mathfrak{n}(u), \psi\rangle_{\partial\RR^n_+} = \int_{\mathbb{R}^{n}_+}  u(x) \cdot\mathrm{d} {\Psi}(x)\mathrm{~d}x - \int_{\mathbb{R}^{n}_+} \mathbf{\delta} u(x) \cdot {\Psi}(x)\mathrm{~d}x.
\end{align}

\textbf{Step 2:} Now, if $p=\infty$, we can perform a similar argument, asking instead for $\psi\in\L^1(\partial\RR^{n}_+,\Lambda^{k-1})$, with $\psi=\Psi_{|_{\partial\mathbb{R}^{n}_+}}$ where $\Psi\in\W^{1,1}({\RR^n_+},\Lambda^{k-1})$ is an arbitrary extension, 
\begin{align*}
    \lvert \kappa_{u}(\Psi) \rvert &\leqslant\lVert (u,\,\delta u) \rVert_{\L^\infty(\mathbb{R}^n_+)} \lVert \Psi \rVert_{{\W}^{1,1}(\mathbb{R}^n_+)}\text{.}
\end{align*}
Since, the trace operator from $\W^{1,1}({\RR^n_+})$ to $\L^1(\partial\RR^n_+)$ is onto, see, \textit{e.g.}, \cite[Theorem~18.13]{Leoni2017},  there exists a unique $\mathfrak{n}(u)\in\L^\infty(\partial\RR^n_+,\Lambda^{k-1})=(\L^1(\partial\RR^n_+,\Lambda^{k-1}))'$ such that \eqref{eq:IBPPartialTraceProof} holds, and with the corresponding estimate.

\textbf{Step 3:} Again, for $u\in\dot{\B}^{s}_{p,q}(\RR^n_+,\Lambda^{k})$ such that $\delta u\in\dot{\B}^{s}_{p,q}(\RR^n_+,\Lambda^{k-1})$, for $\psi\in\S(\partial\RR^{n}_+,\Lambda^{k-1})\subset{\B}^{-s+1-\sfrac{1}{p'}}_{p',q'}(\partial\RR^n_+,\Lambda^{k-1})$ with an extension $\Psi\in\dot{\B}^{-s}_{p',q'}\cap\dot{\B}^{-s+1}_{p',q'}(\RR^n_+,\Lambda^{k-1})$, we define the same way
\begin{align}\label{eq:AbstractIBPPartialTraceProof1}
    \kappa_{u}(\Psi) = \langle u, \d \Psi \rangle_{\RR^n_+} - \langle \delta u, \Psi \rangle_{\RR^n_+}.
\end{align}
Concerning the endpoint spaces, if either $p'$ or $q'$ is finite, we should consider $({\BesSmo}^{-s+1-\sfrac{1}{p'}}_{p',\infty},\dot{\BesSmo}^{-s}_{p',\infty}\cap\dot{\BesSmo}^{-s+1}_{p',\infty})$ instead of $({\B}^{-s+1-\sfrac{1}{p'}}_{p',\infty},\dot{\B}^{-s}_{p',\infty}\cap\dot{\B}^{-s+1}_{p',\infty})$, $({\B}^{-s+1,0}_{\infty,q'},\dot{\B}^{-s,0}_{\infty,q'}\cap\dot{\B}^{-s+1,0}_{\infty,q'})$ instead of $({\B}^{-s+1}_{\infty,q'},\dot{\B}^{-s}_{\infty,q'}\cap\dot{\B}^{-s+1}_{\infty,q'})$ and then $({\BesSmo}^{-s+1,0}_{\infty,\infty},\dot{\BesSmo}^{-s,0}_{\infty,\infty}\cap\dot{\BesSmo}^{-s+1,0}_{\infty,\infty})$ instead of $({\B}^{-s+1}_{\infty,\infty},\dot{\B}^{-s}_{\infty,\infty}\cap\dot{\B}^{-s+1}_{\infty,\infty})$.

By duality, we obtain the estimate
\begin{align*}
    \lvert \kappa_{u}(\Psi) \rvert &\lesssim_{p,s,n}\lVert (u,\,\delta u) \rVert_{\dot{\B}^{s}_{p,q}(\mathbb{R}^n_+)} \left(\lVert \Psi \rVert_{\dot{\B}^{-s}_{p',q'}(\mathbb{R}^n_+)} + \lVert \Psi \rVert_{\dot{\B}^{-s+1}_{p',q'}(\mathbb{R}^n_+)}\right) \text{.}
\end{align*}
Again, $\kappa_u$ is identically zero on $\dot{\B}^{-s}_{p',q'}\cap\dot{\B}^{-s+1}_{p',q',\mathcal{D}}(\RR^n_+)= \overline{\Ccinfty(\RR^n_+)}^{\lVert \cdot \rVert_{\dot{\B}^{-s}_{p',q'}\cap\dot{\B}^{-s+1}_{p',q'}(\RR^n_+)}}$, and similarly for the corresponding endpoint spaces up to the appropriate changes, see Proposition~\ref{prop:DensityIntersectHomSpaces}. Thus, $\kappa_u$ only acts regardless to the choice of the extension $\Psi$ of $\psi$.

Now, as in the previous step, by the ontoness of the trace operator from $\dot{\B}^{-s}_{p',q'}\cap\dot{\B}^{-s+1}_{p',q'}(\mathbb{R}^n_+)$ to ${\B}^{-s+1-\sfrac{1}{p'}}_{p',q'}(\partial\RR^n_+)$\footnote{it also preserves the closure of smooth rapidly decaying functions.} provided by Corollary~\ref{cor:TraceIntersectHomSpacesRn+}, it induces a linear functional $\tilde{\kappa}_{u}$ that can be uniquely represented by an element $\mathfrak{n}(u)\in{\B}^{s-\sfrac{1}{p}}_{p,q}(\partial\RR^n_+,\Lambda^{k-1})$.

The same proof can be achieved for $\dot{\H}^{s,p}(\RR^n_+)$ instead of $\dot{\B}^{s}_{p,p}(\RR^n_+)$, $1<p<\infty$.

\textbf{Step 4:} We prove that, actually, $\tilde{\kappa}_u(\psi)$ only depends on $(\psi_{I})_{I\in\mathcal{I}^{n-1}_{k-1}}$. Hence, we consider to elements $\Psi_1$, $\Psi_2$ such that $((\psi_{1}-\psi_2)_{I})_{I\in\mathcal{I}^{n-1}_{k-1}} =0$. For $\Phi=\Psi_1-\Psi_2$, given $I\in\mathcal{I}^{n-1}_{k-1}$ one has $(\Phi_{I})_{|_{\partial\RR^n_+}}=0$, so we can consider $((\Phi_{j})_{I})_{j\in\NN}\subset \Ccinfty(\RR^{n}_+)$, and $((\Phi_{j})_{I',n})_{j\in\NN}\subset \Ccinfty(\overline{\RR^{n}_+})$, $I'\in\mathcal{I}^{n-1}_{k-2}$, such that $(\Phi_{j})_{j\in\NN}$ converges to $\Phi$. Writing $\d=\d'+\d_n$ $\delta=\delta'+\delta_n$, where $\d'=\nabla'\wedge$ and $\d_n = [\partial_{x_n}]\,\d{x}_{n}\wedge$, and $\delta'= (-\nabla')\iprod$ and $\delta_n = (-\partial_{x_n})\mathfrak{e}_n\iprod$, one writes
\begin{align*}
    \d_n \Phi_j &= \sum_{I'\in\mathcal{I}^{n-1}_{k-1}} (-1)^{k-1} \partial_{x_n} (\Phi_j)_{I'} \, \d x_{I'}\wedge \d x_{n},\\
    \delta_n u &= \sum_{I\in\mathcal{I}^{n}_{k}}  \partial_{x_n} u_{I} \, (-\mathfrak{e}_n)\iprod\d x_{I} = \sum_{I'\in\mathcal{I}^{n-1}_{k-1}} (-1)^{k} \partial_{x_n} u_{I',n} \, \d x_{I'}.
\end{align*}
Where the first equality holds point-wise, and the second one holds in $\mathcal{D}'(\RR^n_+)$.
Consequently, we deduce
\begin{align*}
    \langle u, \d \Phi \rangle_{\RR^n_+} - \langle \delta u, \Phi \rangle_{\RR^n_+} &= \lim_{j\rightarrow \infty } \left( \langle u, \d' \Phi_j \rangle_{\RR^n_+} - \langle \delta' u, \Phi_j \rangle_{\RR^n_+} + \langle u, \delta_n \Phi_j \rangle_{\RR^n_+} - \langle \d_n u, \Phi_j \rangle_{\RR^n_+}\right)\\
    &= \lim_{j\rightarrow \infty } \langle u, \d_n \Phi_j \rangle_{\RR^n_+} - \langle \delta_n u, \Phi_j \rangle_{\RR^n_+} = 0.
\end{align*}
The current step allows in particular to deduce Proposition~\ref{prop:ExtOpDiffFormRn+} and then the density of smooth functions in the nice cases, such as when the third index of Besov spaces $q$ is finite.

\textbf{Step 5:} We show that, actually, $\mathfrak{n}(u) = [(-\mathfrak{e}_n)\iprod u](\cdot,0)$. We assume first that $q<\infty$ while $p\in[1,\infty]$. For $u\in \C^{\infty}_{ub}\cap\dot{\B}^{s}_{p,q}(\overline{\RR^n_+},\Lambda)$ such that $\delta u\in \dot{\B}^{s}_{p,q}({\RR^n_+},\Lambda)$, for $u_t:=u(\cdot,\cdot+t)$ where $t\geqslant0$, one obtains that by smoothness and strong continuity of translations $t\mapsto (-\mathfrak{e}_n)\iprod u(\cdot,t)=\mathfrak{n}(u_t)$  is continuous from $[0,\infty)$ to ${\B}^{s-\sfrac{1}{p}}_{p,q}(\RR^{n-1},\Lambda)$, and one has the estimate
\begin{align*}
    \sup_{t\geqslant 0}\, \lVert [(-\mathfrak{e}_n)\iprod u ](\cdot,t)\rVert_{{\B}^{s-\sfrac{1}{p}}_{p,q}(\RR^{n-1})} \lesssim_{p,s,n}\lVert (u_t,\,\delta u_t) \rVert_{\dot{\B}^{s}_{p,q}(\mathbb{R}^n_+)}\lesssim_{p,s,n}\lVert (u,\,\delta u) \rVert_{\dot{\B}^{s}_{p,q}(\mathbb{R}^n_+)}.
\end{align*}
Therefore, by density provided by Corollary~\ref{cor:ApproxSmoothFuncDivRn+}, one obtains the bounded map
\begin{align*}
 \begin{array}{cllc}
           \{v\in\dot{\B}^{s}_{p,q}({\RR^n_+},\Lambda)\,:\,\delta v\in\dot{\B}^{s}_{p,q}({\RR^n_+},\Lambda)\}& \longrightarrow && \C^{0}_{ub,x_n}(\overline{\RR_+},\B^{s-\frac{1}{p}}_{p,q}(\RR^{n-1},\Lambda))\\
    u&\longmapsto && [x_n\mapsto (-\mathfrak{e}_n)\iprod u (\cdot,x_n)].
    \end{array}
\end{align*}
By uniqueness, one obtains $(-\mathfrak{e}_n)\iprod u(\cdot,0)=\mathfrak{n}(u)$ for all $u\in\dot{\B}^{s}_{p,q}({\RR^n_+},\Lambda)$ such that $\delta u\in\dot{\B}^{s}_{p,q}({\RR^n_+},\Lambda)$. The same goes for $\dot{\B}^{s,0}_{\infty,q}$, $\dot{\BesSmo}^{s}_{p,\infty}$, $\L^1$, or $\dot{\H}^{s,p}$, assuming $1<p<\infty$ for the latter. Furthermore, the density argument shows that $\mathfrak{n}(u)\in\W^{-1,1}(\partial\RR^{n}_+,\Lambda)\subsetneq \mathcal{M}^{-1}(\partial\RR^{n}_+,\Lambda)$, whenever $u,\delta u \in\L^1(\RR^n_+,\Lambda)$.

When $q=\infty$ and $p\in[1,\infty]$, by Corollary~\ref{cor:InterpolationDiffFormRn+}, for all $x_n\geqslant0$, one can apply real interpolation to the map
\begin{align*} \begin{array}{cllc}
           \{v\in\dot{\B}^{s}_{p,r}({\RR^n_+},\Lambda)\,:\,\delta v\in\dot{\B}^{s}_{p,r}({\RR^n_+},\Lambda)\}& \longrightarrow && \B^{s-\sfrac{1}{p}}_{p,r}(\RR^{n-1},\Lambda)\\
    u&\longmapsto && (-\mathfrak{e}_n)\iprod u (\cdot,x_n).
    \end{array}
\end{align*}
where $r\in[1,\infty)$. Thus, we obtain for all $u\in\dot{\B}^{s}_{p,\infty}({\RR^n_+},\Lambda)$ such that $\delta u\in\dot{\B}^{s}_{p,\infty}({\RR^n_+},\Lambda)$, all $x_n\geqslant 0$,
\begin{align*}
    \lVert [(-\mathfrak{e}_n)\iprod u ](\cdot,x_n)\rVert_{{\B}^{s-\sfrac{1}{p}}_{p,\infty}(\RR^{n-1})} \lesssim_{p,s,n}\lVert (u,\,\delta u) \rVert_{\dot{\B}^{s}_{p,\infty}(\mathbb{R}^n_+)}.
\end{align*}
For $u\in\dot{\B}^{s}_{p,\infty}({\RR^n_+},\Lambda)$ such that $\delta u\in\dot{\B}^{s}_{p,\infty}({\RR^n_+},\Lambda)$, by Corollary~\ref{cor:InterpolationDiffFormRn+}, there exist for some $\varepsilon>0$, $v\in\dot{\B}^{s-\varepsilon}_{p,1}({\RR^n_+},\Lambda)$, $w\in \dot{\B}^{s+\varepsilon}_{p,1}({\RR^n_+},\Lambda)$ satisfying $\delta v\in\dot{\B}^{s-\varepsilon}_{p,1}({\RR^n_+},\Lambda)$, $\delta w\in \dot{\B}^{s+\varepsilon}_{p,1}({\RR^n_+},\Lambda)$, and such that $u=v+w$. Hence, by the previous considerations, $$t\mapsto (-\mathfrak{e}_n)\iprod u (\cdot,t)= \mathfrak{n}(u_t) = \mathfrak{n}(v_t) + \mathfrak{n}(w_t)\in \C^0_{ub,x_n}(\overline{\RR}_+,[{\B}^{s-\varepsilon}_{p,1}+{\B}^{s+\varepsilon}_{p,1}]({\RR^{n-1}},\Lambda)).$$ From this, one deduces the weak-$\ast$ continuity on $[0,\infty)$ with values in ${\B}^{s-\sfrac{1}{p}}_{p,\infty}(\RR^{n-1},\Lambda)$, and then the equality $(-\mathfrak{e}_n)\iprod u (\cdot,0)= \mathfrak{n}(u)$ still holds.

Notice that the same arguments from Step 3 to the current Step 5 also apply to the inhomogeneous counterparts obtaining the corresponding map
\begin{align*}
\begin{array}{cllc}
           \{v\in{\B}^{s}_{p,r}({\RR^n_+},\Lambda)\,:\,\delta v\in{\B}^{s}_{p,r}({\RR^n_+},\Lambda)\}& \longrightarrow && \C^0_{ub,x_n}(\overline{\RR_+},\B^{s-\frac{1}{p}}_{p,r}(\RR^{n-1},\Lambda))\\
    u&\longmapsto && (-\mathfrak{e}_n)\iprod u (\cdot,x_n).
    \end{array}
\end{align*}
for all $x_n\geqslant 0$, provided $p\in[1,\infty]$, $r\in[1,\infty)$. Notice that in the case of inhomogeneous function spaces, the proof even simplifies. Hence, the same result also holds for $\B^{s}_{p,\infty}$ and especially $\B^{s}_{\infty,\infty}$.

It remains to deal with the identification in the case of $\L^\infty$. But this follows from the embedding $\mathrm{L}^\infty(\RR^n_+,\Lambda)\hookrightarrow\B^{s}_{\infty,\infty}(\RR^n_+,\Lambda)$, here $-1<s<0$. Indeed, for $u,\delta u\in\mathrm{L}^\infty(\RR^n_+,\Lambda)$, one obtains the equality  $\mathfrak{n}(u_t)=[(-\mathfrak{e}_n)\iprod u_t](\cdot,0)=[(-\mathfrak{e}_n)\iprod u](\cdot,t)$ in $\B^{s}_{\infty,\infty}(\RR^{n-1},\Lambda)$ for all $t\geqslant0$. 

Finally, by uniqueness in the distributional sense, $\mathfrak{n}(u)=[(-\mathfrak{e}_n)\iprod u](\cdot,0)$ in $\mathrm{L}^\infty(\RR^{n-1},\Lambda)= \mathrm{L}^\infty(\partial\RR^{n}_+,\Lambda)$.

\textbf{Step 6:} We provide the argument for the trace space to be homogeneous whenever the condition $\delta u=0 $ is fulfilled. In this case, the identity \eqref{eq:AbstractIBPPartialTraceProof1} becomes
\begin{align}\label{}
    \langle (-\mathfrak{e}_n)\iprod u_{|_{\partial\RR^n_+}}, \psi \rangle_{\partial\RR^n_+} = \langle u, \d \Psi \rangle_{\RR^n_+},
\end{align}
for all $\psi\in\S(\RR^{n-1},\Lambda)$, and $\Psi\in\dot{\B}^{-s+1}_{p',q'}\cap\S(\overline{\RR^n_+},\Lambda)$, such that $\Psi_{|_{\partial\RR^n_+}}=\psi$. Again, the identity above vanishes for $\Psi\in\dot{\B}^{-s+1}_{p',q',\mathcal{D}}(\RR^n_+) = \overline{\Ccinfty(\RR^n_+)}^{\lVert \cdot \rVert_{\dot{\B}^{-s+1}_{p',q'}(\mathbb{R}^n_+)}}$ and hence does not depend on the choice of $\Psi$. So for $\psi\in\S(\RR^{n-1},\Lambda)$, one chooses more specifically $\Psi(x',x_n):=e^{-x_n(-\Delta')^{\sfrac{1}{2}}}\psi(x')$ (component-wise) so that $\Psi\in \dot{\B}^{-s+1}_{p',q'}(\mathbb{R}^n_+,\Lambda)$ by Proposition~\ref{prop:PoissonSemigroup3}, with the estimate
\begin{align*}
    |\langle (-\mathfrak{e}_n)\iprod u_{|_{\partial\RR^n_+}}, \psi \rangle_{\partial\RR^n_+}| \lesssim_{s,n} \lVert u \rVert_{\dot{\B}^{s}_{p,q}(\mathbb{R}^n_+)}\lVert \nabla \Psi \rVert_{\dot{\B}^{-s}_{p',q'}(\mathbb{R}^n_+)}\lesssim_{p,s,n} \lVert u \rVert_{\dot{\B}^{s}_{p,q}(\mathbb{R}^n_+)}\lVert \psi \rVert_{\dot{\B}^{-s+1-\sfrac{1}{p'}}_{p',q'}(\mathbb{R}^{n-1})}.
\end{align*}
By duality, $(-\mathfrak{e}_n)\iprod u_{|_{\partial\RR^n_+}}\in\dot{\B}^{s}_{p,q}(\partial\mathbb{R}^n_+,\Lambda)$, with the desired estimate. The same arguments apply in order to obtain the corresponding result for the remaining function spaces.
\end{proof}

\begin{proposition}\label{prop:DensityIntersectHomSpaces} Let $p,q\in[1,\infty)$, $s\in(-1+\sfrac{1}{p},\sfrac{1}{p})$. The space $\Ccinfty(\RR^{n}_+)$ is a strongly dense subspace of
\begin{enumerate}
    \item $\dot{\H}^{s,p}\cap\dot{\H}^{s+1,p}_{\mathcal{D}}(\RR^n_+)$, $p>1$;
    \item $\dot{\B}^{s}_{p,q}\cap\dot{\B}^{s+1}_{p,q,\mathcal{D}}(\RR^n_+)$;
    \item $\dot{\B}^{s,0}_{\infty,q}\cap\dot{\B}^{s+1,0}_{\infty,q,\mathcal{D}}(\RR^n_+)$;
    \item $\dot{\BesSmo}^{s,0}_{\infty,\infty}\cap\dot{\BesSmo}^{s+1,0}_{\infty,\infty,\mathcal{D}}(\RR^n_+)$.
\end{enumerate}
\end{proposition}

\begin{proof} We give quickly the argument. Let $u$ be an element of $\dot{\B}^{s}_{p,q}\cap\dot{\B}^{s+1}_{p,q,\mathcal{D}}(\RR^n_+)$, then by Propositions~\ref{prop:FundamentalExtby0HomFuncSpaces} and \ref{prop:Ext0DirLip}, $\tilde{u}$ the extension by $0$ belongs to $\dot{\B}^{s}_{p,q,0}\cap\dot{\B}^{s+1}_{p,q,0}(\RR^n_+)\subset\dot{\B}^{s}_{p,q}\cap\dot{\B}^{s+1}_{p,q}(\RR^n)$ with equivalent norms. One can consider a sequence of function $(\varphi_j)_{j\in\NN}\subset\Ccinfty(\RR^n)$ converging towards $\tilde{u}$ in $\dot{\B}^{s}_{p,q}\cap\dot{\B}^{s+1}_{p,q}(\RR^n)$. Now consider the projection operator $\P_0=\I-\E_{\mathcal{N}}^-$, then $(\P_0\varphi_j)_{j\in\NN}\subset\dot{\B}^{s}_{p,q,0}\cap\dot{\B}^{s+1}_{p,q,0}(\RR^n_+)$ by Lemma~\ref{lem:ExtDirNeuRn+} and is a sequence of compactly supported functions that converges to $\tilde{u}$. By translation and mollification one can produce a sequence of smooth compactly supported functions converging towards $\tilde{u}$. Their restriction converges towards $u$ by the definition of function spaces by restriction. 
\end{proof} 

\begin{proposition}\label{prop:ExtOpDiffFormRn+} Let $p,q\in[1,\infty]$, $s\in(-1+\sfrac{1}{p},\sfrac{1}{p})$, $k\in\llb 0,n\rrb$. For all $u\in\dot{\B}^{s}_{p,q}(\RR^n_+,\Lambda^k)$ such that $\delta u \in \dot{\B}^{s}_{p,q}(\RR^n_+,\Lambda^{k-1})$, one has
\begin{align*}
    {\E_{\mathcal{H}_\ast}}u, \delta ({\E_{\mathcal{H}_\ast}}u) \in \dot{\B}^{s}_{p,q}(\RR^n,\Lambda),
\end{align*}
with the estimates
\begin{align*}
    \lVert {\E_{\mathcal{H}_\ast}}u\rVert_{\dot{\B}^{s}_{p,q}(\RR^n)} \lesssim_{s,n}\lVert   u \rVert_{\dot{\B}^{s}_{p,q}(\RR^n_+)}\text{, and }\lVert  \delta ({\E_{\mathcal{H}_\ast}}u ) \rVert_{\dot{\B}^{s}_{p,q}(\RR^n)} \lesssim_{s,n}  \lVert  \delta u \rVert_{\dot{\B}^{s}_{p,q}(\RR^n_+)}\text{.}
\end{align*}
Furthermore, the result still holds if we replace $\dot{\B}^{s}_{p,q}$ by either $\dot{\BesSmo}^s_{p,\infty}$, $\dot{\B}^{s,0}_{\infty,q}$, $\dot{\BesSmo}^{s,0}_{\infty,\infty}$, $\L^1$, $\L^\infty_0$, $\L^\infty$ or even $\dot{\H}^{s,p}$, assuming $p\in(1,\infty)$ in the latter case. A similar result holds with $\d$ instead of $\delta$.
\end{proposition}

\begin{proof} Let $u\in\dot{\B}^{s}_{p,q}(\RR^n_+,\Lambda^k)$ such that $\delta u \in \dot{\B}^{s}_{p,q}(\RR^n_+,\Lambda^{k-1})$. First, by Proposition~\ref{prop:ExtHodgeRn+} notice that $\E_{\mathcal{H}_\ast}u\in\dot{\B}^{s}_{p,q}(\RR^n_+,\Lambda^k)$ and $\E_{\mathcal{H}_\ast}(\delta u)\in\dot{\B}^{s}_{p,q}(\RR^n_+,\Lambda^k)$. Therefore, it is sufficient to show that $\E_{\mathcal{H}_\ast}(\delta u) = \delta \E_{\mathcal{H}_\ast}u$ in $\mathcal{D}'(\RR^n)$.

Let $\Psi\in\Ccinfty(\RR^n,\Lambda^{k-1})$, we write $\d=\d'+\d_n$ where $\d'=\nabla'\wedge$ and $\d_n = [\partial_{x_n}]\,\d{x}_{n}\wedge$, with formal adjoints $\delta'= (-\nabla')\iprod$ and $\delta_n = (-\partial_{x_n})\mathfrak{e}_n\iprod$. We define $\tilde{\Psi}$ by $\tilde{\Psi}_{I} := \Psi_I(\cdot,-\cdot)$ if $I\in\mathcal{I}^{n-1}_{k-1}$, and $\tilde{\Psi}_{I} := - \Psi_I(\cdot,-\cdot)$, if $I=I'\cup\{n\}$ with $I'\in\mathcal{I}^{n-1}_{k-2}$. By Step~4 in the proof of Theorem~\ref{thm:partialtracesDiffFormRn+}, we obtain that $\kappa_u(\Psi_{|_{\partial\RR^n_+}})=-\kappa_{u}(\tilde{\Psi}_{|_{\partial\RR^n_+}})$, and then the following chain of equalities
\begin{align*}
    \langle  \delta(\E_{\mathcal{H}_\ast}u), \d \Psi \rangle_{\RR^n}= \langle \E_{\mathcal{H}_\ast}u, \d \Psi \rangle_{\RR^n} &= \langle u, \d \Psi \rangle_{\RR^n_+} + \langle u_{\mathcal{H}}^{-}, \d \Psi \rangle_{\RR^n_-} \\
    &= \kappa_u(\Psi_{|_{\partial\RR^n_+}}) + \langle \delta u, \Psi \rangle_{\RR^n_+} + \langle u_{\mathcal{H}}^{-}, \d' \Psi \rangle_{\RR^n_-} + \langle u_{\mathcal{H}}^{-}, \d_n \Psi \rangle_{\RR^n_-}\\
    &= \kappa_u(\Psi_{|_{\partial\RR^n_+}}) + \langle \delta u, \Psi \rangle_{\RR^n_+} + \langle u, \d'\tilde{\Psi} \rangle_{\RR^n_+} + \langle u, \d_n\tilde{\Psi}   \rangle_{\RR^n_+}\\
    &= \kappa_u(\Psi_{|_{\partial\RR^n_+}}) + \langle \delta u, \Psi \rangle_{\RR^n_+} + \langle u, \d \tilde{\Psi}  \rangle_{\RR^n_+}\\
    &= \kappa_u(\Psi_{|_{\partial\RR^n_+}}) + \langle \delta u, \Psi \rangle_{\RR^n_+} + \kappa_{u}(\tilde{\Psi}_{|_{\partial\RR^n_+}}) + \langle \delta u, \tilde{\Psi}   \rangle_{\RR^n_+}\\
    &=\langle \delta u, \Psi \rangle_{\RR^n_+} + \langle \delta u, \tilde{\Psi}   \rangle_{\RR^n_+}\\
    &=\langle \delta u, \Psi \rangle_{\RR^n_+} + \langle (\delta u)_{\mathcal{H}}^{-}, {\Psi}   \rangle_{\RR^n_+}\\
    &= \langle \E_{\mathcal{H}_\ast}(\delta u), \Psi \rangle_{\RR^n}.
\end{align*}
Hence, we did obtain the desired the result. 
\end{proof}

\begin{corollary}\label{cor:InterpolationDiffFormRn+}Let $p,p_0,p_1,q,q_0,q_1\in[1,\infty]$, $s_0,s_1,s\in(-1+\sfrac{1}{p},\sfrac{1}{p})$ such that
\begin{align*}
    (s,\text{$\sfrac{1}{p}$})=(1-\theta)(s_0,\text{$\sfrac{1}{p_0}$})+\theta (s_1,\text{$\sfrac{1}{p_1}$}).
\end{align*}
The following assertion hold:
\begin{enumerate}
    \item $(\dot{\D}^{s_0}_{p}(\delta,\RR^{n}_+,\Lambda),\dot{\D}^{s_1}_{p}(\delta,\RR^{n}_+,\Lambda))_{\theta,q}=(\dot{\D}_{p,q_0}^{s_0}(\delta,\RR^{n}_+,\Lambda),\dot{\D}_{p,q_1}^{s_1}(\delta,\RR^{n}_+,\Lambda))_{\theta,q} = \dot{\D}_{p,q}^{s}(\delta,\RR^{n}_+,\Lambda)$, provided $s_0\neq s_1$;
    \item $[\dot{\D}_{p_0,q_0}^{s_0}(\delta,\RR^{n}_+,\Lambda),\dot{\D}_{p_1,q_1}^{s_1}(\delta,\RR^{n}_+,\Lambda)]_{\theta} = \dot{\D}_{p,q}^{s}(\delta,\RR^{n}_+,\Lambda)$, whenever  $\sfrac{1}{q}=\sfrac{(1-\theta)}{q_0}+\sfrac{\theta}{q_1}$ and  $q<\infty$;
    \item $[\dot{\D}^{s_0}_{p}(\delta,\RR^{n}_+,\Lambda),\dot{\D}^{s_1}_{p}(\delta,\RR^{n}_+,\Lambda)]_{\theta} = \dot{\D}^{s}_{p}(\delta,\RR^{n}_+,\Lambda)$, if $1<p_0,p,p_1<\infty$.
\end{enumerate}
The result is also true for $\RR^n$ instead of $\RR^n_+$. A similar result holds with $\d$ instead of $\delta$.
\end{corollary}

\begin{proof} We only deal with the case of Besov spaces, remaining cases can be either obtained similarly, or can be deduced from the embeddings 
\begin{align*}
    \dot{\D}_{p,1}^{s}(\delta,\RR^{n}_+,\Lambda)\hookrightarrow\dot{\D}^{s}_{p}(\delta,\RR^{n}_+,\Lambda)\hookrightarrow\dot{\D}_{p,\infty}^{s}(\delta,\RR^{n}_+,\Lambda).
\end{align*}
We consider the map from $\dot{\D}_{p,q}^{s}(\delta,\RR^{n}_+,\Lambda)$ to $\dot{\B}_{p,q}^{s}\cap\dot{\B}_{p,q}^{s+1}(\RR^{n},\Lambda) \times \dot{\B}_{p,q}^{s}(\RR^{n},\Lambda)$,
\begin{align*}
    \mathfrak{E}\,:\, u\longmapsto ( (\I-\PP_{\RR^n})\E_{\mathcal{H}_\ast}u , \PP_{\RR^n}\E_{\mathcal{H}_\ast}u )
\end{align*}
and from $\dot{\B}_{p,q}^{s}\cap\dot{\B}_{p,q}^{s+1}(\RR^{n},\Lambda) \times \dot{\B}_{p,q}^{s}(\RR^{n},\Lambda)$ to $\dot{\D}_{p,q}^{s}(\delta,\RR^{n}_+,\Lambda)$,
\begin{align*}
    \mathfrak{R}\,:\, (v,w)\longmapsto [(\I-\PP_{\RR^n})v + \PP_{\RR^n}w]_{|_{\RR^n_+}}.
\end{align*}
These maps are well-defined and bounded, thanks to Propositions~\ref{prop:ExtHodgeRn+}~\&~\ref{prop:ExtOpDiffFormRn+}.

One has $\mathfrak{R}\mathfrak{E}=\I$, so by retraction and co-retraction argument, and by cartesian product, it ties to the interpolation of $\dot{\B}_{p,q_0}^{s_0}\cap\dot{\B}_{p,q_0}^{s_0+1}(\RR^{n})$ with $\dot{\B}_{p,q_1}^{s_1}\cap\dot{\B}_{p,q_1}^{s_1+1}(\RR^{n})$, which holds true by Proposition~\ref{prop:InterpHybridBesov}.
\end{proof}

\begin{corollary}\label{cor:ApproxSmoothFuncDivRn+}Let $p,q\in[1,\infty)$, $s\in(-1+\sfrac{1}{p},\sfrac{1}{p})$. The following assertion hold:
\begin{enumerate}
    \item The space $\Ccinfty(\overline{\RR^n_+},\Lambda)$ is universally and strongly dense in $\dot{\D}_{p,q}^{s}(\delta,\RR^{n}_+,\Lambda)$, $\dot{\D}_{\infty,q}^{s,0}(\delta,\RR^{n}_+,\Lambda)$,  $\dot{\mathcal{D}}_{\infty,\infty}^{s,0}(\delta,\RR^{n}_+,\Lambda)$, ${\D}_{1}(\delta,\RR^{n}_+,\Lambda)$ and $\dot{\D}^{s}_p(\delta,\RR^{n}_+,\Lambda)$, assuming $1< p<\infty$ for the latter;
    \item $\C_{ub,h}^\infty(\overline{\RR^n_+},\Lambda)\cap\dot{\B}^{-1}_{\infty,1}(\RR^n_+,\Lambda)$ is universally and strongly dense in $\dot{\D}_{\infty,q}^{s}(\delta,\RR^{n}_+,\Lambda)$,  $\dot{\mathcal{D}}_{\infty,\infty}^{s}(\delta,\RR^{n}_+,\Lambda)$.
\end{enumerate}
A similar result holds for $\d$ instead of $\delta$. The result also holds for inhomogeneous function spaces.
\end{corollary}

\begin{proof} Let $u\in\dot{\D}_{p,q}^{s}(\delta,\RR^{n}_+,\Lambda)$, then $\E_{\mathcal{H}_\ast} u\in\dot{\D}_{p,q}^{s}(\delta,\RR^{n},\Lambda)$ by Propositions~\ref{prop:ExtHodgeRn+}~\&~\ref{prop:ExtOpDiffFormRn+}. Hence, one just considers an approximation $(U_k)_{k\in\NN}\subset \Ccinfty(\RR^n,\Lambda)$ that converges to $\E_{\mathcal{H}_\ast} u$ in $\dot{\D}_{p,q}^{s}(\delta,\RR^{n},\Lambda)$. Therefore, $(u_k)_{k\in\NN}:= (U_k{}_{|_{\RR^n_+}})_{k\in\NN}\in \Ccinfty(\overline{\RR^n_+},\Lambda)$ converges towards $u$ in $\dot{\D}_{p,q}^{s}(\delta,\RR^{n}_+,\Lambda)$ by the definition of function spaces by restriction.

The same arguments apply to the remaining cases.
\end{proof}

\begin{corollary}\label{cor:densityLpLipDiffForms}Let $p\in[1,\infty)$ and $\Omega$ be a special or bounded Lipschitz domain.

The space $\Ccinfty(\overline{\Omega},\Lambda)$ is dense in ${\D}^{p}(\delta,\Omega,\Lambda)$ and $\{v\in\C^0_0(\overline{\Omega},\Lambda)\,:\,\delta v\in\C^0_0(\overline{\Omega},\Lambda)\}$.

Additionally, the space $\C_{ub}^\infty(\overline{\Omega},\Lambda)$ is dense in $\{v\in\C^0_{ub}(\overline{\Omega},\Lambda)\,:\,\delta v\in\C^0_{ub}(\overline{\Omega},\Lambda)\}$ and similarly with $\C_{ub,h}^{\bullet}$ instead $\C_{ub}^{\bullet}$. A similar result holds with $\d$ instead of $\delta$.
\end{corollary}

\begin{proof} \textbf{Step 1:} The $\L^p$-case, $p<\infty$, with $\Omega$ a special Lipschitz domain. Let $u\in{\D}^{p}(\delta,\Omega,\Lambda)$, then for $\Psi\,:\,\RR^n\longmapsto\RR^n$ the globally bi-Lipschitz map such that $\Psi(\RR^n_+)=\Omega$, for $\Psi_{\ast}^{-1}:=(-1)^{k(n-k)} \star \Psi^\ast \star$ on $k$-forms, it holds that
\begin{align*}
    \E_{\mathcal{H}_\ast}\Psi_{\ast}^{-1} u \in{\D}^{p}(\delta,\RR^n,\Lambda),
\end{align*}
with the identity $\delta \E_{\mathcal{H}_\ast}\Psi_{\ast}^{-1} u = \E_{\mathcal{H}_\ast}\Psi_{\ast}^{-1}(\delta u)$. By approximation of $\E_{\mathcal{H}_\ast}\Psi_{\ast}^{-1} u$ by $(V_{k})_{k\in\NN}\subset\Ccinfty(\RR^n,\Lambda)$ provided by Corollary~\ref{cor:ApproxSmoothFuncDivRn+}, it turns out that $(U_k)_{k\in\NN}:=(\Psi_{\ast}V_{k})_{k\in\NN}\subset {\D}^{p}\cap{\D}^{\infty}(\delta, \RR^n,\Lambda)$ has compact support and converges to $\Psi_{\ast}\E_{\mathcal{H}_\ast}\Psi_{\ast}^{-1} u$.  Consequently, $(\tilde{u}_k)_{k\in\NN}:=(U_k{}_{|_{\Omega}})_{k\in\NN}\in {\D}^{p}\cap{\D}^{\infty}(\delta,\Omega,\Lambda)$ converges towards $u$ in ${\D}^{p}(\delta,\Omega,\Lambda)$. One can consider instead the mollification $(u_{k,\varepsilon})_{k\in\NN,\varepsilon>0}:=(\varphi_{\varepsilon
}\ast U_{k}{}_{|_{\Omega}})_{k\in\NN,\varepsilon>0}$ to obtain an approximating sequence in $\C^\infty_{ub}(\overline{\Omega},\Lambda)$. For $\Theta\in\Ccinfty(\RR^n,\RR)$, with $\Theta\geqslant0$, $\Theta_{|_{\B_1(0)}} =1$, and $\Theta_{|_{\B_2(0)^c}} =0$,  one can approximate $u_{k,\varepsilon}$ by
\begin{align*}
    u_{k,\varepsilon}^R:=\Theta(\cdot/R)u_{k,\varepsilon},
\end{align*}
thanks to
\begin{align*}
    \delta ( u_{k,\varepsilon}^R) = \Theta_{R} \, \delta u_{k,\varepsilon} - \frac{1}{R}[\nabla \Theta]_R  \,  \iprod u_{k,\varepsilon},
\end{align*}
and the Dominated Convergence Theorem.

\textbf{Step 2:} The case $\C^0_0$ and $\mathrm{C}^0_{ub}$  with $\Omega$ a special Lipschitz domain. First notice that by strong continuity of translations on $\mathrm{C}_{ub}^0$, for $t>0$, $u_t\,:\,(x',x)\mapsto u(x',x_n+t)$, is such that $(u_t)_{t>0}$ converges strongly towards $u$ as $t$ goes to $0$. Now for $t>0$ being fixed, we set for $x\in\RR^n$

For $\psi\in\Ccinfty((-t,t))$ such that $0\leqslant\psi\leqslant1$, and satisfying $\psi_{|_{[-t/2,t/2]}}=1$, we set
 \begin{equation*}
    \tilde{\psi}_t(x)=\left\{ \begin{array}{rl}
         0 &\text{, if }x\in\Omega^{c}\text{ with } c\mathfrak{d}^\ast(x)\geqslant t\text{,}\\
        \psi(c\mathfrak{d}^\ast(x))&\text{, if }x\in\Omega^{c}\text{ with } 0\leqslant c\mathfrak{d}^\ast(x)\leqslant t\text{,}\\
        1 &\text{, if } x\in \Omega\text{.}
    \end{array}
    \right.,
\end{equation*}
where $0<c<1$ and $\mathfrak{d}^\ast$ are such that $\mathfrak{d}^\ast\sim_{\partial\Omega,n} \d(\cdot,\Omega)$, with $0\leqslant \phi(x')-x_n \leqslant c\mathfrak{d}^\ast(x)$ for all $x\in\Omega^c$, $\mathfrak{d}^\ast \in\C^{0}(\overline{\Omega^c})\cap\C^{\infty}(\overline{\Omega}^c)$ and $\nabla\mathfrak{d}^\ast\in \C^0_{b}(\overline{\Omega}^c,\RR^n)$. See \cite[Chapter~VI,~Section~2.1,~Theorem~2,~\&~~Section~3.2.1,~Lemma~2]{Stein1970} for more details. One can check that $\psi_t\in\C^1_{ub}(\RR^n)$ and $\psi_t{}_{|_{\Omega^c}}(\cdot,\cdot-t)=0$. In this case $\tilde{u}_t:=\tilde{\psi}_t \cdot u_t \in \{v\in\C^0_0(\RR^n,\Lambda)\,:\,\delta v\in\C^0_0(\RR^n,\Lambda)\}$.
Finally, standard techniques such has multiplication by a cut-off $\Theta_R=\Theta(\sfrac{\cdot}{R})$ and mollification yields the result, since, again, one has the formula
\begin{align*}
    \delta ( \Theta_{R} \tilde{u}_t) = \Theta_{R} \, \delta \tilde{u}_t - \frac{1}{R}[\nabla \Theta]_R  \,  \iprod \tilde{u}_t.
\end{align*}
In the case of $\C_{ub}^0$, the cutoff procedure can be skipped.

\textbf{Step 3:} The case of bounded Lipschitz domains reduces to the case of special Lipschitz domains by localization and rotation.
\end{proof}

\paragraph{The case of Lipschitz domains}

Now we prove the counterpart of Theorem~\ref{thm:partialtracesDiffFormRn+} for Lipschitz domains.

\begin{theorem}\label{thm:partialtracesDiffFormLip} Let $p,q\in[1,\infty]$, and $s\in(-1+\sfrac{1}{p},\sfrac{1}{p})$ and $k\geqslant 1$. Let $\Omega$ be either a bounded or special Lipschitz domain.
\begin{enumerate}
    \item For $u\in{\L}^{1}(\Omega,\Lambda^k)$ such that $\delta u \in {\L}^{1}(\Omega,\Lambda^{k-1})$, one can define the partial trace
    \begin{align*}
    \nu\iprod u_{|_{\partial\Omega}} \in \W^{-1,1}(\partial\Omega,\Lambda^{k-1}),
\end{align*}
with the estimate
\begin{align*}
    \lVert  \nu\iprod u_{|_{\partial\Omega}} \rVert_{ \W^{-1,1}(\partial\Omega)} \lesssim_{n,\partial\Omega} \lVert  (u,\,\delta u) \rVert_{{\L}^{1}(\Omega)}.
\end{align*}

\item For $u\in{\L}^{\infty}(\Omega,\Lambda^k)$ such that $\delta u \in {\L}^{\infty}(\Omega,\Lambda^{k-1})$, one can define the partial trace
\begin{align*}
    \nu\iprod u_{|_{\partial\Omega}} \in \L^{\infty}(\partial\Omega,\Lambda^{k-1}),
\end{align*}
with the estimate
\begin{align*}
    \lVert  \nu\iprod u_{|_{\partial\Omega}} \rVert_{{\L}^{\infty}(\partial\Omega)} \lesssim_{n,\partial\Omega} \lVert  (u,\,\delta u) \rVert_{{\L}^{\infty}(\Omega)}.
\end{align*}
    \item For $u\in{\B}^{s}_{p,q}(\Omega,\Lambda^k)$ such that $\delta u\in {\B}^{s}_{p,q}(\Omega,\Lambda^{k-1})$, one can define the partial trace
\begin{align*}
    \nu\iprod u_{|_{\partial\Omega}} \in {\B}^{s-\sfrac{1}{p}}_{p,q}(\partial\Omega,\Lambda^{k-1}),
\end{align*}
with the estimate
\begin{align*}
    \lVert  \nu\iprod u_{|_{\partial\Omega}} \rVert_{{\B}^{s-\sfrac{1}{p}}_{p,q}(\partial\Omega)} \lesssim_{s,n,\partial\Omega} \lVert  (u,\,\delta u) \rVert_{{\B}^{s}_{p,q}(\Omega)}.
\end{align*}
If $p=q\in(1,\infty)$, the result still holds with ${\H}^{s,p}$ instead of ${\B}^{s}_{p,p}$.
\end{enumerate}
Furthermore, a similar result holds with $(\nu\wedge\cdot,\, \d,\, \Lambda^{k+1})$ replacing $(\nu\iprod\cdot,\, \delta,\, \Lambda^{k-1})$.
\end{theorem}

\begin{proof} As in the proof of Theorem~\ref{thm:partialtracesDiffFormRn+}, by the ontoness of the trace operator $$[\cdot]_{|_{\partial\Omega}}\,:\,\B^{-s+1}_{p',q'}(\Omega),\H^{-s+1,p'}(\Omega)\longrightarrow\B^{-s+1-\sfrac{1}{p'}}_{p',q'}(\partial\Omega)$$ and the density of $\Ccinfty(\Omega)$ in $\B^{-s+1}_{p',q',\mathcal{D}}(\Omega),\H^{-s+1,p'}_\mathcal{D}(\Omega)$, both allow to define uniquely a map $\mathfrak{n}_{\partial\Omega}$, satisfying
\begin{align*}
        \big\langle  \mathfrak{n}_{\partial\Omega}(u), \varphi_{|_{\partial\Omega}} \big\rangle_{\partial\Omega}=\big\langle  u, \d \varphi \big\rangle_{\Omega}-\big\langle \delta u, \varphi\big\rangle_{\Omega}.
\end{align*}

From now on, we just have to show the distributional equality
\begin{align}\label{eq:proofequalityLipPartialTrace}
    \mathfrak{n}_{\partial\Omega}(u)=\nu\iprod u_{|_{_{\partial\Omega}}}
\end{align}
In particular, if $p\in[1,\infty)$ smooth functions are dense in by $\D_{p}(\delta,\Omega,\Lambda)$ by Corollary~\ref{cor:densityLpLipDiffForms}. Therefore, \eqref{eq:proofequalityLipPartialTrace} holds in the $\L^p$-case.

By Sobolev embeddings it holds for the $\H^{s,p}$ and $\B^{s}_{p,q}$-cases  $p,q\in[1,\infty]$, $p<\infty$, $0<s<1/p$.

Now, one should stop temporarily this proof and check Proposition~\ref{prop:DiffFormDomainsclosedEtc..}.

If $s<0$, $p,q\in(1,\infty)$, by Proposition~\ref{prop:DiffFormDomainsclosedEtc..} and a property of interpolation scales, one obtains for $\varepsilon>0$ arbitrarily small, that $\D_{p}(\delta,\Omega,\Lambda)=\D_{p}(\delta,\Omega,\Lambda)\cap \D_{p}^{s-\varepsilon}(\delta,\Omega,\Lambda)$ is strongly dense in $\D_{p}^s(\delta,\Omega,\Lambda)$ and $\D_{p,q}^s(\delta,\Omega,\Lambda)$. Since $\B_{p,\infty}^s(\Omega)\hookrightarrow\B_{p,1}^{s-\varepsilon}(\Omega)$, we finally did obtain that \eqref{eq:proofequalityLipPartialTrace} holds in the cases $\H^{s,p}$, $\B^{s}_{p,q}$, $p\in[1,\infty)$, $q\in[1,\infty]$, assuming $s\neq0$. The case $s=0$ follows from the embedding $\B^{0}_{p,q}\hookrightarrow\B^{-\varepsilon}_{p,q}$, for arbitrarily small $0<\varepsilon<1-1/p$. The case $p=\infty$ follows by localization and the previous cases, since for bounded domains one has $\L^{\infty}(\Omega)\hookrightarrow\L^r(\Omega)$, and more generally, for $s\in(-1,0)$, $\B^{s}_{\infty,q}(\Omega)\hookrightarrow\B^{s}_{r,q}(\Omega)$, where $r<\infty$ is large enough, such that $s>-1+\frac{1}{r}$.
\end{proof}

In the case of homogeneous function spaces of differential forms on special Lipschitz domains, we do obtain a similar result with a proof similar in spirit but different (the problem here being the case $s>0$). This is a substantial improvement with regard to  \cite[Theorem~A.1]{Gaudin2023Hodge}.

\begin{theorem}\label{thm:partialtracesDiffFormHomSpeLip} Let $p,q\in[1,\infty]$, and $s\in(-1+\sfrac{1}{p},\sfrac{1}{p})$ and $k\geqslant 1$. Let $\Omega$ be a special Lipschitz domain.

For $u\in\dot{\B}^{s}_{p,q}(\Omega,\Lambda^k)$ such that $\delta u\in \dot{\B}^{s}_{p,q}(\Omega,\Lambda^{k-1})$, one can define the partial trace
\begin{align*}
    \nu\iprod u_{|_{\partial\Omega}} \in {\B}^{s-\sfrac{1}{p}}_{p,q}(\partial\Omega,\Lambda^{k-1}),
\end{align*}
with the estimate
\begin{align*}
    \lVert  \nu\iprod u_{|_{\partial\Omega}} \rVert_{{\B}^{s-\sfrac{1}{p}}_{p,q}(\partial\Omega)} \lesssim_{s,n,\partial\Omega} \lVert  (u,\,\delta u) \rVert_{\dot{{\B}}^{s}_{p,q}(\Omega)}.
\end{align*}
If $p=q\in(1,\infty)$, the result still holds with $\dot{\H}^{s,p}$ instead of $\dot{\B}^{s}_{p,p}$.
\end{theorem}

\begin{proof} \textbf{Step 1:} the case $p\in(1,\infty]$, $-1+\sfrac{1}{p}<s<0$, one has $\dot{{\B}}^{s}_{p,q}(\Omega)\hookrightarrow{{\B}}^{s}_{p,q}(\Omega)$, and consequently $u\in\dot{\B}^{s}_{p,q}(\Omega,\Lambda^k)$ and $\delta u\in \dot{\B}^{s}_{p,q}(\Omega,\Lambda^{k-1})$ implies $u\in{\B}^{s}_{p,q}(\Omega,\Lambda^k)$ and $\delta u\in {\B}^{s}_{p,q}(\Omega,\Lambda^{k-1})$, so that, by Theorem~\ref{thm:partialtracesDiffFormLip},
\begin{align*}
    \lVert  \nu\iprod u_{|_{\partial\Omega}} \rVert_{{\B}^{s-\sfrac{1}{p}}_{p,q}(\partial\Omega)} \lesssim_{s,n,\partial\Omega} \lVert  (u,\,\delta u) \rVert_{{{\B}}^{s}_{p,q}(\Omega)} \lesssim_{s,n,\partial\Omega} \lVert  (u,\,\delta u) \rVert_{\dot{{\B}}^{s}_{p,q}(\Omega)}.
\end{align*}
Similarly, for $-1+\sfrac{1}{p}<s\leqslant 0$, $p<\infty$, $u\in\dot{\H}^{s,p}(\Omega,\Lambda^k)$ and $\delta u\in \dot{\H}^{s,p}(\Omega,\Lambda^{k-1})$ are yielding $\nu\iprod u_{|_{\partial\Omega}}\in {\B}^{s-\sfrac{1}{p}}_{p,p}(\partial\Omega)$.

\textbf{Step 2:} The case $0<s<\sfrac{1}{p}$, $p\in[1,\infty)$, we reconstruct the partial trace and use a density argument provided by interpolation theory. For all $\Psi\in\Ccinfty(\overline{\Omega},\Lambda^{k-1})$, we set
\begin{align*}
    \kappa_{u}(\Psi) := \langle u, \d \Psi \rangle_{\Omega} - \langle \delta u, \Psi \rangle_{\Omega}.
\end{align*}
Here we have,
\begin{align*}
|\kappa_{u}(\Psi)| \leqslant \lVert  (u,\,\delta u) \rVert_{\dot{{\B}}^{s}_{p,q}(\Omega)}\lVert \Psi \rVert_{\dot{{\B}}^{-s}_{p',q'}\cap\dot{{\B}}^{-s+1}_{p',q'}(\Omega)}.
\end{align*}
Since, $0<s<\sfrac{1}{p}$, this implies $\sfrac{1}{p'}<-s+1<1$ --this is really of a paramount importance here--, so that, thanks to \cite[Propositions~4.20,~4.21~\&~Remark~4.22,~\textit{(ii)}]{Gaudin2023Lip} in combination with Proposition~\ref{prop:DumpedPoissonSemigroup2} and \cite[Propositions~3.20,~3.39~\&~Remarks~3.21~\&~3.40]{Gaudin2023Lip}, one has
\begin{itemize}
    \item Density of $\Ccinfty(\Omega,\CC)$ in $\dot{{\B}}^{-s}_{p',q'}\cap\dot{{\B}}^{-s+1}_{p',q',\mathcal{D}}(\Omega)$ (resp. $\dot{{\BesSmo}}^{-s}_{p',\infty}\cap\dot{{{\BesSmo}}}^{-s+1}_{p',\infty,\mathcal{D}}(\Omega)$, $\dot{{\B}}^{-s,0}_{\infty,q'}\cap\dot{{\BesSmo}}^{-s+1,0}_{\infty,q',\mathcal{D}}(\Omega)$, and $\dot{{\BesSmo}}^{-s,0}_{\infty,\infty}\cap\dot{{{\BesSmo}}}^{-s+1,0}_{\infty,\infty,\mathcal{D}}(\Omega)$);
    \item Surjectivity of the trace operator $$[\cdot]_{|_{\partial\Omega}}\,:\, \dot{{\B}}^{-s}_{p',q'}\cap\dot{{\B}}^{-s+1}_{p',q'}(\Omega) \longrightarrow {{\B}}^{-s+1-\sfrac{1}{p'}}_{p',q'}(\partial\Omega),$$
    (resp. from $\dot{{\BesSmo}}^{-s}_{p',\infty}\cap\dot{{\BesSmo}}^{-s+1}_{p',\infty}(\Omega)$ to $\dot{{\BesSmo}}^{-s+1-\sfrac{1}{p'}}_{p',\infty}(\partial\Omega)$,  from $\dot{{\B}}^{-s,0}_{\infty,q'}\cap\dot{{\BesSmo}}^{-s+1,0}_{\infty,q'}(\Omega)$ to ${{\B}}^{-s+1,0}_{\infty,q'}(\partial\Omega)$, and from $\dot{{\BesSmo}}^{-s,0}_{\infty,\infty}\cap\dot{{\BesSmo}}^{-s+1,0}_{\infty,\infty}(\Omega)$ to ${{\BesSmo}}^{-s+1,0}_{\infty,\infty}(\partial\Omega)$).
\end{itemize}
Therefore, there exists a uniquely well defined $\mathfrak{n}_{\partial\Omega}(u)\in {\B}^{s-\sfrac{1}{p}}_{p,q}(\partial\Omega,\Lambda^{k-1})$ such that for all $\Psi\in\Ccinfty(\overline{\Omega},\Lambda^{k-1})$, 
\begin{align*}
    \langle\mathfrak{n}_{\partial\Omega}(u), \Psi_{|_{\partial\Omega}}\rangle_{\partial\Omega} = \langle u, \d \Psi \rangle_{\Omega} - \langle \delta u, \Psi \rangle_{\Omega}.
\end{align*}
By Sobolev embeddings and the case $s<0$, one can deduce that $\mathfrak{n}_{\partial\Omega}(u)=\nu\iprod u_{|_{\partial\Omega}}$.
\end{proof}

\section{Around vector-valued function spaces}\label{App:VectValFuncSpaces}

In order to show some Fubini properties for scalar-valued homogeneous Besov spaces, we need to collect information about their vector-valued counterpart, on $\RR^n$ and $\RR^n_+$, $n\geqslant1$.

Following \cite[Chapters~1~\&~2]{bookBahouriCheminDanchin}, \cite[Sections~1~\&~2]{Gaudin2022}, \cite[Sections~1~\&~2]{Gaudin2023} for homogeneous scalar-valued function spaces, and \cite{SchmeißerSickel2005,ScharfSchmeißerSickel2011,MeyriesVeraar2012,Amann1997VectDistrib,Amann2019BookVolII,HytonenNeervenVeraarWeisbookVolIII2023} for the case of vector-valued inhomogeneous function spaces and the related distribution theory, we introduce
\begin{align*}
    \S'(\RR^n,\X) := \mathcal{L}(\S(\RR^n),\X),
\end{align*}
and the natural subspace of homogeneous distributions
\begin{align*}
    \S'_h(\RR^n,\X):= \Big\{\,u\in\S'(\RR^n,\X) \,:\, \forall \Theta\in\Ccinfty(\RR^n),\,\lVert \Theta(\lambda\mathfrak{D})u\rVert_{\L^\infty(\RR^n,\X)}\xrightarrow[\lambda\rightarrow\infty]{} 0\,\Big\}.
\end{align*}
We consider $(\dot{\Delta}_{j})_{j\in\mathbb{Z}}$ a homogeneous Littlewood-Paley decomposition, for $u\in\S'_h(\RR^n,\X)$, one has the identity
\begin{align*}
    \sum_{j\in\ZZ} \dot{\Delta}_{j}u = u \text{, in }\S'(\RR^n,\X),
\end{align*}
since  the $\S'_h$-condition ensures that one obtains the convergence of the low frequency part $\sum_{j<0} \dot{\Delta}_{j}u$ in $\L^\infty(\RR^n,\X)$ (even in $\C_{ub}^0(\RR^n,\X)$).

\begin{definition}\label{def:HomFuncSpacesRn} For any $p,q\in[1,\infty]$, $s\in\RR$, and $k\in\mathbb{N}$, we define
\begin{itemize}[label={$\bullet$}]
    \item the homogeneous Sobolev (Riesz potential) spaces,
    $$\dot{\H}^{s,p}(\RR^n,\X)=\left\{ u\in\mathcal{S}'_h(\RR^n,\X)\,:\,\left\lVert {u} \right\rVert_{\dot{\H}^{s,p}(\RR^n,\X)}:=\Big\lVert \sum_{j\in\ZZ} (-\Delta)^{s/2}\dot{\Delta}_j{u} \Big\rVert_{{\L}^{p}(\RR^n,\X)}<\infty \right\}\text{;}$$
    \item the standard homogeneous Sobolev spaces,
    $$\dot{\W}^{k,p}(\RR^n,\X)=\left\{ u\in\mathcal{S}'_h(\RR^n,\X) \,:\, \lVert {u} \rVert_{\dot{\W}^{k,p}(\RR^n,\X)}=\lVert \nabla^k{u} \rVert_{{\L}^{p}(\RR^n,\X)}<\infty \right\}\text{;}$$
    \item and the homogeneous Besov spaces,
    $$\dot{\B}^{s}_{p,q}(\RR^n,\X)=\left\{ u\in\mathcal{S}'_h(\RR^n,\X)\,:\,\left\lVert (2^{js}\dot{\Delta}_j{u})_{j\in\ZZ} \right\rVert_{\ell^q(\ZZ,\L^p(\RR^n,\X))}<\infty \right\}\text{.}
    $$
\end{itemize}
These are all normed vector spaces. For all $m\in\mathbb{N}$, we set
\begin{itemize}
    \item $\L^\infty_h(\RR^n,\X):= \L^\infty(\RR^n,\X)\cap \mathcal{S}'_h(\RR^n,\X)$ ($=\dot{\W}^{0,\infty}(\RR^n,\X)$),
    \item $\dot{\C}^{m}_{0}(\RR^n,\X)$ is the space ${\C}^{m}_{0}(\RR^n,\X)$ equipped with the (semi-)norm $\lVert {\cdot} \rVert_{\dot{\W}^{m,\infty}(\RR^n,\X)}$.
    \item $\dot{\C}^{m}_{ub,h}(\RR^n,\X)$ is the space ${\C}^{m}_{ub}(\RR^n,\X)\cap\S'_h(\RR^n,\X)$\footnote{The space of uniformly continuous and bounded functions such that derivatives up to the order $m$ are also uniformly continuous and bounded.} equipped with the (semi-)norm $\lVert {\cdot} \rVert_{\dot{\W}^{m,\infty}(\RR^n,\X)}$.
\end{itemize}

We can consider, for $p,q\in[1,\infty]$, $s\in\RR$, the following spaces
\begin{align*}
    \dot{\mathcal{B}}^{s}_{p,\infty}(\RR^n,\X):=&\Big\{u\in \dot{\B}^{s}_{p,\infty}(\RR^n) \,:\, \lim\limits_{|j|\rightarrow\infty} 2^{js} \lVert\dot{\Delta}_{j}u\rVert_{\L^p(\RR^n,\X)}=0\Big\},\\
    \dot{\B}^{s,0}_{\infty,q}(\RR^n,\X):=&\Big\{u\in \dot{\B}^{s}_{\infty,q}(\RR^n,\X)\,:\, (\dot{\Delta}_ju)_{j\in\mathbb{Z}}\subset\mathrm{C}_0^0(\RR^n,\X)\Big\},\\
    \dot{\mathcal{B}}^{s,0}_{\infty,\infty}(\RR^n,\X):=& \dot{\mathcal{B}}^{s}_{\infty,\infty}(\RR^n,\X)\cap\dot{\B}^{s,0}_{\infty,\infty}(\RR^n,\X).
\end{align*}
\end{definition}

Before we start, we mention that a preliminary study of homogeneous Bessel potential spaces $\dot{\H}^{s,p}(\RR^n,\X)$ on the whole and the half-line has been started in \cite[Section~3]{Gaudin2023} under the assumption that $\X$ has the \textbf{UMD} property. This was necessary to recover the most meaningful properties, as the one obtained in the scalar-valued case. \textbf{Here, we do not make such assumption on the Banach space $\X$}.

The following result holds by a close inspection of the proofs in \cite{Gaudin2023Lip} which shows that most of the arguments remain valid in the vector-valued case. We state the following Meta-Theorem.

\begin{theorem}\label{thm:MetaThmBanachValuedHombesovSpaces}The results Proposition~2.1, Lemma 2.7, Proposition~2.8, Theorem~2.11, and Theorem~2.12 (except the point \textit{(v)}) from \cite[Section~2]{Gaudin2023Lip} remain valid for our definition of vector-valued homogeneous function spaces.
\end{theorem}

We also notice and recall the very well known equivalence of norms which holds for all $p,q\in[1,\infty]$, $s\in\RR$,
\begin{align}\label{eq:EquivNormHomGenInhomogenBesov}
    \lVert u \rVert_{{\B}^{s}_{p,q}(\RR^n,\X)}\sim_{p,s,n}^{\X} \lVert u \rVert_{{\L}^{p}(\RR^n,\X)}+\lVert u \rVert_{\dot{\B}^{s}_{p,q}(\RR^n,\X)}\text{, } u\in{\B}^{s}_{p,q}(\RR^n,\X),
\end{align}
giving the equality of vector spaces ${\B}^{s}_{p,q}(\RR^n,\X)={\L}^{p}(\RR^n,\X)\cap\dot{\B}^{s}_{p,q}(\RR^n,\X)$, whenever $p\in[1,\infty)$, $q\in[1,\infty]$, $s>0$. When $p=\infty$, $s>0$, one has the equality of sets $\dot{\B}^{s}_{\infty,q}(\RR^n,\X) = {\B}^{s}_{\infty,q}(\RR^n,\X)\cap\S'_h(\RR^n,\X)$.

More precisely, in the case of Besov spaces, writing for simplicity $\L^p(\RR^0,\X)=\X$ whenever $p\in[1,\infty]$, one has the following proposition.

\begin{proposition}\label{prop:BanachValuedHomBesovSpaces} Let $p,q\in[1,\infty]$, $s\in\RR$, $n\in\mathbb{N}^\ast$, $\theta\in(0,1)$ and $\X$ be a Banach space. The following properties hold:
\begin{enumerate}
    \item We assume $p,q<\infty$. The space $\mathcal{S}_0(\RR^n,\X)$ is a (universally) dense subspace in $\dot{\B}^{s}_{p,q}(\RR^n,\X)$. The subspace $\Ccinfty(\RR^n,\X)\cap\dot{\B}^{s}_{p,q}(\RR^n,\X)$ is dense in $\dot{\B}^{s}_{p,q}(\RR^n,\X)$. Moreover, under the additional condition $s>-n/p'$, the space $\Ccinfty(\RR^n,\X)$ is also a (universally) dense subspace. A similar result holds for the spaces $\dot{\B}^{s,0}_{\infty,q}(\RR^n,\X)$, $\dot{\BesSmo}^{s}_{p,\infty}(\RR^n,\X)$ and $\dot{\BesSmo}^{s,0}_{\infty,\infty}(\RR^n,\X)$.
    \item $\dot{\B}^{s}_{p,q}(\RR^n,\X)$ is complete whenever \eqref{AssumptionCompletenessExponents} is satisfied.
    \item If $q_0,q_1\in[1,\infty]$, $s_0\neq s_1\in\RR$ and $(1-\theta)s_0+\theta s_1 = s$, then
    \begin{align*}
        \Big(\dot{\B}^{s_0}_{p,q_0}(\RR^n,\X),\dot{\B}^{s_1}_{p,q_1}(\RR^n,\X)\Big)_{\theta,q}=\dot{\B}^{s}_{p,q}(\RR^n,\X).
    \end{align*}
    \item If $p_0,p_1,q_0,q_1\in[1,\infty]$, $s_0,s_1\in\RR$ satisfying \hyperref[AssumptionCompletenessExponents]{$(\mathcal{C}_{s_j,p_j,q_j})$}, $j\in\{0,1\}$, then if $q< \infty$,
    \begin{align*}
        \Big[\dot{\B}^{s_0}_{p_0,q_0}(\RR^n,\X),\dot{\B}^{s_1}_{p_1,q_1}(\RR^n,\X)\Big]_{\theta}=\dot{\B}^{s}_{p,q}(\RR^n,\X)
    \end{align*}
    when $(1-\theta)s_0+\theta s_1 = s$ and $(1-\theta)/r_0+\theta/r_1 = 1/r$ for $r\in\{p,q\}$.
    \item If $p,q\neq1$, and \eqref{AssumptionCompletenessExponents} is satisfied, one has
    \begin{align*}
        \Big(\dot{\B}^{-s}_{p',q'}(\RR^n,\X)\Big)^\ast = \dot{\B}^{s}_{p,q}(\RR^n,\X^\ast) \text{ and } \Big(\dot{\mathcal{B}}^{-s}_{p',\infty}(\RR^n,\X)\Big)^\ast = \dot{\B}^{s}_{p,1}(\RR^n,\X^\ast).
    \end{align*}
    \item We have the following embeddings into Lebesgue spaces:
    \begin{enumerate}
        \item $\dot{\B}^{\frac{n}{p}-\frac{n}{r}}_{p,1}(\RR^n,\X) \hookrightarrow \L^{r}(\RR^n,\X)$, $r\in[p,\infty]$;
        \item $\dot{\B}^{\frac{1}{p}-\frac{1}{r}}_{p,1}(\RR^n,\X) \hookrightarrow \L^{r}(\RR,\L^p(\RR^{n-1},\X))$, $r\in[p,\infty]$.
    \end{enumerate}
    \item If $s\in(0,1)$, the quantity defined for $u\in\L^1_{\text{loc}}(\RR^n,\X)$,
    \begin{align*}
        \sum_{k=1}^n \left\lVert t\mapsto t^{-s}\lVert\tau^{k}_{t}u-u\rVert_{\L^p(\RR^n,\X)} \right\rVert_{\L^q_\ast(\RR_+)}
    \end{align*}
    is an equivalent norm on $\dot{\B}^{s}_{p,q}(\RR^n,\X)$, where $\tau^k_t f (x):= f(x+t\,\mathfrak{e}_k)$.
\end{enumerate}
\end{proposition}

\begin{proof} We only prove \textit{(vi)} and \textit{(vii)}, otherwise the other claimed results are carried over by the Meta-Theorem~\ref{thm:MetaThmBanachValuedHombesovSpaces}. However, for the duality argument, taking into account the vector-valued case, one should modify the proof according to \cite[Theorem~14.4.34]{HytonenNeervenVeraarWeisbookVolIII2023}.

\textbf{Step 1:} Proof of \textit{(vi)}.

\textbf{Step 1.1:} Point \textit{(a)}. Let $u\in\dot{\B}^{n/p}_{p,1}(\RR^n,\X)$, $p\in[1,\infty]$. By Bernstein's inequality
\begin{align*}
    \lVert u\rVert_{\L^\infty(\RR^n,\X)} \leqslant \sum_{j\in\ZZ} \lVert \dot{\Delta}_ju\rVert_{\L^\infty(\RR^n,\X)} \lesssim_{p,s,n}  \sum_{j\in\ZZ} 2^{j\frac{n}{p}}\lVert \dot{\Delta}_ju\rVert_{\L^p(\RR^n,\X)} = \lVert u\rVert_{\dot{\B}^{n/p}_{p,1}(\RR^n,\X)}.
\end{align*}
Similarly, and in an easier way
\begin{align*}
    \lVert u\rVert_{\L^p(\RR^n,\X)} \leqslant \sum_{j\in\ZZ} \lVert \dot{\Delta}_ju\rVert_{\L^p(\RR^n,\X)} = \lVert u\rVert_{\dot{\B}^{0}_{p,1}(\RR^n,\X)}.
\end{align*}

Thus, by complex interpolation provided by point \textit{(iv)}, we obtain \textit{(vi)}-\textit{(a)}.

\textbf{Step 1.2:} Point \textit{(b)}. By the previous point \textit{(vi)}-\textit{(a)}, with $n=1$ and replacing $\X$ by $\L^p(\RR^{n-1},\X)$, one obtains, provided $p\in[1,\infty)$, $r\in[p,\infty]$,
\begin{align*}
    \dot{\B}^{\frac{1}{p}-\frac{1}{r}}_{p,1}(\RR,\L^p(\RR^{n-1},\X)) \hookrightarrow \L^{r}(\RR,\L^p(\RR^{n-1},\X)).
\end{align*}
It remains to show $\dot{\B}^{\frac{1}{p}-\frac{1}{r}}_{p,1}(\RR^n,\X)\hookrightarrow\dot{\B}^{\frac{1}{p}-\frac{1}{r}}_{p,1}(\RR,\L^p(\RR^{n-1},\X))$. This follows by real interpolation, point \textit{(iii)}, for $\theta\in(0,1)$, $m\in\NN$, $s=m\theta$,
\begin{align*}
    \dot{\B}^{s}_{p,q}(\RR^n,\X)=\Big(\dot{\B}^{0}_{p,1}(\RR^n,\X),\dot{\B}^{m}_{p,1}(\RR^n,\X)\Big)_{\theta,q} &\hookrightarrow \Big({\L}^{p}(\RR^n,\X),\dot{\W}^{m,p}(\RR^n,\X)\Big)_{\theta,q}\\
    &\hookrightarrow \Big({\L}^{p}(\RR,{\L}^{p}(\RR^{n-1},\X)),\dot{\W}^{m,p}(\RR,{\L}^{p}(\RR^{n-1},\X))\Big)_{\theta,q}\\
    &\hookrightarrow \Big(\dot{\B}^{0}_{p,\infty}(\RR,{\L}^{p}(\RR^{n-1},\X)),\dot{\B}^{m}_{p,\infty}(\RR,{\L}^{p}(\RR^{n-1},\X))\Big)_{\theta,q}\\
    &= \dot{\B}^{s}_{p,q}(\RR,{\L}^{p}(\RR^{n-1},\X)).
\end{align*}

If $p=\infty$,
\begin{align*}
    \dot{\B}^{0}_{\infty,1}(\RR^n,\X)\hookrightarrow \C^{0}_{ub}(\RR^n,\X) \hookrightarrow \C^{0}_{ub}(\RR,\C^{0}_{ub}(\RR^{n-1},\X)).
\end{align*}

\textbf{Step 2:} We prove \textit{(vii)}. By real interpolation,  assuming $p<\infty$, one obtains for $0<s<1$ 
\begin{align*}
    \dot{\B}^{s}_{p,q}(\RR^n,\X)&= \Big({\L}^{p}(\RR^n,\X),\dot{\W}^{1,p}(\RR^n,\X)\Big)_{s,q}\\
    &\hookrightarrow \bigcap_{1\leqslant k\leqslant n} \Big({\L}^{p}(\RR_{x_k},{\L}^{p}(\RR^{n-1},\X)),\dot{\W}^{1,p}(\RR_{x_k},{\L}^{p}(\RR^{n-1},\X))\Big)_{s,q}\\
    &= \bigcap_{1\leqslant k\leqslant n}\dot{\B}^{s}_{p,q}(\RR_{x_k},{\L}^{p}(\RR^{n-1},\X)).
\end{align*}
So first, we just have to deal with the one dimensional case. But in this case, one can mimic the proof of \cite[Theorem~2.36]{bookBahouriCheminDanchin}, in order to obtain for all $u\in\dot{\B}^{s}_{p,q}(\RR^n,\X)$,
\begin{align}\label{eq:ProofContinuityFiniteDiffHomBesov}
    \left\lVert t\mapsto t^{-s}\lVert\tau_{t}u-u\rVert_{\L^p(\RR,\X)} \right\rVert_{\L^q_\ast(\RR_+)} \lesssim_{p,s,n}^{\X} \lVert u\rVert_{\dot{\B}^{s}_{p,q}(\RR,\X)}.
\end{align}

Following the ideas \cite[Theorem~3.10~\&~Section~5.2]{bookLunardiInterpTheory}, for all $u\in{\B}^{s}_{p,q}(\RR^n,\X)$,
\begin{align*}
    \lVert u\rVert_{{\B}^{s}_{p,q}(\RR^n,\X)} \lesssim_{p,s,n}^{\X} \lVert u\rVert_{\L^p(\RR^n,\X)}+ \sum_{k=1}^n\left\lVert t\mapsto t^{-s}\lVert\tau_{t}u-u\rVert_{\L^p(\RR^n,\X)} \right\rVert_{\L^q_\ast(\RR_+)} 
\end{align*}
By the equivalence of norms \eqref{eq:EquivNormHomGenInhomogenBesov}, one obtains by a dilation argument, for all $u\in{\B}^{s}_{p,q}(\RR^n,\X)$,
\begin{align*}
    \lVert u\rVert_{\dot{\B}^{s}_{p,q}(\RR^n,\X)} \sim_{p,s,n}^{\X} \sum_{k=1}^n \left\lVert t\mapsto t^{-s}\lVert\tau_{t}^k u-u\rVert_{\L^p(\RR^n,\X)} \right\rVert_{\L^q_\ast(\RR_+)}.
\end{align*}
If $q<\infty$, the equivalence follows from the continuity estimate \eqref{eq:ProofContinuityFiniteDiffHomBesov} valid for all $u\in\dot{\B}^{s}_{p,q}(\RR^n,\X)$, and the density of ${\B}^{s}_{p,q}(\RR^n,\X)$ in $\dot{\B}^{s}_{p,q}(\RR^n,\X)$. Now, if $q=\infty$, let $u\in\dot{\B}^{s}_{p,\infty}(\RR^n,\X)$, and consider for all $N\in\mathbb{N}$
\begin{align*}
    u_{N} := \sum_{j\geqslant -N} \dot{\Delta}_ju.
\end{align*}
One has $u_N\in {\B}^{s}_{p,\infty}(\RR^n,\X)$, moreover for almost all $t>0$,
\begin{align*}
    \lVert\tau_{t}^ku_N-u_N\rVert_{\L^p(\RR^n,\X)}\leqslant \lVert\tau_{t}^ku-u\rVert_{\L^p(\RR^n,\X)}\text{, and } 2^{js}\lVert\Delta_j u_N\rVert_{\L^p(\RR^n,\X)}\xrightarrow[N\rightarrow \infty]{} 2^{js}\lVert\Delta_j u\rVert_{\L^p(\RR^n,\X)}.
\end{align*}
By the Fatou lemma for scalar-valued bounded sequences
\begin{align*}
    \lVert u\rVert_{\dot{\B}^{s}_{p,\infty}(\RR^n,\X)}\leqslant \liminf_{N\rightarrow\infty} \lVert u_N\rVert_{\dot{\B}^{s}_{p,\infty}(\RR^n,\X)} &\sim_{p,s,n}^{\X} \sum_{k=1}^n \left\lVert t\mapsto t^{-s}\lVert\tau_{t}^k u_N-u_N\rVert_{\L^p(\RR^n,\X)} \right\rVert_{\L^\infty(\RR_+)}\\
    & \lesssim_{p,s,n}^{\X} \sum_{k=1}^n \left\lVert t\mapsto t^{-s}\lVert\tau_{t}^k u-u\rVert_{\L^p(\RR^n,\X)} \right\rVert_{\L^\infty(\RR_+)}.
\end{align*}
For $p=\infty$, one can reproduce \cite[Theorem~3.10~\&~Section~5.2]{bookLunardiInterpTheory}, giving for all $u\in{\B}^{s}_{\infty,q}(\RR^n,\X)$,
\begin{align*}
    \lVert u\rVert_{{\B}^{s}_{\infty,q}(\RR^n,\X)} \sim_{p,s,n}^{\X} \lVert u\rVert_{\L^\infty(\RR,\X)}+ \sum_{k=1}^n\left\lVert t\mapsto t^{-s}\lVert\tau_{t}u-u\rVert_{\L^\infty(\RR^n,\X)} \right\rVert_{\L^q_\ast(\RR_+)}, 
\end{align*}
and the result holds by a dilation argument, since $\dot{\B}^{s}_{\infty,q}(\RR^n,\X)\subset {\B}^{s}_{\infty,q}(\RR^n,\X)$ (as a set).
\end{proof}

\begin{theorem}\label{thm:fundamentalmultBesovSpacesX} Let $p,q\in[1,\infty]$, $s\in(-1+{\sfrac{1}{p}},{\sfrac{1}{p}})$, $n\in\mathbb{N}^\ast$, and $\X$ be a Banach space. Then for any $u\in\dot{\B}^{s}_{p,q}(\RR^n,\X)$,
\begin{align*}
    \mathbbm{1}_{\RR^n_+}u\in \dot{\B}^{s}_{p,q}(\RR^n,\X)
\end{align*}
with the estimate
\begin{align*}
    \lVert \mathbbm{1}_{\RR^n_+}u\rVert_{\dot{\B}^{s}_{p,q}(\RR^n,\X)} \lesssim_{s,p}^{n,\X} \lVert u\rVert_{\dot{\B}^{s}_{p,q}(\RR^n,\X)}.
\end{align*}
Furthermore, the result remains valid for inhomogeneous Besov spaces.
\end{theorem}

\begin{proof}We define $\chi_0(t) := \mathbbm{1}_{(0,\infty)}(t)$, and identify $\chi_0(x',x_n)= \chi_0(x_n)=\mathbbm{1}_{(0,\infty)}(x_n)= \mathbbm{1}_{\RR^n_+}(x',x_n)$ for $x=(x',x_n)\in\RR^{n-1}\times\RR=\RR^n$. The proof is is an heavily inspired and somewhat simplified version from the scalar-valued case exposed in \cite[Appendix,~Lemma~12]{DanchinMucha2009}.

\medbreak

\textbf{Step 1:} Let $p\in [1,\infty)$, $s\in(0,{\sfrac{1}{p}})$, $u\in\dot{\B}^{s}_{p,1}(\RR^n,\X)$. We will show that 
\begin{align*}
    \lVert \chi_0 u\rVert_{\dot{\B}^{s}_{p,\infty}(\RR^n,\X)} &\lesssim_{s,p}^{n,\X} \lVert  u\rVert_{\dot{\B}^{s}_{p,1}(\RR^n,\X)}.
\end{align*}

\textbf{Step 1.1:} By point \textit{(vii)} of Proposition~\ref{prop:BanachValuedHomBesovSpaces}, by independence of variables and Fubini's theorem
\begin{align*}
    \lVert \chi_0 u\rVert_{\dot{\B}^{s}_{p,\infty}(\RR^n,\X)} &\lesssim_{s,p}^{n,\X} \sum_{k=1}^n \left\lVert t\mapsto t^{-s}\lVert\tau^{k}_{t}(\chi_0 u)-(\chi_0 u)\rVert_{\L^p(\RR^n,\X)} \right\rVert_{\L^\infty(\RR_+)}\\
    &\lesssim_{s,p}^{n,\X} \sum_{k=1}^{n-1} \left\lVert t\mapsto t^{-s}\lVert \chi_0[\tau^{k}_{t} u-u]\rVert_{\L^p(\RR^n,\X)} \right\rVert_{\L^\infty(\RR_+)}\\
    &\qquad + \left\lVert t\mapsto t^{-s}\lVert\tau^{n}_{t}(\chi_0 u)-(\chi_0 u)\rVert_{\L^p(\RR^n,\X)} \right\rVert_{\L^\infty(\RR_+)}.
\end{align*}
Since $\lvert\chi_0\rvert \leqslant 1$, by point \textit{(vii)} of Proposition~\ref{prop:BanachValuedHomBesovSpaces}, we obtain
\begin{align*}
    \lVert \chi_0 u\rVert_{\dot{\B}^{s}_{p,\infty}(\RR^n,\X)} &\lesssim_{s,p}^{n,\X} \lVert u\rVert_{\dot{\B}^{s}_{p,\infty}(\RR^n,\X)} + \left\lVert t\mapsto t^{-s}\lVert (\tau^{n}_{t}\chi_0)(\tau^{n}_{t} u-u)\rVert_{\L^p(\RR^n,\X)} \right\rVert_{\L^\infty(\RR_+)} \\ &\qquad + \left\lVert t\mapsto t^{-s}\lVert(\tau^{n}_{t}\chi_0-\chi_0 )u\rVert_{\L^p(\RR^n,\X)} \right\rVert_{\L^\infty(\RR_+)}
\end{align*}
where we applied the triangle inequality to the decomposition $\tau^{n}_{t}(\chi_0 u)-(\chi_0 u) = (\tau^{n}_{t}\chi_0)(\tau^{n}_{t} u-u) + (\tau^{n}_{t}\chi_0-\chi_0 )u$. Again, since $\rvert\tau^{n}_{t}\chi_0\rvert \leqslant 1$, by point \textit{(vii)} of Proposition~\ref{prop:BanachValuedHomBesovSpaces}, we obtain
\begin{align*}
    \lVert \chi_0 u\rVert_{\dot{\B}^{s}_{p,\infty}(\RR^n,\X)} &\lesssim_{s,p}^{n,\X} \lVert u\rVert_{\dot{\B}^{s}_{p,\infty}(\RR^n,\X)} + \left\lVert t\mapsto t^{-s}\lVert(\tau^{n}_{t}\chi_0-\chi_0 )u\rVert_{\L^p(\RR^n,\X)} \right\rVert_{\L^\infty(\RR_+)}
\end{align*}

\textbf{Step 1.2:} It remains to show that
\begin{align*}
    t^{-s}\lVert(\tau^{n}_{t}\chi_0-\chi_0 )u\rVert_{\L^p(\RR^n,\X)} \lesssim_{s,p}^{n,\X} \lVert  u\rVert_{\dot{\B}^{s}_{p,1}(\RR^n,\X)},\qquad t>0.
\end{align*}
By Fubini's theorem and Hölder's inequality, for $t>0$,
\begin{align*}
    \lVert(\tau^{n}_{t}\chi_0-\chi_0 )u\rVert_{\L^p(\RR^n,\X)} &= \lVert x_n\mapsto (\tau^{n}_{t}\chi_0-\chi_0 )(x_n)u(x_n,\cdot)\rVert_{\L^p(\RR,\L^p(\RR^{n-1},\X))}\\
    & \leqslant \lVert (\tau^{n}_{t}\chi_0-\chi_0 )\rVert_{\L^{\tilde{r}}(\RR)} \lVert  u\rVert_{\L^{r}(\RR,\L^p(\RR^{n-1},\X))}
\end{align*}
with $s=\frac{1}{p}-\frac{1}{r}=\frac{1}{\tilde{r}}$, so that by point \textit{(iv)} of Proposition~\ref{prop:BanachValuedHomBesovSpaces}, one obtains
\begin{align*}
    \lVert(\tau^{n}_{t}\chi_0-\chi_0 )u\rVert_{\L^p(\RR^n,\X)} \lesssim_{s,p}^{n,\X} \lVert \tau_{t}\chi_0-\chi_0 \rVert_{\L^{\Tilde{r}}(\RR)} \lVert  u\rVert_{\dot{\B}^{s}_{p,1}(\RR^n,\X)}.
\end{align*}
Recalling that $\supp (\tau_{t}\chi_0-\chi_0) \subset [-t,0]$ and $\lvert\tau_{t}\chi_0-\chi_0\rvert\leqslant1$, and that $\frac{1}{\Tilde{r}}=s$, it holds
\begin{align*}
    \lVert \tau_{t}\chi_0-\chi_0 \rVert_{\L^{\Tilde{r}}(\RR)} = \lVert \tau_{t}\chi_0-\chi_0 \rVert_{\L^{\Tilde{r}}([-t,0])}= t^s 
\end{align*}
so we can finish the Step 1, with the desired inequality
\begin{align*}
    t^{-s}\lVert (\tau^{n}_{t}\chi_0-\chi_0 )u\rVert_{\L^p(\RR^n,\X)} &\lesssim_{s,p}^{n,\X} \lVert  u\rVert_{\dot{\B}^{s}_{p,1}(\RR^n,\X)},\qquad t>0.
\end{align*}

\textbf{Step 2:} We will show that 
\begin{align*}
    \lVert \chi_0 u\rVert_{\dot{\B}^{s}_{p,q}(\RR^n,\X)} &\lesssim_{s,p}^{n,\X} \lVert  u\rVert_{\dot{\B}^{s}_{p,q}(\RR^n,\X)}.
\end{align*}

We consider $\eta_0,\eta_1 >0$ such that $0<s-\eta_0 <s<s+\eta_1 < {\sfrac{1}{p}}$, since by point \textit{(iii)} of Proposition~\ref{prop:BanachValuedHomBesovSpaces}, we have
\begin{align*}
        \Big(\dot{\B}^{s-\eta_0}_{p,1}(\RR^n,\X),\dot{\B}^{s+\eta_1}_{p,1}(\RR^n,\X)\Big)_{\frac{\eta_0}{\eta_0+\eta_1},q}=\Big(\dot{\B}^{s-\eta_0}_{p,\infty}(\RR^n,\X),\dot{\B}^{s+\eta_1}_{p,\infty}(\RR^n,\X)\Big)_{\frac{\eta_0}{\eta_0+\eta_1},q}=\dot{\B}^{s}_{p,q}(\RR^n,\X)
    \end{align*}

We obtain indeed by interpolation, for all $p\in[1,\infty)$, $q\in[1,\infty]$, $s\in(0,{\sfrac{1}{p}})$, all $u\in \dot{\B}^{s}_{p,q}(\RR^n,\X)$,
\begin{align*}
    \lVert \chi_0 u\rVert_{\dot{\B}^{s}_{p,q}(\RR^n,\X)} &\lesssim_{s,p}^{n,\X} \lVert  u\rVert_{\dot{\B}^{s}_{p,q}(\RR^n,\X)}.
\end{align*}

\textbf{Step 3:}  The case $s\in(-1+{\sfrac{1}{p}},0]$. We apply the result to $\X^\ast$, for $s\in(0,{\sfrac{1}{p}})$, 
\begin{align*}
    v\mapsto \mathbbm{1}_{\RR^n_+}v
\end{align*}
is bounded on $\dot{\B}^{s}_{p,q}(\RR^n,\X^{\ast})$, $p,q\in(1,\infty)$. By duality, the map
\begin{align*}
    u\mapsto \mathbbm{1}_{\RR^n_+}u
\end{align*}
is bounded on $\dot{\B}^{-s}_{p',q'}(\RR^n,\X^{\ast\ast})$, $p',q'\in(1,\infty)$, $-s\in(-1+{\sfrac{1}{p'}},0)$. Duality with the case $p=1$ gives boundedness on the spaces $\dot{\B}^{-s}_{\infty,q'}(\RR^n,\X^{\ast\ast})$, $-s\in(-1,0)$, $q'\in(1,\infty)$.

Now, for $u\in \C^{\infty}_{ub} \cap \dot{\B}^{-s}_{p',q'}(\RR^n,\X)$, for any $x\in\RR^n$, one has $\mathbbm{1}_{\RR^n_+}u(x)\in\X$, so that due to the isometry $\X\hookrightarrow\X^{\ast\ast}$, it holds
\begin{align*}
    \lVert \mathbbm{1}_{\RR^n_+}u \rVert_{\dot{\B}^{-s}_{p',q'}(\RR^n,\X)}=\lVert \mathbbm{1}_{\RR^n_+}u \rVert_{\dot{\B}^{-s}_{p',q'}(\RR^n,\X^{\ast\ast})} \lesssim_{p,s,n,\X}\lVert u \rVert_{\dot{\B}^{-s}_{p',q'}(\RR^n,\X^{\ast\ast})} = \lVert u \rVert_{\dot{\B}^{-s}_{p',q'}(\RR^n,\X)}.
\end{align*}
Since $1<q'<\infty$, by Proposition~\ref{prop:BanachValuedHomBesovSpaces}, $\C^{\infty}_{ub}\cap \dot{\B}^{-s}_{p',q'}(\RR^n,\X)$ is dense in $\dot{\B}^{-s}_{p',q'}(\RR^n,\X)$ whenever $1<p'\leqslant \infty$, so the map extends by density.

Finally, we did obtain that the map
\begin{align*}
    u\mapsto \mathbbm{1}_{\RR^n_+}u
\end{align*}
is bounded on $\dot{\B}^{s}_{p,q}(\RR^n,\X)$, $s\in(-1+{\sfrac{1}{p}},0)\cup(0,{\sfrac{1}{p}})$, $p\in[1,\infty],q\in(1,\infty)$. The remaining cases are obtained by real interpolation.
\end{proof}

For $\Y(\RR^n,\X)\subset\mathcal{S}'(\RR^n,\X)$ a normed function space, we define 
\begin{itemize}
    \item $\Y(\RR^n_+,\X):=\Y(\RR^n,\X)_{|_{\RR^n_+}}$ equipped with the induced quotient norm;
    \item $\Y_0(\RR^n_+,\X):=\{\,u\in\Y(\RR^n,\X)\,:\,\supp u \subset\,\overline{\RR^n_+}\}$ endowed with the $\Y(\RR^n,\X)$-norm.
\end{itemize}

One can follow the Besov spaces part of \cite[Section~3]{Gaudin2023Lip}, and reproduce it to obtain a full vector valued construction of such homogeneous function spaces on half-spaces\footnote{Unfortunately, the corresponding part of \cite[Section~3]{Gaudin2023Lip} for homogeneous Bessel spaces $\dot{\H}^{s,p}$  is \textit{a priori} restricted to the Hilbert-valued case due to the use of square functionals in several proofs. A proper study in the case of UMD Banach-valued case should be possible, but the arguments are no-longer valid}. However, we just choose to highlight here a very small selection of necessary results for our analysis.

\begin{corollary}Let $p,q\in[1,\infty]$ and $s\in(-1+\sfrac{1}{p},{\sfrac{1}{p}})$. One has
\begin{align*}
    \dot{\B}^{s}_{p,q}(\RR^n_+,\X)=\dot{\B}^{s}_{p,q,0}(\RR^n_+,\X),
\end{align*}
and $\Ccinfty(\RR^n_+,\X)$ is dense in $\dot{\B}^{s}_{p,q}(\RR^n_+,\X)$ in
and $\dot{\BesSmo}^{s}_{p,\infty}(\RR^n_+,\X)$,  in $\dot{\B}^{s,0}_{\infty,q}(\RR^n_+,\X)$ and $\dot{\BesSmo}^{s,0}_{\infty,\infty}(\RR^n_+,\X)$, whenever $p,q<\infty$.
\end{corollary}

\begin{proposition}\label{prop:AppendixFiniteDifVectorValuedBesovRN+}Let $p,q\in[1,\infty]$ and $s\in(0,1)$. For all $u\in\L^1_{\text{loc}}(\overline{\RR^n_+},\X)$, one has
\begin{align}\label{eq:estvectorValuedExt}
        \sum_{k=1}^n \Big\lVert t\mapsto t^{-s}\lVert\tau^{k}_{t}[\E_{\mathcal{N}}u]-[\E_{\mathcal{N}}u]\rVert_{\L^p(\RR^n,\X)} \Big\rVert_{\L^q_\ast(\RR_+)} \sim_{p,s,n} \sum_{k=1}^n \Big\lVert t\mapsto t^{-s}\lVert\tau^{k}_{t}u-u\rVert_{\L^p(\RR^n_+,\X)} \Big\rVert_{\L^q_\ast(\RR_+)}. 
\end{align}
In particular, for all $u\in\dot{\B}^{s}_{p,q}(\RR^n_+,\X)$, one has $\E_{\mathcal{N}}u\in\dot{\B}^{s}_{p,q}(\RR^n,\X)$ and
\begin{align*}
    \lVert u\rVert_{\dot{\B}^{s}_{p,q}(\RR^n_+,\X)} \sim_{p,s,n}^\X\lVert \E_{\mathcal{N}}u\rVert_{\dot{\B}^{s}_{p,q}(\RR^n,\X)} \sim_{p,s,n}^\X \sum_{k=1}^n \left\lVert t\mapsto t^{-s}\lVert\tau^{k}_{t}u-u\rVert_{\L^p(\RR^n_+,\X)} \right\rVert_{\L^q_\ast(\RR_+)}. 
\end{align*}
\end{proposition}

\begin{proof} The estimate \eqref{eq:estvectorValuedExt} is elementary. $\E_\mathcal{N}$ is known to map $\L^p(\RR^n_+,\X)$ to $\L^p(\RR^n,\X)$, and $\W^{1,p}(\RR^n_+,\X)$ to $\W^{1,p}(\RR^n,\X)$, for $p\in[1,\infty]$. Hence, by real interpolation, it maps ${\B}^{s}_{p,q}(\RR^n_+,\X)$ to ${\B}^{s}_{p,q}(\RR^n,\X)$ for all $0<s<1$, $p,q\in[1,\infty]$. It also maps $\C^0_0(\overline{\RR^n_+},\X)$ to $\C^0_0(\RR^n,\X)$.

For $p\in[1,\infty)$, for all $u\in\dot{\B}^{s}_{p,q}(\RR^n_+,\X)\subset {\B}^{s}_{p,q}(\RR^n_+,\X) + \C^0_0(\overline{\RR^n_+},\X)$, one has $$\E_{\mathcal{N}}u\in\dot{\B}^{s}_{p,q}(\RR^n_+,\X)\subset {\B}^{s}_{p,q}(\RR^n,\X) + \C^0_0(\RR^n,\X)\subset \S'_h(\RR^n,\X).$$ Consequently, one obtains $\E_{\mathcal{N}}u\in\dot{\B}^{s}_{p,q}(\RR^n,\X)$, and then the result follows from the definition of function spaces by restriction. 

Now if $p=\infty$,  the difficulty is to show that  $\E_{\mathcal{N}}U\in\mathcal{S}'_h(\RR^n,\X)$ for $U\in\dot{\B}^{s}_{\infty,q}(\RR^n,\X)$ but one can proceed as in Lemma~\ref{lem:ExtDirNeuRn+}. 
\end{proof}

We obtain then the standard corollary:

\begin{corollary}\label{cor:InterpVectValuedHomBesov} Let $p,q,q_0,q_1\in[1,\infty]$, $s,s_0,s_1\in(-1+\sfrac{1}{p},1+\sfrac{1}{p})$, such that $s_0<s<s_1$. Then for $s=(1-\theta)s_0+\theta s_1$, $\theta\in(0,1)$, it holds that
\begin{align*}
    (\dot{\B}^{s_0}_{p,q_0}(\RR^n_+,\X),\dot{\B}^{s_1}_{p,q_1}(\RR^n_+,\X))_{\theta,q} =\dot{\B}^{s}_{p,q}(\RR^n_+,\X),
\end{align*}
with equivalence of norms.
\end{corollary}

\begin{proof} We just provide quickly the argument. $\E_{\mathcal{N}}$ is known to be bounded for $s,s_0,s_1\in(-1+\sfrac{1}{p},1)$, taking the derivatives with $\partial_{x_n}\E_{\mathcal{N}}=\E_{\mathcal{D}}\partial_{x_n}$ one obtains boundedness  of $\E_{\mathcal{N}}$  for $s,s_0,s_1\in(-1+\sfrac{1}{p},1+\sfrac{1}{p})$.
Then a retraction and co-retraction argument yields the result.
\end{proof}

\section{Some semigroup estimates in vector-valued Besov spaces}\label{App:BesovSemigrpEst}

Following \cite[Chapter~2]{DanchinHieberMuchaTolk2020}, for $(\D(A),A)$ an injective $\omega$-sectorial operator-- here $\omega\in[0,\frac{\pi}{2})$-- on a Banach space $\X$ satisfying Assumptions \ref{asmpt:homogeneousdomaindef} and \ref{asmpt:homogeneousdomainintersect}, we recall that for any $\theta\in(0,1)$, $q\in[1,\infty]$,
\begin{align*}
    (\X,\D(\mathring{A}))_{\theta,q} =\mathring{\mathcal{D}}_{A}(\theta,q),
\end{align*}
with equivalence norms, having for all $u\in\X+\D(\mathring{A})$,
\begin{align*}
    \lVert u \rVert_{\mathring{\mathcal{D}}_{A}(\theta,q)}^q = \int_{0}^\infty\lVert t^{1-\theta}\mathring{A}e^{-tA}u \rVert_{\X}^q \frac{\d t}{t},
\end{align*}
with the usual change when $q=\infty$.

However, for technical purposes, Danchin \textbf{et al.} did also introduced an \textit{extended scale} allowing to define these spaces for $\theta$ that goes beyond $1$, proving that for all $0<\theta<1$, all $u\in\X+\D(\mathring{A})$
\begin{align*}
    \int_{0}^\infty\lVert t^{1-\theta}\mathring{A}e^{-tA}u \rVert_{\X}^q \frac{\d t}{t} \sim_{\theta,q} \int_{0}^\infty\lVert t^{2-\theta}A\mathring{A}e^{-tA}u \rVert_{\X}^q \frac{\d t}{t}.
\end{align*}
Notice that the right hand-side above still makes sense for $\theta\in[1,2)$. It turns out that we can define for any $0<\theta<2$, $q\in[1,\infty]$
\begin{align*}
    \tilde{\mathcal{D}}_{A}(\theta,q) = \{\,u\in\X+\D(\mathring{A}^2)\,:\,t\mapsto \lVert t^{2-\theta}A\mathring{A}e^{-tA}u \rVert_{\X}\in\L^q_{\ast}(\RR_+)\,\}
\end{align*}
where
\begin{align*}
    &\D(\mathring{A}^2)=\{\,u\in\D(\mathring{A})\,:\,\mathring{A}u\in\D(\mathring{A})\,\} = \{\,u\in\D(\mathring{A})\,:\,\mathring{A}u\in\D({A})\,\},\\
    \text{ with norm }\quad &\lVert u\rVert_{\D(\mathring{A}^2)}:= \lVert A\mathring{A}u\rVert_{\X}.
\end{align*}
It has been proved in \cite[Chapter~2]{DanchinHieberMuchaTolk2020}, that the following holds algebraically and with equivalence of norms
\begin{align*}
     &(\X,\D(\mathring{A}^2))_{\frac{\theta}{2},q} =\mathring{\mathcal{D}}_{A}(\theta,q),\qquad \theta\in(0,1),\\
     \text{ and }& (\X,\D(\mathring{A}^2))_{\frac{\theta}{2},q} =\tilde{\mathcal{D}}_{A}(\theta,q),\qquad \theta\in(0,2).
\end{align*}
Therefore, for $1\leqslant\theta<2$, we set
\begin{align*}
    \mathring{\mathcal{D}}_{A}(\theta,q):=\tilde{\mathcal{D}}_{A}(\theta,q).
\end{align*}

Before we state and prove the next result concerning those spaces, we provide a short remark about $\D(\mathring{A}^2)$. The definition and the use of $\D(\mathring{A}^2)$ and the related extended scale is really here mostly for technical purposes, and considered as such for concrete spaces is a bit unnatural, while sufficient for most purposes. The idea of such an extended scale for us is to have spaces sufficiently ``far away'' in the scale to interpolate with, and then to recover full scale by interpolation when $\theta\in(0,1)$, $q\in[1,\infty]$. To illustrate the weirdness for this definition of $\D(\mathring{A}^2)$, we consider the following setting
\begin{itemize}
    \item $\X=\L^p(\RR^n)$, $1<p<\infty$;
    \item $(\D(A),A)=(\H^{2,p}(\RR^n),-\Delta)$.
\end{itemize}
Therefore, one obtains algebraically and with equivalence of norms
\begin{align*}
    \D(\mathring{A}) =\dot{\H}^{2,p}(\RR^n) 
\end{align*}
but
\begin{align*}
    \D(\mathring{A}^2)=\dot{\H}^{2,p}\cap\dot{\H}^{4,p}(\RR^n),\text{ with }\lVert u \rVert_{\D(\mathring{A}^2)} = \lVert u \rVert_{\dot{\H}^{4,p}(\RR^n)}.
\end{align*}
Hence, defined as such $\D(\mathring{A}^2)$ is not the homogeneous domain of $A^2$ on $\X$ !

When $A$ is invertible recall that $\D(\mathring{A})=\D({A})$ so that $\D(\mathring{A}^2)=\D({A}^2)$.

The goal of the next (not so difficult to prove) Lemma is about to provide a full reiteration results for real interpolation for the family of such spaces $(\mathring{\mathcal{D}}_{A}(\theta,q))_{\theta,q}$. Since such spaces might not be complete\footnote{and recall that one \textbf{SHOULD NOT} consider their completion to preserve compatibility with the framework of "concrete" homogeneous function spaces.}, standard reiteration theorems for real interpolation, such as in \cite[Chapter~3]{BerghLofstrom1976}, do not apply here.

\begin{lemma}\label{lem:ReitThmSmigrpDaPrato}Let $\omega\in[0,\frac{\pi}{2})$ and $(\D(A),A)$ an injective $\omega$-sectorial operator on a Banach space $\X$ such that it satisfies Assumptions \ref{asmpt:homogeneousdomaindef} and \ref{asmpt:homogeneousdomainintersect}. Then the following reiteration results hold for the extended scale
\begin{enumerate}
    \item For all $0<\theta<1$, all $0<s<2$, all $\kappa,q\in[1,\infty]$, it holds algebraically and with equivalence of norms that
    \begin{align*}
        (\mathring{\mathcal{D}}_{A}(s,\kappa),\D(\mathring{A}^2))_{\theta,q} =\mathring{\mathcal{D}}_{A}((1-\theta)s+2\theta,q);
    \end{align*}
    \item For all $0<\theta<1$,  all $q\in[1,\infty]$, it holds algebraically and with equivalence of norms that
    \begin{align*}
        (\D(\mathring{A}),\D(\mathring{A}^2))_{\theta,q} =\mathring{\mathcal{D}}_{A}(1+\theta,q) ;
    \end{align*}
    \item For all $0<\theta<1$, all $0<s<2$, $s\neq1$, all $\kappa,q\in[1,\infty]$, it holds algebraically and with equivalence of norms that
    \begin{align*}
        (\mathring{\mathcal{D}}_{A}(s,\kappa),\D(\mathring{A}))_{\theta,q} =\mathring{\mathcal{D}}_{A}((1-\theta)s+\theta,q);
    \end{align*}
    \item For all $0<\theta<1$, all $0<s_0<s_1<2$, for $s=(1-\theta)s_0+\theta s_1$, all $q_0,q_1,q\in[1,\infty]$, it holds algebraically and with equivalence of norms that
    \begin{align*}
        (\mathring{\mathcal{D}}_{A}(s_0,q_0),\mathring{\mathcal{D}}_{A}(s_1,q_1))_{\theta,q} =\mathring{\mathcal{D}}_{A}(s,q).
    \end{align*}
\end{enumerate}
\end{lemma}

\begin{proof} For point $\textit{(i)}$ it suffices to apply the extremal reiteration property \cite[Lemma~C.1]{Gaudin2023Lip}, so that
\begin{align*}
    (\mathring{\mathcal{D}}_{A}(s,\kappa),\D(\mathring{A}^2))_{\theta,q} &= ((\X,\D(\mathring{A}^2))_{\frac{s}{2},\kappa},\D(\mathring{A}^2))_{\theta,q}\\
    &= (\X,\D(\mathring{A}^2))_{(1-\theta)\frac{s}{2}+\theta,q}\\
    &= \mathring{\mathcal{D}}_{A}((1-\theta)s+2\theta,q).
\end{align*}
Point $\textit{(ii)}$ is a direct corollary from point $\textit{(i)}$, indeed
\begin{align*}
    \mathring{\mathcal{D}}_{A}(1+\theta,q) &= (\mathring{\mathcal{D}}_{A}(1,1),\D(\mathring{A}^2))_{\theta,q}\\
    &\hookrightarrow (\D(\mathring{A}),\D(\mathring{A}^2))_{\theta,q}\\
    &\hookrightarrow (\mathring{\mathcal{D}}_{A}(1,\infty),\D(\mathring{A}^2))_{\theta,q}\\
    &=\mathring{\mathcal{D}}_{A}(1+\theta,q).
\end{align*}
Hence, we did actually obtain
\begin{align*}
        (\D(\mathring{A}),\D(\mathring{A}^2))_{\theta,q} =\mathring{\mathcal{D}}_{A}(1+\theta,q).
\end{align*}
For point \textit{(iii)}, we distinguish two cases $0<s<1$ and $1<s<2$. If $0<s<1$, this follows directly from the extremal reiteration property \cite[Lemma~C.1]{Gaudin2023Lip}, with
\begin{align*}
    (\mathring{\mathcal{D}}_{A}(s,\kappa),\D(\mathring{A}))_{\theta,q} &= ((\X,\D(\mathring{A}))_{{s},\kappa},\D(\mathring{A}))_{\theta,q}\\
    &= (\X,\D(\mathring{A}))_{(1-\theta){s}+\theta,q}\\
    &= \mathring{\mathcal{D}}_{A}((1-\theta)s+\theta,q).
\end{align*}
If $1<s<2$, this is again an application of the extremal reiteration property thanks to point~\textit{(ii)}.

For the last full reiteration theorem, we use point \textit{(i)} and the extremal reiteration property. We set $0<\eta<1$ such that $(1-\eta)s_0+ 2\eta=s_1$, thanks to point \textit{(i)} it holds
\begin{align*}
    (\mathring{\mathcal{D}}_{A}(s_0,q_0),\mathring{\mathcal{D}}_{A}(s_1,q_1))_{\theta,q} &= \Big(\mathring{\mathcal{D}}_{A}(s_0,q_0),(\mathring{\mathcal{D}}_{A}(s_0,q_0),\D(\mathring{A}^2))_{\eta,q_1}\Big)_{\theta,q}\\
    &=(\mathring{\mathcal{D}}_{A}(s_0,q_0),\D(\mathring{A}^2))_{\theta,q}\\
    &=\mathring{\mathcal{D}}_{A}(s,q).
\end{align*}
This finishes the proof.
\end{proof}

\begin{lemma}\label{lem:BesovintimeSemigroup}Let $\omega\in[0,\frac{\pi}{2})$ and $(\D(A),A)$ an injective $\omega$-sectorial operator on a Banach space $\X$ such that it satisfies Assumptions \ref{asmpt:homogeneousdomaindef} and \ref{asmpt:homogeneousdomainintersect}. Let $q,\kappa\in[1,\infty]$ and $-1+{\sfrac{1}{q}}<\alpha<1+{\sfrac{1}{q}}$. For simplicity, we set $\alpha_q:=1+\alpha-\sfrac{1}{q}\in(0,2)$. For all $u_0\in\mathring{\mathcal{D}}_{A}(\alpha_q,\kappa)$,
\begin{align*}
    \lVert t\mapsto \mathring{A}e^{-tA} u_0\rVert_{\dot{\B}^{\alpha}_{q,\kappa}(\RR_+,\X)} \lesssim_{q,\alpha,A,\X} \lVert u_0 \rVert_{\mathring{\mathcal{D}}_{A}(\alpha_q,\kappa)}.
\end{align*}
Furthermore, if $u_0\in{\mathcal{D}}_{A}(\alpha_q,\kappa)$, then
\begin{itemize}
    \item we can drop Assumptions \ref{asmpt:homogeneousdomaindef} and \ref{asmpt:homogeneousdomainintersect};
    \item the result still holds replacing respectively the homogeneous norms $\dot{\B}^{\alpha}_{q,\kappa}(\RR_+,\X)$ and $\mathring{\mathcal{D}}_{A}(\alpha_q,\kappa)$, by the inhomogeneous norms ${\B}^{\alpha}_{q,\kappa}(\RR_+,\X)$ and $\mathcal{D}_{A}(\alpha_q,\kappa)$.
\end{itemize}
\end{lemma}

\begin{proof} We start assuming either $[\kappa=\infty$ and $q<\infty]$ or $[q=\infty$ and $\kappa<\infty]$, while $0<\alpha<1$.

Let $u_0\in\mathring{\mathcal{D}}_{A}(\alpha_q,\kappa)$, one has for $t,\uptau>0$,
\begin{align*}
    t^{-\alpha}[\mathring{A}e^{-(\uptau+t)A} u_0 - \mathring{A}e^{-\uptau A} u_0] = \frac{1}{t^\alpha}\int_{0}^t A\mathring{A}e^{-(\uptau+s) A} u_0 \,\d s.
\end{align*}

\textbf{Step 1:} The case $q=\infty$, $0<\alpha<1$. Let $u\in{\mathcal{D}}_{A}(1+\alpha,\kappa)$,
\begin{align*}
    \lVert t\mapsto \mathring{A}e^{-tA} u_0\rVert_{\dot{\B}^{\alpha}_{\infty,\kappa}(\RR_+,\X)} &\sim_{q,\alpha,A,\X} \left(\int_{0}^{\infty} \left( \sup_{\uptau>0} \left\lVert\frac{1}{t^{\alpha}}\int_{0}^t A\mathring{A}e^{-(\uptau+s) A} u_0 \,\d s\right\rVert_{\X} \right)^{\kappa}  \frac{\d t}{t}\right)^{1/\kappa}\\
    &\leqslant\left(\int_{0}^{\infty} \left( \sup_{\uptau>0} \frac{1}{t^{\alpha}}\int_{0}^t \lVert A\mathring{A}e^{-(\uptau+s) A} u_0 \rVert_{\X} \,\d s \right)^{\kappa}  \frac{\d t}{t}\right)^{1/\kappa}\\
    &\lesssim_{A,\X}\left(\int_{0}^{\infty} \left( \frac{1}{t^{\alpha}}\int_{0}^t \lVert s A\mathring{A}e^{-s A} u_0 \rVert_{\X} \,\frac{\d s}{s} \right)^{\kappa}  \frac{\d t}{t}\right)^{1/\kappa}.\\
    &\qquad\qquad\qquad\qquad\text{(by uniform boundedness of the semigroup $(e^{-\uptau A})_{\uptau>0}$.)}\\
    &\lesssim_{\kappa,\alpha,A,\X}\left( \int_{0}^{\infty} \lVert t^{2-(1+\alpha)}A\mathring{A}e^{-t A} u_0 \rVert_{\X}^\kappa \frac{\d t}{t}\right)^{1/\kappa}.
\end{align*}
The last line above is a consequence of Hardy's inequality \cite[Lemma~6.2.6]{bookHaase2006}.

\textbf{Step 2:} we assume here $q<\infty$, $\kappa=\infty$, $0<\alpha<1$. Since $q$ is finite, one has
\begin{align*}
    \lVert t\mapsto \mathring{A}e^{-tA} u_0\rVert_{\dot{\B}^{\alpha}_{q,\infty}(\RR_+,\X)} &\sim_{q,\alpha,A,\X} \sup_{t>0} \left( \int_{0}^{\infty} \left\lVert\frac{1}{t^{\alpha}}\int_{0}^t A\mathring{A}e^{-(\uptau+s) A} u_0 \,\d s\right\rVert_{\X}^q \d \uptau \right)^{1/q}  \\
    &\leqslant \sup_{t>0}\, t^{-\alpha} \left( \int_{0}^{\infty} \left(\int_{0}^{t} \rVert A\mathring{A}e^{-(s+\uptau) A} u_0\rVert_{\X} \,{\d s}\right)^q \d \uptau \right)^{1/q}  \\
    &\leqslant \sup_{t>0}\, t^{-\alpha} \int_{0}^{t} \left(\int_{s}^{\infty} \rVert A\mathring{A}e^{-\uptau A} u_0\rVert_{\X}^q \,{\d \uptau}\right)^\frac{1}{q} \d s   \\
    &\qquad\qquad\qquad\qquad\qquad\qquad\qquad\qquad\qquad\text{(by Minkowski's inequality.)}\\
    &\leqslant \Bigg(\sup_{t>0} \, t^{-\alpha} \int_{0}^{t} \left(\int_{s}^{\infty}\frac{\d\uptau}{\uptau^{q[2-(1+\alpha-\frac{1}{q})]}}\right)^\frac{1}{q} \,{\d s}\Bigg) \lVert u_0 \rVert_{\mathring{\mathcal{D}}_{A}(\alpha_q,\infty)}
\end{align*}

However, it turns out that for $t>0$, since $0<\alpha<1$,
\begin{align*}
     t^{-\alpha}\int_{0}^{t} \left(\int_{s}^{\infty}\frac{\d\uptau}{\uptau^{q[2-(1+\alpha-\frac{1}{q})]}}\right)^\frac{1}{q} \,{\d s} &=  t^{-\alpha}\int_{0}^{t} \left(\left[ \frac{\uptau^{-q(1-\alpha)}}{-q(1-\alpha)}\right]_s^\infty\right)^\frac{1}{q}\,{\d s}\\ &= \frac{1}{q^\frac{1}{q}}\frac{1}{(1-\alpha)^\frac{1}{q}} t^{-\alpha}\int_{0}^{t} s^{\alpha-1}\,{\d s} \\
     &= \frac{1}{q^\frac{1}{q}}\frac{1}{(1-\alpha)^\frac{1}{q}}\frac{1}{\alpha}.
\end{align*}

\textbf{Step 3:} Now for $q<\infty$, $\kappa=\infty$, we go beyond $1$ assuming $1<\alpha<1+1/q$. Since $q$ is finite, by Step 2, one has
\begin{align*}
    \lVert t\mapsto \mathring{A}e^{-tA} u_0\rVert_{\dot{\B}^{\alpha}_{q,\infty}(\RR_+,\X)} &\sim_{q,\alpha,A,\X} \lVert t\mapsto \partial_t \mathring{A}e^{-tA} u_0\rVert_{\dot{\B}^{\alpha-1}_{q,\infty}(\RR_+,\X)}  \\
    &\sim_{q,\alpha,A,\X} \lVert t\mapsto A e^{-tA} \mathring{A}u_0\rVert_{\dot{\B}^{\alpha-1}_{q,\infty}(\RR_+,\X)}  \\
    &\sim_{q,\alpha,A,\X} \lVert \mathring{A}u_0\rVert_{\mathring{\mathcal{D}}_{A}(\alpha-1/q,\infty)} \lesssim_{q,\alpha,A,\X} \lVert u_0\rVert_{\mathring{\mathcal{D}}_{A}(\alpha_q,\infty)}.
\end{align*}

\textbf{Step 4:} we assume now $q<\infty$, $\kappa=\infty$, $-1+1/q<\alpha<0$. Since $q$ is finite, one has
 one can use similar arguments, writing
\begin{align*}
    \mathring{A}e^{-tA} u_0 = \partial_t\int_{0}^t \mathring{A}e^{-sA} u_0\, \d s
\end{align*}
so that one should consider
\begin{align*}
    \frac{1}{t^{\alpha+1}}\left[ \int_{0}^{t+\uptau} \mathring{A}e^{-sA} u_0 \,\d s- \int_{0}^\uptau \mathring{A}e^{-sA} u_0 \,\d s\right] = \frac{1}{t^{\alpha+1}}\left[ \int_{\uptau}^{t+\uptau} \mathring{A}e^{-sA} u_0 \,\d s\right].
\end{align*}
Therefore,
\begin{align*}
    \lVert t\mapsto \mathring{A}e^{-tA} u_0\rVert_{\dot{\B}^{\alpha}_{q,\infty}(\RR_+,\X)} &= \left\lVert t\mapsto \partial_t\int_{0}^t\mathring{A}e^{-sA} u_0 \,\d s\right\rVert_{\dot{\B}^{\alpha}_{q,\infty}(\RR_+,\X)} \\
    &\lesssim_{\alpha,q} \left\lVert t\mapsto \int_{0}^t\mathring{A}e^{-sA} u_0 \,\d s\right\rVert_{\dot{\B}^{\alpha+1}_{q,\infty}(\RR_+,\X)} \\
    &\lesssim_{\alpha,q,\kappa}  \sup_{t>0}\,\int_{0}^{\infty} \left( \int_{0}^{\infty} \left\lVert\frac{1}{t^{\alpha+1}}\int_{\uptau}^{t+\uptau} \mathring{A}e^{-sA} u_0 \,\d s\right\rVert_{\X}^q \d \uptau \right)^{/q}  
\end{align*}
So that one can apply the exact same arguments from Step 2 in order to deduce
\begin{align*}
    \lVert t\mapsto \mathring{A}e^{-tA} u_0\rVert_{\dot{\B}^{\alpha}_{q,\infty}(\RR_+,\X)} \lesssim_{q,\alpha,A,\X} \lVert u_0\rVert_{\mathring{\mathcal{D}}_{A}(\alpha_q,\infty)}.
\end{align*}

\textbf{Step 5:} The case $q=\infty$, $\kappa<\infty$, $-1<\alpha<0$, admits a proof similar to the one in Step 4, taking advantage of Step 1 instead of Step 2.

\textbf{Step 6:} Only three cases are remaining: $s=0,\sfrac{1}{q},1$ or $\kappa\in[1,\infty]$, if $1< q <\infty$, $s=1$ or $\kappa\in[1,\infty]$  if $q=1$, and $s=0$ or $\kappa=\infty$ if $q=\infty$. Each can be obtained by real interpolation, thanks to  Lemma~\ref{lem:ReitThmSmigrpDaPrato} and Corollary~\ref{cor:InterpVectValuedHomBesov}.
\end{proof}

\section{Some norm estimates for Besov spaces and the Poisson semigroup}\label{App:EquivNormBesov}

\begin{lemma}\label{lem:DumpedPoissonSemigrpEquivNormBesov}Let $p,q\in[1,\infty]$, $s>0$ and $\alpha\in\RR$. Let $\mu\in[0,\pi)$. Then for all $\lambda := e^{i\theta}$, with  $|\theta| < \mu$.
\begin{enumerate}
    \item One has for all $u\in\S'(\RR^{n})$,
    \begin{align*}
        \big\lVert t\mapsto t^s(\lambda \I-\Delta)^\frac{\alpha}{2} e^{-t(\lambda \I-\Delta)^{\sfrac{1}{2}}} u\big\rVert_{\L^q_\ast(\RR_+,{\L}^{p}(\RR^n))} \sim_{n,p,s}^{\alpha,\mu} \lVert u \rVert_{\B^{\alpha-s}_{p,q}(\RR^n)}.
    \end{align*}
    \item We assume moreover that $\alpha\geqslant 0$. One has for all $u\in\S'_h(\RR^{n})$,
    \begin{align*}
        \big\lVert t\mapsto t^s(-\Delta)^\frac{\alpha}{2} e^{-t(-\Delta)^{\sfrac{1}{2}}} u\big\rVert_{\L^q_\ast(\RR_+,{\L}^{p}(\RR^n))} \sim_{n,p,s}^{\alpha} \lVert u \rVert_{\dot{\B}^{\alpha-s}_{p,q}(\RR^n)}.
    \end{align*}
\end{enumerate}
\end{lemma}

\begin{proof}One can reproduce almost \textit{verbatim} the proofs of \cite[Lemma~2.4~\&~Theorem~2.34]{bookBahouriCheminDanchin} for the heat semigroup.
\end{proof}

\begin{proposition}\label{prop:mixedSperateDerivEstBesovSpaces}Let $p,q\in[1,\infty]$, $s\in(-1+{\sfrac{1}{p}},{\sfrac{1}{p}})$.
\begin{enumerate}
    \item If $s>0$, all $u\in\dot{\B}^{s}_{p,q}(\RR^n_+)$,
    \begin{align*}
         \lVert u \rVert_{\dot{\B}^{s}_{p,q}(\RR^n_+)}\sim_{p,s,n} \lVert u \rVert_{\dot{\B}^{s}_{p,q}(\RR_+,\L^p(\RR^{n-1}))}+\lVert u \rVert_{\dot{\B}^{s}_{p,q}(\RR^{n-1},\L^p(\RR_+))};
    \end{align*}
  \item If $s>0$, $1<p<\infty$, for all $u\in\dot{\H}^{s,p}(\RR^n_+)$,
  \begin{align*}
      \lVert u \rVert_{\dot{\H}^{s,p}(\RR^n_+)}\sim_{p,s,n} \lVert (-\partial_{x_n})^s u \rVert_{{\L}^{p}(\RR^{n}_+)}+\lVert (-\Delta')^{\frac{s}{2}}u \rVert_{{\L}^{p}(\RR^{n}_+)};
  \end{align*}
  \item If $s<0$, $1<p<\infty$, for all $u\in\dot{\H}^{s,p}(\RR^n_+)$,
  \begin{align*}
      \lVert u \rVert_{\dot{\H}^{s,p}(\RR^n_+)}\lesssim_{p,s,n} \min \begin{cases}
          \lVert (-\partial_{x_n})^s u \rVert_{{\L}^{p}(\RR^{n}_+)}\\
          \lVert (-\Delta')^{\frac{s}{2}}u \rVert_{{\L}^{p}(\RR^{n}_+)}
      \end{cases}.
  \end{align*}
\end{enumerate}
\end{proposition}

\begin{proof}\textbf{Step 1:} The case $s>0$. By Proposition~\ref{prop:AppendixFiniteDifVectorValuedBesovRN+} and point \textit{(vii)} of Proposition~\ref{prop:BanachValuedHomBesovSpaces}, one obtains
\begin{align*}
         \lVert u \rVert_{\dot{\B}^{s}_{p,q}(\RR^n_+)}&\sim_{p,s,n} \sum_{k=1}^n \left\lVert t\mapsto t^{-s}\lVert\tau^{k}_{t}u-u\rVert_{\L^p(\RR^n_+)} \right\rVert_{\L^q_\ast(\RR_+)}\\
         &\sim_{p,s,n} \left\lVert t\mapsto t^{-s}\lVert\tau^{n}_{t}u-u\rVert_{\L^p(\RR_+,\L^p(\RR^{n-1}))} \right\rVert_{\L^q_\ast(\RR_+)}\\ &\qquad\qquad\qquad\qquad\qquad+ \sum_{k=1}^{n-1} \left\lVert t\mapsto t^{-s}\lVert\tau^{k}_{t}u-u\rVert_{\L^p(\RR^{n-1},\L^p(\RR_+))} \right\rVert_{\L^q_\ast(\RR_+)}\\ 
         &\sim_{p,s,n}\lVert u \rVert_{\dot{\B}^{s}_{p,q}(\RR_+,\L^p(\RR^{n-1}))}+\lVert u \rVert_{\dot{\B}^{s}_{p,q}(\RR^{n-1},\L^p(\RR_+))}.
\end{align*}
\textbf{Step 2:} The case of Sobolev spaces. We start with $s>0$. One applies \cite[Corollary~2.2]{Pruss2002}:

We consider 
\begin{itemize}
    \item $(\D_p(A+B),A+B)=(\H^{2,p}_\mathcal{N}(\RR^n_+),-\Delta_\mathcal{N})$ on $\L^p(\RR^n_+)$, $-\Delta=-\partial_{x_n}^2-\Delta'$;
    \item $(\D_p(A),A)=(\H^{2,p}_{\mathcal{N},{x_n}}(\RR_+,\L^p(\RR^{n-1})),-\partial_{x_n}^2)$
    \item $(\D_p(B),B)=(\L^p_{x_n}(\RR_+,\H^{2,p}(\RR^{n-1})),-\Delta')$
\end{itemize}
are all $0$-sectorial operators, and such that $A$ and $B$ have commuting resolvents, and which admit BIP on the UMD Banach space $\L^p(\RR^n_+) = \L^p(\RR^{n-1},\L^p(\RR_+))= \L^p(\RR_+,\L^p(\RR^{n-1}))$.
This yields for all $0<s<1$, with equivalence of norms
\begin{align*}
    \H^{s,p}(\RR^{n}_+)&=\H^{s,p}_{x_n}(\RR_+,\L^p(\RR^{n-1}))\cap \H^{s,p}(\RR^{n-1},\L^p_{x_n}(\RR_{+}))\\&= \H^{s,p}_{x_n}(\RR_+,\L^p(\RR^{n-1}))\cap \L^p_{x_n}(\RR_{+},\H^{s,p}(\RR^{n-1})).
\end{align*}
From \cite[Proposition~3.3]{Gaudin2022} giving $\L^p\cap\dot{\H}^{s,p}(\RR^n_+) = \H^{s,p}(\RR^n_+)$ with equivalence of norms, and  \cite[Proposition~3.6,~\textit{(vi)}]{Gaudin2023} giving $$\H^{s,p}_{x_n}(\RR_+,\L^p(\RR^{n-1})) = \L^{p}_{x_n}(\RR_+,\L^p(\RR^{n-1}))\cap \dot{\H}^{s,p}_{x_n}(\RR_+,\L^p(\RR^{n-1})),$$ so that for all $u\in\H^{s,p}(\RR^{n}_+)$, it holds
\begin{align*}
    \lVert u\rVert_{\L^p(\RR^n_+)}+ \lVert u\rVert_{\dot{\H}^{s,p}(\RR^n_+)} &\sim_{p,s,n} \lVert u\rVert_{\L^p(\RR^n_+)} + \lVert u\rVert_{\dot{\H}^{s,p}_{x_n}(\RR_+,\L^p(\RR^{n-1}))}+ \lVert u\rVert_{\L^p_{x_n}(\RR_{+},\H^{s,p}(\RR^{n-1}))}\\
    &\sim_{p,s,n} \lVert u\rVert_{\L^p(\RR^n_+)} + \lVert (-\partial_{x_n})^s u\rVert_{\L^p(\RR^n_+)}+ \lVert (-\Delta')^\frac{s}{2}u\rVert_{\L^p(\RR^n_+)}.
\end{align*}
The last estimate is derived from the isomorphism property \cite[Eq.~(3.5)]{Gaudin2023}. A dilation argument and the density of $\H^{s,p}(\RR^n_+)$ in $\dot{\H}^{s,p}(\RR^n_+)$ yield, for all $u\in\dot{\H}^{s,p}(\RR^n_+)$,
\begin{align*}
    \lVert u\rVert_{\dot{\H}^{s,p}(\RR^n_+)} \sim_{p,s,n} \lVert (-\partial_{x_n})^s u\rVert_{\L^p(\RR^n_+)}+ \lVert (-\Delta')^\frac{s}{2}u\rVert_{\L^p(\RR^n_+)}.
\end{align*}
Since $0<s<\frac{1}{p}$ we are using the isomorphism property \cite[eq.~(3.2)]{Gaudin2023}, and one obtains the same result with $(\partial_t)^{s}=(-\partial_t)^{s\,\ast}$, and in particular
\begin{align*}
     \lVert (\partial_{x_n})^s u\rVert_{\L^p(\RR^n_+)} \lesssim_{p,s,n} \lVert u\rVert_{\dot{\H}^{s,p}(\RR^n_+)}\text{, } \forall u\in\dot{\H}^{s,p}(\RR^n_+).
\end{align*}
Now, if $s<0$, let $u\in\dot{\H}^{s,p}(\RR^n_+)$ and $v\in\Ccinfty(\RR^n_+)$ such that $\lVert v\rVert_{\dot{\H}^{-s,p'}(\RR^n_+)}=1$, by duality, it holds
\begin{align*}
    |\langle u,v\rangle_{\RR^n_+}| = |\langle (-\partial_{x_n})^s u, (\partial_{x_n})^{-s}v\rangle_{\RR^n_+}| &\leqslant \lVert (-\partial_{x_n})^s u\rVert_{\L^p(\RR^n_+)} \lVert (\partial_{x_n})^{-s} v\rVert_{\L^{p'}(\RR^n_+)}\\
    &\lesssim_{p,s,n} \lVert (-\partial_{x_n})^s u\rVert_{\L^p(\RR^n_+)} 
\end{align*}
By duality, one can conclude
\begin{align*}
    \lVert u \rVert_{\dot{\H}^{s,p}(\RR^n_+)}\lesssim_{p,s,n} \lVert (-\partial_{x_n})^s u \rVert_{{\L}^{p}(\RR^{n}_+)}.
\end{align*}
A similar argument yields
\begin{align*}
    \lVert u \rVert_{\dot{\H}^{s,p}(\RR^n_+)}\lesssim_{p,s,n} \lVert (-\Delta')^{\frac{s}{2}} u \rVert_{{\L}^{p}(\RR^{n}_+)},
\end{align*}
which ends the proof.
\end{proof}

\begin{proposition}\label{prop:DumpedPoissonSemigroup2} Let $p,q\in[1,\infty]$, $s>-1+{\sfrac{1}{p}}$. For $\mu\in[0,\pi)$, $\lambda=e^{i\theta}$, $|\theta|\leqslant\mu$, we consider
\begin{align*}
    T_\lambda \,:\, f\longmapsto \big[ (x',x_n)\mapsto e^{-x_n(\lambda\I-\Delta')^{\sfrac{1}{2}}}f(x')\big].
\end{align*}
\begin{enumerate}
    \item We assume $1<p<\infty$. For all $f\in \B^{s-{\sfrac{1}{p}}}_{p,p}(\RR^{n-1})$, one has $T_\lambda f \in \dot{\H}^{s,p}(\RR^n_+)$, with the estimate
    \begin{align*}
        \lVert T_\lambda f\rVert_{\dot{\H}^{s,p}(\RR^n_+)} \lesssim_{p,s,n}^{\mu} \lVert  f\rVert_{{\B}^{s-{\sfrac{1}{p}}}_{p,p}(\RR^{n-1})}.
    \end{align*}
    \item  For all $f\in \B^{s-{\sfrac{1}{p}}}_{p,q}(\RR^{n-1})$, one has $T_\lambda f \in \dot{\B}^{s}_{p,q}(\RR^n_+)$, with the estimate
    \begin{align*}
        \lVert T_\lambda f\rVert_{\dot{\B}^{s}_{p,q}(\RR^n_+)} \lesssim_{p,s,n}^{\mu} \lVert  f\rVert_{{\B}^{s-{\sfrac{1}{p}}}_{p,q}(\RR^{n-1})}.
    \end{align*}
    \item If $p<\infty$, for $k\in\mathbb{N}$, for all $f\in\B^{k-{\sfrac{1}{p}}}_{p,p}(\RR^{n-1})$,
    \begin{align*}
        \lVert \nabla^k T_\lambda f\rVert_{{\L}^{p}(\RR^n_+)} \lesssim_{p,s,n}^{\mu} \lVert  f\rVert_{{\B}^{k-{\sfrac{1}{p}}}_{p,p}(\RR^{n-1})}.
    \end{align*}
    When $p=\infty$, the result still holds with ${\B}^{k}_{\infty,1}$ replacing ${\B}^{k-{\sfrac{1}{p}}}_{p,p}$, and $T_\lambda$ maps $\L^\infty(\RR^{n-1})$ to $\L^\infty_h(\RR^{n}_+)$ boundedly.
\end{enumerate}
Furthermore, in the case of Besov spaces for $\mathfrak{B}\in\{\B^{s,0}_{\infty,q},\BesSmo^{s}_{p,\infty},\BesSmo^{s,0}_{\infty,\infty}\}$, $T_\lambda \mathfrak{B}(\RR^{n-1}) \subset \dot{\mathfrak{B}} (\RR^n_+)$.
\end{proposition}

\begin{remark} We do not use directly interpolation to obtain the result, since we want to reach non-positive regularity indices $s<0$ and targeted spaces are the homogeneous function spaces.
\end{remark}

\begin{proof}We write $A_{\lambda}=(\lambda\I-\Delta')^{\sfrac{1}{2}}$.  When $\lambda\notin(-\infty,0]$, we recall that $A_{\lambda}$ has a bounded inverse on $\L^p(\RR^n)$ for all $p\in[1,\infty]$.

\textbf{Step 1:} The case $p=\infty$, $-1<s<0$. For $f\in{\B}^{s}_{\infty,q}(\RR^{n-1})$, one writes $T_\lambda f = \partial_{x_n}T_\lambda (\lambda\I-\Delta')^{-{\sfrac{1}{2}}}f$. It holds that:
\begin{align*}
    \lVert T_\lambda f \rVert_{\dot{\B}^{s}_{\infty,q}(\RR^n_+)} &\leqslant \lVert \partial_{x_n} \E_\mathcal{N}T_\lambda [A_{\lambda}^{-1}f]\rVert_{\dot{\B}^{s}_{\infty,q}(\RR^n)} \\
    &\lesssim_{s,n} \lVert  \E_\mathcal{N}T_\lambda (\lambda\I-\Delta')^{-{\sfrac{1}{2}}}f\rVert_{\dot{\B}^{s+1}_{\infty,q}(\RR^n)}\\
    &\lesssim_{s,n} \sum_{k=1}^n \left\lVert t\mapsto t^{-(s+1)}\lVert\tau^{k}_{t}T_\lambda [A_{\lambda}^{-1}f]-T_\lambda [A_{\lambda}^{-1}f]\rVert_{\L^\infty(\RR^n_+)} \right\rVert_{\L^q_\ast(\RR_+)}\\
    & \qquad \qquad \qquad \text{(by Proposition~\ref{prop:AppendixFiniteDifVectorValuedBesovRN+})}\\ 
    &\lesssim_{s,n,\mu} \sum_{k=1}^{n-1} \left\lVert t\mapsto t^{-(s+1)}\lVert\tau^{k}_{t} [A_{\lambda}^{-1}f]-[A_{\lambda}^{-1}f]\rVert_{\L^\infty(\RR^{n-1})} \right\rVert_{\L^q_\ast(\RR_+)}\\
    & \qquad \qquad \qquad +\left\lVert t\mapsto t^{-(s+1)}\lVert\tau^{n}_{t}T_\lambda [A_{\lambda}^{-1}f]-T_\lambda [A_{\lambda}^{-1}f]\rVert_{\L^\infty(\RR^n_+)} \right\rVert_{\L^q_\ast(\RR_+)}\\
    & \quad \text{(by uniform boundedness of the semigroup $T_\lambda=(e^{-x_n A_\lambda})_{x_n>0}$ on $\L^\infty(\RR^{n-1})$)}\\
    &\lesssim_{s,n,\mu} \lVert (\lambda\I-\Delta')^{-{\sfrac{1}{2}}}f \rVert_{\dot{\B}^{s+1}_{\infty,q}(\RR^{n-1})}\\
    & \qquad \qquad \qquad +\left\lVert t\mapsto t^{-(s+1)}\lVert\tau^{n}_{t}AT_\lambda [A_{\lambda}^{-2}f]-AT_\lambda [A_{\lambda}^{-2}f]\rVert_{\L^\infty(\RR^n_+)} \right\rVert_{\L^q_\ast(\RR_+)}\\
    &    \qquad \qquad \qquad \text{(by point \textit{(vii)} of Proposition~\ref{prop:BanachValuedHomBesovSpaces})}\\
    &\lesssim_{s,n,\mu} \lVert f \rVert_{{\B}^{s}_{\infty,q}(\RR^{n-1})} + \lVert A_{\lambda}^{-2}f\rVert_{\mathring{\mathcal{D}}_{A_\lambda}(2+s,q)}\\
    & \qquad \qquad \qquad \text{(by Proposition~\ref{prop:AppendixFiniteDifVectorValuedBesovRN+} and Lemma~\ref{lem:BesovintimeSemigroup})}\\ 
    &\lesssim_{s,n,\mu} \lVert f \rVert_{{\B}^{s}_{\infty,q}(\RR^{n-1})}.
\end{align*}
The last line is obtained by Lemma~\ref{lem:DumpedPoissonSemigrpEquivNormBesov}.

Now, if $f\in{\B}^{s+1}_{\infty,q}(\RR^{n-1})$, by the previous estimate:
\begin{align*}
    \lVert T_\lambda f \rVert_{\dot{\B}^{s+1}_{\infty,q}(\RR^n_+)} &\leqslant \lVert \E_\mathcal{N}T_\lambda f\rVert_{\dot{\B}^{s+1}_{\infty,q}(\RR^n)} \\
    &\lesssim_{s,n} \lVert  \E_\mathcal{N}T_\lambda (\lambda\I-\Delta')^{-{\sfrac{1}{2}}}[(\lambda\I-\Delta')^{{\sfrac{1}{2}}}f]\rVert_{\dot{\B}^{s+1}_{\infty,q}(\RR^n)}\\
    &\lesssim_{s,n,\mu} \lVert (\lambda\I-\Delta')^{{\sfrac{1}{2}}}f \rVert_{{\B}^{s}_{\infty,q}(\RR^{n-1})}\\
    &\lesssim_{s,n,\mu} \lVert f \rVert_{{\B}^{s+1}_{\infty,q}(\RR^{n-1})}.
\end{align*}
\textbf{Step 2:} The case $p=1$, $0<s<1$. For $f\in{\B}^{s-1}_{1,q}(\RR^{n-1})$, one writes $T_\lambda f = A_{\lambda} T_\lambda [A_{\lambda}^{-1}f]$. By Proposition~\ref{prop:mixedSperateDerivEstBesovSpaces}, it holds that:
\begin{align*}
    \lVert T_\lambda f \rVert_{\dot{\B}^{s}_{1,\infty}(\RR^n_+)} &\sim_{s,n} \lVert T_\lambda f \rVert_{\dot{\B}^{s}_{1,\infty}(\RR_+,\L^1(\RR^{n-1}))}+\lVert  T_\lambda f \rVert_{\dot{\B}^{s}_{1,\infty}(\RR^{n-1},\L^1(\RR_+))}\\
    &\lesssim_{s,n} \lVert A_\lambda T_\lambda [A_\lambda^{-1} f] \rVert_{\dot{\B}^{s}_{1,1}(\RR_+,\L^1(\RR^{n-1}))}+\lVert A_0^{s}  T_\lambda f \rVert_{{\L}^{1}(\RR^{n-1},\L^1(\RR_+))}\\
    &\lesssim_{s,n} \lVert A_\lambda T_\lambda [A_\lambda^{-1} f] \rVert_{\dot{\B}^{s}_{1,1}(\RR_+,\L^1(\RR^{n-1}))}+\lVert A_\lambda^{s} T_\lambda [ A_0^{s} A_\lambda^{-s} f] \rVert_{{\L}^{1}(\RR^{n}_+)}\\
    &\lesssim_{s,n}^{\mu} \lVert [A_\lambda^{-1} f] \rVert_{\mathring{\mathcal{D}}_{A_\lambda}(s,1)}+\lVert A_0^{s} A_\lambda^{-s} f \rVert_{{\B}^{s-1}_{1,1}(\RR^{n-1})}\\
    &\qquad \qquad \qquad \text{(by Lemmas~\ref{lem:DumpedPoissonSemigrpEquivNormBesov} and \ref{lem:BesovintimeSemigroup})}\\
    &\lesssim_{s,n}^{\mu} \lVert f \rVert_{{\B}^{s-1}_{1,1}(\RR^{n-1})},
\end{align*}
where the last inequality follows from Lemma~\ref{lem:DumpedPoissonSemigrpEquivNormBesov} and the boundedness of $A_0^{s} A_\lambda^{-s}$ on $\B^{s-1}_{1,1}(\RR^{n-1})$. Since it holds for all $0<s<1$, by real interpolation, one obtains for all $q\in[1,\infty]$ and all $u\in\B^{s-1}_{1,q}(\RR^{n-1})$,
\begin{align*}
    \lVert T_\lambda f \rVert_{\dot{\B}^{s}_{1,q}(\RR^n_+)} &\lesssim_{s,n}^{\mu} \lVert f \rVert_{{\B}^{s-1}_{1,q}(\RR^{n-1})}.
\end{align*}
Now, if $f\in\B^{s}_{1,q}(\RR^{n-1})$, by the previous estimate
\begin{align*}
    \lVert T_\lambda f \rVert_{\dot{\B}^{s+1}_{1,q}(\RR^n_+)}&\lesssim_{s,n} \lVert \partial_{x_n} T_\lambda f \rVert_{\dot{\B}^{s+1}_{1,q}(\RR^n_+)} + \lVert \nabla' T_\lambda f \rVert_{\dot{\B}^{s+1}_{1,q}(\RR^n_+)}\\
    &\lesssim_{s,n}^{\mu} \lVert A_\lambda f \rVert_{{\B}^{s-1}_{1,q}(\RR^{n-1})} + \lVert \nabla' f \rVert_{{\B}^{s-1}_{1,q}(\RR^{n-1})}\\
    &\lesssim_{s,n}^{\mu} \lVert f \rVert_{{\B}^{s}_{1,q}(\RR^{n-1})}.
\end{align*}

\textbf{Step 3:} The case $p\in(1,\infty)$. Let $f\in\B^{s-{\sfrac{1}{p}}}_{p,p}(\RR^{n-1})$. One can check $(\partial_{x_n})^{s}T_\lambda = A_\lambda^{s}T_\lambda$, hence, by Proposition~\ref{prop:mixedSperateDerivEstBesovSpaces} and Lemma~\ref{lem:DumpedPoissonSemigrpEquivNormBesov}, if $s\leqslant 0$,
\begin{align*}
    \lVert T_\lambda f \rVert_{\dot{\H}^{s,p}(\RR^n_+)} \lesssim_{p,s,n} \lVert A_\lambda^{s}T_\lambda f \rVert_{{\L}^{p}(\RR^n_+)} \sim_{p,s,n}^{\mu} \lVert f \rVert_{{\B}^{s-{\sfrac{1}{p}}}_{p,p}(\RR^{n-1})}.
\end{align*}
If $s>0$,
\begin{align*}
    \lVert T_\lambda f \rVert_{\dot{\H}^{s,p}(\RR^n_+)} &\lesssim_{p,s,n} \lVert A_\lambda^{s}T_\lambda f \rVert_{{\L}^{p}(\RR^n_+)}+\lVert A_0^{s}T_\lambda f \rVert_{{\L}^{p}(\RR^n_+)}\\&\sim_{p,s,n}^{\mu} \lVert f \rVert_{{\B}^{s-{\sfrac{1}{p}}}_{p,p}(\RR^{n-1})} + \lVert  A_0^{s}A_\lambda^{-s}f \rVert_{{\B}^{s-{\sfrac{1}{p}}}_{p,p}(\RR^{n-1})}\\&\lesssim_{p,s,n}\lVert f \rVert_{{\B}^{s-{\sfrac{1}{p}}}_{p,p}(\RR^{n-1})}.
\end{align*}
Now, we assume $s\in(-1+{\sfrac{1}{p}},{\sfrac{1}{p}})$, and we consider $f\in\B^{s+1-{\sfrac{1}{p}}}_{p,p}(\RR^{n-1})$
\begin{align*}
    \lVert T_\lambda f \rVert_{\dot{\H}^{s+1,p}(\RR^n_+)} &\lesssim_{p,s,n} \lVert A_\lambda T_\lambda f \rVert_{\dot{\H}^{s,p}(\RR^n_+)}+\lVert \nabla'T_\lambda f \rVert_{\dot{\H}^{s,p}(\RR^n_+)} \\ &\sim_{p,s,n}^{\mu} \lVert f \rVert_{{\B}^{s+1-{\sfrac{1}{p}}}_{p,p}(\RR^{n-1})} + \lVert  \nabla'f \rVert_{{\B}^{s-{\sfrac{1}{p}}}_{p,p}(\RR^{n-1})}\\&\lesssim_{p,s,n}\lVert f \rVert_{{\B}^{s+1-{\sfrac{1}{p}}}_{p,p}(\RR^{n-1})}.
\end{align*}

Now, one can perform real interpolation so that the result holds also for Besov spaces $p\in(1,\infty)$, $q\in[1,\infty]$, $s\in(-1+{\sfrac{1}{p}},1+{\sfrac{1}{p}})$. Finally, by induction, taking more and more derivatives, the result holds for all $s>-1+{\sfrac{1}{p}}$.

Concerning $\L^\infty$, for $u\in\L^\infty(\RR^{n-1})\subset\B^{s}_{\infty,1}(\RR^{n-1})$, $-1<s<0$, by previous steps, it holds that $T_\lambda u\in\L^\infty(\RR^{n}_+)\cap\dot{\B}^{s}_{\infty,1}(\RR^{n}_+)\subset\L^\infty_h(\RR^{n}_+)$.
\end{proof}

\begin{proposition}\label{prop:PoissonSemigroup3} Let $p,q\in[1,\infty]$, $s\in\RR$. We set
\begin{align*}
    T_0 \,:\, f\longmapsto \big[ (x',x_n)\mapsto e^{-x_n(-\Delta')^{\sfrac{1}{2}}}f(x')\big].
\end{align*}
\begin{enumerate}
    \item We assume $1<p<\infty$. For all $f\in \dot{\B}^{s-{\sfrac{1}{p}}}_{p,p}(\RR^{n-1})$, one has $T_0 f \in \dot{\H}^{s,p}(\RR^n_+)$, with the estimate
    \begin{align*}
        \lVert T_0 f\rVert_{\dot{\H}^{s,p}(\RR^n_+)} \lesssim_{p,s,n} \lVert  f\rVert_{\dot{\B}^{s-{\sfrac{1}{p}}}_{p,p}(\RR^{n-1})}.
    \end{align*}
    \item  For all $f\in \dot{\B}^{s-{\sfrac{1}{p}}}_{p,q}(\RR^{n-1})$, one has $T_0 f \in \dot{\B}^{s}_{p,q}(\RR^n_+)$, with the estimate
    \begin{align*}
        \lVert T_0 f\rVert_{\dot{\B}^{s}_{p,q}(\RR^n_+)} \lesssim_{p,s,n} \lVert  f\rVert_{\dot{\B}^{s-{\sfrac{1}{p}}}_{p,q}(\RR^{n-1})}.
    \end{align*}
    \item If $p<\infty$, for $k\in\mathbb{N}$, for all $f\in\dot{\B}^{k-{\sfrac{1}{p}}}_{p,p}(\RR^{n-1})$,
    \begin{align*}
        \lVert \nabla^k T_0 f\rVert_{{\L}^{p}(\RR^n_+)} \lesssim_{p,s,n} \lVert  f\rVert_{\dot{\B}^{k-{\sfrac{1}{p}}}_{p,p}(\RR^{n-1})}.
    \end{align*}
    When $p=\infty$, the result still holds with $\dot{\B}^{k}_{\infty,1}$ replacing $\dot{\B}^{k-{\sfrac{1}{p}}}_{p,p}$, and $T_0$ maps $\L^\infty_h(\RR^{n-1})$ to $\L^\infty_h(\RR^n_+)$.
\end{enumerate}
Furthermore, in the case of Besov spaces for $\mathfrak{B}\in\{\B^{s,0}_{\infty,q},\BesSmo^{s}_{p,\infty},\BesSmo^{s,0}_{\infty,\infty}\}$, $T_0 \dot{\mathfrak{B}}(\RR^{n-1}) \subset \dot{\mathfrak{B}} (\RR^n_+)$.
\end{proposition}

\begin{proof}We only deal with the case $p=\infty$, otherwise one could consult the very similar result \cite[Proposition~2]{Gaudin2022}, which deals only with the case $p\in(1,\infty)$. It suffices to adapt the proof with the elements presented here to obtain the full result (allowing to get rid of the completeness assumption).

\textbf{Step 1:} We assume first $q<\infty$, $s<0$. Let $f\in\dot{\B}^{s}_{\infty,q}(\RR^{n-1})$, and we set for all $N\in\NN$
\begin{align*}
    f_N = \sum_{|j|\leqslant N} \dot{\Delta}'_jf.
\end{align*}
Notice that $(f_N)_{N\in\NN}\subset\C^{\infty}_{ub,h}(\RR^n)\cap \Big(\bigcap_{\substack{\alpha\in\RR\\ 1\leqslant r\leqslant \infty}}\dot{\B}^{\alpha}_{\infty,r}(\RR^{n-1})\Big)$ converges to $f$ in $\dot{\B}^{s}_{\infty,q}(\RR^{n-1})$, due to the fact that $q<\infty$. Remember that for $s<0$, the homogeneous Besov spaces are always complete. We assume temporarily $s\notin\ZZ$.

So, let $k\in\NN$ such that $0<s+k<1$, it holds
\begin{align*}
     &\lVert T_0 f_N\rVert_{\dot{\B}^{s}_{\infty,q}(\RR^n_+)}= \lVert \partial_{x_n}^k T_0 [(-\Delta')^{-\frac{k}{2}}f_N]\rVert_{\dot{\B}^{s}_{\infty,q}(\RR^n_+)}\\
    &\lesssim_{s,k,n} \lVert \E_\mathcal{N} T_0 [(-\Delta')^{-\frac{k}{2}}f_N]\rVert_{\dot{\B}^{s+k}_{\infty,q}(\RR^n)}\\
    &\lesssim_{s,k,n} \sum_{k=1}^n \left\lVert t\mapsto t^{-(s+k)}\lVert\tau^{k}_{t}T_0 [(-\Delta')^{-\frac{k}{2}}f_N]-T_0 [(-\Delta')^{-\frac{k}{2}}f_N]\rVert_{\L^\infty(\RR^n_+)} \right\rVert_{\L^q_\ast(\RR_+)}\\
    & \qquad \qquad \qquad \text{(by Proposition~\ref{prop:AppendixFiniteDifVectorValuedBesovRN+})}\\ 
    &\lesssim_{s,k,n} \sum_{k=1}^{n-1} \left\lVert t\mapsto t^{-(s+k)}\lVert\tau^{k}_{t} [(-\Delta')^{-\frac{k}{2}}f_N]-[(-\Delta')^{-\frac{k}{2}}f_N]\rVert_{\L^\infty(\RR^{n-1})} \right\rVert_{\L^q_\ast(\RR_+)}\\
    & \qquad \qquad \qquad +\left\lVert t\mapsto t^{-(s+1)}\lVert\tau^{n}_{t}T_\lambda [(-\Delta')^{-\frac{k}{2}}f_N]-T_0 [(-\Delta')^{-\frac{k}{2}}f_N]\rVert_{\L^\infty(\RR^n_+)} \right\rVert_{\L^q_\ast(\RR_+)}\\
    & \qquad \qquad \qquad \text{(by uniform boundedness of the semigroup $T_0=(e^{-x_n (-\Delta')^{\sfrac{1}{2}})})_{x_n>0}$ on $\L^\infty(\RR^{n-1})$)}\\
    &\lesssim_{s,k,n} \lVert (-\Delta')^{-\frac{k}{2}}f_N\rVert_{\dot{\B}^{s+k}_{\infty,q}(\RR^{n-1})}\\
    & \qquad \qquad \qquad +\left\lVert t\mapsto t^{-(s+k)}\lVert\tau^{n}_{t}(-\Delta')^{{\sfrac{1}{2}}}T_\lambda [(-\Delta')^{-\frac{k+1}{2}}f_N]-(-\Delta')^{{\sfrac{1}{2}}}T_0 [(-\Delta')^{-\frac{k+1}{2}}f_N]\rVert_{\L^\infty(\RR^n_+)} \right\rVert_{\L^q_\ast(\RR_+)}\\
    &    \qquad \qquad \qquad \text{(by point \textit{(vii)} of Proposition~\ref{prop:BanachValuedHomBesovSpaces})}\\
    &\lesssim_{s,n,\mu} \lVert f_N \rVert_{\dot{\B}^{s}_{\infty,q}(\RR^{n-1})} + \lVert (-\Delta')^{-\frac{k+1}{2}}f_N\rVert_{\mathring{\mathcal{D}}_{(-\Delta')^{1/2}}(1+s+k,q)}\\
    & \qquad \qquad \qquad \text{(by Proposition~\ref{prop:AppendixFiniteDifVectorValuedBesovRN+} and Lemma~\ref{lem:BesovintimeSemigroup})}\\ 
    &\lesssim_{s,n,\mu} \lVert f_N \rVert_{\dot{\B}^{s}_{\infty,q}(\RR^{n-1})}.
\end{align*}
The last line is obtained by Lemma~\ref{lem:DumpedPoissonSemigrpEquivNormBesov}. Passing to the limit as $N$ goes to infinity yields the result. Real interpolation yields the case $q=\infty$ whenever $s<0$.

\textbf{Step 2:} We show that $T_0$ maps $\L^\infty(\RR^{n-1})$ to $\L^\infty_h(\RR^n_+)$.

By the previous Step 1: one has the boundedness
\begin{align}\label{eq:ProofPoissonSemirgpBesovinfty}
    T_0\,:\,\dot{\B}^{s}_{\infty,q}(\RR^{n-1})\longrightarrow\dot{\B}^{s}_{\infty,q}(\RR^{n}_+),\, s<0,\,q\in[1,\infty].
\end{align}
But by Young's inequality for the convolution one also trivially has the boundedness of 
\begin{align}\label{eq:ProofPoissonSemirgpLinfty}
    T_0\,:\,\L^\infty(\RR^{n-1})\longrightarrow\L^\infty(\RR^{n}_+).
\end{align}

By \cite[Lemma~2.18]{Gaudin2023Lip}, for $u\in\L^\infty_h(\RR^{n-1})$, for all $N\in\NN$, $u_N:=\sum_{j>-N}\dot{\Delta}_ju\in \L^\infty(\RR^{n-1})\cap \dot{\B}^{s}_{\infty,1}(\RR^{n-1})$, so that by \eqref{eq:ProofPoissonSemirgpBesovinfty} and \eqref{eq:ProofPoissonSemirgpLinfty}, one obtains $T_0 u_N \in\L^\infty(\RR^{n}_+)\cap \dot{\B}^{s}_{\infty,1}(\RR^{n}_+)\subset \L^\infty_h(\RR^{n}_+)$, with
\begin{align*}
    \lVert T_0 u -T_0 u_N  \rVert_{\L^\infty(\RR^n_+)}\leqslant \lVert u -u_N  \rVert_{\L^\infty(\RR^{n-1})} \xrightarrow[N\rightarrow+\infty]{}0,
\end{align*}
since $\L^\infty_h(\RR^{n}_+)$ is closed in $\L^\infty(\RR^{n}_+)$, one deduces $T_0 u \in\L^\infty_h(\RR^n_+)$. Thus,
\begin{align}\label{eq:ProofPoissonSemirgpLinftyh}
    T_0\,:\,\L^\infty_h(\RR^{n-1})\longrightarrow\L^\infty_h(\RR^{n}_+).
\end{align}
is well-defined and bounded.

\textbf{Step 3:}  Now, we are going to show that, for all $s\geqslant 0$, one has $T_0( \dot{\B}^{s}_{\infty,q}(\RR^{n-1})) \subset\dot{\B}^{s}_{\infty,q}(\RR^{n}_+)$ with the appropriate control of norms.

We focus on homogeneous Sobolev spaces of integer order. Let $k\in\mathbb{N}^\ast$. Since $\dot{\B}^{k}_{\infty,1}(\RR^{n-1})\hookrightarrow\dot{\W}^{k,\infty}(\RR^{n-1})\subset\L^\infty_h(\RR^{n-1})$, by \eqref{eq:ProofPoissonSemirgpLinftyh} $T_0(\dot{\B}^{k}_{\infty,1}(\RR^{n-1}))\subset \L^\infty_h(\RR^{n}_+)$. For all $u\in\dot{\B}^{k}_{\infty,1}(\RR^{n-1})$, it holds
\begin{align*}
    \lVert T_0u\rVert_{\dot{\W}^{k,\infty}(\RR^{n}_+)}&\sim_{k,n}\lVert \nabla^kT_0u\rVert_{{\L}^{\infty}(\RR^{n}_+)}\\
    &\lesssim_{k,n}\sum_{|\alpha|=k} \lVert T_0 [\partial^{\alpha'}_{x'}(-\Delta')^{\alpha_n/2} u]\rVert_{{\L}^{\infty}(\RR^{n}_+)}\\
    &\lesssim_{k,n} \sum_{|\alpha|=k} \lVert  \partial^{\alpha'}_{x'}(-\Delta')^{\alpha_n/2} u\rVert_{{\L}^{\infty}(\RR^{n-1})}\lesssim_{k,n} \lVert u\rVert_{\dot{\B}^{k}_{\infty,1}(\RR^{n-1})}.
\end{align*}
Thus, we did obtain the well-defined bounded map
\begin{align}\label{eq:ProofPoissonSemirgpHomSobLinfty}
    T_0\,:\,\dot{\B}^{k}_{\infty,1}(\RR^{n-1})\longrightarrow\dot{\W}^{k,\infty}(\RR^{n}_+), \, k\in\NN.
\end{align}
By real interpolation between \eqref{eq:ProofPoissonSemirgpBesovinfty} and \eqref{eq:ProofPoissonSemirgpHomSobLinfty}, one obtains the following well-defined bounded map for all $s<k$:
\begin{align}\label{eq:ProofPoissonSemirgpBesovBesovLinfty}
    T_0\,:\,\dot{\B}^{s}_{\infty,q}(\RR^{n-1})\longrightarrow\dot{\B}^{s}_{\infty,q}(\RR^{n}_+).
\end{align}
Since $k\in\NN^\ast$ is arbitrary, one obtains the boundedness of \eqref{eq:ProofPoissonSemirgpBesovBesovLinfty} for any $s\in\RR$. This finishes the proof.
\end{proof}

\begin{corollary}\label{cor:PoissonSmigrpInhomSpaces} Let $s\in\RR$ and $q\in[1,\infty]$, then for all $u\in\B^{s}_{\infty,q}(\RR^{n-1})$, one has $T_0 u \in\B^{s}_{\infty,q}(\RR^{n}_+)$ with the estimate
\begin{align*}
    \lVert T_0u\rVert_{{\B}^{s}_{\infty,q}(\RR^{n}_+)}\lesssim_{s,n}\lVert u\rVert_{{\B}^{s}_{\infty,q}(\RR^{n-1})}.
\end{align*}
\end{corollary}

\begin{proof}The only case we have to deal with is the case $s\leqslant 0$. The case $s>0$, follows from the $\L^\infty$-case by induction and real interpolation.

First, we assume $s<0$. Consider $u\in{\B}^{s}_{\infty,q}(\RR^{n-1})= \L^\infty(\RR^{n-1})+\dot{\B}^{s}_{\infty,q}(\RR^{n-1})$, we write $u=a+b$,  by the $\L^\infty$-case, and the case of homogeneous function spaces from previous Proposition~\ref{prop:PoissonSemigroup3}, one obtains
\begin{align*}
    \lVert T_0u\rVert_{{\B}^{s}_{\infty,q}(\RR^{n}_+)}&\leqslant \lVert T_0a\rVert_{\L^\infty(\RR^{n}_+)} + \lVert T_0b\rVert_{\dot{\B}^{s}_{\infty,q}(\RR^{n}_+)}\\ 
    &\lesssim_{s,n} \lVert a\rVert_{\L^\infty(\RR^{n-1})} + \lVert b\rVert_{\dot{\B}^{s}_{\infty,q}(\RR^{n-1})}.
\end{align*}
Taking the infimum on such a decomposition yields
\begin{align*}
    \lVert T_0u\rVert_{{\B}^{s}_{\infty,q}(\RR^{n}_+)}\lesssim_{s,n} \lVert u\rVert_{{\B}^{s}_{\infty,q}(\RR^{n-1})}.
\end{align*}
Real interpolation between the case $s<0$ and $s>0$ yields the case $s=0$, which ends the proof.
\end{proof}

\smallskip
\par\noindent
{\bf Acknowledgement}.

\section*{Compliance with Ethical Standards}\label{conflicts}
\smallskip
\par\noindent
{\bf Conflict of Interest}. The author declares that he has no conflict of interest.

\smallskip
\par\noindent
{\bf Data Availability}. Data sharing is not applicable to this article as no datasets were generated or analysed during the current study.

\typeout{}

\end{document}